%% file: Construction_of_SFT.tex
\documentclass{article}

\input{SFT-preamble}

\begin{document}
\maketitle

\begin{abstract}
We construct symplectic field theory in general case completely.
We use Kuranishi theory for the construction.
For the construction of the Kuranishi neighborhood of a holomorphic building
of genus $>0$, we introduce a new space which parametrizes the deformations of
both of the domain curve and the target space.
We also improve the theory of Kuranishi structure and introduce the new notion of
pre-Kuranishi structure and its weakly good coordinate system.
Although the product of good coordinate systems is not
a good coordinate system, weakly good coordinate system is closed with
respect to product, and we can use their product directly for the product of
pre-Kuranishi spaces.
We also explain a new way to prove the smoothness of pre-Kuranishi structure
by using the estimates of the differentials of implicit functions.
We can obtain the estimate of the implicit functions by direct calculations
using appropriate coordinates.
We treat symplectic field theory of Bott-Morse case
by using a triangulation of the space of periodic orbits.
\end{abstract}


\tableofcontents

\input{SFT-01_Introduction_Estimates}
\input{SFT-02_Holomorphic_buildings}
\input{SFT-00_Theory_of_Kuranishi_structure}
\input{SFT-03_Kuranishi_neighborhood}
\input{SFT-04_Smoothness_Embeddings}
\input{SFT-05_Global_Kuranishi_structure}
\input{SFT-06_Fiber_products_for_Y}
\input{SFT-07_Algebras_for_Y}
\input{SFT-08_Fiber_products_for_X}
\input{SFT-09_Algebras_for_X}
\input{SFT-10_Fiber_products_for_homotopy}
\input{SFT-11_Algebras_for_homotopy}
\input{SFT-12_Composition}
\input{SFT-13_Independence}
\input{SFT-14_Circle_action}
\input{SFT-appendix}

\section*{Acknowledgments}
We thank K. Ono, K. Fukaya, M. Akaho, K. Cieliebak and U. Frauenfelder
for useful comments and suggestions.
We also thank H. Hofer for pointing out a lack of a condition in the
definition of grouped multisection which we implicitly used in the proof
of the construction of its extension.

This work was supported by Grant-in-Aid for JSPS Research Fellow.

\input{SFT-bibliography}
Research Institute for Mathematical Sciences, Kyoto University, Kyoto, Japan
\\
suguru@kurims.kyoto-u.ac.jp
\end{document}

%% file: SFT-preamble.tex
\usepackage{amsmath,amsthm,amssymb,enumitem,tikz-cd,subsupscripts,nameref}

\makeatletter
  \def\@thmcountersep{.}
\makeatother

\theoremstyle{plain}
\newtheorem{thm}{Theorem}[section]
\newtheorem{lem}[thm]{Lemma}
\newtheorem{prop}[thm]{Proposition}
\newtheorem{cor}[thm]{Corollary}

\theoremstyle{definition}
\newtheorem{defi}[thm]{Definition}

\theoremstyle{remark}
\newtheorem{rem}[thm]{Remark}
\newtheorem{eg}[thm]{Example}

\newcommand{\tocong}{\stackrel{\cong}{\to}}
\newcommand{\bound}{\overline{\partial}}
\newcommand{\inj}{\hookrightarrow}

\newcommand{\Wedge}{{\textstyle\bigwedge}}

\newcommand{\ab}{\allowbreak}

\newcommand{\hph}{\hphantom}

\newcommand{\id}{\mathrm{id}}
\newcommand{\forget}{\mathop{\mathfrak{forget}}\nolimits}

\newcommand{\vol}{\mathrm{vol}}
\newcommand{\ev}{\mathrm{ev}}
\newcommand{\diam}{\mathrm{diam}}
\newcommand{\St}{\mathrm{St}}

\newcommand{\sign}{\mathop{\mathrm{sign}}\nolimits}
\newcommand{\Ker}{\mathop{\mathrm{Ker}}\nolimits}
\newcommand{\Image}{\mathop{\mathrm{Im}}\nolimits}

\newcommand{\supp}{\mathop{\mathrm{supp}}\nolimits}
\newcommand{\rank}{\mathop{\mathrm{rank}}\nolimits}
\newcommand{\Aut}{\mathop{\mathrm{Aut}}\nolimits}

\newcommand{\gl}{\mathop{\mathrm{gl}}\nolimits}
\newcommand{\Int}{\mathop{\mathrm{Int}}\nolimits}
\newcommand{\CZ}{\mathop{\mathrm{CZ\mathchar`-ind}}\nolimits}

\newcommand{\ind}{\mathop{\mathrm{ind}}\nolimits}
\newcommand{\dist}{\mathop{\mathrm{dist}}\nolimits}
\newcommand{\codim}{\mathop{\mathrm{codim}}\nolimits}

\newcommand{\Av}{\mathop{\mathrm{Av}}\nolimits}
\newcommand{\graph}{\mathop{\mathrm{graph}}\nolimits}

\newcommand{\DM}{\mathrm{D}}

\newcommand{\N}{\mathbb{N}}
\newcommand{\Z}{\mathbb{Z}}
\newcommand{\Q}{\mathbb{Q}}
\newcommand{\R}{\mathbb{R}}
\newcommand{\C}{\mathbb{C}}

\newcommand{\A}{\mathcal{A}}
\newcommand{\B}{\mathcal{B}}
\newcommand{\D}{\mathcal{D}}
\newcommand{\E}{\mathcal{E}}
\newcommand{\F}{\mathcal{F}}
\newcommand{\G}{\mathcal{G}}
\newcommand{\I}{\mathcal{I}}
\newcommand{\J}{\mathcal{J}}
\newcommand{\K}{\mathcal{K}}
\newcommand{\M}{\mathcal{M}}
\renewcommand{\S}{\mathcal{S}}
\newcommand{\U}{\mathcal{U}}
\newcommand{\V}{\mathcal{V}}
\newcommand{\W}{\mathcal{W}}
\newcommand{\X}{\mathcal{X}}
\newcommand{\Y}{\mathcal{Y}}

\allowdisplaybreaks[1]

\title{Construction of general symplectic field theory}
\author{Suguru ISHIKAWA}
\date{\today}

%% file: SFT-01_Introduction_Estimates.tex
%
%
%
%
\section{Introduction}
The aim of this paper is to provide a construction of symplectic field theory (SFT).
SFT is a theory of contact manifolds and symplectic manifolds
with cylindrical ends proposed by Eliashberg, Givental and Hofer in \cite{EGH00}.
It is a generalization of contact homology and Gromov-Witten invariant, and it is
constructed by counting the number of appropriate pseudo-holomorphic
curves in the symplectization of a contact manifold or a symplectic manifold
with cylindrical ends.
In general, we need perturbation to obtain transversality of moduli spaces
of pseudo-holomorphic curves,
and it was a difficult problem to carry out perturbation
with compatibility conditions required for the construction of the algebras.
To give a concrete and transparent proof of the construction,
Hofer, Wysocki and Zehnder
developed the theory of polyfold (\cite{HWZ07I}-\cite{HWZ10II}).
However, they have not yet published a complete proof of
the construction of SFT.
There were various other attempts to overcome this difficulty in special cases.
For example, cylindrical contact homology of some three-dimensional contact manifolds
was constructed by Bao and Honda \cite{BH15} and Hutchings and Nelson \cite{HN15}.
Recently, Contact homology was constructed by Pardon \cite{Pa15}
and Bao and Honda \cite{BH16} independently.
However, the general SFT has not yet been fully constructed.

The main result of this paper is construction of SFT in full generality.
\begin{thm}\label{main result}
For each closed contact manifold $(Y, \xi)$ and each finite subset
$\overline{K}^0 \subset H_\ast(Y, \Q)$, we can define
SFT cohomology $H^\ast_{\mathrm{SFT}}(Y, \xi, \overline{K}^0)$,
rational SFT cohomology $H^\ast_{\mathrm{RSFT}}(Y, \xi, \overline{K}^0)$ and
contact homology $H^\ast_{\mathrm{CH}}(Y, \xi, \overline{K}^0)$
as invariants of $(Y, \xi, \overline{K}^0)$.
\end{thm}
In fact, we construct generating functions
defined in \cite{EGH00} for contact manifolds and symplectic manifolds with
cylindrical ends and prove all of their properties explained in \cite{EGH00}.

We also deal with Bott-Morse case (see Section \ref{asymptotic estimates} for
the definition of the Bott-Morse condition).
Some easy cases of Bott-Morse case was studied by Bourgeois in \cite{Bo02}.
We use the chain complex of triangulation of the space of periodic orbits instead of
Morse chain complex used in \cite{Bo02}.
Constructing SFT by a Bott-Morse contact form, we can calculate
the SFT cohomology of a contact manifold with $S^1$-action generated
by the Reeb vector field. For example, we can prove the following.
\begin{thm}\label{H vanishes}
Assume that $(Y, \xi)$ admits a contact form $\lambda$ whose Reeb flow defines
a locally free $S^1$-action on $Y$.
We also assume that all cycles in $\overline{K}^0$ are $S^1$-invariant.
Let $\overline{P}$ be the space of non-parametrized periodic orbits.
Then $H^\ast_{\mathrm{SFT}}(Y, \xi, \overline{K}^0)$ is the algebra
generated by $H_\ast(\overline{P}; \R)$, $H^\ast_c (\overline{P}; \R)$
and the variables $t_x$ $(x \in \overline{K}^0)$, $\hbar$ with the product defined by the
following commutative relations:
all variables are super-commutative except
\[
[p_c, q_\alpha] =\langle c, \alpha \rangle \hbar
\]
for all $c \in H_\ast(\overline{P}; \R)$ and $\alpha \in H^\ast_c (\overline{P}; \R)$,
where we denote the elements corresponding to $c$ or $\alpha$ by
$p_c$ or $q_\alpha$.
\end{thm}

We use the Kuranishi theory of Fukaya and Ono.
It is one of the general techniques to overcome the transversality problem
and it was first used
in \cite{FO99} for the construction of Gromov-Witten invariant and
Hamiltonian Floer Homology of symplectic manifolds.
We mainly follow the argument of \cite{FO99}.

We explain the new features of this paper briefly.
First we recall the general way to construct a Kuranishi neighborhood of a point in
a moduli space.
For example, consider a point $p = (\hat \Sigma, z, u)$ in the moduli space of
stable curves in a closed symplectic manifold $(M, \omega)$ with a compatible almost
complex structure $J$.
For simplicity, assume that the domain curve $(\hat \Sigma, z)$ is stable
and the automorphism group of $p$ is trivial.
Let $X$ be the deformation space of the domain curve $(\hat \Sigma, z)$.
For each $a = (\hat \Sigma_a, z_a) \in X$,
we construct a approximate solution $u_a$
of $J$-holomorphic equation, and consider the equation as a Fredholm map $F_a$
from $W^{1,p}(\hat \Sigma_a, u_a^\ast TM)$ to $L^p(\hat \Sigma_a,
\Wedge^{0,1} T^\ast \hat \Sigma_a \otimes_\C u_a^\ast TM)$, where $p > 2$.
We construct a finite vector space $E$ and a family of linear maps
$\lambda_a : E \to L^p(\hat \Sigma_a,
\Wedge^{0,1} T^\ast \hat \Sigma_a \otimes_\C u_a^\ast TM)$ which makes
each Fredholm map $F_a^+ = F_a \oplus \lambda_a :
W^{1,p}(\hat \Sigma_a, u_a^\ast TM) \oplus E \to L^p(\hat \Sigma_a,
\Wedge^{0,1} T^\ast \hat \Sigma_a \otimes_\C u_a^\ast TM)$ transverse to
zero.
Define $V = \bigcup_{a \in X} F_a^{-1}(0)$.
Then the zero set of the projection $s : V \to E$ is a neighborhood of $p$.
Roughly speaking, $(V, E, s)$ defines a Kuranishi neighborhood of $p$.

\begin{figure}
\centering
\includegraphics[width= 250pt]{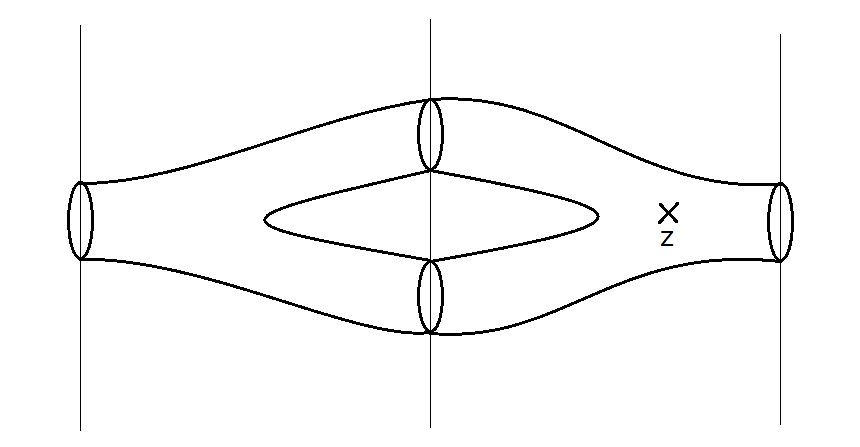}
\caption{$(\Sigma^0, z^0)$}\label{(Sigma0, z0)}
\includegraphics[width= 250pt]{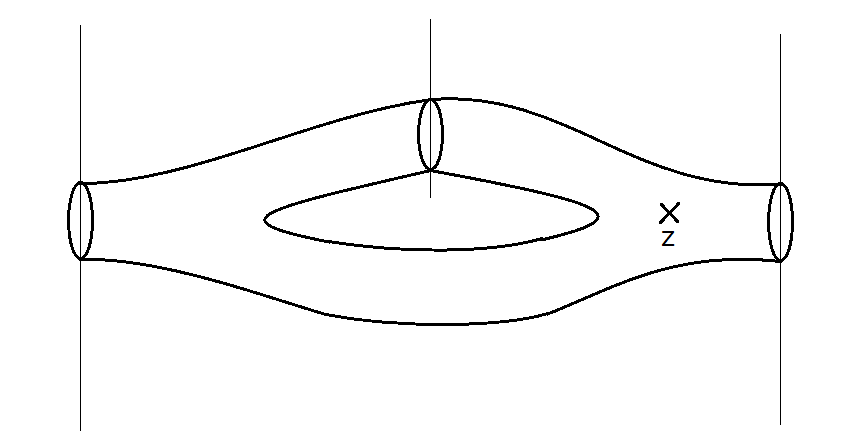}
\caption{$(\Sigma^1, z^1)$}\label{(Sigma1, z1)}
\end{figure}
For the construction of SFT, we count the $J$-holomorphic curves
in the symplectization $Y \times \R$ of a closed contact manifold $Y$.
Hence we consider the case of $M = Y \times \R$.
For example, consider the holomorphic building $(\Sigma^0, z^0, u^0)$
whose domain curve $(\Sigma^0, z^0)$ is as in Figure \ref{(Sigma0, z0)}.
(Holomoprhic buildings are the elements in the compactification of the space of
$J$-holomorphic curves. See Section \ref{space of holomorphic buildings}
for its definition.)
In the neighborhood of its domain curve $(\Sigma^0, z^0)$,
there is a curve like $(\Sigma^1, z^1)$ in Figure \ref{(Sigma1, z1)}.
However, we cannot consider the equation of $J$-holomorphic curves for
the curves like $(\Sigma^1, z^1)$ since they do not have floor structure.
This problem happens because we only consider the deformation of the domain curve
and ignore the deformation of the target space.
Therefore in this case, we need to use not the deformation space $X$
of the domain curve
but the space which parametrizes the deformation of the domain curve and
the deformation of the target space $\R \times Y$ simultaneously.
In Section \ref{construction of Kuranishi}, we define such a parameter space, and
construct an approximate solution and a Fredholm map for each of its points.

For the construction of the counterpart of chain homotopy in SFT, we need to treat
the space of $J$-holomoprhic curves in $1$-parameter family of symplectic manifolds
with cylindrical ends.
For a disjoint curve, we need to use the perturbation induced by the perturbations
for the connected components.
In the case of $1$-parameter family, this implies that the zero set of
the perturbed section for a disjoint curve is the fiber product of those
for the connected components over the parameter space.
However, in general, we cannot make the projections from the zero sets of the perturbed
sections for the connected components to the parameter space submersive,
which implies that the induced section
for the disjoint curve does not satisfy the transversality condition.
To overcome this problem, we use continuous family of perturbations.
(See Section \ref{continuous family of multisections} for its definition.)
It is a technique used in \cite{FOOO10} and \cite{FOOO11}.
Roughly speaking, instead of counting the number of zeros of the perturbed section,
we consider the perturbed section of the product of the moduli space with some
finite vector space and use the average of the number of zeros over the vector space.
If we use the product with appropriate vector spaces, then we can make the projection
from the zero sets to the parameter space submersive.

In the theory of Kuranishi structure, the smoothness of the Kuranishi structre is
one of its difficult part.
If we restrict on $0$- and $1$-dimensional Kuranishi spaces, then
often we do not need to consider the smoothness, but to use continuous
family of multisections, we cannot avoid this problem.
The difficulty is due to the fact that we need to use different Banach spaces for
different domain curves.
If the diffeomorphism type of the domain curve does not change, then the smoothness is
easy to prove since we can use the same Banach space by using diffeomorphisms.
However, if the diffeomorphism type changes, then we cannot identify the Banach spaces.
Hence we need to define artificially the smooth structure and prove the smoothness
of maps in Kuranishi theory (embeddings and evaluation maps).
Fukaya, Oh, Ohta and Ono treated this problem briefly in \cite{FOOO09II}, and
they explained the details of the argument in \cite{FOOO16}.
The key point is the following elementary fact:
If a continuous function $f$ on $\R$ is continuously differentiable on the complement of
a point, and the differential has a limit at this point, then $f$ is continuously
differentiable on the whole of $\R$.
In particular, we can prove the smoothness of $f$
if we check that the norm of its differentials converge to zero at this point.
This implies that it is enough to prove the convergence of the differentials at the strata
where the diffeomorphism type of the domain curve changes.
They proved the convergence
by estimating approximating solutions appearing in Newton's method.
We also prove the smoothness by estimating the limit of the differentials,
but we prove these estimates by
using the estimates of the implicit functions which define the solutions.
Using an appropriate family of identifications of the domain curves,
we estimate the differentials of the implicit functions
by direct calculation (Lemma \ref{estimates of implicit function 0} and
Corollary \ref{estimates of implicit function}).
Once we get the estimates of the implicit functions,
we can prove the estimates of the norm of the differentials of solutions
by Proposition \ref{asymptotic phi} and Corollary \ref{asymptotic estimates of phi}.

Another new feature of this paper is an improvement of the Kuranishi theory.
In the usual Kuranishi theory, the notion of good coordinate system is not
compatible with the product.
Hence usually, for the product space, we need to reconstruct the Kuranishi structures
of the factors from the good coordinate systems and
again construct a good coordinate system from the product of
the new Kuranishi structures.
Furthermore, we need to take care of the order of the product for the product of
more than two factors.
To avoid these complexities, we introduce the new notions of pre-Kuranishi structure and
its weakly good coordinate system.
Roughly speaking, a good coordinate system uses a total order,
but a weakly good coordinate system uses a meet-semilattice.
Similarly to meet-semilattice,
weakly good coordinate system is compatible with product, and we can directly use
their product for the product space.
This simplifies the construction of the algebra.
(See Section \ref{theory of Kuranishi structure}
for pre-Kuranishi sturcture and weakly good coordinate system.)

Finally, we explain about Bott-Morse case.
Bourgeois used Morse function on the space of periodic orbit, but instead,
we triangulate the space of periodic orbit and use the chain complex of the
simplicial complex.
Using this chain complex, we treat the most general case where
bad orbits appear as a subcomplex of the space of periodic orbits.
To construct the algebras by counting intersection numbers with simplices,
we need to use correction terms which correspond to cascades in \cite{Bo02}.
Since the algebra of SFT is more complicated than that of Contact homology,
the correction terms are also complicated.
Hence we need to solve algebraic equations to define appropriate correction terms.
(See Section \ref{algebra for correction}.)

In Bott-Morse case, we need to use the fiber product of pre-Kuranishi spaces
over an orbifold.
For example, we need to consider the fiber products with the diagonal
$\Delta_{\overline{P}}$ in $\overline{P} \times \overline{P}$,
where $\overline{P}$ is the space of non-parametrized periodic orbits.
We treat $\Delta_{\overline{P}}$ not as a suborbifold of $\overline{P} \times \overline{P}$
but as a simplicial complex in $\overline{P} \times \overline{P}$.
(See Definition \ref{def of fiber product of pre-Kuranishi over orbifold} and
Example \ref{fiber prod of pre-Kuranishi with diagonal of orbifold}.)
Although the fiber product of Kuranishi spaces over a manifold was treated before,
this paper is the first which treats the case of orbifold.

We briefly explain the outline of this paper.
First we investigate the local behavior and asymptotic behavior of pseudo-holomorphic
curves in Section \ref{local asymptotic estimates}.
Using them, we define the topology of the moduli space of holomorphic buildings and
prove its topological properties in Section \ref{space of holomorphic buildings}.
Next in Section \ref{theory of Kuranishi structure}, we recall the general theory of
Kuranishi structure and introduce the notions of pre-Kuranishi structure and
its weakly good coordinate system.
In Section \ref{construction of Kuranishi}, we construct a basic pre-Kuranishi structure
of the moduli space of holomorphic buildings.
In Section \ref{fiber prod}, we construct various fiber products of the basic pre-Kuranishi
spaces and construct their compatible multisections.
In this section, we also explain about how to treat the bad orbits.
Defining the orientations of the fiber products, we construct their virtual fundamental
chains, and using them, we construct the algebra.
In Section \ref{case of X} to \ref{composition},
we consider the cases of a symplectic manifold with cylindrical ends,
its $1$-parameter version and the composition of two symplectic cobordisms.
Using them, we prove that the algebras are invariants of contact manifolds in
Section \ref{independence}.
Finally in Section \ref{S^1 action}, we consider the calculation of the SFT cohomology of
contact manifolds with the $S^1$-action generated by the Reeb vector field.

\section{Local estimates and asymptotic estimates}\label{local asymptotic estimates}
Let $(Y,\lambda,J)$ be a triple which consists of
a closed ($2n-1$)-dimensional manifold $Y$ and a contact form $\lambda$,
and a compatible almost complex structure $J$ of $\xi = \Ker \lambda$.
$\xi$ has a symplectic structure given by $d\lambda$, and compatibility of $J$
means $d\lambda(\cdot, J \cdot)$ is a hermitian metric on $\xi$.

We denote by $R_\lambda$ the Reeb vector field of $\lambda$, which is defined by
$\lambda(R_\lambda) =1$ and $i_{R_\lambda}d\lambda = 0$.
We say a loop $\gamma : S^1 \to Y$ is a periodic orbit of period $L = L_\gamma > 0$
if it satisfies $\partial_t \gamma(t) = L R_\lambda(\gamma(t))$.
We note that the period of a periodic orbit $\gamma$ can be expressed as
$L = \int_{S^1} \gamma^\ast \lambda$.
We denote the space of all periodic orbits by $P = P_Y \subset C^\infty(S^1, Y)$.
$S^1 = \R / \Z$ acts on $P$ by $(s \cdot \gamma) (t) = \gamma(t + s)$.
$\overline{P} = P/S^1$ is the space of non-parametrized periodic orbits.
For $L>0$, we denote by $P_L \subset P$ the subspace of periodic orbits with period $L>0$.
It is sometimes convenient to define $P_L$ for $L\leq 0$ by
\[
P_L = \{ \gamma \in C^\infty (S^1, Y); \partial_t \gamma - LR_\lambda(\gamma) = 0 \}.
\]
For example, $P_0 = Y$ is the space of constant loops.
(However, we do not count these loops as periodic orbits.)
 
Let $\hat Y := \R \times Y$ be the symplectization of $Y$.
The coordinate of its $\R$-component is denoted by $\sigma$.
We can extend the complex structure $J$ of $\xi$ to an almost complex structure of
$\hat Y$ by $J (\partial_\sigma) = R_\lambda$, which we still denote by $J$.
The hermitian metric $g$ of $T \hat Y$ is then defined by
$g(\cdot, \cdot) = (d\sigma \wedge \lambda + d\lambda) ( \cdot, J \cdot)$.
In this paper, we construct algebras by counting $J$-holomorphic curves
in manifolds of this type.
Note that if $u : \Sigma \to \hat Y$ is a $J$-holomorphic curve,
then its $\R$-translations $o_{\sigma_0} \circ u : \Sigma \to \hat Y$ are also
$J$-holomorphic, where $o_{\sigma_0} : \R \times Y \to \R \times Y$ ($\sigma_0 \in \R$)
are the translation maps
defined by $o_{\sigma_0}(\sigma, y) = (\sigma + \sigma_0, y)$.

Other symplectic manifolds we consider in this paper are symplectic manifolds
with cylindrical ends. (Sometimes these are called contact ends.)
A symplectic manifold $(X, \omega)$ has cylindrical ends if
there exist contact manifolds $(Y^\pm, \lambda^\pm)$, and
$X$ can be decomposed as
$X = (-\infty, 0] \times Y^- \cup Z \cup [0, \infty) \times Y^+$,
where $Z$ is a compact manifold with boundary $\partial Z = Y^- \coprod Y^+$,
and the symplectic form satisfies
$\omega|_{(-\infty, 0] \times Y^-} = d(e^\sigma \lambda^-)$
and $\omega|_{[0, \infty) \times Y^+} = d(e^\sigma \lambda^+)$.
An almost complex structure $J$ on $X$ is said to be compatible if
$\omega(\cdot, J \cdot)$ is a hermitian metric and the restriction of $J$ on $(-\infty, 0] \times Y^-$ and $[0, \infty) \times Y^+$ are obtained by some complex structures on $\xi^-$ and $\xi^+$ respectively as above.

Two energies of a $J$-holomorphic map $u : (\Sigma, j ) \to (\hat Y, J)$ from
a Riemann surface $(\Sigma, j)$ to $\hat Y$ are defined as follows.
One is
\[
E_{\hat \omega}(u) = \int_{\Sigma} u^\ast d\lambda
\]
and the other is
\[
E_{\lambda}(u) = \sup_{I \subset \R} \frac{1}{|I|}\int_{(\sigma\circ u)^{-1}(I)} u^\ast (d\sigma \wedge \lambda),
\]
where the sup is taken over all intervals $I\subset \R$, and $|I|$ is the length of $I$.
The original energy introduced by Hofer in \cite{Ho93} was
\[
\sup \{\int_\Sigma u^\ast d(\varphi \lambda); \varphi \in C^\infty(\R, [1/2, 1]), \varphi'\geq 0 \}.
\]
This is equivalent to $E_{\hat \omega}(u) + E_{\lambda}(u)$ up to constant factors.

We define the norm of the differential $du(z)$ by
\[
|du(z)|^2 = \frac{|du(z)\zeta|_g^2 + |du(z)j \zeta|_g^2}{|\zeta|_h^2},
\]
where $h$ is a hermitian metric on $\Sigma$ and $\zeta$ is a non-zero vector of $T_z\Sigma$.
This does not depend on $\zeta$ (but depends on $h$).
If $u$ is $J$-holomorphic, then
$\int_\Sigma |du|^2 \vol = \int_\Sigma u^\ast (d\sigma \wedge \lambda + d\lambda)$.
Decomposing the tangent space $T \hat Y$ as $T \hat Y = \R \partial_\sigma \oplus \R R_\lambda \oplus \xi$, we denote the $\xi$-component of $du$ by $d^\xi u$.
Then $E_{\hat \omega}$-norm of $u$ coincides with
$||d^\xi u||_{L^2}^2 = \int_{\Sigma} |d^\xi u|^2 \vol_{\Sigma}$.

The energies of a $J$-holomorphic map $u : (\Sigma, j) \to (X, J)$ are defined as follows.
One is
\[
E_{\hat \omega}(u) = \int_{\Sigma} u^\ast \hat\omega,
\]
where $\hat \omega$ is a (discontinuous) $2$-form defined by
$\hat\omega|_Z = \omega$, $\hat\omega|_{(-\infty,0] \times Y^-}= d\lambda^-$
and $\hat\omega|_{[0, \infty) \times Y^+} = d\lambda^+$.
Note that the integral is invariant by homotopy of $u$ with compact support
(or relative to the boundary $\partial \Sigma$).
The other energy is
\begin{align*}
E_{\lambda}(u) &= \max\biggl\{\sup_{I \subset (-\infty,0]} \frac{1}{|I|}
\int_{u^{-1}(I \times Y^-)} u^\ast (d\sigma \wedge \lambda^-),\\
&\quad \hph{\max\biggl(}
\sup_{I \subset [0, \infty)} \frac{1}{|I|}\int_{u^{-1}(I \times Y^+)}
u^\ast (d\sigma \wedge \lambda^+)\biggr\}.
\end{align*}

\subsection{Local estimates}
The local estimates of $J$-holomorphic curves given in this subsection are not new
and have been already written in various forms.
(See \cite{Ho93} for example.)
However, for the convenience of the subsequent sections,
we state and prove them.

We use the following notation.
For non-negative functions $A$ and $B$,
$A \lesssim B$ means there exists a constant $C>0$ such that $A \leq CB$.
$A \sim B$ means $A \lesssim B$ and $B \lesssim A$.

\begin{lem}\label{L^infty bound}
For any $C_0 > 0$, there exist $\delta >0$ and $C_1>0$ such that
any $J$-holomorphic map $u : B_r(0) \to \hat Y$
($B_r(0) \subset \C$ is a ball with radius $r>0$)
with energies $E_\lambda (u) \leq C_0$ and $E_{\hat \omega}(u) \leq \delta$
satisfies $r |du(0)| \leq C_1$.
\end{lem}
\begin{proof}
If this did not hold, there would exist a constant $C_0 >0$, a sequence $\delta_k \to 0$
and $J$-holomorphic maps $u_k : B_{r_k}(0) \to \hat Y$ such that
$E_\lambda(u_k) \leq C_0$, $E_{\hat \omega} (u_k) \leq \delta_k$
and $r_k |du_k(0)| \to \infty$.
The lemma below implies that we may assume 
$\sup_{B_{r_k}(0)} |du_k(0)| \leq 2 |du_k(0)|$ by changing the center of the ball.
Rescaling the domain if necessary, we may assume $|du_k(0)| =1$.
In this case, the assumption implies $r_k \to \infty$.
Further we may assume $\sigma \circ u_k (0) = 0$ by $\R$-translation.
Then some subsequence of $u_k$ uniformly converges to a $J$-holomorphic map
$u_\infty : \C \to \hat Y$ such that
$|du_\infty(0)| = 1$, $E_\lambda(u_\infty) \leq C_0$ and $E_{\hat \omega} (u_\infty) = 0$.

$E_{\hat \omega} (u_\infty) = 0$ implies that the image of $du_\infty$ is contained
in the integrable subbundle $\R \partial_\sigma \oplus \R R_\lambda \subset T \hat Y$.
Hence the image of $u_\infty$ is contained in one of its leaves.
Each leaf is written as the image of a $J$-holomorphic map
$\Phi : \C \to \hat Y$ given by $\Phi(s + \sqrt{-1}t) = (s, \tilde\gamma(t))$,
where $\tilde \gamma : \R \to Y$ is an integral curve of $R_\lambda$.
Hence $u_\infty$ has a lift $\tilde u_\infty : \C \to \C$ such that
$|d\tilde u_\infty(0)| = 1$ and $u_\infty = \Phi \circ \tilde u_\infty$.

$E_\lambda(u_\infty) \leq C_0$ implies
\begin{align*}
\int_{\tilde u_\infty^{-1}(I \times \R)} \tilde u_\infty^\ast(ds \wedge dt)
&= \int_{\tilde u_\infty^{-1}(I \times \R)} |d\tilde u_\infty|^2 dsdt \\
&\leq C_0 |I| < \infty
\end{align*}
for any interval $I \subset \R$, which is a contradiction since any non-constant
holomorphic function on $\C$ takes all values except at most one value.
\end{proof}

\begin{lem}[\cite{HV92}]
Let $W$ be a complete metric space, and let $\varphi : W \to \R_{\geq 0}$ be
a continuous non-negative function.
For any $x_0 \in W$ and $r_0>0$,
there exist a point $x_1 \in B_{2r_0}(x_0)$ and $0 < r_1 < r_0$ such that
\[
\sup_{B_{r_1}(x_1)} \varphi \leq 2 \varphi(x_1) \text{ and }
r_0 \varphi(x_0) \leq r_1 \varphi(x_1).
\]
\end{lem}


\begin{lem}\label{preannulus}
For any $C_0 >0$, $l \geq 1$ and $\epsilon >0$,
there exist some $\delta >0 $, $A >0$ and $L_0 >0$ such that
any $J$-holomorphic map
$u : [-A, T + A] \times S^1 \to \hat Y$ ($T \geq 0$ is arbitrary) with energies
$E_\lambda (u) \leq C_0$ and $E_{\hat \omega}(u) \leq \delta$ satisfies
\[
\dist_{C^l(S^1, \hat Y)} (o_{-\sigma_s} \circ u(s, \cdot), \bigcup_{|L| \leq L_0} P_L) < \epsilon \text{ for all } s \in  [0,T],
\]
where $\sigma_s = \sigma(u(s,0))$, and we regard $P_L$ as a subset of $C^l(S^1, \hat Y)$ by the embedding
$Y = \{ 0 \} \times Y \inj \hat Y$.
\end{lem}
\begin{proof}
Let $L_0 = 2C_1$ be the double of the constant of Lemma \ref{L^infty bound}.
Note that Lemma \ref{L^infty bound} implies that if $A> \frac{1}{2}$ then
$|du|_{L^\infty([-A + 1/2, T + A - 1/2] \times S^1)} \leq L_0$.

It is enough to prove the claim for $T=0$.
If this lemma did not hold,
there would exist some sequences $A_k \to \infty$ and $\delta_k \to 0$,
some constant $\epsilon >0$, and a sequence of $J$-holomorphic maps
$u_k : [-A_k, A_k] \times S^1 \to \hat Y$ such that
$E_\lambda(u_k) \leq C_0$, $E_{\hat \omega}(u_k) \leq \delta_k$ and
$\dist_{C^l(S^1, \hat Y)} (o_{-\sigma_k} \circ u_k(0, \cdot), \bigcup_{|L| \leq L_0} P_L)
\geq \epsilon$.
We may assume $\sigma_k = \sigma(u_k(0,0)) = 0$.
Then a subsequence of $u_k$ uniformly converges to
a $J$-holomorphic map $u_\infty : \R \times S^1 \to \hat Y$ such that
$E_{\hat \omega}(u_\infty) = 0$ and $|du_\infty|_{L^\infty(\R \times S^1)} \leq L_0$.

We can deduce as follows that there exists some constant $|L|\leq L_0$ and
some periodic orbit $\gamma \in P_L$ such that
$u_\infty(s,t) = (L_\gamma s, \gamma(t))$, which contradicts the assumption
$\dist_{C^l(S^1, \hat Y)} (u_k(0, \cdot) \bigcup_{|L| \leq L_0} P_L) \geq \epsilon$.

As in the proof of Lemma \ref{L^infty bound},
there exists an integral curve $\tilde \gamma : \R \to Y$ such that
the image of $u_\infty$ is contained in the image of the $J$-holomorphic map
$\Phi : \C \to \hat Y$ given by $\Phi(s + \sqrt{-1} t) = (s, \tilde \gamma (t))$.
If $u_\infty$ has a lift $\tilde u_\infty : \R \times S^1 \to \C$, then
$|d \tilde u_\infty|_{L^\infty(\R \times S^1)} < \infty$ implies $u_\infty$ is
a constant map.
(This is the case of $L = 0$.)
If $u_\infty$ does not have such a lift, then there exists $L \neq 0 \in \R$ such that
$u_\infty$ has a lift
\begin{align*}
\tilde u_\infty : \R \times S^1 \to \C/L\sqrt{-1} &\cong \R \times S^1\\
(Ls + \sqrt{-1} Lt) & \leftrightarrow (s, t)
\end{align*}
such that $(\tilde u_\infty)_\ast = 1$ on $\pi_1(\R \times S^1)$.
Since $\tilde u_\infty : \R \times S^1 \to \R\times S^1$ is a $J$-holomorphic map
such that $\sigma \tilde u_\infty(0,0) = 0$, this implies
$\tilde u_\infty(s,t) = (s, t + \theta)$ for some $\theta \in S^1$.
Hence $u_\infty(s,t) = (L s, \gamma(t))$,
where $\gamma(t) = \tilde \gamma(L(t + \theta)) : S^1 \to Y$.
The inequality $|du_\infty|_{L^\infty([0,A]\times S^1)} \leq L_0$ implies $|L| \leq L_0$.
\end{proof}

\begin{cor}\label{first annulus}
For any $C_0 >0$ and $\epsilon >0$,
there exist some $\delta >0 $, $A >0$ and $L_0 >0$ such that
for any $0\leq T\leq \infty$ and any $J$-hoomorphic map
$u : [-A, T + A] \times S^1 \to \hat Y$ with energies
$E_\lambda (u) \leq C_0$ and $E_{\hat \omega}(u) \leq \delta$,
there exists some $|L| \leq L_0$ such that $P_{L} \neq \emptyset$ and
\[
||\partial_t u - L R_\lambda(u) ||_{L^\infty ([0,T] \times S^1)} \leq \epsilon.
\]
\end{cor}


The case of a symplectic manifold $X = (-\infty, 0] \times Y^- \cup Z \cup [0, \infty) \times Y^+$ with cylindrical ends is similar.
\begin{lem}\label{L^infty bound for X}
For any $C_0 > 0$, there exist $\delta >0$ and $C_1>0$ such that
any $J$-holomorphic map $u : B_r(0) \to X$
with energies $E_\lambda (u) \leq C_0$ and $E_{\hat \omega}(u) \leq \delta$
satisfies $r |du(0)| \leq C_1$.
\end{lem}
\begin{proof}
If the claim did not hold, there would exist a constant $C_0 >0$, a sequence $\delta_k \to 0$ and
$J$-holomorphic maps $u_k : B_{r_k}(0) \to X$ such that
$E_\lambda(u_k) \leq C_0$, $E_{\hat \omega} (u_k) \leq \delta_k$ and
$r_k |du_k(0)| \to \infty$.
We may assume 
$\sup_{B_{r_k}(0)} |du_k(0)| \leq 2 |du_k(0)|$.
Rescaling the domain if necessary, we may assume $|du_k(0)| =1$.
In this case, the assumption implies $r_k \to \infty$.

Lemma \ref{L^infty bound} implies there exists a constant $R>0$ such that
every $u_k(B_R(0))$ intersects with $Z$.
Hence some subsequence of $u_k$ uniformly converges to a $J$-holomorphic map
$u_\infty : \C \to X$ such that
$|du_\infty(0)| = 1$, $E_\lambda(u_\infty) \leq C_0$ and $E_{\hat \omega} (u_\infty) = 0$.

Since $du_\infty|_{u_\infty^{-1}(Z)} \equiv 0$,
if the image of $u$ intersects with the interior of $Z$,
unique continuation theorem implies $u_\infty$ is a constant map,
which is a contradiction.
On the other hand, if the image of $u_\infty$ does not intersect with the interior of $Z$,
the same argument as in Lemma \ref{L^infty bound} leads to a contradiction.
\end{proof}

\begin{lem}\label{preannulus for X}
For any $C_0 >0$, $l  \geq 1$ and $\epsilon >0$,
there exist some $\delta >0 $, $A >0$ and $L_0 >0$ such that
any $J$-holomorphic map
$u : [-A, T + A] \times S^1 \to X$ with energies
$E_\lambda (u) \leq C_0$ and $E_{\hat \omega}(u) \leq \delta$ satisfies
\[
\dist_{C^l(S^1, X)} (u(s, \cdot), (-\infty, 0] \times \bigcup_{|L| \leq L_0} (P_{Y^-})_L \cup Z
\cup [0,\infty) \times \bigcup_{|L| \leq L_0} (P_{Y^+})_L) < \epsilon
\]
for all $s \in [0,T]$,
where we regard each point $(\sigma, \gamma) \in (-\infty, 0] \times
\bigcup_{|L| \leq L_0} (P_{Y^-})_L$ as a loop
$(\sigma, \gamma(t)) \in C^l(S^1, (-\infty, 0] \times Y^-) \subset  C^l(S^1, X)$,
each $x\in Z$ as a constant loop in $C^l(S^1, X)$,
and each $(\sigma, \gamma) \in [0, \infty) \times \bigcup_{|L| \leq L_0} (P_{Y^+})_L$
as a loop $(\sigma, \gamma(t)) \in C^l(S^1, [0, \infty) \times Y^+) \subset C^l(S^1, X)$.
\end{lem}
\begin{proof}
Let $L_0 = 2C_1$ be the double of the constant of Lemma \ref{L^infty bound for X}.
Let $A_0>0$ be the constant of Lemma \ref{preannulus} for $\hat Y^\pm$.
Then the claim holds if $u([-A_0, T + A_0] \times S^1)$ does not intersect with the interior of $Z$ by Lemma \ref{preannulus}.

It is enough to prove the claim for $T = 0$.
If it did not hold,
there would exist some sequences $A_k \to \infty$ and $\delta_k \to 0$, some constant
$\epsilon >0$, and a sequence of $J$-holomorphic maps
$u_k : [-A_k, A_k] \times S^1 \to \hat Y$ such that
$E_\lambda(u_k) \leq C_0$, $E_{\hat \omega}(u_k) \leq \delta_k$ and
\[
\dist_{C^l(S^1, X)} (u_k(s, \cdot), (-\infty, 0] \times \bigcup_{|L| \leq L_0} (P_{Y^-})_L \cup
Z \cup [0,\infty) \times \bigcup_{|L| \leq L_0} (P_{Y^+})_L) \geq \epsilon.
\]
Since each $u_k([-A_0, A_0] \times S^1)$ intersects with $Z$,
a subsequence of $u_k$ uniformly converges to
a $J$-holomorphic map $u_\infty : \R \times S^1 \to X$ such that
$E_{\hat \omega}(u_\infty) = 0$ and $|du_\infty|_{L^\infty(\R \times S^1)} \leq L_0$.

Since $du_\infty|_{u_\infty^{-1}(Z)} \equiv 0$,
if the image of $u_\infty$ intersects with the interior of $Z$,
unique continuation theorem implies $u_\infty$ is a constant map, which is a contradiction.
On the other hand, if the image of $u_\infty$ does not intersect with the interior of $Z$,
then the same argument as in Lemma \ref{preannulus} leads to a contradiction.
\end{proof}

\begin{rem}
In the above Lemma,
$(-\infty, 0] \times \bigcup_{0 < |L| \leq L_0} (P_{Y^-})_L$,
$X = (-\infty, 0] \times P_0^- \cup Z \cup [0, \infty) \times P_0^+$ and
$[0, \infty) \times \bigcup_{0 < |L| \leq L_0} (P_{Y^+})_L$ are disjoint closed subsets.
Hence if $\epsilon >0$ is sufficiently small,
then it is independent of $s \in [0, T]$ which of these three
$u|_{\{s\} \times S^1}$ is close to.
\end{rem}

The following lemmas are well known.
See \cite{Gr85} or \cite{Hu97} for example.
\begin{lem}[Removal of Singularities]
Any $J$-holomorphic map $u : D \setminus 0 \to \hat Y$
{\rm(}or $u : D \setminus 0 \to X${\rm)}
with $||du||_{L^2} < \infty$ can be extended uniquely to a $J$-holomorphic map
$u : D \to \hat Y$ {\rm(}or $u : D \to X$ respectively{\rm)}.
\end{lem}
\begin{lem}[Monotonicity Lemma]\label{monotonicity lemma}
There exist some $r_0>0$ and $C>0$ such that
for any compact Riemann surface $\Sigma$ with or without boundary,
any non-constant $J$-holomorphic map $u : \Sigma \to \hat Y$
{\rm(}or $u : \Sigma \to X${\rm)},
any point $z_0\in \Int \Sigma$ and any $0\leq r\leq r_0$, the following holds true.
If $u(\partial \Sigma) \cap B_r(u(z_0)) = \emptyset$ then
\[
||du||_{L^2(u^{-1}(B_r(z_0)))}^2 \geq C r^2.
\]
\end{lem}
\begin{lem}\label{L^infty diam}
For any disc $D_0 \Subset D$, there exist $\delta>0$ and
$C>0$ such that any $J$-holomorphic curve $u: D \to \hat Y$
{\rm(}or $u : D \to X${\rm)}
with $\diam \, u(D) \leq \delta$ satisfies
\[
||du||_{L^\infty(D_0)} \leq C \diam \, u(D).
\]
Similarly, if a $J$-holomorphic curve $u_0 : D \to \hat Y$
{\rm(}or $u_0 : D \to X${\rm)} is given,
then there exist $\delta > 0$ and $C >0$ such that for any $J$-holomorphic curve
$u: D \to \hat Y$ {\rm(}or $u : D \to X$ respectively{\rm)},
if $\dist_{L^\infty(D)}(u, u_0) \leq \delta$ then
\[
||du - du_0||_{L^\infty(D_0)} \leq C \dist_{L^\infty(D)} (u, u_0).
\]
\end{lem}
%

\subsection{Asymptotic estimates}\label{asymptotic estimates}
To obtain asymptotic estimates of the ends of $J$-holomorphic curves, we need to
assume that the contact form satisfies the following condition.
Recall that $P \subset C^\infty(S^1, Y)$ is the space of (parametrized) periodic orbits
of the Reeb flow of $(Y, \lambda)$.
Let $\ev_t : P \to Y$ be the evaluation map at $t \in S$ defined by
$\ev_t \gamma = \gamma(t)$.
\begin{defi}\label{def of Bott-Morse}
For each periodic orbit $\gamma \in P$, we define an $L^2$ self-adjoint operator
$A_\gamma : W^{1,2}(S^1, \gamma^\ast T \hat Y) \to
L^2(S^1, \gamma^\ast T \hat Y)$ by
\[
A_\gamma \xi = J(\gamma) (\nabla_t \xi - L_\gamma \nabla_\xi R_\lambda(\gamma)),
\]
where we regard $\gamma$ as an element of
$C^\infty(S^1, \{0 \} \times Y) \subset C^\infty (S^1, \hat Y)$.
We say $(Y, \lambda)$ satisfies the Bott-Morse condition
(or $(Y, \lambda)$ is Bott-Morse)
if $P \subset C^\infty(S^1 ,Y)$ is a countable union of closed manifolds,
and every operator $A_\gamma$ satisfies
$\Ker A_\gamma = \R \partial_\sigma \oplus T_\gamma P$.
This condition can be stated by using the linearization of the Reeb flow
$\varphi_t^\lambda : Y \to Y$ as
\[
\Ker ((\varphi^{\lambda}_{L_\gamma})_\ast -1 : T_{\gamma(0)} Y \to T_{\gamma(0)} Y)
= T_{\gamma(0)} \ev_0 P_{L_\gamma}
\]
for all periodic orbits $\gamma \in P$.
Note that the Bott-Morse condition implies that
each $P_{\leq L_0} = \coprod_{0< L \leq L_0} P_L$ consists of
finite closed manifolds.
We say $(Y, \lambda)$ satisfies the Morse condition
if it satisfies the Bott-Morse condition and $\overline{P}$ consists of discrete points.
Note that in this case, $\dim \Ker A_\gamma = 2$ for all $\gamma \in P$.
\end{defi}
The above definition of Bott-Morse condition is more natural than that given
in \cite{Bo02} and \cite{BEHWZ03}.
(Their definition assumes another condition.)

In this paper, we always assume $(Y,\lambda)$ is Bott-Morse.
Under this condition,
we can prove more strict estimates on the curves appearing in
Corollary \ref{first annulus}.
\begin{prop}\label{second annulus}
Let $L \in \R$ be a constant such that $P_L \neq \emptyset$.
Then there exist constants $\epsilon > 0$, $\kappa> 0$ and $C>0$ such that
the following holds true.
For any $0 < T \leq \infty$ and any $J$-holomorphic map
$u : [0,T] \times S^1 \to \hat Y$ such that
$||\partial_t u - L R_\lambda(u) ||_{L^\infty([0,T] \times S^1)} \leq \epsilon$,
there exists $(b, \gamma) \in \R \times P_L$ such that
\[
\dist (u(s,t), (Ls + b, \gamma(t)) )
\leq C (e^{-\kappa s} + e^{-\kappa (T-s)})
|| \partial_t u - L R_\lambda(u) ||_{L^\infty([0,T] \times S^1)}
\]
on $[0, T] \times S^1$.
\end{prop}
A similar estimate was proved in \cite{BEHWZ03} under their Bott-Morse condition.
If $T= \infty$ and $L >0$, we say $u$ is positively asymptotic to
a periodic orbit $\gamma \in P_L$.
If $T= \infty$ and $L<0$, we say $u$ is negatively asymptotic to a periodic orbit
$\gamma(-t) \in P_{|L|}$.
In this case, using a biholomorphism $(s,t) \mapsto (-s, -t)$,
we usually consider $u$ as a $J$-holomorphic map
$u : (-\infty, 0]\times S^1 \to \hat Y$ such that
$\lim_{s \to -\infty} u(s,t) = \gamma(-t)$.

This proposition and Corollary \ref{first annulus} imply the following.
\begin{cor}\label{third annulus}
For any constants $C_0 > 0$ and $\epsilon >0$,
there exist $\delta > 0$, $\kappa > 0$, $A > 0$ and $L_0 > 0$
such that the following holds true.
For any $0 \leq T \leq \infty$ and any $J$-holomorphic curve
$u : [-A, T + A] \times S^1 \to \hat Y$ with energies
$E_\lambda(u) \leq C_0$ and $E_{\hat \omega}(u) \leq \delta$,
there exists $L \in \R$ and $(b, \gamma) \in \R \times P_L$ such that
$|L| \leq L_0$ and
\[
\dist ( u(s,t), (Ls + b, \gamma(t)) ) \leq \epsilon ( e^{-\kappa s} + e^{-\kappa (T - s)})
\]
for all $(s, t) \in [0, T] \times S^1$.
\end{cor}

\begin{rem}
The proof below implies that the constant $\kappa > 0$ in Proposition
\ref{second annulus}
can be taken arbitrary close to the minimum of the absolute values of the non-zero
eigenvalues of $A_\gamma$ ($\gamma \in P_L$).
(Instead, we need to take small $\epsilon > 0$.)
Note that in Corollary \ref{third annulus}, $L_0 > 0$ is determined by $C_0 > 0$ and
$\epsilon > 0$, and is independent of $\delta > 0$, $\kappa > 0$, and $A > 0$.
Therefore, also in Corollary \ref{third annulus}, the constant $\kappa > 0$ can be taken
arbitrary close to the minimum of the absolute values of the non-zero eigenvalues of
$A_\gamma$ ($\gamma \in P_{\leq L_0}$).
\end{rem}

To prove the above proposition, we need to rewrite the equation of $J$-holomorphic
curves in a neighborhood of a periodic orbit.

For each coordinate $\phi : B_\epsilon^m(0) \inj P_L$ of $P_L$,
we take a family of open embeddings
$\psi_t : B_\epsilon^m(0) \times B^{2n-1-m}(0) \inj Y$ ($t \in S^1$)
such that $\psi_t(x,0) = \ev_t \phi(x)$ for all $x\in B_\epsilon^m(0)$.
(The existence of such a family is due to the orientability of $Y$.)

First we show that if $\eta : S^1 \to \hat Y$ is a loop such that
$\eta(0) = (\sigma, \psi_0 (x,y))$,
then
\begin{gather}
|y| \lesssim ||\partial_t \eta - L R_\lambda(\eta)||_{L^\infty(S^1)}\label{y estimate}\\
\dist_{C^1(S^1, \hat Y)}(\eta(t), (\sigma, \gamma(t)))
\lesssim ||\partial_t \eta - L R_\lambda(\eta)||_{L^\infty(S^1)}\label{C^1 estimate}
\end{gather}
where $\gamma(t) = \ev_t \phi(x)$.
(\ref{y estimate}) is because
\begin{align*}
|y| &\sim \dist(\pi_Y \circ \eta(0), \ev_0 P_L)\\
&\lesssim \dist (\pi_Y \circ \eta(0), \varphi^\lambda_L (\pi_Y \circ \eta(0)))
\text{ (by the Bott-Morse condition)} \\
&\lesssim ||\partial_t (\varphi^\lambda_{-Lt}(\pi_Y \circ \eta(t)))||_{L^\infty(S^1)}\\
&\lesssim ||\partial_t \eta - L R_\lambda(\eta)||_{L^\infty(S^1)},
\end{align*}
where $\pi_Y : \hat Y = \R \times Y \to Y$ is the projection.
(\ref{C^1 estimate}) is because
\begin{align*}
\dist_{C^1(S^1, \hat Y)}(\eta(t), (\sigma, \gamma(t)))
&\leq \dist_{C^1(S^1, \hat Y)}(\eta(t), (1 \times \varphi^\lambda_{Lt}) \circ \eta(0))\\
&\quad + \dist_{C^1(S^1, \hat Y)}((1 \times \varphi^\lambda_{Lt}) \circ \eta(0),
(\sigma, \gamma(t)))\\
&\sim
\dist_{C^1(S^1, \hat Y)}((1 \times \varphi^\lambda_{-L_t}) \circ \eta(t), \eta(0))\\
&\quad + \dist_{\hat Y}(\eta(0), (\sigma, \gamma(0)))\\
&\lesssim ||\partial_t ((1 \times \varphi^\lambda_{-Lt}) \circ \eta (t))||_{L^\infty(S^1,
\hat Y)} + |y|\\
&\lesssim||\partial_t \eta - L R_\lambda(\eta)||_{L^\infty(S^1)}.
\end{align*}

Define a family of smooth maps $\hat\psi_{s,t} : \R \times B_\epsilon^m(0) \times B^{2n-1-m}(0) \inj \R \times Y$
($(s,t) \in \R \times S^1$) by
$\hat\psi_{s,t} (\sigma, x, y) = (L s + \sigma, \psi_t(x, y))$.

Assume a smooth map $u : I \times S^1 \to \hat Y$ satisfies
$\pi_Y \circ u (I \times \{ t \} ) \subset \mathrm{Im} \psi_t$ for all $t \in S^1$.
Then $u$ can be written as
$u(s,t) = \hat \psi_{s,t} (v(s,t))$, where
$v : I \times S^1 \to \R \times B_\epsilon^m(0) \times B^{2n-1-m}(0)$
is a smooth function.

We regard $N_0 = \R \oplus \R^m \oplus 0^{2n - m - 1} \subset \R^{2n}$ as a subspace of $W^{1,2}(S^1, \R^{2n})$
consisting of constant functions.
Then (\ref{C^1 estimate}) implies that there exists $z^0_s \in N_0$ for each $s\in I$ such that
\begin{equation}
||v|_{\{s\} \times S^1} - z^0_s||_{W^{1,2}(S^1)} \lesssim ||\partial_t u - L R_\lambda(u)||_{L^\infty (\{s\} \times S^1)}\label{pre N_1 estimate}
\end{equation}

The equation
\[
(\partial_s u - L\partial_\sigma) + J(u) (\partial_t u - L R_\lambda(u)) = 0
\]
of $J$-holomorphic curve for $u$
is equivalent to the following equation of $v$.
\begin{align*}
\partial_s v +	 ((\hat\psi_{s,t})_\ast)^{-1}& J(\hat\psi_{s,t}(v)) (\hat\psi_{s,t})_\ast \partial_t v\\
&+ ((\hat\psi_{s,t})_\ast)^{-1} J(\hat\psi_{s,t}(v)) (\partial_t \hat\psi_{s,t}(v) - L R_\lambda(\hat\psi_{s,t}(v))) = 0
\end{align*}
Note that this equation is also $\R$-translation invariant, that is,
if $v$ is a solution of the equation then $v(s,t) + (b,0)$ also satisfies the equation for any $b\in \R$.

We regard the solution $v$ as a map
$v : I \to C^\infty(S^1, \R \times B_\epsilon^m(0) \times B^{2n-1-m}(0))
\ab (\subset C^\infty(S^1, \R^{2n}))$.
Then the above equation has the following form.
\[
\partial_s v + F(v) = 0,
\]
where $F : W^{1,2}(S^1, \R^{2n}) \to L^2(S^1, \R^{2n})$ is a smooth Fredholm map
(more precisely, the domain of $F$ is an open neighborhood of
$0 \in W^{1,2}(S^1, \R^{2n})$)
which satisfies the following conditions:
\begin{itemize}
\item
$F$ maps $W^{k+1,2}(S^1, \R^{2n})$ to $W^{k,2}(S^1, \R^{2n})$ ($k \geq 0$).
\item
$F(v + \sigma) = F(v)$ for any $\sigma \in \R \oplus 0^m \subset N_0$.
\item
For any $z \in \R \times B_\epsilon^m(0) \subset N_0$, $F$ satisfies
$F(z)=0$ and $\Ker DF(z) = N_0$
(This is exactly the Bott-Morse condition.)
\item
There exists a family of inner product $(g_t)_{t \in S}$ of the vector space $\R^{2n}$
which makes the operator $A = DF(0) : W^{1,2}(S^1, \R^{2n}) \to L^2(S^1, \R^{2n})$
$L^2$ self-adjoint.
(In this case, $g_t$ is the pull back of $g$ by $(1 \times \psi_t)_\ast$
at $0 \in \R \times B_\epsilon^m(0) \times B^{2n-1-m}(0)$.)
\end{itemize}

In the following, we denote by $\langle \cdot, \cdot \rangle$ and $|\cdot|$
the inner product and the norm of $L^2(S^1, \R^{2n})$ given by $g_t$ ($t \in S^1$)
respectively.
The norm of $W^{1,2}(S^1, \R^{2n})$ is equivalent to $|v^0| + |A v^1|$.

First note that (\ref{pre N_1 estimate}) implies
\begin{equation}
|A v(s)| \lesssim ||\partial_t u - L R_\lambda(u)||_{L^\infty (\{s \} \times S^1)}.
\label{A v}
\end{equation}

Next we estimate
\begin{align*}
\partial_s^2 \langle A v, A v \rangle
&= 4 \langle A^2 v, A^2 v \rangle
+ 6 \langle A (F(v) - DF(0) v), A^2 v \rangle \\
&\quad + 2|A (F(v) - DF(0) v)|^2
+ 2\langle \pi_1 (DF(v) - DF(0)) F(v), A^2 v \rangle.
\end{align*}
Let $\pi_{\R^m}$ be the second projection of $N_0 = \R \oplus \R^m$.
In the above equation,
\begin{align*}
|A(F(v) - DF(0)v)| &\lesssim (|\pi_{\R^m} v^0| + |A v|) |A^2 v|\\
|\pi_1 (DF(v) - DF(0))F(v)| &\lesssim (|\pi_{\R^m} v^0| + |A v|) |A^2 v|
\end{align*}
because
\begin{align*}
A(F(v) - DF(0)v)
&= A(F(\pi_{\R^m} v^0 + v^1) - F(\pi_{\R^m} v^0) - DF(\pi_{\R^m} v^0)v^1)\\
&\quad +A((DF(\pi_{\R^m} v^0) - DF(0))v^1)\\
&= A\int_0^1\int_0^1 D^2F(\pi_{\R^m} v^0 + \tau_1 \tau_2 v^1) \tau_1 v^1 \cdot v^1
d\tau_1 d\tau_2\\
&\quad + A \int_0^1 D^2F(\tau \pi_{\R^m} v^0) (\pi_{\R^m} v^0) \cdot v^1 d\tau,
\end{align*}
\begin{align*}
&\pi_1 (DF(v) - DF(0))F(v)\\
&= \pi_1 (DF(\pi_{\R^m} v^0 + v^1) - DF(0))(F(v) - F(v^0))\\
&= \pi_1 \int_0^1 D^2F(\tau_1 (\pi_{\R^m} v^0 + v^1)) d\tau_1 (\pi_{\R^m} v^0 + v^1)
\cdot \int_0^1 DF(v^0 + \tau_2 v^1) v^1 d\tau_2,
\end{align*}
and $D^2F$ satisfies
\[
||(D^2F)(v)\xi \cdot \eta||_{W^{k, 2}(S^1, \R^{2n})} \lesssim
\sum_{\substack{i, j \geq 1 \\ i + j = k + 2}}
||\xi||_{W^{i, 2}(S^1, \R^{2n})} ||\eta||_{W^{j, 2}(S^1, \R^{2n})}
\]
for all $k \geq 0$.
(This is because $F$ is a differential operator.)

Therefore, if $||\pi_{\R^m} v^0||_{L^\infty(I, N_0)}$ and
$||A\tilde v||_{L^\infty(I, L^2(S^1, \R^{2n}))}$ are sufficiently small
(this assumption is satisfied if $B_\epsilon^m(0)$ and
$||\partial_t u - L R_\lambda(u)||_{L^\infty(I\times S^1)}$ are sufficiently small),
then there exists $\epsilon \ll 1$ such that
\begin{align*}
\partial_s^2 \langle A v, A v \rangle
&\geq 4 |A^2 v|^2 - C (|\pi_{\R^m} v^0| + |A v|) |A^2 v|^2\\
&\geq (4- \epsilon) |A^2 v|^2\\
&\geq (4- \epsilon ) \kappa_0^2 |A v|^2
\end{align*}
for all $s \in I$, where $\kappa_0 > 0$ is the minimum of the absolute values of
the non-zero eigenvalues of $A$.

Therefore the lemma below (Lemma \ref{absolute annulus}) implies that if $I = [0,T]$ then
\begin{equation}\label{N_1 annulus}
|A v (s)| ^2 \leq ( e^{-\sqrt{4 - \epsilon} \kappa_0 s} + e^{-\sqrt{4 - \epsilon} \kappa_0 (T-s)}) ||A v||^2_{L^\infty(I, L^2(S^1, \R^{2n}))}.
\end{equation}
In particular,
\begin{equation}
||A v||_{L^1(I, L^2(S^1, \R^{2n}))}
\lesssim ||A v||_{L^\infty(I, L^2(S^1, \R^{2n}))} \label{L^1 Av}
\end{equation}
is an estimate uniform with respect to $|I|$.

Since $|\pi_0 F(v(s))| \lesssim |A v^1(s)|$,
the equation $\partial_s v^0 + \pi_0 F (v) = 0$ implies
\begin{equation}
|\partial_s v^0(s)| \lesssim |Av^1(s)|. \label{N_0 norm}
\end{equation}

(\ref{A v}), (\ref{L^1 Av}) and (\ref{N_0 norm}) implies
\begin{align}
||\partial_s v^0||_{L^1(I, N_0)} &\lesssim ||A v||_{L^1(I, L^2(S^1, \R^{2n}))}\notag \\
&\lesssim ||Av||_{L^\infty(I, L^2(S^1, \R^{2n}))} \notag \\
&\lesssim ||\partial_t u - L R_\lambda (u)||_{L^\infty(I \times S^1)}\label{N_0 variation}
\end{align}


Using the above argument, now we prove Proposition \ref{second annulus}.
\begin{proof}[Proof of Proposition \ref{second annulus}]
Suppose $\epsilon>0$ is sufficiently small and that
a $J$-holomorphic map $u : [0,T] \times S^1 \to \hat Y$ satisfies
$||\partial_t u - L R_\lambda(u) ||_{L^\infty([0,T] \times S^1)} \leq \epsilon$.
There exists a coordinate $\phi$ of $P_L$ such that
$\pi_Y u ([0,T] \times \{t\}) $ is contained in the image of $\psi_t$ for all $t \in S^1$
since inequality (\ref{N_0 variation}) implies the variation of $v^0$ on $[0,T]$ is small.
Equalities (\ref{A v}), (\ref{N_1 annulus}) and (\ref{N_0 norm}) imply that
for $z = v^0(T/2) \in N_0$,
\begin{align*}
|v^0 (s) - z| &= \int_{T / 2}^s |\partial_s v^0| |ds|\\
&\lesssim (e^{-\frac{1}{2}\sqrt{4 - \epsilon} \kappa_0 s}
+ e^{-\frac{1}{2}\sqrt{4 - \epsilon} \kappa_0 (T-s)})
|| \partial_t u - L R_\lambda(u) ||_{L^\infty([0,T] \times S^1)}.
\end{align*}
(\ref{A v}) and (\ref{N_1 annulus}) imply
\[
||A v||_{L^2(S^1, \R^{2n})} \lesssim (e^{-\frac{1}{2}\sqrt{4 - \epsilon} \kappa_0 s} + e^{-\frac{1}{2}\sqrt{4 - \epsilon} \kappa_0 (T-s)})|| \partial_t u - L R_\lambda(u) ||_{L^\infty([0,T] \times S^1)}.
\]
Combining the above two inequalities, we see
\[
||v (s) - z||_{W^{1,2}(S^1, \R^{2n})} \lesssim (e^{-\kappa s} + e^{-\kappa (T-s)})
|| \partial_t u - L R_\lambda(u) ||_{L^\infty([0,T] \times S^1)},
\]
where $\kappa = \frac{1}{2}\sqrt{4 - \epsilon} \kappa_0$.
Therefore, if $(b, \gamma) \in \R \times P_L$ corresponds to $z$,
that is, $z = (b,x) \in \R \times B^m_\epsilon(0)$ and $\phi(x) = \gamma \in P_L$,
then
\[
\dist (u(s,t), (Ls + b, \gamma(t)) )
\lesssim (e^{-\kappa s} + e^{-\kappa (T-s)})|| \partial_t u - L R_\lambda(u) ||_{L^\infty([0,T] \times S^1)}.
\]
\end{proof}

\begin{lem}\label{absolute annulus}
If a $C^2$-function $f : [a,b] \to \R$ satisfies $f''(s) \geq \kappa^2 f(s)$ then
\[
f(s) \leq e^{-\kappa (s-a)} f(a)_+ + e^{-\kappa (b-s)} f(b)_+,
\]
where $f(s)_+ = \max (f(s) , 0)$.
\end{lem}
\begin{proof}
Since $g(s) = f(s) - (e^{-\kappa (s-a)} f(a)_+ + e^{-\kappa (b-s)} f(b)_+)$ also satisfies
$g''(s) \geq \kappa^2 g(s)$, we may assume $f(a) \leq 0$ and $f(b)\leq 0$.
If $f$ attained a positive value at some point $s_1$,
then there would exist some $a < s_0 < s_1$ such that
$f(s_0) > 0$ and $f'(s_0) >0$.
However this and the assumption $f''(s) \geq \kappa^2 f(s)$
would imply $f$ is monotone increasing on $s \geq s_0$, 
which contradict the assumption $f(b) \leq 0$.
\end{proof}

The case of a symplectic manifold $X$ with cylindrical ends is covered by
Proposition \ref{second annulus}, Corollary \ref{third annulus},
and the following propositions.

\begin{prop}\label{second annulus for X}
There exist constants $\epsilon > 0$, $\kappa> 0$ and $C>0$ such that
the following holds true.
For any $0 < T \leq \infty$ and any $J$-holomorphic map
$u : [0,T] \times S^1 \to X$ such that
$||\partial_t u||_{L^\infty([0,T] \times S^1)} \leq \epsilon$,
there exists a point $x \in X$ such that
\[
\dist (u(s,t), x ) \leq C (e^{-\kappa s} + e^{-\kappa (T-s)})
||\partial_t u||_{L^\infty([0,T] \times S^1)}.
\]
on $[0, T] \times S^1$.
\end{prop}
The proof of this proposition is the same as that of Proposition \ref{second annulus}.
\begin{cor}\label{third annulus for X}
For any constants $C_0 > 0$ and $\epsilon >0$,
there exist $\delta > 0$, $\kappa > 0$, $A > 0$ and $L_0 > 0$
such that the following holds true.
For any $0 \leq T \leq \infty$ and any $J$-holomorphic curve
$u : [-A, T + A] \times S^1 \to X$ with energies
$E_\lambda(u) \leq C_0$ and $E_{\hat \omega}(u) \leq \delta$,
one of the following two occurs:
\begin{itemize}
\item
There exists a point $x \in X$ such that
\[
\dist ( u(s,t), x ) \leq \epsilon ( e^{-\kappa s} + e^{-\kappa (T - s)})
\]
for all $(s, t) \in [0, T] \times S^1$.
\item
There exists $L \neq 0 \in \R$ and $(b, \gamma) \in \R \times P_L$ such that
$|L| \leq L_0$ and
\[
\dist ( u(s,t), (Ls + b, \gamma(t)) ) \leq \epsilon ( e^{-\kappa s} + e^{-\kappa (T - s)})
\]
for all $(s, t) \in [0, T] \times S^1$.
\end{itemize}
\end{cor}


%% file: SFT-02_Holomorphic_buildings.tex
%
%

\section{The space of holomorphic buildings}\label{space of holomorphic buildings}
In this section, we study the compactification of the space of $J$-holomorphic curves
in the symplectization of a contact manifold or a symplectic manifold with cylindrical ends.
Compactification was studied by Bourgeois, Eliashberg, Hofer, Wysocki and Zehnder
in \cite{BEHWZ03},
and the curves appeared in the compactified space are called holomorphic buildings.

First we recall about holomorphic buildings, and next we explain the topology
of the compactified space.
For the later use, we adopt a different definition of the topology.
This would be the same as that of \cite{BEHWZ03},
but we prove the compactness and Hausdorff property independently.

\subsection{The case of the symplectization}
First we consider holomorphic buildings for the symplectization $\hat Y = \R \times Y$.
The domain curve of a holomorphic building is constructed as follows.

Let $(\check \Sigma, z \cup (\pm \infty_i))$ be a marked semistable curve
or a disjoint union of marked semistable curves.
$z = (z_i)$ and $(\pm \infty_i)$ are sequences of marked points.
See \cite{FO99} for the definition of marked semistable curve.
Assume that an integer $i(\alpha) \in \{1, 2, \dots, k\}$ is attached to each irreducible
component $\check \Sigma_\alpha$ of $\check\Sigma$
(we call this integer the floor of $\check\Sigma_\alpha$) and
\begin{itemize}
\item
the difference of the floors of any adjacent two components is $\leq 1$,
\item
the floor of the component which contains some of the marked points $-\infty_i$ is $1$
(the lowest floor) and
\item
the floor of the component which contains some of the marked points $+\infty_i$ is $k$
(the highest floor).
\end{itemize}

We can construct a new curve from $\check \Sigma$ by oriented blow up.
Oriented blow up is a local deformation defined as follows.
Oriented blow up at $0 \in D = \{ z \in \C; |z|<1\}$ is 
\[
\widetilde{D} = \{ (z,\theta) \in D \times S^1; z = |z|\theta\},
\]
and 
oriented blow up at a nodal point $(0, 0) \in D\cup D = \{ (x,y) \in D\times D; xy =0\}$
by $\varphi \in S^1$ is 
\[
D\widetilde{\cup}_\varphi D = \{ (x, \theta_x, y, \theta_y) \in \widetilde{D} \times
\widetilde{D}; xy= 0, \theta_x \theta_y = \varphi\}.
\]
$S^1 = \{(0, \theta); \theta \in S^1 \} \subset \widetilde{D}$ is called limit circle,
and $S^1 = \{(0, \theta_x, 0, \theta_y); \theta_x \theta_y = \varphi \}
\subset D\widetilde{\cup}_\varphi D$ is called joint circle.
These two circles are collectively called imaginary circles.
The domain curve $(\Sigma, z)$ of a holomorphic building is obtained by oriented blow up
of $(\check \Sigma, z)$ at the points $\pm \infty_i$ and all the nodal points
which join two components with different floors by some $\varphi \in S^1$.
We regard the curve $\Sigma$ as a topological space, and the complement of
its imaginary circles as an open smooth curve with a complex structure.
The topological space $\Sigma$ is compact.
Note that there exists a surjection $\Sigma \to \check \Sigma$ which collapses
the imaginary circles.
For each irreducible component $\check \Sigma_\alpha$ of $\check \Sigma$,
we denote its inverse image by $\Sigma_\alpha \subset \Sigma$ and call it
an irreducible component of $\Sigma$.
We say that the marked curve $(\Sigma, z)$ is connected if $\Sigma$ is connected
as a topological space, that is, if it is constructed from one semistable curve
(not from a disjoint union of several semistable curves).
We emphasize the difference between the notion of irreducible component and
connected component.
For example, two irreducible components of $\Sigma$ connected by a joint circle
are considered to be in the same connected component.
\begin{figure}
\centering
\includegraphics[width= 350pt]{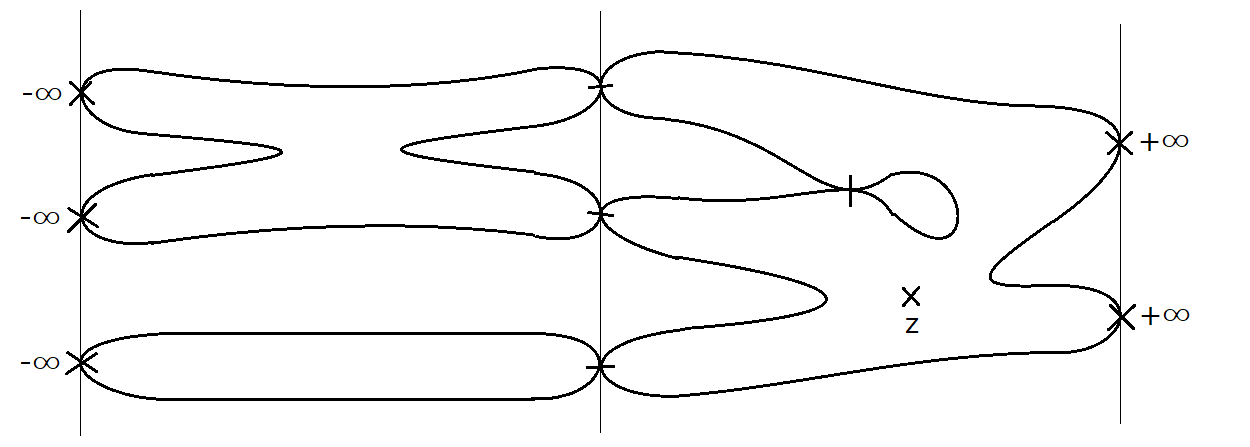}
\caption{$(\hat \Sigma, z \cup (\pm\infty_i))$}\label{(hatSigma, z)}
\end{figure}
\begin{figure}
\centering
\includegraphics[width= 350pt]{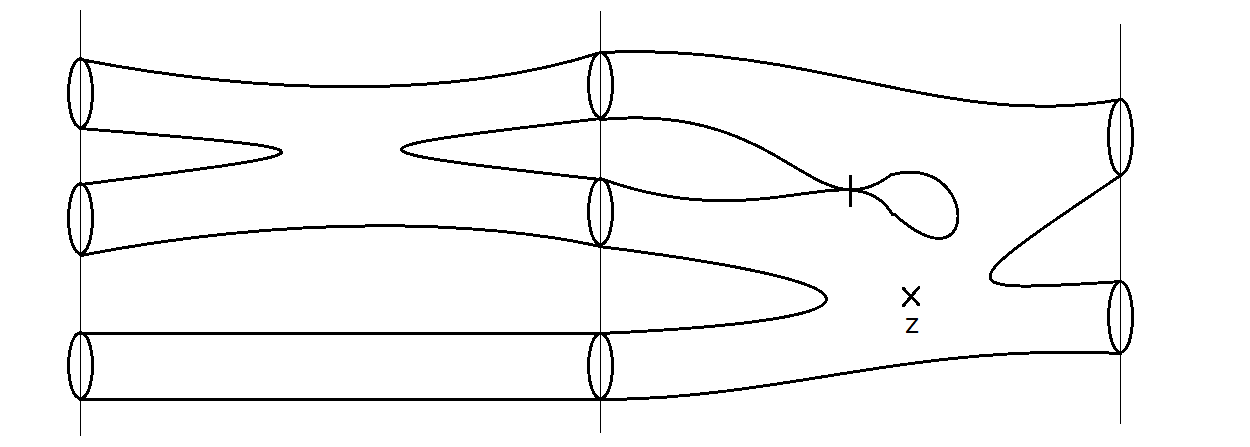}
\caption{$(\Sigma, z)$}\label{(Sigma, z)}
\end{figure}
\begin{defi}
A holomorphic building $(\Sigma, z, u, \phi)$ for $\hat Y$ consists of
\begin{itemize}
\item
a marked curve $(\Sigma, z)$ obtained from some marked semistable curve 
$(\check \Sigma, z \cup (\pm \infty_i))$ (or a union of marked semistable curves)
with a floor structure and some blowing up parameters $\varphi \in S^1$ as above,
\item
a continuous map $u : \Sigma \to (\overline{\R}_1 \cup \overline{\R}_2 \cup \dots
\cup \overline{\R}_k) \times Y$,
where $\overline{\R} = \{-\infty\} \cup \R \cup \{ + \infty\}$ is a compactification of
$\R$ (homeomorphic to an closed interval) and we identify
$+ \infty \in \overline{\R}_i$ and $-\infty \in \overline{\R}_{i+1}$,
and
\item
a family of coordinates $\phi_{\pm\infty_i} : S^1 = \R/ \Z \tocong S^1_{\pm\infty_i}$ of
limit circles, where $S^1_{\pm\infty_i}$ is the limit circle corresponding to
$\pm\infty_i \in \check \Sigma$.
\end{itemize}
which satisfy the following conditions:
\begin{itemize}
\item
$u(\Sigma_\alpha \setminus \coprod_{\text{imaginary circles}} S^1) \subset \R_{i(\alpha)} \times Y$ for each component $\Sigma_\alpha$.
\item
$u|_{\Sigma_\alpha \setminus \coprod S^1} : \Sigma_\alpha \setminus \coprod S^1 \to \R_{i(\alpha)} \times Y$ is $J$-holomorphic.
\item
$E_\lambda(u) <\infty$ and $E_{\hat \omega}(u) <\infty$,
where these energies are defined by
\begin{align*}
E_\lambda(u) &= \max_{1 \leq i \leq k} \sup_{I \subset \R_i} \frac{1}{I}
\int_{(\sigma\circ u)^{-1}(I)} u^\ast (d\sigma \wedge \lambda), \\
E_{\hat \omega}(u) &= \int_{\Sigma} u^\ast d\lambda.
\end{align*}
\item
$u$ is positively asymptotic to a periodic orbit $\gamma_{+\infty_i}
= \pi_Y \circ u \circ \phi_{+\infty_i} \in P$ at each $S^1_{+\infty_i}$,
and negatively asymptotic to a periodic orbit $\gamma_{-\infty_i}
= \pi_Y \circ u \circ \phi_{-\infty_i} \in P$ at each $S^1_{-\infty_i}$.
At every joint circle, $u$ is positively asymptotic to a periodic orbit
on the side of lower floor and negatively asymptotic to the same periodic orbit
on the side of higher floor.
\item
For each component $\check\Sigma_\alpha$, if $u|_{\Sigma_\alpha}$ is a constant map,
then $2g_\alpha + m_\alpha \geq 3$,
where $g_\alpha$ is the genus of $\check\Sigma_\alpha$ and $m_\alpha$ is the sum of
the numbers of marked points and imaginary circles in $\Sigma_\alpha$ and
nodal points which join $\Sigma_\alpha$ with the other components.
\item
An irreducible component is called a trivial cylinder if it is isomorphic to
$\overline{\R} \times S^1$ without any special points such that
the restriction of $u$ on this component is written as
$u(s, t) = (L_\gamma s + b, \gamma(t))$ for some $b\in \R$ and $\gamma \in P$.
The other irreducible components are called nontrivial components.
We assume that for each $i \in \{1,2, \dots, k\}$,
$i$-th floor $u^{-1}(\overline{\R}_i \times Y) \subset \Sigma$ contains
nontrivial components.
(We do not assume the same condition
for each floor of each connected component of $\Sigma$.)
\end{itemize}
We call $k$ the height of $(\Sigma, z, u, \phi)$.
\end{defi}
We say two holomorphic buildings $(\Sigma, z, u, \phi)$ and $(\Sigma', z', u', \phi')$ are
isomorphic if there exist
\begin{itemize}
\item
a biholomorphism $\varphi : \Sigma' \to \Sigma$
(this means $\varphi$ is a homeomorphism which maps each imaginary circle of $\Sigma'$
to a imaginary circle of $\Sigma$ and is biholomorphic on the outside of these circles)
and
\item
an $\R$-translation $\theta : \overline{\R}_1 \cup \overline{\R}_2 \cup \dots \cup
\overline{\R}_k \to \overline{\R}_1 \cup \overline{\R}_2 \cup \dots \cup \overline{\R}_k$
(this means $\theta$ is a map such that $\theta(\overline{\R}_i) \subset \overline{\R}_i$
and $\theta|_{\overline{\R}_i} (s)= s+ a_i$ for some $a_i\in \R$)
\end{itemize}
such that
\begin{itemize}
\item
$\varphi(z'_i) = z_i$ for all $i$,
\item
$u' = (\theta \times 1) \circ u \circ \varphi$, and
\item
$\varphi \circ \phi'_{\pm\infty_i} = \phi_{\pm\infty_i}$ for all $\pm\infty_i$.
\end{itemize}

We denote the space of all connected holomorphic buildings by
$\overline{\M}^0 = \overline{\M}^0(Y, \lambda, J)$,
and the space of all holomorphic buildings without trivial buildings by
$\overline{\M} = \overline{\M}(Y, \lambda, J)$,
where a trivial building in $(\Sigma, z, u, \phi)$ is a connected component of $\Sigma$
which consists of trivial cylinders only.

First we define the topology of $\overline{\M}^0$.
It is enough to define the neighborhoods of each point
$p_0 = (\Sigma_0, z_0, u_0, \phi_0) \in \overline{\M}^0$.
We consider a fibration $(\widetilde{P} \to \widetilde{X}, Z)$ consisting of
some deformations of the domain curve $(\Sigma_0, z_0)$, and
construct a map $\Psi : \widetilde{P} \to \widetilde{P}_0$.
Then the neighborhood of $p_0$ is defined by the set of holomorphic buildings whose
domain curves appear as a fiber $\widetilde{P}_a$ of $\widetilde{P}$ and
which are close to $u_0 \circ \Psi|_{\widetilde{P}_a}$ in $L^\infty$-norm modulo
$\R$-gluings.

Now we explain the details.
First we add marked points $z^+_0$ to $(\Sigma_0, z_0)$ to make
$(\Sigma_0, z_0\cup z_0^+)$ stable,
where $z_0 \cup z_0^+$ is a sequence of marked points obtained by placing the sequence
$z_0^+$ after $z_0$, and stableness of $(\Sigma_0, z_0\cup z_0^+)$ means that
the curve $(\check\Sigma_0, z_0\cup z^+_0 \cup (\pm\infty_i))$ is a stable curve.

The local universal family $(\widetilde{P} \to \widetilde{X}, Z \cup Z^+)$ of
$(\Sigma_0, z_0 \cup z_0^+)$ is defined by the oriented blow up of the local universal
family $(\check P \to \check X, Z \cup Z^+ \cup (Z_{\pm \infty_i}))$ of
the stable curve $(\check\Sigma_0, z_0\cup z^+_0 \cup (\pm\infty_i))$ at
$Z_{\pm \infty_i}$ and the set of nodal points corresponding to the nodal points of
$\check \Sigma_0$ which are blown up in $\Sigma_0$.

Oriented blow up of the local universal family is defined as follows.
For each nodal point of $\check \Sigma_0$, the fibration $\check P \to \check X$ is
locally equivalent to
\begin{align*}
N = D^{m-1} \times D \times D &\to D^{m-1} \times D = \check X,\\
(z, x, y) &\mapsto (z, xy)
\end{align*}
where $(0, 0) \in D^{m-1} \times D = \check X$ is the point corresponding to
the curve $(\check\Sigma_0, z_0\cup z^+_0 \cup (\pm\infty_i))$, and
the nodal point of $\check\Sigma_0$ is $(0, 0, 0) \in N$.
Then the oriented blow up at the set of nodal points $D^{m-1} \times \{(0, 0)\}$ is
defined by
\begin{align*}
\widetilde{N} = D^{m-1} \times \widetilde{D} \times \widetilde{D}
&\to D^{m-1} \times \widetilde{D} = \widetilde{X}.\\
(z, (x, \theta_x), (y, \theta_y)) &\mapsto (z, (xy, \theta_x \theta_y))
\end{align*}
For each marked point $\pm\infty_i$ of
$(\check\Sigma_0, z_0\cup z^+_0 \cup (\pm\infty_i))$, the fibration
$\check P \to \check X$ is locally equivalent to
\begin{align*}
N = D^m \times D &\to D^m = \check X,\\
(z, w) &\mapsto z
\end{align*}
where $0 \in D^m = \check X$ is the point corresponding to the curve
$(\check\Sigma_0, z_0\cup z^+_0 \cup (\pm\infty_i))$, and
$Z_{\pm\infty_i}(z) = (z, 0)$ is the section of marked point corresponding
to the marked point $\pm\infty_i$.
Then the oriented blow up at $Z_{\pm\infty_i}$ is defined by
\begin{align*}
\breve N = D^m \times \widetilde{D} &\to D^m = \widetilde{X}.\\
(z, (w, \theta_w)) &\mapsto z
\end{align*}

We take a discontinuous map $\Psi : \widetilde{P} \to \widetilde{P}_0$ (or a continuous
map which is defined on the complement of some codimension one subset) which satisfies
the following conditions:
\begin{itemize}
\item
$\Psi|_{\widetilde{P}_0} = \id$
\item
For each nodal point of $\Sigma_0$, we fix a neighborhood $\check N \subset
\widetilde{P}$ such that the restriction of the fibration $\widetilde{P} \to \widetilde{X}$
to $\check N$ is equivalent to
\begin{align*}
\check N = A \times D \times D &\to A \times D = \widetilde{X},\\
(a, x, y) &\mapsto (a, xy)
\end{align*}
where $A$ is some complex manifold or its oriented blow up, and $(0, 0) \in A \times D
= \widetilde{X}$ is the point corresponding to the curve $(\Sigma_0, z_0\cup z_0^+)$.
Then the restriction of $\Psi$ to $\check N$ is given by
\[
\Psi (a, x, y) = \begin{cases}
(0, x, 0) \in A \times D \times D \text{ if } |x| \geq |y|\\
(0, 0, y) \in A \times D \times D \text{ if } |y| \geq |x|
\end{cases}.
\]
Note that this is not well defined at the codimension one subset $\{|x| = |y|\}$.
\item
For each joint circle of $\Sigma_0$, we fix its neighborhood $\widetilde{N} \subset
\widetilde{P}$ such that the restriction of the fibration $\widetilde{P} \to \widetilde{X}$
to $\widetilde{N}$ is equivalent to
\begin{align*}
\widetilde{N} = A \times \widetilde{D} \times \widetilde{D} &\to A \times \widetilde{D}
= \widetilde{X},\\
(a, (x, \theta_x), (y, \theta_y)) &\mapsto (a, (xy, \theta_x \theta_y))
\end{align*}
where $(0, 0, e^{2\pi \sqrt{-1} \cdot 0}) \in \widetilde{X}$ is the point corresponding to
the curve $(\Sigma_0, z_0\cup z_0^+)$.
Then the restriction of $\Psi$ to $\widetilde{N}$ is given by
\[
\Psi(a, (x, \theta_x), (y, \theta_y)) = \begin{cases}
(0, (x, \theta_x), (0, \theta_x^{-1})) \text{ if } |x| \geq |y|\\
(0, (0, \theta_y^{-1}), (y, \theta_y)) \text{ if } |y| \geq |x|
\end{cases}.
\]
Note that if we rewrite the above fibration by the isomorphism
$\widetilde{D} \cong [-\infty, 0) \times S^1 \cong (0, \infty] \times S^1$
given by $(e^{2\pi(s + \sqrt{-1} t)}, e^{2\pi\sqrt{-1} t}) \leftrightarrow (s, t)
\leftrightarrow (-s, -t)$
as
\begin{align*}
\widetilde{N} = A \times ((0, \infty] \times S^1) \times ([-\infty, 0) \times S^1)
&\to A \times ((0, \infty] \times S^1) = \widetilde{X},\\
(a, (s_x, t_x), (s_y, t_y)) &\mapsto (a, (s_x - s_y, t_x - t_y))
\end{align*}
then $\Psi|_{\widetilde{N}}$ is expressed as
\[
\Psi(a, (s_x, t_x), (s_y, t_y)) = \begin{cases}
(0, (s_x, t_x), (-\infty, - t_x)) \text{ if } |s_x| \geq |s_y|\\
(0, (+\infty, -t_y), (s_y, t_y)) \text{ if } |s_y| \geq |s_x|
\end{cases}.
\]
\item
For each $+ \infty$-limit circle $S_{+\infty}^1$ of $\Sigma_0$, we fix its neighborhood
$\breve N_{+\infty_i} \subset \widetilde{P}$ such that the restriction of the fibration
$\check P \to \check X$ is locally equivalent to
\begin{align*}
\breve N_{+\infty_i} = A \times ((0, \infty] \times S^1)&\to A = \widetilde{X},\\
(a, s, t) &\mapsto a
\end{align*}
where $0 \in A = \widetilde{X}$ is the point corresponding to the curve
$(\Sigma_0, z_0\cup z_0^+)$.
Then the restriction of $\Psi$ to $\breve{N}_{+\infty_i}$ is given by
\[
\Psi(a, s, t) = (0, s, t).
\]
\item
For each $-\infty$-limit circle $S_{-\infty_i}^1$ of $\Sigma_0$,
we also fix its neighborhood
$\breve N_{-\infty_i} \subset \widetilde{P}$ similarly, and we assume that
the restriction of $\Psi$ to $\breve N_{-\infty_i}$ is given similarly.
\item
$\Psi$ is smooth on the complement $\widetilde{P} \setminus
(\bigcup_{\text{nodal points}} \check N \cup \bigcup_{\text{joint circles}} \widetilde{N}
\cup \bigcup_{\text{limit circles}} \breve N_{\pm\infty_i})$.
\item
$\Psi$ is continuous at the joint $\bigcup \partial \check N \cup \bigcup \partial
\widetilde{N} \cup \bigcup \partial \breve N_{\pm\infty_i}$.
\end{itemize}

A map $\theta : \overline{\R}_1 \sqcup \overline{\R}_2 \sqcup \dots \sqcup
\overline{\R}_k \to \overline{\R}_1 \cup \overline{\R}_2 \cup \dots \cup \overline{\R}_l$
is called an $\R$-gluing if there exist a surjection $\mu : \{1,2,\dots, k\} \to
\{1,2,\dots, l\}$ and constants $c_i \in \R$ ($i = 1,2, \dots, k$) such that
\begin{itemize}
\item
if $i \leq j$ then $\mu(i) \leq \mu(j)$,
\item
$\theta(\overline{\R}_i) = \overline{\R}_{\mu(i)}$, and
\item
$\theta|_{\overline{\R}_i} (s) = s + c_i \ (\in \overline{\R}_{\mu(i)})$.
\end{itemize}
For each $\R$-gluing $\theta$, let
$\theta \times 1
: (\overline{\R}_1 \sqcup \overline{\R}_2 \sqcup \dots \sqcup \overline{\R}_k) \times Y
\to (\overline{\R}_1 \cup \overline{\R}_2 \cup \dots \cup \overline{\R}_l) \times Y$
be the product with the identity map on $Y$.
%

For each neighborhood $U \subset \widetilde{X}$ of $0 \in \widetilde{X}$ and
each constant $\epsilon>0$,
we define a subset $\W_{p_0}(U, \epsilon) = \W_{p_0}(U, \epsilon, \Psi)
\subset \overline{\M}^0$ as follows.
$(\Sigma, z, u, \phi) \in \overline{\M}^0$ belongs to $\W_{p_0}(U, \epsilon)$ if
there exist a point $a\in U$, an isomorphism $(\Sigma, z) \cong (\widetilde{P}_a, Z(a))$
and an $\R$-gluing $\theta$ such that
\begin{equation}
\dist_{L^\infty}(u, (\theta \times 1) \circ u_0 \circ \Psi|_{\widetilde{P}_a}) < \epsilon
\label{L^infty comparison of maps}
\end{equation}
and
\begin{equation}
\dist_{L^\infty(S^1)}(\pi_{S^1}^{\breve N_{\pm\infty_i}} \circ \phi_{\pm\infty_i},
\phi_{0, \pm\infty_i}) < \epsilon,
\label{L^infty comparison of coordinates of S^1}
\end{equation}
where the left hand side of (\ref{L^infty comparison of maps}) is the essential sup of
$\dist(u(z), (\theta \times 1) \circ u_0 \circ \Psi (z))$ over $\widetilde{P}_a$, and
in (\ref{L^infty comparison of coordinates of S^1}),
\[
\pi_{S^1}^{\breve N_{\pm\infty_i}} : \breve N_{\pm\infty_i} \supset \widetilde{X} \times
\{\pm\infty\} \times S^1 \to \{0\} \times \{\pm\infty\} \times S^1
\]
are the projections.
We define a neighborhood of $p_0$ as a subset of $\overline{\M}^0$ which contains
$\W_{p_0}(U, \epsilon)$ for some $U \subset \widetilde{X}$ and $\epsilon >0$.

First we prove that this definition of neighborhood is independent of the choice of
$\check N$, $\widetilde{N}$, $\breve N_{\pm\infty}$ and $\Psi$.
Let $(\check N', \widetilde{N}', \breve N'_{\pm\infty}, \Psi')$ be another choice.
We claim that for any $\epsilon > 0$, there exists a neighborhood
$U \subset \widetilde{X}$ such that for any $a \in U$ and any $\R$-gluing $\theta$,
\begin{equation}
\dist_{L^\infty} ((\theta \times 1) \circ u_0 \circ \Psi'|_{\widetilde{P}_a},
(\theta \times 1) \circ u_0 \circ \Psi|_{\widetilde{P}_a})
< \epsilon + \Delta((\theta \times 1) \circ u_0 \circ \Psi|_{\widetilde{P}_a})
\label{difference}
\end{equation}
and
\begin{equation}
\dist_{L^\infty(U \times \{\pm\infty\} \times S^1)}
(\pi_{S^1}^{\breve N'_{\pm\infty_i}}, \pi_{S^1}^{\breve N_{\pm\infty_i}}) < \epsilon,
\label{S^1 difference}
\end{equation}
where $\Delta((\theta \times 1) \circ u_0 \circ \Psi|_{\widetilde{P}_a})$
is the maximum of the differences of the limits of
$(\theta \times 1) \circ u_0 \circ \Psi|_{\widetilde{P}_a}$ on the both sides
at the discontinuous codimension one subset.

First we prove that these inequalities imply the independence of the choice of
$\check N$, $\widetilde{N}$, $\breve N_{\pm\infty_i}$ and $\Psi$.
For any $(\Sigma, z, u) \in \W_{p_0}(U, \epsilon, \Psi)$,
there exist a point $a\in U$, an isomorphism $(\Sigma, z) \cong (\widetilde{P}_a, Z(a))$
and an $\R$-gluing $\theta$ such that
\[
\dist_{L^\infty}(u, (\theta \times 1) \circ u_0 \circ \Psi|_{\widetilde{P}_a}) < \epsilon
\]
and
\[
\dist_{L^\infty}(\pi_{S^1}^{\breve N_{\pm\infty_i}} \circ \phi_{\pm\infty_i},
\phi_{0, \pm\infty}) < \epsilon.
\]
Since
\[
\Delta((\theta \times 1) \circ u_0 \circ \Psi|_{\widetilde{P}_a})
\leq 2 \dist_{L^\infty}(u, (\theta \times 1) \circ u_0 \circ \Psi|_{\widetilde{P}_a})
< 2 \epsilon,
\]
inequality (\ref{difference}) implies
\[
\dist_{L^\infty} ((\theta \times 1) \circ u_0 \circ \Psi'|_{\widetilde{P}_a},
(\theta \times 1) \circ u_0 \circ \Psi|_{\widetilde{P}_a}) < 3 \epsilon,
\]
hence
\[
\dist_{L^\infty}(u, (\theta \times 1) \circ u_0 \circ \Psi'|_{\widetilde{P}_a}) < 4 \epsilon.
\]
On the other hand, inequality (\ref{S^1 difference}) implies
\[
\dist_{L^\infty}(\pi_{S^1}^{\breve N'_{\pm\infty_i}} \circ \phi_{\pm\infty_i},
\phi_{0, \pm\infty_i})
< 2 \epsilon.
\]
Therefore $\W_{p_0}(U, \epsilon, \Psi) \subset \W_{p_0}(U, 4 \epsilon, \Psi')$,
which implies the independence of $\check N$, $\widetilde{N}$, $\breve N_{\pm\infty_i}$
and $\Psi$.

The above inequalities ((\ref{difference}) and (\ref{S^1 difference})) are proved as follows.
First we need to observe the correspondence of the coordinates of $\check N$ and
$\check N'$.
Since $\{ x = 0 \}$ and $\{ y = 0 \} \subset \check N$ corresponds to $\{ x' = 0 \}$ and
$\{ y' = 0 \} \subset \check N'$ respectively,
\begin{align*}
x' &= C x (1 + O(a, x, y))\\
y' &= C' y (1 + O(a, x, y))
\end{align*}
for some constants $C$ and $C'$.
Hence in the coordinate of $\check N$,
$\Psi'$ is written as
\[
(a, x, y) \mapsto \begin{cases}
(0, x (1 + O(a, x)), 0) \text{ if } |x'| \geq |y'| \\
(0, 0, y (1 + O(a, y))) \text{ if } |y'| \geq |x'|
\end{cases}.
\]
Since $u_0$ is continuous near each nodal point,
there exists a neighborhood $\check N^\circ \subset \check N$ of the nodal point
such that
\[
\dist_{L^\infty(\check N^\circ)} (u_0 \circ \Psi'|_{\check N^\circ},
u_0 \circ \Psi|_{\check N^\circ}) < \epsilon.
\]
Hence
\[
\dist_{L^\infty} ((\theta \times 1) \circ u_0 \circ \Psi'|_{\check N^\circ},
(\theta \times 1) \circ u_0 \circ \Psi|_{\check N^\circ}) < \epsilon
\]
for any $\R$-gluing $\theta$.

Next we consider the neighborhoods $\widetilde{N}$ and $\widetilde{N}'$
of each joint circle of $\Sigma_0$.
As in the case of $\check N$ and $\check N'$,
the correspondence of these two coordinates is
\begin{align*}
(s_x', t_x') &= C + (s_x, t_x) + O(a, e^{-2\pi s_x}, e^{2\pi s_y}) \\
(s_y', t_y') &= C' + (s_y, t_y) + O(a, e^{-2\pi s_x}, e^{2\pi s_y})
\end{align*}
for some constants $C, C' \in \R \times S^1$.
Hence in the coordinate of $\widetilde{N}$,
$\Psi'$ is written as
\begin{align*}
&(a, (s_x, t_x), (s_y, t_y))\\
&\mapsto \begin{cases}
(0, (s_x, t_x) + O(a, e^{-2\pi s_x}), (-\infty, -t_x) + O(a, e^{-2\pi s_x}))
\text{ if } |s_x'| \leq |s_y'| \\
(0, (\infty, -t_y) + O(a, e^{2\pi s_y}), (s_y, t_y) + O(a, e^{2\pi s_y}))
\text{ if } |s_y'| \leq |s_x'|
\end{cases}
\end{align*}
We decompose $\widetilde{N}$ into the following four pieces and prove the inequality
for each piece.
\begin{align*}
A_1 &= \{ |s_x| \leq |s_y| \} \cap \{ |s_x'| \leq |s_y'| \} \\
A_2 &= \{ |s_x| \geq |s_y| \} \cap \{ |s_x'| \geq |s_y'| \} \\
A_3 &= \{ |s_x| \geq |s_y| \} \cap \{ |s_x'| \leq |s_y'| \} \\
A_4 &= \{ |s_x| \leq |s_y| \} \cap \{ |s_x'| \geq |s_y'| \}
\end{align*}

First we consider the pieces $A_1$ and $A_2$.
The above expression of $\Psi'$ implies that
there exists a neighborhood $\widetilde{N}^\circ \subset \widetilde{N}$ of the joint circle
such that
\[
\dist_{L^\infty} (u_0 \circ \Psi'|_{\widetilde{N}^\circ \cap A_i},
u_0 \circ \Psi|_{\widetilde{N}^\circ \cap A_i}) < \epsilon
\]
for $i = 1, 2$.
Hence for any $\R$-gluing $\theta$,
\[
\dist_{L^\infty} ((\theta \times 1) \circ u_0 \circ \Psi'|_{\widetilde{N}^\circ \cap A_i},
(\theta \times 1) \circ u_0 \circ \Psi|_{\widetilde{N}^\circ \cap A_i}) < \epsilon.
\]

Next we consider the piece $A_3$.
For any constant $C > 0$, if $U \subset \widetilde{X}$ is a sufficiently small neighborhood
of $0 \in \widetilde{X}$, then
$|s_x|, |s_x'|, |s_y|, |s_y'| \geq C$ on $\widetilde{N}^\circ|_U \cap A_3$
since $|s_x' - s_x|$ and $|s_y' - s_y|$ are bounded on $\widetilde{N}|_U \cap A_3$.
Applying Corollary \ref{third annulus} to $u_0$, we see that
\begin{align*}
u_0(s_x, t_x) &= (Ls_x + b_x, \gamma(t_y)) + O(e^{-\delta s_x})\\
u_0(s_y, t_y) &= (Ls_y + b_y, \gamma(t_y)) + O(e^{\delta s_y})
\end{align*}
for some $\gamma \in P_L$ and $b_x, b_y \in \R$.
Assume that $U$ is sufficiently small so that
$|O(e^{-\delta s_x})| \leq \epsilon$ and $|O(e^{\delta s_y})| \leq \epsilon$ on
$\widetilde{N}|_U \cap A_3$ in the above equation.
Then for any $z = (a, (s_x, t_x), (s_y, t_y)) \in \widetilde{N}|_U \cap A_3$,
\begin{align*}
&\dist((\theta \times 1) \circ u_0 \circ \Psi'(z),
(\theta \times 1) \circ u_0 \circ \Psi(z))\\
&\leq \dist((\theta \times 1) \circ u_0 \circ \Psi'(z),
(\theta \times 1) \circ u_0 (s_y, t_y)) \\
&\quad +
\dist((\theta \times 1) \circ u_0 (s_y, t_y), (\theta \times 1)(L s_y + b_y, \gamma(t_y)))\\
&\quad + \dist((\theta \times 1) \circ u_0 (s_x, t_x),
(\theta \times 1)(L s_x + b_x, \gamma(t_x)))\\
&\quad + \dist((\theta \times 1)(L s_x + b_x, \gamma(t_x)),
(\theta \times 1)(L s_y + b_y, \gamma(t_y)))\\
&< \epsilon + \epsilon + \epsilon \\
& \quad
+ \dist_{L^\infty(\{|s_x| = |s_y|\})}((\theta \times 1)(L s_x + b_x, \gamma(t_x)),
(\theta \times 1)(L s_y + b_y, \gamma(t_y)))\\
&< 3\epsilon + 2\epsilon + \Delta((\theta \times 1) \circ u_0 \circ \Psi|_{\widetilde{P}_a}).
\end{align*}

Similarly, for any $z = (a, (s_x, t_x), (s_y, t_y)) \in \widetilde{N}|_U \cap A_4$,
\[
\dist ((\theta \times 1) \circ u_0 \circ \Psi'(z), (\theta \times 1) \circ u_0 \circ \Psi(z))
< 5\epsilon + \Delta((\theta \times 1) \circ u_0 \circ \Psi|_{\widetilde{P}_a}).
\]

For each limit circle of $\Sigma_0$,
it is easy to see that there exists a neighborhood $\breve N_{\pm\infty_i}^\circ \subset
\breve N_{\pm\infty_i}$ of the circle such that
\[
\dist_{L^\infty} (u_0 \circ \Psi'|_{\breve N_{\pm\infty_i}^\circ},
u_0 \circ \Psi|_{\breve N_{\pm\infty_i}^\circ}) < \epsilon.
\]

Finally, we consider the complement of the neighborhoods $\check N^\circ$
$\widetilde{N}^\circ$ and $\breve N^\circ_{\pm\infty_i}$.
If $U \subset \widetilde{X}$ is sufficiently small neighborhood of $0 \in \widetilde{X}$,
then the continuity of $\Psi$ and $\Psi'$ on
$\widetilde{P}_U \setminus (\bigcup \check N^\circ \cup \bigcup \widetilde{N}^\circ
\cup \bigcup \breve N_{\pm\infty_i}^\circ)$ implies that
\[
\dist_{L^\infty}(u_0 \circ \Psi'|_{\widetilde{P}_U \setminus (\bigcup \check N^\circ \cup
\bigcup \widetilde{N}^\circ \cup \bigcup \breve N_{\pm\infty_i}^\circ)},
u_0 \circ \Psi|_{\widetilde{P}_U \setminus (\bigcup \check N^\circ \cup
\bigcup \widetilde{N}^\circ \cup \bigcup \breve N_{\pm\infty_i}^\circ)}) < \epsilon.
\]

These estimates prove inequality (\ref{difference}).
Inequality (\ref{S^1 difference}) can be easily checked.

Next we prove the definition of neighborhood does not depend on the choice of
the additional marked points $z^+_0$.
It is enough to compare with another sequence of marked points $z^{++}_0$
which contains $z^+_0$.
We may assume that the local universal family
$(\widetilde{P}^{++} \to \widetilde{X}^{++}, Z \cup Z^{++})$
of $(\Sigma_0, z_0 \cup z^{++}_0)$
has the form $\widetilde{P}^{++} = \widetilde{P} \times D^k$ and
$\widetilde{X}^{++} = \widetilde{X} \times D^k$,
where $D^k$ is the parameter space which determines the value of
$Z^{++} \setminus Z^+$,
and that $Z \cup Z^+$ is independent of $D^k$.
Then we can take $\check N^{++} = \check N \times D^k$,
$\widetilde{N}^{++} = \widetilde{N} \times D^k$ and
$\breve N^{++}_{\pm\infty_i} = \breve N_{\pm\infty_i} \times D^k$
as neighborhoods of nodal points and imaginary circles.
Hence we can take $\Psi^{++} = \Psi \circ \pi_{\widetilde{P}} : \widetilde{P}^{++}
\to \widetilde{P}^{++}_0 = \widetilde{P}_0$,
which implies the definitions of the neighborhood coincide.

Finally, we prove that these neighborhood systems define a topology.
It is enough to prove the following claim:
For each $p_1 \in \W_{p_0}(U, \epsilon)$,
there exists a neighborhood of $p_1$ contained in
$\W_{p_0}(U, \epsilon)$.
This implies not only the well-definedness of the topology but also that
each $\W_{p_0}(U, \epsilon)$ is open.

For each $p_1 = (\Sigma_1, z_1, u_1, \phi_1) \in \W_{p_0}(U, \epsilon)$,
there exist a point $a_1\in U$, an isomorphism
$(\Sigma_1, z_1) \cong (\widetilde{P}_{a_1}, Z(a_1))$
and an $\R$-gluing $\theta_1$ such that
\[
\dist_{L^\infty}(u_1, (\theta_1 \times 1) \circ u_0 \circ \Psi|_{\widetilde{P}_{a_1}})
< \epsilon
\]
and
\[
\dist_{L^\infty(S^1)}(\pi_{S^1}^{\breve N_{\pm\infty_i}} \circ \phi_{1, \pm\infty_i},
\phi_{0, \pm\infty_i}) < \epsilon.
\]
Let $\epsilon_1 > 0$ be a small constant such that
\[
\dist_{L^\infty}(u_1, (\theta_1 \times 1) \circ u_0 \circ \Psi|_{\widetilde{P}_{a_1}})
< \epsilon - 2\epsilon_1
\]
and
\[
\dist_{L^\infty}(\pi_{S^1}^{\breve N_{\pm\infty_i}} \circ \phi_{1, \pm\infty_i},
\phi_{0, \pm\infty_i}) < \epsilon - \epsilon_1.
\]

We use the additional marked points $z_1^+$ of $(\Sigma_1, z_1)$ which correspond to
$Z^+(a_1)$ under the above isomorphism $(\Sigma_1, z_1) \cong
(\widetilde{P}_{a_1}, Z(a_1))$.
Then the local universal family of $(\Sigma_1, z_1 \cup z_1^+)$ is the restriction of
$(\widetilde{P} \to \widetilde{X}, Z \cup Z^+)$ to a neighborhood $U'_1 \subset
\widetilde{X}$ of $a_1$.
Then for the definition of the neighborhoods of $p_1$,
we can take a discontinuous map $\Psi^1 : \widetilde{P}|_{U'_1} \to \widetilde{P}_{a_1}$
which satisfies the following conditions:
\begin{itemize}
\item
$\Psi|_{\widetilde{P}_{a_1}} \circ \Psi^1|_{\check N}
= \Psi|_{\check N} : \check N \to \widetilde{P}_0$
for the neighborhood $\check N$ of each nodal point of $\Sigma_1$.
\item
$\Psi|_{\widetilde{P}_{a_1}} \circ \Psi^1|_{\widetilde{N}}
= \Psi|_{\widetilde{N}} : \widetilde{N} \to \widetilde{P}_0$
for the neighborhood $\widetilde{N}$ of each joint circle of $\Sigma_1$.
\item
On the neighborhood of each limit circle of $\Sigma_1$, $\Psi^1$ is defined by
using the same coordinate of $\breve N_{\pm\infty_i}$ as that for $\Psi$.
\item
Let $\D \subset \widetilde{P}$ be the codimension one subset consisting of nodal points,
imaginary circles and discontinuous points of $\Psi$.
Then, $\Psi^1$ preserves $\D$.
\end{itemize}
Since $u_0$ is continuous on $\widetilde{P}_0 \setminus \D$,
the above assumption of $\Psi^1$ implies that
if $U_1 \subset U'_1$ is sufficiently small, then 
\[
\dist_{L^\infty} (u_0 \circ \Psi|_{\widetilde{P}_{a_1}} \circ \Psi^1|_{\widetilde{P}_{U_1}}, u_0 \circ \Psi|_{\widetilde{P}_{U_1}}) < \epsilon_1.
\]

For any $p = (\Sigma, z, u, \phi) \in \W_{p_1}(U_1, \epsilon_1)$,
there exist a point $a \in U_1$, an isomorphism $(\Sigma, z) \cong (\widetilde{P}_a, Z(a))$
and an $\R$-gluing $\theta$ such that
\[
\dist_{L^\infty}(u, (\theta \times 1) \circ u_1 \circ \Psi^1|_{\widetilde{P}_a}) < \epsilon_1
\]
and
\[
\dist_{L^\infty(S^1)}(\pi_{S^1}^{\breve N_{\pm\infty_i}} \circ \phi_{\pm\infty_i},
\phi_{1, \pm\infty_i}) < \epsilon_1.
\]
Hence
\begin{align*}
&\dist_{L^\infty}(u, (\theta \times 1) \circ (\theta_1 \times 1) \circ u_0 \circ
\Psi|_{\widetilde{P}_a})\\
&\leq \dist_{L^\infty}(u, (\theta \times 1) \circ u_1 \circ \Psi^1|_{\widetilde{P}_a}) \\
&\quad + \dist_{L^\infty}((\theta \times 1) \circ u_1 \circ \Psi^1|_{\widetilde{P}_a},
(\theta \times 1) \circ (\theta_1 \times 1) \circ u_0 \circ \Psi|_{\widetilde{P}_{a_1}}
\circ \Psi^1|_{\widetilde{P}_a})\\
&\quad + \dist_{L^\infty}((\theta \times 1) \circ (\theta_1 \times 1) \circ u_0
\circ \Psi|_{\widetilde{P}_{a_1}} \circ \Psi^1|_{\widetilde{P}_a}, (\theta \times 1)
\circ (\theta_1 \times 1) \circ u_0 \circ \Psi|_{\widetilde{P}_a})\\
&\leq \dist_{L^\infty}(u, (\theta \times 1) \circ u_1 \circ \Psi^1|_{\widetilde{P}_a})
+ \dist_{L^\infty}(u_1, (\theta_1 \times 1) \circ u_0 \circ \Psi|_{\widetilde{P}_{a_1}})\\
&\quad + \dist_{L^\infty}(u_0 \circ \Psi|_{\widetilde{P}_{a_1}} \circ
\Psi^1|_{\widetilde{P}_a}, u_0 \circ \Psi|_{\widetilde{P}_a})\\
&< \epsilon,
\end{align*}
and
\begin{align*}
&\dist_{L^\infty(S^1)}(\pi_{S^1}^{\breve N_{\pm\infty_i}} \circ \phi_{\pm\infty_i},
\phi_{0, \pm\infty_i})\\
&\leq
\dist_{L^\infty(S^1)}(\pi_{S^1}^{\breve N_{\pm\infty_i}|_{U_1}} \circ \phi_{\pm\infty_i}, \phi_{1, \pm\infty_i})
+ \dist_{L^\infty(S^1)}(\pi_{S^1}^{\breve N_{\pm\infty_i}} \circ \phi_{1, \pm\infty_i}, \phi_{0, \pm\infty_i})\\
& < \epsilon.
\end{align*}
These imply $p \in \W_{p_0}(U, \epsilon)$, which proves the claim.

Next we prove the topological properties of
$\overline{\M}^0 = \overline{\M}^0(Y, \lambda, J)$ along the similar lines in the case of
Gromov-Witten theory in \cite{FO99}. $\overline{\M}^0$ is decomposed as
$\overline{\M}^0 = \coprod_{g, \mu, L^-, L^+} \overline{\M}^0_{g, \mu}(L^-,L^+)$,
where $\overline{\M}^0_{g, \mu}(L^-,L^+)$ is the space of holomorphic buildings with
genus $g$ and $\mu$ marked points such that
$\sum L_{\gamma_{-\infty_i}} = L^-$ and $\sum L_{\gamma_{+\infty_i}} = L^+$.
(The genus of blown up curve $\Sigma$ is by definition the genus of $\check \Sigma$.)
First we show that we have a nice way to add marked points to the domain curves.
\begin{lem}\label{nice marked points}
Let $\epsilon > 0$ and $\delta_0 > 0$ be arbitrary small constants,
and let $(\Sigma, z, u, \phi) \in \overline{\M}^0_{g, \mu}(L^-, L^+)$ be an arbitrary
holomorphic building.
We regard $u : \Sigma \to (\overline{\R}_1 \cup \overline{\R}_2 \cup \dots \cup
\overline{\R}_k) \times Y$
not as an equivalence class by $\R$-translations but as a map.
Then there exist a closed subset $I \subset \R_1 \cup \R_2 \cup \dots \cup \R_k$
and additional marked points $z^+$ of $\Sigma$ which satisfy the following conditions:
\begin{itemize}
\item
$I$ is a finite union of intervals in the form $[l, l+1] \subset \R_i$ {\rm(}$l \in \Z${\rm)}.
\item
The length of $I$ and the number of additional marked points are bounded by
some constant determined by $g$, $\mu$, $L^-$, $L^+$, $\epsilon$ and $\delta_0$.
\item
$(\Sigma, z \cup z^+)$ is stable.
\item
There exists a constant $A_1 > 0$ depending only on $g$, $\mu$, $L^-$, $L^+$,
$\epsilon$ and $\delta_0$ such that
if $[-A_1, T + A_1] \times S^1 \subset \Sigma$ does not contain any marked points $z \cup z^+$,
then one of the following two holds true.
\begin{enumerate}[label=\normalfont(\arabic*)]
\item
$u([0, T] \times S^1) \subset I \times Y$ and $\diam \, u([0, T] \times S^1) \leq 20 \epsilon$.
\item
$\sigma \circ u ([0, T] \times S^1)$ is contained in the $\frac{1}{3}$-neighborhood of
the complement of $I \subset \R_1 \cup \R_2 \cup \dots \cup \R_k$,
and there exist $L \in \R$ and $(b, \gamma) \in \R \times P_L$ such that
\[
\dist ( u(s,t), (Ls + b, \gamma(t)) ) \leq \epsilon ( e^{-\kappa s} + e^{-\kappa (T - s)})
\]
on $[0, T] \times S^1$.
\end{enumerate}
In particular,
for any disc $D \subset \Sigma$ such that
$D \setminus 0$ does not contain any marked points,
$\diam \, u(\{ z\in D; |z| \leq e^{-2\pi A_1}\}) \leq 20 \epsilon$.
\item
$\sigma \circ u(z \cup z^+) \subset \R_1 \cup \R_2 \cup \dots \cup \R_k$ is contained
in the $\frac{1}{3}$-neighborhood of $I$.
\item
Each connected component of $u^{-1}(I \times Y)$ either
contains at least one point of $z \cup z^+$
or is contained in the inverse image of the $\frac{1}{3}$-neighborhood of
the complement of $I$ by $\sigma \circ u$.
\item
For the $\frac{1}{3}$-neighborhood $J$ of each connected component of the complement
of $I$,
$E_{\hat \omega}(u|_{u^{-1}(J \times Y)}) \leq \delta_0$.
\end{itemize}
\end{lem}
\begin{proof}
First we see the energy bound:
$E_\lambda(u) \leq L^+$ and $E_{\hat \omega}(u) = L^+ - L^-$
for any $(\Sigma, z, u, \phi) \in \overline{\M}^0_{g, \mu}(L^-,L^+)$.
The former is because for any interval $I \subset \R_i$,
\begin{align*}
\frac{1}{|I|} \int_{u^{-1}(I\times Y)} u^\ast (d\sigma \wedge \lambda)
&= \int u^\ast d\varphi \wedge \lambda\\
&= \int u^\ast d(\varphi \lambda) - \int u^\ast (\varphi d\lambda)\\
&\leq L^+,
\end{align*}
where $\varphi : \R_1 \cup \R_2 \cup \dots \cup \R_k \to \R$ is defined by
\[
\varphi(\sigma) = \int_{-\infty}^\sigma \frac{1}{|I|} 1_I (\sigma') d\sigma' \quad \text{on } \R_i,
\]
$\varphi \equiv 0$ on $\R_j$ ($j < i$)
and $\varphi \equiv 1$ on $\R_j$ ($j > i$), and the last inequality is because
$u^\ast(\varphi d\lambda) \geq 0$ by the equation of $J$-holomorphic curves.
Proof of the latter equation $E_{\hat \omega}(u) = L^+ - L^-$ is straightforward.

Next we prove the number of irreducible components of $\check \Sigma$ is bounded by
some constant depending only on $g$, $\mu$, $L^-$ and $L^+$.
Note that if $E_{\hat\omega} (u|_{\Sigma_\alpha}) > 0$, then
$E_{\hat\omega} (u|_{\Sigma_\alpha}) \geq \min (\sum_i L^+_i - \sum_j L^-_j) \ (>0)$,
where the minimum is taken over all pairs of families of periods $(L^+_i)_i$ and $(L^-_j)_j$
such that $\sum_j L^-_j < \sum_i L^+_i \leq L^+$.
Hence the number of the components $\Sigma_\alpha$ on which $u$ have non-zero
$E_{\hat\omega}$-energies is bounded.

If $E_{\hat\omega} (u|_{\Sigma_\alpha}) = 0$ and $2g_\alpha + m_\alpha < 3$, then
$\Sigma_\alpha$ does not contain any marked points and $(\Sigma_\alpha, u)$ is a trivial cylinder.
We can see it by the following consideration:
\begin{itemize}
\item
If the number of imaginary circles in $\Sigma_\alpha$ were zero,
then $u$ would be a constant map since every closed $J$-holomorphic curve
in $\hat Y$ is a constant map.
However, this contradict the definition of holomorphic building.
\item
The number of imaginary circles in $\Sigma_\alpha$ cannot be one since $E_{\hat\omega}(u|_{\Sigma_\alpha}) = 0$.
\item
If the number of imaginary circles in $\Sigma_\alpha$ is two,
then $g_\alpha = 0$ and $\Sigma_\alpha$ does not contain any marked points or nodal points.
Hence $\Sigma_\alpha \cong \overline{\R} \times S^1$ and $(\Sigma_\alpha, u)$ is a
trivial cylinder.
\end{itemize}
Therefore the number of the nontrivial components $\Sigma_\alpha$ such that
$2 g_\alpha + m_\alpha < 3$ is bounded.

Since the number of the limit circles is bounded, this implies that
the number of the components $\Sigma_\alpha$ such that
$2 g_\alpha + m_\alpha \geq 3$ is also bounded.
This is due to the equality
\[
2g + \mu + \, \text{(the number of the limit circles)} \, -2 = \sum_\alpha (2g_\alpha + m_\alpha -2)
\]
and the fact that trivial cylinders do not contribute to the sum on
the right hand side of the equation.

Therefore the number of the nontrivial components is bounded.
In particular, the height $k$ of the $J$-holomorphic building $(\Sigma, z, u)$ is bounded.
Let $S \subset \Sigma$ be the union of the trivial cylinders of $(\Sigma, z, u)$.
Then each connected component of $S$ consists of at most $(k-1)$ trivial cylinders
and it shares a joint circle with some nontrivial component.
Since the number of the joint circles contained in the nontrivial components is bounded,
it implies that the number of the trivial cylinders is also bounded.
Hence the number of the irreducible components of $\check \Sigma$ is bounded.

Therefore, the number of marked points we need to add to $(\Sigma, z)$ in order to
make $(\Sigma, z \cup z^+)$ stable is bounded.
Assuming that $(\Sigma, z \cup z^+)$ is stable,
we further add marked points $z^{++}$ as follows.
We may assume that $\epsilon < \min(\frac{1}{60}, \frac{1}{24} L_{\min} )$,
where $L_{\min}$ is the minimal period of periodic orbits.
Let $\delta > 0$, $\kappa > 0$, $A > 0$ and $L_0$ be the constant of
Corollary \ref{third annulus} for $C_0 = L^+$ and the given $\epsilon > 0$.
We may assume that $\delta \leq \delta_0$.

First, let $I \subset \R_1 \cup \R_2 \cup \dots \cup \R_k$ be a finite union of intervals
$[l, l+1] \subset \R_i$ ($l \in \Z$)
such that
\begin{itemize}
\item
$E_{\hat\omega} (u|_{u^{-1}(J \times Y)}) \leq \delta$ for the
$\frac{1}{3}$-neighborhood $J$ of each connected component of
the complement of $I \subset \R_1 \cup \R_2 \cup \dots \cup \R_k$, and
\item
$\sigma \circ u(z \cup z^+) \subset I$.
\end{itemize}
We may assume that the length of $I$ is bounded by some constant depending only on
$E_{\hat\omega} (u)$, $\delta$ and the number of marked points $z \cup z^+$.

Let $\bigcup_\alpha B_\alpha^1 \supset I \times Y$ be a finite covering by open balls
with radius $\epsilon$,
where the distance of $\R \times Y$ is given by
$\dist((\sigma, y), (\sigma', y'))^2 = |\sigma - \sigma'|^2 + \dist_Y(y, y')^2$
for some distance $\dist_Y$ of $Y$.
We may assume that the number of open balls is bounded by some constant depending
on the length of $I$ and $\epsilon$.
For each $B_\alpha^1$, let $B_\alpha^2$ be the concentric ball with radius $2\epsilon$.
We may assume that $\sigma(B_\alpha^2) \subset \R_1 \cup \R_2 \cup \dots
\cup \R_k$ is contained in the $\frac{1}{3}$-neighborhood of $I$
since $4\epsilon < \frac{1}{3}$.
Then for each connected component of $u^{-1}(B_\alpha^2)$ which contains some points
of $u^{-1}(B_\alpha^1)$, we choose one of these points in $u^{-1}(B_\alpha^1)$
as an additional marked point.
Then the number of the additional marked points is bounded since
\begin{itemize}
\item
if a connected component $\Omega$ of $u^{-1}(B_\alpha^2) \subset \Sigma$
contains a point $z \in u^{-1}(B_\alpha^1)$,
then $u(\partial \Omega) \cap B_\epsilon(u(z)) = \emptyset$,
hence Lemma \ref{monotonicity lemma} implies
$|du|_{L^2(\Omega)}$ is larger than some positive constant depending on $\epsilon$,
and
\item
the total energy on $u^{-1}(B_\alpha^2)$ is bounded by
$|du|_{L^2(u^{-1}(B_\alpha^2))}^2 \leq E_{\hat \omega}(u) + 4\epsilon E_\lambda (u)$.
\end{itemize}

We rewrite $z^+ \cup z^{++}$ as $z^+$.
We claim that this is the required additional marked points.
The only non-trivial condition is the condition about annuli.
Define $A_1 = (2 A + 2) \cdot \lceil \frac{E_{\hat\omega} (u)}{\delta} \rceil$.
(Recall that $A > 0$ is the constant of Corollary \ref{third annulus}.)

First we claim that for each annulus $[0, A_1] \times S^1 \subset \Sigma$,
there exist $s_0 \in [A, A_1 -A]$, $L\in \R$ and $(b, \gamma) \in \R \times P_L$
such that
\[
\dist( u(s,t), (Ls + b, \gamma(t))) \leq 2\epsilon \quad \text{on } [s_0 -1, s_0 +1] \times S^1.
\]
This is proved as follows.
Decompose $[0, A_1] = \bigcup [(2A + 2) i, (2A + 2) (i + 1)]$
into $\lceil \frac{E_{\hat\omega} (u)}{\delta} \rceil$ pieces of
intervals with length $2A + 2$.
Then one of them $[(2A + 2) i, (2A + 2) (i + 1)]$ satisfies
$E_{\hat\omega} (u|_{[(2A + 2) i, (2A + 2) (i + 1)] \times S^1}) \leq \delta$.
Hence Corollary \ref{third annulus} implies $s_0 = (2A + 2) i + A + 1$ satisfies the above condition.

Now we assume $[-A_1, T + A_1] \times S^1 \subset \Sigma$ does not contain any
marked points and prove that the required condition holds true.
The above claim implies that there exist $s_1 \in [-A_1 + A, -A]$,
$s_2 \in [T + A, T + A_1 -A]$, $L_i\in \R$ and $(b_i, \gamma_i) \in \R \times P_{L_i}$
($i = 1,2$) such that
\begin{equation}
\dist( u(s,t), (L_i s + b_i, \gamma_i (t))) \leq 2\epsilon \quad \text{on } [s_i -1, s_i +1] \times S^1 \label{thin}
\end{equation}
for each $i = 1, 2$.
In particular, this implies $\diam \, \sigma \circ u( \{s_i\} \times S^1) \leq 4 \epsilon$
for each $i = 1, 2$.

For each $z \in [s_1, s_2] \times S^1$ such that $u(z) \in I \times Y$,
there exists some $\alpha$ such that $u(z) \in B_\alpha^1$.
Then $B_\alpha^2$ intersects with $u( \partial ([s_1, s_2] \times S^1) )$ since
$[s_1, s_2] \times S^1$ does not contain any marked points.
(If they did not intersect, then the connected component of $u^{-1}(B_\alpha^2)$
containing $z$ would be contained in $[s_1, s_2] \times S^1$.)
Therefore $u(z)$ is contained in the $3 \epsilon$-neighborhood of
$u(\partial [s_1, s_2] \times S^1)$.

We separate the argument into the following two cases.
\begin{enumerate}[label=\normalfont(\arabic*)]
\item
$\sigma\circ u([s_1, s_2] \times S^1) \subset I$
\item
$\sigma\circ u([s_1, s_2] \times S^1) \not\subset I$
\end{enumerate}
In the first case, $u([s_1, s_2] \times S^1)$ is contained in
the $3 \epsilon$-neighborhood of $u(\{s_1\} \times S^1) \cup u(\{s_2\} \times S^1)$.
Since the diameter of the $3 \epsilon$-neighborhood of each
$\sigma \circ u(\{s_i\} \times S^1)$ is
$\leq 4\epsilon + 2 \cdot 3 \epsilon \leq 10\epsilon$,
it implies $\diam \, \sigma \circ u([s_1, s_2] \times S^1) \leq 20 \epsilon$.
Then $L_i = 0$ ($i = 1, 2$) because if not, (\ref{thin}) implies that
the diameter of $\sigma \circ u ([s_1, s_1 + 1] \times S^1)$ or 
$\sigma \circ u ([s_2 - 1, s_2] \times S^1)$ is $\geq L_{\min} -4\epsilon > 20 \epsilon$.
Therefore (\ref{thin}) implies that
$\diam \, u(\{s_i\} \times S^1) \leq 4 \epsilon$ ($i = 1,2$).
Hence $\diam \, u([s_1, s_2] \times S^1) \leq 20\epsilon$ because
$u([s_1, s_2] \times S^1)$ is contained in
the $3 \epsilon$-neighborhood of $u(\{s_1\} \times S^1) \cup u(\{s_2\} \times S^1)$.

In the second case, $\sigma \circ u ([s_1, s_2] \times S^1)$ is contained in
the $20 \epsilon$-neighborhood of the complement of $I$
because it is covered by the complement of $I$ and the $3 \epsilon$-neighborhood of
$u(\{s_1\} \times S^1) \cup u(\{s_2\} \times S^1)$.
Since $20\epsilon < \frac{1}{3}$, it is contained in the $\frac{1}{3}$-neighborhood of
a connected component of the complement of $I$, which implies
$E_{\hat\omega} (u|_{[s_1, s_2] \times S^1}) \leq \delta$.
Since $[-A, T + A] \subset  [s_1, s_2]$,
Corollary \ref{third annulus} implies
there exists $L \in \R$, $(b, \gamma) \in \R \times P_L$ and $\kappa > 0$ such that
\[
\dist( u(s,t), (Ls + b, \gamma(t))) \leq \epsilon (e^{-\kappa s} + e^{- \kappa(T-s)}) \quad \text{on } [0, T] \times S^1.
\]
\end{proof}

\begin{cor}\label{cor nice marked points}
In Lemma \ref{nice marked points},
we can replace the condition of annuli with the following stronger condition:
\begin{itemize}
\item
If $[-A_1, T + A_1] \times S^1 \subset \Sigma$ does not contain any marked points $z \cup z^+$,
then there exist $L \in \R$ and $(b, \gamma) \in \R \times P_L$ such that
\[
\dist ( u(s,t), (Ls + b, \gamma(t)) ) \leq \epsilon ( e^{-\kappa s} + e^{-\kappa (T - s)})
\]
on $[0, T] \times S^1$.
Furthermore, if $L \neq 0$, then
$\sigma \circ u ([0, T] \times S^1)$ is contained in the $\frac{1}{3}$-neighborhood of
the complement of $I$.
\end{itemize}
\end{cor}
\begin{proof}
This is because if the diameter of $u([-1, T +1] \times S^1)$ is sufficiently small,
then Lemme \ref{L^infty diam} implies $|du|_{L^\infty([0, T] \times S^1)}$ is also small,
and we can apply Proposition \ref{second annulus} on $[0, T] \times S^1$.
\end{proof}

\begin{prop}\label{second countable}
$\overline{\M}^0(Y, \lambda, J)$ is second countable.
\end{prop}
\begin{proof}
It is enough to prove that each $\overline{\M}^0_{g, \mu}(L^-, L^+)$ is second countable.
Basically, this is because Lemma \ref{nice marked points} implies
that $\overline{\M}^0_{g, \mu}(L^-, L^+)$ is covered by a countable family of
open subsets consisting of equicontinuous maps.
To explain the details, first we need a preliminary consideration.

Let $(\Sigma, z, u, \phi) \in \overline{\M}^0_{g}(L^-, L^+)$ be a holomorphic building
with a stable domain curve $(\Sigma, z)$
(the number of marked points may be larger than $\mu$).
Let $(\widetilde{P} \to \widetilde{X}, Z)$ be its local universal family.
Let $R \subset \widetilde{X}$ be the subset of the points whose fibers have
the same number of nodal points and imaginary circles as that of $(\Sigma, z)$.
Take a discontinuous map $\Psi : \widetilde{P} \to \widetilde{P}_0$ as in the definition of
the topology of $\overline{\M}^0(Y, \lambda, J)$. ($0 \in \widetilde{X}$ is the point
whose fiber is isomorphic to $(\Sigma, z)$.)
We may assume that for each $a\in R$,
\[
\Psi|_{\widetilde{P}_a \setminus \coprod S^1} : \widetilde{P}_a \setminus
\coprod_{\text{joint circles}} S^1 \to \widetilde{P}_0 \setminus
\coprod_{\text{joint circles}} S^1
\]
is a homeomorphism.
Hence when we regard $(\widetilde{P} \to \widetilde{X}, Z)$
as the local universal family of $(\widetilde{P}_a, Z(a))$,
we can use $\Psi^a = (\Psi|_{\widetilde{P}_a})^{-1} \circ \Psi : \widetilde{P}
\to \widetilde{P}_a$.

For each open subset $U \subset \widetilde{X}$,
we define $\Xi (U)$ as the set of pairs $(a, u)$
each of which consists of a point $a \in U$ and a holomorphic building
$u : \widetilde{P}_a \to
(\overline{\R}_1 \cup \overline{\R}_2 \cup \dots \cup \overline{\R}_k) \times Y$
which is contained in $\overline{\M}^0_{g}(L^-, L^+)$ and which satisfies the
following condition:
If $[-A_1, T + A_1] \times S^1 \subset \widetilde{P}_a$ does not contain any marked
points $Z(a)$, then there exist $L \in \R$ and $(b, \gamma) \in \R \times P_L$
such that
\[
\dist ( u(s,t), (Ls + b, \gamma(t)) ) \leq \epsilon ( e^{-\kappa s} + e^{-\kappa (T - s)})
\]
on $[0, T] \times S^1$.

Define $\Xi_R(U) = \{ (a,u) \in \Xi(U); a \in U \cap R \}$.
For each $\epsilon >0$ and $(a_0, u_0) \in \Xi_R(U)$,
let $\widetilde{\W}_{(a_0, u_0)} (U, \epsilon)$ be the space of points
$(a,u) \in \Xi (U)$ such that
\[
\dist_{L^\infty}(u, (\theta \times 1) \circ u_0 \circ \Psi^{a_0}|_{\widetilde{P}_a})
< \epsilon
\]
for some $\R$-gluing $\theta$.
First we prove that for any $\epsilon >0$ and any two points
$(a_0, u_0)$, $(a_1, u_1) \in \Xi_R(U)$,
if $(a_0, u_0) \in \widetilde{\W}_{(a_1, u_1)}(U, \epsilon)$,
then $\widetilde{\W}_{(a_1, u_1)}(U, \epsilon) \subset
\widetilde{\W}_{(a_0, u_0)}(U, 2\epsilon)$.

Since $(a_0, u_0) \in \widetilde{\W}_{(a_1, u_1)}(U, \epsilon)$,
there exists an $\R$-translation
$\theta_0 : \overline{\R}_1 \cup \overline{\R}_2 \cup \dots \cup \overline{\R}_k
\to \overline{\R}_1 \cup \overline{\R}_2 \cup \dots \cup \overline{\R}_k$ such that
\[
\dist_{L^\infty} (u_0, (\theta_0 \times 1) \circ u_1 \circ \Psi^{a_1}|_{\widetilde{P}_{a_0}})
< \epsilon.
\]
For any $(a, u) \in \widetilde{\W}_{(a_1, u_1)}(U, \epsilon)$,
there exists an $\R$-gluing $\theta : \overline{\R}_1 \sqcup \overline{\R}_2 \sqcup
\dots \sqcup \overline{\R}_k \to \overline{\R}_1 \cup \overline{\R}_2 \cup \dots \cup
\overline{\R}_l$ such that
\[
\dist_{L^\infty} (u, (\theta \times 1) \circ u_1 \circ \Psi^{a_1}|_{\widetilde{P}_a})
< \epsilon.
\]
Since $a_0$ and $a_1$ are contained in $R$,
$\Psi^{a_0} = (\Psi^{a_1}|_{\widetilde{P}_{a_0}})^{-1} \circ \Psi^{a_1}
: \widetilde{P} \to \widetilde{P}_{a_0}$.
Hence
\begin{align*}
&\dist_{L^\infty} (u, (\theta \circ \theta_0^{-1} \times 1) \circ u_0 \circ
\Psi^{a_0}|_{\widetilde{P_a}})\\
&\leq \dist_{L^\infty} (u, (\theta \times 1) \circ u_1 \circ \Psi^{a_1}|_{\widetilde{P}_a})\\
& \quad + \dist_{L^\infty} ((\theta \times 1) \circ u_1 \circ \Psi^{a_1}|_{\widetilde{P}_a},
(\theta \circ \theta_0^{-1} \times 1) \circ u_0 \circ
(\Psi^{a_1}|_{\widetilde{P}_{a_0}})^{-1} \circ \Psi^{a_1}|_{\widetilde{P}_a})\\
& \leq \dist_{L^\infty} (u, (\theta \times 1) \circ u_1 \circ \Psi^{a_1}|_{\widetilde{P}_a})
+ \dist_{L^\infty} (u_1, (\theta_0 \times 1)^{-1} \circ u_0 \circ
(\Psi^{a_1}|_{\widetilde{P}_{a_0}})^{-1})\\
&< 2\epsilon,
\end{align*}
which implies $\widetilde{\W}_{(a_1, u_1)}(U, \epsilon) \subset
\widetilde{\W}_{(a_0, u_0)}(U, 2\epsilon)$.

We can choose a countable points $(a_i, u_i) \in \Xi_R(U)$
such that $\{\widetilde{\W}_{(a_i, u_i)}(U, \epsilon)\}_i$ covers $\Xi_R(U)$
for any $\epsilon >0$ because the assumption of the holomorphic buildings in $\Xi(U)$
implies the equicontinuity.
Then for any $(a, u) \in \Xi_R(U)$ and $\epsilon >0$,
there exists $(a_i, u_i)$ such that
$(a, u) \in \widetilde{\W}_{(a_i, u_i)}(U, \epsilon) \subset
\widetilde{\W}_{(a, u)}(U, 2\epsilon)$

Let $\{U_j\}$ be a countable open basis of the union of the base spaces of
the universal families and we choose the above $\{(a_i^{(j)}, u_i^{(j)})\}_i$ for each $U_j$.
For each $(a_i^{(j)}, u_i^{(j)})$, we fix a family of coordinates of limit circles
$\phi_{\pm\infty_l} : S^1 \to S^1_{a_i^{(j)}, \pm\infty_l}$.
Let $\mathring{Z}(a_i^{(j)}) \subset Z(a_i^{(j)})$ be an arbitrary subsequence whose
cardinality is $\mu$.
Then
\[
(\widetilde{P}_{a_i^{(j)}}, \mathring{Z}(a_i^{(j)}), u_i^{(j)},
(\phi_{\pm\infty_l} + \frac{c_l}{2^n})) \quad (n \in \N, 1 \leq c_l \leq 2^n)
\]
is a countable family of holomorphic buildings in $\overline{\M}^0_{g, \mu}(L^-, L^+)$.
Let $p_k^{(i, j)}$ ($k \in \N$) be such holomorphic buildings for each $(a_i^{(j)}, u_i^{(j)})$
and all choices of the subsequence $\mathring{Z}(a_i^{(j)}) \subset Z(a_i^{(j)})$.

We claim that $\{\W_{p_k^{(i, j)}}(U_j, 2^{-l}) \}_{i, j, k, l \in \N}$ is a countable basis of
$\overline{\M}^0_{g, \mu}(L^-, L^+)$.
This is proved as follows.
For any $p = (\Sigma, z, u) \in \overline{\M}^0_{g, \mu}(L^-, L^+)$,
we can choose additional marked points $z^+ \subset \Sigma$
as in Corollary \ref{cor nice marked points}.
Let $(\widetilde{P} \to \widetilde{X}, Z \cup Z^+)$ be the local universal family of
$(\Sigma, z \cup z^+)$.
Then for any neighborhood $\mathcal{N}$ of $p$,
there exists $U_j$ and $\epsilon \in \{2^{-l} \}$ such that
$\W_p(U_j, 2\epsilon)$ is contained in $\mathcal{N}$.

Note that we may assume that the point $a \in U_j$ whose fiber is isomorphic to
$(\Sigma, z \cup z^+)$ is contained in $R$.
Choose $(a_i^{(j)}, u_i^{(j)}) \in \Xi_R(U_j)$ such that
\[
(a, u) \in \widetilde{\W}_{(a_i^{(j)}, u_i^{(j)})} (U_j, \epsilon) \subset
\widetilde{\W}_{(a, u)}(U_j, 2\epsilon).
\]
This implies that there exists a holomorphic
building $p_k^{(i, j)} \in \overline{\M}^0_{g, \mu}(L^-, L^+)$ such that
\[
p \in \W_{p_k^{(i, j)}} (U_j, \epsilon) \subset \W_p(U_j, 2 \epsilon) \subset \mathcal{N}.
\]
Therefore $\{\W_{p_k^{(i, j)}}(U_j, 2^{-l}) \}_{i, j, k, l \in \N}$ is a countable basis of
$\overline{\M}^0_{g, \mu}(L^-, L^+)$.
\end{proof}

\begin{prop}\label{compactness}
Each $\overline{\M}^0_{g, \mu}(L^-, L^+)$ is compact.
\end{prop}
\begin{proof}
Since $\overline{\M}^0_{g, \mu}(L^-,L^+)$ is second countable,
in order to prove its compactness, it is enough to prove that any sequence
$p_i = (\Sigma_i, z_i, u_i, \phi_i) \in \overline{\M}^0_{g, \mu}(L^-,L^+)$ contains
a subsequence which converges to a point in $\overline{\M}^0_{g, \mu}(L^-,L^+)$.

Let $I_i \subset \R_1 \cup \R_2 \cup \dots \cup \R_{k_i}$ and $z_i^+ \subset \Sigma_i$
be the pair of closed subset and additional marked points given by Corollary
\ref{cor nice marked points} for sufficiently small $\epsilon > 0$ and $\delta_0 > 0$.
Passing to a subsequence,
we may assume the following:
\begin{itemize}
\item
The number of the additional marked points is independent of $i$.
\item
$(\check \Sigma_i, z_i \cup z_i^+ \cup (\pm \infty_i))$ converges to a stable curve
$(\check \Sigma, z \cup z^+ \cup (\pm\infty_i))$ in the moduli space of
marked stable curves.
\end{itemize}

Let $(\Sigma', z \cup z^+)$ be the oriented blow up of $(\check \Sigma, z \cup z^+)$
at $\pm \infty_i$ and nodal points of $\check \Sigma$ corresponding to joint circles in
$\Sigma_i$ by appropriate $\varphi$'s ($\in S^1$).
Let $(\widetilde{P}' \to \widetilde{X}', Z \cup Z^+)$ be the local universal family of
$(\Sigma', z \cup z^+)$, and let $0 \in \widetilde{X}'$ be the point
whose fiber is isomorphic to $(\Sigma', z \cup z^+)$.
Choosing appropriate $\varphi$'s, we may assume
there exists a sequence $x'_i \in \widetilde{X}'$ converging to $0 \in \widetilde{X}'$
such that $(\widetilde{P}'_{x'_i}, Z(x'_i) \cup Z^+(x'_i)) \cong (\Sigma_i, z_i \cup z_i^+)$.

Let $\Psi' : \widetilde{P}' \to \widetilde{P}'_0$ be the discontinuous map used for the
definition of the topology of $\overline{\M}^0(Y, \lambda, J)$.
We may assume $\Psi'$ maps marked points $Z \cup Z^+$ to $Z(0) \cup Z^+(0)$.
Define a map
\begin{align*}
&v_i = u_i \circ (\Psi'|_{\widetilde{P}'_{x'_i} \setminus \coprod S^1})^{-1} :
\Sigma' (\cong \widetilde{P}_0)
\supset \Psi'(\widetilde{P}'_{x'_i} \setminus \coprod S^1) \\
&\hph{v_i = u_i \circ (\Psi'|_{\widetilde{P}'_{x'_i} \setminus \coprod S^1})^{-1} :
\Sigma' (\cong \widetilde{P}_0) \supset}
\to (\overline{\R}_1 \cup \overline{\R}_2 \cup \dots \cup \overline{\R}_{k_i}) \times Y
\end{align*}
for each $i$.
Let $q^j$ be the new nodal points in $\Sigma'$.
(Namely, neighborhoods of $q^j$ correspond to annuli in $\Sigma_i$.)
Then the annulus condition of Lemma \ref{nice marked points} (or Corollary
\ref{cor nice marked points}) implies that
on any connected compact subset of $\Sigma' \setminus (\coprod S^1 \cup \{q^j\})$,
a subsequence of $v_i$ converges to a $J$-holomorphic map $v_\infty$
if we change each $v_i$ by $\R$-translation.

Let $q \in \Sigma'$ be one of new nodal points.
Recall that the restriction of the fibration $\widetilde{P}' \to \widetilde{X}'$ to the
neighborhood $\check N \subset \widetilde{P}$ of $q$ is equivalent to
\begin{align*}
\check N = A \times D \times D &\to A \times D = \widetilde{X}'\\
(a, x, y) &\mapsto (a, xy)
\end{align*}
and $\Psi'|_{\check N}$ is defined by
\[
\Psi'(a, x, y) =\begin{cases}
(0, x, 0) \quad \text{if } |x| \geq |y|\\
(0, 0, y) \quad \text{if } |y| \geq |x|
\end{cases}
\]
We may assume that $\check N$ does not contain any marked points.
Assume that $x_i' = (a_i, \ab e^{2\pi(-\rho_i + \sqrt{-1} \varphi_i)})$.
Then $\rho_i \to \infty$ as $i \to \infty$,
and $\check N \cap \widetilde{P}'_{x'_i} \cong
[0, \rho_i] \times S^1 \cup [-\rho_i, 0] \times S^1$,
where $\{\rho_i\} \times S^1 \subset [0, \rho_i] \times S^1$ and
$\{-\rho_i\} \times S^1 \subset [-\rho_i, 0] \times S^1$ are identified by
$(\rho_i, t_x) \sim (-\rho_i, t_y)$ if $t_y -t_x = \varphi_i$.
Since $\check N$ does not contain any marked points,
there exist $L \in \R$ and $(b_i, \gamma_i) \in \R \times P_L$ such that
\begin{align*}
\dist (u_i(s,t), (Ls + b_i, \gamma_i(t))) &\leq 2\epsilon e^{-\kappa |s|} \quad \text{on }
[A_1, \rho_i] \times S^1,\\
\dist (u_i(s,t), (Ls + b_i + 2 L \rho_i, \gamma_i(t + \varphi_i))) & \leq 2 \epsilon
e^{-\kappa |s|} \quad
\text{on } [-\rho_i, -A_1] \times S^1.
\end{align*}

We may assume that $L$ is nonnegative and independent of $i$.
If $L = 0$, then $\gamma_i$ is a sequence of points in $Y$,
and its subsequence converges to a point of $Y$.
Hence a subsequence of $v_i$ uniformly converges to a $J$-holomorphic map $v_\infty$
on a neighborhood of this nodal point in $\Sigma'$ if we change each $v_i$ by
$\R$-translation,
where uniform convergence means that the $L^\infty$-distance between
$v_i$ and $v_\infty$ on the intersection of the domain of $v_i$ and the neighborhood of
the nodal point converges to zero.

If $L > 0$,
then a subsequence of $\varphi_i$ converges to some $\varphi \in S^1$.
We blow up these nodal points $q$ of $\Sigma'$ by $\varphi$'s and denote the new curve
by $\Sigma$.
Then it is easy to see that there exists a $J$-holomorphic map $v_\infty$
on a neighborhood $\widetilde{N}^\circ$ of each of these new joint circles to
$(\R \cup \R) \times Y$
such that for each $i$, there exists an $\R$-gluing $\theta_i : \overline{\R} \sqcup
\overline{\R} \to \R$ such that $L^\infty$-distance of $u_i$ and
$(\theta_i \times 1) \circ v_\infty \circ \Psi|_{\widetilde{P}_{x_i}}$ converges to zero
as $i$ goes to infinity,
where $(\widetilde{P} \to \widetilde{X}, Z \cup Z^+)$ is the local universal family
of the blown up curve $(\Sigma, z \cup z^+)$,
$\Psi$ is the discontinuous map for this local universal family,
and $x_i \in \widetilde{X}$ is the point whose fiber is isomorphic to
$(\Sigma_i, z_i \cup z_i^+)$ for each $i$.

Let $\Sigma \setminus \coprod_{\text{imaginary circles}} S^1 = \mathring{\Sigma}^\nu$
be the decomposition into the connected components.
We have already seen that on each closure $\Sigma^\nu
= \overline{\mathring{\Sigma}^\nu}$,
$v_i$ converges to a $J$-holomorphic curve $v_\infty$ if we change each $v_i$ by
$\R$-translation.
(But these $\R$-translations may depend on $\Sigma^\nu$.)

We may assume that each $\Sigma^\nu$ contains some marked points by the following
argument. For each $\Sigma^\nu$ which does not contain any marked points,
we take a holomorphic section $Z' : \widetilde{X} \to \widetilde{P}$ which intersects with
$\Sigma^\nu$,
and we add $z'_i = Z'(x_i) \subset \widetilde{P}_{x_i} \cong \Sigma_i$ as
an additional marked point for each $i$.
Let $I_i^{++}$ be the union of intervals $[k, k+ 1] \subset \R_1 \cup \R_2 \cup \dots
\cup \R_k$ which contains $\sigma \circ u_i (z'_i)$.
We further add marked points $z_i^{++}$ to $\Sigma_i$ as in Lemma
\ref{nice marked points}, that is,
take a finite covering of $I_i^{++} \times Y$ by open balls $B_\alpha^1$
with radius $\epsilon$ and add a marked point for each connected component of
$u_i^{-1}(B_\alpha^2)$ which contains a point of $u^{-1}(B_\alpha^1)$.

We can do the same argument using $(\Sigma_i, z_i \cup z_i^+ \cup z'_i \cup z_i^{++})$
instead of $(\Sigma_i, z_i \cup z_i^+)$, and we get a curve
$(\Sigma^{++}, z \cup z^+ \cup z' \cup z^{++})$ instead of $(\Sigma, z \cup z^+)$.
Then it is clear that $(\Sigma, z \cup z^+ \cup z')$ is obtained by collapsing all unstable
component of $(\Sigma^{++}, z \cup z^+ \cup z')$.
We claim that each connected component of $\Sigma^{++} \setminus
\coprod_{\text{imaginary circles}} S^1$ contains some marked points.
This can be seen as follows.

First we show that every irreducible component $\Sigma^{++}_\alpha$ of $\Sigma^{++}$
which contains at least one imaginary circle and will be collapsed when we forget marked
points $z^{++}$ is a cylinder with at least one additional marked point $z^{++}$ and without
any marked points $z \cup z^+ \cup z'$ or any nodal points.
Such a component $\Sigma^{++}_\alpha$ is either a closed disc or a cylinder,
but the former cannot be occur because if it did, then
$(\Sigma^{++}_\alpha, z \cup z^+)$ would be a closed disc $\C \cup S^1_\infty$
with at most one marked point,
hence the annulus condition for the marked points $z_i \cup z^+_i$ in Lemma
\ref{nice marked points} would imply the diameter of the image of $\Sigma^{++}_\alpha$
by $v_\infty$ is $\leq 2\epsilon$, which is a contradiction.
(We assume that $2\epsilon$ is smaller than the minimal diameter of periodic orbits.)
Hence $\Sigma_\alpha^{++}$ is a cylinder which does not contain any marked points
$z_i \cup z_i^+$ or any nodal points, which implies that $\Sigma^{++}_\alpha$ contains
at least one additional marked point $z^{++}$.

Using this, we can prove each connected component of
$\Sigma^{++} \setminus \coprod S^1$ contains some marked points $z \cup z^+ \cup z'
\cup z^{++}$.
Indeed, if one connected component of $\Sigma^{++} \setminus \coprod S^1$ did not
contain any marked points, then its closure does not collapse to a imaginary circle in
$(\Sigma, z \cup z^+ \cup z')$ and the corresponding component of
$\Sigma \setminus \coprod S^1$ would not contain any marked points
$z \cup z^+ \cup z'$, but this contradicts the choice of $z'$.
Therefore, rewriting $I_i \cup I_i^{++}$ as $I_i$, and $z_i^+ \cup z'_i \cup z^{++}$ as $z^+_i$,
we may assume each connected component $\mathring{\Sigma}^\nu$ of
$\Sigma \setminus \coprod S^1$ contains at least one marked point.

Let $I_i = I_i^1 \cup I_i^2 \cup \dots \cup I_i^{l_i}$ be the decomposition into connected
components for each $i$.
We define an equivalence relation on the set $\{I_i^a\}_{1 \leq a \leq l_i}$
for sufficiently large $i$ as follows
and use these equivalence classes as floors.
Let $\tilde I$ be the $\frac{1}{3}$-neighborhood of $I$ for each interval $I$.
First note the following:
\begin{itemize}
\item
For each $\Sigma^\nu$, $\diam \, \sigma \circ v_i((z\cup z^+) \cap \Sigma^\nu)$ is
bounded uniformly with respect to $i$.
This can be seen by covering a path from one marked point to another by a finite
number of discs in $\Sigma^\nu$ and using the annulus condition for these discs.
Therefore there exists a constant $C>0$ such that if
$v_i((z\cup z^+) \cap \Sigma^\nu)$ intersects with both of
$\tilde I_i^a$ and $\tilde I_i^b$ then $\dist (\tilde I_i^a, \tilde I_i^b) \leq C$.
\item
If $\Sigma^\nu$ and $\Sigma^{\nu'}$ are connected by joint circles in $\Sigma$,
then there exist $a_i \in \{1,2,\dots, l_i\}$ for all large $i$ such that
\begin{itemize}
\item
$v_i((z\cup z^+) \cap \Sigma^\nu)$ intersects with $\tilde I_i^{a_i}$, and
$v_i((z\cup z^+) \cap \Sigma^{\nu'})$ intersects with $\tilde I_i^{a_i + 1}$
(or the condition in which the order of $\nu$ and $\nu'$ is changed is satisfied) and
\item
$\dist (\tilde I_i^{a_i}, \tilde I_i^{a_i +1}) > 2C$ and
$\dist (\tilde I_i^{a_i}, \tilde I_i^{a_i +1}) \to \infty$ as $i \to \infty$.
\end{itemize}
This is because of the asymptotic behavior of $v_i$ on a neighborhood of a joint circle.
\item
For any $\tilde I_i^a$ and $\tilde I_i^b$ ($a < b$),
either of the following two occurs:
\begin{itemize}
\item
There exists $\Sigma^\nu$ such that $v_i((z\cup z^+) \cap \Sigma^\nu)$ intersects with
both of $\tilde I_i^a$ and $\tilde I_i^b$.
\item
There exist $a \leq c < b$ and a pair $\Sigma^\nu$ and $\Sigma^{\nu'}$ connected by
joint circles in $\Sigma$ such that $v_i((z\cup z^+) \cap \Sigma^\nu)$ intersects with
$\tilde I_i^c$ and $v_i((z\cup z^+) \cap \Sigma^{\nu'})$ intersects with $\tilde I_i^{c + 1}$.
\end{itemize}
This is proved as follows.
Since $\Sigma$ is connected, it is easy to see that there exist two marked points
$w_i^a, w_i^b \in z \cup z^+$ such that $v_i(w_i^a) \in \tilde I_i^a$,
$v_i(w_i^b) \in \tilde I_i^b$ and a path $\ell$ in $\Sigma$ from $w_i^a$ to $w_i^b$
such that $v_i(\ell)$ does not intersect with $\tilde I_i^{a-1}$ or $\tilde I_i^{b+1}$.
If $\ell$ intersects with some joint circles then the latter holds and
otherwise the former holds.
\end{itemize}

Therefore, for sufficiently large $i$ and any $\tilde I_i^a$ and $\tilde I_i^b$,
either $\dist (\tilde I_i^a, \tilde I_i^b) \leq C$ or $\dist (\tilde I_i^a, \tilde I_i^b) > 2C$.
Hence we can define the equivalence relation $\sim$ on the set of intervals
$\{\tilde I_i^1, \tilde I_i^2, \dots, \tilde I_i^{l_i}\}$ by
$\tilde I_i^a \sim \tilde I_i^b$ if $\dist (\tilde I_i^a, \tilde I_i^b) \leq C$,
and the set of the equivalent classes has a natural total order.

Fix one large $i$.
Then we can define the floor of each $\Sigma^\nu$ as the equivalence class of
$\tilde I_i^a$ with which $\sigma \circ v_i((z\cup z^+) \cap \Sigma^\nu)$ intersects.
Then for any two components $\Sigma^\nu$ and $\Sigma^{\nu'}$ connected by
some joint circles in $\Sigma$, which of the two has a higher floor is independent of
the choice of $i$ and the difference is one.
Hence we have defined the floor structure of $\Sigma$ independently of $i$.

For each $i$ and floor $j \in \{1,2,\dots l\}$ represented by $I_i^a$,
take one point $b_i^j$ of $I_i^a$.
Define an $\R$-gluing $\theta_i : \overline{\R}_1 \sqcup \overline{\R}_2 \sqcup \dots
\sqcup \overline{\R}_l \to \overline{\R}_1 \cup \overline{\R}_2 \cup \dots \cup
\overline{\R}_{l_i}$ by $\theta_i (0_j) = b_i^j$.
Then it is easy to see that a subsequence of $(\theta_i \times 1)^{-1} \circ u_i \circ
(\Psi|_{\widetilde{P}_{a_i}})^{-1}$
converges to a $J$-holomorphic map
$u_\infty : \Sigma \to
(\overline{\R}_1 \cup \overline{\R}_2 \cup \dots \cup \overline{\R}_l) \times Y$, that is,
\[
\dist_{L^\infty} (u_i, (\theta_i \times 1) \circ u_\infty \circ \Psi|_{\widetilde{P}_{a_i}})
\to 0
\]
as $i \to \infty$.

Finally, passing to a subsequence, we may assume that the sequence
$\pi_{S^1}^{\breve N_{\pm\infty_i}} \circ \phi_{i, \pm\infty_j} : S^1 \to S^1_{\pm\infty_j}$
converges to a family of coordinates $\phi_{\pm\infty_j} : S^1 \to S^1_{\pm\infty_j}$.

The constructed curve $(\Sigma, z, u_\infty, \phi)$ often has unstable components or
floors which consist of trivial cylinders.
Hence we first collapse the unstable components of $(\Sigma, z, u_\infty, \phi)$
(the components $\Sigma_\alpha$ on which $u_\infty$ is constant
and $2g_\alpha + m_\alpha <3$).
Next we collapse all the floors which consist of trivial cylinders.
Then it is clear that $(\Sigma, z_i, u_i, \phi_i)$ converges to this holomorphic building in
the topology of $\overline{\M}^0_{g, \mu}(L^-,L^+)$.
\end{proof}

\begin{prop}\label{Hausdorffness}
$\overline{\M}^0$ is Hausdorff.
\end{prop}
\begin{proof}
The proof is the same as the case of Gromov-Witten theory in \cite{FO99}.
Assuming a sequence $(\Sigma_i, z_i, u_i, \phi_i)\in \overline{\M}^0_{g, \mu}(L^-,L^+)$
converges to two points $(\Sigma, z, u, \phi)$ and $(\Sigma', z', u', \phi')$
in the topology of $\overline{\M}^0_{g, \mu}(L^-,L^+)$,
we prove that these two points coincide.

Let $z^+ \subset \Sigma$ be additional points which make $(\Sigma, z \cup z^+)$ stable,
and let $(\widetilde{P} \to \widetilde{X}, Z \cup Z^+)$ be the local universal family of
$(\Sigma, z \cup z^+)$.
Then by the definition of the topology,
there exists a sequence of points $x_i \to 0 \in \widetilde{X}$ and
a sequence of $\R$-gluings $\theta_i$ such that
$(\Sigma_i, z_i) \cong (\widetilde{P}_{x_i}, Z(x_i))$,
\[
\dist_{L^\infty} (u_i, (\theta_i \times 1) \circ u \circ \Psi|_{\widetilde{P}_{x_i}}) \to 0
\]
and
\[
\dist_{L^\infty(S^1)}(\pi_{S^1}^{\widetilde{N}_{\pm\infty_j}} \circ \phi_{i, \pm\infty_j},
\phi_{\pm\infty_j}) \to 0.
\]
Define additional marked points $z_i^+ = Z^+(x_i) \subset \Sigma_i$.

Similarly, let ${z'}^+ \subset \Sigma'$ be additional points which make
$(\Sigma', z' \cup {z'}^+)$ stable, and let $(\widetilde{P}' \to \widetilde{X}', Z' \cup {Z'}^+)$
be the local universal family of $(\Sigma', z' \cup {z'}^+)$.
Then there exists a sequence of points $x'_i \to 0 \in \widetilde{X}'$
and a sequence of $\R$-gluings $\theta'_i$ such that
$(\Sigma'_i, z'_i) \cong (\widetilde{P}'_{x'_i}, Z'(x'_i))$,
\[
\dist_{L^\infty} (u_i, (\theta'_i \times 1) \circ u \circ \Psi'|_{\widetilde{P}'_{x'_i}}) \to 0
\]
and
\[
\dist_{L^\infty(S^1)}(\pi_{S^1}^{\breve N_{\pm\infty_j}} \circ \phi_{i, \pm\infty_j},
\phi'_{\pm\infty_j}) \to 0.
\]
Define additional marked points ${z'}_i^+ = {Z'}^+(x'_i) \subset \Sigma_i$.

We may assume $\pi_Y \circ u (z^+)$ and $\pi_Y \circ u ({z'}^+)$ do not share any points.
Then $z_i^+$ and ${z'}_i^+$ are disjoint for large $i$.

Starting a holomorphic building $(\Sigma_i, z_i, u_i)$ with additional marked points
$z_i^+ \cup {z'}_i^+$, we further add marked points $z_i^{++}$ by the procedure we
explained in the proof of Proposition \ref{compactness}.
Passing to a subsequence if necessary, there exists a holomorphic building
$(\Sigma'',z \cup z^+ \cup {z'}^+ \cup z^{++}, u'')$ which satisfies the following condition.
Let $(\widetilde{P}'' \to \widetilde{X}'', Z \cup Z^+ \cup {Z'}^+ \cup Z^{++})$
be the local universal family of
$(\Sigma'', z \cup z^+ \cup {z'}^+ \cup z^{++})$.
Then there exists a sequence of points $x''_i \to 0 \in \widetilde{X}''$
and a sequence of $\R$-gluings $\theta''_i$ such that
$(\Sigma_i, z_i \cup z_i^+ \cup {z'}_i^+ \cup z_i^{++}) \cong (\widetilde{P}''_{x''_i},
Z(x''_i) \cup Z^+(x''_i) \cup {Z'}^+(x''_i) \cup Z^{++}(x''_i))$,
\[
\dist_{L^\infty} (u_i, (\theta''_i \times 1) \circ u'' \circ \Psi''|_{\widetilde{P}_{x''_i}}) \to 0
\]
and
\[
\dist_{L^\infty(S^1)}(\pi_{S^1}^{\breve N_{\pm\infty_i}} \circ \phi_{i, \pm\infty_j},
\phi''_{\pm\infty_j}) \to 0.
\]

Since the space of stable curves are Hausdorff and the forgetful map is continuous,
the stabilization of $(\Sigma'', z \cup z^+)$ is $(\Sigma, z \cup z^+)$.
Since the forgetful map $(\widetilde{P}'', \widetilde{X}'') \to (\widetilde{P}, \widetilde{X})$
maps $x''_i$ to $x_i$, $u$ and $\phi$ are the maps induced by $u''$ and $\phi''$.
Since the same is true for $(\Sigma', z', u')$, the two holomorphic buildings
$(\Sigma, z, u)$ and $(\Sigma', z', u')$ coincide.
\end{proof}

We also use the following quotient space
$\widehat{\M}^0(Y, \lambda, J) = \overline{\M}^0(Y, \lambda, J) / \sim$.
This space is obtained by ignoring the coordinates of limit circles and the order of
marked points and limit circles, that is, in $\widehat{\M}^0(Y, \lambda, J)$,
we identify two holomorphic buildings $(\Sigma, z, u, \phi)$ and $(\Sigma', z', u', \phi')$
if there exist a biholomorphism $\varphi : \Sigma' \to \Sigma$ and an $\R$-translation
$\theta$ such that $\varphi(\{z'_i\}) = \{z_i\}$ (that is, $\varphi$ maps $\{z'_i\}$ to
$\{z_i\}$ as a set) and $u' = (\theta \times 1) \circ u \circ \varphi$.
Hence we may write a point of $\widehat{\M}^0(Y, \lambda, J)$ as $(\Sigma, z, u)$,
where $z$ is a set of points of $\Sigma$.
$\widehat{\M}^0$ is also second countable and Hausdorff, and each
$\widehat{\M}^0_{g, \mu}(L^-, L^+)$ is compact because
$\widehat{\M}^0$ is a quotient space of a subspace of $\overline{\M}^0$
by a compact group locally.


Recall that $\overline{\M} = \overline{\M}(Y, \lambda, J)$ is the space of
all (possibly disconnected) holomorphic buildings without trivial buildings.
This space is decomposed by the number of the connected component of the domain
curve.
We can define the topology of each of them similarly and prove the second countability,
compactness and Hausdorff property as $\overline{\M}^0$.
The compactness is stated as follows, where the genus $g$ of
a disconnected holomophic building $(\Sigma, z, u, \phi)$ is defined
by $g = 1 - \frac{1}{2} \chi(\check \Sigma) \in \Z$ ($\chi(\check \Sigma)$ is
the Euler number of the curve $\check \Sigma$).
\begin{prop}\label{compactness for disconnected}
For any $g_0 \in \Z$, $\mu_0 \geq 1$ and $L^+_0 \in \R$,
\[
\bigcup_{\substack{-\infty < g \leq g_0 \\ \mu \leq \mu_0 \\ L^- \leq L^+ \leq L^+_0}}
\overline{\M}_{g, \mu}(L^-, L^+)
\]
is compact.
\end{prop}
\begin{proof}
It is enough to prove that the number of the connected components of the domain
curve of a holomorphic building in the above space is bounded by some constant
depending only on $g_0$, $\mu_0$ and $L^+_0$.
The number of the connected components which have $+\infty$-limit circles is bounded,
and so are the number of the components with marked points.
Since the other components are constant maps, each of them have genus $\geq 2$.
Therefore, the number of them is also bounded.
(Note that the genus of the curve is $g = 1 + \sum_i (g_i - 1)$, where $g_i$ are the
genera of connected components.)
\end{proof}

We define the quotient space  $\widehat{\M} = \overline{\M} / \sim$ similarly.

\subsection{The case of manifolds with cylindrical ends}
\label{holomorphic buildings for X}
Next we consider the holomorphic buildings for a symplectic manifold $X$ with cylindrical
ends.
In this case, floor takes values in $\{ -k_-, -k_- +1, \dots, k_+\}$.
\begin{defi}
A holomorphic building $(\Sigma, z, u, \phi)$ for $X$ consists of
\begin{itemize}
\item
a marked curve $(\Sigma, z)$ which is obtained from
a union of marked semistable curves
$(\check \Sigma, z \cup (\pm \infty_i))$ with a floor structure,
\item
a continuous map $u : \Sigma \to (\overline{\R}_{-k_-} \cup \dots \cup
\overline{\R}_{-1}) \times Y^- \cup \overline{X} \cup (\overline{\R}_1 \cup \dots \cup
\overline{\R}_{k_+}) \times Y^+$ and
\item
a family of coordinates $\phi_{\pm\infty_i} : S^1 \to S^1_{\pm\infty_i}$ of limit circles
\end{itemize}
which satisfy the following conditions:
\begin{itemize}
\item
If $i(\alpha) <0$ then
$u(\Sigma_\alpha \setminus \coprod S^1) \subset \R_{i(\alpha)} \times Y^-$, and
$u|_{\Sigma_\alpha \setminus \coprod S^1} : \Sigma_\alpha \setminus \coprod S^1 \to \R_{i(\alpha)} \times Y$ is $J$-holomorphic
\item
If $i(\alpha) = 0$ then
$u(\Sigma_\alpha \setminus \coprod S^1) \subset X$, and
$u|_{\Sigma_\alpha \setminus \coprod S^1} : \Sigma_\alpha \setminus \coprod S^1 \to X$ is $J$-holomorphic
\item
If $i(\alpha) >0$ then
$u(\Sigma_\alpha \setminus \coprod S^1) \subset \R_{i(\alpha)} \times Y^+$, and
$u|_{\Sigma_\alpha \setminus \coprod S^1} : \Sigma_\alpha \setminus \coprod S^1 \to \R_{i(\alpha)} \times Y^+$ is $J$-holomorphic
\item
$E_\lambda(u) <\infty$ and $E_{\hat \omega}(u) <\infty$,
where these energies are defined by
\begin{align*}
E_\lambda(u) &= \max \biggl\{ \sup_{I \subset \R_{-k_-} \cup \dots \cup \R_{-1} \cup (-\infty, 0]} \frac{1}{|I|} \int_{u^{-1}(I \times Y^-)} u^\ast (d\sigma \wedge \lambda^-),\\
&\hphantom{= \max \biggl\{} \sup_{I \subset [0, \infty) \cup \R_1 \cup \dots \cup \R_{k_+}} \frac{1}{|I|} \int_{u^{-1}(I \times Y^+)} u^\ast (d\sigma \wedge \lambda^+) \biggr\}\\
E_{\hat \omega}(u) &= \int_{u^{-1}(X)} u^\ast \hat \omega
+ \int_{u^{-1}((\overline{\R}_{-k_-} \cup \dots \cup \overline{\R}_{-1}) \times Y^-)} u^\ast d\lambda^-\\
&\quad + \int_{u^{-1}((\overline{\R}_1 \cup \dots \cup \overline{\R}_{k_+}) \times Y^+)} u^\ast d\lambda^+.
\end{align*}
\item
$u$ is positively asymptotic to a periodic orbit $\gamma_{+\infty_i} = \pi_Y \circ u \circ \phi_{+\infty_i} \in P_{Y^+}$ at each $S^1_{+\infty_i}$,
and negatively asymptotic to a periodic orbit $\gamma_{-\infty_i} = \pi_Y \circ u \circ \phi_{-\infty_i} \in P_{Y^-}$ at each $S^1_{-\infty_i}$.
At every joint circle, $u$ is positively asymptotic to a periodic orbit
on the side of lower floor and negatively asymptotic to the same periodic orbit
on the side of higher floor.
\item
For each component $\hat\Sigma_\alpha$, if $u|_{\Sigma_\alpha}$ is a constant map,
then $2g_\alpha + m_\alpha \geq 3$.
\item
For each $i \neq 0$, $i$-th floor $u^{-1}(\overline{\R}_i \times Y^\pm) \subset
\Sigma$ contains nontrivial components.
\end{itemize}
\end{defi}

We say two holomorphic buildings $(\Sigma, z, u, \phi)$ and $(\Sigma', z', u', \phi')$ are
isomorphic if
there exist a biholomorphism $\varphi : \Sigma' \to \Sigma$
and a pair of $\R$-translations $\theta^- : \overline{\R}_{-k_-} \cup \dots \cup
\overline{\R}_{-1} \to \overline{\R}_{-k_-} \cup \dots \cup \overline{\R}_{-1}$ and
$\theta^+ : \overline{\R}_1 \cup \dots \cup \overline{\R}_{k_+}
\to \overline{\R}_1 \cup \dots \cup \overline{\R}_{k_+}$
such that
\begin{itemize}
\item
$\varphi(z'_i) = z_i$ for all $i$,
\item
$u' = (\theta \times 1) \circ u \circ \varphi$,
where
\begin{align*}
&\theta \times 1 : (\overline{\R}_{-k_-} \cup \dots \cup \overline{\R}_{-1})
\times Y^- \cup \overline{X} \cup (\overline{\R}_1 \cup \dots \cup \overline{\R}_{k_+})
\times Y^+ \\
&\hph{\theta \times 1 :}
\to
(\overline{\R}_{-k_-} \cup \dots \cup \overline{\R}_{-1}) \times Y^- \cup \overline{X}
\cup (\overline{\R}_1 \cup \dots \cup \overline{\R}_{k_+}) \times Y^+
\end{align*}
is defined by $(\theta \times 1)|_{\overline{X}} = \id_{\overline{X}}$ and
$(\theta \times 1)|_{\overline{\R_i} \times Y^\pm} = \theta^\pm \times 1$, and
\item
$\varphi \circ \phi'_{\pm\infty_i} = \phi_{\pm\infty_i}$ for all $\pm\infty_i$.
\end{itemize}

Note that the $0$-th floor of a holomorphic building may be empty.
We regard the empty curve, that is, the holomoprphic curve whose domain is
the empty set, as a disconnected holomorphic building for $X$, but we do not regard
it as a connected holomorphic building.
The genus of the empty curve is defined by $1$ ($= 1 - \frac{1}{2} \chi(\emptyset)$).
We denote the space of all holomorphic buildings for $X$ by
$\overline{\M}(X, \omega, J)$, and the space of connected ones
by $\overline{\M}^0(X, \omega, J)$.

The neighborhoods of each point
$p_0 = (\Sigma_0, z_0, u_0, \phi_0) \in \overline{\M}^0(X, \omega, J)$ is defined as follows.
As in the case of $\overline{\M}^0(Y, \lambda, J)$, first we add marked points $z_0^+$ to
$(\Sigma_0, z_0)$ to make $(\Sigma_0, z_0 \cup z_0^+)$ stable.
Let $(\widetilde{P} \to \widetilde{X}, Z \cup Z^+)$ be the local universal family of
$(\Sigma_0, z_0 \cup z_0^+)$.

For a pair $\theta = (\theta^-, \theta^+)$ of $\R$-gluings
$\theta^- : \overline{\R}_{-k_-} \sqcup \dots \sqcup \overline{\R}_0
\to \overline{\R}_{-l_-} \cup \dots \cup \overline{\R}_0$ and
$\theta^+ : \overline{\R}_0 \sqcup \dots \sqcup \overline{\R}_{k_+}
\to \overline{\R}_0 \cup \dots \cup \overline{\R}_{l_+}$,
we define a map
\begin{align*}
&(\theta \times 1) : (\overline{\R}_{-k_-} \sqcup \dots \sqcup \overline{\R}_{-1})
\times Y^- \sqcup \overline{X} \sqcup (\overline{\R}_1 \sqcup \dots \sqcup
\overline{\R}_{k_+}) \times Y^+ \\
&\hph{(\theta \times 1) :}
\to (\overline{\R}_{-l_-} \cup \dots \cup \overline{\R}_{-1}) \times Y^- \cup \overline{X}
\cup (\overline{\R}_1 \cup \dots \cup \overline{\R}_{l_+}) \times Y^+
\end{align*}
by
\begin{itemize}
\item
$(\theta \times 1)|_{\overline{X}} = \id$
\item
$(\theta \times 1)|_{\overline{\R} \times Y^\pm} = \theta^\pm \times 1$
if $\mu(i) \neq 0$.
(Recall $\mu$ is defined by $\theta^\pm(\overline{\R}_i) = \overline{\R}_{\mu(i)}$.)
\item
For each $i<0$ such that $\mu(i) = 0$,
$(\theta \times 1)(\sigma, y) = ( \min(\theta(\sigma), 0), y) \in (-\infty, 0] \times Y^- \subset X$
\item
For each $i>0$ such that $\mu(i) = 0$,
$(\theta \times 1)(\sigma, y) = ( \max(\theta(\sigma), 0), y) \in [0,\infty) \times Y^+ \subset X$
\end{itemize}

For a neighborhood $U \subset \widetilde X$ and $\epsilon>0$,
$\W_{p_0}(U, \epsilon) \subset \overline{\M}^0(X, \omega, J)$ is defined
as follows.
$(\Sigma, z, u) \in \overline{\M}^0(X, \omega, J)$ belongs to $\W_{p_0}(U, \epsilon)$ if
there exist a point $a\in U$, an isomorphism $(\Sigma, z) \cong (\widetilde{P}_a, Z(a))$
and a pair of $\R$-gluings $\theta = (\theta^-, \theta^+)$ such that
\[
\dist_{L^\infty}(u, (\theta \times 1) \circ u_0 \circ \Psi|_{\widetilde{P}_a}) < \epsilon
\]
and
\[
\dist_{L^\infty(S^1)}(\pi_{S^1}^{\breve N_{\pm\infty_i}} \circ \phi_{\pm\infty_i},
\phi_{0, \pm\infty_i}) < \epsilon.
\]
We define a neighborhood of $p_0$ as a subset of $\overline{\M}^0(X, \omega, J)$
which contains some $\W_{p_0}(U, \epsilon)$.
This defines the topology of $\overline{\M}^0(X, \omega, J)$ similarly
to the case of $\overline{\M}^0(Y, \lambda, J)$.

Define a closed two form $\widetilde{\omega}$
on $X = (-\infty,0] \times Y^- \cup Z \cup [0,\infty) \times Y^+$ by
$\widetilde{\omega}|_Z = \omega$,
$\widetilde{\omega}|_{(-\infty, 0] \times Y^-} = d(\varphi \lambda^-)$ and 
$\widetilde{\omega}|_{[0, \infty) \times Y^+} = d(\varphi \lambda^+)$,
where $\varphi :\R \to \R_{\geq 0}$ is a smooth function with compact support
such that $\varphi(0) = 1$.

Then $\overline{\M}^0(X, \omega, J)$ is decomposed as
$\overline{\M}^0(X, \omega, J) = \coprod \overline{\M}^{0, e}_{g, \mu}(L^-, L^+)$,
where $\overline{\M}^{0, e}_{g, \mu}(L^-, L^+)$ consists of holomorphic buildings
$(\Sigma, z, u, \phi)$ with genera $g$ and $\mu$ marked points such that
$\sum L_{\gamma_{-\infty_i}} = L^-$, $\sum L_{\gamma_{+\infty_i}} = L^+$
and $\int_{u^{-1}(X)} u^\ast \widetilde{\omega} = e$.
(This is independent of the choice of the function $\varphi$.)
Note that $(\Sigma, z, u, \phi) \in \overline{\M}^{0, e}_{g, \mu}(L^-, L^+)$ satisfies
\begin{align}
E_\lambda (u) &\leq \max(e + L^+, L^+)
\label{E lambda estimate}\\
E_{\hat \omega} (u) &= e + (L^+ - L^-)
\label{E hat omega estimate}
\end{align}
(\ref{E lambda estimate}) is because
\begin{itemize}
\item
for any interval $I \subset \R_i$ ($i < 0$) or $I \subset (-\infty, 0]$,
\begin{align*}
\frac{1}{|I|}\int_{u^{-1}(I \times Y^-)} u^\ast (d\sigma \wedge \lambda^-)
&= \int u^\ast (d\varphi^- \wedge \lambda^-)\\
&= \int_{u^{-1}(\hat Y^-)} u^\ast d(\varphi^- \lambda^-)
- \int_{u^{-1}(\hat Y^-)} u^\ast(\varphi^- d\lambda^-)\\
&\leq \int_{u^{-1}(\hat Y^-)} u^\ast d(\varphi^- \lambda^-)\\
&= \int_{u^{-1}((-\infty, 0] \times Y^-)} u^\ast d(\varphi \lambda^-)\\
&\leq \int_{u^{-1}((-\infty, 0] \times Y^-)} u^\ast d(\varphi \lambda^-)
+ \int_{u^{-1}(Z)} u^\ast \omega\\
& \quad + \int_{u^{-1}(\hat Y^+)} u^\ast d\lambda^+\\
&= \int_{u^{-1}(X)} u^\ast \widetilde{\omega}
+ \int_{u^{-1}(\hat Y^+)} u^\ast d((1-\varphi)\lambda^+)\\
&= e + L^+,
\end{align*}
where $\hat Y^- = (\overline{\R}_{-k_-} \cup \dots \cup \overline{\R}_{-1}
\cup (-\infty, 0]) \times Y^-$,
$\hat Y^+ = ([0, \infty) \cup \overline{\R}_1 \cup \dots \cup \overline{\R}_{k_+})
\times Y^+$,
and $\varphi^-$ is defined by
\[
\varphi^-(\sigma, y) = \int_{-\infty}^\sigma \frac{1}{|I|} 1_I(\sigma') d\sigma'
\quad \text{on } \R_i \times Y^-
\]
$\varphi^- \equiv 0$ on $\R_j \times Y^-$ ($j < i$),
and $\varphi^- \equiv 1$ on $\R_j \times Y^-$ ($j < i$), $X$, and $\R_j \times Y^+$,
\item
for any interval $I \subset \R_i$ ($i > 0$) or $I \subset [0, \infty)$,
\begin{align*}
\frac{1}{|I|}\int_{u^{-1}(I \times Y^+)} u^\ast (d\sigma \wedge \lambda^+)
&= \int u^\ast (d\varphi^+ \wedge \lambda^+)\\
&= \int_{u^{-1}(\hat Y^+)} u^\ast d(\varphi^+ \lambda^+)
- \int_{u^{-1}(\hat Y^+)} \varphi^+ d\lambda^+\\
&\leq L^+,
\end{align*}
where $\varphi^+$ is defined by
\[
\varphi^+(\sigma) = \int_{-\infty}^\sigma \frac{1}{|I|} 1_{I}(\sigma') d\sigma'
\quad \text{on } \R_i,
\]
$\varphi^+ \equiv 0$ on $\R_j \times Y^-$, $X$, and $\R_j \times Y^+$ for ($j < i$),
and $\varphi^+ \equiv 1$ on $\R_j \times Y^+$ ($j > i$).
\end{itemize}
Proof of (\ref{E hat omega estimate}) is straightforward.

As in the case of $\hat Y$, we have a nice way to add marked points
to the domain curves.
\begin{lem}
Let $\epsilon > 0$ and $\delta_0 > 0$ be arbitrary small constants,
and let $(\Sigma, z, u, \phi) \in \overline{\M}^{0, e}_{g, \mu}(L^-, L^+)$ be an arbitrary
holomorphic building.
Then there exist closed subsets
$I^- \subset \R_{-k_-} \cup \dots \cup \R_{-1} \cup (-\infty, 0]$
and $I^+ \subset [0, \infty) \cup \R_1 \cup \dots \cup \R_{k_+}$,
and additional marked points $z^+$ of $\Sigma$ which satisfy the following conditions:
\begin{itemize}
\item
Both of $I^\pm$ are finite unions of intervals in the form $[l, l+1] \subset \R_i$
{\rm(}$l \in \Z${\rm)}.
\item
The lengths of $I^\pm$ and the number of additional marked points are bounded by
some constant determined by $g$, $\mu$, $L^-$, $L^+$, $e$, $\epsilon$ and $\delta_0$.
\item
$(\Sigma, z \cup z^+)$ is stable.
\item
There exists a constant $A_1 > 0$ depending only on $g$, $\mu$, $L^-$, $L^+$, $e$,
$\epsilon$ and $\delta_0$ such that
if $[-A_1, T + A_1] \times S^1 \subset \Sigma$ does not contain any marked points
$z \cup z^+$, then one of the following two holds true.
\begin{enumerate}[label=\normalfont(\arabic*)]
\item
$u([0, T] \times S^1) \subset I^- \times Y^- \cup Z \cup I^+ \times Y^+$ and
$\diam \, u([0, T] \times S^1) \leq 20 \epsilon$.
\item
$u ([0, T] \times S^1)$ is contained in $J^- \times Y^-$ or $J^+ \times Y^+$,
where $J^-$ is the $\frac{1}{3}$-neighborhood of
the complement of $I^- \subset \R_{-k_-} \cup \dots \cup \R_{-1} \cup (-\infty, 0]$,
and $J^+$ is the $\frac{1}{3}$-neighborhood of
the complement of $I^+ \subset [0, \infty) \cup \R_1 \cup \dots \cup \R_{k_+}$.
Moreover, in this case,
there exist $L \in \R$ and $(b, \gamma) \in \R \times P_L^\pm$ such that
\[
\dist (u(s,t), (Ls + b, \gamma(t))) \leq \epsilon (e^{-\kappa s} + e^{-\kappa (T - s)})
\]
on $[0, T] \times S^1$.
\end{enumerate}
In particular,
for any disc $D \subset \Sigma$ such that
$D \setminus 0$ does not contain any marked points,
$\diam \, u(\{ z\in D; |z| \leq e^{-2\pi A_1}\}) \leq 20 \epsilon$.
\item
$u(z \cup z^+)$ is contained in the $\frac{1}{3}$-neighborhood of
$I^- \times Y^- \cup Z \cup I^+ \times Y^+$.
\item
Each connected component of $u^{-1}(I^- \times Y^- \cup Z \cup I^+ \times Y^+)$ either
contains at least one point of $z \cup z^+$
or is contained in the inverse image of the $\frac{1}{3}$-neighborhood of
the complement of $I^- \times Y^- \cup Z \cup I^+ \times Y^+$ by $u$.
\item
For the $\frac{1}{3}$-neighborhood $\widetilde{J}$ of each connected component of
the complement of $I^- \times Y^- \cup Z \cup I^+ \times Y^+$,
$E_{\hat \omega}(u|_{u^{-1}(\widetilde{J})}) \leq \delta_0$.
\end{itemize}
\end{lem}
\begin{proof}
First we show the number of irreducible components of $\check \Sigma$ is bounded by
some constant depending only on $g$, $\mu$, $L^-$, $L^+$ and $e$.

Define an energy $E'(u)$ by
\begin{align*}
E'(u) &= \int_{u^{-1}(Z)} u^\ast \omega
+ \int_{u^{-1}([-1, 0] \times Y^-)} u^\ast (d\sigma \wedge \lambda^- + d\lambda^-)\\
&\quad + \int_{u^{-1}([0, 1] \times Y^+)} u^\ast (d\sigma \wedge \lambda^+
+ d\lambda^+).
\end{align*}
This is bounded by
$E'(u) \leq E_{\hat \omega}(u) + 2 E_\lambda (u)$.

If $u$ is non-constant on $\Sigma_\alpha$ and $u(\Sigma_\alpha)$ intersects with $Z$,
then Lemma \ref{monotonicity lemma} implies
$E'(u|_{\Sigma_\alpha})$ is larger than some positive constant independent of $u$.
This implies the number of the irreducible components of $\Sigma$ on which $u$ is not
a constant map and whose image by $u$ intersects with $Z$ is bounded.
Hence as in the case of $\hat Y$,
the number of the nontrivial components $\Sigma_\alpha$ such that
$2g_\alpha + m_\alpha <3$ is bounded.
Therefore the height of $(\Sigma, z, u)$ is bounded,
which implies the number of trivial cylinders is also bounded.

Therefore the number of marked points we need to add to $(\Sigma, z)$ in order to
make $(\Sigma, z \cup z^+)$ stable is bounded. Assuming $(\Sigma, z \cup z^+)$ is
stable, we further add marked points $z^{++}$ as follows.
Let $\delta > 0$, $\kappa > 0$, $A > 0$ and $L_0$ be the constant of
Corollary \ref{third annulus for X} for $C_0 = e + 2L^+$ and the given $\epsilon > 0$.
We may assume $\delta \leq \delta_0$.

Let $I^- \subset \R_{-k_-} \cup \dots \cup \R_{-1} \cup (-\infty, 0]$
and $I^+ \subset [0, \infty) \cup \R_1 \cup \dots \cup \R_{k_+}$ be finite unions of
intervals $[l, l+1] \subset \R_i$ such that
\begin{itemize}
\item
$E_{\hat\omega} (u|_{u^{-1}(\widetilde{J})}) \leq \delta$
for the $\frac{1}{3}$-neighborhood $\widetilde{J}$ of each connected component of
the complement of $I^- \times Y^- \cup Z \cup I^+ \times Y^+$, and
\item
$u(z \cup z^+) \subset I^- \times Y^- \cup Z \cup I^+ \times Y^+$.
\end{itemize}
We may assume the lengths of $I^\pm$ are bounded by some constants
depending only on
$E_{\hat \omega}(u)$, $\delta$ and the number of marked points.

Let $\bigcup_\alpha B_\alpha^1 \supset I^- \times Y^- \cup Z \cup I^+ \times Y^+$
be a finite covering by open balls
with radius $\epsilon$.
We may assume the number of open balls is bounded by some constant depending
on the length of $I$ and $\epsilon$.
For each $B_\alpha^1$, let $B_\alpha^2$ be the concentric ball with radius $2\epsilon$.
Using this finite covering, we add marked points $z^{++}$ as in the proof of
Lemma \ref{nice marked points}.
Then by the similar argument, we can easily see that $I^\pm$ and
the additional marked points $z^+ \cup z^{++}$ satisfy the required conditions.
\end{proof}

Using the above lemma, we can prove the following propositions
similarly to the case of $\overline{\M}^0(Y, \lambda, J)$.
\begin{prop}
$\overline{\M}^0(X, \omega, J)$ is second countable.
\end{prop}
\begin{prop}
For any $g$, $\mu$, $e_0$, $L^-$ and $L^+$,
$\bigcup_{e \leq e_0} \overline{\M}^{0, e}_{g, \mu}(L^-, L^+)$ is compact.
\end{prop}
\begin{prop}
$\overline{\M}^0(X, \omega, J)$ is Hausdorff.
\end{prop}

Similarly to the case of symplectization $\hat Y$, we use the quotient space
\[
\widehat{\M}^0(X, \omega, J) = \overline{\M}^0(X, \omega, J) / \sim
\]
obtained by ignoring the coordinates of limit circles and the order of marked points
and limit circles.
We also define the quotient space $\widehat{\M}(X, \omega, J)
= \overline{\M}(X, \omega, J) / \sim$ similarly.
$\overline{\M}(X, \omega, J) / \sim$ and $\widehat{\M}(X, \omega, J)$
are also second countable and Hausdorff.
The compactness is stated as follows.
\begin{prop}
For any $g_0 \in \Z$, $\mu_0 \geq 0$, $L^1_0 \geq 0$ and $L^2_0 \geq 0$,
\[
\bigcup_{\substack{-\infty < g \leq g_0 \\ \mu \leq \mu_0 \\ e + L^+ \leq L^1_0 \\
L^+ \leq L^2_0}} \overline{\M}^e_{g, \mu} (L^-, L^+)
\]
is compact.
\end{prop}
The proof of the above proposition is almost same with that of Proposition
\ref{compactness for disconnected}.
In this case, in order to prove the boundedness of the number of the connected
components, we use the fact that the energy of a non-constant closed
$J$-holomorphic curve in $X$ is bounded below by some positive constant.


%% file: SFT-00_Theory_of_Kuranishi_structure.tex
%
%
\section{Theory of Kuranishi structure}
\label{theory of Kuranishi structure}
We use the theory of Kuranishi structure for the construction of symplectic field theory.
This theory was developed by Fukaya and Ono in \cite{FO99},
and it is a useful tool to perturb the given equation and get algebraic information of the
moduli space of the solutions.
A neighborhood of each point of the moduli space is usually expressed
as the zero set of a Fredholm map between Banach spaces
or the quotient of the zero set by a group actin.
Since this map is not always transverse to zero, in order to get some
algebraic information, we need to perturb the map and make it transverse to zero.
To get information of the moduli space, the most important thing is
perturbing these maps in a compatible way.
The theory of Kuranishi structure is a scheme to manipulate this compatibility problem.
In this section, we recall this theory and introduce new notions of
pre-Kuranishi space and weakly good coordinate system.

First we explain roughly about what is Kuranishi structure and how we construct it.
As mentioned above, a neighborhood of each point of the moduli space
is expressed as the zero set of a Fredholm map between Banach spaces
or its quotient by a group action.
Adding a finite dimensional vector space to the domain of each Fredholm map
and extending the map to this product space,
we first make each map transverse to zero.
Then the original zero set is the zero set of the projection map from the new zero set
to the added finite dimensional space.
This implies that a neighborhood of each point of the moduli space is expressed
as a zero set of a smooth section of an finite dimensional vector bundle or orbibundle.
(This expression is called a Kuranishi neighborhood.)
This expression depends on the choice of the additional finite dimensional space,
but if one additional space is a subspace of another additional space,
then the former orbibundle can be naturally embedded in the latter.
Kuranishi structure is, roughly speaking, a collection of Kuranishi neighborhoods
with the relation of this kind of embeddings.

\subsection{Orbibundle}
\label{section of orbibundle}
In this subsection, we explain about orbibundle.
First we explain our notation of corners of manifold.
For an open subset $V \subset [0, \infty)^n$, we define the interior of
the corner of codimension $k$
\[
\mathring{\partial}^k V = \{ (x_j) \in [0, \infty)^n;
\# \{j; x_j = 0\} = k\}
\]
and the boundary $\partial V = \bigcup_{k \geq 1} \mathring{\partial}^k V$.
For each point $x \in \mathring{\partial}^k V$,
we define the normal space $T^\bot_x \mathring{\partial}^k V
= T_x \R^n / T_x \mathring{\partial}^k V$.
We say a smooth map $\phi$ from an open subset $V \subset [0, \infty)^n$
to another $V' \subset [0, \infty)^{n'}$ is an embedding if
$\phi$ is the restriction of some embedding from an open subset of $\R^n$ to $\R^{n'}$,
$\phi(\mathring{\partial}^k V) \subset \mathring{\partial}^k V'$ for each $k \geq 0$, and
the differential $\phi_{\ast x} : T^\bot_x \mathring{\partial}^k V
\to T^\bot_{\phi(x)} \mathring{\partial}^k V'$ is an isomorphism
for each $k \geq 1$ and $x \in \mathring{\partial}^k V$.
The definition of the diffeomorphism is similar.
Using these definitions, we define manifold with corners and
embedding between two manifolds.

We also explain our definition of
clean intersection.
Let $(M_\alpha)_{\alpha \in A}$ be a locally finite family of
submanifolds of a manifold $M$.
We say $(M_\alpha)_{\alpha \in A}$ intersect cleanly if
for any $p \in M$, there exist a coordinate $U_p \cong \R^m$ of $p \in M$
such that submanifolds $M_\alpha \cap U_p$
coincides with the subspace $\{(t_i) \in \R^m; t_i = 0 \text{ for } i \in I_{\alpha, p}\}$
for some $I_{\alpha, p} \subset \{1, \dots, m\}$
for all $\alpha \in A$ such that $p \in M_\alpha$.
Note that if there exists a smooth function $f_\alpha$ on $M_\alpha$
for each $\alpha \in A$ and they coincide on their intersections,
then we can extend them to a smooth function on $M$ if we shrink $M_\alpha$
to their arbitrary relatively compact subsets.
(In the case of orbifolds,
we say a family of embeddings intersects cleanly
if all of their lifts intersect cleanly.)

First we show the following elementary lemma for the definition of orbifold and
orbibundle.
\begin{lem}\label{effective action}
Let $V$ and $V'$ be connected manifolds (with or without corners).
Suppose a finite group $G'$ acts on $V'$ effectively.
Then the following hold true.
\begin{enumerate}[label=\normalfont(\roman*)]
\item
If two submersions $\phi, \psi : V \to V'$ induce the same map
$\overline{\phi}= \overline{\psi} : V \to V'/G'$, then
there exists a unique $h\in G' $ such that $\psi = h \phi : V\to V'$.
\item
Assume that two embeddings $\phi, \psi : V \to V'$ induce the same map
$\overline{\phi}= \overline{\psi} : V \to V'/G'$.
If $\phi(V) = \psi(V)$, $\phi(V) \subset V'$ is $G' $-invariant and
the $G'$-action on $\phi (V)$ is effective,
then there exists a unique $h\in G' $ such that $\psi = h \phi : V\to V'$.
\end{enumerate}
\end{lem}
\begin{proof}
(i)
First we claim that the differentials $D\phi$ and $D\psi$ induce
the same map $\overline{D\phi} = \overline{D\psi} : TV \to TV'/G'$.
For any curve $l$ in $V$ there exist some $t_j\to 0$ and $h \in G' $
such that $\phi (l(t_j)) = h \psi (l(t_j))$.
Hence $D\phi(l(0)) \dot{l}(0) = h D\psi(l(0)) \dot{l}(0)$, which implies the claim.

Next we show that for any $p\in V$, there exists unique $h_p \in G' $ such that
$D\psi(p) = h_p D\phi(p)$.
Uniqueness is a consequence of the effectiveness of the action. 
The first claim $\overline{D\phi} = \overline{D\psi}$ implies that
\[
T_{\phi(p)} V' = \bigcup _{g \in G'} \{v \in T_{\phi (p)} V';
gv = D\psi(p) D\phi(p)^{-1} v \}.
\]
(For each $v \in T_{\phi (p)} V'$,
$D\psi(p) D\phi(p)^{-1} v \subset T_{\psi(p)} V'$ is an affine space
which is contained in the orbit $G' v$.
Hence it is a point in $G' v$.)
Since the right hand side of this equation is a finite union of subspaces of $T_{\phi(p)} V'$,
one of them coincides with the whole space.
In other words, there exists some $h_p \in G'$ such that $D\psi(p) = h_p D\phi(p)$.

Since the dimension of the space
\[
\{ v \in T_{\phi (p)} V' ; gv = D\psi(p) D\phi(p)^{-1} v \}
\]
is upper semi-continuous with respect to $p$ for each $g \in G'$,
the uniqueness of $h_p \in G'$ and the connectedness of $V$ imply that
$h = h_p \in G'$ does not depend on $p \in V$.
Therefore $\psi = h \phi$ on $V$.

(ii)
Apply (i) to $\phi, \psi : V \to \phi (V)=\psi (V)$.
\end{proof}
\begin{rem}
\label{effective action for composition of embedding and submersion}
In (ii), we assumed that $\phi$ and $\psi$ are embedding.
However, it is clear that the same holds if $\phi$ and $\psi$ are
submersions to their image $\phi(V) = \psi(V)$.
\end{rem}

\begin{defi}
An orbichart $\V = (V, \pi_V, \V)$ consists of 
a connected manifold $V$ (with or without corners), a topological space $\V$ and
a continuous map $\pi_V : V \to \V$ such that
\begin{itemize}
\item
there exists some finite group $G_V $ acting smoothly and effectively on $V$
\item
$\pi_V$ induces a homeomorphism $\pi_V : V/G_V \tocong \V$
\item
if $\mathring{\partial}^k V \neq \emptyset$, then $G_V$ acts effectively
on each connected component of $\mathring{\partial}^k V$.
\end{itemize}
\end{defi}
Lemma \ref{effective action} implies that the image of $G_V$ in $\mathrm{Aut}\,V$ is
$\mathrm{Aut}_{\V} V := \{ g : V \tocong V ; \pi_V g = \pi_V \}$.
We always use $G_V = \mathrm{Aut}_{\V} V$ in this paper.

For a connected open subset $\U \subset \V $, let $U\subset V$ be a connected
component of $\pi_V^{-1}(\U)$.
Then $\U = (U, \pi_V|_{U}, \U) $ is an orbichart, and this does not depend on the choice
of the connected component.
Note that $G_U \subsetneq G_V$ in general.

\begin{defi}
We say a continuous map $\varphi : \V \to \V'$ between two orbicharts is smooth if
there exists a smooth map $\phi : V \to V'$ such that $\pi_{V'} \phi = \varphi \pi_V$
on $V$.
We call $\phi$ a lift of $\varphi$.
\[
\begin{tikzcd}
&V \ar{d}{\pi_V} \ar{r}{\phi} &V' \ar{d}{\pi_{V'}} \\
&\V \ar{r}{\varphi} & \V' 
\end{tikzcd}
\]
\end{defi}
\begin{defi}\label{smooth map}
An embedding $\varphi : \V \to \V'$ is an injective smooth map such that
there exists a connected neighborhood $\U'$ of $\varphi (V) \subset \V'$
which satisfies the following conditions:
The lift $\phi : V \to U'$ of $\varphi : \V \to \U'$ is an embedding,
$\phi (V) \subset U' $ is $G_{U'}$-invariant and
$G_{U'}$ acts on $\phi (V)$ effectively.
If in addition $\phi(V) \subset U'$ is open, then we say $\varphi$ is an open embedding.
(This is the case where $\dim V = \dim V'$.)
\end{defi}
\begin{rem}
In the above definition, we cannot always take $\U' = \V'$ since
$\phi (V) \subset V'$ is not always $G_{V'}$-invariant.
We also note that Lemma \ref{effective action} implies $G_V \cong G_{U'}$.
\end{rem}

\begin{defi}
An orbibundle chart $(\V, \E) = ( (V,\pi_V, \V) , (E, \pi_E , \E) , \tilde{\pi} , \pi)$
consists of 
\begin{itemize}
\item
topological spaces $\V , \E$
\item
a vector bundle $\tilde{\pi} : E \to V$ over a connected manifold $V$
\item
continuous maps $\pi_V : V \to \V, \pi_E : E \to \E, \pi : \E \to \V $
\end{itemize}
which satisfy the following conditions:
\begin{itemize}
\item
There exists a finite group $G$ acting smoothly and effectively on $V$ and $E$
such that
\begin{itemize}
\item
$\tilde{\pi} : E \to V$ is $G$-equivariant
\item
$\pi_V$ and $\pi_E$ are $G$-equivariant, and they induce homeomorphisms
$\pi_V : V/G \to \V$ and $\pi_E : E/G \to \E$
\item
if $\mathring{\partial}^k V \neq \emptyset$, then $G$ acts effectively on
every connected component of $\mathring{\partial}^k V$
\end{itemize}
\item
The following diagram is commutative.
\[
\begin{tikzcd}
&E \ar{d}{\tilde{\pi}} \ar{r}{\pi_E} &\E \ar{d}{\pi} \\
&V \ar{r}{\pi_{V}} & \V 
\end{tikzcd}
\]
\end{itemize}
\end{defi}
Note that we can take $G = \mathrm{Aut}_{\V} V $.
Note also that $\V$ and $ \E$ are orbichart themselves.

For a connected open subspace $\U \subset \V$,
\[
(\U, \E | _{\U}) = ((U, \pi_V|_U, \U),
(E|_U, \pi_E |_{E|_U}, \E | _{\U}), \tilde{\pi} |_{E|_U}, \pi |_{\E |_{\U}})
\]
is also an orbibundle chat.

\begin{defi}
A bundle map $(\varphi , \hat\varphi) : (\V, \E) \to (\V', \E')$ between
two orbibundle charts is a pair of continuous maps
$\varphi : \V \to \V'$ and $\hat\varphi : \E \to \E'$ such that
there exist some smooth bundle map
$(\phi,\hat\phi) : (V,E) \to (V',E')$ which makes the following diagram commutative.
\[
\begin{tikzcd}
&E \ar{r}{\hat\phi} \ar{d} \ar{rrd} &E' \ar{d} \ar{rrd} & & \\
&V \ar{r}{\phi} \ar{rrd} &V'\ar{rrd} &\E \ar{r}[swap]{\hat\varphi} \ar{d} & \E' \ar{d} \\
& &&\V \ar{r}[swap]{\varphi} &\V' 
\end{tikzcd}
\]
\end{defi}
\begin{defi}
We say a bundle map $(\varphi , \hat\varphi) : (\V, \E) \to (\V', \E')$ is
an embedding if $\varphi : \V \to \V' $ is an embedding and
the restriction of $\hat \phi : E \to E' $ to each fiber is injective.
In this case, $\hat\varphi : \E \to \E' $ is also an embedding between two orbicharts.
We say $(\varphi, \hat\varphi)$ is an open embedding if in addition $\varphi$ is
an open embedding and $\hat\phi$ is an isomorphism on each fiber.
\end{defi}


\begin{defi}\label{orbibundle}
An orbibundle $(\V, \E) = (\V, \E, \pi)$ consists of
Hausdorff spaces $\V$, $\E$ and a continuous map $\pi : \E \to \V$
which satisfies the following conditions.
\begin{itemize}
\item
For each $x \in \V$, there exists a neighborhood $\V_x \subset \V$ such that
$(\V_x,\ab \E|_{\pi^{-1}(\V_x)},\ab \pi|_{\E|_{\pi^{-1}(\V_x)}})$ has a structure of
orbibundle chart.
We define $\E_x = \E|_{\pi^{-1}(\V_x)}$.
We always assume $\pi_{V_x}^{-1}(x) \subset V_x$ is one point and
$G_{V_x} = \Aut_{\V_x}V_x$ fixes this point.
We denote this point $\pi_{V_x}^{-1}(x) \subset V_x$ by $x\in V_x$.
We always assume that $x \in \mathring{\partial}^k V_x$ for the largest $k \geq 0$
such that $\mathring{\partial}^k V_x \neq \emptyset$.
\item
For each $y \in \V_x$, if we shrink the neighborhood $\V_y$, the inclusion map
$(\V_y, \E_y) \inj (\V_x, \E_x)$ is an open embedding of orbibundle chart.
\end{itemize}
\end{defi}

\begin{defi}
An embedding $(\varphi, \hat \varphi) : (\V, \E, \pi) \to (\V', \E', \pi')$ of an orbibundle is
a pair of continuous maps $\varphi : \V \to \V'$ and $\hat \varphi : \E \to \E'$
such that
\begin{itemize}
\item
$\pi' \circ \hat \varphi = \varphi \circ \pi : \E \to \V'$
\item
for each $x \in \V$, $(\varphi_x, \hat \varphi_x)
:= (\varphi, \hat \varphi)|_{(\V_x, \E_x)} : (\V_x, \E_x)
\to (\V_{\varphi(x)}, \E_{\varphi(x)})$ is an embedding of an orbibundle chart
if we shrink $\V_x$.
\end{itemize}
\end{defi}

\begin{defi}
Let $(\V, \E) = (\V, \E, \pi)$ be an orbibundle.
A smooth section $s : \V \to \E$ is a continuous map such that
$\pi \circ s = \id_{\V}$ and the restriction of $s$ on each $\V_x$ is
a smooth map between orbicharts $\V_x$ and $\E_x$.
Note that the lift of $s$ on $V_x$ is unique and it is a $G_{V_x}$-equivalent section of
$(V_x, E_x)$.
We also denote this $G_{V_x}$-equivalent section by $s : V_x \to E_x$.
\end{defi}

\begin{defi}
We say a bundle map $(\varphi, \hat \varphi) : (\V, \E) \to (\V', \E')$ between two
orbibundle charts is a submersion if for its lift $(\phi, \hat \phi)$,
$\phi : V \to V'$ is a submersion and
the restriction of $\hat \phi$ to each fiber is an isomorphism.
Note that Lemma \ref{effective action} implies that
there exists a homomorphism $\rho_\phi : G_V \to G_{V'}$ such that
$\phi \circ g = \rho_\phi(g) \circ \phi$.
\end{defi}


Finally we consider fiber product.
Let $(\V, \E)$ be an orbibundle chart and $\varphi : \V \to Y$ be a submersion
to a manifold $Y$.
Then for any submanifold $Z \subset Y$,
$(\varphi^{-1}(Z), \E|_{\varphi^{-1}(Z)})$ is an orbibundle chart
(or a disjoint union of orbibundle charts if $\varphi^{-1}(Z)$ is disconnected).
We note that this satisfies the assumption of the effective group action.
Indeed, the $G_V$-action on $\pi_V^{-1}(\varphi^{-1}(Z))
= \phi^{-1}(Z)$ is effective, where $\phi : V \to Y$ is a lift of $\varphi$.

For the construction of SFT (in particular for Bott-Morse case),
we also need to treat fiber products over orbifolds.
\begin{defi}
\label{def of fiber product of orbibundle over orbifold}
Let $\W = (W, \pi_W, \W)$ be an orbibundle chart and
$K \subset \W$ be an embedded simplicial complex.
We assume that there exists a regular $G_W$-complex
$L \subset W$ (see \cite{Bre72} for regular complex) and an isomorphism
$\psi : L/G_W \cong K$ such that $\psi \circ \pi_L = \pi_W$
on $L \subset W$, where $\pi_L : L \to L/G_W$ is the quotient map.
Let $\varphi$ be a submersion from an orbichart $\V$ to $\W$ and assume that
for any point $p \in V$, the stabilizer $G_p \subset G_V$ of $p$ acts on
a neighborhood of $p$ in $\pi_V^{-1}(\varphi^{-1}(\varphi(\pi_V(p))))$ effectively.
Then $\varphi^{-1}(K) = (\phi^{-1}(L), \pi_V|_{\phi^{-1}(L)}, \varphi^{-1}(K))$
(or its connected components) are not orbicharts in a strict sense because
$\phi^{-1}(L)$ is not a manifold,
but for each connected component $\phi^{-1}(L)_0$ of $\phi^{-1}(L)$,
the group $\{g \in G_V; g\phi^{-1}(L)_0 = \phi^{-1}(L)_0\}$ acts effectively on it.
We regard each connected component of $(\varphi^{-1}(K), \E|_{\varphi^{-1}(K)})$
as an orbibundle chart.
Using this kind of orbibundle charts,
we define the fiber product of an orbibundle with a simplex in an orbibundle.
We say a section (or a multisection) of $(\varphi^{-1}(K), \E|_{\varphi^{-1}(K)})$
is smooth if its lift (or its branches) are the restrictions of some smooth sections
defined on a neighborhood of $\phi^{-1}(L) \subset V$.
\end{defi}

In application, we sometimes need to use the following notion of essential submersion.
In the following definition, we regard usual orbibundle charts $(\V, \E)$ as
a fiber product of $(\V, \E)$ with the $0$-dimensional simplex in a point.
\begin{defi}
\label{def of essential submersion}
Let $(\V_i, \E_i)$ be orbibundle charts for $i = 1,2$, and
let $(\varphi_i^{-1}(K_i), \ab \E_i|_{\varphi_i^{-1}(K_i)})$ be their fiber products
with embedded simplicial complexes $K_i \subset \W_i$
as in Definition \ref{def of fiber product of orbibundle over orbifold}.
As in Definition \ref{def of fiber product of orbibundle over orbifold},
let $\phi_i$ be the lifts of $\varphi_i$, and
let $L_i \subset W_i$ be the regular $G_W$-complex
such that $L_i / G_W \cong K_i$.
An essential submersion $(\varphi, \hat \varphi)$ from
$(\varphi_1^{-1}(K_1), \E_1|_{\varphi_1^{-1}(K_1)})$ to
$(\varphi_2^{-1}(K_2), \E_2|_{\varphi_2^{-1}(K_2)})$ is
a smooth bundle map from $(\V_1, \E_2)$ to $(\V_2, \E_2)$ whose lifts
$(\phi, \hat \phi)$ satisfy the following conditions.
\begin{itemize}
\item
The image $\varphi(\varphi_1^{-1}(K_1))$
is contained in $\varphi_2^{-1}(K_2)$.
\item
For any $k_1, k_2 \geq 0$ and any simplices $s_1$ of $L_1$ and $s_2$ of $L_2$,
\[
A = \phi^{-1}\bigl(\mathring{\partial}^{k_2} V_2 \cap \phi_2^{-1}(\Int s_2)\bigr)
\cap (\mathring{\partial}^{k_1} V_1 \cap \phi_1^{-1}(\Int s_1))
\subset V_1
\]
is a submanifold of $V_1$,
$\phi|_A : A \to \mathring{\partial}^{k_2} V_2 \cap \phi_2^{-1}(\Int s_2)$ is
a submersion.
\item
The restriction of $\hat \phi$ to each fiber is an isomorphism.
\end{itemize}

For fiber products $(\varphi_i^{-1}(K_i), \E_i|_{\varphi_i^{-1}(K_i)})$
of orbibundles $(\V_i, \E_i)$ with simplicial complexes $K_i \subset \W_i$,
an essential submersion $(\varphi, \hat \varphi)$ from
$(\varphi_1^{-1}(K_1), \ab \E_1|_{\varphi_1^{-1}(K_1)})$ to
$(\varphi_2^{-1}(K_2), \E_2|_{\varphi_2^{-1}(K_2)})$ is
a smooth bundle map from a neighborhood $(\mathcal{N}_1, \E_1|_{\mathcal{N}_1})$
of $\varphi_1^{-1}(K_1)$ to $(\V_2, \E_2)$
whose restriction to each orbibundle chart is an essential submersion
in the above sense.
\end{defi}

\begin{eg}
The map $f : [0, \infty)^k \times \R \to [0, \infty) \times \R$ defined by
\[
f((s_i)_{1 \leq i \leq k}, t) = (s_1 \cdots s_k, t)
\]
is an essential submersion if we regard $[0, \infty)^k \times \R$ and
$[0, \infty) \times \R$ as orbibundle charts of rank $0$ without group action.
The map $h : [0, \infty) \times \R \to [0, \infty)^2$ defined by
\[
h(\hat \rho, b) = (\hat \rho, \hat \rho e^b)
\]
is also an essential submersion.
\end{eg}

\begin{rem}
The map $h$ in the above example maps the point $(0, 0)$
of the corner of codimension $1$ to the point $(0, 0)$ of the corner of
codimension $2$.
Note that there does not exist a submersion which maps a point of the corner
of codimension $k$ to a point of the corner of codimension $l > k$.
Therefore the notion of essential submersion is crucial to understand
the structure of moduli space of disconnected holomorphic buildings.
See Section \ref{Kuranishi of disconnected buildings}.
\end{rem}

We also use the following generalization.
In Definition \ref{def of fiber product of orbibundle over orbifold},
we assumed that $L \subset W$ is a regular $G_W$-complex.
Instead, let $L \subset W$ be a $G_W$-invariant Euclidean cell complex.
(We do not assume that $L / G_W \subset W$ is a Euclidean cell complex.)
Then we can similarly define the fiber product
$(\varphi^{-1}(L / G_W), \ab \E|_{\varphi^{-1}(L / G_W)})$.
In this case, we read simplices in Definition \ref{def of essential submersion}
as cells of the Euclidean cell complexes.

For example, for an embedded simplicial complex $K \subset \W$,
the fiber product with $\prod^N K / \mathfrak{S}_N \subset \prod^N \W / \mathfrak{S}$
is defined as the above generalization.
We do not use any subdivision of $\prod^N K$ in this case.

\subsection{Multisections}
In this paper, we use a different definition of multisection.
Perturbed multisection in Definition \ref{sum of section and grouped multisection}
plays the role of multisection in \cite{FO99}.
\begin{defi}
\label{multisection for chart}
A multisection $s = (s^\omega)_{\omega \in \Omega}$
of an orbibundle chart $(\V , \E)$ is a family of
smooth sections $s^\omega : V \to E$ ($\omega \in \Omega$) indexed by
a finite $G_V$-set $\Omega$ such that
$s^{g \omega} = g_{\ast} s^\omega$ for any $\omega \in \Omega$
and $g \in G_V$.
\end{defi}

\begin{defi}\label{grouped multisection for chart}
A grouped multisection $\boldsymbol{\epsilon}
= (\epsilon^\omega)_{\omega \in \coprod_j \Omega_j}$ of an orbibundle chart
$(\V, \E)$ is a multisection of $(\V, \E)$ whose index set $\Omega$ has
a decomposition $\Omega = \coprod_j \Omega_j$ preserved by the action of
$G_V$, that is, for any $g \in G_V$ and $j$,
$g \Omega_j$ coincides with some $\Omega_{j'}$.
We define a family of sections $\epsilon_j = (\epsilon^\omega)_{\omega \in \Omega_j}$
for each $j$,
and we also denote the grouped multisection by
$\boldsymbol{\epsilon} = \{\epsilon_j\}_j$ as a set of such families.
We define the support of each $\epsilon_j$ by
$\supp(\epsilon_j) = \bigcup_{\omega \in \Omega_j} \supp(\epsilon^\omega)
\subset V$.
For a grouped multisection, we also impose the condition
$\supp(\epsilon_j) \neq \emptyset$ for all $j$.
(This is for consistency with the definition of restriction below.)
\end{defi}

\begin{defi}
For a connected open subset $\U \subset \V$, the restriction of a grouped multisection
$\boldsymbol{\epsilon} = (\epsilon^\omega)_{\omega \in \coprod_j \Omega_j}$
of $(\V, \E)$ to $(\U, \E|_{\U})$ is defined by
\[
\boldsymbol{\epsilon}|_{\U}
= (\epsilon^\omega|_U)_{\omega \in \coprod_{j \in I_U} \Omega_j}
\]
where $I_U = \{j ; \supp(\epsilon_j) \cap U \neq \emptyset \}$.
\end{defi}

\begin{eg}
Let $(\V, \E)$ be an orbibundle chart and
let $\epsilon : V \to E$ be a smooth section.
Then its average $\Av \epsilon = (g_{\ast} \epsilon)_{g \in G_V}$ is a multisection.
\end{eg}

\begin{eg}
For finite number of grouped multisections $\boldsymbol{\epsilon}^k$ of $(\V, \E)$,
their union $\coprod_k \boldsymbol{\epsilon}^k$ is also a grouped multisection.
In particular, for finite number of non-zero multisections
$\epsilon_j = (\epsilon^\omega)_{\omega \in \Omega_j}$,
$\boldsymbol{\epsilon} = \{\epsilon_j\}$ is a grouped multisection.
(We cannot always assume that each $\epsilon_j$ is a multisection,
that is, each $\Omega_j$ is not $G_V$-invariant in general.
We need the general case for the induced multisection of the quotient of the product
of the same pre-Kuranishi spaces. See Section \ref{compatible perturbed multisection})
\end{eg}

\begin{defi}
\label{sum of section and grouped multisection}
For a smooth section $s$ and a grouped multisection $\boldsymbol{\epsilon}
= (\epsilon^\omega)_{\omega \in \coprod_j \Omega_j}$ of an orbibundle chart $(\V, \E)$,
their sum is defined by the multisection $s + \boldsymbol{\epsilon}
= (s + \sum_j \epsilon^{\omega_j})_{(\omega_j) \in \prod_j \Omega_j}$
with the product index set $\prod_j \Omega_j$.
We call a multisection of this form a perturbed multisection.
\end{defi}

We will construct a perturbed multisection of a pre-Kuranishi space by the sum
$s + \boldsymbol{\epsilon}$ of the given smooth section $s$ and
a grouped multisection $\boldsymbol{\epsilon}$.
Hence it is enough to define compatibility condition of grouped multisection
$\boldsymbol{\epsilon}$ with embedding instead of the multisection
$s + \boldsymbol{\epsilon}$.

\begin{defi}
\label{(varphi, hat varphi)-relation for chart}
Let $(\varphi, \hat \varphi) : (\mathring{\V}, \mathring{\E}) \to (\V, \E)$
be an embedding between two orbibundle charts.
We say a grouped multisection $\boldsymbol{\mathring{\epsilon}}
= (\mathring{\epsilon}^\omega)_{\omega \in \coprod_j \mathring{\Omega}_j}$ of
$(\mathring{\V}, \mathring{\E})$ and
$\boldsymbol{\epsilon}
= (\epsilon^\omega)_{\omega \in \coprod_j \Omega_j}$ of $(\V, \E)$
are $(\varphi, \hat \varphi)$-related if
there exists an injection $\nu^\phi : \coprod_j \mathring{\Omega}_j
\to \coprod_j \Omega_j$ for each lift $(\phi, \hat \phi)$ of $(\varphi, \hat \varphi)$
and they satisfy the following conditions:
\begin{itemize}
\item
$\nu^\phi$ maps each $\mathring{\Omega}_j$ to some $\Omega_{j'}$
bijectively.
\item
$\epsilon^{\nu^\phi(\omega)} \circ \phi = \hat \phi \circ \mathring{\epsilon}^\omega$
for each $\omega \in \coprod_j \mathring{\Omega}_j$.
\item
$\epsilon^{\nu^\phi(\omega)} = 0$ on a neighborhood of
$\phi(\mathring{V})$ for any $\omega \in \coprod_j \Omega_j \setminus
\nu^\phi(\coprod_j\mathring{\Omega}_j)$.
\item
For any connected open subset $\mathring{U} \subset \mathring{V}$ and $j$,
if $\mathring{\epsilon}^\omega|_{\mathring{U}} = 0$ for all
$\omega \in \mathring{\Omega}_j$,
then $\epsilon^{\nu^\phi(\omega)} = 0$ on a neighborhood of
$\phi(\mathring{U})$ for all $\omega \in \mathring{\Omega}_j$.
\item
$\nu^{g \phi \mathring{g}} = g \circ \nu^\phi \circ \mathring{g}$
for any $g \in G_V$ and $\mathring{g} \in G_{\mathring{V}}$.
\end{itemize}
\end{defi}

\begin{defi}
\label{grouped multisection for orbibundle}
Let $(\V, \E)$ be an orbibundle.
A grouped multisection $\boldsymbol{\epsilon} = (\B, \boldsymbol{\epsilon}_\U,
\nu_{\U_2, \U_1}^\phi)$ of $(\V, \E)$
consists of the following.
$\B = \{\U\}$ is a set of connected open subsets of $\V$ such that
each $(\U, \E|_{\U})$ is an orbibundle chart and
if $\U \in \B$ then every connected open subset of $\U$ is contained in $\B$.
For each $\U \in \B$,
$\boldsymbol{\epsilon}_\U = (\epsilon^\omega_\U)_{\omega \in \coprod_j \Omega_j^U}$
is a grouped multisection of $(\U, \E|_{\U})$.
For each pair $\U_1, \U_2 \in \B$ such that
$\U_1 \subset \U_2$ and a lift $\phi : U_1 \to U_2$
of the inclusion map $\U_1 \inj \U_2$, there exists an injective map
$\nu_{\U_2, \U_1}^\phi : \coprod_j \Omega_j^{U_1} \to \coprod_j \Omega_j^{U_2}$
which satisfy the following conditions:
\begin{itemize}
\item
$\nu_{\U_2, \U_1}^\phi$ maps each $\Omega_j^{U_1}$ to some $\Omega_{j'}^{U_2}$
bijectively.
\item
$\epsilon^{\nu_{\U_2, \U_1}^\phi(\omega)}_{\U_2} \circ \phi
= \hat \phi \circ \epsilon_{\U_1}^\omega$
for any $\omega \in \coprod_j \Omega_j^{U_1}$, where
$\hat \phi$ is the lift of $\hat \varphi$ uniquely determined by $\phi$.)
\item
$\epsilon^{\omega'}_{\U_2} \circ \phi =0$
for any $\omega' \in \coprod_j \Omega_j^{U_2} \setminus
\nu_{\U_2, \U_1}^\phi(\coprod_j \Omega_j^{U_1})$.
\item
$\nu_{\U_2, \U_1}^{g_2 \circ \phi \circ g_1} = g_2 \circ \nu_{\U_2, \U_1}^\phi \circ g_1$
for any $g _1\in G_{U_1}$ and $g_2 \in G_{U_2}$.
\item
$\nu_{\U_3, \U_2}^{\phi_{3,2}} \circ \nu_{\U_2, \U_1}^{\phi_{2,1}}
= \nu_{\U_3, \U_1}^{\phi_{3,2} \circ \phi_{2,1}}$ for any triple $\U_1, \U_2, \U_3 \in \B$
such that $\U_1 \subset \U_2 \subset \U_3$
and lifts $\phi_{2,1} : U_1 \to U_2$ and $\phi_{3,2} : U_2 \to U_3$.
\end{itemize}
\end{defi}

\begin{rem}
We do not define a multisecton of an orbibundle.
(Definition \ref{multisection for chart} is the definition of a multisection of
an orbibundle chart, and Definition \ref{grouped multisection for orbibundle} is
the definition of a grouped multisection of an orbibundle.)
We construct a grouped multisection of an orbibundle,
and for each orbibundle chart, we use the perturbed multisection
$s + \boldsymbol{\epsilon}$ of Definition \ref{sum of section and grouped multisection}.
\end{rem}

\begin{eg}
In general, a grouped multisection $\boldsymbol{\epsilon}
= (\B, \boldsymbol{\epsilon}_\U, \nu_{\U_2, \U_1}^\phi)$ of an orbibundle $(\V, \E)$
does not have a global grouped multisection $\boldsymbol{\epsilon}_\V$.
Namely, even if $(\V, \E)$ itself is an orbibundle chart,
there does not exist a grouped multisection $\boldsymbol{\epsilon}_\V$
of an orbibundle chart $(\V, \E)$
(in the sense of Definition \ref{grouped multisection for chart})
whose restrictions to $\U$ coincide with $\boldsymbol{\epsilon}_\U$
for all $\U \in \B$.
For example, let $f : \R \to \R$ be a periodic smooth function
of period $4\pi$, and let $\chi : \R_{\geq 0} \to \R$ be a smooth function
whose support is contained in $[1/2, 1] \subset \R_{\geq 0}$.
Then $F(r \cos \theta, r \sin \theta) = \chi(r) f(\theta)$ defines
a grouped multisection of the trivial orbibbundle of rank $1$ on $\R^2$
(without group action).
(We define the decomposition of the index sets so that
the indices for the two branches of $F$ constitute one group.)
However, it cannot be represented by a grouped multisection of
the trivial orbibundle chart on $\R^2$.
We also note that on a neighborhood of $(0,0) \in \R^2$,
it is represented by the grouped multisection whose index set
is the empty set.
\end{eg}

\begin{defi}
\label{(varphi, hat varphi)-relation for orbibundles}
For an embedding $(\varphi, \hat \varphi) : (\mathring{\V}, \mathring{\E}) \to
(\V, \E)$ between two orbibundles, we say a grouped multisection
$\boldsymbol{\mathring{\epsilon}} = (\mathring{\B},
\boldsymbol{\mathring{\epsilon}}_{\mathring{\U}},
\mathring{\nu}_{\mathring{\U}_2, \mathring{\U}_1}^\phi)$ of
$(\mathring{\V}, \mathring{\E})$ and
$\boldsymbol{\epsilon} = (\B, \boldsymbol{\epsilon}_\U, \nu_{\U_2, \U_1}^\phi)$ of
$(\V, \E)$ are $(\varphi, \hat \varphi)$-related if
the following conditions hold.
For any $\mathring{\U} \in \mathring{\B}$ and $\U \in \B$ such that
$(\varphi, \hat \varphi)$ defines an embedding
$(\mathring{\U}, \mathring{\E}|_{\mathring{U}}) \to (\U, \E|_{\U})$, and its lift
$(\phi, \hat \phi)$,
there exists an injective map $\nu_{\U, \mathring{\U}}^\phi
: \coprod_j \mathring{\Omega}_j^{\mathring{U}} \to \coprod_j \Omega_j^{U}$
which satisfies the following conditions:
\begin{itemize}
\item
$\nu_{\U, \mathring{\U}}^\phi$ maps each $\mathring{\Omega}_j^{\mathring{U}}$
to some $\Omega_{j'}^{U}$ bijectively.
\item
$\epsilon^{\nu_{\U, \mathring{\U}}^\phi(\omega)}_{\U} \circ \phi
= \hat \phi \circ \mathring{\epsilon}_{\mathring{\U}}^\omega$
for any $\omega \in \coprod_j \mathring{\Omega}_j^{\mathring{U}}$.
\item
$\epsilon^{\omega'}_{\U} = 0$ on a neighborhood of $\phi(\mathring{U})$
for any $\omega' \in \coprod_j \Omega_j^{U} \setminus
\nu_{\U, \mathring{\U}}^\phi(\coprod_j \mathring{\Omega}_j^{\mathring{U}})$.
\item
$\nu_{\U, \mathring{\U}}^{g \circ \phi \circ \mathring{g}}
= g \circ \nu_{\U, \mathring{\U}}^\phi \circ \mathring{g}$
for any $g \in G_U$ and $\mathring{g} \in G_{\mathring{U}}$.
\item
$\nu_{\U_3, \U_2}^{\phi_{3,2}} \circ \nu_{\U_2, \mathring{\U}_2}^\phi \circ
\nu_{\mathring{\U}_2, \mathring{\U}_1}^{\mathring{\phi}_{2,1}}
= \nu_{\U_3, \mathring{\U}_1}^{\phi_{3,2} \circ \phi \circ \mathring{\phi}_{2,1}}$
for any $\mathring{\U}_1 \subset \mathring{\U}_2 \in \mathring{\B}$,
$\U_2 \subset \U_3 \in \B$
such that $(\varphi, \hat \varphi)$ defines an embedding
$(\mathring{\U}_2, \mathring{\E}|_{\mathring{U}_2}) \to (\U_2, \E|_{\U_2})$
and lifts $\mathring{\phi}_{2,1} : \mathring{U}_1 \to \mathring{U}_2$,
$\phi : \mathring{U}_2 \to U_2$ and $\phi_{3,2} : U_2 \to U_3$.
\end{itemize}
More precisely, in the above case,
we say $\boldsymbol{\mathring{\epsilon}}$ and $\boldsymbol{\epsilon}$ are
$(\varphi, \hat \varphi)$-related by
$(\nu_{\U, \mathring{\U}}^\phi)_{(\U, \mathring{\U})}$.
\end{defi}

\begin{defi}
\label{compatibility of (varphi, hat varphi)-relations}
Let $(\varphi^{j,i}, \hat \varphi^{j,i}) : (\V^i, \E^i) \to (\V^j, \E^j)$
be embeddings of orbibundles for $1 \leq i < j \leq 3$ such that
$(\varphi^{3,1}, \hat \varphi^{3,1})
= (\varphi^{3,2}, \hat \varphi^{3,2}) \circ (\varphi^{2,1}, \hat \varphi^{2,1})$.
Let $\boldsymbol{\epsilon}^i = (\B^i, \boldsymbol{\epsilon}^i_{\U^i},
\nu_{\U^i_2, \U^i_1}^\phi)$ be a grouped multisection of $(\V^i, \E^i)$
for each $i = 1,2,3$, and assume that
$\boldsymbol{\epsilon}^i$ and $\boldsymbol{\epsilon}^j$ are
$(\varphi^{j,i}, \hat \varphi^{j,i})$-related by $(\nu_{\U^j, \U^i}^\phi)_{(\U^j, \U^i)}$
for $1 \leq i < j \leq 3$.
We say these relations are compatible if
$\nu_{\U^3, \U^2}^{\phi^{3,2}} \circ \nu_{\U^2, \U^1}^{\phi^{2,1}}
= \nu_{\U^3,\U^1}^{\phi^{3,2} \circ \phi^{2,1}}$ for any
$\U^i \in \B^i$ such that
$(\varphi^{j,i}, \hat \varphi^{j,i})$ defines an embedding
$(\U^i, \E^i|_{\U^i}) \to (\U^j, \E^j|_{\U^j})$
for all $1 \leq i < j \leq 3$, and
lifts $(\phi^{j,i}, \hat \phi^{j,i}) : (U^i, E^i|_{U^i}) \to (U^j, E^j|_{U^j})$
for $(i,j) = (1,2), (2,3)$.
\end{defi}
We always assume the above compatibility condition for compositions of
embeddings.

Next we consider the extension of grouped multisection for
an embedding of an orbibundle.
\begin{lem}
\label{extension of grouped multisection for embedding}
Let $(\varphi, \hat \varphi) : (\mathring{\V}, \mathring{\E}) \to (\V, \E)$
be an embedding between two orbibundles.
For any grouped multisection $\boldsymbol{\mathring{\epsilon}} = (\mathring{\B},
\boldsymbol{\mathring{\epsilon}}_{\mathring{\U}},
\mathring{\nu}_{\mathring{\U}_2, \mathring{\U}_1}^\phi)$ of
$(\mathring{\V}, \mathring{\E})$ and its arbitrary relatively
compact open subset $\mathring{\V}' \Subset \mathring{\V}$,
we can construct a grouped multisection
$\boldsymbol{\epsilon}$ of $(\V, \E)$ which is $(\varphi, \hat \varphi)$-related to
$\boldsymbol{\mathring{\epsilon}}|_{\mathring{\V}'}$.
\end{lem}
\begin{proof}
Let $(\mathring{\V}_\alpha, \mathring{\E}_\alpha)_{\alpha \in \A}$ and
$(\V_\alpha, \E_\alpha)_{\alpha \in \A}$ be finite number of orbibundle charts of
$(\mathring{\V}, \mathring{\E})$ and $(\V, \E)$ respectively such that
$\mathring{\V}_\alpha \in \mathring{\B}$,
$\{\mathring{\V}_\alpha\}_{\alpha \in \A}$ covers the closure of $\mathring{\V}'$,
and $(\varphi, \hat \varphi)$ defines an embedding of
$(\mathring{\V}_\alpha, \mathring{\E}_\alpha)$ to $(\V_\alpha, \E_\alpha)$.
We fix a lift $(\phi_\alpha, \hat \phi_\alpha)$ of this embedding for each $\alpha \in \A$.
Replacing $\V_\alpha$ with a smaller connected open neighborhood of
$\varphi(\mathring{\V}_\alpha)$ if necessary,
we may assume that this lift defines an isomorphism of the automorphism group of
$(\mathring{\V}_\alpha, \mathring{\E}_\alpha)$ and that of $(\V_\alpha, \E_\alpha)$.

Take compact subsets $\mathring{K}_\alpha \subset \mathring{\V}_\alpha$
such that $\bigcup_{\alpha \in \A} \Int \mathring{K}_\alpha \supset \mathring{\V}'$.
We can take finite orbibundle charts
$(\mathring{\V}_\kappa, \mathring{\E}_\kappa)_{\kappa \in \K}$ of
$(\mathring{\V}, \mathring{\E})$ and subsets $A_\kappa \subset \A$ ($\kappa \in \K$)
such that
$\bigcup_{\kappa \in \K} \mathring{\V}_\kappa \Supset \mathring{\V}'$ and
$\mathring{\V}_\kappa \Subset \bigcap_{\alpha \in A_\kappa} \mathring{\V}_\alpha
\setminus \bigcup_{\beta \in \A \setminus A_\kappa} \mathring{K}_\beta$.
For each $\kappa \in \K$, let $(\V_\kappa, \E_\kappa)$ be an orbibundle chart of
$(\V, \E)$ such that
$\V_\kappa \Subset \bigcap_{\alpha \in A_\kappa} \V_\alpha
\setminus \bigcup_{\beta \in \A \setminus A_\kappa} \varphi(\mathring{K}_\beta)$ and
$(\varphi, \hat \varphi)$ defines an embedding of
$(\mathring{\V}_\kappa, \mathring{\E}_\kappa)$ to $(\V_\kappa, \E_\kappa)$.
We fix a lift of this embedding $(\phi_\kappa, \hat \phi_\kappa)$
for each $\kappa \in \K$ and assume that this
lift defines an isomorphism between their automorphism groups.

For each pair $\kappa_1, \kappa_2 \in \K$ such that
$\mathring{\V}_{\kappa_1} \cap \mathring{\V}_{\kappa_2} \neq \emptyset$,
let $\{\mathring{\V}_{\kappa_1, \kappa_2, \gamma}\}_{\gamma}$
be the connected components of
the intersection $\mathring{\V}_{\kappa_1} \cap \mathring{\V}_{\kappa_2}$.
Similarly, for each triple $\kappa_1, \kappa_2, \kappa_3 \in \K$ such that
$\mathring{\V}_{\kappa_1} \cap \mathring{\V}_{\kappa_2} \cap
\mathring{\V}_{\kappa_3} \neq \emptyset$,
let $\{\mathring{\V}_{\kappa_1, \kappa_2, \kappa_3, \gamma}\}_{\gamma}$
be the connected components of
the intersection $\mathring{\V}_{\kappa_1} \cap \mathring{\V}_{\kappa_2}
\cap \mathring{\V}_{\kappa_3}$.
For each $\mathring{\V}_{\kappa_1, \kappa_2, \gamma}$,
let $(\V_{\kappa_1, \kappa_2, \gamma}, \E_{\kappa_1, \kappa_2, \gamma})$ be
an orbibundle chart contained in the intersection
$\V_{\kappa_1} \cap \V_{\kappa_2}$ such that
$(\varphi, \hat \varphi)$ defines an embedding of
$(\mathring{\V}_{\kappa_1, \kappa_2, \gamma},
\mathring{\E}_{\kappa_1, \kappa_2, \gamma})$
to $(\V_{\kappa_1, \kappa_2, \gamma}, \E_{\kappa_1, \kappa_2, \gamma})$
and its lift $(\phi_{\kappa_1, \kappa_2, \gamma}, \hat \phi_{\kappa_1, \kappa_2, \gamma})$
defines an isomorphism between their automorphism groups.
Similarly, for each $\mathring{\V}_{\kappa_1, \kappa_2, \kappa_3, \gamma}$, we define
$(\V_{\kappa_1, \kappa_2, \kappa_3, \gamma},
\E_{\kappa_1, \kappa_2, \kappa_3, \gamma})$
and $(\phi_{\kappa_1, \kappa_2, \kappa_3, \gamma},
\hat \phi_{\kappa_1, \kappa_2, \kappa_3, \gamma})$.
We assume that $\V_{\kappa_1, \kappa_2, \kappa_3, \gamma}$ is contained in
$\V_{\kappa_1, \kappa_2, \gamma_{1,2}} \cap \V_{\kappa_2, \kappa_3, \gamma_{2,3}}
\cap \V_{\kappa_1, \kappa_3, \gamma_{1,3}}$
if $\mathring{\V}_{\kappa_1, \kappa_2, \kappa_3, \gamma}$ is contained in
$\mathring{\V}_{\kappa_1, \kappa_2, \gamma_{1,2}} \cap
\mathring{\V}_{\kappa_2, \kappa_3, \gamma_{2,3}}
\cap \mathring{\V}_{\kappa_1, \kappa_3, \gamma_{1,3}}$.

For each $\kappa \in \K$ and $\alpha \in \A$,
we fix a lift $(\phi_{\alpha, \kappa}^\circ, \hat \phi_{\alpha, \kappa}^\circ)$
of the inclusion map from $(\mathring{\V}_\kappa, \mathring{\E}_\kappa)$ to
$(\mathring{\V}_\alpha, \mathring{\E}_\alpha)$.
Then we can define a lift $(\phi_{\alpha, \kappa}, \hat \phi_{\alpha, \kappa})$
of the inclusion map from $(\V_\kappa, \E_\kappa)$ to
$(\V_\alpha, \E_\alpha)$ by the condition
$(\phi_{\alpha, \kappa}, \hat \phi_{\alpha, \kappa}) \circ
(\phi_\kappa, \hat \phi_\kappa) = (\phi_\alpha, \hat \phi_\alpha)
\circ (\phi_{\alpha, \kappa}^\circ, \hat \phi_{\alpha, \kappa}^\circ)$.
Similarly, we fix a lift
$(\phi_{\alpha, (\kappa_1, \kappa_2, \gamma)}^\circ,
\hat \phi_{\alpha, (\kappa_1, \kappa_2, \gamma)}^\circ)$
of embedding from $(\mathring{\V}_{\kappa_1, \kappa_2, \gamma},
\mathring{\E}_{\kappa_1, \kappa_2, \gamma})$ to
$(\mathring{\V}_\alpha, \mathring{\E}_\alpha)$,
define the lift $(\phi_{\alpha, (\kappa_1, \kappa_2, \gamma)},
\hat \phi_{\alpha, (\kappa_1, \kappa_2, \gamma)})$
of embedding from
$(\V_{\kappa_1, \kappa_2, \gamma}, \E_{\kappa_1, \kappa_2, \gamma})$ to
$(\V_\alpha, \E_\alpha)$,
and so on.

For each $\alpha \in \A$,
we independently construct a grouped multisection
$\boldsymbol{\epsilon}_{\V_\alpha} = (\epsilon_{\V_\alpha}^\omega)_{\omega \in \coprod_j
\mathring{\Omega}_j^{\mathring{V}_\alpha}}$ of $(\V_\alpha, \E_\alpha)$
which is $(\varphi, \hat \varphi)$-related to
$\boldsymbol{\mathring{\epsilon}}_{\mathring{\V}_\alpha}$.
We use the same index set for $\boldsymbol{\epsilon}_{\V_\alpha}$
as that of $\boldsymbol{\mathring{\epsilon}}_{\mathring{\V}_\alpha}$,
and assume that $\nu^{\phi_\alpha} = \id$
in Definition \ref{(varphi, hat varphi)-relation for chart}.
Shrinking $\V_\kappa$, $\V_{\kappa_1, \kappa_2, \gamma}$ and
$\V_{\kappa_1, \kappa_2, \kappa_3, \gamma}$ to smaller neighborhoods of
$\varphi(\mathring{\V}_\kappa)$, $\varphi(\mathring{\V}_{\kappa_1, \kappa_2, \gamma})$
and $\varphi(\mathring{\V}_{\kappa_1, \kappa_2, \kappa_3, \gamma})$ respectively
if necessary, we may assume the following conditions on $\epsilon_{\V_\alpha}^\omega$.
For each $j$, if $\mathring{\epsilon}_{\mathring{\V}_\alpha}^\omega \circ
\phi_{\alpha, \kappa}^\circ = 0$ for all
$\omega \in \mathring{\Omega}_j^{\mathring{V}_\alpha}$, then
$\epsilon_{\V_\alpha}^\omega \circ \phi_{\alpha, \kappa} = 0$
for all $\omega \in \mathring{\Omega}_j^{\mathring{V}_\alpha}$.
Similarly, if $\mathring{\epsilon}_{\mathring{\V}_\alpha}^\omega \circ
\phi_{\alpha, (\kappa_1, \kappa_2, \gamma)}^\circ = 0$
for all $\omega \in \mathring{\Omega}_j^{\mathring{V}_\alpha}$ then
$\epsilon_{\V_\alpha}^\omega \circ \phi_{\alpha, (\kappa_1, \kappa_2, \gamma)} = 0$
for all $\omega \in \mathring{\Omega}_j^{\mathring{V}_\alpha}$,
and if $\mathring{\epsilon}_{\mathring{\V}_\alpha}^\omega \circ
\phi_{\alpha, (\kappa_1, \kappa_2, \kappa_3 \gamma)}^\circ = 0$
for all $\omega \in \mathring{\Omega}_j^{\mathring{V}_\alpha}$ then
$\epsilon_{\V_\alpha}^\omega \circ
\phi_{\alpha, (\kappa_1, \kappa_2, \kappa_3, \gamma)} = 0$
for all $\omega \in \mathring{\Omega}_j^{\mathring{V}_\alpha}$.
We note that these conditions do not depend on the choice of the lifts
$\phi_{\alpha, \kappa}^\circ$, $\phi_{\alpha, (\kappa_1, \kappa_2, \gamma)}^\circ$
or $\phi_{\alpha, (\kappa_1, \kappa_2, \kappa_3, \gamma)}^\circ$.

Since $\{\V_\kappa\}_{\kappa \in \K}$ covers the closure of
$\varphi(\mathring{\V}')$,
we can construct open subsets $\V'_\kappa \Subset \V_\kappa$
such that $\V'_{\kappa_1} \cap \V'_{\kappa_2}$ is contained in
the union of $\V_{\kappa_1, \kappa_2, \gamma}$,
$\V'_{\kappa_1} \cap \V'_{\kappa_2} \cap \V'_{\kappa_3}$ is contained in
the union of $\V_{\kappa_1, \kappa_2, \kappa_3, \gamma}$, and
$\{\V'_\kappa\}_{\kappa \in \K}$ covers the closure of $\varphi(\mathring{\V}')$.

Let $\{\chi_\alpha\}_\alpha$ be a family of smooth functions on $\V$ such that
$\varphi^{-1}(\supp \chi_\alpha) \subset \mathring{K}_\alpha$
and $\sum_\alpha \chi_\alpha \equiv 1$ on $\varphi(\mathring{\V}')$.
We assume that $\V_\kappa \cap \supp \chi_\beta = \emptyset$
for all $\kappa \in \K$ and $\beta \in \A \setminus A_\kappa$.
We also assume that
$\bigcup_{\alpha \in \A} \supp \chi_\alpha \subset \bigcup_{\kappa \in \K} \V'_\kappa$.

Then we define grouped multisections $\boldsymbol{\epsilon}_{\U}$
for all connected open subsets contained in some $\V'_\kappa$
and connected open subsets which do not intersect with
$\bigcup_{\alpha \in \A} \supp \chi_\alpha$.
For the latter, we define $\boldsymbol{\epsilon}_{\U}$
by zero (the grouped multisection whose index set is the empty set).
For the former, we define $\boldsymbol{\epsilon}_{\U}$ as follows.

First we define a grouped multisection $\boldsymbol{\epsilon}_{\V_\kappa}
= (\epsilon_{\V_\kappa}^\omega)_{\omega \in \coprod_j
\mathring{\Omega}^{\mathring{V}_\kappa}_j}$
of $(\V_\kappa, \E_\kappa)$ for each $\kappa \in \K$ by
\[
\epsilon_{\V_\kappa}^\omega = \sum_{\alpha \in A_\kappa}
\chi_\alpha \, \phi_{\alpha, \kappa}^\ast
\epsilon_{\V_\alpha}^{\nu_{\alpha, \kappa}(\omega)},
\]
where $\nu_{\alpha, \kappa} = \nu_{\mathring{\V}_\alpha, \mathring{\V}_\kappa}
^{\phi_{\alpha, \kappa}^\circ}$.

For each connected open subset $\U$ contained in some $\V'_\kappa$,
we fix one of such $\kappa \in \K$, and define its grouped multisection
$\boldsymbol{\epsilon}_\U$
by the restriction of the multisection $\boldsymbol{\epsilon}_{\V_\kappa}$.
Namely, we fix a lift $(\phi_{\kappa, U}, \hat \phi_{\kappa, U})$ of the inclusion map and
define $\boldsymbol{\epsilon}_\U
= (\epsilon_\U^\omega)_{\omega \in \coprod_{j \in I_{\kappa, U}}
\mathring{\Omega}^{\mathring{V}_\kappa}_j}$ by
$\epsilon_\U^\omega = \phi_{\kappa, U}^\ast \epsilon_{\V_\kappa}^\omega$,
where $I_{\kappa, U} = \{j; \phi_{\kappa, U}^\ast \epsilon_{\V_\kappa}^\omega \neq 0$
for some $\omega \in \mathring{\Omega}^{\mathring{\V}_\kappa}_j\}$.
We need to construct $\nu_{\U_2, \U_1}^{\phi_{U_2, U_1}}$ for pairs $\U_1 \subset \U_2$
and lifts $(\phi_{U_2, U_1}, \hat \phi_{U_2, U_1})$ of the inclusion $\U_1 \inj \U_2$.
Assume that the grouped multisections of $\U_1$ and $\U_2$ are
defined by using $\kappa_1$ and $\kappa_2$ respectively.
In particular, $\U_1 \subset \V'_{\kappa_1} \cap \V'_{\kappa_2}$ is contained in
$\V_{\kappa_1, \kappa_2, \gamma}$ for some $\gamma$.
Fix a lift $(\phi_{(\kappa_1, \kappa_2, \gamma), U_1},
\hat \phi_{(\kappa_1, \kappa_2, \gamma), U_1})$ of the inclusion
from $(\U_1, \E|_{\U_1})$ to $(\V_{\kappa_1, \kappa_2, \gamma},
\E_{\kappa_1, \kappa_2, \gamma})$, and define
$g^{\kappa_1}_{(\kappa_1, \kappa_2, \gamma), U_1} \in G_{V_{\kappa_1}}$ by
\[
\phi_{\kappa_1, U_1}
= g^{\kappa_1}_{(\kappa_1, \kappa_2, \gamma), U_1} \circ
\phi_{\kappa_1, (\kappa_1, \kappa_2, \gamma)} \circ
\phi_{(\kappa_1, \kappa_2, \gamma), U_1}.
\]
First we show that $\coprod_{j \in I_{\kappa_1, U_1}}
\mathring{\Omega}^{\mathring{V}_{\kappa_1}}_j$ is contained in the image of
\begin{equation}
\nu_{\mathring{\V}_{\kappa_1}, \mathring{\V}_{\kappa_1, \kappa_2, \gamma}}
^{g^{\kappa_1}_{(\kappa_1, \kappa_2, \gamma), U_1} \circ
\phi_{\kappa_1, (\kappa_1, \kappa_2, \gamma)}^\circ}.
\label{nu g phi}
\end{equation}
For each $j \in I_{\kappa_1, U_1}$,
there exists some $\omega \in \mathring{\Omega}_j^{\mathring{V}_{\kappa_1}}$
such that $\phi_{\kappa_1, U_1}^\ast \epsilon_{\V_{\kappa_1}}^\omega \neq 0$.
Hence
\begin{align*}
&(g^{\kappa_1}_{(\kappa_1, \kappa_2, \gamma), U_1} \circ
\phi_{\kappa_1, (\kappa_1, \kappa_2, \gamma)})^\ast \epsilon_{\V_{\kappa_1}}^\omega \\
&= \sum_{\alpha \in A_{\kappa_1}}
\chi_\alpha (\phi_{\alpha, \kappa_1} \circ
g^{\kappa_1}_{(\kappa_1, \kappa_2, \gamma), U_1} \circ
\phi_{\kappa_1, (\kappa_1, \kappa_2, \gamma)})^\ast
\epsilon_{\V_\alpha}^{\nu_{\alpha, \kappa_1}(\omega)}
\end{align*}
is nonzero.
This implies that
some $(\phi_{\alpha, \kappa_1} \circ
g^{\kappa_1}_{(\kappa_1, \kappa_2, \gamma), U_1} \circ
\phi_{\kappa_1, (\kappa_1, \kappa_2, \gamma)})^\ast
\epsilon_{\V_\alpha}^{\nu_{\alpha, \kappa_1}(\omega)}$ is nonzero.
Since $\phi_{\alpha, \kappa_1} \circ
g^{\kappa_1}_{(\kappa_1, \kappa_2, \gamma), U_1} \circ
\phi_{\kappa_1, (\kappa_1, \kappa_2, \gamma)}$ is a lift of the
open embedding $\V_{\kappa_1, \kappa_2, \gamma} \inj \V_\alpha$,
the assumption of $\boldsymbol{\epsilon}_{\V_\alpha}$ implies that
\[
(\phi_{\alpha, \kappa_1}^\circ \circ
g^{\kappa_1}_{(\kappa_1, \kappa_2, \gamma), U_1} \circ
\phi_{\kappa_1, (\kappa_1, \kappa_2, \gamma)}^\circ)^\ast
\mathring{\epsilon}_{\mathring{\V}_\alpha}^{\nu_{\alpha, \kappa_1}(\omega')}
= (g^{\kappa_1}_{(\kappa_1, \kappa_2, \gamma), U_1} \circ
\phi_{\kappa_1, (\kappa_1, \kappa_2, \gamma)}^\circ)^\ast
\mathring{\epsilon}_{\mathring{\V}_{\kappa_1}}^{\omega'}
\]
is also nonzero for some $\omega' \in \mathring{\Omega}^{\mathring{V}_{\kappa_1}}_j$.
This implies that $\mathring{\Omega}^{\mathring{V}_{\kappa_1}}_j$
is contained in the image of (\ref{nu g phi}).
Hence $\coprod_{j \in I_{\kappa_1, U_1}} \mathring{\Omega}^{\mathring{V}_{\kappa_1}}_j$
is contained in the image of (\ref{nu g phi}).

We define $g^{\kappa_2}_{U_2, U_1} \in G_{\V_{\kappa_2}}$ by
\[
\phi_{\kappa_2, U_2} \circ \phi_{U_2,U_1}
= g_{U_2, U_1}^{\kappa_2} \circ \phi_{\kappa_2, (\kappa_1, \kappa_2, \gamma)} \circ
\phi_{(\kappa_1, \kappa_2, \gamma), U_1},
\]
and define $\nu_{\U_2, \U_1}^{\phi_{U_2, U_1}}
: \coprod_{j \in I_{\kappa_1, U_1}} \mathring{\Omega}^{\mathring{V}_{\kappa_1}}_j
\to \coprod_{j \in I_{\kappa_2, U_2}} \mathring{\Omega}^{\mathring{V}_{\kappa_2}}_j$
by
\[
\nu_{\U_2, \U_1}^{\phi_{U_2, U_1}}
= \nu_{\mathring{\V}_{\kappa_2}, \mathring{\V}_{\kappa_1, \kappa_2, \gamma}}
^{g_{U_2, U_1}^{\kappa_2}
\circ \phi_{\kappa_2, (\kappa_1, \kappa_2, \gamma)}^\circ}
\circ (\nu_{\mathring{\V}_{\kappa_1}, \mathring{\V}_{\kappa_1, \kappa_2, \gamma}}
^{g^{\kappa_1}_{(\kappa_1, \kappa_2, \gamma), U_1} \circ
\phi_{\kappa_1, (\kappa_1, \kappa_2, \gamma)}^\circ})^{-1}.
\]
We need to check that this satisfies the conditions of $\nu_{\U_2, \U_1}^{\phi_{U_2, U_1}}$.

First we check the condition
\begin{equation}
\epsilon^{\nu_{\U_2, \U_1}^{\phi_{U_2, U_1}}(\omega)}_{\U_2} \circ \phi_{U_2, U_1}
= \hat \phi_{U_2, U_1} \circ \epsilon_{\U_1}^\omega
\label{(phi, hat phi)-relation for extension}
\end{equation}
for $\omega \in \coprod_{j \in I_{\kappa_1, U_1}}
\mathring{\Omega}^{\mathring{V}_{\kappa_1}}_j$.
This equation also implies that the image of $\nu_{\U_2, \U_1}^{\phi_{U_2, U_1}}$
is indeed contained in $\coprod_{j \in I_{\kappa_2, U_2}}
\mathring{\Omega}^{\mathring{V}_{\kappa_2}}_j$.
Define $\hat \omega \in \coprod_j
\mathring{\Omega}^{\mathring{V}_{\kappa_1, \kappa_2, \gamma}}_j$
by
\[
\omega = \nu_{\mathring{\V}_{\kappa_1}, \mathring{\V}_{\kappa_1, \kappa_2, \gamma}}
^{g^{\kappa_1}_{(\kappa_1, \kappa_2, \gamma), U_1} \circ
\phi_{\kappa_1, (\kappa_1, \kappa_2, \gamma)}^\circ}(\hat \omega).
\]
Then by definition,
\[
\nu_{\U_2, \U_1}^{\phi_{U_2, U_1}}(\omega)
= \nu_{\mathring{\V}_{\kappa_2}, \mathring{\V}_{\kappa_1, \kappa_2, \gamma}}
^{g_{U_2, U_1}^{\kappa_2}
\circ \phi_{\kappa_2, (\kappa_1, \kappa_2, \gamma)}^\circ}(\hat \omega).
\]
Therefore
\begin{align}
&\phi_{U_2, U_1}^\ast \epsilon^{\nu_{\U_2, \U_1}^{\phi_{U_2, U_1}}(\omega)}_{\U_2} \notag \\
&= (\phi_{\kappa_2, U_2} \circ \phi_{U_2, U_1})^\ast
\epsilon_{\V_{\kappa_2}}^{\nu_{\U_2, \U_1}^{\phi_{U_2, U_1}}(\omega)} \notag \\
&= (g_{U_2, U_1}^{\kappa_2} \circ \phi_{\kappa_2, (\kappa_1, \kappa_2, \gamma)} \circ
\phi_{(\kappa_1, \kappa_2, \gamma), U_1})^\ast
\epsilon_{\V_{\kappa_2}}^{\nu_{\mathring{\V}_{\kappa_2},
\mathring{\V}_{\kappa_1, \kappa_2, \gamma}}^{g_{U_2, U_1}^{\kappa_2}
\circ \phi_{\kappa_2, (\kappa_1, \kappa_2, \gamma)}^\circ}(\hat \omega)} \notag \\
&= (\phi_{\kappa_2, (\kappa_1, \kappa_2, \gamma)} \circ
\phi_{(\kappa_1, \kappa_2, \gamma), U_1})^\ast
\epsilon_{\V_{\kappa_2}}^{\nu_{\mathring{\V}_{\kappa_2},
\mathring{\V}_{\kappa_1, \kappa_2, \gamma}}^{
\phi_{\kappa_2, (\kappa_1, \kappa_2, \gamma)}^\circ}(\hat \omega)} \notag \\
&= \sum_{\alpha \in A_{\kappa_2}}
\chi_\alpha (\phi_{\alpha, \kappa_2} \circ
\phi_{\kappa_2, (\kappa_1, \kappa_2, \gamma)} \circ
\phi_{(\kappa_1, \kappa_2, \gamma), U_1})^\ast
\epsilon_{\V_\alpha}^{\nu_{\mathring{\V}_\alpha,
\mathring{\V}_{\kappa_1, \kappa_2, \gamma}}^{\phi_{\alpha, \kappa_2}^\circ
\circ \phi_{\kappa_2, (\kappa_1, \kappa_2, \gamma)}^\circ}(\hat \omega)}.
\label{phi ast epsilon U2}
\end{align}
Since $\chi_\alpha|_{\U_1} = 0$
for $\alpha \in A_{\kappa_2} \setminus A_{\kappa_1}$,
the terms for $\alpha \in A_{\kappa_2} \setminus A_{\kappa_1}$ are zero.
For each $\alpha \in A_{\kappa_1} \cap A_{\kappa_2}$,
we define $g^\alpha \in G_{V_\alpha}$ by
\[
\phi_{\alpha, \kappa_2}^\circ
\circ \phi_{\kappa_2, (\kappa_1, \kappa_2, \gamma)}^\circ
= g^\alpha \circ \phi_{\alpha, \kappa_1}^\circ \circ
g^{\kappa_1}_{(\kappa_1, \kappa_2, \gamma), U_1} \circ
\phi_{\kappa_1, (\kappa_1, \kappa_2, \gamma)}^\circ.
\]
Then
\begin{align*}
&\phi_{\alpha, \kappa_2} \circ
\phi_{\kappa_2, (\kappa_1, \kappa_2, \gamma)} \circ
\phi_{(\kappa_1, \kappa_2, \gamma), U_1} \\
&= g^\alpha \circ \phi_{\alpha, \kappa_1} \circ
g^{\kappa_1}_{(\kappa_1, \kappa_2, \gamma), U_1} \circ
\phi_{\kappa_1, (\kappa_1, \kappa_2, \gamma)} \circ
\phi_{(\kappa_1, \kappa_2, \gamma), U_1} \\
&= g^\alpha \circ \phi_{\alpha, \kappa_1} \circ \phi_{\kappa_1, U_1}.
\end{align*}
Hence (\ref{phi ast epsilon U2}) is equal to
\begin{align*}
&\sum_{\alpha \in A_{\kappa_1}}
\chi_\alpha (g^\alpha \circ \phi_{\alpha, \kappa_1} \circ \phi_{\kappa_1, U_1})^\ast
\epsilon_{\V_\alpha}^{
\nu_{\mathring{\V}_\alpha, \mathring{\V}_{\kappa_1, \kappa_2, \gamma}}
^{g^\alpha \circ \phi_{\alpha, \kappa_1}^\circ \circ
g^{\kappa_1}_{(\kappa_1, \kappa_2, \gamma), U_1} \circ
\phi_{\kappa_1, (\kappa_1, \kappa_2, \gamma)}^\circ}(\hat \omega)}\\
&= \sum_{\alpha \in A_{\kappa_1}}
\chi_\alpha (\phi_{\alpha, \kappa_1} \circ \phi_{\kappa_1, U_1})^\ast
\epsilon_{\V_\alpha}^{
\nu_{\mathring{\V}_\alpha, \mathring{\V}_{\kappa_1, \kappa_2, \gamma}}
^{\phi_{\alpha, \kappa_1}^\circ \circ
g^{\kappa_1}_{(\kappa_1, \kappa_2, \gamma), U_1} \circ
\phi_{\kappa_1, (\kappa_1, \kappa_2, \gamma)}^\circ}(\hat \omega)} \\
&= \epsilon_{\U_1}^{
\nu_{\mathring{\V}_{\kappa_1}, \mathring{\V}_{\kappa_1, \kappa_2, \gamma}}
^{g^{\kappa_1}_{(\kappa_1, \kappa_2, \gamma), U_1} \circ
\phi_{\kappa_1, (\kappa_1, \kappa_2, \gamma)}^\circ}(\hat \omega)}\\
&= \epsilon_{\U_1}^\omega.
\end{align*}
Hence (\ref{(phi, hat phi)-relation for extension}) holds for
$\omega \in \coprod_{j \in I_{\kappa_1, U_1}}
\mathring{\Omega}^{\mathring{V}_{\kappa_1}}_j$.

Next we check the condition
$\epsilon^{\omega}_{\U_2} \circ \phi_{U_2, U_1} = 0$
for $\omega \in \coprod_{j \in I_{\kappa_2, U_2}}
\mathring{\Omega}^{\mathring{V}_{\kappa_2}}_j$ not contained in the
image of $\nu_{\U_2, \U_1}^{\phi_{U_2, U_1}}$.
If $\omega = \nu_{\mathring{\V}_{\kappa_2}, \mathring{\V}_{\kappa_1, \kappa_2, \gamma}}
^{g_{U_2, U_1}^{\kappa_2}
\circ \phi_{\kappa_2, (\kappa_1, \kappa_2, \gamma)}^\circ}(\hat \omega)$
for some $\hat \omega \in \coprod_j
\mathring{\Omega}^{\mathring{V}_{\kappa_1, \kappa_2, \gamma}}_j$,
then the same argument as above implies that
\[
\phi_{U_2, U_1}^\ast \epsilon^{\omega}_{\U_2}
= \phi_{\kappa_1, U_1}^\ast
\epsilon_{\V_{\kappa_1}}^{
\nu_{\mathring{\V}_{\kappa_1}, \mathring{\V}_{\kappa_1, \kappa_2, \gamma}}
^{g^{\kappa_1}_{(\kappa_1, \kappa_2, \gamma), U_1} \circ
\phi_{\kappa_1, (\kappa_1, \kappa_2, \gamma)}^\circ}(\hat \omega)},
\]
and this is zero because $\nu_{\mathring{\V}_{\kappa_1},
\mathring{\V}_{\kappa_1, \kappa_2, \gamma}}
^{g^{\kappa_1}_{(\kappa_1, \kappa_2, \gamma), U_1} \circ
\phi_{\kappa_1, (\kappa_1, \kappa_2, \gamma)}^\circ}(\hat \omega)
\notin \coprod_{j \in I_{\kappa_1, U_1}}
\mathring{\Omega}^{\mathring{V}_{\kappa_1}}_j$.
If $\omega$ is not contained in the image of
$\nu_{\mathring{\V}_{\kappa_2}, \mathring{\V}_{\kappa_1, \kappa_2, \gamma}}
^{g_{U_2, U_1}^{\kappa_2}
\circ \phi_{\kappa_2, (\kappa_1, \kappa_2, \gamma)}^\circ}$,
then
\[
(\phi_{\alpha, \kappa_2}^\circ \circ
g_{U_2, U_1}^{\kappa_2} \circ \phi_{\kappa_2, (\kappa_1, \kappa_2, \gamma)}^\circ)^\ast
\mathring{\epsilon}_{\mathring{\V}_\alpha}
^{\nu_{\mathring{\V}_\alpha, \mathring{\V}_{\kappa_2}}
^{\phi_{\alpha, \kappa_2}^\circ}(\omega)}
= (g_{U_2, U_1}^{\kappa_2} \circ \phi_{\kappa_2, (\kappa_1, \kappa_2, \gamma)}^\circ)^\ast
\mathring{\epsilon}_{\mathring{\V}_\alpha}^{\omega}
= 0,
\]
and this (and the same equations for the other indices $\omega$ in the same
index group) imply
\[
(\phi_{\alpha, \kappa_2} \circ
g_{U_2, U_1}^{\kappa_2} \circ \phi_{\kappa_2, (\kappa_1, \kappa_2, \gamma)})^\ast
\epsilon_{\V_\alpha}^{\nu_{\mathring{\V}_\alpha, \mathring{\V}_{\kappa_2}}
^{\phi_{\alpha, \kappa_2}^\circ}(\omega)}
= 0
\]
by the assumption of $\boldsymbol{\epsilon}_{\V_\alpha}$.
Hence
\begin{align*}
\phi_{U_2, U_1}^\ast \epsilon_{\U_2}^\omega
&= (g_{U_2, U_1}^{\kappa_2} \circ \phi_{\kappa_2, (\kappa_1, \kappa_2, \gamma)} \circ
\phi_{(\kappa_1, \kappa_2, \gamma), U_1})^\ast
\epsilon_{\V_{\kappa_2}}^\omega \\
&= \sum_{\alpha \in A_{\kappa_2}} \chi_\alpha (\phi_{\alpha, \kappa_2} \circ
g_{U_2, U_1}^{\kappa_2} \circ \phi_{\kappa_2, (\kappa_1, \kappa_2, \gamma)} \circ
\phi_{(\kappa_1, \kappa_2, \gamma), U_1})^\ast
\epsilon_{\V_\alpha}^{\nu_{\mathring{\V}_\alpha, \mathring{\V}_{\kappa_2}}
^{\phi_{\alpha, \kappa_2}^\circ}(\omega)} \\
&= 0.
\end{align*}

Finally we check the condition about composition.
For a triple $\U_1 \subset \U_2 \subset \U_3$ and lifts
$(\phi_{U_2, U_1}, \hat \phi_{U_2, U_1})$, $(\phi_{U_3, U_2}, \hat \phi_{U_3, U_2})$
of the inclusion maps, we prove that
\[
\nu_{\U_2, \U_1}^{\phi_{U_2, U_1}} \circ \nu_{\U_2, \U_1}^{\phi_{U_2, U_1}}
= \nu_{\U_3, \U_1}^{\phi_{U_3, U_2} \circ \phi_{U_2, U_1}}.
\]
Since $\U_1$ is contained in
$\V'_{\kappa_1} \cap \V'_{\kappa_2} \cap \V'_{\kappa_3}$,
$\U_1 \subset \V_{\kappa_1, \kappa_2, \kappa_3, \gamma'}$ for some $\gamma'$.
Fix a lift $(\phi_{(\kappa_1, \kappa_2, \kappa_3, \gamma'), U_1}, \ab
\hat \phi_{(\kappa_1, \kappa_2, \kappa_3, \gamma'), U_1})$ of the inclusion
from $(\U_1, \E|_{\U_1})$ to $(\V_{\kappa_1, \kappa_2, \kappa_3, \gamma'}, \ab
\E_{\kappa_1, \kappa_2, \kappa_3, \gamma'})$, and define
$g^{\kappa_1}_{(\kappa_1, \kappa_2, \kappa_3, \gamma'), U_1} \in G_{\V_{\kappa_1}}$ by
\[
\phi_{\kappa_1, U_1}
= g^{\kappa_1}_{(\kappa_1, \kappa_2, \kappa_3, \gamma'), U_1} \circ
\phi_{\kappa_1, (\kappa_1, \kappa_2, \kappa_3, \gamma')} \circ
\phi_{(\kappa_1, \kappa_2, \kappa_3, \gamma'), U_1}.
\]
By the same argument as above, $\coprod_{j \in I_{\kappa_1, U_1}}
\mathring{\Omega}^{\mathring{V}_{\kappa_1}}_j$ is contained in the image of
\[
\nu_{\mathring{\V}_{\kappa_1}, \mathring{\V}_{\kappa_1, \kappa_2, \kappa_3, \gamma'}}
^{g^{\kappa_1}_{(\kappa_1, \kappa_2, \kappa_3, \gamma'), U_1} \circ
\phi_{\kappa_1, (\kappa_1, \kappa_2, \kappa_3, \gamma')}^\circ}.
\]
Hence it is enough to prove
\begin{align*}
&\nu_{\U_2, \U_1}^{\phi_{U_2, U_1}} \circ \nu_{\U_2, \U_1}^{\phi_{U_2, U_1}} \circ
\nu_{\mathring{\V}_{\kappa_1}, \mathring{\V}_{\kappa_1, \kappa_2, \kappa_3, \gamma'}}
^{g^{\kappa_1}_{(\kappa_1, \kappa_2, \kappa_3, \gamma'), U_1} \circ
\phi_{\kappa_1, (\kappa_1, \kappa_2, \kappa_3, \gamma')}^\circ} \\
&= \nu_{\U_3, \U_1}^{\phi_{U_3, U_2} \circ \phi_{U_2, U_1}} \circ
\nu_{\mathring{\V}_{\kappa_1}, \mathring{\V}_{\kappa_1, \kappa_2, \kappa_3, \gamma'}}
^{g^{\kappa_1}_{(\kappa_1, \kappa_2, \kappa_3, \gamma'), U_1} \circ
\phi_{\kappa_1, (\kappa_1, \kappa_2, \kappa_3, \gamma')}^\circ},
\end{align*}
which follows from the conditions of $\nu$'s of
the grouped multisection $\boldsymbol{\mathring{\epsilon}}$.

By construction, this grouped multisection $\boldsymbol{\epsilon}$
is $(\varphi, \hat \varphi)$-related to
$\boldsymbol{\mathring{\epsilon}}|_{\mathring{\V}'}$.
\end{proof}

We note that we can apply the same argument for more general cases.
Let $(\varphi_i, \hat \varphi_i) : (\mathring{\V}_i, \mathring{\E}_i) \to (\V, \E)$
be embeddings of orbibundles and $\boldsymbol{\mathring{\epsilon}}_i$ be
grouped multisections of $(\mathring{\V}_i, \mathring{\E}_i)$.
Assume that these embeddings intersect cleanly and
the grouped multisections are compatible on the intersections.
Then for any relatively compact subsets $\mathring{\V}'_i \Subset \mathring{\V}_i$,
we can also construct a grouped multisection of $(\V, \E)$
which is $(\varphi_i, \hat \varphi_i)$-related to
$\boldsymbol{\mathring{\epsilon}}_i|_{\mathring{\V}'_i}$ for all $i$.

For a submersion, we can define the pull back of a grouped multisection.
\begin{defi}
Let $\boldsymbol{\epsilon} = (\epsilon^\omega)_{\omega \in \coprod_j \Omega_j}$
be a grouped multisection of an orbibundle chart $(\V, \E)$.
Let $(\varphi, \hat \varphi)$ be a submersion from another orbibundle chart
$(\V', \E')$ to $(\V, \E)$, and $(\phi, \hat \phi)$ be its lift.
Then we can define the pull back $\varphi^\ast \boldsymbol{\epsilon}$ by
$\varphi^\ast \boldsymbol{\epsilon}
= (\phi^\ast \epsilon^\omega)_{\omega \in \coprod_{j \in I} \Omega_j}$,
where $I = \{j; \supp \epsilon_j \cap \phi(V') \neq \emptyset\}$.
We define the $G_{V'}$-action on $\coprod_{j \in I} \Omega_j$ by
the homomorphism $\rho_\phi : G_{V'} \to G_V$ associated to $\phi$.
The pull back of a grouped multisection of an oribibundle by a submersion is
defined by the pull backs for its orbibundle charts.
\end{defi}
We can also define the pull back of a grouped multisection by
an essential submersion.

\subsection{Pre-Kuranishi structure and construction of its perturbed multisection}
We introduce the notion of pre-Kuranishi structure.
This is essentially equivalent to the usual Kuranishi structure in the sense that
we can obtain a Kuranishi structure from a pre-Kuranishi structure and
in application, when we construct a Kuranishi structure,
we usually construct a pre-Kuranishi structure implicitly.
(See Remark \ref{pre-Kuranishi and Kuranishi}.)
However, for a pre-Kuranishi space, we can define weakly good coordinate system,
which is more compatible with product than good coordinate system.
\begin{defi}
Let $X$ be a compact Hausdorff space.
A pre-Kuranishi structure on $X$ consists of the following data
$(\widetilde{X},\ab \mu,\ab (\W_x, \E_x, s_x, \widetilde{\psi}_x),\ab
(\varphi_{x, y}, \hat \varphi_{x, y}))$:
\begin{itemize}
\item
$\widetilde{X}$ is a Hausdorff space,
and $\mu : \widetilde{X} \to X$ is a locally-homeomorphic surjection such that
$\# \mu^{-1}(p)$ $(p \in X)$ is bounded.
\item
Each $(\W_x, \E_x, s_x, \widetilde{\psi}_x)$ is a Kuranishi neighborhood of
$x \in \widetilde{X}$.
Namely, $(\W_x, \E_x)$ is an orbibundle, $s_x : \W_x \to \E_x$ is a smooth section,
and $\widetilde{\psi}_x : s_x^{-1}(0) \inj \widetilde{X}$ is a homeomorphism
onto a neighborhood of $x \in \widetilde{X}$.
We assume that $\psi_x = \mu \circ \widetilde{\psi}_x : s_x^{-1}(0) \inj X$ is also
a homeomorphism onto a neighborhood of $\mu(x)$.
Hence $(\W_x, \E_x, s_x, \psi_x)$ is a Kuranishi neighborhood of $\mu(x)$.
\item
For each $p \in X$, $\mu^{-1}(p)$ has a partial order such that any two elements
$x, y \in \mu^{-1}(p)$ have a unique supremum $x \vee y \in \mu^{-1}(p)$.
Furthermore we assume that $\vee$ is continuous in the following sense:
If $x' \in \widetilde{\psi}_x(s_x^{-1}(0))$, $y' \in \widetilde{\psi}_y(s_y^{-1}(0))$ and
$z' \in \widetilde{\psi}_{x \vee y}(s_{x \vee y}^{-1}(0))$ satisfy
$\mu(x') = \mu(y') = \mu(z')$, then $z' = x' \vee y'$.
Note that this implies the continuity of the partial order, that is,
if $x \geq y$ then $x' \geq y'$ for any $x' \in \widetilde{\psi}_x(s_x^{-1}(0))$ and
$y' \in \widetilde{\psi}_y(s_y^{-1}(0))$ such that $\mu(x') = \mu(y')$.
\item
For each $p \in \psi_x(s_x^{-1}(0))$, define $p_x$ by the unique point in
$\widetilde{\psi}_x(s_x^{-1}(0))$ such that $\mu(p_x) = p$.
We sometimes denote the point $\widetilde{\psi}_x^{-1}(p_x) \in \W_x$
by the same symbol $p_x$.
\item
For any points $x, y \in \widetilde{X}$, if there exists a point
$p \in \psi_x(s_x^{-1}(0)) \cap \psi_y(s_y^{-1}(0))$ such that $p_x \geq p_y$,
then there exists an open neighborhood $\W_{x, y} \subset \W_y$ of
$\psi_y^{-1}(\psi_x(s_x^{-1}(0)))$ and an embedding $(\varphi_{x, y}, \hat\varphi_{x, y}) :
(\W_{x, y}, \E_y|_{\W_{x, y}}) \to (\W_x, \E_x)$ which satisfy the following conditions:
\begin{itemize}
\item
The following diagrams are commutative.
\[
\begin{tikzcd}
\E_y|_{\W_{x, y}} \ar{r}{\hat\varphi_{x, y}} & \E_x\\
\W_{x, y} \ar{u}{s_y} \ar{r}{\varphi_{x, y}} & \W_x \ar{u}{s_x}
\end{tikzcd}
\quad
\begin{tikzcd}
s_y^{-1}(0) \cap \W_{x, y} \ar{r}{\psi_y} \ar[hook]{d}{\varphi_{x, y}}& X\\
s_x^{-1}(0) \ar{ru}[swap]{\psi_x}&
\end{tikzcd}
\]
\item
The vertical differential
\[
d^\bot s_x : \frac{T_{p_x}W_x}{(\phi_{x, y})_\ast T_{p_y}W_y} \tocong
\frac{(E_y)_{p_y}}{\hat\phi_{x, y} (E_x)_{p_x}}
\]
is an isomorphism for each point $p \in \psi_x(s_x^{-1}(0)) \cap \psi_y(s_y^{-1}(0))$,
where $(\phi_{x, y}, \hat \phi_{x, y})$ is a lift of $(\varphi_{x, y}, \hat \varphi_{x, y})$.
(More precisely, $(\phi_{x, y}, \hat \phi_{x, y})$ is a lift of the restriction of
$(\varphi_{x, y}, \hat \varphi_{x, y})$ to orbibundle charts of $(\W_y, \E_y)$
and $(\W_x, \E_x)$ which contains $p_y$ and $p_x$ respectively.)
\item
For $x, y, z \in \widetilde{X}$, if there exists a point
$p \in \psi_x(s_x^{-1}(0)) \cap \psi_y(s_y^{-1}(0)) \cap \psi_z(s_z^{-1}(0))$
such that $p_x \geq p_y \geq p_z$, then
\[
(\varphi_{x, y}, \hat\varphi_{x, y}) \circ (\varphi_{y, z}, \hat\varphi_{y, z})
= (\varphi_{x, z}, \hat\varphi_{x, z})
\]
on some neighborhood $\W_{x, y, z} \subset \W_{x, z} \cap \varphi_{y, z}^{-1}(\W_{x, y})$
of $\psi_z^{-1}(\psi_x(s_x^{-1}(0)) \cap \psi_y(s_y^{-1}(0)))$.
\item
If $p_x > p_y$, then the embedding is not invertible, that is, $\dim \W_x > \dim \W_y$.
\item
(separating condition)\\
For any points $a \in s_y^{-1}(0)$ and $b \in s_x^{-1}(0)$, if $\psi_y(a) \neq \psi_x(b)$,
then there exist some neighborhood $\U_a \subset \W_y$ of $a$ and
$\U_b \subset \W_x$ of $b$ such that $\U_a \cap \varphi_{x, y}^{-1}(\U_b) = \emptyset$.
(This condition is not essential because it always holds true if we replace
$\W_x$, $\W_y$ and $\W_{x, y}$ with their relatively compact subsets.)
\end{itemize}
\end{itemize}

Note that for two points $x, y \in \widetilde{X}$ such that
$\widetilde{\psi}_x(s_x^{-1}(0)) \cap \widetilde{\psi}_y(s_y^{-1}(0)) \neq \emptyset$,
$(\varphi_{x, y}, \hat \varphi_{x, y})$ is an open embedding
since $p_x = p_y$ for any point $p \in \psi_x(s_x^{-1}(0)) \cap \psi_y(s_y^{-1}(0))$.
The Hausdorff space $X$ endowed with a pre-Kuranishi structure is called
a pre-Kuranishi space.
We say $X$ is $n$-dimensional if $\dim \W_x - \dim \E_x = n$ for all $x \in \widetilde{X}$.
For two points $x, y \in \widetilde{X}$, we say $x \unrhd y$ if there exists some point
$p \in \psi_x(s_x^{-1}(0)) \cap \psi_y(s_y^{-1}(0))$ such that $p_x \geq p_y$.
Note that by assumption, this condition is independent of the choice of the point
$p \in \psi_x(s_x^{-1}(0)) \cap \psi_y(s_y^{-1}(0))$.
We also note that this is not a partial order.
Indeed, $x \unrhd y$ and $y \unrhd z$ do not imply $x \unrhd z$ in general.
\end{defi}

\begin{rem}
\label{we can shrink Wx}
We sometimes shrink each $\W_x$ to a smaller open neighborhood $\mathring{\W}_x$ of
$\widetilde{\psi}_x^{-1}(x) \subset \W_x$.
(For example, see Remark \ref{shrink W_x for submersion}.)
Then we replace $\W_{x, y}$ and $\W_{x, y, z}$ with
$\mathring{\W}_{x, y} = \mathring{\W}_y \cap \varphi_{x, y}^{-1}(\mathring{\W}_x)$
and $\mathring{\W}_{x, y, z} = \W_{x, y, z} \cap \mathring{\W}_{x, z} \cap
\varphi_{y, z}^{-1}(\mathring{\W}_{x, y})
= \W_{x, y, z} \cap \mathring{\W}_z \cap \varphi_{y, z}^{-1}(\mathring{\W}_y) \cap
\varphi_{x, z}^{-1}(\mathring{\W}_x)$ respectively.
However, once we construct a weakly good coordinate system, we should not
shrink $\W_x$ any more.
Similarly, we sometimes replace $\widetilde{X}$ to its open subset which satisfies
the same conditions.
\end{rem}

\begin{rem}\label{pre-Kuranishi and Kuranishi}
We can construct a Kuraishi structure from the above pre-Kuranishi structure
as follows.
For a compact subset $\widehat{X} \subset \widetilde{X}$ such that
$\mu(\widehat{X}) = X$, define a compact subset
$\widehat{X}^+ \subset \widetilde{X}$ by
\[
\widehat{X}^+ = \{ x_1 \vee x_2 \vee \dots \vee x_k; x_1, x_2, \dots, x_k \in \widehat{X},
\mu(x_1) = \mu(x_2) = \dots = \mu(x_k)\}.
\]
Then for each $p \in X$, $\mu^{-1}(p) \cap \widehat{X}^+ = \{x_i\}$ has a maximal point
$x_0$.
Take an open subset $\W_p \subset \W_{x_0}$ such that
\[
\psi_{x_0}(s_{x_0}^{-1}(0) \cap \W_p) \cap \mu\biggl(\widehat{X}^+ \setminus
\bigcup_{x_i \in \mu^{-1}(p) \cap \widehat{X}^+} \widetilde{\psi}_{x_i}(s_{x_i}^{-1}(0))\biggr)
= \emptyset.
\]
Then $(\W_p, \E_{x_0}|_{\W_p}, s_{x_0}, \psi_{x_0})$ and restrictions of
$(\varphi_{x, y}, \hat \varphi_{x, y})$ defines a Kuranihsi structure of $X$
in the sense of \cite{FO99}.
\end{rem}


For the construction of each embedding $(\varphi_{x, y}, \hat \varphi_{x, y})
: (\W_{x, y}, \E_y|_{\W_{x, y}}) \to (\W_x, \E_x)$, we usually use the following
easy argument in application.
\begin{lem}
Let $(\W_x, \E_x, s_x, \psi_x)$ and $(\W_y, \E_y, s_y, \psi_y)$ be Kuranishi
neighborhoods of $x$ and $y$ in a compact Hausdorff space $X$ respectively.
Assume that for each point $q \in \psi_y^{-1}(\psi_x(s_x^{-1}(0)))$,
there exist an open neighborhood $\W_y^q \subset \W_y$ of
$\psi_y^{-1}(q)$ and an embedding $(\varphi^q, \hat \varphi^q)
: (\W_y^q, \E_y|_{\W_y^q}) \to (\W_x, \E_x)$.
Assume that each $(\varphi^q, \hat \varphi^q)$ satisfies
$\hat \varphi^q \circ s_y = s_x \circ \varphi^q$ on $\W_y^q$,
$\psi_x \circ \varphi^q = \psi_y$ on $s_y^{-1}(0) \cap \W_y^q$ and
the condition of the vertical differential of $s_x$.
We also assume that for any two points
$q, q' \in \psi_y^{-1}(\psi_x(s_x^{-1}(0)))$,
if $s_y^{-1}(0) \cap \W_y^q \cap \W_y^{q'} \neq \emptyset$,
then $(\varphi^q, \hat \varphi^q)$ and $(\varphi^{q'}, \hat \varphi^{q'})$
coincide on some neighborhood $\W_y^{q, q'}$ of
$s_y^{-1}(0) \cap \W_y^q \cap \W_y^{q'}$.
Then we can construct an open neighborhood $\W_{x, y} \subset \W_y$ of
$\psi_y^{-1}(\psi_x(s_x^{-1}(0)))$ and
an embedding $(\varphi_{x, y}, \hat \varphi_{x, y}) : (\W_{x, y}, \E_y|_{\W_{x, y}})
\to (\W_x, \E_x)$ which coincides with $(\varphi^q, \hat \varphi^q)$
on a neighborhood of $\psi_y^{-1}(q)$ for all
$q \in \psi_y^{-1}(\psi_x(s_x^{-1}(0)))$ and
which satisfies $\hat \varphi_{x, y} \circ s_y = s_x \circ \varphi_{x, y}$
on $\W_{x, y}$, $\psi_x \circ \varphi_{x, y} = \psi_y$ on $s_y^{-1}(0) \cap \W_{x, y}$ and
the condition of the vertical differential of $s_x$.
\end{lem}
\begin{proof}
Let $\psi_y^{-1}(\psi_x(s_x^{-1}(0)) \cap \psi_y(s_y^{-1}(0))) = \bigcup_j K_j$ be
a locally finite covering by compact subsets such that
$K_j \subset \W_y^{q_j}$ for some
$q_j \in \psi_y^{-1}(\psi_x(s_x^{-1}(0)) \cap \psi_y(s_y^{-1}(0)))$ for each $j$.
Let $\mathring{\W}_y^{q_j} \subset \W_y^{q_j}$ be open neighborhoods of $K_j$
such that
$\mathring{\W}_y^{q_j} \cap \mathring{\W}_y^{q'_j} = \emptyset$ if
$K_j \cap K_{j'} = \emptyset$ and
$\mathring{\W}_y^{q_j} \cap \mathring{\W}_y^{q'_j} \subset \W_y^{q_j, q'_j}$
if $K_j \cap K_{j'} \neq \emptyset$.
Define $\W_{x, y} = \bigcup_j \mathring{\W}_y^{q_j}$ and
define $(\varphi_{x, y}, \hat \varphi_{x, y}) : (\W_{x, y}, \E_y|_{\W_{x, y}})
\to (\V_x, \E_x)$ by the union of $(\varphi^q, \hat \varphi^q)|_{\mathring{\W}_y^{q_j}}$.
Then it is well-defined and it is a required embedding.
\end{proof}

Although we may construct a good coordinate system from the Kuranishi structure
obtained in Remark \ref{pre-Kuranishi and Kuranishi} as in \cite{FO99}, in this paper,
we directly construct a good coordinate system from pre-Kuranishi structure.
\begin{defi}
A totally ordered cover of a pre-Kuranishi space $X$ is an open subset
$\Y \subset \widetilde{X}$ such that
$\mu(\Y) = X$ and each fiber $\Y \cap \mu^{-1}(p)$ ($p \in X$) is totally ordered.
\end{defi}
Note that if an open subset $\Y' \subset \Y$ satisfies $\mu(\Y') = X$, then
$\Y'$ is also a totally ordered cover.
The following is our good coordinate system.
\begin{defi}
A good coordinate system of a pre-Kuranishi space $X$ is a family of finite pairs
$(x, \V_x)_{x \in P}$ of points $x \in \widetilde{X}$ and open neighborhoods
$\V_x \subset \W_x$ of $\widetilde{\psi}_x^{-1}(x)$ which satisfies the following conditions.
For two points $x, y \in P$ such that $x \unrhd y$, we define
$\V_{x, y} = \V_y \cap \varphi_{x, y}^{-1}(\V_x)$.
Then $\V_x$ and $\V_{x, y}$ satisfy the following conditions:
\begin{enumerate}[label=$(\arabic*)^{\mathrm{G}}$]
\item
\label{good P totally ordered}
$\bigcup_{x \in P} \widetilde{\psi}_x(\V_x \cap s_x^{-1}(0)) \subset \widetilde{X}$ is
a totally ordered cover.
\item
\label{good (x, y, z)-relation}
For any $x, y, z \in P$, if there exists some point $p \in \psi_x(\V_x \cap s_x^{-1}(0))
\cap \psi_y(\V_y \cap s_y^{-1}(0)) \cap \psi_z(\V_z \cap s_z^{-1}(0))$ such that
$p_x \geq p_y \geq p_z$, then
\begin{equation}
\varphi_{x, y}(\V_{x, y}) \cap \varphi_{x, z}(\V_{x, z}) \subset \varphi_{x, z}(\W_{x, y, z})
\label{good (x, y, z)-overlap}
\end{equation}
and
\begin{equation}
\V_{y, z} \cap \varphi_{y, z}^{-1}(\V_{x, y}) \subset \W_{x, y, z}.
\label{good (x, y, z)-composition}
\end{equation}
\item
\label{good (x, y, z)-no relation}
For any $x, y, z \in P$ such that $\psi_x(\V_x \cap s_x^{-1}(0)) \cap
\psi_y(\V_y \cap s_y^{-1}(0)) \cap \psi_z(\V_z \cap s_z^{-1}(0)) = \emptyset$,
\begin{itemize}
\item
if $x \unrhd y$ and $x \unrhd z$, then
$\varphi_{x, y}(\V_{x, y}) \cap \varphi_{x, z}(\V_{x, z}) = \emptyset$,
\item
if $x \unrhd y$ and $y \unrhd z$, then
$\V_{x, y} \cap \varphi_{y, z}(\V_{y, z}) = \emptyset$, and
\item
if $x \unrhd z$ and $y \unrhd z$, then
$\V_{x, z} \cap \V_{y, z} = \emptyset$.
\end{itemize}
\end{enumerate}
\end{defi}
Condition \ref{good P totally ordered} implies that for any $x, y \in P$,
if $\psi_x(\V_x \cap s_x^{-1}(0)) \cap \psi_y(\V_y \cap s_y^{-1}(0)) \neq \emptyset$
and $\dim \V_x \geq \dim \V_y$, then $x \unrhd y$.
Hence there exists an embedding
$(\varphi_{x, y}, \hat \varphi_{x, y}) : (\V_{x, y}, \E_y|_{\V_{x, y}}) \to
(\V_x, \E_x|_{\V_x})$.
Therefore if we fix a total order $\preceq$ of $P$ such that
$\dim \V_y \leq \dim \V_x$ if $y \leq x$, then
our good coordinate system is essentially the same as that of \cite{FO99}.

We can construct a good coordinate system from a totally ordered cover
as follows.
\begin{lem}\label{good coordinate from totally ordered cover}
Assume that a totally ordered cover $\Y \subset \widetilde{X}$ is given.
Then for any compact subset $\K \subset \Y$,
there exists a good coordinate system $(x, \V_x)_{x \in P}$ such that
$\K \subset \bigcup_{x \in P} \widetilde{\psi}_x(\V_x \cap s_x^{-1}(0)) \subset \Y$.
\end{lem}
\begin{proof}
We may assume that $\mu(\K) = X$.
Choose finite points $P = \{x\} \subset \widetilde{X}$ and compact subsets
$\K_x \subset s_x^{-1}(0)$ such that
$\K \subset \bigcup_{x \in P} \widetilde{\psi}_x(\K_x) \subset \Y$.
We claim that if we choose a sufficiently small open neighborhood $\V_x \subset \W_x$ of
$\K_x$ for each $x \in P$ then $(x, \V_x)_{x \in P}$ is a good coordinate system.

First, it is clear that Condition \ref{good P totally ordered} holds if
$\V_x \subset\W_x$ ($x \in P$) are sufficiently small so that
$\widetilde{\psi}_x(\V_x \cap s_x^{-1}(0)) \subset \Y$.

For Condition \ref{good (x, y, z)-relation}, first we note that for any two points
$x, y \in P$ such that $x \unrhd y$, if we choose sufficiently small neighborhood
$\V_x$ and $\V_y$ of $\K_x$ and $\K_y$ respectively then
$\V_{x, y} = \V_y \cap \varphi_{x, y}^{-1}(\V_x)$ is contained in an arbitrary small
neighborhood of $\K_y \cap \psi_y^{-1}(\psi_x(\K_x))$.
This can be proved as follows.
Let $\V_x^k$ and $\V_y^k$ be decreasing sequences of relatively compact neighborhoods
of $\K_x$ and $\K_y$ such that $\bigcap_k \overline{\V_x^k} = \K_x$ and
$\bigcap_k \overline{\V_y^k} = \K_y$ respectively.
Then $\V_{x, y}^k = \V_y^k \cap \varphi_{x, y}^{-1}(\V_x^k)$ is a decreasing sequence of
relatively compact neighborhoods of $\K_y \cap \psi_y^{-1}(\psi_x(\K_x))$ such that
$\bigcap_k \overline{\V_{x, y}^k} = \K_y \cap \psi_y^{-1}(\psi_x(\K_x))$.
Indeed, for $a \in \bigcap_k \overline{\V_{x, y}^k}$, there exists a sequence
$a_k \in \V_y^k \cap \varphi_{x, y}^{-1}(\V_x^k)$ converging to
$a \in \bigcap_k \overline{\V_y^k} = \K_y \subset s_y^{-1}(0)$.
Taking subsequence, we may assume that $\varphi_{x, y}(a_k) \in \V_x^k$ converges to
some point $b \in \bigcap_k \overline{\V_x^k} = \K_x \subset s_x^{-1}(0)$.
Then the last condition of pre-Kuranishi space (separating condition) implies that
$\psi_y(a) = \psi_x(b)$.
Hence $a$ is contained in $\K_y \cap \psi_y^{-1}(\psi_x(\K_x))$.
Therefore  $\V_{x, y}^k$ is a decreasing sequence of
relatively compact neighborhoods of $\K_y \cap \psi_y^{-1}(\psi_x(\K_x))$ such that
$\bigcap_k \overline{\V_{x, y}^k} = \K_y \cap \psi_y^{-1}(\psi_x(\K_x))$,
which implies that we can make $\V_{x, y}^k$ be an arbitrary small neighborhood of
$\K_y \cap \psi_y^{-1}(\psi_x(\K_x))$.

Consider any triple $x, y, z \in P$ such that there exists some point
$p \in \psi_x(\K_x) \cap
\psi_y(\K_y) \cap \psi_z(\K_z)$ such that $p_x \geq p_y \geq p_z$.
The above argument implies that
if we choose small $\V_x$, $\V_y$ and $\V_z$,
then
\[
\varphi_{x, z}^{-1}(\varphi_{x, y}(\V_{x, y}) \cap \varphi_{x, z}(\V_{x, z}))
= \V_z \cap \varphi_{x, z}^{-1}(\V_x) \cap \varphi_{x, z}^{-1}(\varphi_{x, y}(\V_y))
\]
and
\[
\V_{y, z} \cap \varphi_{y, z}^{-1}(\V_{x, y})
= \V_z \cap \varphi_{y, z}^{-1}(\V_y) \cap \varphi_{y, z}^{-1} (\varphi_{x, y}^{-1}(\V_x))
\]
are contained in an arbitrary small neighborhood
of $\K_z \cap \psi_z^{-1}(\psi_x(\K_x) \cap \psi_y(\K_y))$.
In particular, we may assume that they are contained in $\W_{x, y, z}$.
Then Condition \ref{good (x, y, z)-relation} holds for the triples $(x, y, z)$ such that
$p_x \geq p_y \geq p_z$ for some point $p \in \psi_x(\K_x) \cap
\psi_y(\K_y) \cap \psi_z(\K_z)$.
We may also assume that for any $x, y, z \in P$,
if $\psi_x(\K_x) \cap \psi_y(\K_y) \cap \psi_z(\K_z) = \emptyset$ then
$\psi_x(\V_x \cap s_x^{-1}(0)) \cap \psi_y(\V_y \cap s_y^{-1}(0)) \cap
\psi_z(\V_z \cap s_z^{-1}(0)) = \emptyset$.
Then Condition \ref{good (x, y, z)-relation} holds for all triples $(x, y, z)$.

We can also prove that Condition \ref{good (x, y, z)-no relation} holds if
$\V_x$, $\V_y$ and $\V_z$ are sufficiently small similarly.
Hence we can construct a required good coordinate system.
\end{proof}

We can construct a totally ordered cover by the following lemma.
\begin{lem}\label{totally ordered cover}
Let $\mu : \widetilde{X} \to X$ be a locally homeomorphic and
surjective comtinuous map between Hausdorff spaces.
Assume that for each $p \in X$, $\mu^{-1}(p)$ has a partial order $\leq$ which satisfies
the following conditions:
\begin{itemize}
\item
each $\mu^{-1}(p)$ has a maximum.
\item
$\leq$ is continuous in the following sense:
For any $x, y \in \widetilde{X}$ such that $\mu(x) = \mu(y)$, if $x \leq y$, then
there exist open neighborhoods $U_x \subset \widetilde{X}$ and
$U_y \subset \widetilde{X}$ of $x$ and $y$ respectively such that
$x' \leq y'$ for any $x' \in U_x$ and $y' \in U_y$ such that $\mu(x') = \mu(y')$.
\end{itemize}
We also assume that there exists an integral-valued continuous function
$l : \widetilde{X} \to \Z$ such that $l(x) < l(y)$ if $x < y$.
Then for any compact subset $L \subset X$, there exists an open subset $V \subset
\widetilde{X}$ such that $\mu(V) \supset L$ and each nonempty fiber
$V \cap \mu^{-1}(p)$ $(p \in \mu(V))$ is totally ordered.
\end{lem}

\begin{cor}\label{totally ordered cover for pre-Kuranishi space}
Any pre-Kuranishi space has a totally ordered cover.
\end{cor}
\begin{proof}[Proof of Corollary \ref{totally ordered cover for pre-Kuranishi space}]
Apply Lemma \ref{totally ordered cover} to
$l(x) = \dim \V_x$ ($x \in \widetilde{X}$) and $L = X$.
Then $\Y = V$ is a totally ordered cover.
\end{proof}
\begin{proof}[Proof of Lemma \ref{totally ordered cover}]
For each $p \in L$, let $l(p)$ be the maximal value of $l$ on $\mu^{-1}(p)$.
Define $L_{\leq l} = \{p \in L; l(p) \leq l\}$ for each $l \in \Z$.
Note that it is compact.
For each $l \in \Z$, define an open subset $\widetilde{X}_l = \{x \in \widetilde{X};
l(x) = l\}$.
By the induction in $l$, we construct open subsets $V_l \Subset \widetilde{X}_l$
such that $V_{\leq l} = \bigcup_{k \leq l} V_k$ satisfies
$L_{\leq l} \subset \mu(V_{\leq l})$ and each fiber of $\mu|_{V_{\leq l}} : V_{\leq l} \to X$
is totally ordered.
Then $V = \bigcup_l V_l$ satisfies the conclusion of the claim.

First we consider the minimal $l$ such that $L_{\leq l} \neq \emptyset$.
Note that the restriction of $\mu$ to $\widetilde{X}_l \cap \mu^{-1}(L_{\leq l})$ is
injective.
For each $p \in \widetilde{X}_l \cap \mu^{-1}(L_{\leq l})$, let $U_p \subset \widetilde{X}_l$
be an open neighborhood of $p$ such that $\mu|_{U_p}$ is injective.
Since we can separate $\mu(p)$ and $L_{\leq l} \setminus \mu(U_p)$ by open sets,
there exist an open neighborhood $V_p \Subset U_p$ of $p$ and an open neighborhood
$W_p \subset \widetilde{X}_l$ of $\widetilde{X}_l \cap \mu^{-1}(L_{\leq l}) \setminus U_p$
such that $\mu(V_p) \cap \mu(W_p) = \emptyset$.
Choose finite points $p_i \in \widetilde{X}_l \cap \mu^{-1}(L_{\leq l})$ so that
$V_{p_i}$ covers $\widetilde{X}_l \cap \mu^{-1}(L_{\leq l})$.
Then the restriction of $\mu$ to the open neighborhood $V_l = (\bigcup_i V_{p_i}) \cap
\bigcap_i (U_{p_i} \cup W_{p_i})$ of $\widetilde{X}_l \cap \mu^{-1}(L_{\leq l})$ is injective.
Indeed, if $p \in V_{p_i}$ and $q \in U_{p_i} \cup W_{p_i}$ satisfy $\mu(p) = \mu(q)$,
then $q \notin W_{p_i}$ by the definition of $V_{p_i}$ and $W_{p_i}$.
Hence both of $p$ and $q$ is contained in $U_{p_i}$, which implies that
$p = q$ since $\mu|_{U_{p_i}}$ is injective.
Therefore the restriction of $\mu$ to $V_l$ is injective.

Next we assume that we have already constructed required open subsets
$V_k \Subset \widetilde{X}_k$ for $k < l$.
Namely, we assume that $L_{\leq k} \subset \mu(V_{\leq k})$ for $k < l$
and that each fiber of $\mu|_{V_{\leq l-1}} : V_{\leq l-1} \to X$ is totally ordered.
We construct $V_l \subset \widetilde{X}_l$ as follows.
Since $A_l = \widetilde{X}_l \cap \mu^{-1}(L_{\leq l} \setminus \mu(V_{\leq l-1}))$ consists
of maximums, the restriction of $\mu$ to $A_l$ is injective.
For each $p \in A_l$, let $U_p \subset \widetilde{X}_l$ be an open neighborhood of
$p$ which makes $\mu|_{U_p}$ injective and the following condition hold true:
If $q \in U_p$ and $r \in V_{\leq l-1}$ satisfy $\mu(q) = \mu(r)$, then $q \geq r$.
(This condition holds if $U_p$ is sufficiently small because $p \geq r$
for any $r \in \overline{V_{\leq l-1}}$ such that $\mu(p) = \mu(r)$.)
As in the case of minimal $l$, we define open subsets $V_p$ and $W_p$ for each
$p \in A_l$, and choose finite points $p_i \in A_l$ such that
$V_{p_i}$ covers $A_l$.
Then the restriction of $\mu$ to $V_l = (\bigcup_i V_{p_i}) \cap
\bigcap_i (U_{p_i} \cup W_{p_i})$ is injective, and
if $q \in V_l$ and $r \in V_{\leq l-1}$ satisfy $\mu(q) = \mu(r)$ then $q \geq r$.
Hence this $V_l$ is a required open subset.
\end{proof}

Good coordinate system is enough for the construction of the virtual fundamental chain
of one Kuranishi space, but it is not closed under product operation.
One way which was used before to overcome this problem is that
first we construct a new Kuranishi space from each good coordinate system and
reconstruct a good coordinate system of the product of the new Kuranishi spaces.
(However, this gives rise to another problem about compatibility with the various orders
of product of more than two spaces.)
Instead, we introduce a new notion of weakly good coordinate system,
which is more compatible with product.
This is defined by using the following cover of $X$ instead of a totally ordered cover.
\begin{defi}
A meet-semilattice cover of a pre-Kuranishi space $X$ is an open subset
$\Y \subset \widetilde{X}$ which satisfies $\mu(\Y) = X$ and the following conditions
for each $p \in X$:
\begin{enumerate}[label=$(\arabic*)^{\mathrm{M}}$]
\item
\label{clean intersection for meet-semilattice cover}
For any $x \in \mu^{-1}(p)$,
$\{\varphi_{x, y}(\W_{x, y}); y \leq x, y \in \Y \cap \mu^{-1}(p)\}$
intersect cleanly on a neighborhood of $\widetilde{\psi}_x^{-1}(x) \in \W_x$.
\item
\label{wedge existence for meet-semilattice cover}
For any two points $y, z \in \Y \cap \mu^{-1}(p)$, there exists some point
$w \in \Y \cap \mu^{-1}(p)$ such that
$w \leq y$, $w \leq z$, and
$\varphi_{y \vee z, w}(\W_{y \vee z, w})$ contains the intersection
$\varphi_{y \vee z, y}(\W_{y \vee z, y}) \cap \varphi_{y \vee z, z}(\W_{y \vee z, z})$
in a neighborhood of $\widetilde{\psi}_{y \vee z}^{-1}(y \vee z)$.
\end{enumerate}
(We do not assume that $x \in \Y$ or $y \vee z \in \Y$ in the above conditions.)
\end{defi}
Note that for a meet-semilattice cover $\Y$ and two points
$y, z \in \Y \cap \mu^{-1}(p)$, the point $w \in \Y \cap \mu^{-1}(p)$ which satisfies
Condition \ref{wedge existence for meet-semilattice cover} is unique.
This is easily seen as follows.
If there exist two points $w_1, w_2 \in \Y \cap \mu^{-1}(p)$ satisfying this
condition, then the images of $\varphi_{y \vee z, w_1}$ and $\varphi_{y \vee z, w_2}$
coincides in a neighborhood of $\widetilde{\psi}_{y \vee z}^{-1}(y \vee z)$.
Hence the images of $\varphi_{w_1 \vee w_2, w_1}$ and
$\varphi_{w_1 \vee w_2, w_2}$ also coincide in a neighborhood of
$\widetilde{\psi}_{w_1 \vee w_2}^{-1}(w_1 \vee w_2)$.
Condition \ref{wedge existence for meet-semilattice cover} for
$w_1, w_2 \in \Y \cap \mu^{-1}(p)$ implies that there exists some
$v \in \Y \cap \mu^{-1}(p)$ such that $v \leq w_1$, $v \leq w_2$ and
the image of $\varphi_{w_1 \vee w_2, v}$ coincides with those of
$\varphi_{w_1 \vee w_2, w_1}$ and $\varphi_{w_1 \vee w_2, w_2}$.
Hence $\varphi_{w_1, v}$ and $\varphi_{w_2, v}$ are diffeomorphisms,
which implies $w_1 = w_2 = v$.
We denote the unique point $w$ for a pair $y, z \in \Y \cap \mu^{-1}(p)$
by $y \wedge z$.

We also note that $\wedge$ is continuous, that is,
for any two points $y, z \in \Y \cap \mu^{-1}(p)$, there exist neighborhoods
$\U_y$, $\U_z$ and $\U_{y \wedge z}$ of $y$, $z$ and $y \wedge z$
in $\Y$ respectively such that for any $y' \in \U_y$, $z' \in \U_z$ and
$w' \in \U_{y \wedge z}$, if $\mu(y') = \mu(z') = \mu(w')$ then
$w' = y' \wedge z'$.
\begin{defi}
\label{def of weakly good coordinate system}
A weakly good coordinate system of a pre-Kuranishi space $X$ is a family of finite pairs
$(x, \V_x)_{x \in P}$ of points $x \in \widetilde{X}$ and open neighborhoods
$\V_x \subset \W_x$ of $\widetilde{\psi}_x^{-1}(x)$ which satisfies
the following conditions.
For two points $x, y \in P$ such that $x \unrhd y$, we define
$\V_{x, y} = \V_y \cap \varphi_{x, y}^{-1}(\V_x)$.
Then $\V_x$ and $\V_{x, y}$ satisfy the following conditions:
\begin{enumerate}[label=$(\arabic*)^{\mathrm{W}}$]
\item
\label{P meet-semilattice}
$\bigcup_{x \in P} \psi_x(\V_x \cap s_x^{-1}(0))$ is a meet-semilattice cover of $X$.
\item
\label{(x, y, z)-relation}
For any $x, y, z \in P$, if there exists some point $p \in \psi_x(\V_x \cap s_x^{-1}(0))
\cap \psi_y(\V_y \cap s_y^{-1}(0)) \cap \psi_z(\V_z \cap s_z^{-1}(0))$ such that
$p_x \geq p_y \geq p_z$, then
\begin{equation}
\varphi_{x, y}(\V_{x, y}) \cap \varphi_{x, z}(\V_{x, z}) \subset \varphi_{x, z}(\W_{x, y, z})
\label{(x, y, z)-overlap}
\end{equation}
and
\begin{equation}
\V_{y, z} \cap \varphi_{y, z}^{-1}(\V_{x, y}) \subset \W_{x, y, z}.
\label{(x, y, z)-composition}
\end{equation}
\item
\label{(x, y, z)-no relation}
For any $x, y, z \in P$ such that $\psi_x(\V_x \cap s_x^{-1}(0)) \cap
\psi_y(\V_y \cap s_y^{-1}(0)) \cap \psi_z(\V_z \cap s_z^{-1}(0)) = \emptyset$,
\begin{itemize}
\item
if $x \unrhd y$ and $x \unrhd z$, then
$\varphi_{x, y}(\V_{x, y}) \cap \varphi_{x, z}(\V_{x, z}) = \emptyset$,
\item
if $x \unrhd y$ and $y \unrhd z$, then
$\V_{x, y} \cap \varphi_{y, z}(\V_{y, z}) = \emptyset$, and
\item
if $x \unrhd z$ and $y \unrhd z$, then
$\V_{x, z} \cap \V_{y, z} = \emptyset$.
\end{itemize}
\item
\label{clean intersection for weakly good coordinate system}
For any $x, y_i \in P$ such that $x \unrhd y_i$,
$(\varphi_{x, y_i}(\V_{x, y_i}))_i$ intersect cleanly.
\item
\label{(x, y, z, w)-relation}
For any $x, y, z \in P$, if there exists some point $p \in \psi_x(\V_x \cap s_x^{-1}(0))
\cap \psi_y(\V_y \cap s_y^{-1}(0)) \cap \psi_z(\V_z \cap s_z^{-1}(0))$ such that
$p_x \geq p_y$ and $p_x \geq p_z$,
then 
there exists finite points $w_j \in P$ such that
$y \unrhd w_j$, $z \unrhd w_j$ and
\begin{equation}
\varphi_{x, y}(\V_{x, y}) \cap \varphi_{x, z}(\V_{x, z}) \subset
\bigcup_j \varphi_{x, w_j}(\V_{w_j} \cap \W_{x, y, w_j} \cap \W_{x, z, w_j})
\label{(x, y, z, w)-overlap}
\end{equation}
\end{enumerate}
\end{defi}
Condition \ref{(x, y, z)-relation} and \ref{(x, y, z)-no relation} are the same
with those for good coordinate system.
We also note that in Condition \ref{(x, y, z, w)-relation},
if $p_y \geq p_z$, then (\ref{(x, y, z, w)-overlap}) for $\{w_j\} = \{z\}$
follows from (\ref{(x, y, z)-overlap}).
(We read $\W_{x, z, z}$ as $\W_{x, z}$.)

Similarly to Lemma \ref{good coordinate from totally ordered cover},
we can prove the following.
\begin{lem}
\label{weakly good coordinate from meet-semilattice}
Assume that a meet-semilattice cover $\Y \subset \widetilde{X}$ is given.
Then for any compact subset $\K \subset \Y$,
there exists a weakly good coordinate system $(x, \V_x)_{x \in P}$ such that
$\K \subset \bigcup_{x \in P} \widetilde{\psi}_x(\V_x \cap s_x^{-1}(0)) \subset \Y$.
\end{lem}
\begin{proof}
The proof is similar to Lemma \ref{good coordinate from totally ordered cover}, but
for Condition \ref{(x, y, z, w)-relation}, we need to construct $\V_x$
by the induction in $\dim \W_x$ as follows.
We may assume that $\mu(\K) = X$.
Since we can replace $\K$ with the compact set $\{x_1 \wedge \dots \wedge x_k;
x_i \in \K\}$,
we may also assume that $\K$ is closed under $\wedge$.
Choose finite points $P = \{x\} \subset \widetilde{X}$ and compact subsets
$\K_x \subset s_x^{-1}(0)$ such that $\K = \bigcup_{x \in P} \widetilde{\psi}_x(\K_x)$.
We construct neighborhoods $\V_x$ of $\K_x$ in $\W_x$ by the induction in $\dim \W_x$.

As we saw in the proof of Lemma \ref{good coordinate from totally ordered cover},
the conditions other than \ref{clean intersection for weakly good coordinate system} and
\ref{(x, y, z, w)-relation} hold if each $\V_x$ is sufficiently
small.
First we consider Condition \ref{clean intersection for weakly good coordinate system}.
This condition holds if $\V_x$ and $\V_{y_i}$ are sufficiently small neighborhood of
$\K_x$ and $\K_{y_i}$ respectively because $\widetilde{\psi}_x(\K_x)$ and
$\widetilde{\psi}_{y_i}(\K_{y_i})$ are contained
in a meet-semilattice cover $\Y$.

Next we consider Condition \ref{(x, y, z, w)-relation}.
As in the proof of Lemma \ref{good coordinate from totally ordered cover}, we may
assume that for any triple $x, y, z \in P$,
if $\psi_x(\K_x) \cap \psi_y(\K_y) \cap \psi_z(\K_z) = \emptyset$, then
$\psi_x(\V_x \cap s_x^{-1}(0)) \cap \psi_y(\V_y \cap s_y^{-1}(0)) \cap
\psi_z(\V_z \cap s_z^{-1}(0)) = \emptyset$.
Since the case where $p_y \leq p_z$ or $p_y \geq p_z$ is contained
in Condition \ref{(x, y, z)-relation}, we may assume otherwise.
In particular, the dimension of the intersection
$\varphi_{x, y}(\V_{x, y}) \cap \varphi_{x, z}(\V_{x, z})$ is less than those of
$\W_y$ or $\W_z$.
Let $l \geq 0$ be arbitrary and assume that $\V_w$ for all $w \in P$ such that
$\dim \W_w < l$ are given.
Consider Condition \ref{(x, y, z, w)-relation}
for a triple $x, y, z \in P$ such that $\min(\dim \W_y, \dim \W_z) = l$.
Since $\K$ is closed under $\wedge$, there exists finite points $w_j \in P$
such that $y \unrhd w_j$, $z \unrhd w_j$, and
$\{\varphi_{x, w_j}(\V_{w_j} \cap \W_{x, y, w_j} \cap \W_{x, z, w_j})\}_j$ covers
$\K_x \cap \psi_x^{-1}(\psi_y(\K_y) \cap \psi_z(\K_z))$.
Hence if $\V_x$ for $x \in P$ such that $\dim W_x \geq l$ are sufficiently small
neighborhoods of $\K_x$, then Condition \ref{(x, y, z, w)-relation}
for $x, y, z \in P$ such that $\min(\dim \W_y, \dim \W_z) = l$ holds true.
Therefore we can construct neighborhoods $\V_x$ of $\K_x$ in $\W_x$ which satisfy
the conditions of weakly good coordinate system by the induction in $\dim \W_x$.
\end{proof}

\begin{defi}
Let $(x, \V_x)_{x \in P}$ be a weakly good coordinate system of
a pre-Kuranishi space $X$.
A grouped multisection $\boldsymbol{\epsilon} = (\boldsymbol{\epsilon}_x)_{x \in P}$ of
$(x, \V_x)_{x \in P}$ is a family of grouped multisections $\boldsymbol{\epsilon}_x$ of
orbibundles $(\V_x, \E_x|_{\V_x})$
which satisfies the following compatibility condition:
For any $x, y \in P$, if there exists some $p \in \psi_x(s_x^{-1}(0) \cap \V_x) \cap
\psi_y(s_y^{-1}(0) \cap \V_y)$ such that $p_x \geq p_y$, then
$\boldsymbol{\epsilon}_x$ and $\boldsymbol{\epsilon}_y|_{\V_{x, y}}$ are
$(\varphi_{x, y}, \hat \varphi_{x, y})$-related.
We emphasize that each $\boldsymbol{\epsilon}_x$ is a grouped multisection of
an orbibundle $(\V_x, \E_x|_{\V_x})$, and we do not assume that it is a grouped
multisection of an orbibundle chart.
\end{defi}

The following was proved in \cite{FO99} for the case of good coordinate system.
\begin{lem}
\label{construction of grouped multisection}
For a weakly good coordinate system $(x, \V_x)_{x \in P}$ of a pre-Kuranishi space $X$,
shrinking $\V_x$ slightly if necessary, we can construct a grouped multisection
$(\boldsymbol{\epsilon}_x)_{x \in P}$ which satisfies the following transversality
condition:
For any orbibundle chart $(\V, \E)$ in $(x, \V_x)$,
every branch of the multisection $s_x|_{\V} + \boldsymbol{\epsilon}_x|_{\V}$
is transverse to the zero section of $E$, and
its restriction to each corner of $V$ is also transverse to the zero section.
Furthermore, we can take an arbitrarily $C^\infty$-small grouped multisection.
\end{lem}
\begin{proof}
Fix a total order $\preceq$ of $P$ such that $\dim \V_y \leq \dim \V_x$ if
$y \preceq x$.
We construct the grouped multisection $\boldsymbol{\epsilon}_x$ by the induction in
$x \in P$ with respect to this order.
For the minimum $x \in P$, shrinking $\V_x$ if necessary,
we may assume that the orbibundle $(\V_x, \E_x|_{\V_x})$
is covered by finite number of orbibundle charts $(\V_{x, j}, \E_{x, j})$.
Take a smooth function $\chi_{x, j} \geq 0$ on $\V_{x, j}$ whose support in $\V_x$ is
contained in $\V_{x, j}$ for each $j$ such that $\{\{\chi_{x, j} > 0\}\}_j$ covers $\V_x$.
For each $j$, we take a smooth section $\epsilon^0_{x, j}$ of $E_{x, j} \to V_{x, j}$ and
define a multisection $\epsilon_{x, j}$ of $(\V_{x, j}, \E_x|_{\V_{x, j}})$
by $\epsilon_{x, j} = \Av(\chi_{x, j} \epsilon^0_{x, j})$.
Define a grouped multisection of $(\V_x, \E_x|_{\V_x})$ by the union
$\boldsymbol{\epsilon}_x = \coprod_j \epsilon_{x, j}$.
Sard's theorem implies that we can choose smooth sections $\epsilon^0_{x, j}$
so that every branch of $s_x + \boldsymbol{\epsilon}_x$ is transverse to
the zero section.

Assume that the grouped multisections $\boldsymbol{\epsilon}_y$ for $y \in P$
less than $x \in P$ are given.
We construct the grouped multisection $\boldsymbol{\epsilon}_x$ as follows.
First we check that $(\varphi_{x, y}, \hat \varphi_{x, y})$-relations compatibly
define $\boldsymbol{\epsilon}_x$ on $\bigcup_{y \prec x, y \unlhd x}
\varphi_{x, y}(\V_{x, y}) \subset \V_x$.
Let $y, z \in P$ be two points such that $y, z \prec x$, $y \unlhd x$, $z \unlhd x$
and $\varphi_{x, y}(\V_{x, y}) \cap \varphi_{x, z}(\V_{x, z}) \neq \emptyset$.
Condition \ref{(x, y, z)-no relation} implies that
$\psi_x(\V_x \cap s_x^{-1}(0)) \cap \psi_y(\V_y \cap s_y^{-1}(0)) \cap
\psi_z(\V_z \cap s_z^{-1}(0)) \neq \emptyset$.
$y \unlhd x$ and $z \unlhd x$ implies that $p_x \geq p_y$ and $p_x \geq p_z$
for any point $p \in \psi_x(\V_x \cap s_x^{-1}(0))
\cap \psi_y(\V_y \cap s_y^{-1}(0)) \cap \psi_z(\V_z \cap s_z^{-1}(0))$.
If $p_y \geq p_z$, then Condition \ref{(x, y, z)-relation} implies
\[
\varphi_{x, y}(\V_{x, y}) \cap \varphi_{x, z}(\V_{x, z}) \subset \varphi_{x, z}(\W_{x, y, z}).
\]
Hence the grouped multisection $\boldsymbol{\epsilon}_y$ on
$\varphi_{x, y}^{-1}(\varphi_{x, y}(\V_{x, y}) \cap \varphi_{x, z}(\V_{x, z}))$ is defined
by $(\varphi_{y, z}, \hat \varphi_{y, z})$-relation with $\boldsymbol{\epsilon}_z$.
Therefore $(\varphi_{x, y}, \hat \varphi_{x, y})$-relation and
$(\varphi_{x, z}, \hat \varphi_{x, z})$-relation are compatible on the intersection
$\varphi_{x, y}(\V_{x, y}) \cap \varphi_{x, z}(\V_{x, z})$.

Next we consider the case where $p_y \not \geq p_z$ and $p_y \not \leq p_z$.
Condition \ref{(x, y, z, w)-relation} implies that
there exists finite points $w_j \in P$ such that
$y \unrhd w_j$, $z \unrhd w_j$ and
\[
\varphi_{x, y}(\V_{x, y}) \cap \varphi_{x, z}(\V_{x, z}) \subset
\bigcup_j \varphi_{x, w_j}(\V_{w_j} \cap \W_{x, y, w_j} \cap \W_{x, z, w_j}).
\]
$p_y \not \geq p_z$ and $p_y \not \leq p_z$ imply that
$\dim \V_{w_j} < \min(\dim \V_y, \dim \V_z)$.
In particular, $w_j \prec y, z$.
The above inclusion implies that the grouped multisection $\boldsymbol{\epsilon}_y$
on $\varphi_{x, y}^{-1}(\varphi_{x, y}(\V_{x, y}) \cap \varphi_{x, z}(\V_{x, z}))$ is defined
by $(\varphi_{y, w_j}, \hat \varphi_{y, w_j})$-relations with $\boldsymbol{\epsilon}_{w_j}$,
and $\boldsymbol{\epsilon}_z$ on
$\varphi_{x, z}^{-1}(\varphi_{x, y}(\V_{x, y}) \cap \varphi_{x, z}(\V_{x, z}))$ is defined by
$(\varphi_{z, w_j}, \hat \varphi_{z, w_j})$-relations with $\boldsymbol{\epsilon}_{w_j}$.
Hence $(\varphi_{x, y}, \hat \varphi_{x, y})$-relation and
$(\varphi_{x, z}, \hat \varphi_{x, z})$-relation are compatible on the intersection.
Therefore, $(\varphi_{x, y}, \hat \varphi_{x, y})$-relations compatibly
define $\boldsymbol{\epsilon}_x$ on $\bigcup_{y \prec x, y \unlhd x}
\varphi_{x, y}(\V_{x, y})$.
We also note that Condition \ref{clean intersection for weakly good coordinate system}
implies that $(\varphi_{x, y}(\V_{x, y}))_{y \prec x, y \unlhd x}$ intersect cleanly.

Next we extend $\boldsymbol{\epsilon}_x$ defined on this subset of $\V_x$ to its
neighborhood.
We may shrink $\V_y$ ($y \prec x$) slightly if necessary for the smooth extension.
Using a smooth function $\chi \geq 0$ on $\V_x$ which satisfies $\chi \equiv 1$
on a small neighborhood of this subset and whose support is contained in a slightly larger
neighborhood, we may assume that the support of $\boldsymbol{\epsilon}_x$ is contained
in a small neighborhood of $\bigcup_{y \prec x, y \unlhd x} \varphi_{x, y}(\V_{x, y})$.
By the assumption of pre-Kuranishi space, for any $y \preceq x$,
the vertical differentials
\[
d^\bot s_x : \frac{T_{p_x}V_x}{(\phi_{x, y})_\ast T_{p_y}V_y} \tocong
\frac{(E_y)_{p_y}}{\hat\phi_{x, y} (E_x)_{p_x}}
\]
are isomorphisms for any points $p_x \in s_x^{-1}(0)$ and $p_y \in s_y^{-1}(0)$ such that
$\psi(p_x) = \psi(p_y)$.
Hence if $\boldsymbol{\epsilon}_x$ is sufficiently $C^1$-small, then
the transversality conditions for $\boldsymbol{\epsilon}_y$ ($y \prec x$) imply
that $\boldsymbol{\epsilon}_x$ also satisfies the transversality condition
on a neighborhood of $\bigcup_{y \prec x, y \unlhd x} \varphi_{x, y}(\V_{x, y})$.

On the complement of a neighborhood of $\bigcup_{y \prec x, y \unlhd x}
\varphi_{x, y}(\V_{x, y})$, as in the case of minimal $x \in P$, we take finite number of
orbibundle charts and their grouped multisections,
and add them to $\boldsymbol{\epsilon}_x$ (take the union).
Then the constructed $\boldsymbol{\epsilon}_x$ satisfies the transversality condition
and $(\varphi_{x, y}, \hat \varphi_{x, y})$-relations for all $y \prec x$.
\end{proof}

Next we consider the triangulation of the zero set of the perturbed multisection.
First we explain some notations about simplicial complex.
For a simplicial complex $K$ and its subset $A \subset K$,
we denote by $\St(A, K)$ the minimal subcomplex of $K$ which contains all simplices
intersecting with $A$.
If $K$ is embedded in a space $X$, then for a subset $U \subset X$,
we denote by $K|_U$ the subcomplex consisting of the simplices contained in $U$.
\begin{defi}
For a smooth section $s$ and a grouped multisection
$\boldsymbol{\epsilon} = (\epsilon^\omega)_{\omega \in \coprod_j \Omega_j}$ of
an orbibundle chart $(\V, \E)$,
an embedding of simplicial complex
$K = (K^{(\omega_j)})_{(\omega_j) \in \prod_j \Omega_j}$
to the zero set of $s + \boldsymbol{\epsilon}
= (s + \sum_j \epsilon^{\omega_j})_{(\omega_j) \in \prod_j \Omega_j}$
is a family of embeddings of simplicial complexes
$K^{(\omega_j)} \inj \{s + \sum_j \epsilon^{\omega_j} = 0\}$ such that
$K^{g \cdot (\omega_j)} = g \cdot K^{(\omega_j)}$ for all $g \in G_V$.
For a subset $A \subset \V$, we define
\[
\St(A, K) = (\St(\pi_V^{-1}(A), K^{(\omega_j)}))_{(\omega_j) \in \prod_j \Omega_j}
\]
and
\[
K|_A = (K^{(\omega_j)}|_{\pi_V^{-1}(A)})_{(\omega_j) \in \prod_j \Omega_j}.
\]
For a subset $B \subset \V$, we say $K$ covers $B$ if
each $K^{(\omega_j)}$ contains $\{s + \sum_j \epsilon^{\omega_j} = 0\} \cap B$
for all $(\omega_j) \in \prod_j \Omega_j$.

For a connected open subset $\U \subset \V$, fix a connected component
$U \subset \pi_V^{-1}(\U)$ and regard $\U = (U, \pi_V|_U, \U)$ as an orbichart.
Let $\boldsymbol{\epsilon}|_{\U} = (\epsilon^\omega|_U)_{\omega \in \coprod_{j \in I_U}
\Omega_j}$ be the restriction of the grouped multisection $\boldsymbol{\epsilon}$,
where $I_U = \{j; \supp (\epsilon^\omega)_{\omega \in \Omega_j} \cap
U \neq \emptyset\}$.
Let $K_{\U} = (K_U^{(\omega_j)})_{(\omega_j) \in \prod_{j \in I_U} \Omega_j}$
be an embedding of simplicial complex
to the zero set of $s|_{\U} + \boldsymbol{\epsilon}|_{\U}$.
Choose $g_k \in G_V$ so that $\pi_V^{-1}(\U) = \coprod_k g_k U$.
We say $K$ is equivalent to $K_\U$ if
$K^{(\omega_j)} = \bigcup_k g_k K^{\overline{g_k^{-1} (\omega_j)}}_U$ for all
$(\omega_j) \in \prod_j \Omega_j$,
where $\overline{(\omega_j)} \in \prod_{j \in I_U} \Omega_j$ is the image of
$(\omega_j) \in \prod_j \Omega_j$ by the projection $\prod_j \Omega_j \to
\prod_{j \in I_U} \Omega_j$.
\end{defi}

\begin{defi}
Let $(\boldsymbol{\epsilon}_x)_{x \in P}$ be a grouped multisection of a weakly good
coordinate system $(x, \V_x)_{x \in P}$ of a pre-Kuranishi space $X$
which satisfies the transversality condition in Lemma
\ref{construction of grouped multisection}.
Let $(\U_\tau, \E_x|_{\U_\tau})_{\tau \in T_x}$ be a finite family of orbibundle charts of
each orbibundle $(\V_x, \E_x)$.
Let $\mathring{\U}_\tau \Subset \U_\tau$ be their relatively compact open subsets,
and define $\mathring{\V}_x = \bigcup_{\tau \in T_x} \mathring{\U}_\tau$.
For each $\tau \in T_x$, let
$K_\tau = (K_\tau^{(\omega_j)})_{(\omega_j) \in \prod_j \Omega_{\tau, j}}$ be
an embedding of simplicial complex to the zero set of
$s_x|_{\U_\tau} + \boldsymbol{\epsilon}_x|_{\U_\tau}$.
We say $(\U_\tau, \mathring{\U}_\tau, K_\tau)_{x \in P, \tau \in T_x}$
is a triangulation of the zero set of $(s_x + \boldsymbol{\epsilon}_x)_{x \in P}$ if
the following conditions are satisfied:
\begin{enumerate}[label=$(\arabic*)^{T}$]
\item
\label{each K}
For each $\tau \in T_x$, $K_\tau$ covers $\mathring{\U}_\tau$, and
$K_\tau = \St(\mathring{\U}_\tau, K_\tau)$.
\item
\label{x intersection}
For any $x \in P$ and two indices $\tau, \tau' \in T_x$,
there exists a subset $T_{\tau, \tau'} \subset T_x$ such that
$\mathring{\U}_\tau \cap \mathring{\U}_{\tau'} = \bigcup_{\tau'' \in T_{\tau, \tau'}}
\mathring{\U}_{\tau''}$ and
$\U_{\tau''} \subset \U_\tau \cap \U_{\tau'}$ for all $\tau'' \in T_{\tau, \tau'}$.
\item
\label{K contains K}
For any $\tau, \tau' \in T_x$ such that
$\mathring{\U}_\tau \subset \mathring{\U}_{\tau'}$ and $\U_\tau \subset \U_{\tau'}$,
$\St(\mathring{\U}_\tau, K_{\tau'})$ is equivalent to $K_\tau$.
\item
\label{x y intersection}
For any two points $x, y \in P$ such that $x \unrhd y$ and any $\tau \in T_x$,
$\tau' \in T_y$, there exists a subset $T_{\tau, \tau'} \subset T_y$ such that
$\mathring{\U}_{\tau'} \cap \varphi_{x, y}^{-1}(\mathring{\U}_\tau)
= \bigcup_{\tau'' \in T_{\tau, \tau'}} \mathring{\U}_{\tau''}$ and
$\U_{\tau''} \subset \U_{\tau'} \cap \varphi_{x, y}^{-1}(\U_\tau)$ for all
$\tau'' \in T_{\tau, \tau'}$.
\item
\label{x y same}
For any two points $x, y \in P$ such that $x \unrhd y$ and any $\tau \in T_y$,
if $\U_\tau \subset \V_{x, y}$, then there exists some $\tau' \in T_x$ such that
$\U_\tau = \varphi_{x, y}^{-1}(\U_{\tau'})$ and
$\mathring{\U}_\tau = \varphi_{x, y}^{-1}(\mathring{\U}_{\tau'})$.
Furthermore, we assume that
the automorphism group of $\U_\tau$ and $\U_{\tau'}$ are isomorphic.
\item
\label{K coincides K}
For any any $\tau \in T_y$ and $\tau' \in T_x$ in Condition \ref{x y same},
let $\phi_{\tau', \tau}$ be a lift of $\varphi_{x, y}|_{\U_\tau} : \U_\tau \to \U_{\tau'}$.
Then $K_\tau = (K_\tau^\omega)_{\omega \in \coprod_j \Omega_{\tau, j}}$ and
$K_{\tau'} = (K_{\tau'}^\omega)_{\omega \in \coprod_j \Omega_{\tau', j}}$
satisfy $K_{\tau'}^{\nu^{\phi_{\tau', \tau}}_{\U_{\tau'}, \U_\tau}(\omega)}
= \phi_{\tau', \tau}(K_\tau^\omega)$ for all $\omega \in \coprod_j \Omega_{\tau, j}$.
\item
\label{again weakly good}
$(x, \mathring{\V}_x)_{x \in P}$ is also a weakly good coordinate system.
\end{enumerate}
\end{defi}
We note that in Condition \ref{x y intersection}, $\U_{\tau''}$ is contained in
$\V_{x, y} = \V_y \cap \varphi_{x, y}^{-1}(\V_x)$.
Hence Condition \ref{x y same} implies that there exists some $\tau''' \in T_x$ such that
$\U_{\tau''} = \varphi_{x, y}^{-1}(\U_{\tau'''})$ and
$\mathring{\U}_{\tau''} = \varphi_{x, y}^{-1}(\mathring{\U}_{\tau'''})$.

\begin{lem}
Let $(\boldsymbol{\epsilon}_x)_{x \in P}$ be a grouped multisection of a weakly good
coordinate system $(x, \V_x)_{x \in P}$
which satisfies the transversality condition in Lemma
\ref{construction of grouped multisection}.
Then we can construct a triangulation
$(\U_\tau, \mathring{\U}_\tau, K_\tau)_{x \in P,\tau \in T_x}$ of the zero set of
$(s_x + \boldsymbol{\epsilon}_x)_{x \in P}$.
\end{lem}
\begin{proof}
First we construct open subsets $\mathring{\U}_\tau \Subset \U_\tau \subset \V_x$
($\tau \in T_x$) which satisfy Condition \ref{x intersection}, \ref{x y intersection},
\ref{x y same} and \ref{again weakly good}.
First we take relatively compact open subsets $\mathring{\V}_x \Subset \V_x$
such that $(x, \mathring{\V}_x)_{x \in P}$ is also a weakly good coordinate system.
(We can shrink weakly good coordinate system slightly.)
Let $(\U_\tau, \E_x|_{\U_\tau})_{\tau \in T_x}$ be a family of orbibundle charts of
$(\V_x, \E_x)$ which covers the closure of $\mathring{\V}_x$, and
let $\mathring{\U}_\tau \Subset \U_\tau$ be relatively compact open subsets
such that $\mathring{\V}_x = \bigcup_{\tau \in T_x} \mathring{\U}_\tau$.
We can easily make Condition \ref{x intersection} hold by adding connected components
$\U_{\tau''}$ of $\U_\tau \cap \U_{\tau'}$ which intersect $\mathring{\U}_\tau \cap
\mathring{\U}_{\tau'}$ to $(\U_\tau)_{\tau \in T_x}$ and defining $\mathring{\U}_{\tau''}
= \U_{\tau''} \cap \mathring{\U}_\tau \cap \mathring{\U}_{\tau'}$
for each pair $\tau, \tau' \in T_x$.
Similarly, we can make Condition \ref{x y intersection} and \ref{x y same}
hold by adding appropriate open subsets to $(\U_\tau)_{\tau \in T_x}$ and
$(\mathring{\U}_\tau)_{\tau \in T_x}$.
It is easy to check that these do not break Condition \ref{again weakly good}.
Hence we can construct open subsets $\mathring{\U}_\tau \Subset \U_\tau
\subset \V_x$ ($\tau \in T_x$) which satisfy Condition
\ref{x intersection}, \ref{x y intersection}, \ref{x y same} and \ref{again weakly good}.

We can construct embeddings of simplicial complexes $K_\tau$
($\tau \in \bigcup_{x \in P} T_x$) which satisfy Condition \ref{each K},
\ref{K contains K} and \ref{K coincides K} similarly to the case of usual
triangulation of smooth manifold.
\end{proof}

For the definition of virtual fundamental chain, we need an orientation of the
pre-Kuranishi space $X$ and a strong continuous map from $X$ to a topological space.
\begin{defi}
A strong continuous map $f = (f_x)_{x \in \widetilde{X}}$ from a pre-Kuranishi space $X$
to a topological space $Y$ is a family of continuous maps $f_x : \W_x \to Y$
($x \in \widetilde{X}$) such that
$f_x \circ \varphi_{x, y} = f_y$ on $\W_{x, y}$ for all $x, y \in \widetilde{X}$ such that
$x \unrhd y$.
For a strong continuous map $f = (f_x)_{x \in \widetilde{X}}$, we define
continuous maps $\widetilde{f} : \widetilde{X} \to Y$ and $f : X \to Y$ by
the conditions $\widetilde{f} \circ \widetilde{\psi}_x = f_x$ and
$f \circ \psi_x = f_x$ on $s_x^{-1}(0)$ for all $x \in \widetilde{X}$.
If $Y$ is a smooth manifold and each $f_x$ are smooth, then we call
$f$ a strong smooth map.
\end{defi}
\begin{defi}
We say an orbibundle chart $(\V, \E)$ is orientable if $\det TV \otimes_\R \det E^\ast$
is orientable and the $G_V$-action preserves the orientation.
In this case, an orientation of $(\V, \E)$ is a homotopy type of isomorphism
$\Phi : \det TV \tocong \det E$.
We say an orbibundle is oriented if orientations of its orbibundle charts are given and
they coincide on the intersections.
\end{defi}
\begin{defi}
A pre-Kuranishi space $X$ is oriented if $(\W_x, \E_x)$ are oriented
for all $x \in \widetilde{X}$ and they satisfy the following compatibility condition:
For any $x, y \in \widetilde{X}$ and any point $p \in \psi_x(s_x^{-1}(0)) \cap
\psi_y(s_y^{-1}(0))$ such that $p_x \geq p_y$, let $(\W_{x, p}, \E_{x, p})$ and
$(\W_{y, p}, \E_{y, p})$ be orbibundle charts of $(\W_x, \E_x)$ and $(\W_y, \E_y)$
which contain $p_x$ and $p_y$ respectively such that
$\W_{y, p} \subset \varphi_{x, y}^{-1}(\W_{x, p})$.
Then the condition is that there exists a family of orientations of
$T_{p_y}W_{y, p}$, $T_{p_x}W_{x, p}$, $T_{p_x}W_{x, p} / (\phi_{x, y})_\ast T_{p_y}W_{y, p}$
$E_{x, p}|_{p_x}$, $E_{y, p}|_{p_y}$ and $E_{x, p}|_{p_x} / \hat \phi_{x, y} E_{y, p}|_{p_y}$
which makes the following isomorphisms preserve the orientations:
\begin{gather*}
T_{p_x}W_{x, p} \cong T_{p_y}W_{y, p} \oplus
T_{p_x}W_{x, p} / (\phi_{x, y})_\ast T_{p_y}W_{y, p} \\
E_{x, p}|_{p_x} \cong E_{y, p}|_{p_y} \oplus E_{x, p}|_{p_x} / \hat \phi_{x, y} E_{y, p}|_{p_y} \\
\Phi_y : \det T_{p_y}W_{y, p} \cong \det E_{y, p}|_{p_y} \\
\Phi_x : \det T_{p_x}W_{x, p} \cong \det E_{x, p}|_{p_x} \\
d^\bot s_x : T_{p_x}W_{x, p} / (\phi_{x, y})_\ast T_{p_y}W_{y, p} \cong
E_{x, p}|_{p_x} / \hat \phi_{x, y} E_{y, p}|_{p_y}
\end{gather*}
\end{defi}

\begin{defi}
Let $(\V, \E)$ be an oriented orbibundle chart whose orientation is defined by
$\Phi : \det TV \cong \det E$.
Let $s$ be its smooth section, $\boldsymbol{\epsilon} =
(\epsilon^\omega)_{\omega \in \prod_j \Omega_j}$ be its grouped multisection,
and $K = (K^{(\omega_j)})_{(\omega_j) \in \prod_j \Omega_j}$ be an embedding of
simplicial complex to the zero set of $s + \boldsymbol{\epsilon}$.
For a continuous map $f$ from $\V$ to a topological space $Y$,
we define a singular chain $f_{\#}(K)$ in $Y$ by
\[
f_{\#}(K) = \frac{1}{\# G_V \cdot \prod_j \# \Omega_j}
\sum_{(\omega_j) \in \prod_j \Omega_j} \sum_{\Delta \in (K^{(\omega_j)})^{\text{top}}}
\pm f_{\#}(\Delta),
\]
where the sum $\sum_{\Delta \in (K^{(\omega_j)})^{\text{top}}}$ is taken over
all top-dimensional simplices $\Delta$ of $K^{(\omega_j)}$, and
the sign $\pm$ of each $\Delta$ is defined as follows.
The sign is $+$ if the isomorphism
\[
T_q |\Delta| \oplus E|_q \cong T_q V
\]
given by a split of the exact sequence
\[
0 \to T_q |\Delta| \to T_q V
\xrightarrow{d^\bot (s + \sum_j \epsilon^{(\omega_j)})} E|_q \to 0
\]
preserves the orientations for all $q \in \Delta$, where
the relation of the orientations of $T_q V$ and $E|_q$ are defined by
the isomorphism $\Phi$.
Note that if there exist a connected open subset $\U \subset \V$ and
an embedding of simplicial complex $K_{\U}$ of the zero set of
$s|_\U + \boldsymbol{\epsilon}|_\U$ which is equivalent to $K$,
then $f_{\#}(K) = f_{\#}(K_\U)$.
\end{defi}

Let $f = (f_x)_{x \in \widetilde{X}}$ be a strong continuous map
from an oriented pre-Kuranishi space $X$ to a topological space $Y$.
Assume that a grouped multisection
$\boldsymbol{\epsilon} = (\boldsymbol{\epsilon}_x)_{x \in P}$ of
a weakly good coordinate system $(x, \V_x)_{x \in P}$ of $X$ and
a triangulation $(\U_\tau, \mathring{\U}_\tau, K_\tau)_{x \in P, \tau \in T_x}$ of
the zero set of $(s_x + \boldsymbol{\epsilon}_x)_{x \in P}$ are given.
For $x, y \in P$, $p \in \V_x$ and $q \in \V_y$,
we say $p$ and $q$ are equivalent ($p \sim q$) if there exist some $z \in P$ and
$r \in \V_z$ such that $x \unrhd z$, $y \unrhd z$,
$p = \varphi_{x, z}(r)$ and $q = \varphi_{y, z}(r)$.
This is indeed an equivalence relation because $(x, \V_x)_{x \in P}$ is
a weakly good coordinate system.
Define sets
\[
(s + \boldsymbol{\epsilon})^{-1}(0)|_{\U_\tau}
= \pi_{U_\tau}\biggl(\bigcup_{(\omega_j) \in \prod_j \Omega_j}
\{s_x + \sum_j \epsilon_x^{\omega_j} = 0\} \cap U_\tau\biggr)
\subset \U_\tau
\]
and
\[
(s + \boldsymbol{\epsilon})^{-1}(0) = \bigcup_{x \in P, \tau \in T_x}
(s + \boldsymbol{\epsilon})^{-1}(0)|_{\U_\tau} / \sim.
\]
Let $\pi : (s + \boldsymbol{\epsilon})^{-1}(0)|_{\U_\tau} \inj
(s + \boldsymbol{\epsilon})^{-1}(0)$ be the quotient map.
Then the assumption of $(\U_\tau, \mathring{\U}_\tau, K_\tau)_{x \in P, \tau \in T_x}$
implies that for any subsets $A_1, B_1 \subset \mathring{\U}_{\tau_1}$
and $A_2, B_2 \subset \mathring{\U}_{\tau_2}$ such that $\pi(A_1) = \pi(A_2)$
and $\pi(B_1) = \pi(B_2)$,
the singular chains $f_{\#}(\St(A_1, K_{\tau_1})|_{B_1})$ and
$f_{\#}(\St(A_2, K_{\tau_2})|_{B_2})$ coincide.

Note that if the grouped multisection $(\boldsymbol{\epsilon}_x)_{x \in P}$
is sufficiently small, then $(s + \boldsymbol{\epsilon})^{-1}(0)$ is covered by
$\{\pi(\mathring{\U}_\tau)\}_{\tau \in \bigcup_{x \in P} T_x}$.
(We need to assume that $(\boldsymbol{\epsilon}_x)_{x \in P}$ is sufficiently small.
Otherwise the zeros of the perturbed multisections leak from our open covering.)

Fix an order to the finite set $\bigcup_{x \in P} T_x$, and write it as
$\bigcup_{x \in P} T_x = \{\tau_k\}_{k = 1, 2, \dots}$.
Choose arbitrary subsets $A_k \subset \mathring{\U}_{\tau_k}$ such that
$\bigcup_k \pi(A_k) = (s + \boldsymbol{\epsilon})^{-1}(0)$, and
define $B_k = \U_{\tau_k} \setminus \pi^{-1}(\bigcup_{l < k} \pi(A_l))$.
Then we define the virtual fundamental chain $f_\ast(X)$ by
\[
f_\ast(X) = \sum_k f_{\#}(\St(A_k, K_{\tau_k})|_{B_k}).
\]
This is independent of the order of $\bigcup_{x \in P} T_x$ and
the choice of the subsets $A_k$.

In the case where the dimension of $X$ is zero,
we usually use the trivial strong continuous map to a point.
In this case, we regard the virtual fundamental chain as a rational number.

There is another way to represent the virtual fundamental chain of a pre-Kuranishi
space using differential forms.
For a strong smooth map $f = (f_x)_{x \in \widetilde{X}}$ from $X$ to a manifold
$Y$ and $h = (h_x)_{x \in \widetilde{X}}$ from $X$ to an oriented manifold $Z$,
we represent the virtual fundamental chain as a linear map
$(h_! \circ f^\ast)_X : \Omega(Y) \to \Omega(Z)$ as follows.
If $Z$ is a point, then this map $(h_! \circ f^\ast)_X : \Omega(Y) \to \R$
is the dual representation of the virtual fundamental chain $f_\ast(X)$.
In this case, we often denote the value of this map at $\theta \in \Omega(Y)$
by $\int_{X} f^\ast \theta$.

Let $(x, \V_x)_{x \in P}$ be a weakly good coordinate system of
a pre-Kuranishi space $X$, and let $\beta_x : \V_x \to \R$ be a smooth function
with compact support for each $x \in P$.
Define $\Y = \bigcup_{x \in P} \widetilde{\psi}_x(\V_x \cap s_x^{-1}(0))$.
Note that for any $p \in X$, $\mu^{-1}(p) \cap \Y$ has the unique minimum
$p^\Y_{\min}$ since $\Y$ is a meet-semilattice cover.
We say $(\beta_x)_{x \in P}$ is a partition of unity subordinate to $(x, \V_x)_{x \in P}$
if for any $p \in X$,
$\sum_{x \in P, x \unrhd p^\Y_{\min}} \beta_x \circ \varphi_{x, p^\Y_{\min}} \equiv 1$
on a neighborhood of $\psi_{p^\Y_{\min}}^{-1}(p)$ in $\W_{p^\Y_{\min}}$.

Let $\boldsymbol{\epsilon} = (\boldsymbol{\epsilon}_x)_{x \in P}$ be
a grouped multisection of $(x, \V_x)_{x \in P}$ which satisfies the transversality
condition in Lemma \ref{construction of grouped multisection}.
We assume that the restriction of $h$ to the zero set of each branch of
the multisections $s_x + \boldsymbol{\epsilon}_x$ is submersive.
We can construct such a perturbed multisection if $Z$ is a point.
(In general, we need to use continuous family of multisections.
See Section \ref{continuous family of multisections}.)
We further assume that $\boldsymbol{\epsilon}$ is sufficiently small so that
$\sum_{z \in P, z \unrhd x} \beta_z \circ \varphi_{z, x} = 1$ on
$(s_x + \boldsymbol{\epsilon}_x)^{-1}(0) \cap (\V_x)_{\min}$ for any $x \in P$,
where $(s_x + \boldsymbol{\epsilon}_x)^{-1}(0) \subset \V_x$ is the set of points
at which one of the branchs of the multisection $s_x + \boldsymbol{\epsilon}_x$
takes zero, and $(\V_x)_{\min} \subset \V_x$ is the set of points $q \in \V_x$ such that
there do not exist any $y \in P$ such that $x \unrhd y$, $\dim \W_x > \dim \W_y$
and $q \in \varphi_{x, y}(\V_{x, y})$.
For each $x \in P$, we take finite orbibundle charts $(\V_\tau, \E_\tau)_{\tau \in T_x}$
of $(\V_x, \E_x)$ and smooth functions $\beta_\tau : \V_\tau \to \R$
with compact support such that $\beta_x = \sum_{\tau \in T_x} \beta_\tau$.
Then for each differential form $\theta \in \Omega(Y)$,
$(h_! \circ f^\ast)_X \theta \in \Omega(Z)$ is defined by
\begin{align}
&(h_! \circ f^\ast)_X \theta \notag \\
&= \sum_{x \in P, \tau \in T_x}
\frac{\sum_{(\omega_j) \in \prod_j \Omega_{\tau, j}}
\Bigl(h_x|_{\{s_\tau^{(\omega_j)} = 0\}}\Bigr)_{\textstyle !}\,
\Bigl(\beta_\tau \cdot (f_x|_{V_\tau})^\ast
\theta|_{\{s_\tau^{(\omega_j)} = 0\}}\Bigr)}
{\# G_{V_\tau}
\cdot \prod_j \# \Omega_{\tau, j}},
\label{virtual fundamental chain by forms}
\end{align}
where $s_\tau = s_x|_{\V_\tau}$, $\boldsymbol{\epsilon}_x|_{\V_\tau}
= (\epsilon^\omega)_{\omega \in \coprod_j \Omega_{\tau, j}}$,
$s_\tau^{(\omega_j)} = s_\tau + \sum_j \epsilon_\tau^{(\omega_j)}$, and
$(h_x|_{\{s_\tau^{(\omega_j)} = 0\}})_!$ is integration
along fiber for the fibration
$h_x : \{s_\tau^{(\omega_j)} = 0\} \to Z$.
In our convention, the orientation of the fiber $F$ is defined by
$T_{h_x(p)}Z \oplus T_pF = T_p\{s_\tau^{(\omega_j)} = 0\}$ at each point $p \in F$.
It is easy to check that $(h_! \circ f^\ast)_X \theta$ is independent of the choice of
the partition of unity $(\beta_x)_{x \in P}$ and functions $\beta_\tau$.

If $Z$ is non-orientable, then instead of a compatible family of orientations
$\Phi_x : \det TW_x \cong \det E_x$, we assume that a compatible family of
isomorphisms $\widetilde{\Phi}_x : \det TW_x \cong
h_x^\ast \mathcal{O}_Z \otimes \det E_x$
is given, where $\mathcal{O}_Z = \det TZ$ is the orientation bundle of $Z$.
Then we can define the orientation of the fiber of
each $h_x : \{s_\tau^{(\omega_j)} = 0\} \to Z$ and define
$(h_! \circ f^\ast)_X \theta : \Omega(Y) \to \Omega(Z)$ similarly.

%

\subsection{Compatible perturbed multisections}
\label{compatible perturbed multisection}
In application, we need to construct perturbed multisections of moduli spaces
which respect their algebraic properties.
In this section, we consider the compatibility of the perturbed multisections
of various pre-Kuranishi spaces.
\subsubsection{The boundary of a pre-Kuranishi space}
First, we consider compatibility of the grouped multisection of a pre-Kuranishi space
and the grouped multisection of its boundary.
\begin{defi}
For a pre-Kuranishi space $X$ with corners, we define the boundary
$\partial X \subset X$ by the set of points $p \in X$ such that
for any $x \in \mu^{-1}(p)$, $\widetilde{\psi}_x^{-1}(x)$ is contained in the boundary of
$\W_x$.
(This condition is independent of the choice of $x \in \mu^{-1}(p)$.)
The restriction of the pre-Kuranishi structure of $X$ defines
the pre-Kuranishi structure of $\partial X$.
\end{defi}
For a weakly good coordinate system $(x, \V_x)_{x \in P}$ of $X$,
$(x, \partial \V_x)_{x \in P \cap \mu^{-1}(\partial X)}$ is a weakly good coordinate
system of $\partial X$.
(We note that $\partial \V_x = \emptyset$ for $x \in P \setminus \mu^{-1}(\partial X)$
because $\partial \W_x = \emptyset$ by assumption.)
Conversely, for a weakly good coordinate system
$(x, \V^{\partial X}_x)_{x \in P^{\partial X}}$ of $\partial X$,
we can construct a weakly good coordinate system
$(x, \V_x)_{x \in P}$ of $X$ such that
$(x, \partial \V_x)_{x \in P \cap \mu^{-1}(\partial X)}
= (x, \V^{\partial X}_x)_{x \in P^{\partial X}}$
by the following lemma
and Lemma \ref{weakly good coordinate from meet-semilattice}.
\begin{lem}
\label{extension of meet-semilattice cover}
For any meet-semilattice cover $\Y^{\partial X}$ of $\partial X$,
there exists a meet-semilattice cover $\Y$ of $X$ such that
$\Y \cap \mu^{-1}(\partial X) = \Y^{\partial X}$.
\end{lem}
\begin{proof}
First we construct open neighborhood $\Y^{N(\partial X)}$ of
$\Y^{\partial X} \subset \widetilde{X}$ such that
$\Y^{N(\partial X)} \cap \partial \widetilde{X} = \Y^{\partial X}$ and
$\Y^{N(\partial X)} \subset \widetilde{X}$ satisfies the conditions of
meet-semilattice cover other than the covering condition $\mu(\Y) = X$
as follows.

Let $\Y^{N(\partial X)}_{(0)} \subset \widetilde{X}$ be an open neighborhood of
$\Y^{\partial X}$ such that $\Y^{N(\partial X)}_{(0)} \cap \partial \widetilde{X}
= \Y^{\partial X}$.
We construct a decreasing sequence of open neighborhoods
$\Y^{N(\partial X)}_{(k)} \subset \Y^{N(\partial X)}_{(0)}$ ($k \geq 0$)
of $\Y^{\partial X}$ which satisfies the following conditions,
where $\Y^{N(\partial X)}_{(k), l} = \{x \in \Y^{N(\partial X)}_{(k)}; \dim \W_x = l\}$
for each $l \geq 0$.
\begin{enumerate}[label=(\roman*)]
\item
$\Y^{N(\partial X)}_{(k + 1)} \subset \Y^{N(\partial X)}_{(k)}$ for all $k \geq 0$.
\item
$\Y^{N(\partial X)}_{(k), l} = \Y^{N(\partial X)}_{(l), l}$ for $k \geq l$.
\item
\label{clean intersection for k and (k-1)}
For any $k \geq 1$, $m \geq 2$ and
$y_1, \dots, y_m \in \Y^{N(\partial X)}_{(k)}$ such that $\mu(y_i) = \mu(y_1)$
for all $1 \leq i \leq m$, define $x = y_1 \vee \dots \vee y_m \in \widetilde{X}$.
Then $(\varphi_{x, y_i}(\W_{x, y_i}))_{1 \leq i \leq m}$ intersect cleanly on
a neighborhood of $\widetilde{\psi}_x^{-1}(x) \subset \W_x$.
Furthermore, there exists some $w \in \Y^{N(\partial X)}_{(k-1)}$
such that $\mu(w) = \mu(y_1)$, $w \leq y_i$ for all $1 \leq i \leq m$,
and the intersection of $(\varphi_{x, y_i}(\W_{x, y_i}))_{1 \leq i \leq m}$
is contained in the image of $\varphi_{x, w}$ on a neighborhood of
$\widetilde{\psi}_x^{-1}(x) \subset \W_x$.
\end{enumerate}
Then $\Y^{N(\partial X)} = \bigcup_{l \geq 0} \Y^{N(\partial X)}_{(l), l}
\subset \widetilde{X}$ satisfies the required conditions.

We construct $\Y^{N(\partial X)}_{(k)} \subset \Y^{N(\partial X)}_{(0)}$
by the induction in $k \geq 1$.
Assume that we have already constructed $\Y^{N(\partial X)}_{(k)}$ for $k < k_0$.
Define $\widehat{\Y} = \Y^{N(\partial X)}_{(k_0-1)}$.
Let $\{\U_j\}_j$ be a finite open covering of $X$ such that
each $\widehat{\Y}|_{\overline{\U_j}}$ can be decomposed into its
finite disjoint open subsets
$\coprod_a \widehat{\Y}^{j, a}$ such that
each $\mu|_{\widehat{\Y}^{j, a}}$ is injective.
We say $p \in X$ is a good point if
Condition \ref{clean intersection for k and (k-1)} is satisfied for
$k = k_0$, $\Y^{N(\partial X)}_{(k_0)} = \widehat{\Y}$ and any points
$y_1, \dots, y_m \in \widehat{\Y} \cap \mu^{-1}(p)$.
Note that every point in $\partial X$ is good.

We shrink $\widehat{\Y}$ on each $\overline{\U_j}$ by the induction in $j$
to make every point in $\overline{\U_j}$ a good point by the following argument.

We shrink $\widehat{\Y}|_{\overline{\U_j}}$ by the induction in $m_0 \geq 2$
so that Condition \ref{clean intersection for k and (k-1)} is satisfied
for $k = k_0$, $2 \leq m \leq m_0$ and
any points $y_1, \dots, y_m \in \widehat{\Y}|_{\overline{\U_j}}$
such that $\mu(y_i) = \mu(y_1)$ as follows.
Assume that Condition \ref{clean intersection for k and (k-1)} is satisfied
for $k = k_0$, $2 \leq m \leq m_0 - 1$ and
any points $y_1, \dots, y_m \in \widehat{\Y}|_{\overline{\U_j}}$.
For each sequence $a_1, \dots, a_{m_0}$,
let $C_{a_1, \dots, a_{m_0}} \subset \bigcap_{1 \leq i \leq m_0} \mu(\widehat{\Y}^{j, a})$
be the set of points $p \in \bigcap_{1 \leq i \leq m_0} \mu(\widehat{\Y}^{j, a_i})$ such
that Condition \ref{clean intersection for k and (k-1)} does not hold for
$k = k_0$, $m = m_0$ and
$y_i = (\mu|_{\widehat{\Y}^{j, a_i}})^{-1}(p)$ ($1 \leq i \leq m_0$).
Then $C_{a_1, \dots, a_{m_0}}$ is closed in the relative topology of
$\bigcap_{1 \leq i \leq m_0} \mu(\widehat{\Y}^{j, a})$.
(Its complement is open.)
Furthermore, $C_{a_1, \dots, a_{m_0}}$ does not contain any good point.

Fix a distance of $X$ and
define $(\widehat{\Y}^{j, a_i})^\ast \subset \widehat{\Y}^{j, a_i}$ by the interior of
\begin{align}
\bigl(\widehat{\Y}^{j, a_i} \setminus \mu^{-1}(C_{a_1, \dots, a_{m_0}})\bigr)
\cup \{x \in \widehat{\Y}^{j, a_i}&; \dist(\mu(x), \mu(\widehat{\Y}^{j, a_i}) \cap \partial X)
\notag \\
& < \max_{i' \neq i} \dist(\mu(x), \mu(\widehat{\Y}^{j, a_{i'}}) \cap \partial X)\}.
\label{in order to shrink Y}
\end{align}
(We read $\dist(\mu(x), \emptyset)$ as $\dist(\mu(x), \emptyset) = \infty$.)
We claim that
$(\widehat{\Y}^{j, a_i})^\ast$ contains $\widehat{\Y}^{j, a_i} \cap \partial \widetilde{X}$.
This is proved as follows.
For any point $x \in \widehat{\Y}^{j, a_i} \cap \partial \widetilde{X}$,
if $x \in \bigcap_{i'} \widehat{\Y}^{j, a_{i'}} \cap \partial \widetilde{X}$, then
$x$ is contained in the open subset $\bigcap_{i'} \widehat{\Y}^{j, a_{i'}} \setminus
\mu^{-1}(C_{a_1, \dots, a_{m_0}})$,
which is contained in (\ref{in order to shrink Y}).
If $x \notin \widehat{\Y}^{j, a_{i'}}$ for some $i' \neq i$,
then $x$ is contained in the open subset
\[
\{x \in \widehat{\Y}^{j, a_i}; \dist(\mu(x), \mu(\widehat{\Y}^{j, a_i}) \cap \partial X)
< \max_{i' \neq i} \dist(\mu(x), \mu(\widehat{\Y}^{j, a_{i'}}) \cap \partial X)\},
\]
which is contained in (\ref{in order to shrink Y}).
Hence
$(\widehat{\Y}^{j, a_i})^\ast$ contains $\widehat{\Y}^{j, a_i} \cap \partial \widetilde{X}$,
and we may replace $\widehat{\Y}^{j, a_i}$ with
$(\widehat{\Y}^{j, a_i})^\ast$.
Then $C_{a_1, \dots, a_{m_0}}$ becomes the empty set.

We apply the same argument for all sequences $a_1, \dots, a_{m_0}$.
Then Condition \ref{clean intersection for k and (k-1)} is satisfied
for $k = k_0$, $2 \leq m \leq m_0$ and
any points $y_1, \dots, y_m \in \widehat{\Y}|_{\overline{\U_j}}$
such that $\mu(y_i) = \mu(y_1)$.
Hence the induction in $m_0$ works and
we can shrink $\widehat{\Y}|_{\overline{\U_j}}$ so that
Condition \ref{clean intersection for k and (k-1)} is satisfied
for $k = k_0$, $m \geq 2$ and
any points $y_1, \dots, y_m \in \widehat{\Y}|_{\overline{\U_j}}$
such that $\mu(y_i) = \mu(y_1)$.
Therefore, the induction in $j$ also works, and
we obtain an open neighborhood
$\widehat{\Y} \subset \Y^{N(\partial X)}_{(k_0-1)}$ of $\Y^{\partial X}$
such that Condition \ref{clean intersection for k and (k-1)} holds for
$k = k_0$, $\Y^{N(\partial X)}_{(k_0)} = \widehat{\Y}$ and any points
$y_1, \dots, y_m \in \widehat{\Y}$ such that $\mu(y_i) = \mu(y_1)$.
Then we define $\Y^{N(\partial X)}_{(k_0)} = \widehat{\Y}$.

Note that if $\Y^{N(\partial X)}_{(k), l} = \Y^{N(\partial X)}_{(l), l}$ for
$l \leq k \leq k_0 -1$, then it also holds for $l \leq k = k_0$
by the above construction.
Hence we can construct a decreasing sequence of open neighborhoods
$\Y^{N(\partial X)}_{(k)} \subset \Y^{N(\partial X)}_{(0)}$ ($k \geq 0$)
of $\Y^{\partial X}$ inductively, and
$\Y^{N(\partial X)} = \bigcup_{l \geq 0} \Y^{N(\partial X)}_{(l), l}
\subset \widetilde{X}$ satisfies the conditions of
meet-semilattice cover other than the covering condition $\mu(\Y) = X$.

Take an open subset $\mathring{X} \subset X$ such that
$\mathring{X} \cup \mu(\Y^{N(\partial X)}) =X$ and
$\overline{\mathring{X}} \cap \partial X = \emptyset$.
Let $N_0(\partial X)$ be an open neighborhood of $\partial X \subset X$ such that
$N_0(\partial X) \cap \mathring{X} = \emptyset$.
We construct a totally ordered cover $\mathring{Y}$ of $X$ such that
$\mathring{Y} \subset \mu^{-1}(\mathring{X}) \cup \Y^{N(\partial X)}$.
Then $\Y = \mathring{Y} \cup (\Y^{N(\partial X)} \cap \mu^{-1}(N_0(\partial X)))$
is a required meet-semilattice cover of $X$.
(Since $\mathring{Y} \cap \mu^{-1}(N_0(\partial X))$ is contained in $\Y^{N(\partial X)}$,
$\Y$ satisfies the conditions of meet-semilattice cover.)

We explain the construction of the totally ordered cover $\mathring{Y}$ of $X$.
Applying Lemma \ref{totally ordered cover} to the map
$\mu
: \Y^{N(\partial X)} \cap \mu^{-1}(X \setminus \mathring{X})
\to X \setminus \mathring{X}$ and the opposite partial order
``$\preccurlyeq$''$=$``$\geq$'',
we get a totally ordered cover $\Y^{X \setminus \mathring{X}}$ of
$X \setminus \mathring{X}$ contained in $\Y^{N(\partial X)}$.
Then we can apply Lemma \ref{totally ordered cover} to the map
$\mu
: \mu^{-1}(\mathring{X}) \cup \Y^{X \setminus \mathring{X}}
\to X$
and the partial order $\leq$, and we get a totally ordered cover $\mathring{Y}$ of $X$
such that
$\mathring{Y} \subset \mu^{-1}(\mathring{X}) \cup \Y^{X \setminus \mathring{X}}
\subset \mu^{-1}(\mathring{X}) \cup \Y^{N(\partial X)}$.
Hence we can construct a required meet-semilattice cover.
\end{proof}
By the argument used for the proof of Lemma
\ref{extension of grouped multisection for embedding},
we can extend a grouped multisection of
$(x, \partial \V_x)_{x \in P \cap \mu^{-1}(\partial X)}$ to
a grouped multisection of $(x, \V_x)_{x \in P}$ if we shrink $\V_x$ slightly.
(The only difference of this extension and the extension proved in
Lemma \ref{extension of grouped multisection for embedding}
is whether the rank of the obibundle changes or not,
which has nothing to do with the construction of the extension.)

\subsubsection{Pull back by submersion}
Next we define the pull back of the perturbed multisection
by a submersion from a pre-Kuranishi space to another.
First we define the submersion between pre-Kuranishi spaces.
\begin{defi}
Let $X^k$ ($k = 1,2$) be two pre-Kuranishi spaces with pre-Kuranishi structures
$(\widetilde{X}^k, \mu^k, (\W^k_x, \E^k_x, s^k_x, \widetilde{\psi}^k_x),
(\varphi^k_{x, y}, \hat \varphi^k_{x, y}))$.
A submersion $f = (f, \widetilde{f}, (\varphi^f_x, \hat \varphi^f_x))$
from $X^1$ to $X^2$ consists of
continuous maps $f : X^1 \to X^2$ and $\widetilde{f} : \widetilde{X}^1 \to \widetilde{X}^2$
such that $\widetilde{f} \circ \mu^1 = \mu^2 \circ f$, and
submersions $(\varphi^f_x, \hat \varphi^f_x)$ ($x \in \widetilde{X}^1$) from
$(\W^1_x, \E^1_x)$ to $(\W^2_{\widetilde{f}(x)}, \E^2_{\widetilde{f}(x)})$
which satisfy the following conditions:
\begin{enumerate}[label=$(\arabic*)^{\mathrm{S}}$]
\item
\label{isom of pos for submersion}
For each $p \in X^1$, $\widetilde{f}|_{(\mu^1)^{-1}(p)} : (\mu^1)^{-1}(p) \cong
(\mu^2)^{-1}(f(p))$
is an isomorphism of partially ordered sets.
\item
$s^2_{\widetilde{f}(x)} \circ \varphi^f_x = \hat \varphi^f_x \circ s^1_x$ on $\W^1_x$
and $\widetilde{\psi}^2_{\widetilde{f}(x)} \circ \varphi^f_x
= \widetilde{f} \circ \widetilde{\psi}^1_x$ on $(s^1_x)^{-1}(0)$ for all $x \in \widetilde{X}^1$.
\item
\label{compatibility of submersion with embedding}
For any $x, y \in \widetilde{X}^1$, if $x \unrhd y$, then
$\varphi^f_y(\W^1_{x, y}) \subset \W^2_{\widetilde{f}(x), \widetilde{f}(y)}$ and
\begin{equation}
\label{submersion embedding compatibility equation}
(\varphi^f_x, \hat \varphi^f_x) \circ (\varphi^1_{x, y}, \hat \varphi^1_{x, y})
= (\varphi^2_{\widetilde{f}(x), \widetilde{f}(y)},
\hat \varphi^2_{\widetilde{f}(x), \widetilde{f}(y)})
\circ (\varphi^f_y, \hat \varphi^f_y)
\end{equation}
on $\W^1_{x, y}$.
\end{enumerate}
\end{defi}
\begin{rem}
\label{shrink W_x for submersion}
In the above definition, we assume that each $(\varphi^f_x, \hat \varphi^f_x)$
is defined on the whole of $(\W^1_x, \E^1_x)$.
Hence for the construction of a submersion $f$, we sometimes need
to shrink each $\W^1_x$ as we explained in Remark \ref{we can shrink Wx}.
Similarly, we sometimes need to replace $\widetilde{X}^1$ to
its appropriate open subset.
(See Remark \ref{quotient map is a submersion} for example.)
\end{rem}
\begin{rem}
In Condition \ref{compatibility of submersion with embedding},
we consider the composition of a submersion and an embedding.
In (\ref{submersion embedding compatibility equation}),
the lifts of the right hand side satisfy the same conclusion as that of
Lemma \ref{effective action} as we saw in
Remark \ref{effective action for composition of embedding and submersion}.
Although we cannot apply this lemma to the left hand side a priori,
(\ref{submersion embedding compatibility equation}) implies that
their lifts also satisfy the same property.
\end{rem}
\begin{defi}
\label{compatible weakly good coordinate systems for submersion}
Let $(x_1, \V^1_{x_1})_{x_1 \in P^1}$ and $(x_2, \V^2_{x_2})_{x_2 \in P^2}$ be
weakly good coordinate systems of pre-Kuranishi spaces $X^1$ and $X^2$ respectively.
We say these are compatible with respect to the submersion $f$ if
they satisfy the following conditions:
\begin{enumerate}[label=(\roman*)]
\item
\label{meet-semilattice cover compatibility}
The meet-semilattice covers $\Y^i = \bigcup_{x_i \in P^i}
\widetilde{\psi}^i_{x_i}((s^i_{x_i})^{-1}(0) \cap \V^i_{x_i}))$ ($i = 1,2$) satisfy
$\widetilde{f}(\Y^1) \subset \Y^2$.
\item
\label{V_{x_1} subset V_{x_2}}
For any $x_1 \in P^1$ and $x_2 \in P^2$,
if $\widetilde{f}(x_1) \in \widetilde{\psi}^2_{x_2}((s^2_{x_2})^{-1}(0) \cap \V^2_{x_2}))$,
then $\varphi^f_{x_1}(\V^1_{x_1}) \subset
(\varphi^2_{x_2, \widetilde{f}(x_1)})^{-1}(\V^2_{x_2})$.
\item
\label{composition compatibility for submersion}
For any $y_1, x_1 \in P^1$ and $y_2, x_2 \in P^2$,
if $\widetilde{f}(x_1) \in \widetilde{\psi}^2_{x_2}((s^2_{x_2})^{-1}(0) \cap \V^2_{x_2}))$,
$\widetilde{f}(y_1) \in \widetilde{\psi}^2_{y_2}((s^2_{y_2})^{-1}(0) \cap \V^2_{y_2}))$, and
there exists some $p \in \psi_{x_1}((s^1_{x_1})^{-1}(0) \cap \V^1_{x_1}) \cap
\psi_{y_1}((s^1_{y_1})^{-1}(0) \cap \V^1_{y_1})$ such that $p_{x_1} \geq p_{y_1}$,
then
\[
\varphi_{y_1}^f(\V^1_{x_1, y_1}) \subset \W^2_{x_2, \widetilde{f}(x_1), \widetilde{f}(y_1)}
\cap \W^2_{x_2, y_2, \widetilde{f}(y_1)}.
\]
(Note that the previous condition implies that
$\widetilde{f}(p_{x_1}) \in \widetilde{\psi}^2_{x_2}((s^2_{x_2})^{-1}(0) \cap \V^2_{x_2}))$,
$\widetilde{f}(p_{y_1}) \in \widetilde{\psi}^2_{y_2}((s^2_{y_2})^{-1}(0) \cap \V^2_{y_2}))$ and
$\widetilde{f}(p_{x_1}) \geq \widetilde{f}(p_{y_1})$.)
\end{enumerate}
\end{defi}
If $(x_1, \V^1_{x_1})_{x_1 \in P^1}$ and $(x_2, \V^2_{x_2})_{x_2 \in P^2}$ are
compatible with respect to the submersion $f$,
then for a grouped multisection $\boldsymbol{\epsilon}^2
= (\boldsymbol{\epsilon}^2_{x_2})_{x_2 \in P^2}$ of $(x_2, \V^2_{x_2})_{x_2 \in P^2}$,
we can define a grouped multisection $\boldsymbol{\epsilon}^1
= (\boldsymbol{\epsilon}^1_{x_1})_{x_1 \in P^1}$ of $(x_1, \V^1_{x_1})_{x_1 \in P^1}$
as follows.
For each $x_1 \in P^1$, Condition \ref{meet-semilattice cover compatibility}
implies that there exists some $x_2 \in P^2$ such that
$\widetilde{f}(x_1) \in \widetilde{\psi}^2_{x_2}((s^2_{x_2})^{-1}(0) \cap \V^2_{x_2}))$, and
Condition \ref{V_{x_1} subset V_{x_2}} implies that
$\varphi^f_{x_1}(\V^1_{x_1}) \subset
(\varphi^2_{x_2, \widetilde{f}(x_1)})^{-1}(\V^2_{x_2})$.
We define $\boldsymbol{\epsilon}^1_{x_1}$ by the pull back of
$\boldsymbol{\epsilon}^2_{x_2}$ by
$(\varphi^2_{x_2, \widetilde{f}(x_1)}, \hat \varphi^2_{x_2, \widetilde{f}(x_1)}) \circ
(\varphi^f_{x_1}, \hat \varphi^f_{x_1})$.

We need to check that $\boldsymbol{\epsilon}^1
= (\boldsymbol{\epsilon}^1_{x_1})_{x_1 \in P^1}$ satisfies
$(\varphi^1_{x_1, y_1}, \hat \varphi^1_{x_1, y_1})$-relations.
It is enough to prove the following claim:
If $y_1, x_1 \in P^1$ and $y_2, x_2 \in P^2$ satisfy
$\widetilde{f}(x_1) \in \widetilde{\psi}^2_{x_2}((s^2_{x_2})^{-1}(0) \cap \V^2_{x_2}))$,
$\widetilde{f}(y_1) \in \widetilde{\psi}^2_{y_2}((s^2_{y_2})^{-1}(0) \cap \V^2_{y_2}))$, and
there exists some $p \in \psi_{x_1}((s^1_{x_1})^{-1}(0) \cap \V^1_{x_1}) \cap
\psi_{y_1}((s^1_{y_1})^{-1}(0) \cap \V^1_{y_1})$ such that $p_{x_1} \geq p_{y_1}$,
then the pull backs
$(\varphi^2_{x_2, \widetilde{f}(x_1)} \circ \varphi^f_{x_1})^\ast
\boldsymbol{\epsilon}^2_{x_2}$ and
$(\varphi^2_{y_2, \widetilde{f}(y_1)} \circ \varphi^f_{y_1})^\ast
\boldsymbol{\epsilon}^2_{y_2}|_{\V^1_{x_1, y_1}}$ are
$(\varphi_{x_1, y_1}, \hat \varphi_{x_1, y_1})$-related.
Note that in particular, the case where $y_1 = x_1$ and $y_2 \neq x_2$ implies that
the definition of $\boldsymbol{\epsilon}^1_{x_1}$ does not depend on the
choice of $x_2$.


We can prove this claim as follows.
Condition \ref{compatibility of submersion with embedding} and
\ref{composition compatibility for submersion} imply that
\begin{align*}
&(\varphi^2_{x_2, \widetilde{f}(x_1)}, \hat \varphi^2_{x_2, \widetilde{f}(x_1)}) \circ
(\varphi^f_{x_1}, \hat \varphi^f_{x_1}) \circ
(\varphi^1_{x_1, y_1}, \hat \varphi^1_{x_1, y_1}) \\
&= (\varphi^2_{x_2, \widetilde{f}(x_1)}, \hat \varphi^2_{x_2, \widetilde{f}(x_1)}) \circ
(\varphi^2_{\widetilde{f}(x_1), \widetilde{f}(y_1)},
\hat \varphi^2_{\widetilde{f}(x_1), \widetilde{f}(y_1)}) \circ
(\varphi^f_{y_1}, \hat \varphi^f_{y_1}) \\
&= (\varphi^2_{x_2, \widetilde{f}(y_1)}, \hat \varphi^2_{x_2, \widetilde{f}(y_1)}) \circ
(\varphi^f_{y_1}, \hat \varphi^f_{y_1}) \\
&= (\varphi^2_{x_2, y_2}, \hat \varphi^2_{x_2, y_2}) \circ
(\varphi^2_{y_2, \widetilde{f}(y_1)}, \hat \varphi^2_{y_2, \widetilde{f}(y_1)})
\circ (\varphi^f_{y_1}, \hat \varphi^f_{y_1})
\end{align*}
on $(\V^1_{x_1, y_1}, \E^1_{y_1}|_{\V^1_{x_1, y_1}})$.
This and $(\varphi^2_{x_2, y_2}, \hat \varphi^2_{x_2, y_2})$-relation of
$\boldsymbol{\epsilon}^2_{x_2}$ and $\boldsymbol{\epsilon}^2_{y_2}$ imply
that $(\varphi^2_{x_2, \widetilde{f}(x_1)} \circ \varphi^f_{x_1})^\ast
\boldsymbol{\epsilon}^2_{x_2}$ and
$(\varphi^2_{y_2, \widetilde{f}(y_1)} \circ \varphi^f_{y_1})^\ast
\boldsymbol{\epsilon}^2_{y_2}|_{\V^1_{x_1, y_1}}$ are
$(\varphi_{x_1, y_1}, \hat \varphi_{x_1, y_1})$-related.
Therefore $\boldsymbol{\epsilon}^1
= (\boldsymbol{\epsilon}^1_{x_1})_{x_1 \in P^1}$ is a well-defined grouped multisection
of $(x_1, \V^1_{x_1})_{x_1 \in P^1}$.
We call this grouped multisection the pull back of $\boldsymbol{\epsilon}^2$ by the
submersion $f$, and denote it by $f^\ast \boldsymbol{\epsilon}^2$.


Next we claim that for a meet-semilattice cover $\Y^2 \subset \widetilde{X}^2$
of $X^2$, $\Y^1 = \widetilde{f}^{-1}(\Y^2)$ is a meet-semilattice cover of $X^1$.
To see this, first we note the following fact.
For any $x, y \in \widetilde{X}^1$ such that $x \unrhd y$,
\begin{equation}
\dim \W_x^1- \dim \W_y^1 = \dim \W_{\widetilde{f}(x)}^2 - \dim \W_{\widetilde{f}(y)}^2.
\label{dimension relation for inverse cover}
\end{equation}
This is because $\rank \E_x^1 = \rank \E_{\widetilde{f}(x)}^2$,
$\rank \E_y^1 = \rank \E_{\widetilde{f}(y)}^2$,
$\dim \W_x^1 - \rank \E_x^1 = \dim \W_y - \rank \E_y^1$ and
$\dim \W_{\widetilde{f}(x)}^2 - \rank \E_{\widetilde{f}(x)}^2
= \dim \W_{\widetilde{f}(y)}^2 - \rank \E_{\widetilde{f}(y)}^2$.

Equation (\ref{dimension relation for inverse cover}) and
submersiveness of $\varphi^f_x$ imply that
\begin{align}
(\phi^f_x)_\ast &: T_{\phi^1_{x, y}(z)} W^1_x / (\phi^1_{x, y})_\ast T_z W^1_{x, y} \notag \\
& \quad
\tocong T_{\phi^2_{x, y}(\phi^f_y(z))} W^2_{\widetilde{f}(x)}
/ (\phi^2_{\widetilde{f}(x), \widetilde{f}(y)})_\ast T_{\phi^f_y(z)}
W^2_{\widetilde{f}(x), \widetilde{f}(y)}
\label{transversality for embeddings}
\end{align}
for any $z \in \W^1_{x, y}$, where
$\phi^f_x$, $\phi^f_y$, $\phi^1_{x, y}$ and $\phi^2_{\widetilde{f}(x), \widetilde{f}(y)}$
are lifts of $\varphi^f_x$, $\varphi^f_y$, $\varphi^1_{x, y}$ and
$\varphi^2_{\widetilde{f}(x), \widetilde{f}(y)}$ respectively for
appropriate orbibundle charts of neighborhoods of $z$ and its images
such that $\phi^f_x \circ \phi^1_{x, y} = \phi^2_{\widetilde{f}(x), \widetilde{f}(y)}
\circ \phi^f_y$.
This implies that
for any point $p \in \psi^1_x((s^1_x)^{-1}(0)) \cap \psi^1_y((s^1_y)^{-1}(0))$,
$\varphi^1_{x, y}(\W^1_{x, y})$ coincides with
$(\varphi^f_x)^{-1}(\varphi^2_{\widetilde{f}(x), \widetilde{f}(y)}
(\W^2_{\widetilde{f}(x), \widetilde{f}(y)}))$ on a neighborhood of
$(\psi^1_x)^{-1}(p) \in \W^1_x$.
Therefore, Condition \ref{clean intersection for meet-semilattice cover} for
$\Y^1 = \widetilde{f}^{-1}(\Y^2)$ follows from that for $\Y^2$.

Next we check Condition \ref{wedge existence for meet-semilattice cover} for $\Y^1$.
For any $p \in X^1$ and any two points $y, z \in \Y^1 \cap (\mu^1)^{-1}(p)$,
Condition \ref{isom of pos for submersion} implies that there exists a unique
$w \in \Y^1 \cap (\mu^1)^{-1}(p)$ such that
$\widetilde{f}(w) = \widetilde{f}(y) \wedge \widetilde{f}(z)$.
Then $w$ satisfies the condition of $y \wedge z$.
Hence $\Y^1$ is a meet-semilattice cover of $X^1$.

For a weakly good coordinate system $(x_2, \V^2_{x_2})_{x_2 \in P^2}$ of $X^2$,
define
$\Y^2 = \bigcup_{x_2 \in P^2} \widetilde{\psi}^2_{x_2}(\V^2_{x_2} \cap (s^2_{x_2})^{-1}(0))$.
Then the above argument implies that $\Y^1 = \widetilde{f}^{-1}(\Y^2)$ is
a meet-semilattice cover of $X^1$.
Therefore Lemma \ref{weakly good coordinate from meet-semilattice} implies that
for any compact subset $K \subset \Y^1$, we can construct a weakly good coordinate
system $(x_1, \V^1_{x_1})_{x_1 \in P^1}$ of $X^1$ which is compatible with
$(x_2, \V^2_{x_2})_{x_2 \in P^2}$ and $K \subset
\bigcup_{x_1 \in P^1} \widetilde{\psi}^1_{x_1}(\V^1_{x_1} \cap
(s^1_{x_1})^{-1}(0)) \subset \Y^1$.
(Condition \ref{V_{x_1} subset V_{x_2}} holds if each $\V_{x_1}^1$ is sufficiently small, and
Condition \ref{composition compatibility for submersion} holds
if each $\V_{x_1}^1$ is a sufficiently small neighborhood of
$(s_{x_1}^1)^{-1}(0) \cap \V_{x_1}^1$.)
Then for a grouped multisection $\boldsymbol{\epsilon}^2$ of
$(x_2, \V^2_{x_2})_{x_2 \in P^2}$,
we can define its pull back $f^\ast \boldsymbol{\epsilon}^2$
as a grouped multisection of $(x_1, \V^1_{x_1})_{x_1 \in P^1}$.

\begin{rem}
In the above argument, we can replace submersion with essential submersion
because the pull back of a grouped multisection by an essential submersion
also satisfies the transversality condition in
Lemma \ref{construction of grouped multisection}.
Note that even for an essential submersion,
(\ref{transversality for embeddings}) is an isomorphism
since
\[
T_{\phi^1_{x, y}(z)} W^1_x / (\phi^1_{x, y})_\ast T_z W^1_{x, y}
\cong T_{\phi^1_{x, y}(z)} \mathring{\partial}^k W^1_x
/ (\phi^1_{x, y})_\ast T_z \mathring{\partial}^k W^1_{x, y}
\]
and
\begin{align*}
&T_{\phi^2_{x, y}(w)} W^2_{\widetilde{f}(x)}
/ (\phi^2_{\widetilde{f}(x), \widetilde{f}(y)})_\ast T_w W^2_{\widetilde{f}(x), \widetilde{f}(y)} \\
&\cong T_{\phi^2_{x, y}(w)} \mathring{\partial}^{k'} W^2_{\widetilde{f}(x)}
/ (\phi^2_{\widetilde{f}(x), \widetilde{f}(y)})_\ast T_w \mathring{\partial}^{k'}
W^2_{\widetilde{f}(x), \widetilde{f}(y)}
\end{align*}
for any $z \in \mathring{\partial}^k \W^1_{x, y}$ and
$w \in \mathring{\partial}^{k'} \W^2_{\widetilde{f}(x), \widetilde{f}(y)}$.
(The case of of an essential submersion between fiber products with simplicial complexes
is also similar.
See the next section for the definition of the fiber product of a pre-Kuranishi space
with a simplicial complex.)
\end{rem}

\subsubsection{Product of pre-Kuranishi spaces}
\label{Product of pre-Kuranishi spaces}
Next we define the product of pre-Kuranishi spaces.
The definition of the product of different pre-Kuranishi spaces,
the product of their weakly good coordinate systems and
the product of their grouped multisections are straightforward.
In application, we need to consider the quotient of the product of the same
pre-Kuranishi spaces by the permutation.
In this case, we need to get rid of the products of Kuranishi neighborhoods
which is not compatible with the permutation action.
\begin{defi}
Let $(\widetilde{X}, \mu, (\W_x, \E_x, s_x, \widetilde{\psi}_x),
(\varphi_{x, y}, \hat \varphi_{x, y}))$
be the pre-Kuranishi structure of a compact Hausdorff space $X$.
We assume that $\dim \W_x > 0$ for each $x \in \widetilde{X}$.
Then for each $N \geq 2$, we say a pre-Kuranishi structure
$(\widetilde{X}^{\# N}, \mu^{\# N},
(\W^{\# N}_x, \E^{\# N}_x, s^{\# N}_x, \widetilde{\psi}^{\# N}_x),
(\varphi^{\# N}_{x, y}, \hat \varphi^{\# N}_{x, y}))$ of
$(\prod^N X) / \mathfrak{S}_N$ is compatible
with that of $X$ if the following conditions hold:
\begin{itemize}
\item
$\widetilde{X}^{\# N}$ is an open subset of $\bigl(\prod^N \widetilde{X}\bigr)
/ \mathfrak{S}_N$ defined by
\[
\widetilde{X}^{\# N} = \{(x_i)_{1 \leq i \leq N} \in \bigl(\prod^N \widetilde{X}\bigr)
/ \mathfrak{S}_N; x_i = x_j \text{ if } \mu(x_i) = \mu(x_j) \}.
\]
\item
$\mu^{\# N}$ is the restriction of the product of $\mu$ to $\widetilde{X}^{\# N}$.
\item
For any two elements $x = (x_i), y = (y_i) \in (\mu^{\# N})^{-1}(p)$ in a fiber,
we reorder the sequences so that $\mu(x_i) = \mu(y_i)$ for all $i$.
Then $y \leq x$ if and only if $y_i \leq x_i$ for all $i$.\item
For each $x = (x_i) \in \widetilde{X}^{\# N}$, the following hold.
\begin{itemize}
\item
$\W^{\# N}_x$ is a neighborhood of $x$ in
$\prod_i \W_{x_i} / (\mathfrak{S}_N)_{(x_i)}$,
where $(\mathfrak{S}_N)_{(x_i)} \subset \mathfrak{S}_N$ is the stabilizer of $(x_i)$.
\item
$\E^{\# N}_x$, $s^{\# N}_x$ and $\widetilde{\psi}^{\# N}_x$ are the restriction of the
product of $\E_{x_i}$, $s_{x_i}$ and $\widetilde{\psi}_{x_i}$ respectively to $\W^{\# N}_x$.
\end{itemize}
\item
For any $x = (x_i), y = (y_i) \in \widetilde{X}^{\# N}$ such that $x \unrhd y$,
$(\varphi^{\# N}_{x, y}, \hat \varphi^{\# N}_{x, y})$
are the restrictions of the products of $(\varphi_{x_i, y_i}, \hat \varphi_{x_i, y_i})$
to a neighborhood $\W_{x, y}$ of $\psi_y^{-1}(\psi_x(s_x^{-1}(0)))$.
\end{itemize}
\end{defi}
We note that in the above definition, the action of $(\mathfrak{S}_N)_{(x_i)}$ on
$\prod_i \W_{x_i}$ is effective because of the assumption $\dim \W_{x_i} > 0$.

\begin{defi}
Let $(x, \V_x)_{x \in P}$ and $((x_i), \V^{\# N}_{(x_i)})_{(x_i) \in P^{\# N}}$ be
weakly good coordinate systems of pre-Kuranishi spaces $X$ and
$\prod^N X / \mathfrak{S}_N$ respectively.
We say these are compatible if the following conditions hold.
\begin{enumerate}[label=(\roman*)]
\item
For any $x = (x_i) \in P^{\# N}$, there exist some $\widehat{x}_i \in P$
such that $x_i \in \widetilde{\psi}_{\widehat{x}_i}(s_{x_i}^{-1}(0) \cap \V_{x_i})$ for all $i$.
\item
Let $x = (x_i) \in P^{\# N}$ and $\widehat{x}_i \in P$ be points such that
$x_i \in \widetilde{\psi}_{\widehat{x}_i}(s_{x_i}^{-1}(0) \cap \V_{x_i})$ and
$\widehat{x}_i = \widehat{x}_j$ if $x_i = x_j$.
For such a pair $((x_i), (\widehat{x}_i))$, we impose the condition
$\V^{\# N}_{(x_i)} \subset \prod_i \varphi_{\widehat{x}_i, x_i}^{-1}(\V_{\widehat{x}_i})
/ (\mathfrak{S}_N)_{(x_i)}$.
\item
Let $x = (x_i), y = (y_i) \in P^{\# N}$ and
$\widehat{x}_i, \widehat{y}_i \in P$ be points such that both of $((x_i), (\widehat{x}_i))$
and $((y_i), (\widehat{y}_i))$ are the pairs in the previous condition.
Assume that there exists some
$p = (p_i) \in \psi^{\# N}_x((s_x^{\# N})^{-1}(0) \cap \V^{\# N}_x)
\cap \psi^{\# N}_y((s_y^{\# N})^{-1}(0) \cap \V^{\# N}_y)$ such that $p_x \geq p_y$.
We reorder $p_i$, $x_i$ and $y_i$ so that
\begin{align*}
(p_i) &\in (\prod_i \psi_{x_i})(\pi_x^{-1}((s_x^{\# N})^{-1}(0) \cap \V^{\# N}_x)) \\
& \quad
\cap (\prod_i \psi_{y_i})(\pi_y^{-1}((s_y^{\# N})^{-1}(0) \cap \V^{\# N}_y))
\end{align*}
and $(p_i)_{x_i} \geq (p_i)_{y_i}$,
where $\pi_x : \prod_i \W_{x_i} \to \prod_i \W_{x_i} / (\mathfrak{S}_N)_{(x_i)}$ and
$\pi_y : \prod_i \W_{y_i} \to \prod_i \W_{y_i} / (\mathfrak{S}_N)_{(y_i)}$ are the projections.
Then
\[
\pi_y^{-1}(\V^{\# N}_{x, y}) \subset \prod_i (\W_{\widehat{x}_i, x_i, y_i} \cap
\W_{\widehat{x}_i, \widehat{y}_i, y_i}).
\]
\end{enumerate}
If $(x, \V_x)_{x \in P}$ and $((x_i), \V^{\# N}_{(x_i)})_{(x_i) \in P^{\# N}}$ are
compatible weakly good coordinate systems, then
for a grouped multisection $\boldsymbol{\epsilon} = (\boldsymbol{\epsilon}_x)$ of
$(x, \V_x)_{x \in P}$, we can define a grouped multisection $\boldsymbol{\epsilon}^{\# N}
= (\boldsymbol{\epsilon}_{(x_i)})$ of $((x_i), \V^{\# N}_{(x_i)})_{(x_i) \in P^{\# N}}$
by the restriction of
$\coprod_i \pi_i^\ast(\boldsymbol{\epsilon}_{\widehat{x}_i})$, where
$\widehat{x}_i \in P$ are the points in the above condition, and
$\pi_i^\ast(\boldsymbol{\epsilon}_{\widehat{x}_i})
= (\pi_i^\ast \epsilon^\omega)_{\omega \in \Omega_{x_i, j}}$ is a family of
sections of $\prod_i E_{x_i} \to \prod_i V_{x_i}$ defined by the pull backs of the
sections $(\epsilon^\omega)_{\omega \in \Omega_{x_i, j}}$
by the projection $\pi_i : \prod_{i'} V_{x_{i'}} \to V_{x_i}$.
As in the case of pull back by submersion,
$(\varphi_{x, y}, \hat \varphi_{x, y})$-relations of
$\boldsymbol{\epsilon} = (\boldsymbol{\epsilon}_x)$ and the above conditions
imply $(\varphi_{(x_i), (y_i)}, \hat \varphi_{(x_i), (y_i)})$-relations of
$\boldsymbol{\epsilon}^{\# N} = (\boldsymbol{\epsilon}_{(x_i)})$.
\end{defi}

Note that for a meet-semilattice cover $\Y \subset \widetilde{X}$ of
a pre-Kuranishi space $X$,
$\Y^{\# N} = \prod^N \Y / \mathfrak{S}_N \cap \widetilde{X}^{\# N}$ is
a meet-semilattice cover of $X^N / \mathfrak{S}_N$.
Indeed, for any two points $(x_i), (y_i) \in \Y^{\# N} \cap (\mu^{\# N})^{-1}(p)$,
if we reorder the sequences so that $\mu(x_i) = \mu(y_i)$ for all $i$, then
$(x_i) \wedge (y_i) = (x_i \wedge y_i)$.
For a weakly good coordinate system $(x, \V_x)_{x \in P}$ of $X$,
$\Y = \bigcup_{x \in P} \widetilde{\psi}_x(\V_x \cap s_x^{-1}(0))$ is
a meet-semilattice cover of $X$ by definition.
Hence Lemma \ref{weakly good coordinate from meet-semilattice} implies that
for any compact subset $K \subset \Y^{\# N}$,
we can construct a weakly good coordinate system
$((x_i), \V^{\# N}_{(x_i)})_{(x_i) \in P^{\# N}}$ of $\prod^N X / \mathfrak{S}_N$
which is compatible with $(x, \V_x)_{x \in P}$
and which satisfies
\[
K \subset \bigcup_{(x_i) \in P^{\# N}}
\widetilde{\psi}^{\# N}_{(x_i)}(\V^{\# N}_{(x_i)} \cap (s^{\# N}_{(x_i)})^{-1}(0))
\subset \Y^{\# N}.
\]
Hence a grouped multisection $\boldsymbol{\epsilon} = (\boldsymbol{\epsilon}_x)$ for
$(x, \V_x)_{x \in P}$ defines a grouped multisection $\boldsymbol{\epsilon}^{\# N}
= (\boldsymbol{\epsilon}_{(x_i)})$ for $((x_i), \V^{\# N}_{(x_i)})_{(x_i) \in P^{\# N}}$
as above.

\begin{rem}
\label{quotient of product pre-Kuranishi by subgroup}
For any subgroup $\Gamma \subset \mathfrak{S}_N$, we can
similarly define the compatibility conditions of
a pre-Kuranishi structure of $(\prod^N X) / \Gamma$ with that of $X$.
\end{rem}
\begin{rem}
\label{quotient map is a submersion}
The quotient map
$\prod^N X \to (\prod^N X) / \mathfrak{S}_N$ is a submersion if
we replace $\widetilde{\prod^N X} = \prod^N \widetilde{X}$ with its open subset
$\{(x_i)_{1 \leq i \leq N} \in \prod^N \widetilde{X}; x_i = x_j \text{ if } \mu(x_i) = \mu(x_j) \}$
(and shrink the Kuranishi neighborhoods of $\prod^N X$).
Similarly, for a subgroup $\Gamma \subset \mathfrak{S}$,
we can make the quotient map
$(\prod^N X) / \Gamma \to (\prod^N X) / \mathfrak{S}_N$
a submersion.
\end{rem}

Next we consider fiber product of pre-Kuranishi spaces.
\begin{defi}
\label{def of fiber product of pre-Kuranishi}
Let $f = (f_x)_{x \in \widetilde{X}}$ be a strong continuous map from
a pre-Kuranishi space $X$ to a smooth manifold $Y$ such that
each $f_x : \W_x \to Y$ is a smooth submersion.
Then for a submanifold $Z \subset Y$, the pre-Kuranishi structure
$(\widetilde{X}', \mu', (\W'_x, \E'_x, s'_x, \widetilde{\psi}'_x),
(\varphi'_{x, y}, \hat \varphi'_{x, y}))$ of $f^{-1}(Z) \subset X$
is defined by $\widetilde{X}' = \mu^{-1}(f^{-1}(Z))$, $\mu' = \mu|_{\widetilde{X}'}$,
$\W'_x = f_x^{-1}(Z)$, $\E'_x = \E_x|_{\W'_x}$, $s'_x = s_x|_{\W'_x}$, $\widetilde{\psi}'
= \widetilde{\psi}|_{(s'_x)^{-1}(0)}$,
$\W'_{x, y} = \W'_y \cap\varphi_{x, y}^{-1}(\W'_x)$ and
$(\varphi'_{x, y}, \hat \varphi'_{x, y})
= (\varphi_{x, y}, \hat \varphi_{x, y})|_{\W'_{x, y}}$.
\end{defi}
Similarly to Lemma \ref{extension of meet-semilattice cover},
for a meet-semilattice cover $\Y^Z$ of $f^{-1}(Z)$, we can construct
a meet-semilattice cover $\Y$ of $X$ such that
$\Y \cap \widetilde{f}^{-1}(Z) = \Y^Z$.
Hence for a weakly good coordinate system $(x, \V^Z_x)_{x \in P^Z}$ of $f^{-1}(Z)$,
we can construct a weakly good coordinate system
$(x, \V_x)_{x \in P}$ of $X$ such that $P \cap \widetilde{f}^{-1}(Z) = P^Z$,
$\V_x \cap f_x^{-1}(Z) = \V^Z_x$ for $x \in P^Z$, and
$\V_x \cap f_x^{-1}(Z) = \emptyset$ for $x \in P \setminus P^Z$.
Furthermore, if a grouped multisection of $(x, \V^Z_x)_{x \in P^Z}$ is given, then
we can extend it to a grouped multisection of $(x, \V_x)_{x \in P}$
if we shrink $\V_x$ slightly.

In the above definition of fiber product, $Z$ is a submanifold of a manifold $Y$.
We also consider the case of a simplicial complex in an orbifold.
\begin{defi}\label{def of fiber product of pre-Kuranishi over orbifold}
Let $K \subset \mathcal{Z}$ be an embedded simplicial complex in a smooth orbifold $\mathcal{Z}$.
We assume that for any point $p \in K$, $\St(p, K)$ is contained in an orbichart
$\mathcal{Z}_p = (Z_p, \pi_{Z_p}, \mathcal{Z}_p)$ of $\mathcal{Z}$.
Define $G_p = \Aut_{\mathcal{Z}_p} Z_p$.
We assume that there exists a regular $G_p$-complex
$L \subset Z_p$ and an isomorphism $\varphi : L/G_p \cong \St(p, K)$ such that
$\varphi \circ \pi_L = \pi_{Z_p}$ on $L \subset Z_p$,
where $\pi_L : L \to L/G_p$ is the quotient map.
Let $f = (f_x)_{x \in \widetilde{X}}$ be a strong continuous map from
a pre-Kuranishi space $X$ to $\mathcal{Z}$ such that each $f_x : \W_x \to \mathcal{Z}$
is a smooth submersion.
We assume that for each point $x \in \widetilde{X}$, $G_{W_x}$ acts effectively on
$\pi_{W_x}^{-1}(f_x^{-1}(\widetilde{f}(x))) \subset W_x$.
Then we can define the pre-Kuranishi structure 
$(\widetilde{X}', \mu', (\W'_x, \E'_x, s'_x, \widetilde{\psi}'_x),
(\varphi'_{x, y}, \hat \varphi'_{x, y}))$ of $f^{-1}(K) \subset X$
similarly as in the case of Definition \ref{def of fiber product of pre-Kuranishi},
whose orbibundle charts are defined as in
Definition \ref{def of fiber product of orbibundle over orbifold}.
\end{defi}

We need to modify various definitions for the fiber product
$f^{-1}(K) \subset X$ with simplicial complex as follows.
A meet-semilattice cover of $f^{-1}(K) \subset X$ is
an open subset $\Y \subset \widetilde{X}$ which satisfies
$\mu(\Y) \supset f^{-1}(K)$ and
the conditions of a meet-semilattice cover of $X$
other than the condition $\mu(\Y) = X$.
A weakly good coordinate system of $f^{-1}(K) \subset X$
is a family of finite pairs $(x, \V_x)_{x \in P}$ of points $x \in \widetilde{f}^{-1}(K)$
and open neighborhoods $\V_x \subset \W_x$ of $\widetilde{\psi}_x^{-1}(x)$
which satisfies the following conditions:
\begin{enumerate}[label=(\roman*)]
\item
\label{P and simplex}
For any $x \in P$ and simplex $s$ of $K$,
if $x \notin f^{-1}(s)$ then $f_x(\V_x) \cap s = \emptyset$.
\item
$\bigcup_{x \in P} \psi_x(\V_x \cap s_x^{-1}(0))$ is a meet-semilattice cover of $f^{-1}(K)$
instead of Condition \ref{P meet-semilattice}.
\item
Condition \ref{(x, y, z)-relation} to \ref{(x, y, z, w)-relation}
in Definition \ref{def of weakly good coordinate system} is satisfied.
\item
In Condition \ref{(x, y, z, w)-relation},
if $x, y, z \in P$ is contained in $\widetilde{f}^{-1}(K_0)$ for a subcomplex
$K_0 \subset K$, then
the condition still holds even if we impose the condition
$w_j \in P \cap \widetilde{f}^{-1}(K_0)$.
\end{enumerate}
The first and last condition imply that $(x, \V_x)_{x \in P \cap \widetilde{f}^{-1}(K_0)}$
is also a weakly good coordinate system of $f^{-1}(K_0) \subset X$
for any subcomplex $K_0 \subset K$.

For a weakly good coordinate system $(x, \V_x)_{x \in P}$ of $f^{-1}(K) \subset X$,
similarly to Lemma \ref{construction of grouped multisection},
shrinking $\V_x$ slightly if necessary, we can construct a grouped multisection of
$(\boldsymbol{\epsilon}_x)_{x \in P}$ which satisfies the
following transversality condition:
For any orbibundle chart $(\V, \E)$ in $(x, \V_x)$ and any
lift $\tilde t \subset V$ of a simplex $t$ of $K$,
the restriction of each branch of the multisection
$s_x|_{\V} + \boldsymbol{\epsilon}_x|_{\V}$
to $V \cap f_x^{-1}(\tilde t)$ is transverse to the zero section of $E$,
and the same holds for each corner of $\V$ instead of $\V$.
(In the case of the generalization of the fiber product with
a Euclidean cell complex, we read $\tilde t$ as its cell.)

\begin{rem}
\label{meaning of triangulation}
Even in the case where $K = \Y$ is a triangulation of $\Y$,
the fiber product $f^{-1}(K)$ is meaningful because
the transversality condition imposed on the grouped multisection of
$f^{-1}(K)$ is stronger than that of $X$.
\end{rem}

We note that for the construction of grouped multisection,
we use the same inductive extension of grouped multisections as in the usual case, and
we do not need to use any induction in the dimension of the simplex of $K$.

As to the extension from a subcomplex $K_0 \subset K$ to $K$,
we can prove the following.
First, if a meet-semilattice cover $\Y^0$ of $f^{-1}(K_0)$ is given, then
we can construct a meet-semilattice cover $\Y$ of $f^{-1}(K)$ such that
$\Y|_{N(f^{-1}(K_0))} = \Y^0|_{f^{-1}(K_0)}$ for some neighborhood $N(f^{-1}(K_0))$
of $f^{-1}(K_0) \subset X$.
We can prove this by the argument of Lemma \ref{extension of meet-semilattice cover}.
(In Lemma \ref{extension of meet-semilattice cover},
first we need to extend the given meet-semilattice cover of $\partial X$ to
its neighborhood, but in this case, $\Y^0$ is already defined on a neighborhood of
$f^{-1}(K_0)$.
Hence if we read $\Y^{N(\partial X)}$ in the proof of
Lemma \ref{extension of meet-semilattice cover} as $\Y^0$,
then the last two paragraph of its proof is the proof of the extension in this case.)

Next, if a weakly good coordinate system $(x, \V^{K_0}_x)_{x \in P^0}$ of $f^{-1}(K_0)$
is given, then we can construct a weakly good coordinate system
$(x, \V_x)_{x \in P}$ of $f^{-1}(K)$ such that
$P^0 = P \cap f^{-1}(K_0)$ and for each $x \in P^0$,
$\V_x \subset \V^{K_0}_x$ is a neighborhood of
$\V^{K_0}_x \cap f_x^{-1}(K_0)$.
Furthermore, if a grouped multisection
$(\boldsymbol{\epsilon}^{K_0}_x)_{x \in P \cap \widetilde{f}^{-1}(K_0)}$
of $(x, \V^{K_0}_x)_{x \in P^0}$ is given,
then shrinking $\V_x$ slightly if necessary,
we can construct a grouped multisection
$(\boldsymbol{\epsilon}^K_x)_{x \in P}$ of $(x, \V_x)_{x \in P}$ such that
the restrictions of $\boldsymbol{\epsilon}^{K_0}_x$ and $\boldsymbol{\epsilon}^K_x$
to a neighborhood of $\V_x \cap f_x^{-1}(K_0)$ coincide for all $x \in P^0$.
This is due to Condition \ref{P and simplex} for $(x, \V_x)_{x \in P}$.

\begin{eg}\label{fiber prod of pre-Kuranishi with diagonal of orbifold}
Let $f = (f_x)_{x \in \widetilde{X}}$ be a strong continuous map
from a pre-Kuranishi space $X$ to an orbifold $\Y$ such that
each $f_x$ is a smooth submersion.
Assume that for each point $x \in \widetilde{X}$, the
dimension of $\pi_{W_x}^{-1}(f_x^{-1}(\widetilde{f}(x))) \subset W_x$ is $>0$
if it is not an empty set, and $G_{W_x}$ acts effectively on it.
Then we can define fiber product
$(f \times f)^{-1}(\Delta_\Y / \mathfrak{S}_2) \subset (X \times X) / \mathfrak{S}_2$
by regarding the diagonal
$\Delta_\Y / \mathfrak{S}_2 \subset (\Y \times \Y) / \mathfrak{S}_2$
as an embedded simplicial complex.
The pre-Kuranishi structure depends on the choice of the triangulation of
$\Delta_\Y / \mathfrak{S}_2$ in the following sense.
If we change the triangulation, then the transversality condition imposed on its
grouped multisection also changes.
\end{eg}

\subsubsection{Multi-valued partial submersions}
\label{multi-valued partial submersions}
First consider the following two trivial examples of the construction of
compatible sections.
\begin{eg}
Let $(V, E)$ be a vector bundle and let $V_i \subset V$ be
finite number of submanifolds which intersect cleanly.
Assume that a smooth section $s_i : V_i \to E|_{V_i}$ is given for each $i$, and
that $s_i$ and $s_j$ coincide on the intersection $V_i \cap V_j$ for all $i, j$.
Then we can construct a smooth section $s : V \to E$ whose restriction to
$V_i$ coincides with $s_i$ for all $i$.
\end{eg}
\begin{eg}
Let $N \geq 1$ be an integer.
We denote a decomposition $\{1, \dots, N\} = \coprod_k A_k$ of integers
by $A = (A_k)$.
We say $g \in \mathfrak{S}_N$ preserves the decomposition $A$
if $g$ maps each $A_k$ to some $A_{k'}$ bijectively.
For each decomposition $A$, let $G_A \subset \mathfrak{S}_N$ be
the group of permutations of $\{1, \dots, N\}$ which preserve $A$.
For two decompositions $A = (A_k)$ and $B = (B_k)$,
we say $A \geq B$ if $B$ is a refinement of $A$, that is,
if each $B_k$ is contained in some $A_{k'}$.
Note that $A \geq B$ does not imply $G_A \subset G_B$ in general.
We denote the discrete decomposition $\{1, \dots, N\} = \coprod_{i=1}^N \{i\}$ by $1^N$.
This is the minimum with respect to this partial order of decompositions.

Let $(V, E)$ be a vector bundle such that $\dim V > 0$,
and define an orbibundle charts $(\V_A, \E_A)$ by
$\V_A = (\prod_N V)/G_A$ and $\E_A = (\prod_N E) / G_A$
for each decomposition $A = (A_k)$.
For each pair $A > B$, let $(\varphi^{B, A} ,\hat \varphi^{B, A})$ be
the multi-valued map from $(\V_A, \E_A)$ to $(\V_B, \E_B)$
induced by the identity map of $(\prod_N V, \prod_N E)$.
Namely, $\varphi^{B, A} \subset \V_A \times \V_B$ is the image of
the diagonal set by the quotient map $(\prod_N V) \times (\prod_N V)
\to \V_A \times \V_B$, and
$\hat \varphi^{B, A} \subset \E_A \times \E_B$ is similar.
Let $A_0$ be an arbitrary decomposition, and assume that
a smooth section $s_A$ of $(\V_A, \E_A)$ is given for each $A < A_0$ and
they satisfy $s_B \circ \varphi^{B, A} = \hat \varphi^{B, A} \circ s_A$ for
all pairs $A > B$ such that $A < A_0$.
Then we can construct a smooth section $s_{A_0}$ of $(\V_{A_0}, \E_{A_0})$
such that $s_A \circ \varphi^{A, A_0} = \hat \varphi^{A, A_0} \circ s_{A_0}$ for
all $A < A_0$.
This is because $\varphi^{1^N, A} \circ \varphi^{A, A_0} = \varphi^{1^N, A_0}$
and $\varphi^{1^N , A_0}$ is single-valued.
More directly, $s_{A_0} = (\varphi^{1^N , A_0})^\ast s_{1^N}$ is the required section.
\end{eg}

In application, we need to consider the combination of the above two examples
in the case of pre-Kuranishi spaces.
Although the case of pre-Kuranishi spaces is also essentially nothing more than
the above two trivial examples,
we explain this case in details since it looks complicated
if we do not give precise definitions.
In this section, we define a compatible system of multi-valued partial submersions
and explain how to construct a compatible family of grouped multisections.

\begin{defi}
\label{smooth multi-valued partial map}
Let $\V = (V, \pi_V, \V)$ and $\V'_i = (V'_i, \pi_{V'_i}, \V'_i)$ be
(at most countable number of) orbicharts.
A smooth multi-valued partial map $\varphi : \V \to \coprod_i \V'_i$
is a closed subset $\varphi \subset \V \times \coprod_i \V'_i$ such that
there exists a set of smooth maps $\{\phi : \D(\phi) \to \coprod_i V'_i\}$
which satisfies the following conditions.
We call each $\phi$ a lift of $\varphi$.
\begin{itemize}
\item
The domain $\D(\phi)$ of each lift $\phi$ is a closed connected submanifold of $V$
(with or without corners), and $\phi : \D(\phi) \to \coprod_i V'_i$ is a smooth map.
\item
Assume that the image of $\phi$ is contained in $V'_i$.
Then for any $g \in G_V$ and $g' \in G_{V'_i}$, $g' \circ \phi \circ g
: \D(g' \circ \phi \circ g) = g^{-1} \D(\phi) \to V'_i$ is also a lift of $\varphi$.
\item
$\{\D(\phi)\}$ is locally finite in $V$, and
they intersect cleanly.
\item
For two lifts $\phi$ and $\phi'$, if there exists a point $p \in \D(\phi) \cap \D(\phi')$
such that $\phi(p) = \phi'(p)$,
$T_p \D(\phi) \subset T_p \D(\phi')$ and
$(\phi_\ast)_p = (\phi'_\ast)_p|_{T_p \D(\phi)}$ on $T_p \D(\phi)$,
then $\phi = \phi'$.
\item
$(\pi_V \times \coprod_i \pi_{V'_i})^{-1}(\varphi) \subset V \times \coprod_i V'_i$
coincides with the union of the graphs of the lifts $\{\phi\}$.
\end{itemize}
\end{defi}
\begin{lem}
In the above definition, the set of lifts
$\{\phi : \D(\phi) \to \coprod_i V'_i\}$ is uniquely determined by
a closed subset $\varphi \subset \V \times \coprod_i \V'_i$ if it exists.
\end{lem}
\begin{proof}
Let $D\varphi \subset TV \times \coprod_i TV'_i$ be the
set of tangent vectors which can be represented by
some smooth curves of $V \times \coprod_i V'_i$ contained in
$(\pi_V \times \coprod_i \pi_{V'_i})^{-1}(\varphi)$.
$D\varphi$ coincides with the union of the tangent bundles of
the graphs of the lifts $\phi$.
For each point $p \in V \times \coprod_i V'_i$,
the fiber $(D\varphi)_p$ at $p$ is the union of some subspaces $(E_{p, j})_j$ of
$(TV \times \coprod_i TV'_i)_p$.
We assume that this decomposition is irredundant, that is,
there does not exist two indices $i \neq j$ such that $E_{p, i} \subset E_{p, j}$.
Each $E_{p, j}$ defines a point of the union of the Grassmann bundles
$G(TV \times \coprod_i TV'_i) = \coprod_k G_k(TV \times \coprod_i TV'_i)$
of $TV \times \coprod_i TV'_i$ over all dimensions.
Let $G(D\varphi)$ be the union of these points.
By assumption, $G(D\varphi)$ coincides with
the disjoint union of the closed submanifolds defined by the graphs of the lifts $\phi$.
Hence the lifts $\phi$ are uniquely determined by the connected components of
the submanifold $G(D\varphi)$ of $G(TV \times \coprod_i TV'_i)$.
\end{proof}
\begin{defi}
Let $(\V, \E)$ and $(\V'_i, \E'_i)$ be orbibundle charts.
A multi-valued partial submersion
$(\varphi, \hat \varphi) : (\V, \E) \to \coprod_i (\V'_i, \E'_i)$
is a pair of
smooth multi-valued partial maps $\varphi : \V \to \coprod_i \V'_i$
and $\hat \varphi : \E \to \coprod_i \E'_i$ such that
there is a one-to-one correspondence between their lifts $\{\phi\}$ and
$\{\hat \phi\}$, and each pair $(\phi, \hat \phi)$
is a bundle map $(\D(\phi), \E|_{\D(\phi)}) \to \coprod_i (V'_i, E'_i)$
which satisfies the following conditions:
The underlying map $\phi : \D(\phi) \to \coprod_i V'_i$ is a submersion,
and the restriction of $\hat \phi$ to each fiber is an isomorphism.
We call each pair $(\phi, \hat \phi)$ a lift of $(\varphi, \hat \varphi)$.
\end{defi}
We do not define the composition of two arbitrary multi-valued partial submersions.
Instead, we define the following compatible system.

\begin{defi}
\label{compatible system of multi-valued partial submersions for charts}
Let $(\V^a, \E^a)_{a \in A}$ be a finite family of orbibundle charts whose
index set $A$ has a partial order.
For each pair $a, b \in A$ such that $a > b$,
let $(\varphi^{b, a}, \hat \varphi^{b, a}) : (\V^a, \E^a) \to (\V^b, \E^b)$ be
a multi-valued partial submersion.
We say $((\V^a, \E^a)_{a \in A}, (\varphi^{b, a}, \hat \varphi^{b, a})_{a > b \in A})$
is a compatible system of multi-valued partial submersions if
the following hold.
\begin{itemize}
\item
For any $a \in A$,
$\coprod_{b < a} (\varphi^{b, a}, \hat \varphi^{b, a}) : (\V^a, \E^a) \to
\coprod_{b < a} (\V^b, \E^b)$ is also a multi-valued partial submersion.
This means that the condition of clean intersection of the
domains of the lifts $\{\D(\phi) \subset V^a\}$
holds not only independently for each $b < a$ but also
for the union over all $b < a$.
\item
For any lifts $(\phi_1^{b_1, a}, \hat \phi_1^{b_1, a})$ and
$(\phi_2^{b_2, a}, \hat \phi_2^{b_2, a})$ of
$(\varphi_1^{b_1, a}, \hat \varphi_1^{b_1, a})$ and
$(\varphi_1^{b_1, a}, \hat \varphi_1^{b_1, a})$ respectively,
if $\D(\phi_1^{b_1, a}) \cap \D(\phi_2^{b_2, a}) \neq \emptyset$,
then there exists some $c \in A$ such that $c \leq b_1, b_2$ and
some lifts $(\phi_1^{c, b_1}, \hat \phi_1^{c, b_1})$,
$(\phi_2^{c, b_2}, \hat \phi_2^{c, b_2})$
and $(\phi^{c, a}, \hat \phi^{c, a})$ of
$(\varphi_1^{c, b_1}, \hat \varphi_1^{c, b_1})$,
$(\varphi_2^{c, b_2}, \hat \varphi_2^{c, b_2})$ and
$(\varphi^{c, a}, \hat \varphi^{c, a})$ respectively
which satisfy the following conditions:
\[
\D(\phi_1^{b_1, a}) \cap \D(\phi_2^{b_2, a})
= \D(\phi_1^{c, b_1} \circ \phi_1^{b_1, a})
= \D(\phi_2^{c, b_2} \circ \phi_2^{b_2, a})
= \D(\phi^{c, a}),
\]
and
\[
(\phi_1^{c, b_1}, \hat \phi_1^{c, b_1}) \circ (\phi_1^{b_1, a} \hat \phi_1^{b_1, a})
= (\phi_2^{c, b_2}, \hat \phi_2^{c, b_2}) \circ (\phi_2^{b_2, a}, \hat \phi_2^{b_2, a})
= (\phi^{c, a}, \hat \phi^{c, a}).
\]
(If $c = b_1$, then we read $(\varphi_1^{c, b_1}, \hat \varphi_1^{c, b_1})$
as the identity map.
Hence its lift is an element of the automorphism group $G_{V^c}$.)
\end{itemize}
\end{defi}


\begin{defi}
Let $(\V, \E)$ and $(\V', \E')$ be orbibundles.
A multi-valued partial submersion $(\varphi, \hat \varphi) : (\V, \E) \to (\V', \E')$
is a pair of closed subsets
$\varphi \subset \V \times \V'$ and
$\hat \varphi \subset \E \times \E'$
which satisfies the following conditions.
\begin{itemize}
\item
For any point $x \in \V$, $\varphi(x) \subset \V'$ consists of finite points,
where $\varphi(x)$ is defined by $\varphi \cap (\{x\} \times \V')
= \{x\} \times \varphi(x)$.
\item
For any point $x \in \V$ such that $\varphi(x) \neq \emptyset$,
let $(\V'_i, \E'_i)$ be a finite family of disjoint orbibundle charts of $(\V', \E')$
which covers $\varphi(x)$.
Then there exists some orbibundle chart $(\V_x, \E_x)$ of a neighborhood of
$x \in \V$ such that $(\V'_i, \E'_i)$ covers $\varphi(\V_x)$ and
$(\varphi \cap (\V_x \times \coprod_i \V'_i),
\hat \varphi \cap (\E_x \times \coprod_i \E'_i))$ is a multi-valued partial submersion
from $(\V_x, \E_x)$ to $\coprod_i (\V'_i, \E'_i)$.
\end{itemize}
\end{defi}
Note that for any compact subset $K \subset \V$,
$\varphi(K) = \{x' \in \V'; (x, x') \in \varphi \text{ for some } x \in K\} \subset \V'$
is compact.
This is because in Definition \ref{smooth multi-valued partial map},
we assumed that $\D(\phi)$ is closed for each lift $\phi$.

\begin{defi}
\label{compatible system of multi-valued partial submersions for orbibundles}
Let $(\V^a, \E^a)_{a \in A}$ be a finite family of orbibundles whose
index set $A$ has a partial order.
For each pair $a, b \in A$ such that $a > b$,
let $(\varphi^{b, a}, \hat \varphi^{b, a}) : (\V^a, \E^a) \to (\V^b, \E^b)$ be
a multi-valued partial submersion.
We say $((\V^a, \E^a)_{a \in A}, (\varphi^{b, a}, \hat \varphi^{b, a})_{a > b \in A})$
is a compatible system of multi-valued partial submersions if
the following holds.
\begin{itemize}
\item
For any $a > b > c \in A$,
$\varphi^{c, b} \circ \varphi^{b, a} \subset \varphi^{c, a}$
and $\hat \varphi^{c, b} \circ \hat \varphi^{b, a} \subset \hat \varphi^{c, a}$.
\item
For each point $x_0 \in \V^{a_0}$, we define a partial order of
$\bigcup_{a \leq a_0} \varphi^{a, a_0}(x_0)$ as follows.
For $x \in \varphi^{a, a_0}(x_0)$ and $y \in \varphi^{b, a_0}(x_0)$,
$x > y$ if $a > b$ and $y \in \varphi^{b, a}(x)$.
Then for each $x \in \varphi^{a, a_0}(x_0)$,
there exists an orbibundle chart $(\V^x, \E^x)$ of a neighborhood of $x \in \V^a$
such that the restrictions of $(\varphi^{b, a}, \hat \varphi^{b, a})$
define multi-valued partial submersions $(\varphi^{y, x}, \hat \varphi^{y, x})
: (\V^x, \E^x) \to (\V^y, \E^y)$ for all pairs
$x > y \in \bigcup_{a \leq a_0} \varphi^{a, a_0}(x_0)$, and
\[
((\V^x, \E^x)_{x \in \bigcup_{a \leq a_0} \varphi^{a, a_0}(x_0)},
(\varphi^{y, x}, \hat \varphi^{y, x})_{x > y \in \bigcup_{a \leq a_0} \varphi^{a, a_0}(x_0)})
\]
is a compatible system of multi-valued partial submersions in the sense of
Definition \ref{compatible system of multi-valued partial submersions for charts}.
($\V^x$ may depend on $x_0$.)
\end{itemize}
\end{defi}

We also define $(\varphi, \hat \varphi)$-relation of grouped multisections for
a multi-valued partial submersion $(\varphi, \hat \varphi)$
between orbibundles.
The definition is the same as
Definition \ref{(varphi, hat varphi)-relation for orbibundles}.
For a compatible system of multi-valued partial submersions
$((\V^a, \E^a)_{a \in A}, (\varphi^{b, a}, \hat \varphi^{b, a})_{a > b \in A})$,
we also always assume the compatibility condition similar to
Definition \ref{compatibility of (varphi, hat varphi)-relations}.
The proof of the following lemma is straightforward.
\begin{lem}
\label{extension of grouped multisection for multi-valued partial submersions}
Let $((\V^a, \E^a)_{a \in A}, (\varphi^{b, a}, \hat \varphi^{b, a})_{a > b \in A})$
be a compatible system of multi-valued partial submersions.
Let $a_0 \in A$ be an arbitrary index and assume that
grouped multisections $\boldsymbol{\epsilon}_a$ of $(\V^a, \E^a)$
are given for $a < a_0$.
We also assume that $\boldsymbol{\epsilon}_a$ and $\boldsymbol{\epsilon}_b$
are $(\varphi^{b, a}, \hat \varphi^{b, a})$-related for $a > b \in \A$ such that
$a < a_0$.
Then we can construct a grouped multisection $\boldsymbol{\epsilon}_{a_0}$
of $(\V^{a_0}, \E^{a_0})$ which is $(\varphi^{a, a_0}, \hat \varphi^{a, a_0})$-related
to $\boldsymbol{\epsilon}_a$ for all $a < a_0$.
\end{lem}

\begin{defi}
Let $X^k$ ($k = 1,2$) be two pre-Kuranishi spaces with pre-Kuranishi structures
$(\widetilde{X}^k, \mu^k, (\W^k_x, \E^k_x, s^k_x, \widetilde{\psi}^k_x),
(\varphi^k_{x, y}, \hat \varphi^k_{x, y}))$.
A multi-valued partial submersion
$f = (f, \widetilde{f}, (\W^f_{\widehat{x}| x}, \varphi^f_{\widehat{x}, x},
\hat \varphi^f_{\widehat{x}, x}))$
from $X^1$ to $X^2$ consists of the following.
$f \subset X^1 \times X^2$ and
$\widetilde{f} \subset \widetilde{X}^1 \times \widetilde{X}^2$ are closed subsets
such that $(\mu^1 \times \mu^2)(\widetilde{f}) = f$.
For each $(x, \widehat{x}) \in \widetilde{f}$,
$\W^f_{\widehat{x}| x} \Subset \W^2_{\widehat{x}}$ is an open neighborhood of
$(\widetilde{\psi}^2_{\widehat{x}})^{-1}(\widehat{x})$, and
$(\varphi^f_{\widehat{x}, x}, \hat \varphi^f_{\widehat{x}, x})
: (\W^1_{x}, \E^1_{x}) \to (\W^f_{\widehat{x}| x}, \E^2_{\widehat{x}}|_{\W^f_{\widehat{x}| x}})$
is a multi-valued partial submersion.
We impose the following conditions on them:
\begin{enumerate}[label=$(\arabic*)^{\mathrm{MP}}$]
\item
For any $(p, q) \in f$ and $\widetilde{p} \in (\mu^1)^{-1}(p)$,
there exists a unique $\widetilde{q} \in (\mu^2)^{-1}(q)$
such that $(\widetilde{p}, \widetilde{q}) \in \widetilde{f}$.
Furthermore, this defines an isomorphism
\[
\widetilde{f}|_{(\mu^1)^{-1}(p) \times (\mu^2)^{-1}(q)} : (\mu^1)^{-1}(p) \cong
(\mu^2)^{-1}(q)
\]
between partially ordered sets.
\item
\label{properness of tilde f}
For any compact subset $K \subset \widetilde{X}^2$,
$\widetilde{f} \cap (\widetilde{X}^1 \times K)$ is compact.
\item
For any $x \in \widetilde{X}^1$, closed subsets
$\{\psi^2_{\widehat{x}}((s^2_{\widehat{x}})^{-1}(0) \cap \overline{\W^f_{\widehat{x}| x}})
\subset X^2; \widehat{x} \in \widetilde{f}(x)\}$ are disjoint.
\item
\label{varphi condition for multi-valued partial submersion}
For any $(x, \widehat{x}) \in \widetilde{f}$,
$s^2_{\widehat{x}} \circ \varphi^f_{\widehat{x}, x}
= \hat \varphi^f_{\widehat{x}, x} \circ s^1_{x}$ as a subset
of $\W^1_x \times \E^2_{\widehat{x}}|_{\W^f_{\widehat{x}| x}}$.
Furthermore, for any $x \in \widetilde{X}^1$,
$\coprod_{\widehat{x} \in \widetilde{f}(x)}
(\widetilde{\psi}^2_{\widehat{x}} \circ \varphi^f_{\widehat{x}, x})
= \widetilde{f} \circ \widetilde{\psi}^1_{x}$
as a subset of $(s^1_{x})^{-1}(0) \times \widetilde{X}^2$.
\item
For any $x \in \widetilde{X}^1$ and $(y, \widehat{y}) \in \widetilde{f}$
such that $x \unrhd y$, define
$\widetilde{f}(x)_{(y, \widehat{y})} =
\{\widehat{x} \in \widetilde{f}(x); q_{\widehat{x}} \geq q_{\widehat{y}}
\text{ for some } q \in
\psi_{\widehat{x}}((s^2_{\widehat{x}})^{-1}(0) \cap \W^f_{\widehat{x}| x}) \cap
\psi_{\widehat{y}}((s^2_{\widehat{y}})^{-1}(0) \cap \W^f_{\widehat{y}| y})\}$.
Then the following conditions hold:
\begin{itemize}
\item
$\varphi^f_{\widehat{y}, y}(\W^1_{x, y})$ is contained in
$\bigcup_{\widehat{x} \in \widetilde{f}(x)_{(y, \widehat{y})}}
(\varphi^2_{\widehat{x}, \widehat{y}})^{-1}(\W^f_{\widehat{x}| x})$.
Furthermore,
$\{\varphi^f_{\widehat{y}, y}(\W^1_{x, y}) \cap
(\varphi^2_{\widehat{x}, \widehat{y}})^{-1}(\W^f_{\widehat{x}| x});
\widehat{x} \in \widetilde{f}(x)_{(y, \widehat{y})}\}$ are disjoint.
\item
For any $\widehat{x} \in \widetilde{f}(x)_{(y, \widehat{y})}$,
\begin{equation}
(\varphi^f_{\widehat{x}, x}, \hat \varphi^f_{\widehat{x}, x}) \circ
(\varphi^1_{x, y}, \hat \varphi^1_{x, y})
= (\varphi^2_{\widehat{x}, \widehat{y}}, \hat \varphi^2_{\widehat{x}, \widehat{y}}) \circ
(\varphi^f_{\widehat{y}, y}, \hat \varphi^f_{\widehat{y}, y})
\end{equation}
as pairs of smooth multi-valued partial maps
from $(\W^1_{x, y}, \E^1_y|_{\W^1_{x, y}})$ to
$(\W^f_{\widehat{x}| x}, \E^2_{\widehat{x}}|_{\W^f_{\widehat{x}| x}})$.
\end{itemize}
\end{enumerate}
\end{defi}
In relation to Condition \ref{properness of tilde f},
we note that for a compact subset $K \subset \widetilde{X}^1$,
$\widetilde{f} \cap (K \times \widetilde{X}^2)$ is compact.
This is because each $\varphi^f_{\widehat{x}, x}$ maps compact sets to
compact sets as we noted above.
It implies that for any closed subset $A \subset \widetilde{X}^2$,
$\widetilde{f}^{-1}(A) = \{x \in \widetilde{X}^1; (x, \widehat{x}) \in \widetilde{f}
\text{ for some } \widehat{x} \in A\} \subset \widetilde{X}^1$ is also closed.
In particular, $\widetilde{f}^{-1}(\widetilde{X}^2) \subset \widetilde{X}^1$
is closed.
\begin{defi}
Let $(x_1, \V^1_{x_1})_{x_1 \in P^1}$ and $(x_2, \V^2_{x_2})_{x_2 \in P^2}$ be
weakly good coordinate systems of pre-Kuranishi spaces $X^1$ and $X^2$ respectively.
We say these are compatible with respect to the multi-valued partial submersion $f$ if
they satisfy the following conditions:
\begin{enumerate}[label=$(\arabic*)^{\mathrm{CW}}$]
\item
\label{compatibility of meet-semilattice covers}
The meet-semilattice covers $\Y^i = \bigcup_{x_i \in P^i}
\widetilde{\psi}^i_{x_i}((s^i_{x_i})^{-1}(0) \cap \V^i_{x_i}))$ ($i = 1,2$) satisfy
$\widetilde{f}(\Y^1) \subset \Y^2$.
($\widetilde{f}(\Y^1) \subset \widetilde{X}^2$ is defined by
$\widetilde{f}(\Y^1) = \{\widehat{x} \in \widetilde{X}^2; (x, \widehat{x}) \in \widetilde{f}
\text{ for some } x \in \Y^1\}$.)
\item
For any $x_1 \in P^1$, if
$\widetilde{\psi}^1_{x_1}((s^1_{x_1})^{-1}(0) \cap \V^1_{x_1}) \cap
\widetilde{f}^{-1}(\widetilde{X}^2) \neq \emptyset$,
then $x_1 \in \widetilde{f}^{-1}(\Y^2)$ and
$\widetilde{f}(x_1) \subset \Y^2$.
Furthermore, for any $y_1, x_1 \in P^1$ such that
$y_1 \notin \widetilde{f}^{-1}(\Y^2)$ and $x_1 \in \widetilde{f}^{-1}(\Y^2)$,
if there exists some $p \in \psi^1_{x_1}((s^1_{x_1})^{-1}(0) \cap \V^1_{x_1})
\cap \psi^1_{y_1}((s^1_{y_1})^{-1}(0) \cap \V^1_{y_1})$ such that
$p_{y_1} \leq p_{x_1}$, then $\varphi^1_{x_1, y_1}(\V^1_{x_1, y_1}) \subset \V^1_{x_1}$
does not intersect with the domain of $\varphi^f_{\widehat{x}_1, x_1}$
for any $\widehat{x}_1 \in \widetilde{f}(x_1)$.
\item
\label{(x1, widehat x1, x2)-composition}
For any $x_1 \in P^1$, $\widehat{x}_1 \in \widetilde{X}^2$ and $x_2 \in P^2$,
if $(x_1, \widehat{x}_1) \in \widetilde{f}$ and
$\widehat{x}_1 \in \widetilde{\psi}^2_{x_2}((s^2_{x_2})^{-1}(0) \cap \V^2_{x_2}))$,
then $\varphi^f_{\widehat{x}_1, x_1}(\V^1_{x_1}) \subset
(\varphi^2_{x_2, \widehat{x}_1})^{-1}(\V^2_{x_2})$.
\item
\label{composition compatibility for multi-valued partial submersion}
For any $y_1, x_1 \in P^1$ and $y_2, x_2 \in P^2$,
if there exist some $p \in \psi_{x_1}((s^1_{x_1})^{-1}(0) \cap \V^1_{x_1}) \cap
\psi_{y_1}((s^1_{y_1})^{-1}(0) \cap \V^1_{y_1})$ such that
$p_{x_1} \geq p_{y_1}$, then
the following holds
for any $\widehat{y}_1 \in \widetilde{f}(y_1) \cap
\widetilde{\psi}^2_{y_2}((s^2_{y_2})^{-1}(0) \cap \V^2_{y_2})$ and
$\widehat{x}_1 \in \widetilde{f}(x_1)_{(y, \widehat{y})} \cap
\widetilde{\psi}^2_{x_2}((s^2_{x_2})^{-1}(0) \cap \V^2_{x_2})$.
\[
\varphi^f_{\widehat{y}_1, y_1}(\V^1_{x_1, y_1}) \cap 
(\varphi^2_{\widehat{x}_1, \widehat{y}_1})^{-1}(\W^f_{\widehat{x}_1| x_1})
\subset
\W^2_{x_2, \widehat{x}_1, \widehat{y}_1} \cap
\W^2_{x_2, y_2, \widehat{y}_1}
\]
\end{enumerate}
\end{defi}

\begin{defi}
Let $(x_1, \V^1_{x_1})_{x_1 \in P^1}$ and $(x_2, \V^2_{x_2})_{x_2 \in P^2}$ be
weakly good coordinate systems of pre-Kuranishi spaces $X^1$ and $X^2$ respectively
which are compatible with respect to the multi-valued partial submersion $f$.
Let $(\boldsymbol{\epsilon}^1_{x_1})_{x_1 \in P^1}$ and
$(\boldsymbol{\epsilon}^2_{x_2})_{x_2 \in P^2}$ be grouped multisections of
$(x_1, \V^1_{x_1})_{x_1 \in P^1}$ and $(x_2, \V^2_{x_2})_{x_2 \in P^2}$ respectively.
We say these grouped multisections are compatible with respect to $f$ if
for any $x_1 \in P^1$, $\widehat{x}_1 \in \widetilde{X}^2$ and $x_2 \in P^2$
such that $(x_1, \widehat{x}_1) \in \widetilde{f}$ and
$\widehat{x}_1 \in \widetilde{\psi}^2_{x_2}((s^2_{x_2})^{-1}(0) \cap \V^2_{x_2}))$,
$\boldsymbol{\epsilon}_{x_1}$ and $\boldsymbol{\epsilon}_{x_2}$ are
$(\varphi^2_{x_2, \widehat{x}_1} \circ \varphi^f_{\widehat{x}_1, x_1},
\hat \varphi^2_{x_2, \widehat{x}_1} \circ \hat \varphi^f_{\widehat{x}_1, x_1})$-related.
(Note that the triple $(x_1, \widehat{x}_1, x_2)$ is that
in Condition \ref{(x1, widehat x1, x2)-composition}.
Hence $(\varphi^2_{x_2, \widehat{x}_1} \circ \varphi^f_{\widehat{x}_1, x_1},
\hat \varphi^2_{x_2, \widehat{x}_1} \circ \hat \varphi^f_{\widehat{x}_1, x_1})$ is
a multi-valued partial submersion from $(\V^1_{x_1}, \E^1_{x_1}|_{\V^1_{x_1}})$ to
$(\V^2_{x_2}, \E^2_{x_2}|_{\V^2_{x_2}})$.)
\end{defi}

For a single multi-valued partial submersion $f$, in general,
we cannot construct a weakly good coordinate system
$(x_1, \V^1_{x_1})_{x_1 \in P^1}$ of $X^1$ which is compatible with
a given weakly good coordinate system $(x_2, \V^2_{x_2})_{x_2 \in P^2}$ of $X^2$.
Similarly, even if $(x_1, \V^1_{x_1})_{x_1 \in P^1}$ and $(x_2, \V^2_{x_2})_{x_2 \in P^2}$
are compatible weakly good coordinate systems of $X^1$ and $X^2$ respectively,
we cannot construct a grouped multisection
$(\boldsymbol{\epsilon}^1_{x_1})_{x_1 \in P^1}$ of $(x_1, \V^1_{x_1})_{x_1 \in P^1}$
which is compatible with a given grouped multisection
$(\boldsymbol{\epsilon}^2_{x_2})_{x_2 \in P^2}$ of $(x_2, \V^2_{x_2})_{x_2 \in P^2}$.
What we can say in general is the following.
\begin{lem}
\label{compatible pull backs for multi-valued partial submersion}
Assume that $(x_1, \V^1_{x_1})_{x_1 \in P^1}$ and $(x_2, \V^2_{x_2})_{x_2 \in P^2}$
are compatible weakly good coordinate systems of $X^1$ and $X^2$, and
let $(\boldsymbol{\epsilon}^2_{x_2})_{x_2 \in P^2}$ be
a grouped multisection of $(x_2, \V^2_{x_2})_{x_2 \in P^2}$.
For any $x_1 \in P^1$, $\widehat{x}_1 \in \widetilde{X}^2$
such that $(x_1, \widehat{x}_1) \in \widetilde{f}$,
let $P^2(\widehat{x}_1) \subset P^2$ be the set of points $x_2 \in P^2$ such that
$\widehat{x}_1 \in \widetilde{\psi}^2_{x_2}((s^2_{x_2})^{-1}(0) \cap \V^2_{x_2}))$.
For each $x_2 \in P^2(\widehat{x}_1)$, let
$(\varphi^2_{x_2, \widehat{x}_1})^\ast \boldsymbol{\epsilon}^2_{x_2}$ be the pull back of
$\boldsymbol{\epsilon}^2_{x_2}$ by
\begin{align*}
&(\varphi^2_{x_2, \widehat{x}_1}, \hat \varphi^2_{x_2, \widehat{x}_1})
|_{(\varphi^2_{x_2, \widehat{x}_1})^{-1}(\V^2_{x_2})} \\
&: ((\varphi^2_{x_2, \widehat{x}_1})^{-1}(\V^2_{x_2}),
\E^2_{\widehat{x}_1}|_{(\varphi^2_{x_2, \widehat{x}_1})^{-1}(\V^2_{x_2})})
\to (\V^2_{x_2}, \E^2_{x_2}|_{\V^2_{x_2}}).
\end{align*}
Then its restriction to a neighborhood of $\varphi^f_{\widehat{x}_1, x_1}(\V^1_{x_1})$
does not depend on $x_2 \in P^2(\widehat{x}_1)$.
More precisely, their restrictions to
\[
\bigcap_{x_2, y_2 \in P^2(\widehat{x}_1)} \W^2_{x_2, y_2, \widehat{x}_1}
\cap \bigcap_{x_2 \in P^2(\widehat{x}_1)} (\varphi^2_{x_2, \widehat{x}_1})^{-1}(\V^2_{x_2})
\]
coincides.
\end{lem}
\begin{proof}
First note that
Condition \ref{composition compatibility for multi-valued partial submersion}
for $y_1 = x_1$, $\widehat{y}_1 = \widehat{x}_1$ and
$y_2, x_2 \in P^2(\widehat{x}_1)$ implies that
$\varphi^f_{\widehat{x}_1, x_1}(\V^1_{x_1}) \subset
\bigcap_{x_2, y_2 \in P^2(\widehat{x}_1)} \W^2_{x_2, y_2, \widehat{x}_1}$.
Hence the claim follows from $(\varphi^2_{x_2, y_2}, \hat \varphi^2_{x_2, y_2})$-relation
of $(\boldsymbol{\epsilon}^2_{x_2})_{x_2 \in P^2}$.
\end{proof}

\begin{defi}
Let $(X^\alpha)_{\alpha \in \A}$ be finite number of pre-Kuranishi spaces,
and assume that its index set $\A$ has a partial order.
For each pair $\alpha, \beta \in \A$ such that $\alpha > \beta$,
let $f^{\beta, \alpha} = (f^{\beta, \alpha}, \widetilde{f}^{\beta, \alpha},
(\varphi^{f^{\beta, \alpha}}_{x_2, x_1}, \hat \varphi^{f^{\beta, \alpha}}_{x_2, x_1}))$
be a multi-valued partial submersion from $X^\alpha$ to $X^\beta$.
We say $((X^\alpha)_{\alpha \in \A}, (f^{\beta, \alpha})_{\alpha > \beta})$
is a compatible system of
multi-valued partial submersions if it satisfies the following conditions:
\begin{itemize}
\item
For any triple $\alpha > \beta > \gamma \in \A$,
$f^{\gamma, \beta} \circ f^{\beta, \alpha} \subset f^{\gamma, \alpha}$
as a subset of $X^\alpha \times X^\gamma$
and $\widetilde{f}^{\gamma, \beta} \circ \widetilde{f}^{\beta, \alpha} \subset
\widetilde{f}^{\gamma, \alpha}$ as a subset of $\widetilde{X}^\alpha \times
\widetilde{X}^\gamma$.
\item
For each $x_0 \in \widetilde{X}^{\alpha_0}$,
we define a partial order of
$\bigcup_{\alpha < \alpha_0} \widetilde{f}^{\alpha, \alpha_0}(x_0)$ as follows.
For $x \in \widetilde{f}^{\alpha, \alpha_0}(x_0)$ and $y \in f^{\beta, \alpha_0}(x_0)$,
$x > y$ if $\alpha > \beta$ and $y \in \widetilde{f}^{\beta, \alpha}(x)$.
Then
\[
(\varphi^{f^{\beta, \alpha}}_{y, x}, \hat \varphi^{f^{\beta, \alpha}}_{y, x})
: (\W^\alpha_x, \E^\alpha_x) \to
(\W^{f^{\beta, \alpha}}_{y|x}, \E^\beta|_{\W^{f^{\beta, \alpha}}_{y|x}})
\inj (\W^\beta_y, \E^\beta_y)
\]
for $x \in \widetilde{f}^{\alpha, \alpha_0}(x_0)$ and
$y \in \widetilde{f}^{\beta, \alpha_0}(x_0)$
such that $x > y$ in $\bigcup_{\alpha < \alpha_0} \widetilde{f}^{\alpha, \alpha_0}(x_0)$
constitute a compatible system of multi-valued partial submersions
in the sense of
Definition \ref{compatible system of multi-valued partial submersions for orbibundles}.
%
\end{itemize}
\end{defi}
Note that the above condition implies that
$\bigcup_{\alpha < \alpha_0} \widetilde{f}^{\alpha, \alpha_0}(x_0)$ with the above
partial order has a unique minimum for each $x_0 \in \widetilde{X}^{\alpha_0}$.

\begin{defi}
Let $((X^\alpha)_{\alpha \in \A}, (f^{\beta, \alpha})_{\alpha > \beta})$
be a compatible system of multi-valued partial submersions,
and let $(x_\alpha, \V^\alpha_{x_\alpha})_{x_\alpha \in P^\alpha}$ be a weakly good
coordinate system of $X^\alpha$ for each $\alpha \in \A$.
We say these weakly good coordinate systems
$(x_\alpha, \V^\alpha_{x_\alpha})_{x_\alpha \in P^\alpha}$ are compatible if
for any $\alpha > \beta \in \A$,
$(x_\alpha, \V^\alpha_{x_\alpha})_{x_\alpha \in P^\alpha}$ and
$(x_\beta, \V^\beta_{x_\beta})_{x_\beta \in P^\beta}$ are
compatible with respect to $f^{\beta, \alpha}$.
\end{defi}

For a compatible system of multi-valued partial submersions
$((X^\alpha)_{\alpha \in \A}, \ab (f^{\beta, \alpha})_{\alpha > \beta})$,
we can construct a compatible family of weakly good coordinate systems
$(x_\alpha, \V^\alpha_{x_\alpha})_{x_\alpha \in P^\alpha}$ of $X^\alpha$
for all $\alpha \in \A$ as follows.
First we claim that we can construct meet-semilattice covers
$\Y^\alpha \subset \widetilde{X}^\alpha$ of $X^\alpha$ for all $\alpha \in \A$
such that
\begin{equation}
\Y^\alpha \cap (\widetilde{f}^{\beta, \alpha})^{-1}(\widetilde{X}^\beta)
= (\widetilde{f}^{\beta, \alpha})^{-1}(\Y^\beta)
\label{compatible meet-semilattice covers for multi-valued partial submersions}
\end{equation}
for all pairs $\alpha, \beta \in \A$ such that $\alpha > \beta$.

First we note that if
(\ref{compatible meet-semilattice covers for multi-valued partial submersions})
is satisfied for pairs $\alpha > \beta \in \A$ such that $\alpha < \alpha_0$, then
$\widetilde{f}^{\alpha, \alpha_0}(x_0) \subset \Y^\alpha$
for any $x_0 \in (\widetilde{f}^{\alpha, \alpha_0})^{-1}(\Y^\alpha)$.
This is proved as follows.
$x_0 \in (\widetilde{f}^{\alpha, \alpha_0})^{-1}(\Y^\alpha)$ implies that
$(x_0, x) \in \widetilde{f}^{\alpha, \alpha_0}$ for some $x \in \Y^\alpha$.
Let $x_{\min} \in \widetilde{f}^{\alpha_1, \alpha_0}(x_0)$ be the unique minimum of
$\bigcup_{\alpha < \alpha_0} \widetilde{f}^{\alpha, \alpha_0}(x_0)$.
Then $\{x_{\min}\} = \widetilde{f}^{\alpha_1, \alpha}(x)$,
which implies
\[
x \in \Y^\alpha \cap (\widetilde{f}^{\alpha_1, \alpha})^{-1}(\widetilde{X}^{\alpha_1})
= (\widetilde{f}^{\alpha_1, \alpha})^{-1}(\Y^{\alpha_1}).
\]
Since $\widetilde{f}^{\alpha_1, \alpha}(x) = \{x_{\min}\}$ consists of one point,
this implies that $x_{\min} \in \Y^{\alpha_1}$.
Since any point $x' \in \widetilde{f}^{\alpha, \alpha_0}(x_0)$ satisfies
$(x', x_{\min}) \in \widetilde{f}^{\alpha_1, \alpha}$,
this implies that
$x' \in (\widetilde{f}^{\alpha_1, \alpha})^{-1}(\Y^{\alpha_1}) \subset \Y^\alpha$.
Hence $\widetilde{f}^{\alpha, \alpha_0}(x_0) \subset \Y^\alpha$
for any $x_0 \in (\widetilde{f}^{\alpha, \alpha_0})^{-1}(\Y^\alpha)$.

We also note that this implies that
$(\widetilde{f}^{\alpha, \alpha_0})^{-1}(\Y^\alpha)$ is open
in the relative topology of $(\widetilde{f}^{\alpha, \alpha_0})^{-1}(\widetilde{X}^\alpha)$.

We can construct such meet-semilattice covers
by the same argument as the proof of Lemma \ref{extension of meet-semilattice cover}.
Namely, we can extend
$\bigcup_{\beta < \alpha} (\widetilde{f}^{\beta, \alpha})^{-1}(\Y^\beta)
\subset \widetilde{X}^\alpha$ to the cover of $X^\alpha$
by the induction in $\alpha \in \A$.

Similarly, we can construct decreasing sequences of meet-semilattice covers
$\Y^\alpha_k \subset \widetilde{X}^\alpha$ ($\alpha \in \A$, $k \geq 1$) such that
$\Y^\alpha_{k+1} \Subset \Y^\alpha_k$ and
\[
\Y^\alpha_k \cap (\widetilde{f}^{\beta, \alpha})^{-1}(\widetilde{X}^\beta)
= (\widetilde{f}^{\beta, \alpha})^{-1}(\Y^\beta_k)
\]
for $k \geq 1$ and all pairs $\alpha, \beta \in \A$ such that $\alpha > \beta$.
(This is because Condition \ref{properness of tilde f} implies that
$(\widetilde{f}^{\beta, \alpha})^{-1}(\Y^\beta_{k+1}) \Subset
(\widetilde{f}^{\beta, \alpha})^{-1}(\Y^\beta_k)$.)

We attach an integer $k_\alpha \geq 1$ for each $\alpha \in \A$
so that $k_\alpha > k_\beta$ if $\alpha > \beta$.
Then by the induction in $\alpha \in \A$, we can
construct a compatible family of weakly good coordinate systems
$(x_\alpha, \V^\alpha_{x_\alpha})_{x_\alpha \in P^\alpha}$ of $X^\alpha$
for $\alpha \in \A$ such that
$\Y^\alpha_{k_\alpha + 1} \Subset \bigcup_{x_\alpha \in P^\alpha}
\widetilde{\psi}_{x_\alpha}((s^\alpha_{x_\alpha})^{-1}(0) \cap \V^\alpha_{x_\alpha})
\subset \Y^\alpha_{k_\alpha}$.
(We can apply the argument of
Lemma \ref{weakly good coordinate from meet-semilattice}.
Only non-trivial condition is
Condition \ref{compatibility of meet-semilattice covers},
but it is satisfied if
we construct $(x_\alpha, \V^\alpha_{x_\alpha})_{x_\alpha \in P^\alpha}$ inductively
so that they satisfy
$\Y^\alpha_{k_\alpha + 1} \Subset \bigcup_{x_\alpha \in P^\alpha}
\widetilde{\psi}_{x_\alpha}((s^\alpha_{x_\alpha})^{-1}(0) \cap \V^\alpha_{x_\alpha})
\subset \Y^\alpha_{k_\alpha}$ for all $\alpha \in \A$.
The other conditions hold if each $\V_{x_\alpha}^\alpha$ is sufficiently small
and each $\V_{x_\alpha}^\alpha$ is a sufficiently small neighborhood of
$(s_{x_\alpha}^\alpha)^{-1}(0) \cap \V_{x_\alpha}^\alpha$.)

\begin{defi}
Let $((X^\alpha)_{\alpha \in \A}, (f^{\beta, \alpha})_{\alpha > \beta})$
be a compatible system of multi-valued partial submersions,
and let $(x_\alpha, \V^\alpha_{x_\alpha})_{x_\alpha \in P^\alpha}$ be a
compatible family of weakly good coordinate systems of $X^\alpha$
for $\alpha \in \A$.
For each $\alpha \in \A$,
let $(\boldsymbol{\epsilon}^\alpha_{x_\alpha})_{x_\alpha \in P^\alpha}$ be
a grouped multisection of $(x_\alpha, \V^\alpha_{x_\alpha})_{x_\alpha \in P^\alpha}$.
We say these grouped multisections are compatible if
for any $\alpha > \beta \in \A$,
$(\boldsymbol{\epsilon}^\alpha_{x_\alpha})_{x_\alpha \in P^\alpha}$ and
$(\boldsymbol{\epsilon}^\beta_{x_\beta})_{x_\beta \in P^\beta}$ are
compatible with respect to $f^{\beta, \alpha}$.
\end{defi}

\begin{prop}
Let $((X^\alpha)_{\alpha \in \A}, (f^{\beta, \alpha})_{\alpha > \beta})$
be a compatible system of multi-valued partial submersions,
and let $(x_\alpha, \V^\alpha_{x_\alpha})_{x_\alpha \in P^\alpha}$ be a
compatible family of weakly good coordinate systems of $X^\alpha$
for $\alpha \in \A$.
Let $\alpha_0 \in \A$ be an arbitrary index and
assume that a compatible family of grouped multisections
$(\boldsymbol{\epsilon}^\alpha_{x_\alpha})_{x_\alpha \in P^\alpha}$
for $\alpha < \alpha_0$ are given.
Then we can construct a grouped multisection
$(\boldsymbol{\epsilon}^{\alpha_0}_{x_{\alpha_0}})_{x_{\alpha_0} \in P^{\alpha_0}}$
of $(x_{\alpha_0}, \V^{\alpha_0}_{x_{\alpha_0}})_{x_{\alpha_0} \in P^{\alpha_0}}$
which is compatible with
$(\boldsymbol{\epsilon}^\alpha_{x_\alpha})_{x_\alpha \in P^\alpha}$ $(\alpha < \alpha_0)$.
\end{prop}
\begin{proof}
We saw in Lemma \ref{compatible pull backs for multi-valued partial submersion}
that for any $x_{\alpha_0} \in P^{\alpha_0}$ and $\widehat{x}_{\alpha_0}^{(\alpha, k)}
\in \widetilde{X}^\alpha$ such that
$(x_{\alpha_0}, \widehat{x}_{\alpha_0}^{(\alpha, k)}) \in \widetilde{f}^{\alpha, \alpha_0}$,
the pull back $(\varphi^\alpha_{x_{(\alpha, k)}, \widehat{x}_{\alpha_0}^{(\alpha, k)}})^\ast
\boldsymbol{\epsilon}^\alpha_{x_{(\alpha, k)}}$ does not depend on
$x_{(\alpha, k)} \in P^\alpha(\widehat{x}_{\alpha_0}^{(\alpha, k)})$
on some neighborhood $\U_{\widehat{x}_{\alpha_0}^{(\alpha, k)}| x_{\alpha_0}}$ of
$\varphi^{f^{\alpha, \alpha_0}}_{\widehat{x}_{\alpha_0}^{(\alpha, k)}, x_{\alpha_0}}
(\V^{\alpha_0}_{x_{\alpha_0}})$.
Shrinking $\U_{\widehat{x}_{\alpha_0}^{(\alpha, k)}| x_{\alpha_0}}$ if necessary,
we may assume that
\[
\varphi^{f^{\beta, \alpha}}_{\widehat{x}_{\alpha_0}^{(\beta, k')},
\widehat{x}_{\alpha_0}^{(\alpha, k)}}(\U_{\widehat{x}_{\alpha_0}^{(\alpha, k)}| x_{\alpha_0}})
\subset \U_{\widehat{x}_{\alpha_0}^{(\beta, k')}| x_{\alpha_0}}
\]
for $(\widehat{x}_{\alpha_0}^{(\alpha, k)}, \widehat{x}_{\alpha_0}^{(\beta, k')})
\in \widetilde{f}^{\beta, \alpha}$.
Define a partial order of
$\bigcup_{\alpha \leq \alpha_0} \widetilde{f}^{\alpha, \alpha_0}(x_{\alpha_0})$ by
the condition that
$\widehat{x}_{\alpha_0}^{(\alpha, k)} > \widehat{x}_{\alpha_0}^{(\beta, k')}$ if
$(\widehat{x}_{\alpha_0}^{(\alpha, k)}, \widehat{x}_{\alpha_0}^{(\beta, k')})
\in \widetilde{f}^{\beta, \alpha}$.
Then
\[
((\U_{\widehat{x}_{\alpha_0}^{(\alpha, k)}| x_{\alpha_0}})
_{\widehat{x}_{\alpha_0}^{(\alpha, k)} \in
\bigcup_{\alpha \leq \alpha_0} \widetilde{f}^{\alpha, \alpha_0}(x_{\alpha_0})},
(\varphi^{f^{\beta, \alpha}}_{\widehat{x}_{\alpha_0}^{(\beta, k')},
\widehat{x}_{\alpha_0}^{(\alpha, k)}},
\hat \varphi^{f^{\beta, \alpha}}_{\widehat{x}_{\alpha_0}^{(\beta, k')},
\widehat{x}_{\alpha_0}^{(\alpha, k)}})
_{\widehat{x}_{\alpha_0}^{(\alpha, k)} > \widehat{x}_{\alpha_0}^{(\beta, k')}})
\]
is a compatible system of multi-valued partial submersions in the sense of
Definition \ref{compatible system of multi-valued partial submersions for orbibundles},
where $\U_{x_{\alpha_0}| x_{\alpha_0}} = \V^{\alpha_0}_{x_{\alpha_0}}$.
We claim that $(\varphi^\alpha_{x_{(\alpha, k)}, \widehat{x}_{\alpha_0}^{(\alpha, k)}})^\ast
\boldsymbol{\epsilon}^\alpha_{x_{(\alpha, k)}}$ are compatible
with respect to this system.
This implies that using
Lemma \ref{extension of grouped multisection for multi-valued partial submersions}
and \ref{extension of grouped multisection for embedding},
we can construct a compatible grouped multisection
$(\boldsymbol{\epsilon}^{\alpha_0}_{x_{\alpha_0}})_{x_{\alpha_0} \in P^{\alpha_0}}$
of $(x_{\alpha_0}, \V^{\alpha_0}_{x_{\alpha_0}})_{x_{\alpha_0} \in P^{\alpha_0}}$.

The above claim is proved as follows.
For any $\widehat{x}_{\alpha_0}^{(\alpha, k)}, \widehat{x}_{\alpha_0}^{(\beta, k')}
\in \bigcup_{\alpha \leq \alpha_0} \widetilde{f}^{\alpha, \alpha_0}(x_{\alpha_0})$ such that
$(\widehat{x}_{\alpha_0}^{(\alpha, k)}, \widehat{x}_{\alpha_0}^{(\beta, k')}) \in
\widetilde{f}^{\beta, \alpha}$ and any
$x_{(\alpha, k)} \in P^\alpha(\widehat{x}_{\alpha_0}^{(\alpha, k)})$,
there exists some $\widehat{x}_{(\alpha, k)}^{(\beta, k')} \in
\widetilde{f}^{\beta, \alpha}(x_{(\alpha, k)})$ such that
\begin{equation}
\widehat{x}_{\alpha_0}^{(\beta, k')} \in
\widetilde{\psi}^\beta_{\widehat{x}_{(\alpha, k)}^{(\beta, k')}}
\Bigl(\bigl(s^\beta_{\widehat{x}_{(\alpha, k)}^{(\beta, k')}}\bigr)^{-1}(0) \cap
\varphi^{f^{\beta, \alpha}}_{\widehat{x}_{(\alpha, k)}^{(\beta, k')}, x_{(\alpha, k)}}
(\V^\alpha_{x_{(\alpha, k)}})\Bigr).
\label{(alpha, k) (beta, k') relation}
\end{equation}
This is because Condition \ref{varphi condition for multi-valued partial submersion}
for $f^{\beta, \alpha}$ implies
\[
\coprod_{\widehat{x}_{(\alpha, k)}^{(\beta, k')}
\in \widetilde{f}^{\beta, \alpha}(x_{(\alpha, k)})}
(\widetilde{\psi}^\beta_{\widehat{x}_{(\alpha, k)}^{(\beta, k')}}
\circ \varphi^{f^{\beta, \alpha}}_{\widehat{x}_{(\alpha, k)}^{(\beta, k')}, x_{(\alpha, k)}})
= \widetilde{f}^{\beta, \alpha} \circ \widetilde{\psi}^\alpha_{x_{(\alpha, k)}}
\]
and $\widetilde{\psi}^\alpha_{x_{(\alpha, k)}}
((s^\alpha_{x_{(\alpha, k)}})^{-1}(0) \cap \V^\alpha_{x_{(\alpha, k)}})$
contains $\widehat{x}_{\alpha_0}^{(\alpha, k)}$.
Note that (\ref{(alpha, k) (beta, k') relation}) and
Condition \ref{(x1, widehat x1, x2)-composition} for $f^{\beta, \alpha}$ imply
$P^\beta(\widehat{x}_{(\alpha, k)}^{(\beta, k')}) \subset
P^\beta(\widehat{x}_{\alpha_0}^{(\beta, k')})$.
Choose one point $x_{(\beta, k')} \in P^\beta(\widehat{x}_{(\alpha, k)}^{(\beta, k')})$.
It is enough to prove that
$(\varphi^\alpha_{x_{(\alpha, k)}, \widehat{x}_{\alpha_0}^{(\alpha, k)}})^\ast
\boldsymbol{\epsilon}^\alpha_{x_{(\alpha, k)}}$
on $\U_{\widehat{x}_{\alpha_0}^{(\alpha, k)}| x_{\alpha_0}}$ and
$(\varphi^\beta_{x_{(\beta, k')}, \widehat{x}_{\alpha_0}^{(\beta, k')}})^\ast
\boldsymbol{\epsilon}^\beta_{x_{(\beta, k')}}$
on $\U_{\widehat{x}_{\alpha_0}^{(\beta, k')}| x_{\alpha_0}}$ are
$(\varphi^{f^{\beta, \alpha}}_{\widehat{x}_{\alpha_0}^{(\beta, k')},
\widehat{x}_{\alpha_0}^{(\alpha, k)}},
\hat \varphi^{f^{\beta, \alpha}}_{\widehat{x}_{\alpha_0}^{(\beta, k')},
\widehat{x}_{\alpha_0}^{(\alpha, k)}})$-related,
and this follows from the assumption that
$\boldsymbol{\epsilon}^\alpha_{x_{(\alpha, k)}}$ and
$\boldsymbol{\epsilon}^\beta_{x_{(\beta, k')}}$ are
$(\varphi^\beta_{x_{(\beta, k')}, \widehat{x}_{(\alpha, k)}^{(\beta, k')}} \circ
\varphi^{f^{\beta, \alpha}}_{\widehat{x}_{(\alpha, k)}^{(\beta, k')}, x_{(\alpha, k)}},
\hat \varphi^\beta_{x_{(\beta, k')}, \widehat{x}_{(\alpha, k)}^{(\beta, k')}} \circ
\hat \varphi^{f^{\beta, \alpha}}_{\widehat{x}_{(\alpha, k)}^{(\beta, k')},
x_{(\alpha, k)}})$-related.
\end{proof}

We can treat the case of the fiber products $h^{-1}(K) \subset X$
of pre-Kuranishi spaces $X$ with simplicial complexes $K \subset \W$
by strong smooth submersions $h = (h_x)_{x \in \widetilde{X}}$
and essential submersions between them.
In this case, multi-valued partial essential submersion is defined on
a neighborhood of $h^{-1}(K) \subset X$, and we read the every compatibility
conditions as the conditions on a neighborhood of the fiber products.
Then we can also apply the same argument in this case.

\begin{rem}
Since we need to assume that grouped multisections are sufficiently small,
we cannot treat an infinite family of pre-Kuranishi spaces as a compatible
family of multi-valued partial submersions.
In application, we usually construct invariants for finite subfamilies of
pre-Kuranishi spaces and construct the invariant of the infinite family as a limit.
Hence it is enough to assume that grouped multisections are sufficiently small
for the construction of the invariant of a fixed subfamily and
the proof of its invariance.
\end{rem}

\subsection{Continuous family of multisections}
\label{continuous family of multisections}
First consider the following example.
Let $X$ be a $0$-dimensional pre-Kuranishi space and
$f = (f_x)_{x \in \widetilde{X}}$ be a strong smooth map from $X$ to a manifold $Y$.
We want to construct a perturbed multisection of $X$ such that
the induced perturbed multisection of $(f \times f)^{-1}(\Delta_Y) \subset X \times X$
also satisfies the transversal condition, but
it is impossible unless the perturbed multisection do not take zero or
$\dim Y = 0$.
To treat such a case, we use continuous family of multisections and make
the restriction of $f_x$ to every branch of the perturbed multisection submersive.

Continuous family of multisection were used in \cite{FOOO10} and
\cite{FOOO11}.
We recall its definition in our setting.
\begin{defi}
For an orbibundle chart $(\V, \E)$,
let $D$ be a finite-dimensional oriented open disk, and
consider the pull back bundle $p_V^\ast E$ by the projection $p_V : V \times D \to V$.
Let $s^\omega : V \times D \to p_V^\ast E$ ($\omega \in \Omega$)
be a family of smooth sections indexed by a finite $G_V$-set $\Omega$
such that $s^{g\omega} = g_\ast s^\omega$ for any $\omega \in \Omega$ and
$g \in G_V$.
Let $\alpha$ be a top-dimensional form on $D$ with compact support
such that $\int_D \alpha = 1$.
We call such a triple $\epsilon = (D, (s^\omega)_{\omega \in \Omega}, \alpha)$
a continuous family of multisections of $(\V, \E)$.

We also define the version of grouped multisection similarly as follows.
A continuous family of grouped multisections
$\boldsymbol{\epsilon} = (D_j, (\epsilon^\omega)_{\omega \in \coprod_j \Omega_j},
\alpha_j)$ of an orbibundle chart $(\V, \E)$
consists of the following.
$(D_j)_{j = 1, \dots, k}$ are finite number of finite dimensional oriented open disks,
and for each $j$, $s^\omega : V \times D_j \to p_V^\ast E$ ($\omega \in \Omega_j$) is
a family of smooth sections.
Each $\alpha_j$ is a top-dimensional form $\alpha_j$ on $D_j$ with compact support
which satisfies $\int_{D_j} \alpha_j = 1$.
We assume that there is an $G_V$-action on $\coprod_{1 \leq j \leq k} \Omega_j$
which preserves the decomposition and assume that
if $g \in G_V$ maps $\Omega_j$ to $\Omega_{j'}$, then
$D_j = D_{j'}$ and $\alpha_j = \alpha_{j'}$.
We also assume that the smooth sections satisfy
$s^{g\omega} = (g \times 1_{D_j})_\ast s^\omega$ for any $\omega \in \Omega_j$ and
$g \in G_V$.
For each $j$, we define $\epsilon_j = (D_j, (\epsilon^\omega)_{\omega \in \Omega_j},
\alpha_j)$ and also denote the family of grouped multisections by $\boldsymbol{\epsilon}
= \{\epsilon_j\}$.
\end{defi}

We define the support of each $\epsilon_j$ by
$\supp(\epsilon_j) = \bigcup_{j \in \Omega_j} p_V(\supp(\epsilon^\omega)) \subset V$.
For a connected open subset $\U \subset \V$, the restriction of
a family of grouped multisections $\boldsymbol{\epsilon}
= (D_j, (\epsilon^\omega)_{\omega \in \coprod_j \Omega_j}, \alpha_j)$ of
$(\V, \E)$ to $(\U, \E|_{\U})$ is defined by
\[
\boldsymbol{\epsilon}|_\U = ((D_j)_{j \in I_U},
(\epsilon^\omega|_{U \times D_j})_{\omega \in \coprod_{j \in I_U}  \Omega_j},
(\alpha_j)_{j \in I_U}),
\]
where $I_U = \{j; \supp(\epsilon_j) \cap U \neq \emptyset\}$.
We can similarly define $(\varphi, \hat \varphi)$-relation of continuous families of
grouped multisections for an embedding $(\varphi, \hat \varphi)$ between
orbibundle charts, and pull back of a continuous family of grouped multisections
by a submersion.

For a smooth section $s$ and a continuous family of grouped multisections
$\boldsymbol{\epsilon} = (D_j, (\epsilon^\omega)_{\omega \in \coprod_j \Omega_j},
\alpha_j)$ of an oribibundle chart, we define their sum
by the continuous family of multisections
\[
s + \boldsymbol{\epsilon}
= \bigl(\prod_j D_j, \bigl(s + \sum_j \epsilon^{\omega_j}\bigr)
_{(\omega_j) \in \prod_j \Omega_j},
\alpha_1 \wedge \dots \wedge \alpha_k\bigr).
\]

Let $f = (f_x)_{x \in \widetilde{X}}$ be a strong smooth map from a pre-Kuranishi space
$X$ to a manifold $Y$ such that each $f_x : \V_x \to Y$ is submersive.
Then for a weakly good coordinate system $(x, \V_x)_{x \in P}$ of $X$,
similarly to Lemma \ref{construction of grouped multisection},
shrinking $\V_x$ slightly if necessary,
we can construct a continuous family of grouped multisections
$\boldsymbol{\epsilon} = (\boldsymbol{\epsilon}_x)_{x \in P}$ for $(x, \V_x)_{x \in P}$
which satisfies the following transversality condition:
For any orbibundle chart $(\V, \E)$ in $(\V_x, \E_x)$,
every branch of the multisection $s|_\V + \boldsymbol{\epsilon}_x|_\V$ is
transverse to the zero section, and the restriction of $f_x$ to its zero set is
submersive.
Furthermore, the same holds for the restriction to the corners of $V$.

For a continuous family of perturbed multisections,
it is not suitable to represent the virtual fundamental chain as a singular chain.
Instead, for strong smooth maps $f = (f_x)_{x \in \widetilde{X}}$ from $X$
to a manifold $Y$ and $h = (h_x)_{x \in \widetilde{X}}$ from $X$ to a manifold $Z$,
we represent the virtual fundamental chain as a linear map
$(h_! \circ f^\ast)_X : \Omega(Y) \to \Omega(Z)$.
This map is defined as follows.
As in the usual case, we take a partition of unity $(\beta_x)_{x \in P}$
subordinate to $(x, \V_x)_{x \in P}$, finite number of orbibundle charts
$(\V_\tau, \E_\tau)_{\tau \in T_x}$ of $(\V_x, \E_x)$ and smooth functions
$\beta_\tau : \V_\tau \to \R$ with compact support such that
$\beta_x = \sum_{\tau \in T_x} \beta_\tau$.
Then for each differential form $\theta \in \Omega(Y)$,
$(h_! \circ f^\ast)_X \theta \in \Omega(Z)$ is defined by
\begin{align*}
&(h_! \circ f^\ast)_X \theta \\
&= \sum_{\substack{x \in P \\ \tau \in T_x}}
\frac{\sum_{(\omega_j) \in \prod_j \Omega_{\tau, j}}
\Bigl(h_x|_{\{s_\tau^{(\omega_j)} = 0\}}\Bigr)_{\textstyle !}\,
\Bigl(\beta_\tau \cdot (f_x|_{V_\tau})^\ast
\theta \wedge \alpha_1 \wedge \dots \wedge \alpha_k|_{\{s_\tau^{(\omega_j)} = 0\}}\Bigr)}
{\# G_{V_\tau}
\cdot \prod_j \# \Omega_{\tau, j}}
\end{align*}
instead of Equation (\ref{virtual fundamental chain by forms}).
In this equation, we define the orientation of the fiber of
$h_x|_{\{s_\tau^{(\omega_j)} = 0\}}$ similarly to the usual case
using the following orientation of $V_x \times \prod_j D_j$.
The orientation of $V_x \times \prod_j D_j$ is defined by
$T(V_x \times \prod_j D_j) = (-1)^{(\sum_j \dim D_j) \rank E_x}
T V_x \oplus \bigoplus_j T D_j$.

%% file: SFT-03_Kuranishi_neighborhood.tex
%
%

\section{Construction of pre-Kuranishi structure}\label{construction of Kuranishi}
In this section, we explain the way to construct
a pre-Kuranishi structure of $\widehat{\M} = \widehat{\M}(Y, \lambda, J)$.
This is the basis of the construction of pre-Kuranishi structures of other various spaces
such as fiber products in Section \ref{fiber prod}.

First we explain the Banach spaces we use.
Let $\Sigma_0$ be the domain curve of a holomorphic building $(\Sigma_0, z, u_0) \in
\widehat{\M}$, and let $\{\mu\}$ and $\{\pm\infty_i\}$ be the indices of
its joint circles and limit circles respectively.
Define positive constants $\delta_{0, \mu}$ and $\delta_{0, \pm\infty_i}$ by the
minimal nonzero absolute value of eigenvalues of $A_{\gamma_\mu}$ and
$A_{\gamma_{\pm\infty_i}}$ respectively,
where $\gamma_\mu$ and $\gamma_{\pm\infty_i}$ are the periodic orbits on
the corresponding imaginary circles of $\Sigma$.
(See Definition \ref{def of Bott-Morse} for the definition of the operator $A_\gamma$
for each periodic orbit $\gamma$.)
For a sequence of positive constants
$\delta = ((\delta_\mu)_{\mu}, (\delta_{\pm\infty_i})_{\pm\infty_i})$
such that $\delta_\mu < \delta_{0, \mu}$ and $\delta_{\pm\infty_i}
< \delta_{0, \pm\infty_i}$,
we use the Banach spaces $L_\delta^p(\Sigma)$ and $W_\delta^{1, p}(\Sigma)$
defined as follows.
Fix some coordinate $([0, \infty] \cup [-\infty, 0]) \times S^1$ of a neighborhood
$N_\mu$ of each joint circle $S_\mu^1$ of $\Sigma$,
and some coordinate $[0, \infty] \times S^1$ or $[-\infty, 0] \times S^1$ of
a neighborhood $N_{\pm\infty_i}$ of each limit circle $S^1_{\pm\infty_i}$,
and fix a volume form of $\Sigma$ such that its restriction to these neighborhoods
coincide with the usual Lebesgue measure $ds \wedge dt$.
(On a neighborhood $D\cup D$ of each nodal point,
we use usual volume form of $D$.)
Then $L_\delta^p$-norm of $\xi$ is defined by
\begin{align*}
||\xi||^p_{L_\delta^p} = &\int_{\Sigma \setminus (\coprod N_\mu \sqcup \coprod
N_{\pm\infty_i})} |\xi|^p \vol
+ \sum_\mu \int_{N_\mu} |e^{\delta_\mu |s|}\xi(s,t)|^p ds\wedge dt\\
&+ \sum_{\pm\infty_i}
\int_{\substack{[0, \infty] \times S^1\\ \text{or} \\ [-\infty, 0] \times S^1}}
|e^{\delta_{\pm\infty_i} |s|}\xi(s,t)|^p ds\wedge dt
\end{align*}

The Sobolev space $W^{1, p}_\delta(\Sigma)$ is the space of continuous functions
(or continuous sections) $\xi$ on $\Sigma$ whose $W^{1, p}_\delta$-norms
\begin{align*}
||\xi||^p_{W_\delta^{1, p}} = &\int_{\Sigma \setminus (\coprod N_\mu \sqcup \coprod
N_{\pm\infty_i})} (|\xi|^p + |\nabla \xi|^p) \vol\\
&+ \sum_\mu \int_{N_\mu} (|e^{\delta_\mu |s|}\xi|^p + |e^{\delta_\mu |s|} \partial_s \xi|^p
+ |e^{\delta_\mu |s|} \partial_t \xi|^p) ds\wedge dt\\
&+ \sum_{\pm\infty_i}
\int_{\substack{[0, \infty] \times S^1\\ \text{or} \\ [-\infty, 0] \times S^1}}
(|e^{\delta_{\pm\infty_i} |s|}\xi|^p + |e^{\delta_{\pm\infty_i} |s|} \partial_s \xi|^p
+ |e^{\delta_{\pm\infty_i} |s|} \partial_t \xi|^p) ds\wedge dt
\end{align*}
are finite.
For each holomorphic building $(\Sigma, z, u_0)$,
$\widetilde{W}^{1, p}_\delta(\Sigma, u_0^\ast T \hat Y)$
is the space of continuous sections $\xi$ of $u_0^\ast T \hat Y = \R \oplus
(\pi_Y \circ u_0)^\ast TY$ such that
\[
\xi = \xi_0 + \sum_\mu \beta_\mu v_\mu + \sum_{\pm\infty_i} \beta_{\pm\infty_i}
v_{\pm\infty_i}
\]
for some $\xi_0 \in W^{1, p}_\delta(\Sigma, u_0^\ast T \hat Y)$,
$v_\mu \in \Ker A_{\gamma_\mu}$ and $v_{\pm\infty_i} \in
\Ker A_{\gamma_{\pm\infty_i}}$, where $\beta_\mu$ is a smooth function which is $1$
on a neighborhood of $\mu$-th joint circle and whose support is contained in
its slightly larger neighborhood for each $\mu$, and $\beta_{\pm\infty_i}$ is a smooth
function which is $1$ on a neighborhood of the limit circle $S^1_{\pm\infty_i}$ and whose
support is contained in its slightly large neighborhood for each $\pm\infty_i$.
In the above equation, we regard $v_\mu$ as a section defined on
$([0, \infty] \cup [-\infty, 0] ) \times S^1$ by $v_\mu(s, t) = v_\mu(t)$, where
we fix a trivialization of $u_0^\ast T \hat Y$
on $([0, \infty] \cup [-\infty, 0]) \times S^1$.
The meaning of $v_{\pm\infty_i}$ is similar.
The definition of $\widetilde{W}^{1, p}_\delta(\Sigma, u_0^\ast T \hat Y)$ does not
depend on the choice of $\beta_\mu$ and $\beta_{\pm\infty_i}$.
As a Banach space, we regard $\widetilde{W}^{1, p}_\delta(\Sigma, u_0^\ast T \hat Y)$
as a direct sum of $W^{1, p}_\delta(\Sigma, u_0^\ast T \hat Y)$, $\Ker A_{\gamma_\mu}$
and $\Ker A_{\gamma_{\pm\infty_i}}$.

For a family of deformations of $\Sigma$, we need to use an appropriate family of
norms to obtain a uniform estimate.
This family of norms are used only for the construction of a Kuranishi neighborhood of
a holomorphic building $(\Sigma, z, u)$ and we do not need to assume that the
norm of a curve $\Sigma'$ as a deformation of $\Sigma$ coincides with that
used for the construction of a Kuranishi neighborhood of a holomorphic building
whose domain curve is $\Sigma'$.
Let $\Sigma'$ be a curve obtained from $\Sigma$ by replacing the neighborhood
$([0, \infty] \cup [-\infty, 0]) \times S^1$ of $S^1_\mu$ by $([0, \rho_\mu] \cup
[-\rho_\mu, 0]) \times S^1$ and the neighborhood $D \cup D$ of each nodal point
$q_\nu$ by $\{(x, y) \in D \times D; xy = \zeta_\nu\}$ for some $(\rho_\mu, \zeta_\nu)$.
Then the $L^p_\delta$-norm of $L^p_\delta(\Sigma')$ is defined by
\begin{align*}
||\xi||^p_{L_\delta^p} = &\int_{\Sigma \setminus (\coprod N_\mu \sqcup \coprod
N_{\pm\infty_i})} |\xi|^p \vol
+ \sum_\mu \int_{([0, \rho_\mu] \cup [-\rho_\mu, 0]) \times S^1}
|e^{\delta_\mu |s|}\xi(s,t)|^p ds\wedge dt\\
&+ \sum_{\pm\infty_i}
\int_{\substack{[0, \infty] \times S^1\\ \text{or} \\ [-\infty, 0] \times S^1}}
|e^{\delta_{\pm\infty_i} |s|}\xi(s,t)|^p ds\wedge dt,
\end{align*}
where the volume form on $\{(x, y) \in D \times D; xy = \zeta_\nu\}$ is defined by
$\frac{\sqrt{-1}}{2}dx \wedge d\bar x$ on $\{|x| \geq |y|\}$ and
$\frac{\sqrt{-1}}{2}dy \wedge d\bar y$ on $\{|y| \geq |x|\}$.
The norm of $W^{1, p}_\delta(\Sigma')$ is defined similarly.
The norm of $\widetilde{W}^{1, p}_\delta(\Sigma', u_0^\ast T \hat Y)$ is defined by
\begin{align*}
||\xi||_{\widetilde{W}^{1, p}_\delta}
= \inf \{& ||\xi_0||_{W^{1, p}_\delta(\Sigma')} + \sum_\mu ||v_\mu||_{\Ker A_{\gamma_\mu}}
+ \sum_{\pm\infty_i} ||v_{\pm\infty_i}||_{\Ker A_{\gamma\infty_i}};\\
&\xi = \xi_0 + \sum_\mu \beta_\mu v_\mu
+ \sum_{\pm\infty_i} \beta_{\pm\infty_i} v_{\pm\infty_i},\\
&\xi_0 \in W^{1, p}_\delta(\Sigma', u_0^\ast T \hat Y), v_\mu \in \Ker A_{\gamma_\mu},
v_{\pm\infty_i} \in \Ker A_{\gamma_{\pm\infty_i}}\}.
\end{align*}

In Section \ref{construction of nbds}, we explain the construction of
a Kuranishi neighborhood of a point in $\widehat{\M}(Y, \lambda, J)$, assuming
sufficient data including an additional vector space are given.
To construct a Kuranishi neighborhood by inverse function theorem of Banach spaces,
we need to prove the linearized gluing lemma, which is proved in Section
\ref{linearized gluing}.
In Section \ref{smoothness}, we prove the smoothness of Kuranishi neighborhood, and
in Section \ref{embed}, we consider the embedding of Kuranishi neighborhoods and
prove its smoothness.
In Section \ref{Kuranishi of disconnected buildings}, we consider the relation of
the Kuranishi neighborhoods of $\widehat{\M}^0$ and $\widehat{\M}$.
Finally in Section \ref{global construction},
we construct a global Kuranishi structure of $\widehat{\M}$.

\subsection{Construction of Kuranishi neighborhoods}\label{construction of nbds}
First we explain a way to construct a Kuranishi neighborhood of a point
$p_ 0= (\Sigma_0, z, u_0) \in \widehat{\M}(Y, \lambda, J)$.
The construction is based on the implicit function theorem (or inverse function theorem)
for Banach spaces of functions (or sections of some vector bundles) on deformed curves
of $\Sigma_0$.
Since the Banach space changes if the domain curve changes, we need to apply
the implicit function theorem for each deformed curve.
Using appropriate norms for these Banach spaces, we can apply the implicit function
theorem for them uniformly, and get a Kuranishi neighborhood as an (at least) continuous
fibration over the parameter space of the deformation of the domain curve.
We prove in Section \ref{smoothness} that this fibration is actually smooth in some sense,
and in Section \ref{embed}, we prove the smoothness of the embedding between
two Kuranishi neighborhoods.

We fix an order $z = (z_i)$ of the marked points.
As we have explained, to define a Kuranishi neighborhood,
we need an additional vector space which makes the Fredholm map transverse to zero.
Such an additional vector space is given as the following data
$(p_0^+, S, E^0, \lambda)$:
(These are given in Section \ref{global construction}.)
\begin{figure}
\centering
\includegraphics[width= 350pt]{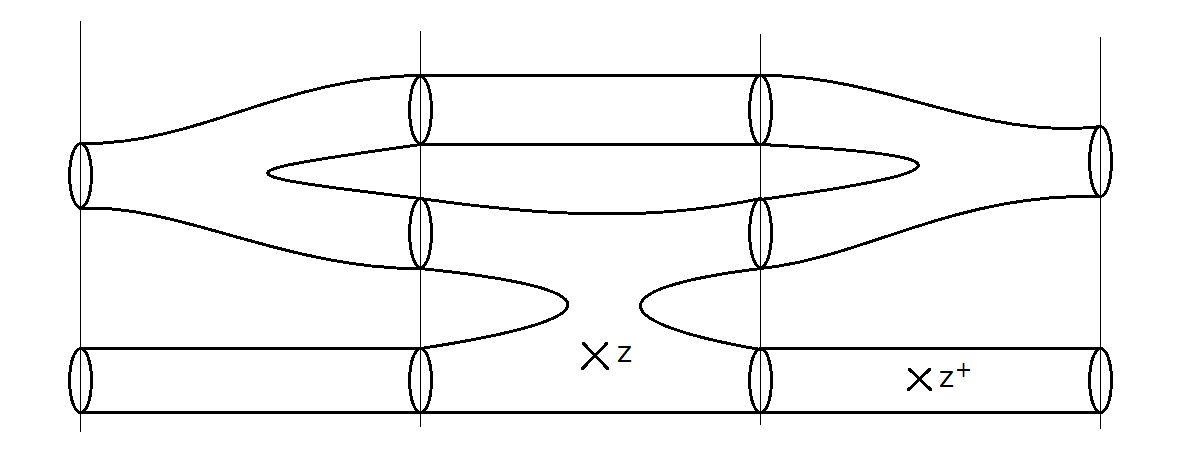}
\caption{$(\Sigma_0, z \cup z^+)$}\label{(Sigma,z,z+)}
\includegraphics[width= 350pt]{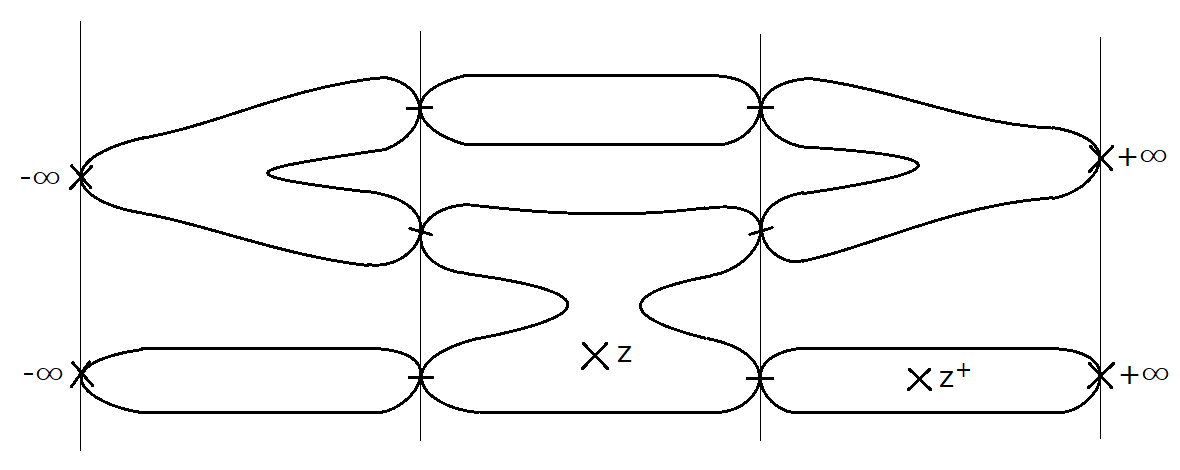}
\caption{$(\check\Sigma_0, z \cup z^+ \cup (\pm\infty_i))$}
\label{(checkSigma,z,z+,pminfty)}
\includegraphics[width= 350pt]{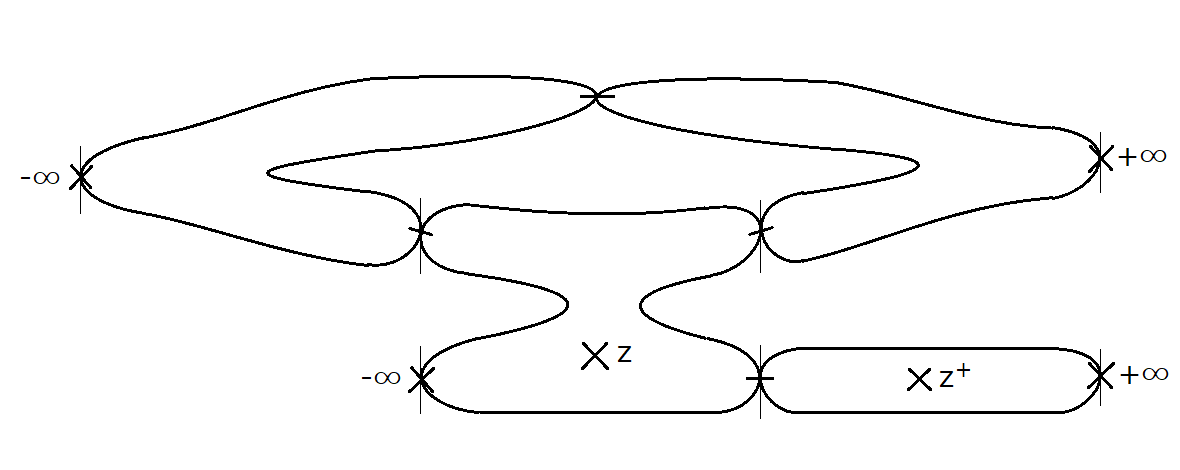}
\caption{$(\hat \Sigma_0, z \cup z^+ \cup (\pm\infty_i))$}\label{(hatSigma,z,z+,pminfty)}
\end{figure}
\begin{itemize}
\item
$p_0^+ = (\Sigma_0, z \cup z^+, u_0)$ is a curve obtained by adding marked points
on the nontrivial components of $\Sigma_0$.
We assume all unstable components of $(\Sigma_0, z \cup z^+)$ are trivial cylinders
of $p_0$.
We assume that $G_0 = \Aut(\Sigma_0, z, u_0)
:= \{g \in \Aut (\Sigma_0); g (\{z_i\}) = \{z_i\}, u_0 \circ g = u_0\}$ preserves $z^+$
as a set, that is, $G_0$ acts on $z^+ = (z_i^+)$ as a symmetric group.
\item
$S \subset Y$ is a finite union of codimension-two submanifolds such that
$\pi_Y \circ u_0$ intersects with $S$ at $z^+$ transversely.
(We do not assume the transversality of the other intersections of $\pi_Y \circ u_0$
with $S$.)
This will be used to kill the excessive dimension of deformation
due to the additional marked points $z^+$.
\item
Let $(\hat \Sigma_0, z \cup z^+ \cup (\pm\infty_i))$ be the stabilization of
$(\check \Sigma_0, z \cup z^+ \cup (\pm\infty_i))$, and
let $(\hat P \to \hat X, Z \cup Z^+ \cup Z_{\pm\infty_i})$ be its local universal family.
Since $G_0$ acts on $\hat \Sigma_0$ preserving $z$, $z^+$ and $\{\pm\infty_i\}$ as
sets, it also acts on $\hat P$ by the universal property of $\hat P$.
Then $E^0$ is a finite dimensional $G_0$-vector space
and $\lambda : E^0 \to C^\infty(\hat P \times Y, \Wedge^{0, 1} V^\ast \hat P
\otimes (\R \partial_\sigma \oplus TY))$ is a $G_0$-equivariant linear map
which satisfies the following conditions:
($V^\ast \hat P$ is the dual of the vertical tangent space $V \hat P \subset T \hat P$
of $\hat P$.)
\begin{itemize}
\item
For each $h \in E^0$, the projection of the support of $\lambda(h)$ to $\hat P$
does not intersect with the nodal points of $\hat P$ or $Z_{\pm\infty_i}$.
\item
Let $E^0 \to C^\infty(\Sigma_0 \times Y, \Wedge^{0, 1} T^\ast\Sigma_0 \otimes
(\R \partial_\sigma \oplus TY))$ be the pullback of $\lambda$ by the composition of the
blowing down $\Sigma_0 \to \check \Sigma_0$ and the forgetful map
$(\check \Sigma_0, z \cup z^+) \tocong (\hat P_0, Z(0) \cup Z^+(0))$.
This pull back is also denoted by $\lambda$.
Then we assume that for a sufficiently small $\delta > 0$, the linear map
\begin{align*}
&D_{p_0}^+ : \widetilde{W}_\delta^{1, p}(\Sigma_0, u_0^\ast T \hat Y) \oplus E^0\\
&\to
L_\delta^p(\Sigma_0, \Wedge^{0, 1}T^\ast \Sigma_0 \otimes u_0^\ast T \hat Y)
\oplus
\bigoplus_{\text{limit circles}} \Ker A_{\gamma_{\pm\infty_i}} / (\R \partial_\sigma
\oplus \R R_\lambda)\\
&\quad \oplus \bigoplus_{z_i} T_{\pi_Y \circ u_0(z_i)} Y
\end{align*}
{\abovedisplayskip=-5pt
\begin{align*}
&(\xi, h) \mapsto (D_{p_0} \xi(z) + \lambda(h)(z, \pi_Y \circ u_0(z)), \\
&\hph{(\xi, h) \mapsto (}
\sum_j \langle\xi|_{S^1_{\pm\infty_i}}, \eta_j^{\pm\infty_i}\rangle
\eta_j^{\pm\infty_i}, \pi_Y \circ \xi (z_i))
\end{align*}}
is surjective,
where $D_{p_0}$ is the linearization of the equation of $J$-holomorphic maps,
that is,
\[
D_{p_0} \xi = \nabla \xi + J(u_0) \nabla \xi j + \nabla_\xi J(u_0) du_0 j,
\]
and $\{\eta_j^{\pm\infty_i}\}_j$ is an orthonormal basis of
the orthogonal complement of $\R \partial_\sigma \oplus \R R_\lambda$ in
$\Ker A_{\gamma_{\pm\infty_i}}$ for each $\pm\infty_i$.
\end{itemize}
\end{itemize}

The above data are given in the form which respects the $\R$-translation invariance.
However, to describe the Kuranishi neighborhood, we further need to fix the following
temporally data $(z^{++}, S', \hat R_i)$ which break the $\R$-translation invariance:
(The Kuranishi neighborhood constructed finally does not depend on these temporally
data. See Section \ref{embed}.)
\begin{itemize}
\item
$z^{++} = (z^{++}_i) \subset \Sigma$ are additional marked points which make
$(\Sigma_0, z \cup z^+ \cup z^{++})$ stable.
We assume $G_0$-action preserves $z^{++}$ as a set.
\item
$S' \subset (\R_1 \cup \dots \cup \R_k) \times Y$ is a codimension-two submanifold
such that $u_0$ intersects with $S'$ at $z^{++}$ transversely.
\item
For each $1 \leq i \leq k$,
let $\hat R_i = (\hat R_{i, l})_{1 \leq l \leq m_i}$ be a family of holomorphic sections
$\hat R_{i, l} : \hat X \to \hat P$ such that
$\sigma_i \circ u_0(\hat R_{i, l}(0)) = 0$, where $\sigma_i$ is the coordinate of $\R_i$,
and $\hat R_i$ is $G_0$-invariant as a family.
($G_0$ may permute $\{ \hat R_{i, l} \}_l$.)
We assume $\hat R_{i, l}$ do not intersect with nodal points or $Z_{\pm\infty_i}$.
Let $(\widetilde{P} \to \widetilde{X}, Z \cup Z^+ \cup Z^{++})$ be the local universal
family of $(\Sigma_0, z \cup z^+ \cup z^{++})$.
Then each $\hat R_{i, l}$ induces a section
$\widetilde{R}_{i, l} : \widetilde{X} \to \widetilde{P}$
which makes following diagram commutative.
\[
\begin{tikzcd}
\widetilde{P} \ar{r}{\forget}& \hat P\\
\widetilde{X} \ar{u}{\widetilde{R}_{i, l}} \ar{r}{\forget}& \hat X \ar{u}{\hat R_{i, l}}
\end{tikzcd}
\]
These families $\widetilde{R}_{i, l}$ are used to kill the $\R$-translations
by imposing the conditions $\sum_l \sigma_i \circ u(\widetilde{R}_{i, l}) = 0$.
The fact that $\widetilde{R}_{i, l}$ are induced by sections $\hat R_{i, l}$ of
$\hat P \to \hat X$ is important to define smooth embeddings in Section \ref{embed}.
\end{itemize}
The pullback $E^0 \to C^\infty(\widetilde{P} \times Y, \Wedge^{0, 1}V^\ast
\widetilde{P} \otimes (\R \partial_\sigma \oplus TY))$
of $\lambda$ by $\widetilde{P} \to \hat P$ is also denoted by $\lambda$.

Using the above data, we construct a Kuranishi neighborhood of $p_0$.
First we explain a convenient way to express curves close to $p_0$.
We separate the domain $\Sigma_0$ into several parts, namely,
neighborhoods of nodal points, neighborhoods of joint circles,
neighborhood of limit circles and the rest.

The local universal family $(\widetilde{P} \to \widetilde{X}, Z \cup Z^+ \cup Z^{++})$
can be described as follows.
Let $N_0 \subset \Sigma_0$ be a neighborhood of nodal points and imaginary circles
such that
\begin{itemize}
\item
$N_0 \cong \coprod_{l_0} (D \cup D) \sqcup \coprod_{l_1} (D \widetilde{\cup} D) \sqcup
\coprod_{l_2} \widetilde{D}$
\item
all marked points and $\widetilde{R}_i(0)$ are contained in $\Sigma_0 \setminus N_0$
\item
the support of $\lambda(h)$ is also contained in $\Sigma_0 \setminus N_0$
for each $h \in E^0$.
\end{itemize}
Let $\J_0$ be a finite dimensional complex manifold which consists of
holomorphic structures of $\check \Sigma_0$ near the original one $j_0$ such that
the restriction of any $j \in \J_0$ to $N_0$ coincides with $j_0$.
If we choose an appropriate $\J_0$, then we may regard $\widetilde{X}$ as
a neighborhood $\widetilde{X} \subset \J_0 \times D^{l_0} \times
\widetilde{D}^{l_1}$ of $(j_0, 0, (0, 0))$,
where $(\zeta_\nu)_{\nu =1}^{l_0} \in D^{l_0}$ are the parameters of deformation of
the neighborhood of nodal points, and
$(\zeta_\mu = \rho_\mu^{2\pi}e^{2\pi \sqrt{-1}\varphi_\mu},
e^{2\pi \sqrt{-1}\varphi_\mu})_{\mu = 1}^{l_1} \in \widetilde{D}^{l_1}$
are the parameters of deformation of the neighborhood of joint circles.
We sometimes denote the parameters $(\zeta_\mu, e^{2\pi \sqrt{-1} \varphi_\mu}) \in
\widetilde{D}$ by $(\rho_\mu, \varphi_\mu) \in [0, 1) \times S^1$.
More precisely, for each $a = (j, (\zeta_\nu)_{1 \leq \nu \leq l_0},
(\rho_\mu, \varphi_\mu)_{1 \leq \mu \leq l_1}) \in \widetilde{X}$,
the fiber $\widetilde{P}_a$ at $a$ has the following form.
\begin{align*}
\widetilde{P}_a =& (\Sigma_0 \setminus N_0)\\
&\cup \coprod_{\nu =1}^{l_0} \{(x, y) \in D \times D; xy = \zeta_\nu\}\\
&\cup \coprod_{\mu = 1}^{l_1} \{((s_x, t_x), (s_y, t_y))
\in [-1, \infty] \times S^1 \times [-\infty, 1] \times S^1;\\
&\quad\quad\quad\quad s_y - s_x = \log \rho_\mu,\ t_y -t_x = \varphi_\mu\}\\
&\cup \coprod_{l_2^-} [-\infty, 0] \times S^1 \cup \coprod_{l_2^+} [0, \infty] \times S^1
\end{align*}
($l_2^\pm$ are the number of $\pm$-limit circles respectively.)
Namely, $\widetilde{P}_a$ is obtained from $\Sigma$ by replacing the neighborhood
$D \cup D$ of the $\nu$-th nodal point with
$\{(x, y) \in D \times D; xy = \zeta_\nu\}$,
and the neighborhood $D \widetilde{\cup} D \cong
([-1, \infty] \cup_{\infty = -\infty} [-\infty, 1]) \times S^1$
of the $\mu$-th joint circle with
\begin{align*}
\widetilde{N}_a^\mu
&= \{((s_x, t_x), (s_y, t_y)) \in [-1, \infty] \times S^1 \times [-\infty, 1] \times S^1;\\
&\quad \quad s_y - s_x = \log \rho_\mu,\, t_y -t_x = \varphi_\mu\}.
\end{align*}
The complex structure of $\widetilde{P}_a$ is defined by $j$ on
$\Sigma_0 \setminus N_0$,
and the usual complex structure of the other parts.
(The complex structure is defined on the complement of the imaginary circles.)
The sections of marked points $Z_i$, $Z^+_i$, and $Z^{++}_i$ are defined by
the constant maps $Z_i \equiv z_i$, $Z^+_i \equiv z^+_i$ and $Z^{++}_i \equiv z^{++}_i$
($\in \Sigma_0 \setminus N_0$).
The above expression of $\widetilde{P}$ can be easily obtained by the local structure
of universal family of stable curves. (See \cite{RS06} for example.)
We identify $\widetilde{N}_a^\mu$ with
\[
([-1, - {\textstyle\frac{1}{2}} \log \rho_\mu]
\cup_{-\frac{1}{2} \log \rho_\mu = \frac{1}{2} \log \rho_\mu}
[{\textstyle\frac{1}{2}} \log \rho_\mu, 1]) \times S^1_\mu
\]
by
\begin{align}
[-1, - {\textstyle\frac{1}{2}} \log \rho_\mu] \times S^1_\mu &\to \widetilde{N}_a^\mu
\notag\\
(s, t) & \mapsto (s_x, t_x) = (s, t - {\textstyle\frac{1}{2}} \chi(s) \varphi_\mu)
\label{left coordinate near joint}\\
[{\textstyle\frac{1}{2}} \log \rho_\mu, 1] \times S^1_\mu &\to \widetilde{N}_a^\mu
\notag\\
(s', t') & \mapsto (s_y, t_y) = (s, t + {\textstyle\frac{1}{2}} \chi(-s) \varphi_\mu)
\label{right coordinate near joint}
\end{align}
where $\chi : \R \to \R_{\geq 0}$ is a smooth function such that
$\chi|_{(-\infty, -1]} \equiv 0$ and $\chi|_{[0, \infty)} \equiv 1$.
Let $j_{\varphi_\mu}$ be the complex structure on $([-1, - \frac{1}{2}\log \rho_\mu]
\cup [\frac{1}{2}\log \rho_\mu, 1]) \times S^1_\mu$ defined by the pull back of the
usual complex structure on $\widetilde{N}_a^\mu$ by the above identification.
We note that $-j_{\varphi_\mu} \partial_t
= \partial_s + \frac{1}{2} \varphi_\mu \chi'(s) \partial_t$
on $[-1, - \frac{1}{2}\log \rho_\mu] \times S^1_\mu$, and
$-j_{\varphi_\mu} \partial_t
= \partial_s + \frac{1}{2} \varphi_\mu \chi'(-s) \partial_t$
on $[\frac{1}{2} \log \rho_\mu, 1] \times S^1_\mu$.
Under this identification, we define the $L_\delta^p$-norms of the function spaces
of $\widetilde{P}_a$ as we explained before this section.

We use a parameter space $\mathring{X}$ which reflects the fact that
the splitting of $\hat Y$ occurs simultaneously with the deformation of
the domain curve.
$\mathring{X} \subset \widetilde{X} \times \prod_{\text{joint circles}} \R_\mu$
is a submanifold defined as follows.
Let $M_i$ be the set of joint circles between the $i$-th floor and the $(i+1)$-th floor.
Then $(a, (b_\mu)_\mu) \in \widetilde{X} \times \prod_{\text{joint circles}} \R_\mu$
belongs to $\mathring{X}$ if
$-L_\mu \log \rho_\mu + b_\mu \in (-\infty, \infty]$ is independent of $\mu \in M_i$ for
each $i = 1,2, \dots, k-1$,
where $L_\mu$ is the period of $\gamma_\mu$.
This implies in particular, whether $\rho_\mu = 0$ or not is independent of
$\mu \in M_i$ for each $i$,
and if $\rho_\mu \neq 0$ then $b_\mu$ is determined by $a \in \widetilde{X}$ and
one of $b_\mu$ for each $i$.
If we use an appropriate smooth structure of $\widetilde{X}$ (see Section
\ref{smoothness}), then $\mathring{X}$ is indeed a smooth submanifold of
$\widetilde{X} \times \prod_{\text{joint circles}} \R_\mu$.
For each $(a, b) \in \mathring{X}$, we define an equivalence relation $\sim_{a, b}$ of
$\overline{\R}_1 \sqcup \overline{\R}_2 \sqcup \dots \sqcup \overline{\R}_k$ by
identifying $s \in \overline{\R}_i$ and $s' \in \overline{\R}_{i+1}$ if $\rho_\mu \neq 0$
and $s - s' = -L_\mu \log \rho_\mu + b_\mu$ for some (and all) $\mu \in M_i$, and
identifying $+\infty \in \overline{\R}_i$ and $-\infty \in \overline{\R}_{i + 1}$
if $\rho_\mu = 0$.
Let $0_i \in (\overline{\R}_1 \sqcup \overline{\R}_2 \sqcup \dots \sqcup
\overline{\R}_k)/ \sim_{a, b}$ be the point corresponds to the zero in $\overline{\R}_i$.
If $\mu \in M_i$ and $\rho_\mu \neq 0$, then $b_\mu$ satisfies
\[
0_{i+1} - 0_i = -L_\mu \log \rho_\mu + b_\mu.
\]
If $\rho_\mu = 0$, the maps $u$ corresponding to the parameter $b_\mu$ will be
related to $b_\mu$ by
\begin{align}
b_\mu &= \lim_{s \to \infty} (\sigma \circ u|_{[0, \infty) \times S_\mu^1} (s, t)
- (0_i + L_\mu s))\notag\\
&\quad - \lim_{s \to -\infty} (\sigma \circ u|_{(-\infty, 0] \times S_\mu^1}(s, t)
- (0_{i+1} + L_\mu s)). \label{asymptotic parameter}
\end{align}
We call $b_\mu$ asymptotic parameters.

\begin{rem}
Before starting to construct a Kuranishi neighborhood, we calculate the virtual dimension
of the Kuranishi neighborhood of $p_0 = (\Sigma_0, z, u_0)
\in \widehat{\M}(Y, \lambda, J)$ and check that it coincides with the expected dimension.
Readers may skip this calculation since we do not use it for the construction of
Kuranishi neighborhood.

First, the dimension of the parameter space $\mathring{X}$ is
$\dim \mathring{X} = \dim \widetilde{X} + (k - 1)$, where $k$ is the height of $p_0$.
For each $(a, b) \in \mathring{X}$, we regard the equation of $J$-holomorphic curves
as a Fredholm map, whose index coincides with that of
the linearization $D_{p_0} : \widetilde{W}_\delta^{1, p}(\Sigma_0, u_0^\ast T \hat Y)
\to L_\delta^p(\Sigma_0, \Wedge^{0, 1}T^\ast \Sigma_0 \otimes u_0^\ast T \hat Y)$.
Since we need to kill the dimension of additional marked points $z^+ \cup z^{++}$
and the dimension ($= k$) corresponding to $\R$-translations,
the virtual dimension $m$ of the Kuranishi neighborhood
(that is, $\dim \V - \dim \E$ of the Kuranishi neighborhood $(\V, \E, s, \psi)$) is
\begin{align*}
m &= \dim \mathring{X} + \ind D_{p_0} - 2 (\# z^+ + \# z^{++}) - k \\
&= (\dim \widetilde{X} - 2 (\# z^+ + \# z^{++})) + \ind D_{p_0} - 1.
\end{align*}

Next we check the relation of the virtual dimension of $p_0$ and
those of its parts.
Assume that we can construct $p_0$ from finite number of
holomorphic buildings $p_\kappa = (\Sigma_\kappa, z_\kappa, u_\kappa)$
and finite number of trivial cylinders by jointing pairs of limit circles to
joint circles and jointing pairs of marked points to nodal points.
(For example, let $\{p_\kappa\}$ be the restrictions of $p_0$ to
the irreducible components $\Sigma_\alpha$ which are not trivial cylinders.)
Let $l_{\text{trivial}}$ be the number of trivial cylinders, and let
$l_{\text{nodal}}$ and $l_{\text{joint}}$ be the number of new nodal points and
new joint circles respectively.
It is easy to check that
\begin{align*}
&(\dim \widetilde{X} - 2 (\# z^+ + \# z^{++})) \\
&= \sum_\kappa (\dim \widetilde{X}_\kappa - 2 (\# z_\kappa^+ + \# z_\kappa^{++}))
- 2 l_{\text{trivial}} + 2 l_{\text{nodal}} + 2 l_{\text{joint}}.
\end{align*}
The index of $D_{p_0}$ and those of $D_{p_\kappa}$ are related by
\[
\ind D_{p_0} = \sum_\kappa \ind D_{p_\kappa} + 2 l_{\text{trivial}}
- \sum_{\star} \dim \Ker A_{\gamma_\mu} - 2n l_{\text{nodal}}
\]
where the sum $\star$ is taken over new joint circles $\{S^1_\mu\}$
and each $\gamma_\mu$ is the periodic orbit on $S^1_\mu$.
The term $- \dim \Ker A_{\gamma_\mu}$ in the above equation is due to
the fact that the Sobolev space $\widetilde{W}_\delta^{1, p}(\Sigma_0, u_0^\ast T \hat Y)$
contains one vector space $\Ker A_{\gamma_\mu}$ for each joint circle
while the direct sum of the Sobolev spaces for $\{p_\kappa\}$ and limit circles
contains a pair of $\Ker A_{\gamma_\mu}$ for each pair of limit circles.
For simplicity, assume Morse condition.
Then the above equations imply
\[
m - \sum_\kappa m_\kappa
= \#\{\kappa\} - 1 + (2 - 2n) l_{\text{nodal}},
\]
where each $m_\kappa$ is the virtual dimension of $p_\kappa$.
For example, this equation implies that if we divide a holomoprhic building into two parts
by a gap of floor, then the virtual dimension of the entire holomorphic building
is larger than the sum of the virtual dimensions of the two by one.
Similarly, the virtual dimension of disjoint holomoprhic building is
larger than the sum of the virtual dimensions of its connected components.
These coincide with the expected relations indeed.
\end{rem}

Now for each $(a, b) \in \mathring{X}$, we construct an approximate solution
$u_{a, b} : \widetilde{P}_a \to (\overline{\R}_1 \sqcup \overline{\R}_2 \sqcup \dots \sqcup
\overline{\R}_k)/ \sim_{a, b} \times Y$
and a map $\Phi_{a, b} : u_{a, b}^\ast T \hat Y \to (\overline{\R}_1 \sqcup
\overline{\R}_2 \sqcup \dots \sqcup \overline{\R}_k)/ \sim_{a, b} \times Y$.
They will satisfy the following conditions:
\begin{itemize}
\item
$u_{a, b}|_{\Sigma_0 \setminus N_0} = u_0|_{\Sigma_0 \setminus N_0}$
\item
The restriction of $\Phi_{a, b}$ to the zero section coincides with $u_{a, b}$,
that is, $\Phi_{a, b}(z, 0) = u_{a, b}(z)$ for all $z \in \widetilde{P}_a$.
\item
The vertical differential of $\Phi_{a, b}$ at the zero section is the identity map of
$u_{a, b}^\ast T \hat Y$.
\item
The restriction of $\Phi_{a, b}$ on $u_0^\ast T \hat Y|_{\Sigma_0 \setminus N_0}$
does not depend on $(a, b) \in \mathring{X}$.
\end{itemize}

First we consider the neighborhood of $\nu$-th nodal point.
Let $\phi^\nu : B_\epsilon^{2n}(0) \to \R \times Y$ be a coordinate centered
at the image of the nodal point by $u_0$.
Define $v_0^\nu : D \cup D \to B_\epsilon^{2n}(0)$ by
\[
u_0|_{(D \cup D)_\nu}(x,y) = \phi^\nu (v_0^\nu (x,y)).
\]
For each $(a, b) \in \mathring{X}$, define a piecewise smooth map $v_{a, b}^\nu :
N_{a, b}^\nu = \{(x,y) \in D \times D; xy = \zeta_\nu \} \to B_\epsilon^{2n}(0)$ by
\[
v_{a, b}^\nu(x,y) = \begin{cases}
v_0^\nu(\frac{r-\sqrt{|\zeta_\nu|}}{1-\sqrt{|\zeta_\nu|}} e^{\sqrt{-1} \theta}, 0)
\text{ if } x = r e^{\sqrt{-1}\theta} \text{ and } r \geq \sqrt{|\zeta_\nu|}\\
v_0^\nu(0, \frac{r-\sqrt{|\zeta_\nu|}}{1-\sqrt{|\zeta_\nu|}} e^{\sqrt{-1} \theta})
\text{ if } y = r e^{\sqrt{-1}\theta} \text{ and } r \geq \sqrt{|\zeta_\nu|}
\end{cases}.
\]
Define piecewise smooth maps $u_{a, b} : N_{a, b}^\nu \to \hat Y$ and
$\Phi_{a, b} : N_{a, b}^\nu \times \R^{2n} \to \hat Y$ by
\begin{align*}
u_{a, b}(x,y) &= \phi^\nu (v_{a, b}^\nu (x,y))\\
\Phi_{a, b}(x, y, \xi) &= \phi^\nu (v_{a, b}^\nu (x,y) + \xi).
\end{align*}
We identify $N_{a, b}^\nu \times \R^{2n}$ and $u_{a, b}^\ast T \hat Y|_{N_{a, b}^\nu}$
by the differential of $\Phi_{a, b}$ at the zero section $N_{a, b}^\nu \times \{0 \}$,
and consider $\Phi_{a, b}$ as a map $u_{a, b}^\ast T \hat Y|_{N_{a, b}^\nu} \to \hat Y$.

Next we consider the neighborhood of $\mu$-th joint circle.
Define $b_\mu^{0, \mathrm{left}},\ab b_\mu^{0, \mathrm{right}},\ab b_\mu^0 \in \R$ by
\begin{align*}
u_0|_{[-1,\infty)_\mu \times S^1}(s,t)
&= (L_\mu s + b_\mu^{0, \mathrm{left}}, \gamma_\mu(t)) + o(1)\\
u_0|_{(-\infty, +1]_\mu \times S^1}(s,t)
&= (L_\mu s + b_\mu^{0, \mathrm{right}}, \gamma_\mu(t)) + o(1)\\
b_\mu^0 &= b_\mu^{0, \mathrm{left}} - b_\mu^{0, \mathrm{right}}.
\end{align*}
Let $\phi^\mu : B^{m_\mu}_\epsilon (0) \to P$ be a coordinate centered at
$\gamma_\mu \in P$ for each $\mu$.
We take a family of open embeddings $\psi^\mu_t : B^{m_\mu}_\epsilon (0) \times
B^{2n-1-m_\mu}(0) \to Y$ ($t \in S^1$) such that
$\psi^\mu_t(x,0) = \ev_t \phi^\mu(x)$ for all $x \in B^{m_\mu}_\epsilon (0)$
as in Section \ref{asymptotic estimates}.
Define families of open embeddings $\hat \psi^{\mu, \mathrm{left}}_{s,t},
\hat \psi^{\mu, \mathrm{right}}_{s,t} : \R \times B^{m_\mu}_\epsilon (0) \times
B^{2n-1-m_\mu}(0) \to \R \times Y$ ($(s,t) \in \R \times S^1$) by
\begin{align*}
\hat \psi^{\mu, \mathrm{left}}_{s,t} (\sigma, (x,y))
&= (L_\mu s + b_\mu^{0, \mathrm{left}} + \sigma, \psi^\mu_t(x,y))\\
\hat \psi^{\mu, \mathrm{right}}_{s,t} (\sigma, (x,y))
&= (L_\mu s + b_\mu^{0, \mathrm{right}} + \sigma, \psi^\mu_t(x,y)).
\end{align*}
Define smooth functions $v_0^{\mu, \mathrm{left}} : [-1, \infty] \times S^1
\to \R \times B^{m_\mu}_\epsilon (0) \times B^{2n-1-m_\mu}(0)$ and
$v_0^{\mu, \mathrm{right}} : [-\infty, +1] \times S^1 \to \R \times
B^{m_\mu}_\epsilon (0) \times B^{2n-1-m_\mu}(0)$ by
\begin{align*}
u_0|_{[-1,\infty]_\mu \times S_\mu^1} (s,t)
&= \hat \psi^{\mu, \mathrm{left}}_{s,t} (v_0^{\mu, \mathrm{left}} (s,t))\\
u_0|_{[-\infty,+1]_\mu \times S_\mu^1} (s,t)
&= \hat \psi^{\mu, \mathrm{right}}_{s,t} (v_0^{\mu, \mathrm{right}} (s,t)).
\end{align*}

For each $\mu$, fix a constant $0 < \kappa_\mu < \delta_{0, \mu}$.
($\delta_{0, \mu}$ is the minimal nonzero absolute value of eigenvalues of
$A_{\gamma_\mu}$.)
Recall that we have identified $\widetilde{N}_a^\mu$ with
$([-1, -\frac{1}{2} \log \rho_\mu] \cup [\frac{1}{2} \log \rho_\mu, 1]) \times S^1_\mu$
by (\ref{left coordinate near joint}) and (\ref{right coordinate near joint}).
For each $(a,b) \in \mathring{X}$ and $\mu$, define
$v_{a,b}^{\mu, \mathrm{left}} : [-1, -\frac{1}{2} \log \rho_\mu] \times S^1
\to \R \times B^{m_\mu}_\epsilon (0) \times B^{2n-1-m_\mu}(0)$ and
$v_{a,b}^{\mu, \mathrm{right}} : [\frac{1}{2} \log \rho_\mu, 1] \times S^1
\to \R \times B^{m_\mu}_\epsilon (0) \times B^{2n-1-m_\mu}(0)$ by
\begin{align*}
v_{a,b}^{\mu, \mathrm{left}} (s,t) &=\begin{cases}
v_0^{\mu, \mathrm{left}} (s, t) &\text{if } s \in [-1, 0]\\
v_0^{\mu, \mathrm{left}} \biggl( {\displaystyle  - \frac{1}{\kappa_\mu} \log
\biggl(\frac{e^{-\kappa_\mu s} - \rho_\mu^{\kappa_\mu/2}}
{1 - \rho_\mu^{\kappa_\mu/2}} \biggr), t }\biggr)
&\text{if } s \in [0, -\frac{1}{2} \log \rho_\mu]
\end{cases},\\
v_{a,b}^{\mu, \mathrm{right}} (s,t) &=\begin{cases}
v_0^{\mu, \mathrm{right}} (s, t) &\text{if } s \in [0, 1]\\
v_0^{\mu, \mathrm{right}} \biggl( {\displaystyle \frac{1}{\kappa_\mu} \log
\biggl(\frac{e^{\kappa_\mu s} - \rho_\mu^{\kappa_\mu/2}}
{1 - \rho_\mu^{\kappa_\mu/2}}\biggr), t }\biggr)
&\text {if } s \in [\frac{1}{2} \log \rho_\mu, 0]
\end{cases}.
\end{align*}
Then piecewise smooth maps $u_{a, b} : ([-1, -\frac{1}{2} \log \rho_\mu] \cup
[\frac{1}{2} \log \rho_\mu, 1]) \times S^1 \to
(\overline{\R}_i \cup \overline{\R}_{i + 1}) / \sim_{a, b} \times Y$ and
$\Phi_{a, b} : ([-1, -\frac{1}{2} \log \rho_\mu] \cup [\frac{1}{2} \log \rho_\mu, 1])
\times S^1 \times \R^{2n} \to (\overline{\R}_i \cup \overline{\R}_{i + 1}) / \sim_{a, b}
\times Y$ are defined by
\begin{equation*}
u_{a, b}(s, t)
 = \begin{cases}
o_{\frac{1}{2}\chi(s) (b_\mu - b_\mu^0)} \circ
\hat \psi^{\mu, \mathrm{left}}_{s, t}
(v_{a,b}^{\mu, \mathrm{left}}(s,t)) \in \overline{\R}_i \times Y\\
\hspace{180pt}
\text{if } s \in [-1, -\frac{1}{2} \log \rho_\mu] \\
o_{-\frac{1}{2}\chi(-s) (b_\mu - b_\mu^0)} \circ
\hat \psi^{\mu, \mathrm{right}}_{s, t}
(v_{a,b}^{\mu, \mathrm{right}}(s,t)) \in \overline{\R}_{i + 1} \times Y\\
\hspace{180pt}
\text{if } s \in [\frac{1}{2} \log \rho_\mu, 1]
\end{cases}
\end{equation*}
and
\begin{equation*}
\Phi_{a, b}(s, t, \xi)
= \begin{cases}
o_{\frac{1}{2}\chi(s) (b_\mu - b_\mu^0)} \circ
\hat \psi^{\mu, \mathrm{left}}_{s, t}
(v_{a,b}^{\mu, \mathrm{left}}(s,t) + \xi) \in \overline{\R}_i \times Y\\
\hspace{170pt}
\text{if } s \in [-1, -\frac{1}{2} \log \rho_\mu] \\
o_{-\frac{1}{2}\chi(-s) (b_\mu - b_\mu^0)} \circ
\hat \psi^{\mu, \mathrm{right}}_{s, t}
(v_{a,b}^{\mu, \mathrm{right}}(s,t) + \xi) \in \overline{\R}_{i + 1} \times Y\\
\hspace{170pt}\text{if } s \in [\frac{1}{2} \log \rho_\mu, 1],
\end{cases}
\end{equation*}
where $o_c : \overline{\R} \times Y \to \overline{\R} \times Y$ is the translation map
of the $\overline{\R}$-factor defined by $o_c(\sigma, y) = (\sigma + c, y)$,
and $\chi : \R \to \R_{\geq 0}$ is a smooth function such that
$\chi|_{(-\infty, 0]} = 0$ and $\chi|_{[0, \infty)} = 1$.
We identify $([-1, -\frac{1}{2} \log \rho_\mu] \cup [\frac{1}{2} \log \rho_\mu, 1])
\times S^1 \times \R^{2n}$ and
$u_{a, b}^\ast T \hat Y|_{([-1, -\frac{1}{2} \log \rho_\mu] \cup
[\frac{1}{2} \log \rho_\mu, 1]) \times S^1}$
by the differential of $\Phi_{a, b}$ at the zero section $([-1, -\frac{1}{2} \log \rho_\mu]
\cup [\frac{1}{2} \log \rho_\mu, 1]) \times S^1 \times \{0 \}$,
and consider $\Phi_{a, b}$ as a map
$u_{a, b}^\ast T \hat Y|_{([-1, -\frac{1}{2} \log \rho_\mu] \cup
[\frac{1}{2} \log \rho_\mu, 1]) \times S^1}
\to (\overline{\R}_i \cup \overline{\R}_{i + 1}) / \sim_{a, b} \times Y$.

Next we consider the neighborhood of each limit circle.
Since this region does not change by $(a, b) \in \mathring{X}$,
we can use $u_{a, b} = u_0$ as an
approximate solution.
Assume this circle is $+ \infty$-limit circle $S_{+ \infty_i}^1$.
(The case of $-\infty$-limit circle is similar.)
Let $\phi^{+\infty_i} : B^{m_{+\infty_i}}_\epsilon (0) \to P$ be a coordinate centered at
$\gamma_{+\infty_i} \in P$, and take a family of open embeddings
$\psi^{+\infty_i}_t : B^{m_{+\infty_i}}_\epsilon (0) \times B^{2n-1-m_{+\infty_i}}(0) \to Y$
($t \in S^1$) such that $\psi^{+\infty_i}_t(x,0) = \ev_t \phi^{+\infty_i}(x)$ for all
$x \in B^{m_{+\infty_i}}_\epsilon (0)$ as in the previous case.
We define a smooth map $v_0^{+\infty_i} : [0, \infty] \times S^1 \to \R \times
B^{m_{+\infty_i}}_\epsilon (0) \times B^{2n-1-m_{+\infty_i}}(0)$ by
\[
u_0|_{[0,\infty] \times S^1} (s,t)
= (1 \times \psi^{+\infty_i}_t) (v_0^{+\infty_i} (s,t)).
\]
Then a smooth map $\Phi_{a, b} : [0, \infty] \times S^1 \times \R^{2n} \to \hat Y$
is defined by
\[
\Phi_{a, b}(s, t, \xi) = (1 \times \psi^{+\infty_i}_t) (v_0^{+\infty_i} (s,t) + \xi).
\]
(This does not depend on $(a, b) \in \mathring{X}$.)

Finally, we consider the rest $\Sigma_0 \setminus N_0$. Since $u_{a, b}|_{\partial N_0}
= u_0|_{\partial N_0}$, we can define a piecewise smooth map
$u_{a, b} : \widetilde{P}_a \to (\overline{\R}_1 \sqcup \overline{\R}_2 \sqcup \dots
\sqcup \overline{\R}_k)/ \sim_{a, b} \times Y$
by $u_{a, b}|_{\Sigma_0 \setminus N_0} = u_0|_{\Sigma_0 \setminus N_0}$.
Note that the restriction of $\Phi_{a, b}$ to $u_0^\ast T \hat Y|_{\partial N_0}$
does not depend on $(a, b)$.
Therefore, we can take a smooth extension
$\Phi : u_0^\ast T \hat Y|_{\Sigma_0 \setminus N_0} \to (\overline{\R}_1 \sqcup
\overline{\R}_2 \sqcup \dots \sqcup \overline{\R}_k)/ \sim_{a, b} \times Y$
which is independent of $(a, b) \in \mathring{X}$ and satisfies the desired conditions,
that is,
\begin{itemize}
\item
the restriction of $\Phi$ to the zero section coincides with $u_0$, and
\item
the vertical differential of $\Phi$ at the zero section is the identity map of
$u_0^\ast T \hat Y$.
\end{itemize}

We will give a differentiable structure to a neighborhood
\begin{align*}
\hat V \subset \bigcup_{(a, b) \in \mathring{X}} \{(a, b)\} \times
\{&(\xi, h) \in \widetilde{W}_\delta^{1, p}(\widetilde{P}_a; u_{a, b}^\ast T \hat Y)
\times E^0;\\
&d(\Phi_{a, b}(\xi)) + J d(\Phi_{a, b}(\xi)) j  + h_{a, b}(z, \Phi_{a, b}(\xi)) = 0\}
\end{align*}
of $(0, b^0, 0,0)$ later, where $h_{a, b}$ is the restriction of $\lambda(h)$ to
$\widehat{P}_a \times Y$.
Then $G_0$ acts on $\hat V$ smoothly, and
a $G_0$-equivariant section $s^0 : \hat V \to \R^k \oplus
\bigoplus_{z_\beta^{++}} \R_i^2$ defined by
\[
s^0(a, b, \xi, h) = (\sigma_i \circ \Phi_{a, b}(\xi)(\widetilde{R}_i(a)),
p' \circ \Phi_{a, b}(\xi)(Z_\beta^{++}(a)))
\]
is a smooth submersion, where each
$\sigma_i \circ \Phi_{a, b}(\xi)(\widetilde{R}_i(a))$ is the abbreviation of
\[
\frac{1}{m_i} \sum_{l=1}^{m_i} \sigma_i \circ \Phi_{a, b}(\xi)(\widetilde{R}_{i, l}(a)),
\]
and $p'$ is a smooth submersion from a neighborhood of $S'$ to $\R^2$ such that
$S' = \{p' = 0\}$.
Let $V = \{s^0 = 0\} \subset \hat V$ be the zero set,
and consider the finite dimensional vector space
$E = E^0 \oplus \bigoplus_{z_\alpha^+} \R_\alpha^2$ as a trivial vector bundled on $V$.
Define a smooth section $s : V \to E$ by
\[
s(a, b, \xi, h) = (h, p \circ \pi_Y \circ \Phi_{a, b}(\xi)(Z_\alpha^+(a))),
\]
where $p$ is a smooth submersion from a neighborhood of $S \subset Y$ to $\R^2$
such that $S = \{p= 0\}$.

Since the zero set of $s$ consists of holomorphic buildings,
we can define a continuous map $\psi : \{s = 0\}/G_0 \to \widehat{\M}^0(Y, \lambda, J)$.
Finally we will prove that this map is a homeomorphism onto a neighborhood of $p$.

Now we start to define a differentiable structure of $\hat V$.
To do so, we express this set as a zero set of a Fredholm map between
Banach spaces.
To define a Fredholm map, first we rewrite the equation of $(\xi, h)$.

Note that the equation
\begin{equation}
d(\Phi_{a, b}(\xi))(z) + J(\Phi_{a, b}(\xi)(z)) d(\Phi_{a, b}(\xi))(z) j_z
+ h_{a, b}(z, \Phi_{a, b}(\xi)(z)) = 0 \label{xi eq}
\end{equation}
is equivalent to the equation of $J$-holomorphic curve on $N_0$ since $h_{a, b}$
vanishes on $N_0 \times Y$.

On $\{x \in D ; |x| \geq \sqrt{|\zeta_\nu|}\} \subset \{(x, y) \in D \times D;
x y = \zeta_\nu\}$ or
$\{y \in D ; |y| \geq \sqrt{|\zeta_\nu|}\} \subset \{(x, y) \in D \times D;
x y = \zeta_\nu\}$, $\Phi_{a, b}(\xi)$ is $J$-holomorphic if and only if
\begin{equation}
\partial_r (v_{a, b}^\nu + \xi) + \frac{1}{r}\widetilde{J}^\nu(v_{a, b}^\nu + \xi)
\partial_\theta (v_{a, b}^\nu + \xi) = 0, \label{node eq}
\end{equation}
where $\widetilde{J}^\nu = (\phi^\nu)^\ast J$ is the pull back of $J$,
and $(r, \theta)$ is the polar coordinate of $x$ or $y$ respectively.

On $[-1, -\frac{1}{2} \log \rho_\mu] \times S^1$,
$\Phi_{a, b}(\xi)$ is $J$-holomoprhic if and only if
\begin{multline*}
\hat\psi_\ast \partial_s (v_{a,b}^{\mu, \mathrm{left}} + \xi)
+ \Bigl(L_\mu + \frac{1}{2} \chi'(s) (b_\mu - b_\mu^0)\Bigr) \partial_\sigma \\
+\Bigl(\frac{1}{2} \varphi_\mu \chi'(s)
+ J(\hat\psi(v_{a,b}^{\mu, \mathrm{left}} + \xi))\Bigr)
((\partial_t \hat\psi) (v_{a,b}^{\mu, \mathrm{left}} + \xi)
+ \hat\psi_\ast \partial_t (v_{a,b}^{\mu, \mathrm{left}} + \xi))= 0
\end{multline*}
since $-j_{\varphi_\mu} \partial_t = \partial_s + \frac{1}{2} \varphi_\mu \chi'(s) \partial_t$.
This can be written as
\begin{multline}
\partial_s (v_{a,b}^{\mu, \mathrm{left}} + \xi)
+ \widetilde{J}^\mu_t(v_{a,b}^{\mu, \mathrm{left}} + \xi)
\partial_t (v_{a,b}^{\mu, \mathrm{left}} + \xi)
+ f^\mu_t(v_{a,b}^{\mu, \mathrm{left}} + \xi)\\
+ \frac{1}{2} (b_\mu -b_\mu^0) \chi'(s) \partial_\sigma
+ \frac{1}{2} \varphi_\mu \chi'(s) (g^\mu_t(v_{a,b}^{\mu, \mathrm{left}} + \xi)
+ \partial_t(v_{a,b}^{\mu, \mathrm{left}} + \xi)) =0, \label{+ joint circle eq}
\end{multline}
where $\widetilde{J}^\mu_t = (1 \times \psi^\mu_t)^\ast J$ and
\begin{align*}
f^\mu_t(\sigma, y) &= (1 \times \psi^\mu_t)_\ast^{-1} J(\psi^\mu_t(y))
(\partial_t \psi^\mu_t(y) - L R_\lambda(y))\\
g^\mu_t(\sigma, y) &= (1 \times \psi^\mu_t)_\ast^{-1} \partial_t \psi_t(y).
\end{align*}
In particular, on $[0, -\frac{1}{2} \log \rho_\mu] \times S^1
\subset [-1, \frac{1}{2} \log \rho_\mu] \times S^1$,
this equation can be written as
\[
\partial_s (v_{a,b}^{\mu, \mathrm{left}} + \xi)
+ \widetilde{J}^\mu_t(v_{a,b}^{\mu, \mathrm{left}} + \xi)
\partial_t (v_{a,b}^{\mu, \mathrm{left}} + \xi)
+ f^\mu_t(v_{a,b}^{\mu, \mathrm{left}} + \xi) = 0.
\]
We note that $f^\mu_t : B^{m_\mu}_\epsilon(0) \times B^{2n-1-m_\mu}(0) \to \R^{2n}$
satisfies $f^\mu_t|_{B^{m_\mu}_\epsilon(0) \times \{0\}} \equiv 0$.

Similarly, on $[\frac{1}{2} \log \rho_\mu, 1] \times S^1$, $\Phi_{a, b}(\xi)$ is
$J$-holomorphic if and only if
\begin{align}
&\partial_s (v_{a,b}^{\mu, \mathrm{right}} \! + \xi)
+ \widetilde{J}^\mu_t(v_{a,b}^{\mu, \mathrm{right}} \! + \xi)
\partial_t (v_{a,b}^{\mu, \mathrm{right}} \! + \xi)
+ f^\mu_t(v_{a,b}^{\mu, \mathrm{right}} \! + \xi)
\notag\\
&+ \frac{1}{2} (b_\mu -b_\mu^0) \chi'(-s) \partial_\sigma
+ \frac{1}{2} \varphi_\mu \chi'(-s) (g^\mu_t(v_{a,b}^{\mu, \mathrm{right}} \! + \xi)
+ \partial_t(v_{a,b}^{\mu, \mathrm{right}} \! + \xi)) \notag\\
& =0. \label{- joint circle eq}
\end{align}

On the neighborhood $[0, \infty] \times S^1$ of the limit circle $S_{+\infty_i}^1$
or on the neighborhood $[-\infty, 0] \times S^1$ of the limit circle $S_{-\infty_i}^1$,
$\Phi_{a, b}(\xi)$ is $J$-holomorphic if and only if
\begin{equation}
\partial_s (v_0^{\pm\infty_i} + \xi)
+ \widetilde{J}^{\pm\infty_i}_t(v_0^{\pm\infty_i} + \xi) \partial_t (v_0^{\pm\infty_i} + \xi)
+ f^{\pm\infty_i}_t (\pi_Y (v_0^{\pm\infty_i} + \xi)) = 0, \label{limit circle eq}
\end{equation}
where $f^{\pm\infty_i}_t : B^m_\epsilon(0) \times B^{2n-1-m}(0) \to \R^{2n}$ ($t\in S^1$)
is a smooth function which satisfies $f^\mu_t|_{B^m_\epsilon(0) \times \{0\}} \equiv 0$.

Now we define a Fredholm map
\begin{align*}
&F^{(a,b)} : \widetilde{W}^{1,p}_\delta( \widetilde{P}_a, u_{a, b}^\ast T \hat Y) \oplus E^0\\
&\to L^p(\Sigma_0 \setminus N_0, \Wedge^{0,1} T^\ast \Sigma_0 \otimes_\C
u_0^\ast T \hat Y)\\
&\quad \oplus \bigoplus_\nu (L^p (\{ x \in D; |x| \geq \sqrt{|\zeta_\nu|}\}, \R^{2n})
\oplus L^p (\{ y \in D; |y| \geq \sqrt{|\zeta_\nu|}\}, \R^{2n}))\\
&\quad \oplus \bigoplus_\mu (L^p_\delta ([-1, -{\textstyle\frac{1}{2}}\log \rho_\mu]
\times S^1, \R^{2n})
\oplus L^p_\delta ([{\textstyle\frac{1}{2}} \log \rho_\mu, +1] \times S^1, \R^{2n}))\\
&\quad \oplus \bigoplus_{+\infty_i} L^p_\delta([0, \infty] \times S^1, \R^{2n})
\oplus \bigoplus_{-\infty_i} L^p_\delta([-\infty, 0] \times S^1, \R^{2n})
\end{align*}
by the left hand sides of the above equations (\ref{xi eq}), (\ref{node eq}),
(\ref{+ joint circle eq}), (\ref{- joint circle eq}) and (\ref{limit circle eq}), that is,
its $L^p(\Sigma_0 \setminus N_0)$-component is defined by
\[
d(\Phi(\xi))(z) + J(\Phi(\xi)(z)) d(\Phi(\xi))(z) j_z + h_{a, b}(z, \Phi(\xi)(z)),
\]
its $L^p (\{ x \in D; |x| \geq \sqrt{|\zeta_\nu|}\})$-component is defined by
\[
\partial_r (v_{a, b}^\nu + \xi) + \frac{1}{r}\widetilde{J}^\nu(v_{a, b}^\nu + \xi)
\partial_\theta (v_{a, b}^\nu + \xi),
\]
and so on.
We abbreviate the range of this Fredholm map as
$L_\delta^p(\widetilde{P}_a, \Wedge^{0, 1} T^\ast \widetilde{P}_a \otimes
u_{a, b}^\ast T \hat Y)$.

We also define a Fredholm map
\[
F^{(a,b) +} : \widetilde{W}^{1,p}_\delta( \widetilde{P}_a, u_{a, b}^\ast T \hat Y) \oplus E^0
\to L_\delta^p(\widetilde{P}_a, \Wedge^{0, 1} T^\ast \widetilde{P}_a \otimes
u_{a, b}^\ast T \hat Y) \oplus \Ker DF^{(0, b^0)}_{(0, 0)}
\]
by
\[
F^{(a,b) +} (\xi, h) = \bigl(F^{(a,b)}(\xi,h),
\sum_i (\langle \xi, \xi_i \rangle_{L^2(\Sigma_0 \setminus N_0)}
+ \langle h, h_i \rangle_{E^0}) \cdot x_i\bigr),
\]
where $\{ x_i = (\xi_i, h_i) \}$ is a orthonormal basis of $\Ker DF^{(0, b^0)}_{(0, 0)}$
with the inner product given by
\[
\langle (\xi, h), (\xi', h') \rangle = \langle \xi, \xi' \rangle_{L^2(\Sigma_0 \setminus N_0)}
+ \langle h, h' \rangle_{E^0}
\]
for some inner product of $E^0$.

In order to apply the implicit function theorem to $F^{(a, b)}$,
or apply inverse function theorem to $F^{(a, b)+}$, we need to check their properties.
First we need to show that $F^{(a, b)}(0, 0)$ is small for any $(a, b) \in \mathring{X}$
sufficiently close to $(0, b^0) \in \mathring{X}$.
(This is equivalent to say that $u_{(a, b)}$ is close to the solution.)
Note that $F^{(a, b)}(0, 0)$ is zero on $\Sigma_0 \setminus N_0$ and
the neighborhoods of limit circles since these regions are independent of
$(a, b) \in \mathring{X}$.
Recall that $\delta_{0, \mu} > 0$ and $\delta_{0, \infty_i}$ are the minimal nonzero
absolute value of eigenvalues of $A_{\gamma_\mu}$ and $A_{\gamma_{\pm\infty_i}}$
respectively.
Assume that a sequence of positive constants
$\delta = ((\delta_\mu)_\mu, (\delta_{\pm\infty_i})_{\pm\infty_i})$ satisfies
$\delta_{\mu} < \delta_{0, \mu}$ and $\delta_{\pm\infty_i} < \delta_{\pm\infty_i}$.
We abbreviate this condition by $\delta < \delta_0$.
For such a sequence of positive constant $\delta$ and a constant $p > 2$,
we use the $L^p_\delta$-norm or $W^{1, p}_\delta$-norm on $\widetilde{P}_a$ as
a deformation of the curve $\Sigma_0$ explained before.
Let $\delta'_0 = ((\delta'_{0, \mu})_\mu, (\delta'_{0, \pm\infty_i})_{\pm\infty_i})$ be
an arbitrary sequence of positive constants such that $\delta < \delta'_0 < \delta_0$.
\begin{lem}\label{estimates of F(0)}
For any $0 < \delta < \delta'_0 < \delta_0$
and $p > 2$, there exists a constant $C>0$ such that
for any $(a,b) \in \mathring{X}$ sufficiently close to $(0, b^0)$,
the following inequalities hold true.
\begin{align*}
||F^{(a,b)}(0,0)|_{[-1, 0] \times S_\mu^1}||_{L^p}
&\leq C (|\varphi_\mu| + |b_\mu - b_\mu^0|)\\
||F^{(a,b)}(0,0)|_{[0, -\frac{1}{2} \log \rho_\mu] \times S_\mu^1}||_{L^p_{\delta_\mu}}
&\leq C \rho_\mu^{\min(\kappa_\mu, \delta'_{0, \mu} - \delta_\mu)/2}
(-\log \rho_\mu)^{1/p} \\
||F^{(a,b)}(0,0)|_{\{(x, y) \in N^\nu_{(a, b)}; |x| \geq \sqrt{|\zeta_\nu|} \}}||_{L^p}
&\leq C |\zeta_\nu|^{1/p}
\end{align*}
\end{lem}
\begin{proof}
First we estimate the $L^p$-norm of $F^{(a,b)}(0,0)|_{[-1, 0] \times S_\mu^1}$.
The equation
\[
0 = F^{(0, b^0)}(0,0)|_{[-1, 0] \times S_\mu^1}
= \partial_s v_0^{\mu, \text{left}}
+ \widetilde{J}^\mu_t(v_0^{\mu, \text{left}}) \partial_t v_0^{\mu, \text{left}}
+ f^\mu_t(v_0^{\mu, \text{left}})
\]
implies
\begin{align*}
&F^{(a, b)}(0, 0)|_{[-1, 0] \times S_\mu^1}\\
&= F^{(a, b)}(0, 0)|_{[-1, 0] \times S_\mu^1} - F^{(0, b^0)}(0,0)|_{[-1, 0] \times S_\mu^1}\\
&= \frac{1}{2} (b_\mu -b_\mu^0) \chi'(s) \partial_\sigma
+ \frac{1}{2} \varphi_\mu \chi'(s) (g^\mu_t(v_{a,b}^{\mu, \mathrm{left}})
+ \partial_t v_{a,b}^{\mu, \mathrm{left}}).
\end{align*}
The first inequality is clear from this equation.

Next we estimate the $L^p_\delta$-norm of
$F^{(a,b)}(0,0)|_{[0, -\frac{1}{2} \log \rho_\mu] \times S^1}$.
Since $\kappa_\mu < \delta_{0, \mu}$, we may assume that
$\delta'_{0, \mu} > \kappa_\mu$.
We omit the subscript $\mu$ of $\kappa_\mu$, $\rho_\mu$ and so on.
First we note that
\begin{align}
&F^{(a,b)}(0,0)|_{[0, -\frac{1}{2} \log \rho] \times S^1}(s, t) \notag\\
&= \biggl(1 + \frac{\rho^{\kappa/2}}{e^{-\kappa s} - \rho^{\kappa/2}}\biggr)
\partial_s v_0^{\mu, \text{left}}(\tilde s, t) \notag \\
&\quad \ + \widetilde{J}^\mu_t(v_0^{\mu, \text{left}}(\tilde s, t))
\partial_t v_0^{\mu, \text{left}}(\tilde s, t) 
+ f_t^\mu(v_0^{\mu, \text{left}}(\tilde s, t)),
\label{F joint circle}
\end{align}
where
\[
\tilde s = -\frac{1}{\kappa} \log
\biggl(\frac{e^{-\kappa s} - \rho^{\kappa/2}}
{1 - \rho^{\kappa/2}}\biggr).
\]
Substituting
\[
\biggl( -\frac{1}{\kappa} \log
\biggl(\frac{e^{-\kappa s} - \rho^{\kappa/2}}
{1 - \rho^{\kappa/2}}\biggr), t \biggr)
\]
for $(s, t)$ in the equation
\[
0 = F^{(0,0)}(0,0)|_{[0,\infty) \times S^1}
= \partial_s v_0^{\mu, \text{left}}
+ \widetilde{J}^\mu_t(v_0^{\mu, \text{left}}) \partial_t v_0^{\mu, \text{left}}
+ f_t^\mu(v_0^{\mu, \text{left}}),
\]
and subtracting it from (\ref{F joint circle}),
we obtain
\begin{align*}
&F^{(a,b)}(0,0)|_{[0, -\frac{1}{2} \log \rho] \times S^1}\\
&= \frac{\rho^{\kappa/2}}{e^{-\kappa s} - \rho^{\kappa/2}}
(\partial_s v_0^{\mu, \text{left}})
\biggl( -\frac{1}{\kappa} \log
\biggl(\frac{e^{-\kappa s} - \rho^{\kappa/2}}
{1 - \rho^{\kappa/2}}\biggr), t \biggr).
\end{align*}
Recall that Proposition \ref{second annulus} implies $|v_0^{\mu, \text{left}}(s,t)|,
|\partial_s v_0^{\mu, \text{left}}(s,t)| \lesssim e^{-\delta'_0 s}$.
Hence
\begin{align*}
&\int_0^{-\frac{1}{2} \log \rho}
\biggl( \frac{\rho^{\kappa/2}}{e^{-\kappa s} - \rho^{\kappa/2}}
\biggl|(\partial_s v_0^{\mu, \text{left}})
\biggl( -\frac{1}{\kappa} \log
\biggl(\frac{e^{-\kappa s} - \rho^{\kappa/2}}
{1 - \rho^{\kappa/2}}\biggr), t \biggr)\biggr| e^{\delta s} \biggr)^p ds\\
&\lesssim \int_0^{-\frac{1}{2} \log \rho}
\biggl( \frac{\rho^{\kappa/2}}{e^{-\kappa s} - \rho^{\kappa/2}}
\biggl(\frac{e^{-\kappa s} - \rho^{\kappa/2}}
{1 - \rho^{\kappa/2}}\biggr)^{\! \delta'_0/\kappa}
e^{\delta s}\biggr)^p ds\\
& = \frac{\rho^{p\kappa/2}}{(1 - \rho^{\kappa/2})^{p\delta'_0/\kappa}}
\int_0^{-\frac{1}{2} \log \rho}
\bigl((e^{-\kappa s} - \rho^{\kappa/2})^{\delta'_0/\kappa - 1} e^{\delta s}\bigr)^p ds\\
& \leq \frac{\rho^{p\kappa/2}}{(1 - \rho^{\kappa/2})^{p\delta'_0/\kappa}}
\int_0^{-\frac{1}{2} \log \rho} e^{-p(\delta'_0 - \delta - \kappa) s} ds\\
& \lesssim \rho^{p\min(\kappa, \delta'_0 - \delta)/2} (-\log \rho).
\end{align*}
This is the proof of the second inequality.

Finally we estimate the $L^p$-norm of
$F^{(a,b)}(0,0)|_{\{x \in D; |x| \geq \sqrt{|\zeta_\nu|}\}}$,
where we denote a point
$(x, y) \in \{(x, y) \in N^\nu_{(a, b)}; |x| \geq \sqrt{|\zeta_\nu|} \}$
by $x \in \{x \in D; |x| \geq \sqrt{|\zeta_\nu|}\}$.
We abbreviate $\zeta_\nu$ to $\zeta$ and define $\rho = \sqrt{|\zeta|}$.
First note that
\begin{align}
&F^{(a,b)}(0,0)|_{\{x \in D; |x| \geq \rho\}}(re^{\sqrt{-1} \theta}) \notag\\
&= \frac{1}{1-\rho} (\partial_r v_0^\nu)
\Bigl(\frac{r-\rho}{1-\rho} e^{\sqrt{-1} \theta}\Bigr) \notag\\
&\quad + \frac{1}{r}\widetilde{J}^\nu
\Bigl(v_0^\nu\Bigl(\frac{r-\rho}{1-\rho} e^{\sqrt{-1} \theta}\Bigr) \Bigr)
\partial_\theta
\Bigl(v_0^\nu\Bigl(\frac{r-\rho}{1-\rho} e^{\sqrt{-1} \theta}\Bigr) \Bigr).
\label{F node}
\end{align}
We also note that $F^{(0,b^0)}(0,0)|_{\{ x \in D\}} = 0$ implies
\begin{equation}
\widetilde{J}^\nu(v_0^\nu(re^{\sqrt{-1} \theta}))
\partial_\theta v_0^\nu(re^{\sqrt{-1} \theta})
= - r \partial_r v_0^\nu(re^{\sqrt{-1} \theta}).
\label{node del theta}
\end{equation}
Substituting $((r - \rho) / (1 - \rho), \theta)$ for
$(r, \theta)$ in (\ref{node del theta}),
and substitute it into (\ref{F node}),
we obtain
\[
F^{(a,b)}(0,0)|_{\{x \in D; |x| \geq \rho\}}
= \frac{\rho}{(1-\rho) r} (\partial_r v_0^\nu)
\Bigl(\frac{r-\rho}{1-\rho} e^{\sqrt{-1} \theta}\Bigr).
\]
Since $|\partial_r v_0^\nu|$ is bounded on $\{ x \in D\}$, this implies
\[
||F^{(a,b)}(0,0)|_{\{x \in D; |x| \geq \rho\}}||_{L^p}
\lesssim \frac{\rho}{1 - \rho}
\biggl(\int_{\rho}^1 r^{-p} rdr \biggr)^{1/p}
\lesssim \rho^{2 / p} = |\zeta|^{1/p}
\]
\end{proof}

Next we need to prove the differential $D F^{(a,b) +}_{(\xi,h)}$ is uniformly invertible for
any $(a,b) \in \mathring{X}$ sufficiently close to $(0, b^0)$ and any
$(\xi,h) \in \widetilde{W}^{1,p}_\delta(\widetilde{P}_a, u_{a, b}^\ast T \hat Y) \oplus E^0$
sufficiently close to $(0,0)$.
Since the assumption of the surjectivity of $D_{p_0}^+$ implies that
$D F^{(0,0) +}_{(0,0)}$ is invertible, the case of $(\xi, h) = (0,0)$ is
Lemma \ref{linearized gluing lemma} in Section \ref{linearized gluing}.
The general case is a consequence of the following lemma,
which can be proved easily by direct calculations.
\begin{lem}\label{D^2F}
For any $\delta < \delta_0$, there exists a constant $C>0$ such that 
for any $(a,b) \in \mathring{X}$ sufficiently close to $(0, b^0)$ and any
$(\xi,h) \in W^{1,p}_\delta( \widetilde{P}_a, (u_{a, b}^\ast T \hat Y)) \oplus E^0$
sufficiently close to $(0,0)$, the following inequalities hold true.
{\belowdisplayskip=0pt
\begin{multline*}
||DF^{(a,b) +}_{(\xi, h)} (\hat \xi, \hat h) - DF^{(a,b) +}_{(0,0)}
(\hat \xi, \hat h)||_{L^p_\delta([-1, -\frac{1}{2} \log \rho_\mu] \times S^1)}\\
\leq C (||\xi||_\infty ||\hat \xi||_{\widetilde{W}^{1,p}_\delta}
+ ||\xi||_{\widetilde{W}^{1,p}_\delta} ||\hat \xi||_\infty)
\end{multline*}
\begin{multline*}
||DF^{(a,b) +}_{(\xi, h)} (\hat \xi, \hat h) - DF^{(a,b) +}_{(0,0)} (\hat \xi, \hat h)||_{L^p
(\{ x \in D; |x| \geq \sqrt{|\zeta^\nu|}\})}\\
\leq C (||\xi||_\infty ||\hat \xi||_{W^{1,p}} + ||\xi||_{W^{1,p}} ||\hat \xi||_\infty)
\end{multline*}
}
\begin{multline*}
||DF^{(a,b) +}_{(\xi, h)} (\hat \xi, \hat h) - DF^{(a,b) +}_{(0,0)} (\hat \xi, \hat h)||_{L^p
(\Sigma_0 \setminus N_0)}\\
\leq C (||\xi||_\infty (||\hat \xi||_{W^{1,p}} + |\hat h|_{E^0})
+ (||\xi||_{W^{1,p}} + |h|_{E^0}) ||\hat \xi||_\infty)
\end{multline*}
\end{lem}

Therefore by the inverse function theorem, there exists some $\epsilon>0$ and $C>0$
such that for any $(a,b) \in \mathring{X}$ sufficiently close to $(0,b^0)$,
there exists a smooth map
\[
\phi^{a,b} : \Ker DF^{(0, b^0)}_{(0, 0)} \supset B_\epsilon(0) \to B_C(0) \subset
\widetilde{W}^{1,p}_\delta (\widetilde{P}_a, u_{a, b}^\ast T \hat Y) \oplus E^0
\]
such that for any $(\xi, h) \in B_C(0)$ and $x \in B_\epsilon(0)$,
\begin{equation}
F^{(a,b) +}(\xi, h) = (0, x) \text{ if and only if } (\xi, h) = \phi^{a,b}(x).
\label{def of phi}
\end{equation}
Note that $\Ker DF^{(0, b^0)}_{(0, 0)}$ does not depend on $p$ or $\delta$.
Although $\epsilon > 0$ may depend on $p$ and $\delta$ since so do the estimates,
$\phi^{a, b}$ does not depend on $p$ or $\delta$ on the intersection of the domains
since $\phi^{a, b}$ is defined by (\ref{def of phi}).

Shrinking $\mathring{X}$, we define $\hat V = \mathring{X} \times B_\epsilon(0)$ and
regard this space as a subspace of
\[
\bigcup_{(a, b) \in \mathring{X}} \{(a, b)\} \times
C^\infty(\widetilde{P}_a, (\overline{\R}_1 \cup
\overline{\R}_2 \cup \dots \cup \overline{\R}_k)
/\!\! \sim_{a, b} \! \times Y) \times E^0
\]
by
\[
(a, b, x) \mapsto (a, b, \Phi_{a, b}(\xi_x), h_x)
\]
where $(\xi_x, h_x) = (\xi_{(a, b, x)}, h_{(a, b, x)}) = \phi^{a, b}(x)$.

Define a map $s^0 : \hat V \to \R^k \oplus \bigoplus_{z_\beta^{++}} \R^2$ by
\begin{equation}
s^0(a, b, x) = (\sigma \circ \Phi_{a, b}(\xi_x)(\widetilde{R}_i(a)),
p' \circ \Phi_{a, b}(\xi_x)(Z_\beta^{++}(a))) \label{s^0}
\end{equation}
as we have already explained.
We will prove in Section \ref{smoothness} that
if we give a nice differentiable structure to the space $\mathring{X}$,
and give the product smooth structure to $\hat V = \mathring{X} \times B_\epsilon(0)$,
then
\begin{align*}
\hat V &\inj \mathring{X} \times C^l(\Sigma_0 \setminus N_0,
(\R_1 \cup \R_2 \cup \dots \cup \R_k) \times Y) \times E^0\\
(a, b, x) &\mapsto (a, b, \Phi_{a, b}(\xi_x)|_{\Sigma_0 \setminus N_0}, h_x)
\end{align*}
is a smooth embedding for any $l \geq 1$.
In particular, $s^0$ is smooth.
Furthermore, the assumption of the surjectivity of $D_{p_0}^+$ implies that
$s^0$ is a submersion.

Define $V = \{s^0 = 0\} \subset \hat V$.
Then the map $s : V \to E = E^0 \oplus \bigoplus_{z_\alpha^+} \R^2$ defined by
\[
s(a, b, x) = (h, p \circ \pi_Y \circ \Phi_{a, b}(\xi_x)(Z_\alpha(a))
\]
is also smooth.
If we fix a family of coordinates $(\phi_{\pm\infty_i})$ of limit circles of $\Sigma_0$,
then the map
\begin{align*}
V &\to \prod_{z_i} Y \times \prod_{\pm\infty_i} P\\
(a, b, x) &\mapsto (\pi_Y \circ \Phi_{a, b}(\xi_x)(z_i), \pi_Y \circ
\Phi_{a, b}(\xi_x)|_{S^1_{\pm\infty_i}}
\circ \phi_{\pm\infty_i})
\end{align*}
is also smooth.
The assumption of the surjectivity of $D_{p_0}^+$ implies that this map is transverse
to the product directions of tangents of the $S^1$-actions on $P$,
that is, its differentials are transverse to $0 \oplus \bigoplus_{\pm\infty_i} \R R_\lambda
\subset \bigoplus_{z_i} TY \oplus \bigoplus_{\pm\infty_i} TS^1$.
It is easy to see that for any $(a, b, x) \in V$ and $\mu \in M_i$ such that
$\kappa_\mu = 0$, the asymptotic parameter $b_\mu$ and the map
$u = \Phi_{a, b}(\xi_x)$ satisfies equation (\ref{asymptotic parameter}).

Assuming the smoothness, we define an orbibundle $(\V, \E)$ by
$\V = V / G_0$ and $\E = E/G_0$, where we regard $E$ as a trivial vector bundle
(with non-trivial $G_0$ action) on $V$.
We denote the induced smooth section $\V \to \E$ by $s$,
and define $\psi : \V \supset s^{-1}(0) \to \widehat{\M}^0$ by
$\psi(a, b, x) = (\widetilde{P}_a, Z(a), \Phi_{a, b}(\xi_x))$.
First we prove that $(\V, \E, s, \psi)$ is a Kuranishi neighborhood
of $p_0 \in \widehat{\M}(Y, \lambda, J)$ if $V$ is sufficiently small,
that is, we prove the following proposition.
\begin{prop}\label{psihomeo}
$\psi : \V \supset s^{-1}(0) \to \widehat{\M}$ is an homeomorphism onto
a neighborhood of $p_0 \in \widehat{\M}$ if $V$ is replaced by a small neighborhood
of $(0, b^0, 0)$.
\end{prop}
To prove this proposition,
first we prove a lemma about additional marked points.
To state the lemma, it is convenient to introduce the space
$\widehat{\M}' = \overline{\M} / \sim$.
This is the quotient space of $\overline{\M}$ obtained by ignoring the coordinates of
limit circles (but respecting the order of the limit circles and the
marked points).
Points of $\widehat{\M}'$ is written as $(\Sigma, z, u, \iota^\pm)$, 
where we regard $z = (z_i)$ as a sequence, and $\iota^\pm$ are bijections from the
index set $\{\pm\infty_i\}$ to the set of $\pm\infty$-limit circles of $\Sigma$.

Let $p_0 = (\Sigma_0, z_0, u_0, \iota^\pm_0) \in \widehat{\M}'(Y, \lambda, J)$ be an
arbitrary curve.
Adding marked points to $p_0$,
we get a curve $p_0^+ = (\Sigma_0, z_0 \cup z_0^+, u_0, \iota^\pm_0)$
such that all nontrivial components of $(\Sigma, z_0 \cup z_0^+)$ are stable.
(Nontrivial components are the components which do not correspond to
the trivial cylinders in $p_0$.)
We assume that there exists a finite union of codimension-two submanifolds
$S \subset Y$ such that $\pi_Y \circ u_0$ intersects $S$ at $z_0^+$ transversely.
(We do not assume the transversality of the other intersections.)
We fix an order of $z_0^+ = (z^+_{0, i})_{i \in \Lambda}$ and regard it as a sequence of
additional marked points.

Assume that $G'_0 = \Aut'(\Sigma_0, z_0, u_0, \iota^\pm_0)
= \{g \in \Aut(\Sigma); u_0 \circ g = u_0, g(z_i) = z_i, g \iota^\pm_0 = \iota^\pm_0\}$
preserves $z_0^+$ as a set.
This implies $G'_0$ acts on the index set $\Lambda$ by
$z^+_{0, g \cdot i} = g^{-1} (z^+_{0, i})$.
Then the following lemma holds.
\begin{lem}\label{additional marked points and lifts}
If $U_0 \subset \widehat{\M}'$ is a sufficiently small neighborhood of $p_0$,
then there exists an open neighborhood $U_0^+ \subset \widehat{\M}'$ of
$p_0^+$ such that the following holds true.
For any $p = (\Sigma, z, u, \iota^\pm) \in U_0$, there exists a sequence of
additional marked points $z^+ = (z^+_i)_{i \in \Lambda} \subset \Sigma$ such that
$p^+ = (\Sigma, z \cup z^+, u, \iota^\pm) \in U_0^+$ and $\pi_Y \circ u(z^+) \subset S$.
Furthermore, for each $p \in U_0$, $G'_0$ acts on the set of such points $\{ p^+ \}$
transitively, where $G'_0$-action is defined by
\[
g \cdot (\Sigma, z \cup (z^+_i)_{i \in \Lambda}, u, \iota^\pm)
= (\Sigma, z \cup (z^+_{g^{-1} \cdot i})_{i \in \Lambda}, u, \iota^\pm).
\]
\end{lem}
We call each $p^+ \in U_0^+$ a lift of $p \in U_0$.
\begin{proof}
Since $\pi_Y \circ u$ intersects with $S$ at $z_i^+$ transversely,
the existence of such a sequence of additional marked points $z^+$ is clear
for each point $p$ in a neighborhood of $p_0$.
(We use the fact that if a $J$-holomorphic curve $v$ on a disc
$D = \{ z \in C; |z| \leq 1\}$ is sufficiently close to a given $J$-holomoprhic curve $u$
on $D$ in $L^\infty$-norm, then $v$ is close to $u$ in $C^\infty$-topology on
$\{ z\in \C; |z| \leq 1/2\}$.)

We need to prove that the $G'_0$-aciton on $\{p^+\}$ is transitive
for any point $p$ sufficiently close to $p_0$.
Suppose not.
Then there exists a sequence $p_k = (\Sigma_k, z_k, u_k, \iota^\pm_k) \in \widehat{\M}'$
converging to $p_0$ and sequences of additional marked points
$z_k^+ = (z_{k, i}^+)$ and ${z'}_k^+ = ({z'}_{k, i}^+) \subset \Sigma_k$ such that
$\pi_Y \circ u_k (z_{k, i}^+), \pi_Y \circ u_k ({z'}_{k, i}^+) \in S$ and 
both of $p_k^+ = (\Sigma_k, z_k \cup z_k^+, u_k, \iota^\pm_k)$ and
${p'}_k^+ = (\Sigma_k, z_k \cup {z'}_k^+, u_k, \iota^\pm_k)$ converge to $p_0^+$,
but there is no $g \in G'_0$ such that ${z'}^+_{k, i} = z^+_{k, g \cdot i}$.

Take additional marked points $z_0^{++}$ of $(\Sigma_0, z_0 \cup z_0^+)$ to make
$(\Sigma_0, z_0 \cup z_0^+ \cup z_0^{++})$ stable, and
let $(\widetilde{P} \to \widetilde{X}, Z \cup Z^+ \cup Z^{++})$ be
the local universal family of $(\Sigma_0, z_0 \cup z_0^+ \cup z_0^{++})$.
Then by the definition of the topology,
there exist sequences $a_k, a'_k \in \widetilde{X}$ converging to $0$ and isomorphisms
\begin{align*}
\varphi_k &: (\widetilde{P}_{a_k}, Z(a_k) \cup Z^+(a_k)) \tocong (\Sigma_k, z_k \cup z_k^+),
\\
\varphi'_k &: (\widetilde{P}_{a'_k}, Z(a'_k) \cup Z^+(a'_k)) \tocong (\Sigma_k, z_k \cup
{z'}_k^+)
\end{align*}
which preserve the order of the limit circles,
and $\R$-gluings $\theta_k, \theta'_k : \overline{\R}_1 \sqcup \overline{\R}_2 \sqcup
\dots \sqcup \overline{\R}_l \to \overline{\R}_1 \cup \overline{\R}_2 \cup \dots
\cup \overline{\R}_{l_k}$ such that
\begin{align*}
\dist_{L^\infty} &(u_k \circ \varphi_k, (\theta_k \times 1) \circ u_0 \circ
\Psi|_{\widetilde{P}_{a_k}}) \to 0,\\
\dist_{L^\infty} &(u_k \circ \varphi'_k, (\theta'_k \times 1) \circ u_0 \circ
\Psi|_{\widetilde{P}_{a'_k}}) \to 0
\end{align*}
as $k \to \infty$.
Hence the biholomoprhisms $\phi_k = (\varphi'_k)^{-1} \circ \varphi_k
: \widetilde{P}_{a_k} \tocong \widetilde{P}_{a'_k}$ satisfy
\[
\dist_{L^\infty}((\theta_k \times 1) \circ u_0 \circ \Psi|_{\widetilde{P}_{a_k}},
(\theta'_k \times 1) \circ u_0 \circ \Psi|_{\widetilde{P}_{a'_k}} \circ \phi'_k)
\to 0 \quad \text{as } k \to \infty.
\]

Note that $\phi_k(Z_i(a_k)) = Z_i(a'_k)$ for all $i$, which implies
$\Psi|_{\widetilde{P}_{a'_k}} \circ \phi_k \circ (\Psi|_{\widetilde{P}_{a_k}})^{-1}$
preserve marked points $Z(0)$.
Similarly, it preserves the order of the limit circles.
Hence a subsequence of
$\Psi|_{\widetilde{P}_{a'_k}} \circ \phi_k \circ (\Psi|_{\widetilde{P}_{a_k}})^{-1}$ converges to
a biholomorphism $g \in \Aut' (\Sigma_0, z_0, u_0, \iota^\pm_0)$ on any compact subset
of the complement of nodal points and imaginary circles.
(First we see that
$\Psi|_{\widetilde{P}_{a'_k}} \circ \phi_k \circ (\Psi|_{\widetilde{P}_{a_k}})^{-1}$
converges to a biholomorphism on all nontrivial components, and then
we see the difference of $\theta_k$ and $\theta'_k$ goes to zero as $k \to \infty$,
which implies $\Psi|_{\widetilde{P}_{a'_k}} \circ \phi_k \circ (\Psi|_{\widetilde{P}_{a_k}})^{-1}$
converges to a biholomoprhism on trivial cylinders.)

Therefore
\[
\dist (\Psi|_{\widetilde{P}_{a'_k}} \circ \phi_k \circ
(\Psi|_{\widetilde{P}_{a_k}})^{-1} (Z_i^+(0)), Z^+_{g^{-1} \cdot i}(0)) \to 0 \quad \text{as }
k \to \infty,
\]
which implies
\[
\dist (\phi_k(Z_i^+(a_k)), Z_{g^{-1} \cdot i}^+(a'_k)) \to 0 \quad \text{as } k \to \infty.
\]
Since $\pi_Y \circ u_k \circ \varphi'_k = \pi_Y \circ u_k \circ \varphi_k \circ \phi_k^{-1}$
intersects transversely with $S$ at both of $\phi_k(Z_i^+(a_k))$ and
$Z_{g^{-1} \cdot i}^+ (a'_k)$,
it implies that $\phi_k(Z_i^+(a_k)) =  Z_{g^{-1} \cdot i}^+(a'_k)$ for all large $k$.
Hence $z_{k, i}^+ = {z'}_{k, g^{-1} \cdot i}^+$, which is a contradiction.
\end{proof}

Next, for the proof of the continuity of $\psi$,
we consider the following approximate solutions centered at each point
$(a, b, x) \in \hat V$.
For $(a', b') \in \mathring{X}$ close to $(a, b)$,
we define $\xi_{(a, b, x)}^{(a', b')} \in
\widetilde{W}^{1, p}_\delta(\widetilde{P}_{a'}, u_{a', b'}^\ast T \hat Y)$ as follows.
We may assume that
$\rho'_\nu \neq 0$ for all $\nu$ such that $\rho_\nu \neq 0$, and that
$\rho'_\mu \neq 0$ for all $\mu$ such that $\rho_\mu \neq 0$.
Recall that $\phi^{a, b}(x) = (\xi_{(a, b, x)}, h_{(a, b, x)})$.

On $\Sigma_0 \setminus N_0$, we define
$\xi_{(a, b, x)}^{(a', b')}|_{\Sigma_0 \setminus N_0}
= \xi_{(a, b, x)}|_{\Sigma_0 \setminus N_0}$.
Similarly, we define $\xi_{(a, b, x)}^{(a', b')}|_{[0, \infty) \times S_{+\infty_i}^1}
= \xi_{(a, b, x)}|_{[0, \infty) \times S_{+\infty_i}^1}$ on $[0, \infty) \times S_{+\infty_i}^1$
and $\xi_{(a, b, x)}^{(a', b')}|_{(-\infty, 0] \times S_{-\infty_i}^1}
= \xi_{(a, b, x)}|_{(-\infty, 0] \times S_{-\infty_i}^1}$ on $(-\infty, 0] \times S_{-\infty_i}^1$.

On $[-1, 0] \times S_\mu^1 \subset[-1, -\frac{1}{2} \log \rho'_\mu] \times S_\mu^1$,
define $\xi_{(a, b, x)}^{(a', b')}(s, t) = \xi_{(a, b, x)}(s, t)$.
On $[0, -\frac{1}{2} \log \rho'_\mu] \times S_\mu^1$,
define $\xi_{(a, b, x)}^{(a', b')}(s', t') = \xi_{(a, b, x)}(s, t')$,
where $s$ is defined by
\begin{equation}
\frac{e^{- \kappa_\mu s'} - (\rho'_\mu)^{\kappa_\mu / 2}}
{1 - (\rho'_\mu)^{\kappa_\mu / 2}}
= \frac{e^{- \kappa_\mu s} - \rho_\mu^{\kappa_\mu / 2}}
{1 - \rho_\mu^{\kappa_\mu / 2}}.
\label{def of function s}
\end{equation}
We define similarly $\xi_{(a, b, x)}^{(a', b')}$ on
$[\frac{1}{2} \log \rho'_\mu, 1] \times S^1_\mu$.

For each $\nu$,
define $\rho_\nu$, $\varphi_\nu$, $\rho'_\nu$ and $\varphi'_\nu$ by
$\zeta_\nu = \rho_\nu^2 e^{2 \sqrt{-1} \varphi_\nu}$
and $\zeta'_\nu = (\rho'_\nu)^2 e^{2 \sqrt{-1} \varphi'_\nu}$.
On $N_{a', b'}^{\nu} = \{(z', w') \in D \times D; z'w' = \zeta'_\mu\}$,
define $\xi_{(a, b, x)}^{(a', b')}(z', w') = \xi_{(a, b, x)}(z, w)$,
where $(z, w) \in \{(z, w) \in D \times D; zw = \zeta_\mu\}$
is defined as follows.
If $|z'| \geq \rho'_\nu$ and $z' = r' e^{\sqrt{-1} \theta'}$
then $z = r e^{\sqrt{-1} \theta}$,
and if $|w'| \geq \rho'_\nu$ and $w' = r' e^{\sqrt{-1} \theta'}$
then $w = r e^{\sqrt{-1} \theta}$,
where $r$ and $\theta$ are defined by
\[
\frac{r' - \rho'_\nu}{1 - \rho'_\nu} = \frac{r - \rho_\nu}{1 - \rho_\nu},
\quad \theta' = \theta + (\varphi'_\nu - \varphi_\nu) \beta_\nu(r),
\quad \beta_\nu(r) = \max\Bigl( \frac{2 \rho_\nu - r}{\rho_\nu}, 0\Bigr).
\]
(If $\rho_\nu = 0$, we define $\beta_\nu = 0$.)

Then $(\xi_{(a, b, x)}^{(a', b')}, h_{(a, b, x)})$ satisfies the following estimates.
\begin{lem}\label{estimates for continuity}
For any $0 < \delta < \delta'_0 < \delta_0$, $p > 2$ and $(a, b, x) \in \hat V$,
there exists a constant $C>0$ such that
for any $(a',b') \in \mathring{X}$ sufficiently close to $(a, b)$,
\begin{align*}
&||F^{(a', b')}(\xi_{(a, b, x)}^{(a', b')}, h_{(a, b, x)})|_{[-1, 0] \times S_\mu^1}||_{L^p}
\leq C (|\varphi'_\mu - \varphi_\mu| + |b'_\mu - b_\mu|), \\
&||F^{(a', b')}(\xi_{(a, b, x)}^{(a', b')}, h_{(a, b, x)})|_{[0, -\frac{1}{2} \log \rho'_\mu]
\times S_\mu^1}||_{L^p_{\delta_\mu}}\\
&\hspace{90pt} \leq
\begin{cases}
C (\rho'_\mu)^{\min (\kappa_\mu, \delta'_{0, \mu} - \delta_\mu)/2}
(-\log \rho'_\mu)^{1/p}, & \text{if } \rho_\mu = 0\\
C |\rho'_\mu - \rho_\mu| & \text{if } \rho_\mu > 0
\end{cases},\\
&||F^{(a', b')}(\xi_{(a, b, x)}^{(a', b')}, h_{(a, b, x)})|_{\{(z, w) \in N^\nu_{(a', b')};
|z| \geq \rho'_\nu \}}||_{L^p} \leq C |\zeta'_\nu - \zeta_\nu|^{1/p}.
\end{align*}
\end{lem}
\begin{proof}
The estimate of $F^{(a', b')}(\xi_{(a, b, x)}^{(a', b')}, h_{(a, b, x)})|_{[-1, 0] \times S_\mu^1}$
is similar to Lemma \ref{estimates of F(0)}.
Since $F^{(a, b)}(\xi_{(a, b, x)}, h_{(a, b, x)})|_{[-1, 0] \times S_\mu^1} = 0$,
\begin{align*}
&F^{(a', b')}(\xi_{(a, b, x)}^{(a', b')}, h_{(a, b, x)})|_{[-1, 0] \times S_\mu^1}\\
&= F^{(a', b')}(\xi_{(a, b, x)}^{(a', b')}, h_{(a, b, x)})|_{[-1, 0] \times S_\mu^1}
- F^{(a, b)}(\xi_{(a, b, x)}, h_{(a, b, x)})|_{[-1, 0] \times S_\mu^1}\\
&= \frac{1}{2} (b'_\mu -b_\mu) \chi'(s) \partial_\sigma\\
&\quad + \frac{1}{2} (\varphi'_\mu - \varphi_\mu) \chi'(s)
(g^\mu_t(v_0^{\mu, \mathrm{left}} + \xi_{(a, b, x)})
+ \partial_t(v_0^{\mu, \mathrm{left}} + \xi_{(a, b, x)})).
\end{align*}
The first inequality follows from this equation.

Next we consider the second inequality.
If $\rho_\mu = 0$, then the proof is similar to Lemma \ref{estimates of F(0)}.
(The proof is obtained by
replacing $v_0^{\nu, \text{left}}$ and $\rho_\mu$ in the proof of
Lemma \ref{estimates of F(0)} with $v_0^{\nu, \text{left}}$ and $\rho'_\mu$
respectively.)
Consider the case of $\rho_\mu > 0$.
We abbreviate the subscript $\mu$.
First note that
\begin{align}
&F^{(a', b')} (\xi_{(a, b, x)}^{(a', b')}, h_{(a, b, x)})|_{[0, -\frac{1}{2} \log \rho'] \times S^1}
(s', t) \notag\\
&= \biggl(1 + \biggl(\frac{1 - \rho^{\kappa/2}}
{1 - (\rho')^{\kappa/2}}
(\rho')^{\kappa/2} - \rho^{\kappa/2}\biggr)
e^{\kappa s}\biggr)
\partial_s \bigl(v_{a,b}^{\mu, \mathrm{left}} + \xi_{(a, b, x)}\bigr)(s, t) \notag\\
&\quad + \widetilde{J}^\mu_t\bigl(
\bigl(v_{a,b}^{\mu, \mathrm{left}} + \xi_{(a, b, x)}\bigr)(s, t)\bigr)
\partial_t \bigl(v_{a,b}^{\mu, \mathrm{left}} + \xi_{(a, b, x)}\bigr)(s, t) \notag\\
&\quad + f^\mu_t \circ \pi_Y\bigl(
\bigl(v_{a,b}^{\mu, \mathrm{left}} + \xi_{(a, b, x)}\bigr)(s, t)\bigr),
\label{F' joint circle}
\end{align}
where $s$ is a function of $s'$ defined by (\ref{def of function s}).
Subtracting the equation
\begin{align*}
0 &= F^{(a, b)} (\xi_{(a, b, x)}, h_{(a, b, x)})|_{[0, -\frac{1}{2} \log \rho] \times S^1}(s,t)\\
&= \partial_s \bigl(v_{a,b}^{\mu, \mathrm{left}} + \xi_{(a, b, x)}\bigr)(s, t)\\
&\quad + \widetilde{J}^\mu_t\bigl(
\bigl(v_{a,b}^{\mu, \mathrm{left}} + \xi_{(a, b, x)}\bigr)(s, t)\bigr)
\partial_t \bigl(v_{a,b}^{\mu, \mathrm{left}} + \xi_{(a, b, x)}\bigr)(s, t) \\
&\quad + f^\mu_t \circ \pi_Y\bigl(
\bigl(v_{a,b}^{\mu, \mathrm{left}} + \xi_{(a, b, x)}\bigr)(s, t)\bigr)
\end{align*}
from (\ref{F' joint circle}), we obtain
\begin{align*}
&F^{(a', b')} (\xi_{(a, b, x)}^{(a', b')}, h_{(a, b, x)})|_{[0, -\frac{1}{2} \log \rho'] \times S^1}
(s', t) \\
&= \biggl(\frac{1 - \rho^{\kappa/2}}
{1 - (\rho')^{\kappa/2}}
(\rho')^{\kappa/2} - \rho^{\kappa/2}\biggr) e^{\kappa s}
\partial_s \bigl(v_{a,b}^{\mu, \mathrm{left}} + \xi_{(a, b, x)}\bigr)(s, t).
\end{align*}
Since $e^{\kappa s} \partial_s \bigl(v_{a,b}^{\mu, \mathrm{left}} + \xi_{(a, b, x)}\bigr)$
is a bounded function,
\[
||F^{(a', b')}(\xi_{(a, b, x)}^{(a', b')}, h_{(a, b, x)})|_{[0, -\frac{1}{2} \log \rho'_\mu]
\times S_\mu^1}||_{L^p_{\delta_\mu}} \lesssim |\rho' - \rho|,
\]
which prove the second inequality.

Finally, we consider the third inequality.
If $\rho_\nu = 0$, then the proof is similar to Lemma \ref{estimates of F(0)}.
We consider the case $\rho_\nu \neq 0$.
First note that
\begin{align}
&F^{(a', b')} (\xi_{(a, b, x)}^{(a', b')}, h_{(a, b, x)})|_{\{z \in D; |z| \geq \rho_\nu\}}
(r' e^{\sqrt{-1} \theta'}) \notag\\
&= \frac{1 - \rho_\nu}{1 - \rho'_\nu}
\Bigl((\partial_r v_{a, b}^\nu)(r e^{\sqrt{-1} \theta'})
+ (\partial_r \xi_{(a, b, x)})(r e^{\sqrt{-1} \theta}) \notag\\
& \quad \hph{\frac{1 - \rho_\nu}{1 - \rho'_\nu} \Bigl(}
+ \frac{\varphi'_\nu - \varphi_\nu}{\rho_\nu}
1_{[\rho_\nu, 2 \rho_\nu]}(r) (\partial_\theta \xi_{(a, b, x)})(r e^{\sqrt{-1} \theta})\Bigr)
\notag\\
&\quad
+ \widetilde{J}^\nu
\bigl(v_{a, b}^\nu(r e^{\sqrt{-1} \theta'}) + \xi_{(a, b, x)}(r e^{\sqrt{-1} \theta})\bigr)
\notag\\
&\quad \quad \cdot
\frac{(\partial_\theta v_{a, b}^\nu)(r e^{\sqrt{-1} \theta'})
+ (\partial_\theta \xi_{(a, b, x)})(r e^{\sqrt{-1} \theta})}{r'}.
\label{F' nodal point}
\end{align}
Subtracting
\begin{align*}
0 &= F^{(a, b)}(\xi_{(a, b, x)}, h_{(a, b, x)})|_{\{z \in D; |z| \geq \rho_\nu\}} \\
&= (\partial_r v_{a, b}^\nu)(r e^{\sqrt{-1} \theta})
+ (\partial_r \xi_{(a, b, x)})(r e^{\sqrt{-1} \theta}) \\
&\quad + \widetilde{J}^\nu
\bigl(v_{a, b}^\nu(r e^{\sqrt{-1} \theta}) + \xi_{(a, b, x)}(r e^{\sqrt{-1} \theta})\bigr)
\frac{(\partial_\theta v_{a, b}^\nu)(r e^{\sqrt{-1} \theta})
+ (\partial_\theta \xi_{(a, b, x)})(r e^{\sqrt{-1} \theta})}{r}
\end{align*}
from (\ref{F' nodal point}), we obtain
\begin{align*}
&F^{(a', b')} (\xi_{(a, b, x)}^{(a', b')}, h_{(a, b, x)})|_{\{z \in D; |z| \geq \rho_\nu\}}
(r' e^{\sqrt{-1} \theta'}) \\
&= ((\partial_r v_{a, b}^\nu)(r e^{\sqrt{-1} \theta'})
- (\partial_r v_{a, b}^\nu)(r e^{\sqrt{-1} \theta})) \\
&\quad + \frac{\rho'_\nu - \rho_\nu}{1 - \rho'_\nu}
\Bigl((\partial_r v_{a, b}^\nu)(r e^{\sqrt{-1} \theta'})
+ (\partial_r \xi_{(a, b, x)})(r e^{\sqrt{-1} \theta})\Bigr) \\
&\quad + \frac{1 - \rho_\nu}{1 - \rho'_\nu}
\cdot \frac{\varphi'_\nu - \varphi_\nu}{\rho_\nu}
1_{[\rho_\nu, 2 \rho_\nu]}(r) (\partial_\theta \xi_{(a, b, x)})(r e^{\sqrt{-1} \theta}) \\
&\quad +
\bigl(\widetilde{J}^\nu
\bigl(v_{a, b}^\nu(r e^{\sqrt{-1} \theta'}) + \xi_{(a, b, x)}(r e^{\sqrt{-1} \theta})\bigr)
- \widetilde{J}^\nu
\bigl(v_{a, b}^\nu(r e^{\sqrt{-1} \theta}) + \xi_{(a, b, x)}(r e^{\sqrt{-1} \theta})\bigr)
\bigr) \\
&\quad \quad \cdot
\frac{(\partial_\theta v_{a, b}^\nu)(r e^{\sqrt{-1} \theta'})
+ (\partial_\theta \xi_{(a, b, x)})(r e^{\sqrt{-1} \theta})}{r'} \\
&\quad + \widetilde{J}^\nu
\bigl(v_{a, b}^\nu(r e^{\sqrt{-1} \theta}) + \xi_{(a, b, x)}(r e^{\sqrt{-1} \theta})\bigr) \\
&\quad \quad \cdot
\frac{(\partial_\theta v_{a, b}^\nu)(r e^{\sqrt{-1} \theta'})
- (\partial_\theta v_{a, b}^\nu)(r e^{\sqrt{-1} \theta})}{r'} \\
&\quad + \Bigl( \frac{1}{r'} - \frac{1}{r} \Bigr)
\widetilde{J}^\nu
\bigl(v_{a, b}^\nu(r e^{\sqrt{-1} \theta}) + \xi_{(a, b, x)}(r e^{\sqrt{-1} \theta})\bigr) \\
&\quad \quad \cdot
\bigl((\partial_\theta v_{a, b}^\nu)(r e^{\sqrt{-1} \theta})
+ (\partial_\theta \xi_{(a, b, x)})(r e^{\sqrt{-1} \theta})\bigr).
\end{align*}
Hence it is easy to check that
\[
||F^{(a', b')} (\xi_{(a, b, x)}^{(a', b')}, h_{(a, b, x)})|_{\{z \in D; |z| \geq \rho_\nu\}}||_{L^p}
\lesssim |\rho'_\nu - \rho_\nu| + |\varphi'_\nu - \varphi_\nu|,
\]
and this inequality implies the claim.
\end{proof}
Note that the $\Ker DF^{(0, b^0)}_{(0, 0)}$-factor of
$F^{(a', b') +}(\xi_{(a, b, x)}^{(a', b')}, h_{(a, b, x)})$ coincides with that of
$F^{(a, b) +}(\xi_{(a, b, x)}, h_{(a, b, x)})$.
Hence the above lemma implies that
\[
\bigl|\bigl|F^{(a', b') +}\bigl(\xi_{(a, b, x)}^{(a', b')}, h_{(a, b, x)}\bigr) - (0, x')\bigr|\bigr|
_{L^p_\delta \oplus \Ker DF^{(0, b^0)}_{(0, 0)}} \to 0
\]
as $(a', b', x') \to (a, b, x)$.
Therefore
\begin{equation}
\bigl|\bigl| \phi^{(a', b')}(x')- \bigl(\xi_{(a, b, x)}^{(a', b')}, h_{(a, b, x)}\bigr)\bigr|\bigr|
_{\widetilde{W}^{1, p}_\delta(\widetilde{P}_{a'}, u_{a', b'}^\ast T \hat Y) \oplus E^0}
\to 0
\label{continuity of phi}
\end{equation}
as $(a', b', x') \to (a, b, x)$.
This implies the continuity of $\psi : \mathcal{V} \supset s^{-1}(0) \to
\widehat{\M}$.

Now we prove Proposition \ref{psihomeo}.
\begin{proof}[Proof of Proposition \ref{psihomeo}]
We have just proved the continuity of $\psi$.
Next we prove the injectivity.
Assume that the image of two points $(a, b, x), (a', b', x') \in \{ s = 0 \}$
($ \subset V$) coincide,
that is, the two holomorphic buildings $(\widetilde{P}_a, Z(a), \Phi_{a, b}(\xi_x))$
and $(\widetilde{P}_{a'}, Z(a'), \Phi_{a', b'}(\xi_{x'}))$ are the same point in $\widehat{\M}$.
We prove that these two points coincide in $V / G_p$.

Since $\widehat{\M}$ is a quotient space of $\widehat{\M}'$,
we may assume that these two holomorphic buildings also coincide in $\widehat{\M}'$
by replacing $(a', b', x')$ with $h \cdot (a', b', x')$ for some $h \in G_0$.

If $V$ is sufficiently small, then Lemma \ref{additional marked points and lifts} implies
that there exist an isomorphism
$\varphi : (\widetilde{P}_a, Z(a)) \tocong (\widetilde{P}_{a'}, Z(a'))$,
an $\R$-translation $\theta$, and $g \in G'_0$ such that
$\Phi_{a', b'}(\xi_{x'}) \circ \varphi = (\theta \times 1) \circ \Phi_{a, b}(\xi_x)$,
$\varphi \circ \iota^\pm = \iota^\pm$ and $\varphi(Z^+_i(a)) = Z^+_{g^{-1} \cdot i}(a')$.
Hence the isomorphism $(\hat P_a, Z(a)) \tocong (\hat P_{a'}, Z(a'))$
induced by $\varphi$ coincides with the restriction of $g : \hat P \to \hat P$.
Therefore $\varphi$ preserves $\widetilde{R}_i$ as a family,
which implies $\theta = \id$.
From this, we can see that $\varphi$ maps $Z^{++}_i(a)$ to $Z^{++}_{g^{-1} \cdot i}(a')$
because these points are contained in the inverse image of $S'$ by
$\Phi_{a, b}(\xi_x)$ and $\Phi_{a', b'}(\xi_{x'})$ respectively.
Hence $\varphi : (\widetilde{P}_a, Z(a)) \tocong (\widetilde{P}_{a'}, Z(a'))$ coincides with
the restriction of $g : \widetilde{P} \to \widetilde{P}$,
which implies $(a', b', x') = g \cdot (a, b, x)$.
Therefore $\psi : \{ s = 0 \}/G_0 \to \widehat{\M}$ is injective.

Finally we prove that the image of $\psi$ contains a neighborhood of $p$.
Assume contrary, that is, assume that there exists a sequence
$(\Sigma_k, z_k, u_k) \in \widehat{\M} \setminus \psi(\{ s = 0 \} /G_p)$ convergent to
$p_0 = (\Sigma_0, z, u_0)$.
We may assume that $(\Sigma_k, z_k, u_k, \iota^\pm_k) \in \widehat{\M}'$ converges to
$p_0 = (\Sigma, z, u_0, \iota^\pm_0) \in \widehat{\M}'$.
Let $(\Sigma_k, z_k \cup z_k^+, u_k)$ be the lift of $(\Sigma_k, z_k, u_k)$ for each $k$.
Then there exist a sequence $a_k \to 0 \in \widetilde{X}$, biholomorphisms
$(\Sigma_k, z_k \cup z_k^+) \cong (\widetilde{P}_{a_k}, Z(a_k) \cup Z^+(a_k))$
and $\R$-gluings $\theta_k$ such that
\[
\dist_{L^\infty} (u_k, (\theta_k \times 1) \circ u_0 \circ \Psi|_{\widetilde{P}_{a_k}}) \to 0
\quad \text{as } k \to \infty.
\]
We may assume that $\theta_k$ maps $0 \in \R_i$ to
$\sigma \circ u_k(\widetilde{R}_i(a_k))$.
Let $z^{++}_k \subset \Sigma_k$ be the points corresponding to $Z^{++}(a_k) \subset
\widetilde{P}_{a_k}$. Changing $a_k$ slightly if necessary,
we may assume $u_k(z^{++}_k) \subset (\theta_k \times 1)(S')$.

Define $(b^k_\mu) \in \prod_{i = 1}^{k - 1} \bigoplus_{\mu \in M_i} \R$ as follows:
\begin{itemize}
\item
For each $\mu$ such that $\rho^k_\mu = 0$,
we define $b_\mu^{k, \mathrm{left}}$ and $b_\mu^{k, \mathrm{right}}$, and $b^k_\mu$ by
\begin{align*}
\sigma \circ u_k|_{[-1,\infty) \times S^1_\mu}(s_x,t_x)
&= \theta_k(0_i) + L_\mu s_x + b_\mu^{k, \mathrm{left}} + O(1)\\
\sigma \circ u_k|_{(-\infty,+1] \times S^1_\mu}(s_y,t_y)
&= \theta_k(0_{i + 1}) + L_\mu s_y + b_\mu^{k, \mathrm{right}} + O(1)\\
b^k_\mu = b_\mu^{k, \mathrm{left}} - b_\mu^{k, \mathrm{right}}.
\end{align*}
\item
For each $\mu$ such that $\rho^k_\mu \neq 0$,
we define $b^k_\mu$ by
\[
\sigma \circ u_k(\widetilde{R}_{i + 1}(a_k)) - \sigma u_k(\widetilde{R}_i(a_k))
= -L_\mu \log \rho_\mu + b^k_\mu.
\]
\end{itemize}
Then $b^k_\mu \to b^0_\mu$ as $k \to \infty$.
(Note that in the former case, the asymptotic estimates of the term
$O(1)$ is uniform with respect to $k$.)

Replacing each map $u_k$ with its appropriate $\R$-translation,
we may assume $\dist_{L^\infty} (u_k, u_{a_k, b_k}) \to 0$.
Then there exists a section $\xi_k$ of $u_{a_k, b_k}^\ast T \hat Y$ for each $k$
such that $||\xi_k||_\infty \to 0$ as $k \to \infty$ and $u_k = \Phi_{a_k, b_k}(\xi_k)$.

To prove $|| \xi_k ||_{\widetilde{W}^{1,p}_\delta} \to 0$, we consider the following equations.
\begin{align*}
F^{(a_k, b^k)+} (\xi_k, 0) =& F^{(a_k, b^k)+}(0, 0) + DF^{(a_k, b^k)+}_{(0, 0)}(\xi_k, 0)\\
&+ \int_0^1\bigl( DF^{(a_k, b^k)+}_{(\lambda \xi_k, 0)} - DF^{(a_k, b^k)+}_{(0, 0)} \bigr)
(\xi_k, 0) d\lambda
\end{align*}
In the above equations,
\begin{align*}
&||F^{(a_k, b^k)+} (\xi_k, 0)||_{L_\delta^p \oplus \Ker DF^{(0, b^0)}_{(0, 0)}}
= ||F^{(a_k, b^k)+} (\xi_k, 0)||_{\Ker DF^{(0, b^0)}_{(0, 0)}} \to 0,\\
&||F^{(a_k, b^k)+} (0, 0)||_{L_\delta^p \oplus \Ker DF^{(0, b^0)}_{(0, 0)}}
\lesssim |a_k| + |b^k - b^0| \to 0,\\
&||DF^{(a_k, b^k)+}_{(0, 0)} (\xi_k, 0)||_{L_\delta^p \oplus \Ker DF^{(0, b^0)}_{(0, 0)}}
\geq \epsilon ||\xi_k||_{\widetilde{W}^{1,p}_\delta} \quad \quad \text{for some } \epsilon>0,
\\
&\Bigl|\Bigl|\int_0^1\bigl( DF^{(a_k, b^k)+}_{(\lambda \xi_k, 0)}
- DF^{(a_k, b^k)+}_{(0, 0)} \bigr)
(\xi_k, 0) d\lambda\Bigr|\Bigr|_{L_\delta^p \oplus \Ker DF^{(0, b^0)}_{(0, 0)}}
\leq ||\xi_k||_\infty ||\xi_k||_{\widetilde{W}^{1,p}_\delta}
\end{align*}
by Lemma \ref{D^2F}.
These imply $|| \xi_k ||_{\widetilde{W}^{1,p}_\delta} \to 0$.
Hence $(\Sigma_k, z_k, u_k)$ is contained in the image of $\psi$ for large $k$,
which is a contradiction.
Therefore the image of $\psi$ contains a neighborhood of $p$.

Since $\{s=0 \} /G_0$ is locally compact and $\widehat{\M}$ is Hausdorff,
$\psi$ is a homeomorphism onto a neighborhood of $p_0$.
\end{proof}

Therefore $(\V, \E, s, \psi)$ is a Kuranishi neighborhood of $p_0$.
In this paper, we sometimes denote the Kuranishi neighborhood by the $5$-tuple
$(V, E, s, \psi, G_0)$.
Sometimes we write a point of $V$ as a $4$-tuple $(\Sigma, z, u, h)$ consisting of
a curve $\Sigma$, its marked points $z$, a map $u$ and a vector $h \in E^0$
which satisfy the equation $du + J(u)du + \lambda(h) = 0$.

\subsection{Linearized gluing lemma}\label{linearized gluing}
In this section, we prove the linearized gluing lemma
(Lemma \ref{linearized gluing lemma} below), which was
used in the previous section to prove the invertibility of $DF^{(a, b) +}_{(0, 0)}$.

Let $\Sigma$ be a domain curve of a holomorphic building,
and let $E \to \Sigma$ be a complex vector bundle of rank $n$.
Assume that on a neighborhood $N_0 \subset \Sigma$ of nodal points and
imaginary circles,
a trivialization $E|_{N_0} \cong N_0 \times \C^n$ is given.
$N_0$ is the union of $D \cup D$, $([0, \infty] \cup [-\infty, 0]) \times S^1$,
$[0, \infty] \times S^1$ and $[-\infty, 0] \times S^1$.

Assume that an elliptic operator $D_0$ on $E$ has the same symbol as $\bound$,
and on the neighborhood $([0, \infty] \cup [-\infty, 0]) \times S^1$ of
each joint circle $S_\mu^1$, $D_0$ has the form
\[
D_0 \xi = \partial_s \xi + J_0 \partial_t \xi + S_\mu(s, t) \xi,
\]
where $S_\mu(s, t) : ([0, \infty] \cup [-\infty, 0]) \times S^1 \to \gl(2n, \R)$ is
a continuous matrix-valued function such that
$S_\mu(t) : = S_\mu(\pm\infty, t) : S^1 \to \gl(2n, \R)$ is a loop of symmetric matrices.
Also on the neighborhood $[0, \infty] \times S_{+\infty_i}^1$ or
$[-\infty, 0] \times S_{-\infty_i}^1$ of each limit circle $S_{\pm\infty_i}^1$,
$D_0$ has the form
\[
D_0 \xi = \partial_s \xi + J_0 \partial_t \xi + S_{\pm\infty_i}(s, t) \xi,
\]
where $S_{\pm\infty_i}$ are continuous matrix-valued functions on
$[0, \infty] \times S_{+\infty_i}^1$ or $[-\infty, 0] \times S_{-\infty_i}^1$ such that
$S_{\pm\infty_i}(t) := S_{\pm\infty_i}(\pm\infty, t) : S^1 \to \gl(2n, \R)$ are loop of
symmetric matrices.

We further assume that there exist a family of positive constants
$\delta_1 = ((\delta_{1, \mu})_\mu, (\delta_{1, \pm\infty_i})_{\pm\infty_i})$
and a constant $C > 0$ such that
\begin{align*}
|S_\mu(s, t) - S_\mu(t)| &\leq C e^{-\delta_{1, \mu} |s|}
\text{ for } s \in [0, \infty] \cup [-\infty, 0]\\
|S_{\pm\infty_i}(s, t) - S_{\pm\infty_i}(t)| &\leq C e^{-\delta_{1, \pm\infty_i} |s|}
\text{ for } s \in [0, \infty] \text{ (or } s\in [-\infty ,0])
\end{align*}

Let $\delta_0 = ((\delta_{0, \mu})_\mu, (\delta_{0, \pm\infty_i})_{\pm\infty_i})$
be the family of positive constants consisting of the minimal
non-zero absolute values of eigenvalues of 
\[
A_\mu = J_0 \partial_t + S_\mu(t) : W^{1,2}(S^1, \R^{2n}) \to L^2(S^1,\R^{2n})
\]
and
\[
A_{\pm\infty_i} = J_0 \partial_t + S_{\pm\infty_i}(t) : W^{1,2}(S^1, \R^{2n}) \to
L^2(S^1,\R^{2n}).
\]
Let $\delta = ((\delta_\mu)_\mu, (\delta_{\pm\infty_i})_{\pm\infty_i})$ be an arbitrary
sequence of constants such that $\delta < \delta_0$ and $\delta < \delta_1$,
and let $2 < p < \infty$ be an arbitrary constant.
We define the $L_\delta^p$-norm on $[0, \infty] \times S^1$ or
$[-\infty, 0] \times S^1$ by $||\xi||_{L_\delta^p} = ||e^{\delta |s|} \xi||_{L^p}$,
using the usual Lebesgue measures of
$[0, \infty) \times S^1$ or $(-\infty, 0] \times S^1$.

Using the trivialization of $E|_{N_0}$, we define the $L_\delta^p$-space by
\begin{align*}
L_\delta^p(\Sigma, \Wedge^{0, 1}T^\ast \Sigma \otimes E)
&= L_\delta^p(\Sigma_0 \setminus N_0, \Wedge^{0, 1}T^\ast \Sigma \otimes E)\\
&\quad \oplus \bigoplus_\nu L^p(D \cup D ,\C^n)\\
&\quad \oplus
\bigoplus_\mu L_{\delta_\mu}^p(([0, \infty] \cup [-\infty, 0]) \times S^1, \C^n)\\
&\quad \oplus \bigoplus_{+\infty_i} L_{\delta_{+\infty_i}}^p([0, \infty] \times S^1, \C^n)\\
&\quad \oplus \bigoplus_{-\infty_i} L_{\delta_{-\infty_i}}^p([-\infty, 0] \times S^1, \C^n).
\end{align*}
We define a Banach space $\widetilde{W}_\delta^{1, p}(\Sigma, E)$ by
\begin{align*}
\widetilde{W}_\delta^{1, p}(\Sigma, E)
=&\ \{\xi = \xi_0 + \sum_\mu \beta_\mu v_\mu
+ \sum_{\pm\infty_i} \beta_{\pm\infty_i} v_{\pm\infty_i} \in C(\Sigma, E);\\
&\quad \quad \xi_0 \in W^{1,p}_\delta (\Sigma, E),
v_\mu \in \Ker A_\mu, v_{\pm\infty_i} \in \Ker A_{\pm\infty_i}\},
\end{align*}
where for each $\mu$, $\beta_\mu$ is a smooth function which is $1$
on some neighborhood of $\mu$-th joint circle and whose support is contained
in its slightly larger neighborhood,
and $\beta_{\pm\infty_i}$ is a similar function for each $\pm\infty_i$.
Each $v_\mu \in \Ker A_\mu$ is regarded as a section
$v_\mu(s, t) = v_\mu(t) : ([0, \infty] \cup [-\infty, 0]) \times S^1 \to \C^n$,
and the meaning of the above $v_{\pm\infty_i}$ is similar for each $\pm\infty_i$.
Then we can regard $D_0$ as a linear operator
$D_0 : \widetilde{W}_\delta^{1, p}(\Sigma, E) \to
L_\delta^p(\Sigma,\ab \Wedge^{0, 1}T^\ast \Sigma \otimes E)$.

For each $(\zeta, r) = (\zeta_\nu, r_\mu)
\in D^{l_0} \times (1, \infty]^{l_1}$,
a new curve $\Sigma_{(\zeta, r)}$ is constructed from $\Sigma$ by
replacing the neighborhood $D \cup D$ of the $\nu$-th nodal point with
\[
N^\nu_{\zeta_\nu} = \{(x, y) \in D \times D; xy = \zeta_\nu\},
\]
and replacing the neighborhood $D \widetilde{\cup} D$ of the $\mu$-th joint circle
with
\[
N^\mu_{r_\mu} = \{((s_x, t_x), (s_y, t_y)) \in [0, \infty] \times S^1 \times
[-\infty, 0] \times S^1;
s_y - s_x = -2r, t_y = t_x\}.
\]
$E$ induces a complex vector bundle on $\Sigma_{(\zeta, r)}$.
(We use the trivialization of $E|_{N_0}$.)
We also denote this vector bundle by $E$.

$L^p$-norm on $N^\nu_{\zeta_\nu}$ is defined by the measure
$\frac{\sqrt{-1}}{2} dx \wedge d\bar x$ on $\{(x, y) \in N^\nu_{\zeta_\nu}; |x| \geq |y|\}$
and the measure $\frac{\sqrt{-1}}{2} dy \wedge d\bar y$ on
$\{(x, y) \in N^\nu_{\zeta_\nu}; |y| \geq |x|\}$.
$L_\delta^p$-norm on $N^\mu_{r_\mu}$ is defined by
\[
||\xi||_{L_\delta^p}
= ||e^{\delta |s|} \xi||_{L^p([0, r_\mu] \times S^1)}
+ ||e^{\delta |s|} \xi||_{L^p([-r_\mu, 0] \times S^1)}.
\]

We define a Banach space $\widetilde{W}_\delta^{1, p}(\Sigma_{(\zeta, r)}, E)$ by
\begin{align*}
\widetilde{W}_\delta^{1, p}(\Sigma_{(\zeta, r)}, E)
=&\ \{\xi = \xi_0 + \sum_\mu \beta_\mu v_\mu
+ \sum_{\pm\infty_i} \beta_{\pm\infty_i} v_{\pm\infty_i} \in C(\Sigma, E);\\
&\quad \quad \xi_0 \in W^{1,p}_\delta (\Sigma_{(\zeta, r)}, E), v_\mu \in
\Ker A_\mu,
v_{\pm\infty_i} \in \Ker A_{\pm\infty_i}\},
\end{align*}
where $\beta_\mu$ and $\beta_{\pm\infty_i}$ are defined by regarding the curve
$\Sigma_{(\zeta, r)}$ as a curve constructed by patching the subsets
$\Sigma_0 \setminus N_0$, $\{x \in D; |x| \geq \sqrt{|\zeta_\nu|}\}$,
$\{y \in D; |y| \geq \sqrt{|\zeta_\nu|}\}$,
$[0, r_\mu] \times S_\mu^1$,
$[-r_\mu, 0] \times S_\mu^1$,
$[0, \infty] \times S_{+\infty_i}^1$ and $[-\infty, 0] \times S_{-\infty_i}^1$ of $\Sigma$.
The norm of $\widetilde{W}_\delta^{1, p}(\Sigma_{(\zeta, \kappa)}, E)$ is defined by
\begin{align*}
||\xi||_{\widetilde{W}_\delta^{1,p}(\Sigma_{(\zeta, \kappa)})}
= \inf\Bigl\{&\Bigl|\Bigl|\xi- \sum_\mu \beta_\mu v_\mu
- \sum_{\pm\infty_i} \beta_{\pm\infty_i} v_{\pm\infty_i}\Bigr|\Bigr|_{W^{1,p}_\delta}
+ \sum_\mu ||v_\mu||_{\Ker A_\mu}\\
&+ \sum_{\pm\infty_i} ||v_{\pm\infty_i}||_{\Ker A_{\pm\infty_i}}; 
v_\mu \in \Ker A_\mu, v_{\pm\infty_i} \in \Ker A_{\pm\infty_i} \Bigr\}.
\end{align*}

Regarding $\Sigma_{(\zeta, r)}$ as the curve constructed by patching
the subsets of $\Sigma$, we define the linear operator
$D_{(\zeta, r)} : \widetilde{W}_\delta^{1, p}(\Sigma_{(\zeta, r)}, E) \to
L_\delta^p(\Sigma_{(\zeta, r)}, \Wedge^{0, 1} T^\ast \Sigma_{(\zeta, r)}
\otimes E)$ from $D_0$.
(The coefficient of the operator is discontinuous in general.)
Let $\lambda : \R^N \to L_\delta^p(\Sigma, \Wedge^{0, 1}T^\ast \Sigma \otimes E)$
be a linear map which makes
\[
D_0 \oplus \lambda : \widetilde{W}_\delta^{1, p}(\Sigma, E) \oplus \R^N \to
L_\delta^p(\Sigma, \Wedge^{0, 1}T^\ast \Sigma \otimes E)
\]
surjective.
We assume the support of $\lambda$ is contained in $\Sigma_0 \setminus N_0$.
Then $\lambda$ induces a map
$\lambda_{(\zeta, r)} : \R^N \to L_\delta^p(\Sigma_{(\zeta, r)}, \Wedge^{0, 1}T^\ast \Sigma_{(\zeta, r)} \otimes E)$.
We prove the surjectivity of
\[
D_{(\zeta, r)} \oplus \lambda_{(\zeta, r)} : \widetilde{W}_\delta^{1, p}(\Sigma_{(\zeta, r)}, E)
\oplus \R^N \to L_\delta^p(\Sigma_{(\zeta, r)}, \Wedge^{0, 1}T^\ast \Sigma_{(\zeta, r)}
\otimes E)
\]
for sufficiently small $(\zeta, r^{-1})$.
Let $\{(\xi_k, h_k)\}$ be a orthonormal basis of $\Ker (D_0 \oplus \lambda)$,
where the inner product of $\Ker (D_0 \oplus \lambda)$ is defined by
\[
\langle (\xi, h), (\xi', h') \rangle
= \langle \xi, \xi' \rangle_{L^2(\Sigma_0 \setminus N_0)} +\langle h, h' \rangle_{\R^N}
\]
\begin{lem}\label{linearized gluing lemma}
There exists some constant $C>0$ such that for any sufficiently small $(\zeta, r^{-1})$,
\begin{align}
&||\xi||_{\widetilde{W}_\delta^{1,p}(\Sigma_{(\zeta, r)})} + |h|_{\R^N} \notag \\
&\ \leq C \bigl( ||D_{(\zeta, r)} \xi + \lambda_{(\zeta, r)} h||
_{L_\delta^p(\Sigma_{(\zeta, r)})}
+ \sum_k |\langle \xi,\xi_k \rangle_{L^2(\Sigma_0 \setminus N_0)}
+ \langle h,h_k \rangle_{\R^N} | \bigr)
\label{inequality of linearized gluing lemma}
\end{align}
\end{lem}
\begin{proof}
We may assume $D_0 = \bound$ on some neighborhood of nodal points, and
$S_\mu(s, t) = s_\mu(t)$ for sufficiently large $|s|$ for all $\mu$ because
the Sobolev embedding $||\xi_0||_{L^\infty} \lesssim ||\xi_0||_{W^{1, p}}$ is uniform
with respect to small $(\zeta, r^{-1})$.

It is enough  prove the inequality for
$\xi \in \widetilde{W}_\delta^{1, p}(\Sigma_{(\zeta, r)}, E) \cap
C^\infty(\Sigma_{(\zeta, r)}, E)$.
We construct a section $\tilde \xi \in \widetilde{W}_\delta^{1, p}(\Sigma, E)$
from $\xi$, and apply the inequality
\begin{align}
&||\tilde \xi||_{\widetilde{W}_\delta^{1,p}(\Sigma)} + |h|_{\R^N} \notag\\
&\ \leq C \bigl( ||D_0 \tilde \xi + \lambda h||_{L_\delta^p(\Sigma)}
+ \sum_k |\langle \tilde \xi,\xi_k \rangle_{L^2(\Sigma_0 \setminus N_0)}
+ \langle h,h_k \rangle_{\R^N}| \bigr)
\label{eq of tilde xi}
\end{align}
followed from the surjectivity of $D_0 \oplus \lambda$ to $(\tilde \xi, h)$.
From this inequality, we will derive the required inequality for $(\xi, h)$.

Define $\tilde \xi|_{\Sigma_0 \setminus N_0} = \xi|_{\Sigma_0 \setminus N_0}$.
We also define $\tilde \xi = \xi$ on the neighborhood of limit circles of $\Sigma$.

Next we consider the neighborhood of the $\nu$-th nodal point.
On $N^\nu_{\zeta_\nu}$, let
\begin{align*}
\xi|_{\{|x| = \sqrt{|\zeta_\nu|}\}}
&= \sum_k a^{(\nu)}_k x^k \in L^2(\{|x|=\sqrt{|\zeta_\nu|}\}, \C^n)\\
\xi|_{\{|y| = \sqrt{|\zeta_\nu|}\}}
&= \sum_k b^{(\nu)}_k y^k \in L^2(\{|y|=\sqrt{|\zeta_\nu|}\}, \C^n)
\end{align*}
be the Fourier expansions.
Note that $b_k^{(\nu)} = a_{-k}^{(\nu)} \zeta^{-k}$.
In particular, $a_0^{(\nu)} = b_0^{(\nu)}$.
Then $\tilde \xi|_{D \cup D}$ is defined by
\begin{align*}
\tilde \xi(x, 0) &= \begin{cases}
\xi(x) - \rho_{\zeta_\nu}(x) \sum_{k < 0} a_k^{(\nu)} x^k
& \text{for } \sqrt{|\zeta_\nu|} \leq |x| \leq 1\\
\sum_{k \geq 0} a_k^{(\nu)} x^k & \text{for } |x| \leq \sqrt{|\zeta_\nu|}
\end{cases}\\
\tilde \xi(0, y) &= \begin{cases}
\xi(y) - \rho_{\zeta_\nu}(y) \sum_{k < 0} b_k^{(\nu)} y^k
& \text{for } \sqrt{|\zeta_\nu|} \leq |y| \leq 1\\
\sum_{k \geq 0} b_k^{(\nu)} y^k & \text{for } |y| \leq \sqrt{|\zeta_\nu|}
\end{cases},
\end{align*}
where $\rho_\zeta$ is defined as follows.
Let $\rho : \R_{\geq 0} \to \R_{\geq 0}$ be a smooth function such that
$\rho|_{[0, 1]} = 1$ and $\supp \rho \subset [0, 2]$, and fix a constant
$0 < \alpha < \frac{1}{2}$.
Then $\rho_\zeta$ is defined by $\rho_\zeta(z) = \rho(\frac{z}{|\zeta|^\alpha})$.
We note that
$\tilde \xi|_{D \cup D} \in W^{2, 2}(D \cup D) \subset W^{1, p}(D \cup D)$
because $\sum_k a^{(\nu)}_k x^k \in C^\infty(\{|x|=\sqrt{|\zeta_\nu|}\}, \C^n)
\subset W^{\frac{3}{2}, 2}(\{|x|=\sqrt{|\zeta_\nu|}\}, \C^n)$ and
$\sum_k b^{(\nu)}_k y^k \in C^\infty(\{|y|=\sqrt{|\zeta_\nu|}\}, \C^n)
\subset W^{\frac{3}{2}, 2}(\{|y|=\sqrt{|\zeta_\nu|}\}, \C^n)$.

Next we consider the neighborhood of the $\mu$-th joint circle.
On $N^\mu_{r_\mu}$, let
\begin{align*}
\xi|_{\{s = r_\mu\} \subset [0, r_\mu] \times S^1}
&= \sum_k a^{(\mu)}_k e^{-\lambda^{(\mu)}_k r_\mu} \phi^{(\mu)}_k(t)
\in L^2(S^1, \R^{2n})\\
\xi|_{\{s = -r_\mu\} \subset [-r_\mu, 0] \times S^1}
&= \sum_k b^{(\mu)}_k e^{\lambda^{(\mu)}_k r_\mu} \phi^{(\mu)}_k(t)
\in L^2(S^1, \R^{2n})
\end{align*}
be expansions by the eigenvectors $\phi^{(\mu)}_k$ of $A_\mu$,
where $\lambda^{(\mu)}_k$ is the eigenvalue corresponding to $\phi^{(\mu)}_k$.
Since $\{s = r_\mu\} \subset [0, r_\mu] \times S^1$ and
$\{s = -r_\mu\} \subset [-r_\mu, 0] \times S^1$ are the same circle,
$b^{(\mu)}_k = e^{-2 \lambda^{(\mu)}_k r_\mu} a^{(\mu)}_k$.
In particular, $b^{(\mu)}_k = a^{(\mu)}_k$ if $\lambda_k=0$.
Then $\tilde \xi|_{([0, \infty] \cup [-\infty, 0]) \times S_\mu^1}$ is defined by
\begin{align*}
\tilde \xi|_{[0, \infty] \times S^1} (s, t) &= \begin{cases}
\xi(s, t) - \chi_{r_\mu}(s) \sum_{\lambda_k^{(\mu)} < 0}
a_k^{(\mu)} e^{- \lambda_k^{(\mu)} s}
\phi_k^{(\mu)}(t) & 0 \leq s \leq r_\mu\\
\sum_{\lambda_k^{(\mu)} \geq 0} a_k^{(\mu)} e^{- \lambda_k^{(\mu)} s} \phi_k^{(\mu)}(t)
& r_\mu \leq s \leq \infty
\end{cases},\\
\tilde \xi|_{[-\infty, 0] \times S^1} (s, t) &= \begin{cases}
\sum_{\lambda_k^{(\mu)} \leq 0} b_k^{(\mu)} e^{- \lambda_k^{(\mu)} s} \phi_k^{(\mu)}(t)
& -\infty \leq s \leq -r_\mu\\
\xi(s, t) - \chi_{r_\mu}(s) \sum_{\lambda_k^{(\mu)} > 0}
b_k^{(\mu)} e^{-\lambda_k^{(\mu)} s}
\phi_k^{(\mu)}(t) & -r_\mu \leq s \leq 0
\end{cases},
\end{align*}
where $\chi : \R_{\geq 0} \to \R_{\geq 0}$ is a smooth function such that
$\chi|_{[0, 1/3]} = 0$ and $\chi|_{[2/3, \infty)} = 1$, and
$\chi_r : \R_{\geq 0} \to \R_{\geq 0}$ is defined by
$\chi_r(s) = \chi(\frac{s}{r})$.
It is easy to see that $\tilde \xi|_{([0, \infty] \cup [-\infty, 0]) \times S_\mu^1} \in
\widetilde{W}^{1, p}_\delta(([0, \infty] \cup [-\infty, 0]) \times S_\mu^1, \C^n)$.

We assume that $\zeta_\nu$ are sufficiently small and $r_\mu$ are sufficiently large
so that
\begin{itemize}
\item
$D_0 = \bound$ on $\{(x, y) \in N^\nu_{\zeta_\nu}; |x| \leq 2 |\zeta_\nu|^\alpha
\text{ or } |y| \leq 2|\zeta_\nu|^\alpha\}$, and
\item
$S_\mu(s, t) = S_\mu(t)$ on $[\frac{1}{3} r_\mu, \infty] \times S_\mu^1 \cup
[-\infty, -\frac{1}{3} r_\mu] \times S_\mu^1$.
\end{itemize}
For each $\mu$, define $v_\mu \in \Ker A_\mu$ by
\[
v_\mu = \sum_{\lambda_k^{(\nu)} = 0} a_k^{(\mu)} \phi_k^{(\mu)}(t).
\]
We also define $v_{\pm\infty_i} \in \Ker A_\mu$ for limit circles $\pm\infty_i$ by
the condition
\[
\xi - \sum_\mu \beta_\mu v_\mu
- \sum_{\pm\infty_i} \beta_{\pm\infty_i} v_{\pm\infty_i}
\in W_\delta^{1,p}(\Sigma_{(\zeta, r)}, E).
\]

We can easily check the following inequalities, where $C > 0$ is some constant and
$0 < \epsilon \leq 1$ is arbitrary.
($C > 0$ does not depend on $\epsilon$.)
\begin{align}
&C \Bigl|\Bigl|\tilde \xi - \sum_\mu \beta_\mu v_\mu
- \sum_{\pm\infty_i} \beta_{\pm\infty_i} v_{\pm\infty_i}\Bigr|\Bigr|
_{W_\delta^{1,p}(\Sigma_0)} \notag \\
&\geq \epsilon \Bigl(\Bigl|\Bigl|\xi - \sum_\mu \beta_\mu v_\mu
- \sum_{\pm\infty_i} \beta_{\pm\infty_i} v_{\pm\infty_i}\Bigr|\Bigr|
_{W_\delta^{1,p}(\Sigma_{(\zeta, r)})} \notag \\
&\quad \hph{\epsilon \Bigl(}
- \sum_\nu \Bigl|\Bigl|\rho_{\zeta_\nu}(x) \sum_{k < 0} a_k^{(\nu)} x^k\Bigr|\Bigr|
_{W^{1,p}(\sqrt{|\zeta_\nu|} \leq |x| \leq 1)} \notag \\
&\quad \hph{\epsilon \Bigl(}
- \sum_\nu \Bigl|\Bigl|\rho_{\zeta_\nu}(y) \sum_{k < 0} b_k^{(\nu)} y^k\Bigr|\Bigr|
_{W^{1,p}(\sqrt{|\zeta_\nu|} \leq |y| \leq 1)} \notag \\
&\quad \hph{\epsilon \Bigl(}
- \sum_\mu \Bigl|\Bigl|\chi_{r_\mu}(s) \sum_{\lambda_k^{(\mu)} < 0} a_k^{(\mu)}
e^{- \lambda_k^{(\mu)} s}
\phi_k^{(\mu)}(t)\Bigr|\Bigr|
_{W_{\delta_\mu}^{1,p}([0, r_\mu] \times S^1)} \notag \\
&\quad \hph{\epsilon \Bigl(}
- \sum_\mu \Bigl|\Bigl|\chi_{r_\mu}(s) \sum_{\lambda_k^{(\mu)} > 0} b_k^{(\mu)}
e^{-\lambda_k^{(\mu)} s}
\phi_k^{(\mu)}(t)\Bigr|\Bigr|
_{W_{\delta_\mu}^{1,p}([-r_\mu, 0] \times S^1)}
\Bigr) \notag \\
&\quad
+ \sum_\nu \Bigl|\Bigl|\sum_{k \geq 0} a_k^{(\nu)} x^k
\Bigr|\Bigr|
_{W^{1,p}(|x| \leq \sqrt{|\zeta_\nu|})}
+ \sum_\nu \Bigl|\Bigl|\sum_{k \geq 0} b_k^{(\nu)} y^k
\Bigr|\Bigr|
_{W^{1,p}(|y| \leq \sqrt{|\zeta_\nu|})} \notag \\
&\quad
+ \sum_\mu \Bigl|\Bigl|\sum_{\lambda_k^{(\mu)} > 0} a_k^{(\mu)}
e^{- \lambda_k^{(\mu)} s} \phi_k^{(\mu)}(t)
\Bigr|\Bigr|
_{W_{\delta_\mu}^{1,p}([r_\mu, \infty] \times S^1)} \notag \\
&\quad
+ \sum_\mu \Bigl|\Bigl|\sum_{\lambda_k^{(\mu)} < 0} b_k^{(\mu)}
e^{- \lambda_k^{(\mu)} s} \phi_k^{(\mu)}(t)
\Bigr|\Bigr|
_{W_{\delta_\mu}^{1,p}([-\infty, -r_\mu] \times S^1)}
\label{tilde xi geq xi}
\end{align}
\begin{align}
&||D_0 \tilde \xi + \lambda h||_{L_\delta^p(\Sigma)} \notag \\
&\leq ||D_{(\zeta, r)} \xi + \lambda_{(\zeta, r)} h||_{L_\delta^p(\Sigma_{(\rho, r)})}
\notag \\
&\quad
+ C\Bigl(\sum_\nu |\zeta_\nu|^{-\alpha}
\Bigl|\Bigl|\sum_{k < 0} a_k^{(\nu)} x^k\Bigr|\Bigr|_{L^p
(|\zeta_\nu|^\alpha \leq |x| \leq 2 |\zeta_\nu|^\alpha)} \notag \\
&\quad \hph{+ C\Bigl( }
+ \sum_\nu |\zeta_\nu|^{-\alpha} \Bigl|\Bigl|\sum_{k < 0} b_k^{(\nu)} y^k\Bigr|\Bigr|_{L^p
(|\zeta_\nu|^\alpha \leq |y| \leq 2 |\zeta_\nu|^\alpha)} \notag \\
&\quad \hph{+ C\Bigl( }
+ \sum_\mu \frac{1}{r_\mu} \Bigl|\Bigl|\sum_{\lambda_k^{(\mu)} < 0} a_k^{(\mu)}
e^{-\lambda_k^{(\mu)} s} \phi_k^{(\mu)}(t)\Bigr|\Bigr|
_{L_{\delta_\mu}^p([0, r_\mu] \times S^1)}
\notag \\
&\quad \hph{+ C\Bigl( }
+ \sum_\mu \frac{1}{r_\mu} \Bigl|\Bigl|\sum_{\lambda_k^{(\mu)} > 0} b_k^{(\mu)}
e^{-\lambda_k^{(\mu)} s} \phi_k^{(\mu)}(t)\Bigr|\Bigr|
_{L_{\delta_\mu}^p([-r_\mu, 0] \times S^1)}
\Bigr)
\label{D tilde xi leq D xi}
\end{align}

Apply inequality (\ref{eq of tilde xi}) to $\tilde \xi \in
\widetilde{W}^{1, p}_\delta(\Sigma_0)$,
and use (\ref{tilde xi geq xi}) for sufficiently small $\epsilon > 0$,
(\ref{D tilde xi leq D xi}) and the following two lemmas
(Lemma \ref{disk cancellation inequality} and \ref{cylinder cancellation inequality}).
Then we can easily see that there exists some constant $C > 0$ such that
for any sufficiently small $(\zeta, r^{-1})$,
\begin{align*}
&\Bigl|\Bigl|\xi - \sum_\mu \beta_\mu v_\mu
- \sum_{\pm\infty_i} \beta_{\pm\infty_i} v_{\pm\infty_i}\Bigr|\Bigr|
_{W_\delta^{1,p}(\Sigma_{(\zeta, r)})} \\
&+ \sum_\mu ||v_\mu||_{\Ker A_\mu}
+ \sum_{\pm\infty_i} ||v_{\pm\infty_i}||_{\Ker A_{\pm\infty_i}}
+ |h|_{\R^N} \\
&\ \leq C \bigl( ||D_{(\zeta, r)} \xi + \lambda_{(\zeta, r)} h||
_{L_\delta^p(\Sigma_{(\zeta, r)})}
+ \sum_k |\langle \xi,\xi_k \rangle_{L^2(\Sigma_0 \setminus N_0)}
+ \langle h,h_k \rangle_{\R^N} | \bigr).
\end{align*}
We do not use any estimates of
$a_k^{(\nu)}$, $b_k^{(\nu)}$, $a_k^{(\mu)}$ or $b_k^{(\mu)}$ by $\xi$.
Lemma \ref{disk cancellation inequality} and \ref{cylinder cancellation inequality}
imply that these terms cancel each other.
(\ref{inequality of linearized gluing lemma}) follows from the above inequality.

\end{proof}

\begin{lem}
\label{disk cancellation inequality}
For any $2<p<\infty$ and $0 < \alpha < \frac{1}{2}$, there exists some $C>0$ such that
for any $\zeta \in D$ and any two sequences $(a_k)_{k \in \Z}$ and $(b_k)_{k \in \Z}$
such that $b_k = a_{-k} \zeta^{-k}$, the following inequalities hold true.
\[
\Bigl|\Bigl|\sum_{k<0}a_k z^k\Bigr|\Bigr|_{W^{1,p}(\sqrt{|\zeta|}\leq |z| \leq 1)}
\leq C \Bigl|\Bigl|\sum_{k\geq0} b_k z^k\Bigr|\Bigr|_{W^{1,p}
(|z|\leq \sqrt{|\zeta|})}
\]
\[
|\zeta|^{-\alpha} \Bigl|\Bigl|\sum_{k<0} a_k z^k\Bigr|\Bigr|_{L^p
(|\zeta|^\alpha \leq |z| \leq 2|\zeta|^\alpha)}
\leq C |\zeta|^{(1-\frac{2}{p})(1-2\alpha)} \Bigl|\Bigl|\sum_{k\geq0} b_k z^k
\Bigr|\Bigr|_{W^{1,p}(|z|\leq \sqrt{|\zeta|})}
\]
\end{lem}
\begin{proof}
Put $f(z) = \sum_{k\geq 0} b_k z^k$.
Then $\sum_{k<0} a_k z^k = f\bigl(\frac{\zeta}{z}\bigr) - b_0$.
Therefore $\frac{d}{dz}(\sum_{k<0} a_k z^k) = f'\bigl(\frac{\zeta}{z}\bigr) \bigl( -
\frac{\zeta}{z^2}\bigr)$.
This implies for any $0 < \alpha \leq \frac{1}{2}$,
\begin{align*}
\int_{|\zeta|^\alpha \leq |z| \leq 1} \biggl|\frac{d}{dz}\Bigl(\sum_{k<0} a_k z^k\Bigr)
\biggr|^p |dz|^2
&= \int_{|\zeta| \leq |w| \leq |\zeta|^{1-\alpha}} |f'(w)|^p
\biggl(\frac{|w|^2}{|\zeta|}\biggr)^{p-2} |dw|^2\\
&\leq
|\zeta|^{(p-2)(1-2\alpha)} \Bigl|\Bigl|\sum_{k\geq 0} b_k z^k\Bigr|\Bigr|_{W^{1,p}
(|z|\leq \sqrt{|\zeta|})}^p.
\end{align*}
On the other hand,
Poincar\'{e}'s inequality on $S^1$ implies
\begin{align*}
\Bigl|\Bigl|\sum_{k<0} a_k z^k\Bigr|\Bigr|_{L^p(|\xi|^\alpha \leq |x| \leq \rho)}^p
&= \int_{|\zeta|^\alpha}^\rho \int_0^{2\pi} \Bigl|\sum_{k<0}a_k r^k e^{\sqrt{-1}k\theta}
\Bigr|^p rdrd\theta\\
&\leq C \int_{|\zeta|^\alpha}^\rho\int_0^{2\pi}
\Bigl|\sum_{k<0}ka_kr^k e^{\sqrt{-1}k\theta}\Bigr|^p rdrd\theta\\
&= C \int_{|\zeta|^\alpha}^\rho\int_0^{2\pi}
\Bigl|\sum_{k<0}ka_kr^{k-1} e^{\sqrt{-1}(k-1)\theta}\Bigr|^p r^{p+1} drd\theta\\
&\leq C \rho^p \biggl|\biggl|\frac{d}{dz}\Bigl(\sum_{k<0} a_k z^k\Bigr)
\biggr|\biggr|_{L^p(|\xi|^\alpha \leq |z|\leq \rho)}^p
\end{align*}
for $\rho = 2|\zeta|^\alpha$ or $1$.
The first of the claimed inequalities is proved by substituting $\alpha = \frac{1}{2}$ and
$\rho = 1$, and the second is proved by substituting $\rho = 2|\zeta|^\alpha$.
\end{proof}

\begin{lem}
\label{cylinder cancellation inequality}
Let $(\phi_k)$ be a family of $W^{1, 2}$-functions on $S^1$.
Let $\delta > 0$ be a positive constant and $(\lambda_k)$ be a sequence of
real numbers such that $\lambda_k < -\delta$.
Then for any $1 < p < \infty$, $r > 0$ and any two sequences $(a_k)$ and $(b_k)$
such that $b_k = e^{-2\lambda_k r} a_k$, the following inequality holds true.
\[
\Bigl|\Bigl|\sum_{\lambda_k<0}a_k e^{-\lambda_k s} \phi_k(t)\Bigr|\Bigr|_{W_\delta^{1,p}
([0,r] \times S^1)}
\leq
\Bigl|\Bigl|\sum_{\lambda_k<0} b_k e^{-\lambda_k s} \phi_k(t)\Bigr|\Bigr|_{W_\delta^{1,p}
((-\infty,-r] \times S^1)}
\]
\end{lem}
\begin{proof}
The $L^p_\delta$-norm is estimated by
\begin{align*}
\int_0^r \Bigl|\sum_{\lambda_k<0}a_k e^{-\lambda_k s} \phi_k(t)\Bigr|^p
e^{p\delta s} dsdt
&= \int_{-2r}^{-r} \Bigl|\sum_{\lambda_k<0} b_k e^{-\lambda_k s} \phi_k(t)\Bigr|^p
e^{p\delta (s + 2r)} dsdt\\
&\leq \int_{-2r}^{-r} \Bigl|\sum_{\lambda_k<0} b_k e^{-\lambda_k s} \phi_k(t)
\Bigr|^p e^{-p\delta s} dsdt.
\end{align*}
Similarly, we can estimate
\[
\int_0^r \Bigl|\partial_s\Bigl(\sum_{\lambda_k<0}a_k e^{-\lambda_k s}
\phi_k(t)\Bigr)\Bigr|^p e^{p\delta s} dsdt
\]
and
\[
\int_0^r \Bigl|\partial_t\Bigl(\sum_{\lambda_k<0}a_k e^{-\lambda_k s}
\phi_k(t)\Bigr)\Bigr|^p e^{p\delta s} dsdt
\]
by the corresponding terms for $b_k$.
\end{proof}

\begin{rem}\label{uniform regularity}
The same argument implies that interior regularity property of $D_{(\zeta, r)}$ is
uniform with respect to small $(\zeta, r)$ on a neighborhood of a nodal point or
a imaginary circle.
\end{rem}


%% file: SFT-04_Smoothness_Embeddings.tex
%
%

\subsection{Smoothness of Kuranishi neighborhoods}
\label{smoothness}
In this section, we prove that if we give $\widetilde{X}$ a stronger differential
structure and give the product differential structure to
$\hat V = \mathring{X} \times B_\epsilon(0) \subset \mathring{X} \times
\Ker DF^{(0, b_0)}_{(0,0)}$, then
\begin{align*}
\hat V &\inj \mathring{X} \times C^l(\Sigma_0 \setminus N_0,
(\R_1 \cup \R_2 \cup \dots \cup \R_k) \times Y) \times E^0\\
(a, b, x) &\mapsto (a, b, \Phi_{a, b}(\xi_x)|_{\Sigma_0 \setminus N_0}, h_x)
\end{align*}
is a smooth embedding for any $l$.
More precisely, we prove that for any $N \geq 1$, we can chose
a stronger differentiable structure of $\widetilde{X}$ such that the map
is of class $C^N$.
Note that we have already proved the continuity of the above map by
(\ref{continuity of phi}).
We also note that once we prove that this is a smooth embedding for $l = 1$,
then it follows that for any $l \geq 1$
and any $\widetilde{N}_0 \supset N_0$,
\[
\hat V \inj \mathring{X} \times C^l(\Sigma_0 \setminus \widetilde{N}_0, (\R_1 \cup \R_2 \cup
\dots \cup \R_k) \times Y) \times E^0
\]
is also a smooth embedding provided that $\widetilde{N}_0$ does not cover any
irreducible components of $\Sigma_0$.

First we explain about the strong differential structure of $\widetilde{X}$.
It is based on the following Lemma.
\begin{lem}
Let $V \subset \C^n$ be an open set and $D \subset \C$ be a disk.
Assume a holomorphic function $f(w, \zeta_1, \dots, \zeta_l) : V \times D^l \to \C$
satisfies $\{(w, \zeta); f(w, \zeta) = 0\} = \bigcup_i \{ \zeta_i = 0\}$.
Define $\varphi_\alpha (r e^{\sqrt{-1} \theta}) = r^\alpha e^{\sqrt{-1} \theta}$
for $\alpha \geq 1$.
Then $\varphi_\alpha^{-1} \circ f(w, \varphi_\alpha (\zeta_1), \dots,
\varphi_\alpha (\zeta_l)) : V \times D^l \to \C$ is of class $C^{\lfloor \alpha \rfloor}$.
(If $\alpha = 2 N + 1$ for some $N \in \Z_{\geq 0}$ then it is real analytic.)
\end{lem}
\begin{proof}
There exists a holomorphic function $g : V \times D^l \to \C \setminus 0$ such that
$f(w,\zeta) = \zeta_1^{k_1} \dots \zeta_l^{k_l} g(w,\zeta)$ for some $k_i \geq 1$.
Then
\[
\varphi_\alpha^{-1}\bigl(f(w,\varphi_\alpha (\zeta_1), \dots,
\varphi_\alpha (\zeta_l))\bigr)
= \zeta_1^{k_1} \dots \zeta_l^{k_l} \varphi_\alpha^{-1}
\bigl(g(w,\varphi_\alpha (\zeta_1), \dots, \varphi_\alpha (\zeta_l))\bigr),
\]
where we have used $\varphi_\alpha (ab) =\varphi_\alpha(a) \varphi_\alpha(b)$.
Since  $\varphi_\alpha$ is of class $C^{\lfloor \alpha \rfloor}$ and
$\varphi_\alpha^{-1} : \C\setminus 0 \to \C\setminus 0$ is real analytic,
$\varphi_\alpha^{-1} \circ f(w, \varphi_\alpha (\zeta_1), \dots,
\varphi_\alpha (\zeta_l))$ is of class $C^{\lfloor \alpha \rfloor}$.
\end{proof}

For any $\alpha \gg 0$ and $\beta \gg 0$,
a new differential structure of $\widetilde{X}$ is defined by
the coordinate
\[
\widetilde{X} \subset \J_0 \times D^{l_0} \times \widetilde{D}^{l_1}
\to \J_0 \times D^{l_0} \times ([0, 1] \times S^1)^{l_1}
\]
{\abovedisplayskip=0pt
\[
(j, (\zeta_\nu = \rho_\nu^2 e^{2\sqrt{-1} \varphi_\nu})_\nu,
(
\rho_\mu^{2\pi} e^{2\pi \sqrt{-1} \varphi_\mu})_\mu) 
\mapsto (j, (\hat \zeta_\nu = \hat \rho_\nu^2 e^{2\sqrt{-1} \varphi_\nu})_\nu,
(\hat \rho_\mu, \varphi_\mu)_\mu)
\]}
defined by $\rho_\nu = \hat \rho_\nu^\alpha$ and
$\rho_\mu = \hat \rho_\mu^{\beta_\mu}$,
where $\beta_\mu = L_\mu^{-1} \beta$ ($L_\mu = L_{\gamma_\mu}$ is
the period of the periodic orbit on $S_\mu^1$).
The above lemma implies that this differential structure is independent of
the local description of the universal family $(\widetilde{P} \to \widetilde{X},
Z \cup Z^+ \cup Z^{++})$ given by a decomposition of $\Sigma_0$
since in any description, $\{\zeta_\nu\}$ consists of the curves which have
$\nu$-th nodal point and is preserved by the coordinate change.
The reason why we use the indeces $\beta_\mu = L_\mu^{-1} \beta$ depending on
$\mu$ is to make $\mathring{X} \subset \widetilde{X} \times \prod_\mu \R_\mu$
a submanifold.
(Recall that $\mathring{X}$ is defined by the condition that
$\rho_\mu^{L_\mu} e^{- b_\mu}$ does not depend on $\mu \in M_i$ for each
$i = 1, 2, \dots, k-1$.)

We fix large constants $\alpha \geq 1$ and $\beta > 0$, and use the differential
structure of $\widetilde{X}$ defined by the same $\alpha$ and $\beta$ for all
Kuranishi neighborhoods of $\widehat{\M}$ for each pre-Kuranishi structure of
$\widehat{\M}$.

Let $\mathring{X} = \coprod_{\Pi, \Pi'} \mathring{X}_{\Pi, \Pi'}$
be the decomposition defined by
\begin{align*}
\mathring{X}_{\Pi, \Pi'}
= \{ (a,b) \in \mathring{X};\, & \rho_\mu \neq 0 \text{ for all } \mu \in M_i
\text{ if and only if } i \in \Pi\\
& \zeta_\nu \neq 0 \text{ if and only if } \nu \in \Pi'\},
\end{align*}
where $\Pi \subset \{1, 2, \dots, k-1\}$ and
$\Pi'$ is a subset of nodal points of $\Sigma_0$.
We prove the differentiability of $\phi$ on each
$\mathring{X}_{\Pi, \Pi'} \times B_\epsilon(0)$ and
investigate its behavior near the boundary.

Fix one point $(a,b) \in \mathring{X}_{\Pi, \Pi'}$ and consider another point
$(\tilde a, \tilde b) \in \mathring{X}_{\Pi, \Pi'}$ close to $(a,b)$.
To investigate the behavior of the differential, we identify $\widetilde{P}_{\tilde a}$
and $\widetilde{P}_a$ by the piecewise smooth map $\Psi$ defined as follows.

On each $[-1, 0] \times S^1_\mu \subset [-1, -\frac{1}{2}\log \rho_\mu] \times
S^1_\mu$,
\begin{align*}
\Psi : [-1, 0] \times S^1_\mu &\to [-1, 0] \times S^1_\mu\\
(s,t) &\mapsto (\tilde s, \tilde t) = (s, t)
\end{align*}
is given by the identity map,
and on each $[0, -\frac{1}{2}\log \rho_\mu] \times S^1_\mu \subset
[-1,-\frac{1}{2}\log \rho_\mu] \times S^1_\mu$,
\begin{align*}
\Psi : [0, - {\textstyle \frac{1}{2}} \log \rho_\mu] \times S^1_\mu
&\to [0, - {\textstyle \frac{1}{2}} \log \tilde \rho_\mu] \times S^1_\mu\\
(s,t) &\mapsto (\tilde s, \tilde t)
\end{align*}
is defined by
\[
\frac{e^{-\kappa_\mu \tilde s} -
\tilde \rho_\mu^{\kappa_\mu/2}}{1 - \tilde \rho_\mu^{\kappa_\mu/2}}
= \frac{e^{-\kappa_\mu s} - \rho_\mu^{\kappa_\mu/2}}{1 - \rho_\mu^{\kappa_\mu/2}},
\quad \tilde t = t.
\]
$\Psi$ on each $[\frac{1}{2} \log \rho_\mu, 1] \times S^1$ is defined similarly.
For simplicity of notation, we denote
$[0, -\frac{1}{2}\log \rho_\mu] \times S^1_\mu \cup
[\frac{1}{2} \log \rho_\mu, 0] \times S^1_\mu$
by $N^\mu_{a, b}$.
Recall the definition of the approximate solutions and note that
$\Psi$ satisfies $v^\mu_{\tilde a, \tilde b} \circ \Psi = v^\mu_{a, b}$ on $N^\mu_{a, b}$.

On each $N^\nu_{a, b}$,
\begin{align*}
N^\nu_{a, b} = \{(x,y)\in \overline{D}\times\overline{D};xy=\zeta_\nu\}
&\to \{(\tilde x,\tilde y)\in \overline{D}\times\overline{D};\tilde x\tilde y=\tilde\zeta_\nu\}
=N^\nu_{\tilde a, \tilde b} \\
(x,y)&\mapsto (\tilde x,\tilde y)
\end{align*}
is defined by
\begin{itemize}
\item
$\tilde x= \tilde r e^{\sqrt{-1}\tilde\theta}$ if $|x| \geq \sqrt{|\zeta_\nu|}$
and $x = r e^{\sqrt{-1}\theta}$
\item
$\tilde y= \tilde r e^{\sqrt{-1}\tilde\theta}$ if $|y| \geq \sqrt{|\zeta_\nu|}$
and $y = r e^{\sqrt{-1}\theta}$
\end{itemize}
where $\tilde r$ and $\tilde \theta$ is defined by
\[
\frac{\tilde r - \tilde \rho_\nu}{1 - \tilde \rho_\nu}
= \frac{r - \rho_\nu}{1 - \rho_\nu}, \quad \quad
\tilde \theta = \theta + \beta_\nu(r) (\tilde \varphi_\nu - \varphi_\nu),
\quad \beta_\nu(r)
=\max \biggl(\frac{2 \rho_\nu - r}{\rho_\nu}, 0 \biggr),
\]
where $\rho_\nu$, $\varphi_\nu$, $\tilde \rho_\nu$ and $\tilde \varphi_\nu$
are defined by $\zeta_\nu = \rho_\nu^2 e^{2 \sqrt{-1} \varphi_\nu}$ and
$\tilde \zeta_\nu = \tilde \rho_\nu^2 e^{2 \sqrt{-1} \tilde \varphi_\nu}$.

On $\Sigma_0 \setminus N_0$,
$\Psi|_{\Sigma_0 \setminus N_0} = \id$.

Then under this identification, we consider $F^{(\tilde a, \tilde b)}$ as a map
\begin{align*}
&F^{(\tilde a,\tilde b)} : \widetilde{W}^{1,p}_\delta( \widetilde{P}_a, u_{a, b}^\ast T \hat Y)
\oplus E^0\\
&\to L^p(\Sigma_0 \setminus N_0, \Wedge^{0,1} T^\ast \Sigma_0 \otimes_\C
u_0^\ast T \hat Y)\\
&\quad \oplus \bigoplus_\mu (L^p_\delta ([-1, - {\textstyle \frac{1}{2}} \log \rho_\mu]
\times S^1, \R^{2n})
\oplus L^p_\delta ([{\textstyle \frac{1}{2}} \log \rho_\mu, 1] \times S^1, \R^{2n}))\\
&\quad \oplus \bigoplus_{+\infty_i} L^p_\delta([0, \infty] \times S^1, \R^{2n})
\oplus \bigoplus_{-\infty_i} L^p_\delta([-\infty, 0] \times S^1, \R^{2n})\\
&\quad \oplus \bigoplus_\nu (L^p (\{ x \in D; |x| \geq \sqrt{|\zeta_\nu|}\}, \R^{2n})
\oplus L^p (\{ y \in D; |y| \geq \sqrt{|\zeta_\nu|}\}, \R^{2n})).
\end{align*}

On each $[-1, 0] \times S^1_\mu \subset [-1, -\frac{1}{2} \log \rho_\mu] \times
S^1_\mu$,
\begin{align*}
F^{(\tilde a,\tilde b)} (\xi, h)
&= 
\partial_s (v_{a,b}^{\mu, \mathrm{left}} + \xi)
+ \widetilde{J}^\mu_t(v_{a,b}^{\mu, \mathrm{left}} + \xi)
\partial_t (v_{a,b}^{\mu, \mathrm{left}} + \xi)
+ f^\mu_t(v_{a,b}^{\mu, \mathrm{left}} + \xi)\\
& \quad + \frac{1}{2} (b_\mu -b_\mu^0) \chi'(s) \partial_\sigma
+ \frac{1}{2} \tilde \varphi_\mu \chi'(s) (g^\mu_t(v_{a,b}^{\mu, \mathrm{left}} + \xi)
+ \partial_t(v_{a,b}^{\mu, \mathrm{left}} + \xi)).
\end{align*}

On each $[0, -\frac{1}{2} \log \rho_\mu] \times S^1_\mu \subset
[-1, -\frac{1}{2} \log \rho_\mu] \times S^1_\mu$,
\begin{align*}
F^{(\tilde a,\tilde b)} (\xi, h)
&= \biggl(1 + \biggl(\frac{1 - \rho_\mu^{\kappa_\mu/2}}
{1 - \tilde \rho_\mu^{\kappa_\mu/2}}
\tilde \rho_\mu^{\kappa_\mu/2} - \rho_\mu^{\kappa_\mu/2}\biggr)
e^{\kappa_\mu s}\biggr)
\partial_s (v_{a,b}^{\mu, \mathrm{left}} + \xi)\\
&\quad + \widetilde{J}^\mu_t(v_{a,b}^{\mu, \mathrm{left}} + \xi)
\partial_t (v_{a,b}^{\mu, \mathrm{left}} + \xi)
+ f^\mu_t(\pi_Y(v_{a,b}^{\mu, \mathrm{left}} + \xi)).
\end{align*}

On each $N^{\nu, \text{left}}_{a, b} = \{x \in D; |x| \geq \sqrt{|\zeta_\nu|}\}
\subset N^\nu_{a, b}$,
\begin{align*}
F^{(\tilde a,\tilde b)} (\xi, h)
=& \frac{1 - \rho_\nu}{1 - \tilde \rho_\nu}
\biggl((\partial_r v_{a, b}^\nu)_{\tilde \varphi_\nu} + \partial_r \xi
+ \frac{\tilde \varphi_\nu - \varphi_\nu}{\rho_\nu}
1_{\{\rho_\nu \leq r_\nu \leq 2 \rho_\nu\}} \partial_\theta \xi
\biggr)\\
&+ \widetilde{J}^\nu((v_{a, b}^\nu)_{\tilde \varphi_\nu} + \xi)
\frac{(\partial_\theta v_{a, b}^\nu)_{\tilde \varphi_\nu} + \partial_\theta \xi}{\tilde r}
\end{align*}
where
$w_{\tilde \varphi_\nu}(r e^{\sqrt{-1} \theta})
= w(r e^{\sqrt{-1}(\theta + \beta_\nu(r) (\tilde \varphi_\nu - \varphi_\nu))})$
for $w = \partial_r v_{a, b}^\nu$, $v_{a, b}^\nu$ or $\partial_\theta v_{a, b}^\nu$.

By the same equations, we can define $F^{(\tilde a,\tilde b)}$ for all
$(\tilde a, \tilde b) \in \widetilde{X}_{\Pi, \Pi'} \times \prod_{\mu} \R$
close to $(a, b)$,
where $\widetilde{X}_{\Pi, \Pi'} \subset \widetilde{X}$ is defined
as $\mathring{X}_{\Pi, \Pi'} \subset \mathring{X}$, that is,
\begin{align*}
\widetilde{X}_{\Pi, \Pi'}
= \{ a \in \widetilde{X};\, & \rho_\mu \neq 0 \text{ if and only if } \mu \in
\bigcup_{i \in \Pi} M_i,\\
& \zeta_\nu \neq 0 \text{ if and only if } \nu \in \Pi'\}.
\end{align*}
In the following lemma, we regard $\tilde a \in \widetilde{X}_{\Pi, \Pi'}$ and
$\tilde b \in \prod_{\mu} \R$ as independent parameters by extending
$F^{(\tilde a,\tilde b)}$ to
$(\tilde a, \tilde b) \in \widetilde{X}_{\Pi, \Pi'} \times \prod_{\mu} \R$
as above and estimate the derivatives at $(a, b)$.
We note that
$\partial_{\tilde \rho_\mu}^k F^{(a, b)}
= \partial_{\tilde \rho_\mu}^k F^{(\tilde a, \tilde b)}|_{(\tilde a, \tilde b) = (a, b)}$
vanishes on the complement of $N_{a, b}^\mu$ for $k > 0$, and
$\partial_{\tilde \rho_\nu}^k \partial_{\tilde \varphi_\nu}^l F^{(a, b)}$ vanishes
on the complement of $N_{a, b}^\nu$ for $(k, l) \neq (0, 0)$.
\begin{lem}\label{estimates of implicit function 0}\ 
\begin{enumerate}[label=\normalfont(\roman*)]
\item
For any $0 < \delta_\mu \leq \delta'_\mu < \kappa_\mu$,
$2 < p < \infty$, $i\in \Pi$, $\mu \in M_i$ and $k > 0$,
there exists some constant $C>0$ and $c_0 > 0$
such that the following hold for $||\xi||_{\widetilde{W}^{1,p}_{\delta'}} \leq c_0$.
\begin{align*}
||(\partial_{\tilde \rho_\mu}^k F^{(a, b)})(\xi, h)||_{L_{\delta_\mu}^p(N_{a, b}^\mu)}
&\leq C \rho_\mu^{(\delta'_\mu - \delta_\mu)/2 - k},\\
||(D \partial_{\tilde \rho_\mu}^k F^{(a, b)})_{(\xi, h)}
(\hat \xi, \hat h)||_{L_{\delta_\mu}^p(N_{a, b}^\mu)}
&\leq C \rho_\mu^{(\delta'_\mu - \delta_\mu)/2 - k}
||\hat \xi||_{\widetilde{W}_{\delta'_\mu}^{1, p}(N_{a, b}^\mu)},\\
D^m \partial_{\tilde \rho_\mu}^k F^{(a, b)} &\equiv 0 \quad (m \geq 2).
\end{align*}
\item
For any $2 < p \leq q < \infty$,
$\nu \in \Pi'$ and any $(m, k, l)$ such that $(k, l) \neq (0, 0)$,
there exists some constants $C > 0$ and $c_0 > 0$ such that
the following holds for $||\xi||_{W^{1, q}(N_{a, b}^\nu)} \leq c_0$.
\begin{align*}
||(D^m \partial_{\tilde \rho_\nu}^k \partial_{\tilde \varphi_\nu}^l F^{(a, b)})_{(\xi, h)}
(\hat \xi^{(m)}, \hat h^{(m)}) (\hat \xi^{(m-1)}, \hat h^{(m-1)}) \dots
(\hat \xi^{(1)}, \hat h^{(1)})||_{L^p(N_{a, b}^\nu)}\\
\leq C \rho_\nu^{2/p - 2/q - k} \prod_{i = 1}^m ||\hat \xi^{(i)}||_{W^{1,q}(N_{a, b}^\nu)}.
\end{align*}
\item
For any $2 < p < \infty$, $0< \delta < \delta_0$,
and any multi-index $(m, k_j, (l_\mu), k_b = (k_{b_\mu}))$,
there exists some constants $C > 0$ and $c_0 > 0$ such that
the following holds for $||\xi||_{\widetilde{W}^{1,p}_\delta} + |h|_E \leq c_0$.
\begin{align*}
||(D^m \partial_{\tilde \jmath}^{k_j} \partial_{(\tilde \varphi_\mu)}^{(l_\mu)}
\partial_{\tilde b}^{k_b} F^{(a, b)})_{(\xi,h)}
(\hat \xi^{(m)}\!, \hat h^{(m)}) (\hat \xi^{(m-1)}\!, \hat h^{(m-1)})\! \dots \!
(\hat \xi^{(1)}\!, \hat h^{(1)})||_{L^p_\delta(\widetilde{P}_a)}\\
\leq C \prod_{i = 1}^m (||\hat \xi^{(i)}||_{\widetilde{W}^{1,p}_\delta(\widetilde{P}_a)}
+ |\hat h^{(i)}|_E).
\end{align*}
\end{enumerate}
(See Appendix \ref{diff notation} for our notation of differential.)
\end{lem}
\begin{rem}\label{independent}
Note that
$\partial_{\tilde \rho_\mu} \partial_{\tilde \rho_{\mu'}}F^{(a, b)} (\xi, h) = 0$
if $\mu \neq \mu'$.
Similarly, the differential of $F$ with respect to two parameters
which correspond to different pieces of the curve vanishes.
Hence the above lemma is enough for the estimate of the differentials of $F$.
\end{rem}
\begin{proof}
(i)
We use a change of variable $\mathring{\rho}_\mu = (\tilde \rho_\mu)^{\kappa_\mu/2}$.
Then on $[0, -\frac{1}{2} \log \rho_\mu] \times S^1$,
\[
\partial_{\mathring{\rho}_\mu}^k F^{(\tilde a,\tilde b)}(\xi, h)|_{(\tilde a,\tilde b) = (a,b)}
= k ! (1 - \rho_\mu^{\kappa_\mu/2})^{-k} e^{\kappa_\mu s}
\partial_s (v_{a, b}^{\mu, \text{lef}} + \xi).
\]
Since $|\partial_s v_0^{\mu, \text{lef}}(s, t)| \lesssim e^{- \delta'_{0, \mu} s}$
on $[0, \infty) \times S^1$ for any $\kappa_\mu < \delta'_{0, \mu} < \delta_{0, \mu}$,
\[
\partial_s v_{a, b}^{\mu, \text{lef}}(s, t)
= \frac{e^{-\kappa_\mu s}}{e^{-\kappa_\mu s} - \rho_\mu^{\kappa_\mu / 2}}
\partial_s v_0^{\mu, \text{lef}}\biggl(- \frac{1}{\kappa_\mu} \log
\biggl(\frac{e^{-\kappa_\mu s} - \rho_\mu^{\kappa_\mu / 2}}
{1 - \rho_\mu^{\kappa_\mu / 2}}\biggr), t
\biggr)
\]
satisfies
\begin{align*}
\bigl|\partial_s v_{a, b}^{\mu, \text{lef}}\bigr| e^{\delta'_\mu s}
&\lesssim \frac{e^{-\kappa_\mu s}}{e^{-\kappa_\mu s} - \rho_\mu^{\kappa_\mu / 2}}
\cdot \biggl(\frac{e^{-\kappa_\mu s} - \rho_\mu^{\kappa_\mu / 2}}
{1 - \rho_\mu^{\kappa_\mu / 2}}\biggr)^{\delta'_{0, \mu} / \kappa_\mu}
e^{\delta'_\mu s} \\
&= e^{-(\kappa_\mu - \delta'_\mu) s} \cdot
\frac{\bigl(e^{-\kappa_\mu s} - \rho_\mu^{\kappa_\mu / 2}\bigr)
^{\delta'_{0, \mu} / \kappa_\mu - 1}}
{\bigl(1 - \rho_\mu^{\kappa_\mu / 2}\bigr)^{\delta'_{0, \mu} / \kappa_\mu}} \\
&\lesssim e^{-(\kappa_\mu - \delta'_\mu) s},
\end{align*}
which implies that
$||\partial_s v_{a, b}^{\mu, \text{lef}}||_{L_{\delta'_\mu}^p([0, -\frac{1}{2} \log \rho_\mu]
\times S^1)} \lesssim 1$.
Hence the assumption
$||\partial_s \xi||_{L_{\delta'_\mu}^p([0, -\frac{1}{2} \log \rho_\mu] \times S^1)}
\lesssim 1$ and
\[
e^{\kappa_\mu s} |\partial_s (v_{a, b}^{\mu, \text{lef}} + \xi)| e^{\delta_\mu s}
\leq \rho_\mu^{(\delta'_\mu - \delta_\mu - \kappa_\mu)/2}
\cdot |\partial_s (v_{a, b}^{\mu, \text{lef}} + \xi)| e^{\delta'_\mu s}
\]
imply
\begin{align*}
||\partial_{\mathring{\rho}_\mu}^k F^{(\tilde a,\tilde b)}(\xi, h)|_{(\tilde a,\tilde b) = (a,b)}||_{
L_{\delta}^p(\widetilde{P}_a)}
&\lesssim \rho_\mu^{(\delta'_\mu - \delta_\mu)/2 - \kappa_\mu/2},\\
||D \partial_{\mathring{\rho}_\mu}^k F^{(\tilde a,\tilde b)}_{(\xi, h)}
(\hat \xi, \hat h)|_{(\tilde a,\tilde b) = (a,b)}||_{L_{\delta}^p(\widetilde{P}_a)}
&\lesssim \rho_\mu^{(\delta'_\mu - \delta_\mu)/2 - \kappa_\mu/2}
||\hat \xi||_{\widetilde{W}_{\delta'_\mu}^{1, p}(N^\mu_{a, b})},\\
D^2 \partial_{\mathring{\rho}_\mu}^k F^{(\tilde a,\tilde b)}_{(\xi, h)} &\equiv 0.
\end{align*}
The claim follows from these inequalities because $\partial_{\tilde \rho_\mu}
= \frac{\kappa_\mu}{2}(\tilde \rho_\mu)^{\kappa_\mu/2 - 1}
\partial_{\mathring{\rho}_\mu}$.

(ii)
We can easily check the following equations by direct calculation
($\natural_{m, l}$ is defined below):
\begin{align*}
\partial_{\tilde \rho_\nu}^k F^{(a,b)}(\xi, h)|_{N_{a, b}^{\nu, \text{left}}}
&= k ! \frac{1}{(1 - \rho_\nu)^k}(\partial_r v_{a, b}^\nu + \partial_r \xi)\\
&\quad + (-1)^k k ! \frac{1}{r^{k + 1}}
\Bigl(\frac{1 - r}{1 - \rho_\nu} \Bigr)^k
\widetilde{J}^\nu(v_{a, b}^\nu + \xi)(\partial_\theta v_{a, b}^\nu + \partial_\theta \xi) \\
&= k ! \frac{1}{(1 - \rho_\nu)^k}(\partial_r v_{a, b}^\nu + \partial_r \xi)\\
&\quad + (-1)^k k ! \frac{1}{r^{k + 1}}
\Bigl(\frac{1 - r}{1 - \rho_\nu} \Bigr)^k \natural_{(0, 0)} \hspace{50pt} (k \geq 1)
\end{align*}
{\abovedisplayskip=0pt
\begin{align*}
\partial_{\tilde \rho_\nu}^k \partial_{\tilde \varphi_\nu}
F^{(a, b)}(\xi, h)|_{N_{a, b}^{\nu, \text{left}}}
&= k ! \frac{1}{(1-\rho_\nu)^k} \Bigl(\beta_\nu \cdot \partial_r \partial_\theta
v_\zeta^\nu + \frac{1_{\{\rho_\nu \leq r_\nu \leq 2 \rho_\nu\}}}{\rho_\nu} 
\partial_\theta \xi \Bigr)\\
&\quad + (-1)^k k ! \frac{\beta_\nu}{r^{k + 1}}
\Bigl(\frac{1 - r}{1 - \rho_\nu}\Bigr)^k
\natural_{(0, 1)} \hspace{50pt} (k \geq 0)
\end{align*}}
{\abovedisplayskip=0pt
\begin{align*}
\partial_{\tilde \rho_\nu}^k \partial_{\tilde \varphi_\nu}^l
F^{(a, b)}(\xi, h)|_{N_{a, b}^{\nu, \text{left}}}
&= k ! \frac{\beta_\nu^l}{(1 - \rho_\nu)^k} \partial_r \partial_\theta^l v_{a, b}^\nu \\
&\quad + (-1)^k k ! \frac{\beta_\nu^l}{r^{k + 1}}
\Bigl(\frac{1 - r}{1 - \rho_\nu}\Bigr)^k \natural_{(0, l)} \quad (k \geq 0, l \geq 2)
\end{align*}}
{\abovedisplayskip=-5pt
\begin{align*}
(D \partial_{\tilde \rho_\nu}^k F^{(a, b)})_{(\xi, h)}(\hat \xi, \hat h)|
_{N_{a, b}^{\nu, \text{left}}}
&= k ! \frac{1}{(1 - \rho_\nu)^k} \partial_r \hat \xi \\
&\quad + (-1)^k k ! \frac{1}{r^{k + 1}}
\Bigl(\frac{1 - r}{1 - \rho_\nu}\Bigr)^k \natural_{(1, 0)} \quad (k \geq 1)
\end{align*}}
{\abovedisplayskip=0pt
\begin{align*}
(D \partial_{\tilde \rho_\nu}^k \partial_{\tilde \varphi_\nu}
F^{(a, b)})_{(\xi, h)}(\hat \xi, \hat h)|_{N_{a, b}^{\nu, \text{left}}}
&= k ! \frac{1}{(1 - \rho_\nu)^k} \frac{1_{\{\rho_\nu \leq r \leq 2\rho_\nu\}}}{\rho_\nu}
\partial_\theta \hat \xi\\
&\quad + (-1)^k k ! \frac{\beta_\nu}{r^{k + 1}}
\Bigl(\frac{1 - r}{1 - \rho_\nu}\Bigr)^k \natural_{(1, 1)} \quad (k \geq 0)
\end{align*}}
{\abovedisplayskip=-5pt
\begin{align*}
&(D^m \partial_{\tilde \rho_\nu}^k \partial_{\tilde \varphi_\nu}^l
F^{(a, b)})_{(\xi, h)}(\hat \xi^{(m)},
\hat h^{(m)}) (\hat \xi^{(m - 1)}, \hat h^{(m - 1)}) \dots
(\hat \xi^{(1)}, \hat h^{(1)})|_{N_{a, b}^{\nu, \text{left}}}\\
&\hspace{175pt}
= (-1)^k k ! \frac{\beta_\nu^l}{r^{k + 1}}
\Bigl(\frac{1 - r}{1 - \rho_\nu}\Bigr)^k \natural_{(m, l)}\\
&\hspace{130pt} (m \geq 2, (k, l) \neq (0, 0) \text{ or } m = 1, k \geq 0, l \geq 2)
\end{align*}}
In the above equations, 
$\natural_{(m,l)}$ is a sum of terms in the following forms:
{\belowdisplayskip=0pt
\begin{multline*}
(D^{\alpha_1 + m -1} \widetilde{J}^\nu)_{v_{a, b}^\nu + \xi} \cdot
\partial_\theta^{j_1}v_{a, b}^\nu \cdot
\partial_\theta^{j_2} v_{a, b}^\nu \dots \partial_\theta^{j_{\alpha_1}} v_{a, b}^\nu
\cdot \hat \xi^{(m)} \cdot \hat \xi^{(m-1)} \dots \hat \xi^{(1)}\\
(\alpha_1 \geq 1, j_1 + j_2 + \dots + j_{\alpha_1} = l + 1)
\end{multline*}}
{\abovedisplayskip=5pt
\begin{multline*}
(D^{\alpha_1 + m} \widetilde{J}^\nu)_{v_{a, b}^\nu + \xi} \cdot
\partial_\theta^{j_1} v_{a, b}^\nu \cdot
\partial_\theta^{j_2} v_{a, b}^\nu \dots \partial_\theta^{j_{\alpha_1}} v_{a, b}^\nu
\cdot \hat \xi^{(m)} \cdot \hat \xi^{(m - 1)} \dots \hat \xi^{(1)} \cdot
\partial_\theta \xi\\
(\alpha_1 \geq 0, j_1 + j_2 + \dots + j_{\alpha_1} + 1 = l +1)
\end{multline*}}
{\abovedisplayskip=-6pt
\begin{multline*}
(D^{\alpha_1 + m -1} \widetilde{J}^\nu)_{v_{a, b}^\nu + \xi} \cdot
\partial_\theta^{j_1} v_{a, b}^\nu
\cdot
\partial_\theta^{j_2} v_{a, b}^\nu \dots \partial_\theta^{j_{\alpha_1}} v_{a, b}^\nu
\cdot \hat \xi^{(m)} \cdot \hat \xi^{(m - 1)} \stackrel{\alpha_2}{\breve \dots}
\hat \xi^{(1)} \cdot \partial_\theta \hat \xi^{(\alpha_2)}\\
(\alpha_1 \geq 1, 1 \leq \alpha_2 \leq m,
j_1 + j_2 + \dots + j_{\alpha_1} + 1 = l + 1)
\end{multline*}
}

We prove the case where $p < q$. The case where $p = q$ is easier.
Define $2 < s < \infty$ by $\frac{1}{s} = \frac{1}{p} - \frac{1}{q}$.
First note that
\[
\biggl( \int_{\rho_\nu \leq |x| \leq 1} |r^{- (l + 1)} \natural_{(m, l)}|^q r dr d\theta
\biggr)^{\textstyle \frac{1}{q}}
\lesssim
||\hat \xi^{(m)}||_{1, q}|| \hat \xi^{(m-1)}||_{1, q} \dots ||\hat \xi^{(1)}||_{1, q}.
\]
This is because
\begin{itemize}
\item
$||r^{-j} \partial_\theta^j v_{a, b}^\nu||_\infty \lesssim
||r^{-j} \partial_\theta^j v_0^\nu||_\infty < \infty$ are uniformly bounded
with respect to small $\zeta_\nu$,
\item
$||\hat \xi^{(i)}||_{L^\infty(N_{a, b}^{\nu, \text{left}})} \lesssim ||\hat \xi^{(i)}||_{W^{1,q}
(N_{a, b}^{\nu, \text{left}})}$ uniformly with respect to small $\zeta_\nu$, and
\item
$||r^{-1} \partial_\theta \hat \xi^{(i)}||_{L^q(N_{a, b}^{\nu, \text{left}})} \lesssim
||\hat \xi^{(i)}||_{W^{1, q}(N_{a, b}^{\nu, \text{left}})}$ uniformly with respect to
small $\zeta_\nu$.
\end{itemize}
We note that
$||r^{-j} \partial_r \partial_\theta^j v_{a, b}^\nu||_\infty
\lesssim ||r^{-j} \partial_r \partial_\theta^j v_0^\nu||_\infty < \infty$.
These imply
\begin{flalign*}
&||\partial_{\tilde \rho_\nu}^k F^{(a, b)}(\xi, h)||_{L^p}&\\
&\lesssim
\biggl(\int_{\rho_\nu \leq |x| \leq 1} (|\partial_r v_{a, b}^\nu|^p + |\partial_r \xi|^p)
r dr d\theta\biggr)^{\textstyle \frac{1}{p}}
\! + \biggl(\int_{\rho_\nu \leq |x| \leq 1} r^{-kp} |r^{-1} \natural_{(0, 0)}|^p
r dr d\theta\biggr)^{\textstyle \frac{1}{p}}\\
&\lesssim
||\partial_r v_{a, b}^\nu||_\infty + ||\partial_r \xi||_{L^q(N_{a, b}^\nu)}\\
&\quad +
\biggl(\int_{\rho_\nu}^{1} r^{-ks + 1} dr\biggr)^{\textstyle \frac{1}{s}}
\biggl(\int_{\rho_\nu \leq |x| \leq 1}
|r^{-1} \natural_{(0, 0)}|^q r dr d\theta\biggr)^{\textstyle \frac{1}{q}}\\
&\lesssim
\rho_\nu^{{\textstyle \frac{2}{s}} - k}
\hspace{50pt} (k \geq 1),
\end{flalign*}
\begin{flalign*}
&||\partial_{\tilde \rho_\nu}^k \partial_{\tilde \varphi_\nu} F^{(a, b)}(\xi, h)||_{L^p}&\\
&\lesssim
\biggl(\int_{\rho_\nu \leq |x| \leq 2\rho_\nu} r^p |r^{-1} \partial_r \partial_\theta
v_{a, b}^\nu|^p
r dr d\theta\biggr)^{\textstyle \frac{1}{p}}\\
&\quad +
\frac{1}{\rho_\nu} \biggl(\int_{\rho_\nu \leq |x| \leq 2\rho_\nu} r^p|r^{-1}
\partial_\theta \xi|^p
r dr d\theta\biggr)^{\textstyle \frac{1}{p}}\\
&\quad +
\biggl( \int_{\rho_\nu \leq |x| \leq 2\rho_\nu} r^{(1 - k)p} |r^{-2} \natural_{(0, 1)}|^p
r dr d\theta \biggr)^{\textstyle \frac{1}{p}}\\
&\lesssim
\biggl(\int_{\rho_\nu}^{2\rho_\nu} r^{p + 1} dr \biggr)^{\textstyle \frac{1}{p}}
|r^{-1} \partial_r \partial_\theta v_{a, b}^\nu|_\infty\\
&\quad +
\frac{1}{\rho_\nu} \biggl(\int_{\rho_\nu}^{2\rho_\nu} r^{s + 1} dr
\biggr)^{\textstyle \frac{1}{s}}
\biggl(\int_{\rho_\nu \leq |x| \leq 2\rho_\nu} |r^{-1} \partial_\theta \xi|^q r dr d\theta
\biggr)^{\textstyle \frac{1}{q}}\\
&\quad +
\biggl(\int_{\rho_\nu}^{2\rho_\nu} r^{(1 - k)s + 1} dr \biggr)^{\textstyle \frac{1}{s}}
\biggl(\int_{\rho_\nu \leq |x| \leq 2\rho_\nu} |r^{-2} \natural_{(0, 1)}|^q r dr d\theta
\biggr)^{\textstyle \frac{1}{q}}\\
&\lesssim
\rho_\nu^{\textstyle \frac{2}{s}} + \rho_\nu^{{\textstyle \frac{2}{s}} + 1 - k}\\
&\lesssim
\rho_\nu^{{\textstyle \frac{2}{s}} - (k - 1)_+} \hspace{50pt} (k \geq 0),
\end{flalign*}
\begin{flalign*}
&||\partial_{\tilde \rho_\nu}^k \partial_{\tilde \varphi_\nu}^l F^{(a, b)}(\xi, h)||_{L^p}&\\
&\lesssim
\biggl(\int_{\rho_\nu \leq |x| \leq 2\rho_\nu} r^{lp} |r^{-l} \partial_r \partial_\theta^l
v_{a, b}^\nu|^p r dr d\theta \biggr)^{\textstyle \frac{1}{p}}\\
&\quad +
\biggl(\int_{\rho_\nu \leq |x| \leq 2\rho_\nu} r^{(l - k)p} |r^{-(l + 1)} \natural_{(0, l)}|^p
r dr d\theta \biggr)^{\textstyle \frac{1}{p}}\\
&\lesssim
\biggl(\int_{\rho_\nu}^{2\rho_\nu} r^{lp + 1} dr \biggr)^{\textstyle \frac{1}{p}}
|r^{-l} \partial_r \partial_\theta^l v_{a, b}^\nu|_\infty\\
&\quad +
\biggl(\int_{\rho_\nu}^{2\rho_\nu} r^{(l - k)s + 1} dr \biggr)^{\textstyle \frac{1}{s}}
\biggl(\int_{\rho_\nu \leq |x| \leq 2\rho_\nu} |r^{-(l + 1)} \natural_{(0, 1)}|^q r dr d\theta
\biggr)^{\textstyle \frac{1}{q}}\\
&\lesssim
\rho_\nu^{{\textstyle \frac{2}{s}} + l - k} \hspace{30pt} (k \geq 0, l \geq 2),
\end{flalign*}
\begin{flalign*}
&||(D \partial_{\tilde \rho_\nu}^k F^{(a, b)})_{(\xi, h)} (\hat \xi, \hat h)||_{L^p}&\\
&\lesssim
\biggl(\int_{\rho_\nu \leq |x| \leq 1} |\partial_r \hat \xi|^p r dr d\theta
\biggr)^{\textstyle \frac{1}{p}}
+ \biggl(\int_{\rho_\nu \leq |x| \leq 1} r^{-kp} |r^{-1} \natural_{(1, 0)}|^p
r dr d\theta \biggr)^{\textstyle \frac{1}{p}}\\
&\lesssim
||\hat \xi||_{1,q}
+ \biggl(\int_{\rho_\nu}^{1} r^{-ks + 1} dr \biggr)^{\textstyle \frac{1}{s}}
\biggl(\int_{\rho_\nu \leq |x| \leq 1} |r^{-1} \natural_{(1, 0)}|^q r dr d\theta
\biggr)^{\textstyle \frac{1}{q}}\\
&\lesssim
\rho_\nu^{{\textstyle \frac{2}{s}} - k} ||\hat \xi||_{1, q} \hspace{20pt} (k \geq 1),
\end{flalign*}
\begin{flalign*}
&||(D \partial_{\tilde \rho_\nu}^k \partial_{\tilde \varphi_\nu} F^{(a, b)})_{(\xi, h)}
(\hat \xi, \hat h)||_{L^p}\\
&\lesssim
\frac{1}{\rho_\nu} \biggl(\int_{\rho_\nu \leq |x| \leq 2\rho_\nu} r^p |r^{-1} \partial_\theta
\hat \xi|^p r dr d\theta\biggr)^{\textstyle \frac{1}{p}}\\
&\quad +
\biggl(\int_{\rho_\nu \leq |x| \leq 2\rho_\nu} r^{-(k - 1)p} |r^{-2} \natural_{(1, 1)}|^p
r dr d\theta
\biggr)^{\textstyle \frac{1}{p}}\\
&\lesssim
\rho_\nu^{-1} \biggl(\int_{\rho_\nu}^{2\rho_\nu} r^{s + 1} dr\biggr)^{\textstyle \frac{1}{s}}
\biggl(\int_{\rho_\nu \leq |x| \leq 2\rho_\nu} |r^{-1} \partial_\theta \hat \xi|^q
r dr d\theta \biggr)^{\textstyle \frac{1}{q}}\\
&\quad +
\biggl(\int_{\rho_\nu}^{2\rho_\nu} r^{-(k - 1)s + 1} dr \biggr)^{\textstyle \frac{1}{s}}
\biggl(\int_{\rho_\nu \leq |x| \leq 2\rho_\nu} |r^{-2} \natural_{(1, 1)}|^q r dr d\theta
\biggr)^{\textstyle \frac{1}{q}}\\
&\lesssim
\rho_\nu^{{\textstyle \frac{2}{s}}-(k - 1)_+} ||\hat \xi||_{1, q} \hspace{20pt} (k \geq 0),
\hspace{200pt}
\end{flalign*}
\begin{flalign*}
&||(D^m \partial_{\tilde \rho_\nu}^k \partial_{\tilde \varphi_\nu}^l F^{(a, b)})_{(\xi, h)}
(\hat \xi^{(m)}, \hat h^{(m)}) (\hat \xi^{(m - 1)}, \hat h^{(m - 1)}) \dots
(\hat \xi^{(1)}, \hat h^{(1)})||_{L^p}&\\
&\lesssim
\biggl( \int_{\rho_\nu \leq |x| \leq 2\rho_\nu} r^{(l - k)p} |r^{-(l + 1)} \natural_{(m, l)}|^p
r dr d\theta \biggr)^{\textstyle \frac{1}{p}}\\
&\lesssim
\biggl(\int_{\rho_\nu}^{2\rho_\nu} r^{(l - k)s + 1} dr \biggr)^{\textstyle \frac{1}{s}}
\biggl(\int_{\rho_\nu \leq |x| \leq 2\rho_\nu} |r^{-(l + 1)} \natural_{(m, l)}|^q r dr d\theta
\biggr)^{\textstyle \frac{1}{q}}\\
&\lesssim
\rho_\nu^{{\textstyle \frac{2}{s}} + (l - k)}
||\hat \xi^{(m)}||_{1, q}||\hat \xi^{(m - 1)}||_{1, q} \dots ||\hat \xi^{(1)}||_{1, q}
\quad (m \geq 2, l > 0 \\
&\hph{\lesssim
\rho_\nu^{{\textstyle \frac{2}{s}} + (l - k)}
||\hat \xi^{(m)}||_{1, q}||\hat \xi^{(m - 1)}||_{1, q} \dots ||\hat \xi^{(1)}||_{1, q}
\quad (}
\text{ or } m = 1, k \geq 0, l \geq 2).
\end{flalign*}
\begin{flalign*}
&||(D^m \partial_{\tilde \rho_\nu}^k F^{(a, b)})_{(\xi, h)}
(\hat \xi^{(m)}, \hat h^{(m)}) (\hat \xi^{(m - 1)}, \hat h^{(m - 1)}) \dots
(\hat \xi^{(1)}, \hat h^{(1)})||_{L^p}&\\
&\lesssim
\biggl( \int_{\rho_\nu \leq |x| \leq 1} r^{(l - k)p} |r^{-(l + 1)} \natural_{(m, l)}|^p
r dr d\theta \biggr)^{\textstyle \frac{1}{p}}\\
&\lesssim
\biggl(\int_{\rho_\nu}^1 r^{-ks + 1} dr \biggr)^{\textstyle \frac{1}{s}}
\biggl(\int_{\rho_\nu \leq |x| \leq 1} |r^{-(l + 1)} \natural_{(m, l)}|^q r dr d\theta
\biggr)^{\textstyle \frac{1}{q}}\\
&\lesssim
\rho_\nu^{{\textstyle \frac{2}{s}} - k}
||\hat \xi^{(m)}||_{1, q}||\hat \xi^{(m - 1)}||_{1, q} \dots ||\hat \xi^{(1)}||_{1, q}
\hspace{50pt} (m \geq 2, k > 0)
\end{flalign*}
These inequalities prove the claim.

(iii)
It is straightforward to prove this case using the equalities
\begin{align*}
\partial_{\tilde \varphi_\mu} F^{(\tilde a, \tilde b)}(\xi, h)|_{(\tilde a, \tilde b) = (a, b)}
&= \frac{1}{2} \chi'(s) (g^\mu_t(v_{a, b}^{\mu, \mathrm{left}} + \xi)
+ \partial_t(v_{a, b}^{\mu, \mathrm{left}} + \xi))\\
\partial_{\tilde \varphi_\mu}^2 F^{(\tilde a,\tilde b)}(\xi, h) &= 0
\end{align*}
and
\begin{align*}
\partial_{\tilde b_\mu} F^{(\tilde a,\tilde b)}(\xi, h)|_{(\tilde a,\tilde b) = (a,b)}
&= \frac{1}{2} \chi'(s) \partial_\sigma\\
\partial_{\tilde b_\mu}^2 F^{(\tilde a,\tilde b)}(\xi, h) &= 0
\end{align*}
on $[-1, 0] \times S^1_\mu$.
\end{proof}

For each $i \in \Pi$, we fix a index $\mu_i \in M_i$.
Then a coordinate of $\mathring{X}_{\Pi, \Pi'}$ is given by
$(j, (b_\mu)_{\mu}, (\rho_{\mu_i})_{i \in \Pi}, (\varphi_\mu)_{\mu},
(\rho_\nu^{2\pi} e^{2\pi\sqrt{-1} \varphi_\nu})_{\nu \in \Pi'})$.
Note that in this coordinate,
$\rho_\mu = \rho_{\mu_i}^{L_{\mu_i} / L_\mu} e^{(b_\mu - b_{\mu_i}) / L_\mu}$
for any $\mu \in M_i$ ($i \in \Pi$).
We rewrite the above lemma in this coordinate and get the following corollary.
(The meaning of $\partial_{\mu_i}$ and $\partial_b$ in the following corollary
are different from that in Lemma \ref{estimates of implicit function 0}.)
\begin{cor}\label{estimates of implicit function}\ 
\begin{enumerate}[label=\normalfont(\roman*)]
\item
For any $0 < \delta_\mu \leq \delta'_\mu < \kappa_\mu$, $2 < p < \infty$,
$i \in \Pi$, $k \neq 0$ and multi-index $k_b = (k_{b_\mu})$,
there exists some constant $C>0$ and
$c_0 > 0$ such that if $||\xi||_{\widetilde{W}^{1,p}_{\delta'}} \leq c_0$, then
\begin{align*}
||(\partial_{\rho_{\mu_i}}^k \partial_b^{k_b} F^{(a, b)})(\xi, h)||
_{L_{\delta}^p(\widetilde{P}_a)}
&\leq C \rho_{\mu_i}^{L_{\mu_i} \tilde \delta_i /2 - k},\\
||(D \partial_{\rho_{\mu_i}}^k \partial_b^{k_b} F^{(a, b)})_{(\xi, h)}
(\hat \xi, \hat h)||_{L_\delta^p(\widetilde{P}_a)}
&\leq C \rho_{\mu_i}^{L_{\mu_i} \tilde \delta_i /2 - k}
||\hat \xi||_{\widetilde{W}_{\delta'}^{1, p}(\bigcup_{\mu \in M_i} N_{a, b}^\mu)},\\
D^m \partial_{\rho_{\mu_i}}^k \partial_b^{k_b} F^{(a, b)}
&\equiv 0 \quad (m \geq 2),
\end{align*}
where $\tilde \delta_i = \min\{(\delta'_\mu - \delta_\mu) / L_\mu;
\mu \in M_i\}$.
\item
For any $2 < p \leq q < \infty$,
$\nu \in \Pi'$ and $(m, k, l)$ such that $(k, l) \neq (0, 0)$,
there exists some constants $C > 0$ and $c_0 > 0$ such that
if $||\xi||_{W^{1, q}(N_{a, b}^\nu)} \leq c_0$, then 
\begin{align*}
||(D^m \partial_{\rho_\nu}^k \partial_{\varphi_\nu}^l F^{(a, b)})_{(\xi, h)}
(\hat \xi^{(m)}, \hat h^{(m)}) (\hat \xi^{(m-1)}, \hat h^{(m-1)}) \dots
(\hat \xi^{(1)}, \hat h^{(1)})||_{L^p(N_{a, b}^\nu)}\\
\leq C \rho_\nu^{2/p - 2/q - k} \prod_{i = 1}^m ||\hat \xi^{(i)}||_{W^{1,q}(N_{a, b}^\nu)}.
\end{align*}
\item
For any $2 < p < \infty$, $0< \delta < \delta_0$,
and multi-index $(m, k_j, (l_\mu), k_b = (k_{b_\mu}))$,
there exists some constants $C > 0$ and $c_0 > 0$ such that
if $||\xi||_{\widetilde{W}^{1,p}_\delta} + |h|_E \leq c_0$, then 
\begin{align*}
||(D^m \partial_j^{k_j} \partial_{(\varphi_\mu)}^{(l_\mu)} \partial_b^{k_b}
F^{(a, b)})_{(\xi,h)}
(\hat \xi^{(m)}\!, \hat h^{(m)}) (\hat \xi^{(m-1)}\!, \hat h^{(m-1)})\! \dots \!
(\hat \xi^{(1)}\!, \hat h^{(1)})||_{L^p_\delta(\widetilde{P}_a)}\\
\leq C \prod_{i = 1}^m (||\hat \xi^{(i)}||_{\widetilde{W}^{1,p}_\delta(\widetilde{P}_a)}
+ |\hat h^{(i)}|_E).
\end{align*}
\end{enumerate}
\end{cor}

Let $U \subset \mathring{X}_{\Pi, \Pi'}$ be a neighborhood of $(a, b)$, and
regard the family of smooth maps
\[
\phi^{\tilde a,\tilde b} : \Ker D_0 \supset B_\epsilon(0)
\to \widetilde{W}^{1,p}_\delta(\widetilde{P}_a, u_{a, b}^\ast T \hat Y) \times E^0
\]
as a map
\begin{equation}
\phi : U \times B_\epsilon(0) \to
\widetilde{W}^{1,p}_\delta(\widetilde{P}_a, u_{a, b}^\ast T \hat Y) \times E^0.
\label{phi 1}
\end{equation}
We estimate the derivative of $\phi$ at $(a, b, x) \in U \times B_\epsilon(0)$.
As we have already mentioned, the domain of $\phi^{a, b}$ or $\epsilon > 0$
may depend on $2 < p < \infty$ and $0 < \delta < \delta_0$.
Hence in the following Proposition, we need to assume that $(a, b, x)$ is
sufficiently close to $(0, b^0, 0)$ for given $p$, $q$, $\delta$ and $\delta'$
to guarantee that $(a, b, x)$ is contained in the domains of various $\phi$.
\begin{prop}\label{asymptotic phi}
For any $2 < p < q$, $0 < \delta_\mu < \delta'_\mu < \kappa_\mu$ and any multi-index
$(k_x, k_j, k_b, (k_{\mu_i})_{i \in \Pi}, (l_\mu)_\mu, (k_\nu)_{\nu \in \Pi'},
(l_\nu)_{\nu \in \Pi'})$, there exists some constant $C > 0$ such that
\begin{multline*}
||\partial_x^{k_x} \partial_j^{k_j} \partial_b^{k_b} \partial_{(\rho_{\mu_i})}^{(k_{\mu_i})}
\partial_{(\varphi_\mu)}^{(l_\mu)} \partial_{(\rho_\nu)}^{(k_\nu)}
\partial_{(\varphi_\nu)}^{(l_\nu)} \phi(a, b, x)||_{\widetilde{W}_\delta^{1,p}
(\widetilde{P}_a, u_{a, b}^\ast T \hat Y) \times E^0}\\
\leq C \prod_{\substack{i \\ k_{\mu_i} \neq 0}}
\rho_{\mu_i}^{L_{\mu_i} \tilde \delta_i /2 - k_{\mu_i}}
\prod_{\substack{\nu \\ (k_\nu, l_\nu) \neq (0,0)}} \rho_\nu^{(2/p - 2/q) - k_\nu}
\end{multline*}
for any $(a, b, x) \in \mathring{X}_{\Pi, \Pi'} \times B_\epsilon(0)$ sufficiently close
to $(0, b^0, 0)$.
Furthermore, if $(k_{\nu_0}, l_{\nu_0}) = (0,0)$ then
\begin{multline*}
||\partial_x^{k_x} \partial_j^{k_j} \partial_b^{k_b} \partial_{(\rho_{\mu_i})}^{(k_{\mu_i})}
\partial_{(\varphi_\mu)}^{(l_\mu)} \partial_{(\rho_\nu)}^{(k_\nu)}
\partial_{(\varphi_\nu)}^{(l_\nu)} \phi(a, b, x)||_{\widetilde{W}^{1,q}
(N_{a, b}^{\nu_0})}\\
\leq C \prod_{\substack{i \\ k_{\mu_i} \neq 0}}
\rho_\mu^{L_{\mu_i} \tilde \delta_i /2 - k_{\mu_i}}
\prod_{\substack{\nu \\ (k_\nu, l_\nu) \neq (0,0)}} \rho_\nu^{(2/p - 2/q) - k_\nu}
\end{multline*}
for any $(a, b, x) \in \mathring{X}_{\Pi, \Pi'} \times B_\epsilon(0)$ sufficiently close
to $(0, b^0, 0)$.
\end{prop}
\begin{proof}
We prove the claim by the induction in $|k_x| + |k_j| + |k_b|
+ |(k_{\mu_i})| + |(l_\mu)| + |(k_\nu)| + |(l_\nu)|$.
The case $(k_x, k_j, k_x, (k_{\mu_i}), (l_\mu), (k_\nu), (l_\nu)) = (0, \dots, 0)$ is obvious.
Differentiating the equation $F^{(\tilde a,\tilde b) +} (\phi(\tilde a, \tilde b, x)) = (0, x)$
of smooth functions on a fixed curve $\widetilde{P}_a$
by $\partial_x^{k_x} \partial_{\tilde \jmath}^{k_j} \partial_{\tilde b}^{k_b}
\partial_{(\tilde \rho_{\mu_i})}^{(k_{\mu_i})} \partial_{(\tilde \varphi_\mu)}^{(l_\mu)}
\partial_{(\tilde \rho_\nu)}^{(k_\nu)} \partial_{(\tilde \varphi_\nu)}^{(l_\nu)}$,
we obtain an equation of the following form.
\begin{align}
&(DF^{(a, b) +})_{\phi(a, b, x)}
\partial_x^{k_x} \partial_j^{k_j} \partial_b^{k_b} \partial_{(\rho_{\mu_i})}^{(k_{\mu_i})}
\partial_{(\varphi_\mu)}^{(l_\mu)} \partial_{(\rho_\nu)}^{(k_\nu)}
\partial_{(\varphi_\nu)}^{(l_\nu)} \phi \notag\\
&+
\sum_{\star_1}
(D^{m} \partial_{\rho_{\nu_0}}^{k_{\nu_0}}\partial_{\varphi_{\nu_0}}^{l_{\nu_0}}
F^{(a, b) +})_{\phi(a, b, x)}
(\hat\xi^{(m)},\hat h^{(m)}) \dots (\hat\xi^{(1)},\hat h^{(1)}) \notag\\
&+
\sum_{\star_2}
(D^{m} \partial_{\rho_{\mu_{i_0}}}^{k_{\mu_{i_0}}} \partial_b^{k'_b}
F^{(a, b) +})_{\phi(a, b, x)}
(\hat\xi^{(m)},\hat h^{(m)}) \dots (\hat\xi^{(1)},\hat h^{(1)})
\notag\\
&+
\sum_{\star_3}
(D^{m} \partial_x^{k'_x} \partial_j^{k'_j} \partial_b^{k'_b}
\partial_{(\rho_\mu)}^{(k'_\mu)} \partial_{(\varphi_\mu)}^{(l'_\mu)}
\partial_{(\rho_\nu)}^{(k'_\nu)} \partial_{(\varphi_\nu)}^{(l'_\nu)}
F^{(a, b) +})_{\phi(a, b, x)}
(\hat\xi^{(m)},\hat h^{(m)}) \dots (\hat\xi^{(1)},\hat h^{(1)}) \notag\\
&=0,\label{good bad}
\end{align}
where each $(\hat\xi^{(l)},\hat h^{(l)})$ is some derivative of $\phi$, and
the sum of the indices of differentials which appear in each term is equal to
$(k_x, k_j, k_b, (k_{\mu_i}), (l_\mu), (k_\nu), (l_\nu))$;
in the sum $\star_1$,
$(k_{\nu_0}, l_{\nu_0}) \neq (0,0)$ and
each $(\hat\xi^{(l)},\hat h^{(l)})$ is some differential of $\phi$
by $\partial_x$, $\partial_j$, $\partial_b$, $\partial_{\rho_{\mu_i}}$,
$\partial_{\varphi_\mu}$, $\partial_{\rho_{\nu}}$ and $\partial_{\varphi_{\nu}}$
except $\partial_{\rho_{\nu_0}}$ and $\partial_{\varphi_{\nu_0}}$;
in the sum $\star_2$, $k_{\mu_{i_0}} \neq 0$, $m = 0, 1$ and
each $(\hat\xi^{(l)},\hat h^{(l)})$ is some differential of $\phi$
by $\partial_x$, $\partial_j$, $\partial_b$, $\partial_{\rho_{\mu_i}}$,
$\partial_{\varphi_{\mu}}$, $\partial_{\rho_\nu}$ and $\partial_{\varphi_\nu}$
except $\partial_{\rho_{\mu_{i_0}}}$;
in the sum $\star_3$, if $(k_\nu, l_\nu) \neq (0,0)$ then
$(k'_\nu, l'_\nu) \neq (k_\nu, l_\nu)$, and if $k_{\mu_i} \neq 0$ then
$k'_{\mu_i} < k_{\mu_i}$.
(As we have noted in Remark \ref{independent}, for example, if $k'_\nu \neq 0$ and
$k'_{\mu_i} \neq 0$, then this term vanishes.)

Corollary \ref{estimates of implicit function} (ii) and the assumption of the induction
(the second inequality) imply that the $L^p$-norm of each term
in the sum $\star_1$ is bounded by
\begin{align*}
&||(D^{m}\partial_{\rho_{\nu_0}}^{k_{\nu_0}}\partial_{\varphi_{\nu_0}}^{l_{\nu_0}}
F^{(a, b) +})_{\phi(a, b, x)}
(\hat\xi^{(m)},\hat h^{(m)}) \dots (\hat\xi^{(1)},\hat
h^{(1)})||_{L^p(N_{a, b}^\nu)}\\
&\lesssim
\rho_{\nu_0}^{2/p - 2/q-k_{\nu_0}} \prod_i ||\hat\xi^{(i)}||_{W^{1,q}(N_{a, b}^{\nu_0})}\\
&\lesssim
\prod_{\substack{i \\ k_{\mu_i} \neq 0}}
\rho_\mu^{L_{\mu_i} \tilde \delta_i /2 - k_{\mu_i}}
\prod_{\substack{\nu \\ (k_{\nu},l_{\nu})\neq (0,0)}} \rho_{\nu}^{2/p - 2/q-k_{\nu}}
\end{align*}
since $\partial_{\rho_{\nu}}$ or $\partial_{\varphi_{\nu}}$ appears in some
$(\hat\xi^{(l)},\hat h^{(l)})$
for each $\nu \neq \nu_0$ such that $(k_{\nu},l_{\nu}) \neq (0,0)$,
$\partial_{\rho_{\mu_i}}$ appears in some $(\hat\xi^{(l)},\hat h^{(l)})$ for each $i$
such that $k_{\mu_i} \neq 0$,
and neither $\partial_{\rho_{\nu_0}}$ nor $\partial_{\varphi_{\nu_0}}$ appears.

Next we consider the sum $\star_2$.
For each $i_0$, define a sequence of positive constants
$\delta'' = ((\delta''_\mu)_\mu, (\delta''_{\pm\infty_i})_{\pm\infty_i})$
by $\delta''_\mu = \delta'_\mu$ for $\mu \in M_{i_0}$, $\delta''_\mu = \delta_\mu$
for $\mu \notin M_{i_0}$, and $\delta''_{\pm\infty_i} = \delta_{\pm\infty_i}$.
Then Corollary \ref{estimates of implicit function} (i) and the assumption of the induction
(the first inequality) imply that the $L_{\delta}^p$-norm of the terms
with $m = 1$ in the sum $\star_2$ is bounded by
\begin{align*}
&||(D \partial_{\rho_{\mu_{i_0}}}^{k_{\mu_{i_0}}} F^{(a, b) +})_{\phi(a, b, x)}
(\hat\xi, \hat h)||_{L_\delta^p(\widetilde{P}_a)}\\
&\lesssim
\rho_{\mu_{i_0}}^{L_{\mu_{i_0}} \tilde \delta_{i_0} /2 - k_{\mu_{i_0}}}
||\hat \xi||_{\widetilde{W}_{\delta''}^{1, p}(\widetilde{P}_a)}\\
&\lesssim
\prod_{\substack{i \\ k_{\mu_i} \neq 0}}
\rho_{\mu_i}^{L_{\mu_i} \tilde \delta_i /2 - k_{\mu}}
\prod_{\substack{\nu \\ (k_\nu,l_\nu)\neq (0,0)}} \rho_\nu^{2/p - 2/q-k_\nu}
\end{align*}
since $\partial_{\rho_\nu}$ or $\partial_{\varphi_\nu}$ appears in some
$(\hat\xi^{(l)},\hat h^{(l)})$
for each $\nu$ such that $(k_\nu,l_\nu)\neq (0,0)$,
$\partial_{\rho_{\mu_i}}$ appears in some $(\hat\xi^{(l)},\hat h^{(l)})$ for each
$i \neq i_0$ such that $k_{\mu_i} \neq 0$,
and $\partial_{\rho_{\mu_{i_0}}}$ does not appear.
If the terms with $m = 0$ appear in the sum $\star_2$, then
$\partial_x^{k_x} \partial_j^{k_j} \partial_b^{k_b} \partial_{(\rho_{\mu_i})}^{(k_{\mu_i})}
\partial_{(\varphi_\mu)}^{(l_\mu)} \partial_{(\rho_\nu)}^{(k_\nu)}
\partial_{(\varphi_\nu)}^{(l_\nu)} = \partial_{\rho_{\mu_{i_0}}}^{k_{\mu_{i_0}}}$, and
Corollary \ref{estimates of implicit function} (i) implies
\[
||\partial_{\rho_{\mu_{i_0}}}^{k_{\mu_{i_0}}} F^{(a, b) +} (\phi(a, b, x))||_{
L_\delta^p(\widetilde{P}_a)} \lesssim
\rho_{\mu_{i_0}}^{L_{\mu_{i_0}} \tilde \delta_{i_0} /2 - k_{\mu_{i_0}}}.
\]

Similarly, Corollary \ref{estimates of implicit function} (i) for $\delta' = \delta$,
(ii) for $q = p$, (iii) and the assumption of the induction (the first inequality) imply
\begin{align*}
&||(D^{m} \partial_x^{k'_x} \partial_j^{k'_j} \partial_b^{k'_b}
\partial_{(\rho_{\mu_i})}^{(k'_{\mu_i})} \partial_{(\varphi_\mu)}^{(l'_\mu)}
\partial_{(\rho_\nu)}^{(k'_\nu)} \partial_{(\varphi_\nu)}^{(l'_\nu)}
F^{(a, b) +})_{\phi(a, b, x)}
(\hat\xi^{(m)}\!,\hat h^{(m)}) \!\dots\! (\hat\xi^{(1)}\!,\hat h^{(1)})||_{L^p_\delta
(\widetilde{P}_a)}\\
&\lesssim
\prod_{\substack{i \\ k'_{\mu_i} \neq 0}} \rho_{\mu_i}^{- k'_{\mu_i}}
\prod_{\substack{\nu \\ (k'_\nu, l'_\nu) \neq (0,0)}} \rho_\nu^{-k'_\nu}
\prod_l (||\hat\xi^{(l)}||_{\widetilde{W}^{1,p}_\delta(\widetilde{P}_a)} + |\hat h^{(l)}|_E)\\
&\lesssim
\prod_{\substack{i \\ k_{\mu_i} \neq 0}}
\rho_\mu^{L_{\mu_i} \tilde \delta_i /2 - k_{\mu_i}}
\prod_{\substack{\nu \\ (k_\nu,l_\nu)\neq (0,0)}} \rho_\nu^{2/p - 2/q-k_\nu}
\end{align*}
since $\partial_{\rho_\nu}$ or $\partial_{\varphi_\nu}$ appears in some
$(\hat\xi^{(l)},\hat h^{(l)})$ for each $\nu$ such that $(k_\nu,l_\nu)\neq (0,0)$,
and $\partial_{\rho_{\mu_i}}$ appears in some $(\hat\xi^{(l)},\hat h^{(l)})$ for each $i$
such that $k_{\mu_i} \neq 0$.

Since $(DF^{(a, b) +}_{\phi(a, b, x)})^{-1}$ is uniformly bounded,
these estimates imply the first inequality of the claim.

Next we prove the second inequality.
If $(k_{\nu_0},l_{\nu_0})=(0,0)$, then the restriction of equation (\ref{good bad})
to $L^p(N_{a, b}^{\nu_0})$ is
\begin{align}
&(DF^{(a, b)})_{\phi(a, b, x)}
\partial_x^{k_x} \partial_j^{k_j} \partial_b^{k_b} \partial_{(\rho_{\mu_i})}^{(k_{\mu_i})}
\partial_{(\varphi_\mu)}^{(l_\mu)} \partial_{(\rho_\nu)}^{(k_\nu)}
\partial_{(\varphi_\nu)}^{(l_\nu)} \phi \notag\\
&+
\sum (D^{m} F^{(a, b)})_{\phi(a, b, x)}(\hat\xi^{(m)},\hat h^{(m)})
\dots (\hat\xi^{(1)},\hat h^{(1)}) =0. \label{near node}
\end{align}
Sobolev embedding and the first inequality of the claim imply
\begin{align*}
&||\partial_x^{k_x} \partial_j^{k_j} \partial_b^{k_b} \partial_{(\rho_{\mu_i})}^{(k_{\mu_i})}
\partial_{(\varphi_\mu)}^{(l_\mu)} \partial_{(\rho_\nu)}^{(k_\nu)}
\partial_{(\varphi_\nu)}^{(l_\nu)} \phi||_{L^q(N_{a, b}^{\nu_0})}\\
&\lesssim
||\partial_x^{k_x} \partial_j^{k_j} \partial_b^{k_b} \partial_{(\rho_{\mu_i})}^{(k_{\mu_i})}
\partial_{(\varphi_\mu)}^{(l_\mu)} \partial_{(\rho_\nu)}^{(k_\nu)}
\partial_{(\varphi_\nu)}^{(l_\nu)} \phi||_{W^{1,p}(N_{a, b}^{\nu_0})}\\
&\lesssim
\prod_{\substack{i \\ k_{\mu_i} \neq 0}}
\rho_\mu^{L_{\mu_i} \tilde \delta_i /2 - k_{\mu_i}}
\prod_{\substack{\nu \\ (k_\nu,l_\nu)\neq (0,0)}} \rho_\nu^{2/p - 2/q-k_\nu}.
\end{align*}
Corollary \ref{estimates of implicit function} (iii) and the assumption of the induction
(the second inequality) imply
\begin{align*}
&||(D^{m} F^{(a, b)})_{\phi(a, b, x)}(\hat\xi^{(m)},\hat h^{(m)})
\dots (\hat\xi^{(1)},\hat h^{(1)})||_{L^q(N_{a, b}^{\nu_0})}\\
&\lesssim
\prod_{\substack{i \\ k_{\mu_i} \neq 0}}
\rho_\mu^{L_{\mu_i} \tilde \delta_i /2 - k_{\mu_i}}
\prod_{\substack{\nu \\ (k_\nu,l_\nu)\neq (0,0)}} \rho_\nu^{2/p - 2/q-k_\nu}.
\end{align*}
Hence $W^{1,q}$ regularity property of the elliptic operator
$(DF^{(a, b)})_{\phi(a, b, x)}$ in (\ref{near node})
implies the second inequality.
(Note that the regularity property of $(DF^{(a, b)})_{\phi(a, b, x)}$
is uniform with respect to small $\zeta_\nu$.
See Remark \ref{uniform regularity}.)
\end{proof}

Next we regard the family of smooth maps
\[
\phi^{a, b} : \Ker D_0 \supset B_\epsilon(0)
\to \widetilde{W}^{1,p}_\delta(\widetilde{P}_a, u_{a, b}^\ast T \hat Y) \times E^0
\]
as a map
\begin{align}
\phi : \mathring{X} \times B_\epsilon(0) &\to
W^{1,p}(\Sigma_0 \setminus N_0, (\R_1 \sqcup \R_2 \sqcup \dots
\sqcup \R_k) \times Y) \times E^0 \notag\\
((a,b), x) &\mapsto (\Phi(z, \xi_{a, b, x}(z))|_{\Sigma_0 \setminus N_0}, h_{a, b, x}),
\label{phi 2}
\end{align}
where $\xi_{a, b, x}$ and $h_{a, b, x}$ is defined by
$\phi^{a,b}(x) = (\xi_{a, b, x}, \xi_{a, b, x})$.
For each $i = 1, \dots, k-1$, fix a index $\mu_i \in M_i$.
Then a coordinate of $\mathring{X} \subset \widetilde{X} \times \prod_\mu \R_\mu$
is given by $(j, (b_\mu), (\rho_{\mu_i})_i, (\varphi_\mu)_\mu,
(\varphi_\nu^{2\pi} e^{2\pi\sqrt{-1} \varphi_\nu})_\nu)$.
For a neighborhood $U \subset \mathring{X}_{\Pi, \Pi'}$ of each point
$(a, b) \in \mathring{X}_{\Pi, \Pi'}$, the restriction of (\ref{phi 2})
to $U \times B_\epsilon(0)$ is the composition of the map (\ref{phi 1})
and the projection
$\widetilde{W}^{1,p}_\delta(\widetilde{P}_a, u_{a, b}^\ast T \hat Y) \times E^0
\to W^{1,p}(\Sigma_0 \setminus N_0, (\R_1 \sqcup \R_2 \sqcup \dots
\sqcup \R_k) \times Y) \times E^0$.
Furthermore, the norm of this projection is uniform with respect to $(a, b)$.
Therefore, Proposition \ref{asymptotic phi} implies
\begin{multline*}
||\partial_x^{k_x} \partial_j^{k_j} \partial_b^{k_b} \partial_{(\rho_{\mu_i})}^{(k_{\mu_i})}
\partial_{(\varphi_\mu)}^{(l_\mu)} \partial_{(\rho_\nu)}^{(k_\nu)}
\partial_{(\varphi_\nu)}^{(l_\nu)} \phi(a, b, x)||_{W^{1,p}(\Sigma_0 \setminus N_0,
(\R_1 \sqcup \R_2 \sqcup \dots \sqcup \R_k) \times Y) \times E^0}\\
\leq C \prod_{\substack{i \\ k_{\mu_i} \neq 0}}
\rho_{\mu_i}^{L_{\mu_i} \tilde \delta_i/2 - k_{\mu_i}}
\prod_{\substack{\nu \\ (k_\nu, l_\nu) \neq (0,0)}} \rho_\nu^{(2/p - 2/q) - k_\nu},
\end{multline*}
where $\tilde \delta_i = \min\{(\delta'_\mu - \delta_\mu) /L_\mu;
\mu \in M_i\}$.
The same estimate holds true for any Sobolev norm $W^{k, p}$
or $C^l$-norm instead of $W^{1, p}$ if we change the constant $C > 0$
because of elliptic regularity.
Since these estimates hold true for
arbitrary $2 < p < q < \infty$ and $0 < \delta \leq \delta' < \delta_0$
such that $0 < \delta_\mu \leq \delta'_\mu < \kappa_\mu$ if we shrink
the domain of $\phi$, the following corollary holds true.
Define $\tilde \delta_{0, i} = \min\{ \kappa_\mu / L_\mu; \mu \in M_i\}$ for each $i$.
We regard $\phi$ as a map
\[
\phi : \mathring{X} \times B_\epsilon(0) \to
C^l(\Sigma_0 \setminus N_0, (\R_1 \sqcup \R_2 \sqcup \dots
\sqcup \R_k) \times Y) \times E^0.
\]
\begin{cor}\label{asymptotic estimates of phi}
For any $l \geq 1$, $0 < \epsilon < 1$,
$0 < \tilde \delta'_{0, i} < \tilde \delta_{0, i}$,
$(\Pi, \Pi')$ and any multi-index $(k_x, \ab k_j, \ab k_b, \ab (k_{\mu_i})_{i \in \Pi}, \ab
(l_\mu)_\mu, (k_\nu)_{\nu \in \Pi'}, (l_\nu)_{\nu \in \Pi'})$,
there exists some constant $C > 0$ such that
\begin{multline*}
||\partial_x^{k_x} \partial_j^{k_j} \partial_b^{k_b} \partial_{(\rho_{\mu_i})}^{(k_{\mu_i})}
\partial_{(\varphi_\mu)}^{(l_\mu)} \partial_{(\rho_\nu)}^{(k_\nu)}
\partial_{(\varphi_\nu)}^{(l_\nu)} \phi(a, b, x)||_{C^l(\Sigma_0 \setminus N_0,
(\R_1 \sqcup \R_2 \sqcup \dots \sqcup \R_k) \times Y) \times E^0}\\
\leq C \prod_{\substack{i \\ k_{\mu_i} \neq 0}}
\rho_\mu^{L_{\mu_i} \tilde \delta'_{0, i}/2 - k_{\mu_i}}
\prod_{\substack{\nu \\ (k_\nu, l_\nu) \neq (0,0)}} \rho_\nu^{\epsilon - k_\nu}
\end{multline*}
for all $(a, b, x) \in \mathring{X}_{\Pi, \Pi'} \times B_\epsilon(0)$
sufficiently close to $(0, b^0, 0)$.
\end{cor}

Recall that we give a strong differential structure to
$\widetilde{X}$ determined by fixed constants $\alpha$ and $\beta$,
and give $\hat V$ the product smooth structure.
(See the beginning of this section.)
\begin{cor}
For any $N$, $\phi$ is of class $C^N$ if $\alpha$ and $\beta$ are sufficiently large.
\end{cor}
\begin{proof}
If we change the coordinate $\rho_\mu$ and $\rho_\nu$ to $\hat \rho_\mu$ and
$\hat \rho_\nu$ respectively by
$\rho_\mu^{L_\mu} = (\hat \rho_\mu)^\beta$ and
$\rho_\nu = (\hat \rho_\nu)^\alpha$,
then the previous corollary implies that for any $l \geq 1$ and $0 < \epsilon < 1$,
\begin{align*}
&\Bigl|\Bigl|\partial_x^{k_x} \partial_j^{k_j} \partial_b^{k_b}
\partial_{(\hat \rho_{\mu_i})}^{(k_{\mu_i})}
\partial_{(\varphi_\mu)}^{(l_\mu)} \partial_{(\hat \rho_\nu)}^{(k_\nu)}
\Bigl(\prod_\nu \frac{1}{\hat \rho_\nu^{l_\nu}}\Bigr)\partial_{(\varphi_\nu)}^{(l_\nu)}
\phi(a, b, x)\Bigr|\Bigr|_{C^l(\Sigma_0 \setminus N_0, (\R_1 \sqcup \R_2 \sqcup
\dots \sqcup \R_k) \times Y) \times E^0}\\
&\lesssim \prod_{\substack{i \\ k_{\mu_i} \neq 0}}
(\hat \rho_{\mu_i})^{\beta \tilde \delta'_{0, i} / 2 - k_{\mu_i}}
\prod_{\substack{\nu \\ (k_\nu, l_\nu) \neq (0,0)}}
(\hat \rho_\nu)^{\epsilon \alpha - (k_\nu + l_\nu)}.
\end{align*}
If $\alpha$ and $\beta$ are sufficiently large, then
$\beta \tilde \delta'_{0, i} / 2 - N > 0$ and
$\epsilon \alpha - N > 0$.

Hence the claim follows from the fact that if a continuous function $f$ on a manifold
$U$ is continuously differentiable on the complement of a submanifold $S \subset U$
and the limit of its differential on $S$ is zero, then $f$ is continuously differentiable
on the entire space $U$.
\end{proof}

Since we use the same coordinates for the neighborhoods of limit circles
of $\widetilde{P}_a$ as those of $\Sigma_0$,
the above argument also implies that for any limit circle $S^1_{\pm\infty_i}$
of $\Sigma_0$,
\begin{align*}
\hat V &\to P\\
(a, b, x) &\mapsto \pi_Y \circ \Phi_{a, b}(\xi_x) \circ \phi_{\pm\infty_i}
\end{align*}
is smooth if we fix a coordinate $\phi_{\pm\infty_i} : S^1 \tocong S^1_{\pm\infty_i}$.
Similarly,
\begin{align*}
\hat V &\to \R\\
(a, b, x) &\to
\lim_{s \to \infty}(\sigma \circ \Phi_{a, b}(\xi_x)|_{[0, \infty) \times S^1_{+\infty_i}}(s, t)
- (0_{k_0} + L_{+\infty_i} s))
\end{align*}
and
\begin{align*}
\hat V &\to \R\\
(a, b, x) &\to
\lim_{s \to -\infty}(\sigma \circ \Phi_{a, b}(\xi_x)|_{(-\infty, 0] \times S^1_{-\infty_i}}(s, t)
- (0_1 + L_{-\infty_i} s))
\end{align*}
are smooth since
\begin{align*}
&\lim_{s \to \infty}(\sigma \circ \Phi_{a, b}(\xi_x)|_{[0, \infty) \times S^1_{+\infty_i}}(s, t)
- (0_{k_0} + L_{+\infty_i} s))\\
& = \pi_{\R_{+\infty_i}} \xi_x
+ \lim_{s \to \infty}(\sigma \circ u_0|_{[0, \infty) \times S^1_{+\infty_i}}(s, t)
- (0_{k_0} + L_{+\infty_i} s))
\end{align*}
and
\begin{align*}
&\lim_{s \to -\infty}(\sigma \circ \Phi_{a, b}(\xi_x)|_{(-\infty, 0] \times S^1_{-\infty_i}}(s, t)
- (0_1 + L_{-\infty_i} s))\\
& = \pi_{\R_{-\infty_i}} \xi_x
+ \lim_{s \to -\infty}(\sigma \circ u_0|_{(-\infty, 0] \times S^1_{-\infty_i}}(s, t)
- (0_1 + L_{-\infty_i} s)),
\end{align*}
where $\pi_{\R_{\pm\infty_i}} : \widetilde{W}^{1, p}(\widetilde{P}_a; u_{a, b}^\ast T \hat Y)
\to \R$ is the projection to $\R\partial_\sigma \subset \Ker A_{\pm\infty_i}$.

\subsection{Embedding of Kuranishi neighborhoods}\label{embed}
In this section, we explain the way to construct an embedding of
a Kuranishi neighborhood $(V_1, E_1, s_1, \psi_1, G_1)$ to
another $(V_2, E_2, s_2, \psi_2, G_2)$.
Assume that $\psi_1(s_1^{-1}(0))$ and $\psi_2(s_2^{-1}(0))$ share a point
$q_0 \in \widehat{\M}$. We also assume that the additional marked points
$z_1^+$ for $(V_1, E_1, s_1, \psi_1, G_1)$ is a subsequence of $z_2^+$ for
$(V_2, E_2, s_2, \psi_2, G_2)$, and $E_1$ is a subspace of $E_2$ at $q_0$.
We do not assume any relationship between the additional temporary data
$(z_1^{++}, S'_1, \hat R^1_i)$ used for the description of $(V_1, E_1, s_1, \psi_1, G_1)$ and
$(z_2^{++} = (z_{2, i}^{++}), S'_2, \hat R^2_i)$ for $(V_2, E_2, s_2, \psi_2, G_2)$.
More precisely, we assume the following conditions:
\begin{itemize}
\item
For each $l = 1, 2$,
a Kuranishi neighborhood $(V_l, E_l, s_l, \psi_l)$
of a point $p_l = (\Sigma_l, z_l, u_l) \in \widehat{\M}(Y, \lambda, J)$ is defined
by the data $(p_l^+ = (\Sigma_l, z_l \cup z_l^+, u_l), S_l, E^0_l, \lambda_l)$
and the additional data $(z_l^{++} = (z_{l, i}^{++}), S'_l, \hat R^l_i)$.
Let $(\hat P_l \to \hat X_l, Z_l \cup Z_l^+ \cup Z_{\pm\infty_i})$ be the local universal
family of the stabilization $(\hat \Sigma_l, z_l \cup z_l^+ \cup (\pm\infty_i))$ of
the blown down curve of $(\Sigma_l, z_l \cup z_l^+)$, and
$(\widetilde{P}_l \to \widetilde{X}_l, Z_l \cup Z_l^+ \cup Z_l^{++})$ be the local universal
family of $(\Sigma_l, z_l \cup z_l^+ \cup z_l^{++})$.
\item
We assume $S_1 \subset S_2$. (We do not assume any correspondence between $S'_1$
and $S'_2$.)
\item
$q_0 = (\Sigma_0, z_0, u_0) \in \widehat{\M}$ is a point in the intersection
$\psi_1(s_1^{-1}(0)) \cap \psi_2(s_2^{-1}(0))$.
Hence there exist $(a^1_0, b^1_0, x^1_0) \in V_1$ and $(a^2_0, b^2_0, x^2_0) \in V_2$
such that $q_0 = \psi_1(a^1_0, b^1_0, x^1_0) = \psi_2(a^2_0, b^2_0, x^2_0)$.
We assume that there exist $\R$-translations
$\theta_1^0 : (\overline{\R}_1 \sqcup \overline{\R}_2 \sqcup \dots \sqcup
\overline{\R}_{k_1})/\sim_{a^1_0, b^1_0}
\to \overline{\R}_1 \cup \overline{\R}_2 \cup \dots \cup
\overline{\R}_{k_0}$,
$\theta_2^0 : (\overline{\R}_1 \sqcup \overline{\R}_2 \sqcup \dots \sqcup
\overline{\R}_{k_2})/\sim_{a^2_0, b^2_0}
\to \overline{\R}_1 \cup \overline{\R}_2 \cup \dots \cup
\overline{\R}_{k_0}$ and
an isomorphism
\[
\Xi_0 : ((\widetilde{P}_1)_{a^1_0}, Z_1(a^1_0)) \tocong ((\widetilde{P}_2)_{a^2_0}, Z_2(a^2_0))
\]
such that $(\theta_2^0 \times 1) \circ \Phi_{a^2_0, b^2_0}(\xi_{x^2_0}) \circ \Xi_0 =
(\theta_1^0 \times 1) \circ \Phi_{a^1_0, b^1_0}(\xi_{x_0^1})$.
\item
$\Xi_0$ maps the marked points $Z^+_1(a^1_0)$ to a subsequence $Z^+_{2|1}(a^2_0)$
of $Z^+_2(a^2_0)$.
\item
We denote by $\hat a^l \in \hat X_l$ the image of $a^l \in \widetilde{X}_l$
by the natural map $\widetilde{X}_l \to \hat X_l$.
Let $\hat U_1 \subset \hat X_1$ and $\hat U_2 \subset \hat X_2$ be small
neighborhoods of $\hat a^1_0$ and $\hat a^2_0$ respectively,
and let $\Theta : \hat P_2|_{\hat U_2} \to \hat P_1|_{\hat U_1}$ be the forgetful map
such that
\begin{itemize}
\item
it maps $Z \cup Z^+_{2|1} \cup Z_{\pm\infty_i}$ to $Z \cup Z^+_1 \cup Z_{\pm\infty_i}$,
\item
its underlying map $\hat U_2 \to \hat U_1$ maps $\hat a^2_0$ to $\hat a^1_0$, and
the isomorphism $\Theta|_{(\hat P_2)|_{\hat a^2_0}} : (\hat P_2)|_{\hat a^2_0} \cong
(\hat P_1)|_{\hat a^1_0}$ coincides with the map induced by $\Xi_0^{-1}$.
\end{itemize}
Let $\Theta^\ast \lambda_1 : E^0_1 \to C^\infty(\hat P_2 \times Y,
\Wedge^{0, 1} V^\ast \hat P_2 \otimes (\R \partial_\sigma \oplus TY))$
be the pull back of $\lambda_1$
by $\Theta$.
Then we assume that $E^0_1$ is embedded in $E^0_2$
as an $\Aut(\Sigma_0, z_0, u_0)$-vector space,
and $\Theta^\ast \lambda_1 = \lambda_2|_{E^0_1}$.
(Note that we may regard $\Aut(\Sigma_0, z_0, u_0)$ as a subgroup of
$\Aut(\Sigma_i, z_i, u_i)$ for each $i = 1, 2$.)
\end{itemize}

Under the above assumption, we define an $\Aut(\Sigma_0, z_0, u_0)$-equivariant
embedding $\phi$ of a neighborhood
$V^0_1$ of $(a^1_0, b^1_0, x^1_0) \in V_1$ to $V_2$ which makes the following diagrams
commutative.
\begin{equation}
\begin{tikzcd}
E_1 \ar[hookrightarrow]{r}{} & E_2\\
V_1^0 \ar{u}{s_1} \ar{r}{\phi} & V_2 \ar{u}{s_2}
\end{tikzcd}
\quad
\begin{tikzcd}
V_1^0 \cap s_1^{-1}(0) \ar{r}{\phi} \ar{dr}{\psi_1} & s_2^{-1}(0) \ar{d}{\psi_2}\\
&\widehat{\M}
\end{tikzcd}
\label{embed commute}
\end{equation}

We regard $V_1$ as a submanifold of
$\mathring{X}_1 \times C^{l_1}(\Sigma_1 \setminus N_1,
(\R_1 \cup \R_2 \cup \dots \cup \R_{k_1}) \times Y) \times E^0_1$
(see Section \ref{smoothness})
and write its point as $(a^1, b^1, u^1, h^1)$, where $(a^1, b^1) \in \mathring{X}_1 \subset
\widetilde{X}_1 \times \prod_\mu \R$, $u^1 \in C^{l_1}(\Sigma_1 \setminus N_1,
(\R_1 \cup \R_2 \cup \dots \cup \R_{k_1}) \times Y)$ and $h^1 \in E^0_1$.
Similarly, we write a point of $V_2$ as $(a^2, b^2, u^2, h^2)
\in \mathring{X}_2 \times C^{l_2}(\Sigma_2 \setminus N_2,
(\R_1 \cup \R_2 \cup \dots \cup \R_{k_1}) \times Y) \times E^0_2$.
We may assume $l_2 \ll l_1$ (since $l_1$, $l_2$ can be taken arbitrary).
The point $q_0$ corresponds to
$(a^1, b^1, u^1, h^1) = (a^1_0, b^1_0, \Phi_{a^1_0, b^1_0}(\xi_{x_0^1}), 0)$
and $(a^2, b^2, u^2, h^2) = (a^2_0, b^2_0, \Phi_{a^2_0, b^2_0}(\xi_{x^2_0}), 0)$.
The embedding $(a^1, b^1, u^1, h^1) \mapsto (a^2, b^2, u^2, h^2)$ is defined
by the following steps.

First, $h^2$ is the image of $h^1$ by the inclusion map $E^0_1 \inj E^0_2$.
This map is obviously smooth.

Next, prior to defining $a^2 \in \widetilde{X}_2$, we define $\hat a^2 \in \hat U_2$
which should be the image of $a^2$ by the natural map $\widetilde{X}_2 \to \hat X_2$.
$\hat a^2 \in \hat U_2$ is the point in the inverse image of $\hat a^1$ by
$\hat U_2 \to \hat U_1$ (the underlying map of $\Theta$) such that
\[
(\pi_Y \circ u^1) \circ (\pi_1|_{(\widetilde{P}_1)_{a^1}})^{-1} \circ
\Theta|_{(\hat P_2)_{\hat a^2}} (Z^+_2(\hat a^2)) \subset S_2,
\]
where $\pi_1 : \widetilde{P}_1 \to \hat P_1$ is the composition of the blow down and
the forgetful map.
Since $\hat U_2 \to \hat U_1$ is a submersion and its fiber is the product of
neighborhoods of the points $Z^+_2(\hat a^2_0) \setminus Z^+_{2|1}(\hat a^2_0)$
in $\Sigma_2 \setminus N_2$, $\hat a^2$ is a smooth function of $(a^1, b^1, u^1, h^1)$.
We denote the sequence of points $(\pi_1|_{(\widetilde{P}_1)_{a^1}})^{-1} \circ
\Theta|_{(\hat P_2)_{\hat a^2}} (Z^+_2(\hat a^2))
\subset (\widetilde{P}_1)_{a^1}$ by $\mathcal{Z}^+_2
= \mathcal{Z}^+_2(a^1, u^1)$.

Define an $\R$-gluing $\theta = \theta_{(a^1, b^1, u^1)}: \overline{\R}_1 \sqcup
\overline{\R}_2 \sqcup \dots \sqcup \overline{\R}_{k_2} \to (\overline{\R}_1 \sqcup
\overline{\R}_2 \sqcup \dots \sqcup \overline{\R}_{k_1}) / \sim_{a^1, b^1}$ by
\[
\theta(0_i) = \sigma \circ u^1 \circ (\pi_1|_{(\widetilde{P}_1)_{a^1}})^{-1} \circ
\Theta|_{(\hat P_2)|_{\hat a^2}} (\hat R^2_i(\hat a^2)),
\]
and let
$\mathcal{Z}^{++}_2 = \mathcal{Z}^{++}_2(a^1, b^1, u^1)
\subset (\widetilde{P}_1)_{a^1}$ be the points near
$\Xi_0^{-1} (Z^{++}_2(\hat a^2_0)) \subset (\widetilde{P}_1)_{a^1_0}
\subset \widetilde{P}_1$ such that
$u^1 (\mathcal{Z}^{++}_2) \subset (\theta \times 1) (S'_2)$.
We assume that $Z^{++}_2(\hat a^2_0) \subset \Sigma_0$ is contained in
$\Sigma_1 \setminus N_1 \subset \Sigma_0$.
Then $\mathcal{Z}^{++}_2$ is a smooth function of $(a^1, b^1, u^1)
\in \widetilde{X}_1 \times
C^{l_1}(\Sigma_1 \setminus N_1, (\R_1 \cup \R_2 \cup \dots \cup \R_{k_1}) \times Y)$.

Let $\widetilde{U}_2 \subset \widetilde{X}_2$ be a neighborhood of $a^2_0$ and let
$a^2 \in \widetilde{U}_2$ be the point over $\hat a^2$ such that
there exists an isomorphism $\Xi_{(a^1, b^1, u^1, h^1)} : (\widetilde{P}_1)_{a^1}
\cong (\widetilde{P}_2)_{a^2}$ which maps
$Z(a^1)$, $\mathcal{Z}^+_2$ and $\mathcal{Z}^{++}_2$ to
$Z(a^2)$, $Z^+_2(a^2)$ and $Z^{++}_2(a^2)$ respectively.
Then $a^2$ is a smooth function of $(a^1, b^1, u^1, h^1)$.
In fact, the map $\Xi : V^0_1 \times_{\widetilde{X}_1} \widetilde{P}_1
\to \widetilde{P}_2$ is smooth because it is the composition of
\begin{itemize}
\item
the map from $V^0_1 \times_{\widetilde{X}_1} \widetilde{P}_1$
to the local universal family $\widetilde{P}_3$ of
$(\Sigma_0, z_0 \cup Z^+_2(a^2_0) \cup Z^{++}_1(a^1_0) \cup Z^{++}_2(a^2_0))$
which maps the marked points $Z$, $\mathcal{Z}^+_2$, $Z^{++}_1$ and
$\mathcal{Z}^{++}_2$ to the corresponding marked points of $\widetilde{P}_3$, and
\item
the forgetful map $\widetilde{P}_3$ to $\widetilde{P}_2$.
\end{itemize}
(We assume that $Z^+_2(a^2_0), Z^{++}_1(a^1_0), Z^{++}_2(a^2_0) \subset \Sigma_0$
are disjoint temporarily.)

We define $u^2 \in C^{l_2}(\Sigma_2 \setminus N_2,
(\R_1 \cup \R_2 \cup \dots \cup \R_{k_2}) \times Y)$ by
\[
u^2 = (\theta_{(a^1, b^1, u^1, h^1)} \times 1)^{-1} \circ u^1 \circ
(\Xi_{(a^1, b^1, u^1, h^1)})^{-1},
\]
where we assume $\Xi_{(a^1, b^1, u^1, h^1)}(N_1) \subset N_2$.
Then this is a smooth function of $(a^1, b^1, u^1, h^1)$ (since $l_2 \ll l_1$).
Then it is easy to see that
$\sigma_i \circ u^2 (\widetilde{R}_i^2(a^2)) = 0$ and
$u^2(Z^{++}_2(a^2)) \subset S'_2$.

Finally, we define the asymptotic parameter $b^2_\mu$.
First we recall the relationship between the parameter $b_\mu^l$ and the map
$u^l$ ($l = 1, 2$).
We denote the index set of joint circles of $\Sigma_l$
between the $j$-th floor and the $(j+1)$-th floor by $M^l_j$.
If $\mu \in M^l_j$ and $\rho_\mu^l \neq 0$, then
\[
b^l_\mu = (\theta_l \circ \sigma \circ u^l (\widetilde{R}^l_{j+1})
- \theta_l \circ \sigma \circ u^l (\widetilde{R}^l_j)) + L_\mu \log \rho^l_\mu,
\]
where $\theta_l : \overline{R}_1 \sqcup \overline{R}_2 \sqcup \dots
\sqcup \overline{R}_{k_1} \to \overline{R}_1 \cup \overline{R}_2 \cup \dots \cup
\overline{R}_k$ is an $\R$-gluing which induces an $\R$-translation
$(\overline{R}_1 \sqcup \overline{R}_2 \sqcup \dots \sqcup \overline{R}_{k_1})
/ \sim_{a^l, b^l} \to \overline{R}_1 \cup \overline{R}_2 \cup \dots \cup \overline{R}_k$.
If $\rho^l_\mu = 0$, then
\begin{align*}
b^l_\mu &= \lim_{s \to \infty} \bigl(\theta_l \circ \sigma \circ
u^l|_{[0, \infty) \times S^1_\mu}(s, t)
- \theta_l \circ \sigma \circ u^l(\widetilde{R}^l_j(a^l)) - L_\mu s \bigr)\\
&\quad - \lim_{s \to -\infty} \bigl(\theta_l \circ \sigma \circ
u^l|_{(-\infty, 0] \times S^1_\mu}(s, t)
- \theta_l \circ \sigma \circ u^l(\widetilde{R}^l_{j + 1}(a^l)) - L_\mu s \bigr)
\end{align*}
Since $u^1$ and $u^2$ represent the same curve, we may assume
$\theta_2 = \theta_1 \circ \theta_{(a^1, b^1, u^1, h^1)}$.

Assume $\mu \in M^2_j$ corresponds to $\iota(\mu) \in M^1_i$.
If $\rho^2_\mu \neq 0$, then
\begin{align}
b^2_\mu &= b^1_{\iota(\mu)} +
\bigl(\theta_1 \circ \sigma \circ u^1 \circ
(\Xi_{(a^1, b^1, u^1, h^1)})^{-1} (\widetilde{R}^2_{i + 1}(a^2))
- \theta_1 \circ \sigma \circ u^1 (\widetilde{R}^1_{j + 1}(a^1))\bigr) \notag\\
&\quad - \bigl(\theta_1 \circ \sigma \circ u^1 \circ
(\Xi_{(a^1, b^1, u^1, h^1)})^{-1}(\widetilde{R}^2_i(a^2))
- \theta_1 \circ \sigma \circ u^1 (\widetilde{R}^1_j(a^1))\bigr) \notag\\
& \quad + L_\mu (-\log \rho^1_{\iota(\mu)} + \log \rho^2_\mu) \notag\\
&= b^1_{\iota(\mu)} +
\bigl( \sigma \circ u^1 \circ (\Xi_{(a^1, b^1, u^1, h^1)})^{-1} (\widetilde{R}^2_{i + 1}(a^2))
- \sigma \circ u^1 (\widetilde{R}^1_{j + 1}(a^1))\bigr) \notag\\
&\quad - \bigl(\sigma \circ u^1 \circ (\Xi_{(a^1, b^1, u^1, h^1)})^{-1}
(\widetilde{R}^2_i(a^2)) - \sigma \circ u^1 (\widetilde{R}^1_j(a^1))\bigr) \notag\\
& \quad + L_\mu (-\log \rho^1_{\iota(\mu)} + \log \rho^2_\mu).
\label{b for nonzero kappa}
\end{align}

If $\rho^2_\mu = 0$, then
\begin{align}
b^2_\mu &= b^1_{\iota(\mu)} \notag\\
&\quad + \lim_{s \to \infty} \bigl(\theta_1 \circ \sigma \circ u^1 \circ
(\Xi_{(a^1, b^1, u^1, h^1)})^{-1}|_{[0, \infty) \times S^1_\mu}(s, t) \notag\\
&\quad \hph{+ \lim_{s \to \infty} \bigl(}
- \theta_1 \circ \sigma \circ u^1|_{[0, \infty) \times S^1_{\iota(\mu)}}(s, t) \bigr)
\notag\\
&\quad - \lim_{s \to -\infty} \bigl(\theta_1 \circ \sigma \circ u^1 \circ
(\Xi_{(a^1, b^1, u^1, h^1)})^{-1}|_{[0, \infty) \times S^1_\mu}(s, t) \notag\\
&\quad \hph{- \lim_{s \to -\infty} \bigl(}
- \theta_1 \circ \sigma \circ u^1|_{[0, \infty) \times S^1_{\iota(\mu)}}(s, t) \bigr)
\notag\\
&\quad + (\theta_1 \circ \sigma \circ u^1 \circ
(\Xi_{(a^1, b^1, u^1, h^1)})^{-1}(\widetilde{R}^2_{j+1}(a^2))
- \theta_1 \circ \sigma \circ u^1(\widetilde{R}^1_{i+1}(a^1))) \notag\\
& \quad - (\theta_1 \circ \sigma \circ u^1 \circ
(\Xi_{(a^1, b^1, u^1, h^1)})^{-1}(\widetilde{R}^2_j(a^2))
- \theta_1 \circ \sigma \circ u^1(\widetilde{R}^1_i(a^1))) \notag\\
&= b^1_{\iota(\mu)} +
\lim_{s \to \infty} L_\mu (p_1 \circ (\Xi_{(a^1, b^1, u^1, h^1)})^{-1}
|_{[0, \infty) \times S^1_\mu}(s, t) - s) \notag\\
&\quad - \lim_{s \to -\infty} L_\mu (p_1 \circ (\Xi_{(a^1, b^1, u^1, h^1)})^{-1}
|_{[0, \infty) \times S^1_\mu}(s, t) -s) \notag\\
&\quad + (\sigma \circ u^1 \circ (\Xi_{(a^1, b^1, u^1, h^1)})^{-1}(\widetilde{R}^2_{j+1}(a^2))
- \sigma \circ u^1(\widetilde{R}^1_{i+1}(a^1))) \notag\\
& \quad - (\sigma \circ u^1 \circ (\Xi_{(a^1, b^1, u^1, h^1)})^{-1}(\widetilde{R}^2_j(a^2))
- \sigma \circ u^1(\widetilde{R}^1_i(a^1))), \label{b for zero kappa}
\end{align}
where $p_1(s, t) = s$ is the projection,
and we have used the asymptotic behavior of $u^1$ near the joint circle
$S^1_{\iota(\mu)}$ for the last equality.

We define $b^2_\mu \in \R$ by (\ref{b for nonzero kappa}) and
(\ref{b for zero kappa}).
It is clear that this is a smooth function of $(a^1, b^1, u^1, h^1)$
at $\rho^2_\mu \neq 0$.
We need to prove the smoothness at $\rho^2_\mu = 0$.
We note that we may assume that if the function $\rho^2_\mu =
\rho^2_\mu(a^1, b^1, u^1, h^1)$ can take zero,
then $\Sigma_0$ has a joint circle corresponding to $\mu$.

To prove the smoothness of $b^2_\mu$, we need to study the map $\Xi$.
First we claim that there exists a smooth function
$f : V^0_1 \to \C^\ast = \C \setminus 0$ such that
\begin{equation}
(\rho^2_\mu)^{2\pi} e^{\sqrt{-1} \varphi_\mu^2}
= (\rho^1_{\iota(\mu)})^{2\pi} e^{\sqrt{-1} \varphi_{\iota(\mu)}^1} f(a^1, b^1, u^1, h^1).
\label{bf eq}
\end{equation}

To prove this claim,
recall that $\Xi$ is the composition of the map $\Xi^1 : V^0_1
\times_{\widetilde{X}_1} \widetilde{P}_1 \to \widetilde{P}_3$
and the forgetful map $\Xi^2 : \widetilde{P}_3 \to \widetilde{P}_2$,
where $(\widetilde{P}_3 \to \widetilde{X}_3, Z_3 \cup Z_3^+ \cup Z_{3, 1}^{++}
\cup Z_{3, 2}^{++})$ is the local universal family of $(\Sigma_0, z_0 \cup Z^+_2(a^2_0)
\cup Z^{++}_1(a^1_0) \cup Z^{++}_2(a^2_0))$.
Since $((\widetilde{P}_3)_0, Z_3(0) \cup Z_3^+(0) \cup Z_{3, 1}^{++}(0))$ is stable and
isomorphic to $((\widetilde{P}_1)_{a^1_0}, Z_1(a^1_0) \cup Z_1^+(a^1_0) \cup
Z_1^{++}(a^1_0))$, we may assume that there exists a neighborhood $U^0_1 \subset
\widetilde{X}_1$ of $a^1_0$ such that $(\widetilde{P}_3 \to \widetilde{X}_3, Z_3 \cup
Z_3^+ \cup Z_{3, 1}^{++} \cup Z_{3, 2}^{++})$ is isomorphic to the product of
$(\widetilde{P}_1|_{U^0_1} \to U^0_1, Z_1 \cup Z_1^+ \cup Z_1^{++})$ and
the parameter space $D^m$ for the marked point $Z_{3, 2}^{++}$.
We can use the coordinate of $\widetilde{X}_3$ defined by the isomorphism
$\widetilde{X}_3 \cong U^0_1 \times D^m \subset \widetilde{X}_1 \times D^m$

Let $S^1_\mu$ be a joint circle of $\Sigma_0 \cong (\widetilde{P}_3)_0$.
Let $S_{\iota_1(\mu)}^1$ and $S_{\iota_2(\mu)}^1$ be the corresponding joint circles
of $\Sigma_1$ and $\Sigma_2$ respectively.
Since the forgetful map $\Xi^2 : \widetilde{P}_3 \to \widetilde{P}_2$ is induced by its
blow down and it is a holomorphic map, there exists a smooth function
$f'_\mu : \widetilde{X}_3 \to \C^\ast$ such that
\[
(\rho^2_{\iota_2(\mu)})^{2\pi} e^{2\pi \sqrt{-1} \varphi^2_{\iota_2(\mu)}}
= (\rho^1_{\iota_1(\mu)})^{2\pi} e^{2\pi \sqrt{-1} \varphi^1_{\iota_1(\mu)}}
\cdot f'_\mu (a^3)
\]
for all $a^3 \in \widetilde{X}_3$, where
$(\rho^1_{\iota_1(\mu)}, \varphi^1_{\iota_1(\mu)})$ is a part of the coordinate
of $a^3 \in \widetilde{X}_3$ under the isomorphism
$\widetilde{X}_3 \cong U^0_1 \times D^m \subset \widetilde{X}_1 \times D^m$, and
$(\rho^2_{\iota_2(\mu)}, \varphi^2_{\iota_2(\mu)})$ is a part of the coordinate of
$\widetilde{X}_2$ at $\Xi^2(a^3)$.

Since the underlying map of $\Xi^1$ is smooth, the claim follows, that is, 
there exists a smooth function
$f : V^0_1 \to \C^\ast = \C \setminus 0$ which satisfies equation (\ref{bf eq}).

Similarly, there exists smooth maps $f_\mu^{\text{left}}, f_\mu^{\text{right}} :
V^0_1 \times_{\widetilde{X}_2} \widetilde{P}_2 \to \C^\ast$ such that
if
\begin{align*}
\Xi_{(a^1, b^1, u^1, h^1)}|_{[0, \infty) \times S^1_{\iota_1(\mu)}}
(s_1^{\text{left}}, t_1^{\text{left}}) &= (s_2^{\text{left}}, t_2^{\text{left}})
\in [0, \infty) \times S^1_{\iota_2(\mu)},\\
\Xi_{(a^1, b^1, u^1, h^1)}|_{(-\infty, 0] \times S^1_{\iota_1(\mu)}}
(s_1^{\text{right}}, t_1^{\text{right}}) &= (s_2^{\text{right}}, t_2^{\text{right}})
\in (-\infty, 0] \times S^1_{\iota_2(\mu)},
\end{align*}
then
\begin{align*}
e^{-2\pi (s_2^{\text{left}} + \sqrt{-1} t_2^{\text{left}})}
&= e^{-2\pi (s_1^{\text{left}} + \sqrt{-1} t_1^{\text{left}})}
\cdot f_\mu^{\text{left}} (s_2^{\text{left}}, t_2^{\text{left}}, a^1, b^1, u^1, h^1),\\
e^{2\pi (s_2^{\text{right}} + \sqrt{-1} t_2^{\text{right}})}
&= e^{2\pi (s_1^{\text{right}} + \sqrt{-1} t_1^{\text{right}})}
\cdot f_\mu^{\text{right}} (s_2^{\text{right}}, t_2^{\text{right}}, a^1, b^1, u^1, h^1).
\end{align*}
Note that $f_\mu$, $f_\mu^{\text{left}}$ and $f_\mu^{\text{right}}$ satisfy
\[
f_\mu(a^1, b^1, u^1, h^1)
= f_\mu^{\text{left}}(s_2^{\text{left}}, t_2^{\text{left}},a^1, b^1, u^1, h^1) 
f_\mu^{\text{right}}(s_2^{\text{right}}, t_2^{\text{right}}, a^1, b^1, u^1, h^1)
\]
if $(s_2^{\text{left}}, t_2^{\text{left}})$ and $(s_2^{\text{right}}, t_2^{\text{right}})$
denote the same point of $(\widetilde{P}_1)_{a^1}$.
In particular, if $\kappa^2_{\iota_2(\mu)} = 0$, then
\[
|f_\mu(a^1, b^1, u^1, h^1)| = \lim_{s \to \infty}
|f_\mu^{\text{left}}(s, t,a^1, b^1, u^1, h^1)|
\lim_{s \to -\infty}
|f_\mu^{\text{right}}(s, t, a^1, b^1, u^1, h^1)|.
\]

We can rewrite the formula of $b^2_{\iota_2(\mu)}$ by using the function $f_\mu$
as follows.
If $\rho^2_{\iota_2(\mu)} \neq 0$, then
\begin{align}
b^2_\mu &= b^1_\mu +
\bigl( \sigma \circ u^1 \circ (\Xi_{(a^1, b^1, u^1, h^1)})^{-1} (\widetilde{R}^2_{i + 1}(a^2))
- \sigma \circ u^1 (\widetilde{R}^1_{j + 1}(a^1))\bigr) \notag\\
&\quad - \bigl(\sigma \circ u^1 \circ (\Xi_{(a^1, b^1, u^1, h^1)})^{-1}
(\widetilde{R}^2_i(a^2)) - \sigma \circ u^1 (\widetilde{R}^1_j(a^1))\bigr)\notag\\
& \quad + L_\mu \cdot \frac{1}{2\pi} \log |f_\mu(a^1, b^1, u^1, h^1)|.
\label{b for any kappa}
\end{align}
If $\rho^2_{\iota_2(\mu)} = 0$, then
\[
L_\mu (p_1 \circ (\Xi_{(a^1, b^1, u^1, h^1)})^{-1}
|_{[0, \infty) \times S^1_\mu}(s, t) - s)
= \frac{1}{2\pi} \log |f_\mu^{\text{left}}(s, t, a^1, b^1, u^1, h^1)|,
\]
and
\[
L_\mu (p_1 \circ (\Xi_{(a^1, b^1, u^1, h^1)})^{-1}
|_{[0, \infty) \times S^1_\mu}(s, t) - s)
= - \frac{1}{2\pi} \log |f_\mu^{\text{right}}(s, t, a^1, b^1, u^1, h^1)|.
\]
Therefore, equation (\ref{b for any kappa}) also holds in this case.
Hence $b^2_{\iota_2(\mu)}$ is a smooth function of $(a^1, b^1, u^1, h^1)$.

Next we prove that the differential of
$\phi : (a^1, b^1, u^1, h^1) \mapsto (a^2, b^2, u^2, h^2)$ is injective everywhere.
It is enough to construct a smooth inverse from an open subset of the submanifold
$\{(a^2, b^2, u^2, h^2) \in V^2; h^2 \in E^0_1\}$ of $V^2$ to $V^1$.
We can construct this map by the same way as $\phi$.
Hence $\phi$ is indeed an embedding.

It is obvious that diagrams (\ref{embed commute}) are commutative for this $\phi$.
Hence it is the required embedding.

So far we have made some assumptions about the position of
the additional marked points $Z^{++}_1(a^1_0)$ or $Z^{++}_2(a^2_0)$.
(For example, we have assumed that $Z^{++}_2(a^2_0) \subset \Sigma_0$ is contained
in $\Sigma_1 \setminus N_1 \subset \Sigma_0$.)
We can remove these assumption because
two Kuranishi neighborhoods of the same point $p_1$ defined by using
the same data $(p_1^+ = (\Sigma_1, z_1 \cup z_1^+, u_1), S_1, E^0_1, \lambda_1)$,
different additional data $(\hat R^1_i, \widetilde{R}^1_i, Z^{++}_{1, i}, S'_1)
\neq (\hat R^2_i, \widetilde{R}^2_i, Z^{++}_{2, i}, S'_2)$ and
different decompositions of the curve $\Sigma_1$ into parts are isomorphic
by the same argument as above.

Note that the smoothness of $\Aut(\Sigma_0, z, u_0)$-action on
a Kuranishi neighborhood of $p_0 = (\Sigma_0, z, u_0)$ also
follows from the above argument because the group action is also a kind of embedding.

\begin{rem}\label{natural projection from hat V to V}
Recall that for a Kuranishi neighborhood $(V, E, s, \psi, G)$,
$V$ is a submanifold of $\hat V = \mathring{X} \times B_\epsilon(0) \subset
\mathring{X} \times C^l(\Sigma_0 \setminus N_0,
(\R_1 \cup \R_2 \cup \dots \cup \R_k) \times Y) \times E^0$.
We can define a natural projection $\hat V \to V : (a', b', u', h') \mapsto (a, b, u, h)$
by a similar way to the above embedding.
It satisfies $h = h'$ and that for each $(a', b', u', h')$,
there exists a biholomorphism $\varphi : \widetilde{P}_{a'} \to \widetilde{P}_a$ and
an $\R$-translation $\theta$ such that
$\varphi$ maps $Z(a')$ and $Z^+(a')$ to $Z(a)$ and $Z^+(a)$ as sequences respectively
and $(\theta \times 1) \circ u \circ \varphi = u'$.
\end{rem}

\subsection{A Kuranishi neighborhood of a disconnected holomorphic building
and those of its connected components}
\label{Kuranishi of disconnected buildings}
In Section \ref{construction of nbds}, we constructed a Kuranishi neighborhood of
a point $p \in \widehat{\M}$ for data $(p^+, S, E^0, \lambda)$ and additional
data $(z^{++}, S', \hat R_j)$.
In Section \ref{embed}, we saw that the Kuranishi neighborhood is determined by
the data $(p^+, S, E^0, \lambda)$ and independent of the
additional data $(z^{++}, S', \hat R_j)$.
To obtain algebraic information of the moduli space,
the data $(p^+, S, E^0, \lambda)$ for a disconnected holomorphic building
should be given by the product of the data
for the connected components, but the additional data $(z^{++}, S', \hat R_j)$ can be
taken independently.
We emphasize that we do not construct a Kuranishi neighborhood from the Kuranishi
neighborhoods of the connected components as a quotient of their product, but
we construct it independently through the same procedure.
Instead, under the above assumption for the data $(p^+, S, E^0, \lambda)$,
we study the relationship between the Kuranishi neighborhood of
a disconnected holomorphic building and those of its connected components.

Let $p'_0 = (\Sigma'_0, z'_0, u'_0) \in \widehat{\M}$ be an arbitrary point and
let $\Sigma'_1, \Sigma'_2, \dots, \Sigma'_N$ be the connected components of
$\Sigma'_0$.
Let $p_i = (\Sigma_i, z_i, u_i)$ be the holomorphic building obtained by collapsing
the floors of $(\Sigma'_i, z'_0|_{\Sigma'_i}, u'_0|_{\Sigma'_i})$ consisting of trivial cylinders.
We call a map $\theta : \overline{\R}_1 \sqcup \overline{\R}_2 \sqcup \dots \sqcup
\overline{\R}_k \to \overline{\R}_1 \cup \overline{\R}_2 \cup \dots \cup
\overline{\R}_l$ an $\R$-compressing if there exist a map $\mu : \{1, 2, \dots, k\} \to
\{\frac{1}{2}, 1, \frac{3}{2}, \dots, l, l + \frac{1}{2}\}$ and constants $c_i \in \R$ ($i \in
\mu^{-1}(\Z)$) such that
\begin{itemize}
\item
the image of $\mu$ contains $\{1, 2, \dots, l\}$,
\item
if $i \leq j$ then $\mu(i) \leq \mu(j)$,
\item
if $\mu(i) \in \Z$ then $\theta(\overline{\R}_i) = \overline{\R}_{\mu(i)}$ and
$\theta|_{\overline{\R}_i}(s) = s + c_i$, and
\item
if $\mu(i) \notin \Z$ then $\theta(\overline{\R}_i) = \{+\infty_{\lfloor \mu(i) \rfloor}\}
= \{-\infty_{\lceil \mu(i) \rceil}\}$.
\end{itemize}

As in the previous section, we assume the following conditions on the Kuranishi
neighborhoods of $p'_0$ and $p_i$:
\begin{itemize}
\item
Kuranishi neighborhoods $(V_i, E_i, s_i, \psi_i, G_i)$
of $p_i = (\Sigma_i, z_i, u_i) \in \widehat{\M}^0$ is defined by the data
$(p_i^+ = (\Sigma_i, z_i \cup z_i^+, u_i), S_i, E^0_i, \lambda_i)$
and the additional data $(z_i^{++} = (z_{i, i}^{++}), S'_i, \hat R^i_j)$.
Let $(\hat P_i \to \hat X_i, Z_i \cup Z_i^+ \cup Z_{\pm\infty_i})$ be the local universal
family of the stabilization $(\hat \Sigma_i, z_i \cup z_i^+ \cup (\pm\infty_i))$ of
the blown down curve of $(\Sigma_i, z_i \cup z_i^+)$, and
$(\widetilde{P}_i \to \widetilde{X}_i, Z_i \cup Z_i^+ \cup Z_i^{++})$ be the local universal
family of $(\Sigma_i, z_i \cup z_i^+ \cup z_i^{++})$.
\item
A Kuranishi neighborhood $(V_0, E_0, s_0, \psi_0, G_0)$
of $p'_0 = (\Sigma'_0, z'_0, u'_0) \in \widehat{\M}$ is defined by the data
$({p'_0}^+ = (\Sigma'_0, z'_0 \cup \bigcup_i z_i^+, u'_0), S_0 = \bigcup_i S_i,
E^0_0 = \bigoplus_i E^0_i, \lambda_0 = \bigoplus_i \lambda_i)$
and the additional data $({z'_0}^{++} = ({z'_{0, j}}^{++}), S'_0, \hat R^0_i)$.
We define ${z'_0}^+ = \bigcup_i z_i^+$.
Let $(\hat P'_0 \to \hat X'_0, Z'_0 \cup {Z'_0}^+ \cup Z_{\pm\infty_i})$ be the local
universal family of the stabilization $(\hat \Sigma'_0, z'_0 \cup {z'_0}^+ \cup (\pm\infty_i))$
of the blown down curve of $(\Sigma'_0, z'_0 \cup {z'_0}^+)$, and
$(\widetilde{P}'_0 \to \widetilde{X}'_0, Z'_0 \cup {Z'_0}^+ \cup {Z'_0}^{++})$
be the local universal family of $(\Sigma'_0, z'_0 \cup {z'_0}^+ \cup {z'_0}^{++})$.
\item
Let $\widetilde{P}'_0 = \coprod_i \widetilde{P}'_i$ be the decomposition into the connected
components corresponding to the decomposition $\Sigma'_0 = \coprod_i \Sigma'_i$.
We define $Z'_i = Z'_0 \cap \widetilde{P}'_i$, ${Z'_i}^+ = {Z'_0}^+ \cap \widetilde{P}'_i$ and
${Z'_i}^{++} = {Z'_0}^{++} \cap \widetilde{P}'_i$.
We assume that for each $i$, there exists a map
\[
\Xi_{0, i} : ((\widetilde{P}'_i)_0, Z'_i(0) \cup {Z'_i}^+(0))
\to ((\widetilde{P}_i)_0, Z_i(0) \cup Z_i^+(0))
\]
which collapses the floors consisting of trivial cylinders and which satisfies
$u_i \circ \Xi_{0, i} = (\theta_0 \times 1) \circ u_0|_{\Sigma'_0}$
on $\Sigma'_i \cong (\widetilde{P}'_i)_0$ for some $\R$-compressing
$\theta_0 : \overline{\R}_1 \sqcup \overline{\R}_2 \sqcup \dots \sqcup
\overline{\R}_{k_0} \to \overline{\R}_1 \cup \overline{\R}_2 \cup \dots \cup
\overline{\R}_{k_i}$.
\end{itemize}

Under theses assumptions, we prove that there exists a natural map
$\phi = (\phi_i) : V_0 \to \prod_i V_i$ which satisfies the following conditions:
\begin{itemize}
\item
The following diagram is commutative.
\[
\begin{tikzcd}
E_0 \ar{r}{} & E_i\\
V_0 \ar{u}{s_0} \ar{r}{\phi_i} & V_i \ar{u}{s_i}
\end{tikzcd}
\]
\item
For any $p \in s_0^{-1}(0)$,
the curve obtained by collapsing trivial floors of the $i$-th connected component of the
curve corresponding to the point $p$ is isomorphic to the curve corresponding to
$\phi_i(p) \in s_i^{-1}(0)$.
\item
For arbitrary integers $l_i \leq 1$ ($i = 1, \dots, N$), let $V^{(l_i)}_i \subset V_i$ be
the submanifold consisting of height-$l_i$ curves.
(This coincides with $\mathring{\partial}^{l_i - 1} V_i$
in Section \ref{section of orbibundle}.)
Then each $\phi^{-1}(\prod_i V_i^{(l_i)}) \subset V_0$ is a union
of the interiors of corners of $V_0$, and
$\phi|_{\phi^{-1}(\prod_i V_i^{(l_i)})} : \phi^{-1}(\prod_i V_i^{(l_i)})
\to \prod_i V_i^{(l_i)}$ is submersive on each of them.
We say that $\phi$ is essentially submersive if it satisfies this condition.
\end{itemize}

Fixing $i_0$, we construct a essential submersion $\phi_{i_0} : V_0 \to V_{i_0}$ as follows.
As in the previous section, we write a point of $V_0$ as $(a^0, b^0, u^0, h^0)$,
where $(a^0, b^0) \in \mathring{X}'_0 \subset \widetilde{X}'_0 \times \prod_\mu \R$,
$u^0 \in C^{l_0}(\Sigma'_0 \setminus N'_0, (\R_1 \cup \R_2 \cup \dots \cup \R_{k_0})
\times Y)$ and $h^0 = (h^0_i) \in E^0_0 = \bigoplus_i E^0_i$.
Similarly, a point of $V_{i_0}$ is written as
$(a^{i_0}, b^{i_0}, u^{i_0}, h^{i_0}) \in \widetilde{X}_{i_0} \times \prod_\mu \R \times
C^{l_{i_0}}(\Sigma_{i_0} \setminus N_{i_0}, (\R_1 \cup \R_2 \cup \dots \cup \R_{k_{i_0}})
\times Y) \times E^0_{i_0}$.
We may assume $l_{i_0} \ll l_0$.
The essential submersion $(a^0, b^0, u^0, h^0) \mapsto (a^{i_0}, b^{i_0}, u^{i_0}, h^{i_0})$
is defined by the following steps as in the previous section.

First, $h^{i_0} \in E^0_{i_0}$ is defined by $h^{i_0} = h^0_{i_0}$.
Next we note that $\hat X'_0 = \prod_i \hat X_i$ is a product
(but $\widetilde{X}'_0$ is not).
Hence we can define $\hat a^{i_0} \in \hat X_{i_0}$ by the $i_0$-th component
of $\hat a^0 = (a^0_i) \in \prod_i \hat X_i$.

Let $\pi_{i_0} : (\widetilde{P}'_{i_0}, Z'_{i_0} \cup {Z'_{i_0}}^+)
\to (\hat P_{i_0}, Z_{i_0} \cup Z_{i_0}^+)$ be the composition of the blow down and the
forgetful map, and define
\[
s_j = \sigma \circ u^0 \circ (\pi_{i_0}|_{(\widetilde{P}'_{i_0})_{\hat a^0}})^{-1}
(\hat R^{i_0}_j(\hat a^{i_0})) \in \R_1 \cup \R_2 \cup \dots \cup \R_{k_0}.
\]
Let $\theta = \theta_{(a^0, b^0, u^0, h^0)} : \overline{\R}_1 \sqcup \overline{\R}_2
\sqcup \dots \sqcup \overline{\R}_{k_0} \to \overline{\R}_1 \cup \overline{\R}_2
\cup \dots \cup \overline{\R}_{k_{i_0}}$ be the $\R$-compressing defined by the
following conditions:
\begin{itemize}
\item
If $s_j \in \R_i$ then $\theta(\overline{\R}_i) = \overline{\R}_j$ and
$\theta|_{\overline{\R}_i}(s) = s - s_j$.
\item
If $\R_i$ does not contain any $s_j$ then $\theta$ maps $\overline{\R}_i$ to some
$\infty$-point. More precisely, if $s_j \in \bigcup_{l < i} \R_l$ and $s_{j+1} \notin
\bigcup_{l \leq i} \R_l$ then $\theta(\overline{\R}_i) = \{+\infty_j\} \subset \overline{\R}_j$.
\end{itemize}
Let $\mathcal{Z}^{++} = \mathcal{Z}^{++}(a^0, b^0, u^0, h^0) \subset
(\widetilde{P}'_{i_0})_{a^0}$ be the sequence of points in a neighborhood of
$\Xi_{0, i_0}^{-1}(Z^{++}_{i_0}(0)) \subset (\widetilde{P}'_{i_0})_0 \subset \widetilde{P}'_{i_0}$
defined by $(\theta \times 1) \circ u^0 (\mathcal{Z}^{++}) \subset S'_{i_0}$.
Let $\Xi : V_0 \times_{\widetilde{X}'_0} \widetilde{P}'_{i_0} \to \widetilde{P}_{i_0}$ be the
natural map which preserves fibers and
which maps $Z'_{i_0}$, ${Z'_{i_0}}^+$ and $\mathcal{Z}^{++}$ to
$Z_{i_0}$, $Z_{i_0}^+$ and $Z_{i_0}^{++}$ respectively.
(The restriction of $\Xi$ to each fiber is the map collapsing trivial floors.)
Let $a^2$ be the image of $(a^0, b^0, u^0, h^0)$ by the underlying map $V_0 \to
\widetilde{X}_{i_0}$.

Define $u^{i_0} \in C^l(\Sigma_{i_0} \setminus N_{i_0}, (\R_1 \sqcup \R_2 \sqcup \dots
\sqcup \R_{k_{i_0}}) \times Y)$ by
\[
u^{i_0} = (\theta_{(a^0, b^0, u^0, h^0)} \times 1) \circ u^0 \circ
(\Xi_{(a^0, b^0, u^0, h^0)})^{-1},
\]
where $\Xi_{(a^0, b^0, u^0, h^0)}$ is the restriction of $\Xi$ to the
fiber at $(a^0, b^0, u^0, h^0) \in V_0$.

Finally, we define the asymptotic parameters $b^{i_0}_\mu$.
We denote the index set of the joint circles of $\Sigma'_0$ between the $j$-th floor
and the $(j+1)$-th floor by $M^0_j$, and
the index set of the joint circles of $\Sigma_{i_0}$ between
the $j$-th floor and the $(j+1)$-th floor by $M^{i_0}_j$.
For each $\mu \in M^{i_0}_{j'}$,
let $S_{\mu_j}^1, S_{\mu_{j+1}}^1, \dots, S_{\mu_{j+m}}^1 \subset \Sigma'_{i_0}
\subset \Sigma'_0$ be the joint circles of $\Sigma'_0$ which collapse to $S_\mu^1$
by $\Sigma'_{i_0} \to \Sigma_{i_0}$, where we assume $\mu_{j+l} \in M^0_{j+l}$.

First we consider the case of $\rho^{i_0}_\mu \neq 0$.
Note that $\rho^0_{\mu_{j+l}} \neq 0$ in this case.
Since $b^{i_0}_\mu$ and $b^0_{\mu_{j+l}}$ should satisfy
\begin{align*}
- L_\mu \log \rho^{i_0}_\mu + b^{i_0}_\mu
&= \sigma \circ u^{i_0}(\widetilde{R}^{i_0}_{j'+1}(a^{i_0}))
- \sigma \circ u^{i_0}(\widetilde{R}^{i_0}_{j'}(a^{i_0}))\\
- L_\mu \log \rho^0_{\mu_{j+l}} + b^0_{\mu_{j+l}}
&= \sigma \circ u^0(\widetilde{R}^0_{j+l+1}(a^0))
- \sigma \circ u^0(\widetilde{R}^0_{j+l}(a^0)),
\end{align*}
we define $b^{i_0}_\mu$ by
\begin{align}
b^{i_0}_\mu &= (b^0_{\mu_j} + b^0_{\mu_{j+1}} + \dots + b^0_{\mu_{j+m}}) \notag\\
&\quad
+ \bigl(\sigma \circ u^0 \circ (\Xi_{(a^0, b^0, u^0, h^0)})^{-1}
(\widetilde{R}^{i_0}_{j'+1}(a^{i_0}))
- \sigma \circ u^0 (\widetilde{R}^0_{j+m+1}(a^0))\bigr) \notag\\
&\quad
- \bigl(\sigma \circ u^0 \circ (\Xi_{(a^0, b^0, u^0, h^0)})^{-1}(\widetilde{R}^{i_0}_{j'}(a^{i_0}))
- \sigma \circ u^0 (\widetilde{R}^0_{j+m}(a^0))\bigr) \notag\\
&\quad
+ L_\mu(-\log \rho^0_{\mu_j} - \dots - \log \rho^0_{\mu_{j+m}}
+ \log \rho^{i_0}_\mu).
\label{b^{i_0} for nonzero kappa}
\end{align}

Next we consider the case of $\rho^{i_0}_\mu = 0$.
Then there exist some $1 \leq c \leq d \leq m$ such that
$\rho^0_{\mu_{j+c}} = 0$, $\rho^0_{\mu_{j+d}} = 0$ and
$\rho^0_{\mu_{j+l}} \neq 0$ for $1 \leq l < c$ and $d < l \leq m$.
Then $b^0_{\mu_{j+l}}$ satisfies
\[
b^0_{\mu_{j+l}} = \bigl(\sigma \circ u^0(\widetilde{R}^0_{j+l+1}(a^0))
- \sigma \circ u^0(\widetilde{R}^0_{j+l}(a^0))\bigr) + L_\mu \log \rho^0_{\mu_{j+l}}
\]
for $1 \leq l < c$ and $d < l \leq m$, and
\begin{align*}
b^0_{\mu_{j+l}}
&= \lim_{s \to \infty}\bigl(\sigma \circ u^0|_{[0, \infty) \times S^1_{\mu_{j+l}}}(s, t)
- \sigma \circ u^0(\widetilde{R}^0_{j+l}(a^0)) - L_\mu s \bigr)\\
& \quad
- \lim_{s \to -\infty}\bigl(\sigma \circ u^0|_{(-\infty, 0] \times S^1_{\mu_{j+l}}}(s, t)
- \sigma \circ u^0(\widetilde{R}^0_{j+l+1}(a^0)) - L_\mu s \bigr)
\end{align*}
for $l = c, d$.
Hence
\begin{align}
&b^0_{\mu_j} + \dots + b^0_{\mu_{j+c}} \notag\\
&= 
\lim_{s \to \infty} \bigl(\sigma \circ u^0|_{[0, \infty) \times S^1_{\mu_{j+c}}}(s, t)
- \sigma \circ u^0(\widetilde{R}^0_j(a^0)) \notag\\
&\quad \hph{\lim_{s \to \infty} \bigl(}
- L_\mu \bigl(s - \log \rho^0_{\mu_j} - \dots
- \log \rho^0_{\mu_{j+c-1}}\bigr)\bigr) \notag\\
& \quad
- \lim_{s \to \infty} \bigl(\sigma \circ u^0|_{(-\infty, 0] \times S^1_{\mu_{j+c}}}(-s, t)
- (\sigma \circ u^0(\widetilde{R}^0_{j+c+1}(a^0)) - L_\mu s)\bigr)
\label{b of 0 to c}
\end{align}
and
\begin{align}
&b^0_{\mu_{j+d}} + \dots + b^0_{\mu_{j+m}} \notag\\
&=
- \lim_{s \to -\infty} \bigl(\sigma \circ u^0|_{(-\infty, 0] \times S^1_{\mu_{j+d}}}(s, t)
- \sigma \circ u^0(\widetilde{R}^0_{j+m+1}(a^0)) \notag\\
&\quad \hph{- \lim_{s \to -\infty} \bigl(}
- L_\mu \bigl(s + \log \rho^0_{\mu_{j+d+1}}
+ \dots + \log \rho^0_{\mu_{j+m}}\bigr) \bigr) \notag\\
& \quad
+ \lim_{s \to \infty} \bigl(\sigma \circ u^0|_{[0, \infty) \times S^1_{\mu_{j+d}}}(s, t)
-(\sigma \circ u^0(\widetilde{R}^0_{j+d}(a^0)) + L_\mu s)\bigr)
\label{b of d to m}
\end{align}

Assume that we use the decomposition of the trivial cylinder of
$(\Sigma'_0, z'_0 \cup {z'_0}^+, u_0)$ between $S^1_{\mu_{j + l}}$ and $S^1_{\mu_{j+l+1}}$
given by
\[
\overline{\R} \times S^1 = (-\infty, 1] \times S^1_{\mu_{j+l}}
\cup
[1, T_{j+l+1} -1] \times S^1 \cup [-1, \infty) \times S^1_{\mu_{j+l+1}}
\]
for the definition of the coordinate of $\widetilde{P}'_0$,
where we identify $\{1\} \times S^1$ and $\{T_{j+l+1} -1\} \times S^1$
with $\{1\} \times S^1_{\mu_{j+l}}$ and $\{-1\} \times S^1_{\mu_{j+l+1}}$ respectively,
and we consider the sections of the additional marked points ${Z'_{i_0}}^{++}$ as functions
to $[1, T_{j+l+1} -1] \times S^1$ instead of deforming the complex structure of
$[1, T_{j+l+1} -1] \times S^1$.
(Other cases can be covered by this case and the embeddings argued in the previous
section.)

First we assume $c < d$.
$u^0$ is trivial on the trivial cylinders between $S^1_{\mu_{j+c}}$ and
$S^1_{\mu_{j+d}}$, and the above assumption on the coordinate of $\widetilde{P}'_0$
implies that the natural coordinate of trivial cylinders and
the coordinates of $[0, \infty) \times S^1_{\mu_{j + l}}$ or
$(-\infty, 0] \times S^1_{\mu_{j + l}}$ coincide up to translation.
Therefore the following equations hold true.
\begin{multline}
\lim_{s \to -\infty} \bigl(\sigma \circ u^0|_{(-\infty, 0] \times S^1_{\mu_{j+c}}}(s, t)
- \sigma \circ u^0(\widetilde{R}^0_{j+c+1}(a^0)) - L_\mu s \bigr) \\
= - \bigl(\sigma \circ u^0(\widetilde{R}^0_{j+c+1}(a^0))
- \sigma \circ u^0|_{(-\infty, 0] \times S^1_{\mu_{j+c}}}(0, t)\bigr) \label{c trivial}
\end{multline}
\begin{multline}
\lim_{s \to \infty} \bigl(\sigma \circ u^0|_{[0, \infty) \times S^1_{\mu_{j+d}}}(s, t)
-\sigma \circ u^0(\widetilde{R}^0_{j+d}(a^0)) - L_\mu s \bigr) \\
= - \bigl(\sigma \circ u^0(\widetilde{R}^0_{j+d}(a^0))
- \sigma \circ u^0|_{[0, \infty) \times S^1_{\mu_{j+d}}}(0, t)\bigr). \label{d trivial}
\end{multline}
Similarly, for any $c < l < d$, whether $\rho^0_{\mu_{j+l}} = 0$ or not,
\begin{align}
b^0_{\mu_{j+l}}
&= \bigl(\sigma \circ u^0(\widetilde{R}^0_{j+l+1}(a^0))
- \sigma \circ u^0|_{(-\infty, 0] \times S^1_{\mu_{j+l}}}(0, t)\bigr) \notag\\
& \quad - \bigl(\sigma \circ u^0(\widetilde{R}^0_{j+l}(a^0))
- \sigma \circ u^0|_{[0, \infty) \times S^1_{\mu_{j+l}}}(0, t)\bigr).
\label{middle trivial}
\end{align}
Therefore equations (\ref{b of 0 to c}) to (\ref{middle trivial}) imply
\begin{align}
&b^0_{\mu_j} + \dots + b^0_{\mu_{j+m}} \notag\\
&= \lim_{s \to \infty} \bigl(\sigma \circ u^0|_{[0, \infty) \times S^1_{\mu_{j+c}}}(s, t)
- \sigma \circ u^0(\widetilde{R}^0_j(a^0)) \notag\\
& \hph{= \lim_{s \to \infty} \bigl(}
- L_\mu \bigl(s - \log \rho^0_{\mu_j} - \dots - \log \rho^0_{\mu_{j+c-1}}\bigr)
\bigr) \notag\\
& \quad
- \lim_{s \to -\infty}\bigl(\sigma \circ u^0|_{(-\infty, 0] \times S^1_{\mu_{j+d}}}(s, t)
- \sigma \circ u^0(\widetilde{R}^0_{j+m+1}(a^0)) \notag\\
& \hph{\quad - \lim_{s \to -\infty}\bigl(}
- L_\mu \bigl(s + \log \rho^0_{\mu_{j+d+1}} + \dots + \log \rho^0_{\mu_{j+m}}\bigr)
\bigr) \notag\\
& \quad
+ \sum_{c \leq l < d} \bigl(\sigma \circ u^0|_{[0, \infty) \times S^1_{\mu_{j+l+1}}}(0, t)
- \sigma \circ u^0|_{(-\infty, 0] \times S^1_{\mu_{j+l}}}(0, t)\bigr). \label{b sum}
\end{align}
It is easy to see that this equation also holds for the case of $c = d$.
The assumption on the coordinate of $\widetilde{P}'_0$ implies that the last terms
of (\ref{b sum}) are
\[
\sigma \circ u^0|_{[0, \infty) \times S^1_{\mu_{j+l+1}}}(0, t)
- \sigma \circ u^0|_{(-\infty, 0] \times S^1_{\mu_{j+l}}}(0, t)
= L_\mu T_{j+l+1}.
\]

Since $b^{i_0}_\mu$ is related to $u^{i_0}$ or $u^0$ by
\begin{align*}
b^{i_0}_\mu &= \lim_{s \to \infty} \bigl(\sigma \circ u^{i_0}|_{[0, \infty) \times S^1_\mu}(s, t)
- \sigma \circ u^{i_0}(\widetilde{R}^{i_0}_{j'}(a^{i_0})) - L_\mu s \bigr)\\
& \quad - \lim_{s \to -\infty} \bigl(\sigma \circ u^{i_0}|_{(-\infty, 0] \times S^1_\mu}(s, t)
- \sigma \circ u^{i_0}(\widetilde{R}^{i_0}_{j'+1}(a^{i_0})) - L_\mu s \bigr)\\
&= \lim_{s \to \infty} \bigl(\sigma \circ u^0 \circ
(\Xi_{(a^0, b^0, u^0, h^0)})^{-1}|_{[0, \infty) \times S^1_\mu}(s, t)\\
&\quad \hph{\lim_{s \to \infty} \bigl(}
- \sigma \circ u^0 \circ (\Xi_{(a^0, b^0, u^0, h^0)})^{-1}(\widetilde{R}^{i_0}_{j'}(a^{i_0}))
- L_\mu s \bigr)\\
& \quad - \lim_{s \to -\infty} \bigl(\sigma \circ u^0 \circ
(\Xi_{(a^0, b^0, u^0, h^0)})^{-1}|_{(-\infty, 0] \times S^1_\mu}(s, t)\\
&\quad \hph{- \lim_{s \to -\infty} \bigl(}
- \sigma \circ u^0 (\Xi_{(a^0, b^0, u^0, h^0)})^{-1}(\widetilde{R}^{i_0}_{j'+1}(a^{i_0}))
- L_\mu s \bigr)
\end{align*}
$b^{i_0}_\mu$ should satisfies
\begin{align}
b^{i_0}_\mu &= (b^0_{\mu_j} + \dots + b^0_{\mu_{j+m}}) \notag\\
& \quad
+ (\sigma \circ u^0 \circ (\Xi_{(a^0, b^0, u^0, h^0)})^{-1}(\widetilde{R}^{i_0}_{j'+1}(a^{i_0}))
- \sigma \circ u^0 (\widetilde{R}^0_{j+m+1}(a^0))) \notag\\
& \quad
- (\sigma \circ u^0 \circ (\Xi_{(a^0, b^0, u^0, h^0)})^{-1}(\widetilde{R}^{i_0}_{j'}(a^{i_0}))
- \sigma \circ u^0(\widetilde{R}^0_j(a^0))) \notag\\
& \quad
+ \lim_{s \to \infty}\bigl(\sigma \circ u^0 \circ
(\Xi_{(a^0, b^0, u^0, h^0)})^{-1}|_{[0, \infty) \times S^1_\mu}(s, t) \notag\\
& \quad \hph{+ \lim_{s \to \infty}\bigl(}
- \sigma \circ u^0|_{[0, \infty) \times S^1_{\mu_{j+c}}}
(s + \log \rho^0_{\mu_j} + \dots + \log \rho^0_{\mu_{j+c-1}}, t)\bigr) \notag\\
& \quad
- \lim_{s \to -\infty}\bigl(\sigma \circ u^0 \circ
(\Xi_{(a^0, b^0, u^0, h^0)})^{-1}|_{(-\infty, 0] \times S^1_\mu}(s, t) \notag\\
&\quad \hph{- \lim_{s \to -\infty}\bigl(}
- \sigma \circ u^0|_{(-\infty, 0] \times S^1_{\mu_{j+d}}}
(s - \log \rho^0_{\mu_{j+d+1}} - \dots - \log \rho^0_{\mu_{j+m}}, t)\bigr) \notag\\
& \quad
- L_\mu \sum_{c < l \leq d}T_{j+l} \notag\\
& =
(b^0_{\mu_j} + \dots + b^0_{\mu_{j+m}}) \notag\\
& \quad
+ (\sigma \circ u^0 \circ (\Xi_{(a^0, b^0, u^0, h^0)})^{-1}(\widetilde{R}^{i_0}_{j'+1}(a^{i_0}))
- \sigma \circ u^0 (\widetilde{R}^0_{j+m+1}(a^0))) \notag\\
& \quad
- (\sigma \circ u^0 \circ (\Xi_{(a^0, b^0, u^0, h^0)})^{-1}(\widetilde{R}^{i_0}_{j'}(a^{i_0}))
- \sigma \circ u^0(\widetilde{R}^0_j(a^0))) \notag\\
& \quad
+ \lim_{s \to \infty}
L_\mu(p_1 \circ (\Xi_{(a^0, b^0, u^0, h^0)})^{-1}|_{[0, \infty) \times S^1_\mu}(s, t)
\notag\\
& \quad \hph{+ \lim_{s \to \infty}L_\mu (}
- s - \log \rho^0_{\mu_j} - \dots - \log \rho^0_{\mu_{j+c-1}}) \notag\\
& \quad
- \lim_{s \to -\infty}
L_\mu(p_1 \circ (\Xi_{(a^0, b^0, u^0, h^0)})^{-1}|_{(-\infty, 0] \times S^1_\mu}(s, t)
\notag\\
&\quad \hph{- \lim_{s \to -\infty}L_\mu (}
- s + \log \rho^0_{\mu_{j+d+1}} + \dots + \log \rho^0_{\mu_{j+m}}) \notag\\
& \quad
- L_\mu \sum_{c < l \leq d}T_{j+l}
\label{b^{i_0} for zero kappa}
\end{align}
We define $b^{i_0}_\mu$ by the above formula (\ref{b^{i_0} for zero kappa}).

It is clear that $b^{i_0}_\mu$ is a smooth function of $(a^0, b^0, u^0, h^0)$
at $\rho^{i_0}_\mu \neq 0$.
We need to prove the smoothness at $\rho^{i_0}_\mu = 0$.
To prove the smoothness, we study the map $\Xi$.
As in the previous section,
we claim that there exists a smooth function $f_\mu : V_0 \to \C^\ast$ such that
\begin{align}
&(\rho^{i_0}_\mu)^{2\pi} e^{2\pi \sqrt{-1} \varphi^{i_0}_\mu)} \notag \\
&= \Bigl(\prod_{l=0}^m (\rho^0_{\mu_{j+l}})^{2\pi}\Bigr)
e^{2\pi (\sum_{l=0}^m
\sqrt{-1} \varphi^0_{\mu_{j+l}} - \sum_{l=1}^m T_{j+l})}
\cdot f_\mu(a^0, b^0, u^0, h^0). \label{bbf eq}
\end{align}
This can be proved as follows.

Let $(\widetilde{P}''_{i_0} \to \widetilde{X}''_{i_0}, Z''_{i_0})$ be the local universal family
of $(\Sigma'_{i_0}, z_{i_0} \cup z_{i_0}^+ \cup {z'_{i_0}}^{++} \cup z_{i_0}^{++})$.
Let ${z'''_{i_0}}^{++}$ be the points in ${z'_{i_0}}^{++} \cup z_{i_0}^{++}$ not contained in
the trivial floors of $(\Sigma'_{i_0}, z_{i_0}, u_0|_{\Sigma'_{i_0}})$,
and let $(\widetilde{P}'''_{i_0} \to \widetilde{X}'''_{i_0}, Z'''_{i_0})$ be the local universal
family of $(\Sigma_{i_0}, z_{i_0} \cup z_{i_0}^+ \cup {z'''_{i_0}}^{++})$.
Since the fiber of the center of $\widetilde{P}''_{i_0}$ is isomorphic to
$\widetilde{P}'_{i_0}$, $\widetilde{P}''_{i_0}$ is isomorphic to the product of
$\widetilde{P}'_{i_0}$ and a parameter space $D^M$ for additional marked points
corresponding to $z_{i_0}^{++}$.
Similarly, $\widetilde{P}'''_{i_0}$ is isomorphic to the product of $\widetilde{P}_{i_0}$
and a parameter space for additional marked points corresponding to
${z'''_{i_0}}^{++} \setminus z_{i_0}^{++}$.

By the assumption of the coordinate of $\widetilde{P}'_0$,
it is easy to see that if $(\rho'''_\mu, \varphi'''_\mu)$ is an appropriately chosen
parameter of $\widetilde{P}'''_{i_0}$ for the deformation of a neighborhood of
the joint circle $S^1_\mu$, then the following holds true under the natural map
$\widetilde{P}''_{i_0} \to \widetilde{P}'''_{i_0}$,
where we use the coordinate of $\widetilde{P}''_{i_0}$ given by
$\widetilde{P}''_{i_0} \cong \widetilde{P}'_{i_0} \times D^M$.
\[
(\rho'''_\mu)^{2\pi} e^{2\pi \sqrt{-1} \varphi'''_{\mu}}
= \Bigl(\prod_{l=0}^m (\rho^0_{\mu_{j+l}})^{2\pi}\Bigr)
e^{2\pi (\sum_{l=0}^m
\sqrt{-1} \varphi^0_{\mu_{j+l}} - \sum_{l=1}^m T_{j+l})}
\]

Since $\widetilde{P}'''_{i_0}$ is isomorphic to the product of $\widetilde{P}_{i_0}$ and
some parameter space, there exists a smooth map $f' : \widetilde{X}'''_{i_0} \to \C^\ast$
such that
\[
(\rho^{i_0}_\mu)^{2\pi} e^{2\pi \sqrt{-1} \varphi^{i_0}_{\mu}}
= (\rho'''_\mu)^{2\pi} e^{2\pi \sqrt{-1} \varphi'''_{\mu}} \cdot f'(a''')
\]
Therefore, there exists a smooth map $f : V_0 \to \C^\ast$ which satisfies
equation (\ref{bbf eq}).

Similarly, there exist smooth maps $f_\mu^{\text{left}}, f_\mu^{\text{right}}
: V_0 \times_{\widetilde{X_{i_0}}} \widetilde{P}_{i_0} \to \C^\ast$ such that
if $\rho^0_{\mu_{j+l}} \neq 0$ for $1 \leq l < c$ and $d < l \leq m$,
and $\rho^0_{\mu_{j+c}} = 0$ and $\rho^0_{\mu_{j+d}} = 0$, and
\begin{align*}
\Xi_{(a^0, b^0, u^0, h^0)}|_{[0, \infty) \times S^1_{\mu_{j+c}}}
(s_0^{\text{left}}, t_0^{\text{left}}) &= (s_{i_0}^{\text{left}}, t_{i_0}^{\text{left}}),\\
\Xi_{(a^0, b^0, u^0, h^0)}|_{(-\infty, 0] \times S^1_{\mu_{j+d}}}
(s_0^{\text{right}}, t_0^{\text{right}}) &= (s_{i_0}^{\text{right}}, t_{i_0}^{\text{right}}),
\end{align*}
then
\begin{align*}
e^{2\pi (s_{i_0}^{\text{left}} + \sqrt{-1} t_{i_0}^{\text{left}})}
&= e^{2\pi (s_0^{\text{left}} + \sqrt{-1} t_0^{\text{left}})}
\Bigl(\prod_{0 \leq l < c} (\rho^0_{\mu_{j+l}})^{2\pi} \Bigr)
e^{2\pi \sqrt{-1} \sum_{0 \leq l < c} \varphi^0_{\mu_{j+l}}}\\
&\quad \cdot e^{-2\pi \sum_{0 < l \leq c} T_{j+l}}
f_\mu^{\text{left}} (s_{i_0}^{\text{left}}, t_{i_0}^{\text{left}}, a^0, b^0, u^0, h^0),
\end{align*}
and
\begin{align*}
e^{2\pi (s_{i_0}^{\text{right}} + \sqrt{-1} t_{i_0}^{\text{right}})}
&= e^{2\pi (s_0^{\text{right}} + \sqrt{-1} t_0^{\text{right}})}
\Bigl(\prod_{d < l \leq m} (\rho^0_{\mu_{j+l}})^{2\pi} \Bigr)
e^{2\pi \sqrt{-1}\sum_{d < l \leq m}
 \varphi^0_{\mu_{j+l}}}\\
& \quad \cdot e^{-2\pi \sum_{d < l \leq m} T_{j+l}}
f_\mu^{\text{right}} (s_{i_0}^{\text{right}}, t_{i_0}^{\text{right}}, a^0, b^0, u^0, h^0).
\end{align*}
Furthermore, $f_\mu$, $f_\mu^{\text{left}}$ and $f_\mu^{\text{right}}$ satisfy
\begin{align*}
f_\mu(a^0, b^0, u^0, h^0) &= \lim_{s_{i_0}^{\text{left}} \to \infty}
f_\mu^{\text{left}}(s_{i_0}^{\text{left}}, t_{i_0}^{\text{left}}, a^0, b^0, u^0, h^0)\\
& \quad \cdot \lim_{s_{i_0}^{\text{right}} \to -\infty}
f_\mu^{\text{right}}(s_{i_0}^{\text{right}}, t_{i_0}^{\text{right}}, a^0, b^0, u^0, h^0).
\end{align*}

Therefore, $b^{i_0}_\mu$ satisfies
\begin{align}
b^{i_0}_\mu
& =
(b^0_{\mu_j} + \dots + b^0_{\mu_{j+m}}) \notag\\
& \quad
+ (\sigma \circ u^0 \circ (\Xi_{(a^0, b^0, u^0, h^0)})^{-1}(\widetilde{R}^{i_0}_{j'+1}(a^{i_0}))
- \sigma \circ u^0 (\widetilde{R}^0_{j+m+1}(a^0))) \notag\\
& \quad
- (\sigma \circ u^0 \circ (\Xi_{(a^0, b^0, u^0, h^0)})^{-1}(\widetilde{R}^{i_0}_{j'}(a^{i_0}))
- \sigma \circ u^0(\widetilde{R}^0_j(a^0))) \notag\\
& \quad
+ L_\mu \cdot \frac{1}{2\pi} \log|f_\mu(a^0, b^0, u^0, h^0)|
- L_\mu \sum_{0 < l \leq m} T_{j+l} \label{b^{i_0} eq}
\end{align}
in both cases.
Hence $b^{i_0}_\mu$ is a smooth function of $(a^0, b^0, u^0, h^0)$.

It is easy to check that the constructed map
$(a^0, b^0, u^0, h^0) \to (a^{i_0}, b^{i_0}, u^{i_0}, h^{i_0})$ is the required essential
submersion.
For example, the essential submersiveness of $\phi : V_0 \to \prod_i V_i$ is seen
as follows.
By the coordinate change (i.e. changing the center $p'$ of the Kuranishi neighborhood),
it is enough to prove that
$\phi|_{\phi^{-1}(\prod_i V_i^{(k_i)})} : \phi^{-1}(\prod_i V_i^{(k_i)}) \to \prod_i V_i^{(k_i)}$
is submersive.
($V_i^{(k_i)} \subset V_i$ are the corners of the highest codimension.)
It is clear that $\phi^{-1}(\prod_i V_i^{(k_i)})$ is a union of corners of $V_0$,
and each of them is defined by
$\{\rho_\mu = 0; \text{ for all } \mu \in \bigcup_{j \in I} M^0_j\}$ for some
$I \subset \{1, 2, \dots, k_0 - 1\}$.
Then $b^0_\mu \in \R$ ($\mu \in \bigcup_{j \in I} M_i$) are independent parameters
in $\mathring{X}$ since $-L_\mu \log \rho^0_\mu + b^0_\mu = \infty$ for all $b^0_\mu$.
For each $\mu \in M^{i_0}_{j'}$, let $S^1_{\mu_j}, S^1_{\mu_{j+1}}, \dots, S^1_{\mu_m}
\subset \Sigma'_{i_0} \subset \Sigma'_0$ be the joint circles which collapse to
$S^1_\mu$ by $\Sigma'_{i_0} \to \Sigma_{i_0}$ as above.
Assume $\mu_j \in M^0_j$.
Then there exists some $j+l \in \{j, j+1, \dots, j+m\}$ such that $j+l \in I$.
Since the derivative of $b^{i_0}_\mu$ by $b^0_{\mu_{j+l}}$ does not vanish by
(\ref{b^{i_0} eq}), it is easy to check
that $\phi|_{\phi^{-1}(\prod_i V_i^{(k_i)})}$ is submersive.

Since $E_0$ is the direct sum of $E_i$ ($1 \leq i \leq N$),
grouped multisections of $(V_i, E_i)$ define a grouped multisection of $(V_0, E_0)$
by the pull back of the product grouped multisection by the essential submersion.
(We assume that the grouped multisection of $(V_i, E_i)$ and $(V_j, E_j)$ coincide
if $p_i = p_j$.)
We note that $\dim V_0 = \dim \prod_{i=1}^N V_i - (N - 1) > \dim \prod_{i=1}^N V_i$
if $N > 1$.


%% file: SFT-05_Global_Kuranishi_structure.tex
%
%
\subsection{Construction of global structure}\label{global construction}
In this section, we construct a global pre-Kuranishi structure of $\widehat{\M}$.
As we explained in Section \ref{embed},
a Kuranishi neighborhood of each point $p \in \widehat{\M}$ is determined
by the data $(z^+, S, E^0, \lambda)$.
Hence construction of a pre-Kuranishi structure of $\widehat{\M}$
is equivalent to constructing a Hausdorff space $\X$ with a locally homeomorphic
surjection $\mu : \X \to \widehat{\M}$ and giving such data for each point of $\X$.
In Section \ref{embed}, $S$ is a codimension-two submanifold of $Y$ and
$z^+$ is a finite subset of the domain curve,
but for the construction of global structure,
it is convenient to use a finite set $\mathcal{S} = \{S\}$ of codimension-two
submanifolds of $Y$ and a finite family $z^+ = (z^S)_{S \in \mathcal{S}}$ of finite subsets
of the domain curve indexed by $\mathcal{S}$ instead.
(We assume that $\pi_Y \circ u$ intersects each $S$ at $z^S$ transversely.)

First we introduce three versions of the space of holomorphic buildings
$\widehat{\M}_{\mathcal{S}}$, $\widehat{\M}_{\mathcal{S}, A}$ and
$\widehat{\M}_{o, \mathcal{S}, A}$.
We will realize the Hausdorff space $\X$ as a subspace of
$\widehat{\M}_{o, \mathcal{S}, A}$.
Let $\mathcal{S} = \{S\}$ be a finite set of codimension-two submanifolds of $Y$.
A point $(\Sigma, z, z^S, z^A, z^o, u)$ of $\widehat{\M}_{o, \mathcal{S}, A}$ consists of
a holomorphic building $(\Sigma, z, u) \in \widehat{\M} = \widehat{\M}(Y, \lambda, J)$,
finite subsets $z^S \subset \Sigma$ ($S \in \mathcal{S}$),
a finite subset $z^A \subset \Sigma$ and
a finite subset $z^o \subset \Sigma$ which satisfy the following conditions:
\begin{itemize}
\item
$\pi_Y \circ u$ intersects $S$ at $z^S$ transversely for each $S \in \mathcal{S}$.
\item
$z^S, z^A, z^o \subset \Sigma$ are disjoint, do not contain any special points
of $(\Sigma, z, u)$ and any points of the imaginary circles of $\Sigma$
and the trivial cylinders of $(\Sigma, z, u)$.
\item
All non-trivial components (i.e. irreducible components other than trivial cylinders) of
$(\Sigma, z, u)$ are stable in $(\Sigma, z, z^S)$.
\end{itemize}
$z^S$ are used to make the domain curve stable,
$z^A$ is used to control the automorphism group of the domain curve, and
$z^o$ is a mark which tells us the additional vector space $E^0$ we used
for the construction of the Kuranishi neighborhood.
In other words, $z^o$ is used to realize the space $\X$ as a subspace of
$\widehat{\M}_{o, \mathcal{S}, A}$.
Two points $(\Sigma, z, z^S, z^A, z^o, u)$ and $(\Sigma', z', (z')^S, (z')^A, (z')^o, u')$ are
the same point if there exists a biholomorphism $\varphi : \Sigma \to \Sigma'$ and
an $\R$-translation $\theta$ such that $\varphi(z) = z'$, $\varphi(z^S) = (z')^S$
for all $S \in \mathcal{S}$, $\varphi(z^A) = (z')^A$,
$\varphi(z^o) = (z')^o$ and $u' \circ \varphi = (\theta \times 1) \circ u$.
The topology of $\widehat{\M}_{o, \mathcal{S}, A}$ is defined as a quotient space of
a subspace of $\overline{\M}(Y, \lambda, J)$
(locally, it is the quotient by the $S^1$-actions on the coordinates of limit circles
and the symmetric group of the sets $z$, $z^S$ ($S \in \mathcal{S}$), $z^A$, $z^o$ and
the set of limit circles).
$\widehat{\M}_{\mathcal{S}, A}$ consists of points $(\Sigma, z, z^S, z^A, u)$,
and $\widehat{\M}_{\mathcal{S}}$ consists of points $(\Sigma, z, z^S, u)$.
We may regard them as the subspaces of $\widehat{\M}_{o, \mathcal{S}, A}$ defined by
$z^o = \emptyset$ and $(z^o, z^S) = (\emptyset, \emptyset)$ respectively.
If $\mathcal{S}' \supset \mathcal{S}$, we regard $\widehat{\M}_{o, \mathcal{S}, A}$
as a subspace of $\widehat{\M}_{o, \mathcal{S}', A}$.
The forgetful map $\forget_{\mathcal{S}, A} : \widehat{\M}_{\mathcal{S}, A} \to
\widehat{\M}$ is defined by forgetting the points $z^S$ and $z^A$.
Similarly, we define
$\forget_{o, \mathcal{S}, A} : \widehat{\M}_{o, \mathcal{S}, A} \to \widehat{\M}$.

For two holomorphic buildings
$p^k = (\Sigma^k, z^k, (z^k)^S, (z^k)^A, (z^k)^o, u^k) \in \widehat{\M}_{o, \mathcal{S}, A}$
($k=1,2$), we say $p^1 \leq p^2$ if $(z^2)^S$, $(z^2)^A$ and $(z^2)^o$
are $\Aut(\Sigma^2, z^2, u^2)$-invariant and
there exists a biholomorphism $\Sigma^1 \cong \Sigma^2$ such that
under this biholomorphism,
$p_1$ is obtained from $p^2$ by forgetting some
subsets of $(z^2)^S$, $(z^2)^A$ and $(z^2)^o$.
(The forgetful map from $p^2$ to $p^1$ does not collapse any components.)

We also define two versions of the space of stable curves
$\overline{\M}^\DM_{\mathcal{S}}$ and $\overline{\M}^\DM_{\mathcal{S}, A}$ as follows.
A point $(\hat \Sigma, z, z^S, z^A)$ of $\overline{\M}^\DM_{\mathcal{S}, A}$
consists of a semistable curve $\hat \Sigma$ and finite disjoint subsets $z$,
$z^S$ ($S \in \mathcal{S}$) and $z^A$ such that
they do not contain any nodal points and the automorphism group of
$(\hat \Sigma, z, z^S)$ is finite.
Similarly, $\overline{\M}^\DM_{\mathcal{S}}$ consists of points
$(\hat \Sigma, z, z^S)$ which satisfy the same conditions.
There is another forgetful map $\forget_u : \widehat{\M}_{\mathcal{S}, A} \to
\overline{\M}^\DM_{\mathcal{S}, A}$ 
defined by forgetting the map $u$, blowing down
joint circles to nodal points, blowing down limit circles and add these points to
marked points $z$, and stabilizing (collapsing all components corresponding to trivial
cylinders).
For example, $\forget_u$ maps $(\Sigma, z, z^S, z^A, u)$ in Figure \ref{(Sigma,z,zS,zA)} to
$(\hat \Sigma, z, z^S, z^A)$ in Figure \ref{(hatSigma,z,zS,zA)}.
\begin{figure}
\centering
\includegraphics[width= 350pt]{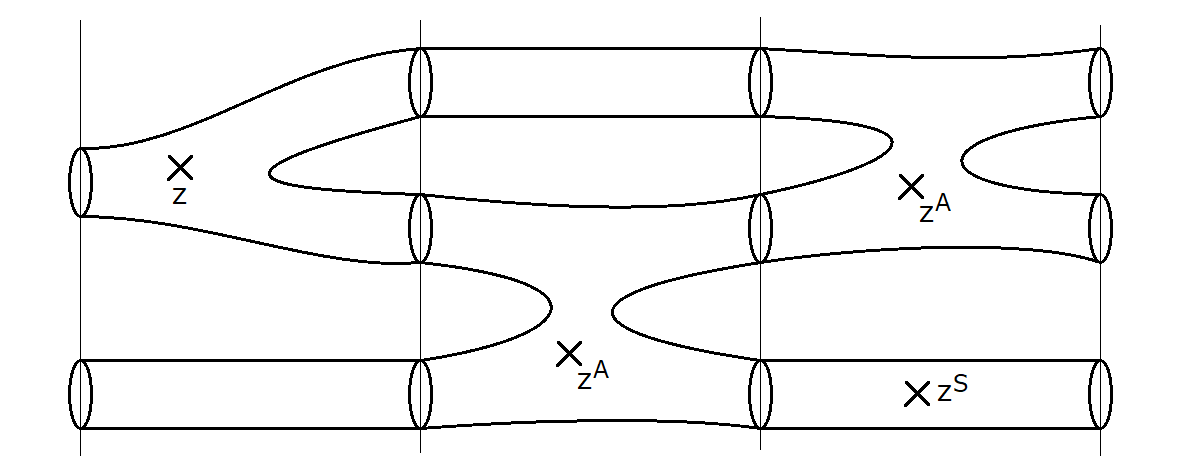}
\caption{$p = (\Sigma, z, z^S, z^A, u)$}\label{(Sigma,z,zS,zA)}
\includegraphics[width= 350pt]{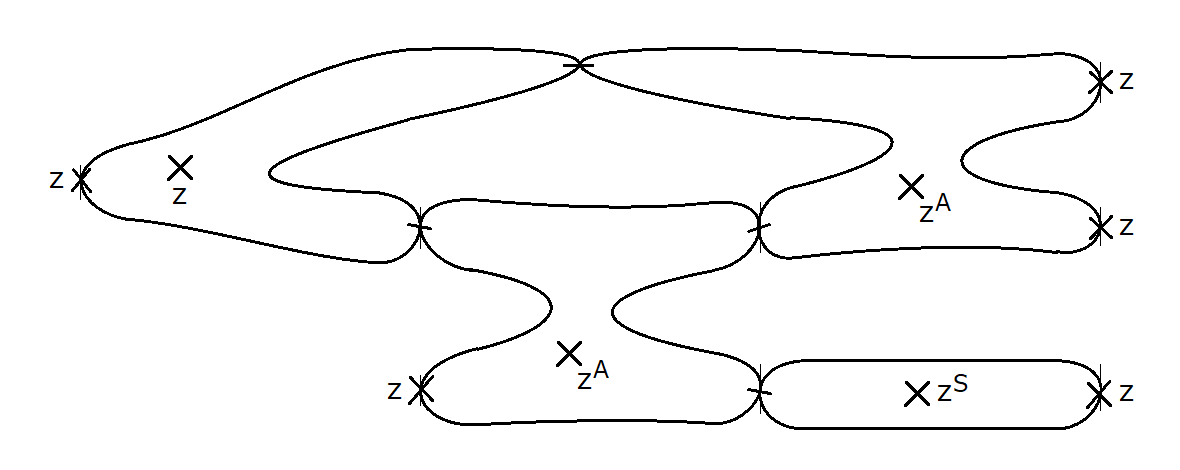}
\caption{$(\hat \Sigma, z, z^S, z^A) = \forget_u(p)$}\label{(hatSigma,z,zS,zA)}
\end{figure}

$(\mathcal{S}, A)$-forgetful map $f$ from
$p \in \overline{\M}^\DM_{\mathcal{S}, A}$
to $q \in \overline{\M}^\DM_{\mathcal{S}, A}$ is a forgetful map
$f : p \to q$ obtained by forgetting some subsets of $z^S$ and $z^A$, and
stabilizing the curve.

Fix an arbitrary large constant $L_{\max} > 0$, and
let $\delta_0 > 0$ be a positive constant such that
$4 \delta_0$ is less than the minimal period of periodic orbits and
$2 \delta_0$ is less than the minimal difference $L^+ - L^-> 0$ of two periods
$L^-, L^+$ of periodic orbits such that $L^- < L^+ \leq L_{\max}$.
Let $\widehat{\M}^{\leq L_{\max}} \subset \widehat{\M}$ be the subspace of
holomorphic buildings the sums of the periods of whose $+\infty$-limit circles
are $\leq L_{\max}$.

We say a holomorphic building $(\Sigma, z, u) \in \widehat{\M}$ is of type
$\theta = (g, k, E_{\hat \omega})$ if the genus of $\Sigma$ is $g$,
the total number of marked points and limit circles is $k$, and
$E_{\hat \omega}(u) = E_{\hat \omega}$.
For each triple $\theta =  (g, k, E_{\hat \omega})$,
we define $\widetilde{e}(\theta) = \widetilde{e}_{\delta_0}(\theta)$ by
$\widetilde{e}_{\delta_0}(\theta) = 5(g-1) + 2k + E_{\hat \omega} / \delta_0$.
For a holomorphic building $p \in \widehat{\M}$ of type $\theta$, we
define $\widetilde{e}(p) = \widetilde{e}(\theta)$.
Then $\widetilde{e}$ satisfies the following:
\begin{itemize}
\item
$\widetilde{e}(p) \geq 1$
for any holomorphic building $p \in \widehat{\M}^{\leq L_{\max}}$.
\item
For any holomorphic building $p \in \widehat{\M}^{\leq L_{\max}}$,
any subset $C_0$ of its nodal points and any subset $C_1$ of the gaps of floors,
replace the nodal points in $C_0$ and the joint circles in the gaps in $C_1$
of $p$ to pairs of marked points and pairs of
limit circles respectively, and let $p'_i$ ($1 \leq i \leq N$) be their non-trivial connected
components. (They are connected holomorphic buildings of height-one without
nodal points.)
Then $\widetilde{e}(p) \geq \sum_i \widetilde{e}(p'_i)$, and
the inequality is strict if $C_0 \neq \emptyset$ or $C_1 \neq \emptyset$.
In particular, $\widetilde{e}(p) > \widetilde{e}(p'_i)$ for all $i$ if $N > 1$.
\end{itemize}
The second property is easy to check.
(We recall that $\widehat{\M}$ consists of holomorphic buildings without trivial buildings,
where a trivial building is a connected component which consists of trivial cylinders only.)
We check the first property.
By the second property, it is enough to check the property for connected height-one
holomorphic buildings.
If the domain curve is stable, then $\widetilde{e} \geq 1$ is clear.
If the domain curve is unstable, then since $E_{\hat \omega} > 0$,
$k$ must be $\geq 1$.
Hence $(g, k)  = (0,1)$ or $(0,2)$.
If $(g, k) = (0,1)$, then $E_{\hat \omega}$ is greater than or equal to the minimal
period of periodic orbit.
If $(g, k) = (0,2)$, then $E_{\hat \omega}$ is greater than or equal to the minimal
difference of two periods of periodic orbits.
Hence in both cases, $\widetilde{e} \geq 1$ by the definition of $\delta_0 > 0$.

For each triple $\theta = (g, k, E_{\hat \omega})$, let
$\widehat{\M}^{\leq L_{\max}}_\theta \subset \widehat{\M}$ be the subspace of
holomorphic buildings of type $\theta$ such that
the sums of the periods of $+\infty$-limit circles are $\leq L_{\max}$.
We also define $\widehat{\M}^{\leq L_{\max}}_{\mathcal{S}, A, \theta}
= \forget_{\mathcal{S}, A}^{-1}(\widehat{\M}^{\leq L_{\max}}_\theta)$.
Note that for any constants $C \geq 0$,
$\widehat{\M}^{\leq L_{\max}}_{\leq C}
= \bigcup_{\widetilde{e}(\theta) \leq C} \widehat{\M}^{\leq L_{\max}}_\theta$
is compact.

For each triple $\theta = (g, k, E_{\hat \omega})$,
let $\overline{\M}^\DM_{\mathcal{S}, A, \theta} \subset
\overline{\M}^\DM_{\mathcal{S}, A}$ be the subspace of stable curves
whose genus and the number of marked points are $g$ and $k$ respectively.
We regard the spaces for different $E_{\hat \omega}$ as disjoint spaces,
and regard the forgetful map $\forget_u$ as a map from
$\widehat{\M}^{\leq L_{\max}}_{\mathcal{S}, A, \theta}$ to
$\overline{\M}^\DM_{\mathcal{S}, A, \theta}$.
For a point $\hat p \in \overline{\M}^\DM_{\mathcal{S}, A, \theta}$,
we define $\theta(\hat p) = \theta$.
For each point $\hat p = (\hat \Sigma, z, z^S, z^A) \in
\overline{\M}^\DM_{\mathcal{S}, A, \theta}$, we define an integer $l(\hat p)$ by
$l(\hat p) = 3(g - 1) + \# z + \sum_{S \in \mathcal{S}} \# z^S + \# z^A$,
where $g$ is the genus of $\hat \Sigma$.
For a point $p \in \widehat{\M}_{\mathcal{S}, A}$, we define $l(p) = l(\forget_u(p))$.
For each $l \geq 0$, let
$\overline{\M}^\DM_{\mathcal{S}, A, \theta, l} \subset
\overline{\M}^\DM_{\mathcal{S}, A, \theta}$ be the subspace of curves $\hat p$
such that $l(\hat p) = l$, and $\widehat{\M}^{\leq L_{\max}}_{\mathcal{S}, A, \theta, l}
\subset \widehat{\M}^{\leq L_{\max}}_{\mathcal{S}, A, \theta}$ be the subspace of
holomorphic buildings $p$ such that $l(p) = l$.

Let $(\hat p, E^0, \lambda)$ be a triple of a stable curve
$\hat p = (\hat \Sigma, z, z^S, z^A, z^o) \in \overline{\M}^\DM_{o, \mathcal{S}, A}$,
a finite-dimensional $\Aut(\hat p)$-vector space $E^0$, and
an $\Aut(\hat p)$-equivariant linear map
\[
\lambda : E^0 \to C^\infty(\hat \Sigma \times Y, \Wedge^{0, 1} T^\ast \hat \Sigma
\otimes_\C (\R \partial_\sigma \oplus T Y)).
\]
We call such a triple $(\hat p, E^0, \lambda)$ a stable curve with perturbation
parameters.
We say two stable curves with perturbation parameters $(\hat p_k, E^0_k, \lambda_k)$
($k =1,2$) are isomorphic if there exist an isomorphism $f : \hat p_1 \to \hat p_2$
and an isomorphism $\hat \phi_f : E^0_1 \to E^0_2$ which is $\Aut(\hat p_1)$-equivariant
with respect to the isomorphism $\rho_f : \Aut(\hat p_1) \to \Aut(\hat p_2)$
associated to $f$, and they make the following diagram commutative.
\[
\begin{tikzcd}
E^0_1 \ar{r}{\lambda_1} \ar{d}{\hat \phi_f} &
C^\infty(\hat \Sigma_1 \times Y, \Wedge^{0, 1} T^\ast \hat \Sigma_1
\otimes_\C (\R \partial_\sigma \oplus T Y)) \\
E^0_2 \ar{r}{\lambda_2} &
C^\infty(\hat \Sigma_2 \times Y, \Wedge^{0, 1} T^\ast \hat \Sigma_2
\otimes_\C (\R \partial_\sigma \oplus T Y))
\ar{u}{f^\ast}
\end{tikzcd}
\]

A holomorphic building with perturbation parameters $(p, E^0, \lambda)$ is a triple
of a holomorphic building $p = (\Sigma, z, z^S, z^A, z^o, u) \in
\widehat{\M}_{o, \mathcal{S}, A}$, a finite-dimensional $\Aut(p)$-vector space $E^0$,
and an $\Aut(p)$-equivariant linear map
\[
E^0 \to C^\infty(\hat \Sigma \times Y, \Wedge^{0, 1} T^\ast \hat \Sigma
\otimes_\C (\R \partial_\sigma \oplus T Y)),
\]
where $(\hat \Sigma, z, z^S, z^A, z^o) = \forget_u(\Sigma, z, z^S, z^A, z^o, u)$
is the stabilization of the blow down of the domain curve.
We say two such triples $(p_k, E^0_k, \lambda_k)$ ($k=1,2$) are isomorphic if
there exist an isomorphism $f : p_1 \to p_2$ and
an $\Aut(p_1)$-equivariant isomorphism $\hat \phi_f : E^0_1 \to E^0_2$ such that
$\lambda_1 = f^\ast \circ \lambda_2 \circ \hat \phi_f$.

For an arbitrary constant $C \geq 0$,
we will construct a space $\X = \X_{\leq C}
= \bigcup_{\widetilde{e}(\theta) \leq C} \X_{\theta}$ consisting of
holomorphic buildings with perturbation parameters
which satisfies the following conditions:
\begin{enumerate}
\item
For each $(p, E^0_p, \lambda_p) \in \X_{\theta}$,
$p$ is contained in $\widehat{\M}^{\leq L_{\max}}_{o, \mathcal{S}, A, \theta}$.
\item
For any $p = (\Sigma, z, u) \in \widehat{\M}^{\leq L_{\max}}_{\theta}$,
there exists $(p^+, E^0, \lambda) \in \X_{\theta}$
such that $p = \forget_{o, \mathcal{S}, A}(p^+)$.
Furthermore, $\forget_{o, \mathcal{S}, A} : \X_{\theta} \to
\widehat{\M}^{\leq L_{\max}}_{\theta}$ is locally homeomorphic.
\item
For any $(p = (\Sigma, z, z^S, z^A, z^o, u), E^0_p, \lambda_p) \in \X$, let
\[
E^0_p \to C^\infty(\Sigma \times Y, \Wedge^{0,1} T^\ast \Sigma \otimes_\C
(\R \partial_\sigma \oplus T Y))
\]
be the pull back of $\lambda_p$ by the forgetful map $p \to \forget_u(p)$ and
also denote it by the same symbol $\lambda_p$.
Then the linear map
\begin{align}
&\widetilde{W}^{1, q}_\delta(\Sigma, u^\ast T\hat Y) \oplus E^0_p \notag\\
&\to L^p_\delta(\Sigma, \Wedge^{0, 1} T^\ast \Sigma \otimes_\C u^\ast T \hat Y)
\oplus \bigoplus_{\text{limit circles}} \Ker A_{\gamma_{\pm\infty_i}}
/ (\R \partial_\sigma \oplus \R R_\lambda) \notag\\
&\quad \oplus \bigoplus_{z_i \in z} T_{\pi_Y \circ u(z_i)} Y \notag\\
&(\xi, h)
\mapsto (D_p \xi + \lambda_p(h)(\cdot, \pi_Y \circ u(\cdot)), \notag \\
&\hph{(\xi, h) \mapsto (}
\sum_j \langle \xi|_{S^1_{\pm\infty_i}}, \eta^{\pm\infty_i}_j \rangle \eta^{\pm\infty_i}_j,
\pi_Y \circ \xi(z_i))
\label{X surjective map}
\end{align}
is surjective, where $D_p$ is the linearization of the equation of $J$-holomorphic maps,
and $\{\eta^{\pm\infty_i}\}_j$ is an orthonormal basis of the orthogonal complement of
$\R \partial_\sigma \oplus \R R_\lambda$ in $\Ker A_{\gamma_{\pm\infty_i}}$
for each $\pm\infty_i$.
\item
\label{X embedding relation}
If two points $(p_k^+ = (\Sigma, z, z^{S, k}, z^{A, k}, z^{o, k}, u), E^0_k, \lambda_k) \in \X$
($k=1,2$) over the same holomorphic building $p = (\Sigma, z, u) \in \widehat{\M}$ satisfy
$z^{S, 2} \supset z^{S, 1}$, $z^{A, 2} \supset z^{A, 1}$ and $z^{o, 2} \supset z^{o, 1}$
(that is, if $p_1^+ \leq p^+_2$),
then $E^0_2$ contains $E^0_1$ as a subspace, and the restriction of $\lambda_2$
coincides with $\lambda_1$.
\item
For any $p = (\Sigma, z, z^S, z^A, z^o, u) \in \X$, $z^S$, $z^A$ and $z^o$ are
$\Aut(\forget_{o, \mathcal{S}, A}(p))$-invariant.
\item
$\X$ is embedded in $\widehat{\M}_{o, \mathcal{S}, A}$.
In fact, we add marked point $z^o$ to distinguish $E^0$ and $\lambda$.
\item
\label{X vee existence}
For any two points
$(p_k^+ = (\Sigma, z, z^{S, k}, z^{A, k}, z^{o, k}, u), E^0_k, \lambda_k) \in \X$
($k=1,2$) for the same holomorphic building $p = (\Sigma, z, u) \in \widehat{\M}$,
there exists some
$(p_3^+ = (\Sigma, z, z^{S, 3}, z^{A, 3}, z^{o, 3}, u), E^0_3, \lambda_3)
\in \X$ such that $z^{S, 3} = z^{S, 1} \cup z^{S, 2}$, $z^{A, 3} = z^{A, 1} \cup z^{A, 2}$
and $z^{o, 3} = z^{o, 1} \cup z^{o, 2}$.
(In the definition of pre-Kuranishi structure, $p_3^+$ will be the unique supremum
$p_1^+ \vee p_2^+$.)
\item
\label{X decomposition into parts}
$\X$ satisfies the following compatibility condition with respect to the decomposition of
a holomorphic building into parts:
For any point $p = (\Sigma, z, z^S, z^A, z^o, u) \in \widehat{\M}_{o, \mathcal{S}, A}$,
replace all nodal points and all joint circles with pairs of marked points and
pairs of limit circles respectively,
and let $p'_i \in \widehat{\M}_{o, \mathcal{S}, A}$ ($i = 1, \dots, k$) be its
connected components other than trivial cylinders.
(Each $p'_i$ is a connected height-one holomorphic building without nodal points.)
Then $(p, E^0, \lambda) \in \X$ for some $E^0$ and $\lambda$
if $z^S$, $z^A$ and $z^o$ are $\Aut((\forget_{o, \mathcal{S}, A}(p'_i))_i)$-invariant
and $(p'_i, E^0_i, \lambda_i) \in \X$ for some $E^0_i$ and $\lambda_i$ for all $i$.
($\Aut((\forget_{o, \mathcal{S}, A}(p'_i))_i)$ is the automorphism group
of $\coprod_i \forget_{o, \mathcal{S}, A}(p'_i)$.)
Furthermore, $E^0$ is isomorphic to the direct sum of $E^0_i$ and
the restriction of $\lambda$ to $E^0$ coincides with the pull back of $E^0_i$ by
the forgetful map.
\end{enumerate}
We define a pre-Kuranishi structure of $\widehat{\M}^{\leq L_{\max}}_{\theta}$ by
$\X_{\theta} \subset
\widehat{\M}^{\leq L_{\max}}_{o, \mathcal{S}, A, \theta}$ and
a locally homeomorphic surjection
$\forget_{o, \mathcal{S}, A} : \X_{\theta} \to \widehat{\M}^{\leq L_{\max}}_{\theta}$.
To define a Kuranishi neighborhood of $\forget_{o, \mathcal{S}, A}(p)$ for
each $(p, E^0, \lambda) \in \X$, we need to extend
$\lambda$ to a local universal family of $\forget_u(p)$.
Hence we also construct a space of stable curves with perturbation parameters
which gives a neighborhood of the domain curves of holomorphic buildings in $\X$
in a sense. (See Lemma \ref{good family of additional vector spaces} for details.)
Condition \ref{X embedding relation} will imply that for any two $p^+_k \in \X$ ($k = 1,2$)
for the same holomorphic building $p \in \widehat{\M}$, if $p^+_1 \leq p^+_2$, then
we can define the embedding of the Kuranishi neighborhood of $p$ defined by
the data associated to $p^+_1$ to that defined by the data associated to $p^+_2$.
Furthermore, Condition \ref{X vee existence} imply the existence of
the unique supremum of any two points in the same fiber of
$\forget_{o, \mathcal{S}, A} : \X_{\theta} \to \widehat{\M}^{\leq L_{\max}}_{\theta}$.

If we ignore the algebraic structure of $\widehat{\M}$ such as the fiber product
structure, then we do not need Condition \ref{X decomposition into parts} and
the construction is easy.
To explain the idea, first we explain this easy version of the construction of $\X$ briefly.
We cover $\widehat{\M}^{\leq L_{\max}}_{\theta}$ by open subsets $\U_i$
($i = 1, \dots, N$) and for each $i$,
choose a family $\mathcal{S}_i =\{S\}$ of codimension-two submanifolds
of $Y$, and add the inverse images $(\pi_Y \circ u)^{-1}(S)$ to each holomorphic buildings
$p = (\Sigma, z, u) \in \U_i$ as the marked points $(z^i)^S$.
If we choose an appropriate family $\mathcal{S}_i$, then the $\pi_Y \circ u$ is transverse
to all submanifolds $S \in \mathcal{S}_i$ and all irreducible components of
$(\Sigma, z, (z^i)^S)$ other than trivial cylinders of $(\Sigma, z, u)$ are stable
for all $p = (\Sigma, z, u) \in \U_i$.
Assume that for each $i$, there exists a local universal family
$(\hat P^i \to \hat X^i, Z^i, (Z^i)^S)$ which contains $\forget_u(\Sigma, z, (z^i)^S, u)$
as fibers for all $p = (\Sigma, z, u) \in \U_i$.
Choose finite dimensional vector space $E^0_i$ and linear map
$\lambda_i : E^0_i \to C^\infty(\hat P^i \times Y, \Wedge^{0,1} V^\ast \hat P^i
\otimes_\C (\R \partial_\sigma \oplus TY))$
which makes the linear map (\ref{X surjective map}) surjective for all $p \in \U_i$.
Then we can define $\X_{\theta}$ by the space of holomorphic buildings with
perturbation parameters $(p^+ = (\Sigma, z, \bigcup_{i \in I} (z^i)^S, u),
\bigoplus_{i \in I} E^0_i, \bigoplus_{i \in I} \lambda_i)$
for holomorphic buildings $p = (\Sigma, z, u) \in \widehat{\M}$ and non empty subsets
$I \subset \{1, \dots, N\}$ such that $p \in \U_i$ for all $i \in I$.
For each point $(p^+ = (\Sigma, z, \bigcup_{i \in I} (z^i)^S, u),
\bigoplus_{i \in I} E^0_i, \bigoplus_{i \in I} \lambda_i)$, we associate the Kuranishi
neighborhood of $p = (\Sigma, z, u)$ defined by the direct sum of the pull backs of
$\lambda_i$ by the $(\mathcal{S}, A)$-forgetful maps for all $i \in I$.
To realize $\X$ as a subspace of $\widehat{\M}_{o, \mathcal{S}}$,
we choose a family of sections $(Z^i)^o = ((Z^i)^o_j)$ of $\hat P^i \to \hat X^i$
for each $i$ and add the union of the values of $(Z^i)^o$ for all $i \in I$ to
each $p^+ = (\Sigma, z, \bigcup_{i \in I} \bigcup_{i \in I} (z^i)^S, u)$
as marked points $z^o$.
This is the outline of the construction of $\X$ in the case where
we ignore the algebraic structure of $\widehat{\M}$.

For Condition \ref{X decomposition into parts}, we need to extend the linear maps
$\lambda$ for decomposable holomorphic buildings
given as the union of those associated for the parts
to their neighborhoods in a compatible way.
To compare the linear maps associated to different points in $\X \subset
\widehat{\M}_{o, \mathcal{S}, A}$ for the same holomorphic building,
we need to assume that these points are related by
$(\mathcal{S}, A)$-forgetful maps.
Hence first we construct the part of the marked points
$z^S$ and $z^A$ which enables us to compare the linear mas $\lambda$.

For any constant $C \geq 0$, we construct a finite set
$\mathcal{S} = \{S\}$ of codimension-two submanifolds of $Y$ and subsets
\[
\V_{\theta, l} \subset \U_{\theta, l} \subset
\widehat{\M}^{\leq L_{\max}}_{\mathcal{S}, A, \theta, l}
\]
and
\[
\U_{\theta, l}^\DM \subset \overline{\M}^\DM_{\mathcal{S}, A, \theta, l}
\]
for all triples $\theta$ such that $\widetilde{e}(\theta) \leq C$ and $l \geq 0$
which satisfy the following conditions, and call a family
$(\mathcal{S}, \V_{\theta, l}, \U_{\theta, l}, \U_{\theta, l}^\DM)$
a domain curve representation of
$\widehat{\M}^{\leq L_{\max}}_{\leq C}$.
\begin{enumerate}[label=$(\arabic*)^{\mathrm{D}}$]
\item
\label{(S, A)-covering}
For any $p \in \widehat{\M}^{\leq L_{\max}}_{\theta}$,
there exist some $l \geq 0$ and $p^+ \in \V_{\theta, l}$
such that $\forget_{\mathcal{S}, A}(p^+) = p$.
\item
\label{D-neighborhood}
The image of $\U_{\theta, l}$ by
$\forget_u$ is contained in $\U_{\theta, l}^\DM$.
Furthermore, there exists an open neighborhood $\W_{\theta, l}
\subset \widehat{\M}^{\leq L_{\max}}_{\mathcal{S}, A, \theta, l}$ of the closure of
$\U_{\theta, l}$ such that
\[
\U_{\theta, l} = \{p \in \W_{\theta, l};
\forget_u(p) \in \U_{\theta, l}^\DM\}.
\]
If $\U_{\theta, l}$ and $\U_{\theta, l}^\DM$ satisfy this condition, then
we say $\U_{\theta, l}^\DM$ is a D-neighborhood of $\U_{\theta, l}$.
\item
\label{open sets and inclusions}
$\V_{\theta, l}$ is open in the relative topology of $\U_{\theta, l}$,
and $\V_{\theta, l} \Subset \U_{\theta, l}$.
\item
\label{l max}
For any $\theta$, there exists some $l_{\theta}^{\max} \geq 0$ such that
$\U_{\theta, l} = \emptyset$ and $\U^\DM_{\theta, l} = \emptyset$ for all
$l > l_{\theta}^{\max}$.
\item
\label{Z^A section}
For each point $\hat p \in \U_{\theta, l}^{\DM}$,
there exist a local universal family
$(\hat P \to \hat X, Z, Z^S)$ of $\forget_A(\hat p)$ and
an $\Aut(\hat p)$-invariant family of smooth sections
$Z^A = (Z^A_j)$ of $\hat P \to \hat X$ such that
\[
\{(\hat P_a, Z(a), Z^S(a), Z^A(a)); a \in \hat X\} / \Aut(\hat p)
\]
is a neighborhood of $\hat p$ in $\U_{\theta, l}^{\DM}$.
(We assume that $Z^A_j$ are disjoint.)
We call $(\hat P \to \hat X, Z, Z^S, Z^A)$ a local representation of a
neighborhood of $\hat p$ in $\U_{\theta, l}^\DM$.
We note that we may regard $Z^A$
as an $\Aut(\hat p)$-equivariant section of $(\prod^{\# z^A} \hat P)_{\hat X}
/ \mathfrak{S}_{\# z^A} \to \hat X$,
where $(\prod^{\# z^A} \hat P)_{\hat X}$ is the fiber product
over $\hat X$, and $\mathfrak{S}_{\# z^A}$ acts on it as permutations.
\item
\label{no collapse for (S, A)-forgetful map}
For any $\theta$, $l \geq l' \geq 0$,
$\hat p \in \U_{\theta, l}^\DM$ and $\hat q \in \U_{\theta, l'}^\DM$,
if there exists an $(\mathcal{S}, A)$-forgetful map $f$ from $\hat p$ to $\hat q$,
then $f$ does not collapse any component of $\hat p$.
(Namely, $\hat p$ is a curve obtained by adding some marked points to $\hat q$.)
\item
\label{Z^A pull back}
Under the same assumption,
let $(\hat P \to \hat X, Z, Z^S, Z^A)$ be the local representation of a neighborhood of
$\hat p$ in $\U_{\theta, l}^\DM$, and
$(\hat P' \to \hat X', Z', (Z')^S, (Z')^A)$ be that of $\hat q$ in $\U_{\theta, l'}^\DM$.
Shrink $\hat X$ and $\hat X'$ if necessary, and let $(\phi, \hat \phi)$ be the unique
forgetful map from $(\hat P \to \hat X, Z, Z^S)$ to $(\hat P' \to \hat X', Z', (Z')^S)$
whose restriction to the central fiber coincides with $f$.
Then the pull back of $(Z')^A$ by $(\phi, \hat \phi)$ is contained in $Z^A$ as a subfamily.
\item
\label{decomposition into parts U DM}
For any $\theta = (g, k, E_{\hat \omega})$,
$\hat p \in \U_{\theta, l}^\DM$ and subset $\mathcal{N}$ of its nodal
points, replace each nodal point in $\mathcal{N}$ with a pair of marked points
(we regard the new marked points as points in the set $z$), and
let $\hat p'_i$ $(1 \leq i \leq N)$ be its connected components or an arbitrary
decomposition into unions of its connected components.
Let $g'_i$ and $k'_i$ be the genus and the number of marked points $z$ of
each $\hat p'_i$ respectively.
Then there exist some $E_{\hat \omega}^i \geq 0$ such that
$E_{\hat \omega} = \sum_i E_{\hat \omega}^i$ and
$\hat p'_i \in \U_{\theta'_i, l(\hat p'_i)}^\DM$ for all $i$,
where $\theta'_i = (g'_i, k'_i, E_{\hat \omega}^i)$.
\item
\label{decomposition into parts U}
$\U_{\theta, l}$ satisfies
the following conditions about decomposition of a holomorphic building into parts.
\begin{itemize}
\item
For any $p \in \U_{\theta, l}$ and any decomposition $p_i$ ($1 \leq k$)
into unions of its connected components, let $p'_i$ be the holomorphic buildings
obtained by collapsing trivial floors (floors consisting of trivial cylinders).
Then $p'_i \in \U_{\theta(p'_i), l(p'_i)}$ for all $i$.
\item
For any $p \in \U_{\theta, l}$ and any gap between floors, let $p_1$ and $p_2$
be the holomorphic buildings obtained by dividing $p$ at this gap.
Then $p'_i \in \U_{\theta(p'_i), l(p'_i)}$ for $i = 1, 2$.
\item
For any $p \in \U_{\theta, l}$ and any subset of its nodal points,
the holomorphic building $p'$ obtained by replacing these nodal points to pairs of
marked points is contained in $\U_{\theta(p'), l(p')}$.
\end{itemize}
\item
\label{decomposition into parts V}
For each $p \in \widehat{\M}^{\leq L_{\max}}_{\mathcal{S}, A, \theta, l}$,
replace all nodal points and joint circles of $p$ to pairs of marked points and
pairs of limit circles respectively
(we regard the new marked points as points in the set $z$), and
let $p'_i$ $(1 \leq i \leq k)$ be their non-trivial connected components.
Then $p \in \V_{\theta, l}$ if and only if $p'_i \in \V_{\theta(p'_i), l(p'_i)}$ for all $i$.
\item
\label{existence of minimum}
For any $p = (\Sigma, z, u) \in \widehat{\M}^{\leq L_{\max}}_{\theta}$
and subsets $(z^k)^S \subset \Sigma$ ($S \in \mathcal{S}$) and $(z^k)^A \subset
\Sigma$ ($k =1,2$), if each $p^k = (\Sigma, z, (z^k)^S, (z^k)^A, u)$ is contained in
$\U_{\theta, l(p^k)}$, then
$p^3 = (\Sigma, z, (z^1)^S \cap (z^2)^S, (z^1)^A \cap (z^2)^A, u)$ is contained in
$\U_{\theta, l(p^3)}$.
Furthermore, $(z^1)^S \cup (z^2)^S$ $(S \in \mathcal{S})$ and $(z^1)^A \cup (z^2)^A$
are disjoint.
(This means that we can define a holomorphic building
$(\Sigma, z, (z^1)^S \cup (z^2)^S, (z^1)^A \cup (z^2)^A, u) \in
\widehat{\M}_{\mathcal{S}, A}$, but we do not assume that it is contained in some
$\U_{\theta, l}$.)
\item
\label{stably unique forgetful map DM}
For any $\hat p = (\hat \Sigma, z, z^S, z^A) \in \U_{\theta, l}^\DM$ and
any subsets $(z^1)^S, (z^2)^S \subset z^S$ $(S \in \mathcal{S})$ and
$(z^1)^A, (z^2)^A \subset z^A$, if
each $\hat p^i = (\hat \Sigma, z, (z^i)^S, (z^i)^A)$ is contained in 
$\U_{\theta, l(\hat p^i)}^\DM$, then
$\hat p^3 = (\hat \Sigma, z, (z^1)^S \cap (z^2)^S, (z^1)^A \cap (z^2)^A)$
is stable and it is also contained in $\U_{\theta, l(\hat p^3)}^\DM$.
\end{enumerate}
Note that Condition \ref{decomposition into parts U DM},
\ref{decomposition into parts U} and \ref{decomposition into parts V} are
conditions about one triple $\theta$ and other triples $\theta'$ such that
$\widetilde{e}(\theta') < \widetilde{e}(\theta)$,
and the others are conditions for each $\theta$.

We will prove the existence of a domain curve representation of
$\widehat{\M}^{\leq L_{\max}}_{\leq C}$ in Lemma
\ref{existence of a domain curve representation}.
First we prove that we can shrink $\U_{\theta, l}$ and $\U_{\theta, l}^\DM$
preserving $\V_{\theta, l}$.
\begin{lem}
\label{shrinking and conditions}
Let $C \geq 0$ be an arbitrary constant, and assume that
a domain curve representation
$(\mathcal{S}, \V_{\theta, l}, \U_{\theta, l}, \U^\DM_{\theta, l})$ of
$\widehat{\M}^{\leq L_{\max}}_{\leq C}$ is given.
Then we can construct open subsets 
\[
\V_{\theta, l} \Subset \mathring{\U}_{\theta, l} \Subset \U_{\theta, l}
\quad (\widetilde{e}(\theta) = C, l \geq 0)
\]
and
\[
\mathring{\U}_{\theta, l}^\DM \Subset \U_{\theta, l}^\DM \quad
(\widetilde{e}(\theta) = C, l \geq 0)
\]
such that if we replace $\U_{\theta, l}$ and $\U_{\theta, l}^\DM$
for $\widetilde{e}(\theta) = C$ in the family
$(\mathcal{S}, \V_{\theta, l}, \U_{\theta, l}, \U^\DM_{\theta, l})$
with $\mathring{\U}_{\theta, l}$ and $\mathring{\U}_{\theta, l}^\DM$
respectively, it still satisfies the conditions of a domain curve representation.
\end{lem}
\begin{proof}
The nontrivial conditions are Condition \ref{existence of minimum} and
\ref{stably unique forgetful map DM}.
We construct $\V_{\theta, l}$, $\U_{\theta, l}$ and $\U^\DM_{\theta, l}$ ($l \geq 0$)
for each triple $\theta$ such that $\widetilde{e}(\theta) = C$.
First we consider Condition \ref{stably unique forgetful map DM}.
Consider the following condition for subsets
$\widehat{B}_{\theta, l}^\DM, B_{\theta, l}^\DM \subset \U_{\theta, l}^\DM$:
\begin{itemize}
\item[\ref{stably unique forgetful map DM}$\:\!\!^+$]
For $l^1, l^2 < l^0$,
$\hat p = (\hat \Sigma, z, z^S, z^A) \in \widehat{B}_{\theta, l^0}^\DM$ and
subsets $(z^1)^S, (z^2)^S \subset z^S$ $(S \in \mathcal{S})$ and
$(z^1)^A, (z^2)^A \subset z^A$, if each $\hat p^i = (\hat \Sigma, z, (z^i)^S, (z^i)^A)$ is
contained in $\widehat{B}_{\theta, l^i}^\DM$ and
$\hat p^3 = (\hat \Sigma, z, (z^1)^S \cap (z^2)^S, (z^1)^A \cap (z^2)^A)$
does not coincide with $\hat p^1$ or $\hat p^2$, then $\hat p^3$ is contained
in $B_{\theta, l^3}^\DM$ for some $l^3 < \min(l^1, l^2)$.
\end{itemize}
Condition \ref{stably unique forgetful map DM} is equivalent to
this condition for $\widehat{B}_{\theta, l}^\DM = B_{\theta, l}^\DM
= \U_{\theta, l}^\DM$.
By the decreasing induction in $l \leq l_{\theta}^{\max}$, we construct
open neighborhoods
$\mathring{\U}_{\theta, l}^\DM \Subset \U_{\theta, l}^\DM$ of
the closure of $\forget_{\mathcal{S}, A}(\V_{\theta, l})$
so that for any $l_0$, Condition \ref{stably unique forgetful map DM}$\:\!\!^+$ holds for
$l^k > l_0$ ($k = 0,1,2$) and
\[
\widehat{B}_{\theta, l}^\DM = \overline{\mathring{\U}_{\theta, l}^\DM}
\quad (l > l_0),
\quad
B_{\theta, l}^\DM = \begin{cases}
\mathring{\U}_{\theta, l}^\DM & l \geq l_0 \\
\U_{\theta, l}^\DM & l < l_0
\end{cases}.
\]

For $l = l_{\theta}^{\max}$, we may choose arbitrary open neighborhood
$\mathring{\U}_{\theta, l}^\DM \Subset \U_{\theta, l}^\DM$ of the closure of
$\forget_{\mathcal{S}, A}(\V_{\theta, l})$.
Assume that $\mathring{\U}_{\theta, l}^\DM$ for $l > l_0$ are given.
Define $K_{\theta, l_0}^\DM \subset \U_{\theta, l_0}^\DM$
by the smallest subset such that the above condition holds for
$l^0, l^1, l^2 > l_0$ and
\[
\widehat{B}_{\theta, l}^\DM = \overline{\mathring{\U}_{\theta, l}^\DM}
\quad (l > l_0),
\quad
B_{\theta, l}^\DM = \begin{cases}
\mathring{\U}_{\theta, l}^\DM & l > l_0\\
K_{\theta, l_0}^\DM & l = l_0\\
\U_{\theta, l}^\DM & l < l_0
\end{cases}.
\]
Namely, $\hat p^3 \in \U_{\theta, l_0}^\DM$ is contained in
$K_{\theta, l_0}^\DM$ if
there exists some $l^0 > l_0$, $\hat p = (\hat \Sigma, z, z^S, z^A) \in
\overline{\mathring{\U}_{\theta, l^0}^\DM}$ and
subsets $(z^i)^S \subset z^S$ $(i = 1,2$, $S \in \mathcal{S})$ and
$(z^i)^A \subset z^A$ such that
each $\hat p^i = (\hat \Sigma, z, (z^i)^S, (z^i)^A)$ is contained in 
$\overline{\mathring{\U}_{\theta, l^i}^\DM}$ for some $l_0 < l^1, l^2 < l^0$
and $\hat p^3$ is isomorphic to
$(\hat \Sigma, z, (z^1)^S \cap (z^2)^S, (z^1)^A \cap (z^2)^A)$.
It is easy to check that this is a compact subset of $\U_{\theta, l_0}^\DM$.
Hence an open subset $\mathring{\U}_{\theta, l_0}^\DM \Subset
\U_{\theta, l_0}^\DM$ such that
$K_{\theta, l_0}^\DM \cup \overline{\forget_{\mathcal{S}, A}(\V_{\theta, l})} \subset
\mathring{\U}_{\theta, l_0}^\DM$ satisfies the required condition.
Therefore we can construct open subsets $\mathring{\U}_{\theta, l}^\DM$
by the decreasing induction in $l \leq l_{\theta}^{\max}$.

Condition \ref{existence of minimum} is also similar.
Namely, we consider the following condition for
subsets $\widehat{B}_{\theta, l} \subset \U_{\theta, l}$
and $B_{\theta, l} \subset \U_{\theta, l}$:
\begin{itemize}
\item[\ref{existence of minimum}$\:\!\!^+$]
For any $0 \leq l^1, l^2 \leq l^{\theta}_{\max}$,
$p = (\Sigma, z, u) \in \widehat{\M}^{\leq L_{\max}}_{\theta}$,
and subsets $(z^k)^S \subset \Sigma$ ($S \in \mathcal{S}$) and $(z^k)^A \subset
\Sigma$ ($k =1,2$), if each $p^k = (\Sigma, z, (z^k)^S, (z^k)^A, u)$ is contained in
$\widehat{B}_{\theta, l^k}$ and
$p^3 = (\Sigma, z, (z^1)^S \cap (z^2)^S, (z^1)^A \cap (z^2)^A, u)$ does not
coincides with $p^1$ or $p^2$,
then $p^3$ is contained in $B_{\theta, l^3}$ for some $l^3 < \min(l^1, l^2)$.
Furthermore, $(z^1)^S \cup (z^2)^S$ ($S \in \mathcal{S}$)
and $(z^1)^A \cup (z^2)^A$ are disjoint.
\end{itemize}
Condition \ref{existence of minimum} is equivalent to this condition for
$\widehat{B}_{\theta, l} = B_{\theta, l} = \U_{\theta, l}$.
We construct open subsets $\mathring{\U}_{\theta, l}$
which make this condition holds for $l^1, l^2 > l_0$ and
\[
\widehat{B}_{\theta, l} = \overline{\mathring{\U}_{\theta, l}} \quad (l > l_0),
\quad
B_l = \begin{cases}
\mathring{\U}_{\theta, l} & l \geq l_0\\
\U_{\theta, l} & l < l_0
\end{cases}
\]
for all $l_0$ by the decreasing induction in $l \leq l_{\theta}^{\max}$.
Then as in the previous case, if $\mathring{\U}_{\theta, l}$ for $l > l_0$ are given,
define $K_{\theta, l_0} \subset \U_{\theta, l_0}$ by the smallest subset
which makes Condition \ref{existence of minimum}$\:\!\!^+$ hold for $l^1, l^2 > l_0$ and
\[
\widehat{B}_{\theta, l} = \overline{\mathring{\U}_{\theta, l}} \quad (l > l_0),
\quad
B_l = \begin{cases}
\mathring{\U}_{\theta, l} & l > l_0\\
K_{\theta, l} & l = l_0 \\
\U_{\theta, l} & l < l_0
\end{cases},
\]
then $K_{\theta, l_0}$ is a compact subset contained in $\U_{\theta, l}$.
Hence any open neighborhood
$\mathring{\U}_{\theta, l} \Subset \U_{\theta, l}$ of
$K_{\theta, l} \cup \overline{\V_{\theta, l}}$ satisfies the condition.
Therefore if we choose appropriate $\mathring{\U}_{\theta, l}^\DM$
and define $\mathring{\U}_{\theta, l}$ by Condition \ref{D-neighborhood}, then
Condition \ref{existence of minimum} also holds for $\mathring{\U}_{\theta, l}$.
\end{proof}

Let $\widehat{\M}^\triangle \subset \widehat{\M}$ and
$\widehat{\M}^\triangle_{\mathcal{S}, A} \subset \widehat{\M}_{\mathcal{S}, A}$
be the subspaces of decomposable holomorphic buildings, that is, the subspaces of
disjoint holomorphic buildings and holomorphic buildings with nodal points or joint circles.
Similarly, let $\overline{\M}^{\DM, \triangle}_{\mathcal{S}, A} \subset
\overline{\M}^\DM_{\mathcal{S}, A}$ be the subspace of disjoint stable curves and
stable curves with nodal points.
\begin{lem}
\label{extension from the decomposable}
Let $C \geq 0$ be an arbitrary constant, and assume that
a domain curve representation
$(\mathcal{S}, \V_{\theta, l}, \U_{\theta, l}, \U^\DM_{\theta, l})$ of
\[
\widehat{\M}^{\leq L_{\max}}_{< C} = \bigcup_{\widetilde{e}(\theta) < C}
\widehat{\M}^{\leq L_{\max}}_\theta
\]
is given.
We also assume that subsets
\[
\V_{\theta, l}^\triangle \subset \U_{\theta, l}^\triangle \subset
\widehat{\M}^{\leq L_{\max}, \triangle}_{\mathcal{S}, A, \theta, l}
\quad (\widetilde{e}(\theta) = C, l \geq 0)
\]
and
\[
\U_{\theta, l}^{\DM, \triangle} \subset
\overline{\M}^{\DM, \triangle}_{\mathcal{S}, A, \theta, l}
\quad (\widetilde{e}(\theta) = C, l \geq 0)
\]
are given and they satisfy the conditions of domain curve representation.
More precisely, they satisfy the conditions obtained by replacing
$\U_{\theta, l}$, $\U_{\theta, l}^\DM$, $\widehat{\M}$ and so on
with the counterparts with superscript $\triangle$.
Condition \ref{Z^A section} for $\U_{e^1_0, e^2_0, l}^{\DM, \triangle}$ is read as follows:
For a local universal family $(\hat P \to \hat X, Z, Z^S)$ of
$\forget_A(\hat p)$, let $\hat X^\triangle \subset \hat X$ be the subset of points whose
fiber are disconnected or whose fiber has nodal points.
Then Condition \ref{Z^A section} for $\U_{\theta, l}^{\DM, \triangle}$ is that
there exists an $\Aut(\hat p)$-invariant family of sections $Z^A$ of
$\hat P|_{\hat X^\triangle} \to \hat X^\triangle$ such that
\[
\{(\hat P_a, Z(a), Z^S(a), Z^A(a)); a \in \hat X^\triangle\} / \Aut(\hat p)
\]
is a neighborhood of $\hat p$ in $\U_{\theta, l}^{\DM, \triangle}$.
In Condition \ref{no collapse for (S, A)-forgetful map}, we read $\U_{\theta, l}^\DM$
$($or $\U_{\theta', l'}^\DM$$)$ for $\widetilde{e}(\theta) = C$
as $\U_{\theta, l}^{\DM, \triangle}$ $($or $\U_{\theta', l'}^{\DM, \triangle}$$)$.

Then we can construct subsets
\[
\V_{\theta, l}^1 \subset \U_{\theta, l}^1 \subset \U_{\theta, l}^2 \subset
\widehat{\M}^{\leq L_{\max}}_{\mathcal{S}, A, \theta, l} \quad
(\widetilde{e}(\theta) = C, l \geq 0)
\]
and
\[
\U_{\theta, l}^{1, \DM} \subset \U_{\theta, l}^{2, \DM} \subset
\overline{\M}^\DM_{\mathcal{S}, A, \theta, l}
\quad (\widetilde{e}(\theta) = C, l \geq 0)
\]
which satisfy the following conditions:
\begin{itemize}
\item
The closure of $\U_{\theta, l}^{1, \DM}$ is contained in
$\U_{\theta, l}^{2, \DM}$.
\item
$\V_{\theta, l}^1 \cap \widehat{\M}_{\mathcal{S}, A}^\triangle
= \V_{\theta, l}^\triangle$,
$\U_{\theta, l}^2 \cap \widehat{\M}_{\mathcal{S}, A}^\triangle \subset
\U_{\theta, l}^\triangle$, and
$\U_{\theta, l}^{2, \DM} \cap \overline{\M}_{\mathcal{S}, A}^{\DM, \triangle}
\subset \U_{\theta}^{\DM, \triangle}$.
\item
For each $k= 1,2$, $\V_{\theta, l} = \V_{\theta, l}^1$,
$\U_{\theta, l} = \U_{\theta, l}^k$ and $\U_{\theta, l}^\DM = \U_{\theta, l}^{k, \DM}$
for $\widetilde{e}(\theta) = C$
and the given $\V_{\theta', l}$, $\U_{\theta', l}$
and $\U_{\theta', l}^\DM$ for $\widetilde{e}(\theta') < C$
satisfy the conditions of domain curve representation of
$\widehat{\M}_{\leq C}$ other than Condition \ref{(S, A)-covering}.
\end{itemize}
\end{lem}
\begin{proof}
We consider each triple $\theta$ such that $\widetilde{e}(\theta) = C$.
By the same argument used for the proof of Lemma \ref{shrinking and conditions},
we see that we can take an open neighborhood
$\mathring{\U}_{\theta, l}^\triangle \Subset \U_{\theta, l}^\triangle$ of 
the closure of $\V_{\theta, l}^\triangle$ and
an open subset $\mathring{\U}_{\theta, l}^{\DM, \triangle} \Subset
\U_{\theta, l}^{\DM, \triangle}$ such that
$\V_{\theta, l}^\triangle$, $\mathring{\U}_{\theta, l}^\triangle$ and
$\mathring{\U}_{\theta, l}^{\DM, \triangle}$ also satisfy the assumption.
We may assume that Condition \ref{stably unique forgetful map DM}$\:\!\!^+$ holds for
$\widehat{B}^\DM_{\theta, l}
= \overline{\mathring{\U}_{\theta, l}^{\DM, \triangle}}$ and
$B^\DM_{\theta, l} = \mathring{\U}_{\theta, l}^{\DM, \triangle}$,
and Condition \ref{existence of minimum}$\:\!\!^+$ holds for
$\widehat{B}_{\theta, l} = \overline{\mathring{\U}^\triangle_{\theta, l}}$ and
$B_{\theta, l} = \mathring{\U}^\triangle_{\theta, l}$.

Choose finite points $\hat p^i \in \U_{\theta, l}^{\DM, \triangle}$ ($i \in I_l$)
for each $l$, let $(\hat P^i \to \hat X^i, Z^i,\ab (Z^i)^S)$ be
a local universal family of each $\forget_A(\hat p^i)$ and
let $(Z^i)^A|_{(\hat X^i)^\triangle}$ be the family of section for each $i \in I_l$
so that
$\overline{\mathring{\U}_{\theta, l}^\triangle}$ are covered by
$\{((\hat X^i)^\triangle, (Z^i)^A|_{(\hat X^i)^\triangle})\}_{i \in I_l}$.
Namely, every stable curve in $\overline{\mathring{\U}_{\theta, l}^\triangle}$
appears as some fiber of the families obtained by adding the marked points
$(Z^i)^A|_{(\hat X^i)^\triangle}$ to the local universal families.
Shrinking each $\hat X^i$ if necessary (keeping the covering condition),
we construct an $\Aut(\hat p^i)$-invariant extension $(Z^i)^A = ((Z^i)^A_j)_j$ of
$(Z^i)^A|_{(\hat X^i)^\triangle}$ to $\hat X^i$ and
an $\Aut(\hat p^i)$-invariant open neighborhood $W^i = \coprod_j W^i_j$ of the value of
$(Z^i)^A = ((Z^i)^A_j)_j$ for each $i$ which satisfy the following condition:
\begin{itemize}
\item[$(\ast)$]
For any $l_0 \geq l$, $i_0 \in I_{l_0}$, $i \in I_l$,
$a \in \hat X^{i_0}$, $b \in \hat X^i$ and an $\mathcal{S}$-forgetful map
\[
\varphi : (\hat P^{i_0}_a, Z^{i_0}(a), (Z^{i_0})^S(a)) \to (\hat P^i_b, Z^i(b), (Z^i)^S(b)),
\]
if $\varphi^{-1}((Z^i)^A(b)) \subset W^{i_0}|_a$
then $\varphi^{-1}((Z^i)^A(b)) \subset (Z^{i_0})^A(a)$.
\end{itemize}
We construct such extensions and neighborhoods of their values by the induction in $l$
as follows.

Let $l = l_{\min}$ be the minimum such that
$\U_{\theta, l}^{\DM, \triangle} \neq \emptyset$.
We fix an order of $I_{l_{\min}}$ and construct $(Z^i)^A$ and $W^i$ by the induction
in $i \in I_{l_{\min}}$.
First for the minimal $i \in I_{l_{\min}}$,
we construct an $\Aut(\hat p^i)$-invariant extension $(Z^i)^A$ of
$(Z^i)^A|_{(\hat X^i)^\triangle}$. Since we may regard $(Z^i)^A$
as an $\Aut(\hat p^i)$-equivariant section of the fiber product
$(\prod^{\# z^A} \hat P^i)_{\hat X^i} / \mathfrak{S}_{\# z^A} \to \hat X^i$,
we can extend $(Z^i)^A|_{(\hat X^i)^\triangle}$ to a neighborhood of
$(\hat X^i)^\triangle$.
Replacing $\hat X^i$ with a small neighborhood of $(\hat X^i)^\triangle$,
we may assume that $(Z^i)^A = (Z^i)^A_j$ is defined on $\hat X^i$.
Let $W^i = \coprod_j \mathring{W}^i_j$ be an $\Aut(\hat p^i)$-invariant neighborhood
of $(Z^i)^A(\hat X^i) = \coprod_j (Z^i)^A_j(\hat X^i)$.

Assume that we have constructed $(Z^i)^A$ and $W^i$ for $i < i_0$
which satisfy Condition $(\ast)$.
First we construct an $\Aut(\hat p^{i_0})$-invariant open neighborhood $W^{i_0}$
of the value of $(Z^{i_0})^A|_{(\hat X^{i_0})^\triangle}$ which satisfies the following
conditions:
\begin{enumerate}[label=\normalfont(\roman*)]
\item
\label{two with the same l}
For any $i < i_0$, $a \in \hat X^{i_0}$, $b_1, b_2 \in \hat X^i$ and isomorphisms
\[
\varphi_k : (\hat P^{i_0}_a, Z^{i_0}(a), (Z^{i_0})^S(a))
\cong (\hat P^i_{b_k}, Z^i(b_k), (Z^i)^S(b_k))
\quad (k = 1, 2),
\]
if $\varphi_k^{-1}((Z^i)^A(b_k)) \subset W^{i_0}|_a$ for $k = 1, 2$,
then there exists some $g \in \Aut(\hat p^i)$ such that
$b_2 = g b_1$ and $\varphi_2 = g \circ \varphi_1$.
In particular, $\varphi_1^{-1}((Z^i)^A(b_1)) = \varphi_2^{-1}((Z^i)^A(b_2))$.
\item
\label{three with the same l}
For any $i_1, i_2 < i_0$, $a \in \hat X^{i_0}$, $b_k \in \hat X^{i_k}$ ($k = 1, 2$) and
isomorphisms
\[
\varphi_k : (\hat P^{i_0}_a, Z^{i_0}(a), (Z^{i_0})^S(a))
\cong (\hat P^{i_k}_{b_k}, Z^{i_k}(b_k), (Z^{i_k})^S(b_k))
\quad (k =1, 2),
\]
if $\varphi_k^{-1}((Z^{i_k})^A(b_k)) \subset W^{i_0}|_a$ for $k = 1, 2$, then the isomorphism
\[
\varphi = \varphi_2 \circ \varphi_1^{-1}
: (\hat P^{i_1}_{a_1}, Z^{i_1}(a_1), (Z^{i_1})^S(a_1))
\cong (\hat P^{i_2}_{a_2}, Z^{i_2}(a_2), (Z^{i_2})^S(a_2))
\]
satisfies $\varphi((Z^{i_1})^A(b_1)) \subset W^{i_2}|_{b_2}$.
Note that by Condition $(\ast)$, this implies that
$(Z^{i_2})^A(b_2) = \varphi((Z^{i_1})^A(b_1))$.
In particular, $\varphi_1^{-1}((Z^{i_1})^A(b_1))$ coincides with
$\varphi_2^{-1}((Z^{i_2})^A(b_2))$.
\end{enumerate}
Note that in Condition \ref{two with the same l} for $a \in (\hat X^{i_0})^\triangle$,
if $W^{i_0}|_a$ is a sufficiently small neighborhood of $(Z^{i_0})^A(a)$, then
the condition $\varphi_k^{-1}((Z^i)^A(b_k)) \subset W^{i_0}|_a$ implies that
$\varphi_k^{-1}((Z^i)^A(b_k)) = (Z^{i_0})^A(a)$.
Hence $g = \varphi_2 \circ \varphi_1^{-1}$ maps $(Z^i)^A(b_1)$ to $(Z^i)^A(b_2)$,
which implies $g \in \Aut(\hat p^i)$.
Therefore Condition \ref{two with the same l} for general $a \in \hat X^{i_0}$ also holds
if $W^{i_0}$ is a sufficiently small neighborhood of the values of
$(Z^{i_0})^A|_{(\hat X^{i_0})^\triangle}$.
Similarly, in Condition \ref{three with the same l} for $a \in (\hat X^{i_0})^\triangle$,
if $W^{i_0}|_a$ is a sufficiently small neighborhood of $(Z^{i_0})^A(a)$, then
the condition $\varphi_k^{-1}((Z^i)^A(b_k)) \subset W^{i_0}|_a$ implies that
$\varphi_k^{-1}((Z^i)^A(b_k)) = (Z^{i_0})^A(a)$.
Hence $\varphi = \varphi_2 \circ \varphi_1^{-1}$ maps $(Z^{i_1})^A(b_1)$ to
$(Z^{i_2})^A(b_2) \subset W^{i_2}|_{b_2}$.
It implies that Condition \ref{three with the same l} for general $a \in \hat X^{i_0}$
holds if $W^{i_0}$ is a sufficiently small neighborhood of the values of
$(Z^{i_0})^A|_{(\hat X^{i_0})^\triangle}$.

For each point $a \in \hat X^{i_0}$ such that there exists some $i < i_0$,
$b \in \hat X^i$ and isomorphism
\[
\varphi : (\hat P^{i_0}_a, Z^{i_0}(a), (Z^{i_0})^S(a)) \cong (\hat P^i_b, Z^i(b), (Z^i)^S(b))
\]
such that $\varphi^{-1}((Z^i)^A(b)) \subset W^{i_0}|_a$, we define
$(Z^{i_0})^A(a) = \varphi^{-1}((Z^i)^A(b))$.
The above conditions of $W^{i_0}$ implies that this definition is independent of the choice
of $i$ and $b \in \hat X^i$ if they exist.
Shrinking $\hat X^i$ ($i < i_0$) slightly if necessary for smooth extension, we extend
$(Z^{i_0})^A$ to a neighborhood of $(\hat X^{i_0})^\triangle \subset \hat X^{i_0}$.
Replacing $\hat X^{i_0}$ to a neighborhood of
$(\hat X^{i_0})^\triangle \subset \hat X^{i_0}$, we assume that $(Z^{i_0})^A$ is defined
on whole of $\hat X^{i_0}$ and its value is contained in $W^{i_0}$.

Next we consider the general $l$.
Assume that we have already constructed the extensions for $l < l_0$.
We fix an order of $I_{l_0}$ and construct $(Z^i)^A$ and $W^i$ by the induction
in $i \in I_{l_0}$.
Assume that $(Z^i)^A$ and $W^i$ for $i < i_0$ are given.
As in the case of minimal $l$,
first we construct an $\Aut(\hat p^i)$-invariant open neighborhood $W^{i_0}$
of the value of $(Z^i)^A|_{(\hat X^{i_0})^\triangle}$ which satisfies
Condition \ref{two with the same l}, \ref{three with the same l} and
the following condition:
\begin{enumerate}[label=\normalfont(\roman*)]
\setcounter{enumi}{2}
\item
\label{stably unique S-forget}
For any $a \in \hat X^{i_0}$, $l_1, l_2 < l_0$, $i_k \in I_{l_k}$ ($k = 1, 2$),
$b_k \in \hat X^{i_k}$ ($k = 1, 2$) and $\mathcal{S}$-forgetful maps
\[
\varphi_k : (\hat P^{i_0}_a, Z^{i_0}(a), (Z^{i_0})^S(a))
\to (\hat P^{i_k}_{b_k}, Z^{i_k}(b_k), (Z^{i_k})^S(b_k))
\quad (k =1,2),
\]
if $\varphi_k^{-1}((Z^{i_k})^A(b_k)) \subset W^{i_0}|_a$ for $k = 1,2$,
then there exist some $l_3 \leq \min(l_1, l_2)$, $i_3 \in I_{l_3}$,
$b_3 \in \hat X^{i_3}$ and $\mathcal{S}$-forgetful maps
\[
\psi_k : (\hat P^{i_k}_{b_k}, Z^{i_k}(b_k), (Z^{i_k})^S(b_k))
\to (\hat P^{i_3}_{b_3}, Z^{i_3}(b_3), (Z^{i_3})^S(b_3))
\quad (k = 1,2)
\]
which satisfy $\psi_1 \circ \varphi_1 = \psi_2 \circ \varphi_2$ and
the following condition:
For any triple $(j, j_1, j_2)$ such that $\varphi_k^{-1}((Z^{i_k})^A_{j_k}(b_k))
\subset W^i_j|_a$ for $k =1,2$, there exists some $j_3$ such that
$\psi_k^{-1}((Z^{i_3})^A_{j_3}(b_3)) \subset W^{i_k}_{j_k}|_{b_k}$ for $k = 1,2$.
By Condition $(\ast)$, this implies that
$(Z^{i_k})^A_{j_k}(b_k) = \psi_k^{-1}((Z^{i_3})^A_{j_3}(b_3))$ for $k = 1,2$.
In particular, $\varphi_1^{-1}((Z^{i_1})^A_{j_1}(b_1)) = \varphi_2^{-1}((Z^{i_2})^A_{j_2}(b_2))$.
\end{enumerate}
Note that in the above condition, if $a \in (\hat X^{i_0})^\triangle$ and
$W^{i_0}|_a$ is a sufficiently small neighborhood of $(Z^{i_0})^A(a)$, then the condition
$\varphi_k^{-1}((Z^{i_k})^A(b_k)) \subset W^{i_0}|_a$ implies that
$\varphi_k^{-1}((Z^{i_k})^A(b_k)) \subset (Z^{i_0})^A(a)$.
Hence Condition \ref{stably unique forgetful map DM} for
$\mathring{\U}_{\theta, l}^{\DM, \triangle}$ implies that
there exist some $l_3 \leq \min(l_1, l_2)$, $i_3 \in I_{l_3}$,
$b_3 \in \hat X^{i_3}$ and $\mathcal{S}$-forgetful maps
\[
\psi_k : (\hat P^{i_k}_{b_k}, Z^{i_k}(b_k), (Z^{i_k})^S(b_k))
\to (\hat P^{i_3}_{b_3}, Z^{i_3}(b_3), (Z^{i_3})^S(b_3))
\quad (k = 1,2)
\]
such that $\psi_1 \circ \varphi_1 = \psi_2 \circ \varphi_2$ and
\[
\varphi_1^{-1}((Z^{i_1})^A(b_1)) \cap \varphi_2^{-1}((Z^{i_2})^A(b_2))
= (\psi_1 \circ \varphi_1)^{-1}((Z^{i_3})^A(b_3)).
\]
Therefore Condition \ref{stably unique S-forget} holds for $a \in (\hat X^{i_0})^\triangle$
if $W^{i_0}$ is a sufficiently small neighborhood of $(Z^{i_0})^A((\hat X^{i_0})^\triangle)$.
Hence it also holds if $a \in \hat X^{i_0}$ is contained in a small neighborhood of
$(\hat X^{i_0})^\triangle$.
Therefore Condition \ref{stably unique S-forget} holds for general $a \in \hat X^{i_0}$
if $W^{i_0}$ is sufficiently small.

The construction of $(Z^{i_0})^A$ is similar to the case of minimal $l$, but in this case,
some part of $(Z^{i_0})^A$ is determined by the pull backs of $(Z^i)^A$ ($i \in I_l, l < l_0$)
as follows.
For $a \in \hat X^{i_0}$, $l < l_0$, $i \in I_l$, $b \in \hat X^i$ and
an $\mathcal{S}$-forgetful map
\[
\varphi : (\hat P^{i_0}_a, Z^{i_0}(a), (Z^{i_0})^S(a)) \to (\hat P^i_b, Z^i(b), (Z^i)^S(b))
\]
such that $\varphi^{-1}((Z^i)^A(b)) \subset W^{i_0}|_a$,
we define $(z^{i_0})^A(a)_{b, \varphi} = \varphi^{-1}((Z^i)^A(b))$.
For each $a \in \hat X^{i_0}$, we define $(z^{i_0})^A(a)$ by the union of
$(z^{i_0})^A(a)_{b, \varphi}$ over the above pairs $(b, \varphi)$.
We need to construct the extension $(Z^{i_0})^A$ which contains $(z^{i_0})^A$
as a subfamily.

Condition \ref{stably unique S-forget} implies that $\varphi^{-1}((Z^i)^A(b)) \cap
W^{i_0}_j|_a$ consists of at most one point for each $j$, and this point is
independent of $(b, \varphi)$ if it exists.
Hence $(z^{i_0})^A(a) \cap W^{i_0}_j|_a$ consists of at most one point for each $j$.
It is clear that $(z^{i_0})^A$ is $\Aut(\hat p^{i_0})$-invariant.
Hence shrinking $\hat X^i$ for $i \in I_l$ ($l < l_0$) and $i \in I_{l_0}$ such that
$i < i_0$ if necessary for smooth extension,
we can construct an $\Aut(\hat p^{i_0})$-invariant extension $(Z^{i_0})^A$ of
$(Z^{i_0})^A|_{(\hat X^{i_0})^\triangle}$ to a neighborhood of $(\hat X^{i_0})^\triangle$
which contains $(z^{i_0})^A$ as a subfamily.
Replacing $\hat X^{i_0}$ with a small neighborhood of $(\hat X^{i_0})^\triangle$,
we get an extension $(Z^{i_0})^A$ on $\hat X^{i_0}$ such that
$(Z^{i_0})^A(\hat X^{i_0}) \subset W^{i_0}$.
Therefore the induction works and we can construct
$\Aut(\hat p^i)$-invariant extensions $(Z^i)^A = ((Z^i)^A_j)_j$ of
$(Z^i)^A|_{(\hat X^i)^\triangle}$ to $\hat X^i$ and
$\Aut(\hat p^i)$-invariant open neighborhoods $W^i = \coprod_j W^i_j$
of the values of $(Z^i)^A = ((Z^i)^A_j)_j$ which satisfy Condition $(\ast)$.

Now we construct $\V_{\theta, l}^1$, $\U_{\theta, l}^k$ and
$\U_{\theta, l}^{k, \DM}$ ($k = 1, 2$) as follows.
First we define $\widehat{\U}_{\theta, l}^\DM \subset
\overline{\M}^\DM_{\mathcal{S}, A, \theta, l}$
by the union of the sets of stable curves
\[
\{(\hat P^i_a, Z^i(a), (Z^i)^S(a), (Z^i)^A(a));
a \in \hat X^i\}
\]
over $i \in I_l$.
We construct $\U_{\theta, l}^{2, \DM}$ as a subset of
$\widehat{\U}_{\theta, l}^\DM$.
For each $l$, let $\W_{\theta, l} \subset
\widehat{\M}^{\leq L_{\max}}_{\mathcal{S}, A, \theta, l}$ be
an open neighborhood of the closure of $\U_{\theta, l}^\triangle$ such that
\[
\U_{\theta, l}^\triangle = \{p \in \W_{\theta, l};
\forget_{\mathcal{S}, A}(p) \in \U_{\theta, l}^{\DM, \triangle}\}.
\]
We construct small neighborhoods $\U_{\theta, l}^{2, \DM} \subset
\widehat{\U}_{\theta, l}^\DM$ of $\mathring{\U}_{\theta, l}^{\DM, \triangle}$
such that they satisfy Condition \ref{no collapse for (S, A)-forgetful map} and
\ref{stably unique forgetful map DM}, and
\[
\U_{\theta, l}^2 = \{p \in \W_{\theta, l};
\forget_{\mathcal{S}, A}(p) \in \U_{\theta, l}^{2, \DM}\}
\]
satisfy Condition \ref{existence of minimum} as follows.
(We also assume that $\U_{\theta, l}^{2, \DM}$ is sufficiently small so that
$\W_{\theta, l}$ is still an open neighborhood of the closure of
$\U_{\theta, l}^2$.)

Since $\overline{\mathring{\U}_{\theta, l}^{\DM, \triangle}}$ ($l \geq 0$)
satisfy Condition \ref{no collapse for (S, A)-forgetful map}, sufficiently small
neighborhoods $\U_{\theta, l}^{2, \DM}$ of $\mathring{\U}_{\theta, l}^{\DM, \triangle}$
also satisfy the same condition.
For Condition \ref{existence of minimum} and \ref{stably unique forgetful map DM},
we construct $\U_{\theta, l}^{2, \DM}$ by the (usual increasing) induction in $l$ so that
for any $l_0$,
Condition \ref{stably unique forgetful map DM}$\:\!\!^+$ holds
for $l^1, l^2 < l^0$ such that
$\min(l^1, l^2) \leq l_0$ and
\[
\widehat{B}_{\theta, l}^\DM = \begin{cases}
\overline{\mathring{\U}_{\theta, l}^{\DM, \triangle}}& l > l_0\\
\overline{\U_{\theta, l}^{2, \DM}} & l \leq l_0
\end{cases},
\quad
B_{\theta, l}^\DM = \U_{\theta, l}^{2, \DM} \quad (l < l_0),
\]
and Condition \ref{existence of minimum}$\:\!\!^+$ holds for $l^1, l^2$ such that
$\min(l^1, l^2) \leq l_0$ and
\[
\widehat{B}_{\theta, l} = \begin{cases}
\overline{\mathring{\U}_{\theta, l}^\triangle}& l > l_0\\
\overline{\U_{\theta, l}^2} & l \leq l_0
\end{cases},
\quad
B_{\theta, l}^\DM = \U_{\theta, l}^2 \quad (l < l_0).
\]
The induction works because of the following reason.
Assume that $\U_{\theta, l}^{2, \DM}$ for $l < l_0$ are given and they satisfy the
above conditions.
We prove that if $\U_{\theta, l_0}^{2, \DM} \subset \widehat{\U}_{\theta, l_0}^\DM$ is
a sufficiently small open neighborhood of
$\overline{\mathring{\U}_{\theta, l_0}^{\DM, \triangle}}$, then
Condition \ref{stably unique forgetful map DM}$\:\!\!^+$ holds for
$l^1, l^2 < l^0$ such that $\min(l^1, l^2) \leq l_0$ and
\[
\widehat{B}_{\theta, l}^\DM = \begin{cases}
\overline{\mathring{\U}_{\theta, l}^{\DM, \triangle}}& l > l_0\\
\overline{\U_{\theta, l}^{2, \DM}} & l \leq l_0
\end{cases},
\quad
B_{\theta, l}^\DM = \U_{\theta, l}^{2, \DM} \quad (l < l_0).
\]
First we consider this condition for $l^1 = l^2 = l_0 < l^0$.
Assume that this condition does not hold for any small open neighborhood
$\U_{\theta, l_0}^{2, \DM}$ of $\overline{\mathring{\U}_{\theta, l_0}^{\DM, \triangle}}$.
Then there exists a sequence of stable curves
$\hat p_k = (\hat \Sigma_k, z_k, z_k^S, z_k^A) \in
\overline{\mathring{\U}_{\theta, l^0}^{\DM, \triangle}}$
and subsets $(z^1_k)^S, (z^2_k)^S \subset z_k^S$ $(S \in \mathcal{S})$ and
$(z^1_k)^A, (z^2_k)^A \subset z_k^A$ such that
two sequences $(\hat p^i_k)_{k \in \N} = (\hat \Sigma_k, z, (z^i_k)^S, (z^i_k)^A)_{k \in \N}$
converge to points in $\overline{\mathring{\U}_{\theta, l_0}^{\DM, \triangle}}$,
but none of $\hat p^3_k
= (\hat \Sigma_k, z, (z^1_k)^S \cap (z^2_k)^S, (z^1_k)^A \cap (z^2_k)^A)$
is not contained in $\bigcup_{l < l_0} \U_{\theta, l}^{2, \DM}$.
Taking a subsequence, assume that $\hat p_k$ converges to a stable curve
$\hat p_\infty = (\hat \Sigma_\infty, z_\infty, z_\infty^S, z_\infty^A) \in
\overline{\mathring{\U}_{\theta, l^0}^{\DM, \triangle}}$.
We may assume that there exists subsets
$(z^1_\infty)^S, \ab (z^2_\infty)^S \subset z_\infty^S$ $(S \in \mathcal{S})$ and
$(z^1_\infty)^A, \ab (z^2_\infty)^A \subset z_\infty^A$ such that
each $(\hat p^i_k)_{k \in \N}$ converges to a stable curve
$\hat p^i_\infty
= (\hat \Sigma_\infty, z, (z^i_\infty)^S, (z^i_\infty)^A)
\in \overline{\mathring{\U}_{\theta, l_0}^{\DM, \triangle}}$ for $i = 1,2$.
(This is because of Condition \ref{no collapse for (S, A)-forgetful map}
for $\overline{\mathring{\U}_{\theta, l}^{\DM, \triangle}}$.)
Hence Condition \ref{stably unique forgetful map DM}$\:\!\!^+$ for
\[
\widehat{B}_{\theta, l}^\DM = \begin{cases}
\overline{\mathring{\U}_{\theta, l}^{\DM, \triangle}}& l \geq l_0\\
\overline{\U_{\theta, l}^{2, \DM}} & l < l_0
\end{cases},
\quad
B_{\theta, l}^\DM = \U_{\theta, l}^{2, \DM} \quad (l < l_0-1),
\]
implies that
$\hat p^3_\infty = (\hat \Sigma_\infty, z, (z^1_\infty)^S \cap (z^2_\infty)^S,
(z^1_\infty)^A \cap (z^2_\infty)^A)$ is contained in
$\bigcup_{l < l_0} \U_{\theta, l}^{2, \DM}$.
Since $(\hat p^3_k)_{k \in \N}$ converges to $\hat p^3_\infty$,
this contradicts to the openness of $\bigcup_{l < l_0} \U_{\theta, l}^{2, \DM}$.
The other cases such as $l^1, l^2 < l^0 = l_0$ are similar.
%
Condition \ref{existence of minimum}$\:\!\!^+$ is also similar.
Hence we can construct open neighborhoods $\U_{\theta, l}^{2, \DM}$
of $\mathring{\U}_{\theta, l}^{\DM, \triangle}$ which satisfy
Condition \ref{no collapse for (S, A)-forgetful map}, \ref{existence of minimum}
and \ref{stably unique forgetful map DM}.

Next we construct $\U_{\theta, l}^{1, \DM}$ by the same way as
$\U_{\theta, l}^{2, \DM}$ under the condition $\U_{\theta, l}^{1, \DM} \Subset
\U_{\theta, l}^{2, \DM}$, and define $\U_{\theta, l}^1$ by
\[
\U_{\theta, l}^1 = \{p \in \W_{\theta, l};
\forget_{\mathcal{S}, A}(p) \in \U_{\theta, l}^{1, \DM}\}.
\]

Finally, we take open subsets
$\V_{\theta, l}^1 \Subset \U_{\theta, l}^1$ such that $\V_{\theta, l}^1 \cap
\widehat{\M}_{\mathcal{S}, A}^\triangle = \V_{\theta, l}^\triangle$.
Then these $\V_{\theta, l}^1$, $\U_{\theta, l}^k$ and $\U_{\theta, l}^{k, \DM}$
($k = 1,2$) are the required subsets.
\end{proof}

Now we explain the construction of a domain curve representation.
\begin{lem}\label{existence of a domain curve representation}
For any constant $C \geq 0$, there exists a domain curve representation of
$\widehat{\M}^{\leq L_{\max}}_{\leq C}$.
\end{lem}
\begin{proof}
First we claim that in general, for a holomorphic building
$p = (\Sigma, z, u) \in \widehat{\M}$ whose domain curve is irreducible and
which has nonzero $E_{\hat \omega}$-energy,
there exist a finite set $\mathcal{S} = \{S\}$ of codimension-two small disks in $Y$
and an open subset $U \subset \widehat{\M}_{\mathcal{S}}$ such that
\begin{itemize}
\item
$p \in \forget_{\mathcal{S}}(U)$,
\item
the restriction of
$\forget_u : \widehat{\M}_{\mathcal{S}} \to \overline{\M}^\DM_{\mathcal{S}}$
to $U$ is injective, and
\item
$\Aut(\forget_u(q)) = \Aut(q)$ for all $q \in U$.
\end{itemize}
This is proved as follows.
Since the critical points of $\pi_Y \circ u$ are discrete,
we can add marked points by the intersections with codimension-two disks in $Y$
to make the domain curve stable.
Namely, we can choose a finite set $\mathcal{S} = \{S\}$ of codimension-two disks
in $Y$ and a point $p^+ = (\Sigma, z, z^S, u) \in \widehat{\M}_{\mathcal{S}}$ such that
$z^S$ is $\Aut(p)$-invariant,
$\forget_{\mathcal{S}} p^+ = p$ and
$\Aut(\forget_u(p^+)) = \Aut(p^+)$ ($= \Aut(p)$).
Let $D_u : \widetilde{W}^{1, p}_\delta(\Sigma, u^\ast T \hat Y) \to
L^p_\delta(\Sigma, \Wedge^{0, 1} T\ast \Sigma \otimes_\C T \hat Y)$ be the
linearization of the equation of $J$-holomorphic curve at $u$ and
consider the linear operator
\begin{align*}
D^+_u : \widetilde{W}^{1, p}_\delta(\Sigma, u^\ast T \hat Y)
&\to L^p_\delta(\Sigma, \Wedge^{0, 1} T^\ast \Sigma \otimes_\C T \hat Y) \\
& \quad \oplus
\bigoplus_{z \in z^S, S \in \mathcal{S}} T_{\pi_Y \circ u(z)} Y/T_{\pi_Y \circ u(z)}S
\oplus \R
\end{align*}
defined by $D^+_u \xi = (D_u \xi, (\pi_Y)_\ast\xi(z), \sigma_\ast \xi (R))$,
where $R \in \Sigma$ is an arbitrary fixed point and $\sigma : \hat Y \to \R$ is
the projection. If it is injective, then
$\forget_u : \widehat{\M}_{\mathcal{S}} \to
\overline{\M}^\DM_{\mathcal{S}}$ is injective on a neighborhood of $p^+$.
We can choose $\mathcal{S} = \{S\}$ which makes $D^+_u$ injective because
for any vector $\xi \in \Ker D_u$ other than $\xi = c \partial_\sigma$ ($c \in \R$ is
a constant), there does not exist any non-empty open subset of $\Sigma$ on which
$\pi_{TY} \xi$ vanishes.
Therefore for any $p \in \widehat{\M}$, we can construct finite number of
disks $\mathcal{S} = \{S\}$ in $Y$ and
an open subset $U \subset \widehat{\M}_{\mathcal{S}}$
which satisfy the above conditions.

We note the following fact:
For each disk $S \in \mathcal{S}$, let $S \times \R^2 \subset Y$ be
its tubular neighborhood. Then for any small $x \in \R^2$,
$\mathcal{S}^x = \{S \times \{x\}; S \in \mathcal{S}\}$ also satisfies
the same condition.
Namely, there exists an open subset $U^x \subset \widehat{\M}_{\mathcal{S}^x}$
such that $p \in \forget_{\mathcal{S}}(U^x)$,
the restriction of
$\forget_u|_{U^x} : U^x \to \overline{\M}^\DM_{\mathcal{S}^x}$ is injective, and
$\Aut(\forget_u(q)) = \Aut(q) = \Aut(\forget_{\mathcal{S}^x}(q))$ for all $q \in U^x$.
We may assume that $\forget_{\mathcal{S}^x}(U^x) = \forget_{\mathcal{S}}(U)$
for all sufficiently small $x \in \R^2$.

We construct
\[
\V_{\theta, l} \subset \U_{\theta, l} \subset
\widehat{\M}^{\leq L_{\max}}_{\mathcal{S}, A, \theta, l}
\]
and
\[
\U_{\theta, l}^\DM \subset \overline{\M}^\DM_{\mathcal{S}, A, \theta, l}
\]
for $\widetilde{e}(\theta) \leq C$ and $l \geq 0$
by the induction in $\widetilde{e}(\theta)$.

For each triple $\theta = (g, k, E_{\hat \omega})$ with minimal $\widetilde{e}(\theta)$,
$\widehat{\M}^{\leq L_{\max}}_{\theta}$ consists of connected height-one
holomorphic buildings without nodal points.
First we consider the case of $E_{\hat \omega} > 0$.
In this case, first we construct a finite set
$\mathcal{S} = \{S\}$ of codimension-two submanifolds of $Y$ and
open subsets $\mathring{U}_{\theta} \subset
\widehat{\M}^{\leq L_{\max}}_{\mathcal{S}, \theta}$ which satisfy
the following conditions:
(Recall that $\widehat{\M}^{\leq L_{\max}}_{\mathcal{S}, \theta}$ is the subspace of
$\widehat{\M}^{\leq L_{\max}}_{\mathcal{S}, A, \theta}$ defined by $z^A = \emptyset$.)
\begin{enumerate}[label=(\roman*)]
\item
\label{circ U injective}
$\forget_u|_{\mathring{U}_{\theta}}
: \mathring{U}_{\theta} \to \overline{\M}^\DM_{\mathcal{S}, \theta}$ is injective,
and $\Aut(\forget_u(p)) = \Aut(p)$ for any $p \in \mathring{U}_{\theta}$.
Furthermore, for any $p, q \in \mathring{U}_{\theta}$, if there exists
an $(\mathcal{S}, A)$-forgetful map from $\forget_u(p)$ to $\forget_u(q)$, then
$p \geq q$.
\item
$\widehat{\M}^{\leq L_{\max}}_{\theta}$ is covered by the image of
$\mathring{U}_{\theta}$ by $\forget_{\mathcal{S}}$.
\item
\label{z^S disjoint circ U}
For any two holomorphic buildings
$p^i = (\Sigma, z, z^{S, i}, u) \in \mathring{U}_{\theta}$ ($i = 1, 2$) such that
$\forget_{\mathcal{S}}(p^1) = \forget_{\mathcal{S}}(p^2) = (\Sigma, z, u)$,
the following holds:
\begin{itemize}
\item
For any two different submanifolds $S_1 \neq S_2 \in \mathcal{S}$,
$z^{S_1, 1}$ and $z^{S_2, 2}$ are disjoint in $\Sigma$.
\item
For any $S \in \mathcal{S}$, if $z^{S, 1} \neq \emptyset$ and
$z^{S, 2} \neq \emptyset$ then $z^{S, 1} = z^{S, 2}$.
\end{itemize}
\end{enumerate}

We can construct such submanifolds $\mathcal{S} = \{S\}$ and
open subsets $\mathring{U}_{\theta}$ as follows.
The claim we proved in the above implies the following:
There exists an open covering $\{V_i\}$ of $\widehat{\M}^{\leq L_{\max}}_{\theta}$,
and for each $V_i$, there exist an infinite family of finite sets
$\mathcal{S}_i^x = \{S^x_i\}_{S_i \in \mathcal{S}_i}$ ($x \in \R^2$) of
codimension two disks of $Y$ and a family of open subsets
$U^x_i \subset \widehat{\M}^{\leq L_{\max}}_{\mathcal{S}^x, \theta}$
such that $\forget_u|_{U^x_i} : U^x_i \to \overline{\M}^\DM_{\mathcal{S}^x}$ is injective,
$\Aut(\forget_u(p)) = \Aut(p) = \Aut(\forget_{\mathcal{S}^x}(p))$ for all $p \in U^x_i$,
and $\forget_{\mathcal{S}^x}(U^x_i) = V_i$.
Furthermore, $\{S^x_i; x \in \R^2\}$ are disjoint for each $i$.
We choose finite numbers $x_i^1, x_i^2, \dots, x_i^{N_i} \in \R^2$
and open subsets $\mathring{U}_i^k \subset U_i^{x_i^k}$ which satisfy
the following conditions,
where we abbreviate $\mathcal{S}_i^{x_i^k}$ as $\mathcal{S}_i^k$:
\begin{itemize}
\item
$\{\forget_{\mathcal{S}_i^k}(\mathring{U}_i^k)\}_{i, k}$ covers
$\widehat{\M}^{\leq L_{\max}}_{\theta}$
\item
For any $p = (\Sigma, z, z^S, u) \in U_i^k$ and $p' = (\Sigma, z, (z')^S, u) \in U_j^l$,
if $\forget_{\mathcal{S}_i^k}(p) = \forget_{\mathcal{S}_j^l}(p')$, then
$z^S$ and $(z')^{S'}$ are disjoint in $\Sigma$ for any two different submanifolds
$S \neq S' \in \coprod_{i, k} \mathcal{S}_i^k$,
\end{itemize}
Then $\mathcal{S} = \coprod_{i, k} \mathcal{S}_i^k$ and
$\mathring{U}_{\theta} = \bigcup_{i, k} \mathring{U}_i^k$ satisfy
Condition \ref{circ U injective} to \ref{z^S disjoint circ U}.

We explain how to choose such numbers $x_i^k \in \R^2$ and open subsets
$\mathring{U}_i^k \subset U_i^{x_i^k}$.
Take open subsets $\mathring{V}_i \Subset V_i$ which cover
$\widehat{\M}^{\leq L_{\max}}_{\theta}$.
We construct $x_i^k \in \R^2$ and $\mathring{U}_i^k \subset U_i^{x_i^k}$
by the induction in $i$ so that $\{\forget_{\mathcal{S}_i^k}(\mathring{U}_i^k)\}_k$
covers $\mathring{V}_i$ for each $i$ as follows.
We assume that $x_i^k \in \R^2$ and $\mathring{U}_i^k \subset U_i^{x_i^k}$ are given
for $i < i_0$, and we construct those for $i = i_0$.
For each $p = (\Sigma, z, u) \in V_{i_0}$, let $A_p \subset \Sigma$ be the subset of the
points which appear in some $z^S$ for $p^ + = (\Sigma, z, z^S, u) \in \mathring{U}_i^k$
($i < i_0$) such that $\forget_{\mathcal{S}_i^k}(p^+) = p$.
Let $N > 0$ be a constant larger than $\# A_p$ for any $p \in V_{i_0}$.
Choose arbitrary points $x_{i_0}^1, \dots, x_{i_0}^N \in \R^2$.
Then for any $p \in V_{i_0}$, there exists at least one $x_{i_0}^k$ such that
for the point $p^+ = (\Sigma, z, z^S, u) \in U_{i_0}^{x_{i_0}^k}$ such that
$\forget_{\mathcal{S}_{i_0}^k}(p^+) = p$,
each $z^S$ is disjoint with $A_p$.
This is because $S^x_{i_0}$ for $x \in \R^2$ are disjoint.
Hence we can construct open subsets $\mathring{U}_{i_0}^k \subset U_{i_0}^{x_{i_0}^k}$
such that their images by $\forget_{\mathcal{S}_{i_0}^k}$ cover $\mathring{V}_i$
and for any $p^+ = (\Sigma, z, z^S, u) \in \mathring{U}_{i_0}^k$ and
$p = \forget_{\mathcal{S}_{i_0}^k}(p^+)$, each $z^S$ is disjoint with $A_p$.
Therefore we can construct a finite set
$\mathcal{S} = \{S\}$ of codimension-two submanifolds of $Y$ and
open subsets $\mathring{U}_{\theta} \subset
\widehat{\M}^{\leq L_{\max}}_{\mathcal{S}, \theta}$ which satisfy Condition
\ref{circ U injective} to \ref{z^S disjoint circ U}.

Let $U_{\theta} \subset \widehat{\M}_{\mathcal{S}}$ be the set of
holomorphic buildings $p = (\Sigma, z, z^S, u)$ such that there exist subsets
$z^{S, i} \subset z^S$ ($i = 1, \dots, k$, $S \in \mathcal{S}$) such that
$p_i = (\Sigma, z, z^{S, i}, u) \in \mathring{U}_{\theta}$ and $z^S = \bigcup_i z^{S, i}$.
Note that the assumption on $\mathcal{S}$ implies that for any finite holomorphic
buildings $p^i = (\Sigma, z, z^{S, i}, u) \in \mathring{U}_{\theta}$ such that
$\forget_{\mathcal{S}}(p^i) = \forget_{\mathcal{S}}(p^1)$,
we can define a holomorphic building $p = (\Sigma, z, z^S, u)$ by
$z^S = \bigcup_i z^{S, i}$.
(Namely, $z^S$ and $z^{S'}$ are disjoint for $S \neq S'$.)
$U_{\theta}$ satisfies the following conditions:
\begin{enumerate}[label=(\roman*)']
\item
$\forget_{\mathcal{S}}|_{U_{\theta}} : U_{\theta}
\to \widehat{\M}^{\leq L_{\max}}_{\theta}$ is surjective.
\item
\label{U DM rep}
$\forget_u|_{U_{\theta}} : U_{\theta} \to \overline{\M}^\DM_{\mathcal{S}, \theta}$
is injective, and $\Aut(\forget_u(p)) = \Aut(p) = \Aut(\forget_{\mathcal{S}}(p))$ for all $p \in U_{\theta}$.
Furthermore, for any $p, q \in U_{\theta}$, if there exists
an $(\mathcal{S}, A)$-forgetful map from $\forget_u(p)$ to $\forget_u(q)$, then
$p \geq q$.
\item
For any $p, q \in U_{\theta}$ such that $\forget_{\mathcal{S}}(p)
= \forget_{\mathcal{S}}(q)$, there exists some $r \in U_{\theta}$ such that
$\forget_{\mathcal{S}}(r) = \forget_{\mathcal{S}}(p)$, $r \geq p$ and $r \geq q$.
(For $p = (\Sigma, z, z^S, u)$ and $q = (\Sigma, z, (z')^S, u)$,
$r = (\Sigma, z, z^S \cup (z')^S, u)$ satisfies these conditions.)
\end{enumerate}
The last condition implies that we can apply Lemma \ref{totally ordered cover} for
$\forget_{\mathcal{S}}|_{U_{\theta}} : U_{\theta}
\to \widehat{\M}^{\leq L_{\max}}_{\theta}$.
It implies that there exists some open subset $V_{\theta} \subset U_{\theta}$
which satisfies the following conditions:
$\forget_{\mathcal{S}}(V_{\theta}) = \widehat{\M}^{\leq L_{\max}}_{\theta}$, and
if $p, q \in V_{\theta}$ satisfy $\forget_{\mathcal{S}}(p) = \forget_{\mathcal{S}}(q)$,
then $p \leq q$ or $p \geq q$.
Take open subsets $\V_{\theta} \Subset \U_{\theta} \Subset V_{\theta}$ such that
$\forget_{\mathcal{S}}(\V_{\theta}) = \widehat{\M}^{\leq L_{\max}}_{\theta}$.
Define $\V_{\theta, l} = \V_{\theta} \cap
\widehat{\M}^{\leq L_{\max}}_{\mathcal{S}, \theta, l}$ and $\U_{\theta, l}
= \U_{\theta} \cap \widehat{\M}^{\leq L_{\max}}_{\mathcal{S}, \theta, l}$.
It is clear that they satisfy Condition \ref{existence of minimum}.

We construct open neighborhoods
$\U_{\theta, l}^\DM \subset \overline{\M}^{\DM}_{\mathcal{S}, \theta, l}$
of $\forget_u(\U_{\theta, l})$,
which satisfy Condition \ref{D-neighborhood} and the following condition stronger than
\ref{stably unique forgetful map DM} as follows.
(In this case, $z^A = \emptyset$.)
\begin{itemize}
\item[\ref{stably unique forgetful map DM}$^{\mathrm{S}}$]
For any 
$l \geq 0$, $\hat p = (\hat \Sigma, z, z^S, z^A) \in \U_{\theta, l}^\DM$ and
any subsets $(z^i)^S \subset z^S$ $(i = 1,2$, $S \in \mathcal{S})$ and
$(z^i)^A \subset z^A$ ($i=1,2$), if each $\hat p^i = (\hat \Sigma, z, (z^i)^S, (z^i)^A)$ is
contained in $\U_{\theta, l(\hat p^i)}^\DM$, then
$(z^1)^S \subset (z^2)^S$ and $(z^1)^A \subset (z^2)^A$, or
$(z^2)^S \subset (z^1)^S$ and $(z^2)^A \subset (z^1)^A$.
\end{itemize}
If we replace $\U_{\theta, l}^\DM$ in this condition with $\forget_u(V_{\theta, l})$,
then it holds.
Indeed, for any $p = (\Sigma, z, z^S, u) \in V_{\theta}$ and subsets
$(z^i)^S \subset z^S$, if $\forget_u((\Sigma, z, (z^i)^S, u)) = \forget_u(p^i)$ for some
$p^i \in V_{\theta} \subset U_{\theta}$, then $p^i \leq p$ by
Condition \ref{U DM rep}, which implies that $p^i = (\Sigma, z, (z^i)^S, u)$.
Since $\forget_{\mathcal{S}}(p^i) = \forget_{\mathcal{S}}(p)$ for $i = 1,2$,
$p^1 \leq p^2$ or $p^2 \leq p^1$ by the property of $V_{\theta}$,
which implies that $(z^1)^S \subset (z^2)^S$ or $(z^2)^S \subset (z^1)^S$.
Hence the above condition holds for $\forget_u(V_{\theta, l})$.
Therefore if $\U_{\theta, l}^\DM \subset \overline{\M}^{\DM}_{\mathcal{S}, \theta, l}$
is a sufficiently small neighborhood of $\forget_u(\U_{\theta, l})$, then
they satisfy the condition.
It is easy to check that the subsets $\V_{\theta, l}$, $\U_{\theta, l}$ and
$\U_{\theta, l}^\DM$ satisfy Condition \ref{(S, A)-covering} to
\ref{stably unique forgetful map DM}.

Next we consider the case of $e^2 = 0$.
In this case, for any $p = (\Sigma, z, u) \in \widehat{\M}^{\leq L_{\max}}_{\theta}$,
the domain curve $(\Sigma, z)$ is already stable.
We take finite points $\hat p^i \in \overline{\M}^\DM_{\theta}$,
a local universal family $(\hat P^i \to \hat X^i, Z^i)$ of each $\hat p^i$
and an $\Aut(\hat p^i)$-invariant family of disjoint smooth sections
$(Z^i)^A = ((Z^i)^A_j)$ of each $\hat P^i \to \hat X^i$ which satisfy the following
conditions:
\begin{itemize}
\item
$\overline{\M}^\DM_{\theta}$ is covered by $\hat X^i$, that is, every stable curve in
$\overline{\M}^\DM_{\theta}$ appears some fiber of the local universal families.
\item
For any $i \neq i'$, $a \in \hat X^i$, $a' \in \hat X^{i'}$ and isomorphism
$\varphi : (\hat P^i_a, Z^i(a)) \cong (\hat P^{i'}_{a'}, Z^{i'}(a'))$, $\varphi((Z^i)^A(a))$ and
$(Z^{i'})^A(a')$ are disjoint.
\end{itemize}
For the construction of the family of smooth sections,
we note that if we take a smooth section $(Z^i)^A_1$ of $\hat P^i \to \hat X^i$
whose values are contained in the open subset of the points of trivial stabilizer,
then the $\Aut(\hat p^i)$-orbit of $(Z^i)^A_1$ is an $\Aut(\hat p^i)$-invariant family of
disjoint smooth sections.

Let $U^\DM \subset \overline{\M}^\DM_{\mathcal{S}, A, \theta}$ be the set of
stable curves $(\hat \Sigma, z, z^A)$ ($z^S = \emptyset$) such that there exist
finite points $a_k \in \hat X^{i_k}$ ($k = 1, \dots, N$) and isomorphisms
$\varphi_k : (\hat P^{i_k}_{a_k}, Z^{i_k}(a_k)) \to (\hat \Sigma, z)$ such that
$z^A = \bigcup_k \varphi_k((Z^{i_k})^A(a_k))$.
Applying Lemma \ref{totally ordered cover} to the forgetful map $\forget_A$
from $U^\DM$ to $\overline{\M}^\DM_{\theta}$, we obtain an open subset
$V^\DM \subset U^\DM$ which satisfies the following conditions:
$\forget_A(V^\DM) = \overline{\M}^\DM_{\theta}$, and
if $\hat p, \hat q \in V^\DM$ satisfy $\forget_A(\hat p) = \forget_A(\hat q)$,
then $\hat p \leq \hat q$ or $\hat p \geq \hat q$.
Take an open subset $\U_{\theta}^\DM \subset V^\DM$ such that
$\forget_{\mathcal{S}, A}(\U_{\theta}^\DM) = \overline{\M}^\DM_{\theta}$ and
define $\U_{\theta, l}^\DM = \U_{\theta}^\DM \cap \overline{\M}^\DM_{\theta, l}$.
Define subsets $\U_{\theta, l} \subset
\widehat{\M}^{\leq L_{\max}}_{\mathcal{S}, A, \theta, l}$ by
$\U_{\theta, l} = \forget_u^{-1}(\U_{\theta, l}^\DM)$, and
take their relatively compact open subsets $\V_{\theta, l} \Subset \U_{\theta, l}$
whose images by $\forget_A$ cover $\widehat{\M}^{\leq L_{\max}}_{\theta}$.
Then they satisfy Condition \ref{(S, A)-covering} to \ref{stably unique forgetful map DM}.

%
For the construction of $\U^\DM_{\theta, l}$ for the general triples $\theta$,
we add the following condition:
\begin{itemize}
\item[$(\natural)$]
For any $\hat p \in \U_{\theta, l}^\DM$ and
$\hat q \in \U_{\theta', l'}^\DM$,
if there exists an $(\mathcal{S}, A)$-forgetful map from $\hat p$ to $\hat q$,
then $E_{\hat \omega} \geq E_{\hat \omega}'$,
where $\theta = (g, k, E_{\hat \omega})$ and $\theta' = (g, k, E_{\hat \omega}')$.
\end{itemize}
We will construct $\U^\DM_{\theta, l}$ so that
they also satisfy this condition.

We consider the general triple $\theta = (g, k, E_{\hat \omega})$.
Assume that $\V_{\theta', l}$, $\U_{\theta', l}$ and $\U_{\theta', l}^\DM$ for $\theta'$
such that $\widetilde{e}(\theta') < \widetilde{e}(\theta)$ are given.
Define $\U_{\theta, l}^{\DM, \triangle} \subset
\overline{\M}_{\mathcal{S}, A, \theta, l}^{\DM, \triangle}$ by the largest subset
which satisfy Condition \ref{decomposition into parts U DM}.
Namely, $\hat p \in \overline{\M}^\DM_{\mathcal{S}, A, \theta, l}$ is contained in
$\U_{\theta, l}^{\DM, \triangle}$ if the following condition holds:
For any set $\mathcal{N}$ of its nodal
points, replace each nodal point in $\mathcal{N}$ with a pair of marked points, and
let $\hat p'_i$ $(1 \leq i \leq N)$ be its connected components or an arbitrary
decomposition into unions of its connected components.
Let $g_i$ and $k_i$ be the genus and the number of marked points of each $\hat p'_i$.
Then there exist some $E_{\hat \omega}^i \geq 0$ such that
$E_{\hat \omega} = \sum_i E_{\hat \omega}^i$ and
$\hat p'_i \in \U_{\theta'_i, l(\hat p'_i)}^\DM$ for all $i$, where
$\theta'_i = (g_i, k_i, E_{\hat \omega}^i)$.
Similarly, we define
$\U_{\theta, l}^\triangle \subset \widehat{\M}_{\mathcal{S}, A, \theta, l}^\triangle$
by the largest subset which satisfy Condition \ref{decomposition into parts U},
and define
$\V_{\theta, l}^\triangle \subset \widehat{\M}_{\mathcal{S}, A, \theta, l}^\triangle$
by Condition \ref{decomposition into parts V}.
Then they satisfy the assumption of Lemma \ref{extension from the decomposable}.

We check that $\U_{\theta, l}^{\DM, \triangle}$ ($l \geq 0$) satisfy Condition
\ref{stably unique forgetful map DM}.
For a stable curve $\hat p = (\hat \Sigma, z, z^S, z^A) \in \U^\DM_{\theta, l}$ and
subsets $(z^1)^S, (z^2)^S \subset z^S$ ($S \in \mathcal{S}$) and
$(z^1)^A, (z^2)^A \subset z^A$, assume that each
$\hat p^j = (\hat \Sigma, z, (z^j)^S, (z^j)^A)$ is contained in
$\U^\DM_{\theta, l(\hat p^j)}$.
We prove that
$\hat p^3 = (\hat \Sigma, z, (z^1)^S \cap (z^2)^S, (z^1)^A \cap (z^2)^A)$
is contained in $\U^\DM_{\theta, l(\hat p^3)}$.
Let $\mathcal{N}$ be an arbitrary set of nodal points of $\hat \Sigma$.
Replace each nodal point of $\hat \Sigma$ in $\mathcal{N}$
with a pair of marked points, and decompose the curve into
arbitrary unions of connected components.
For $\hat p$ and $\hat p^j$ ($j=1,2,3$),
let $\{\hat p'_i\}_{1 \leq i \leq k}$ and $\{(\hat p^j)'_i\}_{1 \leq i \leq k}$ be
the obtained decomposition respectively.
Let $(E_{\hat \omega})_i \geq 0$ be non-negative numbers such that
$E_{\hat \omega} = \sum_i (E_{\hat \omega})_i$ and
$\hat p'_i \in \U_{\theta'_i l(\hat p'_i)}^\DM$ for all $i$, where
$\theta'_i = (g'_i, k'_i, (E_{\hat \omega})_i)$.
Similarly, let $(E_{\hat \omega})^j_i$ ($j =1,2$) be pairs for $(\hat p^j)'_i$.
Since there exists an $(\mathcal{S}, A)$-forgetful map from $\hat p'_i$ to
$(\hat p^j)'_i$, Condition $(\natural)$ implies that
$(E_{\hat \omega})_i \geq (E_{\hat \omega})^j_i$ for all $i$.
Hence $(E_{\hat \omega})_i = (E_{\hat \omega})^j_i$ for all $i$.
Therefore, Condition \ref{stably unique forgetful map DM} for each $\theta'_i$
implies that each $(\hat p^3)'_i$ is contained in $\U^\DM_{\theta'_i, l^3_i}$,
where $l^3_i = l((\hat p^3)'_i)$.
Therefore $\hat p^3$ is contained in $\U^\DM_{\theta, l(\hat p^3)}$.

The other conditions in the assumption of Lemma
\ref{extension from the decomposable}
are easy to check.
Hence there exist subsets
\[
\V_{\theta, l}^1 \subset \U_{\theta, l}^1 \subset \U_{\theta, l}^2 \subset
\widehat{\M}^{\leq L_{\max}}_{\mathcal{S}, A, \theta, l} \quad
(l \geq 0)
\]
and
\[
\U_{\theta, l}^{1, \DM} \subset \U_{\theta, l}^{2, \DM} \subset
\overline{\M}^\DM_{\mathcal{S}, A, \theta, l}
\quad (l \geq 0)
\]
which satisfy the conditions in Lemma \ref{extension from the decomposable}.
Since $\U^{\DM, \triangle}_{\theta, l}$ ($l \geq 0$) satisfy
Condition $(\natural)$ (with the other $\U^\DM_{\theta', l'}$),
we may assume that $\U_{\theta, l}^{2, \DM}$ ($l \geq 0$) also satisfy this condition.

We consider separately the cases where $E_{\hat \omega} > 0$ or not.
First we consider the case where $E_{\hat \omega} > 0$.
Since $\widehat{\M}^{\leq L_{\max}}_{\theta} \setminus \forget_{\mathcal{S}, A}
(\bigcup_l \V^1_{\theta, l})$
consists of connected height-one holomorphic buildings without nodal points,
by the same argument as in the case of triples $\theta$
with minimal $\widetilde{e}(\theta)$,
we obtain a finite set $\mathcal{S}^+ = \{S\}$ of codimension-two submanifolds of $Y$
and open subsets $\mathring{U}_{\theta} \subset
\widehat{\M}^{\leq L_{\max}}_{\mathcal{S}^+, \theta}$
which satisfy the following conditions:
\begin{itemize}
\item
$\forget_u|_{\mathring{U}_{\theta}} : \mathring{U}_{\theta} \to
\overline{\M}^\DM_{\mathcal{S}^+, \theta}$ is injective,
and $\Aut(\forget_u(p)) = \Aut(p)$ for any $p \in \mathring{U}_{\theta}$.
Furthermore, for any $p, q \in \mathring{U}_{\theta}$, if there exists
an $(\mathcal{S}^+, A)$-forgetful map from $\forget_u(p)$ to $\forget_u(q)$, then
$p \geq q$.
\item
$\widehat{\M}^{\leq L_{\max}}_{\theta} \setminus \forget_{\mathcal{S}, A}
(\bigcup_l \V^1_{\theta, l})$ is covered by the image of
$\mathring{U}_{\theta}$ by $\forget_{\mathcal{S}^+}$.
\item
For any two holomorphic buildings
$p^i = (\Sigma, z, z^{S, i}, u) \in \mathring{U}_{\theta}$ ($i = 1, 2$) such that
$\forget_{\mathcal{S}^+}(p^1) = \forget_{\mathcal{S}^+}(p^2) = (\Sigma, z, u)$,
the following holds:
\begin{itemize}
\item
For any two different submanifolds $S_1 \neq S_2 \in \mathcal{S}^+$,
$z^{S_1, 1}$ and $z^{S_2, 2}$ are disjoint in $\Sigma$.
\item
For any $S \in \mathcal{S}^+$, if $z^{S, 1} \neq \emptyset$ and
$z^{S, 2} \neq \emptyset$ then $z^{S, 1} = z^{S, 2}$.
\end{itemize}
\item
$\mathcal{S}^+$ and $\mathcal{S}$ do not share the same submanifolds of $Y$.
\item
For any $p^1 = (\Sigma, z, z^{S, 1}, u) \in \mathring{U}_{\theta}$ and
$p^2 = (\Sigma, z, z^{S, 2}, z^A, u) \in \U_{\theta, l}^2$ such that
$\forget_{\mathcal{S}^+}(p^1) = \forget_{\mathcal{S}, A}(p^2)$,
$z^{S_1, 1}$, $z^{S_2, 2}$ and $z^A$ are disjoint for any $S_1 \in \mathcal{S}^+$ and
$S_2 \in \mathcal{S}$.
\end{itemize}
We add $\mathcal{S}^+$ to $\mathcal{S}$, and denote the union by $\mathcal{S}$
in what follows.

Let $U_{\theta} \subset \widehat{\M}^{\leq L_{\max}}_{\mathcal{S}, \theta}$
be the set of holomorphic buildings $p = (\Sigma, z, z^S, u)$ such that
there exist some subsets $z^{S, i} \subset z^S$ ($i = 1, \dots, k$) such that
$p_i = (\Sigma, z, z^{S, i}, u) \in \mathring{U}_{\theta}$ and
$z^S = \bigcup_i z^{S, i}$.
We note that
\begin{equation}
\widehat{\M}^{\leq L_{\max}}_{\theta} = \forget_{\mathcal{S}}(U_{\theta})
\cup \forget_{\mathcal{S}, A}(\bigcup_l \V^1_{\theta, l}).
\label{U V^1 cover}
\end{equation}
Define
\[
U^3_{\theta} = U_{\theta} \setminus
\forget_{\mathcal{S}}^{-1}(\forget_{\mathcal{S}, A}(\bigcup_l \overline{\U_{\theta, l}^1}))
\]
and apply Lemma \ref{totally ordered cover} for the locally homeomorphic map
$\forget_{\mathcal{S}}|_{U^3_{\theta}} : U^3_{\theta} \to
\widehat{\M}^{\leq L_{\max}}_{\theta}$
and a compact subset $\widehat{\M}^{\leq L_{\max}}_{\theta} \setminus
\forget_{\mathcal{S}, A}(\bigcup_l \U^2_{\theta, l})$.
Then we obtain an open subset $V^3_{\theta} \Subset U^3_{\theta}$ such that
\[
\widehat{\M}^{\leq L_{\max}}_{\theta} = \forget_{\mathcal{S}}(V^3_{\theta})
\cup \forget_{\mathcal{S}, A}(\bigcup_l \U^2_{\theta, l})
\]
and if
$p, q \in V^3_{\theta}$ satisfy $\forget_{\mathcal{S}}(p) = \forget_{\mathcal{S}}(q)$
then $p \leq q$ or $q \leq p$.
We note that
\[
\forget_{\mathcal{S}}(V^3_{\theta}) \cap
\forget_{\mathcal{S}, A}(\bigcup_l \U_{\theta, l}^1)= \emptyset
\]
by the definition of $U^3_{\theta}$.

We define $U^4_{\theta} \subset U_{\theta}$ by the open subset of
holomorphic buildings $p \in U_{\theta}$ such that $p \geq q$
for any $q \in \overline{V^3_{\theta}}$
such that $\forget_{\mathcal{S}}(p) = \forget_{\mathcal{S}}(q)$.
Note that $\forget_{\mathcal{S}}(U^4_{\theta}) = \forget_{\mathcal{S}}(U_{\theta})$
since $U_{\theta}$ is closed under the union of $z^S$.
Hence (\ref{U V^1 cover}) implies that
\begin{equation}
\widehat{\M}^{\leq L_{\max}}_{\theta} = \forget_{\mathcal{S}}(U^4_{\theta})
\cup \forget_{\mathcal{S}, A}(\bigcup_l \V^1_{\theta, l}).
\label{4, 1 cover}
\end{equation}
We also define
$U^2_{\theta} \subset \widehat{\M}^{\leq L_{\max}}_{\mathcal{S}, A, \theta}$
by the set of holomorphic buildings $p = (\Sigma, z, \ab z^S, \ab z^A, \ab u)$
which satisfy the following conditions:
\begin{itemize}
\item
$z^S$ and $z^A$ are $\Aut(\Sigma, z, u)$-invariant.
\item
There exist subsets $z^{S, i} \subset z^S$ and $z^{A, i} \subset z^A$ ($i = 1, \dots, k$)
such that $p_i = (\Sigma, z, z^{S, i}, z^{A, i}, u) \in \U_{\theta, l(p_i)}^2$,
$z^S = \bigcup_i z^{S, i}$ and $z^A = \bigcup_i z^{A, i}$.
\item
$p \geq q$ for any $q \in \bigcup_l \overline{\U^1_{\theta, l}}$ such that
$\forget_{\mathcal{S}, A}(p) = \forget_{\mathcal{S}, A}(q)$.
\end{itemize}
Then $\forget_{\mathcal{S}, A}(U^2_{\theta})
= \forget_{\mathcal{S}, A}(\bigcup_l \U^2_{\theta, l})$, which implies that
\begin{equation}
\widehat{\M}^{\leq L_{\max}}_{\theta} = \forget_{\mathcal{S}}(V^3_{\theta})
\cup \forget_{\mathcal{S}, A}(U^2_\theta)
\label{3, 2 cover}
\end{equation}

Let $U^{2+4}_{\theta} \subset \widehat{\M}^{\leq L_{\max}}_{\mathcal{S}, A, \theta}$
be the set of holomorphic buildings $p = (\Sigma, z, z^S, z^A, u)$ such that
there exist some holomorphic buildings
$p_2 = (\Sigma, z, z^{S, 2}, z^{A, 2}, u) \in U^2_{\theta}$ and
$p_4 = (\Sigma, z, z^{S, 4}, u) \in U^4_{\theta}$ such that
$\forget_{\mathcal{S}, A}(p_i) = \forget_{\mathcal{S}, A}(p)$ ($i = 2,4$),
$z^S = z^{S, 2} \sqcup z^{S, 4}$ and $z^A = z^{A, 2}$.
Then (\ref{4, 1 cover}) and (\ref{3, 2 cover}) imply that
\[
\forget_{\mathcal{S}, A}(U^{2+4}_{\theta})
= \forget_{\mathcal{S}, A}(U^2_{\theta}) \cap
\forget_{\mathcal{S}}(U^4_{\theta})
\]
covers
\begin{equation}
\widehat{\M}^{\leq L_{\max}}_{\theta} \setminus
\bigl(\forget_{\mathcal{S}, A}\bigl(\bigcup_l \V^1_{\theta, l}\bigr) \cup
\forget_{\mathcal{S}}(V^3_{\theta})\bigr).
\label{V^1 V^3 complement}
\end{equation}
Furthermore, $\forget_u$ is injective on $U^{2+4}_{\theta}$.
We apply Lemma \ref{totally ordered cover} to the locally homeomorphic map
$\forget_{\mathcal{S}, A}|_{U^{2+4}_{\theta}}
: U^{2+4}_{\theta} \to \widehat{\M}^{\leq L_{\max}}_{\theta}$ and
the compact subset (\ref{V^1 V^3 complement}).
Then we obtain an open subset $V^{2+4}_{\theta} \subset U^{2+4}_{\theta}$
such that if $p, q \in V^{2+4}_{\theta}$ satisfy
$\forget_{\mathcal{S}, A}(p) = \forget_{\mathcal{S}, A}(q)$ then $p \leq q$ or $q \leq p$,
and
\[
\forget_{\mathcal{S}, A}\bigl(\bigcup_l \V^1_{\theta, l} \cup V^{2+4}_{\theta}
\cup V^3_{\theta}\bigr) = \widehat{\M}_{\theta}.
\]

Take open subsets
$\V_{\theta}^{2+4} \Subset \U_{\theta}^{2+4} \subset V_{\theta}^{2+4}$
and $\V_{\theta}^3 \Subset \U_{\theta}^3 \subset V_{\theta}^3$ such that
\[
\forget_{\mathcal{S}, A}(\bigcup_l \V^1_{\theta, l} \cup \V^{2+4}_{\theta}
\cup \V^3_{\theta}) = \widehat{\M}_{\theta},
\]
and define
$\V_{\theta, l}^k = \V_{\theta}^k \cap \widehat{\M}_{\mathcal{S}, A, \theta, l}$
and $\U_{\theta, l}^k = \U_{\theta}^k \cap \widehat{\M}_{\mathcal{S}, A, \theta, l}$
for $k \in \{2+4, 3\}$ and $l$.
Then
\[
\V_{\theta, l} = \V_{\theta, l}^1 \cup \V_{\theta, l}^{2+4} \cup \V_{\theta, l}^3
\]
and
\[
\U_{\theta, l} = \U_{\theta, l}^1 \cup \U_{\theta, l}^{2+4} \cup \U_{\theta, l}^3
\]
satisfy Condition \ref{(S, A)-covering}, \ref{open sets and inclusions},
\ref{decomposition into parts U}, \ref{decomposition into parts V} and
\ref{existence of minimum}.
Condition \ref{existence of minimum} is due to the following properties of
$\U_{\theta, l}^k$ ($k \in \{1, 2+4, 3\}$).
\begin{itemize}
\item
Each $\U_{\theta, l}^k$ ($k \in \{1, 2+4, 3\}$) satisfies
Condition \ref{existence of minimum}.
\item
$\forget_{\mathcal{S}, A}(\bigcup_l \U_{\theta, l}^1)$ and
$\forget_{\mathcal{S}, A}(\bigcup_l \U_{\theta, l}^3)$ do not intersect.
\item
For any $p \in \U_{\theta, l}^{2+4}$ and $q \in \U_{\theta, l}^1 \cup \U_{\theta, l}^3$,
if $\forget_{\mathcal{S}, A}(p) = \forget_{\mathcal{S}, A}(q)$, then $p \geq q$.
\end{itemize}

Finally, we construct $\U_{\theta, l}^\DM \subset
\overline{\M}^\DM_{\mathcal{S}, A, \theta, l}$ as follows.
First for each point $p = (\Sigma, z, z^S, z^A, u) \in \U_{\theta, l}^{2+4}$,
there exists a decomposition $z^S = z^{S, 2} \sqcup z^{S, 4}$
such that $p_2 = (\Sigma, z, z^{S, 2}, z^A, u) \in U_{\theta}^2$ and
$p_4 = (\Sigma, z, z^{S, 4}, u) \in U_{\theta}^4$ by definition.
By the definition of $U_{\theta}^2$, there exist some $p_{2, i}
= (\Sigma, z, z^{S, 2, i}, z^{A, 2, i}, u) \in \U_{\theta, l(p_{2, i})}^2$
such that $\forget_{\mathcal{S}, A}(p_{2, i}) = \forget_{\mathcal{S}, A}(p_2)$,
$z^{S, 2} = \bigcup_i z^{S, 2, i}$ and $z^A = \bigcup_i z^{A, 2, i}$.
Condition \ref{Z^A section} for $\U_{\theta, l}^{2, \DM}$ implies that
for each $i$, there exist a local universal family
$(\hat P^i \to \hat X^i, Z^i, (Z^i)^S)$ of $\forget_A(\forget_u(p_{2, i}))$ and
an $\Aut(\forget_u(p_{2, i}))$-invariant family of smooth sections
$(Z^i)^A = ((Z^i)^A_j)$ of $\hat P^i \to \hat X^i$ such that
\[
\{(\hat P^i_a, Z^i(a), (Z^i)^S(a), (Z^i)^A(a)); a \in \hat X^i\} / \Aut(\forget_u(p_{2, i}))
\]
is a neighborhood of $\forget_u(p_{2, i})$ in $\U_{\theta, l}^{2, \DM}$.
Let $(\hat P^p \to \hat X^p, Z^p, (Z^p)^S)$ be a local universal family of
$\forget_A(\forget_u(p))$ and define an $\Aut(\forget_u(p))$-invariant family of
sections $(Z^p)^A$ of $\hat P^p \to \hat X^p$ by the union of the pull backs of
$(Z^i)^A$ by the forgetful maps.
We define $\W_{\theta, l}^{2+4, \DM} \subset
\overline{\M}^\DM_{\mathcal{S}, A, \theta, l}$ by the union of
\[
\{(\hat P^p_a, Z^p(a), (Z^p)^S(a), (Z^p)^A(a)); a \in \hat X^p\} / \Aut(\forget_u(p))
\]
over $p \in \U_{\theta, l}^{2+4}$.
Then $\W_{\theta, l}^{2+4, \DM}$ satisfy Condition \ref{Z^A section} instead of
$\U_{\theta, l}^\DM$.

We construct $\U_{\theta, l}^\DM \subset
\overline{\M}^\DM_{\mathcal{S}, A, \theta, l}$ as the union of
$\U_{\theta, l}^{k, \DM}$ ($k \in \{1, 2+4, 3\}$).
We construct $\U_{\theta, l}^{2+4, \DM}$ ($l \geq 0$) as open neighborhoods of
$\forget_{\mathcal{S}, A}(\U_{\theta, l}^{2+4})$ in $\W_{\theta, l}^{2+4, \DM}$
which satisfy Condition \ref{stably unique forgetful map DM}$^{\mathrm{S}}$ and
which are D-neighborhoods of $\U_{\theta, l}^{2+4}$.
$\U_{\theta, l}^{3, \DM}$ are also constructed as open neighborhoods of
$\forget_{\mathcal{S}}(\U_{\theta, l}^3)$ in $\overline{\M}^\DM_{\mathcal{S}, \theta}$
which satisfy Condition \ref{stably unique forgetful map DM}$^{\mathrm{S}}$ and
which are D-neighborhoods of $\U_{\theta, l}^3$.
($\U_{\theta, l}^1$ ($l \geq 0$) have been already constructed.)
Since $\forget_{\mathcal{S}, A}(\U_{\theta, l}^k)$ ($k = 2+4, 3$) satisfy condition
\ref{stably unique forgetful map DM}$^{\mathrm{S}}$,
sufficiently small open neighborhoods satisfy the condition.
Furthermore, if these open neighborhoods are sufficiently small, then they also satisfy
the following condition:
For $\hat p = (\hat \Sigma, z, z^S, z^A) \in \U_{\theta, l(p)}^{2+4, \DM}$
and any subsets $(z^i)^S \subset z^S$ $(i = 1,2$, $S \in \mathcal{S})$ and
$(z^i)^A \subset z^A$, if
each $\hat p^i = (\hat \Sigma, z, (z^i)^S, (z^i)^A)$ is contained in 
$\U_{\theta, l(\hat p^i)}^{1, \DM}$ or $\U_{\theta, l(\hat p^i)}^{3, \DM}$, then
$(z^1)^S \subset (z^2)^S$ and $(z^1)^A \subset (z^2)^A$, or
$(z^2)^S \subset (z^1)^S$ and $(z^2)^A \subset (z^1)^A$.
Then $\U_{\theta, l}^\DM = \bigcup_{k \in \{1, 2+4, 3\}} \U_{\theta, l}^{k, \DM}$
($l \geq 0$) satisfy Condition \ref{stably unique forgetful map DM}.
It is easy to check that $\U_{\theta, l}^\DM$ ($l \geq 0$) satisfy the other conditions.
Hence we can construct the required subsets $\V_{\theta, l}$, $\U_{\theta, l}$ and
$\U_{\theta, l}^\DM$ by the induction.
\end{proof}

Assume that a domain curve representation
$(\mathcal{S}, \V_{\theta, l}, \U_{\theta, l}, \U^\DM_{\theta, l})$ of
$\widehat{\M}^{\leq L_{\max}}_{\leq C}$ is given.
Next we construct spaces $\X_{\theta}$ of holomorphic buildings
with perturbation parameters.
At the same time, we construct spaces $\V_{\alpha, \theta, l}$,
$\U_{\alpha, \theta, l}$ of holomorphic buildings with perturbation parameters and
sets $\U_{\alpha, \theta, l}^\DM$ of stable curves with perturbation parameters
indexed by a finite index set $\A = \{\alpha\}$.
We call a family $(\X_{\theta}, \V_{\alpha, \theta, l}, \U_{\alpha, \theta, l},
\U^\DM_{\alpha, \theta, l})$ Kurainshi data if it satisfies the following conditions:
\begin{enumerate}[label=$(\arabic*)^{\mathrm{K}}$]
\item
We may regard $\V_{\alpha, \theta, l}$ and $\U_{\alpha, \theta, l}$ as subspaces of
$\widehat{\M}^{\leq L_{\max}}_{o, \mathcal{S}, A, \theta, l}$ by
$(p, E^0_p, \lambda_p) \mapsto p$ for each $\alpha$.
Similarly, we may regard $\U_{\alpha, \theta, l}^\DM$ as a subspace of
$\overline{\M}^\DM_{o, \mathcal{S}, A, \theta, l}$ for each $\alpha$.
\item
\label{D-neighborhood o}
There exists an open neighborhood $\W_{\alpha, \theta, l}
\subset \widehat{\M}_{o, \mathcal{S}, A, \theta, l}$ of the closure of
$\U_{\alpha, \theta, l}$
such that
\[
\U_{\alpha, \theta, l} = \{p \in \W_{\alpha, \theta, l};
\forget_u(p) \in \U_{\alpha, \theta, l}^\DM\}
\]
as a space of holomorphic buildings.
Furthermore, for each $p \in \U_{\alpha, \theta, l}$,
the associated vector space $E^0_p$ and $\lambda_p$ are defined by
$E^0_p = E^0_{\forget_u(p)}$ and $\lambda_p = \lambda_{\forget_u(p)}$.
In this case, we say that $\U_{\alpha, \theta, l}^\DM$ is
a D-neighborhood of $\U_{\alpha, \theta, l}$.
\item
\label{forget o relation}
$\forget_o(\V_{\alpha, \theta, l})$, $\forget_o(\U_{\alpha, \theta, l})$ and
$\forget_o(\U_{\alpha, \theta, l}^\DM)$ are contained in
$\V_{\theta, l}$, $\U_{\theta, l}$ and $\U_{\theta, l}^\DM$ respectively.
\item
$\V_{\alpha, \theta, l}$ is open in the relative topology of $\U_{\alpha, \theta, l}$,
and $\V_{\alpha, \theta, l} \Subset \U_{\alpha, \theta, l}$.
\item
The number of marked points $z^o$ of each holomorphic building
in $\U_{\alpha, \theta, l}$ or stable curve in $\U^\DM_{\alpha, \theta, l}$
depends only on $\alpha$.
\item
\label{Z^o section with lambda}
For each point $\hat p \in \U_{\alpha, \theta, l}^\DM$,
let $(\hat P \to \hat X, Z, Z^S, Z^A)$ be the local representation of a
neighborhood of $\forget_o(\hat p)$ in $\U_{\theta, l}^\DM$.
If we shrink $\hat X$ then there exists an $\Aut(\hat p)$-invariant family of smooth
sections $Z^o = (Z^o_j)$ of $\hat P \to \hat X$ such that
\[
\{(\hat P_a, Z(a), Z^S(a), Z^A(a), Z^o(a)); a \in \hat X\} / \Aut(\hat p)
\]
is a neighborhood of $\hat p$ in $\U_{\alpha, \theta, l}^{\DM}$.
Furthermore,
there exists an $\Aut(\hat p)$-equivariant linear map
$\tilde \lambda_{\hat p} : E^0_{\hat p} \to
C^\infty(\hat P \times Y, \Wedge^{0, 1} V^\ast \hat P
\otimes_\C (\R \partial_\sigma \oplus TY))$ which satisfy the following conditions:
\begin{itemize}
\item
For each $h \in E^o_p$, the projection of the support of $\tilde \lambda_{\hat p}(h)$ to
$\hat P$ does not intersect with the nodal points or marked points $Z$.
\item
For any $a \in \hat X$, $\hat q \in \U_{\alpha, \theta, l}^\DM$ and
isomorphism $f : (\hat P_a, Z(a),\ab Z^S(a),\ab Z^A(a), Z^o(a)) \to \hat q$,
there exists an isomorphism $\hat \phi_f : E^0_{\hat p} \to E^0_{\hat q}$ such that
the restriction of $\tilde \lambda_{\hat p}$ to $\hat P_a \times Y$
coincides with $f^\ast \circ \lambda_{\hat q} \circ \hat \phi_f$.
\end{itemize}
We call $(\hat P \to \hat X, Z, Z^S, Z^A, Z^o, E^0_{\hat p}, \tilde \lambda_{\hat p})$
a local representation of a neighborhood of $\hat p$ in $\U_{\alpha, \theta, l}^\DM$.
\item
Any $(\mathcal{S}, A)$-forgetful map from
$\hat p \in \U_{\alpha, \theta, l}^\DM$ to
$\hat q \in \U_{\alpha, \theta, l'}^\DM$
does not collapse any component of $\hat p$.
(This condition follows from Condition \ref{forget o relation} and
Condition \ref{no collapse for (S, A)-forgetful map} of domain curve representation.)
\item
For any $(\mathcal{S}, A)$-forgetful map $f$ from
$\hat p \in \U_{\alpha, \theta, l}^\DM$ to $\hat q \in \U_{\alpha, \theta, l'}^\DM$,
there exists an isomorphism $\hat \phi_f : E^0_{\hat p} \cong E^0_{\hat q}$ such that
$\lambda_{\hat p} = f^\ast \circ \lambda_{\hat q} \circ \hat \phi_f$.
Furthermore, for another $(\mathcal{S}, A)$-forgetful map $h$ from $\hat q$ to
$\hat r \in \U^\DM_{\alpha, \theta, l''}$, $\hat \phi_{h \circ f}$ coincides with
the composition of $\hat \phi_{\hat h}$ and $\hat \phi_{\hat f}$.
\item
\label{(Z^o, E, lambda) pull back relation}
For any $l \geq l'$,
$\hat p \in \U_{\alpha, \theta, l}^\DM$ and $\hat q \in \U_{\alpha, \theta, l'}^\DM$,
if there exists an $(\mathcal{S}, A)$-forgetful map $f$ from $\hat p$ to $\hat q$,
then the following condition holds true:
Let $(\hat P \to \hat X, Z, Z^S, Z^A, Z^o, E^0_{\hat p}, \tilde \lambda_{\hat p})$ be
a local representation of a neighborhood of $\forget_A(\hat p)$
in $\U_{\alpha, \theta, l}^\DM$, and
$(\hat P' \to \hat X', Z', (Z')^S, (Z')^A, (Z')^o, E^0_{\hat q}, \tilde \lambda_{\hat q})$ be
that for $\hat q$.
Shrink $\hat X$ and $\hat X'$ if necessary, and
let $(\phi, \hat \phi)$ be the unique forgetful map
from $(\hat P \to \hat X, Z, Z^S)$ to $(\hat P' \to \hat X', Z', (Z')^S)$
whose restriction to the central fiber coincides with $f$.
Then the pull back of $(Z')^o$ by $(\phi, \hat \phi)$ coincides with $Z^o$,
and $\tilde \lambda_{\hat p}$ coincides with
the pull back of $\tilde \lambda_{\hat q}$ by $(\phi, \hat \phi)$
under the identification $\hat \phi_f : E^0_{\hat p} \cong E^0_{\hat q}$.
\item
\label{decomposition into parts U DM o}
For any $(\hat p, E^0_{\hat p}, \lambda_{\hat p}) \in \U_{\alpha, \theta, l}^\DM$
and any subset $\mathcal{N}$ of the nodal points of $\hat p$,
replace each nodal point in $\mathcal{N}$ with a pair of marked points
(we regard the new marked points as points in the set $z$), and
let $\hat p'_i$ $(1 \leq i \leq k)$ be its connected components or an arbitrary
decomposition into unions of its connected components.
Let $g'_i$ and $k'_i$ be the genus and the number of marked points $z$ of each
$\hat p'_i$ respectively.
Then there exist some $E_{\hat \omega}^i \geq 0$ such that
$E_{\hat \omega} = \sum_i E_{\hat \omega}^i$ and the following holds:
Only one of $\hat p'_i$ contains marked points $z^o$,
the support of $\lambda_{\hat p}(h)$ is contained in this component
for all $h \in E^0_{\hat p}$, and
$(\hat p'_i, E^0_{\hat p}, \lambda_{\hat p})$ is contained in
$\U_{\alpha, \theta'_i, l(\hat p'_i)}^\DM$, where $\theta'_i = (g'_i, k'_i, E_{\hat \omega}^i)$.
Furthermore, the other $\hat p'_i$ are contained in
$\U_{\theta'_i, l(\hat p'_i)}^\DM$.
\item
\label{decomposition into parts U o}
$\U_{\alpha, \theta, l}$ satisfy the following conditions
about decomposition into parts:
\begin{itemize}
\item
For any $p \in \U_{\alpha, \theta, l}$ and any decomposition $p_i$ ($1 \leq k$)
into unions of its connected components, let $p'_i$ be the holomorphic buildings
obtained by collapsing trivial floors (floors consisting of trivial cylinders).
Then only one of $p'_i$ contains marked points $z^o$, and it is contained in
$\U_{\alpha, \theta(p'_i), l(p'_i)}$.
Furthermore, the others are contained in $\U_{\theta(p'_i), l(p'_i)}$.
\item
For any $p \in \U_{\alpha, \theta, l}$ and any gap between floors, let $p_1$ and $p_2$
be the holomorphic buildings obtained by dividing $p$ at this gap.
Then one of $p'_i$ $(i = 1,2)$ is contained in $\U_{\alpha, \theta(p'_i), l(p'_i)}$ and
the other is contained in $\U_{\theta(p'_i), l(p'_i)}$.
\item
For any $p \in \U_{\alpha, \theta, l}$ and any subset of its nodal points,
the holomorphic building $p'$ obtained by replacing these nodal points to pairs of
marked points is contained in $\U_{\alpha, \theta(p'), l(p')}$.
\end{itemize}
\item
\label{decomposition into parts V o}
For each $p \in \widehat{\M}^{\leq L_{\max}}_{o, \mathcal{S}, A, \theta, l}$,
replace all nodal points and joint circles of $p$ to pairs of marked points and
pairs of limit circles respectively
(we regard the new marked points as points in the set $z$), and
let $p'_i$ $(1 \leq i \leq k)$ be the stabilizations of its non-trivial
connected components.
Then $p \in \V_{\alpha, \theta, l}$ if and only if one of $p'_i$ $(1 \leq i \leq k)$ is
contained in $\V_{\alpha, \theta(p'_i), l(p'_i)}$ and the others are contained in
$\V_{\theta(p'_i), l(p'_i)}$.
\item
\label{existence of minimum o}
For any $\alpha \in \A$, $p = (\Sigma, z, u) \in \widehat{\M}^{\leq L_{\max}}_{\theta}$
and subsets $(z^k)^S, (z^k)^A, z^o \subset \Sigma$ ($S \in \mathcal{S}, k =1,2)$,
if each $p^k = (\Sigma, z, (z^k)^S, (z^k)^A, z^o, u)$ is contained in
$\U_{\alpha, \theta, l(p^k)}$, then
$p^3 = (\Sigma, z, (z^1)^S \cap (z^2)^S, (z^1)^A \cap (z^2)^A, z^o, u)$ is contained
in $\U_{\alpha, \theta, l(p^3)}$.
\item
\label{stably unique forgetful map DM o}
For any $\hat p = (\hat \Sigma, z, z^S, z^A, z^o) \in \U_{\alpha, \theta, l}^\DM$ and
any subsets $(z^1)^S, (z^2)^S \subset z^S$ ($S \in \mathcal{S}$) and
$(z^1)^A, (z^2)^A \subset z^A$, if
each $\hat p^i = (\hat \Sigma, z, (z^i)^S, (z^i)^A, z^o)$ is contained in 
$\U_{\alpha, \theta, l(\hat p^i)}^\DM$, then
$\hat p^3 = (\hat \Sigma, z, (z^1)^S \cap (z^2)^S, (z^1)^A \cap (z^2)^A, z^o)$
is also contained in $\U_{\alpha, \theta, l(\hat p^3)}^\DM$.
\item
\label{z^o disjoint}
For any $p = (\Sigma, z, u) \in \widehat{\M}^{\leq L_{\max}}_{\theta}$
and subsets $(z^1)^S, (z^2)^S \subset \Sigma$ ($S \in \mathcal{S}$),
$(z^1)^A, (z^2)^A \subset \Sigma$
and $(z^1)^o, (z^2)^o \subset \Sigma$, if $(z^1)^o \cap (z^2)^o \neq \emptyset$ and
each $p^k = (\Sigma, z, (z^k)^S, (z^k)^A, (z^k)^o, u)$ is contained in
$\U_{\alpha^k, \theta, l(p^k)}$ for some $\alpha^k \in \A$, then
$\alpha^1 = \alpha^2$ and $(z^1)^o = (z^2)^o$.
\item
\label{def of X}
Each $\X_{\theta}$ is determined by $(\V_{\alpha, \theta, l})_{\alpha \in \A, l \geq 0}$
as a subset of $\widehat{\M}^{\leq L_{\max}}_{o, \mathcal{S}, A, \theta}$ as follows.
$p = (\Sigma, z, z^S, z^A, z^o, u) \in
\widehat{\M}^{\leq L_{\max}}_{o, \mathcal{S}, A, \theta}$ is contained in $\X_{\theta}$
if it satisfies the following conditions:
\begin{enumerate}[label=(\alph*)]
\item
$z^S$, $z^A$ and $z^o$ are $\Aut(\forget_{o, \mathcal{S}, A}(p))$-invariant as sets.
Furthermore,
replace all nodal points and joint circles of $p$ to pairs of marked points and
pairs of limit circles respectively
(we regard the new marked points as points in the set $z$), and
let $p'_i$ $(1 \leq i \leq k)$ be the stabilizations of its non-trivial
connected components as in Condition \ref{decomposition into parts V o}.
Then $z^S$, $z^A$ and $z^o$ are $\Aut((\forget_{o, \mathcal{S}, A}(p'_i))_i)$-invariant.
(The latter stronger condition is related to Remark \ref{quotient map is a submersion}
in Section \ref{Product of pre-Kuranishi spaces}.
We need this condition to construct a map corresponding to
the decomposition by gaps of floors,)
\item
\label{(o, S, A)-union}
There exist subsets $z^{S, i} \subset z^S$, $z^{A, i} \subset z^A$, $z^{o, i} \subset z^o$
and indices $\alpha_i \in \A$ $(i = 1, \dots, k)$ such that
$p_i = (\Sigma, z, z^{S, i}, z^{A, i}, z^{o, i}, u) \in \V_{\alpha_i, \theta, l(p_i)}$ for all $i$,
$z^S = \bigcup_i z^{S, i}$, $z^A = \bigcup_i z^{A, i}$ and $z^o = \bigcup_i z^{o, i}$.
\item
The linear map (\ref{X surjective map}) is surjective for the vector space
$E^0_p$ and linear map $\lambda_p$ defined in the next condition.
\end{enumerate}
\item
\label{def of (E, lambda) of X}
For each $p = (\Sigma, z, z^S, z^A, z^o, u) \in \X_{\theta}$,
the associated vector space $E^0_p$ and linear map $\lambda_p$ are defined as follows.
First note that in \ref{(o, S, A)-union} of the above condition,
if $z^{o, i} \cap z^{o, i'} \neq \emptyset$, then $\alpha_i = \alpha_{i'}$ and
$z^{o, i} = z^{o, i'}$ by Condition \ref{z^o disjoint}.
Choose a subset $I \subset \{1, \dots, k\}$ such that $z^o = \coprod_{i \in I} z^{o, i}$,
and fix forgetful maps from $\hat p = \forget_u(p)$ to $\hat p_i = \forget_u(p_i)$
for each $i \in I$.
Then $E^0_p$ is the direct sum of $E^0_{\hat p_i}$, and $\lambda_p$ is the sum of
the pull backs of $\lambda_{\hat p_i}$ by the forgetful map $\hat p \to \hat p_i$.
This definition is independent of the choice of $p_i$ and $I$ by
Condition \ref{(Z^o, E, lambda) pull back relation} and \ref{existence of minimum o}.
\item
\label{(o, S, A)-covering}
For each triple $\theta$, the subspace
$\X_{\theta} \subset \widehat{\M}^{\leq L_{\max}}_{o, \mathcal{S}, A, \theta}$
defined by Condition \ref{def of X} satisfies
$\forget_{o, \mathcal{S}, A}(\X_{\theta}) = \widehat{\M}^{\leq L_{\max}}_{\theta}$.
\end{enumerate}

We can prove the following lemmas similarly to Lemma \ref{shrinking and conditions}
and \ref{extension from the decomposable} respectively.
\begin{lem}
\label{shrinking and conditions o}
Let $(e^1_0, e^2_0)$ be an arbitrary constant, and assume that
Kurainshi data $(\X_{\theta}, \V_{\alpha, \theta, l}, \U_{\alpha, \theta, l},
\U^\DM_{\alpha, \theta, l})$ for $\widehat{\M}_{\leq C}$ are given.
Then we can construct open subsets
\[
\V_{\alpha, e^1_0, e^2_0, l} \Subset \mathring{\U}_{\alpha, e^1_0, e^2_0, l} \Subset
\U_{\alpha, e^1_0, e^2_0, l} \quad
(l \geq 0)
\]
and
\[
\mathring{\U}_{\alpha, e^1_0, e^2_0, l}^\DM \Subset \U_{\alpha, e^1_0, e^2_0, l}^\DM
\quad (l \geq 0)
\]
such that if we replace $\U_{\alpha, e^1_0, e^2_0, l}$ and
$\U_{\alpha, e^1_0, e^2_0, l}^\DM$ in the family
$(\X_{\theta}, \V_{\alpha, \theta, l},\ab \U_{\alpha, \theta, l},\ab
\U^\DM_{\alpha, \theta, l})$ with $\mathring{\U}_{\alpha, e^1_0, e^2_0, l}$ and
$\mathring{\U}_{\alpha, e^1_0, e^2_0, l}^\DM$ respectively,
it still satisfies the conditions of Kuranishi data.
\end{lem}
\begin{proof}
Since only nontrivial conditions are Condition \ref{existence of minimum o} and
\ref{stably unique forgetful map DM o}, we can prove the claim similarly to
Lemma \ref{shrinking and conditions}.
\end{proof}

\begin{lem}
\label{extension from the decomposable o}
Assume that Kurainshi data $(\X_{\theta}, \V_{\alpha, \theta, l},
\U_{\alpha, \theta, l}, \U^\DM_{\alpha, \theta, l})$ of
$\widehat{\M}_{< C}$ are given.
We also assume that spaces
\[
\V_{\alpha, \theta, l}^\triangle \subset \U_{\alpha, \theta, l}^\triangle
\quad (\alpha \in \A, \widetilde{e}(\theta) = C, l \geq 0)
\]
of holomorphic buildings with perturbation parameters and spaces
\[
\U_{\alpha, \theta, l}^{\DM, \triangle} \quad (\alpha \in \A,
\widetilde{e}(\theta) = C, l \geq 0)
\]
of stable curves with perturbation parameters are given and
they satisfy the conditions of Kuranishi data if we replace
$\widehat{\M}$ and $\overline{\M}^\DM$ with $\widehat{\M}^\triangle$ and
$\overline{\M}^{\DM, \triangle}$ respectively.
Then we can construct spaces
\[
\V_{\alpha, \theta, l} \Subset \U_{\alpha, \theta, l}
\quad (\alpha \in \A, \widetilde{e}(\theta) = C, l \geq 0)
\]
of holomorphic buildings with perturbation parameters and spaces
\[
\U_{\alpha, \theta, l}^\DM
\quad (\alpha \in \A, \widetilde{e}(\theta) = C, l \geq 0)
\]
of stable curves with perturbation parameters which satisfy
$\V_{\alpha, \theta, l} \cap \widehat{\M}_{o, \mathcal{S}, A}^\triangle
= \V_{\alpha, \theta, l}^\triangle$,
$\U_{\alpha, \theta, l} \cap \widehat{\M}_{o, \mathcal{S}, A}^\triangle \subset
\U_{\alpha, \theta, l}^\triangle$,
$\U_{\alpha, \theta, l}^\DM \cap \overline{\M}_{o, \mathcal{S}, A}^{\DM, \triangle}
\subset \U_{\alpha, \theta}^{\DM, \triangle}$ and
the conditions of Kuranishi data other than Condition \ref{(o, S, A)-covering}.
\end{lem}
\begin{proof}
For each $\alpha \in \A$, we construct the extensions of $Z^o$
as in the case of $Z^A$ in Lemma \ref{extension from the decomposable}.
In this case, we also construct the extensions of $\lambda$
at the same time by the same induction.
Their construction is also similar.
\end{proof}

Now we explain the construction of Kurainshi data.
\begin{lem}\label{good family of additional vector spaces}
There exist Kurainshi data
$(\X_{\theta}, \V_{\alpha, \theta, l}, \U_{\alpha, \theta, l},
\U^\DM_{\alpha, \theta, l})$ of $\widehat{\M}^{\leq L_{\max}}_{\leq C}$
for any domain curve representation
$(\mathcal{S}, \V_{\theta, l}, \U_{\theta, l}, \U^\DM_{\theta, l})$.
\end{lem}
\begin{proof}
We construct Kurainshi data by the induction in $\widetilde{e}(\theta)$.
For each triple $\theta$ with minimal $\widetilde{e}(\theta)$,
we take finite open subsets $U_\alpha \subset \U_{\theta, l_\alpha}$ and
$U_\alpha^\DM \subset \U_{\theta, l_\alpha}^\DM$ ($\alpha \in \A$)
such that
\begin{itemize}
\item
each $U_\alpha^\DM$ is a D-neighborhood of $U_\alpha$,
\item
$U_\alpha^\DM$ is covered by a local representation
$(\hat P^\alpha \to \hat X^\alpha, Z^\alpha, (Z^\alpha)^S, (Z^\alpha)^A)$
of a neighborhood a point $\hat p^\alpha$ in $\U^\DM_{\theta, l_\alpha}$
for some $l_\alpha$, and
\item
$\{\forget_{o, \mathcal{S}, A}(U_\alpha)\}_{\alpha \in \A}$
covers $\overline{\M}^{\leq L_{\max}}_{\theta}$.
\end{itemize}
For each $\alpha$, we construct an $\Aut(\hat p^\alpha)$ vector space $E^0_\alpha$
and an $\Aut(\hat p^\alpha)$-equivariant linear map
\[
\lambda_\alpha : E^0_\alpha \to C^\infty(\hat P^\alpha \times Y,
\Wedge^{0, 1} V^\ast P^\alpha \otimes_\C (\R \partial_\sigma \oplus TY))
\]
such that for any $p \in U_\alpha$, $E^0_p = E^0_\alpha$ and
the restriction of $\lambda_\alpha$ to the fiber isomorphic to $\forget_u(p)$ make
the linear map (\ref{X surjective map}) surjective.
We assume that for each $h \in E^o_\alpha$,
the projection of the support of $\tilde \lambda_\alpha(h)$
to $\hat P^\alpha$ does not intersect with the nodal points or marked points $Z$.
We also construct an $\Aut(\hat p^\alpha)$-invariant family of section $(Z^\alpha)^o
= ((Z^\alpha)^o_j)$ of $\hat P^\alpha \to \hat X^\alpha$ for each $\alpha \in \A$.
Then we define $\U_{\alpha, \theta}^\DM = \U_{\alpha, \theta, l_\alpha}^\DM$
by the space of stable curves $\hat p \in \forget_o^{-1}(U^\DM_\alpha)$ such that
$\hat p$ is isomorphic to $(\hat P^\alpha_a, Z^\alpha(a), (Z^\alpha)^S(a),
(Z^\alpha)^A(a), (Z^\alpha)^o(a))$
for the point $a \in \hat X^\alpha$ such that
$(\hat P^\alpha_a, Z^\alpha(a), (Z^\alpha)^S(a), (Z^\alpha)^A(a))$ is isomorphic
to $\forget_o(\hat p)$.
For such a stable curve $\hat p$, we define $E^0_{\hat p}$ and
$\lambda_{\hat p}$ by $E^0_{\hat p} = E^0_\alpha$ and the restriction of
$\lambda_\alpha$ respectively, and regard $\U^\DM_{\alpha, \theta, l_\alpha}$
as a space of stable curve with perturbation parameters.
Similarly, we define the space of holomorphic buildings $\U_{\alpha, \theta, l_\alpha}$
by $\U_{\alpha, \theta, l_\alpha} = \forget_o^{-1}(U_\alpha) \cap
\forget_u^{-1}(\U_{\alpha, \theta, l_\alpha}^\DM)$,
and for each $p \in \U_{\alpha, \theta, l_\alpha}$, 
For $l \neq l_\alpha$, we define $\U^\DM_{\alpha, \theta, l} = \emptyset$ and
$\U_{\alpha, \theta, l} = \emptyset$.
It is clear that we can choose the family of sections $(Z^\alpha)^o$ so that
Condition \ref{z^o disjoint} holds.
We take open subsets $\V_{\alpha, \theta, l} \Subset \U_{\alpha, \theta, l}$
such that $\{\forget_{o, \mathcal{S}, A}(V_{\alpha, \theta, l_\alpha})\}_{\alpha \in \A}$
covers $\widehat{\M}^{\leq L_{\max}}_{\theta}$,
and define $\X_{\theta}$ by Condition \ref{def of X} and \ref{def of (E, lambda) of X}.

Next we consider the general triple $\theta$.
We assume that $\V_{\alpha, \theta', l}$, $\U_{\alpha, \theta', l}$,
$\U^\DM_{\alpha, \theta', l}$ and $\X_{\theta'}$
for $\widetilde{e}(\theta') < \widetilde{e}(\theta)$ are already constructed and
construct those for $\theta$.
Define $\U^{\DM, \triangle}_{\alpha, \theta, l} \subset
\overline{\M}^{\DM, \triangle}_{o, \mathcal{S}, A, \theta, l}$ by the largest space
which satisfies Condition \ref{decomposition into parts U DM o},
$\U^\triangle_{\alpha, \theta, l} \subset
\widehat{\M}^{\leq L_{\max}, \triangle}_{o, \mathcal{S}, A, \theta, l}$
by the largest space which satisfies Condition \ref{decomposition into parts U o}, and
$\V^\triangle_{\alpha, \theta, l}$ by Condition \ref{decomposition into parts V o}.
Then they satisfy the assumption of Lemma \ref{extension from the decomposable o}.
Hence we obtain spaces
\[
\V_{\alpha, \theta, l} \Subset \U_{\alpha, \theta, l}
\quad (\alpha \in \A, l \geq 0)
\]
of holomorphic buildings with perturbation parameters and spaces
\[
\U_{\alpha, \theta, l}^\DM
\quad (\alpha \in \A, l \geq 0)
\]
of stable curves with perturbation parameters which satisfy the conclusion of
Lemma \ref{extension from the decomposable o}.

Define $\X^1_{\theta}$ for these spaces
$(\V_{\alpha, \theta, l})_{\alpha \in \A, l \geq 0}$ by Condition \ref{def of X}
and \ref{def of (E, lambda) of X}.
Then its image by $\forget_{o, \mathcal{S}, A}$ contains a neighborhood of
$\widehat{\M}^{\leq L_{\max}, \triangle}_{\theta}$.
For the complement $\widehat{\M}^{\leq L_{\max}}_{\theta} \setminus
\forget_{o, \mathcal{S}, A}(\X^1_{\theta})$,
we use the same argument as in the case of minimal $\widetilde{e}(\theta)$.
Namely, we construct spaces $\U_{\alpha', \theta, l}$ of holomorphic buildings
with perturbation parameters and
spaces $\U^\DM_{\alpha', \theta, l}$ of stable curves with perturbation parameters
indexed by another finite index set $\A' = \{\alpha'\}$ which satisfy the following
conditions:
\begin{itemize}
\item
Each $\U_{\alpha', \theta, l}^\DM$ is a D-neighborhood of
$\U_{\alpha', \theta, l}$.
\item
For each $(p = (\Sigma, z, z^S, z^A, z^o, u), E^0_p, \lambda_p) \in
\U_{\alpha', \theta, l}$, $z^S$, $z^A$ and $z^o$ are $\Aut(\Sigma, z, u)$-invariant.
\item
For any $(p, E^0_p, \lambda_p) \in \U_{\alpha', \theta, l}$,
the linear map (\ref{X surjective map}) is surjective.
\item
$\{\forget_{o, \mathcal{S}, A}(\U_{\alpha', \theta, l})\}_{\alpha' \in \A'}$ covers
$\widehat{\M}^{\leq L_{\max}}_{\theta} \setminus
\forget_{o, \mathcal{S}, A}(\X^1_{\theta})$.
\end{itemize}

Take open subsets $\V_{\alpha', \theta, l} \Subset \U_{\alpha', \theta, l}$
such that $\{\forget_{o, \mathcal{S}, A}(\V_{\alpha', \theta, l})\}_{\alpha' \in \A'}$
covers $\widehat{\M}^{\leq L_{\max}}_{\theta} \setminus
\forget_{o, \mathcal{S}, A}(\X^1_{\theta})$, and
define the space $\X_{\theta}$ of holomorphic buildings
with perturbation parameters for
$(\V_{\alpha, \theta, l})_{\alpha \in \A \cup \A', l\geq 0}$
by Condition \ref{def of X} and \ref{def of (E, lambda) of X}.
Then $(\X_{\theta}, (\V_{\alpha, \theta, l}, \U_{\alpha, \theta, l},
\U^\DM_{\alpha, \theta, l})_{\alpha \in \A \cup \A'})$ is Kurainshi data
of $\widehat{\M}_{\leq C}$.
\end{proof}

For Kuranishi data $(\X_{\theta}, \V_{\alpha, \theta, l}, \U_{\alpha, \theta, l},
\U^\DM_{\alpha, \theta, l})$ of $\widehat{\M}^{\leq L_{\max}}_{\leq C}$,
we define the pre-Kuranishi structure
\[
(\X_{\theta}, \forget_{o, \mathcal{S}, A}, (\W_x, \E_x, s_x, \widetilde{\psi}_x),
(\varphi_{x, y}, \hat \varphi_{x, y}))
\]
of each $\widehat{\M}^{\leq L_{\max}}_{\theta}$ as follows.
For each $p = (\Sigma, z, u) \in \widehat{\M}^{\leq L_{\max}}_{\theta}$ and
two points
\[
p^+_i = (\Sigma, z, (z^i)^S, (z^i)^A, (z^i)^o, u) \in
\X_{\theta} \cap \forget_{o, \mathcal{S}, A}^{-1}(p) \quad (i = 1, 2)
\]
in the same fiber, we define
$p^+_1 \vee p^+_2 \in \X_{\theta}$ by
\[
p^+_1 \vee p^+_2 = (\Sigma, z, (z^1)^S \cup (z^2)^S, (z^1)^A \cup (z^2)^A,
(z^1)^o \cup (z^2)^o, u).
\]

For each point $p = (\Sigma, z, z^S, z^A, z^o, u) \in \X_{\theta}$,
the Kuranishi neighborhood $(\W_p, \E_p, s_p, \widetilde{\psi}_p)$ of
$\forget_{o, \mathcal{S}, A}(p)$ is defined as follows.
By Condition \ref{def of X},
there exist subsets $z^{S, i} \subset z^S$, $z^{A, i} \subset z^A$, $z^{o, i} \subset z^o$
and indices $\alpha_i \in \A$ $(i = 1, \dots, k)$ such that
$p_i = (\Sigma, z, z^{S, i}, z^{A, i}, z^{o, i}, u) \in \V_{\alpha_i, \theta, l(p_i)}$ for all $i$,
$z^S = \bigcup_i z^{S, i}$, $z^A = \bigcup_i z^{A, i}$ and $z^o = \bigcup_i z^{o, i}$.
As in Condition \ref{def of (E, lambda) of X},
choose a subset $I \subset \{1, \dots, k\}$ such that $z^o = \coprod_{i \in I} z^{o, i}$,
and fix forgetful maps $f_i$ from $\hat p = \forget_u(p)$ to $\hat p_i = \forget_u(p_i)$
for each $i \in I$.
Let $(\hat P \to \hat X, Z, Z^S)$ be the local universal family of
$\forget_u(\forget_{o, A}(p))$.
We define an $\Aut(p)$-equivariant linear map $\widetilde{\lambda}_p : E^0_p \to
C^\infty(\hat P \times Y; \Wedge^{0, 1} V^\ast \hat P \otimes_\C
(\R \partial_\sigma \oplus T Y))$ by the sum of the pull backs of $\lambda_{\hat p_i}$
by the forgetful maps from $(\hat P \to \hat X, Z, Z^S)$ to the local universal families
of $\hat p_i$ whose restrictions to the central fiber coincide with $f_i$.
Then the Kuranishi neighborhood $(\W_p, \E_p, s_p, \widetilde{\psi}_p)$ of
$\forget_{o, \mathcal{S}, A}(p)$ is constructed by the argument in
Section \ref{construction of nbds} using this
$\Aut(p)$-equivariant linear map $\widetilde{\lambda}_p$.

For any $p = (\Sigma, z, u) \in \widehat{\M}^{\leq L_{\max}}_{\theta}$
and any two points
\[
p^+_i = (\Sigma, z, (z^i)^S, (z^i)^A, (z^i)^o, u) \in
\X_{\theta} \cap \forget_{o, \mathcal{S}, A}^{-1}(p),
\]
in the same fiber, $p^+_1 \leq p^+_2$ means that $(z^1)^S \subset (z^2)^S$,
$(z^1)^A \subset (z^2)^A$ and $(z^1)^o \subset (z^2)^o$,
which implies that $E^0_{p^+_1}$ is a subspace of $E^0_{p^+_2}$ and
the restriction of $\widetilde{\lambda}_{p^+_2}$ to $E^0_{p^+_1}$ is the pull back of
$\widetilde{\lambda}_{p^+_1}$ by the forgetful map.
Hence the embedding of the Kuranishi neighborhood
$(\W_{p^+_1}, \E_{p^+_1}, s_{p^+_1}, \widetilde{\psi}_{p^+_1})$ to
$(\W_{p^+_2}, \E_{p^+_2}, s_{p^+_2}, \widetilde{\psi}_{p^+_2})$ is defined by the argument
in Section \ref{embed}.
More generally, for any two points $x, y \in \X_{\theta}$,
if there exists some $r \in \psi_x(s_x^{-1}(0)) \cap \psi_y(s_y^{-1}(0))$ such that
$r_x \leq r_y$, where $r_x = \widetilde{\psi}_x^{-1}(\psi_x(r))$ and
$r_y = \widetilde{\psi}_y^{-1}(\psi_y(r))$,
then we can define the embedding of $(\W_x, \E_x, s_x, \widetilde{\psi}_x)$
to $(\W_y, \E_y, s_y, \widetilde{\psi}_y)$ by the argument in that section.

It is straightforward to check that each
\[
(\X_{\theta}, \forget_{o, \mathcal{S}, A}, (\W_x, \E_x, s_x, \widetilde{\psi}_x),
(\varphi_{x, y}, \hat \varphi_{x, y}))
\]
satisfies the other conditions of pre-Kuranishi structure.
Furthermore, they satisfy the compatibility conditions with respect to the fiber product
structure corresponding to the decomposition of holomorphic buildings into parts
and essential submersions corresponding to the decomposition
of holomorphic buildings into their connected components.
(We will consider these compatibility conditions in Section \ref{fiber prod} in details.)

\subsection{Decomposition by floor structure}
\label{floor decomposition}
A holomorphic building in the boundary $\partial \widehat{\M}$ is of height $k > 1$
and it can be decomposed into
the $[1, k_1]$-th floor part and the $[k_1+1, k]$-th floor part for each $1 \leq k_1 < k$.
In this section, we see the relation of the Kuranishi neighborhood of the whole
holomorphic building to those of these two parts.

First we define a space $\widehat{\M}^{\diamond 2}$ as follows.
Its point $((\Sigma^i, z^i, u^i)_{i=1,2}, M^{1,2})$ consists of
two holomorphic buildings $(\Sigma^i, z^i, u^i)$ ($i=1,2$) and
a set $M^{1,2} = \{(S^1_{+\infty_l}, S^1_{-\infty_l})\}$ of pairs of $+\infty$-limit circles
$S^1_{+\infty_l}$ of $(\Sigma^1, z^1, u^1)$ and $-\infty$-limit circles $S^1_{-\infty_l}$
of $(\Sigma^2, z^2, u^2)$ such that
the pairs in $M^{1,2}$ do not share the same limit circles.
Two points $((\Sigma^i, z^i, u^i)_{i=1,2}, M^{1,2})$ and
$(((\Sigma')^i, (z')^i, (u')^i)_{i=1,2}, \ab (M')^{1,2})$ are the same point if there exist
isomorphisms $\varphi^i : \Sigma^i \cong (\Sigma')^i$ and $\R$-translations
$\theta^i$ such that $\varphi^i(z^i) = (z')^i$,
$u^i = (\theta^i \times 1) \circ (u')^i \circ \varphi^i$ and
$(\varphi_1, \varphi_2)$ maps $M^{1,2}$ to $(M')^{1,2}$.
The pre-Kuranishi structure of $\widehat{\M}^{\diamond 2}$ is induced by
that of $(\widehat{\M} \times \widehat{\M}) / \mathfrak{S}_2$
since the only local difference of them is
the automorphism group.

Let $\widehat{\M}^{\diamond 2}_{l_{i, -}, l_{1,2}, l_{i, +}} \subset
\widehat{\M}^{\diamond 2}$ be the subspace of points
$((\Sigma^i, z^i, u^i)_{i=1,2}, M^{1,2})$
such that the number of pairs in $M^{1,2}$ is $l_{1,2}$ and
the number of $\pm\infty$-limit circles of $(\Sigma^i, z^i, u^i)$ which do not appear
in $M^{1,2}$ is $l_{i, \pm}$.
Let
\[
\Psi_{1, 2} : \widehat{\M}^{\diamond 2}_{l_{i, -}, l_{1,2}, l_{i, +}}
\to (\overline{P} \times \overline{P})^{l_{1,2}} / \mathfrak{S}_{l_{1,2}}
\]
be the continuous map which maps each point $((\Sigma^i, z^i, u^i)_{i=1,2}, M^{1,2})$
to the point $(\pi_Y \circ u^1|_{S^1_{+\infty_l}}, \pi_Y \circ u^2|_{S^1_{-\infty_l}})$.
Let $\Delta_{\overline{P}} \subset \overline{P} \times \overline{P}$ be the diagonal.
Since $\Psi_{1,2}$ is realized as a strong smooth map,
$\Psi_{1,2}^{-1}(\Delta_{\overline{P}}^{l_{1,2}} / \mathfrak{S}_{l_{1,2}})$
has a pre-Kuranishi structure.
We study about the map from $\partial \widehat{\M}$ to
$\Psi_{1,2}^{-1}(\Delta_{\overline{P}}^{l_{1,2}} / \mathfrak{S}_{l_{1,2}})$
defined by the decomposition by a gap of floors.
Since the decomposition depends on the choice of the gap,
this map is multi-valued.
%

We study about the relation of the Kuranishi neighborhoods of a point
$(\Sigma, z, u) \in \partial \widehat{\M}$ and
that of one of its image $((\Sigma^i, z^i, u^i)_{i=1,2}, M^{1,2})$ by the multi-valued map
$\partial \widehat{\M} \to
\Psi_{1,2}^{-1}(\Delta_{\overline{P}}^{l_{1,2}} / \mathfrak{S}_{l_{1,2}})$.
Assume that each $(\Sigma^i, z^i, u^i)$ is of height $k_i$.

Let $(V^i, E^i, s^i, \psi^i, G^i)$ be the Kuranishi neighborhood of $(\Sigma^i, z^i, u^i)$
defined by the data $((z^i)^+, S^i, E_i^0, \lambda^i)$ and additional data
$((z^i)^{++}, (S^i)', \hat R^i_j)$ for each $i=1,2$.
We consider the Kuranishi neighborhood
$(V, E, s, \psi, G)$ of $(\Sigma, z, u)$ in $\partial \widehat{\M}$
defined by the data $((z^1)^+ \cup (z^2)^+, S^1 \cup S^2, E_1^0 \oplus E_2^0,
\lambda^1 \oplus \lambda^2)$
and the additional data $((z^1)^{++}, (S^1)' \cup (S^2)', (\hat R^1_j, \hat R^2_j))$.

Fix a coordinate of each joint circle between the $k_1$-th floor and $(k_1 + 1)$-th floor
of $(\Sigma, z, u)$.
These define the coordinates of limit circles of $(\Sigma^1, z^1, u^1)$
and $(\Sigma^2, z^2, u^2)$ which appears in $M^{1,2}$.
Since the curves in each $V^i$ are constructed by patching parts of the curve
$\Sigma^i$, we can define a smooth map
\begin{align*}
\Upsilon:
V^1 \times V^2 &\to \prod_{(S^1_{+\infty_l}, S^1_{-\infty_l}) \in M^{1,2}} (P \times P)\\
((a^1, b^1, u^1, h^1), (a^2, b^2, u^2, h^2)) &\mapsto (\pi_Y \circ u^1|_{S^1_{+\infty_l}},
\pi_Y \circ u^2 \circ \phi|_{S^1_{-\infty_l}})
\end{align*}
by using these coordinates.
Let $I_\epsilon \subset \R$ be a small neighborhood of $0 \in \R$ and
define $I_\epsilon \cdot \Delta_{P} = \{(\gamma, t \cdot \gamma) \in P \times P;
\gamma \in P, t \in I_\epsilon\}$.
For each point in $\Upsilon^{-1}(I_{\epsilon} \cdot \Delta_{P})$,
we can define a (perturbed) holomorphic building by jointing each pair of limit circles
in $M^{1,2}$ by using the coordinates twisted by some $t_l \in I_\epsilon$.
In particular, we can define a continuous map $\psi$ from the zero set of
$(s^1 \oplus s^2)|_{(\Upsilon^{-1}(I_{\epsilon} \cdot \Delta_{P})}$ to $\widehat{\M}$.
Then
$(\Upsilon^{-1}(I_{\epsilon} \cdot \Delta_{P}), E^1 \oplus E^2, s^1 \oplus s^2, \psi, G)$
is isomorphic to a part of the Kuranishi neighborhood
$(V, E, s, \psi, G)$ of $(\Sigma, z, u)$ in $\partial \widehat{\M}$.
($V$ is a boundary of a manifold with corners, and it is a union of manifolds.
More precisely, one of them corresponding to the decomposition at the gap
between the $k_1$-th floor and the $(k_1 + 1)$-th floor
is isomorphic to
$(\Upsilon^{-1}(I_{\epsilon} \cdot \Delta_{P}), E^1 \oplus E^2, s^1 \oplus s^2, \psi, G)$.)
Indeed, we can define a map
\begin{align*}
\Upsilon^{-1}(I_{\epsilon} \cdot \Delta_{P}) &\to V\\
((a^1, b^1, u^1, h^1), (a^2, b^2, u^2, h^2)) &\mapsto (a^0, b^0, u^0, h^0)
\end{align*}
by
$h^0 = (h^1, h^2) \in E^0_1 \oplus E^0_2$,
$a^0 = (a^1, a^2, (0, t_l)_l) \in \widetilde{X} = \widetilde{X}^1 \times \widetilde{X}^2
\times \widetilde{D}^{l_{1,2}}$
($\widetilde{D}^{l_{1,2}}$ is the parameter space for the deformation near
the joint circles between $k_1$-th floor and $(k_1+1)$-th floor),
$u^0 = u^1 \cup u^2$,
$b^0_\mu = b^1_\mu$ for $\mu \in \bigcup_{1 \leq j < k_1} M_j = M^1$,
$b^0_\mu = b^2_\mu$ for $\mu \in \bigcup_{k_1 < j < k_1 + k_2} M_j = M^2$ and
\begin{align*}
b^0_\mu
&= \lim_{s \to \infty} (\sigma \circ u^1|_{[0, \infty) \times S^1_{+\infty_l}}(s, t)
- (0_{k_{i_0}} + L_\mu s))\\
&\quad - \lim_{s \to -\infty} (\sigma \circ u^2|_{(-\infty, 0] \times S^1_{-\infty_l}}(s, t)
- (0_0 + L_\mu s))
\end{align*}
for $\mu = (S^1_{+\infty_l}, S^1_{-\infty_l}) \in M_{k_1} \cong M^{1, 2}$.
As we explained in the last of Section \ref{smoothness}, $b^0_\mu$
($\mu \in M_{k_1}$) are smooth function of $((a^1, b^1, u^1, h^1), (a^2, b^2, u^2, h^2))$.
Hence this map is a diffeomorphism and it defines an isomorphism of
$(\Upsilon^{-1}(I_{\epsilon} \cdot \Delta_{P}), E^1 \oplus E^2, s^1 \oplus s^2, \psi, G)$
and its image in $(V, E, s, \psi, G)$.

The above isomorphism implies that the Kuranishi neighborhood of each point
in $\partial \widehat{\M}$ and those of its image by the map
$\partial \widehat{\M} \to
\Psi_{1,2}^{-1}(\Delta_{\overline{P}}^{l_{1,2}} / \mathfrak{S}_{l_{1,2}})$ are the same
modulo automorphism group.
In particular, the map $\partial \widehat{\M} \to
\Psi_{1,2}^{-1}(\Delta_{\overline{P}}^{l_{1,2}} / \mathfrak{S}_{l_{1,2}})$
is a multi-valued partial submersion between pre-Kuranishi spaces.


%% file: SFT-06_Fiber_products_for_Y.tex
%
%

\section{Fiber products}\label{fiber prod}
The pre-Kuranishi spaces considered in Section \ref{construction of Kuranishi}
are the spaces of holomorphic buildings without any conditions on periodic orbits
on limit circles $S^1_{\pm\infty_i}$.
For the construction of the algebra, we need to use the fiber products of such
Kuranishi spaces with $\overline{P}$ and $Y$.
More precisely, we use the fiber products of $\overline{\M}$ with the lifts of simplices in
$\overline{P}$ to $P$, and we need to perturb the section so that the induced
multisections on the fiber products are independent of the choice of these lifts.

To construct the virtual fundamental chain, we also need to define the orientations of
Kuranishi spaces.
We cannot define the orientations of $\widehat{\M}$ or $\overline{\M}$, but it is enough
to define the orientation of the fiber products we use.
The fiber products with simplices in $Y$ and the lifts of simplices in $\overline{P}$ to $P$
are orientable provided that interiors of these simplices in $\overline{P}$ do not contain
bad orbits.

In the general Bott-Morse case, it is not enough to count the intersection numbers
with simplices in $\overline{P}$, and we need to add correction terms,
which are equivalent to counting cascades in \cite{Bo02}.
This is because the chain which represents the diagonal in Poincar\'e duality is different
from the genuine diagonal in chain level.
These correction terms appear in every Bott-Morse theory if we construct
the algebra by the intersection numbers of the moduli spaces with simplices.
However, since algebraic structure of SFT is more complicated than that of usual Morse
theory, to define the correction terms, we need to solve some algebraic equations.

First we explain the bad orbits in Section \ref{bad orbits}, and
in Section \ref{fiber prod with simpleces}, we explain the fiber products
of $\overline{\M}$ which we use for the construction of the algebra.
In Section \ref{construction of a multisection}, we construct
a family of perturbed multisections of fiber products of
another space $\widehat{\M}^\diamond$
which satisfies appropriate compatibility conditions, and we use
the induced multisections for the fiber products of $\overline{\M}$.
Next in Section \ref{fiber prod and orientation}, we explain the orientations
of the fiber products.
In Section \ref{algebra for correction},
we define the correction terms, and finally in Section \ref{construction of algebra},
we recall the algebra of SFT and explain how to define the algebra by
the virtual fundamental chains of our fiber products.

\subsection{Bad orbits and local coefficients}\label{bad orbits}
Before considering the fiber products of the space of holomorphic buildings,
first we explain about bad orbits.
In Morse case, it is well known that bad orbits should not count as the generators of
the chain complex.
However, in our Bott-Morse case, bad orbits appear as a closed subset of $\overline{P}$.
Hence we need to explain how to treat these bad orbits.

First we define bad orbits.
It is related to orientations of the following $\bound$-operators associated
to periodic orbits.
For each $\gamma \in P \subset C^\infty(S^1, Y)$,
fix one trivialization $\gamma^\ast T \hat Y \cong \C^n$.
Let
\begin{align*}
\mathring{D}^+_\gamma :&\ W_\delta^{1, p}((-\infty, 0] \times S^1 \cup D_\infty, \gamma^\ast T \hat Y \cup \C^n)\\
&\to L_\delta^p((-\infty, 0] \times S^1, \gamma^\ast T \hat Y) \oplus
L^p(D_\infty, \Wedge^{0, 1}T^\ast D_\infty \otimes \C^n)\\
D^+_\gamma :&\ \widetilde{W}_\delta^{1, p}((-\infty, 0] \times S^1 \cup D_\infty, \gamma^\ast T \hat Y \cup \C^n)\\
&\to L_\delta^p((-\infty, 0] \times S^1, \gamma^\ast T \hat Y) \oplus
L^p(D_\infty, \Wedge^{0, 1}T^\ast D_\infty \otimes \C^n)\\
\mathring{D}^-_\gamma :&\ W_\delta^{1, p}(D_0 \cup [0, \infty) \times S^1, \C^n \cup \gamma^\ast T \hat Y)\\
&\to L^p(D_0, \Wedge^{0, 1}T^\ast D_\infty \otimes \C^n) \oplus
L_\delta^p([0, \infty) \times S^1, \gamma^\ast T \hat Y)\\
D^-_\gamma :&\ \widetilde{W}_\delta^{1, p}(D_0 \cup [0, \infty) \times S^1, \C^n \cup \gamma^\ast T \hat Y)\\
&\to L^p(D_0, \Wedge^{0, 1}T^\ast D_\infty \otimes \C^n) \oplus
L_\delta^p([0, \infty) \times S^1, \gamma^\ast T \hat Y)
\end{align*}
be $\bound$-type linear operators such that
\[
\mathring{D}^+_\gamma \xi = D^+_\gamma \xi = \partial_s \xi + J(\gamma)(\nabla_t \xi - L_\gamma \nabla_\xi R_\lambda(\gamma))
\]
on $(-\infty, 0] \times S^1$ and
\[
\mathring{D}^-_\gamma \xi = D^-_\gamma \xi = \partial_s \xi + J(\gamma)(\nabla_t \xi - L_\gamma \nabla_\xi R_\lambda(\gamma))
\]
on $[0, \infty) \times S^1$, where
$D_\infty = \{ z \in \C \cup \{\infty\}; |z| \geq 1\}$, $D_0 = \{ z \in \C; |z| \leq 1\}$,
and we identify $\{0\} \times S^1$ with $\partial D_\infty$ or $\partial D_0$
by $(0, t) \leftrightarrow e^{2\pi \sqrt{-1} t}$.
(The above $\widetilde{W}^{1, p}_\delta$ is defined by
$\widetilde{W}^{1, p}_\delta = W^{1, p}_\delta \oplus \Ker A_\gamma$ as in Section
\ref{construction of Kuranishi})
Adding finite-dimensional complex vector spaces to the domain vector spaces
if necessary, we assume the above operators are surjective.

We consider the orientations of these types of operators, that is, the orientation of their kernels.
Since $\bound$-type operators of each type are connected linearly
(that is, two operators $D$ and $D'$ can be connected by a family of operators
$t D + (1-t) D'$ ($t \in [0, 1]$)),
we can define a consistent orientation of these operators for each type.
Furthermore, changing the trivialization of $\gamma^\ast T \hat Y$ is equivalent to
gluing a $\bound$-operator of a holomorphic bundle on $\C P^1$ to the operators.
Since a $\bound$-operator of a holomorphic bundle has the complex orientation,
an orientation of one $\mathring{D}^+_\gamma$ defines the compatible orientations of
all operators of type $\mathring{D}^+_\gamma$ for each $\gamma \in P$.
Therefore, we can consider an orientation of $\mathring{D}^+_\gamma$ without fixing
particular trivialization of $\gamma^\ast T \hat Y$ or an additional complex vector space.

Let $\S^D$ be the local system of orientation of $\mathring{D}^+_\gamma$ on $P$,
and let $\S^{\lsuperscript{D}{t}}$ be the local system of orientation of
$\mathring{D}^-_\gamma$ on $P$.
We say $\gamma \in \overline{P}$ is a bad orbit if
$\S^D$ is not trivial on $\pi_P^{-1}(\gamma) \subset P$.
Let $\overline{P}^{\text{bad}} \subset \overline{P}$ be the subset of bad orbits.
Similarly, let $\overline{P}^{^t\text{bad}} \subset \overline{P}$ be the set of points
$\gamma \in \overline{P}$ such that $\S^{\lsuperscript{D}{t}}$ is not trivial on
$\pi_P^{-1}(\gamma)$.
By the assumption of $K$, $\overline{P}^{\text{bad}}$ and $\overline{P}^{^t\text{bad}}$ are
subcomplexes of $\overline{P}$.

\begin{rem}
Let $\gamma_0$ be a simple periodic orbit, and $\gamma = \gamma_0^{2^k m}$ be its
$2^k m$-multiple, where $m \geq 1$ is an odd integer.
Then $\gamma$ is a bad orbit if and only if $k \geq 1$ and
$\ind \mathring{D}^+_{\gamma^2_0} - \ind \mathring{D}^+_{\gamma_0}$ is odd.
Similarly, $\gamma$ belongs to $\overline{P}^{^t\text{bad}}$ if and only if $k \geq 1$ and
$\ind D^+_{\gamma^2_0} - \ind D^+_{\gamma_0}$ is odd.
Note that the index of the operator $\mathring{D}^+_\gamma$ is determined by the
Conley Zehnder index of $\gamma$ and $\dim T_\gamma P / TS^1$ as follows.
Fix one trivialization $\gamma^\ast \xi \cong \C^{n-1}$, which induces a trivialization
$\gamma^\ast T \hat Y \cong (\R \partial_\sigma \oplus \R R_\lambda) \oplus
\gamma^\ast \xi \cong \C^n$.
We define the Conley Zehnder index $\CZ \gamma$ of $\gamma$ by
the Conley Zehnder index of the path $\{\varphi^\lambda_t\}_{t \in [0, L_\gamma]}$
of symplectic matrices under the above trivialization of $\gamma^\ast \xi$.
(See \cite{RS93} for the definition of Conley Zehnder index of a path of symplectic
matrices.)
Then it is easy to see that
\[
\ind \mathring{D}^+_\gamma = (n-1) - \CZ \gamma - \frac{1}{2} \dim T\gamma P/ TS^1.
\]
Similarly, the index of the operator $D^+_\gamma$ is
\[
\ind D^+_\gamma = (n+1) - \CZ \gamma + \frac{1}{2} \dim T_\gamma P / TS^1.
\]
\end{rem}
\begin{eg}
We give an example where bad orbits appear as a subcomplex of $\overline{P}$.
This example was given by Bourgeois in \cite{Bo02}.
Let $K = \R^2 / G$ be a Kulein bottle, where $G$ is a group of diffeomorphisms of $\R^2$
generated by $(x, y) \mapsto (x + 1, 1 - y)$ and $(x, y) \mapsto (x, y + 1)$.
We equip $K$ with the flat metric $dx \otimes dx + dy \otimes dy$, and regard its unit
tangent bundle $S(TK) \cong S(T^\ast K)$ as a contact manifold by the Liouville form.
Then the Reeb flow is the geodesic flow.
$\overline{P}_2$ contains a component
\[
\{\gamma_{y}(t) = ((t, y), (1, 0)) : [0, 2]/ \{0, 2\} \to S(TK); y \in [0, 1/2]\},
\]
which is homeomorphic to the interval $[0, 1/2]$.
It contains two multiple orbits $\gamma_0$ and $\gamma_{1/2}$, and the others are
simple.
It is easy to check that the index of the operators $\mathring{D}^+$ for the two are even
and those for $\gamma_0|_{[0, 1]}$ and $\gamma_{1/2}|_{[0, 1]}$ are odd.
Hence these two orbits are bad orbits.
\end{eg}

Let $f : K \to \overline{P}$ be an ordered triangulation.
(``ordered'' means the set of the vertices has a total order.)
For each point $p \in K$, let $d$ be the multiplicity of the periodic orbit corresponding
to $p$.
Then we assume that there exists a regular $\Z/d$-complex $L$
(see \cite{Bre72} for regular complex),
an isomorphism $\varphi : L / (\Z/d) \cong \St(p, K)$ and
a smooth $\Z/d$-equivariant embedding $\tilde f : L \to P$ such that
$f \circ \varphi \circ \pi_L = \pi_P \circ \tilde f : L \to \overline{P}$,
where $\pi_L : L \to L / (\Z/d)$ is the quotient map.
(Note that locally $\pi_P : P \to \overline{P}$ can be written as
$S^1 \times_{\Z / d} W \to W / (\Z/d)$ for some $\Z/d$-manifold $W$.
Hence a $\Z / d$-equivariant triangulation $\check f : L \to W$ defines an embedding
$\tilde f : L \to S^1 \times_{\Z/d} W$ by $\tilde f(x) = [0, \check f (x)]$.)

Let $K^2 \to \overline{P} \times \overline{P}$ be an Euclidean cell decomposition
which is a refinement of $\{s \times t; s, t \in K\}$ and
which contains $\Delta_\ast K = \{\Delta_\ast s; s \in K\}$ and
$\rho_\ast K = \{\partial_{p+1} \dots \partial_n s \times \partial_0 \dots \partial_{p-1} s;
s \in K, 0 \leq p \leq n = \dim s\}$ as subcomplexes.
The chain complex $C_\ast(\overline{P} \times \overline{P})$ is defined by
using this Euclidean cell decomposition as a CW decomposition.
Let $K^0 = (x)$ be a finite sequence of smooth cycles in $Y$.
We denote their cohomologies by $\overline{K}^0 = (\overline{x})$
($\overline{x} \in H_\ast (Y, \Q)$).

We will use the generators $c \theta^D_c$ of the relative chain complex
$C_\ast(\overline{P}, \overline{P}^{\text{bad}}; \S^D \otimes \Q)$ of
ordered simplicial complex, or the generators $(c \theta^D_c)^\ast$ of the cochain complex
with compact support, where $\S^D$ is the induced local system on
$\overline{P} \setminus \overline{P}^{\text{bad}}$.
($(c \theta^D_c)^\ast$ is the cochain which takes one at $c \theta^D_c$ and
which vanishes at the other simplices.)
The $\Z/2$-degree of the above chain complex is defined by
\[
|c \theta^D_c| = \dim c + |\theta^D_c|,
\]
where $|\theta^D_c|$ is the index of the operator $\mathring{D}^+_\gamma$
($\gamma \in |c|$), and
its boundary operator is defined by
\[
\partial (c \theta^D_c) = (\partial c) \theta^D_c.
\]
Note that local system $\S^D$ is not well-defined on $\overline{P}^{\text{bad}}$,
but the above relative chain complex is well-defined.

We construct algebra by counting some intersection numbers with the moduli spaces
and simplices in $\overline{P}$.
Hence we need Poincar\'e duality.
In particular, we need a local system of the orientation of $\overline{P}$.
However, in general, the local orientation of $\overline{P}$ is not well defined.
We treat this as follows.

Let $\S^{\overline{P}}$ be the local system of the orientation of $TP / TS^1$ on $P$,
where $TS^1$ is the tangent of the $S^1$-action on $P$.
We say $\gamma \in \overline{P}$ is a non-orientable point
if $\S^{\overline{P}}$ is not trivial on $\pi_P^{-1}(\gamma)$.
Let $\overline{P}^{\text{no}} \subset \overline{P}$ be the set of non-orientable points.
It is also a subcomplex of $\overline{P}$.
Then $\S^{\overline{P}}$ induces a local system on
$\overline{P} \setminus \overline{P}^{\text{no}}$.

For each top-dimensional simplex $\zeta \in K$
(the top-dimension depends on the connected component of $\overline{P}$),
let $m_\zeta$ be the multiplicity of the periodic orbits in $\Int \zeta$ (it is constant on
$\Int \zeta$).
$m_\zeta$ depends only on the connected component of $P$ containing $|\zeta|$.

Let $\tilde \zeta \inj P$ be a lift of $\zeta$.
Then the orientation of $TP/ TS^1$ defined by the orientation of $\tilde \zeta$
induces a section $\theta^{\overline{P}}_\zeta$ of  $\S^{\overline{P}}$ on $\Int \zeta$.
This section is independent of the choice of the lift $\tilde \zeta$.

We call a chain
\[
[\overline{P}] = \sum_\zeta \frac{1}{m_\zeta} \zeta \theta^{\overline{P}}_\zeta
\in C_{\dim P - 1}(\overline{P}, \overline{P}^{\text{no}}; \S^{\overline{P}} \otimes \Q)
\]
the fundamental cycle of $\overline{P}$, where the sum is taken over all top-dimensional
simplices of $K$.
As usual, this is a cycle in the relative chain complex.

Before considering cap products with the fundamental chain,
we see the relation of the orientations of the operators $\mathring{D}^\pm_\gamma$,
$D^\pm_\gamma$ and that of the tangent space $T_\gamma P / TS^1$.
First recall that
\[
\Ker A_\gamma \cong \R \oplus T_\gamma P
\cong (\R \oplus TS^1) \oplus T_\gamma P / TS^1.
\]
We denote the kernel of a surjective operator $D$ on a curve
(or the kernel of the surjective operator obtained by adding a finite-dimensional complex
vector space to the domain of a non-surjective operator $D$) by $[D]$.
The fiber product
\[
[D^-_\gamma] \underset{\Ker A_\gamma}{\times} [D^+_\gamma]
= [D^-_\gamma] \underset{(\R \oplus TS^1) \oplus T_\gamma P / TS^1}{\times}
[D^+_\gamma]
\]
is equivalent to the kernel of a $\bound$-operator on a complex
vector bundle over $\C P^1$ by gluing.
Hence it has the complex orientation.
The space $[\mathring{D}^-_\gamma] \oplus [\mathring{D}^+_\gamma]$ is a subspace
of the above fiber product,
and its quotient space is isomorphic to $(\R \oplus TS^1) \oplus T_\gamma P / TS^1$.
Therefore, if orientations of $[\mathring{D}^+_\gamma]$ and $T_\gamma P / TS^1$ are
given, we can define the orientation of $[\mathring{D}^-_\gamma]$ so that the orientation
of the above fiber product defined by
\[
[\mathring{D}^-_\gamma] \oplus (\R \oplus TS^1) \oplus T_\gamma P / TS^1 \oplus
[\mathring{D}^+_\gamma]
\]
coincides with the complex orientation.

To define Poincar\'e dual, first we recall the definition of cap product
without local coefficient. (Our definition is a bit different from the usual one.)
For a $p$-cochain $\alpha$ and a simplex $\zeta$ of dimension $n$,
our cap product $\zeta \cap \alpha$ is defined by
\[
\zeta \cap \alpha = \partial_{n- p + 1} \partial_{n - p + 2} \dots \partial_n \zeta \langle \partial_0 \partial_1 \dots \partial_{n - p-1} \zeta, \alpha \rangle.
\]
\begin{rem}
Under this definition, the following equation holds true.
For any $p$-cochain $\alpha$ and $n$-chain $\zeta$,
\[
\partial (\zeta \cap \alpha) = \partial \zeta \cap \alpha + (-1)^{n - p} \zeta \cap \delta \alpha
\]
\end{rem}

For each cochain $\alpha = (c \theta^D_c)^\ast \in
C^\ast (\overline{P}, \overline{P}^{\text{bad}}; \S^D \otimes \Q)$,
we define the chain $[\overline{P}] \cap \alpha \in
C_\ast(\overline{P}, \overline{P}^{^t\text{bad}}; \S^{\lsuperscript{D}{t}} \otimes \Q)$
as follows.
For each top-dimensional simplex $\zeta$ in $K$, let $\tilde \zeta \subset P$ be its lift.
If $c = \partial_0 \partial_1 \dots \partial_{n - p - 1} \zeta$ for some $p$,
then we can extend the orientation $\theta^D_c$ of $\S^D$ on
$\partial_0 \partial_1 \dots \partial_{n - p-1} \tilde \zeta$ to that on $\tilde \zeta$.
Then $\theta^D_c$ and $\theta^{\overline{P}}_{\zeta}$ define the orientation
$\theta^{\lsuperscript{D}{t}}_{c, \zeta}$ of $\S^{\lsuperscript{D}{t}}$ on $\tilde \zeta$
as above.
If $\partial_{n- p + 1} \partial_{n - p + 2} \dots \partial_n \zeta$ is not contained in
$\overline{P}^{^t\text{bad}}$, then $\theta^{\lsuperscript{D}{t}}_{c, \zeta}$ defines
the orientation of $\S^{\lsuperscript{D}{t}}$ on
$\partial_{n- p + 1} \partial_{n - p + 2} \dots \partial_n \zeta$.
We define $[\overline{P}] \cap \alpha$ by the linear combination of the cap products
\[
(\zeta \theta^{\overline{P}}_\zeta) \cap (c \theta^D_c)^\ast
= \theta^{\lsuperscript{D}{t}}_{c, \zeta} (\zeta \cap c^\ast) .
\]

We define the boundary operator of $C_\ast(\overline{P}, \overline{P}^{^t\text{bad}};
\S^{\lsuperscript{D}{t}} \otimes \Q)$ by
\[
\partial (\theta_\eta^{\lsuperscript{D}{t}} \eta) = (-1)^{|\theta_\eta^{\lsuperscript{D}{t}}|}
\theta_\eta^{\lsuperscript{D}{t}} \partial \eta,
\]
where $|\theta_\eta^{\lsuperscript{D}{t}}|$ is the index of the operator
$\mathring{D}_\gamma^-$ ($\gamma \in |\eta|$).
Similarly, the boundary operator of $C_\ast
(\overline{P} \times \overline{P}, \overline{P}^{^t\text{bad}} \times \overline{P} \cup
\overline{P} \times \overline{P}^{\text{bad}};
p_1^\ast \S^{\lsuperscript{D}{t}} \otimes p_2^\ast \S^D \otimes \Q)$
is defined by
\[
\partial (\theta_\eta^{\lsuperscript{D}{t}} \eta \theta_\eta^D)
= (-1)^{|\theta_\eta^{\lsuperscript{D}{t}}|} \theta_\eta^{\lsuperscript{D}{t}}
(\partial \eta) \theta_\eta^D.
\]

Let $\Delta : \overline{P} \to \overline{P} \times \overline{P}$ be the diagonal map.
We define a cycle $\Delta_\ast [\overline{P}]$ of
\[
C_{\dim P - 1}
(\overline{P} \times \overline{P}, \overline{P}^{^t\text{bad}} \times \overline{P} \cup
\overline{P} \times \overline{P}^{\text{bad}};
p_1^\ast \S^{\lsuperscript{D}{t}} \otimes p_2^\ast \S^D \otimes \Q)
\]
by
\[
\Delta_\ast [\overline{P}] = \sum \frac{1}{m_\zeta} \theta^{\lsuperscript{D}{t}}_\zeta
(\Delta_\ast \zeta) \theta^D_\zeta,
\]
where the sum is taken over all top-dimensional simplices of $K$ not contained in
$\overline{P}^{\text{bad}}$,
$\theta^D_\zeta$ is an arbitrary fixed orientation of $p_2^\ast \S^D$ on
$\Int \Delta_\ast \zeta$, and $\theta^{\lsuperscript{D}{t}}_\zeta$ is the orientation of
$p_1^\ast \S^{\lsuperscript{D}{t}}$ defined by $\theta^D_\zeta$ and
$\theta^{\overline{P}}_\zeta$
as above.
This definition is independent of the choice of $\theta^D_\zeta$.

For each simplex $\zeta \in K$ of dimension $n$, we define
a chain $\rho_\ast \zeta$ in $\overline{P} \times \overline{P}$ by
\[
\rho_\ast \zeta = \sum_{0 \leq p \leq n} \partial_{p+1} \dots \partial_n \zeta
\times \partial_0 \dots \partial_{p+1} \zeta.
\]
This corresponds to the image of $\Delta_\ast \zeta$ by Alexander Whitney map
$C_\ast(\overline{P} \times \overline{P})
\to C_\ast(\overline{P}) \otimes C_\ast(\overline{P})$.
(Recall that $K^2$ is not a simplicial somplex but a Euclidean cell complex
which contains $\partial_{p+1} \dots \partial_n \zeta
\times \partial_0 \dots \partial_{p+1} \zeta$.
For the transversality condition, it is convenient not to subdivide these products
because if we subdivide the complex, then we need to make the zero set of
perturbed mulitisection transverse to the new simplices of less dimension.)
We define a cycle $\rho_\ast [\overline{P}]
\in C_{\dim P -1} (\overline{P} \times \overline{P},
\overline{P}^{^t\text{bad}} \times \overline{P} \cup
\overline{P} \times \overline{P}^{\text{bad}};
p_1^\ast \S^{\lsuperscript{D}{t}} \otimes p_2^\ast \S^D \otimes \Q)$
by
\[
\rho_\ast [\overline{P}] = \sum \frac{1}{m_\zeta} \theta^{\lsuperscript{D}{t}}_\zeta
(\rho_\ast \zeta) \theta^D_\zeta.
\]
For later use, we remark that $\rho_\ast [\overline{P}]$ can be written as
\[
\rho_\ast [\overline{P}] = \sum_c ([\overline{P}] \cap (c \theta^D_c)^\ast) \otimes
c \theta^D_c,
\]
where the sum is taken over all simplices $c$ in $K$ which are not contained in
$\overline{P}^{\text{bad}}$.

Let $\epsilon_\ast : C_\ast(\overline{P}) \to C_\ast(\overline{P} \times \overline{P})$
be the natural linear map such that
$\rho_\ast - \Delta_\ast = \partial \circ \epsilon_\ast + \epsilon_\ast \circ \partial$,
and define a chain
$\epsilon_{\overline{P}} \in C_{\dim P} (\overline{P} \times \overline{P},
\overline{P}^{^t\text{bad}} \times \overline{P} \cup
\overline{P} \times \overline{P}^{\text{bad}};
p_1^\ast \S^{\lsuperscript{D}{t}} \otimes p_2^\ast \S^D \otimes \Q)$ by
$\epsilon_{\overline{P}} = \epsilon_\ast [\overline{P}]$.
Then it satisfies
\[
(\rho_\ast - \Delta_\ast) [\overline{P}] = \partial \epsilon_{\overline{P}}.
\]
This chain will be used for the definition of the correction terms.

We define
$\mathring{K}^2 \subset K^2$ by the minimal subcomplex which contains
$\Delta_\ast s$, $\rho_\ast s$ and $\epsilon_\ast s$ for all $s \in K$.


\subsection{Fiber products with simpleces}\label{fiber prod with simpleces}
First we define a Hausdorff space
$\overline{\M}^m_{((l_{i, j}), (l_{i, \pm}), (\mu_i))}$
for each family of non-negative integers
$((l_{i, j})_{1 \leq i < j \leq m}, (l_{i, \pm})_{1 \leq i \leq m}, (\mu_i)_{1 \leq i \leq m})$
as follows.
(We can equip it with a natural pre-Kuranishi structure, but it is not necessary.)
Its point $(\Sigma_i, z_i, u_i, \phi_i)_{1 \leq i \leq m}$ is a sequence of
holomorphic buildings $(\Sigma_i, z_i, u_i, \phi_i) \in \overline{\M}$ with the following
index sets of limit circles.
The index set of $+\infty$-limit circles of $\Sigma_i$ is
\[
\{+\infty^{i, +\infty}_l; 1 \leq l \leq l_{i, +}\} \sqcup
\coprod_{j = i+1}^m \{+\infty^{i, j}_l; 1 \leq l \leq l_{i, j}\}
\]
and the index set of $-\infty$-limit circles of $\Sigma_i$ is
\[
\{-\infty^{i, -\infty}_l; 1 \leq l \leq l_{i, -}\} \sqcup
\coprod_{j = 1}^{i-1} \{-\infty^{i, j}_l; 1 \leq l \leq l_{j, i}\}.
\]
$\overline{\M}^m_{((l_{i, j}), (l_{i, \pm}), (\mu_i))}$
is locally isomorphic to the product $\prod^m \overline{\M}$.
(An isomorphism is determined if we fix a family of bijections between the above
index sets of limit circles and the usual index sets $\{\pm\infty^i_l\}$.)
Note that we respect the indices of limit circles.
Hence even if we change the indices $+\infty^{i, j}_l$ to $+\infty^{i, j}_{g \cdot l}$ and
$-\infty^{j, i}_l$ to $-\infty^{j, i}_{g \cdot l}$ for the same $g \in \mathfrak{S}_{l_{i, j}}$,
we distinguish the obtained curve from the original one
(unless $g$ can be extended to the automorphisms of $(\Sigma_i, z_i, u_i, \phi_i)$).

We define the genus of a point $(\Sigma_i, z_i, u_i, \phi_i)_{1 \leq i \leq m} \in
\overline{\M}^m_{((l_{i, j}), (l_{i, \pm}), (\mu_i))}$ by
\[
g = 1 + \sum_{i = 1}^m (g_i - 1) + \sum_{1 \leq i < j \leq m} l_{i, j},
\]
where $g_i$ is the genus of $\Sigma_i$.
(This is the genus of the curve obtained by gluing joint circles $S^1_{+\infty^{i, j}_l}$ and
$S^1_{-\infty^{j, i}_l}$ for all pairs $(+\infty^{i, j}_l, -\infty^{j, i}_l)$.)

Note that there exists a natural continuous map
\begin{align*}
&\overline{\M}^m_{((l_{i, j}), (l_{i, \pm}), (\mu_i))}
\to \prod_{1 \leq i < j \leq m} (P \times P)^{l_{i, j}}
\times \prod_{1 \leq i \leq m} P^{l_{i, -}}
\times \prod_{1 \leq i \leq m} Y^{\mu_i}
\times \prod_{1 \leq i \leq m} P^{l_{i, +}}\\
&(\Sigma_i, z_i, u_i, \phi_i)_i \mapsto ((\pi_Y \circ u_i \circ \phi_{+\infty^{i, j}_l},
\pi_Y \circ u_i \circ \phi_{-\infty^{j, i}_l}),\\
&\hph{(\Sigma_i, z_i, u_i, \phi_i)_i \mapsto (}
\pi_Y \circ u_i \circ \phi_{+\infty^{i, +\infty}_l},
\pi_Y \circ u_i(z_{i, l}), \pi_Y \circ u_i \circ \phi_{-\infty^{i, -\infty}_l})
\end{align*}

We consider the fiber products with respect to this continuous map.
We consider the following family of sequences
$((\hat \epsilon^{i, j}_l), (\hat c^i_l), (x^i_l), (\hat \eta^i_l))$
of simplices with local coefficients.

$(\hat \epsilon^{i, j}_l = \theta^{\lsuperscript{D}{t}}_{\epsilon^{i, j}_l} \epsilon^{i, j}_l
\theta^D_{\epsilon^{i, j}_l})_{1 \leq i < j \leq m, 1 \leq l \leq l_{i, j}}$
is a sequence of products of
\begin{itemize}
\item
cells $\epsilon^{i, j}_l$ in $\mathring{K}^2$ which are not contained in
$\overline{P}^{^t\text{bad}} \times \overline{P} \cup
\overline{P} \times \overline{P}^{\text{bad}}$, and
\item
orientations $\theta^{\lsuperscript{D}{t}}_{\epsilon^{i, j}_l}$ of
$p_1^\ast \S^{\lsuperscript{D}{t}}$
and $\theta^D_{\epsilon^{i, j}_l}$ of $p_2^\ast \S^D$ on $\Int \epsilon^{i, j}_l$.
\end{itemize}
Take a lift $\tilde \epsilon^{i, j}_l \inj P \times P$ for each $\epsilon^{i, j}_l$, and
define $\breve \epsilon^{i, j}_l = \theta^{\lsuperscript{D}{t}}_{\epsilon^{i, j}_l}
\tilde \epsilon^{i, j}_l \theta^D_{\epsilon^{i, j}_l}$.

$(\hat c^i_l = c^i_l \theta^D_{c^i_l})_{1 \leq l \leq l_{i, -}}$ ($1 \leq i \leq m$)
is a sequence of products of
\begin{itemize}
\item
simplices $c^i_l$ in $K$ which are not contained in $\overline{P}^{\text{bad}}$, and
\item
orientations $\theta^D_{c^i_l}$ of $\S^D$ on $\Int c^i_l$.
\end{itemize}
For each $c^i_l$, we take its lift $\tilde c^i_l \inj P$ and
define $\breve c^i_l = \tilde c^i_l \theta^D_{c^i_l}$.

$(x^i_1, x^i_2, \dots, x^i_{\mu_i})_{i = 1, 2, \dots, m}$
is a sequence of cycles in $K^0$.

$(\hat \eta^i_l = \theta^{\lsuperscript{D}{t}}_{\eta^i_l} \eta^i_l)_{1 \leq l \leq l_{i, +}}$
($1 \leq i \leq m$)
is a sequence of products of
\begin{itemize}
\item
simplices $\eta^i_l$ in $K$ which are not contained in $\overline{P}^{^t\text{bad}}$, and
\item
orientations $\theta_{\eta^i_l}^{\lsuperscript{D}{t}}$ of $\S^{\lsuperscript{D}{t}}$
on $\Int \eta^i_l$.
\end{itemize}
For each $\hat \eta^i_l$, we take its lift $\tilde \eta^i_l \inj P$ and
define $\breve \eta^i_l = \theta^{\lsuperscript{D}{t}}_{\eta^i_l} \tilde \eta^i_l$.


Then for such a family of sequences $((\breve \epsilon^{i, j}_l), (\breve c^i_l),
(x^i_l), (\breve \eta^i_l))$, we define a closed subspace
\begin{align*}
\overline{\M}^m_{((\breve \epsilon^{i, j}_l), (\breve c^i_l),
(x^i_l), (\breve \eta^i_l))}
&\subset \overline{\M}^m_{((l_{i, j}), (l_{i, \pm}), (\mu_i))}
\end{align*}
as the fiber product with
\[
\prod \tilde \epsilon^{i, j}_l \times \prod \tilde c^i_l \times \prod x^i_l \times
\prod \tilde \eta^i_l
\subset \prod (P \times P)^{l_{i, j}}
\times \prod P^{l_{i, -}}
\times \prod Y^{\mu_i}
\times \prod P^{l_{i, +}}.
\]

The pre-Kuranishi structure of the above fiber product is defined as follows.
For a point $(\Sigma_i, z_i, u_i, \phi_i)_{1 \leq i \leq m}
\in \overline{\M}^m_{((\breve \epsilon^{i, j}_l),
(\breve c^i_l), (x^i_l), (\breve \eta^i_l))}$,
let $(V^i, E^i, s^i, \psi^i)$ be the Kuranishi neighborhood of each
$(\Sigma_i, z_i, u_i) \in \widehat{\M}$.
(This is not a Kuranishi neighborhood of $\overline{\M}$ but of $\widehat{\M}$.)
Since the limit circles of each curve in $V^i$ are identified with the limit circles of
$\Sigma_i$ by construction, it is meaningful to say that a coordinate of a limit circle
of a curve corresponding to the point $(a'_i, b'_i, u'_i, h'_i)$ in $V^i$ is close to that of
$(\Sigma_i, z_i, u_i, \phi_i)$.
It is clear that if there exists a family of
coordinates $(\phi'_i)_{+\infty^{i, +\infty}_l}$ close to
$(\phi_i)_{+\infty^{i, +\infty}_l}$ such that
$\pi_Y \circ u'_i \circ (\phi'_i)_{+\infty^{i, +\infty}_l} \in |\tilde \eta^i_l|$ ($\subset P$),
then such a family is unique.
(Furthermore, if the restriction of $\pi_Y \circ u'_i$ to the $+\infty$-limit circle
corresponding to $\eta^i_l$ is contained in $|\eta^i_l|$ ($\subset \overline{P}$),
then there exists a coordinate $(\phi'_i)_{+\infty^{i, +\infty}_l}$ close to $(\phi_i)_{+\infty^{i, +\infty}_l}$
or its rotation by some element of $\Z/d \subset S^1$ such that
$\pi_Y \circ u'_l \circ (\phi'_i)_{+\infty^{i, +\infty}_l} \in |\tilde \eta^i_l|$,
where $d$ is the multiplicity of $\gamma^{i, +}_l$.)
The same is true for the coordinates of the limit circles corresponding to
$c^i_l$ or $\epsilon^{i, j}_l$.

Let
\[
(V^1 \times V^2 \times \dots \times V^m)_{((\breve \epsilon^{i, j}_l), (\breve c^i_l),
(x^i_j), (\breve \eta^i_l))}
\subset V^1 \times V^2 \times \dots \times V^m
\]
be the submanifold consisting of the families of curves which have families of
coordinates of
their limit circles close to that of $(\Sigma_i, z_i, u_i, \phi_i)_{1 \leq i \leq m}$ such that
the periodic orbits on the $\pm\infty$-limit circles are contained in the corresponding
$\tilde \epsilon^{i, j}_l$, $\tilde c^i_j$ and $\tilde \eta^i_l$,
and $\pi_y \circ u'_i$ takes a value in $x^i_l$ at each marked point $z_{i, l}$.
This submanifold can be regarded as a fiber product of
$V^1 \times V^2 \times \dots \times V^m$ with
the product of $(I_\delta \times I_\delta) \cdot \tilde \epsilon^{i, j}_l$,
$I_\delta \cdot \tilde c^i_l$, $x^i_l$ and $I_\delta \cdot \tilde \eta^i_l$,
where $I_\delta \subset S^1$ is a small neighborhood of $0 \in S^1$.

Then a Kuranishi neighborhood of $(\Sigma_i, z_i, u_i, \phi_i)_{1 \leq i \leq m}$ is defined by
this submanifold, the restrictions of the product vector bundle
$E = E^1 \times E^2 \times \dots \times E^m$, its section
$s = s^1 \times s^2 \times \dots \times s^m$,
a finite group $G = \prod_i \Aut(\Sigma_i, z_i, u_i, \phi_i)$, and
the map
\[
\psi : s^{-1}(0)/G \to \overline{\M}^m
_{((\breve \epsilon^{i, j}_l), (\breve c^i_l), (x^i_l), (\breve \eta^i_l))}
\]
induced by the product $\psi^1 \times \psi^2 \times \dots \times \psi^m$ and
the coordinates of the limit circles
defined by the above argument.

Note that the pre-Kuranishi spaces for other lifts of $c^i_l$, $\eta^i_l$ or $\epsilon^{i, j}_l$
are naturally isomorphic to the above pre-Kuranishi space.
We need to construct their perturbed multisections which are independent of the choice
of the lifts.
We construct the perturbed multisections of the above fiber products
as pull backs by submersions to the fiber products of $\widehat{\M}$
in the next section.
Since these submersions forget the coordinates of limit circles,
the pull backs will be independent of the choice of the lifts of simpleces.

\subsection{Construction of a family of multisections}
\label{construction of a multisection}
In this section, we define fiber products of $\widehat{\M}$ and
construct their grouped multisections under appropriate compatibility conditions.
Recall that in Section \ref{global construction},
we construct a pre-Kuranishi structure of
$\widehat{\M}^{\leq L_{\max}}_{\leq C}$ only for fixed constants
$L_{\max} > 0$ and $C >0$.
We cannot treat the whole noncompact space $\widehat{\M}$ as a pre-Kuranishi
space.
Hence we need to read $\widehat{\M}$ below as its compact subset
$\widehat{\M}^{\leq L_{\max}}_{\leq C}$ for some $L_{\max} > 0$ and $C >0$.
We also read the other spaces as their similar compact subsets.

First we define the space $\widehat{\M}^\diamond$.
Its point $((\Sigma^\alpha, z^\alpha, u^\alpha)_{\alpha \in A}, M^{\mathrm{rel}})$
consists of finite number of connected holomorphic buildings
$(\Sigma^\alpha, z^\alpha, u^\alpha)$ and
a set $M^{\mathrm{rel}} = \{(S^1_{+\infty_l}, S^1_{-\infty_l})\}$ of pairs of
their $+\infty$-limit circle $S^1_{+\infty_l}$ and $-\infty$-limit circle $S^1_{-\infty_l}$
which satisfies the following conditions:
\begin{itemize}
\item
Any two pairs in $M^{\mathrm{rel}}$ do not share the same limit circles.
\item
Let $M^{\alpha, \alpha'} \subset M^{\mathrm{rel}}$ be the subset of pairs
$(S^1_{+\infty_l}, S^1_{-\infty_l})$ such that $S^1_{+\infty_l}$ is a $+\infty$-limit circle of
$\Sigma^\alpha$ and $S^1_{-\infty_l}$ is a $-\infty$-limit circle of $\Sigma^{\alpha'}$.
Then there does not exist a sequence
$\alpha_0, \alpha_1, \dots, \alpha_k = \alpha_0 \in A$
such that $M^{\alpha_i, \alpha_{i+1}} \neq \emptyset$ for all $i$.
\end{itemize}
(We emphasize again the meaning of connectedness of a holomorphic building.
The domain curve $\Sigma^\alpha$ of each holomorphic building
$(\Sigma^\alpha, z^\alpha, u^\alpha)$ is a compact space which contains
joint circles and limit circles.
The connectedness of $(\Sigma^\alpha, z^\alpha, u^\alpha)$ means that of
$\Sigma^\alpha$, and it does not imply that $(\Sigma^\alpha, z^\alpha, u^\alpha)$
is of height-one or that its restriction to each floor is connected.
We also note that the height of $(\Sigma^\alpha, z^\alpha, u^\alpha)$
may depend on $\alpha$, and that the total
$((\Sigma^\alpha, z^\alpha, u^\alpha)_{\alpha \in A}, M^{\mathrm{rel}})$
does not have a floor structure.)

Two points $((\Sigma^\alpha, z^\alpha, u^\alpha)_{\alpha \in A}, M^{\mathrm{rel}})$ and
$(((\Sigma')^{\alpha'}, (z')^{\alpha'}, (u')^{\alpha'})_{\alpha' \in A'}, (M')^{\mathrm{rel}})$ are
the same point if there exist a bijection $\nu : A \to A'$,
isomorphisms $\varphi^\alpha : \Sigma^\alpha \to (\Sigma')^{\nu(\alpha)}$,
and $\R$-translations $\theta^\alpha$ such that
$\varphi^\alpha(z^\alpha) = (z')^{\nu(\alpha)}$,
$u^\alpha = (\theta^\alpha \times 1) \circ (u')^{\nu(\alpha)} \circ \varphi^\alpha$
and the family of isomorphisms $\varphi^\alpha$ maps $M^{\mathrm{rel}}$
to $(M')^{\mathrm{rel}}$.
Forgetting $M^{\mathrm{rel}}$ defines a forgetful map from
$\widehat{\M}^\diamond$ to $\bigcup_N \prod^N (\widehat{\M}^0) / \mathfrak{S}_N$,
where $\widehat{\M}^0 \subset \widehat{\M}$ is the space of connected holomorphic
buildings.
Since the only local difference of these two spaces are automorphism group,
$\widehat{\M}^\diamond$ has the natural pre-Kuranishi structure which makes this
forgetful map a submersion.

For subsets $A_1, A_2 \subset A$, we define
$M^{A_1, A_2} = \bigcup_{\alpha_1 \in A_1, \alpha_2 \in A_2} M^{\alpha_1, \alpha_2}$.
We say a point $((\Sigma^\alpha, z, u^\alpha)_{\alpha \in A}, M^{\mathrm{rel}})
\in \widehat{\M}^\diamond$ is disconnected if
there exists a decomposition $A = A_1 \sqcup A_2$ such that
$M^{A_1, A_2} = M^{A_2, A_1} = \emptyset$.
Otherwise we say it is connected.
We denote the space of connected points of $\widehat{\M}^\diamond$ by
$(\widehat{\M}^\diamond)^0$.
Decomposition into connected components defines a map $\widehat{\M}^\diamond
\to \bigcup_N (\prod^N (\widehat{\M}^\diamond)^0) / \mathfrak{S}_N$.


Let $\Upsilon : \widehat{\M}^\diamond \to
(\prod (\overline{P} \times \overline{P})) / \mathfrak{S} \times
\prod \overline{P} / \mathfrak{S} \times \prod Y / \mathfrak{S}$
be the continuous map which maps each point
$((\Sigma^\alpha, z^\alpha, u^\alpha)_{\alpha \in A}, M^{\mathrm{rel}})
\in \widehat{\M}^\diamond$
to
\begin{align*}
((\pi_Y \circ u|_{S^1_{+\infty_l}}, \pi_Y \circ u|_{S^1_{-\infty_l}})
_{(S^1_{+\infty_l}, S^1_{-\infty_l}) \in M^{\mathrm{rel}}},
(\pi_Y \circ u|_{S^1_{\pm\infty}})_{S^1_{\pm\infty} \notin M^{\mathrm{rel}}}, \quad \\
(\pi_Y \circ u(\bigcup_\alpha z^\alpha))),
\end{align*}
where we denote the union of $u^\alpha$ by $u$, and
$S^1_{\pm\infty} \notin M^{\mathrm{rel}}$ means that the limit circle $S^1_{\pm\infty}$
is not contained in any pairs in $M^{\mathrm{rel}}$.
%
%
%
%
It is realized as a strong smooth map.
(The number of the product of the target space depends on the components of
$\widehat{\M}^\diamond$.)
Define the fiber product $(\widehat{\M}^\diamond, \mathring{K}^2, K,K^0) \subset
\widehat{\M}^\diamond$ by $(\widehat{\M}^\diamond, \mathring{K}^2, K,K^0)
= \Upsilon^{-1}(\prod \mathring{K}^2 / \mathfrak{S} \times \prod K / \mathfrak{S}
\times \prod K^0 / \mathfrak{S})$.
(Although $\prod K / \mathfrak{S}$ coincides with
the entire space $\prod P / \mathfrak{S}$,
the fiber product with $\prod K / \mathfrak{S}$ has a meaning.
It implies that we require the stronger transversality condition
for its grouped multisection. See Remark \ref{meaning of triangulation}.
Fiber product with $\prod K^0 / \mathfrak{S} \subset \prod Y / \mathfrak{S}$
is slightly an abuse of notation.
A simplicial complex is a set of simplices, but $K^0$ is a set of smooth cycles, and
we do not assume that they are embedded in $Y$.
The meaning is also that we require that the perturbed multisection
is transverse to the zero section even if we restrict it to the fiber product with
the cycles in $K^0$.)

Define a multi-valued partial submersion $\Xi : \widehat{\M}^\diamond \to
\widehat{\M}^\diamond$ by
\[
\Xi(((\Sigma^\alpha, z^\alpha, u^\alpha)_{\alpha \in A}, M^{\mathrm{rel}}))
= \{((\Sigma^\alpha, z^\alpha, u^\alpha)_{\alpha \in A},
\mathring{M}^{\mathrm{rel}}); \mathring{M}^{\mathrm{rel}} \subsetneq M^{\mathrm{rel}}\}.
\]
(The maps between their Kuranishi neighborhoods are also similarly defined.)
Let $\Xi^\circ$ be the restriction of
$\Xi \subset \widehat{\M}^\diamond \times \widehat{\M}^\diamond$ to
the set of points
\[
(((\Sigma^\alpha, z^\alpha, u^\alpha)_{\alpha \in A}, M^{\mathrm{rel}}),
((\Sigma^\alpha, z^\alpha, u^\alpha)_{\alpha \in A}, \mathring{M}^{\mathrm{rel}}))
\]
such that $((\Sigma^\alpha, z^\alpha, u^\alpha)_{\alpha \in A}, M^{\mathrm{rel}})
\in (\widehat{\M}^\diamond, \mathring{K}^2, K,K^0)$ and
$(\pi_Y \circ u|_{S^1_{+\infty_l}}, \pi_Y \circ u|_{S^1_{-\infty_l}})
\in \rho_\ast K$
for all $(S^1_{+\infty_l}, S^1_{-\infty_l}) \in M^{\mathrm{rel}} \setminus
\mathring{M}^{\mathrm{rel}}$.
Then $\Xi^\circ$ is a multi-valued partial essential submersion
from $(\widehat{\M}^\diamond, \mathring{K}^2, K,K^0)$ to itself.
\begin{rem}
\label{simple restriction of Xi is not essentially submersive}
Note that the restriction of $\Xi$ to $(\widehat{\M}^\diamond, \mathring{K}^2, K,K^0)$
is not a multi-valued partial essential submersion from
$(\widehat{\M}^\diamond, \mathring{K}^2, K,K^0)$ to itself.
This is because $\mathring{K}^2$ is finer than the product $K \times K$
on the outside of $\rho_\ast K$.
This is an important point in Bott-Morse case.
To explain this point, we consider the following easy example.
Let $s_i : \Delta^n \to \R^l$ ($i = 1,2$) be sections of the
trivial vector bundle on the $n$-simplex $\Delta^n$ whose restriction to each of
it faces is transverse to
the zero section.
These transversality conditions of $s_i$ do not imply that
the restriction of $(s_1 + s_2) : \Delta^n \times \Delta^n \to
\R^l \oplus \R^l$ to the diagonal is transverse to the zero section.
However, the restrictions to the cells in
$\rho_\ast \Delta^n
= \{\partial_{p+1} \dots \partial_n \Delta^n \times \partial_0 \dots \partial_{p-1} \Delta^n;
0 \leq p \leq n\}$
are transverse to the zero section.
Hence we use the induced section $s_1 + s_2$ only on $\rho_\ast \Delta^n$.
\end{rem}

We also define a multi-valued partial essential submersion
\[
\Lambda : (\partial \widehat{\M}^\diamond, \mathring{K}^2, K,K^0)
\to (\widehat{\M}^\diamond, \mathring{K}^2, K,K^0)
\]
as follows.
For each point
$((\Sigma^\alpha, z^\alpha, u^\alpha)_{\alpha \in A}, M^{\mathrm{rel}}) \in
(\partial \widehat{\M}^\diamond, \mathring{K}^2, K,K^0)$
and each $\alpha \in A$, let $k_\alpha \geq 1$ be the height of
$(\Sigma^\alpha, z^\alpha, u^\alpha)$.
Let $C = \{(\alpha_i, k_i)\}$ be an arbitrary non-empty set of pairs
$(\alpha_i, k_i) \in A \times \N$ such that $1 \leq k_i < k_{\alpha_i}$.
For all pairs $(\alpha_i, k_i) \in C$,
we replace all joint circles in the gap between the $k_i$-th floor and
the $(k_i + 1)$-th floor of $(\Sigma^{\alpha_i}, z^{\alpha_i}, u^{\alpha_i})$
with pairs of limit circles,
and let $(\Sigma^{\alpha'}, z^{\alpha'}, u^{\alpha'})_{\alpha' \in A'}$ be the
stabilization of the connected components of the new holomorphic buildings.
(Stabilization means we collapse all trivial floors of each connected component.)
$\Lambda : (\partial \widehat{\M}^\diamond, \mathring{K}^2, K,K^0)
\to (\widehat{\M}^\diamond, \mathring{K}^2, K,K^0)$ is the multi-valued partial
submersion which maps each point
$((\Sigma^\alpha, z^\alpha, u^\alpha)_{\alpha \in A}, M^{\mathrm{rel}}) \in
(\partial \widehat{\M}^\diamond, \mathring{K}^2, K,K^0)$
to the above points $(\Sigma^{\alpha'}, z^{\alpha'}, u^{\alpha'})_{\alpha' \in A'}$
for all $C = \{(\alpha_i, k_i)\}$.
(The maps between their Kuranishi neighborhoods are also similarly defined.)

For each point $p = ((\Sigma^\alpha, z^\alpha, u^\alpha)_{\alpha \in A},
M^{\mathrm{rel}}) \in (\widehat{\M}^\diamond, \mathring{K}^2, K,K^0)$,
we define
$\widetilde{e}(p) = \widetilde{e}_{\delta_0}(p)
= \sum_\alpha \widetilde{e}_{\delta_0}(\theta_\alpha)
+ \frac{1}{2} \# M^{\mathrm{rel}}$, where each
$\theta_\alpha$ is the type of $(\Sigma^\alpha, z^\alpha, u^\alpha)$,
and $\# M^{\mathrm{rel}}$ is the number of pairs.
($\# M^{\mathrm{rel}}$ is not the number of limit circles contained in the pairs in
$M^{\mathrm{rel}}$.
Recall that $\widetilde{e}_{\delta_0}(\theta) = 5(g-1) + 2k + E_{\hat \omega} / \delta_0$
for $\theta = (g, k, E_{\hat \omega})$, where
$g$ is the genus, $k$ is the total number of marked points and limit circles,
and $E_{\hat \omega}$ is the $E_{\hat \omega}$-energy.
See Section \ref{global construction} for the definition
and properties of
$\widetilde{e}_{\delta_0}(\theta_\alpha)$.)
Note that the maps $\Xi^\circ$ and $\Lambda$ decrease $\widetilde{e}$.
Hence if we decompose $(\widehat{\M}^\diamond, \mathring{K}^2, K,K^0)$ by
$\widetilde{e}$, then
$\Xi^\circ$ and $\Lambda$ constitute a compatible system of
multi-valued partial submersions.
(See Section \ref{multi-valued partial submersions} for the definition of
a compatible system of multi-valued partial submersions.)

We can construct the perturbed multisections of
$(\widehat{\M}^\diamond, \mathring{K}^2, K,K^0)$
which satisfy the following conditions:
\begin{itemize}
\item
Let $((\widehat{\M}^\diamond)^0, \mathring{K}^2, K, K^0) \subset
(\widehat{\M}^\diamond, \mathring{K}^2, K,K^0)$ be the subset of connected points.
Its grouped multisection induces that of
\[
\bigcup_N (\prod^N ((\widehat{\M}^\diamond)^0, \mathring{K}^2, K, K^0)
/ \mathfrak{S}_N.
\]
Then the grouped multisection of $(\widehat{\M}^\diamond, \mathring{K}^2, K,K^0)$
coincides with its pull back by the submersion
\[
(\widehat{\M}^\diamond, \mathring{K}^2, K,K^0) \to
\bigcup_N (\prod^N ((\widehat{\M}^\diamond)^0, \mathring{K}^2, K, K^0)
/ \mathfrak{S}_N
\]
defined by decomposition into connected components.
\item
The grouped multisections of
$(\widehat{\M}^\diamond, \mathring{K}^2, K,K^0)$
are compatible with respect to the compatible system of
multi-valued partial essential submersions defined by $\Xi^\circ$ and $\Lambda$.
\end{itemize}

Let
$\overline{\M}^m_{((\breve \epsilon^{i, j}_l), (\breve c^i_l), (x^i_l), (\breve \eta^i_l))}
\to (\widehat{\M}^\diamond, \mathring{K}^2, K,K^0)$
be the natural map defined by
decomposing each holomorphic building $(\Sigma_i, z_i, u_i, \phi_i)$
into its connected components (and stabilizing them),
forgetting the order of the sequence of holomorphic buildings,
the coordinates of limit circles and the order of marked points
and limit circles.
$M^{\mathrm{rel}}$ is defined by the set of pairs of limit circles
corresponding to the pairs $(S^1_{+\infty^{i,j}_l}, S^1_{-\infty^{j,i}_l})$.
We define the perturbed multisection of each
$\overline{\M}^m_{((\breve \epsilon^{i, j}_l), (\breve c^i_l), (x^i_l), (\breve \eta^i_l))}$
by the pull back by this natural submersion.

We emphasize the following point.
Although in Section \ref{fiber prod with simpleces},
we only consider the fiber products with simpleces or cells not contained
in $\overline{P}^{\text{bad}}$, $\overline{P}^{^t \text{bad}}$ or 
$\overline{P}^{^t\text{bad}} \times \overline{P} \cup
\overline{P} \times \overline{P}^{\text{bad}}$, in this section,
we construct the perturbed multisections of the fiber products
of $\widehat{\M}^\diamond$ with all simpleces or cells.
We cannot ignore the bad orbits for the construction of the compatible family of
perturbed multisections, but for the construction of the algebra, we only use
the orientable fiber products.

\begin{rem}
In the above construction of grouped multisection, we do not
consider any compatibility condition with respect to
the decompositions at nodal points.
This is just because we do not use this kind of compatibility condition
for the construction of the algebras in this paper.
We can construct a grouped multisection with these compatibility conditions
if we have an appropriate formulation of invariants or algebras
which respect the decompositions at nodal points.

To regard the nodal point as a fiber product over $Y$ by the
maps of evaluation at the marked points, we need to use the
grouped multisection obtained as the fiber product.
As we saw in Remark \ref{simple restriction of Xi is not essentially submersive},
the fiber product of grouped multisections do not satisfy the
transversality condition in general.
There are two ways to solve this problem.
One is the use of continuous family of grouped multisections.
See Section \ref{continuous family of multisections} for details.
(See also Section \ref{fiber product and orientation for homotopy}.
We use this technique for the the case of a $1$-parameter family of
symplectic manifolds with cylindrical ends.)
The other is to use the same argument as the decomposition
by the gaps of floors as above.
In this case, we use the approximation of the diagonal,
and we need to use the correction terms for the difference of
this approximation and the genuine diagonal similar to
those we construct in Section \ref{algebra for correction}.
\end{rem}

\subsection{The orientations of fiber products}\label{fiber prod and orientation}
In this section, we define the orientations of the fiber products
$\overline{\M}^m_{((\breve \epsilon^{i, j}_l), (\breve c^i_l), (x^i_l), (\breve \eta^i_l))}$.
For calculation of orientations, it is convenient to treat these pre-Kuranishi spaces
as fiber products not with $I_\delta \cdot \tilde c^i_l$ but with ``manifold''
$(\R \times S^1) \cdot \tilde c^i_l \times \theta^D_{c^i_l}$.

First we define the orientation of the parameter space
$\mathring{X} \subset \widetilde{X} \times \prod \R_\mu$
used for the construction of the Kuranishi neighborhoods.
For each $i = 1, 2, \dots, k-1$,
we fix one joint circle $S^1_{\mu_i}$ between the $i$-th floor and the $(i+1)$-th
floor.
First we consider the orientation at a point
$(a, b) \in \mathring{X}$ such that $\rho_\mu \neq 0$ for all joint circles $S^1_\mu$.
In this case, we can use $(b_{\mu_i})_i$ and a chart of $\widetilde{X}$ as a chart of $\mathring{X}$
on a neighborhood of this point.
We define the orientation of $\mathring{X}$ by this chart
\[
(b_{\mu_1}, b_{\mu_2}, \dots, b_{\mu_{k-1}}, a)
\in \R \times \R \times \dots \times \R \times \widetilde{X},
\]
where the orientation of each $\R$ is the positive orientation, and the orientation of
$\widetilde{X}$ is the complex orientation
(the orientation induced by the complex orientation of the blow down space).

At a general point $(a, b) \in \mathring{X}$, the orientation of $\mathring{X}$ is defined
as follows.
We can write a point of $\widetilde{X}$ as $a = ((\rho_\mu, \varphi_\mu)_\mu, a^0)$,
where each $(\rho_\mu, \varphi_\mu) \in [0, 1) \times S^1$
is the parameter of the deformation of
the neighborhood of the joint circle $S^1_\mu$, and $a^0$ is the other parameters.
Note that the parameter $a^0$ is a complex parameter.
Then we can use
\[
(-\rho_{\mu_1}, -\rho_{\mu_2}, \dots, -\rho_{\mu_{k-1}},
(b_\mu, \varphi_\mu)_{\mu}, a^0)
\]
as a chart of $\mathring{X}$.
(More precisely, each $\rho_\mu$ should be replaced with its appropriate multiple.
See Section \ref{smoothness}.)
The orientation of $\mathring{X}$ is defined by this chart.
(The order of $\mu$ of $(b_\mu, \varphi_\mu)_{\mu}$ is independent of the orientation
since each $(b_\mu, \varphi_\mu)$ is even dimensional.)

It is easy to see that the two definition of the orientation coincide.
It is also easy to see that the orientation is independent of the choice of $\mu_i$.

To define the orientation of $(V^1 \times V^2 \times \dots \times
V^m)_{((\breve \epsilon^{i, j}_l), (\breve c^i_l), (x^i_j), (\breve \eta^i_l))}$,
we need to see the relations of its tangent space and other various vector spaces.
Recall that $T \hat V^i = T\mathring{X}^i \oplus \Ker D^i$ for each $1 \leq i \leq k$,
where
\[
D^i : \widetilde{W}^{1, p}_\delta (\Sigma_i, u_i^\ast T \hat Y) \oplus E_i^0 \to
L_\delta^p (\Sigma_i, \Wedge^{0, 1} T^\ast \Sigma_i \otimes u_i^\ast T \hat Y)
\]
is the differential $DF^{(0, b^0)}_{(0, 0)}$ in Section \ref{construction of nbds}
for $(\Sigma_i, z_i, u_i)$.
The fiber product
\[
\mathcal{C}^i = [D^{i, -}] \underset{\Ker A^{i, -}}{\times} [D^i]
\underset{\Ker A^{i, +}}{\times} [D^{i, +}]
\]
is equivalent to the kernel of a $\bound$-operator on a complex vector bundle over
a closed semistable curve by gluing, where
\begin{align*}
[D^{i, -}] &= \prod_{1 \leq j < i, 1 \leq l \leq l_{j, i}} [D^-_{\gamma_{-\infty^{i, j}_l}}]
\times \prod_{1 \leq l \leq l_{i, -}} [D^-_{\gamma_{-\infty^{i, -\infty}_l}}],\\
\Ker A^{i, -} &= \prod_{1 \leq j < i, 1 \leq l \leq l_{j, i}}
\Ker A_{\gamma_{-\infty^{i, j}_l}}
\times \prod_{1 \leq l \leq l_{i, -}} \Ker A_{\gamma_{-\infty^{i, -\infty}_l}},\\
\Ker A^{i, +} &= \prod_{i < j \leq m, 1 \leq l \leq l_{i, j}}
\Ker A_{\gamma_{+\infty^{i, j}_l}}
\times \prod_{1 \leq l \leq l_{i, +}} \Ker A_{\gamma_{+\infty^{i, +\infty}_l}},\\
[D^{i, +}] &= \prod_{i < j \leq m, 1 \leq l \leq l_{i, j}} [D^+_{\gamma_{+\infty^{i, j}_l}}]
\times \prod_{1 \leq l \leq l_{i, +}} [D^+_{\gamma_{+\infty^{i, +\infty}_l}}].
\end{align*}
The vector space $[\mathring{D}^{i, -}] \oplus [\mathring{D}^{i, +}]$
is a subspace of $\mathcal{C}^i$, and its quotient space is isomorphic to $[D^i]$,
where
\begin{align*}
&[\mathring{D}^{i, -}] = \prod_{1 \leq j < i, 1 \leq l \leq l_{j, i}}
[\mathring{D}^-_{\gamma_{-\infty^{i, j}_l}}] \times
\prod_{1 \leq l \leq l_{i, -}} [\mathring{D}^-_{\gamma_{-\infty^{i, -\infty}_l}}],\\
&[\mathring{D}^{i, +}] = \prod_{i < j \leq m, 1 \leq l \leq l_{i, j}}
[\mathring{D}^+_{\gamma_{+\infty^{i, j}_l}}] \times
\prod_{1 \leq l \leq l_{i, +}} [\mathring{D}^+_{\gamma_{+\infty^{i, +\infty}_l}}].
\end{align*}
We fix splittings $\mathcal{C}^i \to [\mathring{D}^-_{\gamma_{-\infty^{i, j}_l}}]$,
$\mathcal{C}^i \to [\mathring{D}^-_{\gamma_{-\infty^{i, -\infty}_l}}]$,
$\mathcal{C}^i \to [\mathring{D}^+_{\gamma_{+\infty^{i, j}_l}}]$
and $\mathcal{C}^i \to [\mathring{D}^+_{\gamma_{+\infty^{i, +\infty}_l}}]$.
Then the tangent space of
$(V^1 \times V^2 \times \dots \times V^m)_{((\breve \epsilon^{i, j}_l),
(\breve c^i_l), (x^i_j), (\breve \eta^i_l))}$
is isomorphic to the kernel of the map from
\[
(T\mathring{X}^1 \times \mathcal{C}^1) \oplus (T\mathring{X}^2 \times \mathcal{C}^2) \oplus \dots
\oplus (T\mathring{X}^m \times \mathcal{C}^m)
\]
to the direct sum of the following vector spaces:
\begin{itemize}
\item
$T \hat Y / (\R\partial_\sigma \oplus T x^i_l)$
\item
$[\mathring{D}^-_{\gamma_{-\infty^{i, -\infty}_l}}] \oplus
\bigl(\Ker A_{\gamma_{-\infty^{i, -\infty}_l}} / (\R \oplus T S^1 \oplus T \tilde c^i_l)\bigr)$
\item
$\bigl(\Ker A_{\gamma_{+\infty^{i, +\infty}_l}} / (\R \oplus T S^1 \oplus T \tilde \eta^i_l)\bigr)
\oplus [\mathring{D}^+_{\gamma_{+\infty^{i, +\infty}_l}}]$
\item
$[\mathring{D}^+_{\gamma_{+\infty^{i, j}_l}}] \oplus
\bigl((\Ker A_{\gamma_{+\infty^{i, j}_l}} \oplus \Ker A_{\gamma_{-\infty^{i, j}_l}}) /
(\R \oplus T S^1 \oplus T \tilde \epsilon^{i, j}_l \oplus \R \oplus T S^1)\bigr)
\oplus [\mathring{D}^-_{\gamma_{-\infty^{i, j}_l}}]$
\item
$\R^{k_i} \oplus \bigoplus_{z_{i, \beta}^{++}} \R^2$ (the range of the map $s^0$)
\end{itemize}

Now we explain the definition of the orientations.
For the convenience, we assume all $E^0$ are complex vector spaces and
$\lambda^0$ are complex linear map. (We can always assume this condition.)

First we define the orientation of the vector space
\[
\W^i = T\mathring{X}^i \times \mathcal{C}^i / (\R^{k_i} \oplus \bigoplus_{z_{i, \beta}^{++}} \R^2)
\]
by
\[
\bigl(\R^{k_i} \oplus \bigoplus_{z_{i, \beta}^{++}} \R^2\bigr) \oplus \W^i = T\mathring{X}^i \times \mathcal{C}^i,
\]
where the orientation of $\R^{k_i}$ is the product of the positive orientation of $\R$,
and its order is
\[
(\sigma_1 \circ \Phi_{a, b}(\xi) (R_1), \sigma_2 \circ \Phi_{a, b}(\xi) (R_2), \dots, \sigma_{k_i} \circ \Phi_{a, b}(\xi) (R_{k_i}));
\]
and the orientation of each $\R^2 \cong T^{\bot^Y}\!\! S = T \hat Y / T (\R \times S)$
is the complex orientation defined by the almost complex structure of $\hat Y$.
We note that the dimension of each $\W^i$ is odd.

Next we recall the definition of the orientation of fiber product.
We use only fiber products from right.
Let $f$ be a surjective linear map from an oriented vector space $V$ to another oriented
vector space $W$, and let $A \subset W$ be an oriented subspace.
Define the orientation of $\Ker f$ by $\Ker f \oplus W = V$.
Then the orientation of the fiber product $V_A = \{v \in V; f(v) \in A\}$ is defined by
$V_A = \Ker f \oplus A$.

We deform each curve $\Sigma_i$, preserving a neighborhood of limit circles,
to a curve which consists of inner caps $[-\infty, 0] \times S^1 \cup D_\infty$ of
$-\infty$-limit circles, inner caps $D_0 \cup [0, \infty] \times S^1$ of
$+\infty$-limit circles, and a semistable curve $\mathring{\Sigma}_i$.
Each inner cap is connected to $\mathring{\Sigma}_i$ by a nodal point.
We also deform the linear operator $D^i$ to a $\bound$-type linear operator
which coincide with
\begin{itemize}
\item
$D^+_{\gamma_{-\infty^{i, j}_l}}$ on the inner cap of $S^1_{-\infty^{i, j}_l}$,
\item
$D^+_{\gamma_{-\infty^{i, -\infty}_l}}$ on the inner cap of $S^1_{-\infty^{i, -\infty}_l}$,
\item
$D^-_{\gamma_{+\infty^{i, j}_l}}$ on the inner cap of $S^1_{+\infty^{i, j}_l}$, and
\item
$D^-_{\gamma_{+\infty^{i, +\infty}_l}}$ on the inner cap of $S^1_{+\infty^{i, +\infty}_l}$.
\end{itemize}
We denote the deformed linear operator by $\widetilde{D}^i$
In the definition of $\mathcal{C}$ and $\W^i$, we replace $D^i$ with $\widetilde{D}^i$,
and we get vector spaces $\widetilde{\mathcal{C}}$ and $\widetilde{\W}^i$.

Restriction to the inner cups defines linear maps from
$[\widetilde{D}^i]$ to $[D^+_{\gamma_{-\infty^{i, j}_l}}]$, $[D^+_{\gamma_{-\infty^{i, -\infty}_l}}]$,
$[D^-_{\gamma_{+\infty^{i, j}_l}}]$ and $[D^-_{\gamma_{+\infty^{i, +\infty}_l}}]$.
Therefore there exist linear maps from $\widetilde{\mathcal{C}}^1 \oplus
\widetilde{\mathcal{C}}^2 \oplus \dots \oplus \widetilde{\mathcal{C}}^m$
to the following vector spaces.
\begin{align*}
&[D^-_{\gamma_{-\infty^{i, j}_l}}] \underset{\Ker A_{\gamma_{-\infty^{i, j}_l}}}{\times}
[D^+_{\gamma_{-\infty^{i, j}_l}}],
&&[D^-_{\gamma_{-\infty^{i, -\infty}_l}}] \underset{\Ker A_{\gamma_{-\infty^{i, -\infty}_l}}}{\times}
[D^+_{\gamma_{-\infty^{i, -\infty}_l}}],
\\
&[D^-_{\gamma_{+\infty^{i, j}_l}}] \underset{\Ker A_{\gamma_{+\infty^{i, j}_l}}}{\times}
[D^+_{\gamma_{+\infty^{i, j}_l}}],
&&[D^-_{\gamma_{+\infty^{i, +\infty}_l}}] \underset{\Ker A_{\gamma_{+\infty^{i, +\infty}_l}}}{\times}
[D^+_{\gamma_{+\infty^{i, +\infty}_l}}].
\end{align*}
They define linear maps from $\widetilde{\W}^1 \oplus \widetilde{\W}^2
\oplus \dots \oplus \widetilde{\W}^m$ to the same vector spaces.
The above vector spaces have the complex orientations.

We regard the vector space
\[
\R \oplus TS^1 \oplus T \tilde c^i_l \oplus [\mathring{D}^+_{\gamma_{-\infty^{i, -\infty}_l}}]
\]
as a subspace of
\[
[D^-_{\gamma_{-\infty^{i, -\infty}_l}}] \underset{\Ker A_{\gamma_{-\infty^{i, -\infty}_l}}}{\times}
[D^+_{\gamma_{-\infty^{i, -\infty}_l}}]
\]
by regarding $\R \oplus TS^1 \oplus T \tilde c^i_l \subset \R \oplus TP \cong
\Ker A_{\gamma_{-\infty^{i, -\infty}_l}}$
as a subspace of the above space by a right inverse of the surjection
\[
[D^-_{\gamma_{-\infty^{i, -\infty}_l}}] \underset{\Ker A_{\gamma_{-\infty^{i, -\infty}_l}}}{\times}
[D^+_{\gamma_{-\infty^{i, -\infty}_l}}] \to \Ker A_{\gamma_{-\infty^{i, -\infty}_l}}.
\]
We define its orientation
by the direct sum of the complex orientation of $\R \oplus TS^1 = T(\R \times S^1)$,
the orientation of the ordered simplicial complex $\tilde c^i_l$,
and the orientation $\theta^D_{c^i_l}$ of $[\mathring{D}^+_{\gamma_{-\infty^{i, -\infty}_l}}]$.

We define the orientations of the following spaces similarly.
\[
[\mathring{D}^-_{\gamma_{+\infty^{i, +\infty}_l}}] \oplus \R \oplus TS^1 \oplus T \tilde \eta^i_l
\subset [D^-_{\gamma_{+\infty^{i, +\infty}_l}}] \underset{\Ker A_{\gamma_{+\infty^{i, +\infty}_l}}}{\times}
[D^+_{\gamma_{+\infty^{i, +\infty}_l}}],
\]
\vspace{-12pt}
\begin{align*}
&[\mathring{D}^-_{\gamma_{+\infty^{i, j}_l}}] \oplus \R \oplus T S^1 \oplus T \tilde \epsilon^{i, j}_l \oplus \R \oplus T S^1 \oplus [\mathring{D}^+_{\gamma_{-\infty^{j, i}_l}}]\\
&\subset
([D^-_{\gamma_{+\infty^{i, j}_l}}] \underset{\Ker A_{\gamma_{+\infty^{i, j}_l}}}{\times}
[D^+_{\gamma_{+\infty^{i, j}_l}}]) \\
&\quad \oplus
([D^-_{\gamma_{-\infty^{j, i}_l}}] \underset{\Ker A_{\gamma_{-\infty^{j, i}_l}}}{\times}
[D^+_{\gamma_{-\infty^{j, i}_l}}]).
\end{align*}

For each marked point $z_{i, l}$, the orientation of $T x^i_l \oplus \R \subset
T \hat Y$ is defined by the orientation of the cycle $x^i_l$
and the positive orientation of $\R$.
The orientation of $T \hat Y$ is defined by the complex orientation.
The tangent space of
$(V^1 \times V^2 \times \dots \times V^m)_{((\breve \epsilon^{i, j}_l), (\breve c^i_l),
(x^i_j), (\breve \eta^i_l))}$
is isomorphic to the fiber product
\begin{equation}
(\widetilde{\W}^1 \oplus \widetilde{\W}^2 \oplus \dots \oplus \widetilde{\W}^m)_\star,
\label{ori}
\end{equation}
where
\begin{align*}
\star &= \bigoplus_{(i, j), l}
([\mathring{D}^-_{\gamma_{+\infty^{i, j}_l}}] \oplus \R \oplus T S^1 \oplus
T \tilde \epsilon^{i, j}_l \oplus \R \oplus T S^1 \oplus
[\mathring{D}^+_{\gamma_{-\infty^{j, i}_l}}])\\
&\quad \oplus \bigoplus_{i, l}
(\R \oplus TS^1 \oplus T \tilde c^i_l \oplus [\mathring{D}^+_{\gamma_{-\infty^{i, -\infty}_l}}])
\oplus \bigoplus_{i, l}
(\R\partial_\sigma \oplus T x^i_l)\\
&\quad \oplus \bigoplus_{i, l}
([\mathring{D}^-_{\gamma_{+\infty^{i, +\infty}_l}}] \oplus \R \oplus TS^1 \oplus T \tilde \eta^i_l).
\end{align*}
We give the above space the orientation of fiber product.
(The order of each direct sum $\bigoplus$ is the lexicographic order.)
This orientation and the complex orientation of the obstruction bundle
$E = E^1 \oplus E^2 \oplus \dots \oplus E^m$ define the orientation of the fiber product
as a pre-Kuranishi space.
(The orientation of the zero set of the perturbed multisection on $(V, E)$ is defined by
deleting the vector space $E$ from the tangent space of $V$.)
We note that the parity of the dimension of $\widetilde{\W}^1 \oplus \widetilde{\W}^2
\oplus \dots \oplus \widetilde{\W}^m$ is equal to the parity of $m$
since each $\widetilde{\W}^i$ has odd dimension.

We need to check that this orientation is compatible with the embeddings of Kuranishi
neighborhoods.
First we note that in (\ref{ori}), $\star$ is independent of the Kuranishi neighborhood.
Hence it is enough to compare the orientations of each $\W$ (or $\widetilde{\W}$).

For the convenience of the computation of orientations,
we may change the definition of the map $s^0 : \hat V \to \R^k \oplus
\bigoplus_{z_\beta^{++}} \R^2$ from (\ref{s^0}) to
\[
s^0(a, b, x) = (\sigma_i, p' \circ (\theta_{\sigma_i}^{-1} \times 1) \circ
\Phi_{a, b}(\xi_x)(Z_\beta^{++}(a)))
\]
where $\sigma_i = \sigma \circ \Phi_{a, b}(\xi_x)(\widetilde{R}_i(a))$,
and each $\theta_{\sigma_i} : \R \to \R$ is defined by
$\theta_{\sigma_i}(s) = s + \sigma_i$.
(Note that $V = \{s^0 = 0\}$ does not change.)
Let $C \subset T_{(0, b^0, 0)} \hat V$ be the tangent space of $\R$-translations.
Then the restriction of the differential of $s^0$ gives an isomorphism
$C \cong \R^k \oplus 0 \subset \R^k \oplus \bigoplus_{z_\beta^{++}} \R^2$, and
the differentials of $s$ or the evaluation maps at marked points or limit circles
vanish on $C$.
We note that under this identification $C \cong \R^k$,
the restriction of the differential of $(b_{\mu_1}, \dots, b_{\mu_{k-1}})$ to $C$ is
$(t_1, t_2, \dots, t_k) \mapsto (t_1 - t_2, t_2 - t_3, \dots, t_{k-1} - t_k) : C \cong \R^k \to
\R^{k-1}$.
It is clear that the definition of the orientation does not depend on the
choice of the family of sections $\hat R_j$ or the choice of
the coordinates of the neighborhoods of joint circles of $\Sigma_0$
used for the definition of the asymptotic parameters $b_\mu$.

We consider the situation discussed in Section \ref{embed},
and compare the orientation of $\W_1$ for $(V_1, E_1, s_1, \psi_1)$
and $\W_2$ for $(V_2, E_2, s_2, \psi_2)$ at $q_0$.

First we consider the case where $q_0 = p_1 = p_2$.
We may assume that we have used the same coordinates of the neighborhoods
of joint circles of $\Sigma_1 = \Sigma_2$ for the definition of the asymptotic
parameters $b^1_\mu$ and $b^2_\mu$.
The compatibility of the orientation in this case is essentially
because all the differences have the complex orientations.
To see this more precisely,
it is convenient to consider the following intermediate Kuranishi neighborhood.
Let $(\widetilde{P}_3 \to \widetilde{X}_3, Z)$ be the local universal family of
$(\Sigma_1, z_1 \cup z_2^+ \cup z_1^{++} \cup z_2^{++})$
and let $\lambda_3 : E_3^0 = E_2^0 \to C^\infty(\widetilde{P}_3 \times Y,
\Wedge^{0, 1} V^\ast \widetilde{P}_3 \otimes (\R \partial_\sigma \oplus TY))$
be the pull back of $\lambda_2$.
Let $\widetilde{R}^{3, 1}_j, \widetilde{R}^{3, 2}_j : \widetilde{X}_3 \to \widetilde{P}_3$
be the pull backs of $\widetilde{R}^1_j, \widetilde{R}^2_j$.
Then, using the parameter space $\mathring{X}_3 \subset \widetilde{X}_3 \times
\coprod_{\text{joint circles}} \R_\mu$ and $\lambda_3$,
we can construct a Kuranishi neighborhood of $p_1 = p_2$,
where in this case, for the definition of $V_3 = \{s^0_3 = 0\}$,
we use $s_3^0 : \hat V_3^0 \to \R^k \oplus \bigoplus_{z_{1, \beta^{++}}} \R^2
\oplus \bigoplus_{z_{2, \beta^{++}}} \R^2$
defined by
\begin{align*}
s_0^3(a, b, x) &= (\sigma^1_i, p' \circ (\theta_{\sigma^1_i}^{-1} \times 1) \circ
\Phi_{a, b}(\xi_x)(Z_{1, \beta}^{++}(a)),\\
&\quad\quad p' \circ (\theta_{\sigma^2_i}^{-1} \times 1) \circ
\Phi_{a, b}(\xi_x)(Z_\beta^{++}(a))),
\end{align*}
where $\sigma^l_i = \sigma \circ \Phi_{a, b}(\xi_x)(\widetilde{R}^{3, l}_i(a))$ for $l = 1, 2$.
It is clear that we can define the embeddings from $V_1$ and $V_2$ to $V_3$ as
in Section \ref{embed}.
Then $T_{(0, b^0, 0)} \hat V_3 = T_{(0, b^0, 0)} \hat V_1 \oplus F$,
where $F$ is the tangent of the parameters for the additional marked points
$(z_2^+ \setminus z_1^+) \cup z_2^{++}$.
It is clear that the projection of the restriction of $s_0^3 \times s^3$
gives an isomorphism of complex vector spaces
\[
F \tocong \bigoplus_{z_{2, \beta^{++}}} \R^2 \oplus
\bigoplus_{z^+ \in z_2^+ \setminus z_1^+} \R^2.
\]
(Recall that the complex orientations of $\R^2 \cong T \hat Y / TS'_2$ or
$\R^2 \cong T \hat Y / T(\R \times S_2)$
are defined by the almost complex structure of $\hat Y$.)
This implies that the pair of vector spaces $(TV_3, E_3)$ is isomorphic to
$(TV_1 \oplus F', E_1 \oplus F')$ for some complex vector space $F'$.
Similar condition is satisfied for the embedding $V_2 \inj V_3$.
Therefore, the embedding preserves the orientation.

Next we consider the case where $q_0 = p_2$ and all the data for
the construction of $(V_2, E_2, s_2, \psi_2)$ are the restriction of those for
$(V_1, E_1, s_1, \psi_1)$.
We assume that the $i$-th floor and the $(i+1)$-th floor of $p_1$
are glued into one floor in $p_2$ and the others are not.
We may assume that $(\hat R^2_1, \dots, \hat R^2_{k_2})
= (\hat R^1_1, \stackrel{i}{\check \dots}, \hat R^1_{k_1})$.
Let $C_l \subset T\hat V_l$ be the tangent space of $\R$-translations
for each $l =1, 2$.
Then it is easy to see that $C_1 \cong C_2 \oplus \R$ and the sign of
$b^1_{\mu_i}$ and $\sigma^1_i$
($= \sigma \circ \Phi^1_{a, b}(\xi_x)(\widetilde{R}^1_i(a))$ at $(a, b, x) \in \hat V_1$)
coincide on this $\R$.
Since the orientation of $\W_l$ are defined by subtracting
vector space $(\R^{k_l} \oplus \bigoplus_{z_l^{++}} \R^2)$ from
$T\mathring{X}_l \times \mathcal{C}_l$ and the subtractions are from the left,
this implies that the embedding preserves the orientation.

The general case is covered by the combination of the above two cases.

Furthermore, the orientation is independent of the choice of the lifts of $c^i_l$,
$\eta^i_l$ and $\epsilon^{i, j}_l$ under the natural isomorphism.
Hence we may denote the above fiber product pre-Kuranishi space by
$\overline{\M}^m_{((\hat \epsilon^{i, j}_l), (\hat c^i_l),
(x^i_l), (\hat \eta^i_l))}$.

The algebras of SFT are constructed by the virtual fundamental chains of
the zero-dimensional component of these fiber product pre-Kuranishi spaces, and
the algebraic properties of them are proved by the equation corresponding to
the boundary of the one-dimensional component of the fiber products.
First we study the boundary of $\overline{\M}_{((\hat c_l), (x_l), (\hat \eta_l))}$.
It consists of several parts, and some of them are due to the splitting in $\R$-direction,
and the others are due to the boundaries of the simplices $c_l$ and $\eta_l$.

We consider the former.
For each Kuranishi neighborhood $(V, E, s, \psi)$ of $\hat \M$,
each of these parts corresponds to the subspace
$\{\rho_\mu = 0; \text{ for some (and all) } \ab \mu \in M_i\}$ of $V$.
($M_i$ is the set of the indices of joint circles between the $i$-th floor and the $(i+1)$-th
floor.)
We note that the normal direction is $\kappa_i = - L_\mu \log \rho_\mu + b_\mu$
($\mu \in M_i$) and $\kappa_i^{-1} = 0$ defines the boundary.

As we have seen in Section \ref{construction of a multisection}
(related to the second compatibility condition of the multisection),
the curve corresponding to each zero of the multisection in this boundary is determined
by two curves and a family of diffeomorphisms between some of their limit circles.
Assume that a $+ \infty$-limit circle $S^1_{+\infty_{(1, 2), l}} \subset \Sigma_1$
and a $-\infty$-limit circle $S^1_{-\infty_{(2, 1), l}} \subset \Sigma_2$ are identified
by a diffeomorphism $\phi_{(1, 2), l} : S^1_{+\infty_{(1, 2), l}} \to S^1_{-\infty_{(2, 1), l}}$.
$\phi_{(1, 2), l}$ is determined by a pair of the coordinates $\phi_{+\infty_{(1, 2), l}} :
S^1 \to S^1_{+\infty_{(1, 2), l}}$ of $S^1_{+\infty_{(1, 2), l}}$ and
$\phi_{-\infty_{(2, 1), l}} : S^1 \to S^1_{-\infty_{(2, 1), l}}$ of $S^1_{-\infty_{(2, 1), l}}$
such that $\pi_Y \circ u_1 \circ \phi_{+\infty_{(1, 2), l}} = \pi_Y \circ u_2 \circ
\phi_{-\infty_{(2, 1), l}}$ in $P$.
(Namely, this pair corresponds to
$\phi_{(1, 2), l} = \phi_{-\infty_{(2, 1), l}} \circ \phi_{+\infty_{(1, 2), l}}^{-1}$.)
For any $g \in S^1 \subset \Aut S^1$, $(\phi_{+\infty_{(1, 2), l}} \circ g,
\phi_{-\infty_{(2, 1), l}} \circ g)$ and $(\phi_{+\infty_{(1, 2), l}}, \phi_{-\infty_{(2, 1), l}})$
correspond to the same diffeomorphism.

Assume that $\pi_Y \circ u_1|_{S^1_{\phi_{+\infty_{(1, 2), l}}}}$
($= \pi_Y \circ u_2|_{S^1_{\phi_{-\infty_{(2, 1), l}}}}$) is contained in $\Int \zeta$ for
some top-dimensional simplex $\zeta$ of $\overline{P}$.
(Note that this assumption is satisfied if we restrict to the case of the boundary of the
one-dimensional component of $\overline{\M}_{((\hat c_l), (x_l), (\hat \eta_l))}$.
This is due to the first condition of the multisection in Section
\ref{construction of a multisection}.)
Let $\tilde \zeta \subset P$ be a lift of $\zeta$.
Then we can choose a pair of the coordinates $(\phi_{+\infty_{(1, 2), l}},
\phi_{-\infty_{(2, 1), l}})$ such that $\pi_Y \circ u_1 \circ \phi_{+\infty_{(1, 2), l}} \in
\tilde \zeta$.
For each diffeomorphism $\phi_{(1, 2), l} : S^1_{+\infty_{(1, 2), l}} \to S^1_{-\infty_{(2, 1), l}}$,
the number of such representatives is $m_\zeta$, where $m_\zeta$ is the multiplicity of
the periodic orbits in $\Int \zeta$.
(The number of different diffeomorphisms $S^1_{+\infty_{(1, 2), l}} \to
S^1_{-\infty_{(2, 1), l}}$ is also $m_\zeta$.)

Define a chain $\widetilde{\Delta}_{\overline{P}}$ (not a cycle) in $P \times P$ by
\[
\widetilde{\Delta}_{\overline{P}} = \sum \frac{1}{m_\zeta}
\theta^{\lsuperscript{D}{t}}_{\tilde \zeta}
(\Delta_\ast \tilde \zeta) \theta^D_{\tilde \zeta},
\]
where the sum is taken over all top-dimensional simplices of $K$, including the simplices
contained in $\overline{P}^{\text{bad}}$.
As in the definition of $\Delta_\ast [\overline{P}]$,
$\theta^D_{\tilde \zeta}$ is an arbitrary fixed orientation of $p_2^\ast \S^D$ on
$\Int \Delta_\ast \tilde \zeta$, and $\theta^{\lsuperscript{D}{t}}_{\tilde \zeta}$ is
the orientation of $p_1^\ast \S^{\lsuperscript{D}{t}}$ defined by
$\theta^D_{\tilde \zeta}$ and $\theta^{\overline{P}}_{\tilde \zeta}$.
Then by the above argument, the part of the boundary of the zero-dimensional
component of $\overline{\M}_{((\breve c_l), (x_l), (\breve \eta_l))}$ corresponding to
the splitting in $\R$-direction is the zero-dimensional component of
\[
- \sum (-1)^\ast \overline{\M}^2_{
(e^{\widetilde{\Delta}_{\overline{P}}}, (\breve c^i_l), (x^i_l), (\breve \eta^i_l))}
\]
where the sum is taken over all decompositions
\[
\{\breve c_l\} = \{\breve c^1_l\} \sqcup \{\breve c^2_l\},
\quad \{x_l\} = \{x^1_l\} \sqcup \{x^2_l\},
\quad \{\breve \eta_l\} = \{\breve \eta^1_l\} \sqcup \{\breve \eta^2_l\}
\]
as sets, and
the order of each $(\breve c^i_l)_l$ is defined by the order of $(\breve c_l)_l$.
The orders of $(x^i_l)_l$ or $(\breve \eta^i_l)_l$ are similar.
$\ast$ is the weighted sign of the permutation
\[
\begin{pmatrix}
(\breve c^1_l)_l (\breve c^2_l)_l
& (x^1_l)_l (x^2_l)_l & (\breve \eta^1_l)_l (\breve \eta^2_l)_l\\
(\breve c_l)_l & (x_l)_l & (\breve \eta_l)_l
\end{pmatrix},
\]
where weighted sign is defined as follows.
The weighted sign of the transposition $(a, b)$ is defined by
$\deg a \cdot \deg b \in \Z / 2$, where
the degree is defined by $\deg c \theta^D_c = \deg (c \theta^D_c)^\ast
= \dim c + \dim [\mathring{D}^+_\gamma]$ ($\gamma \in |c|$) and
$\deg s = \codim_Y |s|$.
The weighted sign of a general permutation is defined by the product of the weighted sign
of the transpositions whose product coincides with the permutation.
$e^{\widetilde{\Delta}_{\overline{P}}} = 1 + \widetilde{\Delta}_{\overline{P}}
+ \frac{1}{2} (\widetilde{\Delta}_{\overline{P}}, \widetilde{\Delta}_{\overline{P}}) + \cdots$
is the exponential.

We claim that the virtual fundamental chain of the zero-dimensional component of the
above pre-Kuranishi space does not change if we replace
$\widetilde{\Delta}_{\overline{P}}$
with the sum taken over the top-dimensional simplices $\zeta$ of $K$ not contained
in $\overline{P}^{\text{bad}}$.
This is because if $\pi_Y \circ u_1 \circ \phi_{+\infty_{(1, 2), l}} \in \Int \tilde \zeta$
and $\zeta$ is contained in $\overline{P}^{\text{bad}}$, then the curve obtained by
the pair of coordinates $(\phi_{+\infty_{(1, 2), l}}, \phi_{-\infty_{(2, 1), l}} \circ
g_{1 / m_\zeta})$ instead of $(\phi_{+\infty_{(1, 2), l}}, \phi_{-\infty_{(2, 1), l}})$
($g_{1 / m_\zeta} \in \Aut S^1$ is the translation by $1 / m_\zeta$)
is also a zero of the multisection, but its orientation is opposite.

Similarly, in the parts of the boundary of $\overline{\M}_{((\hat c_l), (x_l), (\hat \eta_l))}$
due to the boundaries of the simplices of $c_l$ and $\eta_l$,
the parts of the boundaries of $c_l$ and $\eta_l$ contained in $\overline{P}^{\text{bad}}$
do not affect the virtual fundamental chain.

Therefore,
\begin{align}
0 &= [\partial \overline{\M}_{((\hat c_l), (x_l), (\hat \eta_l))}]^0 \notag\\
&= - [\overline{\M}_{\partial ((\hat c_l), (x_l), (\hat \eta_l))}]^0
- \sum (-1)^\ast \bigl[\overline{\M}^2_{
(e^{\Delta_\ast [\overline{P}]}, (\hat c^i_l), (x^i_l), (\hat \eta^i_l))}\bigr]^0,
\label{boundary of M}
\end{align}
where
$[\cdot]^0$ denotes the virtual fundamental chain of the zero-dimensional component,
and $\partial((\hat c_l), (x_l), (\hat \eta_l))$ is defined by
\begin{align*}
&\partial((\hat c_l), (x_l), (\hat \eta_l))\\
&= \sum_j (-1)^{\sum_{l < j} |\hat c_l|}
((\hat c_1, \dots, \partial \hat c_j, \dots, \hat c_{l_-}), (x_l), (\hat \eta_l))\\
&\quad + \sum_j (-1)^{\sum_l |\hat c_l| + \sum_l |x_l|^\bot + \sum_{l < j} |\hat \eta_l|}
((\hat c_l), (x_l), (\hat \eta_1, \dots, \partial \hat \eta_j, \dots, \hat \eta_{l_+})),
\end{align*}
where $|x|^\bot = \codim_Y x$.

Similarly, it is easy to see that for any $((\hat c_l), (x_l), (\hat \eta_l))$ and
$(\hat \epsilon^{i, j}_l)$,
\begin{align}
0 &= \sum_{\star_m} (-1)^\ast
\bigl[\partial (\overline{\M}^m_{((\hat \epsilon^{i, j}_l), (\hat c^i_l), (x^i_l), (\hat \eta^i_l))})
\bigr]^0
\notag\\
& = (-1)^m \sum_{\star_m} (-1)^\ast
\bigl[\overline{\M}^m_{\partial ((\hat \epsilon^{i, j}_l), (\hat c^i_l), (x^i_l), (\hat \eta^i_l))}
\bigr]^0 \notag\\
& \quad + \sum_{\substack{1 \leq i_0 \leq m \\ \star_{m+1}}} (-1)^{\ast + i_0}
\bigl[\overline{\M}^{m+1}_{((e^{\Delta_\ast [\overline{P}]})^{i_0, i_0+1} \cup
(\tau_{i_0}\hat \epsilon^{i, j}_l), (\hat c^i_l), (x^i_l), (\hat \eta^i_l))}\bigr]^0,
\label{boundary of MM}
\end{align}
where the sum of $\star_m$ is taken over all decompositions
\[
\{\hat c_l\} = \coprod_i \{\hat c^i_l\}, \quad \{x_l\} = \coprod_i \{x^i_l\},
\quad \{\hat \eta_l\} = \coprod_i \{\hat \eta^i_l\}
\]
as sets, and
$\ast$ is the weighted sign of the permutation
\[
\begin{pmatrix}
(\hat c^1_l)_l & \cdots & (\hat c^m_l)_l
& (x^1_l)_l & \cdots & (x^m_l)_l & (\hat \eta^1_l)_l & \cdots & (\hat \eta^m_l)_l\\
&(\hat c_l)_l&&&(x_l)_l&&&(\hat \eta_l)_l&
\end{pmatrix}.
\]
$\tau_{i_0}\hat \epsilon^{i, j}_l$ is defined by
\[\tau_{i_0} a^{i, j} = \begin{cases}
a^{i+1, j+1} & i_0 < i < j\\
a^{i_0, j+1} + a^{i_0 + 1, j+1} & i = i_0 < j\\
a^{i, j+1} & i < i_0 < j\\
a^{i, i_0} + a^{i, i_0 + 1} &i < j = i_0\\
a^{i, j} &i < j < i_0
\end{cases},
\]
where $a^{i, j}$ means the fiber product with $a$ at a $+\infty$-limit circle of
$i$-th holomorphic building and a $-\infty$-limit circle of $j$-th holomorphic building.
Unfortunately, equation (\ref{boundary of M}) is not the equation for
$\overline{\M}_{((\hat c_l), (x_l), (\hat \eta_l))}$'s in the Bott Morse case
since the second term cannot be written as a function of
$\overline{\M}_{((\hat c_l), (x_l), (\hat \eta_l))}$'s.
(The diagonal $\Delta_\ast [\overline{P}]$ cannot be written as a linear combination of
products of simplices in $K$.)
To obtain a meaningful equation, we add correction terms to
$\overline{\M}_{((\hat c_l), (x_l), (\hat \eta_l))}$ as follows.
The addition of these correction terms are equivalent to count the cascades
in \cite{Bo02}.

Let $(\hat c_l)$ be a family of chains in
$C_\ast(\overline{P}, \overline{P}^{\text{bad}}; \S^D \otimes \Q)$,
let $(x_l)$ be a family of simplices in $K^0$, and let $(\alpha_l)$ be a family of cochains in
$C^\ast(\overline{P}, \overline{P}^{\text{bad}}; \S^D \otimes \Q)$ with compact supports.
Then for such a family $((\hat c_l), (x_l), (\alpha_l))$, we define a pre-Kuranishi space
$\overline{\M}((\hat c_l), (x_l), (\alpha_l))$
(or a linear combination of pre-Kuranishi spaces) by
\begin{align*}
\overline{\M}((\hat c_l), (x_l), (\alpha_l))
&= \overline{\M}_{((\hat c_l), (x_l), ([\overline{P}] \cap \alpha_l))}\\
&\quad + \sum_{m = 2}^\infty \sum_{\star_m} (-1)^\ast
\overline{\M}^m_{(F_m,
(\hat c^i_l), (x^i_l), ([\overline{P}] \cap \alpha^i_l))},
\end{align*}
where the sum of $\star_m$ is taken over all decompositions
\[
\{\hat c_l\} = \coprod_i \{\hat c^i_l\}, \quad \{x_l\} = \coprod_i \{x^i_l\},
\quad \{\alpha_l\} = \coprod_i \{\alpha^i_l\}
\]
as sets, and
the order of each $(\hat c^i_l)_l$ is defined by the order of $(\hat c_l)_l$.
The orders of $(x^i_l)_l$ or $(\alpha^i_l)_l$ are similar.
$\ast$ is the weighted sign of the permutation
\[
\begin{pmatrix}
(\hat c^1_l)_l & \cdots & (\hat c^m_l)_l
& (x^1_l)_l & \cdots & (x^m_l)_l & (\alpha^1_l)_l & \cdots & (\alpha^m_l)_l\\
&(\hat c_l)_l&&&(x_l)_l&&&(\alpha_l)_l&
\end{pmatrix}.
\]
$(F_m)_{m \geq 2}$ is an appropriate family of linear combinations of
\[
((\rho_\ast [\overline{P}])^{i, j}, \dots, (\rho_\ast [\overline{P}])^{i, j},
\epsilon_{\overline{P}}^{i, j}, \dots, \epsilon_{\overline{P}}^{i, j},
(\Delta_\ast [\overline{P}])^{i, j}, \dots, (\Delta_\ast [\overline{P}])^{i, j})_{1 \leq i < j \leq m}
\]
defined in the next section.
The first term
$\overline{\M}_{(\hat c_l), (x_l), ([\overline{P}] \cap \alpha_l)}$
is the main term, and
the second is for the correction of the difference between
$[\overline{P}]$ and $\rho_\ast[\overline{P}]$.

Note that if $(\Sigma_i, z_i, u_i, \phi_i)_{1 \leq i \leq m}$ is in the zero set of the perturbed
multisection of the zero-dimensional component of
$\overline{\M}^m_{(F_m, (\hat c^i_l), (x^i_l), ([\overline{P}] \cap \alpha^i_l))}$,
then each $\Sigma_i$ is connected.
This is because the multisection of $\overline{\M}^m
_{(F_m, (\hat c^i_l), (x^i_l), ([\overline{P}] \cap \alpha^i_l))}$
is the pull back of that of
$(\widehat{\M}^\#, \mathring{K}^2)^{\diamond m}
_{(\bar F_m, (\bar c^i), (\bar x^i), ([\overline{P}] \cap \bar \alpha^i))}$,
and its dimension is $< 0$ if some $\Sigma_i$ is disconnected.
In particular, the genus of each $\Sigma_i$ is $\geq 0$.
Since the total number of $(\rho_\ast [\overline{P}])^{i, j}$,
$\epsilon_{\overline{P}}^{i, j}$ and $(\Delta_\ast [\overline{P}])^{i, j}$ ($1 \leq i < j \leq m$)
contained in each term of $F_m$ is $\geq m - 1$
(in fact, the number of $\epsilon_{\overline{P}}$ is $m - 1$),
the genera of the sequences of curves $(\Sigma_i, z_i, u_i, \phi_i)_{1 \leq i \leq m}$
corresponding to the zeros of the multisection of
the zero-dimensional component of $\overline{\M}^m
_{(F_m, (\hat c^i_l), (x^i_l), ([\overline{P}] \cap \alpha^i_l))}$ are $\geq 0$.

\begin{rem}\label{connected or not}
We do not know whether or not we can choose $(F_m)_{m \geq 2}$ so that
all sequences of holomorphic buildings in the zero-dimensional part of
the fiber products $\overline{\M}^m
_{(F_m, (\hat c^i_l), (x^i_l), ([\overline{P}] \cap \alpha^i_l))}$ are connected
(in the appropriate sense).
However, for the construction of the algebras, it is enough to show that their genera are
$\geq 0$.
\end{rem}

In the next section, we prove that if we choose an appropriate family $(F_m)_{m \geq 2}$,
then the following equation holds true.
\begin{align}
0 &= [\partial \overline{\M}((\hat c_l), (x_l), (\alpha_l))]^0 \notag\\
&=
- [\overline{\M}\bigl(\partial((\hat c_l), (x_l), (\alpha_l))\bigr)]^0 \notag\\
&\quad + \sum_{\blacklozenge} (-1)^{\ast}
\frac{1}{k !} [\overline{\M}((\hat c^1_l), (x^1_l), (\alpha^1_l) \cup (\hat d_1^\ast,
\hat d_2^\ast, \dots, \hat d_k^\ast))]^0\notag\\
&\hspace{70pt}
\times
[\overline{\M}((\hat d_k, \hat d_{k-1}, \dots, \hat d_1) \cup (\hat c^2_l), (x^2_l),
(\alpha^2_l))]^0
\label{boundary formula}
\end{align}
$\partial((\hat c_l), (x_l), (\alpha_l))$ is defined by
\begin{align*}
&\partial((\hat c_l), (x_l), (\alpha_l))\\
&= \sum_j (-1)^{\sum_{l < j} |\hat c_l|}
((\hat c_1, \dots, \partial \hat c_j, \dots, \hat c_{l_-}), (x_l), (\alpha_l))\\
&\quad + \sum_j (-1)^{\sum_l |\hat c_l| + \sum_l |x_l|^\bot + \sum_{l < j} |\alpha_l|}
((\hat c_l), (x_l), (\alpha_1, \dots, \partial \alpha_j, \dots, \alpha_{l_+})),
\end{align*}
where $\partial \alpha$ is defined by
$\partial \alpha = (-1)^{|\alpha|} \delta \alpha = (-1)^{|\alpha|} \alpha \circ \partial$.
The sum $\blacklozenge$ of the last term is taken over
all decompositions
\[
\{\hat c_l\} = \{\hat c^1_l\} \sqcup \{\hat c^2_l\}, \quad \{x_l\} = \{x^1_l\} \sqcup \{x^2_l\},
\quad \{\alpha_l\} = \{\alpha^1_l\} \sqcup \{\alpha^2_l\}
\]
as sets, $k \geq 0$, and all sequences of simplices $d_l$ of $K$ not contained in
$\overline{P}^{\text{bad}}$. (We fix $\theta^D_{d}$ for each simplex $d$ and define $\hat d
= d \theta^D_{d}$.)
The sign $\ast$ of the last term is the weighted sign of the permutation
\[
\begin{pmatrix}
(\hat c^1_l)_l \  (x^1_l)_l \  (\alpha^1_l)_l \  (\hat c^2_l)_l \  (x^2_l)_l \  (\alpha^2_l)_l\\
(\hat c_l)_l  \quad (x_l)_l \quad (\alpha_l)_l
\end{pmatrix}.
\]

For the proof of equation (\ref{boundary formula}), we use the fact
\[
\sum_d ([\overline{P}] \cap (\hat d)^\ast) \otimes \hat d = \rho_\ast [\overline{P}].
\]

\subsection{Construction of the correction terms}
\label{algebra for correction}
In this section, we construct $(F_m)_{\geq 2}$ used for the definition of the correction
terms in $\overline{\M}((\hat c_l), (x_l), (\alpha_l))$, and prove the equation
(\ref{boundary formula}).
For that sake, we consider an algebra modeled on the splitting of holomorphic buildings.

For $m \geq 2$, let $A_m = \bigoplus_{n=0}^{\frac{m(m-1)}{2}} A_m^n$ be the $\Z$-graded
super-commutative algebra with coefficient $\Q$ generated by the variables
$\rho_{(e_i, e_j)}$, $\Delta_{(e_i, e_j)}$ and $\epsilon_{(e_i, e_j)}$ ($1 \leq i < j \leq m$),
where the $\Z$-grading is defined by
$\dim \rho_{(e_i, e_j)} = \dim \Delta_{(e_i, e_j)} = 0$ and $\dim \epsilon_{(e_i, e_j)} = 1$.
$\rho$, $\Delta$ and $\epsilon$ are variables corresponding to
$\rho_\ast [\overline{P}]$, $\Delta_\ast [\overline{P}]$, and
$\epsilon_{\overline{P}}$ respectively.
In particular, the parity of the dimension of a monomial in $A_m$ coincides with
that of the corresponding product of simplices.
(We call $n$ dimension in order to distinguish it from the degree $m$.)
We sometimes use the following notation:
$x_{(\sum_i a_i e_i, \sum_i b_i e_i)} = \sum_{i, j} a_i b_j x_{(e_i, e_j)}$, where $x$ is $\rho$,
$\Delta$ or $\epsilon$.
For $m = 1$, we define $A_1 = \Q$.

For each $m \geq 2$, the differential $\partial' : A_m \to A_m$ is defined by
$\partial' \epsilon_{(a, b)} = (-1)^m (\rho_{(a, b)} - \Delta_{(a, b)})$ and
$\partial' \rho_{(a, b)} = \partial' \Delta_{(a, b)} = 0$.
For $m = 1$, we define $\partial' = 0 : \A_1 \to \A_1$.

We define homomorphisms $\tau_i : A_m \to A_{m+1}$ ($1 \leq i \leq m$, $m \geq 2$) by
$\tau_i(x_{(a, b)}) = x_{(\hat \tau_i(a), \hat \tau_i(b))}$,
where each $\hat \tau_i$ is defined by
\[
\hat \tau_i (e_j) = \begin{cases}
e_j &j < i\\
e_i + e_{i + 1} &j = i\\
e_{j + 1} & j > i
\end{cases}.
\]
For example,
\[
\tau_2(\Delta_{(e_1, e_2)} \epsilon_{(e_2, e_3)})
= (\Delta_{(e_1, e_2)} + \Delta_{(e_1, e_3)}) (\epsilon_{(e_2, e_4)} + \epsilon_{(e_3, e_4)}).
\]
For $m = 1$, we define $\tau_1 = \id_{\Q}$.
For $i > m$, we define $\tau_i = 0 : A_m \to A_{m+1}$.

We also define homomorphisms $\Box : A_m \otimes A_{m'} \to A_{m + m'}$
($m, m' \geq 1$) by
\[
\Box(f \otimes g) = (-1)^{(m - 1)m'} f \cdot
\exp (\rho_{(\sum_{1\leq i \leq m} e_i, \sum_{m + 1 \leq j \leq m + m'} e_j)})
\cdot g^{+ m}
\]
where $g^{+ m}$ is the image of $g$ by the homomorphism $A_{m'} \to A_{m + m'}$
defined by $x_{(e_i, e_j)} \mapsto x_{(e_{i + m}, e_{j + m})}$.
For example, if $m = 2$ and $m' = 2$, then
\begin{align*}
&\Box(\Delta_{(e_1, e_2)} \otimes \rho_{(e_1, e_2)} \epsilon_{(e_1, e_2)})\\
&= \Delta_{(e_1, e_2)} \rho_{(e_3, e_4)} \epsilon_{(e_3, e_4)}
\exp (\rho_{(e_1, e_3)} + \rho_{(e_1, e_4)} + \rho_{(e_2, e_3)} + \rho_{(e_2, e_4)}).
\end{align*}

Define a linear subspace $\Ddot A_m \subset A_m$ as follows.
(It is not an ideal.)
For each $1 \leq i \leq m-1$ and each monomial
\[
f = x^{(1)}_{(a_1, b_1)} x^{(2)}_{(a_2, b_2)} \dots x^{(n)}_{(a_n, b_n)},
\]
(each $x^{(j)}$ is $\rho$, $\Delta$ or $\epsilon$)
such that $(a_j, b_j) \neq (e_i, e_{i + 1})$ for all $1 \leq j \leq n$,
we define a monomial
\[
f^{(e_i, e_{i + 1})} = x^{(1)}_{(a'_1, b'_1)} x^{(2)}_{(a'_2, b'_2)} \dots x^{(n)}_{(a'_n, b'_n)}
\]
by permuting $e_i$ and $e_{i + 1}$ appearing in $\{a_j, b_j\}$.
$\Ddot A_m \subset A_m$ is the subspace spanned by $f + f^{(e_i, e_{i + 1})}$
for all such pairs of $i$ and $f$.

Define $\A_m = A_m / \Ddot A_m$.
It is not an algebra, but the following maps are well defined.
(Namely, the corresponding maps on $A_m$ or $A_m \otimes A_{m'}$
induce the following maps.)
\begin{align*}
\partial' &: \A_m \to \A_m & (m \geq 1)\\
\sum_{i = 1}^m (-1)^i e^{\Delta_{(e_i, e_{i + 1})}} \tau_i &: \A_m \to \A_{m + 1}
& (m \geq 1)\\
\Box  &: \A_m \otimes \A_{m'} \to \A_{m + m'} & (m, m' \geq 1)
\end{align*}
The well-definedness of the first and the third maps are easy to see.
The well-definedness of the second is proved as follows.
If $f \in A_m$ does not contain any $x_{e_{i_0}, e_{i_0 + 1}}$ ($x = \rho, \Delta, \epsilon$),
then
\begin{align*}
&\sum_{i = 1}^m (-1)^i e^{\Delta_{(e_i, e_{i + 1})}} \tau_i (f + f^{(e_{i_0}, e_{i_0 + 1})})\\
& = \sum_{i \neq i_0, i_0 + 1} (-1)^i
\bigl((e^{\Delta_{(e_i, e_{i + 1})}} \tau_i f)
+ (e^{\Delta_{(e_i, e_{i + 1})}} \tau_i f)^{(\hat \tau_i(e_{i_0}), \hat \tau_i(e_{i_0 + 1}))}\bigr)\\
& \quad + (-1)^{i_0} e^{\Delta_{(e_{i_0}, e_{i_0 + 1})}} \tau_{i_0} (f + f^{(e_{i_0}, e_{i_0 + 1})})\\
& \quad + (-1)^{i_0 + 1} e^{\Delta_{(e_{i_0 + 1}, e_{i_0 + 2})}} \tau_{i_0 + 1}
(f + f^{(e_{i_0}, e_{i_0 + 1})}).
\end{align*}
The sum of the last two terms of the right hand side is an element of $\ddot A_{m+1}$
since
\[
\bigl(e^{\Delta_{(e_{i_0}, e_{i_0 + 1})}} \tau_{i_0} f \bigr)^{(e_{i_0 + 1}, e_{i_0 + 2})}
= \bigl(e^{\Delta_{(e_{i_0 + 1}, e_{i_0 + 2})}} \tau_{i_0 + 1}(f^{(e_{i_0}, e_{i_0 + 1})})
\bigr)^{(e_{i_0}, e_{i_0 + 1})}
\]
and
\[
\bigl(e^{\Delta_{(e_{i_0}, e_{i_0 + 1})}} \tau_{i_0}(f^{(e_{i_0}, e_{i_0 + 1})})
\bigr)^{(e_{i_0 + 1}, e_{i_0 + 2})}
= \bigl(e^{\Delta_{(e_{i_0 + 1}, e_{i_0 + 2})}} \tau_{i_0 + 1} f \bigr)^{(e_{i_0}, e_{i_0 + 1})}.
\]
Hence $\sum_{i = 1}^m (-1)^i e^{\Delta_{(e_i, e_{i + 1})}} \tau_i : \A_m \to \A_{m+1}$
is well defined.

Let $\A = (\bigoplus_{m = 1}^\infty \A_m^{m - 1})^\wedge$ be the completion
with respect to the degree $m$.
We also define $\A' = (\bigoplus_{m = 2}^\infty \A_m^{m - 2})^\wedge$.
In this section, we prove that the map $\A \to \A'$ defined by
\[
F \mapsto \partial' F + \sum_i (-1)^i e^{\Delta{(e_i, e_{i + 1})}} \tau_i F + \Box (F \otimes F)
\]
has a zero $F = F_1 + F_2 + \dots \in \A$ such that
$F \equiv 1 \in \A / (\bigoplus_{m = 2}^\infty \A_m^{m - 1})^\wedge \cong \A_1$.

Equation (\ref{boundary formula}) holds for such a zero $F$
if we replace the variables $\rho_{(e_i, e_j)}$, $\epsilon_{(e_i, e_j)}$ and $\Delta_{(e_i, e_j)}$
with $(\rho_\ast [\overline{P}])^{i, j}$, $\epsilon_{\overline{P}}^{i, j}$
and $(\Delta_\ast [\overline{P}])^{i, j}$ respectively.
This can be seen as follows.

Equation (\ref{boundary of MM}) implies that
for any $F_m \in A_m^{m-1}$ and $((\hat c_l), (x_l), (\alpha_l))$,
\begin{align*}
&\sum_{\star_m} (-1)^\ast \partial' \bigl(
\overline{\M}^m_{(F_m,
(\hat c^i_l), (x^i_l), ([\overline{P}] \cap \alpha^i_l))}\bigr)\\
&= - \sum_{\star_m} (-1)^\ast
\overline{\M}^m_{(F_m,
\partial((\hat c^i_l), (x^i_l), ([\overline{P}] \cap \alpha^i_l)))}\\
&\quad + \sum_{\star_m} (-1)^\ast
\overline{\M}^m_{(\partial' F_m,
(\hat c^i_l), (x^i_l), ([\overline{P}] \cap \alpha^i_l))}\\
&\quad + \sum_{\star_{m + 1}} (-1)^\ast
\overline{\M}^{m+1}_{(
\sum_i (-1)^i e^{\Delta{(e_i, e_{i + 1})}} \tau_i F_m,
(\hat c^i_l), (x^i_l), ([\overline{P}] \cap \alpha^i_l))},
\end{align*}
On the other hand,
for any $((\hat c^i_l)_{1 \leq i \leq m + m'}, (x^i_l)_{1 \leq i \leq m + m'},
(\alpha^i_l)_{1 \leq i \leq m + m'})$,
$F_m \in A_m^{m-1}$ and $F_{m'} \in A_{m'}^{m'-1}$,
\begin{align*}
&\sum (-1)^\ast \frac{1}{k !} [
\overline{\M}^m_{(F_m,
(\hat c^i_l)_{i = 1}^m, (x^i_l)_{i = 1}^m,
([\overline{P}] \cap \alpha^i_l)_{i = 1}^m \cup ((\hat d_1^\ast)^m,
\dots, (\hat d_k^\ast)^m))}]^0\\
&\hph{\sum (-1)^\ast \frac{1}{k !}}
\cdot [
\overline{\M}^{m'}_{(F_{m'},
((\hat d_k)^1, \dots, (\hat d_1)^1) \cup
(\hat c^{i+m}_l)_{i = 1}^{m'},
(x^{i+m}_l)_{i = 1}^{m'},
([\overline{P}] \cap \alpha^{i+m}_l)_{i = 1}^{m'})}]^0\\
&= [\overline{\M}^{m+ m'}_{(- \Box
(F_m \otimes F_{m'}),
(\hat c^i_l)_{i = 1}^{m + m'}, (x^i_l)_{i = 1}^{m + m'},
([\overline{P}] \cap \alpha^i_l)_{i = 1}^{m + m'})}]^0,
\end{align*}
where the sum is taken over all $k \geq 0$ and all sequences of simplices
$d_l \in K$ such that $d_l \not \subset \overline{P}^{\text{bad}}$, and
$\ast$ is the weighted sign of the permutation
\[
\begin{pmatrix}
(\hat c^i_l)_{1 \leq i \leq m} \  (x^i_l)_{1 \leq i \leq m} \  (\alpha^i_l)_{1 \leq i \leq m}
\ (\hat c^{i+m}_l)_{1 \leq i \leq m'} \  (x^{i+m}_l)_{1 \leq i \leq m'}
\  (\alpha^{i+m}_l)_{1 \leq i \leq m'}\\
(\hat c^i_l)_{1 \leq i \leq m + m'}  \quad (x^i_l)_{1 \leq i \leq m + m'} \quad
(\alpha^i_l)_{1 \leq i \leq m + m'}
\end{pmatrix}.
\]
These equations imply that equation (\ref{boundary formula}) holds for a zero $F$.
The quotient space $\A_m = A_m / \Ddot A_m$ corresponds to the fact that
we can permute the $i$-th holomorphic building and $(i + 1)$-th holomorphic building
in $(\overline{\M} \times \dots \times \overline{\M})_{((l_{i, j}), (l_{i, \pm}), (\mu_i))}$
if $l_{i, i+1} = 0$.

Note that the homology of $\partial' : \A_m^\ast \to \A_m^\ast$ is zero
at $\ast \neq 0$. This is because that K\"unneth formula implies that
the homology of $A_m^\ast \cong (A_2^\ast)^{\otimes \frac{m(m-1)}{2}}$ is zero
at $\ast \neq 0$, and there exists a splitting $T : \A_m^\ast \to A_m^\ast$.
The splitting $T$ is defined as follows.
For a monomial
\[
f = x^{(1)}_{(a_1, b_1)} x^{(2)}_{(a_2, b_2)} \dots x^{(n)}_{(a_n, b_n)},
\]
we define a subgroup $\mathfrak{S}_f \subset \mathfrak{S}_m$ by
\[
\mathfrak{S}_f = \{\sigma \in \mathfrak{S}_m; \sigma(a_j) < \sigma(b_j) \text{ for all } j\}.
\]
Then $Tf$ is defined by
\[
Tf = \frac{1}{\# \mathfrak{S}_f} \sum_{\sigma \in \mathfrak{S}_f} \sign \sigma \cdot
x^{(1)}_{(\sigma(a_1), \sigma(b_1))} x^{(2)}_{(\sigma(a_2), \sigma(b_2))} \dots
x^{(n)}_{(\sigma(a_n), \sigma(b_n))}.
\]

Starting with $F_1 = 1 \in \A_1$,
we inductively construct $F_{\leq m} = F_1 + \dots + F_m \in \bigoplus_{l = 1}^m
\A_l^{l - 1}$ such that
\begin{equation}
\partial' F_{\leq m} + \sum_i (-1)^i e^{\Delta{(e_i, e_{i + 1})}} \tau_i F_{\leq m-1}
+ \Box (F_{\leq m-1} \otimes F_{\leq m-1}) \equiv 0 \label{F_m eq}
\end{equation}
in $\A' / (\bigoplus_{l = m + 1}^\infty \A_l^{l - 2})^\wedge$.
First we define $F_{\leq 2} = F_1 + F_2 \in \A_1 \oplus \A_2^1$ by
\begin{align*}
F_{\leq 2} = 1 - \frac{1}{k !} \sum_{k = 1}^\infty
(&\underbrace{\epsilon_{(e_1, e_2)} \Delta_{(e_1, e_2)} \cdots \Delta_{(e_1, e_2)}}_k\\
&+
\underbrace{\rho_{(e_1, e_2)} \epsilon_{(e_1, e_2)} \Delta_{(e_1, e_2)} \cdots
\Delta_{(e_1, e_2)}}_k\\
&+ \dots +
\underbrace{\rho_{(e_1, e_2)} \cdots \rho_{(e_1, e_2)} \epsilon_{(e_1, e_2)}}_k)
\end{align*}
It is easy to check that this satisfies equation (\ref{F_m eq}) for $m = 2$.

Next assuming that we have already constructed
$F_{\leq m-1} \in \bigoplus_{l = 1}^m \A_l^{l-1}
$, we need to prove that
there exists a required $F_{\leq m}$ ($m \geq 3$).
Since $\partial'$ is exact at $n \geq 1$,
it is enough to show that
\begin{equation}
\partial' \Bigl(
\sum_i (-1)^i e^{\Delta{(e_i, e_{i + 1})}} \tau_i F_{\leq m-1} + \Box (F_{\leq m-1} \otimes F_{\leq m-1})
\Bigr)
\equiv 0 \label{A closed}
\end{equation}
in $(\bigoplus_{l = 3}^\infty \A_l^{l-3})^\wedge / (\bigoplus_{l = m}^\infty \A_l^{l-3})^\wedge$.

Since $F_{\leq m-1} = F_1 + \dots + F_{m-1}$ satisfies
\[
\partial' F_{\leq m-1} + \sum_i (-1)^i e^{\Delta{(e_i, e_{i + 1})}} \tau_i F_{\leq m-1} + \Box (F_{\leq m-1} \otimes F_{\leq m-1}) \equiv 0
\]
in $\A' / (\bigoplus_{l = m}^\infty \A_l^{l - 2})^\wedge$,
we see that
\begin{align*}
&\partial' \Bigl(
\sum_i (-1)^i e^{\Delta{(e_i, e_{i + 1})}} \tau_i F_{\leq m-1} + \Box (F_{\leq m-1} \otimes F_{\leq m-1})
\Bigr)\\
&=
\sum_i (-1)^{i + 1} e^{\Delta{(e_i, e_{i + 1})}} \tau_i \partial' F_{\leq m-1}
+ \Box \Bigl(\partial' F_{\leq m-1} \otimes \sum_{1 \leq j \leq m - 1} (-1)^j F_j\Bigr)\\
&\quad - \Box (F_{\leq m-1} \otimes \partial' F_{\leq m-1})\\
&=
\sum_i (-1)^i e^{\Delta{(e_i, e_{i + 1})}} \tau_i \Bigl(
\sum_j (-1)^j e^{\Delta{(e_j, e_{j + 1})}} \tau_j F_{\leq m-1} + \Box (F_{\leq m-1} \otimes F_{\leq m-1})
\Bigr)\\
&\quad - \Box \Bigl(\Bigl(\sum_i (-1)^i e^{\Delta{(e_i, e_{i + 1})}} \tau_i F_{\leq m-1}
+ \Box (F_{\leq m-1} \otimes F_{\leq m-1}) \Bigr)\\
&\hph{\quad - \Box \Bigl(\Bigl(}
\otimes \sum_{1 \leq j \leq m - 1} (-1)^j F_j\Bigr)\\
&\quad +\Box \Bigl(F_{\leq m-1} \otimes \Bigl(\sum_i (-1)^i e^{\Delta{(e_i, e_{i + 1})}} \tau_i F_{\leq m-1}
+ \Box (F_{\leq m-1} \otimes F_{\leq m-1})\Bigr)\Bigr).
\end{align*}
By direct calculation, it is easy to see that the following equations hold true.
{\belowdisplayskip= 0pt
\[
\Bigl(\sum_i (-1)^i e^{\Delta{(e_i, e_{i + 1})}} \tau_i\Bigr) \circ
\Bigl(\sum_j (-1)^j e^{\Delta{(e_j, e_{j + 1})}} \tau_j\Bigr) = 0,
\]
}
\begin{multline*}
\sum_i (-1)^i e^{\Delta{(e_i, e_{i + 1})}} \tau_i \Box (f \otimes g)
- \Box \Bigl(\Bigl(\sum_i (-1)^i e^{\Delta{(e_i, e_{i + 1})}} \tau_i f\Bigr) \otimes (-1)^{\deg g} g\Bigr)\\
+ \Box \Bigl(f \otimes \Bigl(\sum_i (-1)^i e^{\Delta{(e_i, e_{i + 1})}} \tau_i g\Bigr)\Bigr) = 0,
\end{multline*}
\[
\Box (f \otimes \Box (g \otimes h)) - \Box ( \Box (f \otimes g) \otimes (-1)^{\deg h} h) = 0.
\]
Therefore
\[
\partial' \Bigl(
\sum_i (-1)^i e^{\Delta{(e_i, e_{i + 1})}} \tau_i F_{\leq m-1} + \Box (F_{\leq m-1} \otimes F_{\leq m-1})
\Bigr)
\equiv 0
\]
in $\bigoplus_{l = 3} \A_l^{l-3} / \bigoplus_{l = m} \A_l^{l-3}$,
and we can construct a required
$F_{\leq m} = F_1 + \dots + F_m \in \bigoplus_{l = 1}^m \A_l^{l-1}$.

\begin{rem}
In fact, we do not need to use $\A_m$, and we can replace $\A_m$ with $A_m$.
However, for the construction of the correction terms for $X$
in Section \ref{correction terms for X}, we need to use a counterpart of $\A_m$.
\end{rem}


%% file: SFT-07_Algebras_for_Y.tex
%
%

\subsection{Construction of the algebras}\label{construction of algebra}
Using the virtual fundamental chains of the $0$-dimensional components of
the pre-Kuranishi spaces in the previous section,
we construct the algebra of symplectic field theory.
We mainly follow the construction explained in \cite{EGH00}.
First we consider general SFT.
We do not consider the $H_2(Y; \Z)$-grading or the $H_1(Y; \Z)$-grading for simplicity.
(See the above paper for these gradings.)

For each simplex $c$ of $K$ not contained in $\overline{P}^{\text{bad}}$,
we fix an orientation $\theta^D_c$ and define $\hat c = c \theta^D_c$.
We use the following variables:
$q_{\hat c^\ast}$ and $p_{\hat c}$ for each simplex $c$ of $K$ not contained in $\overline{P}^{\text{bad}}$,
$t_x$ for each cycle $x$ of $K^0$, and $\hbar$.
The $\Z/2$-degrees of these variables are defined by
$|q_{\hat c^\ast}| = |p_{\hat c}| =  \dim c + \ind \mathring{D}_\gamma^+$
($\gamma \in |c|$), $|t_x| = \codim_Y x$ and $|\hbar| = 0$.
We define the energies of these variables by
$e(q_{\hat c^\ast}) = L_\gamma$ and $e(p_{\hat c}) = - L_\gamma$ for each $c$,
where $\gamma \in |c|$ is an arbitrary periodic orbit and $L_\gamma$ is its period, and
$e(t_x) = e(\hbar) = 0$.

The algebra $\W_Y = \W_{(Y, \lambda, K_Y, \overline{K}_Y^0)}$ is defined as follows.
Its elements are formal series
\[
\sum_{(\hat c_i^\ast), (\hat c'_i)}
f_{(\hat c_i^\ast), (\hat c'_i)}(t, \hbar) q_{\hat c_1^\ast} q_{\hat c_2^\ast} \dots
q_{\hat c_{k_q}^\ast} p_{\hat c'_1} p_{\hat c'_2} \dots p_{\hat c'_{k_p}},
\]
where $f_{(\hat c_i^\ast), (\hat c'_i)}(t, \hbar) \in \R[[t, \hbar]]$ are formal series
of the variables $t_x$ and $\hbar$, and the infinite sum is taken over
all pairs of sequences $(\hat c_i)$ and $(\hat c'_i)$ with the following Novikov condition:
for any $C \geq 0$, the number of the terms with $\sum_i e(p_{\hat c'_i}) \geq -C$
is finite.
(This is equivalent to the condition that for each sequence $(\hat c'_i)$,
all but finite sequences $(\hat c_i^\ast)$ satisfy $f_{(\hat c_i^\ast), (\hat c'_i)} = 0$.)
We sometimes use the following notation:
for a linear combination $\sum_i r_i \hat c_i$, we define
$p_{\sum_i r_i \hat c_i} = \sum_i r_i p_{\hat c_i}$.
We use the similar notation for variables $q$ and $t$.
The associative product $\circ$ of $\W_Y$ is defined by the following commutative
relations: all variables are super-commutative except
\[
[p_{\hat c}, q_\alpha]
= p_{\hat c} \circ q_\alpha - (-1)^{|p_{\hat c}| \cdot |q_\alpha|} q_\alpha \circ p_{\hat c}
= \langle \hat c, \alpha \rangle \hbar.
\]
We often omit the symbol $\circ$ and denote the product $f \circ g$ by $fg$.

For each $\kappa \geq 0$, we define a submodule $\W_Y^{\leq \kappa} \subset \W_Y$
by imposing the condition $\sum_i e(q_{\hat c_i^\ast}) + \sum_i e(p_{\hat c'_i})
\leq \kappa$.
(This condition is stronger than the Novikov condition.)
For each triple $(C_0, C_1, C_2)$,
we define a submodule $I^{\leq \kappa}_{C_0, C_1, C_2} \subset \W_Y^{\leq \kappa}$ by
\begin{align*}
I^{\leq \kappa}_{C_0, C_1, C_2} &= \Bigl\{\sum
a_{(x_i), (\hat c_i^\ast), (\hat c'_i), g}\ t_{x_1} \dots t_{x_{k_t}}
q_{\hat c_1^\ast} \dots q_{\hat c_{k_q}^\ast}
p_{\hat c'_1} \dots p_{\hat c'_{k_p}} \hbar^g \in \W_Y^{\leq \kappa};\\
&\quad \quad a_{(x_i), (\hat c_i^\ast), (\hat c'_i), g} = 0
\text{ for all } ((x_i)_{i = 1}^{k_t}, (\hat c_i^\ast),_{i = 1}^{k_q} (\hat c'_i)_{i = 1}^{k_p}, g)
\text{ such that}\\
&\quad \quad k_t \leq C_0,\, \widetilde{g} \leq C_1 \text{ and }
\sum e(p_{\hat c'_i}) \geq - C_2 \Bigr\},
\end{align*}
where
\[
\widetilde{g} = g + \frac{1}{2}(k_t + k_q + k_p)
- \frac{\sum_i e(q_{\hat c_i^\ast}) + \sum_j e(p_{\hat c'_j})}
{L_{\min}}.
\]
($L_{\min}$ is the minimal period of the periodic orbits of $R_\lambda$.)

We note that
\[
\W_Y \cong \varprojlim_{C_2} \varinjlim_{\kappa} \varprojlim_{C_0, C_1}
\W_Y^{\leq \kappa} / I^{\leq \kappa}_{C_0, C_1, C_2}.
\]
The multiplication of $\W_Y$ defines the maps
\[
\W_Y^{\leq \kappa_1} / I^{\leq \kappa_1}_{C_0, C_1 + \kappa_2 L_{\min}^{-1},
C_2 + \kappa_2}
\times \W_Y^{\leq \kappa_2} / I^{\leq \kappa_2}_{C_0, C_1 + \kappa_1 L_{\min}^{-1}, C_2}
\to \W_Y^{\leq \kappa_1 + \kappa_2} / I^{\leq \kappa_1 + \kappa_2}_{C_0, C_1, C_2}.
\]

Let $(\hbar^{-1} \W_Y^{\leq 0})^+ \subset \hbar^{-1} \W_Y^{\leq 0}$
be the submodule defined by
\begin{align*}
&(\hbar^{-1} \W_Y^{\leq 0})^+ \\
&= \bigl\{\sum
a_{(x_i), (\hat c_i^\ast), (\hat c'_i), g}\ t_{x_1} \dots t_{x_{k_t}}
q_{\hat c_1^\ast} \dots q_{\hat c_{k_q}^\ast}
p_{\hat c'_1} \dots p_{\hat c'_{k_p}} \hbar^g \in \hbar^{-1} \W_Y^{\leq 0};
\widetilde{g} \geq 0\bigr\},
\end{align*}
and $(\hbar^{-1} \W_Y^{\leq 0})^+_{C_0, C_1, C_2} \subset (\hbar^{-1} \W_Y^{\leq 0})^+$
be the submodule defined by
\begin{align*}
&(\hbar^{-1} \W_Y^{\leq 0})^+_{C_0, C_1, C_2} \\
&= \Bigl\{\sum
a_{(x_i), (\hat c_i^\ast), (\hat c'_i), g}\ t_{x_1} \dots t_{x_{k_t}}
q_{\hat c_1^\ast} \dots q_{\hat c_{k_q}^\ast}
p_{\hat c'_1} \dots p_{\hat c'_{k_p}} \hbar^g \in (\hbar^{-1} \W_Y^{\leq 0})^+;\\
&\quad \quad a_{(x_i), (\hat c_i^\ast), (\hat c'_i), g} = 0
\text{ for all } ((x_i)_{i = 1}^{k_t}, (\hat c_i^\ast),_{i = 1}^{k_q} (\hat c'_i)_{i = 1}^{k_p}, g)
\text{ such that}\\
&\quad \quad k_t \leq C_0,\, \widetilde{g} \leq C_1 \text{ and }
\sum e(p_{\hat c'_i}) \geq - C_2 \Bigr\}
\end{align*}
for each triple $(C_0, C_1, C_2)$.

If we fix a triple $(\overline{C}_0, \overline{C}_1, \overline{C}_2)$,
then, choosing a compatible family of perturbations $\B$ of the multisections 
of finite number of pre-Kuranishi spaces and
using their virtual fundamental chains,
we can define the generating function
$\mathcal{H} = \mathcal{H}_{(Y, \lambda, K_Y, K_Y^0, K_Y^2, J, \B)}
= \hbar^{-1} \sum_g \mathcal{H}_g \hbar^g \in
(\hbar^{-1} \W_Y^{\leq 0})^+
/ (\hbar^{-1} \W_Y^{\leq 0})^+_{\overline{C}_0, \overline{C}_1, \overline{C}_2}$ by
\[
\mathcal{H}_g = \sum_{k_q, k_t, k_p \geq 0} \frac{1}{k_q !k_t ! k_p !}
\bigl[\overline{\M}_g(\underbrace{\mathbf{q}, \dots, \mathbf{q}}_{k_q};
\underbrace{\mathbf{t}, \dots, \mathbf{t}}_{k_t};
\underbrace{\mathbf{p}, \dots, \mathbf{p}}_{k_p})\bigr]^0,
\]
where $\mathbf{q} = \sum_c q_{\hat c^\ast} \hat c$,
$\mathbf{t} = \sum_x t_x x$ and
$\mathbf{p} = \sum_c p_{\hat c} \hat c^\ast$ are formal series.
We need to check that
$\mathcal{H}$ is indeed an element of $(\hbar^{-1} \W_Y^{\leq 0})^+$, that is, every
holomorphic building satisfies
\[
\sum_j L_{\gamma_{+\infty_j}} - \sum_i L_{\gamma_{-\infty_i}} \geq 0
\]
and
\[
\tilde g = g + \frac{1}{2}(k_t + k_q + k_p)
+ \frac{\sum_j L_{\gamma_{+\infty_j}} - \sum_i L_{\gamma_{-\infty_i}}}
{L_{\min}} \geq 1,
\]
where $g$ is its genus, $k_t$, $k_q$ and $k_p$ are the numbers of its marked points,
$-\infty$-limit circles, and $+\infty$-limit circles respectively, and
$L_{\gamma_{\pm\infty_i}}$ are the periods of the periodic orbits on its limit circles.
The former is because the left hand side is the $E_{\hat \omega}$-energy.
The latter is proved as follows.
First note that $\tilde g - 1$ is additive with respect to disjoint union or
gluing at limit circles.
Hence it is enough to prove the case of a connected holomorphic building of height one.
Assume that there exists a connected holomorphic building $(\Sigma, z, u, \phi)$
of height one such that $\tilde g < 1$.
Since $\tilde g < 1$ implies $g = 0$ and $k_t \leq 1$, $u$ is not a constant map.
Since $\sigma \circ u$ cannot attain a maximum at the interior, it implies that
$k_p \geq 1$.
Therefore $\tilde g < 1$ implies $k_q = 0$ and $k_p = 1$.
However, this implies
\[
\frac{\sum_j L_{\gamma_{+\infty_j}} - \sum_i L_{\gamma_{-\infty_i}}}
{L_{\min}}
= \frac{L_{\gamma_{+\infty_1}}}
{L_{\min}}
\geq 1,
\]
which contradict the assumption $\tilde g < 1$.
Therefore $\mathcal{H}$ is an element of $(\hbar^{-1} \W_Y^{\leq 0})^+$.
We also note that $\mathcal{H}$ has the odd degree.

Define a differential $\delta : \W_Y \to \W_Y$ by
$\delta q_\alpha = q_{\delta \alpha}$,
$\delta t_x = 0$,
$\delta p_{\hat c} = (-1)^{1 + |\hat c|} p_{\partial \hat c}$ and $\delta \hbar = 0$.
(Note that this is well defined, that is, $\delta [p_{\hat c}, q_\alpha] = 0$.)
Note the following equations:
\[
\sum_c \delta q_{\hat c^\ast} \hat c = \sum_c q_{\hat c^\ast} \partial \hat c, \quad
\sum_c \delta p_{\hat c} \hat c^\ast = \sum_c p_{\hat c} \partial \hat c^\ast.
\]
(Recall that we have defined $\partial \alpha$ by
$\partial \alpha = (-1)^{|\alpha|} \delta \alpha$ for a cochain $\alpha$.)
We also define the differential $\delta$ on $\hbar^{-1} \W_Y$ similarly.
Then equation (\ref{boundary formula}) implies
\begin{equation}
\delta \mathcal{H} - \mathcal{H} \circ \mathcal{H} = 0
\label{main eq}
\end{equation}
in $(\hbar^{-1} \W_Y^{\leq 0})^+
/ (\hbar^{-1} \W_Y^{\leq 0})^+_{\overline{C}_0, \overline{C}_1, \overline{C}_2}$.

For each four-tuple $(\kappa, C_0, C_1, C_2)$ such that
$\overline{C}_0 \geq C_0$, $\overline{C}_1 \geq C_1 + \frac{\kappa}{L_{\min}}$
and $\overline{C}_2 \geq C_2 + \kappa$,
define a linear map $D_Y = D_{(Y, \lambda, K_Y, K_Y^0, K_Y^2, J, \B)}
: \W_Y^{\leq \kappa} / I^{\leq \kappa}_{C_0, C_1, C_2} \to
\W_Y^{\leq \kappa} / I^{\leq \kappa}_{C_0, C_1, C_2}$ by
\[
D_Y f = \delta f - [\mathcal{H}, f].
\]
Then $D_Y$ is a differential, that is,
\begin{align}
D_Y^2 &= 0 \label{D_Y^2}\\
D_Y(fg) &= (D_Y f) g + (-1)^{|f|} f D_Y g \label{D_Y diff alg}
\end{align}
(\ref{D_Y^2}) is a consequence of (\ref{main eq}).
(\ref{D_Y diff alg}) holds if the multiplications are well defined.
Namely, for $f \in \W_Y^{\leq \kappa_1}
/ I^{\leq \kappa_1}_{C_0, C_1 + \kappa_2 L_{\min}^{-1}, C_2 + \kappa_2}$
and $g \in \W_Y^{\leq \kappa_2}
/ I^{\leq \kappa_2}_{C_0, C_1 + \kappa_1 L_{\min}^{-1}, C_2}$,
(\ref{D_Y diff alg}) holds in $\W_Y^{\leq \kappa_1 + \kappa_2}
/ I^{\leq \kappa_1 + \kappa_2}_{C_0, C_1, C_2}$.
We denote the homology of the chain complex
$(\W_Y^{\leq \kappa} / I^{\leq \kappa}_{C_0, C_1, C_2}, D_Y)$
by $H^\ast(\W_Y^{\leq \kappa} / I^{\leq \kappa}_{C_0, C_1, C_2}, D_Y)
= \Ker D_Y / \Image D_Y$.

We will prove that the homology
\[
H^\ast(\W_{(Y, \lambda, K_Y, \overline{K}_Y^0)}^{\leq \kappa}
/ I^{\leq \kappa}_{C_0, C_1, C_2}, D_{(Y, \lambda, K_Y, K_Y^0, K_Y^2, J, \B)})
\]
is independent of the choice of $(K_Y, K_Y^0, K_Y^2, J, \B)$ in Section \ref{independence}
(Lemma \ref{short concordance for Y}).
Therefore we can define the limit
\begin{align*}
&H^\ast(\W_{(Y, \lambda, \overline{K}_Y^0)}, D_{(Y, \lambda, \overline{K}_Y^0)})\\
&= \varprojlim_{C_2} \varinjlim_{\kappa} \varprojlim_{C_0, C_1}
H^\ast(\W_{(Y, \lambda, K_Y, \overline{K}_Y^0)}^{\leq \kappa}
/ I^{\leq \kappa}_{C_0, C_1, C_2}, D_{(Y, \lambda, K_Y, K_Y^0, K_Y^2, J, \B)}).
\end{align*}
(\ref{D_Y diff alg}) implies that this is an algebra.
We will also prove that this is independent of the choice of the contact form $\lambda$
of the contact manifold $(Y, \xi)$ in Section \ref{independence}.

\begin{rem}
We can use the spectral sequence defined by the filtration given by the energy
$\sum e(q_{\hat c_i^\ast}) + \sum e(p_{\hat c'_i})$ for each
$H^\ast(\W_Y^{\leq \kappa} / I^{\leq \kappa}_{C_0, C_1, C_2}, D_Y)$
since $\W_Y^{\leq \kappa} / I^{\leq \kappa}_{C_0, C_1, C_2}$ is finite dimensional.
\end{rem}

Next we briefly explain the construction of rational symplectic field theory.
Define a super-commutative algebra $\mathcal{P}_Y = \mathcal{P}_{(Y, \lambda, K_Y,
\overline{K}^0_Y)}$ by $\mathcal{P}_Y = \W_Y|_{\hbar = 0}$.
It is regarded as a quotient of $\W_Y$.
Its (graded) Poisson structure is defined by
\begin{align*}
\{f, g\} &= (\hbar^{-1} [f, g])|_{\hbar = 0}\\
&= \sum_c \biggl(\frac{\overleftarrow{\partial} f}{\partial p_{\hat c}}
\frac{\overrightarrow{\partial} g}{\partial q_{\hat c^\ast}}
- (-1)^{|f| |g|} \frac{\overleftarrow{\partial} g}{\partial p_{\hat c}}
\frac{\overrightarrow{\partial} f}{\partial q_{\hat c^\ast}}\biggr),
\end{align*}
where $\overrightarrow{\partial}$ and $\overleftarrow{\partial}$ are
differential from left and right respectively.
It is easy to check that it is indeed a Poisson structure, that is, it satisfies the following
equations:
\begin{align*}
\{f, gh\} &= \{f, g\} h + (-1)^{|f| |g|} g\{f, h\},\\
\{g, f\} &= -(-1)^{|f| |g|} \{f, g\},\\
\{\{f, g\}, h\} &= \{f, \{g, h\}\} - (-1)^{|f| |g|} \{g, \{f, h\}\}.
\end{align*}

The differential $\delta : \mathcal{P}_Y \to \mathcal{P}_Y$ is defined similarly
to the case of $\W_Y$.

For each $\kappa \geq 0$, we define a submodule $\mathcal{P}^{\leq \kappa}_Y
\subset \mathcal{P}_Y$ by imposing the condition
$\sum_i e(q_{\hat c_i^\ast}) + \sum_i e(p_{\hat c'_i}) \leq \kappa$.
For each triple $(\kappa, C_0, C_2)$,
We define a submodule $I^{\leq \kappa}_{C_0, C_2} \subset
\mathcal{P}_Y^{\leq \kappa}$ by
\begin{align*}
I^{\leq \kappa}_{C_0, C_2} &= \Bigl\{\sum
a_{(x_i), (\hat c_i^\ast), (\hat c'_i)}\ t_{x_1} \dots t_{x_{k_t}}
q_{\hat c_1^\ast} \dots q_{\hat c_{k_q}^\ast}
p_{\hat c'_1} \dots p_{\hat c'_{k_p}} \in \mathcal{P}_Y^{\leq \kappa};\\
&\quad \quad a_{(x_i), (\hat c_i^\ast), (\hat c'_i)} = 0
\text{ for all } ((x_i)_{i = 1}^{k_t}, (\hat c_i^\ast),_{i = 1}^{k_q} (\hat c'_i)_{i = 1}^{k_p})
\text{ such that}\\
&\quad \quad k_t \leq C_0 \text{ and }
\sum e(p_{\hat c'_i}) \geq - C_2 \Bigr\}.
\end{align*}
In this case the following holds true.
\[
\mathcal{P}_Y \cong \varprojlim_{C_2} \varinjlim_{\kappa} \varprojlim_{C_0}
\mathcal{P}_Y^{\leq \kappa} / I^{\leq \kappa}_{C_0, C_2}.
\]
Note that the Poisson bracket induces the maps
\[
\mathcal{P}_Y^{\leq \kappa_1} / I^{\leq \kappa_1}_{C_0, C_2 + \kappa_2}
\times \mathcal{P}_Y^{\leq \kappa_2} / I^{\leq \kappa_2}_{C_0, C_2 + \kappa_1}
\to \mathcal{P}_Y^{\leq \kappa_1 + \kappa_2} / I^{\leq \kappa_1 + \kappa_2}_{C_0, C_2}.
\]
More generally, the Poisson bracket induces the maps
\begin{align}
&(\mathcal{P}_Y^{\leq \kappa^\circ_1} + I^{\leq \kappa_1}_{C_0, C_2})
/ I^{\leq \kappa_1}_{C_0, C_2 + \kappa^\circ_2} \times
(\mathcal{P}_Y^{\leq \kappa^\circ_2} + I^{\leq \kappa_2}_{C_0, C_2})
/ I^{\leq \kappa_1}_{C_0, C_2 + \kappa^\circ_1} \notag\\
&\hspace{180pt}
\to \mathcal{P}_Y^{\leq \kappa_1 + \kappa_2} / I^{\leq \kappa_1 + \kappa_2}_{C_0, C_2}
\label{Poisson bracket for fiber product}
\end{align}
for $\kappa^\circ_i \leq \kappa_i$ ($i = 1,2$).
Note that for $\kappa^\circ \leq \kappa$ and $C^\circ \leq C$,
$(\mathcal{P}_Y^{\leq \kappa^\circ} + I^{\leq \kappa}_{C_0, C^\circ_2})
/ I^{\leq \kappa}_{C_0, C_2}$ is the fiber product of
$\mathcal{P}_Y^{\leq \kappa^\circ} / I^{\leq \kappa^\circ}_{C_0, C^\circ_2}$
and $\mathcal{P}_Y^{\leq \kappa} / I^{\leq \kappa}_{C_0, C_2}$ over
$\mathcal{P}_Y^{\leq \kappa} / I^{\leq \kappa}_{C_0, C^\circ_2}$.

Equation (\ref{main eq}) implies that $\mathcal{H}_0 \in \mathcal{P}_Y^{\leq 0}
/ I^{\leq 0}_{\overline{C}_0, \overline{C}_2}$ satisfies
\begin{equation}
\delta \mathcal{H}_0 - \frac{1}{2} \{\mathcal{H}_0, \mathcal{H}_0\} = 0.
\label{main eq for rational}
\end{equation}
in $\mathcal{P}_Y^{\leq 0} / I^{\leq 0}_{\overline{C}_0, \overline{C}_2}$.

For each triple $(\kappa, C_0, C_2)$ such that $\overline{C}_0 \geq C_0$,
$\overline{C}_2 \geq C_2 + \kappa$,
define a linear map $d_Y = d_{(Y, \lambda, K_Y, K_Y^0, K_Y^2, J, \B)}
: \mathcal{P}_Y^{\leq \kappa} / I^{\leq \kappa}_{C_0, C_2}
\to \mathcal{P}_Y^{\leq \kappa} / I^{\leq \kappa}_{C_0, C_2}$ by
\[
d_Y f = \delta f - \{\mathcal{H}_0, f\} \ (= D_Y f|_{\hbar = 0}).
\]
Then $d_Y$ satisfies the following.
\begin{align}
d_Y^2 &= 0, \label{d_Y^2}\\
d_Y(fg) &=  (d_Y f) g + (-1)^{|f|} f d_Y g \label{d_Y Leibnitz}\\
d_Y \{f, g\} &= \{d_Y f, g\} + (-1)^{|f|} \{f, d_Y g\}.
\label{d_Y Poisson}
\end{align}
(\ref{d_Y^2}) is due to (\ref{main eq for rational}).
(\ref{d_Y Leibnitz}) and (\ref{d_Y Poisson}) hold if the multiplications or
Poisson brackets are well defined.
We denote the cohomology of the complex
$(\mathcal{P}_Y^{\leq \kappa} / I^{\leq \kappa}_{C_0, C_2}, d_Y)$
by $H^\ast(\mathcal{P}_Y^{\leq \kappa} / I^{\leq \kappa}_{C_0, C_2}, d_Y)$.
We remark that $(\mathcal{P}_Y^{\leq \kappa} / I^{\leq \kappa}_{C_0, C_2}, d_Y)$
can be regarded as a quotient of the chain complex of
general symplectic cohomology by the ideal $(\hbar)$.
As in the case of general SFT, we will define rational SFT cohomology as a limit
\begin{align*}
&H^\ast(\mathcal{P}_{(Y, \lambda, \overline{K}^0_Y)}, d_{(Y, \lambda, \overline{K}^0_Y)}) \\
&= \varprojlim_{C_2} \varinjlim_{\kappa} \varprojlim_{C_0}
H^\ast(\mathcal{P}_{(Y, \lambda, K_Y, \overline{K}^0_Y)}^{\leq \kappa}
/ I^{\leq \kappa}_{C_0, C_2}, d_{(Y, \lambda, K_Y, K_Y^0, K_Y^2, J, \B)}).
\end{align*}

Finally, we consider the construction of contact homology.
We use the super-commutative algebra $\A_Y = \A_{(Y, \lambda, K_Y, \overline{K}^0_Y)}$
defined by $\A_Y = \R[[t]](q)$.
Its elements are written as
\[
\sum_{(\hat c_i)} f_{(\hat c_i)}(t) q_{\hat c_1^\ast} \dots q_{\hat c_{k_q}^\ast},
\]
where $f_{(\hat c_i)}(t) \in \R[[t]]$ are formal series of the variables $t_x$
and the sum is a finite sum.
For each $\kappa \geq 0$, we define a submodule $\A_Y^{\leq \kappa} \subset
\A_Y$ by
\[
\A_Y^{\leq \kappa} = \{ \sum_{(\hat c_i)} f_{(\hat c_i)}(t) q_{\hat c_1^\ast} \dots
q_{\hat c_{k_q}^\ast} \in \A_Y; \sum_i |e(q_{\hat c_i^\ast})| \leq \kappa \text{ if }
f_{(\hat c_i)}(t) \neq 0\}.
\]
For each $C_0 \geq 0$, we also define a submodule $I^{\leq \kappa}_{C_0}
\subset \A^{\leq \kappa}_Y$ by
\[
I^{\leq \kappa}_{C_0}
= \{\sum_{(x_i), (\hat c_i^\ast)} a_{(x_i); (\hat c_i^\ast)} t_{x_1} \dots t_{x_{k_t}}
q_{\hat c_1^\ast} \dots q_{\hat c_{k_q}^\ast} \in \A_Y^{\leq \kappa};
a_{(x_i); (\hat c_i^\ast)} = 0 \text{ for } k_t \leq C_0\}.
\]

Let
\[
\widehat{\mathcal{H}}_0 = \sum_{c}
\frac{\overleftarrow{\partial} \mathcal{H}_0}{\partial p_{\hat c}} \biggr|_{p = 0}
\cdot p_{\hat c} \in \mathcal{P}_Y^{\leq 0} / I^{\leq 0}_{\overline{C}_0, \overline{C}_2}
\]
be the homogeneous component of degree $1$ with respect to the variables $p_{\hat c}$.
Then equation (\ref{main eq for rational}) implies
\begin{equation}
\delta \widehat{\mathcal{H}}_0
- \frac{1}{2} \{\widehat{\mathcal{H}}_0, \widehat{\mathcal{H}}_0\} = 0
\label{main eq for contact}
\end{equation}
in $\mathcal{P}_Y^{\leq 0} / I^{\leq 0}_{\overline{C}_0, \overline{C}_2}$
because $\partial_{q_{\hat c^\ast}} (\mathcal{H}_0|_{p = 0}) = 0$ implies
$\{ \cdot, \mathcal{H}_0|_{p = 0} \} = 0$.
For each pair $(\kappa, C_0)$ such that $\overline{C}_0 \geq C_0$ and
$\overline{C}_2 \geq \kappa$,
define a linear map $\partial_Y = \partial_{(Y, \lambda, K_Y, K_Y^0, K_Y^2, J, \B)}
: \A_Y^{\leq \kappa} / I^{\leq \kappa}_{C_0} \to
\A_Y^{\leq \kappa} / I^{\leq \kappa}_{C_0}$ by
\[
\partial_Y f = \delta f - \{\widehat{\mathcal{H}}_0, f\} \ (= d_Y f|_{p = 0}).
\]
Then it satisfies the following equations:
\begin{align}
\partial_Y^2 &= 0, \label{partial_Y^2} \\
\partial_Y (fg) &= (\partial_Y f) g + (-1)^{|f|} f \partial_Y g.
\label{partial diff alg}
\end{align}
((\ref{partial_Y^2}) is due to (\ref{main eq for contact}) and
((\ref{partial diff alg}) is satisfied if the multiplication is well defined.))

As in the other two theory, we define the contact homologies by their limit
\[
H^\ast(\A_{(Y, \lambda, \overline{K}^0_Y)}, \partial_{(Y, \lambda, \overline{K}^0_Y)})
= \varinjlim_\kappa \varprojlim_{C_0}
H^\ast(\A_{(Y, \lambda, K_Y, \overline{K}^0_Y)}^{\leq \kappa}
/ I^{\leq \kappa}_{C_0}, \partial_{(Y, \lambda, K_Y, K_Y^0, K_Y^2, J, \B)}).
\]


%% file: SFT-08_Fiber_products_for_X.tex
%
%

\section{The case of holomorphic buildings for $X$}
\label{case of X}
In this section, we construct the algebras for a symplectic manifold $X$ with
cylindrical ends.
First we explain the construction of a pre-Kuranishi structure of the space
of holomorphic buildings for $X$ in Section \ref{Kuranishi for X}.
In Section \ref{fiber prod for X},
we construct the perturbed multisections of its fiber products.
We construct the correction terms for $X$ in Section \ref{correction terms for X}
and finally we construct the algebras for $X$ in Section \ref{algebra for X}
\subsection{Construction of pre-Kuranishi spaces for $X$}
\label{Kuranishi for X}
In this section, we construct a pre-Kuranishi structure of
$\widehat{\M}^X = \widehat{\M}(X, \omega, J)$.
The construction is almost the same as the case of the symplectization $\hat Y$.

First we explain the construction of a Kuranishi neighborhood of a point
$p_0 \in \widehat{\M}(X, \omega, J)$.
As in the case of $\hat Y$, we assume the following data
$(p_0^+, S = (S_X, S_{Y^\pm}), E^0, \lambda)$ are given:
\begin{itemize}
\item
$p_0^+ = (\Sigma_0, z \cup z^+, u_0)$ is a curve obtained by adding marked points
on the nontrivial components of $\Sigma_0$.
We assume all unstable components of $(\Sigma_0, z \cup z^+)$ are trivial cylinders
of $p_0$, and $G_0 = \Aut(\Sigma_0, z, u_0)$ preserves $z^+$ as a set.
\item
$S_X \subset X$ and $S_{Y^\pm} \subset Y^\pm$ are finite unions of codimension-two
submanifolds such that
$u_0$ intersects with $(-\infty, 0) \times S_{Y^-} \cup S \cup (0, \infty) \times S_{Y^+}$
at $z^+ \cap \bigcup_{i(\alpha) = 0} \Sigma_\alpha$ transversely,
$\pi_{Y^-} \circ u_0$ intersects with $S_{Y^-}$ at
$z^+ \cap \bigcup_{i(\alpha) < 0} \Sigma_\alpha$ transversely, and
$\pi_{Y^+} \circ u_0$ intersects with $S_{Y^+}$ with at
$z^+ \cap \bigcup_{i(\alpha) > 0} \Sigma_\alpha$ transversely.
\item
Let $(\hat \Sigma_0, z \cup z^+ \cup (\pm\infty_i))$ be the stabilization of
$(\check \Sigma_0, z \cup z^+ \cup (\pm\infty_i))$, and
let $(\hat P \to \hat X, Z \cup Z^+ \cup Z_{\pm\infty_i})$ be its local universal family.
$G_0$ acts on $\hat P$ by the universal property.
Then  we assume $E^0$ is a finite dimensional $G_0$-vector space,
and $\lambda = (\lambda_X, \lambda_{Y^-}, \lambda_{Y^+})$ is a family of
$G_0$-equivariant linear maps
$\lambda_X : E^0 \to C^\infty(\hat P \times X; \Wedge^{0, 1}V^\ast \hat P \otimes TX)$
and $\lambda_{Y^\pm} : E^0 \to
C^\infty(\hat P \times Y^\pm; \Wedge^{0, 1}V^\ast \hat P \otimes (\R \partial_\sigma \oplus TY^\pm))$
which satisfies following conditions:
\begin{itemize}
\item
$\lambda_X$ and $\lambda_{Y^\pm}$ are related by
\begin{align*}
&\lambda_X(h)|_{\hat P \times ((-\infty, -T] \times Y^-)}(z, \sigma, y)
= \lambda_{Y^-}(h)(z, y)\\
&\lambda_X(h)|_{\hat P \times ([T, \infty) \times Y^+)}(z, \sigma, y) = \lambda_{Y^+}(h)(z, y)
\end{align*}
for some $T \geq 0$.
\item
For each $h \in E^0$, the projections of the support of $\lambda_X(h)$ or
$\lambda_{Y^\pm}(h)$ do not intersect with the nodal points of $\hat P$ or $Z_{\pm\infty_i}$.
\item
Let $E^0 \to C^\infty(\Sigma_0, \Wedge^{0, 1} T^\ast \Sigma_0 \otimes u_0^\ast TX)$
be the linear map defined by the pullbacks of $\lambda_X$ and $\lambda_{Y^\pm}$ by
the composition of the blowing down $\Sigma_0 \to \check \Sigma_0$ and the forgetful
map $(\check \Sigma_0, z \cup z^+) \tocong (\hat P_0, Z(0) \cup Z^+(0))$,
where $u_0^\ast TX$ is the vector bundle on $\Sigma_0$ defined by
$(u_0|_{\bigcup_{i(\alpha) = 0} \Sigma_\alpha})^\ast TX$,
$(u_0|_{\bigcup_{i(\alpha) < 0} \Sigma_\alpha})^\ast T \hat Y^-$ and
$(u_0|_{\bigcup_{i(\alpha) > 0} \Sigma_\alpha})^\ast T \hat Y^+$.
Then we assume that for a sufficiently small $\delta > 0$,
the linear map
\begin{align*}
&D_{p_0}^+ : \widetilde{W}_\delta^{1, p}(\Sigma_0, u_0^\ast TX) \oplus E^0\\
&\to
L_\delta^p(\Sigma_0, \Wedge^{0, 1}T^\ast \Sigma \otimes u_0^\ast TX)
\oplus
\bigoplus_{\text{limit circles}} \Ker A_{\gamma_{\pm\infty_i}} / (\R \partial_\sigma \oplus
\R R_\lambda)\\
&\quad \oplus \bigoplus_{z_i \in \bigcup_{i(\alpha) \gtrless 0} \Sigma_\alpha}
T_{\pi_{Y^\pm} \circ u_0(z_i)} Y^\pm
\oplus \bigoplus_{z_i \in \bigcup_{i(\alpha) = 0} \Sigma_\alpha} T_{u_0(z_i)} X\\
&(\xi, h) \mapsto (D_{p_0} \xi + \lambda(h),
\sum_j \langle \xi|_{S^1_{\pm\infty_i}}, \eta^{\pm\infty_i}_j \rangle \eta^{\pm\infty_i}_j,
\pi_{Y^\pm} \circ \xi (z_i), \xi(z_i))
\end{align*}
is surjective,
where $D_{p_0}$ is a linearization of the equation of the $J$-holomorphic maps,
that is,
\[
D_{p_0} \xi = \nabla \xi + J(u_0) \nabla \xi j + \nabla_\xi J(u_0) du_0 j,
\]
and each $\{\eta^{\pm\infty_i}_j\}_j$ is an orthonormal basis of
the complement of $\R \partial_\sigma \oplus \R R_\lambda$ in
$\Ker A_{\gamma_{\pm\infty_i}}$.
\end{itemize}
\end{itemize}

We fix the following temporally data $(z^{++}, S', (\hat R_i)_{i \neq 0})$:
\begin{itemize}
\item
$z^{++} = (z^{++}_i) \subset \Sigma$ are additional marked points which make
$(\Sigma_0, z \cup z^+ \cup z^{++})$ stable.
We assume that $G_0$-action preserves $z^{++}$ as a set.
\item
$S' \subset (\R_{-k_-} \cup \dots \cup \R_{-1}) \times Y^-
\cup (\R_1 \cup \dots \cup \R_k) \times Y^+$ is a codimension-two submanifold
such that $u_0$ intersects with $S'$ at $z^{++}$ transversely.
\item
For each $-k_- \leq i \leq -1$ or $1 \leq i \leq k_+$,
let $\hat R_i = (\hat R_{i, l})_{1 \leq l \leq m_i}$ be a family of holomorphic sections
$\hat R_{i, l} : \hat X \to \hat P$ such that
$\sigma_i \circ u_0(\hat R_{i, l}(0)) = 0$, where $\sigma_i$ is the coordinate of $\R_i$,
and $\hat R_i$ is $G_0$-invariant as a family.
We assume $\hat R_i$ do not intersect with nodal points or $Z_{\pm\infty_i}$.
Let $(\widetilde{P} \to \widetilde{X}, Z \cup Z^+ \cup Z^{++})$ be the local universal family
of $(\Sigma_0, z \cup z^+ \cup z^{++})$.
Then each $\hat R_{i, l}$ induces a section $\widetilde{R}_{i, l} : \widetilde{X} \to \widetilde{P}$
which makes following diagram commutative.
\[
\begin{tikzcd}
\widetilde{P} \ar{r}{\forget}& \hat P\\
\widetilde{X} \ar{u}{\widetilde{R}_{i, l}} \ar{r}{\forget}& \hat X \ar{u}{\hat R_{i, l}}
\end{tikzcd}
\]
We use these families of sections $\widetilde{R}_{i, l}$ to kill the $\R$-translations
by imposing the conditions $\sum_l \sigma_i \circ u(\widetilde{R}_{i, l}) = 0$ ($i \neq 0$).
\end{itemize}
The pullbacks $E^0 \to C^\infty(\widetilde{P} \times X, \Wedge^{0, 1} V^\ast \widetilde{P}
\otimes TX)$ of $\lambda_X$ and $E^0 \to C^\infty(\widetilde{P} \times Y^\pm,
\Wedge^{0, 1} V^\ast \widetilde{P} \otimes (\R \partial_\sigma \oplus TY^\pm))$ of $\lambda_{Y^\pm}$
by $\widetilde{P} \to \hat P$ are also denoted by $\lambda_X$ and
$\lambda_{Y^\pm}$ respectively.

Using the above data, we can construct a smooth Kuranishi neighborhood of $p_0$ as in
Section \ref{construction of nbds}.
The main differences are the following two.
One is that the range of $s^0$ does not contain the factor
to kill the $\R$-translation of $0$-th floor.
The other is about smoothness.
The former does not have any difficulty but the latter do.
The difficulty is that in general, the constants $\alpha$ and $\beta$ for $Y^-$
which determine the differential structure of the parameter space of deformation of
the domain curve and the target space are different from those for $Y^+$.
We explain about this issue in the next section.

The definition of the embedding of a Kuranishi neighborhood to another is also similar.
We can construct an essential submersion from a Kuranishi neighborhood of a
disconnected holomorphic building to the product of those of its connected components.
A holomorphic building for $X$ is also decomposed by its floor structure, and the
relation between the Kuranishi neighborhood of the whole holomorphic building and
the Kuranishi neighborhoods of the parts is similar to the case of $Y$.

Next we consider the construction of a global pre-Kuranishi structure.
Similarly to the case of symplectization $\hat Y$,
we construct a domain curve representation
of the space of holomorphic buildings and Kuranishi data.
The main differences are as follows.

First, instead of a set $\mathcal{S}$ of codimension-two submanifolds of $Y$,
we use a triple $(\mathcal{S}_X, \mathcal{S}_{Y^-}, \mathcal{S}_{Y^+})$ consisting of
sets of codimension two submanifolds of $X$, $Y^-$ and $Y^+$ respectively.
Hence for example, we use the space
$\widehat{\M}^X_{(\mathcal{S}_X, \mathcal{S}_{Y^\pm})}$
of points $(\Sigma, z,
(z^S)_{S \in \mathcal{S}_X \cup \mathcal{S}_{Y^-} \cup \mathcal{S}_{Y^+}}, u)$ which
satisfy the following conditions:
\begin{itemize}
\item
For each $S \in \mathcal{S}_X$, $z^S$ is contained in the $0$-th floor and
$u$ intersects $S$ at $z^S$ transversely,
\item
For each $S \in \mathcal{S}_{Y^-}$, $z^S$ is contained in the union of
the negative floors and the inverse image of $(-\infty, 0] \times Y^- \subset X$ by $u$,
and $u$ intersects $\R \times S$ at $z^S$ transversely
\item
For each $S \in \mathcal{S}_{Y^+}$, $z^S$ is contained in the union of
the positive floors and the inverse image of $[0, \infty) \times Y^+ \subset X$ by $u$,
and $u$ intersects $\R \times S$ at $z^S$ transversely.
\end{itemize}
Similarly, we use a triple $(\lambda_X, \lambda_{Y^-}, \lambda_{Y^+})$ instead of
$\lambda$.

Second, instead of fixing one constant $L_{\max}$, we fix two constants
$L^1_{\max}$ and $L^2_{\max}$ and consider the subspace
$\widehat{\M}^{X, \leq (L^1_{\max}, L^2_{\max})} \subset\widehat{\M}^X$ consisting of
holomorphic buildings $(\Sigma, z, u)$ such that
$e + L^+ \leq L^1_{\max}$ and $L^+ \leq L^2_{\max}$,
where $e = \int_{u^{-1}(X)} u^\ast \widetilde{\omega}$.
(See Section \ref{holomorphic buildings for X} for the definition of $\widetilde{\omega}$
and recall the estimates (\ref{E lambda estimate}) and (\ref{E hat omega estimate}).)
Let $\delta_0 >0$ be a constant which satisfies the following conditions:
\begin{itemize}
\item
$6 \delta_0$ is less than the minimal $E_{\hat \omega}$-energy of
a non-constant $J$-holomorphic sphere in $X$.
\item
$4 \delta_0$ is less than the minimal $E_{\hat \omega}$-energy of
a holomorphic plain in $X$ whose $E_{\lambda}$-energy is
$\leq \max(L^1_{\max}, L^2_{\max})$.
\item
$2 \delta_0$ is less than the minimal $E_{\hat \omega}$-energy of
a holomorphic cylinder in $X$ whose $E_{\lambda}$-energy is
$\leq \max(L^1_{\max}, L^2_{\max})$.
\item
$\delta_0$ is less than the minimal $E_{\hat \omega}$-energy of
a non-constant $J$-holomorphic torus in $X$.
\end{itemize}
Then for a triple $\theta = (g, k, E_{\hat \omega})$, we define
$\widetilde{e}(\theta) = \widetilde{e}_{\delta_0}(\theta)
= 5(g-1) + 2k + E_{\hat \omega} / \delta_0$ as in the
case of symplectization $\hat Y$.
It is easy to check that $\widetilde{e}(p) = \widetilde{e}(\theta) \geq 1$
for any holomorphic building $p \in \widehat{\M}^{X, \leq (L^1_{\max}, L^2_{\max})}$
of type $\theta$.

Assume that domain curve representations
$(\mathcal{S}^{Y^-}, \V_{\theta, l}^{Y^-}, \U_{\theta, l}^{Y^-},
\U_{\theta, l}^{Y^-, \DM})$ and
$(\mathcal{S}^{Y^+}, \V_{\theta, l}^{Y^+}, \U_{\theta, l}^{Y^+},
\U_{\theta, l}^{Y^+, \DM})$
of $\widehat{\M}^{Y^-, \leq L^1_{\max}}_{\leq C}$ and
$\widehat{\M}^{Y^+, \leq L^2_{\max}}_{\leq C}$ are given respectively.
Then we can define a compatible domain curve representation
$(\mathcal{S}^X, \V_{\theta, l}^X, \ab \U_{\theta, l}^X, \ab \U_{\theta, l}^{X, \DM})$
of $\widehat{\M}^{X, \leq (L^1_{\max}, L^2_{\max})}_{\leq C}$ similarly.
To distinguish the negative floors, the $0$-th floor and the positive floors
in the space of domain curves, we add new marked points
$z^{B^-}$, $z^{B^0}$ and $z^{B^+}$ similarly to $z^A$ so that
\begin{itemize}
\item
if a irreducible component contains a marked point in $z^{B^0}$, then
it is contained in the $0$-th floor, and
\item
otherwise, it contains a marked point in either $z^{B^-}$ or $z^{B^+}$ and
in the former case, it is contained in the negative floors, and
in the latter case, it is contained in the positive floors.
\end{itemize}
(The $0$-th floor may contain marked points in $z^{B^-}$ or $z^{B^+}$.)
Hence in this case, we construct $\V^X_{\theta, l}$ and $\U^X_{\theta, l}$
as the subspaces of the space
$\widehat{\M}^X_{(\mathcal{S}_X, \mathcal{S}_{Y^\pm}), A, B^-, B^0, B^+}$
consisting of holomorphic buildings with marked point $z$, $z^A$, $z^{B^-}$, $z^{B^0}$
and $z^{B^+}$.
To formulate the compatibility condition with the domain curve representations
of $(\mathcal{S}^{Y^\pm}, \V_{\theta, l}^{Y^\pm}, \U_{\theta, l}^{Y^\pm},
\U_{\theta, l}^{Y^\pm, \DM})$, first we add marked points $z^{B^\pm}$ to
the curves in $\U_{\theta, l}^{Y^\pm}$ and $\U_{\theta, l}^{Y^\pm, \DM}$ which satisfy
the conditions similar to the marked points $z^A$.
Then the compatibility conditions are formulated as follows instead of Condition
\ref{decomposition into parts U DM}, \ref{decomposition into parts U} and
\ref{decomposition into parts V}.
\begin{enumerate}
\item[\ref{decomposition into parts U DM}$\:\!\!^X$]
For any $\theta = (g, k, E_{\hat \omega})$,
$\hat p \in \U_{\theta, l}^{X, \DM}$ and subset $\mathcal{N}$ of its nodal
points, replace each nodal point in $\mathcal{N}$ with a pair of marked points
(we regard the new marked points as points in the set $z$), and
let $\hat p'_i$ $(1 \leq i \leq N)$ be its connected components or an arbitrary
decomposition into unions of its connected components.
Let $g'_i$ and $k'_i$ be the genus and the number of marked points $z$ of
each $\hat p'_i$ respectively.
Then there exist some $E_{\hat \omega}^i \geq 0$ such that
$E_{\hat \omega} = \sum_i E_{\hat \omega}^i$ and the following hold for
$\theta'_i = (g'_i, k'_i, E_{\hat \omega}^i)$.
\begin{itemize}
\item
$\hat p'_i \in \U_{\theta'_i, l(\hat p'_i)}^{X, \DM}$ if $\hat p'_i$
contains a marked point in $z^{B^0}$.
\item
$\hat p'_i \in \U_{\theta'_i, l(\hat p'_i)}^{Y^\pm, \DM}$ if $\hat p'_i$ does not contain
any marked points in $z^{B^0}$ and it contains a marked point in $z^{B^\pm}$.
\end{itemize} 
\item[\ref{decomposition into parts U}$\:\!\!^X$]
$\U_{\theta, l}^X$ satisfies the following conditions.
\begin{itemize}
\item
For any $p \in \U_{\theta, l}^X$ and any decomposition $p_i$ ($1 \leq k$)
into unions of its connected components, let $p'_i$ be the holomorphic buildings
obtained by collapsing trivial floors.
Then $p'_i \in \U_{\theta(p'_i), l(p'_i)}^X$ for all $i$.
\item
For any $p \in \U_{\theta, l}^X$ and any gap between non-positive floors,
let $p_1$ and $p_2$ be the holomorphic buildings obtained by separating $p$ at this gap.
($p_1$ is the part in the negative floors.)
Then $p'_1 \in \U_{\theta(p'_1), l(p'_1)}^{Y^-}$ and $p'_2 \in \U_{\theta(p'_2), l(p'_2)}^X$.
We also assume the similar condition for the gap between non-negative floors.
\item
For any $p \in \U_{\theta, l}^X$ and any subset of its nodal points,
the holomorphic building $p'$ obtained by replacing these nodal points to pairs of
marked points is contained in $\U_{\theta(p'), l(p')}^X$.
\end{itemize}
\item[\ref{decomposition into parts V}$\:\!\!^X$]
For each $p \in \widehat{\M}^{X, \leq L_{\max}}_{(\mathcal{S}_X, \mathcal{S}_{Y^\pm}),
A, B^-, B^0, B^+, \theta, l}$,
replace all nodal points and joint circles of $p$ to pairs of marked points and
pairs of limit circles respectively
(we regard the new marked points as points in the set $z$), and
let $p'_i$ $(1 \leq i \leq k)$ be their non-trivial connected components.
Then $p \in \V_{\theta, l}^X$ if and only if the following hold:
\begin{itemize}
\item
$p'_i \in \V_{\theta(p'_i), l(p'_i)}^X$ if $p'_i$ contains a marked point in $z^{B^0}$.
\item
$p'_i \in \V_{\theta(p'_i), l(p'_i)}^{Y^\pm}$ if $p'_i$ does not contain a marked point
in $z^{B^0}$ and it contains a marked point in $z^{B^\pm}$.
\end{itemize}
\end{enumerate}

The definition of compatible Kuranishi data for
$\widehat{\M}^{X, \leq (L^1_{\max}, L^2_{\max})}_{\leq C}$
are also similar, and we can construct them by the same argument.
Then the pre-Kuranishi structure of each
$\widehat{\M}^{X, \leq (L^1_{\max}, L^2_{\max})}_\theta$ is defined by these data
as in the case of symplectization.

\subsection{Smoothness of pre-Kuranishi structure in the case of $X$}
Recall that in Section \ref{smoothness},
to obtain a smooth pre-Kuranishi structure of the space of holomorphic buildings
for a contact manifold, we had to use a strong differential structure of
the parameter space of the deformation of a domain curve.
Such a strong differential structure is determined by a fixed pair of large constants
$\alpha$ and $\beta$, and to construct a pre-Kuranishi structure of
the space of holomorphic buildings of higher energy,
we need to choose larger constants in general.
Hence for a cobordism $(X, \omega)$ from $(Y^-, \lambda^-)$ to $(Y^-, \lambda^-)$,
we need to consider the case where
we use different constants $\alpha^\pm$ and $\beta^\pm$ for
the smooth pre-Kuranishi structure of $\widehat{\M}(Y^\pm, \lambda^\pm, J^\pm)$.

The difference of $\beta^\pm$ does not have a difficulty.
We can use the coordinates defined by
$\rho_\mu^{L_\mu} = \hat \rho_\mu^{\beta^-}$ for a joint circle $S_\mu^1$
between non-positive floors and
$\rho_\mu^{L_\mu} = \hat \rho_\mu^{\beta^+}$ for a joint circle $S_\mu^1$
between non-negative floors.
However, for nodal points in the $0$-th floor, there is not such a $\pm$-decomposition.
Hence we need to use a gradation of smooth structures.

We fix a smooth function $\alpha : X \to \R_{>0}$ such that
$\alpha|_{(-\infty, -T] \times Y^-} = \alpha^-$ and $\alpha|_{[T, \infty) \times Y^+}
= \alpha^+$ for some $T \geq 0$.
Roughly speaking, for a nodal point $q_\nu$ of a holomorphic building
$(\Sigma_0, z, u_0) \in \widehat{\M}(X, \omega, J)$, we use the coordinate
defined by $\rho_\nu = \hat \rho_\nu^{\alpha(u_0(q_\nu))}$.

In this section, we explain the precise definition of the smooth structure of
Kuranishi neighborhoods of $\widehat{\M}(X, \omega, J)$,
and prove the smoothness of an embedding between two Kuranishi neighborhoods
or an essential submersion from that of a disconnected holomorphic building
to products of those of its connected components.

Let $(V, E, s, \psi, G)$ be a Kuranishi neighborhood of a point
$(\Sigma_0, z, u_0) \in \widehat{\M}(X, \omega, J)$.
We assume that the height of $(\Sigma_0, z, u_0)$ is $(k_-, k_+)$.
Recall that $V$ is a subset of $\hat V = \mathring{X} \times B_\epsilon(0)$
defined by $V = \{(a, b, x) \in \hat V; s^0(a, b, x) = 0\}$,
where $B_\epsilon(0)$ is a ball in the kernel of a linear operator,
and $s^0 : \hat V \to \R^{k_-} \oplus \R^{k_+} \oplus \bigoplus_{z^{++}_l} \R^2$
is a function on $\hat V$ defined similarly to (\ref{s^0}).
Let $\{\mu\}$ and $\{\nu\}$ be the indices of joint circles and nodal points
of $\Sigma_0$ respectively.
For each $i \in \{-k_-, \dots, -1\}$, let $M_i \subset \{\mu\}$ be the index set of
the joint circles between $i$-th floor and $(i+1)$-th floor,
and for $i \in \{1, \dots, k_+\}$, let $M_i \subset \{\mu\}$ be the index set of
the joint circles between $(i-1)$-th floor and $i$-th floor.
For each pair of subsets $\Pi \subset \{-k_-, \dots, -1\} \cup \{1, \dots, k_+\}$
and $\Pi' \subset \{\nu\}$,
we define $\mathring{X}_{\Pi, \Pi'} \subset \mathring{X}$ by
\begin{align*}
\mathring{X}_{\Pi, \Pi'}
= \{ (a,b) \in \mathring{X};\, & \rho_\mu \neq 0 \text{ for all } \mu \in M_i
\text{ if and only if } i \in \Pi\\
& \zeta_\nu \neq 0 \text{ if and only if } \nu \in \Pi'\},
\end{align*}
\begin{defi}
For any $0 < \epsilon < 1$ and
$\tilde \delta_0 = (\tilde \delta_{0, i})_{i \in \{-k_-, \dots, -1\} \cup \{1, \dots, k_+\}}$,
we say a continuous function $f$ on $\hat V = \mathring{X} \times B_\epsilon(0)$
is $(\epsilon, \tilde \delta_0)$-admissible if for any
$\Pi \subset \{-k_-, \dots, -1\} \cup \{1, \dots, k_+\}$ and $\Pi' \subset \{\nu\}$,
the restriction of $f$ to $\mathring{X}_{\Pi, \Pi'} \times B_\epsilon(0) \subset
\hat V$ is smooth and its differentials satisfy the following estimates
similar to those of $\phi$ given in Corollary \ref{asymptotic estimates of phi}:
For any $l \geq 1$ and any multi-index $(k_x, k_j, k_b, \ab (k_{\mu_i})_{i \in \Pi}, \ab
(l_\mu)_\mu, \ab
(k_\nu)_{\nu \in \Pi'}, \ab (l_\nu)_{\nu \in \Pi'})$,
there exists some constant $C > 0$ such that
\[
|\partial_x^{k_x} \partial_j^{k_j} \partial_b^{k_b} \partial_{(\rho_{\mu_i})}^{(k_{\mu_i})}
\partial_{(\varphi_\mu)}^{(l_\mu)} \partial_{(\rho_\nu)}^{(k_\nu)}
\partial_{(\varphi_\nu)}^{(l_\nu)} f(a, b, x)|
\leq C \prod_{\substack{i \\ k_{\mu_i} \neq 0}}
\rho_\mu^{L_{\mu_i} \tilde \delta_{0, i}/2 - k_{\mu_i}}\!\!
\prod_{\substack{\nu \\ (k_\nu, l_\nu) \neq (0,0)}} \rho_\nu^{\epsilon - k_\nu}
\]
for all $(a, b, x) \in \mathring{X}_{\Pi, \Pi'} \times B_\epsilon(0)$.

We say a continuous function $f$ on $V \subset \hat V$ is
$(\epsilon, \tilde \delta_0)$-admissible if the composition of $f$ and the natural
projection $\hat V \to V$ is $(\epsilon, \tilde \delta_0)$-admissible.
See Remark \ref{natural projection from hat V to V} for the natural projection.
\end{defi}
Corollary \ref{asymptotic estimates of phi} implies that
$\phi : \hat V \to C^l(\Sigma_0 \setminus N_0,
(\R_{-k_-} \sqcup \dots \sqcup \R_{-1}) \times Y^-
\sqcup X \sqcup
(\R_1 \sqcup \dots \sqcup \R_{k_+}) \times Y^+) \times E^0$ is
$(\epsilon, \tilde \delta_0)$-admissible for any $0 < \epsilon < 1$ and
$0 < \tilde \delta_{0, i} < \min\{ \kappa_\mu / L_\mu; \mu \in M_i\}$.

For each $\nu$ such that the $\nu$-th nodal point of $(\Sigma_0, z, u_0)$
is contained in the $0$-th floor,
there exists an $(\epsilon, \tilde \delta_0)$-admissible
function $\alpha_\nu : \hat V \to \R_{> 0}$ such that
$\alpha_\nu(a, b, x) = \alpha(u_{a, b, x}(q_\nu))$ for any $(a, b, x) \in \hat V$ such that
$\rho_\nu = 0$, where $q_\nu$ is the $\nu$-th nodal point of $\widetilde{P}_a$
and $u_{a, b, x} = \Phi_{a, b}(\xi_{a, b, x})$ is the map for $(a, b, x)$.
For example, the composition of the projection $\hat V \to
\{(a, b, x) \in \hat V; \rho_\nu = 0\}$ and the map $\alpha_\nu(a, b, x)
= \alpha(u_{a, b, x}(q_\nu))$ on $\{(a, b, x) \in \hat V; \rho_\nu = 0\}$
satisfies this condition since
the map $u_{a, b, x}(q_\nu) : \{(a, b, x) \in \hat V; \rho_\nu = 0\} \to X$
is $(\epsilon, \tilde \delta_0)$-admissible as well as $\phi$.
We fix such an $(\epsilon, \tilde \delta_0)$-admissible
function $\alpha_\nu : \hat V \to \R_{> 0}$.

For each $\nu$ such that the $\nu$-th nodal point of $(\Sigma_0, z, u_0)$ is
contained in a negative floor or a positive floor, we define $\alpha_\nu$
by $\alpha_\nu = \alpha^-$ or $\alpha_\nu = \alpha^+$ respectively.
For each $\mu$, we define $\beta_\mu$ by $\beta_\mu = L_\mu^{-1} \beta^-$
if $\mu$ is a joint circle between non-positive floors and
$\beta_\mu = L_\mu^{-1} \beta^+$ if $\mu$ is a joint circle between non-negative
floors.
We define a smooth structure of $\hat V  = \mathring{X} \times B_\epsilon(0)$
by the coordinate
{\belowdisplayskip= 0pt
\[
\hat V \subset \J_0 \times D^{l_0} \times \widetilde{D}^{l_1} \times B_\epsilon(0)
\to \J_0 \times D^{l_0} \times ([0, 1] \times S^1)^{l_1} \times B_\epsilon(0)
\]}
{\abovedisplayskip= 0pt
\begin{multline}
(j, (\zeta_\nu = \rho_\nu^2 e^{2\sqrt{-1} \varphi_\nu})_\nu,
(\zeta_\mu = \rho_\mu^{2\pi} e^{2\pi \sqrt{-1} \varphi_\mu})_\mu, x)\\
\mapsto (\hat \jmath, (\hat \zeta_\nu = \hat \rho_\nu^2
e^{2\sqrt{-1} \hat \varphi_\nu})_\nu, (\hat \rho_\mu, \hat \varphi_\mu)_\mu, \hat x) \label{admissible coordinate change}
\end{multline}
}
given by $\rho_\nu = \hat \rho_\nu^{\alpha_\nu}$,
$\rho_\mu = \hat \rho_\mu^{\beta_\mu}$
and $(\hat \jmath, \hat \varphi_\nu, \hat \varphi_\mu, \hat x)
= (j, \varphi_\nu, \varphi_\mu, x)$.

First we prove the smoothness of the map
\[
\phi : \hat V \to
C^l(\Sigma_0 \setminus N_0, (\R_{-k_-} \sqcup \dots \sqcup \R_{-1}) \times Y^-
\sqcup X \sqcup (\R_1 \sqcup \dots \sqcup \R_{k_+}) \times Y^+) \times E^0,
\]
which implies the smoothness of the evaluation maps at the marked points.
This follows from the following lemma.
\begin{lem}\label{admissibility implies smoothness}
For any $(\epsilon, \tilde \delta_0)$-admissible function $f$ on $\hat V$,
\begin{align*}
&\Bigl|\partial_{\hat x}^{k_x} \partial_{\hat \jmath}^{k_j} \partial_{\hat b}^{k_b}
\partial_{(\hat \rho_{\mu_i})}^{(k_{\mu_i})}
\partial_{(\hat \varphi_\mu)}^{(l_\mu)} \partial_{(\hat \rho_\nu)}^{(k_\nu)}
\Bigl(\prod_\nu \frac{1}{\hat \rho_\nu^{l_\nu}}\Bigr)\partial_{(\hat \varphi_\nu)}^{(l_\nu)}
f \Bigr|\\
&\lesssim \prod_{\substack{i \\ k_{\mu_i} \neq 0}}
(\hat \rho_{\mu_i})^{\beta \tilde \delta_{0, i} / 2 - k_{\mu_i}}
\prod_{\substack{\nu \\ (k_\nu, l_\nu) \neq (0,0)}}
(\hat \rho_\nu)^{\epsilon \alpha - (k_\nu + l_\nu)}
(- \log \hat \rho_\nu)^N,
\end{align*}
where $N = |k_x| + |k_j| + |k_b| + |(k_{\mu_i})| + |(l_\mu)| + |(k_\nu)| + |(l_\nu)|$.
\end{lem}
\begin{proof}
It is easy to check that the claim follows from the following estimates of
the differentials of the coordinate change and
the $(\epsilon, \tilde \delta_0)$-admissibility of $f$ and $\alpha$:
\begin{align}
&\Bigl|\partial_x^{k_x} \partial_j^{k_j} \partial_b^{k_b} \partial_{(\rho_{\mu_i})}^{(k_{\mu_i})}
\partial_{(\varphi_\mu)}^{(l_\mu)} \partial_{(\rho_\nu)}^{(k_\nu)}
\partial_{(\varphi_\nu)}^{(l_\nu)}
\Bigl(\frac{\partial \rho_{\nu_1}}{\partial \hat \rho_{\nu_0}}
- \delta^{\nu_0, \nu_1} \alpha_{\nu_1} \hat \rho_{\nu_1}^{\alpha_{\nu_1} - 1}
\Bigr)\Bigr| \notag\\
&\lesssim \rho_{\nu_1} (- \log \hat \rho_{\nu_1}) \prod_{\substack{i \\ k_{\mu_i} \neq 0}}
\rho_{\mu_i}^{L_{\mu_i} \tilde \delta_{0, i} / 2 - k_{\mu_i}}
\prod_{\substack{\nu \\ (k_\nu, l_\nu) \neq (0,0) \\\text{or } \nu = \nu_0}}
\rho_\nu^{\epsilon - k_\nu - \delta^{\nu, \nu_0} \alpha_{\nu_0}^{-1}},
\label{nu nu}
\end{align}
\begin{align}
&\Bigl|\partial_x^{k_x} \partial_j^{k_j} \partial_b^{k_b} \partial_{(\rho_{\mu_i})}^{(k_{\mu_i})}
\partial_{(\varphi_\mu)}^{(l_\mu)} \partial_{(\rho_\nu)}^{(k_\nu)}
\partial_{(\varphi_\nu)}^{(l_\nu)}
\frac{\partial \rho_{\nu_1}}{\partial \hat \rho_{\mu_{i_0}}}\Bigr| \notag\\
&\lesssim \rho_{\nu_1} (- \log \hat \rho_{\nu_1})
\prod_{\substack{i \\ k_{\mu_i} \neq 0\\ \text{or } i = i_0}}
\rho_{\mu_i}^{L_{\mu_i} \tilde \delta_{0, i} / 2 - k_{\mu_i}
- \delta^{i, i_0} \beta_{\mu_{i_0}}^{-1}}
\prod_{\substack{\nu \\ (k_\nu, l_\nu) \neq (0,0)}}
\rho_\nu^{\epsilon - k_\nu},
\label{nu mu}
\end{align}
and
\begin{align}
&|\partial_x^{k_x} \partial_j^{k_j} \partial_b^{k_b} \partial_{(\rho_{\mu_i})}^{(k_{\mu_i})}
\partial_{(\varphi_\mu)}^{(l_\mu)} \partial_{(\rho_\nu)}^{(k_\nu)}
\partial_{(\varphi_\nu)}^{(l_\nu)}g^{\nu_1}| \notag\\
&\lesssim \rho_{\nu_1} (- \log \hat \rho_{\nu_1})
\prod_{\substack{i \\ k_{\mu_i} \neq 0}}
\rho_{\mu_i}^{L_{\mu_i} \tilde \delta_{0, i} / 2 - k_{\mu_i}}
\prod_{\substack{\nu \\ (k_\nu, l_\nu) \neq (0,0)}}
\rho_\nu^{\epsilon - k_\nu}
\label{nu the others}
\end{align}
for
\[
g^{\nu_1} = \frac{\partial \rho_{\nu_1}}{\partial \hat \varphi_{\nu}},
\frac{\partial \rho_{\nu_1}}{\partial \hat \varphi_{\mu}},
\frac{\partial \rho_{\nu_1}}{\partial \hat x},
\frac{\partial \rho_{\nu_1}}{\partial \hat \jmath},
\frac{\partial \rho_{\nu_1}}{\partial \hat b_{\mu}},
\]
where $\delta^{\nu, \nu'}$ and $\delta^{i, i'}$ are the Kronecker deltas.
We sketch the proof of (\ref{nu nu}), (\ref{nu mu}) and (\ref{nu the others}).

Let $A$ be a square-matrix-valued function on
\[
\{(a, b, x) \in \hat V; \rho_\nu \neq 0, \rho_\mu \neq 0 \text{ for all }
\nu \text{ and } \mu\} \times B_\epsilon(0)
\]
defined by
\begin{align*}
&\lsuperscript{(\rho_\nu (-\log \hat \rho_\nu) \partial_{\rho_\nu},\ 
\partial_{\varphi_\nu},\ 
\rho_{\mu_i} \partial_{\rho_{\mu_i}},\ 
\partial_{\varphi_\mu},\ 
\partial_x,\  \partial_j,\  \partial_b)}{t}\\
&= A \cdot
\lsuperscript{
(\alpha_\nu^{-1} \hat \rho_\nu (-\log \hat \rho_\nu) \partial_{\hat \rho_\nu},\ 
\partial_{\hat \varphi_\nu},\ 
\beta_{\mu_i}^{-1} \hat \rho_{\mu_i} \partial_{\hat \rho_{\mu_i}},\ 
\partial_{\hat \varphi_\mu},\ 
\partial_{\hat x},\  \partial_{\hat \jmath},\  \partial_{\hat b})}{t}
\end{align*}
We can easily check the following estimates of the columns of $(A-1)$
corresponding to the vectors $\rho_\nu (-\log \hat \rho_\nu) \partial_{\rho_\nu}$.
It is also easy to check that the other columns of $(A-1)$ are zero.
In the inequalities below,
$(A-1)_{\rho_{\nu_0}, \rho_{\nu_1}}$ is the entry corresponding to
$\rho_{\nu_0} (-\log \hat \rho_{\nu_0}) \partial_{\rho_{\nu_0}}$
and $\rho_{\nu_1} (-\log \hat \rho_{\nu_1}) \partial_{\rho_{\nu_1}}$.
The other entries $(A-1)_{\rho_{\nu_{i_0}}, \rho_{\nu_1}}$ are similar.
In (\ref{A the others}), $\ast$ denotes the other rows:
$\ast = \varphi_{\nu}, \varphi_{\mu}, x, j, b$.
\begin{align}
&\bigl|
\partial_x^{k_x} \partial_j^{k_j} \partial_b^{k_b} \partial_{(\rho_{\mu_i})}^{(k_{\mu_i})}
\partial_{(\varphi_\mu)}^{(l_\mu)} \partial_{(\rho_\nu)}^{(k_\nu)}
\partial_{(\varphi_\nu)}^{(l_\nu)} (A-1)_{\rho_{\nu_0}, \rho_{\nu_1}}
\bigr| \notag\\
&\lesssim (-\log \hat \rho_{\nu_0}) \prod_{\substack{i \\ k_{\mu_i} \neq 0}}
\rho_{\mu_i}^{L_{\mu_i} \tilde \delta_{0, i} / 2 - k_{\mu_i}}
\prod_{\substack{\nu \\ (k_\nu, l_\nu) \neq (0,0) \\\text{or } \nu = \nu_0}}
\rho_\nu^{\epsilon - k_\nu},
\label{A nu}
\end{align}
\begin{align}
&\bigl|
\partial_x^{k_x} \partial_j^{k_j} \partial_b^{k_b} \partial_{(\rho_{\mu_i})}^{(k_{\mu_i})}
\partial_{(\varphi_\mu)}^{(l_\mu)} \partial_{(\rho_\nu)}^{(k_\nu)}
\partial_{(\varphi_\nu)}^{(l_\nu)} (A-1)_{\rho_{\mu_{i_0}}, \rho_{\nu_1}}
\bigr| \notag\\
&\lesssim \prod_{\substack{i \\ k_{\mu_i} \neq 0 \\\text{or } i = i_0}}
\rho_{\mu_i}^{L_{\mu_i} \tilde \delta_{0, i} / 2 - k_{\mu_i}}
\prod_{\substack{\nu \\ (k_\nu, l_\nu) \neq (0,0)}}
\rho_\nu^{\epsilon - k_\nu},
\label{A mu}
\end{align}
and
\begin{align}
&\bigl|
\partial_x^{k_x} \partial_j^{k_j} \partial_b^{k_b} \partial_{(\rho_{\mu_i})}^{(k_{\mu_i})}
\partial_{(\varphi_\mu)}^{(l_\mu)} \partial_{(\rho_\nu)}^{(k_\nu)}
\partial_{(\varphi_\nu)}^{(l_\nu)} (A-1)_{\ast, \rho_{\nu_1}}
\bigr| \notag\\
&\lesssim \prod_{\substack{i \\ k_{\mu_i} \neq 0}}
\rho_{\mu_i}^{L_{\mu_i} \tilde \delta_{0, i} / 2 - k_{\mu_i}}
\prod_{\substack{\nu \\ (k_\nu, l_\nu) \neq (0,0)}}
\rho_\nu^{\epsilon - k_\nu}.
\label{A the others}
\end{align}
These estimates follow from the $(\epsilon, \tilde \delta_0)$-admissibility of
$\alpha_{\nu}$ and the following equations:
\[
(A-1)_{\rho_{\nu_0}, \rho_{\nu_1}}
= (-\log \hat \rho_{\nu_0}) \rho_{\nu_0}
\frac{\partial \alpha_{\nu_1}}{\partial \rho_{\nu_0}},
\]
\[
(A-1)_{\rho_{\mu_{i_0}}, \rho_{\nu_1}}
= \rho_{\mu_{i_0}} \frac{\partial \alpha_{\nu_1}}{\partial \rho_{\mu_{i_0}}}
\]
and
\[
(A-1)_{\ast, \rho_{\nu_1}}
= \partial_\ast \alpha_{\nu_1},
\]
where $\partial_\ast = \partial_{\varphi_{\nu}}, \partial_{\varphi_{\mu}},
\partial_x, \partial_j, \partial_b$.

(\ref{A nu}), (\ref{A mu}) and (\ref{A the others}) imply that
the same inequalities hold for $A^{-1} - 1$ as well as $A-1$.
This is because the derivatives of $A^{-1}$ are polynomials of
$A^{-1} = 1 + (1-A) + (1-A)^2 + \cdots$ and the derivatives of $(A-1)$.
These inequalities are equivalent to (\ref{nu nu}), (\ref{nu mu}) and
(\ref{nu the others}).
\end{proof}

Next we prove the smoothness of the embedding between two Kuranishi neighborhoods.
For this proof, we directly use the admissibility of $\phi$
rather than its smoothness.

The definition of the embedding itself is the same as the case of
$\widehat{\M}(Y, \lambda, J)$.
We assume the similar condition to Section \ref{embed} and use the same notation.
Let $(a^1, b^1, u^1, h^1) \to (a^2, b^2, u^2, h^2)$ be the embedding of 
$V^0_1 \subset V_1$ into $V_2$.
Let $N_{q_0} \subset \{\nu^1\}$ be the set of indices of nodal points of $\Sigma_1$
which remain to be nodal points in $\Sigma_0$, that is, $\rho_{\nu^1} = 0$ at $a_0^1$.
For each $\nu^1 \in N_{q_0}$, let $\iota(\nu^1)$ be the index of the corresponding
nodal point of $\Sigma_2$.
Similarly, let $M_{q_0} \subset \{\mu^1\}$ be the set of indices of joint circles of
$\Sigma_1$ which remain to be joint circles in $\Sigma_0$,
and let $\iota(\mu^1)$ be the index of the corresponding joint circle of $\Sigma_2$
for each $\mu^1 \in M_{q_0}$.
We assume that the maps $\phi$ and $\alpha_\nu$ for the Kuranishi neighborhood
$(V_1, E_1, s_1, \psi_1, G_1)$ are $(\epsilon, \tilde \delta_0)$-admissible,
those for $(V_2, E_2, s_2, \psi_2, G_2)$ are $(\epsilon, \tilde \delta'_0)$-admissible,
and $\tilde \delta_{0, i} \leq \tilde \delta'_{0, j}$ if
the joint circles of $\Sigma_1$ which belong to $M^1_i$
remain to be joint circles in $\Sigma_0$ and they correspond to those which
belong to $M^2_j$.
First we check the following:
\begin{itemize}
\item[($\clubsuit$)]
$\zeta^2_{\iota(\nu^1)} / \zeta^1_{\nu^1}$ for $\nu^1 \in N_{q_0}$
and $\rho^2_{\iota(\mu^1)} / \rho^1_{\mu^1}$ for $\mu^1 \in M_{q_0}$
are $(\epsilon, \tilde \delta_0)$-admissible and bounded away from zero
(i.e. the continuous extensions do not take zero on $V^0_1$).
$\zeta^2_{\nu^2}$ for $\nu^2 \notin \iota(N_{q_0})$,
$\rho^2_{\mu^2}$ for $\mu^2 \notin \iota(M_{q_0})$ and
$(\varphi^2_\mu, j^2, b^2_\mu, u^2, h^2)$ are $(\epsilon, \tilde \delta_0)$-admissible.
\end{itemize}

The $(\epsilon, \tilde \delta_0)$-admissibility of $h^2$ is clear.
$\mathcal{Z}_2^+ \in \Sigma_1 \setminus N_1$ is $(\epsilon, \tilde \delta_0)$-admissible
because of the $(\epsilon, \tilde \delta_0)$-admissibility of $u^1 \in
C^{l_1}(\Sigma_1 \setminus N_0, (\R_{-k_-} \sqcup \dots \sqcup \R_{-1}) \times Y^-
\sqcup X \sqcup (\R_1 \sqcup \dots \sqcup \R_{k_+}) \times Y^+)$.
Hence $\hat a^2 \in \hat U_2$ is also $(\epsilon, \tilde \delta_0)$-admissible.
Therefore in the definition of $\theta$,
$\sigma \circ u^1 \circ (\pi_1|_{(\widetilde{P}_1)_{a^1}})^{-1} \circ
\Theta|_{(\hat P_2)|_{\hat a^2}} (\hat R^2_i(\hat a^2)) \in \R_{i'}$ is
$(\epsilon, \tilde \delta_0)$-admissible.
($i'$ is the floor of $\Sigma_1$ which corresponds to the $i$-th floor of $\Sigma_2$
in $\Sigma_0$.)
Together with the $(\epsilon, \tilde \delta_0)$-admissibility of $u^1$,
it implies that $\mathcal{Z}_2^{++} \in \Sigma_1 \setminus N_1$ is also
$(\epsilon, \tilde \delta_0)$-admissible.
The $(\epsilon, \tilde \delta_0)$-admissibility of $\mathcal{Z}_2^+$ and
$\mathcal{Z}_2^{++}$ implies the $(\epsilon, \tilde \delta_0)$-admissibility of
$a^2 \in \widetilde{U}_2$.
Furthermore, it implies that
$\zeta^2_{\iota(\nu^1)} / \zeta^1_{\nu^1}$ for $\nu^1 \in N_{q_0}$
and $\rho^2_{\iota(\mu^1)} / \rho^1_{\mu^1}$ for $\mu^1 \in M_{q_0}$
are $(\epsilon, \tilde \delta_0)$-admissible and bounded away from zero.
The $(\epsilon, \tilde \delta_0)$-admissibility of $a^2$ and $u^1$ implies
$u^2 \in
C^{l_2}(\Sigma_2 \setminus N_2, (\R_{-k_-} \sqcup \dots \sqcup \R_{-1}) \times Y^-
\sqcup X \sqcup (\R_1 \sqcup \dots \sqcup \R_{k_+}) \times Y^+)$ is
$(\epsilon, \tilde \delta_0)$-admissible.
Finally, $b^2_\mu$ are also $(\epsilon, \tilde \delta_0)$-admissible because
the function $f_\mu$ in (\ref{b for any kappa}) is $(\epsilon, \tilde \delta_0)$-admissible.

It is easy to check that ($\clubsuit$) and the $(\epsilon, \tilde \delta'_0)$-admissibility
of $\alpha^2_{\nu^2}$ on $V_2$ imply that $\alpha^2_{\nu^2}$ are
$(\epsilon, \tilde \delta_0)$-admissible as functions on $V^0_1$.
Note that for any $\nu^1 \in N_{q_0}$, $\alpha^2_{\iota(\nu^1)} = \alpha^1_{\nu^1}$
on $\{(a^1, b^1, u^1, h^1) \in V^0_1; \rho^1_{\nu^1} = 0\}$ by definition.
Therefore their $(\epsilon, \tilde \delta'_0)$-admissibility implies that
for any multi-index $(k_{x^1}, k_{j^1}, k_{b^1}, k_{\mu^1_i}, \ab l_{\mu^1},
\ab k_{\nu^1}, \ab l_{\nu^1})$
such that $(k_{\nu^1_0}, l_{\nu^1_0}) = (0, 0)$,
\begin{align}
&\bigl|\partial_{x^1}^{k_{x^1}} \partial_{j^1}^{k_{j^1}} \partial_{b^1}^{k_{b^1}}
\partial_{(\rho_{\mu^1_i})}^{(k_{\mu^1_i})}
\partial_{(\varphi_{\mu^1})}^{(l_{\mu^1})} \partial_{(\rho_{\nu^1})}^{(k_{\nu^1})}
\partial_{(\varphi_{\nu^1})}^{(l_{\nu^1})}
(\alpha^2_{\iota(\nu^1_0)} - \alpha^1_{\nu^1_0})\bigr| \notag\\
&\lesssim (\rho^1_{\nu^1_0})^\epsilon \cdot \prod_{\substack{i \\ k_{\mu_i} \neq 0}}
\rho_{\mu_i}^{L_{\mu_i} \tilde \delta_{0, i} / 2 - k_{\mu_i}}
\prod_{\substack{\nu \\ (k_\nu, l_\nu) \neq (0,0)}}
\rho_\nu^{\epsilon - k_\nu}.
\label{difference of alpha}
\end{align}

Now we prove the smoothness of the embedding.
For any $\nu^1 \in N_{q_0}$,
($\clubsuit$), (\ref{difference of alpha}) and the $(\epsilon, \tilde \delta_0)$-admissibility
of $\alpha^1_{\nu^1}$ and $\alpha^2_{\iota(\nu^1)}$ imply that
\[
\hat \rho^2_{\iota(\nu^1)} / \hat \rho^1_{\nu^1}
= (\rho^1_{\nu^1})^{(\alpha^2_{\iota(\nu^1)})^{-1} - (\alpha^1_{\nu^1})^{-1}}
\cdot \bigl(\rho^2_{\iota(\nu^1)} / \rho^1_{\nu^1}\bigr)^{(\alpha^2_{\iota(\nu^1)})^{-1}}
\]
is $(\epsilon, \tilde \delta_0)$-admissible.
Assume that $\mu^1_i \in M^1_i$ and that
$\mu^2_\ast = \iota(\mu^1_i) \in M^2_{i'}$.
Then ($\clubsuit$) implies
\[
\hat \rho^2_{\mu^2_{i'}} / \hat \rho^1_{\mu^1_i}
= (\rho^2_{\iota(\mu^1_i)} / \rho^1_{\mu^1_i})^{(\beta_{\mu^1_i})^{-1}}
\cdot e^{(b^2_{\mu^2_{i'}} - b^2_{\mu_\ast}) \beta}
\]
is also $(\epsilon, \tilde \delta_0)$-admissible.
Therefore, $(a^2, b^2, u^2, h^2) \in V_2$ is an $(\epsilon, \tilde \delta_0)$-admissible
function of $(a^1, b^1, u^1, h^1) \in V^0_1$ if the differential structure of $V_2$
is defined by $(\alpha^2_{\nu^2}, \beta^\pm)$.
Hence Lemma \ref{admissibility implies smoothness} implies the smoothness of
the embedding.

We can similarly prove the smoothness of the essential submersion
from a Kuranishi neighborhood of a disconnected holomorphic building
to the product of those of its connected components.

\subsection{Fiber products and multisections}\label{fiber prod for X}
Let $K_{Y^\pm} \inj \overline{P}_{Y^\pm}$ be triangulations, and let $K^0_{Y^\pm}$ be
finite sets of smooth cycles in $Y^\pm$.
Assume that a finite sequence $K^0_X = (x)$ of smooth cycles
with closed supports in $X$ is given which satisfies the following conditions.
For each cycle $x \in K^0_X$, $\supp x \cap (-\infty, 0] \times Y^-$ is empty set or
there exists some cycle $y \in K^0_{Y^-}$ such that
$x |_{(-\infty, 0] \times Y^-} = (-\infty, 0] \times y$.
Similarly, for each cycle $x \in K^0_X$, $\supp x \cap [0, \infty) \times Y^+$ is
empty set or there exists some cycle $y \in K^0_{Y^+}$ such that
$x |_{[0, \infty) \times Y^+} = [0, \infty) \times y$.
Further we assume that these relations give bijections
\begin{align*}
\mu_- &: \{x \in K^0_X; \supp x \cap (-\infty, 0] \times Y^- \neq \emptyset\}
\to K^0_{Y^-} \text{ and}\\
\mu_+ &: \{x \in K^0_X; \supp x \cap [0, \infty) \times Y^+ \neq \emptyset\}
\to K^0_{Y^+}.
\end{align*}

First we explain about the construction of the multisections of the fiber products.
Assume that the multisection of
$(\widehat{\M}_{Y^\pm}^\diamond, \mathring{K}_{Y^\pm}^2, K_{Y^\pm}, K_{Y^\pm}^0)$
is given for each $Y^\pm$.
We define a space
$\widehat{\M}_X^\diamond$
as follows.
Its point
\[
((\Sigma^\alpha, z^\alpha, u^\alpha)_{\alpha \in A^- \sqcup A^0 \sqcup A^+},
M^{\mathrm{rel}})
\]
consists of holomorphic buildings
$(\Sigma^\alpha, z^\alpha, u^\alpha)_{\alpha \in A^-}$ for $Y^-$,
$(\Sigma^\alpha, z^\alpha, u^\alpha)_{\alpha \in A^0}$ for $X$,
$(\Sigma^\alpha, z^\alpha, u^\alpha)_{\alpha \in A^+}$ for $Y^+$, and
a set $M^{\mathrm{rel}} = \{(S^1_{+\infty_l}, S^1_{-\infty_l})\}$ of pairs of
limit circles which satisfy the following conditions:
\begin{itemize}
\item
Any two pairs in $M^{\mathrm{rel}}$ do not share the same limit circle.
\item
For each pair $\alpha_1, \alpha_2 \in A = A^- \sqcup A^0 \sqcup A^+$,
let $M^{\alpha_1, \alpha_2} \subset M^{\mathrm{rel}}$ be the subset of pairs
$(S^1_{+\infty_l}, S^1_{-\infty_l})$ such that $S^1_{+\infty_l}$ is a $+\infty$-limit circle
of $\Sigma^{\alpha_1}$ and $S^1_{-\infty_l}$ is a $-\infty$-limit circle of
$\Sigma^{\alpha_2}$.
Then there does not exists any sequence
$\alpha_0, \alpha_1, \dots, \alpha_k = \alpha_0 \in A$ such that
$M^{\alpha_i, \alpha_{i+1}} \neq \emptyset$ for all $i = 0, 1, \dots, k-1$.
\item
For subsets $A_1, A_2 \subset A$, define
$M^{A_1, A_2} = \bigcup_{\alpha_1 \in A_1, \alpha_2 \in A_2} M^{\alpha_1, \alpha_2}$.
Then $M^{\mathrm{rel}}$ is the union of
$M^{\mathrm{rel}, \leq 0} = M^{A^-, A^- \sqcup A^0}$
and $M^{\mathrm{rel}, \geq 0} = M^{A^0 \sqcup A^+, A^+}$.
\end{itemize}

We regard $\widehat{\M}^\diamond_{Y^-}$ and $\widehat{\M}^\diamond_{Y^+}$
as subspaces of $\widehat{\M}^\diamond_X$ consisting of points
such that $A^0 = A^+ = \emptyset$ and $A^- = A^0 = \emptyset$ respectively.

We say a point
$((\Sigma^\alpha, z^\alpha, u^\alpha)_{\alpha \in A^- \sqcup A^0 \sqcup A^+},
M^{\mathrm{rel}}) \in \widehat{\M}^\diamond_X$ is disconnected if
there exists a decomposition $A^- \sqcup A^0 \sqcup A^+ = A_1 \sqcup A_2$ such that
$M^{A_1, A_2} = M^{A_2, A_1} = \emptyset$.
Otherwise we say it is connected.
We denote the space of connected points of $\widehat{\M}^\diamond_X$ by
$(\widehat{\M}^\diamond_X)^0$.
Decomposition into connected components defines a map $\widehat{\M}^\diamond_X
\to \bigcup_N (\prod^N (\widehat{\M}^\diamond_X)^0) / \mathfrak{S}_N$.

Let
\begin{align*}
\Upsilon : \widehat{\M}^\diamond_X &\to
\prod (\overline{P}_{Y^-} \times \overline{P}_{Y^-}) / \mathfrak{S}
\times \prod (\overline{P}_{Y^+} \times \overline{P}_{Y^+}) / \mathfrak{S} \\
&\quad \quad
\times \prod \overline{P}_{Y^-} / \mathfrak{S}
\times \prod \overline{P}_{Y^+} / \mathfrak{S}
\times \prod (Y^- \cup X \cup Y^+) / \mathfrak{S}.
\end{align*}
be the continuous map which maps a point
$((\Sigma^\alpha, z^\alpha, u^\alpha)_{\alpha \in A^- \sqcup A^0 \sqcup A^+},
M^{\mathrm{rel}})$ to the following point.
Its $\prod (\overline{P}_{Y^-} \times \overline{P}_{Y^-}) / \mathfrak{S}$-factor is
\[
(\pi_{Y^-} \circ u|_{S^1_{+\infty_l}}, \pi_{Y^-} \circ u|_{S^1_{+\infty_l}})
_{(S^1_{+\infty_l}, S^1_{-\infty_l}) \in M^{\mathrm{rel}, \leq 0}},
\]
its $\prod (\overline{P}_{Y^+} \times \overline{P}_{Y^+}) / \mathfrak{S}$-factor is
\[
(\pi_{Y^+} \circ u|_{S^1_{+\infty_l}}, \pi_{Y^+} \circ u|_{S^1_{+\infty_l}})
_{(S^1_{+\infty_l}, S^1_{-\infty_l}) \in M^{\mathrm{rel}, \geq 0}},
\]
its $\prod \overline{P}_{Y^-} / \mathfrak{S}
\times \prod \overline{P}_{Y^+} / \mathfrak{S}$-factor is
$(\pi_{Y^\pm} \circ u|_{S^1_{\pm\infty}})_{S^1_{\pm\infty} \notin M^{\mathrm{rel}}}$,
and its $\prod (Y^- \cup X \cup Y^+) / \mathfrak{S}$-factor is
$(u(z))_{z \in \bigcup_\alpha z^\alpha}$, where if $z \in z^\alpha$ is contained in
a negative (or positive) floor of $(\Sigma^\alpha, z^\alpha, u^\alpha)$,
then we read $u(z)$ as $\pi_{Y^-} \circ u(z)$
(or $\pi_{Y^-} \circ u(z)$ respectively).
$\Upsilon$ is realized as a strong smooth map.
Define the fiber product
\[
(\widehat{\M}^\diamond_X, (\mathring{K}_{Y^-}^2, \mathring{K}_{Y^+}^2),
(K_{Y^-}, K_{Y^+}), (K^0_{Y^-}, K^0_X, K^0_{Y^+})) \subset
\widehat{\M}^\diamond_X
\]
by
\begin{align*}
&(\widehat{\M}^\diamond_X, (\mathring{K}_{Y^-}^2, \mathring{K}_{Y^+}^2),
(K_{Y^-}, K_{Y^+}), (K^0_{Y^-}, K^0_X, K^0_{Y^+})) \\
&= \Upsilon^{-1}\Bigl(\prod \mathring{K}_{Y^-}^2 / \mathfrak{S} \times
\prod \mathring{K}_{Y^+}^2 / \mathfrak{S} \times
\prod K_{Y^-} / \mathfrak{S} \times \prod K_{Y^+} / \mathfrak{S} \\
&\hph{= \Upsilon^{-1}\Bigl(}
\times \prod (K^0_{Y^-} \cup K^0_X \cup K^0_{Y^+}) / \mathfrak{S}\Bigr).
\end{align*}
We abbreviate the above fiber product
as $(\widehat{\M}^\diamond_X, \mathring{K}^2_X, K_X, K^0_X)$.
(This is simply an abbreviation, and $\mathring{K}^2_X$ or $K_X$
are not defined.)

Define a multi-valued partial submersion
$\Xi : \widehat{\M}^\diamond_X \to \widehat{\M}^\diamond_X$ by
\begin{align*}
&\Xi(((\Sigma^\alpha, z^\alpha, u^\alpha)_{\alpha \in A^- \sqcup A^0 \sqcup A^+},
M^{\mathrm{rel}})) \\
&= \{((\Sigma^\alpha, z^\alpha, u^\alpha)_{\alpha \in A^- \sqcup A^0 \sqcup A^+},
\mathring{M}^{\mathrm{rel}}); \mathring{M}^{\mathrm{rel}} \subsetneq M^{\mathrm{rel}}\},
\end{align*}
and let $\Xi^\circ$ be the restriction of
$\Xi \subset \widehat{\M}^\diamond_X \times \widehat{\M}^\diamond_X$
to the set of points
\[
(((\Sigma^\alpha, z^\alpha, u^\alpha)_{\alpha \in A^- \sqcup A^0 \sqcup A^+},
M^{\mathrm{rel}}),
((\Sigma^\alpha, z^\alpha, u^\alpha)_{\alpha \in A^- \sqcup A^0 \sqcup A^+},
\mathring{M}^{\mathrm{rel}}))
\]
such that
$((\Sigma^\alpha, z^\alpha, u^\alpha)_{\alpha \in A^- \sqcup A^0 \sqcup A^+},
M^{\mathrm{rel}})
\in (\widehat{\M}^\diamond_X, \mathring{K}^2_X, K_X, K^0_X)$
and $(\pi_Y \circ u|_{S^1_{+\infty_l}}, \pi_Y \circ u|_{S^1_{-\infty_l}})
\in \rho_\ast K_{Y^\pm}$
for all $(S^1_{+\infty_l}, S^1_{-\infty_l}) \in M^{\mathrm{rel}} \setminus
\mathring{M}^{\mathrm{rel}}$.
Then $\Xi^\circ$ is a multi-valued partial essential submersion
from $(\widehat{\M}^\diamond_X, \mathring{K}^2_X, K_X, K^0_X)$ to itself.

We also define a multi-valued partial essential submersion
\[
\Lambda : (\partial \widehat{\M}^\diamond_X, \mathring{K}^2_X, K_X, K^0_X)
\to (\widehat{\M}^\diamond_X, \mathring{K}^2_X, K_X, K^0_X)
\]
similarly to the case of $Y$.

Similarly to the case of symplectization, for each point
\[
p = ((\Sigma^\alpha, z^\alpha, u^\alpha)
_{\alpha \in A^- \sqcup A^0 \sqcup A^+}, M^{\mathrm{rel}})
\in (\widehat{\M}^\diamond_X, \mathring{K}^2_X, K_X, K^0_X),
\]
we define
$\widetilde{e}(p) = \widetilde{e}_{\delta_0}(p)
= \sum_\alpha \widehat{e}_{\delta_0}(\theta_\alpha)
+ \frac{1}{2} \# M^{\mathrm{rel}}$, where each
$\theta_\alpha$ is the type of $(\Sigma^\alpha, z^\alpha, u^\alpha)$.
Then the maps $\Xi^\circ$ and $\Lambda$ decrease $\widetilde{e}$.
Hence if we decompose $(\widehat{\M}^\diamond_X, \mathring{K}^2_X, K_X, K^0_X)$ by
$\widetilde{e}$, then
$\Xi^\circ$ and $\Lambda$ constitute a compatible system of
multi-valued partial submersions.

We can construct the perturbed multisections of
$(\widehat{\M}^\diamond_X, \mathring{K}^2_X, K_X, K^0_X)$
which satisfy the following conditions:
\begin{itemize}
\item
The restrictions of the grouped multisection of
$(\widehat{\M}^\diamond_X, \mathring{K}^2_X, K_X, K^0_X)$ to
$(\widehat{\M}_{Y^\pm}^\diamond, \mathring{K}_{Y^\pm}^2, K_{Y^\pm}, K_{Y^\pm}^0)$
coincide with the given grouped multisections.
\item
Let $((\widehat{\M}^\diamond_X)^0, \mathring{K}^2_X, K_X, K^0_X) \subset
(\widehat{\M}^\diamond_X, \mathring{K}^2_X, K_X, K^0_X)$
be the subset of connected points.
Its grouped multisection induces that of
\[
\bigcup_N (\prod^N (\widehat{\M}^\diamond_X)^0, \mathring{K}^2_X, K_X, K^0_X)
/ \mathfrak{S}_N.
\]
Then the grouped multisection of
$(\widehat{\M}^\diamond_X, \mathring{K}^2_X, K_X, K^0_X)$
coincides with its pull back by the submersion
\[
(\widehat{\M}^\diamond_X, \mathring{K}^2_X, K_X, K^0_X) \to
\bigcup_N (\prod^N (\widehat{\M}^\diamond_X)^0, \mathring{K}^2_X, K_X, K^0_X)
/ \mathfrak{S}_N
\]
defined by decomposition into connected components.
\item
The grouped multisections of
$(\widehat{\M}^\diamond_X, \mathring{K}^2_X, K_X, K^0_X)$
are compatible with respect to the compatible system of
multi-valued partial essential submersions defined by $\Xi^\circ$ and $\Lambda$.
\end{itemize}


Next we define the fiber products we use for the construction of the algebra.
As in Section \ref{fiber prod and orientation},
let $((\hat \epsilon^{i, j}_l), (\hat c^i_l), (x^i_l), (\eta^i_l))$ be sequences of simplices
with local coefficients such that
\begin{itemize}
\item
$\hat \epsilon^{i, j}_l = \theta^{\lsuperscript{D}{t}}_{\epsilon^{i, j}_l} \epsilon^{i, j}_l
\theta^D_{\epsilon^{i, j}_l}$ ($-m_- \leq i < j \leq 0$)
are products of simplices $\epsilon^{i, j}_l$ in $\mathring{K}^2_{Y^-}$not contained in
$\overline{P}_{Y^-}^{^t\text{bad}} \times \overline{P}_{Y^-} \cup
\overline{P}_{Y^-} \times \overline{P}_{Y^-}^{\text{bad}}$ and orientations
$\theta^{\lsuperscript{D}{t}}_{\epsilon^{i, j}_l}$ of $p_1^\ast S_{Y^-}^{\lsuperscript{D}{t}}$
and $\theta^D_{\epsilon^{i, j}_l}$ of $p_2^\ast S_{Y^-}^{\lsuperscript{D}{t}}$
on $\Int \epsilon^{i, j}_l$,
\item
$\hat \epsilon^{i, j}_l = \theta^{\lsuperscript{D}{t}}_{\epsilon^{i, j}_l} \epsilon^{i, j}_l
\theta^D_{\epsilon^{i, j}_l}$ ($0 \leq i < j \leq m_+$)
are products of simplices $\epsilon^{i, j}_l$ in $\mathring{K}^2_{Y^+}$ not contained in
$\overline{P}_{Y^+}^{^t\text{bad}} \times \overline{P}_{Y^+} \cup
\overline{P}_{Y^+} \times \overline{P}_{Y^+}^{\text{bad}}$ and orientations
$\theta^{\lsuperscript{D}{t}}_{\epsilon^{i, j}_l}$ of $p_1^\ast S_{Y^+}^{\lsuperscript{D}{t}}$
and $\theta^D_{\epsilon^{i, j}_l}$ of $p_2^\ast S_{Y^+}^{\lsuperscript{D}{t}}$
on $\Int \epsilon^{i, j}_l$,
\item
$\hat c^i_l = c^i_l \theta^D_{c^i_l}$ ($-m_- \leq i \leq 0$) are products of simplices
$c^i_l$ in $K_{Y^-}$ not contained in $\overline{P}_{Y^-}^\text{bad}$ and
orientations $\theta^D_{c^i_l}$ of $\S_{Y^-}^D$ on $\Int c^i_l$,
\item
$\hat \eta^i_l = \theta^{\lsuperscript{D}{t}}_{\eta^i_l} \eta^i_l$ are products of simplices
$\eta^i_l$ in $K_{Y^+}$ not contained in $\overline{P}_{Y^+}^{^t\text{bad}}$ and
orientations $\theta^{\lsuperscript{D}{t}}_{\eta^i_l}$ of $\S_{Y^+}^{\lsuperscript{D}{t}}$
on $\Int \eta^i_l$,
\item
$x^i_l$ ($-m_- \leq i < 0$) are cycles in $K^0_{Y^-}$,
\item
$x^0_l$ are cycles in $K^0_X$, and
\item
$x^i_l$ ($0 < i \leq m_+$) are cycles in $K^0_{Y^+}$.
\end{itemize}
Take lifts $\tilde \epsilon^{i, j}_l$, $\tilde c^i_l$ and $\tilde \eta^i_l$, and define
$\breve \epsilon^{i, j}_l = \theta^{\lsuperscript{D}{t}}_{\epsilon^{i, j}_l} \tilde \epsilon^{i, j}_l
\theta^D_{\epsilon^{i, j}_l}$, $\breve c^i_l = \tilde c^i_l \theta^D_{c^i_l}$ and
$\breve \eta^i_l = \theta^{\lsuperscript{D}{t}}_{\eta^i_l} \tilde \eta^i_l$ as in Section
\ref{fiber prod and orientation}.

For such a sequence, the pre-Kuranishi space
$\overline{\M}^{(m_-, X, m_+)}
_{((\breve \epsilon^{i, j}_l), (\breve c^i_l), (x^i_l), (\breve \eta^i_l))}$
is defined similarly.
Its grouped multisection is defined by the pull back of that of
$(\widehat{\M}^\diamond_X, \mathring{K}^2_X, K_X, K^0_X)$ by the natural map
\[
\overline{\M}^{(m_-, X, m_+)}
_{((\breve \epsilon^{i, j}_l), (\breve c^i_l), (x^i_l), (\breve \eta^i_l))} \to
(\widehat{\M}^\diamond_X, \mathring{K}^2_X, K_X, K^0_X).
\]
The definition of its orientation is almost the same with the case of $\hat Y$.
The only difference is that we define the orientation of
$\W^0 = TX^0 \times \mathcal{C}^0 / \R^{k_- + k_+} \oplus
\bigoplus_{z^{++}_{0, \beta}} \R^2$
by
\[
(-1)^{k_+} (\R^{k_- + k_+} \oplus \bigoplus_{z^{++}_{0, \beta}} \R^2) \oplus \W^0
= TX^0 \times \mathcal{C}^0
\]
if the range of the holomorphic building corresponding to the center of the Kuranishi
neighborhood is $(\overline{\R}_{-k_-} \cup \dots \cup \overline{\R}_{-1}) \times Y^-
\cup \overline{X} \cup (\overline{\R}_1 \cup \dots \cup \overline{\R}_{k_+}) \times Y^+$.
It is easy to check that this is well defined
and independent of the choice of the lifts of $\eta^i_j$, $c^i_j$ and $\epsilon^i_j$
under the natural isomorphism.
Hence we may denote the above pre-Kuranishi space by
$\overline{\M}^{(m_-, X, m_+)}_{((\hat \epsilon^{i, j}_l), (\hat c^i_l), (x^i_l), (\hat \eta^i_l))}$.
Similarly to equation (\ref{boundary of MM}),
it is easy to see that for any $((\hat c_l), (x_l), (\hat \eta_l))$ and
$(\hat \epsilon^{i, j}_l)$,
\begin{align}
0 &= \sum_{\star_{m_-, m_+}} (-1)^\ast
[\partial \overline{\M}^{(m_-, X, m_+)}_{((\hat \epsilon^{i, j}_l), (\hat c^i_l), (x^i_l),
(\hat \eta^i_l))}]^0 \notag\\
&= \sum_{\star_{m_-, m_+}} (-1)^{\ast + m_- + m_+}
[\overline{\M}^{(m_-, X, m_+)}_{\partial ((\hat \epsilon^{i, j}_l), (\hat c^i_l), (x^i_l),
(\hat \eta^i_l))}]^0 \notag\\
&\quad + \sum_{\substack{-m_- \leq i_0 \leq 0 \\ \star_{m_- +1, m_+}}}
(-1)^{\ast + m_- + 1 + i_0} [\overline{\M}^{(m_- +1, X, m_+)}
_{((e^{\Delta_\ast [\overline{P}_{Y^-}]})^{i_0 -1, i_0} \cup
(\tau^-_{i_0}\hat \epsilon^{i, j}_l), (\hat c^i_l),
(x^i_l), (\hat \eta^i_l))}]^0 \notag\\
&\quad + \sum_{\substack{0 \leq i_0 \leq m_+ \\ \star_{m_-, m_+ +1}}}
(-1)^{\ast + m_- + i_0} [\overline{\M}^{(m_-, X, m_+ +1)}
_{((e^{\Delta_\ast [\overline{P}_{Y^-}]})^{i_0, i_0 +1} \cup (\tau^+_{i_0} \hat \epsilon^{i, j}_l),
(\hat c^i_l), (x^i_l), (\hat \eta^i_l))}]^0, \label{boundary of MMX}
\end{align}
where the sum $\star_{m_-, m_+}$ is taken over all decompositions
\[
\{\hat c_l\} = \coprod_{-m_- \leq i \leq 0} \{\hat c^i_l\},
\quad \{\hat \eta_l\} = \coprod_{0 \leq i \leq m_+} \{\hat \eta^i_l\}
\]
as sets
and all decompositions
\[
\{x_l\} = \coprod_{-m_- \leq i \leq m_+} \{x^i_l\}
\]
such that $x^i_l \in K_{Y^-}^0$ for $-m_- \leq i <0$ and $x^i_l \in K_{Y^+}^0$
for $0 < i \leq m_+$.
(We identify $x \in K_X^0$ with $\mu_-(x) \in K_{Y^-}^0$ and $\mu_+(x) \in K_{Y^+}^0$
in the above decomposition.)
The sign $\ast$ is the weighted sign of the permutation
\[
\begin{pmatrix}
(\hat c^1_l)_l & \cdots & (\hat c^m_l)_l
& (x^1_l)_l & \cdots & (x^m_l)_l & (\hat \eta^1_l)_l & \cdots & (\hat \eta^m_l)_l\\
&(\hat c_l)_l&&&(x_l)_l&&&(\hat \eta_l)_l&
\end{pmatrix}.
\]
The definition of $\tau^\pm_{i_0} \hat \epsilon_l^{i, j}$ are similar to
$\tau_{i_0} \hat \epsilon^{i, j}_l$ in Section \ref{fiber prod and orientation}.
See the next section for the precise definition.

Let $((\hat c_l), (x_l), (\alpha_l))$ be a triple of
\begin{itemize}
\item
a sequence of chains $\hat c_l$ in
$C_\ast(\overline{P}_{Y^-}, \overline{P}_{Y^-}^{\text{bad}}; \S_{Y^-}^D \otimes \Q)$
\item
a sequence of cycles $x_l$ in $K_X^0$, and
\item
a sequence of cochains $\alpha_l$ with compact supports in
$C^\ast(\overline{P}_{Y^+}, \overline{P}_{Y^+}^{\text{bad}}; \S_{Y^+}^D \otimes \Q)$,
\end{itemize}
For such a triple $((\hat c_l), (x_l), (\alpha_l))$, we define a pre-Kuranishi space
(or a linear combination of pre-Kuranishi spaces)
$\overline{\M}^X((\hat c_l), (x_l), (\alpha_l))$ by
\begin{align*}
&\overline{\M}^X((\hat c_l), (x_l), (\alpha_l))\\
& = \sum_{m_-, m_+ \geq 0}\sum_{\star_{m_-, m_+}}(-1)^\ast \overline{\M}^{(m_-, X, m_+)}
_{((\widetilde{G}^+_{m_+}, \widetilde{G}^-_{-m_-}), (\hat c^i_l), (x^i_l),
([\overline{P}_{Y^+}] \cap \alpha^i_l))}
\end{align*}
where
$\widetilde{G}^\pm = \widetilde{G}^\pm_0 + \widetilde{G}^\pm_{\pm1} +
\widetilde{G}^\pm_{\pm2} + \cdots = \Theta^\pm(e^{\otimes G^\pm})$ are appropriate
linear combinations of
\begin{multline*}
((\kappa\rho\Delta_\ast [\overline{P}_{Y^\pm}])^{i, j}, \dots,
(\kappa\rho\Delta_\ast [\overline{P}_{Y^\pm}])^{i, j},
\epsilon_{\overline{P}_{Y^\pm}}^{i, j}, \dots, \epsilon_{\overline{P}_{Y^\pm}}^{i, j},\\
(\Delta_\ast [\overline{P}_{Y^\pm}])^{i, j}, \dots, (\Delta_\ast [\overline{P}_{Y^\pm}])^{i, j})
\end{multline*}
defined in the next section.
(Pay attention to the order of $(\widetilde{G}^+_{m_+}, \widetilde{G}^-_{-m_-})$.
This is equivalent to $(-1)^{m_- m_+}(\widetilde{G}^-_{-m_-}, \widetilde{G}^+_{m_+})$.)
The sum $\star_{m_-, m_+}$ is taken over all decompositions 
\[
\{\hat c_l\} = \coprod_{-m_- \leq i \leq 0} \{\hat c^i_l\},
\quad \{\alpha_l\} = \coprod_{0 \leq i \leq m_+} \{\alpha^i_l\}
\]
as sets
and all decompositions
\[
\{x_l\} = \coprod_{-m_- \leq i \leq m_+} \{x^i_l\}
\]
such that $x^i_l \in K_{Y^-}^0$ for $-m_- \leq i <0$ and $x^i_l \in K_{Y^+}^0$
for $0 < i \leq m_+$.
The sign $\ast$ is the weighted sign of the permutation
\[
\begin{pmatrix}
(\hat c^1_l)_l & \cdots & (\hat c^m_l)_l
& (x^1_l)_l & \cdots & (x^m_l)_l & (\alpha^1_l)_l & \cdots & (\alpha^m_l)_l\\
&(\hat c_l)_l&&&(x_l)_l&&&(\alpha_l)_l&
\end{pmatrix}.
\]

We note $\widetilde{G}^\pm_0 = 1$.
Hence the main term is
\[
\overline{\M}^X_{((\hat c_l), (x_l), ([\overline{P}_{Y^+}] \cap \alpha_l))}.
\]

The following equation holds true.
\begin{align}
0 &= [\partial\overline{\M}^X((\hat c_l), (x_l), (\alpha_l))]^0\notag\\
&= [\overline{\M}^X(\partial((\hat c_l), (x_l), (\alpha_l)))]^0\notag\\
&\quad - \sum_{\star_-} (-1)^{\ast_-} \frac{1}{k !}
[\overline{\M}^{Y^-}((\hat c^-_l), (x^-_l), (\hat d_1^\ast, \hat d_2^\ast, \dots,
\hat d_k^\ast))]^0\notag\\
&\hph{\quad - \sum_{\star_-} (-1)^{\ast_-} \frac{1}{k !}}
\times [\overline{\M}^X((\hat d_k, \hat d_{k - 1}, \dots, \hat d_1) \cup (\hat c^0_l),
(x^0_l), (\alpha_l))]^0\notag\\
&\quad + \sum_{\star_+} (-1)^{\ast_+} \frac{1}{k !}
[\overline{\M}^X((\hat c_l), (x^0_l), (\alpha^0_l) \cup (\hat d_1^\ast, \hat d_2^\ast,
\dots, \hat d_k^\ast))]^0\notag\\
&\hph{\quad + \sum_{\star_+} (-1)^{\ast_+} \frac{1}{k !}}
\times [\overline{\M}^{Y^+}((\hat d_k, \hat d_{k - 1}, \dots, \hat d_1), (x^+_l), (\alpha^+_l))]^0
\label{boundary formula for X}
\end{align}
where the sum $\star_-$ is taken over $k \geq 0$, all simplices $d_l$ of $K_{Y^-}$
not contained in $\overline{P}_{Y^-}^{\text{bad}}$,
and all decompositions
\[
\{\hat c_l\} = \{\hat c^-_l\} \sqcup \{\hat c^0_l\}, \quad
\{x_l\} = \{x^-_l\} \sqcup \{x^0_l\}
\]
such that $x^-_l \in K_{Y^-}^0$.
The sign $\ast_-$ is the weighted sign of the permutation
\[
\begin{pmatrix}
(\hat c^-_l) \ (x^-_l) \ (\hat c^0_l) \ (x^0_l)\\
(\hat c_l) \quad (x_l)
\end{pmatrix}.
\]
The sum $\star_+$ is taken over $k \geq 0$, all simplices $d_l$ of $K_{Y^+}$ not contained
in $\overline{P}_{Y^+}^{\text{bad}}$, and all decompositions
\[
\{s_l\} = \{x^0_l\} \sqcup \{x^+_l\}, \quad \{\alpha_l\} = \{\alpha^0_l\} \sqcup \{\alpha^+_l\}
\]
such that $x^+_l \in K_{Y^+}^0$.
$\ast_+$ is the weighted sign of the permutation
\[
\begin{pmatrix}
(x^0_l) \ (\alpha^0_l) \ (x^+_l) \ (\alpha^+_l)\\
(x_l) \quad (\alpha_l)
\end{pmatrix}.
\]

To construct the algebra, we need to use the space of irreducible sequences of
holomorphic buildings.
Let $f^\pm_a$ be monomials of the form
\begin{multline*}
((\rho_\ast [\overline{P}_{Y^\pm}])^{i, j}, \dots,
(\rho_\ast [\overline{P}_{Y^\pm}])^{i, j},
\epsilon_{\overline{P}_{Y^\pm}}^{i, j}, \dots, \epsilon_{\overline{P}_{Y^\pm}}^{i, j},\\
(\Delta_\ast [\overline{P}_{Y^\pm}])^{i, j}, \dots, (\Delta_\ast [\overline{P}_{Y^\pm}])^{i, j}
)_{\substack{0 \leq i < j \leq m_{f_a^\pm}\\ (\text{or} -m_{f_a^\pm} \leq i < j \leq 0)}}
\end{multline*}
such that $m_\pm = \sum m_{f_a^\pm}$.
Then we define the space of irreducible sequences of holomorphic buildings
\begin{align*}
&(\overline{\M}^{(m_-, X, m_+)}
)_{(f^+_1 \otimes \dots \otimes f^+_{n^+}, f^-_1 \otimes \dots \otimes f^-_{n^-}),
(\breve c^i_l), (x^i_l), (\breve \eta^i_l))}^0\\
&\subset \overline{\M}^{(m_-, X, m_+)}
_{(\Theta^+(f^+_1 \otimes \dots \otimes f^+_{n^+}),
\Theta^-(f^-_1 \otimes \dots \otimes f^-_{n^-})), (\breve c^i_l), (x^i_l), (\breve \eta^i_l))}
\end{align*}
as follows. ($\Theta^\pm$ is defined by (\ref{Theta^pm}) in the next section.)
First we consider the case of $(n_-, n_+) \neq (0, 0)$.
A sequence of holomorphic buildings
\[
(\Sigma_i, s_i, u_i, \phi_i)_{-m_- \leq i \leq m_+} \in
\overline{\M}^{(m_-, X, m_+)}_{(\Theta^+(f^+_1 \otimes \dots \otimes f^+_{n^+}),
\Theta^-(f^-_1 \otimes \dots \otimes f^-_{n^-})), (\breve c^i_l), (x^i_l), (\breve \eta^i_l))}
\]
is contained in the above space if
\begin{itemize}
\item
each connected component of $\Sigma_0$ ($\Sigma_0$ is the $\overline{\M}^X$-factor
of $(\Sigma_i)_{-m_- \leq i \leq m_+}$. It is not necessarily of height one.)
concerns at least one monomial $f^\pm_i$,
that is, it contains at least one limit circle corresponding to
a variable in $f^\pm_i$,
and
\item
for any decomposition $\{f^-_1, \dots f^-_{n^-}, f^+_1, \dots, f^+_{n^+}\} = A \sqcup B$,
there exists a connected component of $\Sigma_0$ which concerns both of
some $f \in A$ and some $g \in B$.
\end{itemize}
If $(n_-, n_+) = (0, 0)$, then a holomorphic building $(\Sigma, z, u, \phi) \in
\overline{\M}^X_{((\breve c^i_l), (x^i_l), (\breve \eta^i_l))}$
is irreducible if it is connected.

First we note that all irreducible sequences of holomorphic buildings corresponding to
the zeros of the multisection of the $0$-dimensional component of the above
Kuranishi space have genera $\geq 0$ if each $f_a^\pm$ is contained
in $\mathring{\B}_{m_{f_a^\pm}}^\pm$
(this is also defined in the next section),
that is, if the number of variables in each $f_a^\pm$ each of which defines a relation of
the periodic orbit on one $+\infty$-limit circle of $\Sigma_i$ ($i \neq 0$) and
the periodic orbit on one $-\infty$-limit circle of $\Sigma_j$ ($j \neq 0$) is
$\geq m_{f_a^\pm} -1$.
This is because each factor $\Sigma_i$ except the $\overline{\M}^X$-factor is
connected by the dimensional reason.

We also note that for any sequence of holomorphic buildings
$(\Sigma_i, z_i, u_i, \phi_i)_i$ in
\[
\overline{\M}^{(m_-, X, m_+)}_{(\Theta^+(f^+_1 \otimes \dots \otimes f^+_{n^+}),
\Theta^-(f^-_1 \otimes f^-_2 \otimes \dots \otimes f^-_{n^-})), (\breve c^i_l), (x^i_l),
(\breve \eta^i_l))},
\]
we can decompose the set $\{f^\pm_a\}$ into sets $A_j$ such that
for any $j \neq j'$, there does not exist a connected component of $\Sigma_0$
which concerns both of some $f \in A_j$ and some $g \in A_j$,
and each $A_j$ cannot be decomposed further.
Hence each sequence of holomorphic buildings corresponding to a zero of the multisection
of the $0$-dimensional component can be decomposed into
ireducible sequences of holomorphic buildings contained
in the factors corresponding to $A_j$ and the connected holomorphic
buildings with height one.
For each connected holomorphic buildings with height one,
we add an empty set to $\{A_j\}$, and call $\{A_j\}$ as the irreducible decomposition
of $\{f^\pm_a\}$ corresponding to $(\Sigma_i, z_i, u_i, \phi_i)_i$.

For each triple $((\hat c_l), (x_l), (\alpha_l))$, we define a pre-Kuranishi space
\begin{align*}
&(\overline{\M}^X)^0((\hat c_l), (x_l), (\alpha_l))\\
& = \sum_{m_-, m_+ \geq 0}\sum_{\star_{m_-, m_+}}(-1)^\ast
(\overline{\M}^{(m_-, X, m_+)})^0_{((\widehat{G}^+_{m_+}, \widehat{G}^-_{-m_-}),
(\hat c^i_l), (x^i_l), ([\overline{P}_{Y^+}] \cap \alpha^i_l))}
\end{align*}
where
$e^{\otimes G^\pm} = \widehat{G}^\pm_0 + \widehat{G}^\pm_1 + \widehat{G}^\pm_2
+ \dots \in (\bigoplus_{m = 0}^\infty \bigotimes_{\sum l_i = m}(\B_{l_i}^+)^{l_i})^\wedge$.
Note that since $G^\pm$ is contained in
$(\bigoplus_{l = 1}^\infty (\mathring{\B}_l^+)^l)^\wedge$, the genera of the zero of the
multisection of the zero-dimensional component of
$(\overline{\M}^X)^0((\hat c_l), (x_l), (\alpha_l))$ are $\geq 0$.

The irreducible decomposition implies the following equation.
\begin{equation}
[\overline{\M}^X((\hat c_l), (x_l), (\alpha_l))]^0
= \sum (-1)^\ast \frac{1}{k !} \prod_{i = 1}^k
[(\overline{\M}^X)^0((\hat c^i_l), (x^i_l), (\alpha^i_l))]^0,
\label{irreducible decomposition}
\end{equation}
where the sum is taken over all $k \geq 0$ and all decompositions
\[
\{\hat c_l\} = \coprod_{i = 1}^k \{\hat c^i_l\}, \quad \{x_l\} = \coprod_{i = 1}^k \{x^i_l\}, \quad
\{\alpha_l\} = \coprod_{i = 1}^k \{\alpha^i_l\}
\]
as sets.
The sign $\ast$ is the weighted sign of the permutation
\[
\begin{pmatrix}
(c^1_l)\ (x^1_l)\ (\alpha^1_l) \dots (c^k_l)\ (x^k_l)\ (\alpha^k_l)\\
(c_l) \quad (x_l) \quad (\alpha_l)
\end{pmatrix}
\]
If $((\hat c_l), (x_l), (\alpha_l)) = (\emptyset, \emptyset, \emptyset)$, then the term
$\prod_{i = 1}^0 [(\overline{\M}^X)^0((\hat c^i_l), (x^i_l), (\alpha^i_l))]^0$ corresponding to
$k = 0$ on the right hand side of equation (\ref{irreducible decomposition}) is
defined by $1$, and otherwise it is defined by zero.
It corresponds to the number of the empty curve.
Equation (\ref{irreducible decomposition}) is proved as follows.
We write $\overline{\M}^X((\hat c_l), (x_l), (\alpha_l))$ as
\begin{align*}
&\overline{\M}^X((\hat c_l), (x_l), (\alpha_l)) \\
&= \sum_{N^-, N^+ \geq 0} (-1)^\ast \overline{\M}^{(\ast, X, \ast)}
_{((\Theta^+((G^+)^{\otimes N^+}), \Theta^-((G^-)^{\otimes N^-})), (\hat c^i_l), (x^i_l),
([\overline{P}_{Y^+}] \cap \alpha^i_l))},
\end{align*}
where we omit $m_\pm$ because they differ according to the variables in
$\Theta^\pm((G^+)^{\otimes N^\pm})$.
For each point $(\Sigma_i, z_i, u_i, \phi_i)_i$ in
\[
\overline{\M}^{(\ast, X, \ast)}
_{((\Theta^+((G^+)^{\otimes N^+}), \Theta^-((G^-)^{\otimes N^-})), (\hat c^i_l), (x^i_l),
([\overline{P}_{Y^+}] \cap \alpha^i_l))},
\]
we decompose the sequence of holomorphic buildings $(\Sigma_i, z_i, u_i, \phi_i)_i$
into irreducible sequences $(\Sigma^j_i, z^j_i, u^j_i, \phi^j_i)^j_i$ ($j = 1, \dots, k$).
Then it corresponds to a point in
$\prod^k (\overline{\M}^{(\ast, X, \ast)})^0_{((e^{\otimes G^+}, e^{\otimes G^-}),
\ast, \ast, \ast)} / \mathfrak{S}_k$,
and the group $\mathfrak{S}_k$ of permutation corresponds to the coefficient
$1 / k!$ in the right hand side of Equation (\ref{irreducible decomposition}).
Next we count the number of points $(\Sigma_i, z_i, u_i, \phi_i)_i$ in
\[
\overline{\M}^{(\ast, X, \ast)}
_{((\Theta^+((G^+)^{\otimes N^+}), \Theta^-((G^-)^{\otimes N^-})), \ast, \ast, \ast)}
\]
corresponding to a given
irreducible sequences $(\Sigma^j_i, z^j_i, u^j_i, \phi^j_i)^j_i$ ($j = 1, \dots, k$).
Define $n^\pm_j \geq 0$ by the condition that
for each $j$, $(\Sigma^j_i, z^j_i, u^j_i, \phi^j_i)_i$ is contained in
$(\overline{\M}^{(\ast, X, \ast)})^0_{(((G^+)^{\otimes n^+_j}, (G^-)^{\otimes n^-_j}),
\ast, \ast, \ast)}$.
Consider the irreducible decomposition $\{A_j\}_{j=1,\dots, k}$ of
$\{G^+_{(1)}, \dots, G^+_{(N^+)}, G^-_{(1)}, \dots, G^-_{(N^-)}\}$
($G^\pm_{(i)} = G^\pm$) corresponding to the point $(\Sigma_i, z_i, u_i, \phi_i)_i$,
and let $A_j = A^+_j \sqcup A^-_j$ be the decomposition into the sets
$A^\pm_j$ consisting of $G^\pm$.
Then $\# A^\pm_j = n^\pm_j$.
Conversely, for any decomposition $\{A_j\}_{j=1,\dots, k}$ of
$\{G^+_{(1)}, \dots, G^+_{(N^+)}, G^-_{(1)}, \dots, G^-_{(N^-)}\}$ such that
$\# A^\pm_j = n^\pm_j$, there exists a unique point in $\overline{\M}^{(\ast, X, \ast)}
_{((\Theta^+((G^+)^{\otimes N^+}), \Theta^-((G^-)^{\otimes N^-})), \ast, \ast, \ast)}$
corresponding to $(\Sigma^j_i, z^j_i, u^j_i, \phi^j_i)^j_i$ ($j = 1, \dots, k$) and
the decomposition $\{A_j\}_{j=1,\dots, k}$.
The number of such decompositions of the set is
\[
\frac{N^+ !}{n^+_1 ! n^+_2 ! \dots n^+_k !} \cdot \frac{N^- !}{n^-_1 ! n^-_2 ! \dots n^-_k !},
\]
and it coincides with the ratios of the product of the coefficients
$1 / N^\pm !$ of $(G^\pm)^{\otimes N^\pm}$ in $e^{\otimes G^\pm}$
on the left hand side of Equation (\ref{irreducible decomposition})
to the product of the coefficients
$1 / n^\pm_j !$ of $(G^\pm)^{\otimes n^\pm_j}$ in $e^{\otimes G^\pm}$
on the right hand side.
Hence Equation (\ref{irreducible decomposition}) holds true.

\begin{rem}
As in Remark \ref{connected or not}, we do not know whether or not we can choose
$G^\pm$ so that all irreducible sequences of holomorphic buildings
in the zero-dimensional
component of $(\overline{\M}^X)^0((\hat c_l), (x_l), (\alpha_l))$ are connected.
However, for the construction of the algebra in Section \ref{algebra for X},
it is enough to observe their genera are $\geq 0$.
\end{rem}

\subsection{Construction of the correction terms}\label{correction terms for X}
In this section, we construct $(G_{\pm m}^\pm)_{m \geq 1}$ used for the definition of
the correction terms of $\overline{\M}^X((\hat c_l), (x_l), (\alpha_l))$.
As in the case of the construction of $(F_m)_{m \geq 2}$, we consider algebras modeled
on the splitting of holomorphic buildings.
We construct $(G_m^+)_{m \geq 1}$ and $(G_{-m}^-)_{m \geq 1}$ independently.
First we construct $(G_m^+)_{m \geq 1}$.

For $m \geq 1$, let $B^+_m = \bigoplus_{n = 0}^{\frac{m(m + 1)}{2}} (B^+_m)^n$ be the
$\Z$-graded super-commutative algebra with coefficient $\R$ generated by variables
$\rho_{(e_i, e_j)}$, $\Delta_{(e_i, e_j)}$ and $\epsilon_{(e_i, e_j)}$ ($0 \leq i < j \leq m$).
The $\Z$-grading is defined by $\dim \rho_{(e_i, e_j)} = \dim \Delta_{(e_i, e_j)} = 0$ and
$\dim \epsilon_{(e_i, e_j)} = 1$.
For $m = 0$, we define $B^+_0 = \R$.

For each $m \geq 1$, the differential $\partial' : B^+_m \to B^+_m$ is defined by
$\partial' \epsilon_{(a, b)} = (-1)^m (\rho_{(a, b)} - \Delta_{(a, b)})$ and $\partial' \rho_{(a, b)}
= \partial' \Delta_{(a, b)} = 0$.
For $m = 0$, we define $\partial' = 0 : B_0^+ \to B_0^+$.

Homomorphisms $\tau^+_i : B^+_m \to B^+_{m + 1}$ ($0 \leq i \leq m$, $m \geq 1$)
are defined by
$\tau^+_i(x_{(a, b)}) = x_{(\hat \tau_i(a), \hat \tau_i(b))}$,
where $x$ is $\rho$, $\epsilon$ or $\Delta$, and each $\hat \tau^+_i$ is defined by
\[
\hat \tau^+_i (e_j) = \begin{cases}
e_j &j < i\\
e_i + e_{i + 1} &j = i\\
e_{j + 1} & j > i
\end{cases}.
\]
For $m = 0$, we define $\tau^+_0 = \id_{\R}$.
For $i > m$, we define $\tau^+_i = 0 : B^+_m \to B^+_{m + 1}$.

We define homomorphisms $\Diamond^+ : B^+_m \otimes A_{m'} \to B^+_{m+ m'}$
($m \geq 0$, $m' \geq 1$) by
\[
\Diamond^+ (f \otimes g) = (-1)^{1 + m m'} f \cdot
\exp(\rho_{(\sum_{0 \leq i \leq m} e_i, \sum_{m + 1 \leq j \leq m + m'} e_j)}) \cdot
g^{+ m}.
\]

We define homomorphisms $\Theta^+ : \bigotimes_{i = 1}^n B^+_{m_i} \to B^+_{\sum m_i}$
by
\begin{equation}
\Theta^+ (f_1 \otimes f_2 \otimes \dots \otimes f_n)
= f_1^{+\sum_{i = 2}^n m_i} \cdot f_2^{+\sum_{i = 3}^n m_i} \cdots f_n,
\label{Theta^pm}
\end{equation}
where $f^{+k}$ is defined by
\[
e^{+k}_j = \begin{cases}
e_0 &j = 0\\
e_{j + k} & j \neq 0
\end{cases}.
\]
For $n = 0$, we define $\Theta^+ = \id_\R : \R \to \R$.

Define a linear subspace $\Ddot B^+_m \subset B^+_m$ as follows.
For each $1 \leq i \leq m-1$ ($i \neq 0$) and each monomial
\[
f = x^{(1)}_{(a_1, b_1)} x^{(2)}_{(a_2, b_2)} \dots x^{(n)}_{(a_n, b_n)},
\]
such that $(a_j, b_j) \neq (e_i, e_{i+1})$, we define a monomial
\[
f^{(e_i, e_{i + 1})} = x^{(1)}_{(a'_1, b'_1)} x^{(2)}_{(a'_2, b'_2)} \dots x^{(n)}_{(a'_n, b'_n)}
\]
by permuting $e_i$ and $e_{i + 1}$ in $\{a_j, b_j\}$.
Then $\Ddot B^+_m \subset B^+_m$ is the subspace spanned by
$f + f^{(e_i, e_{i + 1})}$ for all such pair $i$ and $f$.

Define $\B^+_m = B^+_m / \Ddot B^+_m$.
This is not an algebra but the following maps are well defined.
\begin{align*}
\partial' &: \B^+_m \to \B^+_m & (m \geq 0)\\
\sum_{i \geq 1} (-1)^i e^{\Delta_{(e_i, e_{i + 1})}} \tau^+_i &: \B_m \to \B_{m + 1} &
(m \geq 0)\\
e^{\Delta_{(e_0, e_1)}} \tau^+_0 &: \B_m \to \B_{m + 1} & (m \geq 0)\\
\Diamond^+ &: \B^+_m \otimes \A_{m'} \to \B^+_{m + m'} & (m \geq 0, m' \geq 1)\\
\Theta^+ &: \bigotimes_{i = 1}^n \B^+_{m_i} \to \B^+_{\sum m_i} & (n \geq 0, m_i \geq 0)
\end{align*}

Further we define $\mathring{\B}^+_m \subset \B^+_m$ as follows.
First we define a new degree $\deg'$ by
\[
\deg' x_{(e_i, e_j)} = \begin{cases}
0 & i = 0\\
1 & i \geq 1
\end{cases}.
\]
Let $\mathring{B}^+_m \subset B^+_m$ be the ideal generated by monomials with
$\deg' \geq m - 1$ and define $\mathring{\B}^+_m = \mathring{B}^+_m
/ (\Ddot B^+_m \cap \mathring{B}^+_m)$.
It is easy to see that the homology of $((\mathring{\B}^+_m)^\ast, \partial')$ is zero
at $\ast \neq 0$.
($\ast$ is the dimension.)

Let $F^+ \in \A$ be a zero obtained in Section \ref{algebra for correction}.
We prove that there exists some $G^+ = G^+_1 +G^+_2 + \dots \in
(\bigoplus_{l = 1}^\infty (\mathring{\B}^+_l)^l)^\wedge$ such that
\begin{multline}
\partial' (\Theta^+(e^{\otimes G^+})) + \sum_i (-1)^i e^{\Delta{(e_i, e_{i + 1})}} \tau^+_i
\Theta^+(e^{\otimes G^+})\\
+ \Diamond^+ (\Theta^+(e^{\otimes G^+}) \otimes F^+) = 0,
\label{+G eq}
\end{multline}
where $e^{\otimes G^+} = 1 + G^+ + \frac{1}{2!} G^+ \otimes G^+
+ \frac{1}{3!} G^+ \otimes G^+ \otimes G^+ + \cdots$.
We inductively construct $G^+_{\leq m}= G^+_1 + G^+_2 + \cdots + G^+_m \in
\bigoplus_{l = 1}^m (\mathring{\B}_l^+)^l$ such that
\begin{multline}
\partial' (\Theta^+(e^{\otimes G^+_{\leq m}})) + \sum_i (-1)^i e^{\Delta{(e_i, e_{i + 1})}} \tau^+_i
\Theta^+(e^{\otimes G^+_{\leq m-1}})\\
 + \Diamond^+ (\Theta^+(e^{\otimes G^+_{\leq m-1}}) \otimes
F^+) \equiv 0 \label{G^+_m eq}
\end{multline}
in $(\bigoplus_{l = 1}^\infty (\B^+_l)^{l-1})^\wedge /
(\bigoplus_{l = m + 1}^\infty (\B^+_l)^{l - 1})^\wedge$.

First we define $G^+_{\leq 1} = G^+_1 \in (\mathring{\B}_1^+)^1$ by
\begin{align*}
G^+_1 = - \sum_{k = 1}^\infty \frac{1}{k !} (
&\underbrace{\epsilon_{(e_0, e_1)} \Delta_{(e_0, e_1)} \dots \Delta_{(e_0, e_1)}}_k\\
&+ \underbrace{\rho_{(e_0, e_1)} \epsilon_{(e_0, e_1)} \Delta_{(e_0, e_1)} \dots
\Delta_{(e_0, e_1)}}_k\\
&+ \dots + \underbrace{\rho_{(e_0, e_1)} \dots \rho_{(e_0, e_1)} \epsilon_{(e_0, e_1)}}_k).
\end{align*}
Then it is easy to check that this satisfies equation (\ref{G^+_m eq}).

Next assuming we have constructed $G^+_{\leq m - 1}$,
we prove there exists a required $G^+_{\leq m}$ ($m \geq 2$).
It is enough to show that
\begin{multline}
\partial' (\Theta^+(e^{\otimes G^+_{\leq m-1}})) + \sum_i (-1)^i e^{\Delta{(e_i, e_{i + 1})}} \tau^+_i
\Theta^+(e^{\otimes G^+_{\leq m-1}})\\
 + \Diamond^+ (\Theta^+(e^{\otimes G^+_{\leq m-1}}) \otimes
F^+) \equiv 0 \label{B mathring}
\end{multline}
in $(\bigoplus_{l = 1}^\infty (\B^+_l)^{l-1})^\wedge /((\bigoplus_{l = m + 1}^\infty
(\B^+_l)^{l - 1})^\wedge \oplus \bigoplus_{l = 1}^\infty (\mathring{\B}^+_l)^{l-1})$
and
\begin{equation}
\partial' \Bigl( \sum_i (-1)^i e^{\Delta{(e_i, e_{i + 1})}} \tau^+_i
\Theta^+(e^{\otimes G^+_{\leq m-1}}) + \Diamond^+ (\Theta^+(e^{\otimes G^+_{\leq m-1}})
\otimes F^+) \Bigr) \equiv 0 \label{B closed}
\end{equation}
in $(\bigoplus_{l = 2}^\infty (\B^+_l)^{l-2})^\wedge /
(\bigoplus_{l = m + 1}^\infty (\B^+_l)^{l-2})^\wedge$.

First we prove equation (\ref{B mathring}).
For the proof, we use the following maps
$\mathring{\tau}^+_0$ and $\mathring{\Diamond}^+$.

The linear map $\mathring{\tau}^+_0 : \bigotimes_{i = 1}^n B^+_{m_i} \to B^+_{\sum m_i}$ is
defined as follows.
Let $f_i \in B_{m_i}^+$ ($1 \leq i \leq n$) be monomials, and consider each term of
\[
\tau^+_0 \Theta^+(f_1 \otimes f_2 \otimes \dots \otimes f_n).
\]
In each term, some of $e_0$'s appearing in $f_i$ are changed to $e_1$ since
$\hat \tau_0^+$ maps $e_0$ to $e_0 + e_1$.
$\mathring{\tau}^+_0(f_1 \otimes f_2 \otimes \dots f_n)$ is defined by the sum of
the terms appearing in $\tau^+_0 \Theta^+(f_1 \otimes f_2 \otimes \dots \otimes f_n)$
such that each $f_i$ has at least one $e_0$ which is changed to $e_1$.
Then it induces a linear map $\mathring{\tau}^+_0 : \bigotimes_{i = 1}^n \B^+_{m_i}
\to \B^+_{\sum m_i}$.
For $n = 0$, we defined $\mathring{\tau}^+_0 = \id_\R : \R \to \R$.
For example, if $m_1 = 2$ and $m_2 = 1$, then
\begin{align*}
&\mathring{\tau}^+_0 (\Delta_{(e_0, e_1)} \epsilon_{(e_0, e_2)} \epsilon_{(e_1, e_2)} \otimes
\epsilon_{(e_0, e_1)})\\
&= (\Delta_{(e_0, e_3)} \epsilon_{(e_0, e_4)} + \Delta_{(e_0, e_3)} \epsilon_{(e_1, e_4)}
+ \Delta_{(e_1, e_3)} \epsilon_{(e_0, e_4)}) \epsilon_{(e_3, e_4)} \epsilon_{(e_0, e_2)}.
\end{align*}

The linear map $\mathring{\Diamond}^+ : (\bigotimes_{i = 1}^n B^+_{m_i}) \otimes A_{m'} \to
B^+_{\sum m_i + m'}$ is defined as follows.
Put $m = \sum_{i = 1}^n m_i$ and let
\[
\mathring{\exp}(\rho_{(\sum_{0 \leq i \leq m} e_i, \sum_{m + 1 \leq j \leq m + m'} e_j)})
\]
be the sum of all terms in
$\exp(\rho_{(\sum_{0 \leq i \leq m} e_i, \sum_{m + 1 \leq j \leq m + m'} e_j)})$
in which at least one $e_k$ appears for each $1 \leq l \leq n$ such that
$\sum_{1 \leq a < l} m_a + 1 \leq k \leq \sum_{1 \leq a \leq l} m_a$.
For example, if $n = 2$ and $m_1 = m_2 = m' = 2$, then
\begin{align*}
&\mathring{\exp}(\rho_{(\sum_{0 \leq i \leq 4} e_i, \sum_{j = 5, 6} e_j)})\\
& = \exp(\rho_{(\sum_{0 \leq i \leq 4} e_i, \sum_{j = 5, 6} e_j)})
- \exp(\rho_{(\sum_{i = 0, 3, 4} e_i, \sum_{j = 5, 6} e_j)})\\
& \quad - \exp(\rho_{(\sum_{i = 0, 1, 2} e_i, \sum_{j = 5, 6} e_j)})
+ \exp(\rho_{(e_0, \sum_{j = 5, 6} e_j)}).
\end{align*}
$\mathring{\Diamond}^+$ is defined by
\begin{align*}
\mathring{\Diamond}^+(f_1 \otimes f_2 \otimes \dots \otimes f_n \otimes g)
&=  (-1)^{1 + (\sum m_i) m'} \Theta^+(f_1 \otimes f_2 \otimes \dots \otimes f_n)\\
&\quad \cdot \mathring{\exp}(\rho_{(\sum_{0 \leq i \leq m} e_i, \sum_{m + 1 \leq j \leq m + m'} e_j)})
g^{+ \sum m_i}.
\end{align*}
It also induces a linear map
$\mathring{\Diamond}^+ : (\bigotimes_{i = 1}^n \B^+_{m_i}) \otimes \A_{m'} \to
\B^+_{\sum m_i + m'}$.

We can easily check the following equations for any
$G \in (\bigoplus_{l = 1}^\infty (\mathring{\B}^+_l)^l)^\wedge$.
\begin{equation}
\partial' \Theta^+\biggl(\frac{1}{k !} G^{\otimes k}\biggr) = \Theta^+\biggl(\frac{1}{(k-1) !}
G^{\otimes (k-1)} \otimes \partial' G\biggr)\label{Geq 1}
\end{equation}
\begin{align}
&\sum_{i \geq 1} (-1)^i e^{\Delta_{(e_i, e_{i+1})}} \tau^+_i \Theta^+
\biggl(\frac{1}{k !} G^{\otimes k}\biggr) \notag\\
&\hspace{60pt}
= \Theta^+\biggl(\frac{1}{(k-1) !} G^{\otimes (k-1)} \otimes
\sum_{i \geq 1} (-1)^i e^{\Delta_{(e_i, e_{i+1})}} \tau^+_i G\biggr)\label{Geq 2}
\end{align}
\begin{align}
&e^{\Delta_{(e_0, e_1)}} \tau^+_0 \Theta^+\biggl(\frac{1}{k !} G^{\otimes k}\biggr) \notag\\
&\hspace{60pt}
= \sum_{l = 0}^k  \Theta^+\biggl(\frac{1}{(k-l) ! l !}G^{\otimes (k-l)} \otimes
(e^{\Delta_{(e_0, e_1)}} \mathring{\tau}^+_0(G^{\otimes l}))\biggr)\label{Geq 3}
\end{align}
\begin{align}
&\Diamond^+\biggl(\Theta^+\biggl(\frac{1}{k !} G^{\otimes k}\biggr) \otimes F^+\biggr)
\notag\\
&\hspace{70pt}
=
\sum_{l = 0}^k  \Theta^+\biggl(\frac{1}{(k-l) ! l !}G^{\otimes (k-l)} \otimes
\mathring{\Diamond}^+ (G^{\otimes l} \otimes F^+)\biggr)\label{Geq 4}
\end{align}
Furthermore, it is easy to see that
\begin{equation}
\Theta^+(f_1 \otimes \dots \otimes f_k \otimes
\Theta^+(f_{k + 1} \otimes \dots \otimes f_n))
= \Theta^+(f_1 \otimes \dots f_n).\label{Theta eq}
\end{equation}

The assumption of the induction implies
\begin{align*}
R^{(m-1)} =& \ 
\partial' (\Theta^+(e^{\otimes G^+_{\leq m-1}}))
+ \sum_i (-1)^i e^{\Delta{(e_i, e_{i + 1})}}
\tau^+_i \Theta^+(e^{\otimes G^+_{\leq m-1}})\\
&+ \Diamond^+ (\Theta^+(e^{\otimes G^+_{\leq m-1}}) \otimes F^+)
\end{align*}
is zero in $(\bigoplus_{l = 1}^\infty (\B^+_l)^{l-1})^\wedge /
(\bigoplus_{l = m}^\infty (\B^+_l)^{l-1})^\wedge$.
Hence
\[
\frac{(-1)^l}{l !} \Theta^+((G^+_{\leq m-1})^{\otimes l} \otimes R^{(m-1)})
\equiv 0
\]
in $(\bigoplus_{l = 1}^\infty (\B^+_l)^{l-1})^\wedge /
(\bigoplus_{l = m + 1}^\infty (\B^+_l)^{l-1})^\wedge$ for all $l \geq 1$.
Therefore, for the proof of (\ref{B mathring}), it is enough to prove that
\begin{equation}
\sum_{l \geq 0} \frac{(-1)^l}{l !} \Theta^+((G^+_{\leq m-1})^{\otimes l} \otimes R^{(m-1)})
\equiv 0 \label{sum eq}
\end{equation}
in $(\bigoplus_{l = 1}^\infty (\B^+_l)^{l-1})^\wedge /(
(\bigoplus_{l = m + 1}^\infty (\B^+_l)^{l-1})^\wedge \oplus
\bigoplus_{l=1}^m (\mathring{\B}^+_l)^{l-1})$.
Equations (\ref{Geq 1}) to (\ref{Theta eq}) imply that the left hand side of
(\ref{sum eq}) is equal to
the sum of the following terms:
\begin{equation}
\sum_{l \geq 0} \frac{(-1)^l}{l !} \Theta^+((G^+_{\leq m-1})^{\otimes l} \otimes
\partial' (\Theta^+(e^{\otimes G^+_{\leq m-1}})))
= \partial' G^+_{\leq m-1} \label{term 1}
\end{equation}
\begin{align}
&\sum_{l \geq 0} \frac{(-1)^l}{l !} \Theta^+\Bigl((G^+_{\leq m-1})^{\otimes l} \otimes
\sum_{i \geq 1} (-1)^i e^{\Delta{(e_i, e_{i + 1})}}
\tau^+_i \Theta^+(e^{\otimes G^+_{\leq m-1}})\Bigr) \notag\\
&= \sum_{i \geq 1} (-1)^i e^{\Delta{(e_i, e_{i + 1})}}
\tau^+_i G^+_{\leq m-1} \label{term 2}
\end{align}
\begin{align}
&\sum_{l \geq 0} \frac{(-1)^l}{l !} \Theta^+((G^+_{\leq m-1})^{\otimes l} \otimes
e^{\Delta{(e_0, e_1)}} \tau^+_0 \Theta^+(e^{\otimes G^+_{\leq m-1}})) \notag\\
&= e^{\Delta{(e_0, e_1)}} \mathring{\tau}^+_0 \Theta^+(e^{\otimes G^+_{\leq m-1}})
\label{term 3}
\end{align}
\begin{align}
&\sum_{l \geq 0} \frac{(-1)^l}{l !} \Theta^+((G^+_{\leq m-1})^{\otimes l} \otimes
\Diamond^+ (\Theta^+(e^{\otimes G^+_{\leq m-1}}) \otimes F^+)) \notag\\
&= \mathring{\Diamond}^+ (\Theta^+(e^{\otimes G^+_{\leq m-1}}) \otimes F^+)
\label{term 4}
\end{align}
Terms (\ref{term 1}), (\ref{term 3}), (\ref{term 4}) and
$(e^{\Delta_{(e_i, e_{i + 1})}} - 1) \tau^+_i G^+_{\leq m-1}$ ($i \geq 1$)
are contained in
$(\bigoplus_l (\mathring{\B}^+_l)^{l-1})^\wedge$,
and $\tau^+_i G^+_{\leq m-1} \equiv 0$ in
$\bigoplus_{l = 1}^\infty (\B_l^+)^{l-1}$ for $i > 0$.
(In general, $\tau_i^+ f$ is contained in $\ddot B_{m + 1}^+$ for any $f \in B_m^+$
and $i > 0$.)
These prove equation (\ref{sum eq}).
Therefore we can construct $G^+_{\leq m}$ inductively.

As with equation (\ref{A closed}) in Section \ref{algebra for correction},
equation (\ref{B closed}) is proved as follows.
Put $\widetilde{G}^+ = \Theta^+(e^{\otimes G^+_{\leq m-1}})$.
The left hand side of (\ref{B closed}) is
\begin{align*}
&\sum_i (-1)^{i+1} e^{\Delta_{(e_i, e_{i + 1})}} \tau^+_i \partial' \widetilde{G}^+
+ \Diamond^+ \Bigl(\partial' \widetilde{G}^+ \otimes \sum_j (-1)^j F^+_j \Bigr)
+ \Diamond^+(\widetilde{G}^+ \otimes \partial' F^+)\\
&= \sum_i (-1)^i e^{\Delta_{(e_i, e_{i + 1})}} \tau^+_i \Bigl(\sum_j (-1)^j
e^{\Delta_{(e_j, e_{j + 1})}}
\tau^+_j \widetilde{G}^+ + \Diamond^+(\widetilde{G}^+ \otimes F^+)\Bigr)\\
&\quad - \Diamond^+\Bigl(\Bigl(\sum_i (-1)^i e^{\Delta_{(e_i, e_{i + 1})}}
\tau^+_i \widetilde{G}^+ + \Diamond^+(\widetilde{G}^+ \otimes F^+)\Bigr)
\otimes \sum_j (-1)^j F^+_j\Bigr)\\
&\quad - \Diamond^+\Bigl(\widetilde{G}^+ \otimes \Bigl(\sum_i (-1)^i
e^{\Delta_{(e_i, e_{i + 1})}} \tau_i F^+ + \Box(F^+ \otimes F^+)\Bigr)\Bigr)
\end{align*}
and this is zero because
{\belowdisplayskip= 0pt
\[
\Bigl(\sum_i (-1)^i e^{\Delta_{(e_i, e_{i + 1})}} \tau^+_i \Bigr) \circ \Bigl(\sum_j (-1)^j
e^{\Delta_{(e_j, e_{j + 1})}} \tau^+_j \Bigr) = 0,
\]
}
\begin{align*}
\sum_i (-1)^i e^{\Delta_{(e_i, e_{i + 1})}} \tau^+_i \Diamond^+(f \otimes g)
- \Diamond^+\Bigl(\sum_i (-1)^i e^{\Delta_{(e_i, e_{i + 1})}} \tau^+_i f \otimes
(-1)^{\deg g} g\Bigr)\\
- \Diamond^+\Bigl(f \otimes \sum_i (-1)^i e^{\Delta_{(e_i, e_{i + 1})}} \tau^+_i g \Bigr)
= 0,
\end{align*}
\[
\Diamond^+(f \otimes \Box(g \otimes h)) + \Diamond^+(\Diamond^+(f \otimes g)
\otimes (-1)^{\deg h} h) = 0.
\]

Next, we construct $(G_{-m}^-)_{m \geq 1}$.
For $m \geq 1$, let $B^-_{-m}$ be the $\Z$-graded super-commutative algebra
with coefficient $\R$ generated by variables 
$\rho_{(e_i, e_j)}$, $\Delta_{(e_i, e_j)}$ and $\epsilon_{(e_i, e_j)}$ ($-m \leq i < j \leq 0$).
The $\Z$-grading is defined by $\dim \rho_{(e_i, e_j)} = \dim \Delta_{(e_i, e_j)} = 0$ and
$\dim \epsilon_{(e_i, e_j)} = 1$.
For $m = 0$, we define $B^-_0 = \R$.

For each $m \geq 1$, the differential $\partial' : B^-_{-m} \to B^-_{-m}$ is defined by
$\partial' \epsilon_{(a, b)} = (-1)^m (\rho_{(a, b)} - \Delta_{(a, b)})$ and $\partial' \rho_{(a, b)}
= \partial' \Delta_{(a, b)} = 0$.
Homomorphisms $\tau^-_i : B^-_{-m} \to B^-_{- m - 1}$ ($-m \leq i \leq 0$, $m \geq 1$)
are defined by
$\tau^-_i(x_{(a, b)}) = x_{(\hat \tau_i(a), \hat \tau_i(b))}$,
where $\hat \tau^-_i$ is defined by
\[
\hat \tau^-_i (e_j) = \begin{cases}
e_{j - 1} &j < i\\
e_{i - 1} + e_i &j = i\\
e_j & j > i
\end{cases}.
\]
For $m = 0$, we define $\tau^-_0 = \id_{\R}$.
For $i < -m$, we define $\tau^-_i = 0 : B^-_{-m} \to B^-_{-m - 1}$.
We define $\tilde\tau^-_i = (-1)^{m + 1 + i} \tau^-_i : B^-_{-m} \to B^-_{-m - 1}$.

We define homomorphisms $\Diamond^- : A_m \otimes B^-_{-m'} \to B^-_{- m - m'}$
($m \geq 1$, $m' \geq 0$) by
\[
\Diamond^- (f \otimes g) = (-1)^{(m - 1)m'} f^{-m'} \cdot
\exp(\rho_{(\sum_{-m -m' \leq i \leq -m' - 1} e_i, \sum_{-m' \leq j \leq 0} e_j)}) \cdot g.
\]

We define homomorphisms $\Theta^- : \bigotimes_{i = 1}^n B^-_{-m_i} \to B^-_{\sum -m_i}$
by
\[
\Theta^- (f_1 \otimes f_2 \otimes \dots \otimes f_n)
= f_1 \cdot f_2^{-m_1} \cdots f_n^{-\sum_{i=1}^{n - 1} m_i},
\]
where $f^{-k}$ is defined by
\[
e^{-k}_j = \begin{cases}
e_0 &j = 0\\
e_{j - k} & j \neq 0
\end{cases}.
\]

Define $\B^-_m$ and $\mathring{\B}^-_m \subset \B^-_m$ similarly.
In this case, $\deg'$ is defined by
\[
\deg' x_{(e_i, e_j)} = \begin{cases}
0 & j = 0\\
1 & j \leq -1
\end{cases}.
\]

Let $F^- \in \A$ be a zero in Section \ref{algebra for correction}.
(We do not need to assume $F^- = F^+$.)
As in the case of $G^+$, we can construct
$G^- = G^-_{-1} + G^-_{-2} + \dots \in (\bigoplus_{l = 1} (\B^-_{-l})^l)^\wedge$ such that
\begin{align}
&\partial' (\Theta^-(e^{\otimes G^-})) + \sum_i e^{\Delta{(e_{i - 1}, e_i)}} \tilde\tau^-_i
\Theta^-(e^{\otimes G^-}) \notag\\
&\hspace{120pt}
+ \Diamond^- (F^- \otimes \Theta^-(e^{\otimes G^-})) = 0. \label{-G eq}
\end{align}
Note that
\begin{align*}
G^-_1 = \sum_{k=1}^\infty \frac{1}{k !}(
&\underbrace{\epsilon_{(e_{-1}, e_0)} \Delta_{(e_{-1}, e_0)} \dots \Delta_{(e_{-1}, e_0)}}_k\\
&+ \underbrace{\rho_{(e_{-1}, e_0)} \epsilon_{(e_{-1}, e_0)} \Delta_{(e_{-1}, e_0)} \dots
\Delta_{(e_{-1}, e_0)}}_k\\
&+ \dots + \underbrace{\rho_{(e_{-1}, e_0)} \dots
\rho_{(e_{-1}, e_0)} \epsilon_{(e_{-1}, e_0)}}_k).
\end{align*}

Equation (\ref{boundary formula for X}) is satisfied for the solutions
$G^+$ of (\ref{+G eq}) and $G^-$ of (\ref{-G eq}) because (\ref{boundary of MMX})
implies
\begin{align*}
&\sum_{\star_{m_-, m_+}} (-1)^\ast \partial' \overline{\M}^{(m_-, X, m_+)}
_{((\widetilde{G}^+_{m_+}, \widetilde{G}^-_{-m_-}), (\hat c^i_l), (x^i_l), (\hat \eta^i_l))}\\
&= \sum_{\star_{m_-, m_+}} (-1)^\ast \overline{\M}^{(m_-, X, m_+)}
_{((\widetilde{G}^+_{m_+}, \widetilde{G}^-_{-m_-}), \partial((\hat c^i_l), (x^i_l),
(\hat \eta^i_l)))}\\
& \quad + \sum_{\star_{m_-, m_+}} (-1)^{\ast + m_-} \overline{\M}^{(m_-, X, m_+)}
_{((\partial' \widetilde{G}^+_{m_+}, \widetilde{G}^-_{-m_-}), (\hat c^i_l), (x^i_l),
(\hat \eta^i_l))}\\
& \quad + \sum_{\star_{m_-, m_+ +1}}
(-1)^{\ast + m_-} \overline{\M}^{(m_-, X, m_+ +1)}_{((\sum_{i = 0}^{m_+}
(-1)^i e^{\Delta_{(e_i, e_{i+1})}}\widetilde{G}^+_{m_+},
\widetilde{G}^-_{-m_-}),
(\hat c^i_l), (x^i_l), (\hat \eta^i_l))}\\
& \quad + \sum_{\star_{m_-, m_+}} (-1)^\ast \overline{\M}^{(m_-, X, m_+)}
_{((\widetilde{G}^+_{m_+}, \partial' \widetilde{G}^-_{-m_-}), (\hat c^i_l), (x^i_l),
(\hat \eta^i_l))}\\
& \quad + \sum_{\star_{m_- + 1, m_+}}
(-1)^\ast \overline{\M}^{(m_- +1, X, m_+)}_{((\widetilde{G}^+_{m_+},
\sum_{i = -m_-}^0 e^{\Delta_{(e_{i-1}, e_i)}} \tilde \tau^-_i \widetilde{G}^-_{-m_-}),
(\hat c^i_l), (x^i_l), (\hat \eta^i_l))}
\end{align*}
and the following equations hold true.
For any $m$, $(m_-, m_+)$ and $((\hat c^i_l)_{i = -m - m_-}^{m_+},\ab
(x^i_l)_{i = -m - m_-}^{m_+},\ab (\hat \eta^i_l)_{i = - m_-}^{m_+})$,
\begin{align*}
&\sum \frac{1}{k!} \Bigl[(\overline{\M}_{Y^-})^m_{(F^-_m, (\hat c^{i -m_-}_l)_{i = -m}^{- 1},
(x^{i - m_-}_l)_{i = -m}^{- 1},
([\overline{P}_{Y^-}] \cap \hat d_1^\ast, \dots, [\overline{P}_{Y^-}] \cap
\hat d_k^\ast))}\Bigr]^0\\
& \hph{\sum \frac{1}{k!}} \cdot
\bigl[\overline{\M}^{(m_-, X, m_+)}_{((\widetilde{G}^+_{m_+}, \widetilde{G}^-_{-m_-}),
(\hat d_k, \dots, \hat d_1) \cup (\hat c^i_l), (x^i_l), (\hat \eta^i_l))}\bigr]^0\\
& = \bigl[\overline{\M}^{(m_- + m, X, m_+)}_{((\widetilde{G}^+_{m_+}, \Diamond^-(F^-_m
\otimes \widetilde{G}^-_{-m_-})), (\hat c^i_l), (x^i_l), (\hat \eta^i_l))}\bigr]^0
\end{align*}
where the sum is taken over all $k \geq 0$ and all simplices $d_l$ of $K_{Y^-}$
not contained in $\overline{P}_{Y^-}^{\text{bad}}$, and
for any $m$, $(m_-, m_+)$ and $((\hat c^i_l)_{i = -m_-}^{m_+},\ab (x^i_l)_{i = -m_-}^{m_+ + m},
\ab (\hat \eta^i_l)_{i = -m_-}^{m_+ + m})$,
\begin{align*}
&\sum \frac{1}{k!} \bigl[\overline{\M}^{(m_-, X, m_+)}_{((\widetilde{G}^+_{m_+},
\widetilde{G}^-_{-m_-}), (\hat c^i_l), (x^i_l), (\hat \eta^i_l) \cup
([\overline{P}_{Y^-}] \cap \hat d_1^\ast, \dots, [\overline{P}_{Y^-}] \cap
\hat d_k^\ast))}\bigr]^0\\
& \hph{\sum \frac{1}{k!}} \cdot\Bigl[(\overline{\M}_{Y^+})^m
_{(F^+_m, (\hat d_k, \dots, \hat d_1),
(x^{i + m_+}_l)_{i = 1}^m, (\hat \eta^{i+ m_+}_l)_{i = 1}^m)}\Bigr]^0\\
&= (-1)^{1 + m_-}\bigl[\overline{\M}^{(m_-, X, m_+ + m)}
_{((\Diamond^+(\widetilde{G}^+_{m_+} \otimes F^+_m), \widetilde{G}^-_{-m_-}),
(\hat c^i_l), (x^i_l), (\hat \eta^i_l))}\bigr]^0
\end{align*}
where the sum is taken over all $k \geq 0$ and all simplices $d_l$ of $K_{Y^+}$
not contained in $\overline{P}_{Y^+}^{\text{bad}}$


%% file: SFT-09_Algebras_for_X.tex
%
%

\subsection{Construction of the algebras}\label{algebra for X}
In this section, we construct the algebra for $X$.
It gives a kind of chain map between the algebras for $Y^-$ and $Y^+$
in the sense of SFT.
We follow the argument of \cite{EGH00}.

First we consider the case of general SFT.
We define a super-commutative algebra $\D_X = \D_{(X, \omega, Y^\pm, \lambda^\pm,
K_{Y^\pm}, \overline{K}_X^0)}$ as follows.
Its elements are formal series
\[
\sum_{(\hat c_i^\ast), (\hat c'_i), e}
f_{(\hat c_i^\ast), (\hat c'_i), e}(t, \hbar) q^-_{\hat c_1^\ast}
\dots q^-_{\hat c_{k_q}^\ast} p^+_{\hat c'_1} \dots p^+_{\hat c'_{k_p}} T^e,
\]
where $f_{(\hat c_i^\ast), (\hat c'_i), e}(t, \hbar) \in \R[[t, \hbar]]$ is a formal series
of the variables $t_x$ ($x \in K_X^0$) and $\hbar$,
and the infinite sum is taken over all sequences $((\hat c_i), (\hat c'_i))$ consisting of
simplices $\hat c_i$ of $K_{Y^-}$ not contained in $\overline{P}^{\text{bad}}_{Y^-}$ and
simplices $(\hat c'_i)$ of $K_{Y^+}$ not contained in $\overline{P}^{\text{bad}}_{Y^+}$,
and $e \in \tilde \omega H_2(\overline{X}, \partial \overline{X}; \Z)$
($\cong H_2(\overline{X}, \partial \overline{X}; \Z) / \Ker \tilde \omega$)
with the following Novikov condition:
for any $C \geq 0$, the number of the non-zero terms with
$\sum_j e(p^+_{\hat c'_j}) \geq - C$
and $e + \sum_j e(p^+_{\hat c'_j}) \geq - C$ is finite.
The product is defined so that all variables are super-commutative,
where $\Z / 2$-degree is similar to the case of $\W_Y$ except
$|t_x| = \codim_X x$ and $|T^e| = 0$.
We also define a submodule $\D_X^{\leq \kappa} \subset \D_X$ for each
$\kappa \geq 0$ by the condition
$\sum_i e(q^-_{\hat c_i^\ast}) + e + e(p^+_{\hat c'_i}) \leq \kappa$.

To define differentials on quotients of $\D_X$, we use a bigger
super-commutative algebra
$\D\D_X = \D\D_{(X, \omega, Y^\pm, \lambda^\pm, K_{Y^\pm}, \overline{K}^0_X,
\overline{K}^0_{Y^\pm}, \mu^\pm)}$.
Its elements are formal series
\[
\sum_{(\hat c_i^\ast), (\hat c'_i), e}
f_{(\hat c_i^\ast), (\hat c'_i), e}(t, \hbar) q^-_{\hat c_1^\ast}
\dots q^-_{\hat c_{k_q}^\ast} p^+_{\hat c'_1} \dots p^+_{\hat c'_{k_p}} T^e,
\]
where in this case, $f_{(\hat c_i^\ast), (\hat c'_i), e}(t, \hbar) \in
\R[[\hbar]] [\hbar^{-1}] [[t]]$, namely,
the coefficient of each monomial of the $t$-variables in
$f_{(\hat c_i^\ast), (\hat c'_i), e}(t, \hbar)$ is arrowed to have a pole of finite degree
at $\hbar = 0$.
(The degrees do not need to be bounded.)
For each $\kappa \geq 0$, we define a submodule $\D\D_X^{\leq \kappa} \subset
\D\D_X$ by the condition
$\sum_i e(q^-_{\hat c_i^\ast}) + e + e(p^+_{\hat c'_i}) \leq \kappa$.
For each positive constant $\delta > 0$,
we also define a submodule $\D\D_X^{\leq \kappa, \delta} \subset
\D\D_X^{\leq \kappa}$ by the condition
\begin{equation}
\widetilde{g}_\delta := g + \frac{1}{2}(k_t + k_q + k_p)
- \frac{\sum_i e(q^-_{\hat c_i^\ast}) + e + e(p^+_{\hat c'_i})}{\delta}
\geq - \frac{\kappa}{\delta}. \label{tilde gC}
\end{equation}
Note that $\D_X^{\leq \kappa} \subset \D\D_X^{\leq \kappa, \delta}$
and $\D\D_X^{\leq \kappa, \delta} \subset \D\D_X^{\leq \kappa, \delta'}$
for $\delta \geq \delta'$.

Define submodules $\widetilde{J}^{\leq \kappa, \delta}_{C_0, C_1, C_2}
= \widetilde{J}^{\leq \kappa, \delta}_{X, C_0, C_1, C_2} \subset
\D\D_X^{\leq \kappa, \delta}$ by
\begin{align*}
\widetilde{J}^{\leq \kappa, \delta}_{C_0, C_1, C_2} &= \Bigl\{\sum
a_{(x_i), (\hat c_i^\ast), (\hat c'_i), g, e} t_{x_1} \dots t_{x_{k_t}} q^-_{\hat c_1^\ast}
\dots q^-_{\hat c_{k_q}^\ast} p^+_{\hat c'_1} \dots p^+_{\hat c'_{k_p}} \hbar^g T^e
\in \D\D_X^{\leq \kappa, \delta};\\
&\quad \quad a_{(x_i), (\hat c_i^\ast), (\hat c'_i), g, e} = 0 \text{ for all }
((x_i)_{i = 1}^{k_t}, (\hat c_i^\ast),_{i = 1}^{k_q} (\hat c'_i)_{i = 1}^{k_p}, g, e)
\text{ such that}\\
&\quad \quad
k_t \leq C_0,\ \widetilde{g}_\delta \leq C_1,\ \sum e(p^+_{\hat c'_i}) \geq - C_2
\text{ and } e + \sum e(p^+_{\hat c'_i}) \geq - C_2 \Bigr\}.
\end{align*}
Note that these are ideals if $\kappa = 0$.
Note also that
$\widetilde{J}^{\leq \kappa, \delta}_{C_0, C_1 + \kappa((\delta')^{-1} - \delta^{-1}), C_2}
\subset \widetilde{J}^{\leq \kappa, \delta'}_{C_0, C_1, C_2}$ for $\delta \geq \delta'$,
which implies that we have a natural map
\[
\D\D_X^{\leq \kappa, \delta}
/ \widetilde{J}^{\leq \kappa, \delta}_{C_0, C_1 + \kappa((\delta')^{-1} - \delta^{-1}), C_2}
\to \D\D_X^{\leq \kappa, \delta'}
/ \widetilde{J}^{\leq \kappa, \delta'}_{C_0, C_1, C_2}.
\]
We also define submodules $J^{\leq \kappa, \delta}_{C_0, C_1, C_2} \subset
\D_X^{\leq \kappa}$ by $J^{\leq \kappa, \delta}_{C_0, C_1, C_2}
= \widetilde{J}^{\leq \kappa, \delta}_{C_0, C_1, C_2} \cap \D_X^{\leq \kappa}$.

Let $(\hbar^{-1} \D_X^{\leq 0})^{\star, \delta} \subset \hbar^{-1} \D_X^{\leq 0}$
be a submodule defined by the following conditions:
\begin{itemize}
\item
$\widetilde{g}_\delta$ is nonnegative.
(Hence $(\hbar^{-1} \D_X^{\leq 0})^{\star, \delta} \subset \D\D_X^{\leq 0, \delta}$.)
\item
The constant term is zero.
\end{itemize}
We also define submodules $J^{\star, \delta}_{C_0, C_1, C_2}
= J^{\star, \delta}_{X, C_0, C_1, C_2} \subset
(\hbar^{-1} \D_X^{\leq 0})^{\star, \delta}$ by
\begin{align*}
&J^{\star, \delta}_{C_0, C_1, C_2} \\
&= \Bigl\{\sum
a_{(x_i), (\hat c_i^\ast), (\hat c'_i), g, e} t_{x_1} \dots t_{x_{k_t}} q^-_{\hat c_1^\ast}
\dots q^-_{\hat c_{k_q}^\ast} p^+_{\hat c'_1} \dots p^+_{\hat c'_{k_p}} \hbar^g T^e
\in (\hbar^{-1} \D_X^{\leq 0})^{\star, \delta};\\
&\quad \quad a_{(x_i), (\hat c_i^\ast), (\hat c'_i), g, e} = 0 \text{ for all }
((x_i)_{i = 1}^{k_t}, (\hat c_i^\ast),_{i = 1}^{k_q} (\hat c'_i)_{i = 1}^{k_p}, g, e)
\text{ such that}\\
&\quad \quad
k_t \leq C_0,\ \widetilde{g}_\delta \leq C_1,\ \sum e(p^+_{\hat c'_i}) \geq - C_2
\text{ and } e + \sum e(p^+_{\hat c'_i}) \geq - C_2 \Bigr\}.
\end{align*}

We say that $\delta > 0$ is admissible for $C_2$ if
\begin{itemize}
\item
$\delta \leq L^\pm_{\min}$, where $L^\pm_{\min}$ is the minimal period of
periodic orbits in $(Y^\pm, \lambda^\pm)$, and
\item
$\delta \leq E_{\hat \omega}(u)$
for any non-constant holomorphic building $(\Sigma, u)$ for $X$
of genus $0$ and height $1$ such that the number of the limit circle is $\leq 1$
and the period of the periodic orbit on the circle is $\leq C_2$ (if it exists).
\end{itemize}

If we fix a triple $(\overline{C}_0, \overline{C}_1, \overline{C}_2)$ and an admissible
constant $\delta$ for $\overline{C}_2$, then, choosing a compatible family of
perturbations $\B_X$ of the multisections of finite number of pre-Kuranishi spaces
(these also need to be compatible with $\B_{Y^\pm}$)
and using their virtual fundamental chains, we can define the generating functions
\begin{align*}
\F &= \hbar^{-1} \sum_{g \geq 0} \F_g \hbar^g
\in (\hbar^{-1} \D_X^{\leq 0})^{\star, \delta}
/ J^{\star, \delta}_{\overline{C}_0, \overline{C}_1, \overline{C}_2}\\
\widetilde{\F} &= \hbar^{-1} \sum_{g \in \Z} \widetilde{\F}_g \hbar^g
\in \D\D_X^{\leq 0, \delta}
/ \widetilde{J}^{\leq 0, \delta}_{\overline{C}_0, \overline{C}_1, \overline{C}_2}
\end{align*}
by
\[
\F_g = \sum_{k_q, k_t, k_p \geq 0, e} \frac{1}{k_q ! k_t ! k_p !} [ (\overline{\M}^X_{g, e})^0(
\underbrace{\mathbf{q}, \dots, \mathbf{q}}_{k_q},
\underbrace{\mathbf{t}, \dots, \mathbf{t}}_{k_t},
\underbrace{\mathbf{p}, \dots, \mathbf{p}}_{k_p})]^0 T^{-e}
\]
and
\[
\widetilde{\F}_g = \sum_{k_q, k_t, k_p \geq 0, e}
\frac{1}{k_q ! k_t ! k_p !} [ (\overline{\M}^X_{g, e})(
\underbrace{\mathbf{q}, \dots, \mathbf{q}}_{k_q},
\underbrace{\mathbf{t}, \dots, \mathbf{t}}_{k_t},
\underbrace{\mathbf{p}, \dots, \mathbf{p}}_{k_p})]^0 T^{-e},
\]
where $\mathbf{q} = \sum_c q_{\hat c^\ast} \hat c$, $\mathbf{t} = \sum_x t_x x$ and
$\mathbf{p} = \sum_c p_{\hat c} \hat c^\ast$ are formal series.
Sometimes we explicitly indicate the dependence of $\F$ to various data as
\[
\F = \F_{(X, \omega, Y^\pm, \lambda^\pm, K_{Y^\pm}, K_X^0, K_{Y^\pm}^0,
\mu^\pm, K_{Y^\pm}^2, J, \B_X)}.
\]

$\widetilde{\F}$ indeed satisfies the condition of
$\D\D_X^{\leq 0, \delta}
/ \widetilde{J}^{\leq 0, \delta}_{\overline{C}_0, \overline{C}_1, \overline{C}_2}$, that is,
$\widetilde{g}_\delta \geq 0$ for all terms such that
$k_t \leq \overline{C}_0$, $\sum e(p^+_{\hat c'_i}) \geq - \overline{C}_2$
and $e + \sum e(p^+_{\hat c'_i}) \geq - \overline{C}_2$.
It is enough to see that every holomorphic building $(\Sigma, z, u, \phi)$ such that
$\sum_j L_{\gamma_{+\infty_j}} \leq \overline{C}_2$
and $e + \sum_j L_{\gamma_{+\infty_j}} \leq \overline{C}_2$ satisfies
\begin{equation}
\widetilde{g}_\delta = g + \frac{1}{2}(k_t + k_q + k_p)
+ \frac{E_{\hat \omega}(u)}
{\delta} \geq 1,
\label{tilde gC ineq}
\end{equation}
where $g$ is its genus, $k_t$, $k_q$ and $k_p$ are the numbers of its marked points,
$-\infty$-limit circles, and $+\infty$-limit circles respectively.
($L_{\gamma_{\pm\infty_i}}$ are the periods of the periodic orbits on its limit circles.)
First note that $\widetilde{g}_\delta - 1$ is additive with respect to
disjoint union of holomorphic buildings, and that if a holomorphic building
$(\Sigma', z', u', \phi')$ for $Y^-$ or $Y^+$ is glued to a holomorphic
building for $X$, then $\widetilde{g}$ is changed by
more than or equal to the corresponding $\widetilde{g}$ of $(\Sigma', z', u', \phi')$
since $\delta \leq L^\pm_{\min}$.
Therefore, it is enough to show inequality (\ref{tilde gC ineq}) for a connected
holomorphic building of height one.
Assume contrary, that is, assume that there exists a holomorphic building
$(\Sigma, z, u, \phi)$ of height one such that $\tilde g < 1$.
Since $\tilde g < 1$ implies $g = 0$ and $k_t \leq 1$, $u$ is not a constant map.
Note that the period of the periodic orbits on its $-\infty$-limit circle is
$\leq e + \sum_j L_{\gamma_{+\infty_j}} \leq \overline{C}_2$ (if it exists) by
(\ref{E hat omega estimate}), and
$\tilde g < 1$ implies that the number of the limit circles ($= k_q + k_p$) is $\leq 1$.
Therefore, the assumption of $\delta$ implies that $\delta \leq E_{\hat \omega}(u)$,
which contradicts to the assumption $\tilde g < 1$.
Hence $\tilde F$ satisfies the condition $\widetilde{g}_\delta \geq 0$.

$\F$ also satisfies the condition $\widetilde{g}_\delta \geq 0$.
Furthermore, the degree of $\F$ is even because of the dimension of pre-Kuranishi
spaces, and $\F$ does not contain constant term because there does not exist
any holomorphic buildings of genus $g = 1$ without marked points or limit circles
whose $E_{\hat \omega}$-energy is zero.

It is easy to check that for any
$\G \in (\hbar^{-1} \D_X^{\leq 0})^{\star, \delta}
/ J^{\star, \delta}_{\overline{C}_0, \overline{C}_1, \overline{C}_2}$ of even degree
and any formal series $P(x) \in \R[[x]]$,
$P(\G) \in \D\D_X^{\leq 0, \delta}
/ \widetilde{J}^{\leq 0, \delta}_{\overline{C}_0, \overline{C}_1, \overline{C}_2}$ is well
defined.
Equation (\ref{irreducible decomposition}) implies that $\widetilde{\F} = e^{\F}$
in $\D\D_X^{\leq 0, \delta}
/ \widetilde{J}^{\leq 0, \delta}_{\overline{C}_0, \overline{C}_1, \overline{C}_2}$.

$\D\D_X$ has a structure of a left $D$-module
over $\W_{Y^-}$ as follows.
For each variable $p_{\hat c}$ ($c \in K^-$),
we define a differential operator on $\D\D_X$ by
\[
\overrightarrow{p_{\hat c}}
= \hbar \overrightarrow{\frac{\partial}{\partial q_{\hat c^\ast}}}.
\]
Then each
\[
f = \sum f_{(\hat c_i^\ast), (\hat c'_i)}(t, \hbar) q_{\hat c_1^\ast} q_{\hat c_2^\ast} \dots
q_{\hat c_k^\ast} p_{\hat c'_1} p_{\hat c'_2} \dots p_{\hat c'_l} \in \W_{Y^-}
\]
acts on $\D\D_X$ as a differential operator
\[
\overrightarrow{f}
= \sum f_{(\hat c_i^\ast), (\hat c'_i)}(\tilde t, \hbar)
q^-_{\hat c_1^\ast} q^-_{\hat c_2^\ast}
\dots q^-_{\hat c_k^\ast} \overrightarrow{p_{\hat c'_1}} \overrightarrow{p_{\hat c'_2}}
\dots \overrightarrow{p_{\hat c'_l}},
\]
where we replace each variable $t_x$ ($x \in K_{Y^-}^0$) with
$\tilde t_x = t_{(\mu_-)^{-1}(x)}$.
($\mu_-$ is the bijection defined in Section \ref{fiber prod for X}.)

Similarly, $\D\D_X$ has a structure of a right $D$-module
over $\W_{Y^+}$.
In this case, each variable $q_{\hat c^\ast}$ ($c \in K^+$)
defines a differential operator
\[
\overleftarrow{q_{\hat c^\ast}} = \hbar \overleftarrow{\frac{\partial}{\partial p_{\hat c}}}
\]
from right,
and each
\[
f = \sum f_{(\hat c_i^\ast), (\hat c'_i)}(t, \hbar) q_{\hat c_1^\ast} q_{\hat c_2^\ast} \dots
q_{\hat c_k^\ast} p_{\hat c'_1} p_{\hat c'_2} \dots p_{\hat c'_l} \in \W_{Y^+}
\]
acts on $\D\D_X$ as a differential operator
\[
\overleftarrow{f} = \sum f_{(\hat c_i^\ast), (\hat c'_i)}(\tilde t, \hbar)
\overleftarrow{q_{\hat c_1^\ast}}
\overleftarrow{q_{\hat c_2^\ast}} \dots \overleftarrow{q_{\hat c_k^\ast}}
p^+_{\hat c'_1} p^+_{\hat c'_2} \dots p^+_{\hat c'_l},
\]
where we replace each variable $t_x$ ($x \in K_{Y^+}^0$) with
$\tilde t_x = t_{(\mu_+)^{-1}(x)}$.

These $D$-module structures
\begin{align*}
\W_{Y^-} \times \D\D_X &\to \D\D_X,\\
\D\D_X \times \W_{Y^+} &\to \D\D_X
\end{align*}
induce the following maps:
{\belowdisplayskip=0pt
\begin{multline*}
\W_{Y^-}^{\leq \kappa_1}
/ I^{\leq \kappa_1}_{C_0, C_1 + \kappa_2 \delta^{-1}
+ \kappa_1 (\delta^{-1} - L_{\min}^{-1}), C_2 + \kappa_2} \times
\D\D_X^{\leq \kappa_2, \delta}
/ \widetilde{J}^{\leq \kappa_2, \delta}_{C_0, C_1 + \kappa_1 \delta^{-1}, C_2}\\
\to \D\D_X^{\leq \kappa_1 + \kappa_2, \delta}
/ \widetilde{J}^{\leq \kappa_1 + \kappa_2, \delta}_{C_0, C_1, C_2},
\end{multline*}
}
{\abovedisplayskip=0pt
\begin{multline*}
\D\D_X^{\leq \kappa_1, \delta}
/ \widetilde{J}^{\leq \kappa_1, \delta}_{C_0, C_1 + \kappa_2 \delta^{-1}, C_2 + \kappa_2}
\times \W_{Y^+}^{\leq \kappa_2} / I^{\leq \kappa_2}_{C_0, C_1 + \kappa_1 \delta^{-1}
+ \kappa_2 (\delta^{-1} - L_{\min}^{-1}), C_2}\\
\to \D\D_X^{\leq \kappa_1 + \kappa_2, \delta}
/ \widetilde{J}^{\leq \kappa_1 + \kappa_2, \delta}_{C_0, C_1, C_2}.
\end{multline*}
}

Assume that a generating function $\mathcal{H}_{Y^\pm} \in (\hbar^{-1} \W_Y^{\leq 0})^+
/ (\hbar^{-1} \W_Y^{\leq 0})^+_{\overline{C}_0, \overline{C}_1, \overline{C}_2}$
are defined and that $\overline{C}_0 \geq C_0$, $\overline{C}_1 \geq C_1
+ \kappa \delta^{-1}$ and $\overline{C}_2 \geq C_2 + \kappa$.
Then they define a linear map $\widehat{D}_X : \D\D_X^{\leq \kappa, \delta}
/ \widetilde{J}^{\kappa, \delta}_{C_0, C_1, C_2}
\to \D\D_X^{\leq \kappa, \delta}
/ \widetilde{J}^{\kappa, \delta}_{C_0, C_1, C_2}$ by
\[
\widehat{D}_X f = \delta f - \overrightarrow{\mathcal{H}_{Y^-}} f
+ (-1)^{|f|} f \overleftarrow{\mathcal{H}_{Y^+}}.
\]
Equations (\ref{main eq}) for $\mathcal{H}_{Y^\pm}$ imply that
$\widehat{D}_X$ is a differential of $\D\D_X^{\leq \kappa, \delta}
/ \widetilde{J}^{\kappa, \delta}_{C_0, C_1, C_2}$.

Equation (\ref{boundary formula for X}) implies $\widetilde{\F} = e^{\F}$ satisfies
\begin{equation}
\widehat{D}_X e^\F = 0 \label{main eq for X}
\end{equation}
in $\D\D_X^{\leq 0, \delta}
/ \widetilde{J}^{\leq 0, \delta}_{\overline{C}_0, \overline{C}_1, \overline{C}_2}$.

Define maps
\begin{multline*}
\W_{Y^-}^{\leq \kappa_1}
/ I^{\leq \kappa_1}_{C_0, C_1 + \kappa_2 \delta^{-1}
+ \kappa_1(\delta^{-1} - L_{\min}^{-1}), C_2 + \kappa_2} \times
\D_X^{\leq \kappa_2, \delta}
/ J^{\leq \kappa_2, \delta}_{C_0, C_1 + \kappa_1 \delta^{-1}, C_2 + \kappa_1}\\
\to \D_X^{\leq \kappa_1 + \kappa_2, \delta}
/ J^{\leq \kappa_1 + \kappa_2, \delta}_{C_0, C_1, C_2}
\end{multline*}
by
\[
(f, g) \mapsto f \underset{\F}{\overrightarrow{\ast}} g
= e^{- \F} \overrightarrow{f} (e^{\F} g)
\]
for $\kappa_1, \kappa_2, C_0, C_1, C_2$ such that
$\overline{C}_0 \geq C_0$, $\overline{C}_1 \geq C_1
+ (\kappa_1 + \kappa_2) \delta^{-1}$, $\overline{C}_2 \geq C_2 + \kappa_1 + \kappa_2$.
This family of maps defines a left module-like structure, that is,
the associativity law is satisfied if it is well defined.

Similarly, we define maps 
\begin{multline*}
\D_X^{\leq \kappa_1, \delta}
/ J^{\leq \kappa_1, \delta}_{C_0, C_1 + \kappa_2 \delta^{-1}, C_2 + \kappa_2} \times
\W_{Y^+}^{\leq \kappa_2}
/ I^{\leq \kappa_2}_{C_0, C_1 + \kappa_1 \delta^{-1}
+ \kappa_2(\delta^{-1} - L_{\min}^{-1}), C_2 + \kappa_1}
\\
\to \D_X^{\leq \kappa_1 + \kappa_2, \delta}
/ J^{\leq \kappa_1 + \kappa_2, \delta}_{C_0, C_1, C_2}
\end{multline*}
by
\[
(g, f) \mapsto g \underset{\F}{\overleftarrow{\ast}} f
= (g e^{\F}) \overleftarrow{f} e^{-\F}
\]
for $\kappa_1, \kappa_2, C_0, C_1, C_2$ such that
$\overline{C}_0 \geq C_0$, $\overline{C}_1 \geq C_1
+ (\kappa_1 + \kappa_2) \delta^{-1}$, $\overline{C}_2 \geq C_2 + \kappa_1 + \kappa_2$.
This family of maps defines a right module-like structure.
Note that these module-like structures are a bimodule structure, that is,
\[
(f \underset{\F}{\overrightarrow{\ast}} g) \underset{\F}{\overleftarrow{\ast}} h
= f \underset{\F}{\overrightarrow{\ast}} (g \underset{\F}{\overleftarrow{\ast}} h)
\]
for all $f \in \W_{Y^-}$, $g \in \D_X$ and $h \in \W_{Y^+}$.

Define a linear map $D_{\F} : \D_X^{\leq \kappa} / J^{\leq \kappa, \delta}_{C_0, C_1, C_2}
\to \D_X^{\leq \kappa} / J^{\leq \kappa, \delta}_{C_0, C_1, C_2}$ by
\begin{align*}
D_{\F} f &= e^{-\F} [\widehat{D}_X, f] (e^{\F}) \\
&= e^{-\F}\widehat{D}_X(f e^{\F}). &(\text{by } (\ref{main eq for X}))
\end{align*}
Then it satisfies the following:
\begin{itemize}
\item
$D_{\F}$ is a differential, that is, $D_{\F}^2 = 0$.
\item
For any $f \in \W_{Y^-}^{\leq \kappa_1} / I^{\leq \kappa_1}_{C_0, C'_1, C_2 + \kappa_2}$
and $g \in \D_X^{\leq \kappa_2} / J^{\leq \kappa_2, \delta}_{C_0, C''_1, C_2}$,
\begin{equation}
D_{\F}(f \underset{\F}{\overrightarrow{\ast}} g)
= (D_{Y^-}f) \underset{\F}{\overrightarrow{\ast}} g
+ (-1)^{|f|} f \underset{\F}{\overrightarrow{\ast}} D_{\F} (g)
\label{D_F left Leibnitz}
\end{equation}
in $\D_X^{\leq \kappa_1 + \kappa_2} / J^{\leq \kappa_1 + \kappa_2, \delta}
_{C_0, C_1, C_2}$, where $C'_1 = C_1 + \kappa_1(\delta^{-1} - L_{\min}^{-1})
+ \kappa_2 \delta^{-1}$ and $C''_1 = C_1 + \kappa_1 \delta^{-1}$.
\item
For any $g \in \D_X^{\leq \kappa_2} / J^{\leq \kappa_2, \delta}_{C_0, C''_1,
C_2 + \kappa_1}$ and
$f \in \W_{Y^+}^{\leq \kappa_1} / I^{\leq \kappa_1}_{C_0, C'_1, C_2}$,
\begin{equation}
D_{\F}(g \underset{\F}{\overleftarrow{\ast}} f)
= D_{\F}(g) \underset{\F}{\overleftarrow{\ast}} f
+ (-1)^{|g|} g \underset{\F}{\overleftarrow{\ast}} (D_{Y^+} f)
\label{D_F right Leibnitz}
\end{equation}
in $\D_X^{\leq \kappa_1 + \kappa_2} / J^{\leq \kappa_1 + \kappa_2, \delta}
_{C_0, C_1, C_2}$.
\end{itemize}
They imply that the family of cohomology groups
$H^\ast(\D_X^{\leq \kappa} / J^{\leq \kappa, \delta}_{C_0, C_1, C_2}, D_{\F})$ has a
$(H^\ast(\W_{Y^-}^{\leq \kappa} / I^{\leq \kappa}_{C_0, C_1, C_2}, D_{Y^-}),
H^\ast(\W_{Y^+}^{\leq \kappa} / I^{\leq \kappa}_{C_0, C_1, C_2}, D_{Y^+}))$-bimodule-like
structure.
Sometimes we denote the linear maps $D_{\F}$ for
the generating function
$\F = \F_{(X, \omega, Y^\pm, \lambda^\pm, K_{Y^\pm}, K^0_X,
K^0_{Y^\pm}, \mu^\pm, K^2_{Y^\pm}, J, \B_X)}$ by
\[
D_{(X, \omega, Y^\pm, \lambda^\pm, K_{Y^\pm}, K^0_X,
K^0_{Y^\pm}, \mu^\pm, K^2_{Y^\pm}, J, \B_X)}.
\]

By definition, $D_{\F} 1 = 0$.
Therefore the linear maps
\begin{align*}
i_\F^-(f) &= f \underset{\F}{\overrightarrow{\ast}} 1 = e^{-\F} \overrightarrow{f} e^{\F}
: \W_{Y^-}^{\leq \kappa}
/ I^{\leq \kappa}_{C_0, C_1 + \kappa(\delta^{-1} - L_{\min}^{-1}), C_2}
\to \D_X^{\leq \kappa} / J^{\leq \kappa, \delta}_{C_0, C_1, C_2},\\
i_\F^+(f) &= 1 \underset{\F}{\overleftarrow{\ast}} f = e^{\F} \overleftarrow{f} e^{-\F}
: \W_{Y^+}^{\leq \kappa}
/ I^{\leq \kappa}_{C_0, C_1 + \kappa(\delta^{-1} - L_{\min}^{-1}), C_2}
\to \D_X^{\leq \kappa} / J^{\leq \kappa, \delta}_{C_0, C_1, C_2}
\end{align*}
induce homomorphisms
\[
i_\F^\pm : H^\ast(\W_{Y^\pm}^{\leq \kappa}
/ I^{\leq \kappa}_{C_0, C_1 + \kappa(\delta^{-1} - L_{\min}^{-1}), C_2}, D_{Y^\pm})
\to H^\ast(\D_X^{\leq \kappa} / J^{\leq \kappa, \delta}_{C_0, C_1, C_2}, \D_\F).
\]
This pair of homomorphisms $i_\F^\pm$ is the chain map in the sense of general SFT.

Next we consider rational SFT.
Define $\mathcal{L}_X = \mathcal{L}_{(X, \omega, Y^\pm, \lambda^\pm,
K_{Y^\pm}, \overline{K}_X^0)} = \D_X|_{\hbar = 0}$ as a quotient super-commutative
algebra of $\D_X$.
We also use a bigger super-commutative algebra $\widehat{\mathcal{L}}_X$.
Its elements are formal series
\[
\sum
f_{(\hat c_i), (\hat c'_i), (\hat c''_i), (\hat c'''_i), e}(t)
q_{\hat c_1^\ast}^- \dots q_{\hat c_{k^-_q}^\ast}^-
q_{(\hat c'_1)^\ast}^+ \dots q_{(\hat c'_{k^+_q})^\ast}^+
p_{\hat c''_1}^- \dots p_{\hat c''_{k^-_p}}^-
p_{\hat c'''_1}^+ \dots p_{\hat c'''_{k^+_p}}^+ T^e,
\]
where each $f_{(\hat c_i), (\hat c'_i), (\hat c''_i), (\hat c'''_i), e}(t) \in \R[[t]]$ is a formal
series of the variables $t_x$ ($x \in K^0_X$) and
the infinite sum is taken over all sequences
$((\hat c_i), (\hat c'_i), (\hat c''_i), (\hat c'''_i), e)$
consisting of the simplices $(\hat c_i)$ of $K_{Y^-}$ not contained in
$\overline{P}^{\text{bad}}_{Y^-}$,
$(\hat c'_i)$ of $K_{Y^+}$ not contained in $\overline{P}^{\text{bad}}_{Y^+}$,
$(\hat c''_i)$ of $K_{Y^-}$ not contained in $\overline{P}^{\text{bad}}_{Y^-}$,
$(\hat c'''_i)$ of $K_{Y^+}$ not contained in $\overline{P}^{\text{bad}}_{Y^+}$ and
$e \in \tilde \omega H_2(\overline{X}, \partial \overline{X})$.
We impose the following Novikov condition on the infinite sum:
for any $C > 0$, the number of the non-zero terms with
$\sum_j e(p^-_{\hat c''_j}) + \sum_j e(p^+_{\hat c'''_j}) \geq -C$
and $e + \sum_j e(p^-_{\hat c''_j}) + \sum_j e(p^+_{\hat c'''_j}) \geq -C$ is finite.

The Poisson structure of $\widehat{\mathcal{L}}_X$ is defined by
\begin{align*}
\{f, g\} &= \sum_{c \in K_{Y^-}} \biggl(\frac{\overleftarrow{\partial} f}{\partial p_{\hat c}^-}
\frac{\overrightarrow{\partial} g}{\partial q_{\hat c^\ast}^-}
- (-1)^{|f| |g|} \frac{\overleftarrow{\partial} g}{\partial p_{\hat c}^-}
\frac{\overrightarrow{\partial} f}{\partial q_{\hat c^\ast}^-}\biggr)\\
&\quad - \sum_{c \in K_{Y^+}} \biggl(\frac{\overleftarrow{\partial} f}{\partial p_{\hat c}^+}
\frac{\overrightarrow{\partial} g}{\partial q_{\hat c^\ast}^+}
- (-1)^{|f| |g|} \frac{\overleftarrow{\partial} g}{\partial p_{\hat c}^+}
\frac{\overrightarrow{\partial} f}{\partial q_{\hat c^\ast}^+}\biggr).
\end{align*}
We regard $\mathcal{P}_{Y^-}$ and $\mathcal{P}_{Y^+}$ as subspaces of
$\widehat{\mathcal{L}}_X$ by
$q_{\hat c^\ast} \mapsto q^-_{\hat c^\ast}$, $p_{\hat c} \mapsto p^-_{\hat c}$
and $q_{\hat c^\ast} \mapsto q^+_{\hat c^\ast}$, $p_{\hat c} \mapsto p^+_{\hat c}$
respectively.
Then the inclusions $\mathcal{P}_{Y^-} \inj \widehat{\mathcal{L}}_X$ and
$\mathcal{P}_{Y^+} \inj \widehat{\mathcal{L}}_X$
are a Poisson map and an anti-Poisson map respectively.

For each even element $g \in \mathcal{L}_X$,
define a map
\[
f \mapsto f|_g : \widehat{\mathcal{L}}_X \to \mathcal{L}_X
\]
by the evaluation map given by $p_{\hat c}^-
= \frac{\overrightarrow{\partial} g}{\partial q_{\hat c^\ast}^-}$ and
$q_{\hat c^\ast}^+ = \frac{\overleftarrow{\partial} g}{\partial p_{\hat c}^+}$.

For each $\kappa \geq 0$, we define submodules
$\mathcal{L}_X^{\leq \kappa} \subset \mathcal{L}_X$
and $\widehat{\mathcal{L}}_X^{\leq \kappa} \subset \widehat{\mathcal{L}}_X$
by the conditions
$\sum_i e(q_{\hat c_i^\ast}^-) + e + \sum_i e(p_{\hat c'_i}^+) \leq \kappa$
and $\sum_i e(q_{\hat c_i^\ast}^-) + \sum_i e(q_{(\hat c'_i)^\ast}^+) + e
+ \sum_i e(p_{\hat c''_i}^-) + \sum_i e(p_{\hat c''_i}^+) \leq \kappa$ respectively.
Define submodules
$J^{\leq \kappa}_{C_0, C_2} \subset \mathcal{L}_X^{\leq \kappa}$
and $\widetilde{J}^{\leq \kappa}_{C_0, C_2} \subset
\widehat{\mathcal{L}}_X^{\leq \kappa}$ by
\begin{align*}
J^{\leq \kappa}_{C_0, C_2} &= \Bigl\{\sum
a_{(x_i), (\hat c_i^\ast), (\hat c'_i), e} t_{x_1} \dots t_{x_{k_t}} q^-_{\hat c_1^\ast}
\dots q^-_{\hat c_{k_q}^\ast} p^+_{\hat c'_1} \dots p^+_{\hat c'_{k_p}} T^e
\in \mathcal{L}_X^{\leq \kappa};\\
&\quad \quad a_{(x_i), (\hat c_i^\ast), (\hat c'_i), e} = 0 \text{ for all }
((x_i)_{i = 1}^{k_t}, (\hat c_i^\ast),_{i = 1}^{k_q} (\hat c'_i)_{i = 1}^{k_p}, e)
\text{ such that}\\
&\quad \quad
k_t \leq C_0,\ \sum e(p^+_{\hat c'_i}) \geq - C_2
\text{ and } e + \sum e(p^+_{\hat c'_i}) \geq - C_2 \Bigr\}
\end{align*}
and
\begin{align*}
\widetilde{J}^{\leq \kappa}_{C_0, C_2} &= \Bigl\{\sum
a_{(x_i), (\hat c_i), (\hat c'_i), (\hat c''_i), (\hat c'''_i), e} t_{x_1} \dots t_{x_{k_t}}
q^-_{\hat c_1^\ast} \dots q^-_{\hat c_{k_q}^\ast}
q^+_{(\hat c'_1)^\ast} \dots q^+_{(\hat c'_{k_q})^\ast} \\
&\quad \hph{\Bigl\{\sum
a_{(x_i), (\hat c_i), (\hat c'_i), (\hat c''_i), (\hat c'''_i), e} t_{x_1} \dots t_{x_{k_t}}}
p^-_{\hat c''_1} \dots p^-_{\hat c''_{k_p}}
p^+_{\hat c'''_1} \dots p^+_{\hat c'''_{k_p}} T^e
\in \widehat{\mathcal{L}}_X^{\leq \kappa};\\
&\quad \quad a_{\alpha} = 0 \text{ for all }
\alpha = ((x_i)_{i = 1}^{k_t}, (\hat c_i)_{i = 1}^{k^-_q}, (\hat c'_i)_{i = 1}^{k^+_q},
(\hat c''_i)_{i = 1}^{k^-_p}, (\hat c'''_i)_{i = 1}^{k^+_p}, e)\\
&\quad \quad \text{such that }
k_t \leq C_0,\ \sum e(p^+_{\hat c''_i}) + \sum e(p^+_{\hat c'''_i}) \geq - C_2
\text{ and}\\
&\quad \quad \hph{\text{such that }
k_t \leq C_0,\ }
e + \sum e(p^+_{\hat c''_i}) + \sum e(p^+_{\hat c'''_i}) \geq - C_2 \Bigr\}.
\end{align*}
%
%
%
%
%
%

First we note that $\mathbf{h} = \mathcal{H}_{Y^-, 0} - \mathcal{H}_{Y^+, 0}
\in \mathcal{L}_X^{\leq 0} / J^{\leq 0}_{\overline{C}_0, \overline{C}_2}$ satisfies
\begin{equation}
\delta \mathbf{h} - \frac{1}{2} \{\mathbf{h}, \mathbf{h}\} = 0.
\label{eq for h}
\end{equation}
For each triple $(\kappa, C_0, C_2)$ such that $\overline{C}_0 \geq C_0$ and
$\overline{C}_2 \geq C_2 + \kappa$,
we define a linear map
$\widehat{d}_X : \widehat{\mathcal{L}}^{\leq \kappa}_X
/ \widetilde{J}^{\leq \kappa}_{C_0, C_2}
\to \widehat{\mathcal{L}}^{\leq \kappa}_X
/ \widetilde{J}^{\leq \kappa}_{C_0, C_2}$ by
\[
\widehat{d}_X f = \delta f - \{\mathbf{h}, f\}.
\]
Then (\ref{eq for h}) implies that $\widehat{d}_X^2 = 0$.
$\widehat{d}_X$ also satisfies
\begin{align}
\widehat{d}_X(fg) = (\widehat{d}_X f) g + (-1)^{|f|} f \widehat{d}_X g
\label{hat d_X Leibnitz}\\
\widehat{d}_X\{f, g\} = \{\widehat{d}_X f, g\} + (-1)^{|f|} \{f, \widehat{d}_X g\}
\label{hat d_X bracket}
\end{align}
if the multiplications or Poisson brackets are well defined.

We use the genus zero part $\F_0 \in \mathcal{L}_X^{\leq 0}
/ J^{\leq 0}_{\overline{C}_0, \overline{C}_2}$ of the generating function.
Equation (\ref{main eq for X}) implies that
\begin{equation}
\delta \F_0 - \mathbf{h}|_{\F_0} = 0
\label{main eq for X rational}
\end{equation}
in $\mathcal{L}_X^{\leq 0} / J^{\leq 0}_{\overline{C}_0, \overline{C}_2}$.

For each triple $(\kappa, C_0, C_2)$ such that $\overline{C}_0 \geq C_0$ and
$\overline{C}_2 \geq C_2 + \kappa$,
define linear maps
$d_{\F_0} : \mathcal{L}^{\leq \kappa}_X / J^{\leq \kappa}_{C_0, C_2}
\to \mathcal{L}^{\leq \kappa}_X / J^{\leq \kappa}_{C_0, C_2}$
and $i_{\F_0}^\pm : \mathcal{P}^{\leq \kappa}_{Y^\pm} / I^{\leq \kappa}_{C_0, C_2}
\to \mathcal{L}^{\leq \kappa}_X / J^{\leq \kappa}_{C_0, C_2}$ by
\begin{align*}
d_{\F_0} f &= (\widehat{d}_X f)|_{\F_0} \\
&= \delta f - \{\mathbf{h}, f\}|_{\F_0}
\ (= (D_{\F} f)|_{\hbar = 0})
\end{align*}
and
\[
i_{\F_0}^\pm (f) = f|_{\F_0}.
\]
We claim that (\ref{main eq for X rational}) and (\ref{eq for h}) imply
that $d_{\F_0}$ is a differential (i.e. $d_{\F_0}^2 = 0$) and $i_{\F_0}^\pm$
are chain maps.
For its proof,
it is convenient to introduce a linear map
\[
f \mapsto \widetilde{f} : \widehat{\mathcal{L}}_X \to \widehat{\mathcal{L}}_X
\]
defined by
\begin{align*}
\widetilde{f} &= (k^-_p + k^+_q - 1) f \\
&= \bigl(\sum p^-_{\hat c} \overrightarrow{\partial}_{p^-_{\hat c}}
+ \sum q^+_{\hat c^\ast} \overrightarrow{\partial}_{q^+_{\hat c^\ast}} - 1\bigr) f
\end{align*}
for each monomial
\[
f = t_{x_1} \dots t_{x_{k_t}} q_{\hat c_1^\ast}^- \dots q_{\hat c_{k^-_q}^\ast}^-
q_{(\hat c'_1)^\ast}^+ \dots q_{(\hat c'_{k^+_q})^\ast}^+
p_{\hat c''_1}^- \dots p_{\hat c''_{k^-_p}}^-
p_{\hat c'''_1}^+ \dots p_{\hat c'''_{k^+_p}}^+ T^e.
\]
\begin{lem}
\begin{enumerate}[label=\normalfont (\roman*)]
\item
For any $f, g \in \widehat{\mathcal{L}}_X$,
\begin{equation}
\widetilde{\{f, g\}} = \{ \widetilde{f}, g \} + \{f, \widetilde{g}\}.
\label{tilde{}}
\end{equation}
\item
For all $g \in \mathcal{L}_X^{\text{even}}$ and
$f \in \widehat{\mathcal{L}}_X$,
\begin{equation}
\{f, g\}|_g = \widetilde{f}|_g + f|_g
\label{{}ev}
\end{equation}
and
\begin{equation}
\delta(f|_g) = (\delta f)|_g - \{\delta g, f\}|_g.
\label{delta ev}
\end{equation}
\item
For all $g \in \mathcal{L}_X^{\text{even}}$ and $f, h \in \widehat{\mathcal{L}}_X$,
\begin{equation}
\{h, f|_g\}|_g = \{h, \{f, g\}\}|_g - \{h, \widetilde{f}\}|_g.
\label{{}evev}
\end{equation}
In particular,
\begin{equation}
\{h|_g, f\}|_g + \{h, f|_g\}|_g = \{h, f\}|_g.
\label{evev}
\end{equation}
\end{enumerate}
\end{lem}
\begin{proof}
(\ref{tilde{}}) and (\ref{{}ev}) are easy.
(\ref{delta ev}) is proved as follows.
First note that if we regard each side as an operator $A$ for $f$ then it satisfies
\[
A(f_1 f_2) = A(f_1) f_2|_g + (-1)^{|f_1|} f_1|_g A(f_2).
\]
Hence we may assume that $f$ is some variable $q^-$, $q^+$, $p^-$ or $p^+$.

If $f$ is $q^-$ or $p^+$, then it satisfies (\ref{delta ev}) since
$\delta(f|_g) = (\delta f)|_g = \delta f$ and $\{\delta g, f\}|_g = 0$.

Next we consider the case of $f = q^+_{\hat c^\ast}$.
Define $a_{c', c} \in \Q$ by $\partial \hat c' = \sum_c a_{c', c} \hat c$.
Then (\ref{delta ev}) is equivalent to
\[
\delta\biggl(\frac{\overleftarrow{\partial} g}{\partial p_{\hat c}^+}\biggr)
- \frac{\overleftarrow{\partial} (\delta g)}{\partial p_{\hat c}^+}
= (-1)^{|g|} \sum_{c'} a_{c', c} \frac{\overleftarrow{\partial} g}{\partial p_{\hat c'}^+}.
\]
($(-1)^{|g|} = 1$ for $g \in \mathcal{L}_X^{\text{even}}$.)
We prove that this equation holds for all $g \in \widehat{\mathcal{L}}_X$ as follows.
If we regard each side as an operator $B$ for $g$ then it satisfies
\[
B(g_1 g_2) = (-1)^{|\hat c| |g_2|} B(g_1) g_2 + (-1)^{|g_1|} g_1 B(g_2).
\]
Hence it is enough to prove the equation for the case where $g$ is
some variable $q^-$, $q^+$, $p^-$ or $p^+$ and it can be easily checked.

Finally, if $f = p^-_{\hat c}$ then (\ref{delta ev}) is equivalent to
\[
\delta\biggl(\frac{\overrightarrow{\partial} g}{\partial q_{\hat c^\ast}^-}\biggr)
- (-1)^{|\hat c|} \frac{\overrightarrow{\partial} (\delta g)}{\partial q_{\hat c^\ast}^-}
= (-1)^{1 + |\hat c|} \sum_{c'} a_{c, c'}
\frac{\overrightarrow{\partial} g}{\partial q_{(\hat c')^\ast}^-},
\]
and it can be proved similarly.

(\ref{{}evev}) is proved as follows.
If we regard each side as an operator $C$ for $f$ then it satisfies
\[
C(f_1 f_2) = C(f_1) f_2|_g + (-1)^{|h| |f_1|} f_1|_g C(f_2).
\]
Hence it is enough to prove (\ref{{}evev}) for the case where $f$ is
some variable $q^-$, $q^+$, $p^-$ or $p^+$ and it can be easily checked.
(\ref{evev}) is a corollary of (\ref{{}evev}).
\end{proof}

Now we prove the following proposition.
The fourth claim is used to define a Poisson structure of rational SFT cohomology of
$(Y, \lambda)$ in Section \ref{independence}.
\begin{prop}\label{properties of d_{F_0}}
\begin{enumerate}[label=\normalfont (\roman*)]
\item
$d_{\F_0}$ is a differential, that is, $d_{\F_0}^2 = 0$.
\item
$f \mapsto f|_{\F_0}
: (\widehat{\mathcal{L}}_X^{\leq \kappa} / \widetilde{J}^{\leq \kappa}_{C_0, C_2},
\widehat{d}_X)
\to (\mathcal{L}_X^{\leq \kappa} / J^{\leq \kappa}_{C_0, C_2}, d_{\F_0})$ is a chain map,
that is, $d_{\F_0}(f|_{\F_0}) = (\widehat{d}_X f)|_{\F_0}$ for all
$f \in \widehat{\mathcal{L}}_X^{\leq \kappa} / \widetilde{J}^{\leq \kappa}_{C_0, C_2}$.
In particular, $i_{\F_0}^\pm$ are chain maps, that is, $d_{\F_0} \circ i_{\F_0}^\pm
= i_{\F_0}^\pm \circ d_{Y^\pm}$.
\item
For any $f \in \mathcal{P}^{\leq \kappa_1}_{Y^\pm}
/ I^{\leq \kappa_1}_{C_0, C_2 + \kappa_2}$
and $g \in \mathcal{L}^{\leq \kappa_2}_X / J^{\leq \kappa_2}_{C_0, C_2}$,
\[
d_{\F_0}(i_{\F_0}^\pm(f) g)
= i_{\F_0}^\pm(d_{Y^\pm} f) g + (-1)^{|f|} i_{\F_0}^\pm(f) d_{\F_0} g
\]
in $\mathcal{L}^{\leq \kappa_1 + \kappa_2}_X / J^{\leq \kappa_1 + \kappa_2}_{C_0, C_2}$.
\item
Assume that $f \in \widehat{\mathcal{L}}_X^{\leq \kappa_1}
/ \widetilde{J}^{\leq \kappa_1}_{C_0, C_2 + \kappa_2}$,
$g \in \widehat{\mathcal{L}}_X^{\leq \kappa_2}
/ \widetilde{J}^{\leq \kappa_2}_{C_0, C_2 + \kappa_1}$,
$a \in \mathcal{L}_X^{\leq \kappa_1} / J^{\leq \kappa_1}_{C_0, C_2 + \kappa_2}$
and $b \in \mathcal{L}_X^{\leq \kappa_2} / J^{\leq \kappa_2}_{C_0, C_2 + \kappa_1}$
satisfy $f|_{\F_0} = d_{\F_0} a$ and $g|_{\F_0} = d_{\F_0} b$.
Then
\begin{align}
&\{f, g\}|_{\F_0}
+ (-1)^{|f|} \bigl(\{\widehat{d}_X f, b\}|_{\F_0}
- \{a, \widehat{d}_X g\}|_{\F_0}\bigr)\notag\\
&= d_{\F_0} \bigl( \{a, g\}|_{\F_0} + (-1)^{|f|} \{f, b\}|_{\F_0}
+ \{a, \{\mathbf{h}, b\}\}|_{\F_0} \bigr) \label{eq for rational independence}
\end{align}
in $\mathcal{L}_X^{\leq \kappa_1 + \kappa_2}
/ J^{\leq \kappa_1 + \kappa_2}_{C_0, C_2}$.
In particular, if in addition $\widehat{d}_X f = 0$ and $\widehat{d}_X g = 0$, then
$\{f, g\}|_{\F_0}$ is exact.
\end{enumerate}
\end{prop}
\begin{proof}
First note that (\ref{main eq for X rational}), (\ref{delta ev}) and (\ref{evev}) imply that
for any $f \in \widehat{\mathcal{L}}_X^{\leq \kappa}
/ \widetilde{J}^{\leq \kappa}_{C_0, C_2}$,
\begin{align}
d_{\F_0}(f|_{\F_0})
&= (\delta f)|_{\F_0} - \{ \delta \mathcal{F}_0, f\}|_{\F_0}
-\{\mathbf{h}, f|_{\F_0}\}|_{\F_0} \notag\\
&= (\delta f)|_{\F_0} - \{ \mathbf{h}|_{\F_0}, f\}|_{\F_0}
-\{\mathbf{h}, f|_{\F_0}\}|_{\F_0} \notag\\
&= (\delta f)|_{\F_0} - \{\mathbf{h}, f\}|_{\F_0} \notag\\
&= (\widehat{d}_X f)|_{\F_0}.
\label{d ev}
\end{align}
(i) is because (\ref{d ev}) implies
\[
d_{\F_0}^2 f = d_{\F_0} ((\widehat{d}_X f)|_{\F_0}) = (\widehat{d}_X^2 f)|_{\F_0} = 0.
\]
(ii) is due to (\ref{d ev}).
(iii) is because for any $f \in \widehat{\mathcal{L}}_X^{\leq \kappa_1}
/ \widetilde{J}^{\leq \kappa_1}_{C_0, C_2 + \kappa_1}$
and $g \in \mathcal{L}^{\leq \kappa_2}_X / J^{\leq \kappa_2}_{C_0, C_2}$,
\begin{align*}
d_{\F_0}(f|_{\F_0} \cdot g)
&= d_{\F_0}((f g)|_{\F_0}) & (\text{since } g = g|_{\F_0})\\
&= (\widehat{d}_X (fg))|_{\F_0} & (\text{by } (\ref{d ev}))\\
&= (\widehat{d}_X f)|_{\F_0} \cdot g + (-1)^{|f|} f|_{\F_0} \cdot d_{\F_0} g
& (\text{by } (\ref{hat d_X Leibnitz}))
\end{align*}
in $\mathcal{L}^{\leq \kappa_1 + \kappa_2}_X / J^{\leq \kappa_1 + \kappa_2}_{C_0, C_2}$.

(iv) is because
\begin{align*}
d_{\F_0}(\{a, g\}|_{\F_0})
&= (\widehat{d}_X \{a, g\})|_{\F_0} \\
&= \{\widehat{d}_X a, g\}|_{\F_0} + (-1)^{|a|}\{a, \widehat{d}_X g\}|_{\F_0} \\
&= \{(\widehat{d}_X a)|_{\F_0}, g\}|_{\F_0}
+ \{\widehat{d}_X a, g|_{\F_0}\}|_{\F_0}
+ (-1)^{|a|}\{a, \widehat{d}_X g\}|_{\F_0} \\
&= \{f|_{\F_0}, g\}|_{\F_0}
+ \{\widehat{d}_X a, (\widehat{d}_X b)|_{\F_0}\}|_{\F_0}
+ (-1)^{|a|}\{a, \widehat{d}_X g\}|_{\F_0},
\end{align*}
\begin{align*}
(-1)^{|f|} d_{\F_0}(\{f, b\}|_{\F_0})
&= (-1)^{|f|} (\widehat{d}_X \{f, b\})|_{\F_0} \\
&= \{f, \widehat{d}_X b\}|_{\F_0} + (-1)^{|f|} \{\widehat{d}_X f, b\}|_{\F_0} \\
&= \{f, (\widehat{d}_X b)|_{\F_0}\}|_{\F_0} + \{f|_{\F_0}, \widehat{d}_X b\}|_{\F_0}
+ (-1)^{|f|} \{\widehat{d}_X f, b\}|_{\F_0} \\
&= \{f, g|_{\F_0}\}|_{\F_0} + \{(\widehat{d}_X a)|_{\F_0}, \widehat{d}_X b\}|_{\F_0}
+ (-1)^{|f|} \{\widehat{d}_X f, b\}|_{\F_0}
\end{align*}
and
\begin{align*}
d_{\F_0}(\{a, \{\mathbf{h}, b\}\}|_{\F_0})
&= - d_{\F_0}(\{a, \widehat{d}_X b\}|_{\F_0}) \\
&= - (\widehat{d}_X \{a, \widehat{d}_X b\})|_{\F_0} \\
&= - \{\widehat{d}_X a, \widehat{d}_X b\}|_{\F_0} \\
&= - \{(\widehat{d}_X a)|_{\F_0}, \widehat{d}_X b\}|_{\F_0}
- \{\widehat{d}_X a, (\widehat{d}_X b)|_{\F_0}\}|_{\F_0}.
\end{align*}
\end{proof}

Finally we consider the case of contact homology.
Assume that $(X, \omega)$ is an exact cobordism, that is,
$\omega = d \theta$ for some 1-form $\theta$ on $X$ such that
\[
\theta|_{(-\infty, 0] \times Y^-} = e^\sigma \lambda^-\quad \text{and} \quad
\theta|_{[0, \infty) \times Y^+} = e^\sigma \lambda^+.
\]
Further we assume that the domains of $\mu_\pm$ are the whole of $K_X^0$ and
$\mu_\pm : K_X^0 \to K_{Y^\pm}^0$ are bijections.
Define
\[
\widehat{\mathcal{F}}_0 = \sum_{c}
\frac{\overleftarrow{\partial} \mathcal{F}_0}{\partial p^+_{\hat c}} \biggr|_{p^+ = 0}
\cdot p^+_{\hat c} \in \mathcal{L}^{\leq 0}_X / J^{\leq 0}_{\overline{C}_0, \overline{C}_2}.
\]

Exactness of $(X, \omega)$ implies
\[
\frac{\overrightarrow{\partial}}{\partial q_{\hat c^\ast}^-} (\mathcal{F}_0|_{p = 0}) = 0.
\]
Hence equation (\ref{main eq for X rational}) implies
\begin{equation}
\delta \widehat{\mathcal{F}}_0 - \widehat{\mathbf{h}}|_{\widehat{\mathcal{F}}_0} = 0
\label{main eq for X contact}
\end{equation}
in $\mathcal{L}^{\leq 0}_X / J^{\leq 0}_{\overline{C}_0, \overline{C}_2}$,
where $\widehat{\mathbf{h}}
= \widehat{\mathcal{H}}_{Y^-, 0} - \widehat{\mathcal{H}}_{Y^+, 0}$.
For each pair $(\kappa, C_0)$ such that $\overline{C}_0 \geq C_0$ and
$\overline{C}_2 \geq \kappa$, we define a homomorphism
$\Psi_{\widehat{\mathcal{F}}_0} :
\A_{Y^+}^{\leq \kappa} / I^{\leq \kappa}_{C_0}
\to \A_{Y^-}^{\leq \kappa} / I^{\leq \kappa}_{C_0}$
by the evaluation
\[
\Psi_{\widehat{\mathcal{F}}_0}(f) = f|_{\widehat{\mathcal{F}}_0}
= f\Big|_{q^+_{\hat c^\ast}
= \bigl(\frac{\overleftarrow{\partial} \F_0}{\partial p^+_{\hat c}}\bigr|_{p^+ = 0}\bigr)}.
\]
Then (\ref{main eq for X contact}) implies that this is a chain map, that is,
$\partial_{Y^-} \circ \Psi_{\widehat{\mathcal{F}}_0}
= \Psi_{\widehat{\mathcal{F}}_0} \circ \partial_{Y^+}$,
where we identify each $t_x$ ($x \in K_{Y^+}^0$) with
$t_{\mu_- \circ \mu_+^{-1}(x)}$.
Therefore it induces a homomorphism $(\Psi_{\widehat{\mathcal{F}}_0})_\ast :
H^\ast (\A_{Y^+}^{\leq \kappa} / I^{\leq \kappa}_{C_0}, \partial_{Y^+})
\to H^\ast(\A_{Y^-}^{\leq \kappa} / I^{\leq \kappa}_{C_0}, \partial_{Y^-})$.

\subsection{Algebras with further energy conditions}
\label{algebras with further energy conditions}
Assume that $Z$ contains contact manifolds $(Y_i, \lambda_i)$ ($1 \leq i \leq m$) and
that for each $i = 1, 2, \dots, m$, there is a pair of symplectic cobordisms
$Z_i^-$ and $Z_i^+$ such that $Z = Z_i^- \cup_{Y_i} Z_i^+$.
We assume that the pull back of the symplectic form $\omega$ to $Y_i$ is $d\lambda_i$.
Then we can construct the algebras which respect these decompositions as follows.
(We need these algebras for the definition of the composition of generating
functions in Section \ref{composition of generating functions}.)

Let $((-\epsilon, \epsilon) \times Y_i, d(e^\sigma \lambda_i)) \inj (Z, \omega)$
be a neighborhood of each $Y_i$ and define a closed two form $\widetilde{\omega}_{Y_i}$
on $X$ by
$\widetilde{\omega}_{Y_i} = \omega$ on $Z_i^+$,
$\widetilde{\omega}_{Y_i} = d(\varphi \lambda^+)$ on $[0, \infty) \times Y^+$,
$\widetilde{\omega}_{Y_i} = d(\varphi \lambda_i)$ on $(-\epsilon, 0] \times Y_i$,
and $\widetilde{\omega}_{Y_i} = 0$
on $(-\infty, 0] \times Y^- \cup (Z_i^- \setminus (-\epsilon, 0] \times Y_i)$,
where $\varphi : \R \to \R_{\geq 0}$ is a smooth function with compact support
such that $\varphi(0) = 1$ and $\varphi|_{(-\infty, -\epsilon]} \equiv 0$.

For a holomorphic building $(\Sigma, z, u) \in \widehat{\M}(X, \omega, J)$,
define $e = \int u^\ast \widetilde{\omega}$ and
$e_{Y_i} = \int u^\ast \widetilde{\omega}_{Y_i}$.
Then these satisfy
\[
e + \sum_{+\infty\text{-limit circles}} L_{\gamma_{+\infty_j}} \geq
\sum_{-\infty\text{-limit circles}} L_{\gamma_{-\infty_j}}
\]
and
\[
e_{Y_i} + \sum_{+\infty\text{-limit circles}} L_{\gamma_{+\infty_j}} \geq 0.
\]
The former is due to (\ref{E hat omega estimate}), and the latter is because
\begin{align*}
&e_{Y_i} + \sum_{+\infty\text{-limit circles}} L_{\gamma_{+\infty_j}}\\
& = \int_{u^{-1}((-\epsilon, 0] \times Y^-)} u^\ast d(\varphi \lambda_i)
+ \int_{u^{-1}(Z^+_i)} u^\ast \omega
+ \int_{u^{-1}(([0, \infty) \cup \R_1 \cup \dots \cup \R_{k_+}) \times Y^+)}
u^\ast d\lambda^+\\
&\geq 0,
\end{align*}
where
\begin{align*}
&\int_{u^{-1}((-\epsilon, 0] \times Y^-)} u^\ast d(\varphi \lambda_i)\\
&= \int_{u^{-1}((-\epsilon, 0] \times Y^-)} u^\ast (d\varphi \wedge \lambda_i)
+ \int_{u^{-1}((-\epsilon, 0] \times Y^-)} u^\ast (\varphi d \lambda_i)\\
&\geq 0
\end{align*}
since we may assume $\partial_\sigma \varphi \geq 0$ on $(-\epsilon, 0]$.

Define a super-commutative algebra $\D_{X, (Y_i)}$ as follows.
Its elements are formal series
\[
\sum_{(\hat c_i^\ast), (\hat c'_i), A}
f_{(\hat c_i^\ast), (\hat c'_i), A}(t, \hbar) t_{x_1} \dots t_{x_{k_t}} q^-_{\hat c_1^\ast}
\dots q^-_{\hat c_{k_q}^\ast} p^+_{\hat c'_1} \dots p^+_{\hat c'_{k_p}} T^A,
\]
where
$f_{(\hat c_i^\ast), (\hat c'_i), A}(t, \hbar) \in \R[[t, \hbar]]$ are
formal series of the variables $t_x$ ($x \in K^0_X$) and $\hbar$, and
the infinite sum is taken over all sequences $((\hat c_i), (\hat c'_i))$
as in the usual case and
$A \in H_2(\overline{X}, \partial \overline{X}; \Z) \big/
(\Ker e \cap \bigcap_i \Ker e_{Y_i})$,
where $e(A) = \widetilde{\omega}A$ and $e_{Y_i}(A) = \widetilde{\omega}_{Y_i} A$.
We impose the following Novikov condition on the elements of $\D_{X, (Y_i)}$:
for any $C > 0$, the number of the non-zero terms with
$\sum_j e(p^+_{\hat c'_j}) \geq - C$,
$e(A) + \sum_j e(p^+_{\hat c'_j}) \geq - C$ and
$e_{Y_i}(A) + \sum_j e(p^+_{\hat c'_j}) \geq - C$ is finite.

We also define a bigger super-commutative algebra $\D\D_{X, (Y_i)}$.
Its elements are formal series
\[
\sum_{(\hat c_i^\ast), (\hat c'_i), A}
f_{(\hat c_i^\ast), (\hat c'_i), A}(t, \hbar) t_{x_1} \dots t_{x_{k_t}} q^-_{\hat c_1^\ast}
\dots q^-_{\hat c_{k_q}^\ast} p^+_{\hat c'_1} \dots p^+_{\hat c'_{k_p}} T^A,
\]
similar to those of $\D_{X, (Y_i)}$ except that
$f_{(\hat c_i^\ast), (\hat c'_i), A}(t, \hbar) \in \R[[\hbar]][\hbar^{-1}][[t]]$.

The submodules $\D_{X, (Y_i)}^{\leq \kappa} \subset \D_{X, (Y_i)}$
and $\D\D_{X, (Y_i)}^{\leq \kappa} \subset \D\D_{X, (Y_i)}$ are defined by
the conditions
$\sum_i e(q^-_{\hat c_i^\ast}) + e(A) + \sum_j e(p^+_{\hat c'_j}) \leq \kappa$
and $e_{Y_i}(A) + \sum_j e(p^+_{\hat c'_j}) \leq \kappa$.
For each positive constant $\delta > 0$, we define a submodule
$\D_{X, (Y_i)}^{\leq \kappa, \delta} \subset \D_{X, (Y_i)}^{\leq \kappa}$ by the condition
$\widetilde{g}_\delta \geq - \kappa / \delta$.
The submodules $\widetilde{J}^{\leq \kappa, \delta}_{C_0, C_1, C_2} \subset
\D\D_{X, (Y_i)}^{\leq \kappa, \delta}$ are defined by
\begin{align*}
&\widetilde{J}^{\leq \kappa, \delta}_{C_0, C_1, C_2} \\
&= \Bigl\{\sum
a_{(x_i), (\hat c_i^\ast), (\hat c'_i), g, e} t_{x_1} \dots t_{x_{k_t}} q^-_{\hat c_1^\ast}
\dots q^-_{\hat c_{k_q}^\ast} p^+_{\hat c'_1} \dots p^+_{\hat c'_{k_p}} \hbar^g T^A
\in \D\D_{X, (Y_i)}^{\leq \kappa, \delta};\\
&\quad \quad a_{(x_i), (\hat c_i^\ast), (\hat c'_i), g, A} = 0 \text{ for all }
((x_i)_{i = 1}^{k_t}, (\hat c_i^\ast),_{i = 1}^{k_q} (\hat c'_i)_{i = 1}^{k_p}, g, A)
\text{ such that}\\
&\quad \quad
k_t \leq C_0,\ \widetilde{g}_\delta \leq C_1,\ \sum e(p^+_{\hat c'_i}) \geq - C_2,\ 
e(A) + \sum e(p^+_{\hat c'_i}) \geq - C_2 \\
& \quad \quad \text{and }
e_{Y_j}(A) + \sum e(p^+_{\hat c'_i}) \geq - C_2 \text{ for all } j
\Bigr\},
\end{align*}
and submodules $J^{\leq \kappa, \delta}_{C_0, C_1, C_2} \subset
\D_{X, (Y_i)}^{\leq \kappa}$ are defined by
$J^{\leq \kappa, \delta}_{C_0, C_1, C_2} = \D_{X, (Y_i)} \cap
\widetilde{J}^{\leq \kappa, \delta}_{C_0, C_1, C_2}$.

The same argument are valid for $\D_{X, (Y_i)}$ and $\D\D_{X, (Y_i)}$.
Namely, a compatible finite family of virtual fundamental chains defines
generating functions $\F \in (\hbar^{-1} \D_{X, (Y_i)}^{\leq 0})^{\star, \delta}
/ J^{\star, \delta}_{\overline{C}_0, \overline{C}_1, \overline{C}_2}$ and
$\widetilde{\F} \in \D\D_{X, (Y_i)}^{\leq 0, \delta}
/ \widetilde{J}^{\leq 0, \delta}_{\overline{C}_0, \overline{C}_1, \overline{C}_2}$,
and they define differentials
$D_\F : \D_{X, (Y_i)}^{\leq \kappa} /J^{\leq \kappa, \delta}_{C_0, C_1, C_2}
\to \D_{X, (Y_i)}^{\leq \kappa} /J^{\leq \kappa, \delta}_{C_0, C_1, C_2}$.

The rational version is similarly defined.
Let $\mathcal{L}_{X, (Y_i)} = \D_{X, (Y_i)}|_{\hbar = 0}$ be the quotient super-commutative
algebra.
Elements of the Poisson space $\widehat{\mathcal{L}}_{X, (Y_i)}$ are formal series
\[
\sum
f_{\alpha}(t)
q_{\hat c_1^\ast}^- \dots q_{\hat c_{k^-_q}^\ast}^-
q_{(\hat c'_1)^\ast}^+ \dots q_{(\hat c'_{k^+_q})^\ast}^+
p_{\hat c''_1}^- \dots p_{\hat c''_{k^-_p}}^-
p_{\hat c'''_1}^+ \dots p_{\hat c'''_{k^+_p}}^+ T^A,
\]
as in the usual $\widehat{\mathcal{L}}_X$ with Novikov condition, that is,
for any $C > 0$, the number of the non-zero terms with
$\sum_j e(p^-_{\hat c''_j}) + \sum_j e(p^+_{\hat c'''_j}) \geq - C$,
$e(A) + \sum_j e(p^-_{\hat c''_j}) + \sum_j e(p^+_{\hat c'''_j}) \geq - C$
and $e_{Y_i}(A) + \sum_j e(p_{\hat c''_j}^-) + \sum_j e(p_{\hat c'''_j}^+) \geq - C$
is finite.
Then we can define a differential $d_{\F_0} : \mathcal{L}_{X, (Y_i)}^{\leq \kappa}
/ J^{\leq \kappa}_{C_0, C_2} \to \mathcal{L}_{X, (Y_i)}^{\leq \kappa}
/ J^{\leq \kappa}_{C_0, C_2}$ as in the usual case.

We note that we do not need to consider the case of contact homology since
in this case, we consider only exact cobordisms.

%% file: SFT-10_Fiber_products_for_homotopy.tex
%
%

\section{The case of homotopy}
\label{Homotopy}
In this section, we prove that two generating functions for $(X, \omega)$ and $K_X^0$
defined by using different almost complex structures and perturbations are homotopic
in the sense of \cite{EGH00}.
Furthermore, we prove that its homotopy type does not change if we change
the symplectic form $\omega$ on $Z$ by an exact form, and $K_X^0$ by
boundaries.
\subsection{Fiber products and their orientations}
\label{fiber product and orientation for homotopy}
Let $(X^\tau, \omega^\tau)_{\tau \in I = [0, 1]}$ be a family of symplectic manifolds
with cylindrical ends such that
the manifold $X^\tau = X = (-\infty, 0] \times Y^- \cup Z \cup [0, \infty) \times Y^+$ is
independent of $\tau$ and the symplectic forms have the form
$\omega^\tau = \omega^0 + d \theta^\tau$ for some
one-forms $\theta^\tau$ whose supports are contained in $Z$.
Let $J^\tau$ be a family of $\omega^\tau$-compatible almost complex structures
whose restriction to $(-\infty, 0] \times Y^-$ and $[0, \infty) \times Y^+$ are
independent of $\tau$ and obtained by some complex structures of $\xi^\pm$.
For each $i = 0, 1$, let $K_{X^i}^0$ be a finite set of smooth cycles with closed
support in $X$ with bijections
\begin{align*}
\mu^i_- &: \{x \in K_{X^i}^0; \supp x \cap (-\infty, 0] \times Y^- \neq \emptyset\}
\to K_{Y^-}^0\\
\mu^i_+ &: \{x \in K_{X^i}^0; \supp x \cap [0, \infty) \times Y^+ \neq \emptyset\}
\to K_{Y^+}^0
\end{align*}
such that $x|_{(-\infty, 0] \times Y^-} = (-\infty, 0] \times \mu^i_-(x)$
and $x_{[0, \infty) \times Y^+} = [0, \infty) \times \mu^i_+(x)$.
Assume that a finite set $K_{X^I}^0 = \{(x^\tau)_{\tau \in I}\}$ of
$C^\infty(I, \R)$-linear combinations of smooth cycles with closed supports in $X$
is given which satisfies the following conditions:
\begin{itemize}
\item
$\{x^0\} = K_{X^0}^0$ and $\{x^1\} = K_{X^1}^0$.
Furthermore, we assume that $x^\tau$ is constant near $\tau = 0,1$.
\item
$\frac{d}{d\tau} x^\tau$ are boundaries of some $C^\infty(I, \R)$-linear combinations of
smooth chains $(y^\tau)_{\tau \in I}$ in $X$ whose supports are contained in $Z$
for each $(x^\tau)_{\tau \in I}$.
In particular, $x^\tau$ is independent of $\tau$ on the complement of $Z$.
\item
There exist bijections
\begin{align*}
\mu_- &: \{(x^\tau)_{\tau \in I} \in K_{X^I}^0; \supp x^\tau \cap (-\infty, 0] \times Y^-
\neq \emptyset\} \to K_{Y^-}^0\\
\mu_+ &: \{(x^\tau)_{\tau \in I} \in K_{X^I}^0; \supp x^\tau \cap [0, \infty) \times Y^+
\neq \emptyset\} \to K_{Y^+}^0
\end{align*}
such that
$x^\tau|_{(-\infty, 0] \times Y^-} = (-\infty, 0] \times \mu_-((x^\tau)_{\tau \in I})$,
$x^\tau_{[0, \infty) \times Y^+} = [0, \infty) \times \mu_+((x^\tau)_{\tau \in I})$ and
$\mu_\pm((x^\tau)_{\tau \in I}) = \mu^i_\pm(x^i)$ for $i = 0, 1$.
\end{itemize}

As with $\widehat{\M}(X, \omega, J)$ in Section \ref{Kuranishi for X},
we can construct a pre-Kuranishi structure of
\[
\widehat{\M}_{X^I} = \bigcup_{\tau \in I} \widehat{\M}(X^\tau, \omega^\tau, J^\tau).
\]
There exists a natural strong smooth map from $\widehat{\M}_{X^I}$ to
$I$ which maps $\widehat{\M}_{X^\tau} = \widehat{\M}(X^\tau, \omega^\tau, J^\tau)$
to $\tau \in I$, and the associated maps from the Kuranishi neighborhoods to
$I$ are submersive.

Assume that grouped multisections of the fiber products
$(\widehat{\M}_{Y^\pm}^\diamond, \mathring{K}_{Y^\pm}^2, \ab K_{Y^\pm},
\ab K_{Y^\pm}^0)$ and
$(\widehat{\M}^\diamond_{X^i}, \mathring{K}^2_{X^i}, K_{X^i}, K^0_{X^i})$ for $i = 0,1$
are given and that they satisfy the compatibility conditions.
We define $\widehat{\M}_{X^I}^\diamond$ by
$\widehat{\M}_{X^I}^\diamond = \bigcup_{\tau \in I} \widehat{\M}_{X^\tau}^\diamond$.
Note that this is not locally the product of $\widehat{\M}_{X^I}$ but
its fiber product over $I$.
Its element is written as $(\tau,
(\Sigma^\alpha, z^\alpha, u^\alpha)_{\alpha \in A^- \sqcup A^0 \sqcup A^+},
M^{\mathrm{rel}})$, and
all of $(\Sigma^\alpha, z^\alpha, u^\alpha)$ are holomorphic buildings
for the same $X^\tau$.
Note that $\widehat{\M}_{X^I}$ contains $I \times \widehat{\M}_{Y^\pm}$.

Similarly to the case of $X$ in Section \ref{Kuranishi for X},
we define the fiber product
\[
(\widehat{\M}^\diamond_{X^I}, (\mathring{K}_{Y^-}^2, \mathring{K}_{Y^+}^2),
(K_{Y^-}, K_{Y^+}), (K^0_{Y^-}, K^0_{X^I}, \partial_\tau K^0_{X^I}, K^0_{Y^+})) \subset
\widehat{\M}^\diamond_{X^I}.
\]
In this case, we need to explain the meaning of fiber product with
$K^0_{X^I}$ or $\partial_\tau K^0_{X^I}$.
Let $\boldsymbol{K}^0_{X^I}$ be the finite set of smooth chains
which appear in $x^\tau$ or $y^\tau$.
Assume that $y^\tau = 0$ on $[0, \epsilon] \cup [1 - \epsilon, 1]$
for all $(x^\tau)_{\tau \in I} \in K^0_{X^I}$ for some $\epsilon > 0$.
Then the meaning of fiber product with $K^0_{X^I}$ and $\partial_\tau K^0_{X^I}$
is that we assume the transversality condition of the grouped multisection
with respect to $K^0_{X^0}$, $\boldsymbol{K}^0_{X^I}$ and $K^0_{X^1}$
on $[0, \epsilon]$, $[\epsilon, 1 - \epsilon]$ and $[1 - \epsilon, 1]$ respectively.
We abbreviate the above fiber product as
$(\widehat{\M}^\diamond_{X^I}, \mathring{K}^2_{X^I}, K_{X^I}, K^0_{X^I})$.

We define multi-valued partial essential submersions
$\Xi^\circ$ and $\Lambda$ from
$(\widehat{\M}^\diamond_{X^I}, \ab \mathring{K}^2_{X^I}, K_{X^I}, K^0_{X^I})$
to itself similarly.
We need to construct the grouped multisections of
$(\widehat{\M}^\diamond_{X^I}, \mathring{K}^2_{X^I}, K_{X^I}, K^0_{X^I})$
which satisfy the similar compatibility conditions.
Notice that for a disconnected holomorphic building of $\widehat{\M}_{X^I}$,
the perturbed multisection induced by the product of the perturbed multisections
transverse to the zero sections for the connected components is not always
transverse to the zero section.
This is because we need to use the same factor $I$ for all connected components.
In other words, it is not a product but a fiber product with respect to $I$.
To overcome this problem, we use a continuous family of grouped multisections.
We can construct the continuous families of grouped multisections of
$(\widehat{\M}^\diamond_{X^I}, \mathring{K}^2_{X^I}, K_{X^I}, K^0_{X^I})$
which satisfy the following conditions:

\begin{itemize}
\item
The restrictions of the natural map
$(\widehat{\M}^\diamond_{X^I}, \mathring{K}^2_{X^I}, K_{X^I}, K^0_{X^I}) \to I$
to the fiber products of the zero sets of the perturbed multisection
are essentially submersive, that is,
even if we restrict to the fiber product of the zero set with the simplices, it is
submersive.
\item
The restrictions of the grouped multisection of
$(\widehat{\M}^\diamond_{X^I}, \mathring{K}^2_{X^I}, K_{X^I}, K^0_{X^I})$ to
$I \times
(\widehat{\M}_{Y^\pm}^\diamond, \mathring{K}_{Y^\pm}^2, K_{Y^\pm}, K_{Y^\pm}^0)$
coincide with the pull backs of the grouped multisections of
$(\widehat{\M}_{Y^\pm}^\diamond, \mathring{K}_{Y^\pm}^2, K_{Y^\pm}, K_{Y^\pm}^0)$
by the projection.
\item
The restrictions of the grouped multisection of
$(\widehat{\M}^\diamond_{X^I}, \mathring{K}^2_{X^I}, K_{X^I}, K^0_{X^I})$ to
$(\widehat{\M}^\diamond_{X^i}, \mathring{K}^2_{X^i}, K_{X^i}, K^0_{X^i})$
for $i = 0,1$ coincide with the given grouped multisections.
\item
Let $((\widehat{\M}^\diamond_{X^I})^0, \mathring{K}^2_{X^I}, K_{X^I}, K^0_{X^I}) \subset
(\widehat{\M}^\diamond_{X^I}, \mathring{K}^2_{X^I}, K_{X^I}, K^0_{X^I})$
be the subset of connected points.
Its continuous family of grouped multisections induces that of
\[
\bigcup_N (\prod^N (\widehat{\M}^\diamond_{X^I})^0, \mathring{K}^2_{X^I}, K_{X^I},
K^0_{X^I})_I / \mathfrak{S}_N,
\]
where each $(\prod^N (\widehat{\M}^\diamond_{X^I})^0, \mathring{K}^2_{X^I}, K_{X^I},
K^0_{X^I})_I$ is the fiber product over $I$.
This is because of the first condition.
Then the grouped multisection of
$(\widehat{\M}^\diamond_{X^I}, \mathring{K}^2_{X^I}, K_{X^I}, K^0_{X^I})$
coincides with its pull back by the submersion
\[
(\widehat{\M}^\diamond_{X^I}, \mathring{K}^2_{X^I}, K_{X^I}, K^0_{X^I}) \to
\bigcup_N (\prod^N (\widehat{\M}^\diamond_{X^I})^0, \mathring{K}^2_{X^I}, K_{X^I},
K^0_{X^I})_I / \mathfrak{S}_N
\]
defined by decomposition into connected components.
\item
The grouped multisections of
$(\widehat{\M}^\diamond_{X^I}, \mathring{K}^2_{X^I}, K_{X^I}, K^0_{X^I})$
are compatible with respect to the compatible system of
multi-valued partial essential submersions defined by $\Xi^\circ$ and $\Lambda$.
\end{itemize}

The grouped multisection of each
$\overline{\M}^{(m_-, X^I, m_+)}_{((\hat \epsilon^{i, j}_l), (\hat c^i_l), (x^i_l), (\hat \eta^i_l))}$
is defined by the pull back of that of
$(\widehat{\M}^\diamond_{X^I}, \mathring{K}^2_{X^I}, K_{X^I}, K^0_{X^I})$
by the natural essential submersion.

The definition of the orientation of
$\overline{\M}^{(m_-, X^I, m_+)}_{((\hat \epsilon^{i, j}_l), (\hat c^i_l), (x^i_l), (\hat \eta^i_l))}$
is almost the same with the case of $X$.
The only difference is that it is defined by
\[
(TI \oplus \W^{-m_-} \oplus \dots \oplus \W^{m_+})_{\star}
\]
instead of (\ref{ori}).

Assume that two pairs of solutions $(G^\pm)^0$, $(G^\pm)^1$ of (\ref{+G eq})
and (\ref{-G eq}) in Section \ref{correction terms for X} are given.
Then we can construct a smooth family of the solutions $(G^\pm)^\tau$ ($\tau \in I$)
of the equations which coincide with the given solution at $\tau = 0, 1$.

For a triple of sequences $((\hat c_l), (x_l), (\alpha_l))$, we define a pre-Kuranishi space
(or a $C^\infty(I, \R)$-linear combination of pre-Kuranishi spaces)
$\overline{\M}^{X^I}((\hat c_l), (x_l), (\alpha_l))$ by
\begin{align*}
&\overline{\M}^{X^I}((\hat c_l), (x_l), (\alpha_l))\\
& = \sum_{m_-, m_+ \geq 0}\sum_{\star_{m_-, m_+}}(-1)^\ast \,
\overline{\M}^{(m_-, X^I, m_+)}
_{(((\widetilde{G}^+_{m_+})^I, (\widetilde{G}^-_{-m_-})^I), (\hat c^i_l), (x^i_l),
([\overline{P}_{Y^+}] \cap \alpha^i_l))}
\end{align*}
where
$(\widetilde{G}^\pm)^\tau = (\widetilde{G}^\pm_0)^\tau
+ (\widetilde{G}^\pm_{\pm1})^\tau +
(\widetilde{G}^\pm_{\pm2})^\tau + \cdots = \Theta^\pm(e^{\otimes (G^\pm)^\tau})$.
The sum and the sign $\ast$ are the same as those of the non-parametrized case.

Let
\begin{align*}
\bigl[\overline{\M}^{X^I \!, e}_g((\hat c_l), (x_l), (\alpha_l))\bigr]
&= (f_{0, g}^e)^\tau((\hat c_l), (x_l), (\alpha_l))
\oplus (f_{1, g}^e)^\tau((\hat c_l), (x_l), (\alpha_l)) d\tau\\
\bigl[(\overline{\M}^{X^I\!,  e}_g)^0((\hat c_l), (x_l), (\alpha_l))\bigr]
&= (h_{0, g}^e)^\tau((\hat c_l), (x_l), (\alpha_l))
\oplus (h_{1, g}^e)^\tau((\hat c_l), (x_l), (\alpha_l)) d\tau
\end{align*}
be the counterparts of virtual fundamental chains, where $f_{j, g}^e((\hat c_l), (x_l), (\alpha_l))$
and $h_{j, g}^e((\hat c_l), (x_l), (\alpha_l))$ ($j = 0, 1$) are smooth functions of
$\tau \in I = [0, 1]$.

Let $(H^\pm)^\tau = (H^\pm)^\tau_2 + (H^\pm)^\tau_3 + \cdots$ be
an appropriate $C^\infty(I, \R)$-linear combination of
\begin{multline*}
((\rho_\ast [\overline{P}_{Y^\pm}])^{i, j}, \dots, (\rho_\ast [\overline{P}_{Y^\pm}])^{i, j},
\epsilon_{\overline{P}_{Y^\pm}}^{i, j}, \dots, \epsilon_{\overline{P}_{Y^\pm}}^{i, j},\\
(\Delta_\ast [\overline{P}_{Y^\pm}])^{i, j}, \dots, (\Delta_\ast [\overline{P}_{Y^\pm}])^{i, j})
\end{multline*}
defined in the next section, and define
$(\mathring{f}_{0, g}^e)^\tau((\hat c_l), (x_l), (\alpha_l))$
by the $\Omega^0(I)$ part of the virtual fundamental chain of the $(g, e)$-part of
\[
\sum_{\substack{m_- \geq 0\\ m_+ \geq 0}}\sum_{\star_{m_-, m_+}}
(-1)^\ast\bigl(\overline{\M}^{(m_-, X^I, m_+)}
_{\kappa_1} + \overline{\M}^{(m_-, X^I, m_+)}
_{\kappa_2}\bigr),
\]
where
\[
\kappa_1 = (\Theta^+(e^{\otimes (G^+)^\tau})_{m_+},
\Theta^-(e^{\otimes (G^-)^\tau} \otimes (H^-)^\tau)_{-m_-}, (\hat c^i_l), (x^i_l), ([\overline{P}_{Y^+}] \cap \alpha^i_l))
\]
and
\[
\kappa_2 = (\Theta^+(e^{\otimes (G^+)^\tau} \otimes (H^+)^\tau)_{m_+},
\Theta^-((-1)^{m_-}e^{\otimes (G^-)^\tau})_{-m_-},
(\hat c^i_l), (x^i_l), ([\overline{P}_{Y^+}] \cap \alpha^i_l)).
\]
We also define $(\mathring{h}_{0, g}^e)^\tau((\hat c_l), (x_l), (\alpha_l))$ by
the $\Omega^0(I)$ part of the virtual fundamental chain of its irreducible part.

Define $\hat f_{1, g}^e((\hat c_l), (x_l), (\alpha_l))$ and
$\hat h_{1, g}^e((\hat c_l), (x_l), (\alpha_l))$ by
\begin{align*}
&\hat f_{1, g}^e((\hat c_l), (x_l), (\alpha_l))\\
&= - f_{1, g}^e((\hat c_l), (x_l), (\alpha_l)) + \mathring{f}_{0, g}^e((\hat c_l), (x_l), (\alpha_l))\\
&\quad + \sum_j (-1)^{\sum |\hat c_l| + \sum_{i < j} |x_i|}
f_{0, g}^e((\hat c_l), (x_1, x_2, \dots, y_j, \dots, x_{k_t}), (\alpha_l))
\end{align*}
and
\begin{align*}
&\hat h_{1, g}^e((\hat c_l), (x_l), (\alpha_l))\\
&= - h_{1, g}^e((\hat c_l), (x_l), (\alpha_l))
+ \mathring{h}_{0, g}^e((\hat c_l), (x_l), (\alpha_l))\\
&\quad + \sum_j (-1)^{\sum |\hat c_l| + \sum_{i < j} |x_i|}
h_{0, g}^e((\hat c_l), (x_1, x_2, \dots, y_j, \dots, x_{k_t}), (\alpha_l)).
\end{align*}
The second terms $\mathring{f}_{0, g}^e((\hat c_l), (x_l), (\alpha_l))$ and
$\mathring{h}_{0, g}^e((\hat c_l), (x_l), (\alpha_l))$ correspond to the differential of
$(G^\pm)^\tau$, and the third terms corresponds to the differential of $x^\tau$.
There terms are added to make equation (\ref{boundary formula for X^I}) below
hold true.

Then $f_{0, g}^e$, $h_{0, g}^e$, $\hat f_{1, g}^e$ and $\hat h_{1, g}^e$ satisfy
the following equations.
\begin{equation}
f_{0, g}^e((\hat c_l), (x_l), (\alpha_l)) = \sum_{\star_0} (-1)^{\ast_0} \frac{1}{k !} \prod_{i = 1}^k
h_{0, g_i}^{e_i}((\hat c^i_l), (x^i_l), (\alpha^i_l)) \label{f_0 equation}
\end{equation}
\begin{equation}
\hat f_{1, g}^e((\hat c_l), (x_l), (\alpha_l)) = \sum_{\star_1} (-1)^{\ast_1}
f_{0, g_0}^{e_0}((\hat c^0_l), (x^0_l), (\alpha^0_l))
\, \hat h_{1, g_1}^{e_1}((\hat c^1_l), (x^1_l), (\alpha^1_l))
\label{f_1 equation}
\end{equation}
\begin{align}
&df_{0, g}^e((\hat c_l), (x_l), (\alpha_l)) \notag\\
&= \hat f_{1, g}^e(\partial((\hat c_l), (x_l), (\alpha_l))) d\tau \notag\\
&\quad - \sum_{\star'_-} (-1)^{\ast'_-} \frac{1}{k !}
[\overline{\M}^{Y^-}_{g_{-}}((\hat c^{-}_l), (x^{-}_l), (\hat d_1^\ast, \hat d_2^\ast, \dots,
\hat d_k^\ast))]^0 \notag\\
&\hph{\quad + \sum_{\star_-} (-1)^{\ast_-} \frac{1}{k !}}
\cdot \hat f_{1, g_0}^{e_0}((\hat d_k, \hat d_{k - 1}, \dots, \hat d_1) \cup (\hat c^0_l),
(x^0_l), (\alpha_l)) d\tau \notag\\
&\quad - \sum_{\star'_+} (-1)^{\ast'_+} \frac{1}{k !}
\hat f_{1, g_0}^e((\hat c_l), (x^0_l), (\alpha^0_l) \cup (\hat d_1^\ast, \hat d_2^\ast, \dots,
\hat d_k^\ast)) \notag\\
&\hph{\quad + \sum_{\star_+} (-1)^{\ast_+} \frac{1}{k !}}
\cdot [\overline{\M}^{Y^+}_{g_+}((\hat d_k, \hat d_{k - 1}, \dots, \hat d_1), (x^+_l), (\alpha^+_l))]^0 d\tau
\label{boundary formula for X^I}
\end{align}
\begin{align}
0 &= f_0(\partial((\hat c_l), (x_l), (\alpha_l)))\notag\\
&\quad - \sum_{\star_-} (-1)^{\ast_-} \frac{1}{k !}
[\overline{\M}^{Y^-}((\hat c^-_l), (x^-_l), (\hat d_1^\ast, \hat d_2^\ast, \dots, \hat d_k^\ast))]^0\notag\\
&\hph{\quad + \sum_{\star_+} (-1)^{\ast_+} \frac{1}{k !}}
\times f_0((\hat d_k, \hat d_{k - 1}, \dots, \hat d_1) \cup (\hat c^0_l), (x^0_l), (\alpha_l))\notag\\
&\quad + \sum_{\star_+} (-1)^{\ast_+} \frac{1}{k !}
f_0((\hat c_l), (x^0_l), (\alpha^0_l) \cup (\hat d_1^\ast, \hat d_2^\ast, \dots, \hat d_k^\ast))\notag\\
&\hph{\quad + \sum_{\star_+} (-1)^{\ast_+} \frac{1}{k !}}
\times [\overline{\M}^{Y^+}((\hat d_k, \hat d_{k - 1}, \dots, \hat d_1), (x^+_l), (\alpha^+_l))]^0
\label{boundary formula for each tau}
\end{align}
The sum $\star_0$ is taken over all $k \geq 0$, all decompositions
$g - 1 = \sum_{i = 1}^k (g_i - 1)$,
$e = \sum_{i = 1}^k e_i$ and all decompositions
\[
\{\hat c_l\} = \coprod_{i = 1}^k \{\hat c^i_l\}, \quad
\{x_l\} = \coprod_{i = 1}^k \{x^i_l\}, \quad
\{\alpha_l\} = \coprod_{i = 1}^k \{\alpha^i_l\}
\]
as sets.
The sign $\ast_0$ is the weighted sign of the permutation
\[
\begin{pmatrix}
(c^1_l) \ (x^1_l) \ (\alpha^1_l) \dots (c^k_l) \ (x^i_l) \ (\alpha^i_l)\\
(c_l) \quad (x_l) \quad (\alpha_l)
\end{pmatrix}.
\]
The sum $\star_1$ is taken over all decompositions
$g - 1= (g_0 - 1) + (g_1 - 1)$, $e = e_0 + e_1$
and all decompositions
\[
\{c_l\} = \{c^0_l\} \sqcup \{c^1_l\}, \quad \{x_l\} = \{x^0_l\} \sqcup \{x^1_l\}, \quad
\{\alpha_l\} = \{\alpha^0_l\} \sqcup \{\alpha^1_l\}
\]
as sets, and the sign $\ast_1$ is the weighted sign of the permutation
\[
\begin{pmatrix}
(c^0_l) \ (x^0_l) \ (\alpha^0_l) \ (c^1_l) \ (x^1_l) \ (\alpha^1_l)\\
(c_l) \ (x_l) \ (\alpha_l)
\end{pmatrix}.
\]
The sum $\star'_-$ is taken over $k \geq 0$, all simplices $d_l$ of $K_{Y^-}$ not contained
in $\overline{P}_{Y^-}^{\text{bad}}$,
all decompositions
\[
\{\hat c_l\} = \{\hat c^-_l\} \sqcup \{\hat c^0_l\}, \quad
\{x_l\} = \{x^-_l\} \sqcup \{x^0_l\}
\]
such that $x^-_l \in K^0_{Y^-}$, and all pairs
$(g_-, g_0)$ such that $g = g_- + g_0 + k -1$.
The sign $\ast'_-$ is the weighted sign of the permutation
\[
\begin{pmatrix}
(\hat c^-_l) \ (x^-_l) \ (\hat c^0_l) \ (x^0_l)\\
(\hat c_l) \quad (x_l)
\end{pmatrix}.
\]
The sum $\star'_+$ is taken over $k \geq 0$, all simplices $d_l$ of $K_{Y^+}$ not contained
in $\overline{P}_{Y^+}^{\text{bad}}$, and all decompositions
\[
\{x_l\} = \{x^0_l\} \sqcup \{x^+_l\}, \quad \{\alpha_l\} = \{\alpha^0_l\} \sqcup \{\alpha^+_l\}
\]
such that $x^+_l \in K^0_{Y^+}$, and all pairs
$(g_0, g_+)$ such that $g = g_0 + g_+ + k - 1$.
The sign $\ast'_+$ is the weighted sign of the permutation
\[
\begin{pmatrix}
(x^0_l) \ (\alpha^0_l) \ (x^+_l) \ (\alpha^+_l)\\
(x_l) \quad (\alpha_l)
\end{pmatrix}.
\]
Equation (\ref{boundary formula for each tau}) is a counterpart of equation
(\ref{boundary formula for X}), and the meaning of the sums and the signs are
the same.

As with equation (\ref{irreducible decomposition}),
(\ref{f_0 equation}) and (\ref{f_1 equation}) are due to the irreducible decomposition.
(\ref{boundary formula for X^I}) is proved in the next section.
(\ref{boundary formula for each tau}) is same as (\ref{boundary formula for X}).

\subsection{Construction of $H^\tau$}
In this section, we construct smooth families
$(H^+)^\tau = (H^+_2)^\tau + (H^+_3)^\tau + \dots \in (\bigoplus_{m \geq 2}
(\mathring{\B}_m^+)^{m+1})^\wedge$
and $(H^-)^\tau = (H^-_{-2})^\tau + (H^-_{-3})^\tau + \dots \in (\bigoplus_{m \geq 2}
(\mathring{\B}_{-m}^-)^{m+1})^\wedge$ such that
\begin{multline}
\partial \Theta^+(e^{\otimes G^+} \otimes H^+)
+ \sum_{i \geq 0} (-1)^i e^{\Delta_{(e_i, e_{i+1})}} \tau^+_i
\Theta^+(e^{\otimes G^+} \otimes H^+)\\
- \Diamond^+ \Bigl(\Theta^+(e^{\otimes G^+} \otimes H^+) \otimes
\sum_{j \geq 1} (-1)^j F^+_j\Bigr)
- \Theta^+\Bigl(e^{\otimes G^+} \otimes \frac{d}{d\tau}G^+\Bigr) = 0, \label{H^+ eq}
\end{multline}
\begin{multline}
\partial \Theta^-(e^{\otimes G^-} \otimes H^-)
+ \sum_{i \leq 0} e^{\Delta_{(e_{i-1}, e_i)}} \tilde \tau^-_i
\Theta^-(e^{\otimes G^-} \otimes H^-)\\
+ \Diamond^-(F^- \otimes \Theta^-(e^{\otimes G^-} \otimes H^-))
- \Theta^-(e^{\otimes G^-} \otimes \frac{d}{d\tau}G^-) = 0, \label{H^- eq}
\end{multline}
and prove (\ref{boundary formula for X^I}) for these $(H^+)^\tau$ and $(H^-)^\tau$.

First we construct $(H^+)^\tau$.
We inductively construct smooth families $(H^+)^\tau_{\leq m} = (H^+_2)^\tau
+ (H^+_3)^\tau + \dots + (H^+_m)^\tau \in
\bigoplus_{l = 2}^m (\mathring{\B}_l^+)^{l+1}$ such that
\begin{multline}
\partial \Theta^+(e^{\otimes G^+} \otimes H^+_{\leq m})
+ \sum_{i \geq 0} (-1)^i e^{\Delta_{(e_i, e_{i+1})}} \tau^+_i
\Theta^+(e^{\otimes G^+} \otimes H^+_{\leq m})\\
- \Diamond^+ \Bigl(\Theta^+(e^{\otimes G^+} \otimes H^+_{\leq m}) \otimes
\sum_{j \geq 1} (-1)^j F^+_j\Bigr)
- \Theta^+\Bigl(e^{\otimes G^+} \otimes \frac{d}{d\tau}G^+\Bigr) \equiv 0
\end{multline}
in $(\bigoplus_{l = 1}^\infty (\B_l^+)^l)^\wedge
/(\bigoplus_{l = m+1}^\infty (\B_l^+)^l)^\wedge$

First note that $\frac{d}{d\tau}G^+ \in (\bigoplus_{l \geq 2}
(\mathring{\B}_l^+)^l)^\wedge$
since $G^+_1$ is independent of $\tau$.
Hence we do not need $H_1^+$-part.

Assume we have already constructed $H^+_{\leq m-1} = H^+_2 + H^+_3 + \dots
+ H^+_{m-1} \in (\bigoplus_{l = 2}^{m-1} (\mathring{\B}_l^+)^{l+1})^\wedge$.
Then it is enough to prove that
\begin{multline}
\partial \Theta^+(e^{\otimes G^+} \otimes H^+_{\leq m-1})
+ \sum_{i \geq 0} (-1)^i e^{\Delta_{(e_i, e_{i+1})}} \tau^+_i
\Theta^+(e^{\otimes G^+} \otimes H^+_{\leq m-1})\\
- \Diamond^+ \Bigl(\Theta^+(e^{\otimes G^+} \otimes H^+_{\leq m-1}) \otimes
\sum_{j \geq 1} (-1)^j F^+_j\Bigr)
- \Theta^+ \Bigl(e^{\otimes G^+} \otimes \frac{d}{d\tau}G^+\Bigr) \equiv 0
\label{H tau mathring}
\end{multline}
in $(\bigoplus_{l = 1}^\infty (\B_l^+)^l)^\wedge
/((\bigoplus_{l = m+1}^\infty (\B_l^+)^l)^\wedge \oplus
\bigoplus_{l = 1}^\infty (\mathring{\B}_l^+)^l)	$ and
\begin{align}
&\partial\Bigl(\sum_{i \geq 0} (-1)^i e^{\Delta_{(e_i, e_{i+1})}} \tau^+_i
\Theta^+(e^{\otimes G^+} \otimes H^+_{\leq m-1})\notag\\
&\quad - \Diamond^+\Bigl(\Theta^+(e^{\otimes G^+} \otimes H^+_{\leq m-1}) \otimes
\sum_{j \geq 1} (-1)^j F^+_j\Bigr)
- \Theta^+\Bigl(e^{\otimes G^+} \otimes \frac{d}{d\tau}G^+\Bigr)\Bigr)\notag\\
& \equiv 0
\label{H tau closed}
\end{align}
in $(\bigoplus_{l = 1}^\infty (\B_l^+)^{l-1})^\wedge
/(\bigoplus_{l = m+1}^\infty (\B_l^+)^{l-1})^\wedge$

The latter equation is proved by an argument similar to those for (\ref{A closed}) or
(\ref{B closed}).
The former can be proved in a similar way to equation (\ref{B mathring}) by using
the following equations.
\begin{align}
&\partial \Theta^+\biggl(\frac{1}{k !} G^{\otimes k} \otimes H \biggr)\notag\\
&= \Theta^+\biggl(\frac{1}{k !} G^{\otimes k} \otimes \partial H \biggr)
+ \Theta^+\biggl(\frac{1}{(k-1) !} G^{\otimes (k-1)} \otimes \partial G \otimes
\sum_j (-1)^j H_j \biggr)
\end{align}
\begin{align}
&\sum_{i \geq 1} (-1)^i e^{\Delta_{(e_i, e_{i+1})}} \tau^+_i \Theta^+
\biggl(\frac{1}{k !} G^{\otimes k} \otimes H \biggr)\notag\\
&= \Theta^+\biggl(\frac{1}{k !} G^{\otimes k} \otimes
\sum_{i \geq 1} (-1)^i e^{\Delta_{(e_i, e_{i+1})}} \tau^+_i H \biggr)\notag\\
&\quad + \Theta^+\biggl(\frac{1}{(k-1) !} G^{\otimes (k-1)} \otimes
\sum_{i \geq 1} (-1)^i e^{\Delta_{(e_i, e_{i+1})}} \tau^+_i G \otimes \sum_j (-1)^j H_j\biggr)
\end{align}
\begin{align}
&e^{\Delta_{(e_0, e_1)}} \tau^+_0 \Theta^+\biggl(\frac{1}{k !} G^{\otimes k} \otimes
H \biggr)\notag\\
&= \sum_{l = 0}^k  \Theta^+\biggl(\frac{1}{(k-l) ! l !}G^{\otimes (k-l)} \otimes
(e^{\Delta_{(e_0, e_1)}} \mathring{\tau}^+_0(G^{\otimes l} \otimes H))\biggr)\notag\\
&\quad + \sum_{l = 0}^k  \Theta^+\biggl(\frac{1}{(k-l) ! l !}G^{\otimes (k-l)} \otimes
(e^{\Delta_{(e_0, e_1)}} \mathring{\tau}^+_0(G^{\otimes l})) \otimes \sum_j (-1)^j H_j\biggr)
\end{align}
\begin{align}
&\Diamond^+\biggl(\Theta^+\biggl(\frac{1}{k !} G^{\otimes k} \otimes H \biggr)
\otimes \sum_j (-1)^j F^+_j\biggr)\notag\\
&= \sum_{l = 0}^k  \Theta^+\biggl(\frac{1}{(k-l) ! l !}G^{\otimes (k-l)} \otimes
\mathring{\Diamond}^+ (G^{\otimes l} \otimes H \otimes
\sum_j (-1)^j F^+_j)\biggr)\notag\\
&\quad - \sum_{l = 0}^k  \Theta^+\biggl(\frac{1}{(k-l) ! l !}G^{\otimes (k-l)} \otimes
\mathring{\Diamond}^+ (G^{\otimes l} \otimes F^+) \otimes \sum_j (-1)^j H_j\biggr)
\end{align}

Similarly, we can inductively construct a smooth family
$(H^-)^\tau_{\geq -m} = (H^-_{-2})^\tau + (H^-_{-3})^\tau + \dots
+ (H^-_{-m}) \in (\bigoplus_{l = 2}^m
(\mathring{\B}_{-l}^-)^{l+1})^\wedge$ such that
\begin{multline}
\partial \Theta^-(e^{\otimes G^-} \otimes H^-_{\geq -m})
+ \sum_{i \leq 0} e^{\Delta_{(e_{i-1}, e_i)}} \tilde \tau^-_i
\Theta^-(e^{\otimes G^-} \otimes H^-_{\geq -m})\\
+ \Diamond^-(F^- \otimes \Theta^-(e^{\otimes G^-} \otimes H^-_{\geq -m}))
- \Theta^-(e^{\otimes G^-} \otimes \frac{d}{d\tau}G^-) \equiv 0
\end{multline}
in $(\bigoplus_{l = 1}^\infty (\B_{-l}^-)^l)^\wedge
/(\bigoplus_{l = m+1}^\infty (\B_{-l}^-)^l)^\wedge$,
and we obtain a required solution $(H^-)^\tau \in (\bigoplus_{m \geq 2}
(\mathring{\B}_{-m}^-)^{m+1})^\wedge$.

Now we prove (\ref{boundary formula for X^I}) for these $H^\pm$.
In what follows, we omit the subscripts $g$ or $e$ for the simplification of
notation.
We abbreviate
\[
\sum (-1)^\ast ((\hat c^i_l), (x^i_l), ([\overline{P}_{Y^+}] \cap \alpha^i_l))
\]
as $((\hat c_l), (x_l), (\alpha_l))$,
where $(-1)^\ast$ is the weighted sign of the permutation corresponding to
$((\hat c^i_l), (x^i_l), (\alpha^i_l))$, and the sum is taken over all decomposition of
$(\hat c_l)$, $(x_l)$ and $(\alpha_l)$.
The sums below are taken over all $(m_-, m_+)$
(and all $m$, $k$, all sequences of simplices
$(\hat d_l)_{l = 1}^k$ of $K_{Y^\pm}$ not contained in $\overline{P}_{Y^\pm}^{\text{bad}}$,
and all decomposition $((\hat c_l), (x_l), (\alpha_l)) = ((\hat c'_l), (x'_l), (\alpha'_l)) \sqcup
((\hat c''_l), (x''_l), (\alpha''_l))$
if they appear).
We abbreviate $\Theta^-(e^{\otimes G^-} \otimes H^-)$ and
$\Theta^+(e^{\otimes G^+} \otimes H^+)$ to $\widetilde{H}^-$ and $\widetilde{H}^+$
respectively.
It is easy to check that
\begin{align}
&d f_0((\hat c_l), (x_l), (\alpha_l)) \notag\\
& = \sum  d\bigl[\overline{\M}^{(m_-, X^I, m_+)}
_{((\widetilde{G}^+_{m_+}, \widetilde{G}^-_{-m_-}), (\hat c_l), (x_l),
(\alpha_l))}\bigr]^0 \notag\\
&= \sum   \bigl[\overline{\M}^{(m_-, X^I, m_+)}
_{(\partial_\tau(\widetilde{G}^+_{m_+}, \widetilde{G}^-_{-m_-}), (\hat c_l), (x_l),
(\alpha_l))}
\bigr]^0 \notag\\
& \quad + \sum  \bigl[\overline{\M}^{(m_-, X^I, m_+)}
_{((\widetilde{G}^+_{m_+}, \widetilde{G}^-_{-m_-}), (\hat c_l), \partial_\tau (x_l),
(\alpha_l))}\bigr]^0 \notag\\
& \quad - \sum  \bigl[\overline{\M}^{(m_-, X^I, m_+)}
_{((\widetilde{G}^+_{m_+}, \widetilde{G}^-_{-m_-}), \partial((\hat c_l), (x_l),
(\alpha_l)))}\bigr]^1 \notag\\
& \quad - \sum  (-1)^{m_-} \bigl[\overline{\M}^{(m_-, X^I, m_+)}
_{((\partial' \widetilde{G}^+_{m_+}, \widetilde{G}^-_{-m_-}), (\hat c_l), (x_l),
(\alpha_l))}\bigr]^1 \notag\\
& \quad - \sum  \bigl[\overline{\M}^{(m_-, X^I, m_+)}
_{((\widetilde{G}^+_{m_+}, \partial'\widetilde{G}^-_{-m_-}), (\hat c_l), (x_l),
(\alpha_l))}\bigr]^1 \notag\\
& \quad -\sum  (-1)^{m_-} \bigl[\overline{\M}^{(m_-, X^I, m_+ + 1)}
_{((\sum_i (-1)^i e^{\Delta_\ast [\overline{P}_{Y^+}]^{i, i+1}} \tau^+_i \widetilde{G}^+_{m_+},
\widetilde{G}^-_{-m_-}), (\hat c_l), (x_l),
(\alpha_l))}\bigr]^1 \notag\\
& \quad - \sum  \bigl[\overline{\M}^{(m_- +1, X^I, m_+)}
_{((\widetilde{G}^+_{m_+}, \sum_i e^{\Delta_\ast [\overline{P}_{Y^-}]^{i-1, i}}
\tilde \tau^-_i \widetilde{G}^-_{-m_-}), (\hat c_l), (x_l),
(\alpha_l))}\bigr]^1,
\label{df_0}
\end{align}
(\ref{+G eq}) implies
\begin{align}
&- \sum  (-1)^{m_-}\bigl(\bigl[\overline{\M}^{(m_-, X^I, m_+)}
_{((\partial' \widetilde{G}^+_{m_+}, \widetilde{G}^-_{-m_-}), (\hat c_l), (x_l),
(\alpha_l))}\bigr]^1 \notag\\
& \hph{- \sum  (-1)^{m_-}}
+ \bigl[\overline{\M}^{(m_-, X^I, m_+)}
_{((\sum_i (-1)^i e^{\Delta_\ast [\overline{P}_{Y^+}]^{i, i+1}} \tau^+_i \widetilde{G}^+_{m_+},
\widetilde{G}^-_{-m_-}), (\hat c_l), (x_l),
(\alpha_l))}\bigr]^1\bigr) \notag\\
&= \sum  (-1)^{m_-} \bigl(\bigl[\overline{\M}^{(m_-, X^I, m_+ + m)}
_{((\Diamond^+(\widetilde{G}^+_{m_+} \otimes F^+_m), \widetilde{G}^-_{-m_-}),
(\hat c_l), (x_l), (\alpha_l))}\bigr]^1 \bigr) \notag\\
&= \sum  \frac{1}{k!}
\bigl[\overline{\M}^{(m_-, X^I, m_+)}
_{((\widetilde{G}^+_{m_+}, \widetilde{G}^-_{-m_-}), (\hat c'_l), (x'_l),
(\alpha'_l) \cup
(\hat d_l^\ast)_{l= 1}^k)}\bigr]^1 \notag\\
& \quad \ \hph{\sum  \frac{1}{k!}} \cdot
\bigl[(\overline{\M}_{Y^+})^m
_{(F^+_m, (\hat d_l)_{l=k}^1, (x''_l),
(\alpha''_l))} \bigr]^0,
\label{G^+ 2}
\end{align}
and (\ref{-G eq}) implies
\begin{align}
& - \sum  \bigl(\bigl[\overline{\M}^{(m_-, X^I, m_+)}
_{((\widetilde{G}^+_{m_+}, \partial' \widetilde{G}^-_{-m_-}), (\hat c_l), (x_l),
(\alpha_l))}\bigr]^1 \notag\\
& \hph{- \sum  \bigl(}
+ \bigl[\overline{\M}^{(m_-, X^I, m_+)}
_{((\widetilde{G}^+_{m_+}, \sum_i e^{\Delta_\ast [\overline{P}_{Y^-}]^{i-1, i}} \tilde \tau^-_i
\widetilde{G}^-_{-m_-}), (\hat c_l), (x_l),
(\alpha_l))}\bigr]^1 \bigr) \notag\\
&= \sum  \bigl[\overline{\M}^{(m_- + m, X^I, m_+)}
_{((\widetilde{G}^+_{m_+}, \Diamond^-(F^-_m \otimes \widetilde{G}^-_{-m_-})),
(\hat c_l), (x_l),
(\alpha_l))}\bigr]^1 \notag\\
&= \sum  \frac{1}{k!} \bigl[(\overline{\M}_{Y^-})^m
_{(F^-_m, (\hat c'_l), (x'_l),
(\hat d_l^\ast)_{l=1}^k)} \bigr]^0 \notag\\
&\quad \ \hph{\sum  \frac{1}{k!}} \cdot
\bigl[\overline{\M}^{(m_-, X^I, m_+)}
_{((\widetilde{G}^+_{m_+}, \widetilde{G}^-_{-m_-}), (\hat d_l)_{l=k}^1
\cup (\hat c''_l), (x''_l),
(\alpha_l))}\bigr]^1.
\label{G^- 2}
\end{align}
(\ref{H^- eq}) implies
\begin{align}
&\sum  \bigl[\overline{\M}^{(m_-, X^I, m_+)}
_{((\widetilde{G}^+_{m_+}, \partial_\tau \widetilde{G}^-_{-m_-}), (\hat c_l), (x_l),
(\alpha_l))}\bigr]^0 \notag\\
& = \sum  \bigl[\overline{\M}^{(m_-, X^I, m_+)}
_{((\widetilde{G}^+_{m_+}, \partial' \widetilde{H}^-_{-m_-}),
(\hat c_l), (x_l),
(\alpha_l))}\bigr]^0 \notag\\
& \quad + \sum  \bigl[\overline{\M}^{(m_-, X^I, m_+)}
_{((\widetilde{G}^+_{m_+},
\sum_i e^{\Delta_\ast [\overline{P}_{Y^-}]^{i-1, i}} \tilde \tau^-_i
\widetilde{H}^-_{-m_-}), (\hat c_l), (x_l),
(\alpha_l))}\bigr]^0 \notag\\
& \quad + \sum  \bigl[\overline{\M}^{(m_-, X^I, m_+)}
_{((\widetilde{G}^+_{m_+},
\Diamond^-(F^-_m \otimes \widetilde{H}^-_{-m_-})),
(\hat c_l), (x_l),
(\alpha_l))}\bigr]^0.
\label{partial tau - first}
\end{align}
(\ref{boundary of MMX}) implies
\begin{align}
0 &= \sum  \bigl[\partial \overline{\M}^{(m_-, X^\tau, m_+)}
_{((\widetilde{G}^+_{m_+}, \widetilde{H}^-_{-m_-}),
(\hat c_l), (x_l),
(\alpha_l))}\bigr]^0 \notag\\
&= \sum (-1)^{m_-} \bigl( \bigl[\overline{\M}^{(m_-, X^I, m_+)}
_{((\partial' \widetilde{G}^+_{m_+}, \widetilde{H}^-_{-m_-}),
(\hat c_l), (x_l),
(\alpha_l))}\bigr]^0 \notag\\
&\quad \hph{\sum (-1)^{m_-} \bigl(}
+ \bigl[\overline{\M}^{(m_-, X^I, m_+ +1)}
_{((\sum_i (-1)^i e^{\Delta_\ast [\overline{P}_{Y^+}]^{i, i+1}} \tau^+_i \widetilde{G}^+_{m_+},
\widetilde{H}^-_{-m_-}), (\hat c_l), (x_l),
(\alpha_l))}\bigr]^0
\bigr) \notag\\
&\quad + \sum  \bigl(\bigl[\overline{\M}^{(m_-, X^I, m_+)}
_{((\widetilde{G}^+_{m_+}, \partial' \widetilde{H}^-_{-m_-}),
(\hat c_l), (x_l),
(\alpha_l))}\bigr]^0 \notag\\
& \quad \hph{+ \sum  \bigl(}
+ \bigl[\overline{\M}^{(m_- + 1, X^I, m_+)}
_{((\widetilde{G}^+_{m_+},
\sum_i e^{\Delta_\ast [\overline{P}_{Y^-}]^{i-1, i}} \tilde \tau^-_i
\widetilde{H}^-_{-m_-}), (\hat c_l), (x_l),
(\alpha_l))}\bigr]^0
\bigr) \notag\\
&\quad - \sum  \bigl[\overline{\M}^{(m_-, X^I, m_+)}
_{((\widetilde{G}^+_{m_+}, \widetilde{H}^-_{-m_-}),
\partial((\hat c_l), (x_l),
(\alpha_l)))}\bigr]^0.
\end{align}
(\ref{+G eq}) implies
\begin{align}
&\sum (-1)^{m_-} \bigl( \bigl[\overline{\M}^{(m_-, X^I, m_+)}
_{((\partial' \widetilde{G}^+_{m_+}, \widetilde{H}^-_{-m_-}),
(\hat c_l), (x_l),
(\alpha_l))}\bigr]^0 \notag\\
&\hph{\sum (-1)^{m_-} \bigl(}
+ \bigl[\overline{\M}^{(m_-, X^I, m_+ +1)}
_{((\sum_i (-1)^i e^{\Delta_\ast [\overline{P}_{Y^+}]^{i, i+1}} \tau^+_i \widetilde{G}^+_{m_+},
\widetilde{H}^-_{-m_-}), (\hat c_l), (x_l),
(\alpha_l))}\bigr]^0
\bigr) \notag\\
&= -\sum (-1)^{m_-} \bigl[\overline{\M}^{(m_-, X^I, m_+ + m)}
_{((\Diamond^+(\widetilde{G}^+_{m_+} \otimes F^+_m), \widetilde{H}^-_{-m_-}), (\hat c_l), (x_l),
(\alpha_l))}\bigr]^0
\end{align}
It is easy to check the following equations.
\begin{align}
&\sum (-1)^{m_-} \bigl[\overline{\M}^{(m_-, X^I, m_+ + m)}
_{((\Diamond^+(\widetilde{G}^+_{m_+} \otimes F^+_m),
\widetilde{H}^-_{-m_-}),
(\hat c_l), (x_l),
(\alpha_l))}\bigr]^0 \notag\\
&= - \sum  \frac{1}{k!} \bigl[\overline{\M}^{(m_-, X^I, m_+)}
_{((\widetilde{G}^+_{m_+}, \widetilde{H}^-_{-m_-}),
(\hat c_l), (x'_l),
(\alpha'_l) \cup (
\hat d_l^\ast)_{l=1}^k)}\bigr]^0 \notag\\
& \quad \hph{- \sum  \frac{1}{k!}} \cdot
\bigl[(\overline{\M}_{Y^+})^m_{(F^+_m, (\hat d_l)_{l=k}^1, (x''_l), (\alpha''_l))} \bigr]^0,
\end{align}
\begin{align}
&\sum  \bigl[\overline{\M}^{(m_- + m, X^I, m_+)}
_{((\widetilde{G}^+_{m_+},
\Diamond^-(F^-_m \otimes \widetilde{H}^-_{-m_-})),
(\hat c_l), (x_l),
(\alpha_l))}\bigr]^0 \notag\\
&= - \sum  \frac{1}{k!}
\bigl[(\overline{\M}_{Y^-})^m
_{(F^-_m, (\hat c'_l), (x'_l), (\hat d_l^\ast)_{l=1}^k)} \bigr]^0 \notag\\
& \quad \hph{- \sum  \frac{1}{k!}} \cdot
\bigl[\overline{\M}^{(m_-, X^I, m_+)}
_{((\widetilde{G}^+_{m_+}, \widetilde{H}^-_{-m_-}), (\hat d_l)_{l=k}^1,
(\hat c''_l), (x''_l),
(\alpha_l))}\bigr]^0.
\label{partial tau - end}
\end{align}
(\ref{partial tau - first}) to (\ref{partial tau - end}) imply
\begin{align}
&\sum  \bigl[\overline{\M}^{(m_-, X^I, m_+)}
_{((\widetilde{G}^+_{m_+}, \partial_\tau \widetilde{G}^-_{-m_-}), (\hat c_l), (x_l),
(\alpha_l))}\bigr]^0 \notag\\
&= \sum (-1)^{m_-} \bigl[\overline{\M}^{(m_-, X^I, m_+ + m)}
_{((\Diamond^+(\widetilde{G}^+_{m_+} \otimes F^+_m), \widetilde{H}^-_{-m_-}), (\hat c_l), (x_l),
(\alpha_l))}\bigr]^0 \notag\\
&\quad + \sum  \bigl[\overline{\M}^{(m_-, X^I, m_+)}
_{((\widetilde{G}^+_{m_+}, \widetilde{H}^-_{-m_-}),
\partial((\hat c_l), (x_l), (\alpha_l)))}\bigr]^0 \notag\\
&\quad + \sum  \bigl[\overline{\M}^{(m_-, X^I, m_+)}
_{((\widetilde{G}^+_{m_+},
\Diamond^-(F^-_m \otimes \widetilde{H}^-_{-m_-})), (\hat c_l), (x_l),
(\alpha_l))}\bigr]^0 \notag\\
&= \sum  \bigl[\overline{\M}^{(m_-, X^I, m_+)}
_{((\widetilde{G}^+_{m_+}, \widetilde{H}^-_{-m_-}),
\partial((\hat c_l), (x_l), (\alpha_l)))}\bigr]^0 \notag\\
&\quad - \sum  \frac{1}{k!} \bigl[\overline{\M}^{(m_-, X^I, m_+)}
_{((\widetilde{G}^+_{m_+}, \widetilde{H}^-_{-m_-}),
(\hat c_l), (x'_l), (\alpha'_l) \cup (\hat d_l^\ast)_{l=1}^k)}\bigr]^0 \notag\\
& \quad \hph{- \sum  \frac{1}{k!}} \cdot
\bigl[(\overline{\M}_{Y^+})^m_{(F^+_m, (\hat d_l)_{l=k}^1, (x''_l),
(\alpha''_l))} \bigr]^0 \notag\\
&\quad - \sum  \frac{1}{k!}
\bigl[(\overline{\M}_{Y^-})^m_{(F^-_m, (\hat c'_l), (x'_l),
(\hat d_l^\ast)_{l=1}^k)} \bigr]^0 \notag\\
& \quad \hph{- \sum  \frac{1}{k!}} \cdot
\bigl[\overline{\M}^{(m_-, X^I, m_+)}
_{((\widetilde{G}^+_{m_+}, \widetilde{H}^-_{-m_-}), (\hat d_l)_{l=k}^1,
(\hat c''_l), (x''_l), (\alpha_l))}\bigr]^0.
\label{partial tau -}
\end{align}
Similarly,
\begin{align}
&\sum  \bigl[\overline{\M}^{(m_-, X^I, m_+)}
_{((\partial_\tau \widetilde{G}^+_{m_+}, \widetilde{G}^-_{-m_-}), (\hat c_l), (x_l),
(\alpha_l))}\bigr]^0 \notag\\
&= \sum (-1)^{m_-} \bigl[\overline{\M}^{(m_-, X^I, m_+)}
_{((\widetilde{H}^+_{m_+}, \widetilde{G}^-_{-m_-}),
\partial((\hat c_l), (x_l),
(\alpha_l)))}\bigr]^0 \notag\\
&\quad -\sum (-1)^{m_-} \frac{1}{k!}
\bigl[\overline{\M}^{(m_-, X^I, m_+)}
_{((\widetilde{H}^+_{m_+}, \widetilde{G}^-_{-m_-}),
(\hat c_l), (x'_l),
(\alpha'_l) \cup
(\hat d_l^\ast)_{l=1}^k)}\bigr]^0 \notag\\
& \quad \hph{-\sum (-1)^{m_-} \frac{1}{k!}}
\cdot [\bigl[(\overline{\M}_{Y^+})^m_{(F^+_m, (\hat d_l)_{l=k}^1, (x''_l),
(\alpha''_l))} \bigr]^0 \notag\\
&\quad - \sum (-1)^{m_-} \frac{1}{k!}
\bigl[(\overline{\M}_{Y^-})^m_{(F^-_m, (\hat c'_l), (x'_l),
(\hat d_l^\ast)_{l=1}^k)} \bigr]^0
\notag\\
& \quad \hph{- \sum (-1)^{m_-} \frac{1}{k!}}
\cdot \bigl[\overline{\M}^{(m_-, X^I, m_+)}
_{((\widetilde{H}^+_{m_+}, \widetilde{G}^-_{-m_-}),
(\hat d_l)_{l=k}^1,
(\hat c''_l), (x''_l),
(\alpha_l))}\bigr]^0.
\label{partial tau +}
\end{align}
(\ref{boundary of MMX}) implies
\begin{align}
&\sum \bigl[\partial \overline{\M}^{(m_-, X^\tau, m_+)}
_{((\widetilde{G}^+_{m_+}, \widetilde{G}^-_{-m_-}), (\hat c_l), (x_1, \dots, y_r, \dots, x_{k_t}),
(\alpha_l))}\bigr]^0 \notag\\
&= \sum \bigl[\overline{\M}^{(m_-, X^I, m_+)}
_{((\widetilde{G}^+_{m_+}, \widetilde{G}^-_{-m_-}), \partial((\hat c_l),
(x_1, \dots, y_r, \dots, x_{k_t}),
(\alpha_l)))}\bigr]^0 \notag\\
&\quad + \sum (-1)^{m_-} \bigl(\bigl[\overline{\M}^{(m_-, X^I, m_+)}
_{((\partial' \widetilde{G}^+_{m_+}, \widetilde{G}^-_{-m_-}), (\hat c_l),
(x_1, \dots, y_r, \dots, x_{k_t}),
(\alpha_l))}\bigr]^0 \notag\\
& \quad \hph{+ \sum (-1)^{m_-} \bigl(}
+ \bigl[\overline{\M}^{(m_-, X^I, m_+)}
_{((\sum_i (-1)^i e^{\Delta_\ast [\overline{P}]^{i, i+1}} \tau^+_i \widetilde{G}^+_{m_+},
\widetilde{G}^-_{-m_-}), (\hat c_l), (x_1, \dots,  y_r, \dots, x_{k_t}),
(\alpha_l))}\bigr]^0 \bigr) \notag\\
& \quad + \sum \bigl(\bigl[\overline{\M}^{(m_-, X^I, m_+)}
_{((\widetilde{G}^+_{m_+}, \partial' \widetilde{G}^-_{-m_-}), (\hat c_l),
(x_1, \dots, y_r, \dots, x_{k_t}),
(\alpha_l))}\bigr]^0 \notag\\
& \quad \hph{+ \sum \bigl(}
+ \bigl[\overline{\M}^{(m_-, X^I, m_+)}
_{((\widetilde{G}^+_{m_+}, \sum_i e^{\Delta_\ast [\overline{P}_{Y^-}]^{i-1, i}} \tilde \tau^-_i
\widetilde{G}^-_{-m_-}), (\hat c_l), (x_1, \dots, y_r, \dots, x_{k_t}),
(\alpha_l))}\bigr]^0 \bigr) \notag\\
&= \sum \bigl[\overline{\M}^{(m_-, X^I, m_+)}
_{((\widetilde{G}^+_{m_+}, \widetilde{G}^-_{-m_-}), (\hat c_l),
\partial (x_1, \dots, y_r, \dots, x_{k_t}),
(\alpha_l))}\bigr]^0 \notag\\
& \quad - \sum (-1)^{m_-}\bigl[\overline{\M}^{(m_-, X^I, m_+ + m)}
_{((\Diamond^+(\widetilde{G}^+_{m_+} \otimes F^+_m), \widetilde{G}^-_{-m_-}), (\hat c_l),
(x_1, \dots, y_r, \dots, x_{k_t}),
(\alpha_l))}\bigr]^0 \notag\\
& \quad - \sum \bigl[\overline{\M}^{(m_-, X^I, m_+)}
_{((\widetilde{G}^+_{m_+}, \Diamond^-(F^-_m \otimes \widetilde{G}^-_{-m_-})), (\hat c_l),
(x_1, \dots, y_r, \dots, x_{k_t}),
(\alpha_l))}\bigr]^0.
\end{align}
Hence
\begin{align}
&\sum (-1)^{\sum |c_l| + \sum_{l < r}|x_l|}\bigl[\overline{\M}^{(m_-, X^I, m_+)}
_{((\widetilde{G}^+_{m_+}, \widetilde{G}^-_{-m_-}), \partial((\hat c_l),
(x_1, \dots, y_r, \dots, x_{k_t}),
(\alpha_l)))}\bigr]^0 \notag\\
&= - \sum (-1)^{\sum |c_l| + \sum_{l < r'}|x'_l|}\bigl[\overline{\M}^{(m_-, X^I, m_+)}
_{((\widetilde{G}^+_{m_+}, \widetilde{G}^-_{-m_-}), (\hat c_l),
(x'_1, \dots, y_{r'}, \dots,  x'_{k'_t}),
(\alpha'_l) \cup
(\hat d_l^\ast)_{l=1}^k)}\bigr]^0 \notag\\
& \quad \hph{- \sum (-1)^{\sum |c_l| + \sum_{l < r'}|x'_l|}}
\cdot \bigl[(\overline{\M}_{Y^+})^m_{(F^+_m, (\hat d_l)_{l_=k}^1, (x''_l),
(\alpha''_l))} \bigr]^0 \notag\\
&\quad - \sum \bigl[(\overline{\M}_{Y^-})^m_{(F^-_m, (\hat c'_l), (x'_l),
(\hat d_l^\ast)_{l=1}^k)} \bigr]^0 \notag\\
& \quad \hph{- \sum}
\cdot (-1)^{\sum |c''_l| + \sum_{l < r''}|x''_l|}\bigl[\overline{\M}^{(m_-, X^I, m_+)}
_{((\widetilde{G}^+_{m_+}, \widetilde{G}^-_{-m_-}), (\hat c''_l),
(x''_1, \dots, y_{r''}, \dots, x''_{k''_t}), (\alpha_l))}\bigr]^0
\label{partial tau x}
\end{align}

(\ref{df_0}), (\ref{G^+ 2}), (\ref{G^- 2}), (\ref{partial tau -}), (\ref{partial tau +})
and (\ref{partial tau x})
imply (\ref{boundary formula for X^I}).


%% file: SFT-11_Algebras_for_homotopy.tex
%
%

\subsection{Construction of homotopies}\label{algebra for homotopy}
We define families of generating functions
$\F^\tau, \K^\tau \in (\hbar^{-1} \D_X^{\leq 0})^{\star, \delta}
/ J^{\star, \delta}_{\overline{C}_0, \overline{C}_1, \overline{C}_2}$
and $\widetilde{\F}^\tau, \widetilde{\K}^\tau \in \D\D_X^{\leq 0, \delta}
/ \widetilde{J}^{\leq 0, \delta}_{\overline{C}_0, \overline{C}_1, \overline{C}_2}$
by
\begin{align*}
\F^\tau &= \hbar^{-1} \sum \frac{1}{k_q ! k_t ! k_p !} (h_{0, g}^e)^\tau
(\underbrace{\mathbf{q}, \dots, \mathbf{q}}_{k_q},
\underbrace{\mathbf{t}, \dots, \mathbf{t}}_{k_t},
\underbrace{\mathbf{p}, \dots, \mathbf{p}}_{k_p}) \hbar^g T^e\\
\widetilde{\F}^\tau &= \hbar^{-1} \sum \frac{1}{k_q ! k_t ! k_p !} (f_{0, g}^e)^\tau
(\underbrace{\mathbf{q}, \dots, \mathbf{q}}_{k_q},
\underbrace{\mathbf{t}, \dots, \mathbf{t}}_{k_t},
\underbrace{\mathbf{p}, \dots, \mathbf{p}}_{k_p}) \hbar^g T^e\\
\K^\tau &= \hbar^{-1} \sum \frac{1}{k_q ! k_t ! k_p !} (\hat h_{1, g}^e)^\tau
(\underbrace{\mathbf{q}, \dots, \mathbf{q}}_{k_q},
\underbrace{\mathbf{t}, \dots, \mathbf{t}}_{k_t},
\underbrace{\mathbf{p}, \dots, \mathbf{p}}_{k_p}) \hbar^g T^e\\
\widetilde{\K}^\tau &= \hbar^{-1} \sum \frac{1}{k_q ! k_t ! k_p !} (\hat f_{1, g}^e)^\tau
(\underbrace{\mathbf{q}, \dots, \mathbf{q}}_{k_q},
\underbrace{\mathbf{t}, \dots, \mathbf{t}}_{k_t},
\underbrace{\mathbf{p}, \dots, \mathbf{p}}_{k_p}) \hbar^g T^e
\end{align*}

Then (\ref{f_0 equation}), (\ref{f_1 equation}),
(\ref{boundary formula for X^I}) and (\ref{boundary formula for each tau})
imply the following equations.
\begin{gather*}
\widetilde{\F}^\tau = e^{\F^\tau}\\
\widetilde{\K}^\tau = \K^\tau \widetilde{\F}^\tau\\
\frac{d}{d\tau} \widetilde{\F}^\tau = \delta \widetilde{\K}^\tau
- \overrightarrow{\mathcal{H}^-} \widetilde{\K}^\tau
- \widetilde{\K}^\tau \overleftarrow{\mathcal{H}^+}\\
\widehat{D}_X(\widetilde{\F}^\tau) = 0
\end{gather*}
Therefore, the following equation holds true in $\D\D_X^{\leq 0, \delta}
/ \widetilde{J}^{\leq 0, \delta}_{\overline{C}_0, \overline{C}_1, \overline{C}_2}$.
\begin{align}
\frac{d}{d\tau} (e^{\F^\tau}) &= \widehat{D}_X (\K^\tau e^{\F^\tau}) \notag\\
&= [\widehat{D}_X, \K^\tau] (e^{\F^\tau}) \label{main eq for X^I}
\end{align}
Namely, the family of functions $\F^\tau$ is a homotopy in the sense
of \cite{EGH00}.
\begin{defi}
One-parameter family of functions $F^\tau \in (\hbar^{-1} \D_X^{\leq 0})^{\star, \delta}
/ J^{\star, \delta}_{\overline{C}_0, \overline{C}_1, \overline{C}_2}$ ($\tau \in [0, 1]$)
of even degree is said to be a homotopy
if (\ref{main eq for X}) holds for all $\F = \F^\tau$ and
there exists a family of functions $\K^\tau \in (\hbar^{-1} \D_X^{\leq 0})^{\star, \delta}
/ J^{\star, \delta}_{\overline{C}_0, \overline{C}_1, \overline{C}_2}$ of odd degree
which makes equation (\ref{main eq for X^I}) holds for all $\tau \in [0, 1]$.
\end{defi}
\begin{rem}
If (\ref{main eq for X}) holds for some $\F^\tau$ and (\ref{main eq for X^I}) is satisfied
for all $\tau \in [0, 1]$, then (\ref{main eq for X}) holds for all $\F^\tau$.
\end{rem}
\begin{rem}
(\ref{main eq for X^I}) is equivalent to
\begin{equation}
\frac{d}{d\tau} \F^\tau = D_{\F^\tau}(\K^\tau)
\label{main eq for X F}
\end{equation}
in $(\hbar^{-1} \D_X^{\leq 0})^{\star, \delta}
/ J^{\star, \delta}_{\overline{C}_0, \overline{C}_1, \overline{C}_2}$.
\end{rem}

First we consider the case of general SFT.
As in \cite{EGH00}, we define flows by linear differential equations.
For each four-tuple $(\kappa, C_0, C_1, C_2)$ such that
$\overline{C}_0 \geq C_0$, $\overline{C}_1 \geq C_1
+ \kappa \delta^{-1}$ and $\overline{C}_2 \geq C_2 + \kappa$,
we define a flow $\Phi^\tau : \D\D_X^{\leq \kappa, \delta}
/ \widetilde{J}^{\kappa, \delta}_{C_0, C_1, C_2}
\to \D\D_X^{\leq \kappa, \delta}
/ \widetilde{J}^{\kappa, \delta}_{C_0, C_1, C_2}$ by
\[
\frac{d}{d\tau} \Phi^\tau(f) = [\widehat{D}_X, \K^\tau] \Phi^\tau(f), \quad \Phi^0 = \id,
\]
and define $T^\tau : \D_X^{\leq \kappa} / J^{\leq \kappa, \delta}_{C_0, C_1, C_2}
\to \D_X^{\leq \kappa} / J^{\leq \kappa, \delta}_{C_0, C_1, C_2}$ by
\[
T^\tau(f) = e^{-\F^\tau} \Phi^\tau (e^{\F^0} f).
\]
$\Phi^\tau$ is well defined because it is defined by a linear differential equation
on a finite dimensional vector space.
$T^\tau$ is well defined because it is also defined by
\[
\frac{d}{d\tau} T^\tau(f)
= [[D_{\F^\tau}, \K^\tau], T^\tau(f)] (1),
\quad T^0 = \id.
\]
Some of the following were proved in \cite{EGH00}, and some are straightforward,
but we prove all of them for the convenience of the reader.
\begin{lem}\label{homotopy and T}
\begin{enumerate}[label = \normalfont (\roman*)]\ 
\item
$T^\tau$ is a chain map from
$(\D_X^{\leq \kappa} / J^{\leq \kappa, \delta}_{C_0, C_1, C_2}, D_{\F^0})$ to
$(\D_X^{\leq \kappa} / J^{\leq \kappa, \delta}_{C_0, C_1, C_2}, D_{\F^\tau})$
for each $\tau$, that is,
$D_{\F^\tau} \circ T^\tau = T^\tau \circ D_{\F^0}$.
(This is equivalent to $\widehat{D}_X \circ \Phi^\tau = \Phi^\tau \circ \widehat{D}_X$.)
Furthermore, up to chain homotopy, it is determined by $(\F^{\tau'})_{\tau' \in [0, \tau]}$
and independent of the choice of the family $(\K^{\tau'})_{\tau' \in [0, \tau]}$
which satisfies equation (\ref{main eq for X^I}).
\item
If a smooth family of generating functions $\F^{\tau, \sigma}
\in (\hbar^{-1} \D_X^{\leq 0})^{\star, \delta}
/ J^{\star, \delta}_{\overline{C}_0, \overline{C}_1, \overline{C}_2}$
$((\tau, \sigma) \in [0, 1] \times [0, 1])$
satisfies $\F^{0, \sigma} \equiv \F^{0, 0}$ and the one-parameter family
$(\F^{\tau, \sigma})_{\tau \in [0, 1]}$ is a homotopy for each $\sigma \in [0, 1]$,
then the one-parameter family $(\F^{\tau, \sigma})_{\sigma \in [0, 1]}$ is also
a homotopy for each $\tau \in [0, 1]$.
\item
Further assume that the above family of generating functions satisfies
$\F^{1, \sigma} = \F^{0, 0}$ and $\F^{\tau, 0} \equiv \F^{0, 0}$.
Let $T^\tau : \D_X^{\leq \kappa} / J^{\leq \kappa, \delta}_{C_0, C_1, C_2}
\to \D_X^{\leq \kappa} / J^{\leq \kappa, \delta}_{C_0, C_1, C_2}$
be the flow defined by the one-parameter homotopy $(\F^{\tau, 1})_{\tau \in [0, 1]}$.
Then $T^1$ is equal to the identity map up to chain homotopy.
In other words, if a loop homotopy $(\F^{\tau, 1})_{\tau \in S^1}$ is contractible
in the space of loop homotopies with the base point $\F^{0, 1}$,
then the chain map $T^1$ is the identity map up to chain homotopy.
Hence for a general one-parameter homotopy $(\F^\tau)_{\tau \in [0, 1]}$,
the end $T^1$ of the family of the chain maps $(T^\tau)_{\tau \in [0, 1]}$
is determined up to chain homotopy by the homotopy type of the homotopy
$(\F^\tau)_{\tau \in [0, 1]}$ relative to the end points.
\item
There exists a family of linear maps
$A^{\pm, \tau} : \W^{\leq \kappa}_{Y^\pm}
/ I^{\leq \kappa}_{C_0, C_1 + \kappa(\delta^{-1} - L^{-1}_{\min}), C_2}
\to \D^{\leq \kappa}_X / J^{\leq \kappa, \delta}_{C_0, C_1, C_2}$
such that
\begin{equation}
i_{\F^0}^\pm - (T^\tau)^{-1} \circ i_{\F^\tau}^\pm = D_{\F^0} \circ A^{\pm, \tau}
+ A^{\pm, \tau} \circ D_{Y^\pm},
\label{A chain homotopy}
\end{equation}
that is, the following diagrams are commutative up to chain homotopy.
\[
\begin{tikzcd}
(\D_X^{\leq \kappa} / J^{\leq \kappa, \delta}_{C_0, C_1, C_2}
, D_{\F^0})\ar{r}{T^\tau}
& (\D_X^{\leq \kappa} / J^{\leq \kappa, \delta}_{C_0, C_1, C_2}, D_{\F^\tau})\\
(\W_{Y^\pm}^{\leq \kappa}
/ I^{\leq \kappa}_{C_0, C_1 + \kappa(\delta^{-1} - L^{-1}_{\min}), C_2}, D_{Y^\pm})
\ar{u}{i_{\F^0}^\pm} \ar{ur}[swap]{i_{\F^\tau}^\pm}&
\end{tikzcd}
\]
\end{enumerate}
\end{lem}
\begin{rem}
In the following proof, we need to take care of the degree with respect to $\hbar$.
Multiplication of $\F$ or $\K$ may decrease the degree at most by one,
but super-commutator $[ , ]$ increase the degree at least by one.
Hence in order to see that the linear maps defined below are well-defined,
we need to check that the number of super-commutators are greater than or equal to
the number of multiplications of $\F$ or $\K$.
\end{rem}
\begin{proof}
(i)
To prove that $T^\tau$ is a chain map,
it is enough to see that $(\Phi^\tau)^{-1} \widehat{D}_X \Phi^\tau
: \D\D_X^{\leq \kappa, \delta} / J^{\leq \kappa, \delta}_{C_0, C_1, C_2}
\to \D\D_X^{\leq \kappa, \delta} / J^{\leq \kappa, \delta}_{C_0, C_1, C_2}$ is
independent of $\tau \in [0, 1]$. This can be proved by
\begin{align*}
\frac{d}{d\tau} (\Phi^\tau)^{-1} \widehat{D}_X \Phi^\tau(f)
&= - (\Phi^\tau)^{-1} [\widehat{D}_X, \K^\tau] \widehat{D}_X \Phi^\tau(f)
+ (\Phi^\tau)^{-1} \widehat{D}_X [\widehat{D}_X, \K^\tau] \Phi^\tau(f)\\
&= 0.
\end{align*}
The latter claim is proved as follows.
If $\F^\tau$ and $\mathring{\K}^\tau$ also satisfy equation (\ref{main eq for X^I}),
then $\G^\tau = \mathring{\K}^\tau - \K^\tau \in \Ker D_{\F^\tau}$.
Let $\mathring{T}^\tau$ be the flow defined by $\F^\tau$ and $\mathring{\K}^\tau$.
Then
\begin{align*}
\frac{d}{d\tau} (\mathring{T}^\tau)^{-1} T^\tau (f)
&= - (\mathring{T}^\tau)^{-1} [[D_{\F^\tau}, \mathring{\K}^\tau], T^\tau(f)] (1)
+ (\mathring{T}^\tau)^{-1} [[D_{\F^\tau}, \K^\tau], T^\tau(f)] (1)\\
&= - (\mathring{T}^\tau)^{-1} [D_{\F^\tau}, \G^\tau] T^\tau(f)\\
&= - [D_{\F^0}, (\mathring{T}^\tau)^{-1} \G^\tau T^\tau] f.
\end{align*}
Therefore $(\mathring{T}^\tau)^{-1} T^\tau
: (\D_X^{\leq \kappa} / J^{\leq \kappa, \delta}_{C_0, C_1, C_2}, D_{\F^0})
\to (\D_X^{\leq \kappa} / J^{\leq \kappa, \delta}_{C_0, C_1, C_2}, D_{\F^0})$
are chain homotopic to the identity map for all $\tau$,
which implies the claim.

(ii)
Let $T^{\tau, \sigma} : (\hbar^{-1} \D_X^{\leq 0})^{\star, \delta}
/ J^{\star, \delta}_{\overline{C}_0, \overline{C}_1, \overline{C}_2}
\to (\hbar^{-1} \D_X^{\leq 0})^{\star, \delta}
/ J^{\star, \delta}_{\overline{C}_0, \overline{C}_1, \overline{C}_2}$ be the flow
defined similarly for each pair of one-parameter families
$(\F^{\tau, \sigma})_{\tau \in [0, 1]}$
and $(\K^{\tau, \sigma})_{\tau \in [0, 1]}$ satisfying (\ref{main eq for X F}).
Namely, they are defined by
\[
\partial_\tau T^{\tau, \sigma}(f)
= [[D_{\F^{\tau, \sigma}}, \K^{\tau, \sigma}], T^{\tau, \sigma}(f)] (1),
\quad T^{0, \sigma} = \id.
\]
Similarly to (i), each $T^{\tau, \sigma}$ is a chain map from
$((\hbar^{-1} \D_X^{\leq 0})^{\star, \delta}
/ J^{\star, \delta}_{\overline{C}_0, \overline{C}_1, \overline{C}_2}, D_{\F^{0, 0}})$
to $((\hbar^{-1} \D_X^{\leq 0})^{\star, \delta}
/ J^{\star, \delta}_{\overline{C}_0, \overline{C}_1, \overline{C}_2}, D_{\F^{\tau, \sigma}})$.
Hence it is enough to show that there exists a family of functions
$\mathring{\K}^{\tau, \sigma} \in (\hbar^{-1} \D_X^{\leq 0})^{\star, \delta}
/ J^{\star, \delta}_{\overline{C}_0, \overline{C}_1, \overline{C}_2}$ of odd degree
satisfying the following equations.
\begin{equation}
(T^{\tau, \sigma})^{-1} \partial_\sigma \F^{\tau, \sigma}
= D_{\F^{0, 0}}\bigl((T^{\tau, \sigma})^{-1} \mathring{\K}^{\tau, \sigma}\bigr)
\label{circle K condition}
\end{equation}
This is proved by the following calculations.
\begin{align}
&\partial_\tau
\bigl((T^{\tau, \sigma})^{-1} \partial_\sigma \F^{\tau, \sigma} \bigr) \notag\\
&= (T^{\tau, \sigma})^{-1}
\partial_\sigma \partial_\tau \F^{\tau, \sigma}
- (T^{\tau, \sigma})^{-1} \partial_\tau T^{\tau, \sigma} (T^{\tau, \sigma})^{-1}
\partial_\sigma \F^{\tau, \sigma} \notag\\
&= (T^{\tau, \sigma})^{-1}
\partial_\sigma (D_{\F^{\tau, \sigma}}\K^{\tau, \sigma})
- (T^{\tau, \sigma})^{-1} [[D_{\F^{\tau, \sigma}}, \K^{\tau, \sigma}],
\partial_\sigma \F^{\tau, \sigma}](1) \notag\\
& =(T^{\tau, \sigma})^{-1}
D_{\F^{\tau, \sigma}} \partial_\sigma \K^{\tau, \sigma}
+ (T^{\tau, \sigma})^{-1}
[D_{\F^{\tau, \sigma}}, \partial_\sigma\F^{\tau, \sigma}](\K^{\tau, \sigma}) \notag\\
&\quad - (T^{\tau, \sigma})^{-1}
[[D_{\F^{\tau, \sigma}}, \partial_\sigma\F^{\tau, \sigma}], \K^{\tau, \sigma}](1)
- (T^{\tau, \sigma})^{-1}
[D_{\F^{\tau, \sigma}}, [\K^{\tau, \sigma}, \partial_\sigma\F^{\tau, \sigma}]](1) \notag\\
&= D_{\F^{0, 0}}((T^{\tau, \sigma})^{-1} \partial_\sigma\K^{\tau, \sigma})
- (T^{\tau, \sigma})^{-1}
\K^{\tau, \sigma} D_{\F^{\tau, \sigma}}(\partial_\sigma\F^{\tau, \sigma}) \notag\\
&\quad - (T^{\tau, \sigma})^{-1}
[D_{\F^{\tau, \sigma}}, [\K^{\tau, \sigma}, \partial_\sigma\F^{\tau, \sigma}]](1) \notag\\
&= D_{\F^{0, 0}}((T^{\tau, \sigma})^{-1} \partial_\sigma\K^{\tau, \sigma})
\label{diff of sigma homotopy}
\end{align}
In the last equality, we have used the following facts:
\begin{itemize}
\item
$D_{\F^{\tau, \sigma}}(\partial_\sigma\F^{\tau, \sigma}) = 0$
because
\begin{align*}
D_{\F^{\tau, \sigma}}(\partial_\sigma\F^{\tau, \sigma})
&= e^{-\F^{\tau, \sigma}}
\widehat{D}_X(\partial_\sigma\F^{\tau, \sigma} e^{\F^{\tau, \sigma}})\\
&= e^{-\F^{\tau, \sigma}} \partial_\sigma \widehat{D}_X(e^{\F^{\tau, \sigma}})\\
&= 0.
\end{align*}
\item
$[\K^{\tau, \sigma}, \partial_\sigma\F^{\tau, \sigma}] = 0$ because multiplication
in $\D\D_X$ is super-commutative.
\end{itemize}
(\ref{diff of sigma homotopy}) implies that
\begin{equation}
\mathring{\K}^{\tau, \sigma}
:= T^{\tau, \sigma} \int_0^\tau
(T^{{\tau'}, \sigma})^{-1} \partial_\sigma\K^{{\tau'}, \sigma} d\tau'
\in (\hbar^{-1} \D_X^{\leq 0})^{\star, \delta}
/ J^{\star, \delta}_{\overline{C}_0, \overline{C}_1, \overline{C}_2}
\label{circle K}
\end{equation}
satisfies equation (\ref{circle K condition}).

(iii)
Let $\K^{\tau, \sigma}, \mathring{\K}^{\tau, \sigma}
\in (\hbar^{-1} \D_X^{\leq 0})^{\star, \delta}
/ J^{\star, \delta}_{\overline{C}_0, \overline{C}_1, \overline{C}_2}$ be families of functions
satisfying
\begin{align}
\partial_\tau\F^{\tau, \sigma} &= D_{\F^{\tau, \sigma}}(\K^{\tau, \sigma}),\\
\partial_\sigma\F^{\tau, \sigma} &= D_{\F^{\tau, \sigma}}(\mathring{\K}^{\tau, \sigma}).
\end{align}
Let $T^{\tau, \sigma}, \mathring{T}^{\tau, \sigma}
: \D_X^{\leq \kappa} / J^{\leq \kappa, \delta}_{C_0, C_1, C_2}
\to \D_X^{\leq \kappa} / J^{\leq \kappa, \delta}_{C_0, C_1, C_2}$ be flows defined by
\begin{align}
\partial_\tau T^{\tau, \sigma}(f)
&= [[D_{\F^{\tau, \sigma}}, \K^{\tau, \sigma}], T^{\tau, \sigma}(f)] (1),
\quad T^{0, \sigma} = \id,\\
\partial_\sigma \mathring{T}^{\tau, \sigma}(f)
&= [[D_{\F^{\tau, \sigma}}, \mathring{\K}^{\tau, \sigma}],
\mathring{T}^{\tau, \sigma}(f)] (1),
\quad \mathring{T}^{\tau, 0} = \id.
\end{align}
Since $\mathring{T}^{1, \sigma} = \id$, it is enough to prove that the chain maps
$(\mathring{T}^{\tau, \sigma})^{-1} T^{\tau, \sigma}$ from
$(\D_X^{\leq \kappa} / J^{\leq \kappa, \delta}_{C_0, C_1, C_2}, D_{\F^{0, 0}})$ to itself are
equal to the identity map up to chain homotopy for all $\sigma, \tau \in [0, 1]$.
The latter claim of (i) implies that we may assume that
$\K^{\tau, 0} = 0$ (since $\F^{\tau, 0} = \F^{0, 0}$)
and that the family $\mathring{\K}^{\tau, \sigma}$ is defined by (\ref{circle K}).
(\ref{circle K}) implies that the following equation holds in
$(\hbar^{-1} \D_X^{\leq 0})^{\star, \delta}
/ J^{\star, \delta}_{\overline{C}_0, \overline{C}_1, \overline{C}_2}$.
\begin{equation}
\partial_\sigma\K^{\tau, \sigma} - \partial_\tau\mathring{\K}^{\tau, \sigma}
+ [[D_{\F^\sigma_\tau}, \K^{\tau, \sigma}], \mathring{\K}^{\tau, \sigma}](1) = 0
\label{exactness of K}
\end{equation}

Define linear maps $S^{\tau, \sigma}, U^{\tau, \sigma}
: \D_X^{\leq \kappa} / J^{\leq \kappa, \delta}_{C_0, C_1, C_2}
\to \D_X^{\leq \kappa} / J^{\leq \kappa, \delta}_{C_0, C_1, C_2}$
($\tau, \sigma \in [0, 1]$) by
\begin{equation}
S^{\tau, \sigma}(f)
= (T^{\tau, \sigma})^{-1} \bigl(-\partial_\tau \mathring{T}^{\tau, \sigma}(f) +
[[D_{\F^{\tau, \sigma}}, \K^{\tau, \sigma}], \mathring{T}^{\tau, \sigma}(f)](1)\bigr)
\end{equation}
\begin{equation}
U^{\tau, \sigma}(f) = (\mathring{T}^{\tau, \sigma})^{-1}
[[[D_{\F^{\tau, \sigma}}, \K^{\tau, \sigma}], \mathring{K}^{\tau, \sigma}],
T^{\tau, \sigma} (f)](1)
\end{equation}
Then the following hold true.
\begin{enumerate}[label=\normalfont(\alph*)]
\item
Each $S^{\tau, \sigma}$ is a chain map from
$(\D_X^{\leq \kappa} / J^{\leq \kappa, \delta}_{C_0, C_1, C_2}, D_{\F^{0, 0}})$ to itself.
\item
$\partial_\sigma (\mathring{T}^{\tau, \sigma})^{-1} T^{\tau, \sigma}(f) \in
\D_X^{\leq \kappa} / J^{\leq \kappa, \delta}_{C_0, C_1, C_2}$
satisfies the following differential equation in variable $\tau$
for any $f \in \D_X^{\leq \kappa} / J^{\leq \kappa, \delta}_{C_0, C_1, C_2}$
and $\sigma \in [0, 1]$:
\begin{align}
&\partial_\tau\partial_\sigma (\mathring{T}^{\tau, \sigma})^{-1} T^{\tau, \sigma}(f)
= S^{\tau, \sigma} \partial_\sigma
\bigl((\mathring{T}^{\tau, \sigma})^{-1} T^{\tau, \sigma}(f)\bigr)\notag\\
&\hph{\partial_\tau\partial_\sigma (\mathring{T}^{\tau, \sigma})^{-1} T^{\tau, \sigma}(f) =}
+ D_{\F^{0, 0}} U(f) + U(D_{\F^{0, 0}} f), \label{eq of circle T T}\\
&\partial_\sigma (\mathring{T}^{\tau, \sigma})^{-1} T^{\tau, \sigma}(f)\bigr|_{\tau = 0} = 0.
\label{initial circle T T}
\end{align}
\end{enumerate}

First we prove the claim assuming the above two.
(a) implies that we can regard (\ref{eq of circle T T}) and (\ref{initial circle T T})
as equations of one-parameter families
$\partial_\sigma((\mathring{T}^{\tau, \sigma})^{-1} T^{\tau, \sigma}(\cdot))$
in the quotient space of chain maps from 
$(\D_X^{\leq \kappa} / J^{\leq \kappa, \delta}_{C_0, C_1, C_2}, D_{\F^{0, 0}})$ to itself
modulo null homotopies.
Then they become a linear differential equation with the trivial initial condition,
which implies that $\partial_\sigma \bigl((\mathring{T}^{\tau, \sigma})^{-1}
T^{\tau, \sigma}(\cdot)\bigr)$ is a family of null homotopies from
$(\D_X^{\leq \kappa} / J^{\leq \kappa, \delta}_{C_0, C_1, C_2}, D_{\F^{0, 0}})$ to itself.
Hence their integrations $(\mathring{T}^{\tau, \sigma})^{-1} T^{\tau, \sigma}(\cdot)$
are chain homotopic to $(\mathring{T}^{\tau, 0})^{-1} T^{\tau, 0}(\cdot) = \id$.

Now we prove the above two claims.
First we check (a).
By direct calculations, we see
\begin{align*}
D_{\F^{0, 0}} S^{\tau, \sigma}(f)
&= (T^{\tau, \sigma})^{-1} D_{\F^{\tau, \sigma}}
\bigl(-\partial_\tau \mathring{T}^{\tau, \sigma}(f)
+ [[D_{\F^{\tau, \sigma}}, \K^{\tau, \sigma}], \mathring{T}^{\tau, \sigma}(f)](1)\bigr) \\
&= (T^{\tau, \sigma})^{-1}
\bigl(-\partial_\tau(D_{\F^{\tau, \sigma}} \mathring{T}^{\tau, \sigma}(f))
+ [D_{\F^{\tau, \sigma}}, D_{\F^{\tau, \sigma}}(\K^{\tau, \sigma})]
(\mathring{T}^{\tau, \sigma}(f)) \\
&\hph{= (T^{\tau, \sigma})^{-1}
\bigl(}
+ D_{\F^{\tau, \sigma}}
([[D_{\F^{\tau, \sigma}}, \K^{\tau, \sigma}], \mathring{T}^{\tau, \sigma}(f)](1))
\bigr)
\\
&= (T^{\tau, \sigma})^{-1}
\bigl(-\partial_\tau(D_{\F^{\tau, \sigma}} \mathring{T}^{\tau, \sigma}(f))
- [D_{\F^{\tau, \sigma}}(\mathring{T}^{\tau, \sigma}(f)), D_{\F^{\tau, \sigma}}]
(\K^{\tau, \sigma})
\bigr)
\end{align*}
and
\begin{align*}
S^{\tau, \sigma} D_{\F^{0, 0}}(f)
&= (T^{\tau, \sigma})^{-1} \bigl(-\partial_\tau \mathring{T}^{\tau, \sigma}(D_{\F^{0, 0}}f) +
[[D_{\F^{\tau, \sigma}}, \K^{\tau, \sigma}],
D_{\F^{\tau, \sigma}}(\mathring{T}^{\tau, \sigma}(f))](1)\bigr) \\
&= (T^{\tau, \sigma})^{-1}
\bigl(-\partial_\tau(\mathring{T}^{\tau, \sigma} D_{\F^{0, 0}}(f))
- [D_{\F^{\tau, \sigma}}(\mathring{T}^{\tau, \sigma}(f)), D_{\F^{\tau, \sigma}}]
(\K^{\tau, \sigma})
\bigr)
\end{align*}
Hence $D_{\F^{0, 0}} S^{\tau, \sigma} = S^{\tau, \sigma} D_{\F^{0, 0}}$.

Next we prove (b).
This is also proved by direct calculation.
The key is equation (\ref{exactness of K}).
First we separate
\begin{align}
&\partial_\tau\partial_\sigma((\mathring{T}^{\tau, \sigma})^{-1}T^{\tau, \sigma}(f))
\notag\\
&= (\partial_\tau\partial_\sigma (\mathring{T}^{\tau, \sigma})^{-1}) T^{\tau, \sigma}(f)
+ (\mathring{T}^{\tau, \sigma})^{-1} \partial_\tau\partial_\sigma T^{\tau, \sigma}(f)
\notag\\
&\quad
+ (\partial_\sigma(\mathring{T}^{\tau, \sigma})^{-1}) \partial_\tau T^{\tau, \sigma}(f)
+ (\partial_\tau(\mathring{T}^{\tau, \sigma})^{-1}) \partial_\sigma T^{\tau, \sigma}(f)
\end{align}
into four parts and calculate each of them as follows.
\begin{align}
&(\partial_\tau\partial_\sigma (\mathring{T}^{\tau, \sigma})^{-1}) T^{\tau, \sigma}(f)
\notag\\
&= - \partial_\tau\bigl((\mathring{T}^{\tau, \sigma})^{-1}
[[D_{\F^{\tau, \sigma}}, \mathring{\K}^{\tau, \sigma}], \cdot](1)\bigr)(T^{\tau, \sigma}(f))
\notag\\
&= - (\mathring{T}^{\tau, \sigma})^{-1}
[[D_{\F^{\tau, \sigma}}, \partial_\tau \mathring{\K}^{\tau, \sigma}], T^{\tau, \sigma}(f)](1)
\notag \\
& \quad - (\mathring{T}^{\tau, \sigma})^{-1}
[[[D_{\F^{\tau, \sigma}}, D_{\F^{\tau, \sigma}}(\K^{\tau, \sigma})], \mathring{\K}^{\tau, \sigma}],
T^{\tau, \sigma}(f)](1) \notag\\
&\quad
+ (\mathring{T}^{\tau, \sigma})^{-1} \partial_\tau\mathring{T}^{\tau, \sigma}
(\mathring{T}^{\tau, \sigma})^{-1} [[D_{\F^{\tau, \sigma}}, \mathring{\K}^{\tau, \sigma}],
T^{\tau, \sigma}(f)](1)
\label{00 first}
\end{align}
\begin{align}
&(\mathring{T}^{\tau, \sigma})^{-1} \partial_\tau\partial_\sigma T^{\tau, \sigma}(f)
\notag\\
&= (\mathring{T}^{\tau, \sigma})^{-1}[[D_{\F^{\tau, \sigma}}, \K^{\tau, \sigma}],
\partial_\sigma T^{\tau, \sigma}(f)](1) \notag\\
& \quad
+ (\mathring{T}^{\tau, \sigma})^{-1}
[[D_{\F^{\tau, \sigma}}, \partial_\sigma \K^{\tau, \sigma}],
T^{\tau, \sigma}(f)](1) \notag\\
&\quad
+ (\mathring{T}^{\tau, \sigma})^{-1}
[[[D_{\F^{\tau, \sigma}}, D_{\F^{\tau, \sigma}}(\mathring{\K}^{\tau, \sigma})],
\K^{\tau, \sigma}],
T^{\tau, \sigma}(f)](1)
\label{00 second}
\end{align}
\begin{align}
&(\partial_\sigma(\mathring{T}^{\tau, \sigma})^{-1}) \partial_\tau T^{\tau, \sigma}(f)
\notag\\
&= - (\mathring{T}^{\tau, \sigma})^{-1} [[D_{\F^{\tau, \sigma}},
\mathring{\K}^{\tau, \sigma}],
[[D_{\F^{\tau, \sigma}}, \K^{\tau, \sigma}], T^{\tau, \sigma}(f)](1)](1)
\label{00 third}
\end{align}
\begin{align}
(\partial_\tau(\mathring{T}^{\tau, \sigma})^{-1}) \partial_\sigma T^{\tau, \sigma}(f)
= - (\mathring{T}^{\tau, \sigma})^{-1} \partial_\tau \mathring{T}^{\tau, \sigma}
(\mathring{T}^{\tau, \sigma})^{-1} \partial_\sigma T^{\tau, \sigma}(f)
\label{00 fourth}
\end{align}
We also calculate the following two.
\begin{align}
&- S^{\tau, \sigma} \partial_\sigma
\bigl((\mathring{T}^{\tau, \sigma})^{-1} T^{\tau, \sigma}(f)\bigr) \notag\\
&= (\mathring{T}^{\tau, \sigma})^{-1} \partial_\tau \mathring{T}^{\tau, \sigma}
(\mathring{T}^{\tau, \sigma})^{-1} \partial_\sigma T^{\tau, \sigma}(f) \notag\\
&\quad
- (\mathring{T}^{\tau, \sigma})^{-1} \partial_\tau \mathring{T}^{\tau, \sigma}
(\mathring{T}^{\tau, \sigma})^{-1} [[D_{\F^{\tau, \sigma}}, \mathring{K}^{\tau, \sigma}],
T^{\tau, \sigma}(f)](1) \notag\\
& \quad
- (\mathring{T}^{\tau, \sigma})^{-1} [[D_{\F^{\tau, \sigma}}, \K^{\tau, \sigma}],
\partial_\sigma T^{\tau, \sigma}(f)](1) \notag\\
&\quad
+ (\mathring{T}^{\tau, \sigma})^{-1} [[D_{\F^{\tau, \sigma}}, \K^{\tau, \sigma}],
[[D_{\F^{\tau, \sigma}}, \mathring{\K}^{\tau, \sigma}], T^{\tau, \sigma}(f)](1)](1)
\label{00 fifth}
\end{align}
\begin{align}
&-(D_{\F^{0, 0}} U(f) + U(D_{\F^{0, 0}} f)) \notag\\
&= -(\mathring{T}^{\tau, \sigma})^{-1}D_{\F^{\tau, \sigma}}(
[[[D_{\F^{\tau, \sigma}}, \K^{\tau, \sigma}], \mathring{K}^{\tau, \sigma}],
T^{\tau, \sigma} (f)](1)) \notag\\
& \quad
- (\mathring{T}^{\tau, \sigma})^{-1}
[[[D_{\F^{\tau, \sigma}}, \K^{\tau, \sigma}], \mathring{K}^{\tau, \sigma}],
D_{\F^{\tau, \sigma}}(T^{\tau, \sigma} (f))](1)
\notag\\
&= - (\mathring{T}^{\tau, \sigma})^{-1} [D_{\F^{\tau, \sigma}}, [[D_{\F^{\tau, \sigma}},
\K^{\tau, \sigma}], \mathring{\K}^{\tau, \sigma}]] T^{\tau, \sigma}(f) \notag\\
&\quad
+ (\mathring{T}^{\tau, \sigma})^{-1}
[D_{\F^{\tau, \sigma}}, [[D_{\F^{\tau, \sigma}}, \K^{\tau, \sigma}],
\mathring{\K}^{\tau, \sigma}](1)] T^{\tau, \sigma}(f) \notag\\
&= - (\mathring{T}^{\tau, \sigma})^{-1} [[D_{\F^{\tau, \sigma}}, \K^{\tau, \sigma}],
[D_{\F^{\tau, \sigma}}, \mathring{\K}^{\tau, \sigma}]] T^{\tau, \sigma}(f) \notag\\
&\quad
+ (\mathring{T}^{\tau, \sigma})^{-1}
[D_{\F^{\tau, \sigma}}, [[D_{\F^{\tau, \sigma}}, \K^{\tau, \sigma}],
\mathring{\K}^{\tau, \sigma}](1)] T^{\tau, \sigma}(f)
\label{00 sixth}
\end{align}
We need to show that the sum of (\ref{00 first}) to (\ref{00 sixth}) is zero.
The sum of the third term of (\ref{00 second}) and the second term of (\ref{00 first}) is
\begin{align}
&(\mathring{T}^{\tau, \sigma})^{-1}
[[[D_{\F^{\tau, \sigma}}, D_{\F^{\tau, \sigma}}(\mathring{\K}^{\tau, \sigma})],
\K^{\tau, \sigma}],
T^{\tau, \sigma}(f)](1) \notag\\
&- (\mathring{T}^{\tau, \sigma})^{-1}
[[[D_{\F^{\tau, \sigma}}, D_{\F^{\tau, \sigma}}(\K^{\tau, \sigma})],
\mathring{\K}^{\tau, \sigma}],
T^{\tau, \sigma}(f)](1) \notag\\
&= (\mathring{T}^{\tau, \sigma})^{-1}
([[D_{\F^{\tau, \sigma}}, D_{\F^{\tau, \sigma}}(\mathring{\K}^{\tau, \sigma})],
\K^{\tau, \sigma}] \notag\\
& \hph{= (\mathring{T}^{\tau, \sigma})^{-1}
(}
- [[D_{\F^{\tau, \sigma}}, D_{\F^{\tau, \sigma}}(\K^{\tau, \sigma})],
\mathring{\K}^{\tau, \sigma}]) T^{\tau, \sigma}(f) \notag\\
&\quad - (\mathring{T}^{\tau, \sigma})^{-1}
\bigl(T^{\tau, \sigma}(f)
([[D_{\F^{\tau, \sigma}}, D_{\F^{\tau, \sigma}}(\mathring{\K}^{\tau, \sigma})],
\K^{\tau, \sigma}] \notag\\
&\quad \hph{- (\mathring{T}^{\tau, \sigma})^{-1}
\bigl(T^{\tau, \sigma}(f)
(}
- [[D_{\F^{\tau, \sigma}}, D_{\F^{\tau, \sigma}}(\K^{\tau, \sigma})],
\mathring{\K}^{\tau, \sigma}])(1) \bigr) \notag\\
&= (\mathring{T}^{\tau, \sigma})^{-1}
(- [D_{\F^{\tau, \sigma}}(\mathring{\K}^{\tau, \sigma}), [D_{\F^{\tau, \sigma}},
\K^{\tau, \sigma}]] \notag\\
&\quad \hph{(\mathring{T}^{\tau, \sigma})^{-1}
(}
+ [D_{\F^{\tau, \sigma}}(\K^{\tau, \sigma}), [D_{\F^{\tau, \sigma}},
\mathring{\K}^{\tau, \sigma}]]) T^{\tau, \sigma}(f) \notag\\
&\quad - (\mathring{T}^{\tau, \sigma})^{-1} T^{\tau, \sigma}(f)
([D_{\F^{\tau, \sigma}}, D_{\F^{\tau, \sigma}}(\mathring{\K}^{\tau, \sigma})]
(\K^{\tau, \sigma}) \notag\\
&\quad \hph{- (\mathring{T}^{\tau, \sigma})^{-1} T^{\tau, \sigma}(f)
(}
- [D_{\F^{\tau, \sigma}}, D_{\F^{\tau, \sigma}}(\K^{\tau, \sigma})]
(\mathring{\K}^{\tau, \sigma})) \notag\\
&= (\mathring{T}^{\tau, \sigma})^{-1}
([D_{\F^{\tau, \sigma}}(\K^{\tau, \sigma}), [D_{\F^{\tau, \sigma}},
\mathring{\K}^{\tau, \sigma}]] \notag\\
&\quad \hph{(\mathring{T}^{\tau, \sigma})^{-1}
(}
- [D_{\F^{\tau, \sigma}}(\mathring{\K}^{\tau, \sigma}), [D_{\F^{\tau, \sigma}},
\K^{\tau, \sigma}]] ) T^{\tau, \sigma}(f) \notag\\
&\quad - (\mathring{T}^{\tau, \sigma})^{-1} T^{\tau, \sigma}(f)
D_{\F^{\tau, \sigma}}[[D_{\F^{\tau, \sigma}},
\K^{\tau, \sigma}], \mathring{\K}^{\tau, \sigma}](1)
\end{align}
The sum of the fourth term of (\ref{00 fifth}) and (\ref{00 third}) is
\begin{align}
&(\mathring{T}^{\tau, \sigma})^{-1} [[D_{\F^{\tau, \sigma}}, \K^{\tau, \sigma}],
[[D_{\F^{\tau, \sigma}}, \mathring{\K}^{\tau, \sigma}], T^{\tau, \sigma}(f)](1)](1) \notag\\
&- (\mathring{T}^{\tau, \sigma})^{-1} [[D_{\F^{\tau, \sigma}}, \mathring{\K}^{\tau, \sigma}],
[[D_{\F^{\tau, \sigma}}, \K^{\tau, \sigma}], T^{\tau, \sigma}(f)](1)](1) \notag\\
&= (\mathring{T}^{\tau, \sigma})^{-1} [[D_{\F^{\tau, \sigma}}, \K^{\tau, \sigma}],
[D_{\F^{\tau, \sigma}}, \mathring{\K}^{\tau, \sigma}](T^{\tau, \sigma}(f))](1) \notag\\
&\quad
- (\mathring{T}^{\tau, \sigma})^{-1} [[D_{\F^{\tau, \sigma}}, \K^{\tau, \sigma}],
D_{\F^{\tau, \sigma}}(\mathring{\K}^{\tau, \sigma}) T^{\tau, \sigma}(f)](1) \notag\\
&\quad
- (\mathring{T}^{\tau, \sigma})^{-1} [[D_{\F^{\tau, \sigma}}, \mathring{\K}^{\tau, \sigma}],
[D_{\F^{\tau, \sigma}}, \K^{\tau, \sigma}](T^{\tau, \sigma}(f))](1) \notag\\
&\quad
+ (\mathring{T}^{\tau, \sigma})^{-1} [[D_{\F^{\tau, \sigma}}, \mathring{\K}^{\tau, \sigma}],
D_{\F^{\tau, \sigma}}(\K^{\tau, \sigma}) T^{\tau, \sigma}(f)](1) \notag\\
&= (\mathring{T}^{\tau, \sigma})^{-1} [[D_{\F^{\tau, \sigma}}, \K^{\tau, \sigma}],
[D_{\F^{\tau, \sigma}}, \mathring{\K}^{\tau, \sigma}]] T^{\tau, \sigma}(f) \notag\\
&\quad - (\mathring{T}^{\tau, \sigma})^{-1}
\bigl([D_{\F^{\tau, \sigma}}(\K^{\tau, \sigma}), [D_{\F^{\tau, \sigma}},
\mathring{\K}^{\tau, \sigma}]] \notag\\
&\quad - [D_{\F^{\tau, \sigma}}(\mathring{\K}^{\tau, \sigma}), [D_{\F^{\tau, \sigma}},
\K^{\tau, \sigma}]] \bigr) T^{\tau, \sigma}(f)
\end{align}
The sum of the second term of (\ref{00 second}) and the first term of (\ref{00 first}) is
\begin{align}
& (\mathring{T}^{\tau, \sigma})^{-1}
[[D_{\F^{\tau, \sigma}}, \partial_\sigma \K^{\tau, \sigma}],
T^{\tau, \sigma}(f)](1) \notag\\
& \quad
- (\mathring{T}^{\tau, \sigma})^{-1}
[[D_{\F^{\tau, \sigma}}, \partial_\tau \mathring{\K}^{\tau, \sigma}], T^{\tau, \sigma}(f)](1)
\notag\\
&= (\mathring{T}^{\tau, \sigma})^{-1}
[[D_{\F^{\tau, \sigma}}, \partial_\sigma \K^{\tau, \sigma}
- \partial_\tau \mathring{\K}^{\tau, \sigma}],
T^{\tau, \sigma}(f)](1)
\notag\\
&= - (\mathring{T}^{\tau, \sigma})^{-1}
[[D_{\F^{\tau, \sigma}}, [[D_{\F^{\tau, \sigma}}, \K^{\tau, \sigma}],
\mathring{\K}^{\tau, \sigma}](1)], T^{\tau, \sigma}(f)](1)
\quad (\text{by (\ref{exactness of K})})
\notag\\
&= - (\mathring{T}^{\tau, \sigma})^{-1}
[D_{\F^{\tau, \sigma}}, [[D_{\F^{\tau, \sigma}}, \K^{\tau, \sigma}],
\mathring{\K}^{\tau, \sigma}](1)]T^{\tau, \sigma}(f) \notag\\
&\quad
+ (\mathring{T}^{\tau, \sigma})^{-1} T^{\tau, \sigma}(f)
D_{\F^{\tau, \sigma}}([[D_{\F^{\tau, \sigma}}, \K^{\tau, \sigma}],
\mathring{\K}^{\tau, \sigma}](1))
\end{align}
Therefore the sum of (\ref{00 first}) to (\ref{00 sixth}) is zero.

(iv)
We prove the existence of $A^{-, \tau}$.
Since
\begin{align*}
&\frac{d}{d \tau} (\mathring{T}^\tau)^{-1} \circ i_{\F^\tau}^-(f) \\
&= (\mathring{T}^\tau)^{-1} [f \underset{\F^\tau}{\overrightarrow{\ast}},
D_{\F^\tau}(\K^\tau)](1)
- (\mathring{T}^\tau)^{-1}
[[D_{\F^\tau}, \K^\tau], i_{\F^\tau}^-(f)](1) \\
&= - (\mathring{T}^\tau)^{-1} [[D_{\F^\tau}, \K^\tau],
f \underset{\F^\tau}{\overrightarrow{\ast}}](1) \\
&= - (T^\tau)^{-1}
[D_{\F^\tau}, [\K^\tau, f \underset{\F^\tau}{\overrightarrow{\ast}}]](1)
- (T^\tau)^{-1}
[\K^\tau, [D_{\F^\tau}, f \underset{\F^\tau}{\overrightarrow{\ast}}]](1) \\
&= - (T^\tau)^{-1} D_{\F^\tau}
\bigl([\K^\tau, f \underset{\F^\tau}{\overrightarrow{\ast}}](1)\bigr)
- (T^\tau)^{-1}
[\K^\tau, (D_{Y^-} f) \underset{\F^\tau}{\overrightarrow{\ast}}](1),
\end{align*}
(\ref{A chain homotopy}) is satisfied for
$A^{-, \tau} : \W^{\leq \kappa}_{Y^-}
/ I^{\leq \kappa}_{C_0, C_1 + \kappa(\delta^{-1} - L^{-1}_{\min}), C_2}
\to \D^{\leq \kappa}_X / J^{\leq \kappa, \delta}_{C_0, C_1, C_2}$ defined by
\[
A^{-, \tau}(h) = - \int_0^\tau (T^{\tau'})^{-1}
[\K^{\tau'}, h \underset{\F^{\tau'}}{\overrightarrow{\ast}}](1) d\tau'.
\]
We can similarly construct a family of chain homotopies
$A^{+, \tau}$.
\end{proof}

Next we consider the case of rational SFT.
Equation (\ref{main eq for X^I}) implies that the two families of
functions $\F^\tau_0, \K^\tau_0 \in
\mathcal{L}_X^{\leq 0} / J^{\leq 0}_{\overline{C}_0, \overline{C}_2}$ satisfy
\[
\frac{d}{d\tau} \F^\tau_0 = \delta \K^\tau_0
- \{\mathbf{h}, \K^\tau_0\}|_{\F^\tau_0}
\ (= d_{\F^\tau_0} \K^\tau_0)
\]
in $\mathcal{L}_X^{\leq 0}
/ J^{\leq 0}_{\overline{C}_0, \overline{C}_2}$, where
$\mathbf{h} = \mathcal{H}_{Y^-, 0} - \mathcal{H}_{Y^+, 0}$.
Namely, the family of functions $\F^\tau_0$ is a homotopy in the following sense.
\begin{defi}
One-parameter family of functions $F^\tau_0 \in \mathcal{L}_X^{\leq 0}
/ J^{\leq 0}_{\overline{C}_0, \overline{C}_2}$ ($\tau \in [0, 1]$)
of even degree is said to be a homotopy
if (\ref{main eq for X rational}) holds for all $\F_0 = \F^\tau_0$ and
there exists a family of functions $\K^\tau_0 \in \mathcal{L}_X^{\leq 0}
/ J^{\leq 0}_{\overline{C}_0, \overline{C}_2}$ of odd degree
which makes the following equation hold for all $\tau \in [0, 1]$.
\begin{equation}
\frac{d}{d\tau} \F^\tau_0 = d_{\F^\tau_0} \K^\tau_0 \label{main eq for X^I rational}
\end{equation}
\end{defi}
For each triple $(\kappa, C_0, C_2)$ such that $\overline{C}_0 \geq C_0$ and
$\overline{C}_2 \geq C_2 + \kappa$,
define a flow $T_0^\tau : \mathcal{L}^{\leq \kappa}_X / J^{\leq \kappa}_{C_0, C_2}
\to \mathcal{L}^{\leq \kappa}_X / J^{\leq \kappa}_{C_0, C_2}$ by
\[
\frac{d}{d\tau} T^\tau_0(f) = -\{\{\mathbf{h}, \K^\tau_0\},
T^\tau_0(f)\}|_{\F^\tau_0}.
\]
(It is related to the flow $T^\tau$ by $T^\tau_0(f) = T^\tau(f)|_{\hbar = 0}$.)
Then the following hold true as in the case of general SFT.
\begin{lem}
\begin{enumerate}[label = \normalfont (\roman*)]\ 
\item
$T^\tau_0$ is a chain map from
$(\mathcal{L}^{\leq \kappa}_X / J^{\leq \kappa}_{C_0, C_2}, d_{\F^0_0})$ to
$(\mathcal{L}^{\leq \kappa}_X / J^{\leq \kappa}_{C_0, C_2}, d_{\F^\tau_0})$
for each $\tau$.
Furthermore, up to chain homotopy, it is determined by
$(\F^{\tau'}_0)_{\tau' \in [0, \tau]}$
and independent of the choice of the family $(\K^{\tau'}_0)_{\tau' \in [0, \tau]}$
which satisfies equation $(\ref{main eq for X^I rational})$.
\item
If a smooth family of generating functions $\F^{\tau, \sigma}_0
\in \mathcal{L}_X^{\leq 0} / J^{\leq 0}_{\overline{C}_0, \overline{C}_2}$
$((\tau, \sigma) \in [0, 1] \times [0, 1])$
satisfies $\F^{0, \sigma}_0 \equiv \F^{0, 0}_0$ and the one-parameter family
$(\F^{\tau, \sigma}_0)_{\tau \in [0, 1]}$ is a homotopy for each $\sigma \in [0, 1]$,
then the one-parameter family $(\F^{\tau, \sigma}_0)_{\sigma \in [0, 1]}$ is also
a homotopy for each $\tau \in [0, 1]$.
\item
Further assume that the above family of generating functions satisfies
$\F^{1, \sigma}_0 = \F^{0, 0}_0$ and $\F^{\tau, 0}_0 \equiv \F^{0, 0}_0$.
Let $T^\tau_0 : \mathcal{L}^{\leq \kappa}_X / J^{\leq \kappa}_{C_0, C_2}
\to \mathcal{L}^{\leq \kappa}_X / J^{\leq \kappa}_{C_0, C_2}$
be the flow defined by the one-parameter homotopy $(\F^{\tau, 1}_0)_{\tau \in [0, 1]}$.
Then $T^1_0$ is equal to the identity map up to chain homotopy.
In other words, if a loop homotopy $(\F^{\tau, 1}_0)_{\tau \in S^1}$ is contractible
in the space of loop homotopies with the base point $\F^{0, 1}_0$,
then the chain map $T^1_0$ is the identity map up to chain homotopy.
Hence for a general one-parameter homotopy $(\F^\tau_0)_{\tau \in [0, 1]}$,
the end $T^1_0$ of the family of the chain maps $(T^\tau_0)_{\tau \in [0, 1]}$
is determined up to chain homotopy by the homotopy type of the homotopy
$(\F^\tau_0)_{\tau \in [0, 1]}$ relative to the end points.
\item
There exists a family of linear maps
$A^{\pm, \tau}_0 : \mathcal{P}_Y^{\leq \kappa} / I^{\leq \kappa}_{C_0, C_2}
\to \mathcal{L}^{\leq \kappa}_X / J^{\leq \kappa}_{C_0, C_2}$
such that
\begin{equation}
i_{\F^0_0}^\pm - (T^\tau_0)^{-1} \circ i_{\F^\tau_0}^\pm = d_{\F^0_0} \circ A^{\pm, \tau}_0
+ A^{\pm, \tau}_0 \circ d_{Y^\pm},
\label{A chain homotopy rational}
\end{equation}
that is, the following diagrams are commutative up to chain homotopy.
\[
\begin{tikzcd}
(\mathcal{L}^{\leq \kappa}_X / J^{\leq \kappa}_{C_0, C_2}, d_{\F^0_0}) \ar{rr}{T^\tau_0}
& & (\mathcal{L}^{\leq \kappa}_X / J^{\leq \kappa}_{C_0, C_2}, d_{\F^\tau_0})\\
&(\mathcal{P}_Y^{\leq \kappa} / I^{\leq \kappa}_{C_0, C_2}, d_{Y^\pm})
\ar{ul}{i_{\F^0_0}^\pm} \ar{ur}[swap]{i_{\F^\tau_0}^\pm}&
\end{tikzcd}
\]
\end{enumerate}
\end{lem}

Finally, we consider the case of contact homology.
Define
\[
\widehat{\mathcal{K}}^\tau_0 :=
\sum \frac{\overleftarrow{\partial} \K^\tau_0}{\partial p_{\hat c}^+} \bigg|_{p^+=0} \cdot
p^+_{\hat c}
\in \mathcal{L}_X^{\leq 0} / J^{\leq 0}_{\overline{C}_0, \overline{C}_2}.
\]
Then $\widehat{\mathcal{F}}^\tau_0$ and $\widehat{\mathcal{K}}^\tau_0$ satisfy
\begin{gather*}
\frac{d}{d\tau} \widehat{\mathcal{F}}^\tau_0 =
\delta \widehat{\mathcal{K}}^\tau_0
- \{\widehat{\mathbf{h}},
\widehat{\mathcal{K}}^\tau_0\} |_{\widehat{\mathcal{F}}^\tau_0},\\
\delta \widehat{\mathcal{F}}^\tau_0
= \widehat{\mathbf{h}}
|_{\widehat{\mathcal{F}}^\tau_0},
\end{gather*}
where $\widehat{\mathbf{h}}
= \widehat{\mathcal{H}}_{Y^-, 0} - \widehat{\mathcal{H}}_{Y^+, 0}$.
For pairs $(\kappa, C_0)$ such that $\overline{C}_0 \geq C_0$ and
$\overline{C}_2 \geq \kappa$, define linear maps
$\Delta^\tau : \A^{\leq \kappa}_{Y^+} / I^{\leq \kappa}_{C_0}
\to \A^{\leq \kappa}_{Y^-} / I^{\leq \kappa}_{C_0}$ by
\[
\Delta^\tau(f) = - \int_0^\tau \{\widehat{\mathcal{K}}^s_0, f\}|
_{\widehat{\mathcal{F}}^s_0} ds.
\]
Then the above equations imply that these are chain homotopies from $\Psi^\tau
= \Psi_{\widehat{\mathcal{F}}^\tau_0}$
to $\Psi^0$, that is,
\[
\Psi^\tau - \Psi^0 = \partial_{Y^-} \circ \Delta^\tau + \Delta^\tau \circ \partial_{Y^+}.
\]
In fact, $\widehat{\mathcal{F}}^\tau_0$ and $\widehat{\mathcal{K}}^\tau_0$
give a DGA homotopy in the sense that the following maps satisfy
the conditions of DGA homomorphism.
\begin{align}
\Psi : (\A^{\leq \kappa}_{Y^+} / I^{\leq \kappa}_{C_0}, \partial_{Y^+})
&\to (\Omega^\ast(I), d) \otimes
(\A^{\leq \kappa}_{Y^-} / I^{\leq \kappa}_{C_0}, \partial_{Y^-}) \label{DGA homotopy?}\\
f &\mapsto f|_{\widehat{\mathcal{F}}^\tau_0}
- d \tau \otimes \{\widehat{\mathcal{K}}^\tau_0, f\}|_{\widehat{\mathcal{F}}^\tau_0}
\notag
\end{align}
More precisely, $\Psi : \A^{\leq \kappa}_{Y^+} / I^{\leq \kappa}_{C_0} \to
\Omega^\ast_{C^N}(I) \otimes \A^{\leq \kappa}_{Y^-} / I^{\leq \kappa}_{C_0}$
is a linear map which satisfies
\[
(d \otimes 1 + (-1)^\ast \otimes \partial_{Y^-}) \Psi(f)
= \Psi(\partial_{Y^+} f)
\]
for $f \in \A^{\leq \kappa}_{Y^+} / I^{\leq \kappa}_{C_0}$, and
\[
\Psi(fg) = \Psi (f) \Psi(g)
\]
in $\Omega^\ast_{C^N}(I) \otimes \A^{\leq \kappa_1 + \kappa_2}_{Y^-}
/ I^{\leq \kappa_1 + \kappa_2}_{C_0}$
for all $f \in \A^{\leq \kappa_1}_{Y^+} / I^{\leq \kappa_1}_{C_0}$ and
$g \in \A^{\leq \kappa_2}_{Y^+} / I^{\leq \kappa_2}_{C_0}$
if $\overline{C}_0 \geq C_0$ and $\overline{C}_2 \geq \kappa_1 + \kappa_2$.
($\widehat{\mathcal{F}}^\tau_0$ is not of class $C^\infty$, but
$(\Omega^\ast(I), d)$ is a DGA of differential forms of class $C^\infty$.
Hence (\ref{DGA homotopy?}) is not strictly a DGA homomorphism.)


%% file: SFT-12_Composition.tex
%
%

\section{Composition}
\label{composition}
Let $X^- = (-\infty, 0] \times Y^- \cup Z^- \cup [0, \infty) \times Y^0$ and
$X^+ = (-\infty, 0] \times Y^0 \cup Z^+ \cup [0, \infty) \times Y^+$ be two
symplectic manifolds with cylindrical ends.
We regard them as symplectic cobordisms.
Then their composition $X = X^- \# X^+$ is defined by
\[
X = (-\infty, 0] \times Y^- \cup Z^- \cup Z^+ \cup [0, \infty) \times Y^+.
\]
Let $K_{X}^0$ be the set of cycles consisting of
\begin{itemize}
\item
cycles $x$ in $K_{X^-}^0$ such that $\supp x \cap [0, \infty) \times Y^0
= \emptyset$,
\item
cycles $x$ in $K_{X^+}^0$ such that $\supp x \cap (-\infty, 0] \times Y^0
= \emptyset$ and
\item
the cycles $x = x^- \# x^+$ obtained by the sums of the restrictions of
cycles $x^-$ in $K_{X^-}^0$ to $(-\infty, 0] \times Y^- \cup Z^-$
and the restrictions of cycles $x^+$ in $K_{X^+}^0$ to
$Z^+ \cup [0, \infty) \times Y^+$
corresponding to the same cycles $y$ in $K_{Y^0}^0$.
\end{itemize}
In this section, we prove that the composition of symplectic cobordisms corresponds to
the composition of the algebras.
First in Section \ref{composition of generating functions}, we recall the composition
of generating functions, and in Section \ref{composition of cobordisms},
we prove that the generating function of $X$ is homotopic to the composition of the
generating functions of $X^-$ and $X^+$.
In Section \ref{correction terms for composition}, we construct the correction terms
needed for Section \ref{composition of cobordisms}.

\subsection{Composition of generating functions}
\label{composition of generating functions}
In this section, we recall the definition of the composition of generating functions of
$X^-$ and $X^+$ and its linearizations defined in \cite{EGH00}.

First we consider the case of general SFT.
The composition map $\star : \D\D_{X^-} \otimes \D\D_{X^+}
\to \D\D_{X, Y^0}$ is defined by
\[
f \star g = (\overrightarrow{f} g)|_{q^0_{\hat c^\ast} = 0 \text{ for } c \in K_{Y^0}},
\]
where $\overrightarrow{f}$ is the differential operator obtained from $f$ by replacing
the variables $p^0_{\hat c}$ ($c \in K_{Y^0}$) with
$\hbar \overrightarrow{\frac{\partial}{\partial q^0_{\hat c^\ast}}}$,
and we replace the variables $t_{x^-}$ in $f$ and $t_{x^+}$ in $g$ with $t_{x^- \# t_x+}$.
(We denote the two variables corresponding to each simplex $c$ of $K_{Y^0}$
by $p^0_{\hat c}$ and $q^0_{\hat c^\ast}$.)
In the above definition, we regard $A \in \tilde \omega_{X^-}
H_2(\overline{X^-}, \partial \overline{X^-}; \Z) \cong
H_2(\overline{X^-}, \partial \overline{X^-}; \Z) / \Ker \tilde \omega_{X^-}$
of the variables $T^A$ appearing in $f$
and $B \in \tilde \omega_{X^+}
H_2(\overline{X^+}, \partial \overline{X^+}; \Z) \cong
H_2(\overline{X^+}, \partial \overline{X^+}; \Z) / \Ker \tilde \omega_{X^+}$ of the
variable $T^B$ in $g$
as elements of $H_2(\overline{X}, \partial{X}; \Z) \big/ (\Ker e \cap \Ker e_{Y^0})$
by the isomorphism
\begin{align*}
&H_2(\overline{X}, \partial \overline{X}) \big/ (\Ker e \cap \Ker e_{Y^0})\\
& \cong H_2(\overline{X}, Y^0 \cup \partial \overline{X}) \big/
(\Ker e \cap \Ker e_{Y^0})\\
& \cong H_2(Z^-, \partial Z^-) / \Ker \tilde \omega_{X^-} \oplus
H_2(Z^+, \partial Z^+) / \Ker \tilde \omega_{X^+}\\
& \cong H_2(\overline{X^-}, \partial \overline{X^-}) / \Ker \tilde \omega_{X^-} \oplus
H_2(\overline{X^+}, \partial \overline{X^+}) / \Ker \tilde \omega_{X^+}.
\end{align*}
(In the above equation, we use $\tilde \omega_{X^+} = e_{Y^0}$ and $\tilde \omega_{X^-}
+ \tilde \omega_{X^+} = e$.)

Note that the above composition map induces maps
\begin{multline*}
\star : \D\D_{X^-}^{\leq \kappa_1, \delta}
/ \widetilde{J}^{\leq \kappa_1, \delta}_{C_0, C_1 + \kappa_2 \delta^{-1}, C_2 + \kappa_2}
\otimes \D\D_{X^+}^{\leq \kappa_2, \delta}
/ \widetilde{J}^{\leq \kappa_2, \delta}_{C_0, C_1 + \kappa_1 \delta^{-1}, C_2}\\
\to \D\D_X^{\leq \kappa_1 + \kappa_2, \delta}
/ \widetilde{J}^{\leq \kappa_1 + \kappa_2, \delta}_{C_0, C_1, C_2}.
\end{multline*}
The composition
\[
\F^- \Diamond \F^+ \in (\hbar^{-1} \D_X^{\leq 0})^{\star, \delta}
/ J^{\star, \delta}_{\overline{C}_0, \overline{C}_1, \overline{C}_2}
\]
of generating functions
$\F^\pm \in (\hbar^{-1} \D_{X^\pm}^{\leq 0})^{\star, \delta}
/ J^{\star, \delta}_{\overline{C}_0, \overline{C}_1, \overline{C}_2}$ of $X^\pm$
are defined by
\[
e^{\F^- \Diamond \F^+} = e^{\F^-} \star e^{\F^+}
\]
in $\D\D_X^{\leq 0, \delta}
/ \widetilde{J}^{\leq 0, \delta}_{\overline{C}_0, \overline{C}_1, \overline{C}_2}$.
Then equations (\ref{main eq for X}) for $\F^\pm$ imply
\[
\widehat{D}_X (e^{\F^- \Diamond \F^+}) = 0.
\]
in $\D\D_X^{\leq 0, \delta}
/ \widetilde{J}^{\leq 0, \delta}_{\overline{C}_0, \overline{C}_1, \overline{C}_2}$.
In fact, any
$f^- \in \D\D_{X^-}^{\leq \kappa_1, \delta}
/ \widetilde{J}^{\leq \kappa_1, \delta}_{C_0, C_1 + \kappa_2 \delta^{-1}, C_2 + \kappa_2}$
and $f^+ \in \D\D_{X^+}^{\leq \kappa_2, \delta}
/ \widetilde{J}^{\leq \kappa_2, \delta}_{C_0, C_1 + \kappa_1 \delta^{-1}, C_2}$
satisfy
\begin{equation}
\widehat{D}_{X} (f^- \star f^+)
= (\widehat{D}_{X^-} f^-) \star f^+ + (-1)^{|f^-|} f^- \star (\widehat{D}_{X^+} f^+)
\label{star eq}
\end{equation}
in $\D\D_X^{\leq \kappa_1 + \kappa_2, \delta}
/ \widetilde{J}^{\leq \kappa_1 + \kappa_2, \delta}_{C_0, C_1, C_2}$.
More generally,
if $X^\pm$ contains contact manifolds $(Y^\pm_i, \lambda^\pm_i)$ as in
Section \ref{algebras with further energy conditions}, then
for $\F^\pm \in (\hbar^{-1} \D_{X^\pm, (Y^\pm_i)}^{\leq 0})^{\star, \delta}
/ J^{\star, \delta}_{\overline{C}_0, \overline{C}_1, \overline{C}_2}$,
we can define the composition
$\F^- \Diamond \F^+ \in (\hbar^{-1} \D_{X, (Y^-_i, Y^0, Y^+_i)}^{\leq 0})^{\star, \delta}
/ J^{\star, \delta}_{\overline{C}_0, \overline{C}_1, \overline{C}_2}$.

Define linear maps
$T_{\F^-}(\cdot \Diamond \F^+)
: \D^{\leq \kappa}_{X^-} / J^{\leq \kappa, \delta}_{C_0, C_1, C_2} \to
\D^{\leq \kappa}_{X, Y^0} / J^{\leq \kappa, \delta}_{C_0, C_1, C_2}$
and $T_{\F^+}(\F^- \Diamond \cdot) : \D^{\leq \kappa}_{X^+}
/ J^{\leq \kappa, \delta}_{C_0, C_1, C_2} \to
\D^{\leq \kappa}_{X, Y^0} / J^{\leq \kappa, \delta}_{C_0, C_1, C_2}$ by
\begin{align*}
T_{\F^-}( \cdot \Diamond \F^+)(f)
&= e^{- \F^- \Diamond \F^+}((f e^{\F^-}) \star e^{\F^+})\\
T_{\F^+}(\F^- \Diamond \cdot)(f)
&= e^{- \F^- \Diamond \F^+}(e^{\F^-} \star (f e^{\F^+})).
\end{align*}
(These are the linearizations of the composition map.)
We also define a map
\begin{align*}
&T^2_{\F^-, \F^+} : \D^{\leq \kappa_1}_{X^-}
/ J^{\leq \kappa_1, \delta}_{C_0, C_1 + \kappa_2 \delta^{-1}, C_2 + \kappa_2}
\otimes  \D^{\leq \kappa_2}_{X^+}
/ J^{\leq \kappa_2, \delta}_{C_0, C_1 + \kappa_1 \delta^{-1}, C_2} \\
&\hspace{170pt}
\to \D^{\leq \kappa_1 + \kappa_2}_{X, Y^0}
/ J^{\leq \kappa_1 + \kappa_2, \delta}_{C_0, C_1, C_2}
\end{align*}
by
\[
T^2_{\F^-, \F^+}(f \otimes g) = e^{- \F^- \Diamond \F^+} ((f e^{\F^-}) \star (g e^{\F^+})).
\]
Note that $T_{\F^-}(\cdot \Diamond \F^+)(f) = T^2_{\F^-, \F^+}(f \otimes 1)$
and $T_{\F^+}(\F^- \Diamond \cdot)(f) = T^2_{\F^-, \F^+}(1 \otimes f)$.
Some of the following properties of these maps were proved in \cite{EGH00}.
\begin{lem}\label{linearized T}
The linearizations of the composition map satisfy the following.
\begin{enumerate}[label=\normalfont(\roman*)]
\item
They are chain maps, that is,
\begin{align*}
T_{\F^-}(\cdot \Diamond \F^+) \circ D_{\F^-}
&= D_{\F^- \Diamond \F^+} \circ T_{\F^-}(\cdot \Diamond \F^+),\\
T_{\F^+}(\F^- \Diamond \cdot) \circ D_{\F^+}
&= D_{\F^- \Diamond \F^+} \circ T_{\F^+}(\F^- \Diamond \cdot).
\end{align*}
More generally,
\[
T^2_{\F^-, \F^+} \circ (D_{\F^-} \otimes 1 + (-1)^\ast \otimes D_{\F^+})
= D_{\F^- \Diamond \F^+} \circ T^2_{\F^-, \F^+}.
\]
\item
They satisfy the following compatibility conditions with $i_{\F^\pm}^\pm$ and
$i_{\F^- \Diamond \F^+}^\pm$.
{\belowdisplayskip=0pt
\begin{multline*}
T_{\F^-}(\cdot \Diamond \F^+) \circ i_{\F^-}^- = i_{\F^- \Diamond \F^+}^-\\
: \W_{Y^-}^{\leq \kappa}
/ I^{\leq \kappa}_{C_0, C_1 + \kappa(\delta^{-1} - L_{\min}^{-1}), C_2}
\to \D_{X, Y^0}^{\leq \kappa, \delta} / J^{\leq \kappa, \delta}_{C_0, C_1, C_2},
\end{multline*}
}{\abovedisplayskip=0pt
\begin{multline*}
T_{\F^+}(\F^- \Diamond \cdot) \circ i_{\F^+}^+ = i_{\F^- \Diamond \F^+}^+\\
: \W_{Y^+}^{\leq \kappa}
/ I^{\leq \kappa}_{C_0, C_1 + \kappa(\delta^{-1} - L_{\min}^{-1}), C_2}
\to \D_{X, Y^0}^{\leq \kappa, \delta} / J^{\leq \kappa, \delta}_{C_0, C_1, C_2}.
\end{multline*}
}
More generally, they are compatible with the multiplication as follows.
For any $f \in \W_{Y^-}^{\leq \kappa_1} / I^{\leq \kappa}_{C_0, C'_1, C_2 + \kappa_1}$
and $g \in \D_{X^-}^{\leq \kappa_2} / J^{\leq \kappa_2, \delta}_{C_0, C_1
+ \kappa_1 \delta^{-1}, C_2}$,
\[
T_{\F^-}(\cdot \Diamond \F^+) (f \underset{\F^-}{\overrightarrow{\ast}} g)
= f \underset{\F^- \Diamond \F^+}{\overrightarrow{\ast}}
(T_{\F^-}(\cdot \Diamond \F^+) g)
\]
in $\D_{X, Y^0}^{\leq \kappa_1 + \kappa_2, \delta}
/ J^{\leq \kappa_1 + \kappa_2, \delta}_{C_0, C_1, C_2}$,
where $C'_1 = C_1 + \kappa_1 (\delta^{-1} - L_{\min}^{-1}) + \kappa_2 \delta^{-1}$.
For any $g \in \D_{X^+}^{\leq \kappa_2} / J^{\leq \kappa_2, \delta}_{C_0, C_1
+ \kappa_1 \delta^{-1}, C_2 + \kappa_1}$ and
$f \in \W_{Y^+}^{\leq \kappa_1} / I^{\leq \kappa_1}_{C_0, C'_1, C_2}$,
\[
T_{\F^+}(\F^- \Diamond \cdot) (g \underset{\F^+}{\overleftarrow{\ast}} f)
= (T_{\F^+}(\F^- \Diamond \cdot) g)
\underset{\F^- \Diamond \F^+}{\overleftarrow{\ast}} f
\]
in $\D_{X, Y^0}^{\leq \kappa_1 + \kappa_2, \delta}
/ J^{\leq \kappa_1 + \kappa_2, \delta}_{C_0, C_1, C_2}$.
\item
They satisfy the following compatibility condition with $i^{\pm}_{\F^\mp}$.
\begin{align*}
T_{\F^-}(\cdot \Diamond \F^+) \circ i^+_{\F^-}
&= T_{\F^+}(\F^- \Diamond \cdot) \circ i^-_{\F^+}\\
&: \W_{Y^0}^{\leq \kappa}
/ I^{\leq \kappa}_{C_0, C_1 + \kappa(\delta^{-1} - L_{\min}^{-1}), C_2}
\to \D_{X, Y^0}^{\leq \kappa, \delta} / J^{\leq \kappa, \delta}_{C_0, C_1, C_2}
\end{align*}
More generally, they are compatible with the multiplication:
\[
T^2_{\F^-, \F^+}((g \underset{\F^-}{\overleftarrow{\ast}} f) \otimes h)
= T^2_{\F^-, \F^+}(g \otimes (f \underset{\F^+}{\overrightarrow{\ast}} h))
\]
in $\D^{\leq \kappa_1 + \kappa_2 + \kappa_3}_{X, Y^0}
/ J^{\leq \kappa_1 + \kappa_2 + \kappa_3, \delta}_{C_0, C_1, C_2}$
for any
\begin{align*}
f &\in \W_{Y^0}^{\leq \kappa_1}
/ I^{\leq \kappa_1}_{C_0, C_1 + \kappa_1 (\delta^{-1} - L_{\min}^{-1})
+ (\kappa_2 + \kappa_3) \delta^{-1},
C_2 + \kappa_3}, \\
g &\in \D^{\leq \kappa_2}_{X^-} / J^{\leq \kappa_2, \delta}_{C_0, C_1
+ (\kappa_1 + \kappa_3) \delta^{-1}, C_2 + \kappa_1 + \kappa_3}, \\
h &\in  \D^{\leq \kappa_3}_{X^+} / J^{\leq \kappa_3, \delta}_{C_0, C_1
+ (\kappa_1 + \kappa_2) \delta^{-1}, C_2}.
\end{align*}
\item
Let $X^i$ ($i = 1, 2, 3$) be symplectic cobordisms from $Y^{i-1}$ to $Y^i$,
and let $\F^i$ be a generating function for each $X^i$.
Then
\begin{align*}
T_{\F^1 \Diamond \F^2}(\cdot \Diamond \F^3) \circ T_{\F^1}(\cdot \Diamond \F^2)
&= T_{\F^1}(\cdot \Diamond (\F^2 \Diamond \F^3)),\\
T_{\F^2 \Diamond \F^3}(\F^1 \Diamond \cdot) \circ T_{\F^3}(\F^2 \Diamond \cdot)
&= T_{\F^3}((\F^1 \Diamond \F^2) \Diamond \cdot).
\end{align*}
More generally,
\[
T^2_{\F^1 \Diamond \F^2, \F^3} \circ (T^2_{\F^1, \F^2} \otimes 1)
= T^2_{\F^1, \F^2 \Diamond \F^3} \circ (1 \otimes T^2_{\F^2, \F^3}).
\]
\item
Let $(\F^{\pm, \tau}, \K^{\pm, \tau})$ be homotopies of generating functions for
$X^\pm$.
Then
\[
(\F ^\tau = \F^{-, \tau} \Diamond \F^{+, \tau},\ \K^\tau
= T^2_{\F^{-, \tau}, \F^{+, \tau}}
(\K^{-, \tau} \otimes 1 + 1 \otimes \K^{+, \tau}))
\]
is a homotopy of generating functions of $X$.
Furthermore, there exist families of linear maps $A^{\pm, \tau}
: \D_{X^\pm}^{\leq \kappa} / J^{\leq \kappa, \delta}_{C_0, C_1, C_2}
\to \D_{X, Y^0}^{\leq \kappa} / J^{\leq \kappa, \delta}_{C_0, C_1, C_2}$
such that
\begin{multline*}
(T^\tau)^{-1} \circ T_{\F^{-, \tau}}(\cdot \Diamond \F^{+, \tau}) \circ T^{-, \tau}
- T_{\F^{-, 0}}(\cdot \Diamond \F^{+, 0})\\
= D_{\F^{-, 0} \Diamond \F^{+, 0}} \circ A^{-, \tau} + A^{-, \tau} \circ D_{\F^{-, 0}},
\end{multline*}
\begin{multline*}
(T^\tau)^{-1} \circ T_{\F^{+, \tau}}(\F^{-, \tau} \Diamond \cdot) \circ T^{+, \tau}
- T_{\F^{+, 0}}(\F^{-, 0} \Diamond \cdot)\\
= D_{\F^{-, 0} \Diamond \F^{+, 0}} \circ A^{+, \tau} + A^{+, \tau} \circ D_{\F^{+, 0}},
\end{multline*}
where $T^{\pm, \tau}$ and $T^\tau$ are the flows for the homotopies
$(\F^{\pm, \tau}, \K^{\pm, \tau})$ and $(\F^\tau, \K^\tau)$ respectively.
Namely, the following diagrams are commutative up to chain homotopy.
\[
\begin{tikzcd}[column sep= huge]
(\D_{X^-}^{\leq \kappa} / J^{\leq \kappa, \delta}_{C_0, C_1, C_2}, D_{\F^{-, 0}})
\ar{r}{T_{\F^{-, 0}}(\cdot \Diamond \F^{+, 0})}
\ar{d}{T^{-, \tau}}
& (\D_{X, Y^0}^{\leq \kappa} / J^{\leq \kappa, \delta}_{C_0, C_1, C_2},
D_{\F^{-, 0} \Diamond \F^{+, 0}}) \ar{d}{T^{\tau}}\\
(\D_{X^-}^{\leq \kappa} / J^{\leq \kappa, \delta}_{C_0, C_1, C_2}, D_{\F^{-, \tau}})
\ar{r}{T_{\F^{-, \tau}}(\cdot \Diamond \F^{+, \tau})}
& (\D_{X, Y^0}^{\leq \kappa} / J^{\leq \kappa, \delta}_{C_0, C_1, C_2},
D_{\F^{-, \tau} \Diamond \F^{+, \tau}})
\end{tikzcd}
\]
\[
\begin{tikzcd}[column sep= huge]
(\D_{X^+}^{\leq \kappa} / J^{\leq \kappa, \delta}_{C_0, C_1, C_2}, D_{\F^{+, 0}})
\ar{r}{T_{\F^{+, 0}}(\F^{-, 0} \Diamond \cdot)}
\ar{d}{T^{+, \tau}}
& (\D_{X, Y^0}^{\leq \kappa} / J^{\leq \kappa, \delta}_{C_0, C_1, C_2},
D_{\F^{-, 0} \Diamond \F^{+, 0}}) \ar{d}{T^{\tau}}\\
(\D_{X^+}^{\leq \kappa} / J^{\leq \kappa, \delta}_{C_0, C_1, C_2}, D_{\F^{+, \tau}})
\ar{r}{T_{\F^{+, \tau}}(\F^{-, \tau} \Diamond \cdot)}
& (\D_{X, Y^0}^{\leq \kappa} / J^{\leq \kappa, \delta}_{C_0, C_1, C_2},
D_{\F^{-, \tau} \Diamond \F^{+, \tau}})
\end{tikzcd}
\]
More generally, there exists a family of linear maps
\begin{multline*}
A^\tau
: \D_{X^-}^{\leq \kappa_1}
/ J^{\leq \kappa_1, \delta}_{C_0, C_1 + \kappa_2 \delta^{-1}, C_2 + \kappa_2}
\otimes \D_{X^+}^{\leq \kappa_2}
/ J^{\leq \kappa_2, \delta}_{C_0, C_1 + \kappa_1 \delta^{-1}, C_2}\\
\to \D_{X, Y^0}^{\leq \kappa_1 + \kappa_2, \delta}
/ J^{\leq \kappa_1 + \kappa_2, \delta}_{C_0, C_1, C_2}
\end{multline*}
such that
\begin{multline*}
(T^\tau)^{-1} \circ T^2_{\F^{-, \tau}, \F^{+, \tau}} \circ (T^{-, \tau} \otimes T^{+, \tau})
- T^2_{\F^{-, 0}, \F^{+, 0}} \\
= D_{\F^{-, 0} \Diamond \F^{+, 0}} \circ A^\tau
+ A^\tau \circ (D_{\F^{-, 0}} \otimes 1 + (-1)^\ast \otimes D_{\F^{+, 0}}).
\end{multline*}
\end{enumerate}
\end{lem}
\begin{proof}
(i) is due to (\ref{star eq}).
(ii), (iii) and (iv) are straightforward.
(v) is proved as follows.
Using (\ref{star eq}), we can easily check that $(\F^\tau, \K^\tau)$ is a homotopy.
We construct $A^\tau$.
For any $f \in \D_{X^-}^{\leq \kappa_1}
/ J^{\leq \kappa_1, \delta}_{C_0, C_1 - \kappa_2 \delta^{-1}, C_2 + \kappa_2}$
and $g \in \D_{X^+}^{\leq \kappa_2}
/ J^{\leq \kappa_2, \delta}_{C_0, C_1 - \kappa_1 \delta^{-1}, C_2}$,
\begin{align*}
&\frac{d}{d \tau} (T^\tau)^{-1} \circ T^2_{\F^{-, \tau}, \F^{+, \tau}} \circ
(T^{-, \tau} \otimes T^{+, \tau})(f \otimes g)\\
&= (T^\tau)^{-1} \circ T^2_{\F^{-, \tau}, \F^{+, \tau}}
([[D_{\F^{-, \tau}}, \K^{-, \tau}], T^{-, \tau}(f)](1) \otimes T^{+, \tau}(g))\\
& \quad + (T^\tau)^{-1} \circ T^2_{\F^{-, \tau}, \F^{+, \tau}}
(T^{-, \tau}(f) \otimes [[D_{\F^{+, \tau}}, \K^{+, \tau}], T^{+, \tau}(g)](1))\\
&\quad + (T^\tau)^{-1} T^2_{\F^{-, \tau}, \F^{+, \tau}}
(D_{\F^{-, \tau}}(\K^{-, \tau}) T^{-, \tau}(f) \otimes T^{+, \tau}(g)) \\
&\quad + (T^\tau)^{-1} T^2_{\F^{-, \tau}, \F^{+, \tau}}
(T^{-, \tau}(f) \otimes D_{\F^{+, \tau}}(\K^{+, \tau}) T^{+, \tau}(g))\\
&\quad - (T^\tau)^{-1} D_{\F^\tau}(\K^\tau)
T^2_{\F^{-, \tau}, \F^{+, \tau}}(T^{-, \tau}(f) \otimes T^{+, \tau}(g))\\
&\quad - (T^\tau)^{-1} [[D_{\F^\tau}, \K^\tau],
T^2_{\F^{-, \tau}, \F^{+, \tau}}(T^\tau(f) \otimes T^{+, \tau}(g))](1)\\
&= (T^\tau)^{-1} T^2_{\F^{-, \tau}, \F^{+, \tau}}
([D_{\F^{-, \tau}}, \K^{-, \tau}] T^{-, \tau}(f) \otimes T^{+, \tau}(g)) \\
&\quad + (T^\tau)^{-1} T^2_{\F^{-, \tau}, \F^{+, \tau}}
(T^{-, \tau}(f) \otimes [D_{\F^{+, \tau}}, \K^{+, \tau}] T^{+, \tau}(g)) \\
&\quad - (T^\tau)^{-1}
[D_{\F^\tau}, \K^\tau](T^2_{\F^{-, \tau}, \F^{+, \tau}}
(T^{-, \tau}(f) \otimes T^{+, \tau}(g))\\
&= D_{\F^0} S^\tau(f) + S^\tau(D_{\F^{-, 0}}f \otimes g
+ (-1)^{|f|} f \otimes D_{\F^{+, 0}}g)\\
\end{align*}
where
\begin{align*}
S^\tau(f \otimes g)
&= (T^\tau)^{-1}\bigl( T^2_{\F^-, \F^+}
(\K^{-, \tau} T^{-, \tau}(f) \otimes T^{+, \tau}(g))\\
&\quad \hph{(T^\tau)^{-1}\bigl(}
+ (-1)^{|f|} T^2_{\F^-, \F^+}(T^{-, \tau}(f) \otimes \K^{+, \tau} T^{+, \tau}(g))\\
&\quad \hph{(T^\tau)^{-1}\bigl(}
- \K^\tau T^2_{\F^-, \F^+}(T^{-, \tau}(f) \otimes T^{+, \tau}(g)) \bigr).
\end{align*}
Therefore,
\[
A^\tau(f \otimes g) = \int_0^\tau S^{\tau'}(f \otimes g) d \tau'
\]
is a required family of linear maps.
\end{proof}

Next we consider rational SFT.
The composition $\F^-_0 \sharp \F^+_0 \in \mathcal{L}^{\leq 0}_{X, Y^0}
/ J^{\leq 0}_{\overline{C}_0, \overline{C}_2}$ of
generating functions $\F^\pm \in \mathcal{L}^{\leq 0}_{X^\pm}
/ J^{\leq 0}_{\overline{C}_0, \overline{C}_2}$ is defined by
\begin{align*}
\F^-_0 \sharp \F^+_0 &= ((\hbar^{-1} \F^-_0 \Diamond \hbar^{-1} \F^+)
\cdot \hbar)|_{\hbar = 0}\\
&= (\F^- \Diamond \F^+)_0.
\end{align*}
We define linear maps
\begin{align*}
T_{\F^-_0}(\cdot \sharp \F^+_0) &: \mathcal{L}^{\leq \kappa}_{X^-}
/ J^{\leq \kappa}_{C_0, C_2}
\to \mathcal{L}^{\leq \kappa}_{X, Y^0} / J^{\leq \kappa}_{C_0, C_2}, \\
T_{\F^+_0}(\F^-_0 \sharp \cdot) &: \mathcal{L}^{\leq \kappa}_{X^+}
/ J^{\leq \kappa}_{C_0, C_2}
\to \mathcal{L}^{\leq \kappa}_{X, Y^0} / J^{\leq \kappa}_{C_0, C_2}, \\
(T_0)^2_{\F^-_0, \F^+_0} &: \mathcal{L}^{\leq \kappa_1}_{X^-}
/ J^{\leq \kappa_1}_{C_0, C_2 + \kappa_2} \otimes \mathcal{L}^{\leq \kappa_2}_{X^+}
/ J^{\leq \kappa_2}_{C_0, C_2} \to
\mathcal{L}^{\leq \kappa_1 + \kappa_2}_{X, Y^0}
/ J^{\leq \kappa_1 + \kappa_2}_{C_0, C_2}
\end{align*}
by
\begin{align*}
T_{\F^-_0}(\cdot \sharp \F^+_0)(f)
&= T_{\hbar^{-1}\F^-_0}(\cdot \Diamond \hbar^{-1}\F^+_0)(f)|_{\hbar = 0},\\
T_{\F^+_0}(\F^-_0 \sharp \cdot)(f)
&= T_{\hbar^{-1}\F^+_0}(\hbar^{-1}\F^-_0 \Diamond \cdot)(f)|_{\hbar = 0},\\
(T_0)^2_{\F^-_0, \F^+_0}(f \otimes g)
&= T^2_{\hbar^{-1} \F^-_0, \hbar^{-1} \F^+_0}(f \otimes g)|_{\hbar = 0}.
\end{align*}
Then they satisfy the counterpart of Lemma \ref{linearized T}.

Finally we consider the case of contact homology.
Note that
\[
\widehat{\mathcal{F}^-_0 \sharp \mathcal{F}^+_0} =
\sum_{c}
\frac{\overleftarrow{\partial} (\mathcal{F}^-_0 \sharp \mathcal{F}^+_0)}
{\partial p^+_{\hat c}} \biggr|_{p^+ = 0} \cdot p^+_{\hat c}
= \widehat{\mathcal{F}}^+_0\Bigr|_{q^0_{\hat c^\ast} =
\frac{\overleftarrow{\partial} \widehat{\mathcal{F}}^-_0}{\partial p^0_{\hat c}}}.
\]
This implies that the composition $\Psi_{\widehat{\mathcal{F}}^-_0} \circ
\Psi_{\widehat{\mathcal{F}}^+_0} : \A_{Y^+}^{\leq \kappa} / I^{\leq \kappa}_{C_0} \to
\A_{Y^-}^{\leq \kappa} / I^{\leq \kappa}_{C_0}$ coincides with the chain map
defined by $\widehat{\F^-_0 \sharp \F^+_0}$.

\subsection{Composition of cobordisms}\label{composition of cobordisms}
In this section, we construct a homotopy between the generating function of $X$
and the composition of the generating functions of $X^-$ and $X^+$.

For each $0 \leq T < \infty$, a new manifold $X^T$ is defined by
\[
X^T = (-\infty, 0] \times Y^- \cup Z^- \cup ([0, T]_{0^-} \cup [-T, 0]_{0^+}) \times Y^0
\cup Z^+ \cup [0, \infty) \times Y^+,
\]
where we identify $T \in [0, T]_{0^-}$ with $-T \in [-T, 0]_{0^+}$.

First we define a holomorphic building for $X^{[0, \infty]}$.
\begin{defi}
A holomorphic building $(T, \Sigma, z, u, \phi)$ for $X^{[0, \infty]}$
consists of the following:
\begin{itemize}
\item
$0 \leq T \leq \infty$
\item
A marked curve $(\Sigma, z)$ which is obtained from some union of
marked semistable curves $(\check \Sigma, z \cup (\pm \infty_i))$
with a floor structure.
In this case, floor takes values in $\{-k_-, \dots, -1, 0, 1, \dots, k_+\}$
if $0 \leq T < \infty$,
and $\{-k_-, \dots, -1,\ab 0^-,\ab 0_1, \dots, 0_l, 0^+, 1, \dots, k_+\}$ ($l \geq 0$)
if $T = \infty$.
\item
If $T < \infty$, then $u$ is a continuous map
$u : \Sigma \to (\overline{\R}_{-k_-} \cup \dots \cup \overline{\R}_{-1}) \times Y^-
\cup X^T \cup (\overline{\R}_1 \cup \dots \cup \overline{\R}_{k_+}) \times Y^+$,
and if $T = \infty$, then $u$ is a continuous map
$u : \Sigma \to (\overline{\R}_{-k_-} \cup \dots \cup \overline{\R}_{-1}) \times Y^-
\cup X^- \cup (\overline{\R}_{0_1} \cup \dots \cup \overline{\R}_{0_l})
\times Y^0 \cup X^+ \cup (\overline{\R}_{1} \cup \dots \cup
\overline{\R}_{k_+}) \times Y^+$.
\item
$\phi_{\pm\infty_i} : S^1 \to S^1_{\pm\infty_i}$ is a family of coordinates of limit circles.
\end{itemize}
We assume the following conditions:
If $T < \infty$, then $(\Sigma, z, u, \phi)$ is a holomorphic building for $X^T$.
In this case, the energies $E_\lambda(u)$ and $E_{\hat \omega}(u)$ are defined by
\begin{align*}
E_\lambda(u) &= \max \biggl\{
\sup_{I \subset \R_{-k_-} \cup \dots \cup \R_{-1} \cup (-\infty, 0]}
\frac{1}{|I|} \int_{u^{-1}(I \times Y^-)} u^\ast (d\sigma \wedge \lambda^-),\\
&\hphantom{= \max \biggl\{}
\sup_{I \subset [0, T] \cup [-T, 0]} \frac{1}{|I|} \int_{u^{-1}(I \times Y^0)}
u^\ast (d\sigma \wedge \lambda^0),\\
&\hphantom{= \max \biggl\{}
\sup_{I \subset [0, \infty) \cup \R_1 \cup \dots \cup \R_{k_+}}
\frac{1}{|I|} \int_{u^{-1}(I \times Y^+)} u^\ast (d\sigma \wedge \lambda^+) \biggr\},\\
E_{\hat \omega}(u) &= \int_{u^{-1}(X^T)} u^\ast \hat \omega^T
+ \int_{u^{-1}((\overline{\R}_{-k_-} \cup \dots \cup \overline{\R}_{-1}) \times Y^-)} u^\ast d\lambda^-\\
&\quad + \int_{u^{-1}((\overline{\R}_1 \cup \dots \cup \overline{\R}_{k_+}) \times Y^+)} u^\ast d\lambda^+,
\end{align*}
where $\hat \omega^T$ is defined by
$\hat \omega^T|_{Z^\pm} = \omega^\pm$,
$\hat \omega^T|_{(-\infty, 0] \times Y^-} = d\lambda^-$,
$\hat \omega^T|_{([0, T] \cup [-T, 0]) \times Y^0} = d\lambda^0$, and
$\hat \omega^T |_{[0, \infty) \times Y^+} = d\lambda^+$.

If $T = \infty$, then we assume that $(\Sigma, z, u, \phi)$ satisfies the following conditions:
\begin{itemize}
\item
If $i(\alpha) <0^-$ then
$u(\Sigma_\alpha \setminus \coprod S^1) \subset \R_{i(\alpha)} \times Y^-$, and
$u|_{\Sigma_\alpha \setminus \coprod S^1} : \Sigma_\alpha \setminus \coprod S^1 \to \R_{i(\alpha)} \times Y^-$ is $J$-holomorphic.
\item
If $i(\alpha) = 0^-$ then
$u(\Sigma_\alpha \setminus \coprod S^1) \subset X^-$, and
$u|_{\Sigma_\alpha \setminus \coprod S^1} : \Sigma_\alpha \setminus \coprod S^1 \to X^-$
is $J$-holomorphic.
\item
If $0_1 \leq i(\alpha) \leq 0_l$ then
$u(\Sigma_\alpha \setminus \coprod S^1) \subset \R_{i(\alpha)} \times Y^0$, and
$u|_{\Sigma_\alpha \setminus \coprod S^1} : \Sigma_\alpha \setminus \coprod S^1 \to \R_{i(\alpha)} \times Y^0$ is $J$-holomorphic.
\item
If $i(\alpha) = 0^+$ then
$u(\Sigma_\alpha \setminus \coprod S^1) \subset X^+$, and
$u|_{\Sigma_\alpha \setminus \coprod S^1} : \Sigma_\alpha \setminus \coprod S^1 \to X^+$
is $J$-holomorphic.
\item
If $i(\alpha) >0^+$ then
$u(\Sigma_\alpha \setminus \coprod S^1) \subset \R_{i(\alpha)} \times Y^+$, and
$u|_{\Sigma_\alpha \setminus \coprod S^1} : \Sigma_\alpha \setminus \coprod S^1 \to \R_{i(\alpha)} \times Y^+$ is $J$-holomorphic.
\item
The energies $E_\lambda(u) <\infty$ and $E_{\hat \omega}(u) <\infty$ are finite
which are defined by
\begin{align*}
E_\lambda(u) &= \max \biggl\{
\sup_{I \subset \R_{-k_-} \cup \dots \cup \R_{-1} \cup (-\infty, 0]}
\frac{1}{|I|} \int_{u^{-1}(I \times Y^-)} u^\ast (d\sigma \wedge \lambda^-),\\
&\hphantom{= \max \biggl\{}
\sup_{I \subset [0, \infty)_{0^-} \cup \R_{0_1} \cup \dots \cup \R_{0_l} \cup (-\infty, 0]_{0^+}}
\frac{1}{|I|} \int_{u^{-1}(I \times Y^0)}
u^\ast (d\sigma \wedge \lambda^0),\\
&\hphantom{= \max \biggl\{}
\sup_{I \subset [0, \infty) \cup \R_1 \cup \dots \cup \R_{k_+}}
\frac{1}{|I|} \int_{u^{-1}(I \times Y^+)} u^\ast (d\sigma \wedge \lambda^+) \biggr\}, \\
E_{\hat \omega}(u) &= \int_{u^{-1}(X^-)} u^\ast \hat \omega^-
+ \int_{u^{-1}(X^+)} u^\ast \hat \omega^+\\
&\quad + \int_{u^{-1}((\overline{\R}_{-k_-} \cup \dots \cup \overline{\R}_{-1}) \times Y^-)} u^\ast d\lambda^-
+ \int_{u^{-1}((\overline{\R}_{0_1} \cup \dots \cup \overline{\R}_{0_l}) \times Y^0)}
u^\ast d\lambda^0\\
&\quad + \int_{u^{-1}((\overline{\R}_1 \cup \dots \cup \overline{\R}_{k_+}) \times Y^+)} u^\ast d\lambda^+.
\end{align*}
\item
$u$ is positively asymptotic to a periodic orbit $\gamma_{+\infty_i} = \pi_Y \circ u \circ \phi_{+\infty_i} \in P_{Y^+}$ at each $S^1_{+\infty_i}$,
and negatively asymptotic to a periodic orbit $\gamma_{-\infty_i} = \pi_Y \circ u \circ \phi_{-\infty_i} \in P_{Y^-}$ at each $S^1_{-\infty_i}$.
At every joint circle, $u$ is positively asymptotic to a periodic orbit
on the side of lower floor and negatively asymptotic to the same periodic orbit
on the side of higher floor.
\item
For each component $\hat\Sigma_\alpha$, if $u|_{\Sigma_\alpha}$ is a constant map,
then $2g_\alpha + m_\alpha \geq 3$.
\item
For each $i \neq 0^\pm$, the $i$-th floor $u^{-1}(\overline{\R}_i \times Y^\pm) \subset
\Sigma$ (or $u^{-1}(\overline{\R}_i \times Y^0) \subset \Sigma$) contains
nontrivial components.
\end{itemize}
\end{defi}

We denote the space of holomorphic buildings for $X^{[0, \infty]}$ by
$\overline{\M}_{X^{[0, \infty]}}$.
We define $\widehat{\M}_{X^{[0, \infty]}}$ similarly.
Kuranishi neighborhoods of $\widehat{\M}_{X^{[0, \infty]}}$ are defined in a similar way to
those of $\widehat{\M}_X$.

Let $\beta_0 > 0$ be the constant used for the definition of the strong differential
structure of the space of deformation of the domain curves in the case of $Y^0$
in order to construct smooth Kuranishi neighborhoods.
(Recall that in Section \ref{smoothness}, we use the coordinate change
$\rho_\mu = \hat \rho_\mu^{L_\mu^{-1} \beta}$ to define the strong differential
structure. We denote this $\beta$ in the case of $Y^0$ by $\beta_0$.)
We define the differential structure of $[0, \infty]$ by
the coordinate $\varphi : [0, \infty] \to [0, 1]$ given by
$\varphi(T) = \exp(-2T / \beta_0)$.
Then it is easy to check that the natural map
$\widehat{\M}_{X^{[0, \infty]}} \to [0, \infty]$ is a strong smooth map.
The associated maps from the Kuranishi neighborhoods to $[0, \infty]$ is
not submersive in general, but if we regard $[0, \infty]$ as
a $1$-dimensional simplex, then these maps are essentially submersive.
(Locally, at the corner corresponding to the splitting of $Y^0$-floors,
these maps looks like
$(t_1, \dots, t_k) \to t_1 \cdots t_k : [0, 1)^k \to [0, 1)$.)


First we define a space $\widehat{\M}^\diamond_{X^{[0, \infty]}}$.
Its point $(T, (\Sigma^\alpha, z^\alpha, u^\alpha)_{\alpha \in A_Y^- \sqcup A_X^0
\sqcup A_Y^+}, M^{\mathrm{rel}})$ consists of $0 \leq T \leq \infty$,
holomorphic buildings
$(\Sigma^\alpha, z^\alpha, u^\alpha)_{\alpha \in A_Y^-}$ for $Y^-$,
$(\Sigma^\alpha, \ab z^\alpha, \ab u^\alpha)_{\alpha \in A_X^0}$ for $X^T$
(that is, $(T, \Sigma^\alpha, z^\alpha, u^\alpha)_{\alpha \in A_X^0}$ are
holomorphic buildings for $X^{[0, \infty]}$),
$(\Sigma^\alpha, z^\alpha, u^\alpha)_{\alpha \in A_Y^+}$ for $Y^+$, and
a set $M^{\mathrm{rel}} = \{(S^1_{+\infty_l}, S^1_{-\infty_l})\}$ of pairs of
limit circles which satisfy the following conditions:
\begin{itemize}
\item
Any two pairs in $M^{\mathrm{rel}}$ do not share the same limit circle.
\item
For each pair $\alpha_1, \alpha_2 \in A = A_Y^- \sqcup A_X^0 \sqcup A_Y^+$,
let $M^{\alpha_1, \alpha_2} \subset M^{\mathrm{rel}}$ be the subset of pairs
$(S^1_{+\infty_l}, S^1_{-\infty_l})$ such that $S^1_{+\infty_l}$ is a $+\infty$-limit circle
of $\Sigma^{\alpha_1}$ and $S^1_{-\infty_l}$ is a $-\infty$-limit circle of
$\Sigma^{\alpha_2}$.
Then there does not exists any sequence
$\alpha_0, \alpha_1, \dots, \alpha_k = \alpha_0 \in A$ such that
$M^{\alpha_i, \alpha_{i+1}} \neq \emptyset$ for all $i = 0, 1, \dots, k-1$.
\item
For subsets $A_1, A_2 \subset A$, define
$M^{(A_1, A_2)} = \bigcup_{\alpha_1 \in A_1, \alpha_2 \in A_2} M^{\alpha_1, \alpha_2}$.
Then $M^{\mathrm{rel}}$ is the union of
$M^{\mathrm{rel}, \leq 0} = M^{(A_Y^-, A_Y^- \sqcup A_X^0)}$
and $M^{\mathrm{rel}, \geq 0} = M^{(A_X^0 \sqcup A_Y^+, A_Y^+)}$.
\end{itemize}
This is not a usual pre-Kuranishi structure but is a fiber product with
the diagonal in the product of $[0, \infty]$ by the essential submersion.

We also define a space $\widehat{\M}^\diamond_{X^-, X^+}$ as follows.
Its point
\[
((\Sigma^\alpha, z^\alpha, u^\alpha)_{\alpha \in A_Y^- \sqcup A_X^-
\sqcup A_Y^0 \sqcup A_X^+ \sqcup A_Y^+}, M^{\mathrm{rel}})
\]
consists of holomorphic buildings
$(\Sigma^\alpha, z^\alpha, u^\alpha)_{\alpha \in A_Y^-}$ for $Y^-$,
$(\Sigma^\alpha, z^\alpha, u^\alpha)_{\alpha \in A_X^-}$ for $X^-$,
$(\Sigma^\alpha, z^\alpha, u^\alpha)_{\alpha \in A_Y^0}$ for $Y^0$,
$(\Sigma^\alpha, z^\alpha, u^\alpha)_{\alpha \in A_X^+}$ for $X^+$,
$(\Sigma^\alpha, z^\alpha, u^\alpha)_{\alpha \in A_Y^+}$ for $Y^+$, and
a set $M^{\mathrm{rel}} = \{(S^1_{+\infty_l}, S^1_{-\infty_l})\}$ of pairs of
limit circles which satisfy the following conditions:
\begin{itemize}
\item
Any two pairs in $M^{\mathrm{rel}}$ do not share the same limit circle.
\item
For each pair $\alpha_1, \alpha_2 \in A = A_Y^- \sqcup A_X^- \sqcup A_Y^0
\sqcup A_X^+ \sqcup A_Y^+$,
let $M^{\alpha_1, \alpha_2} \subset M^{\mathrm{rel}}$ be the subset of pairs
$(S^1_{+\infty_l}, S^1_{-\infty_l})$ such that $S^1_{+\infty_l}$ is a $+\infty$-limit circle
of $\Sigma^{\alpha_1}$ and $S^1_{-\infty_l}$ is a $-\infty$-limit circle of
$\Sigma^{\alpha_2}$.
Then there does not exists any sequence
$\alpha_0, \alpha_1, \dots, \alpha_k = \alpha_0 \in A$ such that
$M^{\alpha_i, \alpha_{i+1}} \neq \emptyset$ for all $i = 0, 1, \dots, k-1$.
\item
For subsets $A_1, A_2 \subset A$, define
$M^{(A_1, A_2)} = \bigcup_{\alpha_1 \in A_1, \alpha_2 \in A_2} M^{\alpha_1, \alpha_2}$.
Then $M^{\mathrm{rel}}$ is the union of
$M^{\mathrm{rel}, -} = M^{(A_Y^-, A_Y^- \sqcup A_X^-)}$,
$M^{\mathrm{rel}, 0} = M^{(A_X^- \sqcup A_Y^0, A_Y^0 \sqcup A_X^+)}$
and $M^{\mathrm{rel}, +} = M^{(A_X^+ \sqcup A_Y^+, A_Y^+)}$.
\end{itemize}
The definition of the connected points of $\widehat{\M}^\diamond_{X^{[0, \infty]}}$ and
$\widehat{\M}^\diamond_{X^-, X^+}$ is similarly to the case of $X$ and $Y$.
$\widehat{\M}^\diamond_{X^{[0, \infty]}}$ contains
$[0, \infty] \times \widehat{\M}^\diamond_{Y^\pm}$, and
$\widehat{\M}^\diamond_{X^-, X^+}$
contains $\widehat{\M}^\diamond_{Y^\pm}$, $\widehat{\M}^\diamond_{Y^0}$ and
$\widehat{\M}^\diamond_{X^\pm}$.

We define 
\begin{align*}
&(\widehat{\M}^\diamond_{X^{[0, \infty]}}, (\mathring{K}_{Y^-}^2,
\mathring{K}_{Y^0}^2, \mathring{K}_{Y^+}^2),
(K_{Y^-}, K_{Y^+}), \\
& \hspace{50pt}
(K^0_{Y^-}, K^0_{X^-}, K^0_{Y^0}, K^0_{X^+}, K^0_{Y^+})) \subset
\widehat{\M}^\diamond_{X^{[0, \infty]}}
\end{align*}
and
\begin{align*}
&(\widehat{\M}^\diamond_{X^-, X^+}, (\mathring{K}_{Y^-}^2,
\mathring{K}_{Y^0}^2, \mathring{K}_{Y^+}^2),
(K_{Y^-}, K_{Y^0}, K_{Y^+}), \\
& \hspace{50pt}
(K^0_{Y^-}, K^0_{X^-}, K^0_{Y^0}, K^0_{X^+}, K^0_{Y^+})) \subset
\widehat{\M}^\diamond_{X^-, X^+}
\end{align*}
similarly.
We abbreviate the above two spaces by
$(\widehat{\M}^\diamond_{X^{[0, \infty]}}, \mathring{K}_{X^{[0, \infty]}}, \ab
K_{X^{[0, \infty]}}, \ab K^0_{X^{[0, \infty]}})$ and
$(\widehat{\M}^\diamond_{X^-, X^+}, \mathring{K}_{X^-, X^+},
K_{X^-, X^+}, K^0_{X^-, X^+})$ respectively.

Let $(\widehat{\M}^\diamond_{X^\infty}, \mathring{K}_{X^\infty},
K_{X^\infty}, K^0_{X^\infty})$ be the fiber product of
$(\widehat{\M}^\diamond_{X^{[0, \infty]}}, \mathring{K}_{X^{[0, \infty]}}, \ab
K_{X^{[0, \infty]}}, \ab K^0_{X^{[0, \infty]}})$ with $\infty \in [0, \infty]$.
Then both of $(\widehat{\M}^\diamond_{X^{[0, \infty]}}, \mathring{K}_{X^{[0, \infty]}}, \ab
K_{X^{[0, \infty]}}, \ab K^0_{X^{[0, \infty]}})$ and
$(\widehat{\M}^\diamond_{X^-, X^+}, \mathring{K}_{X^-, X^+},
K_{X^-, X^+}, K^0_{X^-, X^+})$ contains
$(\widehat{\M}^\diamond_{X^\infty}, \mathring{K}_{X^\infty}, \ab
K_{X^\infty}, \ab K^0_{X^\infty})$.

Similarly to the usual case, we can define multi-valued partial essential submersions
$\Xi^\circ$ and $\Lambda$ for both of
$(\widehat{\M}^\diamond_{X^{[0, \infty]}}, \mathring{K}_{X^{[0, \infty]}}, \ab
K_{X^{[0, \infty]}}, \ab K^0_{X^{[0, \infty]}})$ and
$(\widehat{\M}^\diamond_{X^-, X^+}, \ab \mathring{K}_{X^-, X^+}, \ab
K_{X^-, X^+}, \ab K^0_{X^-, X^+})$.

We can construct a continuous family of grouped multisections of
$(\widehat{\M}^\diamond_{X^{[0, \infty]}}, \ab \mathring{K}_{X^{[0, \infty]}}, \ab
K_{X^{[0, \infty]}}, \ab K^0_{X^{[0, \infty]}})$ and
a usual grouped multisection of
$(\widehat{\M}^\diamond_{X^-, X^+}, \ab \mathring{K}_{X^-, X^+}, \ab
K_{X^-, X^+}, \ab K^0_{X^-, X^+})$
which satisfy the following conditions.
\begin{itemize}
\item
On a neighborhood of $T \in \{0\} \cup \{\infty\} \subset [0, \infty]$,
the continuous family of grouped multisections of
$(\widehat{\M}^\diamond_{X^{[0, \infty]}}, \ab \mathring{K}_{X^{[0, \infty]}}, \ab
K_{X^{[0, \infty]}}, \ab K^0_{X^{[0, \infty]}})$ is a usual grouped multisection.
\item
The restrictions of the natural map
$(\widehat{\M}^\diamond_{X^{[0, \infty]}}, \ab \mathring{K}_{X^{[0, \infty]}}, \ab
K_{X^{[0, \infty]}}, \ab K^0_{X^{[0, \infty]}}) \to [0, \infty]$
to the fiber products of the zero sets of the perturbed multisection
are essentially submersive.
\item
The restrictions of the continuous family of
grouped multisections of
$(\widehat{\M}^\diamond_{X^{[0, \infty]}}, \ab \mathring{K}_{X^{[0, \infty]}}, \ab
K_{X^{[0, \infty]}}, \ab K^0_{X^{[0, \infty]}})$ to
$[0, \infty] \times
(\widehat{\M}_{Y^\pm}^\diamond, \ab \mathring{K}_{Y^\pm}^2, \ab
K_{Y^\pm}, K_{Y^\pm}^0)$
coincide with the pull backs of the grouped multisections of
$(\widehat{\M}_{Y^\pm}^\diamond, \mathring{K}_{Y^\pm}^2, K_{Y^\pm}, K_{Y^\pm}^0)$
by the projection.
\item
The restrictions of the grouped multisection of
$(\widehat{\M}^\diamond_{X^-, X^+}, \ab \mathring{K}_{X^-, X^+}, \ab
K_{X^-, X^+}, \ab K^0_{X^-, X^+})$ to
$(\widehat{\M}_{Y^\pm}^\diamond, \ab \mathring{K}_{Y^\pm}^2, \ab
K_{Y^\pm}, K_{Y^\pm}^0)$,
$(\widehat{\M}_{Y^0}^\diamond, \ab \mathring{K}_{Y^0}^2, \ab
K_{Y^0}, K_{Y^0}^0)$ and
$(\widehat{\M}^\diamond_{X^\pm}, \ab \mathring{K}^2_{X^\pm}, \ab
K_{X^\pm}, K^0_{X^\pm})$
coincide with the given grouped multisections.
\item
The restrictions of the 
grouped multisections of
$(\widehat{\M}^\diamond_{X^{[0, \infty]}}, \ab \mathring{K}_{X^{[0, \infty]}}, \ab
K_{X^{[0, \infty]}}, \ab K^0_{X^{[0, \infty]}})$ and
$(\widehat{\M}^\diamond_{X^-, X^+}, \ab \mathring{K}_{X^-, X^+}, \ab
K_{X^-, X^+}, \ab K^0_{X^-, X^+})$ to
$(\widehat{\M}^\diamond_{X^\infty}, \ab \mathring{K}_{X^\infty}, \ab
K_{X^\infty}, \ab K^0_{X^\infty})$ coincide.
\item
Let
\begin{align*}
&((\widehat{\M}^\diamond_{X^{[0, \infty]}})^0, \mathring{K}^2_{X^{[0, \infty]}},
K_{X^{[0, \infty]}}, K^0_{X^{[0, \infty]}}) \\
&\quad \subset
(\widehat{\M}^\diamond_{X^{[0, \infty]}}, \mathring{K}^2_{X^{[0, \infty]}}, \ab
K_{X^{[0, \infty]}}, \ab K^0_{X^{[0, \infty]}})
\end{align*}
be the subset of connected points.
Its continuous family of grouped multisections induces that of
\[
\bigcup_N (\prod^N (\widehat{\M}^\diamond_{X^{[0, \infty]}})^0,
\mathring{K}^2_{X^{[0, \infty]}}, K_{X^{[0, \infty]}},
K^0_{X^{[0, \infty]}})_{[0, \infty]} / \mathfrak{S}_N,
\]
where each $(\prod^N (\widehat{\M}^\diamond_{X^{[0, \infty]}})^0,
\mathring{K}^2_{X^{[0, \infty]}}, K_{X^{[0, \infty]}},
K^0_{X^{[0, \infty]}})_{[0, \infty]}$ is the fiber product over ${[0, \infty]}$.
Then the grouped multisection of
$(\widehat{\M}^\diamond_{X^{[0, \infty]}}, \mathring{K}^2_{X^{[0, \infty]}},
\ab K_{X^{[0, \infty]}}, \ab K^0_{X^{[0, \infty]}})$
coincides with its pull back by the essential submersion
\begin{align*}
&(\widehat{\M}^\diamond_{X^{[0, \infty]}}, \mathring{K}^2_{X^{[0, \infty]}},
K_{X^{[0, \infty]}}, K^0_{X^{[0, \infty]}}) \\
&\to
\bigcup_N (\prod^N (\widehat{\M}^\diamond_{X^{[0, \infty]}})^0,
\mathring{K}^2_{X^{[0, \infty]}}, K_{X^{[0, \infty]}},
K^0_{X^{[0, \infty]}})_{[0, \infty]} / \mathfrak{S}_N
\end{align*}
defined by decomposition into connected components.
\item
Let
\begin{align*}
&((\widehat{\M}^\diamond_{X^-, X^+})^0, \mathring{K}^2_{X^-, X^+},
K_{X^-, X^+}, K^0_{X^-, X^+}) \\
&\quad \subset
(\widehat{\M}^\diamond_{X^-, X^+}, \mathring{K}^2_{X^-, X^+},
\ab K_{X^-, X^+}, \ab K^0_{X^-, X^+})
\end{align*}
be the subset of connected points.
Its grouped multisection induces that of
\[
\bigcup_N (\prod^N (\widehat{\M}^\diamond_{X^-, X^+})^0,
\mathring{K}^2_{X^-, X^+}, K_{X^-, X^+},
K^0_{X^-, X^+}) / \mathfrak{S}_N.
\]
Then the grouped multisection of
$(\widehat{\M}^\diamond_{X^-, X^+}, \mathring{K}^2_{X^-, X^+},
K_{X^-, X^+}, K^0_{X^-, X^+})$
coincides with its pull back by the submersion
\begin{align*}
&(\widehat{\M}^\diamond_{X^-, X^+}, \mathring{K}^2_{X^-, X^+},
K_{X^-, X^+}, K^0_{X^-, X^+}) \\
&\to
\bigcup_N (\prod^N (\widehat{\M}^\diamond_{X^-, X^+})^0,
\mathring{K}^2_{X^-, X^+}, K_{X^-, X^+},
K^0_{X^-, X^+}) / \mathfrak{S}_N
\end{align*}
defined by decomposition into connected components.
\item
The (continuous families of) grouped multisections of
$(\widehat{\M}^\diamond_{X^{[0, \infty]}}, \ab \mathring{K}_{X^{[0, \infty]}}, \ab
K_{X^{[0, \infty]}}, \ab K^0_{X^{[0, \infty]}})$ and
$(\widehat{\M}^\diamond_{X^-, X^+}, \ab \mathring{K}_{X^-, X^+}, \ab
K_{X^-, X^+}, \ab K^0_{X^-, X^+})$ are compatible with respect to
the compatible systems of
multi-valued partial essential submersions defined by $\Xi^\circ$ and $\Lambda$.
\end{itemize}

We define pre-Kuranishi spaces
$\overline{\M}^{(m_-, X^{[0, \infty]}, m_+)}
_{(\hat \epsilon^{i, j}_l, \hat c^i_l, x^i_l, \hat \eta^i_l)}$ and
$\overline{\M}^{(m_-, X^-, m, X^+, m_+)}
_{(\hat \epsilon^{i, j}_l, \hat c^i_l, x^i_l, \hat \eta^i_l)}$ similarly
to $\overline{\M}^{(m_-, X^I, m_+)}_{((\hat \epsilon^{i, j}_l), (\hat c^i_l), (x^i_l),
(\hat \eta^i_l))}$ and
$\overline{\M}^{(m_-, X, m_+)}_{((\hat \epsilon^{i, j}_l), (\hat c^i_l), (x^i_l), (\hat \eta^i_l))}$,
and define their grouped multisections (or a continuous family of grouped multisections
for the former)
by the pull back by the natural maps to
$(\widehat{\M}^\diamond_{X^{[0, \infty]}}, \mathring{K}_{X^{[0, \infty]}}, \ab
K_{X^{[0, \infty]}}, \ab K^0_{X^{[0, \infty]}})$ and
$(\widehat{\M}^\diamond_{X^-, X^+}, \ab \mathring{K}_{X^-, X^+}, \ab
K_{X^-, X^+}, \ab K^0_{X^-, X^+})$ respectively.


For a triple $((\hat c_l), (x_l), (\alpha_l))$, we define a pre-Kuranishi space
(or a linear combination of pre-Kuranishi spaces)
$\overline{\M}^{X^{[0, \infty]}}((\hat c_l), (x_l), (\alpha_l))$ by
\[
\overline{\M}^{X^{[0, \infty]}}((\hat c_l), (x_l), (\alpha_l))
= \sum_{\star} (-1)^\ast 
\overline{\M}^{(m_-, X^{[0, \infty]}, m_+)}
_{(\Theta^+(e^{\otimes G_{X^+}^+}), \Theta^-(e^{\otimes G_{X^-}^-}),
(\hat c^i_l), (x^i_l), ([\overline{P}] \cap \alpha^i_l))},
\]
where
($G_{X^\pm}^\pm$ are the solutions of (\ref{+G eq}) and (\ref{-G eq})
used for the definition of the generating functions of $X^\pm$)
and the sum $\star$ is taken over all decompositions
\[
\{\hat c_l\} = \coprod_{-m_- \leq i \leq 0} \{\hat c^i_l\}, \quad
\{x_l\} = \coprod_{-m_- \leq i \leq m_+} \{x^i_l\}, \quad
\{\alpha_l\} = \coprod_{0 \leq i \leq m_+} \{\alpha^i_l\}
\]
as sets, and the sign $\ast$ is the weighted sign of the permutation
\[
\begin{pmatrix}
(\hat c^{-m_-}_l) \ (x^{-m_-}_l) \dots (x^{m_+}_l) \ (\alpha^{m_+}_l)\\
(\hat c_l) \quad (x_l) \quad (\alpha_l)
\end{pmatrix}.
\]

Similarly, we define its subspace of irreducible sequences of holomorphic buildings by
\[
\bigl(\overline{\M}^{X^{[0, \infty]}}\bigr)^0((\hat c_l), (x_l), (\alpha_l))
= \sum_{\star} (-1)^\ast
\Bigl(\overline{\M}^{(m_-, X^{[0, \infty]}, m_+)}_{((e^{\otimes G_{X^+}^+}),
(e^{\otimes G_{X^-}^-}),
(\hat c^i_l), (x^i_l), ([\overline{P}] \cap \alpha^i_l))}\Bigr)^0.
\]

Let
\[
\bigl[\overline{\M}^{X^{[0, \infty]}\!, e}_g((\hat c_l), (x_l), (\alpha_l))\bigr]
= (f_{0, g}^e)^\tau((\hat c_l), (x_l), (\alpha_l))
\oplus (f_{1, g}^e)^\tau((\hat c_l), (x_l), (\alpha_l)) d\tau
\]
and
\[
\bigl[(\overline{\M}^{X^{[0, \infty]}\!, e}_g)^0((\hat c_l), (x_l), (\alpha_l))\bigr]
= (h_{0, g}^e)^\tau((\hat c_l), (x_l), (\alpha_l))
\oplus (h_{1, g}^e)^\tau((\hat c_l), (x_l), (\alpha_l)) d\tau
\]
be the counterparts of the virtual fundamental chains, where
$(f_{j, g}^e)((\hat c_l), (x_l), (\alpha_l))$ and $(h_{j, g}^e)((\hat c_l), (x_l), (\alpha_l))$
are smooth functions of $\tau \in [0, \infty] \cong [0, 1]$.
Note that since we assume that the continuous family of grouped multisections
is a usual grouped multisection on a neighborhood of
$\tau \in \{0\} \cup \{\infty\} \subset [0, \infty]$, the zero set of the
perturbed multisections of the $1$-dimensional parts of the above
pre-Kuranishi spaces do not intersect with the corner of codimension $\geq 2$.
In particular, on the zero set, the strong continuous map to
$[0, \infty] \cong [0, 1]$ is submersive.
Hence on a neighborhood of $\tau \in \{0\} \cup \{\infty\} \subset [0, \infty]$,
$(f_{0, g}^e)^\tau$ and $(h_{0, g}^e)^\tau$ are constant as functions of $\tau$,
and $(f_{1, g}^e)^\tau$ and $(h_{1, g}^e)^\tau$ are zero.

Then $f_{0, g}^e$, $h_{0, g}^e$, $\hat f_{1, g}^e = f_{1, g}^e$ and
$\hat h_{1, g}^e = h_{1, g}^e$ satisfy (\ref{f_0 equation}), (\ref{f_1 equation}),
and (\ref{boundary formula for X^I}).
Furthermore, $(f_{0, g}^e)^{\tau = \infty}$ coincides with the $(g, e)$ part of
\[
\sum_{\star} (-1)^\ast
\overline{\M}^{(m_-, X^-, 0, X^+, m_+)}
_{(\Theta^+(e^{\otimes G_{X^+}^+}),
(e^{(\Delta_\ast[\overline{P}])^{0^-, 0^+}}), \Theta^-(e^{\otimes G_{X^-}^-}),
(\hat c^i_l), (x^i_l), ([\overline{P}] \cap \alpha^i_l))}.
\]

If $(Y^0, \lambda^0)$ satisfies Morse condition (i.e. if $P$ is a union of circles),
then the above pre-Kuranishi space is enough for the construction of
a homotopy from the generating function of $X^0$ to the composition
of the generating functions of $X^-$ and $X^+$.
However, in general, we need another parametrized pre-Kuranishi space.

Let $G^\theta$ ($\theta \in [0, 1]$) be an appropriate $C^\infty(I, \R)$-linear
combination of
\[
((\rho_\ast [\overline{P}])^{i, j}, \dots,
(\rho_\ast [\overline{P}])^{i, j},
\epsilon_{\overline{P}}^{i, j}, \dots, \epsilon_{\overline{P}}^{i, j},
(\Delta_\ast [\overline{P}])^{i, j}, \dots, (\Delta_\ast [\overline{P}])^{i, j})_{(i, j)}
\]
defined in the next section.
For each family $((\hat c_l), (x_l), (\alpha_l))$, we define
a $C^\infty(I, \R)$-linear combination of pre-Kuranishi spaces
$\overline{\M}^{X^\infty, \theta \in [0, 1]}((\hat c_l), (x_l), (\alpha_l))$ by
\begin{align*}
&\overline{\M}^{X^\infty, \theta \in [0, 1]}((\hat c_l), (x_l), (\alpha_l)) \\
&= \sum_{\star} (-1)^\ast
\overline{\M}^{(m_-, X^-, m, X^+, m_+)}
_{(\Theta^+(e^{\otimes G_{X^-}^+}), \Theta(e^{\otimes G^\theta}),
\Theta^-(e^{\otimes G_{X^+}^-}), (\hat c^i_l), (x^i_l), ([\overline{P}] \cap \alpha^i_l))}.
\end{align*}
Similarly, we define a $C^\infty(I, \R)$-linear combination of Kuranishi spaces
of irreducible sequences of holomorphic buildings
$(\overline{\M}^{X^\infty, \theta \in [0, 1]})^0((\hat c_l), (x_l), (\alpha_l))$ by
\begin{align*}
&(\overline{\M}^{X^\infty, \theta \in [0, 1]})^0((\hat c_l), (x_l), (\alpha_l))\\
&= \sum_{\star} (-1)^\ast
\Bigl(\overline{\M}^{(m_-, X^-, m, X^+, m_+)}_{((e^{\otimes G_{X^-}^+}), (e^{\otimes G^\theta}),
(e^{\otimes G_{X^+}^-}), (\hat c^i_l), (x^i_l), ([\overline{P}] \cap \alpha^i_l))}\Bigr)^0
\end{align*}
where in this case, the irreducibility is defined as follows.
First we consider the case of $(n_-, n, n_+) \neq (0, 0, 0)$.
A sequence of holomorphic buildings
\[
(\Sigma_i, z_i, u_i, \phi_i)_{i \in
\{-m_-, \dots, -1, 0^-, 1, \dots, m, 0^+, 1, \dots, m_+\}}
\]
in
\[
\overline{\M}^{(m_-, X^-, m, X^+, m_+)}
_{(\Theta^+(f_1^+ \otimes \dots \otimes f_{n_+}^+),
\Theta(f_1 \otimes \dots \otimes f_n),
\Theta^-(f_1^- \otimes \dots \otimes f_{n_-}^-),
(\hat c^i_l), (x^i_l), ([\overline{P}] \cap \alpha^i_l))}
\]
is said to be irreducible if
\begin{itemize}
\item
each connected component of $\Sigma_{0^-}$ and $\Sigma_{0^+}$ concerns at least one
monomial in $\{f_i^\pm, f_i\}$, and
\item
for any decomposition $\{f_i^\pm, f_i\} = A \sqcup B$, there exists some connected
component of $\Sigma_{0^-}$ or $\Sigma_{0^+}$ which concerns both of some $f \in A$
and some $g \in B$.
\end{itemize}
If $(n_-, n, n_+) = (0, 0, 0)$, then a point
$((\Sigma_{0^-}, z_{0^-}, u_{0^-}, \phi_{0^-}), (\Sigma_{0^+}, z_{0^+}, u_{0^+}, \phi_{0^+}))$
is irreducible if one of $(\Sigma_{0^\pm}, z_{0^\pm}, u_{0^\pm}, \phi_{0^\pm})$ is
connected and the other is the empty curve.

Let
\begin{align*}
\bigl[\overline{\M}^{X^\infty, \theta \in [0, 1], e}_g((\hat c_l), (x_l), (\alpha_l))\bigr]
= (\mathring{f}_{0, g}^e)^\theta((\hat c_l), (x_l), (\alpha_l))
\end{align*}
and
\begin{align*}
\bigl[(\overline{\M}^{X^\infty, \theta \in [0, 1], e}_g)^0((\hat c_l), (x_l), (\alpha_l))\bigr]
= (\mathring{h}_{0, g}^e)^\theta((\hat c_l), (x_l), (\alpha_l))
\end{align*}
be the virtual fundamental chains.
Equation (\ref{G^0 eq}) in the next section implies that
$(\mathring{f}_{0, g}^e)^{\theta = 0}$ coincides with $(f_{0, g}^e)^{\tau = \infty}$.

Let $H^\theta = H^\theta_1 + H^\theta_2 + \dots$ be an appropriate
$C^\infty(I, \R)$-linear combination of
\[
((\rho_\ast [\overline{P}])^{i, j}, \dots,
(\rho_\ast [\overline{P}])^{i, j},
\epsilon_{\overline{P}}^{i, j}, \dots, \epsilon_{\overline{P}}^{i, j},
(\Delta_\ast [\overline{P}])^{i, j}, \dots, (\Delta_\ast [\overline{P}])^{i, j})_{(i, j)}
\]
defined in the next section, and define
$(\Ddot f_{0, g}^e)^\theta((\hat c_l), (x_l), (\alpha_l))$ by the virtual fundamental chain of
\[
\sum_{\star} (-1)^\ast
\overline{\M}^{(m_-, X^-, m, X^+, m_+)}
_{(\Theta^+(e^{\otimes G_{X^+}^+}),
\Theta(e^{\otimes G^\theta} \otimes H^\theta),
\sum_{m_-}(-1)^{m_-}\Theta^-(e^{\otimes G_{X^-}^-})_{-m_-}, (\hat c^i_l), (x^i_l),
([\overline{P}] \cap \alpha^i_l))}.
\]
We also define $(\Ddot h_{0, g}^e)^\theta((\hat c_l), (x_l), (\alpha_l))$ by
the virtual fundamental chains of its irreducible part.

Then $f_{0, g}^e = \mathring{f}_{0, g}^e$, $h_{0, g}^e = \mathring{h}_{0, g}^e$,
$\hat f_{1, g}^e = \Ddot f_{0, g}^e$ and $\hat h_{1, g}^e = \Ddot h_{0, g}^e$ also satisfy
(\ref{f_0 equation}), (\ref{f_1 equation}), and (\ref{boundary formula for X^I}).

Define the following families of generating functions.
\begin{align*}
\F^\tau &= \hbar^{-1} \sum \frac{1}{k_q ! k_t ! k_p !} (h_{0, g}^e)^\tau
(\underbrace{\mathbf{q}, \dots, \mathbf{q}}_{k_q},
\underbrace{\mathbf{t}, \dots, \mathbf{t}}_{k_t},
\underbrace{\mathbf{p}, \dots, \mathbf{p}}_{k_p}) \hbar^g T^e\\
\widetilde{\F}^\tau &= \hbar^{-1} \sum \frac{1}{k_q ! k_t ! k_p !} (f_{0, g}^e)^\tau
(\underbrace{\mathbf{q}, \dots, \mathbf{q}}_{k_q},
\underbrace{\mathbf{t}, \dots, \mathbf{t}}_{k_t},
\underbrace{\mathbf{p}, \dots, \mathbf{p}}_{k_p}) \hbar^g T^e\\
\K^\tau_g &= \hbar^{-1} \sum \frac{1}{k_q ! k_t ! k_p !} (h_{1, g}^e)^\tau
(\underbrace{\mathbf{q}, \dots, \mathbf{q}}_{k_q},
\underbrace{\mathbf{t}, \dots, \mathbf{t}}_{k_t},
\underbrace{\mathbf{p}, \dots, \mathbf{p}}_{k_p}) \hbar^g T^e\\
\widetilde{\K}^\tau_g &= \hbar^{-1} \sum \frac{1}{k_q ! k_t ! k_p !} (f_{1, g}^e)^\tau
(\underbrace{\mathbf{q}, \dots, \mathbf{q}}_{k_q},
\underbrace{\mathbf{t}, \dots, \mathbf{t}}_{k_t},
\underbrace{\mathbf{p}, \dots, \mathbf{p}}_{k_p}) \hbar^g T^e\\
\F^\theta &= \hbar^{-1} \sum \frac{1}{k_q ! k_t ! k_p !} (\mathring{h}_{0, g}^e)^\theta
(\underbrace{\mathbf{q}, \dots, \mathbf{q}}_{k_q},
\underbrace{\mathbf{t}, \dots, \mathbf{t}}_{k_t},
\underbrace{\mathbf{p}, \dots, \mathbf{p}}_{k_p}) \hbar^g T^e\\
\widetilde{\F}^\theta &= \hbar^{-1} \sum \frac{1}{k_q ! k_t ! k_p !}
(\mathring{f}_{0, g}^e)^\theta
(\underbrace{\mathbf{q}, \dots, \mathbf{q}}_{k_q},
\underbrace{\mathbf{t}, \dots, \mathbf{t}}_{k_t},
\underbrace{\mathbf{p}, \dots, \mathbf{p}}_{k_p}) \hbar^g T^e\\
\K^\theta_g &= \hbar^{-1} \sum \frac{1}{k_q ! k_t ! k_p !} (\Ddot h_{0, g}^e)^\theta
(\underbrace{\mathbf{q}, \dots, \mathbf{q}}_{k_q},
\underbrace{\mathbf{t}, \dots, \mathbf{t}}_{k_t},
\underbrace{\mathbf{p}, \dots, \mathbf{p}}_{k_p}) \hbar^g T^e\\
\widetilde{\K}^\theta_g &= \hbar^{-1} \sum \frac{1}{k_q ! k_t ! k_p !}
(\Ddot f_{0, g}^e)^\theta
(\underbrace{\mathbf{q}, \dots, \mathbf{q}}_{k_q},
\underbrace{\mathbf{t}, \dots, \mathbf{t}}_{k_t},
\underbrace{\mathbf{p}, \dots, \mathbf{p}}_{k_p}) \hbar^g T^e
\end{align*}

Then it is easy to see that the concatenation of the homotopies $\F^\tau$ and
$\F^\theta$ defined by the above generating functions
gives a homotopy from the generating function of $X^0$ to the composition of
the generating functions of $X^-$ and $X^+$.

\subsection{Construction of the correction terms}\label{correction terms for composition}
\subsubsection{Constuction of $G^\theta$}
For $m \geq 1$, let $C_m = \bigoplus_{n = 0}^{\frac{m(m+1)}{2}} C_m^n$ be the $\Z$-graded
super-commutative algebra with coefficient $\R$ generated by variables
$\rho_{(e_i, e_j)}$, $\Delta_{(e_i, e_j)}$ and $\epsilon_{(e_i, e_j)}$ ($0 \leq i < j \leq m$).
The $\Z$-grading is defined by $\dim \rho_{(e_i, e_j)} = \dim \Delta_{(e_i, e_j)} = 0$ and
$\dim \epsilon_{(e_i, e_j)} = 1$.

For each $m \geq 1$, the differential $\partial' : C_m^n \to C_m^{n-1}$ is defined by
$\partial' \epsilon_{(a, b)} = (-1)^{m-1} (\rho_{(a, b)} - \Delta_{(a, b)})$ and
$\partial' \rho_{(a, b)} = \partial' \Delta_{(a, b)} = 0$.
Homomorphisms $\tau_i : C_m \to C_{m + 1}$ ($0 \leq i \leq m$) are defined by
$\tau_i (x_{(a, b)}) = x_{(\hat \tau_i(a), \hat \tau_i(b))}$,
where each $\hat \tau_i$ is defined by
\[
\hat \tau_i(e_j) = \begin{cases}
e_j & j < i\\
e_i + e_{i + 1} & j = i\\
e_{j + 1} & j > i
\end{cases}.
\]

Define homomorphism $\Theta : \bigotimes_{i = 1}^n C_{m_i} \to C_{1 + \sum_{i=1}^n
(m_i - 1)}$ by
\[
\Theta(f_1 \otimes f_2 \otimes \dots f_n) = f_1^{+\sum_{i = 2}^n (m_i - 1)} f_2^{+\sum_{i = 3}^n
(m_i - 1)} \dots f_n,
\]
where each $f_a^{+\sum_{i = a + 1}^n (m_i - 1)}$ is defined by
\[
e_j^{+\sum_{i = a + 1}^n (m_i - 1)} = \begin{cases}
e_0 & j = 0\\
e_{j + \sum_{i = a + 1}^n (m_i - 1)} & j \neq 0, m_a\\
e_{1 + \sum_{i=1}^n (m_i - 1)} & j = m_a
\end{cases}.
\]

We also define $\boxminus : B_m^+ \otimes B_{m'}^- \to C_{m + m' + 1}$ by
\[
\boxminus(f \otimes g) = (-1)^{m m'} f \cdot \exp(\rho_{(\sum_{0 \leq i \leq m} e_i,
\sum_{m + 1 \leq j \leq m + m' + 1} e_j)}) \cdot g^{+(m + m' + 1)}.
\]

We define a linear subspace $\Ddot C_m \subset C_m$ as follows.
For each $1 \leq i \leq m - 2$ and each monomial
\[
f = x^{(1)}_{(a_1, b_1)} x^{(2)}_{(a_2, b_2)} \dots x^{(n)}_{(a_n, b_n)}
\]
such that $(a_j, b_j) \neq (i, i + 1)$, we define a monomial
\[
f^{(e_i, e_{i + 1})} = x^{(1)}_{(a'_1, b'_1)} x^{(2)}_{(a'_2, b'_2)} \dots x^{(n)}_{(a'_n, b'_n)}
\]
by permuting $i$ and $i + 1$ of $\{a_j, b_j\}$.
Then $\Ddot C_m \subset C_m$ is the subspace spanned by $f + f^{(e_i, e_{i + 1})}$ for all
such pair $i$ and $f$.

Define $\mathcal{C}_m = C_m / \Ddot C_m$.
Then the following maps are well defined.
\begin{align*}
\partial' &: \mathcal{C}_m \to \mathcal{C}_m \\
\sum_{0 < i < \max} (-1)^i e^{\Delta_{(e_i, e_{i + 1})}}
\tau_i &: \mathcal{C}_m \to \mathcal{C}_{m + 1} \quad\quad ({\max} = m) \\
e^{\Delta_{(e_0, e_1)}} \tau_0 &: \mathcal{C}_m \to \mathcal{C}_{m + 1} \\
(-1)^{\max} e^{\Delta_{(e_{\max}, e_{\max + 1})}} \tau_{\max} &: \mathcal{C}_m \to \mathcal{C}_{m + 1}\\
\Theta &: \otimes_{i = 1}^n \mathcal{C}_{m_i} \to \mathcal{C}_{1 + \sum_i (m_i - 1)}\\
\boxminus &: \B_m^+ \otimes \B_{m'}^- \to \mathcal{C}_{m + m' + 1}
\end{align*}

Further we define $\mathring{\mathcal{C}}_m \subset \mathcal{C}_m$ as follows.
We define a new degree $\deg'$ by
\[
\deg' x_{(e_i, e_j)} = \begin{cases}
0 & \text{if } i = 0 \text{ or } j = m\\
1 & \text{otherwise}
\end{cases}
\]
For $m \geq 2$, let $\mathring{C}_m \subset C_m$ be the subspace spanned by monomials
with $\deg' \geq m - 2$ which do not contain variables $\rho_{(e_0, e_m)}$,
$\Delta_{(e_0, e_m)}$ or $\epsilon_{(e_0, e_m)}$.
Define $\mathring{\mathcal{C}}_m = \mathring{C}_m
/ (\Ddot{C}_m \cap \mathring{C}_m) \subset \mathcal{C}_m$.

In this section, we prove that there exists a smooth family
$G^\theta = G^\theta_1 + G^\theta_2 + \dots
\in (\bigoplus_{m = 1}^\infty \mathring{C}_m^{m-1})^\wedge$ ($\theta \in [0, 1]$)
which satisfies the following equations.
\begin{gather}
\partial' (\Theta(e^{\otimes G^\theta})) + \sum_{i \geq 0} e^{\Delta_{(e_i, e_{i+1})}} \tau_i
\Theta(e^{\otimes G^\theta}) = 0\\
G^\theta_1 = (1 - \theta) \Delta_{(e_0, e_1)} + \theta \rho_{(e_0, e_1)}\\
G^0 = G^0_1= \Delta_{(e_0, e_1)} \label{G^0 eq}\\
\boxminus(e^{\otimes G^+_{X^-}} \otimes e^{\otimes G^-_{X^+}}) = \Theta(e^{\otimes G^1})
\label{G^1 eq}
\end{gather}
In the previous section, we replace $\rho_{(e_i, e_j)}$, $\Delta_{(e_i, e_j)}$ and
$\epsilon_{(e_i, e_j)}$ in $G^\theta$ with $(\rho_\ast [\overline{P_{Y^0}}])^{0_i, 0_j}$,
$(\Delta_\ast [\overline{P_{Y^0}}])^{0_i, 0_j}$ and $(\epsilon_{\overline{P_{Y^0}}})^{0_i, 0_j}$
respectively, where $0_0$ and $0_{\max}$ should be read as $0^-$ and $0^+$ respectively.

First we note that the last two equations define $G^0$ and $G^1$.
We inductively construct $G^\theta_{\leq m} = G^\theta_1 + \dots + G^\theta_m \in
\bigoplus_{l = 1}^m \mathring{\mathcal{C}}_l^{l-1}
$ such that
\begin{equation}
\partial' (\Theta(e^{\otimes G^\theta_{\leq m}})) + \sum_{i \geq 0} e^{\Delta_{(e_i, e_{i+1})}} \tau_i
\Theta(e^{\otimes G^\theta_{\leq m-1}}) \equiv 0 \label{G theta eq}
\end{equation}
in $\bigoplus_{l=2}^\infty \mathcal{C}_l^{l-2} / \bigoplus_{l=m+1}^\infty \mathcal{C}_l^{l-2}$.

First we define $G^\theta_2 \in \mathring{\mathcal{C}}_2^1$ by
\begin{align*}
G^\theta_2 &= e^{(1-\theta)(\Delta_{(e_0, e_1)} + \Delta_{(e_1, e_2)})}\\
&\quad \cdot \biggl(-\sum_{k \geq 1} \frac{\theta^k}{k !} e^{\theta\rho_{(e_1, e_2)}}
(\underbrace{\epsilon_{(e_0, e_1)} \Delta_{(e_0, e_1)} \dots \Delta_{(e_0, e_1)}}_k\\
&\hph{\quad \cdot \biggl(-\sum_{k \geq 1} \frac{\theta^k}{k !} e^{\theta\rho_{(e_1, e_2)}}(}
+ \underbrace{\rho_{(e_0, e_1)} \epsilon_{(e_0, e_1)} \Delta_{(e_0, e_1)} \dots
\Delta_{(e_0, e_1)}}_k\\
&\hph{\quad \cdot \biggl(-\sum_{k \geq 1} \frac{\theta^k}{k !} e^{\theta\rho_{(e_1, e_2)}}(}
+ \dots + \underbrace{\rho_{(e_0, e_1)} \dots \rho_{(e_0, e_1)} \epsilon_{(e_0, e_1)}}_k)\\
&\quad \quad + \sum_{k \geq 1} \frac{\theta^k}{k !} e^{\theta\rho_{(e_0, e_1)}}
(\underbrace{\epsilon_{(e_1, e_2)} \Delta_{(e_1, e_2)} \dots \Delta_{(e_1, e_2)}}_k\\
&\hph{\quad \quad + \sum_{k \geq 1} \frac{\theta^k}{k !} e^{\theta\rho_{(e_0, e_1)}}(}
+ \underbrace{\rho_{(e_1, e_2)} \epsilon_{(e_1, e_2)} \Delta_{(e_1, e_2)} \dots
\Delta_{(e_1, e_2)}}_k\\
&\hph{\quad \quad + \sum_{k \geq 1} \frac{\theta^k}{k !} e^{\theta\rho_{(e_0, e_1)}}(}
+ \dots + \underbrace{\rho_{(e_1, e_2)} \dots \rho_{(e_1, e_2)} \epsilon_{(e_1, e_2)}}_k)\biggr).
\end{align*}
Then it is easy to see that this satisfies equation (\ref{G theta eq}) for $m = 2$.

Next assuming we have constructed $G^\theta_{\leq m-1}$, we prove there exists
a required family $G^\theta_m$. It is enough to show that
\begin{equation}
\Theta\Bigl(\Bigl(\partial' \Theta(e^{\otimes G^\theta_{\leq m-1}}) + \sum (-1)^i
e^{\Delta_{(e_i, e_{i+1})}} \tau_i \Theta(e^{\otimes G^\theta_{\leq m-1}})\Bigr)
\otimes e^{-G^\theta_1}\Bigr) \equiv 0 \label{G theta mathring}
\end{equation}
in $\bigoplus_{l=2}^\infty \mathcal{C}_l^{l-2} / (\bigoplus_{l = m+1}^\infty \mathcal{C}_l^{l-2}
\oplus \bigoplus_{l = 2}^\infty \mathring{\mathcal{C}}_l^{l-2})$ and
\begin{equation}
\partial'\Bigl(\sum_{0 \leq i \leq \max} (-1)^i e^{(e_i, e_{i+1})} \tau_i
\Theta (e^{\otimes G^\theta_{\leq m-1}})\Bigr) \equiv 0 \label{G theta closed}
\end{equation}
in $\bigoplus_{l = 3}^\infty \mathcal{C}_l^{l-3} / \bigoplus_{l = m+1}^\infty \mathcal{C}_l^{l-3}$.

The latter is proved by an argument similar to that for equation (\ref{B mathring}).
We can prove the former similarly to equation (\ref{B closed}) using the following
equations.
\begin{align}
&\partial' \Theta\Bigl(\frac{1}{k !}(G^\theta_{\geq m-1} - G^\theta_1)^{\otimes k}\Bigr)
\notag\\
&= \Theta\Bigl(\frac{1}{(k-1) !} (G^\theta_{\geq m-1} - G^\theta_1)^{\otimes (k-1)} \otimes
\partial' (G^\theta_{\geq m-1} - G^\theta_1)\Bigr)
\end{align}
\begin{align}
&\sum_{0 < i < \max} (-1)^i e^{\Delta_{(e_i, e_{i+1})}} \tau_i
\Theta\Bigl(\frac{1}{k !} (G^\theta_{\geq m-1} - G^\theta_1)^{\otimes k}\Bigr)\notag\\
&= \Theta\Bigl(\frac{1}{(k-1) !} (G^\theta_{\geq m-1} - G^\theta_1)^{\otimes (k-1)}
\notag\\
&\hph{= \Theta\Bigl(}
\otimes \sum_{0 < i < \max} (-1)^i e^{\Delta_{(e_i, e_{i+1})}}
\tau_i (G^\theta_{\geq m-1} - G^\theta_1)
\Bigr)
\end{align}
\begin{align}
&e^{\Delta_{(e_0, e_1)}} \tau_0
\Theta\Bigl(\frac{1}{k !}(G^\theta_{\geq m-1})^{\otimes k}\Bigr)
\notag\\
&= \sum_{l_1 + l_2 + l_3 = k} \Theta\Bigl(\frac{1}{l_1 ! l_2 ! l_3 !}
(G^\theta_{\geq m-1} - G^\theta_1)^{\otimes l_1} \otimes
\mathring{\tau}_0((G^\theta_{\geq m-1})^{\otimes l_2}) \notag \\
&\hph{= \sum_{l_1 + l_2 + l_3 = k} \Theta\Bigl(}
\otimes
(G^\theta_1)^{\otimes l_3} \Bigr)
\end{align}
\begin{align}
&(-1)^{\max} e^{\Delta_{(e_{\max}, e_{\max +1})}} \tau_{\max}
\Theta\Bigl(\frac{1}{k !}(G^\theta_{\geq m-1})^{\otimes k}\Bigr)\notag\\
&= \sum_{l_1 + l_2 + l_3 = k} \Theta\Bigl(\frac{1}{l_1 ! l_2 ! l_3 !}
(G^\theta_{\geq m-1} - G^\theta_1)^{\otimes l_1} \otimes
(-1)^{\max}\mathring{\tau}_{\max}((G^\theta_{\geq m-1})^{\otimes l_2})
\notag\\
&\hph{= \sum_{l_1 + l_2 + l_3 = k} \Theta\Bigl(}
\otimes (G^\theta_1)^{\otimes l_3}\Bigr)
\end{align}
In the above equations, $\mathring{\tau}_0$ and $\mathring{\tau}_{\max}$ are defined
in a similar way to $\mathring{\tau}^+_0$ in Section \ref{correction terms for X}.

Therefore we can inductively construct a required family
$G^\theta_{\leq m} = G^\theta_1 + \dots + G^\theta_m \in
\bigoplus_{l = 1}^m \mathring{\mathcal{C}}_l^{l-1} \cong \bigoplus_{l = 1}^\infty
\mathring{\mathcal{C}}_l^{l-1} / \bigoplus_{l = m+1}^\infty \mathring{\mathcal{C}}_l^{l-1}$.

\subsubsection{Construction of $H^\theta$}
Next we construct a smooth family
$H^\theta = H^\theta_1 + H^\theta_2 + \dots
\in (\bigoplus_{m=1}^\infty \mathring{\mathcal{C}}_m^m)^\wedge$ which satisfies
the following equation.
\[
\partial' \Theta(e^{\otimes G} \otimes H)
+ \sum_{i \geq 0} (-1)^i e^{\Delta_{(e_i, e_j)}} \tau_i \Theta(e^{\otimes G} \otimes H)
+ \Theta\Bigl(e^{\otimes G} \otimes \frac{d}{d\theta} G\Bigr) = 0.
\]

We inductively construct $H^\theta_{\leq m} = H^\theta_1 + \dots + H^\theta_m \in
\bigoplus_{l = 1}^m \mathring{\mathcal{C}}_l^l$ such that
\begin{multline}
\partial' \Theta(e^{\otimes G} \otimes H_{\leq m})
+ \sum_{i \geq 0} (-1)^i e^{\Delta_{(e_i, e_j)}} \tau_i
\Theta(e^{\otimes G} \otimes H_{\leq m}) \\
+ \Theta\Bigl(e^{\otimes G} \otimes \frac{d}{d\theta} G\Bigr) \equiv 0 \label{H eq}
\end{multline}
in $(\bigoplus_{l=1}^\infty \mathcal{C}_l^{l-1})^\wedge /
(\bigoplus_{l=m+1}^\infty \mathcal{C}_l^{l-1})^\wedge$

Since $\frac{d}{d\theta} G_1^\theta = \rho_{(e_0, e_1)} - \Delta_{(e_0, e_1)}$,
$H_1^\theta = - \epsilon_{(e_0, e_1)}$ satisfies equation (\ref{H eq}) for $m = 1$.

Assuming we have already constructed $H_{\leq m-1}$,
we prove that there exists a required family $H_m^\theta$.
It is enough to show that
\begin{multline*}
\Bigl(\partial' \Theta(e^{\otimes G} \otimes H_{\leq m-1})
+ \sum_{i \geq 0} (-1)^i e^{\Delta_{(e_i, e_j)}} \tau_i
\Theta(e^{\otimes G} \otimes H_{\leq m-1})\\
+ \Theta\Bigl(e^{\otimes G} \otimes \frac{d}{d\theta} G\Bigr)\Bigr) \otimes 
e^{-\otimes G_1} \equiv 0
\end{multline*}
in $(\bigoplus_{l=1}^\infty \mathcal{C}_l^{l-1})^\wedge /
((\bigoplus_{l=m+1}^\infty \mathcal{C}_l^{l-1})^\wedge
\oplus \bigoplus_{l = 1}^\infty \mathring{\mathcal{C}}_l^{l-1})$
and
\[
\partial'\Bigl(\sum_{i \geq 0} (-1)^i e^{\Delta_{(e_i, e_j)}}
\tau_i \Theta(e^{\otimes G} \otimes H_{\leq m-1})
+ \Theta\Bigl(e^{\otimes G} \otimes \frac{d}{d\theta} G\Bigr)\Bigr) \equiv 0
\]
in $(\bigoplus_{l=2}^\infty \mathcal{C}_l^{l-2})^\wedge /
(\bigoplus_{l=m+1}^\infty \mathcal{C}_l^{l-2})^\wedge$.
The former can be proved by a similar argument to
those for (\ref{H tau mathring}) or (\ref{G theta mathring}),
and the latter can be proved similarly to (\ref{A closed}), (\ref{B closed}),
(\ref{H tau closed}) or (\ref{G theta closed}).
Therefore, we can inductively construct a required family $H^\theta \in
(\bigoplus_{m=1}^\infty \mathring{\mathcal{C}}_m^m)^\wedge$.


%% file: SFT-13_Independence.tex
%
%

\section{Independence}\label{independence}
Let $(Y, \xi)$ be a contact manifold and let $\overline{K}_Y^0 \subset H_\ast(Y, \Z)$
be a finite subset (or a finite sequence).
We have seen that if we fix 
a contact form $\lambda$, a triangulation $K_Y$ of $\overline{P}_Y$,
a Euclidean cell complex $K^2_Y$, a representative $K^0_Y$ of $\overline{K}_Y^0$,
and a complex structure $J$ of $\ker \lambda$,
construct a family of pre-Kuranishi spaces and choose a compatible family of
perturbed multisections, then
we obtain chain complexes
\begin{gather*}
(\W^{\leq \kappa}_{(Y, \lambda, K_Y, \overline{K}_Y^0)} / I^{\leq \kappa}_{C_0, C_1, C_2},
D_{(Y, \lambda, K_Y, K_Y^0, K_Y^2, J, \B)}), \\
(\mathcal{P}^{\leq \kappa}_{(Y, \lambda, K_Y, \overline{K}_Y^0)}
/ I^{\leq \kappa}_{C_0, C_2},
d_{(Y, \lambda, K_Y, K_Y^0, K_Y^2, J, \B)}),\\
(\A^{\leq \kappa}_{(Y, \lambda, K_Y, \overline{K}_Y^0)} / I^{\leq \kappa}_{C_0},
\partial_{(Y, \lambda, K_Y, K_Y^0, K_Y^2, J, \B)}),
\end{gather*}
where $\B$ denotes the other choices for the construction of
the pre-Kuranishi structure and the perturbed multisections.
The aim of this section is to construct SFT cohomologies of a contact manifold
by the limits of the cohomologies of the above chain complexes and
to prove that they are invariants of $(Y, \xi, \overline{K}_Y^0)$.
We also construct SFT cohomologies of a symplectic cobordism as limits.

First we note that for any constant $a > 0$,
chain complexes for $(Y, a \lambda)$ can be constructed by using
the same $(K_Y, K^2_Y, K^0_Y, J, \B)$ as those of $(Y, \lambda)$.
Then the chain complex
$(\W^{\leq a \kappa}_{(Y, a \lambda, K_Y, K_Y^0)} / I^{\leq a \kappa}_{C_0, C_1, a C_2},
D_{(Y, \lambda, K_Y, K_Y^0, K_Y^2, J, \B)})$ is naturally isomorphic to
$(\W^{\leq \kappa}_{(Y, \lambda, K_Y, K_Y^0)} / I^{\leq \kappa}_{C_0, C_1, C_2},
D_{(Y, \lambda, K_Y, K_Y^0, K_Y^2, J, \B)})$.
The cases of the other two chain complexes are similar.

Let $\mathcal{C}_{Y^\pm}
= (Y^\pm (= Y), \lambda^\pm, K_{Y^\pm}, K_{Y^\pm}^0, K_{Y^\pm}^2, J^\pm, \B_{Y^\pm})$
be two choices to define the above chain complexes.
A concordance $\mathcal{C}_X
= (X, \omega, Y^\pm, \lambda^\pm, K_{Y^\pm},\ab K_X^0,\ab K^0_{Y^\pm},\ab
\mu^\pm, K^2_{Y^\pm}, J, \B_X)$
from $\mathcal{C}_{Y^-}$ to $\mathcal{C}_{Y^+}$
consists of
\begin{itemize}
\item
a cobordism $(X, \omega)$ from $(Y^-, \lambda^-)$ to $(Y^+, \lambda^+)$ of the form
$X = (-\infty, 0] \times Y^- \cup [0, T_0] \times Y \cup [0, \infty) \times Y^+$
for $T_0 \geq 0$ and $\omega|_{[0, T_0] \times Y} = d(f \lambda^-)$ for some
smooth function $f : [0, T_0] \times Y \to \R_{>0}$ such that
$f \lambda^-|_{\{0\} \times Y} = \lambda^-$ and
$f \lambda^-|_{\{T_0\} \times Y} = \lambda^+$,
\item
a sequence $K_X^0$ of smooth cycles in $X$ with closed support and
bijections $\mu^\pm : K_X^0 \to K_{Y^\pm}^0$ such that
for some $T \geq 0$,
$x|_{(-\infty, -T] \times Y^-} = (-\infty, -T] \times \mu^-(x)$
and $x|_{[T, \infty) \times Y^+} = [T, \infty) \times \mu^+(x)$,
\item
an $\omega$-compatible almost complex structure $J$ of $X$ whose restrictions to
$(-\infty, -T] \times Y^-$ and $[T, \infty) \times Y^+$ coincide with those induced
by $J^-$ and $J^+$ respectively for some $T \geq 0$ and
\item
a pre-Kuranishi structure of $\widehat{\M}(X, \omega, J)$ and a family of multisections
of its fiber products compatible with $\B_{Y^\pm}$, which is denoted by $\B_X$.
\end{itemize}
We note that for the algebra of SFT of $X$,
$\delta = \min(L_{Y^-, \min}, L_{Y^+, \min})$ is admissible for any $C_2 \geq 0$.
(We can define the generating function $\F$ for $X$ as an element of
$(\hbar^{-1} \D_X^{\leq 0})^{\star, \delta}
/ J^{\star, \delta}_{\overline{C}_0, \overline{C}_1, \overline{C}_2}$.)

We say a concordance
$\mathcal{C}_X$ is trivial if
$(Y^-, \lambda^-, K_{Y^-}, K_{Y^-}^0, K_{Y^-}^2, J^-, \B_{Y^-})
= (Y^+, a \lambda^+, K_{Y^+}, K_{Y^+}^0, K_{Y^+}^2, J^+, \B_{Y^+})$ for some
$a > 0$.
A short concordance $\mathcal{C}_X$ is a concordance such that
$(Y^-, \lambda^-) = (Y^+, \lambda^+)$ and $T_0 = 0$, that is,
$X = (-\infty, 0] \times Y^- \cup [0, \infty) \times Y^+$.

First we prove the following.
\begin{lem}\label{trivial generating function}
For a trivial concordance $\mathcal{C}_X$,
the generating function $\F$
is homotopic to the trivial generating function
\[
\F^{\text{\it tri}} = \hbar^{-1} \sum_c q^-_{\hat c^\ast} p^+_{\hat c},
\]
where the sum is taken over all simplices in $K_{Y^+}$ not contained in
$\overline{P}_{Y^+}^{\text{bad}}$.
\end{lem}
First we consider the case of a trivial short concordance.
We denote the same
$(Y^\pm, \lambda^\pm, K_{Y^\pm}, K_{Y^\pm}^0, K_{Y^\pm}^2, J^\pm, \B_{Y^\pm})$
by $(Y, \lambda, K_Y, K_Y^0, K_Y^2, J, \B_Y)$,
and regard the symplectization $(X, \omega) = (Y \times \R, d(e^\sigma \lambda))$
as a trivial short concordance.

For each pair $(\hat c = c \theta^D_c, \hat \eta = \theta^{^t\! D}_\eta \eta)$,
let $\overline{\M}^{E_{\hat \omega}=0}_{g = 0, \# \pm\infty = 1}(\hat c, \hat \eta)
\subset (\overline{\M}^X)^0_{(\hat c, \emptyset, \hat \eta)}$
be the component which consists of connected holomorphic buildings of genera $g = 0$
with one limit circle for each end and
without marked points whose $E_{\hat \omega}$-energies are zero.
(Namely, these are trivial cylinders in the $0$-th floor.)
\begin{lem}
The chain map
$\varphi : C_\ast(\overline{P}_Y, \overline{P}_Y^{\text{\normalfont bad}};
\S^D \otimes \Q)
\to C_\ast(\overline{P}_Y, \overline{P}_Y^{\text{\normalfont bad}}; \S^D \otimes \Q)$
defined by
\[
\varphi(\hat c) = \sum_{c'} \bigl[\overline{\M}^{E_{\hat \omega}=0}_{g = 0, \# \pm\infty = 1}(\hat c,
[\overline{P}] \cap (\hat c')^\ast)\bigr]^0 \hat c'
\]
is chain homotopic to the identity, where
the sum is taken over all simplices $c'$ in $K_Y^0$ not contained in
$\overline{P}_Y^{\text{\normalfont bad}}$.
\end{lem}
\begin{proof}
It is easy to check that $\varphi$ is indeed a chain map.
Therefore it is enough to show that
\begin{equation}
\bigl[\overline{\M}^{E_{\hat \omega}=0}_{g = 0, \# \pm\infty = 1}(x, [\overline{P}] \cap \alpha)\bigr]^0
= \langle x, \alpha \rangle \label{equation of intersection theory}
\end{equation}
for any cycle
$x \in C_\ast(\overline{P}_Y, \overline{P}_Y^{\text{bad}}; \S^D \otimes \Q)$
and any cocycle $\alpha \in C^\ast(\overline{P}_Y, \overline{P}_Y^{\text{bad}};
\S^D \otimes \Q)$.
Recall that for the fundamental chain
\[
[\overline{P}] = \sum_\zeta \frac{1}{m_\zeta} \zeta \theta^{\overline{P}}_\zeta
\in C_{\dim P - 1}(\overline{P}, \overline{P}^{\text{no}}; \S^{\overline{P}} \otimes \Q),
\]
$\rho_\ast [\overline{P}] \in C_{\dim P -1} (\overline{P} \times \overline{P},
\overline{P}^{^t\text{bad}} \times \overline{P} \cup
\overline{P} \times \overline{P}^{\text{bad}};
p_1^\ast \S^{\lsuperscript{D}{t}} \otimes p_2^\ast \S^D \otimes \Q)$ is defined by
\[
\rho_\ast [\overline{P}] = \sum_\zeta \frac{1}{m_\zeta} \theta^{\lsuperscript{D}{t}}_\zeta
(\rho_\ast \zeta) \theta^D_\zeta,
\]
where
\[
\rho_\ast \zeta = \sum_{0 \leq p \leq n} \partial_{p+1} \dots \partial_n \zeta
\times \partial_0 \dots \partial_{p+1} \zeta.
\]
Note that by definition,
the left hand side of Equation (\ref{equation of intersection theory}) coincides with
\begin{equation}
\Bigl\langle \sum_{\zeta, p}
\Bigl[\overline{\M}^{E_{\hat \omega}=0}_{g = 0, \# \pm\infty = 1}(x,
\theta^{\lsuperscript{D}{t}}_\zeta \partial_{p+1} \dots \partial_n \zeta)\Bigr]^0
\partial_0 \dots \partial_{p+1} \zeta \theta^D_\zeta, \alpha \Bigr\rangle.
\label{expansion of cap product}
\end{equation}
We rewrite $\sum_{\zeta, p} [\overline{\M}^{E_{\hat \omega}=0}_{g = 0, \# \pm\infty = 1}(x,
\theta^{\lsuperscript{D}{t}}_\zeta \partial_{p+1} \dots \partial_n \zeta)]^0
\partial_0 \dots \partial_{p+1} \zeta \theta^D_\zeta$ as
the virtual fundamental chain of the fiber product of
$\overline{\M}^{E_{\hat \omega}=0}_{g = 0, \# \pm\infty = 1}(x, \cdot)$ with
$\rho_\ast [\overline{P}]$,
and prove that this is homologous to the virtual fundamental chain of
the fiber product with $\Delta_\ast [\overline{P}]$.
More precisely, we construct these fiber products as follows.

For simplices $c \subset \overline{P}$ and
$\eta \subset \overline{P} \times \overline{P}$,
we define $\widehat{\M}^{E_{\hat \omega}=0}_{g = 0, \# \pm\infty = 1}(c, \eta)$
by the inverse image of $c \times \Delta_{\overline{P}} \subset
\overline{P} \times (\overline{P} \times \overline{P})$ by the map
\[
(\ev_{-\infty} \times \ev_{+\infty}) \times \pi_1
: \widehat{\M}^{E_{\hat \omega}=0}_{g = 0, \# \pm\infty = 1} \times \eta
\to (\overline{P} \times \overline{P}) \times \overline{P},
\]
where $\pi_1 : \overline{P} \times \overline{P} \to \overline{P}$ is the first projection.
Similarly, for simplices with local coefficients $\hat c = c \theta^D_c$ and
$\hat \eta = \theta^{^t\! D}_\eta \eta \theta^D_\eta$,
we define $\overline{\M}^{E_{\hat \omega}=0}_{g = 0, \# \pm\infty = 1}(\hat c, \hat \eta)$
by choosing lifts $\tilde c \subset P$ and
$\tilde \eta \subset \overline{P} \times \overline{P}$ of $c$ and $\eta$ respectively.
Its orientation is defined by using $\theta^D_c$ and $\theta^{^t\! D}_\eta$.

For each
\[
\widehat{\M}^{E_{\hat \omega}=0}_{g = 0, \# \pm\infty = 1}
(c, \partial_{p+1} \dots \partial_n \zeta \times \partial_0 \dots \partial_{p+1} \zeta)
\subset \widehat{\M}^{E_{\hat \omega}=0}_{g = 0, \# \pm\infty = 1}
(x, \rho_\ast [\overline{P}]),
\]
we use the perturbed multisection
defined by the pull back by the submersion to
$\widehat{\M}^{E_{\hat \omega}=0}_{g = 0, \# \pm\infty = 1}
(c, \partial_{p+1} \dots \partial_n \zeta)$.
Then (\ref{expansion of cap product}) coincides with
\begin{equation}
\bigl\langle (\pi_2)_\ast
\bigl(\overline{\M}^{E_{\hat \omega}=0}_{g = 0}(x, \rho_\ast [\overline{P}])\bigr),
\alpha\bigr\rangle,
\label{another expression of cap product}
\end{equation}
where $\pi_2$ is the strong smooth map defined by
the second projection $\overline{P} \times \overline{P} \to \overline{P}$.
In (\ref{another expression of cap product}), we can replace $\alpha$ with
a closed form $\tilde \alpha$ (with local coefficient) which represents
$\alpha \in H^\ast(\overline{P}_Y, \overline{P}_Y^{\text{bad}}; \S^D \otimes \Q)$
and rewrite (\ref{another expression of cap product}) as
\begin{equation}
\int_{\overline{\M}^{E_{\hat \omega}=0}_{g = 0}(x, \rho_\ast [\overline{P}])}
\pi_2^\ast \tilde \alpha.
\label{int on rho ast}
\end{equation}
Since $\overline{\M}^{E_{\hat \omega}=0}_{g = 0}(x, \rho_\ast [\overline{P}])$ and
$\overline{\M}^{E_{\hat \omega}=0}_{g = 0}(x, \Delta_\ast [\overline{P}])$ are
cobordant by $\overline{\M}^{E_{\hat \omega}=0}_{g = 0}(x, \epsilon_\ast [\overline{P}])$,
(\ref{int on rho ast}) coincides with
\begin{equation}
\int_{\overline{\M}^{E_{\hat \omega}=0}_{g = 0}(x, \Delta_\ast [\overline{P}])}
\pi_2^\ast \tilde \alpha.
\label{int on Delta ast}
\end{equation}
For a simplex $c \subset \overline{P}$,
let $\widehat{\M}^{E_{\hat \omega}=0}_{g = 0}(c, \cdot)$ be the space
defined by the fiber product with $c$ on the $-\infty$-side. (For $+\infty$-side,
we do not take fiber product.)
Note that there exists a submersion from
$\widehat{\M}^{E_{\hat \omega}=0}_{g = 0}(c, \cdot)$ to
$\widehat{\M}^{E_{\hat \omega}=0}_{g = 0}(c, \Delta_\ast [\overline{P}])$.
In fact, the only difference is that for the construction of a perturbed multisection
of the latter, we need to make the zero set transverse to all simplices in $\overline{P}$.
Define $\overset{\,\!_-\,\!_\wedge}{\M}\,\!_{g = 0}^{E_{\hat \omega}=0}(x, \cdot)$
by the space of holomorphic buildings with $S^1$-coordinates
only on $-\infty$-limit circle.
Then (\ref{int on Delta ast}) coincides with
\begin{equation}
\int_{\overset{\,\!_-\,\!_\wedge}{\M}\,\!_{g = 0}^{E_{\hat \omega}=0}(x, \cdot)}
\pi_{+\infty}^\ast \tilde \alpha.
\label{int on cdot}
\end{equation}
Since we do not need perturbation for
$\overset{\,\!_-\,\!_\wedge}{\M}\,\!_{g = 0}^{E_{\hat \omega}=0}(x, \cdot)$,
(\ref{int on cdot}) coincides with $\langle x, \alpha \rangle$.
\end{proof}

\begin{proof}[Proof of Lemma \ref{trivial generating function}]
First we prove the case of trivial short concordance.
For each $A \geq 0$, define an ideal $\I^\delta_A \subset \D\D_X^{\leq 0, \delta}$ by
\begin{align*}
\I^\delta_A = \{&\sum
a_{(x_i), (\hat c_i^\ast), (\hat c'_i), g} t_{x_1} \dots t_{x_{k_t}} q^-_{\hat c_1^\ast}
\dots q^-_{\hat c_{k_q}^\ast} p^+_{\hat c'_1} \dots p^+_{\hat c'_{k_p}} \hbar^g
\in \D\D_X^{\leq 0, \delta};\\
&a_{(x_i), (\hat c_i^\ast), (\hat c'_i), g} = 0 \text{ if }
\widetilde{g}_\delta \leq A\},
\end{align*}
and define $\I_A^{\star, \delta}
= \I^\delta_A \cap (\hbar^{-1} \D_X^{\leq 0})^{\star, \delta}$.
Then the generating function $\F$ satisfies
\[
\F \equiv \hbar^{-1} \sum_{\hat c, \hat c'} [\overline{\M}^{E_{\hat \omega}=0}_{g = 0, \# \pm\infty = 1}
(\hat c, [\overline{P}] \cap (\hat c')^\ast)]^0 q^-_{\hat c^\ast} p^+_{\hat c'}
\]
in $(\hbar^{-1} \D_X^{\leq 0})^{\star, \delta}
/ (J^{\star, \delta}_{\overline{C}_0, \overline{C}_1, \overline{C}_2} + \I_0^{\star, \delta})$

Let
\begin{align*}
C_\ast(\overline{P}_Y, \overline{P}_Y^{\text{bad}}; \S^D \otimes \Q)
&\to C_\ast(\overline{P}_Y, \overline{P}_Y^{\text{bad}};\ab \S^D \otimes \Q)\\
\hat c &\mapsto \sum_{\hat c'} a_{\hat c, \hat c'} \hat c'
\end{align*}
be the chain homotopy from $\varphi$ to $\id$ given in the above lemma,
that is, the family $a_{\hat c, \hat c'}$ satisfies
\[
\hat c - \sum_{c'} [\overline{\M}^{E_{\hat \omega}=0}_{g = 0, \# \pm\infty = 1}(\hat c,
[\overline{P}] \cap (\hat c')^\ast)]^0 \hat c'
= \sum_{\hat c'} a_{\hat c, \hat c'} \partial \hat c'
+ \sum_{\hat c'} a_{\partial \hat c, \hat c'} \hat c'
\]
for any $\hat c$.
Define
\[
\K = \hbar^{-1} \sum_{c, c'} a_{\hat c, \hat c'} q_{\hat c^\ast} p_{\hat c'}
\in (\hbar^{-1} \D_X^{\leq 0})^{\star, \delta}
/ J^{\star, \delta}_{\overline{C}_0, \overline{C}_1, \overline{C}_2}.
\]
Then $e^{[\widehat{D}_X, \tau \K]} e^{\F}$ ($\tau \in [0, 1]$) is a homotopy from
$\F$ to a generating function $\F^1$ which satisfies
\begin{equation}
\F^1 \equiv \hbar^{-1} \sum_{\hat c} q^-_{\hat c^\ast} p^+_{\hat c} \label{F trivial}
\end{equation}
in $(\hbar^{-1} \D_X^{\leq 0})^{\star, \delta}
/ (J^{\star, \delta}_{\overline{C}_0, \overline{C}_1, \overline{C}_2} + \I_0^{\star, \delta})$.
Hence we may assume $\F$ also satisfies the above equation
in $(\hbar^{-1} \D_X^{\leq 0})^{\star, \delta}
/ (J^{\star, \delta}_{\overline{C}_0, \overline{C}_1, \overline{C}_2} + \I_0^{\star, \delta})$.

We claim that there exists $\G \in (\hbar^{-1} \D_X^{\leq 0})^{\star, \delta}
/ J^{\star, \delta}_{\overline{C}_0, \overline{C}_1, \overline{C}_2}$ such that
\[
e^{\F} \star e^{\G} = e^{\F^{\text{\it tri}}}
\]
in $\D\D_X^{\leq 0, \delta}
/ \widetilde{J}^{\leq 0, \delta}_{\overline{C}_0, \overline{C}_1, \overline{C}_2}$, that is,
$\F \Diamond \G = \F^{\text{\it tri}}$ in
$(\hbar^{-1} \D_X^{\leq 0})^{\star, \delta}
/ J^{\star, \delta}_{\overline{C}_0, \overline{C}_1, \overline{C}_2}$.
This can be proved as follows.
Let $0 = A_0 < A_1 < A_2 < \dots $ be all constants $A$ such that
$\bigcap_{\epsilon > 0} (\widetilde{J}^{\leq 0, \delta}_{\overline{C}_0, \overline{C}_1,
\overline{C}_2} + \I^\delta_{A - \epsilon}) \supsetneq
\widetilde{J}^{\leq 0, \delta}_{\overline{C}_0, \overline{C}_1, \overline{C}_2} + \I^\delta_A$.
Since $\F$ satisfies equation (\ref{F trivial}),
it is easy to construct $\G_{\leq m} = \G_0 + \G_1 + \dots + \G_m \in
(\hbar^{-1} \D_X^{\leq 0})^{\star, \delta}
/ (J^{\star, \delta}_{\overline{C}_0, \overline{C}_1, \overline{C}_2}
+ \I^{\star, \delta}_{A_m})$ inductively such that
$\G_{\leq m} \equiv \G_{\leq m-1}$ in
$(\hbar^{-1} \D_X^{\leq 0})^{\star, \delta}
/ (J^{\star, \delta}_{\overline{C}_0, \overline{C}_1, \overline{C}_2}
+ \I^{\star, \delta}_{A_{m-1}})$ and
\[
e^{\F} \star e^{\G_{\leq m}} \equiv e^{\F^{\text{\it tri}}}
\]
in $\D\D_X^{\leq 0, \delta}
/ (\widetilde{J}^{\leq 0, \delta}_{\overline{C}_0, \overline{C}_1, \overline{C}_2}
+ \I^\delta_{A_m})$.
Therefore we can construct a required $\G$.


Since the composition of $X$ and $X$ is isomorphic to $X$,
$\F \Diamond \F$ is homotopic to $\F$.
Hence $\F \Diamond \F \Diamond \G$ is homotopic to $\F \Diamond \G$.
Therefore, any generating function $\F$ ($= \F \Diamond \F^{\text{\it tri}}$) of $X$
is homotopic to $\F^{\text{\it tri}}$.
(All generating functions
$\F$, $\F \Diamond \F$, $\F \Diamond \F \Diamond \G$ and $\F \Diamond \G$
are elements of $(\hbar^{-1} \D_X^{\leq 0})^{\star, \delta}
/ J^{\star, \delta}_{\overline{C}_0, \overline{C}_1, \overline{C}_2}$.)

Finally we consider the case of general trivial concordance.
Since $\omega|_{[0, T_0] \times Y} = d(f \lambda^-)$ for some
smooth function $f : [0, T_0] \times Y \to \R_{>0}$ such that
$f \lambda^-|_{\{0\} \times Y} = \lambda^-$ and
$f \lambda^-|_{\{T_0\} \times Y} = \lambda^+$,
$(X, \omega)$ is isomorphic to the trivial short concordance
$((-\infty, 0] \cup [0, \infty)) \times Y^+$ of $Y^+$ by
\begin{align*}
[0, T_0] \times Y &\inj (-\infty, 0] \times Y^+ \\
(\sigma, y) &\mapsto (\log f(\sigma), y)
\end{align*}
and
\begin{align*}
(-\infty, 0] \times Y^- &\inj (-\infty, 0] \times Y^+ \\
(\sigma, y) &\mapsto (\sigma + \log a, y).
\end{align*}
We can construct the generating function for $(X, \omega)$ by
the same data as those for the trivial short concordance of $Y^+$.
Then it is easy to check that this generating function is also homotopic to
the trivial generating function.
\end{proof}

Let $(Y^\pm, \lambda^\pm)$ be two arbitrary contact manifolds and
$(X, \omega)$ be an arbitrary cobordism from $(Y^-, \lambda^-)$ to $(Y^+, \lambda)$.
We assume that the generating functions for $Y^-$ and $Y^+$ are defined by
$\mathcal{C}_{Y^-} = (Y^-, \lambda^-, K_{Y^-}, K_{Y^-}^0, K_{Y^-}^2, J^-,
\B_{Y^-})$ and $\mathcal{C}_{Y^+}
= (Y^+, \lambda^+, K_{Y^+}, K_{Y^+}^0, K_{Y^+}^2, J^+, \B_{Y^+})$ respectively,
and that the generating function $\F_X \in (\hbar^{-1} \D_X^{\leq 0})^{\star, \delta}
/ J^{\star, \delta}_{\overline{C}_0, \overline{C}_1, \overline{C}_2}$ for $X$
are defined by the data
$\mathcal{C}_X = (X, \omega, Y^\pm, \lambda^\pm, K_{Y^\pm},\ab
K_X^0,\ab K^0_{Y^\pm},\ab \mu^\pm, K^2_{Y^\pm}, J, \B_X)$ compatible with
$\mathcal{C}_{Y^-}$ and $\mathcal{C}_{Y^+}$.
The argument in Section \ref{Homotopy} implies that the homotopy type of
$\F_X$ does not depend on
the choice of $\mathcal{C}_X$ if we fix $\mathcal{C}_{Y^\pm}$.
We denote the cohomology
$H^\ast(\D_X^{\leq \kappa} / J^{\leq \kappa, \delta}_{C_0, C_1, C_2}, D_{\F})$
for $\mathcal{C}_X$ by
$H^\ast(\D_{\mathcal{C}_X}^{\leq \kappa} / J^{\leq \kappa, \delta}_{C_0, C_1, C_2},
D_{\mathcal{C}_X})$.
Then this implies that cohomologies
$H^\ast(\D_{\mathcal{C}_X}^{\leq \kappa} / J^{\leq \kappa, \delta}_{C_0, C_1, C_2},
D_{\mathcal{C}_X})$
for $\mathcal{C}_X$ compatible with a fixed pair
$(\mathcal{C}_{Y^-}, \mathcal{C}_{Y^+})$ (and with the same $\mu^\pm$)
are naturally isomorphic.
Namely, for every pair $(\mathcal{C}_X, \mathcal{C}'_X)$,
there exists a unique isomorphism
\[
T_{\mathcal{C}'_X, \mathcal{C}_X}
: H^\ast(\D_{\mathcal{C}_X}^{\leq \kappa} / J^{\leq \kappa, \delta}_{C_0, C_1, C_2},
D_{\mathcal{C}_X})
\to H^\ast(\D_{\mathcal{C}'_X}^{\leq \kappa} / J^{\leq \kappa, \delta}_{C_0, C_1, C_2},
D_{\mathcal{C}'_X}),
\]
and these isomorphisms satisfy
$T_{\mathcal{C}_X, \mathcal{C}_X} = \id$ and
$T_{\mathcal{C}''_X, \mathcal{C}'_X} \circ T_{\mathcal{C}'_X, \mathcal{C}_X}
= T_{\mathcal{C}''_X, \mathcal{C}_X}$.

Similarly, cohomologies $H^\ast(\mathcal{L}_{\mathcal{C}_X}^{\leq \kappa}
/ J^{\leq \kappa}_{C_0, C_2}, d_{\mathcal{C}_X})
= H^\ast(\mathcal{L}_{X}^{\leq \kappa}
/ J^{\leq \kappa}_{C_0, C_2}, d_{\F_0})$
for $\mathcal{C}_X$ compatible with a fixed pair
$(\mathcal{C}_{Y^-}, \mathcal{C}_{Y^+})$ (and with the same $\mu^\pm$)
are naturally isomorphic.

Next we compare two SFT cohomologies of $X$ compatible with different pairs
$\mathcal{C}_{Y^\pm}$ for $(Y^\pm, \xi^\pm)$.
First we treat the case where we do not change the contact forms $\lambda^\pm$.
(To treat the general case, we cannot fix a filtration and need to take the limit
with respect to the filtration.)
\begin{lem}\label{short concordance for X}
Let $\mathcal{C}_X$ be a cobordism from $\mathcal{C}_{Y^-}$ to $\mathcal{C}_{Y^+}$,
and let
\begin{multline*}
\mathcal{C}_{X_1}
= (X_1, \omega_1, (Y^+, Y^+), (\lambda^+, \lambda^+), (K_{Y^+}, K_{Y^+_1}), K_X^0,
(K^0_{Y^+}, K^0_{Y^+_1}),\\
\mu^\pm, (K^2_{Y^+}, K^2_{Y^+_1}), J, \B_X)
\end{multline*}
be a short concordance from
$\mathcal{C}_{Y^+}$
to $\mathcal{C}_{Y^+_1}
= (Y^+, \lambda^+, K_{Y^+_1}, K_{Y^+_1}^0, K_{Y^+_1}^2, J^+_1, \B_{Y^+_1})$.
Then
\[
T_{\F_X}(\cdot \Diamond \F_{X_1})
: (\D^{\leq \kappa}_X / J^{\leq \kappa, \delta}_{C_0, C_1, C_2}, D_{\F_X})
\to (\D^{\leq \kappa}_{X \# X_1} / J^{\leq \kappa, \delta}_{C_0, C_1, C_2},
D_{\F_X \Diamond \F_{X_1}})
\]
and
\[
T_{(\F_X)_0}(\cdot \sharp (\F_{X_1})_0)
: (\mathcal{L}^{\leq \kappa}_X / J^{\leq \kappa}_{C_0, C_2}, d_{(\F_X)_0})
\to (\mathcal{L}^{\leq \kappa}_{X \# X_1} / J^{\leq \kappa}_{C_0, C_2},
d_{(\F_X) \sharp (\F_{X_1})_0})
\]
are chain homotopy equivalences.
\end{lem}
\begin{proof}
First we consider the case of general SFT.
Let
\begin{multline*}
\mathcal{C}_{X_2}
= (X_2, \omega_2, (Y^+, Y^+), (\lambda^+, \lambda^+), (K_{Y^+_1}, K_{Y^+}), K_X^0,
(K^0_{Y^+_1}, K^0_{Y^+}),\\
\mu^\pm, (K^2_{Y^+_1}, K^2_{Y^+}), J, \B_X)
\end{multline*}
be a short concordance from
$\mathcal{C}_{Y^+_1}$
to $\mathcal{C}_{Y^+}$.
Since $X_1 \# X_2$ is a trivial short concordance, its generating function
$\F_{X_1} \Diamond \F_{X_2}$ is homotopic to the trivial generating function
$\F^{\text{\it tri}}$.
Hence Lemma \ref{linearized T} (v) implies that
\begin{align*}
T_{\F_X}(\cdot \Diamond (\F_{X_1} \Diamond \F_{X_2}))
&: (\D^{\leq \kappa}_X / J^{\leq \kappa, \delta}_{C_0, C_1, C_2}, D_{\F_X}) \\
&\quad \to (\D^{\leq \kappa}_{X \# X_1 \# X_2} / J^{\leq \kappa, \delta}_{C_0, C_1, C_2},
D_{\F_X \Diamond \F_{X_1} \Diamond \F_{X_2}})
\end{align*}
is a chain homotopy equivalence.
By Lemma \ref{linearized T} (iv), this map coincides with the composition
$T_{\F_X \Diamond \F_{X_1}}(\cdot \Diamond \F_{X_2}) \circ
T_{\F_X}(\cdot \Diamond \F_{X_1})$.
Hence $T_{\F_X}(\cdot \Diamond \F_{X_1})$ has a left homotopy inverse and
$T_{\F_X \Diamond \F_{X_1}}(\cdot \Diamond \F_{X_2})$ has a right homotopy inverse.
Since we can apply the above argument for a cobordism $X \# X_1$ and
a short concordance $\mathcal{C}_{X_2}$, 
$T_{\F_X \Diamond \F_{X_1}}(\cdot \Diamond \F_{X_2})$ has a right inverse.
Hence $T_{\F_X}(\cdot \Diamond \F_{X_1})$ is a chain homotopy equivalence.

The case of rational SFT is similar.
\end{proof}

The above lemma and the counterpart of a short concordance from
$\mathcal{C}_{Y^-_1}$ to $\mathcal{C}_{Y^-}$ imply that homologies
$H^\ast(\D_{\mathcal{C}_X}^{\leq \kappa} / J^{\leq \kappa, \delta}_{C_0, C_1, C_2},
D_{\mathcal{C}_X})$ and
$H^\ast(\mathcal{L}_{\mathcal{C}_X}^{\leq \kappa}
/ J^{\leq \kappa}_{C_0, C_2}, d_{\mathcal{C}_X})$
for $\mathcal{C}_X$ compatible with a fixed pair
$((Y^-, \lambda^-, \overline{K}^0_{Y^-}),
(Y^+, \lambda^+, \overline{K}^0_{Y^+}))$ (and with the same $\mu^\pm$)
are naturally isomorphic respectively.
(The naturality is due to Lemma \ref{linearized T} (iv), (v) (or its rational version) and
Lemma \ref{trivial generating function}.)
Therefore for any cobordism $(X, \omega)$ between two strict contact manifolds
$(Y^\pm, \lambda^\pm)$ and any $(\overline{K}^0_X, \overline{K}^0_{Y^\pm},
\mu^\pm)$,
we can define the limits
\begin{align}
&H^\ast_{\mathrm{SFT}}(X, \omega, Y^\pm, \lambda^\pm, \overline{K}^0_X,
\overline{K}^0_{Y^\pm}, \mu^\pm) \notag\\
&= \varprojlim_{C_2} \varinjlim_{\kappa, \delta} \varprojlim_{C_0, C_1}
H^\ast(\D_{\mathcal{C}_X}^{\leq \kappa} / J^{\leq \kappa, \delta}_{C_0, C_1, C_2},
D_{\mathcal{C}_X}).
\label{limit general X}
\end{align}
and
\begin{equation}
H^\ast_{\mathrm{RSFT}}(X, \omega, Y^\pm, \lambda^\pm, \overline{K}^0_X,
\overline{K}^0_{Y^\pm}, \mu^\pm)
= \varprojlim_{C_2} \varinjlim_{\kappa} \varprojlim_{C_0}
H^\ast(\mathcal{L}_{\mathcal{C}_X}^{\leq \kappa} / J^{\leq \kappa}_{C_0, C_2},
d_{\mathcal{C}_X}).
\label{limit rational X}
\end{equation}
We sometimes abbreviate these limits as $H^\ast(\D_X, D_X)$ and
$H^\ast(\mathcal{L}_X, d_X)$ respectively.
We will prove that these cohomology groups do not depend on the choice of
the contact forms of $(Y^\pm, \xi^\pm)$ later.
It is easy to check that for a pair of composable cobordism $(X, \omega)$ and
$(X', \omega')$, the limit of the linearizations of the composition maps define maps
\[
T_{\F_X}(\cdot \Diamond \F_{X'}) : H^\ast(\D_X, D_X)
\to H^\ast(\D_{X \# X'}, D_{X \# X'}),
\]
\[
T_{\F_{X'}}(\F_X \Diamond \cdot) : H^\ast(\D_{X'}, D_{X'})
\to H^\ast(\D_{X \# X'}, D_{X \# X'}),
\]
\[
T_{(\F_X)_0}(\cdot \sharp (\F_{X'})_0) : H^\ast(\mathcal{L}_X, d_X)
\to H^\ast(\mathcal{L}_{X \# X'}, d_{X \# X'}),
\]
and
\[
T_{(\F_{X'})_0}((\F_X)_0 \sharp \cdot) : H^\ast(\mathcal{L}_{X'}, d_{X'})
\to H^\ast(\mathcal{L}_{X \# X'}, d_{X \# X'}).
\]

Next we consider the SFT cohomologies of a contact manifold.
First we compare two cohomology groups defined by the same contact form with
different other data.
\begin{lem}\label{short concordance for Y}
For a short concordance $\mathcal{C}_X$ from $\mathcal{C}_{Y^-}$
to $\mathcal{C}_{Y^+}$, the linear maps
\[
i_{\F_X}^\pm : (\W^{\leq \kappa}_{Y^\pm} / I^{\leq \kappa}_{C_0, C_1, C_2}, D_{Y^\pm})
\to (\D^{\leq \kappa, L_{\min}}_X / J^{\leq \kappa, L_{\min}}_{C_0, C_1, C_2}, D_{\F_X})
\]
and
\[
i_{(\F_X)_0}^\pm
: (\mathcal{P}^{\leq \kappa}_{Y^\pm} / I^{\leq \kappa}_{C_0, C_2}, d_{Y^\pm})
\to (\mathcal{L}^{\leq \kappa}_X / J^{\leq \kappa}_{C_0, C_2}, d_{(\F_X)_0})
\]
are chain homotopy equivalences,
and the compositions of the induced maps
\[
A = i_{\F_X}^- \circ (i_{\F_X}^+)^{-1}
: H^\ast(\W^{\leq \kappa}_{Y^+} / I^{\leq \kappa}_{C_0, C_1, C_2}, D_{Y^+})
\to H^\ast(\W^{\leq \kappa}_{Y^-} / I^{\leq \kappa}_{C_0, C_1, C_2}, D_{Y^-})
\]
and
\[
A^0 = i_{(\F_X)_0}^- \circ (i_{(\F_X)_0}^+)^{-1}
: H^\ast(\mathcal{P}^{\leq \kappa}_{Y^+} / I^{\leq \kappa}_{C_0, C_2}, d_{Y^+})
\to H^\ast(\mathcal{P}^{\leq \kappa}_{Y^-} / I^{\leq \kappa}_{C_0, C_2}, d_{Y^-})
\]
do not depend on the short concordance $\mathcal{C}_X$.
\end{lem}
\begin{proof}
We prove the case of general SFT. The case of rational SFT is similar.
First we consider the case of a trivial short concordance.
Note that for the trivial generating function $\F^{\text{\it tri}}$,
$i_{\F^{\text{\it tri}}}^\pm$ coincide with the identity map under the natural
identification $\W^{\leq \kappa}_{Y^\pm} / I^{\leq \kappa}_{C_0, C_1, C_2}
\cong \D^{\leq \kappa, L_{\min}}_X / J^{\leq \kappa, L_{\min}}_{C_0, C_1, C_2}$
which maps $q_{\hat c^\ast}$ and $p_{\hat c}$ to $q^-_{\hat c^\ast}$ and
$p^+_{\hat c}$ respectively.
Since the generating function $\F_X$ is homotopic to the trivial generating function,
Lemma \ref{homotopy and T} (iv) implies that
$i_{\F_X}^\pm$ is chain homotopic to the composition of $i_{\F^{\text{\it tri}}}^\pm$
and the isomorphism defined by the homotopy.
Hence $i_{\F_X}^\pm$ are also chain homotopy equivalence.

Next we consider the general case.
Let $\mathcal{C}_{X'}$ be a short concordance from $\mathcal{C}_{Y^+}$ to
$\mathcal{C}_{Y^-}$.
Then $T_{\F_X}(\cdot \Diamond \F_{X'}) \circ i_{\F_X}^-
= i_{\F_X \Diamond \F_{X'}}^-
: (\W^{\leq \kappa}_{Y^-} / I^{\leq \kappa}_{C_0, C_1, C_2}, D_{Y^\pm})
\to (\D^{\leq \kappa}_{X \# X'} / J^{\leq \kappa, L_{\min}}_{C_0, C_1, C_2},
D_{\F_X \Diamond \F_{X'}})$
are chain homotopy equivalence since $\F_X \Diamond \F_{X'}$ is homotopic to the
generating function of a trivial short concordance.
Since $T_{\F_X}(\cdot \Diamond \F_{X'})$ is also a chain homotopy equivalence by
Lemma \ref{short concordance for X}, so is $i_{\F_X}^-$.
Similarly, $i_{\F_X}^+$ is also a chain homotopy equivalence.

Finally we check the independence of $A = i_{\F_X}^- \circ (i_{\F_X}^+)^{-1}$.
For any two short concordances $\mathcal{C}_X$ and $\mathcal{C}'_X$,
there exists a smooth family of exact cobordisms
$(X^\tau, \omega^\tau)_{\tau \in I}$ such that $(X^0, \omega^0) = (X, \omega)$
and $(X^1, \omega^1) = (X', \omega')$, and
we can construct a homotopy from $\F_X$ to $\F_{X'}$.
This implies that there exists an isomorphism
$T : (\D_X / J_{C_0, C_1}^\delta, D_\F) \to
(\D_{X'} / J_{C_0, C_1}^\delta, D_{\F'})$
such that $i_{\F_{X'}}^\pm$ coincides with $T \circ i_{\F_X}^\pm$ up to chain homotopy.
Hence $A$ does not depend on the choice of the short concordance.
\end{proof}
We denote the isomorphisms $A$ and $A^0$ in the above lemma by
$A_{\mathcal{C}_{Y^-}, \mathcal{C}_{Y^+}}$ and
$A_{\mathcal{C}_{Y^-}, \mathcal{C}_{Y^+}}^0$ respectively.
The above lemma implies that if we fix $(Y, \lambda, \overline{K}_Y^0)$,
the cohomologies
\[
H^\ast(\W_{(Y, \lambda, K_Y, \overline{K}_Y^0)}^{\leq \kappa}
/ I^{\leq \kappa}_{C_0, C_1, C_2}, D_{(Y, \lambda, K_Y, K_Y^0, K_Y^2, J, \B)})
\]
and
\[
H^\ast(\mathcal{P}_{(Y, \lambda, K_Y, \overline{K}_Y^0)}^{\leq \kappa}
/ I^{\leq \kappa}_{C_0, C_2}, d_{(Y, \lambda, K_Y, K_Y^0, K_Y^2, J, \B)})
\]
defined by various data $\mathcal{C}_Y = (Y, \lambda, K_Y, K_Y^0, K_Y^2, J, \B)$
of the same $(Y, \lambda, \overline{K}_Y^0)$ are naturally isomorphic respectively, and
the isomorphisms are given by the above $A_{\mathcal{C}_{Y^-}, \mathcal{C}_{Y^+}}$
and $A_{\mathcal{C}_{Y^-}, \mathcal{C}_{Y^+}}^0$.
The naturality of isomorphisms $A_{\mathcal{C}_{Y^-}, \mathcal{C}_{Y^+}}$
is proved as follows.
$A_{\mathcal{C}_Y, \mathcal{C}_Y} = \id$ is due to Lemma \ref{trivial generating function}.
$A_{\mathcal{C}_{Y''}, \mathcal{C}_{Y'}} \circ A_{\mathcal{C}_{Y'}, \mathcal{C}_Y}
= A_{\mathcal{C}_{Y''}, \mathcal{C}_Y}$ is because the following diagram is commutative
by Lemma \ref{linearized T},
where $\mathcal{C}_{X^-}$ is a short concordance from
$\mathcal{C}_{Y''}$ to $\mathcal{C}_{Y'}$ with a generating function $\F^-$
and $\mathcal{C}_{X^+}$ is a short concordance from
$\mathcal{C}_{Y'}$ to $\mathcal{C}_Y$ with a generating function $\F^+$,
and we abbreviate $H^\ast(\W_Y^{\leq \kappa} / I^{\leq \kappa}_{C_0, C_1, C_2}, D_Y)$
or $H^\ast(\D_{\mathcal{C}_{X}}^{\leq \kappa} / J^{\leq \kappa, L_{\min}}_{C_0, C_1, C_2})$
by $H^\ast(\mathcal{C}_Y)$ or $H^\ast(\mathcal{C}_X)$ respectively.
\[
\begin{tikzcd}
H^\ast(\mathcal{C}_{Y''})
\ar{r}{i^-_{\F^-}}
\ar{ddrr}[swap]{i^-_{\F^- \Diamond \F^+}}
&H^\ast(\mathcal{C}_{X^-})
\ar{ddr}{\hspace{-5pt}T_{\F^-}(\cdot \Diamond \F^+)}
&\ar{l}[swap]{i^+_{\F^-}}
H^\ast(\mathcal{C}_{Y'})
\ar{r}{i^-_{\F^+}}
&H^\ast(\mathcal{C}_{X^+})
\ar{ddl}[swap, near end]{T_{\F^+}(\F^- \Diamond \cdot) \hspace{-5pt}}
&\ar{l}[swap]{i^+_{\F^+}}
H^\ast(\mathcal{C}_Y)
\ar{ddll}{i^+_{\F^- \Diamond \F^+}}\\
\\
&&H^\ast(\mathcal{C}_{X^- \# X^+})&&
\end{tikzcd}
\]
Therefore we can define the limits
\begin{align}
&H^\ast_{\mathrm{SFT}}(Y, \lambda, \overline{K}^0_Y) \notag\\
&= \varprojlim_{C_2} \varinjlim_{\kappa} \varprojlim_{C_0, C_1}
H^\ast(\W_{(Y, \lambda, K_Y, \overline{K}_Y^0)}^{\leq \kappa}
/ I^{\leq \kappa}_{C_0, C_1, C_2}, D_{(Y, \lambda, K_Y, K_Y^0, K_Y^2, J, \B)}).
\label{limit general Y}
\end{align}
and
\begin{align}
&H^\ast_{\mathrm{RSFT}}(Y, \lambda, \overline{K}^0_Y) \notag\\
&= \varprojlim_{C_2} \varinjlim_{\kappa} \varprojlim_{C_0}
H^\ast(\mathcal{P}_{(Y, \lambda, K_Y, \overline{K}_Y^0)}^{\leq \kappa}
/ I^{\leq \kappa}_{C_0, C_2}, d_{(Y, \lambda, K_Y, K_Y^0, K_Y^2, J, \B)}).
\label{limit rational Y}
\end{align}
We sometimes abbreviate these limits as $H^\ast(\W_Y, D_Y)$ and
$H^\ast(\mathcal{P}_Y, d_Y)$ respectively.

For any cobordism $(X, \omega)$ from $(Y^-, \lambda^-)$ to
$(Y^+, \lambda^+)$, we can define
\begin{equation}
i_{X}^\pm : H^\ast_{\mathrm{SFT}}(Y^\pm, \lambda^\pm, \overline{K}_{Y^\pm}^0)
\to H^\ast_{\mathrm{SFT}}(X, \omega, Y^\pm, \lambda^\pm, \overline{K}^0_X,
\overline{K}^0_{Y^\pm}, \mu^\pm)
\label{well def of i_X^pm}
\end{equation}
and
\begin{equation}
i_{X, 0}^\pm : H^\ast_{\mathrm{RSFT}}(Y^\pm, \lambda^\pm, \overline{K}^0_{Y^\pm})
\to H^\ast_{\mathrm{RSFT}}(X, \omega, Y^\pm, \lambda^\pm, \overline{K}^0_X,
\overline{K}^0_{Y^\pm}, \mu^\pm)
\label{well def of i_{X, 0}^pm}
\end{equation}
by the limits of $i_{\F}^\pm$ and $i_{\F_0}^\pm$ respectively.
For example, the well-definedness of $i_X^+$ is due to the following fact:
Let $\mathcal{C}_{Y^\pm}$ and $\mathcal{C}'_{Y^\pm}$ be two data for
$(Y^\pm, \lambda^\pm)$,
$\mathcal{C}_{X}$ be a cobordism from $\mathcal{C}_{Y^-}$ to $\mathcal{C}_{Y^+}$,
and $\mathcal{C}_{X'}$ be a cobordism from $\mathcal{C}'_{Y^-}$ to $\mathcal{C}'_{Y^+}$.
Assume that both of $\mathcal{C}_{X}$ and $\mathcal{C}_{X'}$ are data of
the same cobordism $(X, \omega)$.
Then the following diagram is commutative by Lemma \ref{linearized T},
where $\mathcal{C}_{X_0}$ is a short concordance from $\mathcal{C}_{Y^+}$ to
$\mathcal{C}'_{Y^+}$, $\mathcal{C}_{X_1}$ is a short concordance from
$\mathcal{C}_{Y^-}$ to $\mathcal{C}'_{Y^-}$, and
$T : H^\ast(\mathcal{C}_{X \# X_0}) \to H^\ast(\mathcal{C}_{X_1 \# X'})$
is the isomorphism for a homotopy from $\F_{X} \Diamond \F_{X_0}$ to
$\F_{X_1} \Diamond \F_{X'}$.
The left column is the natural isomorphism for the SFT cohomology for
$(Y^+, \lambda^+)$,
and the right column is the natural isomorphism for the SFT cohomology for
$(X, \omega)$.
Therefore the compatibility of these isomorphisms and the maps $i^+_{\F_X}$,
$i^+_{\F_{X'}}$ implies the well-definedness of $i^+_{X}$.
\[
\begin{tikzcd}[column sep = huge]
H^\ast(\mathcal{C}_{Y^+})
\ar{rr}{i^+_{\F_X}}
\ar{d}[swap]{i^-_{\F_{X_0}}}
&&H^\ast(\mathcal{C}_{X})
\ar{d}{T_{\F_X}(\cdot \Diamond \F_{X_0})}\\
H^\ast(\mathcal{C}_{X_0})
\ar{rr}{T_{\F_{X_0}}(\F_X \Diamond \cdot)}
&&H^\ast(\mathcal{C}_{X \# X_0})
\ar{d}{T}\\
H^\ast(\mathcal{C}'_{Y^+})
\ar{u}{i^+_{\F_{X_0}}}
\ar{urr}{i^+_{\F_{X} \Diamond \F_{X_0}}\hspace{-10pt}}
\ar{rr}[swap]{\hspace{30pt} i^+_{\F_{X_1} \Diamond \F_{X'}}}
\ar{drr}[swap]{i^+_{\F_{X'}}}
&&H^\ast(\mathcal{C}_{X_1 \# X'})\\
&&H^\ast(\mathcal{C}_{X'})
\ar{u}[swap]{T_{\F_{X'}}(\F_{X_1} \Diamond \cdot)}
\end{tikzcd}
\]

\begin{prop}\label{algebraic structures of limits}
The limits of SFT cohomologies satisfies the following.
\begin{enumerate}[label= \normalfont (\roman*)]
\item
$H^\ast_{\mathrm{SFT}}(Y, \lambda, \overline{K}^0)$ inherits a structure of algebra.
\item
$H^\ast_{\mathrm{RSFT}}(Y, \lambda, \overline{K}^0)$ inherits a structure
of Poisson algebra.
\item
$H^\ast_{\mathrm{SFT}}(X, \omega, Y^\pm, \lambda^\pm, \overline{K}^0_X,
\overline{K}^0_{Y^\pm}, \mu^\pm)$ has a structure of
$H^\ast_{\mathrm{SFT}}(Y^\pm, \lambda^\pm, \overline{K}_{Y^\pm}^0)$-bimodule.
\item
$H^\ast_{\mathrm{RSFT}}(X, \omega, Y^\pm, \lambda^\pm, \overline{K}^0_X,
\overline{K}^0_{Y^\pm}, \mu^\pm)$ has a structure of
$H^\ast_{\mathrm{RSFT}}(Y^\pm, \lambda^\pm, \overline{K}_{Y^\pm}^0)$-bimodule.
\end{enumerate}
\end{prop}
\begin{rem}
$H^\ast_{\mathrm{RSFT}}(X, \omega, Y^\pm, \lambda^\pm, \overline{K}^0_X,
\overline{K}^0_{Y^\pm}, \mu^\pm)$ does not have a structure of Poisson module over
$H^\ast_{\mathrm{RSFT}}(Y^\pm, \lambda^\pm, \overline{K}_{Y^\pm}^0)$.
\end{rem}
\begin{proof}
First we consider (i) and (iii).
(\ref{D_Y diff alg}) implies that the multiplication of $\W_Y$ induces maps
\begin{align}
&H^\ast(\W_Y^{\leq \kappa_1} / I^{\leq \kappa_1}_{C_0, C_1 + \kappa_2 L_{\min}^{-1},
C_2 + \kappa_2}, D_Y) \times
H^\ast(\W_Y^{\leq \kappa_2} / I^{\leq \kappa_2}_{C_0, C_1 + \kappa_1 L_{\min}^{-1}, C_2},
D_Y) \notag\\
&\hph{H^\ast(\W_Y^{\leq \kappa_1} / I^{\leq \kappa_1}_{C_0, C_1 + \kappa_2 L_{\min}^{-1},
C_2 + \kappa_2}, D_Y)}
\to H^\ast(\W_Y^{\leq \kappa_1 + \kappa_2}
/ I^{\leq \kappa_1 + \kappa_2}_{C_0, C_1, C_2}, D_Y).
\label{multiplication of homology of W_Y with filtration}
\end{align}
Similarly, (\ref{D_F left Leibnitz}) and (\ref{D_F right Leibnitz}) imply that for
any cobordism $(X, \omega)$ from $(Y^-, \lambda^-)$ to $(Y^+, \lambda^+)$, the
$\W_{Y^\pm}$-bimodule structure of $\D_X$ induces maps
\begin{multline*}
H^\ast(\W_{Y^-}^{\leq \kappa_1} / I^{\leq \kappa_1}_{C_0, C'_1, C_2 + \kappa_2}, D_Y)
\times H^\ast(\D_X^{\leq \kappa_2} / J^{\leq \kappa_2, \delta}_{C_0, C_1
+ \kappa_1 \delta^{-1}, C_2}, D_{\F_X})\\
\to H^\ast(\D_X^{\leq \kappa_1 + \kappa_2} / J^{\leq \kappa_1 + \kappa_2, \delta}
_{C_0, C_1, C_2}, D_{\F_X}),
\end{multline*}
where $C'_1 = C_1 + \kappa_1 (\delta^{-1} - L_{\min}^{-1}) + \kappa_2 L_{\min}^{-1}$,
and
\begin{multline*}
H^\ast(\D_X^{\leq \kappa_1} / J^{\leq \kappa_1, \delta}_{C_0, C_1
+ \kappa_2 \delta^{-1}, C_2 + \kappa_1}, D_{\F_X}) \times
H^\ast(\W_{Y^+}^{\leq \kappa_2} / I^{\leq \kappa_2}_{C_0, C''_1, C_2}, D_{Y^+})\\
\to H^\ast(\D_X^{\leq \kappa_1 +\kappa_2} / J^{\leq \kappa_1 + \kappa_2, \delta}
_{C_0, C_1, C_2}, D_{\F_X}),
\end{multline*}
where $C''_1 = C_1 + \kappa_1 (\delta^{-1} - L_{\min}^{-1})$.
These multiplications satisfy the associativity condition.
Therefore the map $A = i_X^- \circ (i_X^+)^{-1}$ in Lemma
\ref{short concordance for Y} preserves the multiplication.
Namely, for any
$f \in H^\ast(\W^{\leq \kappa_1}_{Y^+}
/ I^{\leq \kappa}_{C_0, C_1 + \kappa_2 L_{\min}^{-1}, C_2 + \kappa_2}, D_{Y^+})$ and
$g \in H^\ast(\W^{\leq \kappa_2}_{Y^+}
/ I^{\leq \kappa}_{C_0, C_1 + \kappa_1 L_{\min}^{-1}, C_2}, D_{Y^+})$,
\[
i_{\F}^+(fg)
= 1 \underset{\F}{\overleftarrow{\ast}} f \underset{\F}{\overleftarrow{\ast}} g
= A(f) \underset{\F}{\overrightarrow{\ast}} A(g) \underset{\F}{\overrightarrow{\ast}} 1
= i_{\F}^-(A(f) A(g))
\]
in $H^\ast(\D^{\leq \kappa_1 + \kappa_2, L_{\min}}_X
/ J^{\leq \kappa_1 + \kappa_2, L_{\min}}_{C_0, C_1, C_2}, D_{\F_X})$.
Hence (\ref{multiplication of homology of W_Y with filtration}) depends only on
the triple $(Y, \lambda, \overline{K}^0)$, and it does not depend on the other choices
of $\mathcal{C}_Y$.

(\ref{multiplication of homology of W_Y with filtration}) induces the multiplication of
the limit $H^\ast_{\mathrm{SFT}}(Y, \lambda, \overline{K}^0)$ as follows.
First (\ref{multiplication of homology of W_Y with filtration}) induces
\begin{multline*}
\varinjlim_{\kappa_1} \varprojlim_{C_0, C_1}
H^\ast(\W_Y^{\leq \kappa_1} / I^{\leq \kappa_1}_{C_0, C_1, C_2 + \kappa_2}, D_Y)
\times \varprojlim_{C_0, C_1}
H^\ast(\W_Y^{\leq \kappa_2} / I^{\leq \kappa_2}_{C_0, C_1, C_2}, D_Y) \\
\to \varinjlim_{\kappa} \varprojlim_{C_0, C_1}
H^\ast(\W_Y^{\leq \kappa} / I^{\leq \kappa}_{C_0, C_1, C_2}, D_Y),
\end{multline*}
and this induces
\begin{multline*}
\varprojlim_{C'_2} \varinjlim_{\kappa_1} \varprojlim_{C_0, C_1}
H^\ast(\W_Y^{\leq \kappa_1} / I^{\leq \kappa_1}_{C_0, C_1, C'_2}, D_Y)
\times \varprojlim_{C_0, C_1}
H^\ast(\W_Y^{\leq \kappa_2} / I^{\leq \kappa_2}_{C_0, C_1, C_2}, D_Y) \\
\to \varinjlim_{\kappa} \varprojlim_{C_0, C_1}
H^\ast(\W_Y^{\leq \kappa} / I^{\leq \kappa}_{C_0, C_1, C_2}, D_Y).
\end{multline*}
Then this induces
\begin{align*}
\varprojlim_{C'_2} \varinjlim_{\kappa_1} \varprojlim_{C_0, C_1}
H^\ast(\W_Y^{\leq \kappa_1} / I^{\leq \kappa_1}_{C_0, C_1, C'_2}, D_Y)
\times \varinjlim_{\kappa_2} \varprojlim_{C_0, C_1}
H^\ast(\W_Y^{\leq \kappa_2} / I^{\leq \kappa_2}_{C_0, C_1, C_2}, D_Y) \ \\
\to \varinjlim_{\kappa} \varprojlim_{C_0, C_1}
H^\ast(\W_Y^{\leq \kappa} / I^{\leq \kappa}_{C_0, C_1, C_2}, D_Y),
\end{align*}
and finally this induces the multiplication of the limit.

(iii) also follows from the above argument and a similar argument to the proof of
well-definedness of (\ref{well def of i_X^pm}).
($i_{\F_X}^\pm$ are special case of multiplication.)

Next we consider (ii) and (iv).
A similar argument implies
$H^\ast_{\mathrm{RSFT}}(Y, \lambda, \overline{K}^0)$ inherits a structure of algebra
and that
$H^\ast_{\mathrm{RSFT}}(X, \omega, Y^\pm, \lambda^\pm, \overline{K}^0_X,
\overline{K}^0_{Y^\pm}, \mu^\pm)$ has a structure of
$H^\ast_{\mathrm{RSFT}}(Y^\pm, \lambda^\pm, \overline{K}_{Y^\pm}^0)$-bimodule.
We need to prove that $H^\ast_{\mathrm{RSFT}}(Y, \lambda, \overline{K}^0)$ inherits
a Poisson structure.
First we prove that the map $A^0$ in Lemma \ref{short concordance for Y}
preserves the Poisson structure.
Namely, we prove that for any $f \in H^\ast(\mathcal{P}^{\leq \kappa_1}_{Y^+}
/ I^{\leq \kappa_1}_{C_0, C_2 + \kappa_2}, d_{Y^+})$ and
$g \in H^\ast(\mathcal{P}^{\leq \kappa_2}_{Y^+}
/ I^{\leq \kappa_2}_{C_0, C_2 + \kappa_1}, d_{Y^+})$,
\begin{equation}
A^0(\{f, g\}) = \{A^0(f), A^0(g)\} \label{A^0 Poisson}
\end{equation}
in $H^\ast(\mathcal{P}^{\leq \kappa_1 + \kappa_2}_{Y^-}
/ I^{\leq \kappa_1 + \kappa_2}_{C_0, C_2}, d_{Y^-})$.
We denote the subspace of cycles of a chain complex $(C^\ast, d)$ by $Z(C^\ast, d)$.
Assume that $f^\pm \in Z(\mathcal{P}^{\leq \kappa_1}_{Y^\pm}
/ I^{\leq \kappa_1}_{C_0, C_2 + \kappa_2}, d_{Y^\pm})$ and
$g^\pm \in Z(\mathcal{P}^{\leq \kappa_2}_{Y^\pm}
/ I^{\leq \kappa_2}_{C_0, C_2 + \kappa_1}, d_{Y^\pm})$ satisfy
\begin{align*}
i_{(\F_X)_0}^- f^- - i_{(\F_X)_0}^+ f^+ &= d_{\F_0} a, \\
i_{(\F_X)_0}^- g^- - i_{(\F_X)_0}^+ g^+ &= d_{\F_0} b
\end{align*}
for some $a, b \in \mathcal{L}^{\leq \kappa_1 + \kappa_2}_X
/ J^{\leq \kappa_1 + \kappa_2}_{C_0, C_2}$.
Note that $\{f^-, g^+\} = \{f^+, g^-\} = 0$.
Then
\begin{align*}
i_{(\F_X)_0}^-(\{f^-, g^-\}) - i_{(\F_X)_0}^+(\{f^+, g^+\})
&= (\{f^-, g^-\} - \{f^+, g^+\})|_{\F_0} \\
&= \{f^- - f^+, g^- - g^+\}|_{\F_0}
\end{align*}
is exact in $(\mathcal{L}^{\leq \kappa_1 + \kappa_2}_X
/ J^{\leq \kappa_1 + \kappa_2}_{C_0, C_2}, d_{(\F_X)_0})$
by Proposition \ref{properties of d_{F_0}} (iv) since $\widehat{d}_X(f^- - f^+) = 0$
and $\widehat{d}_X(g^- - g^+) = 0$.
This proves equation (\ref{A^0 Poisson}).
Therefore $A^0$ in Lemma \ref{short concordance for Y} preserves
the Poisson structure.

Recall that for $\kappa^\circ \leq \kappa$ and $C^\circ \leq C$,
$(\mathcal{P}_Y^{\leq \kappa^\circ} + I^{\leq \kappa}_{C_0, C^\circ_2})
/ I^{\leq \kappa}_{C_0, C_2}$ is the fiber product of
$\mathcal{P}_Y^{\leq \kappa^\circ} / I^{\leq \kappa^\circ}_{C_0, C^\circ_2}$
and $\mathcal{P}_Y^{\leq \kappa} / I^{\leq \kappa}_{C_0, C_2}$ over
$\mathcal{P}_Y^{\leq \kappa} / I^{\leq \kappa}_{C_0, C^\circ_2}$,
and the Poisson bracket of $\mathcal{P}_Y$
induces (\ref{Poisson bracket for fiber product}).
Its homology
$H^\ast((\mathcal{P}_Y^{\leq \kappa^\circ} + I^{\leq \kappa}_{C_0, C^\circ_2})
/ I^{\leq \kappa}_{C_0, C_2}, d_Y)$ is also well-defined.
Since $H^\ast$ preserves fiber product structure, it is the fiber product of
$H^\ast(\mathcal{P}_Y^{\leq \kappa^\circ} / I^{\leq \kappa^\circ}_{C_0, C^\circ_2}, d_Y)$
and $H^\ast(\mathcal{P}_Y^{\leq \kappa} / I^{\leq \kappa}_{C_0, C_2}, d_Y)$ over
$H^\ast(\mathcal{P}_Y^{\leq \kappa} / I^{\leq \kappa}_{C_0, C^\circ_2}, \ab d_Y)$.
Furthermore, since fiber product commutes with limits,
$H^\ast_{\mathrm{RSFT}}(Y, \lambda, \overline{K}^0_Y)$ is isomorphic to
\[
\varprojlim_{C_2^\circ} \varinjlim_{\kappa^\circ}
\varprojlim_{C_2} \varinjlim_{\kappa} \varprojlim_{C_0}
H^\ast((\mathcal{P}_Y^{\leq \kappa^\circ} + I^{\leq \kappa}_{C_0, C^\circ_2})
/ I^{\leq \kappa}_{C_0, C_2}, d_Y).
\]
First, (\ref{Poisson bracket for fiber product}) induces the map
\begin{align*}
H^\ast((\mathcal{P}_Y^{\leq \kappa^\circ_1} + I^{\leq \kappa_1}_{C_0, C_2})
/ I^{\leq \kappa_1}_{C_0, C'_2}, d_Y) \times
H^\ast((\mathcal{P}_Y^{\leq \kappa^\circ_2} + I^{\leq \kappa_2}_{C_0, C_2})
/ I^{\leq \kappa_1}_{C_0, C''_2}, d_Y) \\
\to H^\ast(\mathcal{P}_Y^{\leq \kappa_1 + \kappa_2}
/ I^{\leq \kappa_1 + \kappa_2}_{C_0, C_2}, d_Y)
\end{align*}
for $C'_2 \geq C_2 + \kappa^\circ_2$ and $C''_2 \geq C_2 + \kappa^\circ_1$,
and then this induces the map
\begin{align*}
&\varprojlim_{C'_2} \varinjlim_{\kappa_1} \varprojlim_{C_0}
H^\ast((\mathcal{P}_Y^{\leq \kappa^\circ_1} + I^{\leq \kappa_1}_{C_0, C_2})
/ I^{\leq \kappa_1}_{C_0, C'_2}, d_Y) \\
&\quad
\times
\varprojlim_{C''_2} \varinjlim_{\kappa_2} \varprojlim_{C_0}
H^\ast((\mathcal{P}_Y^{\leq \kappa^\circ_2} + I^{\leq \kappa_2}_{C_0, C_2})
/ I^{\leq \kappa_1}_{C_0, C''_2}, d_Y) \\
&\to \varinjlim_{\kappa} \varprojlim_{C_0}
H^\ast(\mathcal{P}_Y^{\leq \kappa}
/ I^{\leq \kappa}_{C_0, C_2}, d_Y).
\end{align*}
Finally, this induces the map
\begin{align*}
&\varprojlim_{C_2} \varinjlim_{\kappa_1^\circ}
\varprojlim_{C'_2} \varinjlim_{\kappa_1} \varprojlim_{C_0}
H^\ast((\mathcal{P}_Y^{\leq \kappa^\circ_1} + I^{\leq \kappa_1}_{C_0, C_2})
/ I^{\leq \kappa_1}_{C_0, C'_2}, d_Y) \\
&\quad
\times \varprojlim_{C_2} \varinjlim_{\kappa_2^\circ}
\varprojlim_{C''_2} \varinjlim_{\kappa_2} \varprojlim_{C_0}
H^\ast((\mathcal{P}_Y^{\leq \kappa^\circ_2} + I^{\leq \kappa_2}_{C_0, C_2})
/ I^{\leq \kappa_1}_{C_0, C''_2}, d_Y) \\
&\to \varprojlim_{C_2} \varinjlim_{\kappa} \varprojlim_{C_0}
H^\ast(\mathcal{P}_Y^{\leq \kappa}
/ I^{\leq \kappa}_{C_0, C_2}, d_Y),
\end{align*}
which is the Poisson bracket of $H^\ast_{\mathrm{RSFT}}(Y, \lambda, \overline{K}^0_Y)$.
\end{proof}

Next we show that SFT cohomologies of cobordisms $(X, \omega)$ from $(Y^-, \xi^-)$
to $(Y^+, \xi^+)$ does not depend on the choice of the contact structure of
$(Y^\pm, \xi^\pm)$.
\begin{prop}\label{general concordance for X}
Let $(X, \omega)$ be a cobordism from $(Y^-, \lambda^-)$ to $(Y^+, \lambda^+)$,
and let $(X_1, \omega_1)$ be a (general) concordance from
$(Y^+, \lambda^+)$ to $(Y^+, \lambda_1^+)$.
($\lambda^+$ and $\lambda_1^+$ are contact forms for the same contact structure
$\xi^+$.)
Then
\[
T_{\F_X}(\cdot \Diamond \F_{X_1})
: H^\ast(\D_X, D_X) \to H^\ast(\D_{X \# X_1}, D_{X \# X_1})
\]
and
\[
T_{(\F_X)_0}(\cdot \sharp (\F_{X_1})_0)
: H^\ast(\mathcal{L}_X, d_X) \to H^\ast(\mathcal{L}_{X \# X_1}, d_{X \# X_1})
\]
are isomorphisms.
\end{prop}
\begin{proof}
We consider the case of general SFT. The case of rational SFT is similar.
By the argument similar to Lemma \ref{short concordance for X},
it is enough to prove for the case where $(X_1, \omega_1)$ is a trivial concordance.
In this case, Lemma \ref{trivial generating function} implies that
$T_{\F_X}(\cdot \Diamond \F_{X_1}) : H^\ast(\D_X, D_X) \to H^\ast(\D_X, D_X)$
is the limit of the inclusion-quotient maps similar to those given by the filtration.
Hence this limit is the identity map.
Therefore $T_{\F_X}(\cdot \Diamond \F_{X_1})$ is an isomorphism.
\end{proof}
This Proposition implies that SFT cohomologies
\[
H^\ast_{\mathrm{SFT}}(X, \omega, Y^\pm, \lambda^\pm, \overline{K}^0_X,
\overline{K}^0_{Y^\pm}, \mu^\pm)
\]
and
\[
H^\ast_{\mathrm{RSFT}}(X, \omega, Y^\pm, \lambda^\pm, \overline{K}^0_X,
\overline{K}^0_{Y^\pm}, \mu^\pm)
\]
defined by (\ref{limit general X}) and (\ref{limit rational X}) for cobordisms
$(X, \omega)$ from $(Y^-, \xi^-)$ to $(Y^+, \lambda^+)$ with different contact forms
are naturally isomorphic respectively.
We denote these isomorphic cohomology groups by
\[
H^\ast_{\mathrm{SFT}}(X, \omega, Y^\pm, \xi^\pm, \overline{K}^0_X,
\overline{K}^0_{Y^\pm}, \mu^\pm)
\]
and
\[
H^\ast_{\mathrm{RSFT}}(X, \omega, Y^\pm, \xi^\pm, \overline{K}^0_X,
\overline{K}^0_{Y^\pm}, \mu^\pm)
\]
respectively.

Finally we show that SFT cohomologies of $(Y, \xi)$ do not depend on the choice of
the contact structure $\lambda$.
\begin{prop}\label{general concordance for Y}
For any concordance $(X, \omega)$ from $(Y^-, \lambda^-)$ to $(Y^+, \lambda^+)$,
the homomorphisms
\[
i_X^\pm : H^\ast (\W_{Y^\pm} , D_{Y^\pm}) \to H^\ast(\D_X, D_X)
\]
and
\[
i_{X, 0}^\pm : H^\ast (\mathcal{P}_{Y^\pm} , d_{Y^\pm}) \to H^\ast(\mathcal{L}_X, d_X)
\]
are isomorphisms of modules, and the composition
\[
A = i_X^- \circ (i_X^+)^{-1} : H^\ast(\W_{Y^+}, D_{Y^+})
\to H^\ast(\W_{Y^-}, D_{Y^-})
\]
is an isomorphism of algebras, and
the composition
\[
A^0 = i_{X, 0}^- \circ (i_{X, 0}^+)^{-1} : H^\ast(\mathcal{P}_{Y^+}, d_{Y^+})
\to H^\ast(\mathcal{P}_{Y^-}, d_{Y^-})
\]
is an isomorphism of Poisson algebras.
Furthermore, $A$ and $A^0$ do not depend on the concordance $(X, \omega)$.
\end{prop}
\begin{proof}
The proof of the first claim is similar to that of Lemma \ref{short concordance for Y}.
First we consider the case of trivial concordance.
As in the proof of Proposition \ref{general concordance for X},
in this case, Lemma \ref{trivial generating function} implies that
$i_X^\pm : H^\ast (\W_{Y^\pm} , D_{Y^\pm}) \to H^\ast(\D_X, D_X)$ is the limit of
the inclusion-quotient maps similar to those given by the filtration.
Hence the limit is an isomorphism.
In general case, there exists a concordance $(X', \omega')$ from
$(Y^+, \lambda^+)$ to $(Y^-, \lambda^-)$ such that $X \# X'$ is a trivial concordance.
Then $T_{\F_X}(\cdot \Diamond \F_{X'}) \circ i_X^- = i_{X \# X'}^-
: H^\ast (\W_{Y^-} , D_{Y^-}) \to H^\ast(\D_{X \# X'}, D_{X \# X'})$ is an isomorphism.
Since $T_{\F_X}(\cdot \Diamond \F_{X'})$ is also an isomorphism by Proposition
\ref{general concordance for X}, this implies that $i_X^-$ is an isomorphism.
The cases of $i_X^+$ or $i_{X, 0}^\pm$ are similar.

By an argument similar to that of Proposition \ref{algebraic structures of limits},
we can prove that $A = i_X^- \circ (i_X^+)^{-1}$ is an algebra homomorphism,
and $A^0 = i_{X, 0}^- \circ (i_{X, 0}^+)^{-1}$ is an isomorphism of Poisson algebras.

The independence of $A = i_X^- \circ (i_X^+)^{-1}$ and
$A^0 = i_{X, 0}^- \circ (i_{X, 0}^+)^{-1}$ are similar to
Lemma \ref{short concordance for Y}.
\end{proof}
This proposition implies that SFT cohomologies
$H^\ast(\W_{(Y, \lambda, \overline{K}_Y^0)}, D_{(Y, \lambda, \overline{K}_Y^0)})$
and
$H^\ast(\mathcal{P}_{(Y, \lambda, \overline{K}_Y^0)}, d_{(Y, \lambda, \overline{K}_Y^0)})$
defined by (\ref{limit general Y}) and (\ref{limit rational Y}) for different strict
contact manifolds $(Y, \lambda)$ of the same contact structure $\xi$ are
naturally isomorphic respectively.
Hence we denote these cohomology groups by
$H^\ast_{\mathrm{SFT}}(Y, \xi, \overline{K}^0)$ and
$H^\ast_{\mathrm{RSFT}}(Y, \xi, \overline{K}^0)$ respectively.

It is easy to check that
$H^\ast_{\mathrm{SFT}}(X, \omega, Y^\pm, \xi^\pm, \overline{K}^0_X,
\overline{K}^0_{Y^\pm}, \mu^\pm)$ has a structure of
$H^\ast_{\mathrm{SFT}}(Y^\pm, \xi^\pm, \overline{K}_{Y^\pm}^0)$-bimodule,
and
$H^\ast_{\mathrm{RSFT}}(X, \omega, Y^\pm, \xi^\pm,\ab \overline{K}^0_X,
\overline{K}^0_{Y^\pm}, \mu^\pm)$ has a structure of
$H^\ast_{\mathrm{RSFT}}(Y^\pm, \xi^\pm, \overline{K}_{Y^\pm}^0)$-bimodule.

Finally we consider the case of contact homology.
This case is more standard.
\begin{prop}\label{short concordance for contact homology}
For any short concordance $\mathcal{C}_X$ from $\mathcal{C}_{Y^-}$
to $\mathcal{C}_{Y^+}$, the homomorphism
\[
\Psi_{(\widehat{\mathcal{F}}_X)_0} :
H^\ast (\A_{Y^+}^{\leq \kappa} / I^{\leq \kappa}_{C_0}, \partial_{Y^+}) \to
H^\ast (\A_{Y^-}^{\leq \kappa} / I^{\leq \kappa}_{C_0}, \partial_{Y^-})
\]
is an isomorphism.
Furthermore, it does not depend on the short concordance $\mathcal{C}_X$.
\end{prop}
\begin{proof}
If $\mathcal{C}_X$ is a trivial short concordance, then the claim follows from
Lemma \ref{trivial generating function}.
For a general short concordance $\mathcal{C}_X$, let $\mathcal{C}_{X'}$ be
a short concordance from $\mathcal{C}_{Y^+}$ to $\mathcal{C}_{Y^-}$.
Then since the compositions of these two short concordance are
trivial short concordances, $\Psi_{(\widehat{\mathcal{F}}_{X'})_0} \circ
\Psi_{(\widehat{\mathcal{F}}_X)_0}$ and $\Psi_{(\widehat{\mathcal{F}}_X)_0} \circ
\Psi_{(\widehat{\mathcal{F}}_{X'})_0}$ are isomorphisms.
Hence $\Psi_{(\widehat{\mathcal{F}}_X)_0}$ is also an isomorphism.
\end{proof}
Therefore we can define the limit
\[
H^\ast_{\mathrm{CH}}(Y, \lambda, \overline{K}^0)
= \varinjlim_{\kappa} \varprojlim_{C_0}
H^\ast(\A_{(Y, \lambda, K_Y, \overline{K}_Y^0)}^{\leq \kappa}
/ I^{\leq \kappa}_{C_0}, \partial_{(Y, \lambda, K_Y, K_Y^0, K_Y^2, J, \B)}).
\]
We sometimes abbreviate this limit as $H^\ast(\A_Y, \partial_Y)$.
For any exact cobordism $(X, \omega)$ from $(Y^-, \lambda^-)$ and
$(Y^+, \lambda^+)$, we can define
\[
\Psi_X : H^\ast_{\mathrm{CH}}(Y^+, \lambda^+, \overline{K}^0_{Y^+})
\to H^\ast_{\mathrm{CH}}(Y^-, \lambda^-, \overline{K}^0_{Y^-})
\]
by the limit of $\Psi_{(\widehat{\mathcal{F}}_X)_0}$.
We can easily prove the following.
\begin{prop}
For any concordance $(X, \omega)$ from $(Y^-, \lambda^-)$ to $(Y^+, \lambda^+)$,
the homomorphism
\[
\Psi_X : H^\ast_{\mathrm{CH}}(Y^+, \lambda^+, \overline{K}^0_{Y^+})
\to H^\ast_{\mathrm{CH}}(Y^-, \lambda^-, \overline{K}^0_{Y^-})
\]
is an isomorphism of algebras.
Furthermore, it does not depend on the concordance $(X, \omega)$.
\end{prop}
We denote the isomorphism class of contact homology by
$H^\ast_{\mathrm{CH}}(Y, \xi, \overline{K}^0)$.


%% file: SFT-14_Circle_action.tex
%
%

\section{SFT of a contact manifold with the $S^1$-action induced by the Reeb flow}
\label{S^1 action}
The arguments in \cite{EGH00} or \cite{Bo02} are easily adapted to our construction of
SFT.
In this section, we demonstrate how to calculate the SFT cohomology of pre-quantization
spaces, or more generally, contact manifolds with the locally free $S^1$-action
generated by the Reeb vector field.

Let $(Y, \lambda)$ be a closed contact manifold and assume that there exists a constant
$L > 0$ such that $\varphi^\lambda_L = \id$.
Then $S^1 = \R / L\Z$ acts on $Y$ by $t \cdot y = \varphi^\lambda_t(y)$.
We consider the SFT of such a contact manifold.
We may assume $L = 1$.
First we consider the case where every cycle in $K^0$ is invariant by this action.
In this case, we can calculated the SFT cohomology by the following proposition.
\begin{prop}\label{H = 0}
All periodic orbits are good, and the local systems $\S^D$ and $\S^{\overline{P}}$ are
trivial on $\overline{P}$.
Furthermore, we can construct the virtual fundamental chains which make
$\mathcal{H} = 0$.
\end{prop}
Theorem \ref{H vanishes} is a corollary of this proposition.

Let $J$ be an $S^1$-invariant $d\lambda$-compatible complex structure of
$\xi = \Ker \lambda$.
First we prove the claim about the local systems $\S^D$ and $\S^{\overline{P}}$.
For each $l \geq 1$, let $Y^{l^{-1} \Z / \Z} = \ev_0 P_{l^{-1}} \subset Y$ be the
fixed manifold of the subgroup $l^{-1} \Z / \Z \subset S^1$.
Then $l^{-1} \Z / \Z$ acts on each fiber of $\xi|_{Y^{l^{-1} \Z / \Z}}$.
Since this is a unitary action, we can decompose this complex vector bundle by
the eigenvalues:
\[
\xi|_{Y^{l^{-1} \Z / \Z}} = W_0 \oplus W_1 \oplus \dots \oplus W_{l-1},
\]
where $(\varphi^\lambda_{l^{-1}})_\ast$ acts on each $W_k$ by $e^{2\pi \sqrt{-1} k / l}$.
Then for each point $y \in Y^{l^{-1} \Z / \Z}$, we can define a unitary trivialization of
$\xi$ on the periodic orbit $\gamma(t) = \varphi^\lambda_{l^{-1}t}(y)$ by
\[
(\varphi^\lambda_{l^{-1}t})_\ast \circ
\bigl(\bigoplus_k e^{-2\pi \sqrt{-1} kt / l} 1_{W_k}\bigr):\\
\xi_{\gamma(0)} = W_0 \oplus W_1 \oplus \dots
\oplus W_{l-1} \tocong \xi_{\gamma(t)}
\]
if we fix a unitary basis of each $W_k$.
Under this trivialization, $(\varphi^\lambda_{l^{-1}t})_\ast$ are given by the diagonal
matrices
\[
\bigoplus_k e^{2\pi \sqrt{-1} kt / l} 1_{W_k}.
\]
Hence the linear operator $\mathring{D}^+_\gamma$ is complex linear.
In particular, its kernel has the complex orientation.
Therefore $\S^D$ is a trivial local system on $P_{l^{-1}}$.
Similarly, $\S^D$ is trivial on $P_{k / l}$ for each $k / l$
since $\ev_0 P_{k / l} = \ev_0 P_{1 / l}$ if $k$ and $l$ are coprime.
Hence there are no bad orbits and the induced local system on
$\overline{P}$ is also trivial.
Similarly, $\overline{P}$ does not contain any non-orientable points,
and $\S^{\overline{P}}$ is trivial on $\overline{P}$.

Next we construct required virtual fundamental chains.
$\widehat{\M} = \widehat{\M}(Y, \lambda, J)$ has a locally free $S^1$-action
defined by $t \cdot (\Sigma, z, u)
= (\Sigma, z, (1 \times \varphi^\lambda_t) \circ u)$.
We will construct a pre-Kuranishi structure of the quotient space
$\widehat{\M}^\bullet = \widehat{\M} / S^1$ which induces a pre-Kuranishi structure of
$\widehat{\M}$.
Since the evaluation maps to $\overline{P}$ or $Y / S^1$ are well-defined on
$\widehat{\M}^\bullet$, we can define its fiber products
$((\widehat{\M}^\bullet)^\diamond, \mathring{K}^2, K, K^0 / S^1)$
and multi-valued partial submersions $\Xi^\circ$ and $\Lambda$ similarly.
We can construct a grouped multisection of
$((\widehat{\M}^\bullet)^\diamond, \mathring{K}^2, K, K^0 / S^1)$ which satisfies
the similar compatibility conditions.
Then we define the grouped multisection of
$(\widehat{\M}^\diamond, \mathring{K}^2, K,K^0)$ by its pull back.
Since the $S^1$-action is locally free, it makes the virtual fundamental chains of
the zero-dimensional fiber products used for the definition of $\mathcal{H}$ vanish.
Therefore, it is enough to define a required pre-Kuranishi structure of
$\widehat{\M}^\bullet$.

First we explain the construction of a Kuranishi neighborhood of a point $\bar p_0
= (\Sigma_0, z, u_0) \in \widehat{\M}^\bullet$.
Define a finite group $G^+_0$ by
\[
G^+_0 = \{(g, t) \in \Aut (\Sigma_0) \times S^1; g (\{z_i\}) = \{z_i\},
u_0 \circ g = (1 \times \varphi^\lambda_t) \circ u_0\}.
\]
We also define a group $G^{++}_0 \subset \Aut (\Sigma_0) \times S^1$ by
$G^{++}_0 = G^+_0 \cdot S^1$.
We assume that the following data
$(\bar p_0^+, S, (\mathcal{O}_a, \mathcal{N}_a, E^0_a, \lambda_a)_{a \in A})$
are given instead of the data in the usual case.

$\bar p_0^+ = (\Sigma_0, z, z^+, u_0) \in \widehat{\M}^\bullet$ is, as in the usual case,
a curve obtained by adding marked points on the nontrivial components of $\Sigma_0$.
We assume that all unstable components of $(\Sigma_0, z, z^+)$ are trivial cylinders
of $\bar p_0$ and $G^+_0$ preserves $z^+ = \{z^+_i\}$ as a set.

$S \subset Y$ is a finite union of $S^1$-invariant codimension-two submanifolds
such that $\pi_Y \circ u_0$ intersects with $S$ at $z^+$ transversely.
We can take such an $S^1$-invariant submanifold for the following reason.
Choosing appropriate additional marked points $z^+$, we assume that
the differential $d^\xi u_0$ does not vanish at $z^+$.
Let $l^{-1} \Z / \Z$ be the stabilizer of the point $y = u_0(z_i^+)$.
Then an $S^1$-equivariant tubular neighborhood of the orbit $S^1 \cdot y$ is
isomorphic to $\R / \Z \times_{l^{-1} \Z / \Z} \xi_y$.
Since the $l^{-1} \Z / \Z$-action on $\xi_y$ is unitary and commutative,
$\xi_y$ can be decomposed into irreducible representations of complex dimension one.
Therefore there exists an $l^{-1} \Z / \Z$-invariant subspace $\xi^0_y \subset \xi_y$ of
complex codimension one such that $\Image d^\xi u_0 (z_i^+) \pitchfork \xi^0_y$.
Then $\pi_Y \circ u_0$ intersects with the $S^1$-invariant submanifold
$S = \R / \Z \times_{l^{-1} \Z / \Z} \xi^0_y$ transversely at $z_i^+$.

For the construction of the global pre-Kuranishi structure, we used an infinite family of
disjoint submanifolds $\{S^x\}_{x \in \R^2}$.
(See the proof of Lemma \ref{existence of a domain curve representation}.)
It was constructed as constant sections of the trivial tubular neighborhood of $S$.
To construct such a family of $S^1$-invariant submanifolds,
it is enough to make the $l^{-1} \Z / \Z$-action on $\xi_y / \xi^0_y \cong
\Image d^\xi u_0 (z_i^+)$ trivial.
In particular, it is enough to choose $z^+$ so that the stabilizer $l^{-1} \Z / \Z$ of
each $\pi_Y \circ u_0(z^+_i)$ is locally minimal in the image of $\pi_Y \circ u_0$.

Let $(\hat P \to \hat X, Z, Z^+, Z_{\pm\infty})$ be the local universal family of
the stabilization $(\hat \Sigma, z, z^+, \pm\infty)$ of the blow down curve of
$(\Sigma, z, z^+)$
We need an additional vector space $E^0$ and a linear map $\lambda$.
If we can take a $G^+_0$-equivariant linear map
$\lambda : E^0 \to C^\infty(\hat P \times Y, \Wedge^{0, 1} V^\ast \hat P
\otimes (\R \partial_\sigma \oplus TY))$
which is $S^1$-invariant, that is,
\[
\lambda(h)(z, \varphi^\lambda_t(y)) = (1 \otimes (\varphi^\lambda_t)_\ast)
\lambda(h)(z, y)
\]
for all $t \in S^1$, and which makes the linear operator $D_{p_0}^+$ defined in
Section \ref{construction of nbds} surjective, then it is easy to construct
a Kuranishi neighborhood of $\bar p_0 \in \widehat{\M}^\bullet$
which is independent of the choice of the representative $p_0 \in \widehat{\M}$.
However, since the $S^1$-action on $Y$ is not necessarily free, we cannot construct
such a $G^+_0$-equivariant linear map in general.
Instead, we take the following data
$(\mathcal{O}_a, \mathcal{N}_a, E^0_a, \lambda_a, I_a)_{a \in A}$:
\begin{itemize}
\item
$A$ is a finite index set.
\item
For each $a \in A$,
$\mathcal{O}_a \subset Y$ is an $S^1$-orbit,
$\mathcal{N}_a \subset Y$ is its $S^1$-invariant tubular neighborhood,
and $\pi_{\mathcal{N}_a} : \mathcal{N}_a \to \mathcal{O}_a$ is its $S^1$-equivariant
projection.
\item
Let $\pi_{\widetilde{\mathcal{O}}_a} : \widetilde{\mathcal{O}}_a \to \mathcal{O}_a$ be
the covering space of $\mathcal{O}_a$ such that the $S^1$-action lifts to
$\widetilde{\mathcal{O}}_a$ as a free (and transitive) action.
Then $\pi_a : E^0_a \to \widetilde{\mathcal{O}}_a$ is a finite dimensional
$G^{++}_0$-vector bundle.
(The action of $G^{++}_0 \subset \Aut (\Sigma_0) \times S^1$ on
$\widetilde{\mathcal{O}}_a$ is defined by the projection $G^{++}_0 \to S^1$.)
\item
Define $\widetilde{\mathcal{N}}_a = \mathcal{N}_a \times_{\mathcal{O}_a}
\widetilde{\mathcal{O}}_a$ and let
$\widetilde{\pi}_{\mathcal{N}_a} : \widetilde{\mathcal{N}}_a \to
\widetilde{\mathcal{O}}_a$ be the projection.
Define $\pi_{\hat P \times \widetilde{\mathcal{N}}_a} : \hat P \times
\widetilde{\mathcal{N}}_a \to \widetilde{\mathcal{O}}_a$ by
$\pi_{\hat P \times \widetilde{\mathcal{N}}_a}(z, y)
= \widetilde{\pi}_{\mathcal{N}_a}(y)$.
Then
$\lambda_a : \pi_{\hat P \times \widetilde{\mathcal{N}}_a}^\ast E^0_a \to
\Wedge^{0,1} V^\ast \hat P \otimes_{\C} (\R \partial_\sigma \oplus TY)
|_{\hat P \times \widetilde{\mathcal{N}}_a}$
is a $G^{++}_0$-equivariant bundle map with compact support
$\supp \lambda_a \subset \hat P \times \widetilde{\mathcal{N}}_a$.
\item
$I_a \subset S^1$ is a union of finite number of intervals
which is invariant by the $G^+_0$-action.
\end{itemize}
We impose the following conditions on them:
\begin{enumerate}[label=(\arabic*)]
\item
The projection of $\supp \lambda_a \subset \hat P \times \widetilde{\mathcal{N}}_a$ to
$\hat P$ does not intersect with the nodal points of $\hat P$ or $Z_{\pm\infty}$.
\item
\label{condition of I}
There exists a simply connected neighborhood $\mathcal{I}_a \subset S^1$ of $0$ and
a finite subgroup $\Gamma \subset S^1$ such that
$I_a = \mathcal{I}_a + \Gamma$,
$\mathcal{I}_a = - \mathcal{I}_a$ and
$(\mathcal{I}_a + \mathcal{I}_a) \cap \Gamma = \{0\}$.
(Namely, the intervals in $I_a$ have the same length, and
the intervals in the complement $S^1 \setminus I_a$ also have the same length.
Furthermore, the former is smaller than the latter.)
\item
\label{support lambda small}
Let $p_0 = (\Sigma_0, z, u_0) \in \widehat{\M}$ be a representative of
$\bar p_0 \in \widehat{\M}^\bullet$.
Then there exists a point $x_a \in \widetilde{\mathcal{O}}_a$ such that
\[
\supp \lambda_a \cap (1 \times \pi_{\widetilde{\mathcal{N}}_a})^{-1}
\graph (\pi_Y \circ u_0) \subset \pi_{\hat P \times \widetilde{\mathcal{N}}_a}^{-1}
(I_a \cdot x_a).
\]
\item
\label{transverse condition of lambda}
Let $E^0_{a, x_a}$ be the vector space of locally $S^1$-invariant sections of
$E^0_a|_{I_a \cdot x_a}$.
(A locally $S^1$-invariant section is a section which is $S^1$-invariant on each connected
component of $I_a \cdot x_a$.
Namely, if we trivialize $E^0_a|_{I_a \cdot x_a}$ by the $S^1$-action,
then it is a locally constant section.)
Note that the $G^+_0$-action on $E^0_a$ induces a $G^+_0$-action on $E^0_{a, x_a}$.
Define $\widetilde{\mathcal{N}}_{a, x_a} = \widetilde{\pi}_{\mathcal{N}_a}^{-1}(I_a \cdot x_a)
\subset \widetilde{\mathcal{N}}_a$ and $\mathcal{N}_{a, x_a}
= \pi_{\widetilde{\mathcal{N}}_a}(\widetilde{\mathcal{N}}_{a, x_a}) \subset \mathcal{N}_a$,
where $\pi_{\widetilde{\mathcal{N}}_a} : \widetilde{\mathcal{N}}_a \to \mathcal{N}_a$ is
the projection.
Define a $G^+_0$-equivariant linear map
\[
\lambda_{a, x_a} : E^0_{a, x_a} \to C^\infty(\hat P \times \mathcal{N}_{a, x_a},
\Wedge^{0,1} V^\ast \hat P \otimes_{\C} (\R \partial_\sigma \oplus TY)).
\]
by
\[
\lambda_{a, x_a}(h)(z, y) = \sum_{\tilde y \in \widetilde{\mathcal{N}}_a,
\pi_{\widetilde{\mathcal{N}}_a}(\tilde y) = y} \lambda_a(h(z, \tilde y)).
\]
Let $E^0$ and $\lambda$ be the direct sums of $E^0_{a, x_a}$ and $\lambda_{a, x_a}$
over $a \in A$ respectively.
Then the liner map
\begin{align*}
&D_{p_0}^+ : \widetilde{W}_\delta^{1, p}(\Sigma_0, u_0^\ast T \hat Y) \oplus E^0\\
&\to
L_\delta^p(\Sigma_0, \Wedge^{0, 1}T^\ast \Sigma_0 \otimes u_0^\ast T \hat Y)
\oplus
\bigoplus_{\text{limit circles}} \Ker A_{\gamma_{\pm\infty_i}} / (\R \partial_\sigma
\oplus \R R_\lambda)\\
&\quad \oplus \bigoplus_{z_i} T_{\pi_Y \circ u_0(z_i)} Y\\
&(\xi, h) \mapsto (D_{p_0} \xi(z) + \lambda(h)(z, \pi_Y \circ u_0(z)),
\sum_j \langle\xi|_{S^1_{\pm\infty_i}}, \eta_j^{\pm\infty_i}\rangle
\eta_j^{\pm\infty_i}, \pi_Y \circ \xi (z_i))
\end{align*}
is surjective,
where $D_{p_0}$ is the linearization of the equation of $J$-holomorphic maps,
and $\{\eta_j^{\pm\infty_i}\}_j$ is an orthonormal basis of
the orthogonal complement of $\R \partial_\sigma \oplus \R R_\lambda$ in
$\Ker A_{\gamma_{\pm\infty_i}}$ for each $\pm\infty_i$.
\end{enumerate}

We can construct such data $(\mathcal{O}_a, \mathcal{N}_a, E^0_a, \lambda_a)_{a \in A}$
as follows.
First we explain the construction of $E^0_a$
for each $S^1$-orbit $\mathcal{O}_a \subset Y$.
Define a map $\pi_{\widetilde{\mathcal{O}}^+_a} : \widetilde{\mathcal{O}}^+_a
= G^{++}_0 \times_{S^1} \widetilde{\mathcal{O}}_a \to \widetilde{\mathcal{O}}_a$ by
$\pi_{\widetilde{\mathcal{O}}^+_a}(g, t, x) = t \cdot x$.
Let $\widehat{E}^0_a \to \widetilde{\mathcal{O}}^+_a$ be the pull back of
$(\R \partial_\sigma \oplus TY)|_{\mathcal{O}_a}$ by $\pi_{\widetilde{\mathcal{O}}_a} \circ
\pi_{\widetilde{\mathcal{O}}^+_a} : \widetilde{\mathcal{O}}^+_a \to \mathcal{O}_a$,
and define a $G^{++}_0$-vector bundle $\pi_a : E^0_a \to \widetilde{\mathcal{O}}_a$ by
$E^0_a|_x = \bigoplus_{y \in \pi_{\widetilde{\mathcal{O}}^+_a}^{-1}(x)} \widehat{E}^0_a|_y$.
Define maps $\pi_{\widetilde{\mathcal{N}}^+_a} : \widetilde{\mathcal{N}}^+_a
= G^{++}_0 \times_{S^1} \widetilde{\mathcal{N}}_a \to \widetilde{\mathcal{N}}_a$
and $\pi_{\hat P \times \widetilde{\mathcal{N}}^+_a} : \hat P \times
\widetilde{\mathcal{N}}^+_a \to \widetilde{\mathcal{O}}^+_a$ by
$\pi_{\widetilde{\mathcal{N}}^+_a}(g, t, x) = t \cdot x$ and
$\pi_{\hat P \times \widetilde{\mathcal{N}}^+_a}(z, g, t, x)
= (g, t, \widetilde{\pi}_{\mathcal{N}_a}(x))$ respectively.
We note that
\[
(\pi_{\hat P \times \widetilde{\mathcal{N}}_a}^\ast E^0_a)|_\gamma
= \bigoplus_{\delta \in (1 \times \pi_{\widetilde{\mathcal{N}}^+_a})^{-1}(\gamma)}
(\pi_{\hat P \times \widetilde{\mathcal{N}}^+_a}^\ast \widehat{E}^0_a)|_\delta.
\]
We construct $\lambda_a$ as follows.
Take a $G^{++}_0$-invariant section $\rho_a$ of the pull back of
$\Wedge^{0,1} V^\ast \hat P$ to $\hat P \times \widetilde{\mathcal{N}}^+_a$
such that the projection of its support to $\hat P$ is contained in a small neighborhood of
some $G^+_0$-orbit.
Since $\pi_{\hat P \times \widetilde{\mathcal{N}}^+_a}^\ast \widehat{E}^0_a$ is a pull back
of $(\R \partial_\sigma \oplus TY)|_{\mathcal{O}_a}$,
$\rho_a$ defines a linear map $\pi_{\hat P \times \widetilde{\mathcal{N}}^+_a}^\ast
\widehat{E}^0_a \to \Wedge^{0,1} V^\ast \hat P \otimes_{\C}
(\R \partial_\sigma \oplus TY)$, which defines the $G^{++}_0$-linear map
$\lambda_a : \pi_{\hat P \times \widetilde{\mathcal{N}}_a}^\ast E^0_a \to
\Wedge^{0,1} V^\ast \hat P \otimes_{\C} (\R \partial_\sigma \oplus TY)$.
If the support of $\rho_a$ is sufficiently small, then there exists a union  of intervals
$I_a \subset S^1$ which satisfies Condition \ref{condition of I} and
\ref{support lambda small}.
Since the $G^{++}_0$-action on $\hat P \times \widetilde{\mathcal{N}}^+_a$ is free,
if we choose appropriate $\mathcal{O}_a$ and $\rho_a$
($a \in A$), then Condition \ref{transverse condition of lambda} also holds true.

Using the above data, we construct the Kuranishi neighborhood of
$\bar p_0 \in \widehat{\M}^\bullet$ as follows.
As in the usual case, we fix a temporally data $(z^{++}, S', \hat R_i)$,
where in this case, we assume that they are $G^+_0$-invariant.
In addition, we take a $G^+_0$-invariant family of sections
$\hat R_{S^1} = (\hat R_{S^1, l})$ of $\hat P \to \hat X$ and
a codimension-one submanifold $S_{S^1} \subset Y$ transverse to the Reeb vector field
such that $\pi_Y \circ u_0(\widetilde{R}_{S^1, l}(0)) \in S_{S^1}$ for all $l$,
where $\widetilde{R}_{S^1, l}$ is the section of $\widetilde{P} \to \widetilde{X}$
induced by $\hat R_{S^1, l}$.
Define a function $p_{S_{S^1}}$ on a small neighborhood of $S_{S^1}$
by $y \in \varphi^\lambda_{p_{S_{S^1}}(y)}(S_{S^1})$ and $|p_{S_{S^1}}(y)| \ll 1$.
These data are used to kill the $S^1$-action.

As in the usual case, we define a smooth manifold $\hat V = X \times B_\epsilon(0)$
and define a smooth map $s^0 : \hat V \to \R^k \oplus \bigoplus_{z_\beta^{++}} \R^2$.
In addition, we define a smooth map $s^1 : \hat V \to \R$ by
\[
s^1(a, b, x) = \frac{1}{m_i} \sum_{l=1}^{m_i} p_{S_{S^1}} \circ \pi_Y \circ
\Phi_{a, b}(\xi_x)(\widetilde{R}_{S~1, l}(a)),
\]
and define $\mathring{V} = \{s^0 = 0, s^1 = 0\} \subset \hat V$.
As in the usual case, we define a smooth map
$s : \mathring{V} \to E := E^0_{p_0} \oplus \bigoplus_{z_\alpha^+} \R^2$.
It is easy to see that the natural map
$\bar \psi : s^{-1}(0) / G^+_0 \to \widehat{\M}^\bullet$
is a homeomorphism onto a neighborhood of $\bar p_0 \in \widehat{\M}^\bullet$.
Hence $(\mathring{V}, E, s, \bar \psi, G^+_0)$
define the Kuranishi neighborhood of $\bar p_0 \in \widehat{\M}^\bullet$.
We note that the Kuranishi neighborhood is independent of the choice of
$x_a \in \widetilde{\mathcal{O}}_a$ because of Condition \ref{condition of I}.

The definition of the embedding of Kuranishi neighborhoods are similar to the usual
one explained in Section \ref{embed}.
A global pre-Kuranishi structure of $\widehat{\M}^\bullet$ is defined similarly,
and it induces a pre-Kuranishi structure of $\widehat{\M}$.
As we noted, we can define the fiber products
$((\widehat{\M}^\bullet)^\diamond, \mathring{K}^2, K, K^0 / S^1)$
and its multi-valued partial submersions $\Xi^\circ$ and $\Lambda$ similarly.
We can construct its grouped multisection satisfying the conditions similar
to those for $(\widehat{\M}^\diamond, \mathring{K}^2, K,K^0 / S^1)$.
Then its pull back defines the grouped multisection of
$(\widehat{\M}^\diamond, \mathring{K}^2, K,K^0)$.
As we have explained, the virtual fundamental chains defined by these grouped
multisections are the required ones.
Therefore Proposition \ref{H = 0} holds true.

Next we consider the case where $K^0$ contains cycles which are not invariant by
the $S^1$-action.
We assume that the $S^1$-action is free, that is, we only consider the case of
a pre-quantization space of some closed symplectic manifold.
We show that some terms of the rational part $\mathcal{H}_0$ of the generating function
are calculated by the Gromov-Witten invariants of the closed symplectic manifold.
The following argument is an adaptation of that given in \cite{EGH00} and \cite{Bo02}.

Let $(M, \omega)$ be a closed symplectic manifold of dimension $2(n-1)$ with an
integral cohomology class $[\omega] \in H^\ast (M; \Z)$.
Let $\pi_M : Y \to M$ be a principal $U(1)$-bundle with first Chern class
$c_1(Y) = [\omega]$, and $\alpha$ be a connection form such that $\pi_M^\ast \omega
= -\frac{1}{2\pi\sqrt{-1}} d\alpha$.
Then $\lambda = -\frac{1}{2\pi\sqrt{-1}} \alpha$ is a contact form of $Y$ such that
$d\lambda = \pi_M^\ast \omega$.
Note that the Reeb flow of the pre-quantization space $(Y, \lambda)$ is opposite
to the usual $U(1)$-action on $Y$.
Since $\overline{P} = \bigcup_{k = 1}^\infty \overline{P}_k$ and $\overline{P}_k \cong M$,
a smooth triangulation of $M$ defines a triangulation $K$ of $\overline{P}$.

Let $J$ be an $\omega$-compatible almost complex structure on $M$.
It induces a complex structure of $\xi = \Ker \lambda \cong \pi_M^\ast TM$, which
we also denote by $J$.
Then as an almost complex manifold,
$\hat Y = \R \times Y$ is isomorphic to
\[
Y \underset{U(1)}{\times} (\C \setminus 0)
= Y \underset{U(1)}{\times} (\C P^1 \setminus \{0, \infty\})
\]
by $(\sigma, y) \mapsto [y, e^{-2\pi \sigma}]$,
where the almost complex structure of $\mathcal{L} = Y \times_{U(1)} \C$ is defined by
$T_{[y, z]} \mathcal{L} \cong \xi_y \oplus T_z \C$.
The almost complex structure of $Y \times_{U(1)} \C P^1$ is similar.

Holomorphic buildings for $(Y, \lambda, J)$ and stable maps in $(M, J)$ are related
as follows.
For a holomorphic building $(\Sigma, z, u) \in \widehat{\M}(Y, \lambda, J)$ of height $k$,
a $J$-holomorphic map
\[
\check u : \check \Sigma \to Y \underset{U(1)}{\times} (
\underbrace{\C P^1 \lrsubscripts{\cup}{0}{\infty} \C P^1 \lrsubscripts{\cup}{0}{\infty}
\dots \lrsubscripts{\cup}{0}{\infty} \C P^1}_k)
\]
is defined by $\check u|_{\Sigma \setminus \coprod S^1}
= u|_{\Sigma \setminus \coprod S^1}$
(and removal of singularity), where
$(\check \Sigma, z, \pm\infty)$ is the blow down curve of $(\Sigma, z)$,
$\coprod S^1 \subset \Sigma$ is the union of imaginary cirlces in $\Sigma$
and we regard $\Sigma \setminus \coprod S^1$ as a subset of $\check \Sigma$.

Let $\hat \pi_M : Y \times_{U(1)} (
\C P^1 \lrsubscripts{\cup}{0}{\infty}
\dots \lrsubscripts{\cup}{0}{\infty} \C P^1) \to M$ be the projection.
Then $\bar u = \hat \pi_M \circ \check u : \check \Sigma \to M$ is a
$J$-holomorphic map, and the restriction of $\check u$ to the $i$-th floor component
$\check \Sigma_i \subset \check \Sigma$ can be regarded as a meromorphic section of
$\bar u^\ast \mathcal{L}$ on $\check \Sigma_i$.
Then each zero of $\check u$ with degree $k$ corresponds to a $+\infty$-limit circle of
$\Sigma$, and the asymptotic periodic orbit of $u$ on this circle has multiplicity $k$.
Similarly, each pole of $\check u$ with degree $k$ corresponds to a $-\infty$-limit circle
of $\Sigma$, and the asymptotic periodic orbit of $u$ on this circle has multiplicity $k$.

Let $(\hat \Sigma, z, \pm\infty)$ be the curve obtained by collapsing the irreducible
component of $(\check \Sigma, z, \pm\infty)$ corresponding to the trivial cylinders of
$(\Sigma, z, u)$.
Note that $\bar u$ induces a stable map $(\hat \Sigma, z \cup \{\pm\infty_i\}, \hat u)$
of $(M, J)$ since $\bar u$ is constant on each irreducible component of $\check \Sigma$
corresponding to a trivial cylinder of $(\Sigma, z, u)$.
Then the $E_{\hat \omega}$-energy of $(\Sigma, z, u)$ is
\begin{equation}
E_{\hat \omega}(u) =  E(\hat u) := \int_{\hat \Sigma} \hat u^\ast \omega
= \sum_{+\infty_i} k_{\gamma_{+\infty_i}} - \sum_{-\infty_i} k_{\gamma_{-\infty_i}},
\label{energy of stable map}
\end{equation}
where $k_{\gamma_{\pm\infty_i}}$ is the multiplicity of $\gamma_{\pm\infty_i}$,
which is equivalent to the degree of the corresponding zero or pole of $\hat u$.

Conversely, let $(\check \Sigma, z \cup \{\pm\infty_i\})$ be a semistable curve
of genus $g = 0$ with a floor structure and
$\bar u : \check \Sigma \to M$ be a $J$-holomorphic map.
We assume that an integer $k_{\pm\infty_i} \geq 1$ is attached to
each marked point $\pm\infty_i$, and an integer $k_\mu \geq 1$ to
each nodal point $p_\mu$ which joints two components with different floors.
We assume that these integers satisfy the energy condition
for each component of $\check \Sigma$.
Namely, we assume that the sum of $k_{+\infty_i}$ and $k_\mu$ corresponding to
the zeros on the component is larger than the sum of $k_{-\infty}$ and $k_\mu$
corresponding to the poles on the component.
Then there exists a $J$-holomorphic map
\[
\check u : \check \Sigma \to Y \underset{U(1)}{\times} (
\C P^1 \lrsubscripts{\cup}{0}{\infty} \C P^1 \lrsubscripts{\cup}{0}{\infty}
\dots \lrsubscripts{\cup}{0}{\infty} \C P^1)
\]
which is obtained by patching meromorphic sections of $\bar u^\ast \mathcal{L}$ on
$\check \Sigma_i$
such that each $+\infty_i$ is a zero of degree $k_{+\infty_i}$,
each $-\infty_i$ is a pole of degree $k_{-\infty_i}$, and each nodal point $p_\mu$ is
a pole on the component of higher floor and a zero on the component of the lower floor
of degree $k_\mu$.
Furthermore, $\check u$ is unique modulo $\C^\ast$-valued holomorphic functions
on $\coprod_i \check \Sigma_i$.
(The uniqueness is true for $g \geq 1$ but the existence is not always true for
$g \geq 1$.)
Let $(\Sigma, z)$ be the curve obtained by the oriented blow up of $(\check \Sigma, z)$
at $\pm\infty_i$ and $p_\mu$ with appropriate $\varphi_\mu \in S^1$.
Then $\check u$ defines a holomorphic building $(\Sigma, z, u) \in \widehat{\M}$.
(There are $k_\mu$ choices of $\varphi_\mu \in S^1$ for each $\mu$.)

Assume that all cycles in $K^0$ except one cycle $y$ are $S^1$-invariant.
We show that if we use an appropriate virtual fundamental cycles then
$\mathcal{H}_0 \in \W_Y|_{g = 0} / (t_y^2)$ is calculated by
the rational Gromov-Witten invariants of $(M, \omega)$.

First we recall the definition of Gromov-Witten invariants.
Since the pre-Kuranishi spaces used for its definition do not have boundary of codimension
one, usually we do not need any compatibility conditions of the virtual fundamental
chains of them for construction.
However, in order to use the induced grouped multisection of the pre-Kuranishi spaces
for the definition of SFT cohomology of $(Y, \lambda)$,
we need some compatibility conditions.

We need the following data $(\hat p_0^+, \hat S, \hat E^0, \hat \lambda)$
to define a Kuranishi neighborhood of $\hat p_0 = (\hat \Sigma_0, z, \hat u_0)
\in \overline{\M}(M, J)$:
\begin{itemize}
\item
$\hat p_0^+ = (\hat \Sigma_0, z \cup z^+, \hat u_0)$ is a curve obtained
by adding marked points to make $(\hat \Sigma_0, z \cup z^+)$ stable.
We assume that $G_0 = \Aut(\hat \Sigma_0, z, \hat u_0)$ preserves $z^+$ as a set.
\item
$\hat S \subset M$ is a finite union of codimension-two submanifolds such that
$u_0$ intersects with $\hat S$ at $z^+$ transversely.
\item
Let $(\hat P \to \hat X, Z \cup Z^+)$ be the local universal family of
$(\hat \Sigma_0, z \cup z^+)$.
Then $\hat E^0$ is a finite dimensional $G_0$-vector space
and $\hat \lambda : \hat E^0 \to C^\infty(\hat P \times M;
\Wedge^{0, 1} V^\ast \hat P \otimes TM)$
is a $G_0$-equivariant linear map which satisfies following conditions:
\begin{itemize}
\item
For each $h \in \hat E^0$, the projection of the support of $\hat \lambda(h)$ to
$\hat P$ does not intersect with the nodal points or marked points $Z$.
(It may intersect with $Z^+$.)
\item
The linear map
\begin{align*}
\widehat{D}_{\hat p_0}^+ : \widetilde{W}^{1, p}(\hat \Sigma_0, \hat u_0^\ast TM)
\oplus \hat E^0
&\to L_\delta^p(\hat \Sigma_0, \Wedge^{0, 1}T^\ast \hat \Sigma_0 \otimes
\hat u_0^\ast T M)\\
&\quad \oplus \bigoplus_{z_i} T_{\hat u_0(z_i)} M\\
(\xi, h) &\mapsto (\widehat{D}_{\hat p_0} \xi + \hat \lambda(h), \xi (z_i))
\end{align*}
is surjective,
where $\widehat{D}_{\hat p_0}$ is a linearization of the equation of
the $J$-holomorphic maps, that is,
\[
\widehat{D}_{\hat p_0} \xi = \nabla \xi + J(\hat u_0) \nabla \xi j
+ \nabla_\xi J(u_0) d\hat u_0 j.
\]
\end{itemize}
\end{itemize}

Using the above data, we can construct a Kuranishi neighborhood of
$\hat p_0 \in \overline{\M}(M, J)$ similarly.
A global pre-Kuranishi structure of $\overline{\M}(M, J)$ is also constructed similarly.
Define its fiber products $\overline{\M}(M, J)^m
_{(\epsilon^{i, j}_l, c^i_l, \mathring{x}^i_l, \eta^i_l)}$
for all sequences $(\epsilon^{i, j}_l, c^i_l, \mathring{x}^i_l, \eta^i_l)$ consisting of
$\epsilon^{i, j}_l \in K^2$, $c^i_l \in K$,
$\mathring{x}^i_l \in \{ x / S^1; x \in K^0 \setminus \{y\}\} \cup \{\pi_M(y)\}$ and
$\eta^i_l \in K$,
where we regard each $x / S^1$ ($x \in K^0 \setminus \{y\}\}$)
as a cycle of dimension $\dim x - 1$ defined by the map $x / S^1 \to \pi_M(x)$,
and $\pi_M(y)$ as a cycle of dimension $\dim y$ defined by the map $y \to \pi_M(y)$.
We also define the fiber product $(\overline{\M}(M, J)^\diamond, \mathring{K}^2, K,
K^0 / S^1)$
similarly to $(\widehat{\M}^\diamond, \mathring{K}^2, K,K^0)$, and construct its grouped
multisection satisfying the compatibility conditions.
Then using the induced grouped multisection of the fiber products $\overline{\M}(M, J)^m
_{(\epsilon^{i, j}_l, c^i_l, \mathring{x}^i_l, \eta^i_l)}$, we can construct
the Gromov-Witten invariant of $(M, \omega)$.

Now we consider the pre-Kuranishi structure of $\widehat{\M}_{g = 0}(Y, \lambda, J)$.
The Kuranishi neighborhood of a point $p_0 \in \widehat{\M}_{g = 0}(Y, \lambda, J)$
is defined by the data $(p_0^+, S, E^0, \lambda)$ obtained from the data
$(\hat p_0^+, \hat S, \hat E^0, \hat \lambda)$
for the stable curve $\hat p_0 = (\hat \Sigma, z, \hat u)$,
where
$p_0^+ = (\Sigma, z \cup z^+, u) \in \widehat{\M}(Y, \lambda, J)$ is a curve obtained
by adding the marked points $z^+$ to $(\Sigma, z)$ corresponding to the
additional marked points of $\hat p_0^+$, $S$ and $E^0$ are defined by
$S = \pi_M^{-1}(\hat S)$ and $E^0 = \hat E^0$, and
$\lambda : E^0 \to C^\infty(\hat P \times Y; \Wedge^{0, 1} V^\ast \hat P
\otimes (\R \partial_\sigma \oplus TY))$ is the map defined by the pull back of $\hat \lambda$ and
the isomorphism
\[
T_{(\sigma, y)} \hat Y = \xi_y \oplus (\R \partial_\sigma \oplus \R R_\lambda(y))
\cong (\pi_M^\ast TM)_y \oplus \C.
\]
Then the linear operator $D_{p_0}^+$ is not necessarily surjective,
but if we replace all vector spaces
$T_{\pi_Y \circ u_0(z_i)} Y$ except one in the range of $D_{p_0}^+$ with 
$T_{\pi_M \circ \pi_Y \circ u_0(z_i)} M$, then it becomes surjective.
Hence we can define the generating function $\mathcal{H}_0$ modulo $(t_y^2)$ using
the grouped multisections of the fiber products of $\widehat{\M}$ induced by those of
the corresponding fiber products of $\overline{\M}(M, j)$.
Then it is easy to see that
\begin{align*}
&[\overline{\M}^Y((\hat c_l), (y, x_l), (\alpha_l))]^0\\
&= (-1)^{\sum |\hat c_l|} \prod k_{\hat c_l} \cdot
[\overline{\M}(M, J)_{((\hat c_l), (\pi_M(y), (x_l / S^1)), ([M] \cap \alpha_l))}]^0,
\end{align*}
where $x_l$ are cycles in $M$, each $k_{\hat c_l}$ is the multiplicity of the
periodic orbits in $c_l$.
Note that in the left hand side of the above equation, the correction terms vanish
because they correspond to linear combinations of fiber products of several pre-Kuranishi
spaces, and for each fiber product, at least one factor has a locally free $S^1$-action.
(See \cite{EGH00} or \cite{Bo02} for more sophisticated expression of the above
equation.)


%% file: SFT-appendix.tex
%
%
%
\appendix
\section{Notation of differential}\label{diff notation}
We use the following notation in Section \ref{smoothness}.
\begin{defi}
Let $X$, $Y$ be real Banach spaces
(or finite dimensional vector spaces).
A continuous map $A : X \to Y$ is said to be differentiable at $x\in X$
if there exists a bounded operator
$DA_x :X\to Y$
such that for any $\epsilon >0$ there exists some constant $\delta>0$
such that
$||A(x+v) - A(x) + DA_x \cdot v ||_Y \leq \epsilon ||v||_X$
for any $||v||_X\leq \delta$.
We call $DA_x$ the differential of $A$ at $x\in X$.
$A$ is said to be of class $C^1$ if it is differentiable at every point of $X$ and
$DA : X\to L(X,Y)$ is continuous.
$A$ is said to be of class $C^k$ if it is of class $C^1$ and $DA$ is of class $C^{k-1}$.
Define $D^kA=D(D^{k-1}A) : X \to L(X,L(X,\dots, L(X,Y)\dots ))$
inductively.
Hence
\begin{align*}
&D^kA_x\cdot v^k\cdot v^{k-1} \cdot \dots \cdot v^1\\
&= \frac{\partial^k}
{\partial t^1\partial t^2 \dots \partial t^k}A(x+t^1v^1+t^2v^2+\dots
+t^kv^k)\Bigr|_{t^1=t^2=\dots=t^k=0} \in Y
\end{align*}
for any $v^1,v^2,\dots,v^k\in X$.
\end{defi}
